\documentclass[12pt,openany,twoside]{book}

\usepackage[english]{babel}

\usepackage{amssymb}
\usepackage{amsmath}
\usepackage{amsthm}
\usepackage{amsfonts}
\usepackage{mathrsfs}
\usepackage{color}
\usepackage[perpage,symbol*]{footmisc}
\usepackage{esvect}
\usepackage{latexsym}

\usepackage{anysize}

\usepackage{txfonts}

\usepackage{indentfirst}
\usepackage{soul}
\usepackage[normalem]{ulem}
\usepackage{cancel}

\usepackage{imakeidx}

\makeindex[name=words,title=Subject Index,columns=2,
intoc=true, options=-s index_stylenodots.ist]
\makeindex[name=symbols,title=Symbol Index,columns=2,
intoc=true, options=-s index_stylenodots.ist]

\usepackage{multicol}    

\usepackage{harpoon}
\newcommand{\myvec}[1]

\usepackage[colorlinks=true,
linkcolor=blue,
citecolor=red,
urlcolor=magenta,
]{hyperref}

\allowdisplaybreaks

\newtheorem{theorem}{Theorem}[section]
\newtheorem{lemma}[theorem]{Lemma}
\newtheorem{corollary}[theorem]{Corollary}
\newtheorem{proposition}[theorem]{Proposition}
\newtheorem{example}[theorem]{Example}
\theoremstyle{definition}
\newtheorem{remark}[theorem]{Remark}
\newtheorem{definition}[theorem]{Definition}
\renewcommand{\appendix}{\par
\setcounter{section}{0}%
\setcounter{subsection}{0}%
\setcounter{subsubsection}{0}%
\gdef\thesection{\@Alph\c@section}%
\gdef\thesubsection{\@Alph\c@section.\@arabic\c@subsection}%
\gdef\theHsection{\@Alph\c@section.}%
\gdef\theHsubsection{\@Alph\c@section.\@arabic\c@subsection}%
\csname appendixmore\endcsname
}

\usepackage{tikz}
\usepackage{wrapfig}

\let\oldmarginpar\marginpar
\renewcommand{\marginpar}[2][rectangle,draw,text width= 2cm,rounded corners]{
\oldmarginpar{
\scriptsize \tikz \node at (0,0) [#1]{#2};}
}
\numberwithin{section}{chapter}
\numberwithin{equation}{section}

\begin{document}

\title{\vspace{-6cm}\bf\Large Real-Variable
Theory of Hardy--Lorentz Spaces on
Quasi-Ultrametric Spaces of Homogeneous Type
with Reverse-Doubling Property}
\author{Chenfeng Zhu, Ryan Alvarado,
Xianjie Yan, \\ Dachun Yang and Wen Yuan}
\date{\color{red}{\today}}
\maketitle

\noindent Chenfeng Zhu

\medskip

\noindent Laboratory of Mathematics and Complex Systems
(Ministry of Education of China),
School of Mathematical Sciences,
Beijing Normal University, Beijing 100875,
The People's Republic of China

\smallskip

\noindent \textit{Current address:} School of Mathematical Sciences, Zhejiang
University of Technology, Hangzhou 310023,
The People's Republic of China

\smallskip

\noindent{\it E-mail address:} \texttt{chenfengzhu@zjut.edu.cn}

\bigskip

\noindent Ryan Alvarado

\medskip

\noindent Department of Mathematics, Amherst College, Amherst, MA, USA

\smallskip

\noindent{\it E-mail address:} \texttt{rjalvarado@amherst.edu}

\bigskip

\noindent Xianjie Yan

\medskip

\noindent Institute of Contemporary Mathematics,
School of Mathematics and Statistics, Henan University,
Kaifeng 475004,
The People's Republic of China

\smallskip

\noindent{\it E-mail address:} \texttt{xianjieyan@henu.edu.cn}

\bigskip

\noindent Dachun Yang (Corresponding author) and Wen Yuan

\medskip

\noindent Laboratory of Mathematics and Complex Systems (Ministry of Education of China),
School of Mathematical Sciences, Institute for Advanced Study, Beijing Normal University, Beijing 100875,
The People's Republic of China

\smallskip

\noindent {\it E-mails}: \texttt{dcyang@bnu.edu.cn} (D. Yang)

\noindent\phantom{{\it E-mails:}} \texttt{wenyuan@bnu.edu.cn} (W. Yuan)

\vspace{5cm}


\medskip


\frontmatter


\markboth{\scriptsize\rm\sc Preface}{\scriptsize\rm\sc Preface}
\chapter*{Preface}
\addcontentsline{toc}{chapter}{Preface}

Function spaces and their applications to the boundedness
of operators have always been one
of the core topics in harmonic analysis, with
Lebesgue spaces $L^p(\mathbb{R}^n)$ for $p\in(0,\infty]$
serving as one of the most fundamental and widely used examples.
As a refinement of the Lebesgue scale, the Lorentz spaces $L^{p,q}(\mathbb{R}^n)$
with $q\in(0,\infty]$ were introduced and studied by Lorentz
in the early 1950s \cite{Lorentz1950,Lorentz1951}.
Lorentz spaces naturally contain Lebesgue spaces and weak-Lebesgue spaces as
special cases, with
$$
L^{p,p}(\mathbb{R}^n)=L^p(\mathbb{R}^n)\
\mathrm{and}\ L^{p,\infty}(\mathbb{R}^n)=\mathrm{weak}-L^p(\mathbb{R}^n),
$$
and they also arise as
intermediate spaces between $L^p(\mathbb{R}^n)$ and $L^{\infty}(\mathbb{R}^n)$
via the real interpolation
method (see, for example, \cite{Calderon,LP1964}).
It is well known that,
due to their finer structural properties,
Lorentz spaces play an important role in addressing various
critical or endpoint problems across different areas of analysis,
including harmonic analysis (see, for example,
\cite{mtt2003,osttw2012,st2001,tw2001})
and partial differential equations
(see, for example, \cite{lly2015,mr2015,P2015}).
For further studies on Lorentz spaces defined over $\mathbb{R}^n$,
such as their dual spaces, we refer to the works of
Hunt \cite{Hunt}, Cwikel \cite{Cwikel1975},
Cwikel and Fefferman \cite{CF8081,CF84},
as well as \cite{bs1988,eks2008,kv2014,m1948,s1979,sw1971}.
This paradigm of refining classical scales extends naturally
to the theory of Hardy spaces, which also plays an important
role in harmonic analysis (see, for example,
\cite{cgn2017,c1974,Cowe77,fs1972,g1979,ky2014,L1978,lyh2022,lz2015,yyz2021})
and partial differential equations and
has been extensively studied
(see, for example, \cite{B14,bdt2021,clms1993,G2014M,
Muller1994,shyy2017,Semmes1994,st1989}).
Recall that the classical Hardy space $H^p(\mathbb{R}^n)$ with
$p\in(0,1]$ was originally introduced by Stein and Weiss
\cite{sw60} and further studied by Fefferman and
Stein \cite{fs1972}. In particular, when $p\in(0,1]$,
in studying the boundedness
of certain classes of singular integral operators, $H^p(\mathbb{R}^n)$
serves as a good substitute for $L^p(\mathbb{R}^n)$,
where, in the endpoint case,
$H^p(\mathbb{R}^n)$ is often replaced
by the weak Hardy space $H^{p,\infty}(\mathbb{R}^n)$.
Indeed, Alvarez and Milman \cite{am1986} proved that a
convolutional $\delta$-type Calder\'on--Zygmund operator $T$
with $\delta\in(0,1]$ and $T^*1=0$,
where $T^*$ denotes the adjoint operator of $T$ on $L^2(\mathbb{R}^n)$,
is bounded on $H^p(\mathbb{R}^n)$ for any given $p\in(\frac{n}{n+\delta},1]$.
However, in the endpoint case $p=\frac{n}{n+\delta}$,
$T$ may fail to be bounded on
$H^{\frac{n}{n+\delta}}(\mathbb{R}^n)$.
Liu \cite{Liu1991} later showed that in this case $T$ is bounded
from $H^{\frac{n}{n+\delta}}(\mathbb{R}^n)$
to $H^{\frac{n}{n+\delta},\infty}(\mathbb{R}^n)$.
Moreover, in order to find the largest function space from which
the Hilbert transform is bounded
into the weak Lebesgue space $L^{1,\infty}(\mathbb{R}^n)$,
Fefferman and Soria \cite{fs86} introduced
the weak Hardy space $H^{1,\infty}(\mathbb{R}^n)$ and,
using its $\infty$-atomic characterization,
proved that a convolutional Calder\'on--Zygmund operator
whose kernel satisfies the Dini-type condition is bounded from
$H^{1,\infty}(\mathbb{R}^n)$ to
$L^{1,\infty}(\mathbb{R}^n)$.
For more investigations on Hardy spaces and their variants,
we refer to the well-known monographs
\cite{Duo,fs1982,gr1985,G2014M,L1995,Pag,st1989,T1986}
and the articles by
Bonami et al. \cite{bcklyy2019,bgk2021},
Bui et al. \cite{bdk2020,BTAXT1,BTAXT2,BTAXT3},
Cruz-Uribe et al. \cite{cmn2019-1,cmn2019-2,cmn2020,cn2021,cw2014},
Grafakos et al. \cite{gnns2019,gnns2021},
Ho \cite{h2012,h2017,h2019},
Nakai and Sawano \cite{EY,ns2014},
and Sawano et al. \cite{s2013,shyy2017}.
We also refer to \cite{Alvarez1994,Alvarez1998,dls2006,
dlx2007,h20171,lyyy2025-2,syy1,syy21,yyyz16}
for further results on weak Hardy spaces
and to \cite{bow05,bow08,dl2003,dlx2007,Grafakos1992,HHLC,lyz2024,n2017,zyyw2021}
for the boundedness of operators on (weak) Hardy spaces.

In 1974, Fefferman et al. \cite{frs74} introduced
the classical Hardy--Lorentz spaces
$H^{p,q}(\mathbb{R}^n)$,
with $p\in(0,\infty)$ and $q\in(0,\infty]$,
which include the Hardy space $H^{p}(\mathbb{R}^n)$
and the weak Hardy space
$H^{p,\infty}(\mathbb{R}^n)$ as special cases and which also
naturally appear as the intermediate spaces
of $H^p(\mathbb{R}^n)$ and $L^\infty(\mathbb{R}^n)$
under real interpolation.
Parilov \cite{Parilov2006} later
obtained atomic characterizations of
$H^{1,q}(\mathbb{R}^n)$ for $q\in(1,\infty)$,
and in 2007 Abu-Shammala and Torchinsky \cite{at07}
investigated $H^{p,q}(\mathbb{R}^n)$ with $p\in(0,1]$ and $q\in(0,\infty)$.
In particular, they established their $\infty$-atomic characterizations,
real interpolation properties in the parameter $q$,
and boundedness results for
Calder\'on--Zygmund and related operators on these spaces.
Since then, various extensions of the classical
Hardy--Lorentz spaces have been well developed
in different settings.
In particular,
Almeida and Caetano \cite{ac2010}
introduced generalized Hardy spaces
which include the spaces $H^{p,q}(\mathbb{R}^n)$ from \cite{at07}
as special cases, and
obtained their maximal characterizations,
real interpolation results with the function parameters,
and boundedness results for certain classical operators acting on these spaces.
More recently, Jiao et al. \cite{jzzw2019} introduced the
variable Hardy--Lorentz spaces $H^{p(\cdot),q}(\mathbb{R}^n)$, where
$p(\cdot)$ is a nonnegative Lebesgue measurable function on $\mathbb{R}^n$
and $q\in(0,\infty]$. They established their various
real-variable characterizations involving
atoms and Littlewood--Paley functions,
and further obtained results concerning real interpolation, duality,
and boundedness properties of some singular integral operators
on these spaces.
In the anisotropic setting,
Liu et al. \cite{lyy2016,lyy2018} introduced the anisotropic
Hardy--Lorentz spaces
$H^{p,q}_A(\mathbb{R}^n)$
with $p\in(0,1]$ and $q\in(0,\infty]$, where
$A$ is a general expansive matrix on $\mathbb{R}^n$.
They established  several real-variable characterizations,
respectively, in terms of (finite) atoms, molecules,
maximal functions, and Littlewood--Paley functions,
and applied these to characterize the dual space of $H^{p,q}_A(\mathbb{R}^n)$
and obtain the boundedness of Calder\'on--Zygmund operators from
$H^{p,q}_A(\mathbb{R}^n)$ to $L^{p,q}(\mathbb{R}^n)$,
including the corresponding critical case.
Very recently, Liu et al. \cite{lwyy2019,lyy17} extended the
results in \cite{lyy2016,lyy2018} to anisotropic
variable Hardy--Lorentz spaces $H^{p(\cdot),q}_A(\mathbb{R}^n)$,
which include the variable exponent spaces $H^{p(\cdot),q}(\mathbb{R}^n)$
from \cite{jzzw2019} as special cases.
Independently, Almeida et al. \cite{abr2017} investigated
the anisotropic variable Hardy--Lorentz spaces
$H^{p(\cdot),q(\cdot)}(\mathbb{R}^n,A)$, where $p(\cdot)$ and
$q(\cdot)$ are nonnegative Lebesgue measurable functions on $(0,\infty)$
and $A$ is a general expansive matrix on $\mathbb{R}^n$.
In \cite{abr2017}, they obtained equivalent characterizations of
$H^{p(\cdot),q(\cdot)}(\mathbb{R}^n,A)$ using maximal
functions and atoms.
Notably, the two anisotropic variable Hardy--Lorentz space theories
in \cite{abr2017} and
\cite{lwyy2019,lyy17} are not directly comparable,
since in the former the variable
exponents $p(\cdot)$ are only defined on $(0,\infty)$
rather than on $\mathbb{R}^n$.
It is worth noting that anisotropic Euclidean spaces serve as
concrete and important models of spaces of homogeneous type studied
in this monograph. For research conducted in the
anisotropic Euclidean setting, we refer to the articles by
Bownik et al. \cite{b2003,b2005,b2007,b2008,bh2006,bll2023,bw2023,lby2014},
Cleanthous et al. \cite{cgn2017,cgn2019-1,cgn2019-2},
and Ho \cite{h2003}.
In particular, a differential characterization
for anisotropic Hardy spaces using a
parabolic differential equation of Calder\'on and
Torchinsky has been established in \cite{bw2023}.
For further results on Hardy--Lorentz spaces on $\mathbb{R}^n$,
we refer to the articles by
Jiao et al. \cite{jwxy2021,wzzj2018,zsj2019} and
Weisz et al. \cite{hlw2025,jcwyy2025,jwyyz2023,w2018}
and also to more related articles
\cite{h2014,lyyy2024,lyyy2025-1,llh2023,lql2020,lw2026,wyy2020,wyy2023}.

One fundamental development in function space theory came
with the consideration of function spaces defined
on environments much more general than
the Euclidean ambient.
In particular, in the 1970s, Coifman and Weiss \cite{CW1,Cowe77}
introduced the notion of a space of homogeneous type,
which is a triplet
$(X,\rho,\mu),$
where $(X,\rho)$ is a quasi-metric space and
$\mu$ is a doubling measure, as a natural setting for
the study of function spaces and the boundedness of
operators arising in harmonic analysis. In this context,
Coifman and Weiss \cite{Cowe77} introduced the atomic Hardy space
$H^{p,q}_{\rm CW}(X)$\index[symbols]{H@$H^{p,q}_{\rm CW}(X)$} with
$p\in(0,1]$ and $q\in(p,\infty]\cap[1,\infty]$,
proved that $H^{p,q}_{\rm CW}(X)$ is independent of the choice of $q$
(and therefore can be denoted simply by $H^{p}_{\rm CW}(X)$),
and further showed that the dual of $H^{p}_{\rm CW}(X)$
is the Lipschitz space
$\mathcal{L}_{\frac{1}{p}-1}(X)$
when $p\in(0,1)$, or the John--Nirenberg class of functions of bounded
mean oscillation $\mathrm{BMO\,}(X)$\index[symbols]{B@$\mathrm{BMO\,}(X)$} when $p=1$.
Given the natural limitations of such a general setting,
the spaces $H^{p}_{\rm CW}(X)$ are only expected to coincide with
their Euclidean counterparts when $p$ is sufficiently close to 1,
since higher-order vanishing moments involving polynomials are
typically unavailable on $X$. Nevertheless, a rich theory of Hardy spaces
has since developed.
Indeed, spaces of homogeneous type have become the most natural
and general underlying spaces for developing the real-variable
theory of various function spaces and their applications
to the boundedness of operators.
We refer to the monographs \cite{AlMi2015,DH,Hei,MiMiMiMo13,sfh20a},
the articles by Alvarado et al. \cite{AGH20,zwartestale,AWYY21,AYY24,AYY22},
Georgiadis et al. \cite{gk2020,gk2023,gkp2025,gn2018},
Nakai et al. \cite{kn2020,n06,n2010,ntyyz2026,sgn2018,tnyyz2025},
and Sawano et al. \cite{ss2021,sst18},
and more related articles  \cite{kmyz2024,lyyy2025-2,MaSe79i,w2025,yy2025}
for further studies about function spaces and their
applications to the boundedness of operators on spaces of homogeneous type.

In \cite{Cowe77}, Coifman and Weiss also established
the radial maximal function characterization of $H^{1}_{\rm CW}(X)$,
provided that an additional geometric assumption is imposed on the
space $(X,\rho,\mu)$. This led them to ask what extent
additional geometric conditions are necessary for the validity of
such characterizations. This question has since inspired a substantial body of work
aimed at establishing various real-variable
characterizations of the atomic Hardy
spaces on $X$ under minimal assumptions on the underlying ambient space.
Early progress in this vein was made in certain subclasses of spaces of
homogeneous type. More specifically, although the measure of any dilated
ball $\lambda B$ with $\lambda\in[1,\infty)$ has a quantitative
upper bound estimate in a space of homogeneous type,
it was recognized that there are situations where
a more precise lower bound estimate of the measure of $\lambda B$ is needed,
beyond merely comparing the relative sizes of $\lambda B$ and $B$.
With the Euclidean space $\mathbb{R}^n$ as the
canonical example, this need naturally led to the study of $n$-Ahlfors-regular
spaces, which are spaces of homogeneous type where the measure of a ball
is required to be exactly comparable to its radius raised to the power $n$.
In the setting of $1$-Ahlfors-regular spaces,
Mac\'ias and Segovia \cite{MaSe79ii} characterized $H^p_{\rm CW}(X)$
via the grand maximal function
when $p\in(0,1]$ is sufficiently close to 1,
Li \cite{Li} obtained a new maximal function
characterization of $H^p_{\rm CW}(X)$ using the space
of test functions introduced in \cite{HS}, and
Duong and Yan \cite{DY} characterized Hardy spaces via
the Lusin area function associated with certain semigroup.
Since then, Ahlfors-regular spaces have continued to attract considerable attention
in the study of the real-variable theory of function spaces;
see, for example, \cite{AlMi2015,DY,Li,MaSe79i,MaSe79ii}.
Further insight into this question came from realizing that
a robust function space theory can still be developed in settings more general
than Ahlfors-regular spaces. This led
to the study of RD-spaces, which are spaces of homogeneous type
where the equipped measure additionally satisfies
a reverse doubling condition. Such spaces provide
a natural framework for studying the
real-variable theory of function spaces and the boundedness
of operators. For further studies about function spaces and their
applications to the boundedness of operators on RD-spaces,
we refer to \cite{LLD,GLY,hmy06,hmy08,KYZ1,KYZ2,YZ1,YZ2,YZ3,YZ}.

Beyond geometric considerations, analytic tools such as
Calder\'on reproducing formulae  play an equally fundamental role
in function space theory,
providing a powerful tool for establishing
characterizations of function spaces and proving boundedness
results for operators.
Recall that the classical Calder\'on reproducing formulae on
$\mathbb{R}^n$ gain their prototypes from Calder\'on \cite{Calderon},
in which problems of the intermediate space and the
interpolation were studied via a reproducing method.
In 1974, Heideman \cite{Hei} gave a detailed treatment of the
homogeneous and the inhomogeneous
continuous Calder\'on reproducing formulae
and used them to study duality and  fractional
integrals on Lipschitz spaces (see also \cite{CT1,CT2}).
Since then, Calder\'on reproducing formulae and their variants
have become indispensable in the study of a wide range
of function spaces, including Hardy spaces
(see, for example, \cite{BTAXT1,BTAXT2,LLD,YZ3}),
Besov spaces and Triebel--Lizorkin spaces
(see, for example, \cite{MB1,MB2,MB3,Peetre,Y20051,Y20052,YZ}),
and also more general scale of function spaces
(see, for example, \cite{DW1,DW2,YSY}).
These tools also have further applications in
proving the boundedness of Calder\'on--Zygmund operators
(see, for example, \cite{LDC1,LDC2,LJ}).
For more applications of Calder\'on reproducing formulae,
particularly in the study of Musielak--Orlicz Hardy spaces and
variable exponent Hardy spaces, we refer to
\cite{XJB,hyy2025,YJD,EY,yyyz16,ylk2017,ZSY}.

Despite their success in the Euclidean setting,
Calder\'on reproducing formulae were initially unavailable
in the setting of a space of homogeneous type $X$, leading
much of the work on function spaces and their applications
to the operator boundedness
in this context, including \cite{DY,Li,MaSe79ii} mentioned above,
to rely on additional geometrical assumptions on the
underlying space, such as the reverse doubling
property of the measure
(see, for example, \cite{LLD,GLY,sst18,ZSY}).
Recently, a breakthrough on the analysis over $X$
without any additional geometrical assumptions came with
the work of Auscher and Hyt\"onen \cite{AH1,ah2015}, who,
using the system of random dyadic cubes
established by Hyt\"onen and Kairema \cite{HK}, constructed
an orthonormal wavelet basis of $L^2(X)$ with exponential decay.
Later, Han et al. \cite{HLW} established the wavelet reproducing
formulae, which hold in both the sense of spaces of
test functions and distributions.
Motivated by these advances, He et al. \cite{HLYY} introduced an
approximation of the identity with exponential decay
and obtained new Calder\'on reproducing formulae on $X$.
Together, these advances overcame initial obstacles in
function space theory and laid the foundation for developing the real-variable
theory of function spaces and for studying the
operator boundedness in general
spaces of homogeneous type.

Very recently,
He et al. \cite{ZYJLDW,HYY} completely answered
the aforementioned question of
Coifman and Weiss by developing a real-variable theory
of Hardy spaces and their localized versions,
including various equivalent characterizations and the boundedness
of sublinear operators, valid for $p\in(\frac{n}{n+\theta},1]$,
on spaces of homogeneous type with upper
dimension $n\in(0,\infty)$,
where $\theta\in(0,1)$ is the H\"older continuity exponent of the
orthonormal wavelet basis from \cite{AH1}.
Fu et al. \cite{fmy19} and Zhou et al. \cite{zhy} further
generalized the corresponding results in \cite{ZYJLDW},
respectively, to Musielak--Orlicz Hardy spaces and
Hardy--Lorentz spaces. Moreover,
Yan et al. \cite{YHYY2,YHYY1} introduced
the Hardy-type spaces $H_Y(X)$,
associated with ball quasi-Banach
function spaces $Y(X)$ on $X$,
and developed a complete real-variable theory of these spaces.
We also refer to \cite{AWYY21,FY,HHHLP,HHL,HWYY,syy1,syy21,w2025,whyy21,y2024}
for further related studies.

One issue that arises in the consideration of Hardy spaces
$H^p(X)$ in the setting of a space of homogeneous type $X$
is that unless $p\in(0,1]$ is sufficiently near to 1,
$H^p(X)$ spaces become trivial, in the sense that $H^p(X)$ only contains
constant elements. At the heart of the matter
is the fact that spaces of H\"older continuous functions reduce
to exactly constant
functions if the smoothness parameter is too large.
In the Euclidean setting, for example, H\"older spaces
of any order $\alpha>1$ are trivial; however,
this upper smoothness bound need not be $1$
in an arbitrary quasi-metric space.
Therefore, a central question---one that implicitly
underpins the work of Coifman and Weiss \cite{Cowe77}---is to determine
the range of $p$, depending on the geometry of $X$,
for which a satisfactory theory of Hardy spaces exists.
Coifman and Weiss \cite{Cowe77} mentioned
an unspecified range of $p$ for which their construction works.
The first significant
attempt to clarify this range of $p$ in a quantitative
fashion is due to Mac\'ias and Segovia \cite{MaSe79ii},
who, in the setting of $1$-Ahlfors-regular spaces,
obtained a grand maximal function
characterization for the Coifman--Weiss atomic Hardy spaces $H^{p}_{\rm CW}(X)$ for
$p\in(p_{\rho},1]$, where
\begin{equation*}
p_{\rho}:=\frac{1}{1+[\log_2(A_\rho(2A_\rho+1))]^{-1}}
\end{equation*}
and $A_\rho$ is the best constant usable in the
quasi-triangle inequality [see \eqref{2234-intro} below]
satisfied by the quasi-metric $\rho$.
The formula for $p_{\rho}$ is due
to the metrization result that Mac\'ias and Segovia developed in \cite{MaSe79i},
which has continued to play a very influential role in the analysis on $X$.
Despite \cite{MaSe79ii} being a cornerstone result in the theory of Hardy spaces
on spaces of homogeneous type and having been widely
cited (see, for example, \cite{DH,HS,H2001,Pag}),
the value of $p_{\rho}$ is far from optimal. For example, in the $1$-dimensional
Euclidean setting, where $A_\rho=1$,
we have $p_{\rho}=\frac{1}{1+(\log_23)^{-1}}$, which is strictly
larger than the expected threshold $\frac{1}{2}$, and therefore the work in
\cite{MaSe79ii} does not provide a genuine generalization of the classical
theory in the Euclidean setting.

For many years, the limitations of the theory of Mac\'ias and Segovia,
in particular the nonoptimality of the value of $p_{\rho}$,
have been implicitly accepted as perhaps an inevitable consequence
of working within such a general framework.
An important advance toward addressing these limitations was provided
by Mitrea et al. \cite{MiMiMiMo13} who introduced
and studied the index
\begin{align}\label{index-intro}
\mathrm{ind\,}(X,\rho)
:=\sup_{\varrho\sim\rho}\left(\log_2C_{\varrho}\right)^{-1},
\end{align}
where the supremum is taken over all quasi-metrics $\varrho$
on $X$ that are \emph{bi-Lipschitzly}\footnote{For any given
quasi-ultrametric spaces $(X_1,\rho_1)$
and $(X_2,\rho_2)$, a mapping $\Phi:(X_1,\rho_1)\to(X_2,\rho_2)$
is said to be \emph{bi-Lipschitz}\index[words]{bi-Lipschitz}
if there exist $0<C_1\leq C_2<\infty$
such that, for any $x,y\in X_1$,
$$C_1\rho_1(x,y)\leq\rho_2(\Phi(x),
\Phi(y))\leq C_2\rho_1(x,y).$$}
equivalent to $\rho$, and, for any quasi-metric
$\rho$ on $X$, $C_{\rho}$ is the best constant $C\in(0,\infty)$ for which
the following quasi-ultrametric inequality
$\rho(x,y)\leq C\max\left\{\rho(x,z),\,\rho(z,y)\right\}$
holds for all $x,y,z\in X$.
In a sense, this index captures the nature of the ``best" quasi-metric
on $X$ that is bi-Lipschitz equivalent to $\rho$
and captures the maximal regularity that any
approximation of the identity on $X$ can achieve.
Its value also reflects
certain geometric information about the underlying space; for example,
${\rm ind}\,(\mathbb{R}^n,|\cdot-\cdot|)=1$ and
${\rm ind}\,([0,1]^n,|\cdot-\cdot|)=1$,
where $|\cdot-\cdot|$ denotes the Euclidean distance, whereas
${\rm ind}\,(X,\rho)\geq 1$ whenever there exists a genuine metric on $X$ that is
pointwise equivalent to $\rho$; ${\rm ind}\,(X,\rho)<1$ if
$(X,\rho)$ cannot be bi-Lipschitz embedded into some ${\mathbb{R}}^n$ with
$n\in{\mathbb{N}}$; and ${\rm ind}\,(X,\rho)=\infty$ if there exists an ultrametric
on $X$ [see Remark~\ref{remark2.2}(iii)] which is pointwise equivalent to $\rho$.
It is instructive to note that, although the quasi-ultrametric inequality
is equivalent to the more standard quasi-triangle inequality
\begin{align}\label{2234-intro}
\rho(x,y)\leq A_\rho\left[\rho(x,z)+\rho(z,y)\right],
\end{align}
as we will demonstrate throughout this monograph,
it is the constant $C_\rho$, rather than $A_\rho$, that will prove
to be of utmost importance with regards to the sharpness of the main results in
this work. Given our use of the quasi-ultrametric inequality
and the fact that we will allow quasi-metrics to be
merely quasi-symmetric, we choose to adopt the terminology
quasi-ultrametric throughout this monograph, to avoid any confusion
with the standard term quasi-metric used in the literature.

Based on the index \eqref{index-intro} and a new metrization theorem,
which sharpens the work of Mac\'ias and
Segovia \cite{MaSe79i}, Mitrea et al. \cite{MiMiMiMo13}
constructed a discrete approximation of the identity of any order
$$
0<\varepsilon_0<\min\left\{n+1,\,\mathrm{ind\,}(X,\rho)\right\}
$$
on $n$-Ahlfors-regular spaces
with $n\in(0,\infty)$, under the additional assumption
that all singletons have measure zero.
Building on \cite{MiMiMiMo13}, Alvarado and Mitrea \cite{AlMi2015}
constructed an approximation of the identity
in general $n$-Ahlfors-regular spaces
and were successful in further extending the range of order to
$$
0<\varepsilon_0\preceq\mathrm{ind\,}(X,\rho),
$$
where the symbol $\varepsilon_0\preceq\mathrm{ind\,}(X,\rho)$
means that $\varepsilon_0\leq\mathrm{ind\,}(X,\rho)$ and
$\varepsilon_0=\mathrm{ind\,}(X,\rho)$ is only permissible
when the supremum in \eqref{index-intro} is attained.
By developing fundamental tools such as
a sharp Lebesgue differentiation theorem,
a maximally smooth approximation of the identity,
and a Calder\'on--Zygmund decomposition,
Alvarado and Mitrea \cite{AlMi2015} obtained
various real-variable characterizations of $H^p(X)$
in terms of atoms, molecules, and maximal functions,
valid for the sharp range
\begin{align}\label{rangep}
p\in\left(\frac{n}{n+\mathrm{ind\,}(X,\rho)},1\right].
\end{align}
They also developed a sharp theory of Besov and
Triebel--Lizorkin spaces on $X$ in the same article.
Notably, the range of $p$ in \eqref{rangep} represents
a marked improvement over the work in \cite{MaSe79ii}. Indeed,
when the underlying space $X$ is $\mathbb{R}^n$, \eqref{rangep} reduces to
the expected range $(\frac{n}{n+1},1]$, and so, in contrast to \cite{MaSe79ii},
the results in \cite{AlMi2015} can be regarded as a genuine generalization of the
classical Hardy space theory from $\mathbb{R}^n$ to  $n$-Ahlfors regular
spaces. It is also worth noting that,
when $X$ is an $n$-Ahlfors regular space,
the range of $p$ in \eqref{rangep} is, in general, larger than the range
$p\in(\frac{n}{n+\theta},1]$ appearing in \cite{ZYJLDW}.

In order to set the stage for the main results in this monograph,
let us briefly recall the idea of the construction
of the approximation of the identity
on spaces of homogeneous type $(X,\rho,\mu)$,
which is originally due to Coifman;
see also \cite[p.\,40]{DaJoSe85}.
More precisely, choose a smooth function $h:[0,\infty)\to[0,\infty)$
such that $h(x)=1$ for any $x\in[0,\frac{1}{2}]$
and $h(x)=0$ for any $x\in[2,\infty)$.
Let $T_k$, for any $k\in\mathbb{Z}$,
be the integral operator with kernel $2^{k}h(2^{k}\rho(x,y))$.
Let $M_k$ and $W_k$, for any $k\in\mathbb{Z}$, be the operators
of multiplication by $\frac{1}{T_k(1)}$ and $(T_k(\frac{1}{T_k(1)}))^{-1}$,
respectively. For any $k\in\mathbb{Z}$, let
\begin{align}\label{1518}
\mathcal{S}_k:=M_kT_kW_kT_kM_k.
\end{align}
Then, using the regularity property of $\rho$,
which stems from the metrization
theorem in \cite{MiMiMiMo13}
(see, for example, Theorem~\ref{1000}(iii) below),
together with the doubling property satisfied by
the measure $\mu$ (see \eqref{doublingcond} below),
we can prove that the kernel of $\mathcal{S}_k$
satisfies the desired estimates for an approximation of the identity
(see, for instance, \cite[p.\,40]{DaJoSe85}).
Under the additional assumption that $\mu$ is Borel regular
(or, more generally, Borel-semiregular),
we can use the Lebesgue differentiation theorem
(see, for instance, \cite[Theorem~1.8]{H2001}
and \cite[Theorem~3.14]{AlMi2015})
to conclude that, for any $p\in[1,\infty)$,
\begin{align}\label{1519}
\lim_{k\to\infty}\mathcal{S}_k=\mathcal{I}
\end{align}
in $L^p(X)$,
where $\mathcal{I}$ is the identity operator on $L^p(X)$.
Moreover, if $\mu(X)=\infty$, then
\begin{align}\label{1520}
\lim_{k\to-\infty}\mathcal{S}_k=0
\end{align}
in $L^p(X)$ for any $p\in(1,\infty)$.
As for applications,
such an approximation of the identity
can be used, for example, to establish the Littlewood--Paley theory
for $L^p(X)$ with $p\in(1,\infty)$ and the $T(1)$ theorem on
$(X,\rho,\mu)$; see, for example,
\cite{DaJoSe85}.

Based on the aforementioned approximation of
the identity $\{\mathcal{S}_k\}_{k\in\mathbb{Z}}$
in \eqref{1518},
one can establish the Calder\'on reproducing formula
by using Coifman's idea of the decomposition of the identity on $L^2(X)$.
For any $k\in\mathbb{Z}$, let
$$\mathcal{D}_k:=\mathcal{S}_k-\mathcal{S}_{k-1}.$$
Then, by \eqref{1519} and \eqref{1520}, we conclude that
\begin{align*}
\sum_{k=-\infty}^\infty\mathcal{D}_k=\lim_{k\to\infty}\mathcal{S}_k
-\lim_{k\to-\infty}\mathcal{S}_k=\mathcal{I}
\end{align*}
in $L^2(X)$
and hence,
for any fixed integer $N\in\mathbb{N}$, we can write
\begin{align*}
\mathcal{I}=\left(\sum_{k=-\infty}^\infty\mathcal{D}_k\right)
\left(\sum_{l=-\infty}^\infty\mathcal{D}_l\right)
=\sum_{k=-\infty}^\infty\mathcal{D}_k^N\mathcal{D}_k
+\sum_{k=-\infty}^\infty\sum_{|l|>N}\mathcal{D}_{k+l}\mathcal{D}_k
=:\mathcal{T}_N+\mathcal{R}_N
\end{align*}
in $L^2(X)$, where $\mathcal{D}_k^N:=\sum_{|l|\leq N}\mathcal{D}_{k+l}$
for any $k\in\mathbb{Z}$.
If $N$ is sufficiently large,
then we can show that the operator $\mathcal{R}_N$
is bounded on both $L^2(X)$ and the
space $\mathring{\mathcal{G}}^\varepsilon_0(\beta,\gamma)$
of test functions (see Definitions~\ref{testfunction} and~\ref{2207} below)
with operator norms at most 1, where $0<\beta,\gamma<\varepsilon<\infty$
and $\varepsilon$ is the order of $\{\mathcal{S}_k\}_{k\in\mathbb{Z}}$.
In this case, $\mathcal{T}_N$ is invertible and
$\mathcal{T}_N^{-1}$ is bounded on
both $L^2(X)$ and $\mathring{\mathcal{G}}^\varepsilon_0(\beta,\gamma)$.
Formally, we have
\begin{align*}
f=\mathcal{T}_N^{-1}\mathcal{T}_Nf
=\mathcal{T}_N^{-1}\left(\sum_{k=-\infty}^\infty
\mathcal{D}_k^N\mathcal{D}_k\right)f
=\sum_{k=-\infty}^\infty\mathcal{T}_N^{-1}
\mathcal{D}_k^N\mathcal{D}_kf
=\sum_{k=-\infty}^\infty\widetilde{\mathcal{D}}_k
\mathcal{D}_kf
\end{align*}
after setting
$\widetilde{\mathcal{D}}_k:=\mathcal{T}_N^{-1}\mathcal{D}_k^N$
for any $k\in\mathbb{Z}$.
We then obtain the following homogeneous continuous Calder\'on
reproducing formula
\begin{align}\label{1532}
f=\sum_{k=-\infty}^\infty\widetilde{\mathcal{D}}_k\mathcal{D}_kf,
\end{align}
where the series in \eqref{1532}
converges in $\mathring{\mathcal{G}}^\varepsilon_0(\beta,\gamma)$
for any $f\in\mathring{\mathcal{G}}^\varepsilon_0(\beta,\gamma)$,
in $L^p(X)$ with $p\in(1,\infty)$ for any $f\in L^p(X)$,
and in the space $(\mathring{\mathcal{G}}^\varepsilon_0(\beta,\gamma))'$
of distributions (see Definition~\ref{1558} below)
for any $f\in(\mathring{\mathcal{G}}^\varepsilon_0(\beta,\gamma))'$.
There are also corresponding inhomogeneous or discrete versions of \eqref{1532}.

In this monograph, we begin by constructing
a maximally smooth approximation of the
identity of order $0<\varepsilon\preceq\mathrm{ind\,}(X,\rho)$
on general quasi-ultrametric spaces of homogeneous type $(X,\rho,\mu)$,
which we then use to prove the boundedness of Calder\'on--Zygmund
operators on spaces of test functions and thereby establish
sharp Calder\'on reproducing formulae
on ultra-RD-spaces (see Definition~\ref{RD}).
As applications, we also obtain the Littlewood--Paley
function characterizations of Hardy spaces and Triebel--Lizorkin
spaces on ultra-RD-spaces.
Finally, we establish a complete real-variable theory including
equivalent characterizations and boundedness of operators
on Hardy--Lorentz spaces $H^{p,q}(X)$ in
quasi-ultrametric spaces of homogeneous type.
It is remarkable that, in this monograph, we only assume that $\mu$
is Borel-semiregular, that singletons may have
strictly positive measure, and that $\rho$ is a quasi-ultrametric, which
provides a more general theoretic setting than the space of
homogeneous type in the sense of Coifman and Weiss.
Moreover, we point out that the maximal order of approximation
of the identity and the regularity exponent of Calder\'on
reproducing formulae obtained in this monograph are sharp,
since they are based on the sharp regularity of the quasi-ultrametric.
Given the sharpness of the tools constructed in this monograph,
we are able to develop the real-variable
theory of $H^{p,q}(X)$ for sharp ranges of $p$ and $q$,
under minimal assumptions on the space $X$.
In addition, all results presented in this monograph are new and appear
in print for the first time.

To be precise, the whole monograph consists of five chapters.

As a foundational step, in Chapter~\ref{Chapter1},
we construct a maximally smooth approximation of the identity
on quasi-ultrametric spaces of homogeneous type.
We begin by laying out several
necessary conventions regarding notation and recalling the
basic notions of the underlying spaces. We then
introduce a new type of approximation of
the identity with bounded support on quasi-ultrametric
spaces of homogeneous type.
By an argument similar to that used in the proof of
\cite[Theorem 3.22]{AlMi2015} and the properties of
the regularized quasi-ultrametric,
we prove the existence of such approximation of the identity,
which is the main result of this chapter and extends
\cite[Theorem 3.22]{AlMi2015} from Ahlfors-regular quasi-ultrametric
spaces to quasi-ultrametric spaces of homogeneous type.

Building on this foundation,
Chapter~\ref{Chapter2} establishes sharp Calder\'on
reproducing formulae on ultra-RD-spaces $X$. Specifically,
we first recall the notions of test functions and
distributions on $X$.
By successfully constructing maximally smooth cut-off functions,
we prove the boundedness on spaces of test functions
of Calder\'on--Zygmund operators
whose kernels satisfy certain size and regularity conditions.
To obtain sharp homogeneous continuous Calder\'on reproducing formulae,
we first divide the identity operator into a main operator
$\mathcal{T}_N$ and a remainder operator
$\mathcal{R}_N$ for any fixed $N\in\mathbb{N}$.
Then, by using some properties of the composition of
two homogeneous approximations of the identity,
we prove that the operator norms of $\mathcal{R}_N$ on both $L^2(X)$
and spaces of test functions can be made sufficiently small if $N$
is sufficiently large, which ensures that $\mathcal{T}_N^{-1}$
exists and is bounded on both $L^2(X)$ and spaces of test functions.
From this and the homogeneous approximation of the identity,
we establish sharp homogeneous continuous Calder\'on reproducing formulae on $X$.
Moreover, by an argument similar to that used in the proof of the
homogeneous case and by the inhomogeneous approximation of the identity,
we also establish sharp inhomogeneous continuous Calder\'on
reproducing formulae on $X$. On the other hand,
to establish sharp homogeneous discrete Calder\'on reproducing formulae,
we first split the $k$-level dyadic cubes into a sum of dyadic
subcubes in level $k+j$ and then divide the identity operator
into the Riemann sum operator $\mathfrak{S}_N$
and the remainder operator $\mathfrak{R}_N$.
By proving that the operator norms of $\mathfrak{R}_N$ on both $L^2(X)$
and spaces of test functions can be made sufficiently small if both $N$ and $j$
are sufficiently large, we conclude that $\mathfrak{S}_N^{-1}$ exists and
is bounded on both $L^2(X)$ and spaces of test functions,
which, combined with the homogeneous approximation of the identity,
further implies sharp homogeneous discrete Calder\'on reproducing
formulae on $X$.
We also obtain sharp inhomogeneous discrete Calder\'on
reproducing formulae $X$ by a similar argument.
All of these sharp Calder\'on reproducing formulae are established
in the sense of spaces of test functions, distributions,
and $L^p(X)$ with $p\in(1,\infty)$.
The range of parameters for the
spaces of test functions and distributions
in which Calder\'on reproducing formulae hold
is sharp, given the sharp regularity of both
Calder\'on--Zygmund operators and
approximations of the identity.

As a first application,
Chapter~\ref{Chapter3} characterizes
the Lusin-area-function-based Hardy space
on ultra-RD-spaces $(X,\rho,\mu)$
via atoms and various Littlewood--Paley functions.
We begin by defining the Hardy space $H^p(X)$
via the Lusin area function
and prove, using the sharp homogeneous
discrete Calder\'on reproducing formula and the Fefferman--Stein
vector-valued maximal inequality on $X$,
that this definition of $H^p(X)$ is independent of the choice of the
approximation of the identity.
By utilizing the sharp homogeneous continuous Calder\'on reproducing formula,
we obtain the atomic characterization of $H^p(X)$.
From this and a key lemma concerning the Littlewood--Paley
auxiliary function, we further deduce the Littlewood--Paley
$g$-function and the $g_\lambda^*$-function characterizations of $H^p(X)$
for the best-known range of $\lambda$.
In particular, in the special case when $X$
is an $n$-Ahlfors-regular space as in \cite{AlMi2015},
the main results in this chapter give brand new characterizations
of the Hardy space studied in \cite{AlMi2015} via
the Lusin area function,
the Littlewood--Paley $g$-function,
and the Littlewood--Paley $g_\lambda^*$-function,
for the sharp range of $p\in(\frac{n}{n+\mathrm{ind\,}(X,\rho)},1]$.

Chapter~\ref{Chapter4}
extends this analysis to the more general Triebel--Lizorkin scales
via characterizing
Triebel--Lizorkin spaces on ultra-RD-spaces $(X,\rho,\mu)$
with upper dimension $n\in(0,\infty)$
by the Lusin area function and the
Littlewood--Paley $g_\lambda^*$-function. To this end,
we begin by recalling the notion of homogeneous
Triebel--Lizorkin spaces $\dot{F}^{p,q}_s(X)$ with the sharp ranges
$s\in(-{\mathrm{ind\,}}(X,\rho),{\mathrm{ind\,}}(X,\rho))$ and
$$
\max\left\{\frac{n}{n+{\mathrm{ind\,}}(X,\rho)},\,
\frac{n}{n+{\mathrm{ind\,}}(X,\rho)+s}\right\}<p,q\leq\infty.
$$
Then we establish the homogeneous
Plancherel--P\'olya inequality
by using the homogeneous discrete Calder\'on reproducing formula,
an almost orthogonal estimate,
and the Fefferman--Stein
vector-valued maximal inequality.
Using this inequality,
we establish the Lusin area function
characterization of $\dot{F}^{p,q}_s(X)$.
From this and an important change-of-angle formula for a
variant of the Lusin area function, we further deduce the
Littlewood--Paley $g_\lambda^*$-function characterization
of $\dot{F}^{p,q}_s(X)$ for the best-known range of $\lambda$.
When $s=0$ and $q=2$, these characterizations
of $\dot{F}^{p,2}_0(X)$ reduce to the corresponding ones of
$H^p(X)$ obtained in Chapter~\ref{Chapter3}.
We also extend all of the corresponding results
to the inhomogeneous case.
In particular, in the special case when $X$
is an $n$-Ahlfors-regular space as in \cite{AlMi2015},
the main results of this chapter give completely new characterizations
of the Triebel--Lizorkin spaces studied in \cite{AlMi2015} via
the Lusin area function and the Littlewood--Paley $g_\lambda^*$-function,
for the same optimal ranges of $s$, $p$, and $q$ as above.

Finally, armed with these tools,
Chapter~\ref{ChapterHL} develops a complete real-variable
theory of Hardy--Lorentz spaces on ultra-RD-spaces $(X,\rho,\mu)$
with upper
dimension $n\in(0,\infty)$. To this end, we first introduce
the Hardy--Lorentz space $H^{p,q}_*(X)$ via the grand maximal function
with the sharp ranges
$$
p\in\left(\frac{n}{n+{\mathrm{ind\,}}(X,\rho)},\infty\right)
\ \ \text{and}\ \
q\in(0,\infty],
$$
and show that it coincides with both the radial maximal-function-based
and non-tangential maximal-function-based Hardy--Lorentz spaces.
When $p\in(1,\infty)$ and $q\in(0,\infty]$,
we show that $H^{p,q}_*(X)$ coincides with the Lorentz
space $L^{p,q}(X)$, and, for
\begin{align}\label{rangepq}
p\in\left(\frac{n}{n+{\mathrm{ind\,}}(X,\rho)},1\right]
\ \ \text{and}\ \
q\in(0,\infty],
\end{align}
we obtain various real-variable characterizations of $H^{p,q}_*(X)$
in terms of (finite) atoms, molecules,
and various Littlewood--Paley functions.
Based on these characterizations, we establish a duality theorem
between Hardy--Lorentz spaces and
Campanato--Lorentz spaces,
prove a real interpolation theorem for Hardy--Lorentz spaces,
and derive boundedness results for
Calder\'on--Zygmund operators on
Hardy--Lorentz spaces, including the critical cases.
It should be emphasized that many of the main results in
this chapter are actually established in the more general setting of
quasi-ultrametric spaces of homogeneous type.

A salient feature of the work undertaken in Chapter~\ref{ChapterHL}
is the range of $p$ in \eqref{rangepq}, which captures,
in an optimal and quantitative manner,
the relationship between the geometry of the underlying space $X$
and the range of $p$ for which a satisfactory theory of
$H^{p,q}_\ast(X)$ exists. In particular, we recover the familiar
condition $p\in(\frac{n}{n+1},1]$ in the case where $X=\mathbb{R}^n$.
More strikingly, there exist ultra-RD spaces
where the range of $p$ in \eqref{rangepq}
is strictly larger than $(\frac{n}{n+1},1]$.
For example, in any ultrametric space (such as
the four-corner planar Cantor set equipped with the $1$-dimensional Hausdorff
measure and the standard Euclidean distance),
the range of $p$ in \eqref{rangepq} reduces to $(0,1]$, so that
a complete real-variable theory of
Hardy--Lorentz spaces is available for any $p\in(0,1]$ and $q\in(0,\infty]$.
This is a remarkable feature, since these full ranges $p\in(0,1]$ and $q\in(0,\infty]$
cannot be obtained in any of the related results currently available
in the literature, including recent works such as \cite{YHYY2,YHYY1,zhy}
for Hardy--Lorentz spaces and \cite{ZYJLDW}
for Hardy spaces [namely $H^{p,p}_\ast(X)=H^{p}_\ast(X)$ spaces],
where the techniques employed therein do not allow for having $p\leq\frac{n}{n+1}$.

Throughout this monograph, we always let
$\mathbb{N}:=\{1,2,\ldots\}$,\index[symbols]{N@$\mathbb{N}$}
$\mathbb{Z}_+:=\mathbb{N}\cup\{0\}$,\index[symbols]{Z@$\mathbb{Z}_+$}
and $\mathbb{R}_+:=[0,\infty)$.\index[symbols]{R@$\mathbb{R}_+$}
We denote by $C$ a \emph{positive constant} which is independent
of the main parameters involved, but may vary from line to line.
We use $C_{\alpha,\dots}$ to denote a positive constant depending
on the indicated parameters $\alpha,\, \dots$.
The symbol $f\lesssim g$ means $f\le Cg$
and, if $f\lesssim g\lesssim f$, we then write $f\sim g$.
If $f\le Cg$ and $g=h$ or $g\le h$,
we then write $f\lesssim g=h$ or $f\lesssim g\le h$.
For any $r\in\mathbb{R}$, let $r_+:=\max\{0,\,r\}$\index[symbols]{R@$r_+$}
and, for any $a,b\in\mathbb{R}$,
$a\land b:=\min\{a,\,b\}$\index[symbols]{L@$\land$}.
For any $\alpha\in\mathbb{R}$,
we denote by $\lceil\alpha\rceil$ the smallest integer
not less than $\alpha$. Given sets $A$ and $B$,
the notation $A \subsetneqq B$\index[symbols]{S@$\subsetneqq$}
will be used to indicate that $A$ is a proper subset of $B$.
For any index $p\in[1,\infty]$,
we denote by $p'$ its \emph{conjugate index}; that is,
$\frac{1}{p}+\frac{1}{p'}=1$.
Let $(X,\rho)$ be a quasi-ultrametric space and $\mu$
be a nonnegative measure on $X$.
For any $\mu$-measurable function $f:X\to\mathbb{C}$
and for any $\mu$-measurable set $E\subset X$
with $\mu(E)\in(0,\infty)$, $m_E(f)$,
the average of $f$ on $E$, is defined by setting
\begin{align}\label{mBf}
m_E(f):=f_E:=\frac{1}{\mu(E)}\int_Ef(x)\,d\mu(x).
\end{align}
If $E$ is a subset of $X$, we denote by ${\mathbf{1}}_E$ its
\emph{characteristic function},
by $E^\complement$ the set $X\setminus E$,
and by $\#E$ the \emph{cardinality} of $E$.
Let $(\mathscr{X}_1,\|\cdot\|_{\mathscr{X}_1})$
and $(\mathscr{X}_2,\|\cdot\|_{\mathscr{X}_2})$
be two (quasi-)normed vector spaces.
Then we use $\mathscr{X}_1\hookrightarrow\mathscr{X}_2$
\index[symbols]{X@$\mathscr{X}_1\hookrightarrow\mathscr{X}_2$} to
denote that $\mathscr{X}_1\subset \mathscr{X}_2$ and there exists
a positive constant $C$ such that,
for any $f\in\mathscr{X}_1$, $\|f\|_{\mathscr{X}_2}
\leq C\|f\|_{\mathscr{X}_1}$.
Finally, in all proofs, we consistently
retain the symbols introduced in the
original theorem (or related statement).

This project is partially supported by
the National Natural Science Foundation of China
(Grant Nos. 12431006, 12501129, 12371093, and 12301112),
the Beijing Natural Science Foundation (Grant No. 1262011),
the Fundamental Research Funds
for the Central Universities (Grant No. 2253200028),
and the China Postdoctoral Science Foundation (Grant No. 2025M783094).

\vspace{0.8cm}

\noindent Beijing, The People's Republic of China                     \hfill Chenfeng Zhu

\noindent Amherst, Massachusetts, United States of America            \hfill Ryan Alvarado

\noindent Kaifeng, The People's Republic of China                     \hfill Xianjie Yan

\noindent Beijing, The People's Republic of China                     \hfill Dachun Yang

\noindent Beijing, The People's Republic of China                     \hfill Wen Yuan

\medskip

\noindent \today

\newpage

\markboth{\scriptsize\rm\sc Contents}{\scriptsize\rm\sc Contents}
\tableofcontents

\mainmatter

\chapter[A Maximally Smooth
Approximation of the Identity\\ on
Quasi-Ultrametric Spaces
of Homogeneous Type]{A Maximally Smooth
Approximation \\ of the Identity on
Quasi-Ultrametric Spaces
of Homogeneous Type}
\label{Chapter1}
\markboth{\scriptsize\rm\sc An
approximation of the identity on
quasi-ultrametric spaces
of homogeneous type}
{\scriptsize\rm\sc
An approximation of the identity on
quasi-ultrametric spaces
of homogeneous type}

The main target of this chapter is to establish an
approximation of the identity on quasi-ultrametric
spaces of homogeneous type that exhibits a
maximal degree of smoothness (measured on the H\"older scale)
that the underlying ambient space can permit.
To this end, in Section \ref{ss2.1},
we first recall some definitions and
basic properties of
quasi-ultrametric and
quasi-ultrametric spaces.
Then we present the notion of quasi-ultrametric measure spaces,
which includes spaces such as
quasi-ultrametric spaces of homogeneous type,
ultra-RD-spaces, and Ahlfors-regular quasi-ultrametric spaces,
and list some basic assumptions on the underlying
quasi-ultrametric measure space.
Next, we recall the lower smoothness index of quasi-ultrametric spaces,
which plays an important role throughout this monograph.
In Section \ref{ss2.2}, we introduce a maximally smooth
approximation of the identity with bounded support.
Then we prove its existence and further obtain some useful
results related to this approximation of the identity
on quasi-ultrametric spaces of homogeneous type.

\section{Underlying Quasi-Ultrametric Measure Spaces}
\label{ss2.1}

In this section, we recall some definitions and
basic properties of
quasi-ultrametric and
quasi-ultrametric spaces,
and then present the notion of quasi-ultrametric measure spaces
and some basic assumptions that we will use throughout this monograph.
To this end, we first begin with recalling the
definition of quasi-ultrametrics;
see, for instance, \cite[pp.\,34--35]{AlMi2015}.

\begin{definition}\label{D2.1}
Let $X$ be a nonempty set. A function
$\rho:X\times X\to\mathbb{R}_+$\index[symbols]{R@$\rho$} is
called a \emph{quasi-ultrametric}\index[words]{quasi-ultrametric}
if there exist two constants $C_0,C_1\in(0,\infty)$ such that,
for any $x,y,z\in X$,
\begin{enumerate}
\item[\rm(i)] $\rho(x,y)=0$ if and only if $x=y$;

\item[\rm(ii)] (the \emph{quasi-symmetry condition})
$\rho(y,x)\le C_0\rho(x,y)$;

\item[\rm(iii)] (the \emph{quasi-ultrametric condition})
$\rho(x,y)\le C_1
\max\{\rho(x,z),\,\rho(z,y)\}$.
\end{enumerate}	
\end{definition}

If $X$ has cardinality at least $2$,
then it is easy to check that $C_0,C_1\in[1,\infty)$.
Throughout this monograph,
let $C_\rho$\index[symbols]{C@$C_\rho$}
and $\widetilde{C}_\rho$\index[symbols]{C@$\widetilde{C}_\rho$}
be the least constants,
respectively, in (iii) and (ii) of Definition~\ref{D2.1}; that is,
\begin{align}
\label{Crho}
C_\rho:=\sup_{\genfrac{}{}{0pt}{}{x,y,z\in X}{\text{not all equal}}}
\frac{\rho(x,y)}
{\max\left\{\rho(x,z),\,\rho(z,y)\right\}}
\in[1,\infty)
\end{align}
and
\begin{align}
\label{Crhotilde}
\widetilde{C}_\rho:=\sup_{\genfrac{}{}{0pt}{}{x,y\in X}{x\ne y}}
\frac{\rho(y,x)}{\rho(x,y)}\in[1,\infty).
\end{align}

\begin{remark}\label{remark2.2}
\begin{enumerate}
\item
Let $X:=\mathbb{R}^n$ with $n\in\mathbb{N}$ and
let $\rho:=|\cdot-\cdot|$ be the standard
Euclidean distance on $\mathbb{R}^n$.
Then we easily conclude that $C_\rho=2$ and $\widetilde{C}_\rho=1$.

\item
$\rho$ is said to
be \emph{symmetric}\index[words]{symmetric}
if $\widetilde{C}_\rho=1$.

\item Recall that a \emph{metric}\footnote{A function
$d:X\times X\to\mathbb{R}_+$
is call a \emph{metric} if, for any $x,y,z\in X$,
$d$ satisfies the following three conditions:
(i) $d(x,y)=0$ if and only if $x=y$,
(ii) $d(x,y)=d(y,x)$, (iii) $d(x,y)\leq
d(x,z)+d(z,y)$.} $d$\index[words]{metric} on $X$ is referred to as an
\emph{ultrametric}\index[words]{ultrametric}
provided that, for any
$x,y,z\in X$,
$$
d(x,y)\leq\max\left\{d(x,z),\,d(z,y)\right\}
$$
holds. Hence,
if $C_\rho=\widetilde{C}_\rho=1$, then a quasi-ultrametric
$\rho$ in this case becomes an
\emph{ultrametric} on $X$ and, from this perspective, it is natural
to refer to the inequality in
Definition~\ref{D2.1}(iii) as the \emph{quasi-ultrametric inequality}.
\index[words]{quasi-ultrametric inequality}
For more properties of ultrametrics,
we refer to \cite[Proposition~1.3.28]{cmn2019}.

\item
By the obvious inequalities that
$\frac{a+b}{2}\leq\max\{a,\,b\}\leq a+b$
for any $a,b\in\mathbb{R}_+$, we easily
conclude that the quasi-ultrametric inequality in
Definition~\ref{D2.1}(iii) is equivalent to
the \emph{quasi-triangle inequality};\index[words]{quasi-triangle inequality}
that is, for any $x,y,z\in X$,
\begin{align}\label{2234}
\rho(x,y)\leq A_\rho\left[\rho(x,z)+\rho(z,y)\right],
\end{align}
where $A_\rho\in[1,\infty)$ is
independent of $x$, $y$, and $z$.
Although the quasi-ultrametric inequality in Definition~\ref{D2.1}(iii)
and the quasi-triangle inequality \eqref{2234}
are equivalent, the constant
$A_\rho$ fails to fully capture the various
subtleties within the very class of metrics and, as we
will demonstrate throughout this monograph,
it is the constant $C_\rho$,
rather than $A_\rho$, that will prove to be
of utmost importance with regards to the sharpness
of the main results in this work.
\end{enumerate}
\end{remark}

In what follows, we denote by
$\mathfrak{Q}(X)$\index[symbols]{Q@$\mathfrak{Q}(X)$} the set of
all quasi-ultrametrics on $X$.
A pair $(X,\rho)$ is called a
\emph{quasi-ultrametric space}\index[words]{quasi-ultrametric space}
if $X$ is a set of cardinality at least 2 and
$\rho\in\mathfrak{Q}(X)$. Recall that,
for any given nonempty set $\mathscr{X}$,
two functions
$f,g:\mathscr{X}\to[0,\infty]$
are said to be \emph{equivalent}, denoted by $f\sim g$\index[symbols]{S@$\sim$},
if there exists a constant $C\in[1,\infty)$
such that $C^{-1}f\leq g\leq Cf$
pointwise on $\mathscr{X}$.
With this notation, it is easy to see that, if $\rho\in\mathfrak{Q}(X)$
and $\varrho:X\times X\to\mathbb{R}_+$ is
any function such that $\varrho\sim\rho$,
then $\varrho\in\mathfrak{Q}(X)$, as well.
Thus, there is a natural equivalence relation $\sim$
on $\mathfrak{Q}(X)$ and we call each
equivalence class
$\mathbf{q}\in\mathfrak{Q}(X)/\sim$\index[symbols]{Q@$\mathbf{q}$}
a \emph{quasi-ultrametric space
structure}\index[words]{quasi-ultrametric space structure}
on $X$.
If $X$ is a set of cardinality at least 2 and
$\mathbf{q}\in\mathfrak{Q}(X)/\sim$,
then the pair $(X,\mathbf{q})$\index[symbols]{X@$(X,\mathbf{q})$}
shall also be referred to as a
\emph{quasi-ultrametric space}\index[words]{quasi-ultrametric space}.

Next, we make some conventions on the notation
used throughout this monograph.
Let $(X,\mathbf{q})$ be a
quasi-ultrametric space
and $\rho\in\mathbf{q}$.
The \emph{$\rho$-ball}\index[symbols]{R@$\rho$-ball}
$B_\rho(x,r)$\index[symbols]{B@$B_\rho(x,r)$}
centered at $x\in X$
with radius $r\in(0,\infty)$ is naturally
defined by setting
\begin{align*}
B_\rho(x,r):=\left\{y\in X:\rho(x,y)<r\right\}.
\end{align*}
Since we do not assume that $\rho$ is symmetric, the reader is alerted
to the fact that care must be taken when discussing the
membership of a point to any $\rho$-ball. For any $\lambda\in(0,\infty)$
and $\rho$-ball $B_\rho\subset X$,
we denote by $\lambda B_\rho$ the ball that
is concentric with $B_\rho$ and radius $\lambda$
times that of $B_\rho$. That is,
if $B_\rho=B_\rho(x,r)$ for some $x\in X$ and $r\in(0,\infty)$, then
$\lambda B_\rho=B_\rho(x,\lambda r)$.
The quasi-ultrametric $\rho$ naturally induces a topology on $X$,
denoted by $\tau_\rho$,\index[symbols]{T@$\tau_\rho$}
where a set $\mathcal{O}\subset X$ is open in $\tau_\rho$
if and only if, for any $x\in\mathcal{O}$,
there exists $r\in(0,\infty)$ such that $B_\rho(x,r)\subset\mathcal{O}$.
The reader is also alerted to the fact that,
in contrast to the setting of metric spaces,
there is no guarantee that $\rho$-balls belong to $\tau_\rho$.
Note that, for any given $\rho,\rho'\in\mathbf{q}$,
there exists a positive constant $c\in[1,\infty)$,
depending only on $\rho$ and $\rho'$,
such that, for any $x\in X$ and $r\in(0,\infty)$,
\begin{align}\label{1557-xx}
B_{\rho'}\left(x,c^{-1}r\right)\subset B_\rho(x,r)\subset B_{\rho'}(x,cr),
\end{align}
from which it follows that $\tau_{\rho'}=\tau_\rho$.
Therefore, it is meaningful to define
the topology on $(X,\mathbf{q})$ as
$\tau_{\mathbf{q}}:=\tau_\rho$\index[symbols]{T@$\tau_{\mathbf{q}}$}.
Note that, for any $\alpha\in(0,\infty)$,
we have $\rho^\alpha\in\mathfrak{Q}(X)$, where, for any $x,y\in X$,
$\rho^\alpha(x,y):=[\rho(x,y)]^\alpha$. Hence,
\begin{align}\label{1557}
\mathbf{q}^\alpha:=\left\{\rho^\alpha:\rho\in\mathbf{q}\right\}
\end{align}
is a quasi-ultrametric space structure on $X$,
and it follows from definitions
that, for any $x\in X$ and $r\in(0,\infty)$,
\begin{align}\label{935}
B_\rho(x,r)=B_{\rho^\alpha}(x,r^\alpha),
\end{align}
which further implies that
$\tau_{\mathbf{q}^\alpha}=\tau_{\rho^\alpha}
=\tau_\rho=\tau_{\mathbf{q}}$. This
is remarkable since $\rho$ and $\rho^\alpha$
are not pointwise equivalent, in general.
For any nonempty set $E\subset X$,
the \emph{diameter $\mathrm{diam}_{\rho}(E)$
of $E$ (with respect to $\rho$)} is
defined by setting
\index[symbols]{D@$\mathrm{diam}_\rho(E)$}
\begin{align}
\label{diam}
\mathrm{diam}_{\rho}(E):=
\sup\left\{\rho(x,y):x,y\in E\right\}.
\end{align}
Clearly, for any given $\rho,\rho'\in\mathbf{q}$,
we have
$$
\text{diam}_{\rho'}(E)\sim\text{diam}_{\rho}(E)
$$
with the positive equivalence constants depending only on $\rho$ and $\rho'$.
Therefore, we sometimes simply write $\mathrm{diam\,}(E)$
in place of $\text{diam}_{\rho}(E)$ if the particular
choice of quasi-ultrametric is not important.
A set
$E\subset X$ is said to be
\emph{bounded}\index[words]{bounded}
if $E$ is contained in a $\rho$-ball
for some (hence all) $\rho\in\mathbf{q}$.
In other words, a set $E\subset X$ is bounded,
relative to the quasi-ultrametric space structure
$\mathbf{q}$ on $X$, \emph{if and only if} $\text{diam}_{\rho}(E)<\infty$
for some (hence all) $\rho\in \mathbf{q}$.
For any nonempty sets $F_0,F_1\subset X$,
the \emph{$\rho$-distance between $F_0$ and $F_1$}
is defined by setting
\index[symbols]{D@$\mathrm{dist}_\rho(F_0,F_1)$}
\begin{align*}
\mathrm{dist}_\rho(F_0,F_1)
:=\inf\left\{\rho(x,y):x\in F_0\text{ and }y\in F_1\right\}.
\end{align*}

Now, we recall the definition of the lower
smoothness index of a quasi-ultrametric space,
which was originally introduced by Mitrea et al.
in \cite[Definition 4.26]{MiMiMiMo13}.
This index plays a pivotal role in identifying the optimal range of $p$
for which there exists a rich theory of Hardy spaces $H^p(X)$
(see, for example, \cite{AlMi2015}).
More specifically, the index permits us to
determine just how far $p$ can be below 1
while still having a maximal function characterization
of the atomic Hardy spaces introduced in \cite[p.\,592]{Cowe77}.
In this monograph, the notion of indices is going to play a fundamental
role in the formulation and proofs of many of our main results.
For more information about this
lower smoothness index,
we refer to \cite[Section 4.7]{MiMiMiMo13};
see also \cite[Section 2.5]{AlMi2015}.

\begin{definition}
Let $(X,\mathbf{q})$ be a quasi-ultrametric space and $\rho\in\mathbf{q}$.
The \emph{lower smoothness index}\index[words]{lower smoothness index}
${\mathrm{ind\,}}(X,\mathbf{q})$\index[symbols]{I@${\mathrm{ind\,}}(X,\mathbf{q})$}
[resp. ${\mathrm{ind\,}}(X,\rho)$]\index[symbols]{I@${\mathrm{ind\,}}(X,\rho)$}
of $(X,\mathbf{q})$ [resp. $(X,\rho)$]
is defined by setting
\begin{align}\label{index}
{\mathrm{ind\,}}(X,\mathbf{q})
:=\sup\left\{(\log_2C_\varrho)^{-1}:
\varrho\in\mathbf{q}\right\}
\end{align}
\begin{align}\label{1553}
\left[\text{resp. }{\mathrm{ind\,}}(X,\rho):=
\sup_{\varrho\sim\rho}\,(\log_2C_\varrho)^{-1}\right],
\end{align}
where, for any $\varrho\in\mathbf{q}$,
$C_\varrho$ is the same as in \eqref{Crho}.
\end{definition}

\begin{remark}
\label{rmk390}
\begin{enumerate}
\item
By the fact that $C_\varrho\in[1,\infty)$
for any $\varrho\in\mathbf{q}$,
we find that ${\mathrm{ind\,}}(X,\mathbf{q})\in(0,\infty]$.
In particular, if there exists an ultrametric
$\rho\in\mathbf{q}$, then
${\mathrm{ind\,}}(X,\mathbf{q})=\infty$;
if there exists a metric
$\rho\in\mathbf{q}$, then
${\mathrm{ind\,}}(X,\mathbf{q})\in[1,\infty]$.
If the cardinality of $X$ is finite, then, for any $\rho\in\mathfrak{Q}(X)$,
$\mathrm{ind\,}(X,\rho)=\infty$.
\item
Given an arbitrary quasi-ultrametric space $(X,\mathbf{q})$, the supremum
in \eqref{index} which defines
${\mathrm{ind\,}}(X,\mathbf{q})$  may or may not be attained;
see \cite[p.\,65, Example 1]{AlMi2015}
for the case that ${\mathrm{ind\,}}(X,\mathbf{q})$
is attained and \cite[Theorem 1.1]{AlMi2015}
for the case that ${\mathrm{ind\,}}(X,\mathbf{q})$
is not attained.
\item
For the remainder of this monograph, we write
$0<\varepsilon\preceq{\mathrm{ind\,}}(X,\mathbf{q})$\index[symbols]{P@$\preceq$}
to mean that $\varepsilon\in(0,\infty)$ and
$\varepsilon\leq{\mathrm{ind\,}}(X,\mathbf{q})$, where the equality
$\varepsilon={\mathrm{ind\,}}(X,\mathbf{q})$ is
only \emph{permissible} when the supremum
in \eqref{index} is attained.
\item
By \cite[Theorem~4.43]{MiMiMiMo13},
we find that, for any quasi-ultrametric space $(X,\mathbf{q})$, one has
\begin{align*}
\mathrm{ind\,}(X,\mathbf{q})=\frac{1}{\mathrm{Cv\,}(X,\mathbf{q})},
\end{align*}
where $\mathrm{Cv\,}(X,\mathbf{q})$ is Assouad's \emph{convexity index},
introduced by Assouad in \cite[Definition~3]{a1979},
of $(X,\mathbf{q})$, which
is defined by setting
\begin{align*}
\mathrm{Cv\,}(X,\mathbf{q}):=\inf\Big\{&p\in(0,\infty):
\text{there exists }\rho\in\mathbf{q}\text{ such that }\\
&\rho^\frac{1}{p}
\text{ is equivalent to a metric on }X\Big\}.
\end{align*}
\item By \cite[Corollary~4.46]{MiMiMiMo13},
we have the following characterization of
$\mathrm{ind\,}(X,\rho)$,
\begin{align*}
\mathrm{ind\,}(X,\rho)
&=\sup\left\{\alpha\in(0,\infty):
\inf_{x,y\in X,\,x\neq y}\frac{\rho_\alpha(x,y)}{\rho(x,y)}>0\right\}\\
&=\sup\left\{\alpha\in(0,\infty):\rho_\alpha\sim\rho\right\},
\end{align*}
where $\rho_\alpha$, for any $\alpha\in(0,\infty)$,
is defined in Definition~\ref{treatment}(i)
below.
In particular, for any $0<\alpha\preceq\mathrm{ind\,}(X,\rho)$,
$\rho_\alpha\sim\rho$.
\item From \cite[Corollary~4.60]{MiMiMiMo13},
it follows that, for any
quasi-ultrametric space $(X,\mathbf{q})$, it holds that
\begin{align}\label{twj-348}
\mathrm{ind\,}(X,\mathbf{q})\leq\mathrm{ind}_H(X,\mathbf{q}),
\end{align}
where\index[symbols]{I@$\mathrm{ind}_H(X,\mathbf{q})$}
$\mathrm{ind}_H(X,\mathbf{q})$ is
the \emph{H\"older index}
of $(X,\mathbf{q})$ that is defined by setting\index[words]{H\"older index}
\begin{align}\label{1550}
\mathrm{ind}_H(X,\mathbf{q}):=\sup\left\{\beta\in(0,\infty):
\dot{C}^\beta(X,\mathbf{q})\neq\mathbb{C}\right\}\subset(0,\infty],
\end{align}
and $\dot{C}^\beta(X,\mathbf{q})$ with $\beta\in(0,\infty)$
is the (homogeneous) H\"older space defined in Definition~\ref{1454} below.
By definition, for any $\beta\in(0,\infty)$
satisfying $\beta>\mathrm{ind}_H(X,\mathbf{q})$,
we have $\dot{C}^\beta(X,\mathbf{q})$
only contains constant functions.
On the other hand, as a consequence of Corollary~\ref{Lem3.3},
if $0<\beta\preceq\mathrm{ind\,}(X,\mathbf{q})$,
then $\dot{C}^\beta(X,\mathbf{q})$ is nontrivial,
in the sense that $\dot{C}^\beta(X,\mathbf{q})$
contains nonconstant functions;
for more studies of the H\"older index $\mathrm{ind}_H(X,\mathbf{q})$,
we refer to \cite[Section~2.5]{AlMi2015} and
\cite[Section~4.7]{MiMiMiMo13}.
\end{enumerate}
\end{remark}

The following lemma is a part of \cite[Proposition 2.20]{AlMi2015}.

\begin{lemma}\label{11111}
Let $(X,\mathbf{q})$ be
a quasi-ultrametric space.
Then, for any $\alpha\in(0,\infty)$,
$$
{\mathrm{ind\,}}(X,\mathbf{q}^\alpha)
=\frac{1}{\alpha}{\mathrm{ind\,}}(X,\mathbf{q})
\ \ \text{and}\ \
\mathrm{ind}_H(X,\mathbf{q}^\alpha)
=\frac{1}{\alpha}\mathrm{ind}_H(X,\mathbf{q}),
$$
where
$\mathbf{q}^\alpha$ is the same as in \eqref{1557}
and
${\mathrm{ind\,}}(X,\mathbf{q}^\alpha)$
and $\mathrm{ind}_H(X,\mathbf{q}^\alpha)$
are the same as, respectively, in \eqref{index} and \eqref{1550}
with $\mathbf{q}$ replaced by $\mathbf{q}^\alpha$.
\end{lemma}

Next, we provide a few concrete examples to show how the value of the
lower smoothness varies across different
quasi-ultrametric spaces.
For more examples, we refer to
\cite[pp.\,65--69]{AlMi2015}.

\begin{example}
Let $n\in\mathbb{N}$, $\alpha\in(0,\infty)$,
and $|\cdot-\cdot|$ be the standard Euclidean distance on
$\mathbb{R}^n$. Then
$$
{\mathrm{ind\,}}(\mathbb{R}^n,|\cdot-\cdot|^\alpha)=\frac{1}{\alpha}
={\mathrm{ind\,}}([0,1]^n,|\cdot-\cdot|^\alpha),
$$
and
$$
\mathrm{ind}_H(\mathbb{R}^n,|\cdot-\cdot|^\alpha)=\frac{1}{\alpha}
=\mathrm{ind}_H([0,1]^n,|\cdot-\cdot|^\alpha),
$$
which is a direct corollary of both
Remark~\ref{remark2.2}(i) and Lemma~\ref{11111}.
In particular,
$$
{\mathrm{ind\,}}(\mathbb{R}^n,|\cdot-\cdot|)=1
={\mathrm{ind\,}}([0,1]^n,|\cdot-\cdot|).
$$
and
$$
\mathrm{ind}_H(\mathbb{R}^n,|\cdot-\cdot|)=1
=\mathrm{ind}_H([0,1]^n,|\cdot-\cdot|).
$$
More generally, if $(V,\|\cdot\|)$ is a nontrivial
normed vector space and if $\mathbf{q}$
stands for the quasi-ultrametric space
structure induced by the norm $\|\cdot\|$,
then, for any convex subset $U$ of $V$ of cardinality at least 2,
$$
{\mathrm{ind\,}}(V,\mathbf{q})=1=\mathrm{ind}_H(V,\mathbf{q}).
$$
\end{example}

\begin{example}
For any given $p,\alpha\in(0,\infty)$,
\begin{align*}
{\mathrm{ind\,}}(L^p(\mathbb{R}^n),\|\cdot-\cdot\|_{L^p(\mathbb{R}^n)}^\alpha)
=\frac{\min\,\{1,\,p\}}{\alpha}
={\mathrm{ind\,}}(\ell^p(\mathbb{N}),\|\cdot-\cdot\|_{\ell^p(\mathbb{N})}^\alpha)
\end{align*}
and
$$
\mathrm{ind}_H(L^p(\mathbb{R}^n),\|\cdot-\cdot\|_{L^p(\mathbb{R}^n)}^\alpha)
=\frac{\min\,\{1,\,p\}}{\alpha}
=\mathrm{ind}_H(\ell^p(\mathbb{N}),\|\cdot-\cdot\|_{\ell^p(\mathbb{N})}^\alpha).
$$
To verify the above first equality,
it suffices to show
\begin{align}\label{1400}
{\mathrm{ind\,}}(L^p(\mathbb{R}^n),\|\cdot-\cdot\|_{L^p(\mathbb{R}^n)})
=\min\,\{1,\,p\}
\end{align}
due to Lemma~\ref{11111}.
To this end,
we consider the following two
cases for $p$.

\emph{Case (1)} $p\in[1,\infty)$.
In this case, by the triangle inequality of $\|\cdot\|_{L^p(\mathbb{R}^n)}$,
we find that, for any $f,g,h\in L^p(\mathbb{R}^n)$,
\begin{align*}
\|f-g\|_{L^p(\mathbb{R}^n)}
&\leq\|f-h\|_{L^p(\mathbb{R}^n)}+\|h-g\|_{L^p(\mathbb{R}^n)}\\
&\leq2\max\left\{\|f-h\|_{L^p(\mathbb{R}^n)},\,\|h-g\|_{L^p(\mathbb{R}^n)}\right\}.
\end{align*}
Moreover, if we take $g:=-f\in L^p(\mathbb{R}^n)$,
then we obtain the following equality
\begin{align*}
\|f-g\|_{L^p(\mathbb{R}^n)}=2\|f\|_{L^p(\mathbb{R}^n)}=
2\max\left\{\|f-0\|_{L^p(\mathbb{R}^n)},\,\|0-g\|_{L^p(\mathbb{R}^n)}\right\}.
\end{align*}
Thus, $C_{\|\cdot-\cdot\|_{L^p(\mathbb{R}^n)}}=2$ and \eqref{1400} holds.

\emph{Case (2)} $p\in(0,1)$. In this case, using Lemma~\ref{11111},
we have
\begin{align*}
{\mathrm{ind\,}}(L^p(\mathbb{R}^n),\|\cdot-\cdot\|_{L^p(\mathbb{R}^n)})
=p\,{\mathrm{ind\,}}(L^p(\mathbb{R}^n),\|\cdot-\cdot\|_{L^p(\mathbb{R}^n)}^p)
\end{align*}
From this, an argument similar to that used in the proof in Case (1),
and the triangle inequality of $\|\cdot\|_{L^p(\mathbb{R}^n)}^p$,
we further deduce that
\begin{align*}
{\mathrm{ind\,}}(L^p(\mathbb{R}^n),\|\cdot-\cdot\|_{L^p(\mathbb{R}^n)}^p)=1
\end{align*}
and hence \eqref{1400} holds.
\end{example}

We also give an example $X$ of an ultrametric space
on which we have ${\mathrm{ind\,}}(X,\rho)=\infty$
and $\mathrm{ind}_H(X,\rho)=\infty$.

\begin{example}
Let $X:=[0,1)$.
For any $x,y\in X$, if $x=y$ then
let $\rho(x,y):=0$, and if $x\neq y$ then let
$\rho(x,y)$ be the length of the smallest dyadic interval such that
\begin{align*}
\left[\frac{k}{2^n},\frac{k+1}{2^n}\right)\supset\{x,y\},
\end{align*}
where $n\in\mathbb{Z}_+$ and $k\in\mathbb{Z_+}\cap[0,2^n-1]$.
Then $(X,\rho)$ is an ultrametric space with
\begin{align}\label{125}
C_\rho=1
\ \ \text{and}\ \
{\mathrm{ind\,}}(X,\rho)=\infty=\mathrm{ind}_H(X,\rho).
\end{align}
Indeed, by the definition of $\rho$, we find that,
for any $x,y\in X$ with $x\neq y$, there exist
$n\in\mathbb{Z}_+$ and $k\in\mathbb{Z_+}\cap[0,2^n-1]$ such that
\begin{align*}
\frac{k}{2^n}\leq x<\frac{2k+1}{2^{n+1}}\leq y<\frac{k+1}{2^{n}},
\end{align*}
where the second and the third inequalities follow from the minimality of
the dyadic interval $[\frac{k}{2^n},\frac{k+1}{2^n})$.
Then
\begin{align*}
\frac{\rho(x,y)}{\max\{\rho(x,z),\,\rho(z,y)\}}
=\frac{\rho(x,y)}{\rho(x,z)}=1\ \text{if}\ z\in\left[\frac{2k+1}{2^{n+1}},y\right],
\end{align*}
\begin{align*}
\frac{\rho(x,y)}{\max\{\rho(x,z),\,\rho(z,y)\}}
=\frac{\rho(x,y)}{\rho(z,y)}=1\ \text{if}\ z\in\left[x,\frac{2k+1}{2^{n+1}}\right),
\end{align*}
and
\begin{align*}
\frac{\rho(x,y)}{\max\{\rho(x,z),\,\rho(z,y)\}}
\leq1\ \text{if}\ z\in[0,x)\cup(y,1),
\end{align*}
which, combined with \eqref{1553},
further implies that \eqref{125} holds.
\end{example}

The following example highlights the fact that the
inequality in \eqref{twj-348} can be strict.
The reader is referred to
Comment~4.38 and Remark~4.49 in \cite{MiMiMiMo13} for more details.

\begin{example}
Let $a$, $b$, $c$, and $d$ be four real numbers with the property that $a<b<c<d$.
Then ${\mathrm{ind\,}}([a,b]\cup[c,d],|\cdot-\cdot|)=1$
and $\mathrm{ind}_H([a,b]\cup[c,d],|\cdot-\cdot|)=\infty$.
\end{example}

Recall that a topological space $(X,\tau)$ is said to be \emph{metrizable}
if there exists a metric $d$ on $X$ satisfying that the topology
induced by $d$ on $X$ coincides with $\tau$.
There are several metrization results in various branches of
mathematics, such as the Mac\'ias--Segovia
metrization theorem (see Theorem~\ref{1000} below)
for quasi-metric spaces (harmonic analysis),
the Alexandroff--Urysohn metrization
theorem (see, for instance, \cite[Theorem~1.1]{MiMiMiMo13}
and \cite[Theorem 34.1]{m2000})
for uniform spaces (topology),
and the Aoki--Rolewicz theorem
(see, for instance, \cite[Theorem 5]{a42} and \cite{r1957})
for quasi-normed vector spaces (functional analysis).
Here, we focus on the contribution of Mac\'ias and Segovia \cite{MaSe79i}
which has been hailed as a \emph{cornerstone of the analysis on quasi-metric spaces}.
We record it as follows; see \cite[Theorem 2]{MaSe79i}.

\begin{theorem}[\bf Mac\'ias--Segovia]\label{1000}
Let $(X,\rho)$ be a quasi-metric space; that is, $X$ is a nonempty set
and $\rho:X\times X\to\mathbb{R}_+$ is a
quasi-metric [a function satisfies
Definition~\ref{D2.1}(i), Definition~\ref{D2.1}(ii) with $C_0=1$,
and \eqref{2234}].
Then there exists another quasi-metric $\varrho$ on $X$
satisfying that
\begin{enumerate}
\item[\rm (i)]
$\rho\sim\varrho$;
\item[\rm (ii)]
$\varrho^\alpha$ is a metric on $X$ that induces the same topology on $X$
as $\rho$, where
\begin{align*}
\alpha:=\frac{1}{\log_2[A_\rho(2A_\rho+1)]}\in(0,1).
\end{align*}
In particular, this topology on $X$ is metrizable;
\item[\rm (iii)]
$\varrho$ satisfies the following H\"older-type
regularity condition of order $\alpha$; that is,
for any $x,y,z\in X$,
\begin{align*}
\left|\varrho(x,z)-\varrho(y,z)\right|
\leq\frac{1}{\alpha}\left[\varrho(x,y)\right]^\alpha
\max\left\{\left[\varrho(x,z)\right]^{1-\alpha},\,
\left[\varrho(y,z)\right]^{1-\alpha}\right\}.
\end{align*}
\end{enumerate}
\end{theorem}

Until now, Theorem~\ref{1000}
has still played a key role in analysis on spaces of homogeneous type
because the natural setting for analysis in this
context is the setting of quasi-metric spaces.
Via Theorem~\ref{1000},
a considerable portion of the well-established results in the area of analysis
on metric spaces may be quite successfully developed in the more general setting
of quasi-metric spaces.
Indeed, quasi-metric spaces are metrizable,
but it is a rather subtle matter to associate metrics,
inducing the same topology,
in a way that brings out the quantitative features of the quasi-metric space
in question in an optimal manner.
The seminal work on this topic done by R. Mac\'ias and C. Segovia
has been very influential in the study of spaces of homogeneous type; see,
for example, the monographs by Christ \cite{Christ1990},
by Stein \cite{Pag}, by Triebel \cite{Triebel2006},
by Heinonen \cite{H2001}, by Han and Sawyer \cite{HS},
by David and Semmes \cite{DS}, and by Deng and Han \cite{DH}.

Next, we recall a sharp metrization result of Mitrea et al. \cite{MiMiMiMo13},
which is a sharpened version of Theorem~\ref{1000}.
To this end, we begin with recalling some
notions related to the quasi-ultrametric $\rho$,
which is useful to construct the metric we needed;
see, for instance,
\cite[pp.\,36--37]{AlMi2015} and \cite[Theorem 3.46]{MiMiMiMo13}.

\begin{definition}
\label{treatment}
Let $X$ be a nonempty set and $\rho:X\times X\to\mathbb{R}_+$ be
any function.
\begin{enumerate}
\item
Given $\alpha\in(0,\infty]$,
$\rho_\alpha:X\times X\to\mathbb{R}_+$
is defined by setting, for any $x,y\in X$,
when $\alpha\in(0,\infty)$,\index[symbols]{R@$\rho_\alpha$}
\begin{align}
\label{2.1}
\rho_\alpha(x,y):=
&\inf\Bigg\{\bigg[\sum_{j=1}^{N}
\Big\{\rho(\xi_j,\xi_{j+1})\Big\}^\alpha
\bigg]^{\frac{1}{\alpha}}:
\text{there exist an }N\in\mathbb{N}\
\text{and}\nonumber\\
&\qquad
\{\xi_j\}_{j=1}^{N+1}\text{ in }X\
\text{such that }
\xi_1=x\text{ and }\xi_{N+1}=y\Bigg\}
\end{align}
and, when $\alpha=\infty$,
\begin{align}\label{2.2}
\rho_\infty\index[symbols]{R@$\rho_\infty$}
\left(x,y\right):=
&\inf\Bigg\{\max_{1\leq j\leq N}
\rho(\xi_j,\xi_{j+1}):
\text{there exist an }N\in\mathbb{N}\ \text{and}
\nonumber\\
&\qquad
\{\xi_j\}_{j=1}^{N+1}\text{ in }X\
\text{such that }
\xi_1=x\text{ and }\xi_{N+1}=y\Bigg\}
\end{align}
with the infima in \eqref{2.1}
and \eqref{2.2} taken over
all $N\in\mathbb{N}$ and $\{\xi_j\}_{j=1}^{N+1}$ in $X$.
\item
The \emph{symmetrization}\index[words]{symmetrization}
$\rho_{\mathrm{sym}}:X\times X\to\mathbb{R}_+$ of $\rho$
is defined by setting, for any $x,y\in X$,
$$
\rho_{\text{sym}}(x,y):=
\max\left\{\rho(x,y),\,\rho(y,x)\right\}.\index[symbols]{R@$\rho_{\textup{sym}}$}
$$
\end{enumerate}
\end{definition}

\begin{remark}
\label{byw-45}
Let the notation be the same as in Definition~\ref{treatment}(i).
Then
\begin{enumerate}
\item
$\rho_\alpha$ is well defined and $\rho_\alpha\leq\rho$ pointwise on $X\times X$.

\item
For any given finite $\beta\in(0,\alpha]$,
$\rho_\alpha$ satisfies the triangle inequality; that is,
for any $x,y,z\in X$,
\begin{align*}
\left[\rho_\alpha(x,y)\right]^\beta\leq
\left[\rho_\alpha(x,z)\right]^\beta
+\left[\rho_\alpha(z,y)\right]^\beta.
\end{align*}

\item
For any $x,y,z\in X$,
\begin{align*}
\rho_\alpha(x,y)\leq2^\frac{1}{\alpha}
\max\left\{\rho_\alpha(x,z),\,\rho_\alpha(z,y)\right\}.
\end{align*}
\end{enumerate}
We refer to \cite[Theorem~3.46(1)]{MiMiMiMo13}
for more properties of $\rho_\alpha$.
As we will see, Remark~\ref{byw-45}(ii) gives a direct approach to
construct a metric via any quasi-ultrametric.
Moreover, Remark~\ref{byw-45}(iii) with $\alpha:=\infty$
allows us to construct an ultrametric via a
quasi-ultrametric.
\end{remark}

\begin{remark}
\label{1.1.7}
It is easy to see that $\rho_{\text{sym}}$ is symmetric;
that is, for any $x,y\in X$,
\begin{align*}
\rho_{\text{sym}}(x,y)=\rho_{\text{sym}}(y,x).
\end{align*}
Furthermore, if $\rho\in\mathfrak{Q}(X)$, then
$\rho_{\text{sym}}\in\mathfrak{Q}(X)$,
$C_{\rho_{\text{sym}}}\leq C_{\rho}$,
$\widetilde{C}_{\rho_{\text{sym}}}=1$,
and
\begin{align*}
\rho\leq\rho_{\text{sym}}
\leq\widetilde{C}_{\rho}\,\rho
\end{align*}
pointwise on $X\times X$,
where $C_{\rho}$ and $\widetilde{C}_{\rho}$ are
the same as, respectively,
in \eqref{Crho} and \eqref{Crhotilde}.
\end{remark}

The following sharp metrization lemma is a part of
\cite[Theorem 3.46]{MiMiMiMo13}
(see also \cite[Theorem 2.1]{AlMi2015}),
which plays a key role in this monograph.

\begin{lemma}\label{L2.3}
Let $(X,\mathbf{q})$ be a quasi-ultrametric
space and $\rho\in\mathbf{q}$.
In this context, set
$$
\rho_{\#}:=(\rho_{\text{sym}})_{\alpha_\rho},\index[symbols]{R@$\rho_{\#}$}
$$
where
\begin{align}\label{2219}
\alpha_\rho:=\frac{1}{\log_2C_\rho}\in(0,\infty]
\end{align}
with $C_\rho\in[1,\infty)$ as in \eqref{Crho}.
Then the following assertions hold.
\begin{enumerate}
\item [\rm (i)]
$\rho_{\#}\in\mathbf{q}$ and satisfies that
$C_{\rho_\#}\leq C_\rho$ and $\widetilde{C}_{\rho_\#}=1$.
If $C_\rho=1$,
then $\rho_\#=(\rho_{\mathrm{sym}})_\infty=\rho_{\mathrm{sym}}$.

\item [\rm (ii)]
For any $x,y\in X$,
$$
(C_\rho)^{-2}\rho(x,y)
\leq\rho_\#(x,y)
\leq\widetilde{C}_\rho\rho(x,y).
$$

\item [\rm (iii)]
For any finite $\alpha\in(0,\alpha_\rho]$,
$\rho_\#$ satisfies the \emph{local H\"older-type regularity
condition of order $\alpha$};\index[words]{local H\"older-type regularity
condition}
that is, for any $x,y,z\in X$,
\begin{equation}\label{2.3}
\left|\rho_\#(x,y)-\rho_\#(x,z)\right|
\leq\frac{1}{\alpha}\left[\rho_\#(y,z)\right]^{\alpha}
\max\left\{\left[\rho_\#(x,y)\right]^{1-\alpha},\,
\left[\rho_\#(x,z)\right]^{1-\alpha}\right\}
\end{equation}
with the understanding that, when $\alpha\in[1,\infty)$, one also
imposes the condition that $x\notin\{y,z\}$. In particular,
the function $\rho_\#:X\times X\to\mathbb{R}_+$ is continuous
in the natural product topology
$\tau_\mathbf{q}\times\tau_\mathbf{q}$ and hence
all $\rho_\#$-balls are open in the topology $\tau_\mathbf{q}$.

\item[\rm (iv)]
For any finite $\alpha\in(0,\alpha_\rho]$,
$(\rho_\#)^{\alpha}$ is a metric on $X$; that is,
$(\rho_\#)^{\alpha}$ is a nondegenerate and symmetric function
on $X\times X$ satisfying, for any $x,y,z\in X$,
\begin{align*}
\left[\rho_\#(x,y)\right]^{\alpha}
\leq\left[\rho_\#(x,z)\right]^{\alpha}
+\left[\rho_\#(z,y)\right]^{\alpha}.
\end{align*}
Moreover, $(\rho_\#)^{\alpha}$ induces the same topology on $X$ as $\rho$
and hence the topological space $(X,\tau_\rho)$ is metrizable.
\end{enumerate}
\end{lemma}

\begin{remark}
\begin{enumerate}
\item[\rm (i)]
By \cite[Theorem~3.46(14)]{MiMiMiMo13},
we find that,
if we replace the quasi-ultrametric condition in Definition~\ref{D2.1}(iii)
by the following condition that
there exist $p\in(0,\infty)$ and $C_2\in[1,\infty)$
such that, for any $x,y,z\in X$,
\begin{align*}
\rho(x,y)\leq C_2\left\{\left[\rho(x,z)
\right]^p+\left[\rho(z,y)\right]^p\right\}^\frac{1}{p},
\end{align*}
then the same conclusion of Lemma~\ref{L2.3} holds if in place of
\eqref{2219} one defines $\alpha_p$ to be
\begin{align*}
\alpha_p:=\frac{p}{1+p\log_2C_2}\in(0,\infty],
\end{align*}
which, when $p=1$, is obviously greater
than $\alpha$ in Theorem~\ref{1000} and hence
Lemma~\ref{L2.3} is indeed an improved version of Theorem~\ref{1000}.

\item[\rm (ii)]
As shown in \cite[p.\,144]{MiMiMiMo13},
it is possible to further refine a large number of basic
results in the area of analysis on quasi-metric spaces, whose
proofs make use of Theorem~\ref{1000} instead of
the sharp metrization theorem \cite[Theorem 3.46]{MiMiMiMo13};
see \cite[Section 4]{MiMiMiMo13} for more information,
such as extensions and separation properties of H\"older functions,
approximations to the identity and Hardy spaces on
Ahlfors-regular quasi-metric spaces,
and Kuratowski's and Fr\'echet's embedding theorems.

\item[\rm (iii)]
It is worth mentioning that
Lemma~\ref{L2.3} gives a concrete formula,
namely $d:=[(\rho_{\mathrm{sym}})_{\alpha_\rho}]^\beta$ for
any given finite $\beta\in(0,\alpha_\rho]$, for the metric $d$ yielding the
same topology as that induced by the original quasi-ultrametric,
which, in particular, makes it transparent how this metric is
related to the original quasi-ultrametric.

\item[\rm (iv)]
The H\"older-type regularity result in Lemma~\ref{L2.3}(iii)
is \emph{sharp} in the following precise sense.
Given $(X,\mathbf{q})$ and $\rho\in\mathbf{q}$,
if there exists $\varrho\in\mathbf{q}$ satisfying that
there exist two constants $\alpha,C\in(0,\infty)$
such that, for any $x,y,z\in X$ (and also $x\notin\{y,\,z\}$ if $\alpha\in[1,\infty)$),
\begin{align*}
\left|\varrho(x,y)-\varrho(x,z)\right|
\leq C\left[\varrho(y,z)\right]^{\alpha}
\max\left\{\left[\varrho(x,y)\right]^{1-\alpha},\,
\left[\varrho(x,z)\right]^{1-\alpha}\right\},
\end{align*}
then necessarily
\begin{align*}
\alpha\leq\frac{1}{\log_2C_\rho};
\end{align*}
see also \cite[Theorem~3.46(13)]{MiMiMiMo13}.

\item[\rm (v)]
In Lemma~\ref{L2.3}, if $\alpha_\rho\in[1,\infty]$,
then, using Lemma~\ref{L2.3}(iv) and taking $\alpha:=1$,
it is easy to see that $\rho_\#$ is actually a metric on $X$.
\end{enumerate}
\end{remark}

In light of Lemma~\ref{L2.3}, we record the following definition
for a regularized quasi-ultrametric;
see \cite[Theorem~3.26]{MiMiMiMo13} and \cite[Theorem~2.1]{AlMi2015}.

\begin{definition}\label{D2.3}
Let $(X,\mathbf{q})$ be a
quasi-ultrametric space, $\rho\in\mathbf{q}$,
and $\alpha_\rho:=(\log_2C_\rho)^{-1}\in(0,\infty]$,
where $C_\rho$ is the same as in \eqref{Crho}. The quasi-ultrametric
$\rho_{\#}:=(\rho_{\text{sym}})_{\alpha_\rho}$\index[symbols]{R@$\rho_{\#}$}
is called the
\emph{regularized quasi-ultrametric (of $\rho$)}.\index[words]{regularized quasi-ultrametric}
More generally, $\varrho\in\mathbf{q}$ is referred to as
\emph{a regularized quasi-ultrametric}
if $\varrho$ is the regularized quasi-ultrametric
of some $\varrho'\in\mathbf{q}$.
\end{definition}

\begin{remark}\label{RmkD2.3}
Let the notation be as in Definition~\ref{D2.3}.
\begin{enumerate}
\item[\rm(i)]
If $\rho$ is a metric on $X$ and $\alpha_\rho=1$ (namely $C_\rho=2$),
then $\rho_{\#}:=(\rho_{\text{sym}})_{\alpha_\rho}=\rho$. In particular,
$\rho_{\#}=\rho$ if $\rho$ is the standard Euclidean distance on $\mathbb{R}^n$.
Moreover, if $\rho$ is a regularized quasi-ultrametric then $\rho_{\#}=\rho$.

\item[\rm(ii)]
By Lemma~\ref{L2.3}(i), we find that the regularized quasi-ultrametric of
any quasi-ultrametric
satisfying Definition~\ref{D2.1}(iii)
with $C_1:=1$ is precisely its symmetrization.
In particular,
the regularized quasi-ultrametric of any ultrametric
is exactly itself.

\item[\rm(iii)]
From Remark~\ref{rmk390}(iii), \eqref{index}, and Lemma~\ref{L2.3},
we infer that,
for any $0<\varepsilon\preceq{\mathrm{ind\,}}(X,\mathbf{q})$,
there exists a regularized quasi-ultrametric $\rho\in\mathbf{q}$ satisfying
$\varepsilon\leq(\log_2C_\rho)^{-1}$.

\item[\rm(iv)]
The advantage of the regularized quasi-ultrametric
is that it guarantees the measurability
of certain functions (which are defined via
such a regularized quasi-ultrametric),
such as the Hardy--Littlewood maximal function (see Lemma~\ref{M}),
the Lusin area function (see Lemma~\ref{1955}),
and the non-tangential maximal function (see Proposition~\ref{measurability}).
\end{enumerate}
\end{remark}

In order to extend certain classical results
of function spaces over Euclidean spaces to more general settings,
we need to introduce some measure theoretic notions, beginning with
the definition of quasi-ultrametric
spaces of homogeneous type in the sense of
Coifman and Weiss \cite{CW1,Cowe77}.

\begin{definition}\label{Homogenous}
A triplet $(X,\rho,\mu)$ is called a
\textit{quasi-ultrametric space of homogeneous type}
\index[words]{quasi-ultrametric space of homogeneous type}
if $(X,\rho)$ is a quasi-ultrametric space and
$\mu$ is a nonnegative measure on $X$ that satisfies
the following \emph{doubling condition}\index[words]{doubling condition}
(with respect to $\rho$):\ all $\rho$-balls
are $\mu$-measurable and there exists a constant
$C_\mu\in[1,\infty)$\index[symbols]{C@$C_\mu$} such that,
for any $x\in X$ and $r\in(0,\infty)$,
\begin{equation}
\label{doublingcond}
0<\mu\big(B_\rho(x,2r)\big)\le C_\mu\mu\big(B_\rho(x,r)\big)<\infty.
\end{equation}
If $\rho$ is a metric, then $(X,\rho,\mu)$
is called a \emph{doubling metric measure space}
\index[words]{doubling metric measure space}.
Finally, the triplet $(X,\mathbf{q},\mu)$\index[symbols]{X@$(X,\mathbf{q},\mu)$}
is called
a \textit{quasi-ultrametric space of homogeneous type} if $(X,\rho,\mu)$ is a
quasi-ultrametric space of homogeneous type for some $\rho\in\mathbf{q}$.
\end{definition}

\begin{remark}
\label{rm1112}
Let $(X,\rho,\mu)$ be a
quasi-ultrametric space of homogeneous type.
\begin{enumerate}
\item
It follows from the doubling
condition \eqref{doublingcond} that,
for any
$\lambda\in\left[1,\infty\right)$, $x\in X$, and $r\in(0,\infty)$,
\begin{equation}
\label{upperdoub}
\mu\big(B_\rho(x,\lambda r)\big)\le
C_\mu\lambda^n\mu\big(B_\rho(x,r)\big),
\end{equation}
where $n:=\log_2C_\mu$\index[symbols]{N@$n$} is called
the \emph{upper dimension}\index[words]{upper dimension} of $X$.
Moreover,
if the cardinality of the set $X$ is at least $2$,
then, by \eqref{doublingcond} again,
we find that $C_\mu\in(1,\infty)$ and hence $n\in(0,\infty)$;
see \cite[p.\,72]{AlMi2015}.

\item
It is known that, under the assumption that
$\mu(B_\rho(x,r))\in(0,\infty)$ for any $x\in X$ and $r\in(0,\infty)$,
$\mathrm{diam}_\rho(X)=\infty$
if and only if $\mu(X)=\infty$
(see, for instance, \cite[Lemma~5.1]{ny1997} and
\cite[p.\,7]{YHYY1}).

\item
We refer to \cite[pp.\,588--590]{Cowe77}
for several examples of quasi-ultrametric spaces of homogeneous type.
\end{enumerate}
\end{remark}

\begin{remark}\label{Rmkdoub}
If $(X,\mathbf{q},\mu)$ is a quasi-ultrametric space of homogeneous type, then
the measure $\mu$ satisfies the doubling
condition \eqref{doublingcond} with respect to
any other $\rho'\in\mathbf{q}$
for which all $\rho'$-balls are $\mu$-measurable; see Proposition~\ref{L2.10} below.
\end{remark}

The following lemma is precisely
\cite[Lemma~2.1(i)]{k2014};
see also \cite[Theorem~1]{MaSe79i}.

\begin{lemma}\label{mu-r}
Let $(X,\rho,\mu)$ be a quasi-ultrametric space
of homogeneous type.
Then, for any $x\in X$,
$\mu(\{x\})\in(0,\infty)$
if and only if there exists $r\in(0,\infty)$ such that
$\{x\}=B_\rho(x,r)$.
\end{lemma}

Next, we discuss the notion of an ultra-RD-space
which is a special case of the quasi-ultrametric space
of homogeneous type, where, in addition to the doubling condition,
the measure also satisfies a reverse doubling condition. As a preamble, we
first introduce some notation.
Let $(X,\mathbf{q},\mu)$ be a quasi-ultrametric space
of homogeneous type. For any $\rho\in\mathbf{q}$ and
$x\in X$, define
\begin{align}
\label{rrho}
r_\rho(x)\index[symbols]{R@$r_\rho(x)$}:=
\inf\left\{\tau\in(0,\infty):
B_\rho(x,\tau)\supsetneqq\{x\}\right\}
\end{align}
and
\begin{align}
\label{Rrho}
R_\rho(x)\index[symbols]{R@$R_\rho(x)$}:=	
\begin{cases}
\sup\left\{\tau\in(0,\infty):
B_\rho(x,\tau)\subsetneqq X\right\}\
&\text{if}\ \mu(X)<\infty,\\
\infty\ &\text{if}\ \mu(X)=\infty.
\end{cases}
\end{align}
Some properties of $r_\rho$ and $R_\rho$ are
recorded in the following two remarks.

\begin{remark}
\label{1032}
Let $(X,\mathbf{q},\mu)$ be a quasi-ultrametric space of
homogeneous type and $\rho\in\mathbf{q}$.
Using \eqref{rrho} and \eqref{Rrho}, it is straightforward to check that the
following properties hold for any $x\in X$:
\begin{enumerate}
\item  $r_\rho(x)\in\left[0,\infty\right)$ and
$R_\rho(x)\in\left(0,\infty\right]$ are well defined.

\item $r_\rho(x)\leq R_\rho(x)$.

\item If $r_\rho(x)>0$, then, for any $r\in(0,r_\rho(x)]$,
$B_\rho(x,r)=\{x\}$.

\item If $R_\rho(x)<\infty$, then
$X\setminus B_\rho(x,R_\rho(x))=\{y\in X:\rho(x,y)=R_\rho(x)\}$.

\item If $\varrho\in\mathbf{q}$, more specifically,
if there exist two constants $C_1,C_2\in(0,\infty)$
such that $C_1\varrho\leq\rho\leq C_2\varrho$
pointwise on $X\times X$, then, for any $x\in X$,
\begin{align*}
C_1r_{\varrho}(x)\leq r_\rho(x)\leq C_2r_{\varrho}(x)\
\text{and}\
C_1R_{\varrho}(x)\leq R_\rho(x)\leq C_2R_{\varrho}(x).
\end{align*}
\end{enumerate}
\end{remark}

\begin{remark}
\label{1031}
Let $(X,\mathbf{q},\mu)$ be a quasi-ultrametric space of
homogeneous type and $\rho\in\mathbf{q}$.
By Remark~\ref{rm1112}(ii), \eqref{Rrho},
and \eqref{diam}, it is easy to see that
the following statements are equivalent:
\begin{enumerate}
\item
$\mu(X)=\infty$.

\item
$\inf_{x\in X}R_\rho(x)=\infty$.

\item
$\sup_{x\in X}R_\rho(x)=\infty$.

\item
$\mathrm{diam}_\rho(X)=\infty$.
\end{enumerate}
Moreover, if
$\mathrm{diam}_\rho(X)<\infty$,
by \cite[Proposition 2.12(1)]{AlMi2015}, we find that
\begin{align*}
\left(C_\rho\widetilde{C}_\rho\right)^{-1}\mathrm{diam}_\rho(X)
\leq\inf_{x\in X}R_\rho(x)
\leq\sup_{x\in X}R_\rho(x)\leq\mathrm{diam}_\rho(X).
\end{align*}
\end{remark}

The stage has now been set to introduce the notion of an ultra-RD-space
in the following definition.

\begin{definition}\label{RD}
A triplet $(X,\rho,\mu)$ is called a
$(\kappa,n)$-\emph{ultra-RD-space}
(or simply \emph{ultra-RD-space})\index[words]{ultra-RD-space}
if $(X,\rho,\mu)$ is a quasi-ultrametric space of homogeneous type with upper
dimension $n\in(0,\infty)$ in \eqref{upperdoub}
and if the measure $\mu$ satisfies the following additional
\emph{reverse doubling condition}\index[words]{reverse doubling condition}
(with respect to $\rho$):
there exist four constants $C_1,c_1\in(0,1]$,
$c_2\in[1,\infty)$, and $\kappa\in(0,n]$\index[symbols]{K@$\kappa$}
such that, for any $x\in X$,
$r\in(0,\infty)\cap[c_1r_\rho(x),c_2R_\rho(x)]$,
and $\lambda\in(0,\infty)\cap[1,\frac{c_2R_\rho(x)}{r}]$,
\begin{align}\label{RDcondition}
C_1\lambda^\kappa\mu\big(B_\rho(x,r)\big)
\le\mu\big(B_\rho(x,\lambda r)\big),
\end{align}
where $r_\rho(x)$ and $R_\rho(x)$ are as, respectively,
in \eqref{rrho} and \eqref{Rrho}.
A triplet $(X,\mathbf{q},\mu)$ is
called a $(\kappa,n)$-\emph{ultra-RD-space}
if $(X,\rho,\mu)$ is a $(\kappa,n)$-ultra-RD-space
for some $\rho\in\mathbf{q}$.
\end{definition}

\begin{remark}
\begin{enumerate}
\item[{\rm(i)}]
In 2004, Nagel and Stein \cite{ns2004} developed the product theory on
Carnot--Carath\'eodory spaces with a quasi-ultrametric $\rho$
and a measure $\mu$ satisfying
both the doubling and the reverse doubling conditions
\eqref{doublingcond} and \eqref{RDcondition}.
This reverse doubling property
was further extended by Han et al. \cite{hmy06,hmy08} into
metric measure spaces.
Such quasi-ultrametric spaces of homogeneous type
having the reverse-doubling property are
called RD-spaces, which were originally
introduced in \cite{hmy08}.
Recall that an \emph{RD-space}  in \cite{hmy08} is a doubling metric measure space
$(X,\rho,\mu)$, where the measure $\mu$ satisfies
the following additional condition:
there exist two constants $\widetilde{C}\in(0,1]$ and $\kappa\in(0,n]$
with $n$ in \eqref{upperdoub} such that,
for any $x\in X$, $r\in(0,\frac{\mathrm{diam}_\rho(X)}{2})$,
and $\lambda\in[1,\frac{\mathrm{diam}_\rho(X)}{2r})$,
$$\widetilde{C}\lambda^\kappa\mu
(B_\rho(x,r))\leq\mu(B_\rho(x,\lambda r)).$$
By Remark~\ref{1031},
we find that these RD-spaces
are special cases of ultra-RD-spaces in Definition~\ref{RD}
where $\rho$ is a metric and $r_\rho(x)=0$ for any $x\in X$.
For further studies about RD-spaces,
we refer to \cite{hmy06,hmy08,YZ}.
\item[\rm{(ii)}]
We point out that the condition that $r\ge c_1r_\rho(x)$
in Definition~\ref{RD}
is necessary if $X$ has the property that the measure of a singleton is
strictly positive.
However, as it was shown in \cite[p.\,259, Theorem 1]{MaSe79i},
there can only be at most countably many such points.

\item[\rm{(iii)}]
There are many examples of ultra-RD-spaces,
such as Euclidean
spaces with $A_\infty$-weights (Muckenhoupt's class),
Ahlfors-regular spaces, Lie groups of polynomial growth,
and Carnot--Carath\'eodory spaces with doubling measure.

\item[\rm{(iv)}]
There also exist several nontrivial examples of quasi-ultrametric spaces
of homogeneous type, which are not ultra-RD-spaces;
see, for instance, \cite[p.\,267]{AH1}
and \cite[Proposition~2.3]{wyyTLinfty}.
\end{enumerate}
\end{remark}

\begin{remark}
If $(X.\mathbf{q},\mu)$ is an ultra-RD-space, then
the measure $\mu$ satisfies the reverse doubling
condition \eqref{RDcondition} with respect to
any other quasi-ultrametric $\rho'\in\mathbf{q}$
for which all $\rho'$-balls are $\mu$-measurable;
see Proposition~\ref{RDsharp} below.
\end{remark}

Now, we recall the notion of
Ahlfors-regular quasi-ultrametric spaces
(see, for instance, \cite[Definition 2.11]{AlMi2015}).

\begin{definition}\label{defAR}
Let $n\in(0,\infty)$.
A triplet $(X,\rho,\mu)$ is called
an \emph{$n$-Ahlfors-regular quasi-ultrametric space}\index[words]{Ahlfors-regular
quasi-ultrametric space}
if $(X,\rho)$ is a quasi-ultrametric space
and $\mu$ is a nonnegative measure on $X$ with the property that
all $\rho$-balls are $\mu$-measurable and there exist
$c_1\in(0,1]$, $c_2\in[1,\infty)$, and
$C_1,C_2\in(0,\infty)$
such that the following
\emph{$n$-Ahlfors-regular condition}\index[words]{Ahlfors-regular condition}
holds:
for any $x\in X$ and $r\in(0,\infty)$
satisfying $r\in[c_1r_{\rho}(x),c_2R_{\rho}(x)]$,
\begin{align}\label{AR}
C_1r^n\leq\mu\big(B_{\rho}(x,r)\big)\leq C_2r^n,
\end{align}
where $r_{\rho}(x)$ and $R_{\rho}(x)$
are the same as, respectively, in \eqref{rrho} and \eqref{Rrho}.
A triplet $(X,\mathbf{q},\mu)$ is
called an \emph{$n$-Ahlfors-regular quasi-ultrametric space}
if $(X,\rho,\mu)$ is an $n$-Ahlfors-regular quasi-ultrametric space
for some $\rho\in\mathbf{q}$.
\end{definition}

\begin{remark}
\begin{enumerate}
\item
Recall that a triplet $(X,\rho,\mu)$ is called
an \emph{$n$-Ahlfors-regular metric space},
with $n\in(0,\infty)$, if $\rho$ is a metric and
there exists a constant $C\in[1,\infty)$
such that, for any $x\in X$ and
$r\in(0,\mathrm{diam}_\rho(X))$, $C^{-1}r^n\leq\mu(B_\rho(x,r))\leq Cr^n$.
By Remark~\ref{1031}, we find that
$n$-Ahlfors-regular metric spaces
are special cases of the $n$-Ahlfors-regular quasi-ultrametric
spaces in Definition~\ref{defAR}
when $\rho$ is a metric
and $r_\rho(x)=0$ for any $x\in X$.
\item
It is clear from definitions that the $n$-Ahlfors-regularity condition
\eqref{AR} implies the reverse doubling
condition \eqref{RDcondition} with $\kappa=n$.
Moreover, it can be shown that \eqref{AR}
also implies that the doubling condition
\eqref{upperdoub} holds;
see, for instance, \cite[Proposition~2.12]{AlMi2015}. Therefore,
an $n$-Ahlfors-regular quasi-ultrametric space
is an $(n,n)$-ultra-RD-space and hence
a quasi-ultrametric space
of homogeneous type with upper dimension $n$.
\item
Assume that $(X,\rho,\mu)$ is a quasi-ultrametric
space of homogeneous type with the property that all balls
are open.
Let the \emph{function $\delta(x,y)$} be defined by setting, for any $x,y\in X$,
\begin{align*}
\delta(x,y):=
\begin{cases}
\inf\left\{\mu(B):B\text{ is a ball
containing $x$ and $y$}\right\}\ &\text{if }x\neq y,\\
0&\text{if }x=y.
\end{cases}
\end{align*}
Then \cite[Theorem~3]{MaSe79i}
implies that $\delta$ is a quasi-ultrametric on $X$
satisfying that the topologies induced on $X$, respectively,
by $\rho$ and $\delta$ coincide.
Moreover, $(X,\delta,\mu)$
is a $1$-Ahlfors-regular quasi-ultrametric space.
\end{enumerate}
\end{remark}

\begin{remark}
If $(X,\mathbf{q},\mu)$ is an $n$-Ahlfors-regular
quasi-ultrametric space with $n\in(0,\infty)$, then
the measure $\mu$ satisfies the $n$-Ahlfors-regularity
condition \eqref{AR} with respect to
any other quasi-ultrametric $\rho'\in\mathbf{q}$
for which all $\rho'$-balls are $\mu$-measurable;
see Corollary~\ref{ADRequiv} below.
\end{remark}

Historically, assuming a certain degree of regularity
for the underlying measure has proven extremely useful
in establishing a number of fundamental results used
in the analysis on spaces of homogeneous type, such as
the density of smooth functions in Lebesgue spaces and
the Lebesgue differentiation theorem. It turns out that
the sharp regularity condition on the underlying measure
is the notion of Borel-semiregularity introduced in
\cite[Definition~3.9]{AlMi2015}.
We now take a moment to record this concept here as
it will be of central importance throughout this monograph.
To this end, in what follows, for an arbitrary
set $X$ equipped with a topology $\tau$\index[symbols]{T@$\tau$},
we denote by $Borel_\tau(X)$\index[symbols]{B@$Borel_\tau(X)$}
the smallest $\sigma$-algebra on $X$
containing $\tau$.
We now recall several related notions of measure regularity in the
following definition (see, for instance,
\cite[Definitions 2.9 and 3.9]{AlMi2015} and \cite[p.\,3]{H2001}).

\begin{definition}\label{Defmeasure}
Let $X$ be a set equipped with
a topology $\tau$,
and let $\mathfrak{M}$ be a $\sigma$-algebra of subsets of $X$.
Assume that $\mu:\mathfrak{M}\to[0,\infty]$ is a nonnegative measure.
\begin{enumerate}
\item
$\mu$ is called a \emph{Borel measure}\index[words]{Borel measure}
on $(X,\tau)$
if $Borel_\tau(X)\subset\mathfrak{M}$.
\item
$\mu$ is called a
\emph{Borel-semiregular measure}\index[words]{Borel-semiregular measure}
on $(X,\tau)$
if $\mu$ is a Borel measure on $X$
and satisfies that, for any $E\in\mathfrak{M}$
with $\mu(E)<\infty$,
there exists $B\in Borel_\tau(X)$
such that $\mu(E\Delta B)=0$,
where $E\Delta B\index[symbols]{E@$E\Delta B$}
:=(E\setminus B)\cup(B\setminus E)$.
\item
$\mu$ is called
a \emph{Borel-regular measure}\index[words]{Borel-regular measure}
on $(X,\tau)$
if $\mu$ is a Borel measure on $X$
and satisfies that, for any $E\in\mathfrak{M}$,
there exists $B\in Borel_\tau(X)$
such that $E\subset B$ and $\mu(E)=\mu(B)$.
\end{enumerate}
\end{definition}

\begin{remark}\label{Rmk502}
In regards to Definition~\ref{Defmeasure},
we have the following observations.
\begin{enumerate}
\item
Note that Zhou et al. \cite{zhy} assumed that $\mu$ is a Borel
regular measure and satisfies that the measure
of every singleton is zero.
Moreover, by (ii) and (iii) of Definition~\ref{Defmeasure},
it is easy to see that,
if $\mu$ is a Borel-regular measure on $(X,\tau)$,
then $\mu$ is a Borel-semiregular measure on $(X,\tau)$.
Thus,
all the results of this monograph are done under
some weaker assumptions of $\mu$ than those in \cite{zhy}.
\item
By Whitney's decomposition theorem for open sets
(see \cite[Theorem~6.2]{AlMiMi13},
\cite[Theorem~2.4]{AlMi2015}, and Lemma~\ref{Wtd}
in this monograph; see also Remark~\ref{940}(iii) in this regard)
every open subset of a quasi-ultrametric
space of homogeneous type $(X,\mathbf{q},\mu)$
can be decomposed into a countable union of $\rho$-balls,
each of which is $\mu$-measurable, since $\rho\in\mathbf{q}$
in Definition~\ref{Homogenous} is chosen
so that all $\rho$-balls are $\mu$-measurable.
Consequently, every doubling measure is automatically a Borel measure.
In particular, if $\varrho\in\mathbf{q}$ is any
regularized quasi-ultrametric as in Definition~\ref{D2.3},
then, all $\varrho$-balls are open, and thus
we find that all $\varrho$-balls are $\mu$-measurable.

\item
We stress here that, for
$\mu:\mathfrak{M}\to[0,\infty]$ to be a Borel measure, we only require
that $Borel_\tau(X)$ is a subset of $\mathfrak{M}$
and not necessarily that
$Borel_\tau(X)=\mathfrak{M}$. In the latter
case, the measure $\mu$ would automatically be Borel-regular.
\end{enumerate}
\end{remark}

We conclude this section by recalling the notation for volume functions
on a quasi-ultrametric space of homogeneous type and recording
some of their properties. To this end, let $(X,\rho,\mu)$ be
a quasi-ultrametric space of homogeneous type. Throughout this monograph,
for any $x,y\in X$ with $x\neq y$ and for any $r\in(0,\infty)$,
define
\begin{align*}
(V_\rho)_r(x):=\mu\big(B_\rho(x,r)\big)
\end{align*}
and
\begin{align*}
V_\rho(x,y):=(V_\rho)_{\rho(x,y)}(x)=\mu\big(B_\rho(x,\rho(x,y))\big).
\end{align*}

For a quasi-ultrametric space of homogeneous type,
we have the following proposition.

\begin{proposition}
\label{L2.10}
Let $(X,\mathbf{q},\mu)$ be a quasi-ultrametric space of homogeneous type
with upper dimension $n\in(0,\infty)$ in \eqref{upperdoub},
where $\mu$ satisfies the doubling condition \eqref{doublingcond}
with respect to $\rho\in\mathbf{q}$.
\begin{enumerate}
\item[{\rm (i)}]
For any given $\rho'\in\mathbf{q}$
that has the property that all $\rho'$-balls are $\mu$-measurable,
the triplet $(X,\rho',\mu)$ is also
a quasi-ultrametric space of homogeneous type.
Moreover,
there exists a positive constant $C$
(depending only on $\rho$, $\rho'$, and $\mu$) such that,
for any $\lambda\in[1,\infty)$, $x\in X$, and $r\in(0,\infty)$,
$$
(V_{\rho'})_{\lambda r}(x)
\le C\lambda^n
(V_{\rho'})_{r}(x).
$$

\item[{\rm (ii)}]
For any given $\gamma\in(0,\infty)$,
the triplet $(X,\rho^\gamma,\mu)$ is a
quasi-ultrametric space of homogeneous
type.
Moreover, for any $\lambda\in[1,\infty)$, $x\in X$, and $r\in(0,\infty)$,
$$
(V_{\rho^\gamma})_{\lambda r}(x)
\le C_\mu\lambda^\frac{n}{\gamma}
(V_{\rho^\gamma})_{r}(x),
$$
where $C_\mu\in[1,\infty)$ is the same as in
\eqref{upperdoub}.
\end{enumerate}
\end{proposition}

\begin{proof}
We first prove (i).
Suppose that $\rho'\in\mathbf{q}$ has the property
that all $\rho'$-balls are $\mu$-measurable.
Since $\rho,\rho'\in\mathbf{q}$,
it follows that there exists a constant $C\in[1,\infty)$
such that, for any $x,y\in X$,
\begin{align*}
C^{-1}\rho'(x,y)\leq\rho(x,y)\leq C\rho'(x,y).
\end{align*}
By this, we find that,
for any $\lambda\in[1,\infty)$, $x\in X$, and $r\in(0,\infty)$,
$$
B_{\rho'}(x,\lambda r)
\subset B_\rho(x,{C}\lambda r)
\ \ \text{and}\ \
B_\rho\left(x,C^{-1}r\right)\subset B_{\rho'}(x,r),
$$
which, combined with \eqref{upperdoub},
further implies that
\begin{align*}
(V_{\rho'})_{\lambda r}(x)
\leq(V_{\rho})_{{C}\lambda r}(x)
\lesssim\lambda^n(V_{\rho})_{C^{-1}r}(x)
\leq\lambda^n(V_{\rho'})_{r}(x).
\end{align*}
This finishes the proof of (i).

To show (ii), observe that, for any $x\in X$, $r\in(0,\infty)$,
and $\gamma\in(0,\infty)$,
$$
B_{\rho^\gamma}(x,r)=B_\rho\left(x,r^\frac{1}{\gamma}\right),
$$
which, combined with \eqref{upperdoub},
further implies that, for any $\lambda\in[1,\infty)$,
\begin{align*}
(V_{\rho^\gamma})_{\lambda r}(x)
=(V_{\rho})_{(\lambda r)^\frac{1}{\gamma}}(x)
\leq C_\mu\lambda^\frac{n}{\gamma}
(V_{\rho})_{r^\frac{1}{\gamma}}(x)
=C_\mu\lambda^\frac{n}{\gamma}
(V_{\rho^\gamma})_{r}(x).
\end{align*}
This finishes the proof of (ii)
and hence Proposition~\ref{L2.10}.
\end{proof}

Combining Proposition~\ref{L2.10}(i), Lemma~\ref{L2.3}(i),
and Remarks~\ref{Rmk502}(ii) and~\ref{1.1.7},
we immediately obtain the following conclusion.

\begin{corollary}
\label{2025}
Let $(X,\mathbf{q},\mu)$ be a quasi-ultrametric space of homogeneous type
and fix $\rho\in\mathbf{q}$.
Let $\rho_{\#}$ be the regularized quasi-ultrametric of $\rho$
as in Definition~\ref{D2.3}
and $\rho_{\textup{sym}}$ the symmetrization of $\rho$
as in Definition~\ref{treatment}(ii).
Then $(X,\rho_{\#},\mu)$ is a quasi-ultrametric space of homogeneous type.
Moreover, if all $\rho$-balls are $\mu$-measurable,
then $(X,\rho_{\textup{sym}},\mu)$
is also a quasi-ultrametric space of homogeneous type.
\end{corollary}

Some basic estimates are stated in the following lemma,
which is a generalization of \cite[Lemma 2.1]{hmy06}
(see also \cite[Lemma 2.4]{HLYY})
from doubling metric measure spaces
to quasi-ultrametric spaces of homogeneous type and
whose proof remains true for
quasi-ultrametric spaces of homogeneous type.
Indeed, recall that, in contrast to
doubling metric measure spaces,
for the quasi-ultrametric space of homogeneous type in this monograph,
we have no restriction on $\mu(\{x\})$ for any $x\in X$,
which means that it could happen that either
$\mu(\{x\})=0$ or $\mu(\{x\})>0$.
If, for some $x\in X$, $\mu(\{x\})>0$,
then, by \cite[Theorem~1]{MaSe79i},
we find that $r_\rho(x)>0$.
As such, it is easy to verify, using the techniques in
\cite[Lemma 2.1]{hmy06}, that
the following lemma also holds
for such $x\in X$; we omit the details.

\begin{lemma}\label{A3.2}
Let $(X,\rho,\mu)$ be a quasi-ultrametric space
of homogeneous type, where $\rho$ is a regularized
quasi-ultrametric as in Definition~\ref{D2.3}.
Then there exists a
positive constant $C$ such that
the following assertions hold.
\begin{enumerate}
\item[\textup{(i)}]
For any $x\in X$ and $\alpha,\delta\in(0,\infty)$,
$$
\int_{\rho(x,y)\leq\delta}
\frac{\left[\rho(x,y)\right]^\alpha}
{{V_\rho(x,y)}}\,d\mu(y)
\leq C\delta^\alpha,
$$
where ${V_\rho(x,x)}:=\mu(\{x\})$,
and
$$
\int_{\rho(x,y)\ge\delta}\frac{1}{V_\rho(x,y)}
\left[\frac{\delta}{\rho(x,y)}\right]^\alpha\,d\mu(y)\leq C.
$$
\item[\textup{(ii)}]
For any $x\in X$ and $\alpha,\delta\in(0,\infty)$,
$$
\int_X\frac{1}{(V_\rho)_\delta(x)+V_\rho(x,y)}
\left[\frac{\delta}{\delta+\rho(x,y)}\right]^\alpha\,d\mu(y)
\leq C.
$$
\item[\textup{(iii)}]
For any $x,x',x_1\in X$, $r\in(0,\infty)$,
and $\theta\in(0,1)$ satisfying
$\rho(x,x')\leq\theta[r+\rho(x,x_1)]$,
$$
r+\rho(x_1,x')\leq C\left[r+\rho(x_1,x)\right],
$$
$$
(V_\rho)_r(x_1)+V_\rho(x_1,x')
\leq C\left[(V_\rho)_r(x_1)+V_\rho(x_1,x)\right],
$$
and
$$
(V_\rho)_r(x')+V_\rho(x',x_1)
\leq C\left[(V_\rho)_r(x)+V_\rho(x,x_1)\right].
$$
\item[\textup{(iv)}]
For any $f\in L^1_{\mathrm{loc}}(X)$, $x\in X$,
and $\alpha,\delta\in(0,\infty)$,
$$
\int_{\rho(x,y)\geq\delta}
\frac{1}{V_\rho(x,y)}
\left[\frac{\delta}
{\rho(x,y)}\right]^\alpha|f(y)|\,d\mu(y)
\leq C\mathcal{M}_\rho f(x).
$$
\item[\textup{(v)}]
For any $x\in X$ and $\alpha,\delta,r\in(0,\infty)$,
$$
\int_{\rho(y,x)\ge\delta}
\frac{1}{(V_\rho)_r(x)+V_\rho(y,x)}
\frac{1}{[r+\rho(y,x)]^\alpha}\,d\mu(y)
\leq C\frac{1}{(r+\delta)^\alpha}.
$$
\item[\textup{(vi)}]
For any $x,y\in X$ and $r\in(0,\infty)$,
\begin{align*}
&(V_\rho)_r(x)+(V_\rho)_r(y)+V_\rho(x,y)\\
&\quad\sim (V_\rho)_r(y)+V_\rho(y,x)\sim(V_\rho)_r(x)+V_\rho(x,y)\\
&\quad\sim (V_\rho)_{r+\rho(x,y)}(x)\sim(V_\rho)_{r+\rho(x,y)}(y),
\end{align*}
where the positive equivalence constants are independent of $x$, $y$, and $r$.
\item[\rm(vii)]
For any given $\alpha\in(0,\infty)$
and for any $r\in(0,\infty)$ and
$x,y\in X$ satisfying $\rho(x,y)\leq\alpha r$,
\begin{align*}
(V_\rho)_r(x)
\sim(V_\rho)_{\alpha r}(x)
\sim(V_\rho)_{\alpha r}(y)
\sim(V_\rho)_r(y),
\end{align*}
where the positive equivalence constants
are independent of $x$, $y$, and $r$.
In particular,
$V_\rho(x,y)\sim V_\rho(y,x)$.
\end{enumerate}
\end{lemma}

\begin{remark}
By Lemma~\ref{L2.3}, it is straightforward to
verify that, if $\theta\in(C_\rho^2\widetilde{C}_\rho,\infty)$,
then, for any $x\in X$ and $r\in(0,\infty)$,
$$
\overline{B_\rho(x,r)}\subset B_\rho(x,\theta r),
$$
where $\overline{B_\rho(x,r)}$ denotes the closure of $B_\rho(x,r)$ in $X$.
Therefore, if $\mu$ is a doubling measure on $X$, then we have
$\mu(\overline{B_\rho(x,r)})\sim\mu({B_\rho(x,r)})$,
where the positive equivalence constants are independent of
both $x$ and $r$,
and hence the first inequality in Lemma~\ref{A3.2}(i)
is equivalent to
$$
\int_{\rho(x,y)\leq\delta}
\frac{\left[\rho(x,y)\right]^\alpha}
{\overline{V_\rho(x,y)}}\,d\mu(y)
\lesssim\delta^\alpha,
$$
where, for any $x,y\in X$, the \emph{symbol
$\overline{V_\rho(x,y)}$} denotes the
measure of the closure of $B_\rho(x,\rho(x,y)))$ in $X$
with the understanding that $B_\rho(x,\rho(x,x))=\varnothing$
and $\overline{B_\rho(x,\rho(x,x))}=\{x\}$.
\end{remark}

Next, we present some properties of the reverse doubling
condition and ultra-RD-spaces that will be helpful in what follows.
The following lemma shows that
the reverse doubling condition of the equipped measure can
be altered to accommodate a larger range of radii.

\begin{lemma}\label{kst.1}
Let $(X,\rho,\mu)$ be a $(\kappa,n)$-ultra-RD-space
with $0<\kappa\leq n<\infty$.
Let $\theta_1\in(0,1]$ and $\theta_2\in[1,\infty)$.
Then, for any $x\in X$,
$r\in(0,\infty)\cap[\theta_1c_1r_\rho(x),\theta_2c_2R_\rho(x)]$,
and $\lambda\in(0,\infty)\cap[1,\frac{\theta_2c_2R_\rho(x)}{r}]$,
\begin{align}\label{vow1}
C_1\theta_1^\kappa\theta_2^{-\kappa}\lambda^\kappa (V_\rho)_r(x)
\le (V_\rho)_{\lambda r}(x),
\end{align}
where $C_1,c_1\in(0,1]$ and $c_2\in[1,\infty)$
are the same positive constants as in Definition~\ref{RD}. 	
\end{lemma}	

\begin{proof}
We only prove \eqref{vow1} in the case both $\theta_1=1$
and $\theta_2\in[1,\infty)$ since the proof of the other case
is similar.
Fix $x\in X$.
If $R_\rho(x)=\infty$,
then the proof is trivial.
Thus, we may assume that $R_\rho(x)\in(0,\infty)$ and
we consider the following three cases.

\emph{Case (1)}
$r\leq c_2R_\rho(x)$
and $\lambda\leq \frac{c_2R_\rho(x)}{r}$.
In this case, we are done given that $(X,\rho,\mu)$ is
a $(\kappa,n)$-ultra-RD-space.

\emph{Case (2)}
$r\leq c_2R_\rho(x)$
and $\frac{c_2R_\rho(x)}{r}\leq\lambda\leq\frac{\theta_2c_2R_\rho(x)}{r}$.
In this case, by \eqref{RDcondition}, we obtain
\begin{align*}
C_1\theta_2^{-\kappa}\lambda^\kappa (V_\rho)_r(x)
&\le C_1\left[\frac{c_2R_\rho(x)}{r}\right]^\kappa
(V_\rho)_r(x)
\leq(V_\rho)_{c_2R_\rho(x)}(x)
\leq(V_\rho)_{\lambda r}(x).
\end{align*}

\emph{Case (3)}
$c_2R_\rho(x)<r\leq\theta_2c_2R_\rho(x)$
and $\lambda\leq\frac{\theta_2c_2R_\rho(x)}{r}$.
In this case, we have
\begin{align}
\label{1222-1}
\lambda\leq\frac{\theta_2c_2R_\rho(x)}{r}<\theta_2.
\end{align}
Moreover, from Remark~\ref{1032}(iv),
we infer that $B_\rho(x,r)=X$,
which, combined with \eqref{1222-1}
and $C_1\in(0,1]$, further implies that
\begin{align*}
C_1\theta_2^{-\kappa}\lambda^\kappa (V_\rho)_r(x)
\leq(V_\rho)_r(x)=\mu(X)=(V_\rho)_{\lambda r}(x).
\end{align*}

To summarize, we finish the proof of \eqref{vow1}
and hence Lemma~\ref{kst.1}.
\end{proof}

For the reverse doubling condition \eqref{RDcondition},
we also have the following proposition.

\begin{proposition}
\label{RDsharp}
Let $(X,\mathbf{q},\mu)$ be a $(\kappa,n)$-ultra-RD-space
with $0<\kappa\leq n<\infty$.
Assume that $\rho'\in\mathbf{q}$ is any quasi-ultrametric with the property
that all $\rho'$-balls are $\mu$-measurable.
Then there exists a positive constant $C$
such that, for any $x\in X$,
$r\in(0,\infty)\cap[c_1r_{\rho'}(x),c_2R_{\rho'}(x)]$,
and $\lambda\in(0,\infty)\cap[1,\frac{c_2R_{\rho'}(x)}{r}]$,
\begin{align}
\label{12.25.x1}
C\lambda^\kappa(V_{\rho'})_r(x)
\le(V_{\rho'})_{\lambda r}(x),
\end{align}
where $r_{\rho'}(x)$ and $R_{\rho'}(x)$
are as, respectively, in \eqref{rrho} and \eqref{Rrho}
with $\rho$ replaced by $\rho'$
and where $c_1\in(0,1]$ and $c_2\in[1,\infty)$ are
the same positive constants as in Definition~\ref{RD}.
In particular,
$(X,\rho',\mu)$ is a $(\kappa,n)$-ultra-RD-space.
\end{proposition}

\begin{proof}
By Proposition~\ref{L2.10}(i), we find that
$(X,\rho',\mu)$ is a quasi-ultrametric space of
homogeneous type with upper dimension $n$.
Now, we verify \eqref{12.25.x1}.
Since $(X,\mathbf{q},\mu)$ is a $(\kappa,n)$-ultra-RD-space,
it follows that there
exists $\rho\in\mathbf{q}$ such that $(X,\rho,\mu)$
is a $(\kappa,n)$-ultra-RD-space.
Given that $\rho,\rho'\in\mathbf{q}$,
there exists a constant $c\in[1,\infty)$ such that,
for any $x,y\in X$,
\begin{align*}
c^{-1}\rho'(x,y)\leq\rho(x,y)\leq c\rho'(x,y).
\end{align*}
By this and Remark~\ref{1032}(v),
we find that, for any $x\in X$,
$r\in(0,\infty)\cap[c_1r_{\rho'}(x),c_2R_{\rho'}(x)]$,
and $\lambda\in(0,\infty)\cap[1,\frac{c_2R_{\rho'}(x)}{r}]$,
$B_{\rho'}(x,r)
\subset B_\rho(x,cr)$,
$B_\rho(x,c^{-1}\lambda r)\subset B_{\rho'}(x,\lambda r)$,
$$
c_1r_{\rho}(x)
\leq c_1cr_{\rho'}(x)
\leq cr
\leq c_2cR_{\rho'}(x)
\leq c_2c^2R_{\rho}(x),
$$
and
$\lambda\leq\frac{c_2c^2R_{\rho}(x)}{cr}$,
which, combined with Lemma~\ref{kst.1} in the case
both $\theta_1:=1$ and $\theta_2:=c^2$
and with \eqref{upperdoub},
further implies that
\begin{align*}
\lambda^\kappa(V_{\rho'})_{r}(x)
\leq\lambda^\kappa
(V_\rho)_{cr}(x)
\lesssim(V_\rho)_{c\lambda r}(x)
\lesssim(V_\rho)_{c^{-1}\lambda r}(x)
\leq(V_{\rho'})_{\lambda r}(x).
\end{align*}
This shows \eqref{12.25.x1} and hence
finishes the proof of Proposition~\ref{RDsharp}.
\end{proof}

Combining Corollary~\ref{2025}, Proposition~\ref{RDsharp},
Lemma~\ref{L2.3}(i), and Remark~\ref{1.1.7},
we immediately obtain the following conclusion;
we omit the details.

\begin{corollary}
Let $(X,\mathbf{q},\mu)$ be a $(\kappa,n)$-ultra-RD-space
with $0<\kappa\leq n<\infty$ and fix $\rho\in\mathbf{q}$.
Let $\rho_\#\in\mathbf{q}$ be the regularized
quasi-ultrametric of $\rho$ as in Definition~\ref{D2.3} and
$\rho_{\textup{sym}}\in\mathbf{q}$ be the symmetrization
of $\rho$ as in Definition~\ref{treatment}(ii).
Then $(X,\rho_{\#},\mu)$ is a $(\kappa,n)$-ultra-RD-space.
Moreover, if all $\rho$-balls are $\mu$-measurable,
then $(X,\rho_{\textup{sym}},\mu)$
is also a $(\kappa,n)$-ultra-RD-space.
\end{corollary}

The following corollary on
Ahlfors-regular quasi-ultrametric spaces
is a direct consequence of Propositions~\ref{L2.10} and~\ref{RDsharp}
(see also \cite[Proposition~2.12]{AlMi2015});
we omit the details.

\begin{corollary}
\label{ADRequiv}
Let $(X,\mathbf{q},\mu)$ be an $n$-Ahlfors-regular quasi-ultrametric space
with $n\in(0,\infty)$ and fix $\rho\in\mathbf{q}$.
Let $\rho_{\#}\in\mathbf{q}$ be the regularized quasi-ultrametric of $\rho$
as in Definition~\ref{D2.3} and $\rho_{\textup{sym}}$ the symmetrization
of $\rho$ as in Definition~\ref{treatment}(ii).
Then $(X,\rho_{\#},\mu)$
is an $n$-Ahlfors-regular quasi-ultrametric space. Moreover, if all $\rho$-balls
are $\mu$-measurable, then $(X,\rho_{\textup{sym}},\mu)$
is also an $n$-Ahlfors-regular
quasi-ultrametric space. Furthermore,
for any given $\rho'\in\mathbf{q}$, if all $\rho'$-balls
are $\mu$-measurable, then $(X,\rho',\mu)$ is
an $n$-Ahlfors-regular quasi-ultrametric space.
\end{corollary}

Next, we recall the concepts of geometrically doubling
quasi-ultrametric spaces and systems of dyadic cubes
on geometrically doubling
quasi-ultrametric spaces.

\begin{definition}\label{2248}
A quasi-ultrametric space $(X,\rho)$ is said to be
\emph{geometrically doubling (with respect to $\rho$)}
\index[words]{geometrically doubling} if
there exists a constant $N\in\mathbb{N}$,
called the \emph{geometric doubling constant of
$(X,\rho)$}\index[words]{geometrically doubling constant of
$(X,\rho)$}, such that any $\rho$-ball
of radius $r\in(0,\infty)$ in $X$ can be covered by at most $N$ $\rho$-balls
in $X$ of radii $\frac{r}{2}$.
Moreover, a quasi-ultrametric space $(X,\mathbf{q})$ is called
a \emph{geometrically doubling quasi-ultrametric space} if $(X,\rho)$ is
a geometrically doubling quasi-ultrametric space for some
$\rho\in\mathbf{q}$ and, in this case, the constant $N\in\mathbb{N}$ is
called the \emph{geometric doubling constant
of $(X,\mathbf{q})$ with respect to $\rho$}.
\end{definition}

\begin{remark}\label{940}
\begin{enumerate}
\item[(i)]
A quasi-ultrametric space $(X,\mathbf{q})$ is
a geometrically doubling quasi-ultrametric space
if and only if, for any $\rho\in\mathbf{q}$ and $\theta\in(0,1)$,
there exists $N_{(\theta)}\in\mathbb{N}$, depending on $\theta$,
such that any $\rho$-ball of radius $r\in(0,\infty)$ in $X$
can be covered by at most $N_{(\theta)}$ $\rho$-balls
in $X$ of radii $\theta r$; see also \cite[(2.34)]{AlMi2015}.
This ensures that the definition
of $(X,\mathbf{q})$ as a
geometrically doubling quasi-ultrametric space
in Definition~\ref{2248} is meaningful.

\item[(ii)]
By \cite[(2.35)]{AlMi2015}, we find that,
if $(X,\mathbf{q})$ is a
geometrically doubling quasi-ultrametric space,
then the topological space $(X,\tau_\mathbf{q})$ is separable,
which is a consequence of the geometrically doubling property.
\item[(iii)] Any quasi-ultrametric space of homogeneous type $(X,\mathbf{q},\mu)$
is a geometrically doubling quasi-ultrametric space;
see \cite[Proposition~3.28]{AlMi2015} and \cite{CW1}.
\end{enumerate}
\end{remark}

We now recall the following construction
given by Hyt\"onen and Kairema in \cite[Theorem~2.2]{HK},
which provides an analogue of
the grid of Euclidean dyadic cubes
on geometrically doubling quasi-ultrametric spaces;
see also Christ \cite[Theorem 11]{C} and Hofmann et al.
\cite[Proposition~2.11]{HMMM}.

\begin{lemma}\label{DyadicSys}
Let $(X,\mathbf{q})$ be a
geometrically doubling quasi-ultrametric space
and $\rho\in\mathbf{q}$.
Then there exist a collection
$\{Q^k_\alpha
\subset X:k\in\mathbb{Z},\ \alpha\in I_k\}$
of open subsets,
where $I_k$ is some countable index set,
and constants $\delta\in(0,1)$ and
$C_l,C_r\in(0,\infty)$ such that\index[symbols]{C@$C_l$}\index[symbols]{C@$C_r$}
\begin{enumerate}
\item[\textup{(i)}]
for any given $k\in\mathbb{Z}$,
$X=\bigcup_{\alpha\in I_k} Q^k_\alpha$;
$Q^k_\alpha\cap Q^k_\beta=\varnothing$
if $\alpha,\beta\in I_k$
with $\alpha\neq\beta$;
\item[\textup{(ii)}]
for any given $k,l\in\mathbb{Z}$
with $l\ge k$ and for any $\alpha\in I_k$ and $\beta\in I_l$,
either $Q^l_\beta\subset Q^k_\alpha$
or $Q^l_\beta\cap Q^k_\alpha=\varnothing$;
\item[\textup{(iii)}]
for any given $k,l\in\mathbb{Z}$ such that
$l<k$ and for any $\alpha\in I_k$,
there exists a unique $\beta\in I_l$ such that
$Q^k_\alpha\subset Q^l_\beta$;
\item[\textup{(iv)}]
for any $k\in\mathbb{Z}$
and $\alpha\in I_k$,
$\mathrm{diam}_\rho(Q^k_\alpha)<C_r\delta^k$;
\item[\textup{(v)}]
for any $k\in\mathbb{Z}$
and $\alpha\in I_k$,
$Q^k_\alpha$
contains some ball
$B_\rho(z^k_\alpha,C_l\delta^k)$,
where $z^k_\alpha\in X$.
\end{enumerate}
\end{lemma}

\begin{remark}\label{remarkcube}
Indeed, we can think of $Q^k_\alpha$
as a \emph{dyadic cube}\index[words]{dyadic cube} with diameter rough $\delta^k$
and center $z^k_\alpha$. By \cite[pp.23--24]{HMMM}, it is always
possible to set $\delta=2^{-1}$ in the construction given by
Christ in \cite{C}. As such, in what follows,
we always suppose that $\delta=2^{-1}$ to simplify our presentation;
otherwise, we would replace $2^{-k}$
in the definition of the approximation
of the identity by $\delta^k$
as well as make some other necessary changes.
\end{remark}

\begin{remark}\label{jcubes}
In what follows, we always fix a sufficiently
large integer $j\in\mathbb{N}$ such that
\begin{align}\label{2cubej}
2^{-j}C_r<\frac{1}{3},
\end{align}
where $C_r$ is the same as in Lemma~\ref{DyadicSys}(iv).
For any given $k\in\mathbb{Z}$
and $\tau\in I_k$, we let the \emph{symbol $N(k,\tau)$}\index[symbols]{N@$N(k,\tau)$}
be the cardinality of the
set $\{\tau'\in I_{k+j}:Q^{k+j}_{\tau '}\subset Q^k_\tau\}$ and
we arrange the set of all
cubes $Q^{k+j}_{\tau '}\subset Q^k_\tau$ with $\tau'\in I_{k+j}$
as $\{Q^{k,v}_\tau\}_{v=1}^{N(k,\tau)}$\index[symbols]{Q@$Q^{k,v}_\tau$}.
For any $k\in\mathbb{Z}$,
$\tau\in I_k$, and $v\in\{1,...,N(k,\tau)\}$,
we denote by the \emph{symbol $z^{k,v}_\tau$}
the ``center" of $Q^{k,v}_\tau$
and by $y^{k,v}_\tau$ any point in $Q^{k,v}_\tau$.
It follows from both (iv) and (v) of Lemma~\ref{DyadicSys}
and the geometrically doubling property in Remark~\ref{940}(i) that,
for any $k\in\mathbb{Z}$ and $\tau\in I_k$,
$$
N(k,\tau)\leq N_{(C_l2^{-j}/C_r)}\ \ \left[\text{so, in particular},
N(k,\tau)\in\mathbb{N}\right],
$$
where
$N_{(C_l2^{-j}/C_r)}$ is the same positive constant as in
Remark~\ref{940}(i) with $\theta$ therein
replaced by $\frac{C_l2^{-j}}{C_r}$.
Moreover, by Lemma~\ref{DyadicSys}(i), we find that,
for any $k\in\mathbb{Z}$ and $\tau\in I_k$,
\begin{align}
\label{un-3iu}
\bigcup_{v=1}^{N(k,\tau)}Q^{k,v}_\tau=
\bigcup_{\tau'\in I_{k+j},\ Q^{k+j}_{\tau'}\subset Q_\tau^k
}Q^{k+j}_{\tau'}=Q_\tau^k.
\end{align}
\end{remark}

\section{A Maximally Smooth
Approximation of the Identity}
\label{ss2.2}

In this section, let
$(X,\mathbf{q},\mu)$ be
a quasi-ultrametric space of homogeneous type.
We aim to build a maximally smooth
approximation of the identity
with bounded support on $(X,\mathbf{q},\mu)$
and show some properties of
this brand of approximation of the identity.
As a preamble, we first record some definitions.
The \emph{Lebesgue space}\index[words]{Lebesgue space}
$L^p(X)$,\index[symbols]{L@$L^p(X)$} with $p\in(0,\infty]$,
is defined to be the set of
all $\mu$-measurable functions $f:X\to\mathbb{C}$ such that
\begin{align*}
\|f\|_{L^p(X)}:=\left\{
\begin{aligned}
&\left[\int_X|f(x)|^p\,d\mu(x)\right]^\frac{1}{p}<\infty&&\text{if }p\in(0,\infty),
\\
&\mathop\mathrm{\,ess\,sup\,}_{x\in X}|f(x)|<\infty&&\text{if }p=\infty.
\end{aligned}\right.
\end{align*}
The \emph{weak Lebesgue space}\index[words]{weak Lebesgue space}
$L^{p,\infty}(X)$,\index[symbols]{L@$L^{p,\infty}(X)$} with $p\in(0,\infty]$,
is defined to be the set of
all $\mu$-measurable functions $f:X\to\mathbb{C}$ such that
\begin{align*}
\|f\|_{L^{p,\infty}(X)}:=\left\{
\begin{aligned}
&\sup_{\lambda\in(0,\infty)}
\lambda\left[\mu\left(\left\{x\in X:|f(x)|>\lambda
\right\}\right)\right]^\frac{1}{p}<\infty&&\text{if }p\in(0,\infty),
\\
&\|f\|_{L^\infty(X)}<\infty&&\text{if }p=\infty.
\end{aligned}\right.
\end{align*}
The choice to denote the weak Lebesgue space
by $L^{p,\infty}(X)$ is due to the fact that the
weak Lebesgue space may be viewed as a particular example
of a more general scale of function spaces,
namely \emph{Lorentz spaces} $L^{p,q}(X)$, which we will
introduce in Section~\ref{section5.1}. In that setting,
the identity $L^{p,\infty}(X)=$ weak-$L^p(X)$ will become clear.
For now, we simply adopt this notation with the understanding
that weak Lebesgue spaces are a distinct object in their own right,
independent of the Lorentz spaces to be discussed later.
For $p\in(0,\infty)$, the \emph{symbol $L_\mathrm{loc}^p(X)$}\index[symbols]{L@$L_\mathrm{loc}^p(X)$} denotes
the space of all $p$-locally integrable
functions\index[words]{locally integrable function}
on $X$. That is, $L_\mathrm{loc}^p(X)$ is the set of
all $\mu$-measurable functions $f:X\to\mathbb{C}$
such that, for any point $x\in X$,
there exists a neighborhood of $x$ where $f$ is
$p$-integrable; see, for instance, \cite[p.\,5]{H2001}.
In what follows, $1$-locally integrable functions will simply be referred to
as \emph{locally integrable}.

We now introduce the following definition of the approximation of the identity
on quasi-ultrametric spaces of homogeneous type.

\begin{definition}
\label{DefATI}
Let $(X,\mathbf{q},\mu)$ be
a quasi-ultrametric space of homogeneous type.
A family $\{\mathcal{S}_t\}_{t\in(0,\infty)}$
\index[symbols]{S@$\{\mathcal{S}_t\}_{t\in(0,\infty)}$}
of integral operators
with kernels $\{S_t\}_{t\in(0,\infty)}$ mapping $X\times X$ to $\mathbb{R}$,
defined by setting, for any $t\in(0,\infty)$,
$f\in L_{\mathrm{loc}}^1(X)$, and $x\in X$,
\begin{align*}
\mathcal{S}_tf(x)
:=\int_{X}S_t(x,y)f(y)\,d\mu(y),
\end{align*}
is called an
\emph{approximation of the identity}
of order $\varepsilon\in(0,\infty)$
\index[words]{approximation of the identity}
on $(X,\mathbf{q},\mu)$
if there exist $\theta,C,C'\in(0,\infty)$
and $\rho\in\mathbf{q}$, where all $\rho$-balls are
$\mu$-measurable, such that, for any $t\in(0,\infty)$,
the following properties hold:

\begin{enumerate}
\item
For any $x,y\in X$,
\begin{align*}
0\leq S_t(x,y)\leq C\frac{1}{(V_\rho)_{t^\theta}(x)};
\end{align*}
moreover, if $\rho(x,y)\ge C't^\theta$, then $S_t(x,y)=0$;
\item
For any $x,y,x'\in X$,
$$
\left|S_t(x,y)-S_t(x',y)\right|
\leq C\left[\frac{\rho(x,x')}{t^\theta}\right]^\varepsilon
\left[\frac{1}{(V_\rho)_{t^\theta}(x)}
+\frac{1}{(V_\rho)_{t^\theta}(x')}\right];
$$
\item
For any $x,y,x',y'\in X$,
\begin{align*}
&\left|\left[S_t(x,y)-S_t(x',y)\right]
-\left[S_t(x,y')-S_t(x',y')\right]\right|\\
&\quad\leq C
\left[\frac{\rho(x,x')\rho(y,y')}{t^{2\theta}}\right]^\varepsilon
\left[\frac{1}{(V_\rho)_{t^\theta}(x)}
+\frac{1}{(V_\rho)_{t^\theta}(x')}\right];
\end{align*}
\item
For any $x,y\in X$,
$S_t(x,y)=S_t(y,x)$;
\item
For any $x\in X$,
$$
\int_{X}S_t(x,z)\,d\mu(z)=\int_{X}S_t(z,x)\,d\mu(z)=1.
$$
\end{enumerate}
\end{definition}

\begin{remark}\label{Rm}
Let the notation be as in Definition~\ref{DefATI}.
\begin{enumerate}
\item Definition~\ref{DefATI}(i) implies that,
for any $t\in(0,\infty)$, $x,y\in X$ with $x\neq y$,
and $\gamma\in\mathbb{R}_+$,
\begin{align}\label{1302}
0\leq S_t(x,y)
\lesssim\frac{1}{V_\rho(x,y)+(V_\rho)_{t^\theta}(x)+(V_\rho)_{t^\theta}(y)}
\left[\frac{{t^\theta}}{{t^\theta}+\rho(x,y)}\right]^\gamma,
\end{align}
where the implicit positive constant is independent of
$t$, $x$, and $y$, but depends on $\gamma$.
Indeed, if $t\in(0,\infty)$ and $x,y\in X$ satisfy $\rho(x,y)\ge C't^\theta$,
then Definition~\ref{DefATI}(i) implies that
\eqref{1302} holds. If $\rho(x,y)<C't^\theta$, then,
by Definition~\ref{DefATI}(i) and both
(vi) and (vii) of Lemma~\ref{A3.2},
we find that
\begin{align*}
0&\leq S_t(x,y)\leq C\frac{1}{(V_\rho)_{t^\theta}(x)}
\sim\frac{1}{V_\rho(x,y)+(V_\rho)_{t^\theta}(x)+(V_\rho)_{t^\theta}(y)}\\
&\leq(1+C')^\gamma\frac{1}{V_\rho(x,y)+
(V_\rho)_{t^\theta}(x)+(V_\rho)_{t^\theta}(y)}
\left[\frac{{t^\theta}}{{t^\theta}+\rho(x,y)}\right]^\gamma.
\end{align*}
Thus, \eqref{1302} holds.

\item
Both (iii) and (iv) of Definition~\ref{DefATI} imply that,
for any $t\in(0,\infty)$ and
$x,y,x',y'\in X$,
\begin{align*}
&\left|\left[S_t(x,y)-S_t(x',y)\right]
-\left[S_t(x,y')-S_t(x',y')\right]\right|\\
&\quad\lesssim\left[\frac{\rho(x,x')
\rho(y,y')}{t^{2\theta}}\right]^\varepsilon
\left[\frac{1}{(V_\rho)_{t^\theta}(y)}
+\frac{1}{(V_\rho)_{t^\theta}(y')}\right],
\end{align*}
where the implicit positive constant is independent of
$t$, $x$, $x'$, $y$, and $y'$.
\item[\textup{(iii)}]
Recall that, for any given $\rho_1,\rho_2\in\mathbf{q}$,
there exists a constant $C\in[1,\infty)$
such that, for any $x,y\in X$,
$$
C^{-1}\rho_1(x,y)\leq\rho_2(x,y)\leq C\rho_1(x,y),
$$
which further implies that, for any $r\in(0,\infty)$,
\begin{align*}
C^{-1}B_{\rho_2}(x,r)\subset B_{\rho_1}(x,r)
\subset CB_{\rho_2}(x,r).
\end{align*}
Therefore, if the balls with respect to
$\rho_1$ and $\rho_2$ are $\mu$-measurable, then, by \eqref{upperdoub},
we find that, for
any $x\in X$ and $r\in(0,\infty)$,
\begin{align*}
(V_{\rho_2})_r(x)
\lesssim(V_{\rho_2})_{C^{-1}r}(x)
\leq(V_{\rho_1})_{r}(x)
\leq(V_{\rho_2})_{Cr}(x)
\lesssim(V_{\rho_2})_r(x).
\end{align*}
Combining this and Lemma~\ref{L2.3}(ii),
we deduce that, if $\{\mathcal{S}_{t}\}_{t\in(0,\infty)}$ is
an approximation of the identity of
order $\varepsilon\in(0,\infty)$,
then its integral kernels still satisfy (i)--(v)
of Definition~\ref{DefATI} with $\rho$
replaced by any $\rho'\in\mathbf{q}$
having the property that all $\rho'$-balls
are $\mu$-measurable. In particular,
the kernels of $\{\mathcal{S}_{t}\}_{t\in(0,\infty)}$ still satisfy
(i)--(v) of Definition~\ref{DefATI} with $\rho$ replaced by its regularized
quasi-ultrametric $\rho_\#$. As such, without loss of generality, we
can always assume that $\rho$
is a regularized quasi-ultrametric as in Definition~\ref{DefATI}
in the remainder of this monograph.
\item[\textup{(iv)}]
By possibly increasing the constants $C,C'\in(0,\infty)$
in Definition~\ref{DefATI}, we may assume that $C,C'\in[1,\infty)$.
\item[\rm(v)]
In order to defray our notation and simplify our
presentation, we always assume that $\theta=1$ in Definition~\ref{DefATI}
throughout this monograph.
\end{enumerate}
\end{remark}

Now, we present the main result of this section
in the following theorem,
which shows the existence of the above approximation
of the identity with its order in a maximal range.
To prove it,
we borrow some ideas from
the proof of \cite[Theorem 3.22]{AlMi2015};
see also the proof of \cite[Theorem 4.93]{MiMiMiMo13}.
The reader is reminded of the significance of the
symbol $\preceq$ in Remark~\ref{rmk390}(iii).

\begin{theorem}\label{ThmATI}
Let $(X,\mathbf{q},\mu)$ be
a quasi-ultrametric space of homogeneous type.
Then, for any $\varepsilon_0\in(0,\infty)$
satisfying
$0<\varepsilon_0\preceq{\mathrm{ind\,}}(X,\mathbf{q})$
with ${\mathrm{ind\,}}(X,\mathbf{q})$ as in \eqref{index},
there exists an approximation of the identity
$\{\mathcal{S}_t\}_{t\in(0,\infty)}$
of any order $\varepsilon\in(0,\varepsilon_0]$.
\end{theorem}

\begin{proof}
We divide the proof into the following two steps.

\emph{Step 1.} In this step,
we prove the existence of the
approximation of the identity $\{\mathcal{S}_t\}_{t\in(0,\infty)}$
(with $\theta=1$ in Definition~\ref{DefATI})
of any order $\varepsilon\in(0,\varepsilon_0]$
in the case both $\varepsilon_0\in(0,1)$
and $\varepsilon_0\preceq{\mathrm{ind\,}}(X,\mathbf{q})$.
To this end, we follow Coifman's
idea (see the discussion on \cite[pp.\,16--17 and p.\,40]{DaJoSe85})
with the goal of monitoring the maximal amount
of H\"older regularity for the integral
kernels on quasi-ultrametric spaces of homogeneous type;
see also \cite[Theorem~(1.13)]{mst1992}.
By \eqref{index},
we find that there exists $\rho\in\mathbf{q}$
satisfying that
$\varepsilon_0\in(0,(\log_2C_\rho)^{-1}]$.
Let $\rho_\#\in\mathbf{q}$ be the
regularized quasi-ultrametric of $\rho$
as in Definition~\ref{D2.3}.
From this and \eqref{2.3}, we infer that,
for any $\varepsilon\in(0,\varepsilon_0]$ and
$x,y,z\in X$,
\begin{equation}
\label{rhohold}	
\left|\rho_\#(x,y)-\rho_\#(x,z)\right|
\leq\frac{1}{\varepsilon}
\left[\rho_\#(y,z)\right]^\varepsilon
\max\left\{\left[\rho_\#(x,y)\right]^{1-\varepsilon},
\,
\left[\rho_\#(x,z)\right]^{1-\varepsilon}\right\}.
\end{equation}
Let $h\in C^1(\mathbb{R})$
(where, generally speaking,
$C^k(\mathbb{R}^d)$\index[symbols]{C@$C^k(\mathbb{R}^d)$}
with $k\in\mathbb{N}\cup\{\infty\}$ denotes
the class of \emph{$k$-fold continuously differentiable functions}
on $\mathbb{R}^d$)
satisfy that $0\leq h\leq1$
pointwise on $\mathbb{R}$,
$h\equiv1$ on $[-\frac{1}{2},\frac{1}{2}]$,
and $h\equiv0$ on $\mathbb{R}\setminus(-2,2)$.
For any $t\in(0,\infty)$,
the integral operator $T_t$ is defined by setting,
for any $f\in L_{\mathrm{loc}}^1(X)$ and $x\in X$,
$$
T_tf(x):=\int_{X}h\left(t^{-1}\rho_{\#}(x,y)\right)
f(y)\,d\mu(y).
$$
We claim that, for any $t\in(0,\infty)$ and $x\in X$,
\begin{equation}\label{2.15}
T_t1(x)\sim(V_{\rho_\#})_t(x)
\end{equation}
and
\begin{equation}\label{2.15x}
T_t\left(\frac{1}{T_t1}\right)(x)\sim 1,
\end{equation}
where the positive equivalence constants
are independent of both $x$ and $t$.
Indeed, by the properties of $h$ and \eqref{upperdoub},
we find that, for any $t\in(0,\infty)$ and $x\in X$,
\begin{align}
\label{2.14}
T_t1(x)
=\int_{B_{\rho_{\#}}(x,2t)}
h\left(t^{-1}\rho_{\#}(x,y)\right)\,
d\mu(y)
\leq(V_{\rho_\#})_{2t}(x)
\lesssim(V_{\rho_\#})_t(x).
\end{align}	
From the properties of $h$ and \eqref{upperdoub},
we deduce that, for any $t\in(0,\infty)$ and $x\in X$,
\begin{align*}
T_t1(x)
\ge\int_{B_{\rho_{\#}}(x,\frac{t}{2})}
h\left(t^{-1}\rho_{\#}(x,y)\right)\,
d\mu(y)
=(V_{\rho_\#})_{\frac{t}{2}}(x)
\gtrsim(V_{\rho_\#})_t(x).
\end{align*}	
Combining this with \eqref{2.14},
we obtain \eqref{2.15}.

Next, we prove \eqref{2.15x}.
By Lemma~\ref{A3.2}(vii), we conclude that,
for any given $t\in(0,\infty)$ and $z\in X$ and
for any $y\in B_{\rho_\#}(z,2t)$,
$$
(V_{\rho_\#})_t(z)\sim(V_{\rho_\#})_t(y),
$$
where the positive equivalence constants
are independent of $y$, $z$, and $t$.
Combining this, \eqref{2.15},
and the properties of $h$,
we find that, for any $t\in(0,\infty)$ and $z\in X$,
\begin{align*}
T_t\left(\frac{1}{T_t1}\right)(z)
&=\int_{X}h\left(t^{-1}\rho_{\#}(z,y)\right)
\frac{1}{T_t1(y)}\,d\mu(y)
\\
&\sim\int_{X}\frac{h(t^{-1}\rho_{\#}(z,y))}
{(V_{\rho_\#})_t(y)}\,d\mu(y)
\sim\int_{B_{\rho_{\#}}(z,2t)}
\frac{h(t^{-1}\rho_{\#}(z,y))}
{(V_{\rho_\#})_t(z)}\,d\mu(y)
\sim 1,
\end{align*}
where the positive equivalence constants are
independent of both $z$ and $t$.
This finishes the proof of \eqref{2.15x}
and hence the above claim.

Now, for any $t\in(0,\infty)$ and $x,y\in X$, define
$$
S_t(x,y)
:=\frac{1}{T_t1(x)T_t1(y)}
\int_{X}\frac{h(t^{-1}\rho_{\#}(x,z))h(t^{-1}\rho_{\#}(z,y))}
{T_t(\frac{1}{T_t1})(z)}\,d\mu(z).
$$
From the above claim, we infer that $S_t$ is well defined
and $S_t(x,y)\ge0$ for any $t\in(0,\infty)$ and $x,y\in X$.
Next, we turn to show that $S_t$ satisfies
(i)--(v) of Definition~\ref{DefATI} with $\rho$ replaced by $\rho_\#$.

Fix $t\in(0,\infty)$.
We first show that $S_t$ satisfies Definition~\ref{DefATI}(i).
By \eqref{2.15}, \eqref{2.15x}, the properties of $h$,
and \eqref{doublingcond},
we find that, for any $x,y\in X$,
\begin{align}
\label{v.w1}
S_t(x,y)
&\sim\frac{1}{(V_{\rho_\#})_t(x)(V_{\rho_\#})_t(y)}
\int_{X}h\left(t^{-1}\rho_{\#}(x,z)\right)
h\left(t^{-1}\rho_{\#}(z,y)\right)\,d\mu(z)
\nonumber\\
&\leq\frac{(V_{\rho_\#})_{2t}(y)}{(V_{\rho_\#})_t(x)(V_{\rho_\#})_t(y)}
\lesssim\frac{1}{(V_{\rho_\#})_t(x)},
\end{align}
where the implicit positive
constants are independent of $x$, $y$, and $t$.
In addition, from Definition~\ref{D2.1}(iii),
it follows that,
for any given $x,y\in X$ satisfying $\rho_\#(x,y)\ge 2C_{\rho_{\#}}t$
and for any $z\in X$,
$$
2C_{\rho_{\#}}t\leq\rho_\#(x,y)\leq C_{\rho_{\#}}
\max\left\{\rho_{\#}(x,z),\,\rho_{\#}(y,z)\right\},
$$
which further implies that
$\max\{{\rho_{\#}(x,z),\,\rho_{\#}(y,z)}\}\ge 2t$
and hence
$$
h\left(t^{-1}\rho_{\#}(x,z)\right)
h\left(t^{-1}\rho_{\#}(z,y)\right)
\equiv0.
$$
Thus, for any $x,y\in X$ satisfying
$\rho_\#(x,y)\ge 2C_{\rho_{\#}}t$,
$S_t(x,y)=0$.
Combining this and \eqref{v.w1}, we conclude that $S_t$
satisfies Definition~\ref{DefATI}(i).	

Next, we show that $S_t$ satisfies both
(iv) and (v) of Definition~\ref{DefATI}.
By the symmetry of $\rho_\#$ and
the definition of $S_t(x,y)$, it is easy to check that,
for any $x,y\in X$,
$S_t(x,y)=S_t(y,x)$.
Moreover, from Tonelli's theorem, we deduce that, for any $x\in X$,
\begin{align*}
\int_{X}S_t(x,y)\,d\mu(y)
&=\int_{X}\frac{h(t^{-1}\rho_{\#}(x,z))}
{T_t1(x)T_t{({\frac{1}{T_t1}})}(z)}
\left[\int_{X}
\frac{h(t^{-1}\rho_{\#}(z,y))}
{T_t1(y)}\,d\mu(y)\right]\,d\mu(z)\\
&=\int_{X}\frac{h(t^{-1}\rho_{\#}(x,z))}
{T_t1(x)T_t({\frac{1}{T_t1}})(z)}
{T_t{\left({\frac{1}{T_t1}}\right)}(z)}\,d\mu(z)=1.
\end{align*}
Thus, $S_t$ satisfies both
(iv) and (v) of  Definition~\ref{DefATI}.

Now, we prove that
$S_t$ satisfies Definition~\ref{DefATI}(ii) and has the following
property (which will be used later in this proof):\ for any $x,x',y\in X$,
\begin{align}\label{2y3new}
\left|S_t(x,y)-S_t(x',y)\right|
\lesssim
\left[\frac{\rho_{\#}(x,x')}{t}\right]^\varepsilon
\frac{1}{(V_{\rho_\#})_t(y)}.
\end{align}
To this end,
let $x,x',y\in X$. We divide the proof into the following three cases.

\emph{Case (1)}
$\rho_\#(x,x')\ge t$.
In this case, by Definition~\ref{DefATI}(i), we conclude that
\begin{align}\label{4.12y2}
\left|S_t(x,y)-S_t(x',y)\right|
&\leq\left|S_t(x,y)\right|+\left|S_t(x',y)\right|
\lesssim
\frac{1}{(V_{\rho_\#})_t(x)}
+\frac{1}{(V_{\rho_\#})_t(x')}
\nonumber\\
&\leq\left[\frac{\rho_{\#}(x,x')}{t}\right]^\varepsilon
\left[\frac{1}{(V_{\rho_\#})_t(x)}
+\frac{1}{(V_{\rho_\#})_t(x')}\right].
\end{align}
This is the desired estimate in Definition~\ref{DefATI}(ii).
Moreover, from both (i) and (iv) of Definition~\ref{DefATI},
we infer that
\begin{align}\label{2y1}
\left|S_t(x,y)-S_t(x',y)\right|
&\lesssim
\frac{1}{(V_{\rho_\#})_t(y)}
\leq\left[\frac{\rho_{\#}(x,x')}{t}\right]^\varepsilon
\frac{1}{(V_{\rho_\#})_t(y)},
\end{align}
which is \eqref{2y3new}.

\emph{Case (2)}
$\rho_\#(x,x')<t$ and $\rho_\#(x,y)\ge2C^2_{\rho_\#}t$.
In this case, by Definition~\ref{D2.1}(iii), we obtain
$\rho_\#(x',y)\ge2C_{\rho_\#}t,$
which, combined with Definition~\ref{DefATI}(i), further implies that
\begin{align}\label{2y2}
S_t(x,y)-S_t(x',y)=0.
\end{align}
Hence, the desired estimates in Definition~\ref{DefATI}(ii) and \eqref{2y3new}
are trivially satisfied.

\emph{Case (3)}
$\rho_\#(x,x')<t$ and $\rho_\#(x,y)<2C^2_{\rho_\#}t$.
In this case, we write
\begin{align}
\label{4.12z}
&S_t(x,y)-S_t(x',y)
\nonumber\\
&\quad=\frac{1}{T_t1(y)}
\left[\frac{1}{T_t1(x)}-\frac{1}{T_t1(x')}\right]
\int_{X}\frac{h(t^{-1}\rho_{\#}(x,z))
h(t^{-1}\rho_{\#}(z,y))}
{T_t(\frac{1}{T_t1})(z)}\,d\mu(z)
\nonumber\\
&\qquad+\frac{1}{T_t1(y)T_t1(x')}
\int_{X}\frac{[h(t^{-1}\rho_{\#}(x,z))
-h(t^{-1}\rho_{\#}(x',z))]
h(t^{-1}\rho_{\#}(z,y))}
{T_t(\frac{1}{T_t1})(z)}\,d\mu(z)
\nonumber\\
&\quad=:I(x,y,x')+II(x,y,x').
\end{align}
To estimate $I(x,y,x')$, let $D:=B_{\rho_\#}(x,2t)\cup B_{\rho_\#}(x',2t)$ and
observe that
\begin{align}\label{1444}
D\subset\left[B_{\rho_{\#}}(x,2C_{\rho_{\#}}t)
\cap B_{\rho_{\#}}(x',2C_{\rho_{\#}}t)\right].
\end{align}
Then, in the current scenario, we have
\begin{align*}
\left|\frac{1}{T_t1(x)}-\frac{1}{T_t1(x')}\right|
&=\left|\frac{T_t1(x')-T_t1(x)}
{T_t1(x)T_t1(x')}\right|
\\
&\lesssim\frac{\int_{X}{|h(t^{-1}\rho_{\#}(x,z))
-h(t^{-1}\rho_{\#}(x',z))|}\,d\mu(z)}
{(V_{\rho_\#})_t(x)
(V_{\rho_\#})_t(x')}\quad\text{by \eqref{2.15}}
\\
&\leq\frac{\int_{D}\|h'\|_{\infty}
|\rho_{\#}(x,z)-\rho_{\#}(x',z)|\,d\mu(z)}
{t(V_{\rho_\#})_t(x)(V_{\rho_\#})_t(x')}
\quad\text{by the mean value theorem}
\\
&\lesssim
\frac{[\rho_{\#}(x,x')]^\varepsilon
\int_{D}\Big\{[\rho_{\#}(x,z)]^{1-\varepsilon}
+[\rho_{\#}(z,x')]^{1-\varepsilon}\Big\}\,d\mu(z)}
{t(V_{\rho_\#})_t(x)(V_{\rho_\#})_t(x')}
\quad\text{by \eqref{rhohold}},
\end{align*}
which further implies that
\begin{align}\label{2.19}
&\left|\frac{1}{T_t1(x)}-\frac{1}{T_t1(x')}\right|
\nonumber\\
&\quad\lesssim\frac{
[\rho_{\#}(x,x')]^\varepsilon}
{t(V_{\rho_\#})_t(x)(V_{\rho_\#})_t(x')}
\left\{\int_{B_{\rho_{\#}(x,2C_{\rho_{\#}}t)}}
\left[\rho_{\#}(x,z)\right]^{1-\varepsilon}\,d\mu(z)\right.
\nonumber\\
&\qquad\left.
+\int_{B_{\rho_{\#}(x',2C_{\rho_{\#}}t)}}
\left[\rho_{\#}(z,x')\right]^{1-\varepsilon}\,d\mu(z)\right\}
\quad\text{by \eqref{1444}}
\nonumber\\
&\quad\lesssim\left[\frac{\rho_{\#}(x,x')}{t}\right]^\varepsilon\frac{
\mu(B_{\rho_{\#}}
(x',2C_{\rho_{\#}}t))+\mu(B_{\rho_{\#}}(x,2C_{\rho_{\#}}t))}
{(V_{\rho_\#})_t(x)(V_{\rho_\#})_t(x')}
\quad\text{since $\varepsilon\leq\varepsilon_0<1$}
\nonumber\\
&\quad\lesssim
\left[\frac{\rho_{\#}(x,x')}{t}\right]^\varepsilon
\left[\frac{1}{(V_{\rho_\#})_t(x)}
+\frac{1}{(V_{\rho_\#})_t(x')}\right]
\quad\text{by \eqref{upperdoub}}.
\end{align}
As such,
by \eqref{2.15}, \eqref{2.19},
\eqref{2.15x}, and the properties of $h$,
we conclude that
\begin{align}
\label{4.12y1}
\left|I(x,y,x')\right|
&\lesssim\left[\frac{\rho_{\#}(x,x')}{t}\right]^\varepsilon
\left[\frac{1}{(V_{\rho_\#})_t(x)}
+\frac{1}{(V_{\rho_\#})_t(x')}\right].
\end{align}
As concerns $II(x,y,x')$, from Lemma~\ref{A3.2}(vii),
it follows that
$(V_{\rho_\#})_t(x)
\sim(V_{\rho_\#})_t(y)$,
which, together with an argument similar to
that used in the proof of \eqref{2.19}, implies that
\begin{align}
\label{942}
\left|II(x,y,x')\right|
&\lesssim\frac{1}{
(V_{\rho_\#})_t(y)(V_{\rho_\#})_t(x')}
\int_{X}\left|h\left(t^{-1}\rho_{\#}(x,z)\right)
-h\left(t^{-1}\rho_{\#}(x',z)\right)\right|
\nonumber\\
&\quad\times
h\left(t^{-1}\rho_{\#}(z,y)\right)\,d\mu(z)
\quad\text{by \eqref{2.15} and \eqref{2.15x}}
\nonumber\\	
&\leq\frac{1}{(V_{\rho_\#})_t(x)
(V_{\rho_\#})_t(x')}
\int_{X}\left|h\left(t^{-1}\rho_{\#}(x,z)\right)
-h\left(t^{-1}\rho_{\#}(x',z)\right)\right|\,d\mu(z)
\nonumber\\
&\lesssim\left[\frac{\rho_{\#}(x,x')}{t}\right]^\varepsilon
\left[\frac{1}{(V_{\rho_\#})_t(x)}
+\frac{1}{(V_{\rho_\#})_t(x')}\right].
\end{align}
This is the desired estimate in Definition~\ref{DefATI}(ii).
Moreover,
by Definition~\ref{D2.1}(iii), we obtain
$\rho_\#(x',y)\leq2C^3_{\rho_\#}t$.
From this and Lemma~\ref{A3.2}(vii),
we deduce that
$$(V_{\rho_\#})_t(x)\sim
(V_{\rho_\#})_t(y)\sim(V_{\rho_\#})_t(x'),$$
which, combined with \eqref{4.12z},
\eqref{4.12y1}, and \eqref{942},
further implies that
\begin{align}\label{2y3}
\left|S_t(x,y)-S_t(x',y)\right|
\lesssim\left[\frac{\rho_{\#}(x,x')}{t}\right]^\varepsilon
\frac{1}{(V_{\rho_\#})_t(y)},
\end{align}
which is \eqref{2y3new}.

Putting \eqref{4.12z},
\eqref{4.12y1}, \eqref{4.12y2}, \eqref{942},
and \eqref{2y2} together implies that $S_t$
satisfies Definition~\ref{DefATI}(ii).
Moreover, by \eqref{2y1}, \eqref{2y2},
and \eqref{2y3}, we conclude that,
for any $x,x',y\in X$,
$S_t$ satisfies \eqref{2y3new}.

Finally, we show that $S_t$
satisfies Definition~\ref{DefATI}(iii).
To this end,
let $x,x',y,y'\in X$. We divide the proof into the following four cases.

\emph{Case (1)}
$\rho_{\#}(x,x')\leq 2C_{\rho_{\#}}^2t$
and $\rho_{\#}(y,y')\leq 2C_{\rho_{\#}}^2t$.
In this case, we write
\begin{align*}
S_t(x,y)-S_t(x',y)-\left[S_t(x,y')-S_t(x',y')\right]
=:\sum_{j=1}^4I_j(x,y,x',y'),
\end{align*}
where
\begin{align*}
I_1(x,y,x',y'):&=\left[\frac{1}{T_t1(x)}-\frac{1}{T_t1(x')}\right]
\left[\frac{1}{T_t1(y)}-\frac{1}{T_t1(y')}\right]\\
&\quad\times\int_X
\frac{h(t^{-1}\rho_{\#}(x,z))h(t^{-1}\rho_{\#}(z,y))}
{{T_t(\frac{1}{T_t1})(z)}}\,d\mu(z),
\end{align*}
\begin{align*}
I_2(x,y,x',y'):&=\frac{1}{T_t1(y')}\left[\frac{1}{T_t1(x)}-\frac{1}{T_t1(x')}\right]
\\
&\quad\times\int_X
\frac{h(t^{-1}\rho_{\#}(x,z))
[h(t^{-1}\rho_{\#}(z,y))
-h(t^{-1}\rho_{\#}(z,y'))]}
{{T_t(\frac{1}{T_t1})(z)}}\,d\mu(z),
\end{align*}
\begin{align*}
I_3(x,y,x',y'):&=\frac{1}{T_t1(x')}
\left[\frac{1}{T_t1(y)}-\frac{1}{T_t1(y')}\right]\\
&\quad\times\int_X\frac{[h(t^{-1}\rho_{\#}(x,z))
-h(t^{-1}\rho_{\#}(x',z))]
h(t^{-1}\rho_{\#}(z,y))}
{{T_t(\frac{1}{T_t1})(z)}}\,d\mu(z),
\end{align*}
and
\begin{align*}
&I_4(x,y,x',y')\\
&\quad:=\frac{1}{T_t1(x')}\frac{1}{T_t1(y')}\\
&\qquad\times\int_X\frac{
[h(t^{-1}\rho_{\#}(x,z))-h(t^{-1}\rho_{\#}(x',z))]
[h(t^{-1}\rho_{\#}(z,y))-h(t^{-1}\rho_{\#}(z,y'))]}
{T_t(\frac{1}{T_t1})(z)}\,d\mu(z).
\end{align*}
For $I_1(x,y,x',y')$, from \eqref{upperdoub},
an argument similar to that used in the estimation of
\eqref{2.19}, \eqref{2.15x}, the properties of $h$,
and Lemma~\ref{A3.2}(vii),
we infer that
\begin{align*}
\left|I_1(x,y,x',y')\right|
&\lesssim
\left[\frac{\rho_{\#}(x,x')}{t}\right]^\varepsilon
\left[\frac{1}{(V_{\rho_\#})_t(x)}
+\frac{1}{(V_{\rho_\#})_t(x')}\right]\\
&\quad\times
\left[\frac{\rho_{\#}(y,y')}{t}\right]^\varepsilon
\left[\frac{1}{(V_{\rho_\#})_t(y)}
+\frac{1}{(V_{\rho_\#})_t(y')}\right]
(V_{\rho_\#})_{2t}(y)\\
&\lesssim
\left[\frac{\rho_{\#}(x,x')\rho_{\#}(y,y')}{t^2}\right]^\varepsilon
\left[\frac{1}{(V_{\rho_\#})_t(x)}
+\frac{1}{(V_{\rho_\#})_t(x')}\right].
\end{align*}
For $I_2(x,y,x',y')$, by \eqref{upperdoub},
an argument similar to that used in the proof of
\eqref{2.19}, \eqref{2.15},
\eqref{2.15x}, the properties of $h$, and Lemma~\ref{A3.2}(vii),
we find that
\begin{align*}
\left|I_2(x,y,x',y')\right|
&\lesssim\left[\frac{\rho_{\#}(x,x')}{t}\right]^\varepsilon
\left[\frac{1}{(V_{\rho_\#})_t(x)}
+\frac{1}{(V_{\rho_\#})_t(x')}\right]
\frac{1}{(V_{\rho_\#})_t(y')}\\
&\quad\times\left[\frac{\rho_{\#}(y,y')}{t}\right]^\varepsilon
\left[(V_{\rho_\#})_t(y)
+(V_{\rho_\#})_t(y')\right]\\
&\lesssim\left[\frac{\rho_{\#}(x,x')\rho_{\#}(y,y')}{t^2}\right]^\varepsilon
\left[\frac{1}{(V_{\rho_\#})_t(x)}
+\frac{1}{(V_{\rho_\#})_t(x')}\right].
\end{align*}
The proof of the estimate for $I_3(x,y,x',y')$ is similar to $I_2(x,y,x',y')$
and we omit the details.
For $I_4(x,y,x',y')$, from \eqref{2.15}, \eqref{2.15x}, the properties of $h$,
an argument similar to that used
in the proof of \eqref{2.19}, and Lemma~\ref{A3.2}(vii),
we deduce that
\begin{align}\label{4130}
\left|I_4(x,y,x',y')\right|
&\lesssim
\frac{1}{(V_{\rho_\#})_t(y')(V_{\rho_\#})_t(x')}
\int_X\left|h\left(t^{-1}
\rho_{\#}(x,z)\right)-h\left(t^{-1}\rho_{\#}(x',z)\right)\right|
\nonumber\\
&\quad\times
\left|h\left(t^{-1}\rho_{\#}(z,y)\right)-h\left(t^{-1}
\rho_{\#}(z,y')\right)\right|\,d\mu(z)
\nonumber\\
&\leq\frac{1}{(V_{\rho_\#})_t(y')(V_{\rho_\#})_t(x')}
(t\varepsilon)^{-2}
\left[\rho_{\#}(x,x')\right]^\varepsilon
\left[\rho_{\#}(y,y')\right]^\varepsilon
\nonumber\\
&\quad\times
\int_{D_1\cap D_2}
\left\{\left[\rho_{\#}(x,z)\right]^{1-\varepsilon}+
\left[\rho_{\#}(x',z)\right]^{1-\varepsilon}\right\}
\nonumber\\
&\quad\times
\left\{\left[\rho_{\#}(y,z)\right]^{1-\varepsilon}+
\left[\rho_{\#}(y',z)\right]
^{1-\varepsilon}\right\}\,d\mu(z),
\end{align}
where
\begin{align*}
D_1:=\left[B_{\rho_\#}(x,2t)\cup B_{\rho_\#}(x',2t)\right]
\subset B_{\rho_\#}(x,2C^3_{\rho_\#}t)
\end{align*}
and
\begin{align*}
D_2:=\left[B_{\rho_\#}(y,2t)\cup B_{\rho_\#}(y',2t)\right]
\subset B_{\rho_\#}(y,2C^3_{\rho_\#}t).
\end{align*}
If $x,y\in X$ satisfy
$\rho_\#(x,y)\ge2C^4_{\rho_\#}t$, then
\begin{align*}
\left[D_1\cap D_2\right]\subset
\left[B_{\rho_\#}(x,2C^3_{\rho_\#}t)
\cap B_{\rho_\#}(y,2C^3_{\rho_\#}t)\right]=\varnothing
\end{align*}
and hence
\begin{align*}
\left|I_4(x,y,x',y')\right|=0.
\end{align*}
If $x,y\in X$ satisfy $\rho_\#(x,y)<2C^4_{\rho_\#}t$, then,
by Definition~\ref{D2.1}(iii) and the fact that we are currently assuming
$\rho_\#(x,x')\leq2C^2_{\rho_\#}t$
and $\rho_\#(y,y')\leq2C^2_{\rho_\#}t$, we conclude that,
for any $z\in B_{\rho_\#}(x,2C^3_{\rho_\#}t)$,
\begin{align}\label{4131}
\rho_\#(y,z)<2C^5_{\rho_\#}t, \
\rho_\#(y',z)<2C^6_{\rho_\#}t,\
\mbox{and}\
\rho_\#(x',z)<2C^4_{\rho_\#}t.
\end{align}
From this and Lemma~\ref{A3.2}(vii),
we infer that, for any $z\in B_{\rho_\#}(x,2C^3_{\rho_\#}t)$,
\begin{align*}
(V_{\rho_\#})_t(x)
\sim(V_{\rho_\#})_t(x')
\sim(V_{\rho_\#})_t(z)
\sim(V_{\rho_\#})_t(y)
\sim(V_{\rho_\#})_t(y'),
\end{align*}
which, combined with both \eqref{4130}
and \eqref{4131}, further implies that
\begin{align*}
\left|I_4(x,y,x',y')\right|
\lesssim\left[\frac{\rho_{\#}(x,x')\rho_{\#}(y,y')}{t^2}\right]^\varepsilon
\left[\frac{1}{(V_{\rho_\#})_t(x)}
+\frac{1}{(V_{\rho_\#})_t(x')}\right].
\end{align*}
By the estimates of $I_1(x,y,x',y')$,
$I_2(x,y,x',y')$, $I_3(x,y,x',y')$, and $I_4(x,y,x',y')$,
we obtain the desired estimate in Case (1).

\emph{Case (2)}
$\rho_{\#}(x,x')>2C_{\rho_{\#}}^2t$
and $\rho_{\#}(y,y')>2C_{\rho_{\#}}^2t$.
In this case,
\begin{align*}
&\left|\left[S_t(x,y)-S_t(x',y)\right]
-\left[S_t(x,y')-S_t(x',y')\right]\right|\\
&\quad\leq\left|S_t(x,y)\right|+
\left|S_t(x',y)\right|+\left|S_t(x,y')\right|
+\left|S_t(x',y')\right|\\
&\quad\lesssim\frac{1}{(V_{\rho_\#})_t(x)}
+\frac{1}{(V_{\rho_\#})_t(x')}
\quad\text{by Definition~\ref{DefATI}(i)}\\
&\quad\lesssim\left[\frac{\rho_{\#}(x,x')
\rho_{\#}(y,y')}{t^2}\right]^\varepsilon
\left[\frac{1}{(V_{\rho_\#})_t(x)}
+\frac{1}{(V_{\rho_\#})_t(x')}\right].
\end{align*}
This is the desired estimate.

\emph{Case (3)}
$\rho_{\#}(x,x')>2C_{\rho_{\#}}^2t$
and $\rho_{\#}(y,y')\leq2C_{\rho_{\#}}^2t$.
In this case,
\begin{align*}
&\left|\left[S_t(x,y)-S_t(x',y)\right]
-\left[S_t(x,y')-S_t(x',y')\right]\right|\\
&\quad\leq\left|S_t(y,x)-S_t(y',x)\right|+
\left|S_t(y,x')-S_t(y',x')\right|
\quad\text{by Definition~\ref{DefATI}(iv)}\\
&\quad\lesssim
\left[\frac{\rho_{\#}(y,y')}{t}\right]^\varepsilon
\left[\frac{1}{(V_{\rho_\#})_t(x)}
+\frac{1}{(V_{\rho_\#})_t(x')}\right]
\quad\text{by \eqref{2y3new}}
\\
&\quad\lesssim\left[\frac{\rho_{\#}(x,x')\rho_{\#}(y,y')}{t^2}\right]^\varepsilon
\left[\frac{1}{(V_{\rho_\#})_t(x)}
+\frac{1}{(V_{\rho_\#})_t(x')}\right].
\end{align*}
This is the desired estimate.

\emph{Case (4)}
$\rho_{\#}(x,x')\leq2C_{\rho_{\#}}^2t$
and $\rho_{\#}(y,y')>2C_{\rho_{\#}}^2t$.
In this case,
from an argument
similar to that used in the proof of
Case (3), it is easy to verify that
the desired estimate in the present case holds.
To summarize, we finish the proof of Step 1.

\emph{Step 2.}
In this step, we aim to prove the existence of
$\{\mathcal{S}_t\}_{t\in(0,\infty)}$
of any order $\varepsilon\in(0,\varepsilon_0]$
in the case $0<\varepsilon_0\preceq{\mathrm{ind\,}}(X,\mathbf{q})$.
To this end,
note that it is easy to see from Definition~\ref{DefATI} and \eqref{935}
that, if $\{\mathcal{S}_t\}_{t\in(0,\infty)}$
is an approximation of the identity
of order $\eta\in(0,\infty)$ with $\theta=1$
on $(X,\mathbf{q},\mu)$,
then, for any $\gamma\in(0,\infty)$,
$\{\mathcal{S}_{t}\}_{t\in(0,\infty)}$
is an approximation of the identity of order
$\frac{\eta}{\gamma}$ with $\theta=\gamma$ on $(X,\mathbf{q}^\gamma,\mu)$,
where $\mathbf{q}^\gamma:=\{\rho^\gamma:\rho\in\mathbf{q}\}$.
Now, choose $\gamma\in(0,\infty)$
sufficiently large such that $\varepsilon_0\in(0,\gamma)$.
By this, Lemma~\ref{11111}, and Step~1, we conclude that
there exists an approximation of the identity
$\{{\mathcal{S}}_{t}\}_{t\in(0,\infty)}$
of any order $\varepsilon\in(0,\frac{\varepsilon_0}{\gamma}]$ with $\theta=1$
on $(X,\mathbf{q}^\gamma,\mu)$,
which, combined with the above claim, further implies that the family
$\{{\mathcal{S}}_{t}\}_{t\in(0,\infty)}$
is an approximation of the identity with $\theta=\frac{1}{\gamma}$
of order $\varepsilon\gamma$
for any $\varepsilon\in(0,\frac{\varepsilon_0}{\gamma}]$
on $(X,(\mathbf{q}^\gamma)^\frac{1}{\gamma},\mu)$; that is,
$\{{\mathcal{S}}_{t}\}_{t\in(0,\infty)}$
is an approximation of the identity
of any order $\varepsilon\in(0,\varepsilon_0]$ with $\theta=\frac{1}{\gamma}$
on $(X,\mathbf{q},\mu)$.
This finishes the proof of Step 2
and hence Theorem~\ref{ThmATI}.
\end{proof}

\begin{remark}
\label{RemThmATI}
\begin{enumerate}
\item[\textup{(i)}]
From the proof of \eqref{2y3new}
and Definition~\ref{DefATI}, it follows that,
for any $\{\mathcal{S}_t\}_{t\in(0,\infty)}$
as in Definition~\ref{DefATI},
$t\in(0,\infty)$, and $x,y,x'\in X$,
\begin{align*}
\left|S_t(x,y)-S_t(x',y)\right|
\lesssim\left[\frac{\rho(x,x')}{t}\right]^\varepsilon
\frac{1}{(V_\rho)_t(y)}.
\end{align*}
\item[\textup{(ii)}]
The proof of Theorem~\ref{ThmATI} relies heavily on the
properties of the quasi-ultrametric as defined in Definition~\ref{D2.1}
and on the quasi-ultrametric space structure $\mathbf{q}$.
More precisely, the reason why the order $\varepsilon$
of the approximation of the identity can be extended to such a
wide range is because of Lemma~\ref{L2.3}(iii) and
Lemma~\ref{11111}. On the one hand, Lemma~\ref{L2.3}(iii)
makes it possible that $\varepsilon$ can be bigger than 1
and even belong to $(0,\infty)$ if we choose some special $\mathbf{q}$.
On the other hand, Lemma~\ref{11111} helps
us to extend the range of $\varepsilon\in(0,1]$
to $0<\varepsilon\preceq{\mathrm{ind\,}}(X,\mathbf{q})$
whenever ${\mathrm{ind\,}}(X,\mathbf{q})\in(1,\infty]$.
However, unlike the constant $C_\rho$ in \eqref{Crho}, the constant $A_{\rho}$ in
the quasi-triangle inequality
as in Remark~\ref{remark2.2}(iv) does not have the property
that $(A_{\rho})^\alpha=A_{\rho^\alpha}$
for any $\alpha\in(0,\infty)$ and hence a conclusion similar to
Lemma~\ref{11111} for quasi-metrics does not hold.
From this perspective, it is crucial for us to work with the
quasi-ultrametric inequality instead of the quasi-metric
inequality in order to achieve the full strength
of our results.
\item[\textup{(iii)}]
Recall that, if $(X,\rho,\mu)$ is a doubling metric measure space,
then there always exists
an approximation of the identity with any order $\varepsilon\in(0,1]$
(see, for instance, \cite[Theorem 2.6]{hmy08}).
By Theorem~\ref{ThmATI}, we find that,
if $\mathbf{q}$ contains a metric,
then ${\mathrm{ind\,}}(X,\mathbf{q})\in[1,\infty]$
and hence, for any given $0<\varepsilon_0
\preceq{\mathrm{ind\,}}(X,\mathbf{q})$, there always exists
an approximation of the identity of any order $\varepsilon\in(0,\varepsilon_0]$
and that, if $\mathbf{q}$ contains an ultrametric,
which is a special case of metrics, then ${\mathrm{ind\,}}(X,\mathbf{q})=\infty$
and hence, for any given $\varepsilon_0\in(0,\infty)$, there always exists
an approximation of the identity of any order $\varepsilon\in(0,\varepsilon_0]$.
Thus, Theorem~\ref{ThmATI}
is sharper than \cite[Theorem 2.6]{hmy08}.
\item[\textup{(iv)}]
By the proof of Theorem~\ref{ThmATI}, we find that
Theorem~\ref{ThmATI} holds in a more general measure
theoretic setting in which the measure
of a singleton is allowed to be strictly positive.
\item[\textup{(v)}]
We point out that \cite[Definition 3.21 and Theorem 3.22]{AlMi2015}
are, respectively, special cases of Definition~\ref{DefATI} and
Theorem~\ref{ThmATI}
when $(X,\mathbf{q},\mu)$ is an Ahlfors-regular quasi-ultrametric space.
\end{enumerate}
\end{remark}

The following proposition shows that, for any
$t\in(0,\infty)$, $\mathcal{S}_tf$
is a $\mu$-measurable function when $f$ is ``nice enough".

\begin{proposition}\label{1507}
Let $(X,\mathbf{q},\mu)$ be a quasi-ultrametric space of homogeneous type.
Assume that $\{\mathcal{S}_t\}_{t\in(0,\infty)}$ is
an approximation of the identity
of order $\varepsilon\in(0,\infty)$.
Then, for any $f\in L_{\mathrm{loc}}^1(X)$ and $t\in(0,\infty)$,
the function $\mathcal{S}_tf$ is lower semicontinuous.
Moreover, $\mathcal{S}_tf$ is $\mu$-measurable.
\end{proposition}

\begin{proof}
Let $\rho\in\mathbf{q}$ be the same as in Definition~\ref{DefATI}.
Let $f\in L_{\mathrm{loc}}^1(X)$ and $t\in(0,\infty)$.
To show that $\mathcal{S}_tf$ is lower semicontinuous,
it suffices to prove that,
for any given $a\in\mathbb{R}$ and $x_0\in X$
satisfying that $\mathcal{S}_tf(x_0)>a$,
there exists $\delta\in(0,\infty)$ such that
$\mathcal{S}_tf(x)>a$
for any $x\in B_\rho(x_0,\delta)$.
To this end, using Definition~\ref{D2.1}(iii),
we find that,
for any $x\in B_\rho(x_0,t)$,
$$
B_\rho(x_0,C't)\cup
B_\rho(x,C't)\subset B_\rho\left(x_0,\max\{1,\,C'\}C_\rho t\right),
$$
where $C'$ and $C_\rho$ are the same
constants as, respectively, in Definition~\ref{DefATI}(i) and \eqref{Crho}.
From this, both (i) and (ii) of Definition~\ref{DefATI},
and \eqref{upperdoub},
we deduce that, for any $x\in B_\rho(x_0,t)$,
\begin{align*}
&\left|\mathcal{S}_tf(x_0)-\mathcal{S}_tf(x)\right|\\
&\quad\leq\int_{B_\rho(x_0,\max\{1,\,C'\}C_\rho t)}
\left|S_t(x_0,y)-S_t(x,y)\right||f(y)|\,d\mu(y)\\
&\quad\leq C\left[\frac{\rho(x_0,x)}{t}\right]^\varepsilon
\left[\frac{1}{(V_\rho)_t(x_0)}+\frac{1}{(V_\rho)_t(x)}\right]
\int_{B_\rho(x_0,\max\{1,\,C'\}C_\rho t)}|f(y)|\,d\mu(y)\\
&\quad\leq\frac{C[1+C_\mu(C_\rho)^n]}{(V_\rho)_t(x_0)}
\left[\frac{\rho(x_0,x)}{t}\right]^\varepsilon
\int_{B_\rho(x_0,\max\{1,\,C'\}C_\rho t)}|f(y)|\,d\mu(y),
\end{align*}
where $C$ and $C_\mu$ are the same positive constants as, respectively, in
Definition~\ref{DefATI}(ii) and \eqref{doublingcond}.
Choose $\delta\in(0,\infty)$ such that
\begin{align*}
&\frac{C[1+C_\mu(C_\rho)^n]}{(V_\rho)_t(x_0)}
\left(\frac{\delta}{t}\right)^\varepsilon
\int_{B_\rho(x_0,\max\{1,\,C'\}C_\rho t)}
|f(y)|\,d\mu(y)<\mathcal{S}_tf(x_0)-a.
\end{align*}
Then we have, for any $x\in B_\rho(x_0,\delta)$,
\begin{align*}
\mathcal{S}_tf(x)&\ge\mathcal{S}_tf(x_0)-
\left|\mathcal{S}_tf(x_0)-\mathcal{S}_tf(x)\right|
>\mathcal{S}_tf(x_0)-[\mathcal{S}_tf(x_0)-a]=a,
\end{align*}
which completes the proof of that
$\mathcal{S}_tf$ is lower semicontinuous.

Next, we show that
$\mathcal{S}_tf$ is a $\mu$-measurable function.
If $\lambda\in\mathbb{R}$ is such that
$$
\{x\in X:\mathcal{S}_tf(x)>\lambda\}=X
\ \text{or}\
\{x\in X:\mathcal{S}_tf(x)>\lambda\}=\varnothing,
$$
then there is nothing to prove.
Now, assume that $\lambda\in\mathbb{R}$ is such that
$\{x\in X:\mathcal{S}_tf(x)>\lambda\}$
is a proper and nonempty subset of $X$.
Moreover, by the above proof, we find that
$\{x\in X:\mathcal{S}_tf(x)>\lambda\}$ is open.
Since $(X,\mathbf{q},\mu)$ is a quasi-ultrametric
space of homogeneous type,
it follows from \cite[Proposition~3.28]{AlMi2015}
[see also Remark~\ref{940}(iii) below] that
$(X,\mathbf{q})$ is also a geometrically doubling
quasi-ultrametric space as in Definition~\ref{2248}.
By this and the Whitney decomposition
theorem for open sets \cite[Theorem~2.4]{AlMi2015}
[see also Lemma~\ref{Wtd}(i) below],
we find that the open set
$\{x\in X:\mathcal{S}_tf(x)>\lambda\}$
can be decomposed into a countable union of $\rho$-balls,
which are obviously $\mu$-measurable because the choice of $\rho\in\mathbf{q}$ is
such that all $\rho$-balls are $\mu$-measurable.
This finishes the proof of Proposition~\ref{1507}.
\end{proof}

In the next main result of this section (Theorem~\ref{corATI}),
we collect some important properties
that are satisfied by any approximation of the identity
defined in Definition~\ref{DefATI}. To state this result, we first need
to recall the definition of the H\"older space.
Throughout this monograph, given a function $f:X\to\mathbb{C}$,
the \emph{support of $f$} is defined by setting
\begin{align*}
{\mathrm{supp\,}}f\index[symbols]{S@${\mathrm{supp\,}}f$}:=
\left\{x\in X:f(x)\neq 0\right\}.\index[words]{support of a funcion}
\end{align*}

\begin{definition}\label{1454}
Let $(X,\mathbf{q})$ be a quasi-ultrametric
space, $\rho\in\mathbf{q}$, and $\beta\in(0,\infty)$.
The \emph{H\"older seminorm}\footnote{
Let $X$ be a vector space over $\mathbb{C}$.
A function $\|\cdot\|:X\to\mathbb{R}_+$
is called a \emph{seminorm} if,
for any $x,y\in X$, the following
three conditions hold:
(i) $x=0$ implies $\|x\|=0$,
(ii) for any $\lambda\in\mathbb{C}$,
$\|\lambda x\|=|\lambda|\|x\|$,
(iii) $\|x+y\|\leq\|x\|+\|y\|$.}\index[words]{H\"older seminorm}
$\|\cdot\|_{\dot{C}^\beta(X,\,\rho)}$ of order $\beta$,
related to the quasi-ultrametric $\rho$, is defined by setting,
for any given complex-valued function $f$ defined on $X$,
\begin{align}
\label{C(X)}
{\lVert f\rVert}_{\dot{C}^\beta(X,\,\rho)}
:=\sup_{\genfrac{}{}{0pt}{}{x,y\in X}{x\ne y}}
\frac{|f(x)-f(y)|}{[\rho(x,y)]^\beta}.
\end{align}
The \emph{homogeneous H\"older
space}\index[words]{homogeneous H\"older space}
$\dot{C}^\beta(X,\mathbf{q})$\index[symbols]{C@$\dot{C}^\beta(X,\mathbf{q})$}
is defined by setting
\begin{align*}
\dot{C}^\beta(X,\mathbf{q})
&:=\left\{f:\|f\|_{\dot{C}^\beta(X,\,\rho)}<\infty
\text{ for some }\rho\in\mathbf{q}\right\}
\\
&\ =\left\{f:\|f\|_{\dot{C}^\beta(X,\,\rho)}<\infty
\text{ for all }\rho\in\mathbf{q}\right\}.
\end{align*}
The \emph{inhomogeneous H\"older
space}\index[words]{inhomogeneous H\"older space}
$C^\beta(X,\mathbf{q})$\index[symbols]{C@$C^\beta(X,\mathbf{q})$}
is defined by setting
\begin{align*}
C^\beta(X,\mathbf{q})
&:=\left\{f:\|f\|_{C^\beta(X,\rho)}<\infty
\text{ for some }\rho\in\mathbf{q}\right\}
\\
&\ =\left\{f:\|f\|_{C^\beta(X,\rho)}<\infty
\text{ for all }\rho\in\mathbf{q}\right\},
\end{align*}
where
$\|\cdot\|_{C^\beta(X,\rho)}:
=\|\cdot\|_{L^\infty(X)}
+\|\cdot\|_{\dot{C}^\beta(X,\,\rho)}$.
The \emph{H\"older space with bounded support},
$C^\beta_{\mathrm{b}}(X,\mathbf{q})$,
\index[symbols]{C@$C^\beta_{\mathrm{b}}(X,\mathbf{q})$}
is defined by setting
\begin{align*}
C^\beta_{\mathrm{b}}(X,\mathbf{q})
:=\left\{f\in C^\beta(X,\mathbf{q}):
f\text{ has bounded support}\right\}.\index[words]{H\"older space with bounded support}
\end{align*}
In this monograph, we always simply
write $\dot{C}^\beta(X,\mathbf{q})$,
$C^\beta(X,\mathbf{q})$,
and $C_{\mathrm{b}}^\beta(X,\mathbf{q})$
as $\dot{C}^\beta(X)$,
$C^\beta(X)$,
and $C_{\mathrm{b}}^\beta(X)$,
respectively.
\end{definition}

For the approximation of the identity defined in
Definition~\ref{DefATI}, we have the following conclusions.

\begin{theorem}\label{corATI}
Let $(X,\mathbf{q},\mu)$ be a
quasi-ultrametric space of homogeneous type.
Assume that $\{\mathcal{S}_t\}_{t\in(0,\infty)}$ is
an approximation of the identity
of order $\varepsilon\in(0,\infty)$.
Then the following statements hold:
\begin{enumerate}
\item [\textup{(i)}]
For any $p\in [1,\infty]$,
$\sup_{t\in(0,\infty)}\|\mathcal{S}_t\|_{L^p(X)
\to L^p(X)}<\infty$.
\item [\textup{(ii)}]
$\sup_{t\in(0,\infty)}
\|\mathcal{S}_t\|_{\dot{C}^\varepsilon(X)
\to\dot{C}^\varepsilon(X)}<\infty$.
\item [\textup{(iii)}]
For any $g\in\dot{C}^\varepsilon(X)$
and $\alpha\in(0,\varepsilon)$,
$\lim_{t\to 0^+}\|\mathcal{S}_tg-g\|_{\dot{C}^\alpha(X)}=0$.
\item [\textup{(iv)}]
$\sup_{t\in(0,\infty)}
\|\mathcal{S}_t\|_{{C}^\varepsilon(X)
\to{C}^\varepsilon(X)}<\infty$.
\item [\textup{(v)}]
For any $g\in {C}^\varepsilon(X)$
and $\alpha\in(0,\varepsilon)$,
$\lim_{t\to 0^+}\|\mathcal{S}_tg-g\|_{C^\alpha(X)}=0$.
\item [\textup{(vi)}]
There exists a constant $C\in(0,\infty)$ such that,
for any $t\in(0,\infty)$
and $f\in L^\infty(X)$,
\begin{align*}
\left\|\mathcal{S}_tf\right\|_
{\dot{C}^\varepsilon(X)}
\leq C t^{-\varepsilon}
\left\|f\right\|_{L^\infty(X)}.
\end{align*}
\item [\textup{(vii)}]
If $\mathrm{diam\,}(X)=\infty$, then,
for any given $p\in(1,\infty)$ and
for any $f\in L^p(X)$,
\begin{align*}
\lim_{t\to\infty}\left\|\mathcal{S}_tf\right\|_{L^p(X)}=0.
\end{align*}
\item [\textup{(viii)}]
Let $\rho_0\in\mathbf{q}$,
$x_0\in X$, and $r_0\in(0,\infty)$.
If $f\in L_{\mathrm{loc}}^1(X)$ is such that
${\mathrm{supp\,}}f\subset B_{\rho_0}(x_0,r_0)$,
then there exists a positive constant $C$, independent of
$f$, $x_0$, and $r_0$, such that,
for any $t\in(0,\infty)$,
${\mathrm{supp\,}}\mathcal{S}_tf\subset
B_{\rho_0}(x_0,C(r_0+t))$.
\item [\textup{(ix)}]
Let $p\in[1,\infty)$.
Then
\begin{equation}\label{2.12}
\lim_{t\to 0^+}\left\|\mathcal{S}_tf-f\right\|_{L^p(X)}=0
\end{equation}
for any $f\in L^p(X)$
if and only if $\mu$ is a
Borel-semiregular measure on $X$.
\end{enumerate}
\end{theorem}

\begin{remark}
\begin{enumerate}
\item[\textup{(i)}]
In the proof of Theorem~\ref{corATI}(i),
we find that the supremum in Theorem~\ref{corATI}(i)
is dominated by a constant independent of $p\in[1,\infty]$.

\item[\textup{(ii)}]
The Borel-semiregular assumption on $\mu$ in Theorem~\ref{corATI}(ix)
is weaker than that of \cite[Proposition~2.7]{hmy08}
in which $\mu$ is
assumed to be Borel-regular.
Indeed,
Theorem~\ref{corATI}(ix) is sharp
due to Lemma~\ref{LDT}.

\item[\textup{(iii)}]
If $\mu$ is a Borel measure on $X$,
then Theorem~\ref{corATI} holds
with the same proof except \eqref{2.12}.
In this case, \eqref{2.12} holds only for any $p\in[1,\infty]$
and $f$ belonging to the closure of
$\dot{C}^\varepsilon(X)\cap L^p(X)$
in $L^p(X)$.
\end{enumerate}
\end{remark}

To prove Theorem~\ref{corATI}, we need some technical lemmas.
Let $(X,\mathbf{q},\mu)$ be a quasi-ultrametric space of homogeneous
type and suppose that $\rho\in\mathbf{q}$ has the property that
all $\rho$-balls are $\mu$-measurable. In this context, recall that, for any
$f\in L^1_\mathrm{loc}(X)$, the
\emph{Hardy--Littlewood maximal function}
$\mathcal{M}_\rho f$\index[symbols]{M@$\mathcal{M}_\rho$}
of $f$ (with respect to $\rho$) is defined by setting,
for any $x\in X$,
\begin{equation}\label{Mf}
\mathcal{M}_\rho f(x)
:=\sup_{r\in(0,\infty)}
\frac{1}{(V_{\rho})_r(x)}\int_{B_\rho(x,r)}
|f(y)|\,d\mu(y).
\end{equation}
It is easy to see from \eqref{upperdoub}
that, if $\rho'\in\mathbf{q}$
is any other quasi-ultrametric with the property
that all $\rho'$-balls are $\mu$-measurable,
then, for any $f\in L^1_\mathrm{loc}(X)$ and $x\in X$,
\begin{align}\label{eq-538}
\mathcal{M}_{\rho'} f(x)\sim\mathcal{M}_\rho f(x),
\end{align}
where the positive equivalence constants
are independent of $f$ and $x$. Although
equivalent quasi-ultrametrics induce Hardy--Littlewood
maximal functions which are pointwise comparable,
the particular choice of the quasi-ultrametric
becomes an important issue when establishing
the measurability of function $\mathcal{M}_\rho f$.

The following lemma is a part of \cite[Theorem 3.7]{AlMi2015}.

\begin{lemma}\label{M}
Let $(X,\mathbf{q},\mu)$
be a quasi-ultrametric space of homogeneous type.
Assume that
$\rho\in\mathbf{q}$ is a
regularized quasi-ultrametric
as in Definition~\ref{D2.3}.
Then, for any $f\in L^1_\mathrm{loc}(X)$, the function
$\mathcal{M}_{\rho} f:X\to[0,\infty]$ is well defined
and $\mu$-measurable. Moreover,
for any $p\in(1,\infty]$,
$$
\mathcal{M}_{\rho}:L^p(X)\to L^p(X)
$$
and
$$
\mathcal{M}_{\rho}:L^1(X)\to L^{1,\infty}(X)
$$
are well defined and bounded.
\end{lemma}

Let $(X,\mathbf{q},\mu)$ be a quasi-ultrametric space of homogeneous
type and suppose that $\rho\in\mathbf{q}$ has the property that
all $\rho$-balls are $\mu$-measurable. The maximal
operator defined in \eqref{Mf} is often referred to
as the \emph{centered} Hardy--Littlewood maximal operator.
A closely related version of this is the so-called
\emph{uncentered Hardy--Littlewood maximal operator},
which is defined by setting, for any $f\in L^1_{\rm loc}(X)$ and $x\in X$,
\begin{align*}
\widetilde{\mathcal{M}}_\rho f(x)
:=\sup_{B\ni x}
\frac{1}{\mu(B)}\int_{B}
|f(y)|\,d\mu(y),
\end{align*}
where the supremum is taken over all $\rho$-balls $B\subset X$ containing the point $x$.
Note that, for any $f\in L^1_\mathrm{loc}(X)$, the function
$\widetilde{\mathcal{M}}_\rho f:X\to[0,\infty]$ is $\mu$-measurable
since $\mu$ is a Borel measure [see Remark~\ref{Rmk502}(ii)] and sets of the form
$\{x\in X:\widetilde{\mathcal{M}}_\rho f(x)>\lambda\}$, $\lambda\in(0,\infty)$,
are open in the topology $\tau_\rho$. Moreover,
there exists a positive constant $C$, depending only
on the doubling constant of $\mu$, such that, for any $x\in X$,
\begin{equation}\label{HLmaxequiv}
({\mathcal{M}_{\rho}}f)(x)\leq
(\widetilde{\mathcal{M}}_{\rho}f)(x)\leq C
({\mathcal{M}_{\rho}}f)(x).
\end{equation}
Finally, it follows from \eqref{upperdoub}
that, if $\rho'\in\mathbf{q}$
is any other quasi-ultrametric with the property
that all $\rho'$-balls are $\mu$-measurable,
\begin{align*}
\widetilde{\mathcal{M}}_{\rho'} f(x)\sim\widetilde{\mathcal{M}}_\rho f(x),
\end{align*}
where the positive equivalence constants
are independent of $f$ and $x$.

The following lemma is a sharp variant of the Lebesgue
differentiation theorem that was proved
in \cite[Theorem~3.14]{AlMi2015}.
\begin{lemma}\label{LDT}
Let $(X,\mathbf{q},\mu)$ be
a quasi-ultrametric space of homogeneous type and assume that $\rho\in\mathbf{q}$
is a
regularized quasi-ultrametric as in Definition~\ref{D2.3}.
Then the following assertions are equivalent:
\begin{enumerate}
\item[\textup{(i)}]
$\mu$ is Borel-semiregular.
\item[\textup{(ii)}]
For any $f\in L^1_{\mathrm{loc}}(X)$ and
for $\mu$-almost every $x\in X$,
$$
\lim_{r\to 0^+}\frac{1}{(V_{\rho})_r(x)}
\int_{B_{\rho}(x,r)}
\left|f(y)-f(x)\right|\,d\mu
(y)=0.
$$
\item[\textup{(iii)}]
For any $p\in(0,\infty)$ and $0<\beta\preceq{\mathrm{ind\,}}(X,\mathbf{q})$,
$C^\beta_{\mathrm{b}}(X)$ continuously embeds into $L^p(X)$
and is dense in $L^p(X)$.
\item[\textup{(iv)}]
For any $p\in(0,\infty)$,
$C_{\mathrm{b}}(X)\cap L^p(X)$\index[symbols]{C@$C_{\mathrm{b}}(X)$}
continuously embeds into $L^p(X)$
and is dense in $L^p(X)$,
where the \emph{symbol $C_{\mathrm{b}}(X)$} denotes
the space of all continuous functions with bounded support.
\end{enumerate}
\end{lemma}

Now, we are ready to show Theorem~\ref{corATI}.

\begin{proof}[Proof of Theorem~\ref{corATI}]
Fix $\rho\in\mathbf{q}$.
By Remark~\ref{Rm}(iii),
we find that the kernels of $\{\mathcal{S}_t\}_{t\in(0,\infty)}$ satisfy
(i)--(v) in Definition~\ref{DefATI} with the regularized quasi-ultrametric
of $\rho$ defined in Definition~\ref{D2.3}. Thus, without loss of generality,
we may assume that $\rho$ is a regularized quasi-ultrametric.
In addition, using Proposition~\ref{1507},
we conclude that, for any $t\in(0,\infty)$,
$\mathcal{S}_tf$ is integrable.

Moving on, we first prove (i).
By Definition~\ref{DefATI}(i)
and \eqref{upperdoub},
we find that,
for any $t\in(0,\infty)$, $f\in L^\infty(X)$, and $x\in X$,
$\mathcal{S}_tf(x)$ is well defined and
\begin{align}
\label{atimax}
\left|\mathcal{S}_tf(x)\right|
&\leq\int_{B_{\rho(x,C't)}}
S_t(x,y)\left|f(y)\right|\,d\mu(y)
\nonumber\\
&\lesssim\frac{(V_\rho)_{C't}(x)}{(V_\rho)_{t}(x)}
\frac{1}{(V_\rho)_{C't}(x)}
\int_{B_\rho(x,C't)}\left|f(y)\right|\,d\mu(y)
\lesssim\mathcal{M}_\rho f(x),
\end{align}
where $C'\in[1,\infty)$ is the same positive constant as in Definition~\ref{DefATI}
and
$\mathcal{M}_\rho f$ is
the Hardy--Littlewood maximal function of $f$
defined in \eqref{Mf}.
Note that the implicit positive constants in
\eqref{atimax} are independent of $t$.
From \eqref{atimax}
and Lemma~\ref{M}, we deduce that (i) holds when $p=\infty$.
On the other hand, by (i) and (v) of Definition~\ref{DefATI}
and Tonelli's theorem, we conclude that,
for any $t\in(0,\infty)$ and $f\in L^1(X)$,
\begin{align}
\label{StL1}
\left\|\mathcal{S}_tf\right\|_{L^1(X)}
&\leq\int_{X}\int_{X}
S_t(x,y)\left|f(y)\right|\,
d\mu(y)\,d\mu(x)
\nonumber\\
&=\int_{X}\left[\int_{X}S_t(x,y)\,d\mu(x)\right]
\left|f(y)\right|\,d\mu(y)
=\|f\|_{L^1(X)},
\end{align}
which implies that (i) holds when $p=1$.
By using the Riesz--Thorin
interpolation theorem (see, for instance,
\cite[Theorem~1.19]{Duo})
between $p=1$ and $p=\infty$
of operator $\mathcal{S}_t$ in $L^p(X)$,
we further find that (i) holds for any $p\in(1,\infty)$.

Next, we turn to prove (ii). Let
$f\in\dot{C}^\varepsilon(X)$ and $t\in(0,\infty)$.
On the one hand, by Definition~\ref{D2.1}(iii), we obtain, for any
$x,y\in X$ satisfying $\rho(x,y)<C't$,
\begin{equation}
\label{zj.48}
\left[B_\rho(x,C't)\cup B_\rho(y,C't)\right]\subset B_\rho(x,C_\rho C't).
\end{equation}
From this and both (i) and (v) of Definition~\ref{DefATI},
we infer that, for any
$x,y\in X$ satisfying $\rho(x,y)<C't$,
\begin{align}
\label{2.25}
&\left|\mathcal{S}_tf(x)-\mathcal{S}_tf(y)\right|\nonumber\\
&\quad\leq\int_{B_\rho(x,C_\rho C't)}
\left|S_t(x,z)-S_t(y,z)\right|
\left|f(z)-f(x)\right|\,d\mu(z)
\nonumber\\
&\quad\lesssim\int_{B_\rho(x,C_\rho C't)}
\frac{[\rho(x,y)]^\varepsilon}
{t^{\varepsilon}(V_\rho)_t(z)}
\|f\|_{\dot{C}^\varepsilon(X)}
\left[\rho(x,z)\right]^\varepsilon\,d\mu(z)
\quad\text{by \eqref{2y3new} and \eqref{C(X)}}
\nonumber\\
&\quad\sim\frac{[\rho(x,y)]^\varepsilon}
{t^{\varepsilon}(V_\rho)_t(x)}
\|f\|_{\dot{C}^\varepsilon(X)}\int_{B_\rho(x,C_\rho C't)}
\left[\rho(x,z)\right]^\varepsilon\,d\mu(z)
\quad\text{by Lemma~\ref{A3.2}(vii)}
\nonumber\\
&\quad\lesssim\frac{(V_\rho)_{C_\rho C't}(x)}{(V_\rho)_t(x)}
\left[\rho(x,y)\right]^\varepsilon
\|f\|_{\dot{C}^\varepsilon(X)}
\lesssim\left[\rho(x,y)\right]^\varepsilon
\|f\|_{\dot{C}^\varepsilon(X)}
\quad\text{by \eqref{upperdoub}}.
\end{align}
On the other hand,
by both (i) and (v) of Definition~\ref{DefATI}
and \eqref{C(X)},
we conclude that,
for any $x\in X$,
\begin{align}\label{2.27}
\left|\mathcal{S}_tf(x)-f(x)\right|
&\leq\int_{B_\rho(x,C't)}
S_t(x,z)\left[\rho(x,z)\right]^\varepsilon
\|f\|_{\dot{C}^\varepsilon(X)}\,d\mu(z)
\nonumber\\
&\lesssim t^\varepsilon\|f\|_{\dot{C}^\varepsilon(X)}
\int_{B_\rho(x,C't)}
S_t(x,z)\,d\mu(z)
\leq t^\varepsilon
\|f\|_{\dot{C}^\varepsilon(X)}.
\end{align}
From this,
we deduce that, for any $x,y\in X$ satisfying $\rho(x,y)\ge C't$,
\begin{align*}
\left|\mathcal{S}_tf(x)-\mathcal{S}_tf(y)\right|
&\leq\left|\mathcal{S}_tf(x)-f(x)\right|
+\left|f(x)-f(y)\right|
+\left|f(y)-\mathcal{S}_tf(y)\right|
\nonumber\\
&\lesssim t^\varepsilon
\|f\|_{\dot{C}^\varepsilon(X)}
+\left[\rho(x,y)\right]^\varepsilon
\|f\|_{\dot{C}^\varepsilon(X)}
\lesssim\left[\rho(x,y)\right]^\varepsilon
\|f\|_{\dot{C}^\varepsilon(X)},
\end{align*}
which, combined with \eqref{2.25},
further implies that, for any $x,y\in X$,
$$|\mathcal{S}_tf(x)-\mathcal{S}_tf(y)|\lesssim[\rho(x,y)]^\varepsilon
\|f\|_{\dot{C}^\varepsilon(X)},$$
where the implicit positive constant is independent of $t$.
This finishes the proof of (ii).

Now, we show (iii).
Let $\alpha\in(0,\varepsilon)$,
$g\in\dot{C}^\varepsilon(X)$, and $t\in(0,\infty)$.
For any $x,y\in X$ satisfying $\rho(x,y)<C't$,
\begin{align}\label{2.28}
&\left|\left(\mathcal{S}_tg-g\right)(x)
-\left(\mathcal{S}_tg-g\right)(y)\right|\nonumber\\
&\quad\leq\left|\mathcal{S}_tg(x)-\mathcal{S}_tg(y)\right|
+\left|g(x)-g(y)\right|
\nonumber\\
&\quad\lesssim\left[\rho(x,y)\right]^\varepsilon
\|g\|_{\dot{C}^\varepsilon(X)}
\lesssim\left[\rho(x,y)\right]^\alpha
t^{\varepsilon-\alpha}\|g\|_{\dot{C}^\varepsilon(X)}
\quad\text{by \eqref{2.25} and \eqref{C(X)}}.
\end{align}
On the other hand,
for any $x,y\in X$ satisfying $\rho(x,y)\ge C't$,
\begin{align}\label{2.29}
&\left|\left(\mathcal{S}_tg-g\right)(x)
-\left(\mathcal{S}_tg-g\right)(y)\right|\nonumber\\
&\quad\leq\left|\mathcal{S}_tg(x)-g(x)\right|
+|\mathcal{S}_tg(y)-g(y)|
\nonumber\\
&\quad\lesssim t^\varepsilon
\|g\|_{\dot{C}^\varepsilon(X)}
\lesssim\left[\rho(x,y)\right]^\alpha
t^{\varepsilon-\alpha}
\|g\|_{\dot{C}^\varepsilon(X)}
\quad\text{by \eqref{2.27}}.
\end{align}
Combining \eqref{2.28} and \eqref{2.29},
we conclude that, for any $\alpha\in(0,\varepsilon)$,
$g\in\dot{C}^\varepsilon(X)$, and $t\in(0,\infty)$,
$$
\left\|\mathcal{S}_tg-g\right\|_{\dot{C}^\alpha(X)}
\lesssim t^{\varepsilon-\alpha}
\|g\|_{\dot{C}^\varepsilon(X)},
$$
which readily yields (iii) as $t\to 0^+$.

Note that (iv) is a direct consequence of
both (i) with $p=\infty$ and (ii).
We then turn to prove (v).
By \eqref{2.27},
we find that,
for any $t\in(0,\infty)$ and $g\in\dot{C}^\varepsilon(X)$,
\begin{align*}
\left\|\mathcal{S}_tg-g\right\|_{L^\infty(X)}
&\leq\sup_{x\in X}
\left|\left(\mathcal{S}_tg-g\right)(x)\right|
\lesssim t^\varepsilon\|g\|_{\dot{C}^\varepsilon(X)},
\end{align*}
which, combined with (iii), further implies (v).

Next, we show (vi).
Let $t\in(0,\infty)$, $f\in L^\infty(X)$, and $x,y\in X$.
If $\rho(x,y)<C't$, then, from Lemma~\ref{A3.2}(vii), we infer that
$(V_\rho)_t(x)\sim(V_\rho)_t(y)$,
which further implies that
\begin{align}
\label{2.22}
&\left|\mathcal{S}_tf(x)-\mathcal{S}_tf(y)\right|\nonumber\\
&\quad\leq\int_{B_\rho(x,C_\rho C't)}
\left|S_t(x,z)-S_t(y,z)\right|
\left|f(z)\right|\,d\mu(z)
\quad\text{by Definition~\ref{DefATI}(i) and \eqref{zj.48}}
\nonumber\\
&\quad\lesssim\left[\frac{\rho(x,y)}{t}\right]^\varepsilon
\left[\frac{1}{(V_\rho)_t(x)}
+\frac{1}{(V_\rho)_t(y)}\right]\int_{B_\rho(x,C_\rho C't)}
\left|f(z)\right|\,d\mu(z)
\quad\text{by Definition~\ref{DefATI}(ii)}
\nonumber\\
&\quad\lesssim\left[\frac{\rho(x,y)}{t}\right]^\varepsilon
\frac{(V_\rho)_{C_\rho C't}(x)}
{(V_\rho)_t(x)}\|f\|_{L^\infty(X)}
\lesssim\left[\frac{\rho(x,y)}{t}\right]^\varepsilon
\|f\|_{L^\infty(X)}
\quad\text{by \eqref{upperdoub}},
\end{align}
where the implicit positive constants are independent of
both $t$ and $f$.
On the other hand,
if $\rho(x,y)\ge C't$, then, by (i) and (v) of Definition~\ref{DefATI},
we find that
\begin{align*}
\left|\mathcal{S}_tf(x)
-\mathcal{S}_tf(y)\right|
&\lesssim\left[\frac{\rho(x,y)}{t}\right]^\varepsilon
\left[\left|\mathcal{S}_tf(x)\right|
+\left|\mathcal{S}_tf(y)\right|\right]
\lesssim\left[\frac{\rho(x,y)}{t}\right]^\varepsilon
\|f\|_{L^\infty(X)},
\end{align*}
where the implicit positive constants are independent of
both $t$ and $f$.
Combining this and \eqref{2.22},
we further conclude that
\begin{align*}
\left\|\mathcal{S}_tf\right\|_{\dot{C}^\varepsilon(X)}
&=\sup_{x,y\in X,\,x\neq y}
\frac{|\mathcal{S}_tf(x)-\mathcal{S}_tf(y)|}
{[\rho(x,y)]^\varepsilon}
\lesssim t^{-\varepsilon}\|f\|_{L^\infty(X)},
\end{align*}
which proves (vi).

Now, we turn to prove (vii).
From $\mathrm{diam}_\rho(X)=\infty$,
Remark~\ref{rm1112}(ii), and \eqref{upperdoub},
we deduce that, for any $x\in X$,
\begin{align}
\label{2230}
(V_\rho)_{t}(x)\to\infty
\end{align}
as $t\to\infty$.
Moreover, by Definition~\ref{DefATI}(i), \eqref{upperdoub},
and H\"older's inequality,
we conclude that, for any $t\in(0,\infty)$,
$f\in L^p(X)$ with $p\in(1,\infty)$,
and $x\in X$,
\begin{align*}
\left|\mathcal{S}_tf(x)\right|
&\leq\int_{X}S_t(x,y)\left|f(y)\right|\,d\mu(y)
\lesssim\frac{(V_\rho)_{Ct}(x)}{(V_\rho)_{t}(x)}\frac{1}{(V_\rho)_{Ct}(x)}
\int_{B_\rho(x,Ct)}
\left|f(y)\right|\,d\mu(y)
\\
&\lesssim\frac{1}{(V_\rho)_{Ct}(x)}
\int_{B_\rho(x,Ct)}
\left|f(y)\right|\,d\mu(y)
\leq{\left[\frac{1}
{(V_\rho)_{Ct}(x)}\right]}^{\frac{1}{p}}
\|f\|_{L^p(X)},
\end{align*}
which, combined with \eqref{2230}, further implies that
$|\mathcal{S}_tf(x)|\to 0$
as $t\to\infty$.
Combining this, \eqref{atimax}, Lemma~\ref{M}
[here, we need $p\in(1,\infty)$],
and the Lebesgue dominated convergence theorem
(see, for instance, \cite[Theorem~1.34]{r1987}),
we obtain (vii).

Next, we prove (viii).
Since $\rho_0\in\mathbf{q}$,
it follows that there exists
a constant $c\in[1,\infty)$ such that, for any $x,y\in X$,
$c^{-1}\rho_0(x,y)\leq\rho(x,y)\leq c\rho_0(x,y)$.
Given that the regularized quasi-ultrametric
$(\rho_0)_{\#}\in\mathbf{q}$,
we may assume, without loss of generality, that $\rho_0=(\rho_0)_{\#}$.
From Definition~\ref{DefATI}(i),
we infer that, for any given $t\in(0,\infty)$
and for any $f\in L_{\mathrm{loc}}^1(X)$ and $x\in X$,
$$
\mathcal{S}_tf(x)
=\int_{B_{\rho_0}(x_0,r_0)\cap B_\rho(x,C't)}
S_t(x,y)f(y)\,d\mu(y).
$$
Observe that, if $\mathcal{S}_tf(x)\neq 0$,
then $B_{\rho_0}(x_0,r_0)\cap B_\rho(x,C't)\neq\varnothing.$
Thus, we can choose
$y\in B_{\rho_0}(x_0,r_0)\cap B_\rho(x,C't)$ such that
\begin{align*}
\rho(x_0,x)&\leq C_\rho
\max\left\{\rho(x_0,y),\,\rho(y,x)\right\}
\nonumber\\
&\leq C_\rho\max\left\{c\rho_0(x_0,y),\,\rho(y,x)\right\}
<C_\rho(cr_0+Ct),
\end{align*}
which further implies (viii), given that $C',C_\rho,c\geq1$.

There remains to prove (ix).
We first show the sufficiency.
Assume that $\mu$ is Borel-semiregular.
We consider the following two cases for $p$.

\emph{Case (1)} $p\in(1,\infty)$. In this case,
by (i) and (v) of Definition~\ref{DefATI} and
\eqref{upperdoub}, we find that,
for any $t\in(0,\infty)$, $f\in L^p(X)$, and $x\in X$,
\begin{align*}
\left|\mathcal{S}_tf(x)-f(x)\right|
&\lesssim\frac{1}{(V_\rho)_t(x)}
\int_{B_\rho(x,C't)}\left|f(y)-f(x)\right|\,d\mu(y)
\\
&\lesssim\frac{1}{(V_\rho)_{C't}(x)}
\int_{B_\rho(x,C't)}\left|f(y)-f(x)\right|\,d\mu(y),
\end{align*}
which, combined with Lemma~\ref{LDT}(ii),
further implies that
\begin{align*}
\lim_{t\to 0^+}\mathcal{S}_tf(x)
=f(x).
\end{align*}
From this, \eqref{atimax}, Lemma~\ref{M}
(here, we need $p\in(1,\infty)$),
and the Lebesgue dominated convergence theorem,
we deduce that \eqref{2.12} holds in this case.

\emph{Case (2)} $p=1$. In this case, let $f\in L^1(X)$,
$\delta\in(0,\infty)$, and
$\beta\in(0,\log_2C_\rho)^{-1}]\cap(0,\infty)$.
By the assumption that $\mu$ is Borel-semiregular
and by Lemma~\ref{LDT}(iii), we find that there exists
$g\in C^\beta_\mathrm{b}(X)$ such that
\begin{align}\label{choiceg}
\|f-g\|_{L^1(X)}<\delta.
\end{align}
For such $g$,
from \eqref{2.27}, we infer that, for any $t\in(0,\infty)$,
\begin{align}\label{12.17.x1}
\left\|\mathcal{S}_tg-g\right\|_{L^\infty(X)}
\lesssim t^\beta\|g\|_{\dot{C}^\beta(X)},
\end{align}
where the implicit positive constant
is independent of both $g$ and $t$.
Moreover, since $g$ vanishes
outside of a $\rho$-bounded subset of $X$,
it follows from (viii) that
there exists a bounded subset $B\subset X$
outside of which both $S_tg$ and $g$
vanish for any $t\in(0,1]$.
By this, the fact that $\mu$
is finite on bounded sets, and \eqref{12.17.x1},
we conclude that
\begin{align}\label{2.32}
\left\|\mathcal{S}_tg-g\right\|_{L^1(X)}
&=\int_{B}
\left|\mathcal{S}_tg(x)-g(x)\right|\,d\mu(x)
\nonumber\\
&\leq\mu(B)
\left\|\mathcal{S}_tg-g\right\|_{L^\infty(X)}
\lesssim t^\beta\mu(B)\|g\|_{\dot{C}^\beta(X)}\to 0
\end{align}
as $t\to 0^+$.
In addition, from \eqref{StL1}
and \eqref{choiceg}, we deduce that
\begin{align*}
\left\|\mathcal{S}_tf-f\right\|_{L^1(X)}
&\leq \left\|\mathcal{S}_t(f-g)\right\|_{L^1(X)}
+\left\|\mathcal{S}_tg-g\right\|_{L^1(X)}
+\left\|g-f\right\|_{L^1(X)}\\
&\leq 2\|f-g\|_{L^1(X)}
+\left\|\mathcal{S}_tg-g\right\|_{L^1(X)}
\lesssim\delta+
\left\|\mathcal{S}_tg-g\right\|_{L^1(X)},
\end{align*}
which, together with \eqref{2.32}, further implies that
$$
\limsup_{t\to 0^+}
\left\|\mathcal{S}_tf-f\right\|_{L^1(X)}
\lesssim\delta.
$$
By this and the arbitrariness of $\delta$,
we conclude that	
$\lim_{t\to 0^+}
\left\|\mathcal{S}_tf-f\right\|_{L^1(X)}=0$
which implies that \eqref{2.12} holds in this case.
Combining the above two cases,
we complete the proof of the sufficiency of (ix).

Next, we show the necessity of (ix).
Assume that $p\in[1,\infty)$ and \eqref{2.12}
holds for any $f\in L^p(X)$.
Let $f\in L^p(X)$. In light of Lemma~\ref{LDT},
the goal is to use \eqref{2.12}
to approximate $f$ arbitrarily
well in $L^p(X)$ with functions
from $C_{\mathrm{b}}(X)$.
To this end,
since simple functions are dense in $L^p(X)$,
without loss of generality, we may assume that
$f$ is a characteristic function, ${\mathbf{1}_E}$,
for some $\mu$-measurable set $E\subset X$
with $\mu(E)<\infty$.
We can also assume that $E$ is actually bounded.
Indeed, if we take any point $x_0\in X$,
then, for any $j\in\mathbb{N}$,
$E_j:=E\cap B_{\rho}(x_0,j)$
is bounded and $\mu$-measurable
and satisfies that
$\mathbf{1}_{E_j}\to\mathbf{1}_E$
in $L^p(X)$ as $j\to\infty$.
Since $E$ is bounded, it follows
that ${\mathbf{1}_E}$ has bounded support.
By this and (viii), we find that,
for any $t\in(0,\infty)$, $\mathcal{S}_t\mathbf{1}_E$
also has bounded support.
From this, (vi),
and the fact that ${\mathbf{1}_E}\in L^\infty(X)$,
we infer that, for any $t\in(0,\infty)$,
$\mathcal{S}_t\mathbf{1}_E\in C^\varepsilon_\mathrm{b}(X)
\subset C_\mathrm{b}(X)$.
By \eqref{2.12}, we conclude that
$\mathcal{S}_t\mathbf{1}_E\to\mathbf{1}_E$ in $L^p(X)$
as $t\to0^+$ and, consequently,
${\mathbf{1}_E}$ is approximated arbitrarily well
in $L^p(X)$ by functions from
$C_{\mathrm{b}}(X)$,
and we ultimately deduce that
$C_{\mathrm{b}}(X)$ is dense in $L^p(X)$.
From this and Lemma~\ref{LDT},
we further infer that $\mu$ is Borel-semiregular,
which completes the proof of the necessity and hence (ix).
This finishes the proof of Theorem~\ref{corATI}.
\end{proof}


\chapter[Sharp Calder\'on Reproducing Formulae on Ultra-RD-Spaces]
{Sharp Calder\'on Reproducing Formulae on Ultra-RD-Spaces}\label{Chapter2}
\markboth{\scriptsize\rm\sc Sharp Calder\'on Reproducing Formulae on
Ultra-RD-Spaces}
{\scriptsize\rm\sc Sharp Calder\'on Reproducing Formulae on
Ultra-RD-Spaces}

The main goal of this chapter is to use
the approximation of the identity from Definition~\ref{DefATI}
to construct sharp homogeneous continuous/discrete
Calder\'on reproducing formulae (as well as their inhomogeneous counterparts)
on the ultra-RD-space $(X,\mathbf{q},\mu)$,
where $\mu$ is assumed to be Borel-semiregular on $X$.
Indeed, this is an improvement over works in \cite{hmy08} and \cite{HLYY}
in the sense that these Calder\'on reproducing formulae
obtained in this section are sharp with respect to the regularity exponent
in the setting of ultra-RD-spaces.
To this end, in Section~\ref{section2.1}, we first recall
the definitions of spaces of test functions
and distributions on
$(X,\mathbf{q},\mu)$ and
present some basic estimates.
In Section~\ref{SS3.1}, we first recall the definitions of
Calder\'on--Zygmund kernels and Calder\'on--Zygmund operators,
and then establish the boundedness of Calder\'on--Zygmund
operators on spaces of test functions with (or without) cancellation.
Finally, we establish the sharp homogeneous continuous
Calder\'on reproducing formula in Section~\ref{SS3.2},
the sharp inhomogeneous continuous
Calder\'on reproducing formula in Section~\ref{S3.3},
the sharp homogeneous discrete
Calder\'on reproducing formula in Section~\ref{SS3.4},
and the sharp inhomogeneous discrete
Calder\'on reproducing formula in Section~\ref{SS3.5},
all of which hold, respectively, in spaces of test functions,
spaces of distributions, and $L^p(X)$ with $p\in(1,\infty)$.

\section{Test Function and Distribution Spaces}
\label{section2.1}

In this section, we first
recall the notions of test functions and distributions on
quasi-ultrametric spaces of homogeneous type,
and then state some related estimates that will be frequently used
in the sequel.
To this end, we first recall
the definition of spaces of test functions,
which was originally introduced by Han et al. \cite[Definition 2.8]{hmy08};
see also \cite[Definition 2.2]{hmy06}.

\begin{definition}\label{testfunction}
Let $(X,\mathbf{q},\mu)$ be a quasi-ultrametric space of homogeneous type,
where $\mu$ satisfies the doubling condition \eqref{doublingcond}
with respect to $\rho\in\mathbf{q}$. Fix
$x_1\in X$, $r\in(0,\infty)$,
$\beta\in(0,\infty)$, and $\gamma\in(0,\infty)$.
A $\mu$-measurable function $f:X\to\mathbb{C}$ is called a
\emph{test function of type $(x_1,r,\beta,\gamma)$},
\index[words]{test function of type $(x_1,r,\beta,\gamma)$}
denoted by $f\in{\mathcal G}(x_1,r,\beta,\gamma)$,
\index[symbols]{G@${\mathcal G}(x_1,r,\beta,\gamma)$}
if there exists a
positive constant $C$ such that
\begin{enumerate}
\item (the \emph{size condition})\index[words]{size condition}
For any $x\in X$,
$$
\left|f(x)\right|
\le C\frac{1}{(V_\rho)_r(x_1)+V_\rho(x_1,x)}
\left[\frac{r}{r+\rho(x_1,x)}\right]^\gamma.
$$
\item (the \emph{regularity condition})\index[words]{regularity condition}
For any $x,y\in X$ satisfying
$\rho(x,y)\leq\frac{r+\rho(x_1,x)}{2C_\rho}$,
\begin{equation*}
\left|f(x)-f(y)\right|
\le C\left[\frac{\rho(x,y)}{r+\rho(x_1,x)}\right]^\beta
\frac{1}{(V_\rho)_r(x_1)+V_\rho(x_1,x)}
\left[\frac{r}{r+\rho(x_1,x)}\right]^\gamma.
\end{equation*}
\end{enumerate}
For any $f\in{\mathcal G}(x_1,r,\beta,\gamma)$, define its norm
$$
\|f\|_{{\mathcal G}\left(x_1,r,\beta,\gamma\right)}
:=\inf\left\{C\in(0,\infty):C
\textup{\ satisfies (i) and (ii)}\right\}.
$$
Moreover, define
$$
\mathring{\mathcal{G}}(x_1,r,\beta,\gamma)
:=\left\{f\in{\mathcal G}(x_1,r,\beta,\gamma):
\int_Xf(x)\,d\mu(x)=0\right\}\index[symbols]{G@$\mathring{\mathcal{G}}(x_1,r,\beta,\gamma)$}
$$
equipped with the norm
$\|\cdot\|_{\mathring{\mathcal{G}}(x_1,r,\beta,\gamma)}
:=\|\cdot\|_{{\mathcal G}(x_1,r,\beta,\gamma)}$.
\end{definition}

\begin{remark}
\label{Rmtest}
Let the notation be as in Definition~\ref{testfunction}.
\begin{enumerate}
\item[\textup{(i)}]
Fix $x_1\in X$.
It can be seen that,
for any $x\in X$ and $r\in(0,\infty)$,
\begin{equation*}
{\mathcal G}(x,r,\beta,\gamma)={\mathcal G}(x_1,1,\beta,\gamma)
\end{equation*}
with equivalent norms,
but the positive equivalence constants
depend on both $x$ and $r$.
In what follows, we simply write
$$
{\mathcal G}(\beta,\gamma)
:={\mathcal G}(x_1,1,\beta,\gamma)
\ \ \text{and}\ \
\mathring{\mathcal{G}}(\beta,\gamma)
:=\mathring{\mathcal{G}}(x_1,1,\beta,\gamma).
\index[symbols]{G@${\mathcal G}(\beta,\gamma)$}\index[symbols]{G@$\mathring{\mathcal{G}}(\beta,\gamma)$}
$$

\item[\textup{(ii)}]
If $\rho'\in\mathbf{q}$ is any other quasi-ultrametric
with the property that all $\rho'$-balls are $\mu$-measurable,
then the volume functions defined by $\rho'$ and $\rho$ are equivalent
and hence the spaces ${\mathcal G}(\beta,\gamma)$,
constructed with respect to $\rho$ or $\rho'$,
are the same (with equivalent norms).
In this sense, ${\mathcal G}(\beta,\gamma)$ is independent
of the choice of $\rho\in\mathbf{q}$.
Throughout this monograph, we always assume that
${\mathcal G}(\beta,\gamma)$ is constructed by using a
regularized quasi-ultrametric belonging to $\mathbf{q}$.

\item[\textup{(iii)}]
${\mathcal G}(\beta,\gamma)$ is a Banach space
with respect to the norm in ${\mathcal G}(\beta,\gamma)$;
see, for instance, \cite[p.\,19]{hmy08}.

\item[\textup{(iv)}]
The spaces ${\mathcal G}(\beta,\gamma)$ are nested
with respect to $\beta$ and $\gamma$ in the sense that, for any
$0<\beta_1\leq\beta_2<\infty$ and $0<\gamma_1\leq\gamma_2<\infty$,
$$
{\mathcal G}(\beta_2,\gamma_2)\subset{\mathcal G}(\beta_1,\gamma_1)
$$
with appropriate control of the norms.

\item[\textup{(v)}]
It is easy to see from definitions that,
for any $0<\varepsilon\leq\beta<\infty$ and $0<\gamma<\infty$,
$$
\mathcal{G}(\beta,\gamma)\subset C^\varepsilon(X,\mathbf{q}).
$$
In particular,
if $\mathrm{ind}_H(X,\mathbf{q})<\beta<\infty$, where $\mathrm{ind}_H(X,\mathbf{q})$
denotes the H\"older index of $(X,\mathbf{q})$, then,
by Remark~\ref{rmk390}(vi), we find that,
for any $\gamma\in(0,\infty)$,
$\mathcal{G}(\beta,\gamma)$ only contains constant functions.
See also Lemma~\ref{CGL} below for the relationship between
the spaces $\mathcal{G}(\beta,\gamma)$ and $C^\beta_\mathrm{b}(X)$.

\item[\rm(vi)]
As pointed out in \cite[p.\,18]{hmy08},
Definition~\ref{testfunction} gives a quantified meaning to the
notion of a sufficiently ``smooth" function on quasi-ultrametric spaces
of homogeneous type.
To be more precise, Definition~\ref{testfunction}(i)
means that such test functions are essentially supported in $B_\rho(x_1,r)$
in the sense that it decays of sufficiently high order (measured by
$\gamma$) at infinity
and Definition~\ref{testfunction}(ii) means that
such test functions are H\"older continuous of order $\beta$
(at the right scale $r$).

\item[\rm(vii)]
We point out that, if we replace
$\frac{1}{(V_\rho)_r(x_1)+V_\rho(x_1,x)}$ in both
(i) and (ii) of Definition~\ref{testfunction}
by $\frac{1}{(V_\rho)_{r+\rho(x_1,x)}(x_1)}$
or $\frac{1}{(V_\rho)_{r+\rho(x_1,x)}(x)}$,
then, by Lemma~\ref{A3.2}(vi), we obtain the
same space of test functions,
but with different norms that are equivalent.
\item[\rm(viii)]
By Lemma~\ref{CGL}, we find that
${\mathcal G}(x_1,r,\beta,\gamma)\subset L^1(X)$,
and hence the integral intervening in the
definition of $\mathring{\mathcal{G}}(x_1,r,\beta,\gamma)$
is well defined.
\end{enumerate}
\end{remark}

The following lemma shows that the space
of test functions ${\mathcal G}(\beta,\gamma)$
is non-trivial whenever
$0<\beta\preceq{\mathrm{ind\,}}(X,\mathbf{q})$
and $\gamma\in(0,\infty)$.

\begin{lemma}
\label{ATItest}
Let $(X,\mathbf{q},\mu)$ be a quasi-ultrametric space of homogeneous type,
where $\mu$ satisfies \eqref{doublingcond}
with respect to $\rho\in\mathbf{q}$.
Assume that $0<\beta\preceq{\mathrm{ind\,}}(X,\mathbf{q})$
with ${\mathrm{ind\,}}(X,\mathbf{q})$ as in \eqref{index},
$\gamma\in(0,\infty)$,
and $\{\mathcal{S}_{2^{-k}}\}_{k\in\mathbb{Z}}$
is an approximation of the identity of order $\beta$
as in Definition~\ref{DefATI}.
Then the kernels $\{{S}_{2^{-k}}\}_{k\in\mathbb{Z}}$ of
$\{\mathcal{S}_{2^{-k}}\}_{k\in\mathbb{Z}}$ satisfy,
for any given $k\in\mathbb{Z}$ and $y\in X$,
$$
S_{2^{-k}}(\cdot,y),S_{2^{-k}}(y,\cdot)\in{\mathcal G}(y,2^{-k},\beta,\gamma)
$$
with their norms independent of both $y$ and $k$.
\end{lemma}

\begin{proof}
According to Remark~\ref{Rmkdoub}, we may assume
that $\rho\in\mathbf{q}$ is a regularized quasi-ultrametric.
In particular, suppose that $\rho$ is symmetric.
Note that the existence of such
an approximation of the identity
is guaranteed by Theorem~\ref{ThmATI}.
Fix $y\in X$ and $k\in\mathbb{Z}$.
By Remarks~\ref{Rm}(iii)
and \ref{Rmtest}(ii), without loss of generality,
we may assume that both
$\mathcal{S}_{2^{-k}}$
and ${\mathcal G}(y,2^{-k},\beta,\gamma)$ are related to $\rho$.
Let $C'$ be the same positive constant as in Definition~\ref{DefATI} and let $C:=C'$.
The size condition [see Definition~\ref{testfunction}(i)]
of $S_{2^{-k}}(\cdot,y)$ follows directly from Remark~\ref{Rm}(i)
with both $t:=2^{-k}$ and $\theta:=1$.

Next, we prove that
$S_{2^{-k}}(\cdot,y)$ satisfies
the regularity condition in
Definition~\ref{testfunction}(ii).
Let $x_1,x_2\in X$
satisfy $\rho(x_1,x_2)\leq\frac{2^{-k}+\rho(y,x_1)}{2C_\rho}$.
We divide the proof into three cases.

\emph{Case (1)}
$\rho(y,x_1)\leq C_12^{-k}$ with
some $C_1\in(\max\{1,\,CC_\rho\},\infty)$.
In this case, $\rho(x_1,x_2)\leq\frac{(1+C_1)2^{-k}}{2C_\rho}$ and
Definition~\ref{D2.1}(iii) implies that
$$
\rho(y,x_2)\leq
C_\rho\max\left\{\rho(y,x_1),\,\rho(x_1,x_2)\right\}\lesssim 2^{-k}.
$$
By this, \eqref{upperdoub}, and Lemma~\ref{A3.2}(vii), we find that
$$
V_\rho(y,x_1)\leq(V_\rho)_{C_12^{-k}}(y)
\lesssim (V_\rho)_{2^{-k}}(y)
\sim (V_\rho)_{2^{-k}}(x_1)
\sim (V_\rho)_{2^{-k}}(x_2).
$$
From this, Definition~\ref{DefATI}(ii),
and $\rho(y,x_1)\leq C_12^{-k}$,
we deduce that
\begin{align*}
&\left|S_{2^{-k}}(x_1,y)-S_{2^{-k}}(x_2,y)\right|\\
&\quad\lesssim\left[\frac{\rho(x_1,x_2)}{2^{-k}}\right]^\beta
\left[\frac{1}{(V_\rho)_{2^{-k}}(x_1)}
+\frac{1}{(V_\rho)_{2^{-k}}(x_2)}\right]\\
&\quad\sim
\left[\frac{\rho(x_1,x_2)}{2^{-k}+\rho(y,x_1)}\right]^{\beta}
\frac{1}{(V_\rho)_{2^{-k}}(y)+V_\rho(y,x_1)}
\left[\frac{2^{-k}}{2^{-k}+\rho(y,x_1)}\right]^{\gamma}.
\end{align*}
This is the desired result.

\emph{Case (2)}
$\rho(y,x_1)>C_12^{-k}$
and $\rho(y,x_2)>\frac{C_12^{-k}}{C_\rho}$.
In this case, by Definition~\ref{DefATI}(i),
we conclude that
$$
S_{2^{-k}}(x_1,y)-S_{2^{-k}}(x_2,y)=0.
$$
Hence, the desired result is trivially satisfied.

\emph{Case (3)}
$\rho(y,x_1)>C_12^{-k}$
and $\rho(y,x_2)\leq \frac{C_12^{-k}}{C_\rho}$.
In this case,
from Definition~\ref{D2.1}(iii) and
the fact that
$$
C_\rho\rho(x_1,x_2)\leq\frac{2^{-k}+\rho(y,x_1)}{2}
<\frac{(C_1^{-1}+1)\rho(y,x_1)}{2}<\rho(y,x_1),
$$
we infer that
\begin{align*}
C_12^{-k}<\rho(y,x_1)
\leq C_\rho\max\{\rho(y,x_2),\,\rho(x_2,x_1)\}
=C_\rho\rho(y,x_2)
\leq C_12^{-k},
\end{align*}
which is a contradiction.

In summary, we have shown that $S_{2^{-k}}(\cdot,y)$
satisfies Definition~\ref{testfunction}(ii)
and hence $$
S_{2^{-k}}(\cdot,y)\in{\mathcal G}(y,2^{-k},\beta,\gamma).
$$
By this and the symmetry of $S_{2^{-k}}$ [see Definition~\ref{DefATI}(iv)],
we find that
$S_{2^{-k}}(y,\cdot)\in{\mathcal G}(y,2^{-k},\beta,\gamma)$.
This finishes the proof of Lemma~\ref{ATItest}.
\end{proof}

A well-known \emph{drawback} of the space
${\mathcal G}(\beta,\gamma)$ of test functions
is that, even in the case where
$X=\mathbb{R}^n$ and $\rho$ is the standard Euclidean distance on $\mathbb{R}^n$,
${\mathcal G}(\beta_2,\gamma)$
is not dense in ${\mathcal G}(\beta_1,\gamma)$ with
$0<\beta_1<\beta_2<\infty$ and $\gamma\in(0,\infty)$,
which will bring some inconvenience.
To overcome this defect, we now recall the following definition
of spaces of test functions; see also \cite[p.\,19]{hmy08}.

\begin{definition}\label{2207}
Let $(X,\mathbf{q},\mu)$ be a quasi-ultrametric space of homogeneous type,
where $\mu$ satisfies the doubling condition \eqref{doublingcond}
with respect to $\rho\in\mathbf{q}$. Assume that
$0<\beta,\gamma\leq\varepsilon\preceq{\mathrm{ind\,}}(X,\mathbf{q})$
with ${\mathrm{ind\,}}(X,\mathbf{q})$
as in \eqref{index}.
Let ${\mathcal G}_0^\varepsilon(\beta,\gamma)$
\index[symbols]{G@${\mathcal G}_0^\varepsilon(\beta,\gamma)$}
[resp. $\mathring{{\mathcal G}}_0^\varepsilon(\beta,\gamma)$]
\index[symbols]{G@$\mathring{{\mathcal G}}_0^\varepsilon(\beta,\gamma)$}
be the completion of ${\mathcal G}(\varepsilon,\varepsilon)$
[resp. $\mathring{\mathcal{G}}(\varepsilon,\varepsilon)$]
in ${\mathcal G}(\beta,\gamma)$.
For any $f\in{\mathcal G}_0^\varepsilon(\beta,\gamma)$
[resp. $f\in\mathring{\mathcal G}^\varepsilon_0(\beta,\gamma)$],
let
$\|f\|_{{\mathcal G}_0^\varepsilon(\beta,\gamma)}:=\|f\|_{{\mathcal G}(\beta,\gamma)}$
[resp.
$\|f\|_{\mathring{\mathcal G}^\varepsilon_0(\beta,\gamma)}
:=\|f\|_{{\mathcal G}(\beta,\gamma)}$].
\end{definition}

\begin{remark}
Let the notation be as in Definition~\ref{2207}.
\begin{enumerate}
\item
Obviously, ${\mathcal G}_0^\varepsilon(\varepsilon,\varepsilon)
={\mathcal G}(\varepsilon,\varepsilon)$
and $\mathring{{\mathcal G}}_0^\varepsilon(\varepsilon,\varepsilon)
=\mathring{{\mathcal G}}(\varepsilon,\varepsilon)$.
\item
$f\in{\mathcal G}_0^\varepsilon(\beta,\gamma)$
[resp. $f\in\mathring{\mathcal G}_0^\varepsilon(\beta,\gamma)$] \emph{if and only if}
$f\in\mathcal{G}(\beta,\gamma)$ [resp. $f\in\mathring{\mathcal G}(\beta,\gamma)$] and
there exists a sequence $\{f_j\}_{j\in\mathbb{N}}$
in ${\mathcal G}(\varepsilon,\varepsilon)$
[resp. $\{f_j\}_{j\in\mathbb{N}}$ in $\mathring{\mathcal{G}}(\varepsilon,\varepsilon)$]
such that
$\|f_j-f\|_{{\mathcal G}(\beta,\gamma)}\to0$
as $j\to\infty$.
\item
Both ${\mathcal G}_0^\varepsilon(\beta,\gamma)$
and $\mathring{\mathcal G}^\varepsilon_0(\beta,\gamma)$
are Banach spaces.
\item
For any given $\beta_1,\gamma_1,\beta_2,\gamma_2\in(0,\varepsilon]$,
we find that both $\mathcal{G}^\varepsilon_0(\beta_1,\gamma_1)$
and $\mathcal{G}^\varepsilon_0(\beta_2,\gamma_2)$ have the same
dense subset $\mathcal{G}(\varepsilon,\varepsilon)$.
\end{enumerate}
\end{remark}

The following result is a quantitative
version of the classical Urysohn's lemma,
which was originally proved in
\cite[Theorem 4.12]{MiMiMiMo13}
and subsequently generalized in \cite[Theorem 4.1]{AlMiMi13}.

\begin{lemma}
\label{lemUrysohn}
Let $(X,\rho)$ be a quasi-ultrametric space.
Assume that $\varepsilon\in(0,(\log_2C_\rho)^{-1}]\cap(0,\infty)$
with $C_\rho\in[1,\infty)$ as in \eqref{Crho}.
Suppose that $F_0,F_1\subset X$ are two
nonempty sets with the property that
$\mathrm{dist}_\rho(F_0,F_1)>0$.
Then there exist a constant $\theta_\rho\in(0,\infty)$
and a function $\varphi\in C^\varepsilon(X)$ such that
\begin{enumerate}
\item[\textup{(i)}]
For any $x\in X$, $\varphi(x)\in[0,1]$.

\item[\textup{(ii)}]
For any $x\in F_0$, $\varphi(x)=0$.

\item[\textup{(iii)}]
For any $x\in F_1$, $\varphi(x)=1$.

\item[\textup{(iv)}]
$\|\varphi\|_{\dot{C}^\varepsilon(X)}
\leq\theta_\rho[\mathrm{dist}_\rho(F_0,F_1)]^{-\varepsilon}$.
\end{enumerate}
\end{lemma}

The following conclusion is a direct corollary of
Lemma~\ref{lemUrysohn}.

\begin{corollary}
\label{Lem3.3}
Let $(X,\rho)$ be a quasi-ultrametric space and
$\varepsilon\in(0,(\log_2C_\rho)^{-1}]\cap(0,\infty)$.
Let $x\in X$ and $r\in(0,\infty)$.
Then, for any $A\in[C_\rho,\infty)$,
there exists a function $f\in C_{\rm b}^\varepsilon(X)$ such that
\begin{enumerate}
\item[\textup{(i)}]
For any $y\in X$,
\begin{equation*}
\mathbf{1}_{B_\rho(x,r)}(y)\leq f(y)
\leq\mathbf{1}_{B_\rho(x,Ar)}(y).
\end{equation*}

\item[\textup{(ii)}]
$\|f\|_{\dot{C}^\varepsilon(X)}\leq Cr^{-\varepsilon}$,
where $C\in(0,\infty)$ is independent of both $x$ and $r$.
\end{enumerate}
\end{corollary}

\begin{remark}
If $\rho$ is a metric, then $(\log_2C_\rho)^{-1}\in[1,\infty]$
and Corollary~\ref{Lem3.3} holds for any $\beta\in(0,(\log_2C_\rho)^{-1}]$.
If $\rho$ is an ultrametric,
then $(\log_2C_\rho)^{-1}=\infty$ and
Corollary~\ref{Lem3.3} holds for any $\beta\in(0,\infty)$.
However, He et al. \cite[Lemma 3.3]{HLYY}
only showed the existence of such smooth cut-off functions
as in Lemma~\ref{Lem3.3} for any $\varepsilon\in(0,1)$.
In this sense,
Corollary~\ref{Lem3.3} is sharper than \cite[Lemma 3.3]{HLYY}.
\end{remark}

The following lemma gives the relation among $C^\beta_\mathrm{b}(X)$,
${\mathcal G}(x_1,r,\beta,\gamma)$, and $L^p(X)$
with $p\in[1,\infty)$.

\begin{lemma}
\label{CGL}
Let $(X,\mathbf{q},\mu)$ be a quasi-ultrametric
space of homogeneous type. Then, for any
$x_1\in X$, $r\in(0,\infty)$, $\beta,\gamma\in(0,\infty)$,
and $p\in[1,\infty)$,
\begin{align*}
C^\beta_\mathrm{b}(X)\subset{\mathcal G}(x_1,r,\beta,\gamma)\subset L^p(X).
\end{align*}
\end{lemma}

\begin{proof}
Let $\rho\in\mathbf{q}$ be a regularized quasi-ultrametric
as in Definition~\ref{D2.3}.
From Definition~\ref{testfunction}(i) and Lemma~\ref{A3.2}(ii),
we deduce that,
for any $f\in{\mathcal G}(x_1,r,\beta,\gamma)$,
\begin{align*}
\int_X|f(x)|^p\,d\mu(x)&\leq
\|f\|_{{\mathcal G}(x_1,r,\beta,\gamma)}
\int_X\frac{1}{(V_\rho)_r(x_1)+V_\rho(x_1,x)}\left[\frac{r}{r+\rho(x_1,x)}\right]^\gamma
\\
&\quad\times\left\{\frac{1}{(V_\rho)_r(x_1)+V_\rho(x_1,x)}
\left[\frac{r}{r+\rho(x_1,x)}\right]^\gamma\right\}^{p-1}\,d\mu(x)
\\
&\lesssim\|f\|_{{\mathcal G}(x_1,r,\beta,\gamma)}\left[\frac{1}{(V_\rho)_r(x_1)}\right]^{p-1}
<\infty,
\end{align*}
which further implies that $f\in L^p(X)$
and hence ${\mathcal G}(x_1,r,\beta,\gamma)\subset L^p(X)$.

Thus, to complete the proof of the present lemma,
it suffices to show that
$$
C^\beta_\mathrm{b}(X)\subset{\mathcal G}(x_1,r,\beta,\gamma).
$$
Let $f\in C^\beta_\mathrm{b}(X)$. By Remark~\ref{Rmtest}(i),
without loss of generality, we may assume that
${\mathrm{supp\,}} f\subset B_\rho(x_1,r)$ and
that there exist $C\in(1,\infty)$ and
$y\in B_\rho(x_1,Cr)$
such that $f(y)=0$.
We first prove that $f$ satisfies
Definition~\ref{testfunction}(i).
Note that, for any $x\in B_\rho(x_1,r)$,
we have $\rho(x_1,x)<r$ and $V_\rho(x_1,x)\leq(V_\rho)_{r}(x_1)$,
which further imply that
\begin{align*}
\frac{1}{(V_\rho)_r(x_1)+V_\rho(x_1,x)}
\left[\frac{r}{r+\rho(x_1,x)}\right]^\gamma
\sim\frac{1}{(V_\rho)_r(x_1)},
\end{align*}
where the positive equivalence constants
depend only on $\gamma$.
From this and
$f(y)=0$,
we infer that, for any $x\in B_\rho(x_1,r)$,
\begin{align*}
|f(x)|&\leq|f(x)-f(y)|
\leq\|f\|_{\dot{C}^\beta(X)}\left[\rho(x,y)\right]^\beta
\lesssim\|f\|_{\dot{C}^\beta(X)}r^\beta
\\
&\sim\|f\|_{\dot{C}^\beta(X)}r^\beta(V_\rho)_r(x_1)
\frac{1}{(V_\rho)_r(x_1)+V_\rho(x_1,x)}
\left[\frac{r}{r+\rho(x_1,x)}\right]^\gamma,
\end{align*}
where the implicit positive constants
depend only on $C$, $C_\rho$, and $\gamma$.
This, combined with ${\mathrm{supp\,}} f\subset B_\rho(x_1,r)$,
implies that $f$ satisfies the desired size condition.

Next, we show that $f$ satisfies
Definition~\ref{testfunction}(ii).
Let $x,x'\in X$ satisfy
$\rho(x,x')\leq\frac{r+\rho(x_1,x)}{2C_\rho}$.
We divide the proof into three cases.

\emph{Case (1)} $\rho(x_1,x)\ge r$ and $\rho(x_1,x')\ge r$.
In this case, since ${\mathrm{supp\,}} f\subset B_\rho(x_1,r)$,
it follows that
$|f(x)-f(x')|=0$.

\emph{Case (2)} $\rho(x_1,x)\ge r$ and $\rho(x_1,x')<r$.
In this case, by Definition~\ref{D2.1}(iii) and
$\rho(x,x')\leq\frac{r+\rho(x_1,x)}{2C_\rho}$, we conclude that
\begin{align*}
\rho(x_1,x)
&\leq C_\rho\max\{\rho(x_1,x'),\,
\rho(x',x)\}
\leq\max\left\{C_\rho\rho(x_1,x'),\,
\frac{r+\rho(x_1,x)}{2}\right\}
\\
&\leq C_\rho r+\frac{r+\rho(x_1,x)}{2},
\end{align*}
which further implies that $\rho(x_1,x)\leq(2C_\rho+1)r$.
From this and \eqref{upperdoub},
we deduce that $V_\rho(x_1,x)\lesssim (V_\rho)_r(x_1)$ and hence
\begin{align}\label{2052}
\left[\frac{1}{r+\rho(x_1,x)}\right]^\beta
\frac{1}{(V_\rho)_r(x_1)+V_\rho(x_1,x)}
\left[\frac{r}{r+\rho(x_1,x)}\right]^\gamma\sim
\frac{1}{r^\beta(V_\rho)_r(x_1)},
\end{align}
where the positive equivalence constants are
independent of $x_1$ and $r$.
Consequently, we obtain
\begin{align}\label{2049}
\left|f(x)-f(x')\right|
&\leq\|f\|_{\dot{C}^\beta(X)}\left[\rho(x,x')\right]^\beta
\quad\text{by \eqref{C(X)}}\nonumber\\
&\sim\|f\|_{\dot{C}^\beta(X)}r^\beta(V_\rho)_r(x_1)
\left[\frac{\rho(x,x')}{r+\rho(x_1,x)}\right]^\beta\nonumber\\
&\quad\times\frac{1}{(V_\rho)_r(x_1)+V_\rho(x_1,x)}
\left[\frac{r}{r+\rho(x_1,x)}\right]^\gamma
\quad\text{by \eqref{2052}},
\end{align}
where the positive equivalence constants are independent of $x_1$ and $r$.

\emph{Case (3)} $\rho(x_1,x)<r$.
In this case, we have $V_\rho(x_1,x)\leq(V_\rho)_r(x_1)$
and hence \eqref{2052} still holds in this case.
Applying this and an argument similar to that used in the estimation
of \eqref{2049},
we find that \eqref{2049} still holds in this case.

Altogether, we have verified that $f$ satisfies
Definition~\ref{testfunction}(ii) and hence
$$f\in{\mathcal G}(x_1,r,\beta,\gamma)\ \mathrm{and}\
\|f\|_{{\mathcal G}(x_1,r,\beta,\gamma)}\lesssim
\|f\|_{\dot{C}^\beta(X)}r^\beta(V_\rho)_r(x_1),$$
where the implicit positive constant is
independent of $x_1$ and $r$. This finishes the proof
of Lemma~\ref{CGL}.
\end{proof}

For the density among $C^\beta_\mathrm{b}(X)$,
${\mathcal G}(x_1,r,\beta,\gamma)$, and $L^p(X)$ with $p\in(1,\infty)$,
we have the following lemma.

\begin{lemma}
\label{density}
Let $(X,\mathbf{q},\mu)$ be a quasi-ultrametric space of homogeneous type,
where $\mu$ is a Borel-semiregular measure on $X$.
Then, for any
$x_1\in X$, $r\in(0,\infty)$,
$0<\beta\preceq{\mathrm{ind\,}}(X,\mathbf{q})$,
and $\gamma\in(0,\infty)$,
\begin{enumerate}
\item[\textup{(i)}]
Both $C^\beta_\mathrm{b}(X)$ and ${\mathcal G}(x_1,r,\beta,\gamma)$ are
dense in $L^p(X)$ for any $p\in[1,\infty)$.

\item[\textup{(ii)}]
If $\mu(X)=\infty$, then
$\mathring{C}^\beta_{\mathrm{b}}(X)$ and
$\mathring{{\mathcal G}}(x_1,r,\beta,\gamma)$ are
both dense in $L^p(X)$ for any $p\in(1,\infty)$,
where, for any $0<\beta\preceq{\mathrm{ind\,}}(X,\mathbf{q})$,
\begin{align*}
\mathring{C}^\beta_{\mathrm{b}}(X):=\left\{f\in C^\beta_\mathrm{b}(X):
\int_Xf(x)\,d\mu(x)=0\right\}.
\end{align*}
\end{enumerate}
\end{lemma}

\begin{proof}
From Remark~\ref{RmkD2.3}(iii), it follows that
there exists a regularized quasi-ultrametric
$\rho\in\mathbf{q}$ as in Definition~\ref{D2.3}
satisfying $\beta\leq(\log_2C_{\rho})^{-1}$.
By Lemmas~\ref{LDT} and~\ref{CGL}, we easily conclude that (i) holds.
From this, we infer that, to prove (ii),
it suffices to show that every function in $C^\beta_\mathrm{b}(X)$ can be
approximated in the $L^p$ norm by functions in $\mathring{C}^\beta_{\mathrm{b}}(X)$.
Let $f\in C^\beta_\mathrm{b}(X)$,
where ${\mathrm{supp\,}} f\subset B_\rho(x_0,r_0)$ with $x_0\in X$
and $r_0\in(0,\infty)$.
If $\int_Xf(x)\,d\mu(x)=0$ then there is nothing to prove. Thus,
we may assume that $\int_Xf(x)\,d\mu(x)\neq0$.
By Corollary~\ref{Lem3.3}, we find that,
for any given $n\in\mathbb{N}$, there exists
a function $h_n$ such that, for any $y\in X$,
\begin{align*}
\mathbf{1}_{B_\rho(x_0,n)}(y)\leq h_n(y)
\leq\mathbf{1}_{B_\rho(x_0,An)}(y)
\end{align*}
and $\|h_n\|_{\dot{C}^{\beta}(X)}
\lesssim n^{-\beta}$, where $A:=C_\rho\in[1,\infty)$.
Then, for any $n\in\mathbb{N}$ and $y\in X$, it is meaningful to define
\begin{align*}
\widetilde{h}_n(y):=\frac{\int_Xf(x)\,d\mu(x)}{\int_Xh_n(x)\,d\mu(x)}
h_n(y)
\end{align*}
and
\begin{align*}
g_n(y):=f(y)-\widetilde{h}_n(y).
\end{align*}

Now, we claim that, for any $n\in\mathbb{N}$,
$g_n\in\mathring{C}^\beta_\mathrm{b}(X)$.
Indeed,
from the definitions of $\widetilde{h}_n$
and $g_n$, it is easy to verify that, for any $n\in\mathbb{N}$,
\begin{align*}
\int_Xg_n(x)\,d\mu(x)=0.
\end{align*}
Since both $f$ and $h_n$ have bounded supports,
it follows that, for any $n\in\mathbb{N}$,
$g_n$ has a bounded support.
By \eqref{C(X)},
the fact that $f,h_n\in\dot{C}^\beta(X)$,
and the definition of $\widetilde{h}_n$,
we conclude that,
for any $n\in\mathbb{N}$,
\begin{align*}
\left\|g_n\right\|_{\dot{C}^\beta(X)}
&\leq\sup_{x,y\in X,\,x\neq y}
\frac{|f(x)-f(y)|}{[\rho(x,y)]^\beta}
+\sup_{x,y\in X,\,x\neq y}
\frac{|\widetilde{h}_n(x)-\widetilde{h}_n(y)|}{[\rho(x,y)]^\beta}
\\
&=\left\|f\right\|_{\dot{C}^\beta(X)}
+\left\|\widetilde{h}_n\right\|_{\dot{C}^\beta(X)}
\sim\left\|f\right\|_{\dot{C}^\beta(X)}
+\left\|h_n\right\|_{\dot{C}^\beta(X)}
<\infty.
\end{align*}
This finishes the proof of the claim.

Next, we show that
\begin{align}\label{1449}
\left\|f-g_n\right\|_{L^p(X)}
=\left\|\widetilde{h}_n\right\|_{L^p(X)}\to0
\end{align}
as $n\to\infty$.
From the definition of $h_n$,
$A\in[1,\infty)$, \eqref{upperdoub},
$p\in(1,\infty)$, and $\mu(X)=\infty$,
we deduce that
\begin{align*}
\left\|\widetilde{h}_n\right\|_{L^p(X)}
&=
\left|\frac{\int_Xf(x)\,d\mu(x)}{\int_Xh_n(x)\,d\mu(x)}\right|
\left[\int_X\left|h_n(y)\right|^p\,d\mu(y)\right]^\frac{1}{p}
\\
&\leq\frac{\left\|f\right\|_{L^1(X)}}{(V_\rho)_n(x_0)}
\left[(V_\rho)_{An}(x_0)\right]^\frac{1}{p}
\lesssim\left\|f\right\|_{L^1(X)}
\left[(V_\rho)_{n}(x_0)\right]^{\frac{1}{p}-1}\to0
\end{align*}
as $n\to\infty$,
which is exactly \eqref{1449}.
This finishes the proof of Lemma~\ref{density}.
\end{proof}

For any $\beta,\gamma\in(0,\infty)$,
let
$$
{\mathcal{G}}_{\mathrm{b}}(\beta,\gamma):=
\left\{f\in\mathcal{G}(\beta,\gamma):f\
\text{has bounded support}\right\}
$$
and
$$
\mathring{{\mathcal{G}}}_{\mathrm{b}}(\beta,\gamma):=
\left\{f\in\mathring{\mathcal{G}}(\beta,\gamma):f\
\text{has bounded support}\right\}.
$$
Note that ${\mathcal{G}}_{\mathrm{b}}(\beta,\gamma)
=C_{\mathrm{b}}^\beta(X,\mathbf{q})$
as sets.

\begin{lemma}\label{dens}
Let $(X,\mathbf{q},\mu)$ be a quasi-ultrametric space
of homogeneous type. Assume that $0<\varepsilon\preceq\mathrm{ind\,}(X,\mathbf{q})$.
Then, for any $\beta,\gamma\in(0,\varepsilon)$,
${\mathcal{G}}_{\mathrm{b}}(\varepsilon,\varepsilon)$
[resp. $\mathring{\mathcal{G}}_{\mathrm{b}}(\varepsilon,\varepsilon)$]
is dense in ${{\mathcal G}}_0^\varepsilon(\beta,\gamma)$
[resp. $\mathring{{\mathcal G}}_0^\varepsilon(\beta,\gamma)$].
\end{lemma}

\begin{proof}
Let $\beta,\gamma\in(0,\varepsilon)$.
We only prove that ${\mathcal{G}}_{\mathrm{b}}(\varepsilon,\varepsilon)$
is dense in ${{\mathcal G}}_0^\varepsilon(\beta,\gamma)$.
By
${{\mathcal G}}_0^\varepsilon(\beta,\gamma)
=\overline{\mathcal{G}(\varepsilon,\varepsilon)}^{\mathcal{G}(\beta,\gamma)}$,
it suffices to show that,
for any $f\in\mathcal{G}(\varepsilon,\varepsilon)$,
there exists a sequence $\{f_j\}_{j\in\mathbb{N}}$ in
${\mathcal{G}}_{\mathrm{b}}(\varepsilon,\varepsilon)$
such that
\begin{align}\label{1923}
\lim_{j\to\infty}\left\|f-f_j\right\|_{\mathcal{G}(\beta,\gamma)}=0.
\end{align}
To this end, fix a point $x_1\in X$.
From Remark~\ref{RmkD2.3}(iii), we infer that
there exists a regularized
quasi-ultrametric $\rho\in\mathbf{q}$
that satisfies $\varepsilon\leq(\log_2C_\rho)^{-1}$.
By Corollary~\ref{Lem3.3}, we find that there exists a sequence
$\{\varphi_j\}_{j\in\mathbb{N}}$ of smooth
cut-off functions such that, for any $j\in\mathbb{N}$,
\begin{align}\label{2399}
\mathbf{1}_{B_\rho(x_1,j)}\leq\varphi_j\leq\mathbf{1}_{B_\rho(x_1,Aj)},
\end{align}
where $A\in[C_\rho,\infty)$ is the same positive constant as in Corollary~\ref{Lem3.3}
and $\|\varphi_j\|_{\dot{C}^\varepsilon(X)}\lesssim j^{-\varepsilon}$.
For any $j\in\mathbb{N}$, let $f_j:=f\varphi_j$.

To prove \eqref{1923}, we first estimate the size condition of $f-f_j$.
From \eqref{2399} and Definition~\ref{testfunction}(i)
with $\gamma$ therein replaced by $\varepsilon$,
we deduce that, for any $x\in X$,
\begin{align}\label{size19}
\left|f(x)-f_j(x)\right|
&\leq|f(x)|\mathbf{1}_{[B_\rho(x_1,j)]^\complement}(x)
\nonumber\\
&\leq\frac{\|f\|_{\mathcal{G}(\varepsilon,\varepsilon)}}{(V_\rho)_1(x_1)+V_\rho(x_1,x)}
\left[\frac{1}{1+\rho(x_1,x)}\right]^\gamma\left[\frac{1}{1+\rho(x_1,x)}\right]^{\varepsilon-\gamma}
\mathbf{1}_{[B_\rho(x_1,j)]^\complement}(x)\nonumber\\
&<j^{\gamma-\varepsilon}
\frac{\|f\|_{\mathcal{G}(\varepsilon,\varepsilon)}}{(V_\rho)_1(x_1)+V_\rho(x_1,x)}
\left[\frac{1}{1+\rho(x_1,x)}\right]^\gamma.
\end{align}
Now, we estimate the regularity of $f-f_j$.
Assume that $x,x'\in X$ satisfy $\rho(x,x')\leq\frac{1+\rho(x_1,x)}{2C_\rho}$.
Then we write
\begin{align}\label{2019II}
&\left|\left[f(x)-f_j(x)\right]-
\left[f(x')-f_j(x')\right]\right|\nonumber\\
&\quad=\left|f(x)\left[1-\varphi_j(x)\right]
-f(x)\left[1-\varphi_j(x')\right]\right.\nonumber\\
&\qquad\left.+f(x)\left[1-\varphi_j(x')\right]
-f(x')\left[1-\varphi_j(x')\right]\right|\nonumber\\
&\quad\leq|f(x)|\left|\varphi_j(x)-\varphi_j(x')\right|
+\left|f(x)-f(x')\right|\left|1-\varphi_j(x')\right|\nonumber\\
&\quad=:\mathrm{I}+\mathrm{II}.
\end{align}

To estimate $\mathrm{I}$, we consider the following
three cases for $\rho(x_1,x)$.

\emph{Case (1)} $\rho(x_1,x)\ge C_\rho Aj$.
In this case, by Definition~\ref{D2.1}(iii)
and $\rho(x,x')\leq\frac{1+\rho(x_1,x)}{2C_\rho}$
(keep in mind that we assume that $\rho$ is a regularized
quasi-ultrametric and hence $\rho$ is symmetric),
we obtain
$$
\rho(x_1,x)\leq C_\rho
\max\left\{\rho(x_1,x'),\,\rho(x',x)\right\}
=C_\rho\rho(x_1,x'),
$$
which further implies that $\rho(x_1,x')\ge Aj$.
From this and \eqref{2399}, we infer that
$$
\varphi_j(x)=0=\varphi_j(x')
$$
and hence $\mathrm{I}=0$.

\emph{Case (2)} $\rho(x_1,x)<\frac{j}{C_\rho}$.
In this case, by Definition~\ref{D2.1}(iii)
and $\rho(x,x')\leq\frac{1+\rho(x_1,x)}{2C_\rho}$, we find that
\begin{align*}
\rho(x_1,x')&\leq C_\rho
\max\left\{\rho(x_1,x),\,\rho(x,x')\right\}\\
&\leq C_\rho\max\left\{\rho(x_1,x),\,
\frac{1+\rho(x_1,x)}{2C_\rho}\right\},
\end{align*}
which, given the assumption $\rho(x_1,x)<\frac{j}{C_\rho}$,
implies that $\rho(x_1,x')<j$.
From this and \eqref{2399}, it follows that
$$
\varphi_j(x)=1=\varphi_j(x')
$$
and hence $\mathrm{I}=0$.

\emph{Case (3)} $\frac{j}{C_\rho}\leq\rho(x_1,x)<C_\rho Aj$.
In this case, by Definition~\ref{testfunction}(i),
the fact that
$\|\varphi_j\|_{\dot{C}^\varepsilon(X)}\lesssim j^{-\varepsilon}$,
and $\beta,\gamma\in(0,\varepsilon)$,
we conclude that
\begin{align*}
\mathrm{I}&\lesssim\|f\|_{\mathcal{G}(\varepsilon,\varepsilon)}
\frac{1}{(V_\rho)_1(x_1)+V_\rho(x_1,x)}
\left[\frac{1}{1+\rho(x_1,x)}\right]^\varepsilon
\left[\frac{\rho(x,x')}{j}\right]^\varepsilon\\
&\sim\|f\|_{\mathcal{G}(\varepsilon,\varepsilon)}
\frac{1}{(V_\rho)_1(x_1)+V_\rho(x_1,x)}
\left[\frac{1}{1+\rho(x_1,x)}\right]^\gamma
\left[\frac{1}{1+\rho(x_1,x)}\right]^{\varepsilon-\gamma}
\left[\frac{\rho(x,x')}{1+\rho(x_1,x)}\right]^\varepsilon\\
&\lesssim j^{\gamma-\varepsilon}
\|f\|_{\mathcal{G}(\varepsilon,\varepsilon)}
\frac{1}{(V_\rho)_1(x_1)+V_\rho(x_1,x)}
\left[\frac{1}{1+\rho(x_1,x)}\right]^\gamma
\left[\frac{\rho(x,x')}{1+\rho(x_1,x)}\right]^\beta.
\end{align*}
Combining the above three cases, we find that
\begin{align}\label{2019I}
\mathrm{I}\lesssim
j^{\gamma-\varepsilon}
\|f\|_{\mathcal{G}(\varepsilon,\varepsilon)}
\frac{1}{(V_\rho)_1(x_1)+V_\rho(x_1,x)}
\left[\frac{1}{1+\rho(x_1,x)}\right]^\gamma
\left[\frac{\rho(x,x')}{1+\rho(x_1,x)}\right]^\beta.
\end{align}

To estimate $\mathrm{II}$, from Definition~\ref{testfunction}(ii)
and \eqref{2399}, we deduce that
\begin{align*}
\mathrm{II}
&\leq\|f\|_{\mathcal{G}(\varepsilon,\varepsilon)}
\left[\frac{\rho(x,x')}{1+\rho(x_1,x)}\right]^\beta
\frac{1}{(V_\rho)_1(x_1)+V_\rho(x_1,x)}
\left[\frac{1}{1+\rho(x_1,x)}\right]^\gamma\\
&\quad\times\left[\frac{1}{1+\rho(x_1,x)}\right]^{\varepsilon-\gamma}
\mathbf{1}_{[B_\rho(x_1,j)]^\complement}(x'),
\end{align*}
which, together with $\gamma\in(0,\varepsilon)$
and the estimate that
\begin{align*}
j&<\rho(x_1,x')\leq C_\rho
\max\{\rho(x_1,x),\,\rho(x,x')\}\\
&\leq C_\rho\max\left\{\rho(x_1,x),\,
\frac{1+\rho(x_1,x)}{2C_\rho}\right\}
<C_\rho\left[1+\rho(x_1,x)\right],
\end{align*}
implies that
\begin{align*}
\mathrm{II}
&\lesssim j^{\gamma-\varepsilon}
\|f\|_{\mathcal{G}(\varepsilon,\varepsilon)}
\left[\frac{\rho(x,x')}{1+\rho(x_1,x)}\right]^\beta
\frac{1}{(V_\rho)_1(x_1)+V_\rho(x_1,x)}
\left[\frac{1}{1+\rho(x_1,x)}\right]^\gamma.
\end{align*}
This, combined with \eqref{2019II}, \eqref{2019I},
and \eqref{size19},
further implies that
\begin{align*}
&\left|\left[f(x)-f_j(x)\right]-
\left[f(x')-f_j(x')\right]\right|\nonumber\\
&\quad\lesssim
j^{\gamma-\varepsilon}
\|f\|_{\mathcal{G}(\varepsilon,\varepsilon)}
\left[\frac{\rho(x,x')}{1+\rho(x_1,x)}\right]^\beta
\frac{1}{(V_\rho)_1(x_1)+V_\rho(x_1,x)}
\left[\frac{1}{1+\rho(x_1,x)}\right]^\gamma
\end{align*}
and \eqref{1923} holds.
Moreover, note that an argument similar to that used in the proof
of \eqref{1923} implies that $\{f_j\}_{j\in\mathbb{N}}$ in
${\mathcal{G}}_{\mathrm{b}}(\varepsilon,\varepsilon)$,
which completes the proof of Lemma~\ref{dens}.
\end{proof}

Finally, we recall the concept of spaces of distributions on
quasi-ultrametric spaces of homogeneous type;
see \cite[p.\,19]{hmy08}. To this end, we introduce some notation:
Given a vector space $Y$
(over $\mathbb{C}$ or similarly over $\mathbb{R}$)
and a quasi-norm $\|\cdot\|_Y$ defined on $Y$,
the \emph{topological dual} $Y'$ of $Y$ is defined by setting
\begin{align*}
Y':\index[symbols]{Y@$Y'$}
=\left\{f:(Y,\tau_{\|\cdot\|_Y})\to\mathbb{C}
\text{ is continuous and linear}\right\},
\end{align*}
where $\tau_{\|\cdot\|_Y}$\index[symbols]{T@$\tau_{\|\cdot\|_Y}$}
is the topology induced
by the quasi-norm $\|\cdot\|_Y$ on $Y$.
We can  now state the following definition.

\begin{definition}\label{1558}
Let $(X,\mathbf{q},\mu)$ be a quasi-ultrametric space of homogeneous type,
where $\mu$ satisfies the doubling condition \eqref{doublingcond}
with respect to $\rho\in\mathbf{q}$. Assume that
$0<\beta,\gamma\leq\varepsilon\preceq{\mathrm{ind\,}}(X,\mathbf{q})$
with ${\mathrm{ind\,}}(X,\mathbf{q})$
as in \eqref{index}.
The \emph{dual space}
$({\mathcal G}^\varepsilon_0(\beta,\gamma))'$
\index[symbols]{G@$({\mathcal G}^\varepsilon_0(\beta,\gamma))'$}
[resp. $(\mathring{{\mathcal G}}_0^\varepsilon
(\beta,\gamma))'$]\index[symbols]{G@$(\mathring{{\mathcal G}}_0^\varepsilon(\beta,\gamma))'$}
is defined to be the set of
all continuous linear functionals on ${\mathcal G}^\varepsilon_0(\beta,\gamma)$
[resp. $\mathring{{\mathcal G}}_0^\varepsilon(\beta,\gamma)$]
equipped with the weak-$*$ topology,
which is also called the
\emph{space of distributions}.\index[words]{space of distributions}
\end{definition}

It is well understood in the classical Euclidean setting that function
$f\in L^r(\mathbb{R}^n)$ with $r\in[1,\infty]$
induces a distribution on ${\mathscr{D}}(\mathbb{R}^n):=
C^\infty_\mathrm{c}(\mathbb{R}^n)$ (denoted by $\Lambda_f$)
via integrating $f$ against
any function from ${\mathscr{D}}(\mathbb{R}^n)$
over the entire $\mathbb{R}^n$. Indeed, this association is injective
and we may unambiguously identify such a distribution $\Lambda_f$ with
the function $f$ itself. As the following proposition asserts,
this continues to remain valid in quasi-ultrametric spaces of homogeneous type.

\begin{proposition}
\label{distfncttype}
Let $(X,\mathbf{q},\mu)$ be a quasi-ultrametric space of homogeneous type,
$r\in[1,\infty]$, and $\varepsilon,\beta,\gamma\in(0,\infty)$ be such that
$0<\beta,\gamma<\varepsilon
\preceq\mathrm{ind\,}(X,\mathbf{q})$.
Then, for any $f\in L^r(X)$, $f\in({\mathcal{G}}^\varepsilon_0(\beta,\gamma))'$
in the sense that the linear functional
$\Lambda_f:{\mathcal{G}}^\varepsilon_0(\beta,\gamma)\to\mathbb{C}$,
defined by setting,
for any $\varphi\in{\mathcal{G}}^\varepsilon_0(\beta,\gamma)$,
\begin{align}\label{1503}
\langle\Lambda_f,\varphi\rangle:=\int_Xf(x)\varphi(x)\,d\mu(x),
\end{align}
belongs to $({\mathcal{G}}^\varepsilon_0(\beta,\gamma))'$.
Moreover, if $\mu$ is Borel-semiregular, then
$\Lambda_f=0$ in $({\mathcal{G}}^\varepsilon_0(\beta,\gamma))'$
if and only if $f=0$ pointwise
$\mu$-almost everywhere on $X$. Consequently, the association of $f\in L^r(X)$
to the distribution $\Lambda_f\in({\mathcal{G}}^\varepsilon_0(\beta,\gamma))'$
is injective, and hence
$L^r(X)\subset({\mathcal{G}}^\varepsilon_0(\beta,\gamma))'$.
\end{proposition}

\begin{remark}
In the context of Proposition~\ref{distfncttype},
since the association of $f\in L^r(X)$ to the distribution
$\Lambda_f\in({\mathcal{G}}^\varepsilon_0(\beta,\gamma))'$ is
injective whenever $\mu$ is Borel-semiregular,
we simply write $f$ in place of $\Lambda_f$ in what follows.
\end{remark}

\begin{proof}[Proof of Proposition~\ref{distfncttype}]
Let $\rho\in\mathbf{q}$ be a regularized quasi-ultrametric.
Assume that $f\in L^r(X)$ and $\phi\in{\mathcal{G}}^\varepsilon_0(\beta,\gamma)
\subset{\mathcal{G}}(\beta,\gamma)$
and recall that, in Remark~\ref{Rmtest}(i),
we set ${\mathcal{G}}(\beta,\gamma):={\mathcal{G}}(x_1,1,\beta,\gamma)$,
where $x_1\in X$ is some fixed point. On the one hand, it is easy to see from
Definition~\ref{testfunction} that,
for any $x\in B_\rho(x_1,1)$,
\begin{align}\label{1803}
\phi\in {\mathcal{G}}(x,1,\beta,\gamma)
\ \ \text{and}\ \
\|\phi\|_{{\mathcal{G}}(x,1,\beta,\gamma)}
\lesssim\|\phi\|_{{\mathcal{G}}(x_1,1,\beta,\gamma)},
\end{align}
where the implicit positive
constant is independent of $f$, $\phi$, and $x$. On the other hand, by
Lemma~\ref{M}, we find that
$\mathcal{M}_{\rho} f$ belongs to $L^r(X)$ if $r>1$ and to $L^{1,\infty}(X)$ if $r=1$,
from which we can deduce that
there exists a point $x\in B_\rho(x_1,1)$
such that
$\mathcal{M}_\rho f(x)<\infty$.
Using this, we find that
\begin{align}
\label{alk-34}
\left|\left\langle\Lambda_f,\varphi\right\rangle\right|
&\leq\int_{\rho(x,y)<1}|f(y)\phi(y)|\,d\mu(y)
+\int_{\rho(x,y)\geq1}|f(y)\phi(y)|\,d\mu(y)
\nonumber\\
&\leq\|\phi\|_{{\mathcal{G}}(x,1,\beta,\gamma)}
\frac{1}{(V_\rho)_1(x)}\int_{\rho(x,y)<1}|f(y)|\,d\mu(y)
\nonumber\\
&\quad
+\|\phi\|_{{\mathcal{G}}(x,1,\beta,\gamma)}
\int_{\rho(x,y)\geq1}|f(y)|\frac{1}{V_\rho(x,y)}
\left[\frac{1}{\rho(x,y)}\right]^\gamma\,d\mu(y)
\nonumber\\
&\qquad\ \text{by Definition~\ref{testfunction}(i)}
\nonumber\\
&\lesssim\|\phi\|_{{\mathcal{G}}(x,1,\beta,\gamma)}\mathcal{M}_\rho f(x)
\quad\text{by Lemma~\ref{A3.2}(iv)}
\nonumber\\
&\lesssim\|\phi\|_{{\mathcal{G}}(x_1,1,\beta,\gamma)}\mathcal{M}_\rho f(x)<\infty\quad\text{by \eqref{1803}}.
\end{align}
Hence, $\Lambda_f\in({\mathcal{G}}^\varepsilon_0(\beta,\gamma))'$, as wanted.

As concerns the last claim in the statement of this proposition,
it is clear
from the formula of $\Lambda_f$ that, if $f=0$ pointwise
$\mu$-almost everywhere on $X$, then $\Lambda_f=0$
in $({\mathcal{G}}^\varepsilon_0(\beta,\gamma))'$.
To see that the converse is also true, suppose that
$\Lambda_f=0$ in $({\mathcal{G}}^\varepsilon_0(\beta,\gamma))'$ for some $f\in L^r(X)$.
We claim that, for any $u\in L^\infty(X,\mu)$ with bounded support,
\begin{equation}\label{DKJ-092-X}
\int_Xf(x)u(x)\,d\mu(x)=0.
\end{equation}
To show \eqref{DKJ-092-X}, let $u\in L^\infty(X,\mu)$
be any function that has bounded support
and suppose that
$\{\mathcal{S}_{2^{-k}}\}_{k\in\mathbb{Z}}$ is
an approximation of the identity of order
$\varepsilon$ as in Definition~\ref{DefHATI}
with $t$ replaced by $2^{-k}$.
By (vi) and (viii) of Theorem~\ref{corATI} and by Lemma~\ref{CGL},
we conclude that, for any $k\in\mathbb{Z}$,
\begin{align*}
{\mathcal{S}}_{2^{-k}}u\in\dot{C}^\varepsilon_{\mathrm{b}}(X)\subset
{\mathcal{G}}(\varepsilon,\varepsilon)\subset{\mathcal{G}}^\varepsilon_0(\beta,\gamma),
\end{align*}
which, together with the assumption that $\Lambda_f=0$
in $({\mathcal{G}}^\varepsilon_0(\beta,\gamma))'$,
implies that, for any $k\in\mathbb{Z}$,
\begin{align*}
\int_Xf(x){\mathcal{S}}_{2^{-k}}u(x)\,d\mu(x)=0.
\end{align*}
Since $u\in L^\infty(X,\mu)$ has bounded support, we find
that $u\in L^1(X,\mu)$ and hence it follows from Theorem~\ref{corATI}(ix)
(keeping in mind that we are now assuming that $\mu$ is Borel-semiregular)
that there exists a numerical sequence
$\{k_j\}_{j\in\mathbb{N}}$ in $\mathbb{N}$
such that $\lim_{j\to\infty}k_j=\infty$ and, for $\mu$-almost every $x\in X$,
\begin{equation}\label{WWW.111}
\lim_{j\to\infty}{\mathcal{S}}_{2^{-k_j}}u(x)=u(x).
\end{equation}
Moreover, by Theorem~\ref{corATI}(viii),
we conclude that there exist $x_0\in X$
and $r_0\in(0,\infty)$ such that, for any $j\in\mathbb{N}$,
\begin{equation}\label{DBR-854}
{\rm supp}\,{\mathcal{S}}_{2^{-k_j}}u\subset
B_{\rho}(x_0,r_0).
\end{equation}
On the other hand, Theorem~\ref{corATI}(i) (used here with $p=\infty$)
implies that, for any $j\in\mathbb{N}$,
\begin{equation*}
\|{\mathcal{S}}_{2^{-k_j}}u\|_{L^\infty(X)}\lesssim
\|u\|_{L^\infty(X)},
\end{equation*}
where the implicit positive constant is independent of $u$, which,
combined with \eqref{DBR-854}, gives
\begin{equation*}
|f{\mathcal{S}}_{2^{-k_j}}u|\lesssim
\|u\|_{L^\infty(X)}|f|{\mathbf{1}}_{B_{\rho}(x_0,r_0)}
\in L^1(X).
\end{equation*}
From this, \eqref{WWW.111}, and the Lebesgue
dominated convergence theorem, we infer that
\begin{equation*}
0=\lim_{j\to\infty}\int_Xf(x){\mathcal{S}}_{2^{-k_j}}u(x)\,d\mu(x)=
\int_Xfu\,d\mu,
\end{equation*}
which completes the proof of \eqref{DKJ-092-X}.

Next, fix $x_\ast\in X$ and, for each $j\in\mathbb{N}$, define
$A_j:=\{x\in B_{\rho}(x_\ast,j): f(x)>0\}$.
Then
\begin{align}\label{tijn444aa}
\bigcup\limits_{j=1}^{\infty}A_j=A_\ast,
\end{align}
where $A_\ast:=\{x\in X:f(x)\neq0\}$.
By design, each $A_j$ is $\mu$-measurable with finite
$\mu$-measure and hence each ${\mathbf{1}}_{A_j}\in L^{\infty}(X,\mu)$
has bounded support in $X$.
As such, by \eqref{DKJ-092-X} for any $j\in\mathbb{N}$,
\begin{align*}
\int_{A_j}f(x)\,d\mu(x)=0.
\end{align*}
Since $f>0$ pointwise on $A_j$, we have $\mu(A_j)=0$.
From this and \eqref{tijn444aa},
we deduce that $\mu(A_\ast)=0$. That is,
$f(x)\leq0$ for $\mu$-almost every $x\in X$.
In order to complete the proof of the present proposition, introduce the set
$\widetilde{A}_\ast:=\{x\in X: -f(x)>0\}$.
Then using a similar reasoning as above with the sets $A_j$ replaced with
\begin{align*}
\widetilde{A}_j:=\left\{x\in B_{\rho}(x_\ast,j):f(x)<0\right\}
\end{align*}
for any $j\in\mathbb{N}$
implies $\mu(\widetilde{A}_\ast)=0$. Hence,
$f(x)\geq0$ for $\mu$-almost every $x\in X$ and hence we find
that $f=0$ pointwise $\mu$-almost everywhere on $X$,
as desired. This finishes the proof of Proposition~\ref{distfncttype}.
\end{proof}

\section[Boundedness of Calder\'on--Zygmund
Operators\\ on Spaces of Test Functions]{Boundedness of Calder\'on--Zygmund
Operators \\ on Spaces of Test Functions}
\label{SS3.1}

In this section, we aim to establish
the boundedness of Calder\'on--Zygmund
operators on spaces of test functions.
To this end,
we first recall the definition of the
Calder\'on--Zygmund kernel as follows.

\begin{definition}
\label{CZk}
Let $(X,\mathbf{q},\mu)$ be a quasi-ultrametric
space of homogeneous type. Assume that
$\rho\in\mathbf{q}$ and $0<\varepsilon\preceq{\mathrm{ind\,}}(X,\mathbf{q})$
with ${\mathrm{ind\,}}(X,\mathbf{q})$ as in \eqref{index}.
A function
$$
K:(X\times X)\setminus\{(x,x):x\in X\}\to\mathbb{C}
$$
is called an
\emph{$\varepsilon$-Calder\'on--Zygmund kernel}
\index[words]{Calder\'on--Zygmund kernel}
(with respect to $\rho$) if all $\rho$-balls
are $\mu$-measurable and there exist two
constants $C_\ast\in(0,1)$ and $C_T\in(0,\infty)$  such that
\begin{enumerate}
\item For any $x,y\in X$ with $x\neq y$,
\begin{equation}\label{CZk1}
|K(x,y)|\leq C_T\frac{1}{V_\rho(x,y)}.
\end{equation}

\item
For any $x,x',y\in X$ satisfying $\rho(x,x')\leq C_\ast\rho(x,y)$
with $x\neq y$,
\begin{align}\label{CZk2}
\left|K(x,y)-K(x',y)\right|
+\left|K(y,x)-K(y,x')\right|
\leq C_T\left[\frac{\rho(x,x')}{\rho(x,y)}
\right]^\varepsilon\frac{1}{V_\rho(x,y)}.
\end{align}
\end{enumerate}
\end{definition}

\begin{remark}
\begin{enumerate}
\item[\rm(i)]
Let $0<\varepsilon\preceq{\mathrm{ind\,}}(X,\mathbf{q})$
and $K$ be an $\varepsilon$-Calder\'on--Zygmund
kernel with respect to $\rho$ as in Definition~\ref{CZk}.
Then $K$ is also an $\varepsilon$-Calder\'on--Zygmund
kernel with respect to
any other quasi-ultrametric $\rho'\in\mathbf{q}$
for which all $\rho'$-balls are $\mu$-measurable.
In particular, for the remainder of this monograph,
without loss of generality, we may always assume
that $\rho$ in Definition~\ref{CZk} is a
regularized quasi-ultrametric as in Definition~\ref{D2.3}.

\item[\rm(ii)]
In order to simplify the presentation of the
calculations in this chapter, we assume that
$C_\ast=(2C_\rho)^{-1}$ in Definition~\ref{CZk},
where $C_\rho$ is the same positive constant as in \eqref{Crho}.
All of the proofs that use
this assumption can be suitably modified for a general constant $C_\ast\in(0,1)$.
\end{enumerate}
\end{remark}

Now, we recall the definition of the Calder\'on--Zygmund operator as follows.

\begin{definition}
\label{CZO}
Let $(X,\mathbf{q},\mu)$ be a quasi-ultrametric space of homogeneous type.
Assume that
$\rho\in\mathbf{q}$ and $0<\varepsilon\preceq{\mathrm{ind\,}}(X,\mathbf{q})$
with ${\mathrm{ind\,}}(X,\mathbf{q})$ as in \eqref{index}.
A linear operator
$$T:C^\varepsilon_{\mathrm{b}}(X)\to
\left(C^\varepsilon_{\mathrm{b}}(X)\right)'$$
is called an
\emph{$\varepsilon$-Calder\'on--Zygmund operator}
\index[words]{Calder\'on--Zygmund operator}
(with respect to $\rho$)
if $T$ can be extended to a bounded linear operator on $L^2(X)$
and if there exists an $\varepsilon$-Calder\'on--Zygmund kernel $K$
(with respect to $\rho$) as in Definition~\ref{CZk}
such that, for any $f\in L^2(X)$ and $x\notin\overline{{\mathrm{supp\,}} f}$,
$$
Tf(x)=\int_XK(x,y)f(y)\,d\mu(y),
$$
where $(C^\varepsilon_{\mathrm{b}}(X))'$ denotes the \emph{dual
space} of $C^\varepsilon_{\mathrm{b}}(X)$ equipped with the weak-$*$ topology
and $\overline{{\mathrm{supp\,}} f}$ denotes the closure of
${\mathrm{supp\,}} f:=\{x\in X:f(x)\neq 0\}$ in $X$.
\end{definition}

\begin{remark}
If $T$ is an $\varepsilon$-Calder\'on--Zygmund operator with respect to $\rho$ as in Definition~\ref{CZO}
for some $0<\varepsilon\preceq{\mathrm{ind\,}}(X,\mathbf{q})$,
then $T$ is also an $\varepsilon$-Calder\'on--Zygmund operator with respect to
any other quasi-ultrametric $\rho'\in\mathbf{q}$
for which all $\rho'$-balls are $\mu$-measurable.
In particular, for the remainder of this monograph,
without loss of generality, we may always assume
that $\rho$ in Definition~\ref{CZO} is a
regularized quasi-ultrametric as in Definition~\ref{D2.3}.
\end{remark}

Next, we establish the boundedness
of Calder\'on--Zygmund operators on
$\mathring{\mathcal{G}}(x_1,r,\beta,\gamma)$ which is a sharpening of
\cite[Theorem~3.1]{HLYY}.

\begin{theorem}
\label{CZOonGO}
Let $(X,\mathbf{q},\mu)$ be a
quasi-ultrametric space of homogeneous type.
Assume that
$\rho\in\mathbf{q}$ has the property that all $\rho$-balls are $\mu$-measurable and
$0<\beta,\gamma<\varepsilon\preceq{\mathrm{ind\,}}(X,\mathbf{q})$
with ${\mathrm{ind\,}}(X,\mathbf{q})$ as in \eqref{index}.
Assume that $T$ is an
$\varepsilon$-Calder\'on--Zygmund operator
(with respect to $\rho$) and its kernel $K$
satisfies the following additional conditions:
\begin{enumerate}
\item[\textup{(i)}]  there exists a
positive constant $C_T$ such that, for any
$x,x',y,y'\in X$ satisfying
$\rho(x,x')\leq(2C_\rho)^{-2}\rho(x,y)$,
$\rho(y,y')\leq(2C_\rho)^{-2}\rho(x,y)$, and
$x\neq y$,
\begin{align}
\label{CZO1}
&\left|\left[K(x,y)-K(x',y)\right]
-\left[K(x,y')-K(x',y')\right]\right|
\nonumber\\
&\quad\leq C_T\left[\frac{\rho(x,x')}{\rho(x,y)}\right]^\varepsilon
\left[\frac{\rho(y,y')}{\rho(x,y)}\right]^\varepsilon
\frac{1}{V_\rho(x,y)},
\end{align}
where $C_\rho$ is the same positive constant as in \eqref{Crho};
\item[\textup{(ii)}]
for any $f\in C^\beta(X)$ and $x\in X$,
$$
Tf(x)=\int_XK(x,y)f(y)\,d\mu(y);
$$
\item[\textup{(iii)}]
there exists a constant $c_0\in\mathbb{C}$ such that,
for any $x\in X$, $K(x,\cdot)$ is integrable and
\begin{equation}\label{CZO2}
\int_XK(x,y)\,d\mu(y)=c_0.
\end{equation}
\end{enumerate}
Then there exists a positive constant
$C$ such that,
for any $f\in\mathring{\mathcal{G}}(x_1,r,\beta,\gamma)$
with $x_1\in X$ and $r\in(0,\infty)$,
$$
\|Tf\|_{{\mathcal G}(x_1,r,\beta,\gamma)}\leq C
\left[C_T+\|T\|_{L^2(X)\to L^2(X)}+|c_0|\right]
\|f\|_{{\mathcal G}(x_1,r,\beta,\gamma)}.
$$
\end{theorem}

To prove Theorem~\ref{CZOonGO}, we need
some technical lemmas.
The following lemma is a sharpening of \cite[Lemma~3.4]{HLYY}
(see also \cite[(3.24)]{HLW})
in the setting of quasi-ultrametric spaces of
homogeneous type,
which is deduced from an argument similar to
that used in the proof of \cite[(3.24)]{HLW};
we present the details
here for the convenience of the reader.

\begin{lemma}
\label{Lem3.4L}
Let $(X,\mathbf{q},\mu)$ be a quasi-ultrametric space
of homogeneous type,
$\rho\in\mathbf{q}$ have the property that
all $\rho$-balls are $\mu$-measurable,
$0<\varepsilon\preceq{\mathrm{ind\,}}(X,\mathbf{q})$,
and $T$ be an $\varepsilon$-Calder\'on--Zygmund
operator with its kernel $K$
satisfying the following additional conditions:
\begin{enumerate}
\item[\textup{(i)}]
For any $\beta\in(0,\varepsilon)$, $f\in C^\beta(X)$, and $x\in X$,
$$
Tf(x)=\int_XK(x,y)f(y)\,d\mu(y).
$$

\item[\textup{(ii)}]
There exists a constant $c_0\in\mathbb{C}$ such that,
for any $x\in X$, $K(x,\cdot)$ is integrable and
\begin{align}\label{intKxy=0}
\int_XK(x,y)\,d\mu(y)=c_0.
\end{align}
\end{enumerate}
Assume that $x_0\in X$, $t\in(0,\infty)$,
and a function $f\in C_{\rm b}^\varepsilon(X)$ satisfies both
$|f|\leq\mathbf{1}_{B_\rho(x_0,t)}$
and $\|f\|_{\dot{C}^\varepsilon(X)}\leq t^{-\varepsilon}$.
Then there exists a positive constant $C$,
independent of $f$, $x_0$, $t$, and $T$, such that,
for any $x\in X$,
\begin{align}\label{Tf(y)}
\left|Tf(x)\right|
\leq C\left[C_T+\|T\|_{L^2(X)\to L^2(X)}+|c_0|\right].
\end{align}
\end{lemma}

\begin{proof}
Without loss of generality,
we may assume that $\rho\in\mathbf{q}$ is a regularized quasi-ultrametric
as in Definition~\ref{D2.3}.
We first claim that there exists a positive constant $K_1$,
independent of $f$, $x_0$, $t$, and $T$, such that
\begin{align}\label{TfLinfty}
\left\|Tf\right\|_{L^\infty(X)}
\leq K_1\left[C_T+\|T\|_{L^2(X)\to L^2(X)}+|c_0|\right].
\end{align}
Let $h\in C_{\rm b}^\varepsilon(X)$ be a smooth cut-off function satisfying
the conclusions of Corollary~\ref{Lem3.3} with
\begin{align*}
h\equiv1\text{ on }B_\rho(x_0,2C_\rho t)
\ \ \text{and}\ \
h\equiv0\text{ on }[B_\rho(x_0,4C_\rho^2t)]^\complement.
\end{align*}
Then $f=fh\in C_{\rm b}^\varepsilon(X)\subset C^\beta(X)$
and, for any $\phi\in C_{\rm b}^\varepsilon(X)$,
\begin{align}\label{2027}
\left\langle Tf,\phi\right\rangle
&=\int_{\rho(x_0,x)\leq4C_\rho^3t}
\left\{\int_{B_\rho(x_0,4C_\rho^2t)}K(x,y)\left[f(y)-f(x)\right]
h(y)\,d\mu(y)\right\}\phi(x)\,d\mu(x)
\nonumber\\
&\quad+\int_{\rho(x_0,x)>4C_\rho^3t}
\left\{\int_{B_\rho(x_0,4C_\rho^2t)}K(x,y)\left[f(y)-f(x)\right]
h(y)\,d\mu(y)\right\}\phi(x)\,d\mu(x)
\nonumber\\
&\quad+\int_{B_\rho(x_0,t)}\left[\int_XK(x,y)h(y)\,d\mu(y)
\right]f(x)\phi(x)\,d\mu(x)
\nonumber\\
&=:\mathrm{I}+\mathrm{II}+\mathrm{III}.
\end{align}
To estimate $\mathrm{I}$,
from the fact that
$B_\rho(x_0,4C_\rho^2t)\subset B_\rho(x,4C_\rho^4t)$
for any given $x\in B_\rho(x_0,4C_\rho^3t)$,
we deduce that
\begin{align}\label{2028}
\left|\mathrm{I}\right|
&\leq\int_{\rho(x_0,x)\leq4C_\rho^3t}\left\{\int_{B_\rho(x_0,4C_\rho^2t)}
C_T\frac{1}{V_\rho(x,y)}\left[\frac{\rho(x,y)}{t}\right]^\varepsilon\,d\mu(y)\right\}
|\phi(x)|\,d\mu(x)\nonumber\\
&\qquad\ \text{by \eqref{CZk1}}
\nonumber\\
&\leq\int_{\rho(x_0,x)\leq4C_\rho^3t}\left\{\int_{B_\rho(x,4C_\rho^4t)}
C_T\frac{1}{V_\rho(x,y)}\left[\frac{\rho(x,y)}{t}\right]^\varepsilon\,d\mu(y)\right\}
|\phi(x)|\,d\mu(x)
\nonumber\\
&\lesssim C_T\|\phi\|_{L^1(X)}
\quad\text{by Lemma~\ref{A3.2}(i)}.
\end{align}
To estimate $\mathrm{II}$,
by Definition~\ref{D2.1}(iii),
we find that, for any $y\in B_\rho(x_0,4C_\rho^2t)$
and $x\in[B_\rho(x_0,4C_\rho^3 t)]^\complement$,
it holds $\rho(x,y)\ge\frac{\rho(x_0,x)}{C_\rho}$,
which, together with \eqref{upperdoub},
further implies that $V_\rho(x_0,x)\lesssim V_\rho(x,y)$.
From this, the fact that $f(x)=0$ in this case,
the assumption $|f|\leq\mathbf{1}_{B_\rho(x_0,t)}$,
and \eqref{CZk1},
we infer that
\begin{align}\label{2029}
\left|\mathrm{II}\right|
&\leq\int_{\rho(x,x_0)>4C_\rho^3t}\left[
\int_{B_\rho(x_0,4C_\rho^2t)}C_T\frac{1}{V_\rho(x,y)}\,d\mu(y)
\right]\left|\phi(x)\right|\,d\mu(x)
\nonumber\\
&\lesssim C_T\frac{\mu(B_\rho(x_0,t))}{V_\rho(x_0,x)}\|\phi\|_{L^1(X)}
\leq C_T\|\phi\|_{L^1(X)}.
\end{align}
To estimate $\mathrm{III}$, we first show that
\begin{align}\label{2020}
\|Th\|_{L^\infty(B_\rho(x_0,t))}\lesssim\|T\|_{L^2(X)\to L^2(X)}+|c_0|+C_T.
\end{align}
To prove this, we use Meyer's idea in \cite{m1985}.
Let $g:=1-h\in C^\beta(X)$ and suppose that
$\psi\in C^\varepsilon_{\rm b}(X)$ satisfies
both $\mathrm{supp\,}\psi\subset B_\rho(x_0,t)$
and $\int_X\psi(x)\,d\mu(x)=0$.
Then, by \eqref{intKxy=0} and $\mathrm{supp\,}g
\subset[B_\rho(x_0,2C_\rho t)]^\complement$,
we conclude that
\begin{align*}
&\left|\int_{B_\rho(x_0,t)}Th(x)\psi(x)\,d\mu(x)\right|
\\
&\quad=\left|\int_{B_\rho(x_0,t)}T1(x)\psi(x)\,d\mu(x)-
\int_{B_\rho(x_0,t)}Tg(x)\psi(x)\,d\mu(x)\right|
\\
&\quad\leq|c_0|\,\|\psi\|_{L^1(B_\rho(x_0,t))}
+\left|\int_{B_\rho(x_0,t)}Tg(x)\psi(x)\,d\mu(x)\right|
\\
&\quad=|c_0|\,\|\psi\|_{L^1(B_\rho(x_0,t))}
+\left|\int_{B_\rho(x_0,t)}\left[Tg(x)-Tg(x_0)\right]\psi(x)\,d\mu(x)\right|,
\end{align*}
which, combined with \eqref{CZk2} and Lemma~\ref{A3.2}(i),
further implies that
\begin{align*}
&\left|\int_{B_\rho(x_0,t)}Th(x)\psi(x)\,d\mu(x)\right|
\\
&\quad\leq|c_0|\,\|\psi\|_{L^1(B_\rho(x_0,t))}
\\
&\qquad+\int_{B_\rho(x_0,t)}\left\{\int_{[B_\rho(x_0,2C_\rho t)]^\complement}
C_T\left[\frac{\rho(x_0,x)}{\rho(x_0,y)}\right]^\varepsilon
\frac{1}{V_\rho(x_0,y)}\,d\mu(y)\right\}
|\psi(x)|\,d\mu(x)
\\
&\quad\lesssim\left(|c_0|+C_T\right)\|\psi\|_{L^1(B_\rho(x_0,t))}.
\end{align*}
From this, the Hahn--Banach theorem,
and $(L^1(B_\rho(x_0,t)))'=L^\infty(B_\rho(x_0,t))$,
we deduce that $Th(x)=\Lambda+H(x)$ for any $x\in B_\rho(x_0,t)$,
where $\Lambda$ is a constant and
\begin{align}\label{2018}
\|H\|_{L^\infty(B_\rho(x_0,t))}\lesssim|c_0|+C_T.
\end{align}
To estimate $\Lambda$, choose a nonnegative
function $\phi_1\in C_{\rm b}^\varepsilon(X)$
with $\mathrm{supp\,}\phi_1\subset B_\rho(x_0,t)$,
$\|\phi_1\|_{L^\infty(X)}\leq1$,
$\|\phi_1\|_{\dot{C}^\varepsilon(X)}\leq t^{-\varepsilon}$,
and
$$
\int_X\phi_1(x)\,d\mu(x)=\widetilde{C}\mu\left(B_\rho(x_0,t)\right),
$$
where the positive constant $\widetilde{C}$
is independent of $f$, $x_0$, $t$, and $T$.
Then, using H\"older's inequality,
we find that
\begin{align*}
|\Lambda|\widetilde{C}\mu(B_\rho(x_0,t))
&=\int_X|\Lambda|\phi_1(x)\,d\mu(x)\\
&\leq\int_X|Th(x)|\phi_1(x)\,d\mu(x)+\int_XH(x)\phi_1(x)\,d\mu(x)\\
&\leq\left\|Th\right\|_{L^2(X)}\|\phi_1\|_{L^2(X)}
+\|H\|_{L^\infty(X)}\widetilde{C}\mu(B_\rho(x_0,t)).
\end{align*}
From this, the boundedness of $T$ on $L^2(X)$, and \eqref{2018},
it follows that
\begin{align*}
|\Lambda|\widetilde{C}\mu(B_\rho(x_0,t))
&\leq\|T\|_{L^2(X)\to L^2(X)}\left\|h\right\|_{L^2(X)}\|\phi_1\|_{L^2(X)}
+\|H\|_{L^\infty(X)}\widetilde{C}\mu(B_\rho(x_0,t))\\
&\lesssim\left[\|T\|_{L^2(X)\to L^2(X)}+|c_0|+C_T\right]\mu(B_\rho(x_0,t)),
\end{align*}
which, together with \eqref{2018} again,
further implies that
\begin{align*}
\|Th\|_{L^\infty(B_\rho(x_0,t))}\lesssim\|T\|_{L^2(X)\to L^2(X)}+|c_0|+C_T
\end{align*}
and hence the proof of \eqref{2020} is complete.

By \eqref{2020}, we conclude that
\begin{align*}
\left|\mathrm{III}\right|
&=\left|\int_{B_\rho(x_0,t)}Th(x)f(x)\phi(x)\,d\mu(x)\right|\\
&\lesssim\left[\|T\|_{L^2(X)\to L^2(X)}+|c_0|+C_T\right]\|\phi\|_{L^1(X)}.
\end{align*}
Combining this, \eqref{2027}, \eqref{2028}, and \eqref{2029},
we obtain, for any $\phi\in C_{\rm b}^\varepsilon(X)$,
\begin{align}\label{2033}
\left|\left\langle Tf,\phi\right\rangle\right|
\lesssim\left[\|T\|_{L^2(X)\to L^2(X)}+|c_0|+C_T\right]\|\phi\|_{L^1(X)},
\end{align}
which, together with Lemma~\ref{LDT}(iii) and a standard density argument,
further implies that \eqref{2033} holds for any $\phi\in L^1(X)$.
Hence, $Tf\in L^\infty(X)$ with
\begin{align*}
\|Tf\|_{L^\infty(X)}\lesssim\left[\|T\|_{L^2(X)\to L^2(X)}+|c_0|+C_T\right].
\end{align*}
This finishes the proof of \eqref{TfLinfty} and hence
the aforementioned claim.

From \eqref{TfLinfty}, we infer that there exists a set $X_1\subset X$
satisfying that $\mu(X\setminus X_1)=0$ and that
\eqref{Tf(y)} holds for any $x\in X_1$.
To complete the proof of the present lemma, it remains to
show that, for any $x\in X\setminus X_1$,
\eqref{Tf(y)} holds. To this end,
let $x\in X\setminus X_1$. Then $\mu(\{x\})\leq\mu(X\setminus X_1)=0$.
By this and Lemma~\ref{mu-r}, we find that
$r_\rho(x)=0$.
From this, we deduce that there exists a sequence
$\{x_n\}_{n\in\mathbb{N}}$ in $B_\rho(x,t)$
of points such that
\begin{align}\label{1624}
\lim_{n\to\infty}\rho(x,x_n)=0.
\end{align}
Indeed, we can choose $\{x_n\}_{n\in\mathbb{N}}$ in $B_\rho(x,t)\cap X_1$
such that \eqref{1624} holds.
Otherwise, there exists $r\in(0,t)$ such that
$B_\rho(x,r)\cap X_1=\varnothing$
and hence $\mu(B_\rho(x,r))\leq\mu(X\setminus X_1)=0$,
which contradicts the assumption that
$\mu(B_\rho(y,r))\in(0,\infty)$ for any $y\in X$ and $r\in(0,\infty)$.
By \eqref{intKxy=0}, we conclude that, for any $n\in\mathbb{N}$,
\begin{align}\label{3.41}
\left|Tf(x)-Tf(x_n)\right|
&=\left|\int_X\left[K(x,y)-K(x_n,y)\right]
\left[f(y)-f(x_n)\right]\,d\mu(y)\right|\nonumber\\
&\leq I_1(x,x_n)+I_2(x,x_n)+I_3(x,x_n)+I_4(x,x_n),
\end{align}
where,
\begin{align*}
I_1(x,x_n):=\int_{\rho(x_n,y)\ge2C_\rho\rho(x,x_n)}
\left|K(x,y)-K(x_n,y)\right|\left[\left|f(y)\right|
+\left|f(x_n)\right|\right]\,d\mu(y),
\end{align*}
\begin{align*}
I_2(x,x_n):=\int_{\rho(x_n,y)<2C_\rho\rho(x,x_n)}
\left|K(x_n,y)\right|
\left|f(y)-f(x_n)\right|\,d\mu(y),
\end{align*}
\begin{align*}
I_3(x,x_n):=\int_{\rho(x_n,y)<2C_\rho\rho(x,x_n)}
\left|K(x,y)\right|
\left|f(y)-f(x)\right|\,d\mu(y),
\end{align*}
and
\begin{align*}
I_4(x,x_n):=\left|f(x)-f(x_n)\right|\int_{\rho(x_n,y)<2C_\rho\rho(x,x_n)}
\left|K(x,y)\right|\,d\mu(y).
\end{align*}
For $I_1(x,x_n)$, from \eqref{CZk2}, $|f|\leq1$,
and Lemma~\ref{A3.2}(i), we infer that
\begin{align}\label{3.42}
I_1(x,x_n)\lesssim\int_{\rho(x_n,y)\ge2C_\rho\rho(x,x_n)}
\left[\frac{\rho(x,x_n)}{\rho(x_n,y)}\right]^\varepsilon
\frac{1}{V_\rho(x_n,y)}\,d\mu(y)\lesssim1,
\end{align}
where the implicit positive constants depend on $C_T$.
For $I_2(x,x_n)$,
by \eqref{CZk1},
$\|f\|_{\dot{C}^\varepsilon(X)}\leq t^{-\varepsilon}$,
$\rho(x,x_n)<t$,
and Lemma~\ref{A3.2}(i), we find that
\begin{align}\label{3.44}
I_2(x,x_n)
&\lesssim\int_{\rho(x_n,y)<2C_\rho\rho(x,x_n)}
\frac{1}{V_\rho(x_n,y)}
\left[\frac{\rho(x_n,y)}{t}\right]^\varepsilon\,d\mu(y)\nonumber\\
&\leq\int_{\rho(x_n,y)<2C_\rho t}
\frac{1}{V_\rho(x_n,y)}
\left[\frac{\rho(x_n,y)}{t}\right]^\varepsilon\,d\mu(y)\lesssim1,
\end{align}
where the implicit positive constants depend on both $C_T$ and $C_\rho$.
From an argument similar to that used in
the estimation of \eqref{3.44}, we deduce that
\begin{align}\label{3.45}
I_3(x,x_n)\lesssim1,
\end{align}
where the implicit positive constant only depends on $C_T$ and $C_\rho$.
For $I_4(x,x_n)$, note that, \eqref{intKxy=0} implies that,
for any given $x\in X$,
\begin{align*}
\int_X\left|K(x,y)\right|\,d\mu(y)<\infty.
\end{align*}
By this, \eqref{1624}, and \eqref{C(X)} with
$\|f\|_{\dot{C}^\varepsilon(X)}\leq t^{-\varepsilon}$, we conclude that
\begin{align*}
I_4(x,x_n)\leq\left[\frac{\rho(x,x_n)}{t}\right]^\varepsilon
\int_X\left|K(x,y)\right|\,d\mu(y)\to0
\end{align*}
as $n\to\infty$, which, further implies that
there exists $N_{(x)}\in\mathbb{N}$
such that $I_4(x,x_{N_{(x)}})\leq C_T$.
From this, \eqref{3.41}, \eqref{3.42},
\eqref{3.44}, \eqref{3.45}, and \eqref{Tf(y)}
(keeping in mind that $x_{N_{(x)}}\in X_1$),
it follows that there exists a positive constant $K_2$
such that
\begin{align*}
\left|Tf(x)\right|
&\leq\left|Tf(x)-Tf(x_{N_{(x)}})\right|+\left|Tf(x_{N_{(x)}})\right|
\leq K_2C_T+\left|Tf(x_{N_{(x)}})\right|\\
&\leq(K_1+K_2)\left[C_T+\|T\|_{L^2(X)\to L^2(X)}
+|c_0|\right],
\end{align*}
which further implies that \eqref{Tf(y)}
holds for any $x\in X\setminus X_1$
and hence holds for any $x\in X$.
This finishes the proof of Lemma~\ref{Lem3.4L}.
\end{proof}

Now, we prove Theorem~\ref{CZOonGO}.
Although a version of Theorem~\ref{CZOonGO}
was shown in \cite[Theorem~3.1]{HLYY}
on spaces of homogeneous type (for a smaller range of parameters),
we still include its proof here with some additional details for
the sake of completeness.

\begin{proof}[Proof of Theorem~\ref{CZOonGO}]
Let $A:=2C_\rho$, where $C_\rho$ is the same as in \eqref{Crho}.
Without loss of generality, we may assume that
$\rho\in\mathbf{q}$ is a regularized quasi-ultrametric
as in Definition~\ref{D2.3}.
By Remark~\ref{RmkD2.3}(iii),
we can also assume  that $\rho\in\mathbf{q}$
satisfies $\varepsilon\leq(\log_2C_\rho)^{-1}$.
Moreover, we may assume
$\|f\|_{{\mathcal G}(x_1,r,\beta,\gamma)}=1$
and
$$C_T+\|T\|_{L^2(X)\to L^2(X)}+|c_0|\leq1.$$
From Lemma~\ref{Lem3.4L}, it follows that there exists a positive
constant $C$ such that, for any $g\in\dot{C}^\varepsilon(X)$
as in Lemma~\ref{Lem3.4L} and for any $x\in X$,
\begin{align}
\label{Lem3.4}
|Tg(x)|\leq C\left[C_T+\|T\|_{L^2(X)\to L^2(X)}
+|c_0|\right].
\end{align}
To verify Theorem~\ref{CZOonGO},
we divide the proof into the following two steps.

\emph{Step 1.} In this step, we prove that,
for any $x\in X$,
\begin{align}\label{Tfs}
|Tf(x)|\lesssim\frac{1}{(V_\rho)_r(x_1)+V_\rho(x_1,x)}
\left[\frac{r}{r+\rho(x_1,x)}\right]^\gamma.
\end{align}
To this end, let $x\in X$.
We consider the following two cases.

\emph{Case (1)} $\rho(x_1,x)<Ar$.
In this case, note that
$\rho(x_1,x)+r\sim r$, which, combined with
Lemma~\ref{A3.2}(vii),
further implies that
\begin{align*}
V_\rho(x_1,x)+(V_\rho)_r(x_1)\sim (V_\rho)_{r+\rho(x_1,x)}(x_1)
\sim (V_\rho)_r(x_1)\sim (V_\rho)_r(x).
\end{align*}
Thus, to show \eqref{Tfs},
it suffices to prove
\begin{align*}
|Tf(x)|\lesssim\frac{1}{(V_\rho)_r(x)}.
\end{align*}
By Corollary~\ref{Lem3.3}, we conclude that there exists
a function $\theta\in{C}^\varepsilon(X)$ satisfying
\begin{align}\label{theta}
\mathbf{1}_{B_\rho(x_1,A^2r)}\leq\theta
\leq\mathbf{1}_{B_\rho(x_1,A^3r)}
\end{align}
and $\|\theta\|_{\dot{C}^\varepsilon(X)}
\lesssim r^{-\varepsilon}$. Let $\psi:=1-\theta$.
From the fact that $f\in{\mathcal G}(x_1,r,\beta,\gamma)
\subset C^\beta(X)$ and (ii) of the present theorem,
we deduce that
\begin{align}\label{Tf=123}
Tf(x)
&=\int_X K(x,y)[f(y)-f(x)]\theta(y)\,d\mu(y)
\nonumber\\
&\quad+\int_X K(x,y)f(y)\psi(y)\,d\mu(y)
+f(x)\int_X K(x,y)\theta(y)\,d\mu(y)
\nonumber\\
&=:Z_1(x)+Z_2(x)+Z_3(x).
\end{align}

To deal with $Z_1(x)$, since $\rho(x_1,x)<Ar$,
by \eqref{theta}, Definition~\ref{D2.1}(iii),
and the fact that regularized quasi-ultrametrics are symmetric,
we find that, for any $y\in X$ satisfying $\theta(y)\neq0$,
\begin{align*}
\rho(x,y)\leq C_\rho\max\left\{\rho(x,x_1),\,\rho(x_1,y)\right\}
<C_\rho A^3r.
\end{align*}
From this, \eqref{theta}, \eqref{CZk1}, and
both (i) and (ii) of Definition~\ref{testfunction},
we infer that
\begin{align*}
|Z_1(x)|
&\lesssim\int_{\rho(x,y)\leq\frac{r+\rho(x_1,x)}{2C_\rho}}
\frac{1}{V_\rho(x,y)}\left[\frac{\rho(x,y)}{r+\rho(x_1,x)}\right]^\beta
\\
&\quad\times\frac{1}{(V_\rho)_r(x_1)+V_\rho(x_1,x)}
\left[\frac{r}{r+\rho(x_1,x)}\right]^\gamma
\,d\mu(y)
\\
&\quad+\int_{\frac{r+\rho(x_1,x)}{2C_\rho}\leq\rho(x,y)
<C_\rho A^3r}\frac{1}{V_\rho(x,y)}
\left\{\frac{1}{(V_\rho)_r(x_1)+V_\rho(x_1,y)}
\left[\frac{r}{r+\rho(x_1,y)}\right]^\gamma
\right.
\\
&\quad\left.+\frac{1}{(V_\rho)_r(x_1)+V_\rho(x_1,x)}
\left[\frac{r}{r+\rho(x_1,x)}\right]^\gamma\right\}\,d\mu(y),
\end{align*}
which, together with Lemma~\ref{A3.2}(i)
and \eqref{upperdoub},
further implies that
\begin{align}\label{Z_1}
|Z_1(x)|
&\lesssim\frac{1}{(V_\rho)_r(x_1)}
\int_{\rho(x,y)\leq\frac{r+\rho(x_1,x)}{2C_\rho}}
\frac{1}{V_\rho(x,y)}\left[\frac{\rho(x,y)}{r+\rho(x_1,x)}\right]^\beta
\,d\mu(y)
\nonumber\\
&\quad+\frac{1}{(V_\rho)_r(x_1)}
\int_{\frac{r+\rho(x_1,x)}{2C_\rho}\leq\rho(x,y)
<C_\rho A^3r}\frac{1}{V_\rho(x,y)}\,d\mu(y)
\nonumber\\
&\lesssim\frac{1}{(V_\rho)_r(x_1)}\left[1+
\int_{\rho(x,y)<C_\rho A^3r}\frac{1}{(V_\rho)_\frac{r}{2C_\rho}(x)}
\,d\mu(y)\right]
\sim\frac{1}{(V_\rho)_r(x_1)}.
\end{align}

To deal with $Z_2(x)$, note that, for any $y\in X$ satisfying
$\psi(y)\neq0$, we have $y\in[B_\rho(x_1,A^2r)]^\complement$,
which, combined with Definition~\ref{D2.1}(iii)
and the fact that $\rho(x_1,x)<Ar<\frac{A^2r}{C_\rho}$, further implies that
\begin{align*}
\frac{A^2r}{C_\rho}\leq\frac{\rho(x_1,y)}{C_\rho}\leq
\max\left\{\rho(x_1,x),\,
\rho(x,y)\right\}=\rho(x,y).
\end{align*}
This, the definition of $\psi$,
the size condition of $f$, and (ii) and (vii)
of Lemma~\ref{A3.2} imply that
\begin{align}\label{Z_2}
|Z_2(x)|&\lesssim\int_{\rho(x,y)\ge\frac{A^2r}{C_\rho}}
\frac{1}{V_\rho(x,y)}|f(y)|\,d\mu(y)
\nonumber\\
&\lesssim\frac{1}{(V_\rho)_r(x)}\int_X
\frac{1}{(V_\rho)_r(x_1)+V_\rho(x_1,y)}
\left[\frac{r}{r+\rho(x_1,y)}\right]^\gamma\,d\mu(y)
\nonumber\\
&\lesssim\frac{1}{(V_\rho)_r(x)}
\sim\frac{1}{(V_\rho)_r(x_1)}.
\end{align}

To deal with $Z_3(x)$, by the size condition of $f$
and \eqref{Lem3.4}, we find that
\begin{align*}
|Z_3(x)|\lesssim\frac{1}{(V_\rho)_r(x_1)+V_\rho(x_1,x)}
\left[\frac{r}{r+\rho(x_1,x)}\right]^\gamma
|T\theta(x)|\lesssim\frac{1}{(V_\rho)_r(x_1)}.
\end{align*}
From this, \eqref{Tf=123}, \eqref{Z_1}, and \eqref{Z_2},
we deduce that
\begin{align*}
|Tf(x)|\leq|Z_1(x)|+|Z_2(x)|+|Z_3(x)|\lesssim\frac{1}{(V_\rho)_r(x_1)}.
\end{align*}
Thus, we finish the proof of \eqref{Tfs} in
the case that $\rho(x_1,x)<Ar$.

\emph{Case (2)} $\rho(x_1,x)\ge Ar$.
In this case, according to Corollary~\ref{Lem3.3},
we choose $u_1,u_2\in{C}^\varepsilon(X)$ such that
\begin{align}\label{u1}
\mathbf{1}_{B_\rho(x,A^{-3}\rho(x_1,x))}
\leq u_1\leq\mathbf{1}_{B_\rho(x,A^{-2}\rho(x_1,x))},
\end{align}
$$
\|u_1\|_{\dot{C}^\varepsilon(X)}\lesssim[\rho(x_1,x)]^{-\varepsilon},
$$
$$
\mathbf{1}_{B_\rho(x_1,A^{-3}\rho(x_1,x))}
\leq u_2\leq\mathbf{1}_{B_\rho(x_1,A^{-2}\rho(x_1,x))},
$$
and
$$
\|u_2\|_{\dot{C}^\varepsilon(X)}\lesssim[\rho(x_1,x)]^{-\varepsilon}.
$$
Let $u_3:=1-u_1-u_2$ and $f=fu_1+fu_2+fu_3=:f_1+f_2+f_3$.
To show \eqref{Tfs} holds in this case,
we claim that, for any $y,y'\in X$,
\begin{align}\label{f1size}
|f_1(y)|\lesssim\frac{1}{V_\rho(x_1,x)}
\left[\frac{r}{\rho(x_1,x)}\right]^\gamma,
\end{align}
\begin{align}\label{f1r}
|f_1(y)-f_1(y')|\lesssim\left[\frac{\rho(y,y')}
{\rho(x_1,x)}\right]^\beta\frac{1}{V_\rho(x_1,x)}
\left[\frac{r}{\rho(x_1,x)}\right]^\gamma,
\end{align}
\begin{align}\label{f3size}
|f_3(y)|\lesssim\frac{1}{(V_\rho)_r(x_1)+V_\rho(x_1,y)}
\left[\frac{r}{r+\rho(x_1,y)}\right]^\gamma
\mathbf{1}_{\{z\in X:\rho(x,z)\ge A^{-3}\rho(x_1,x)\}}(y),
\end{align}
\begin{align}\label{intf3}
\int_X|f_3(y)|\,d\mu(y)\lesssim\left[
\frac{r}{\rho(x_1,x)}\right]^\gamma,
\end{align}
and
\begin{align}\label{intf2}
\left|\int_Xf_2(y)\,d\mu(y)\right|
\lesssim\left[
\frac{r}{\rho(x_1,x)}\right]^\gamma.
\end{align}

We first prove \eqref{f1size}. Indeed,
if $y\in{\mathrm{supp\,}} u_1\subset B_\rho(x,A^{-2}\rho(x_1,x))$,
then, by Definition~\ref{D2.1}(iii),
we conclude that
\begin{align*}
Ar\leq\rho(x_1,x)\leq C_\rho
\max\left\{\rho(x_1,y),\,
\rho(y,x)\right\}=C_\rho\rho(x_1,y)
\end{align*}
and
\begin{align*}
\rho(x_1,y)\leq C_\rho\max\left\{\rho(x_1,x),\,\rho(x,y)\right\}
=C_\rho\rho(x_1,x),
\end{align*}
which, together with the assumption $\rho(x_1,x)\ge Ar$,
further implies that
\begin{align}\label{1533}
r+\rho(x_1,y)\sim\rho(x_1,x).
\end{align}
From this,
we infer that, for any $y\in X$,
\begin{align*}
|f_1(y)|
&\leq|f(y)|\mathbf{1}_{B_\rho(x,A^{-2}\rho(x_1,x))}(y)
\quad\text{by the definition of $f_1$ and \eqref{u1}}
\\
&\leq\frac{1}{(V_\rho)_r(x_1)+V_\rho(x_1,y)}
\left[\frac{r}{r+\rho(x_1,y)}\right]^\gamma
\mathbf{1}_{B_\rho(x,A^{-2}\rho(x_1,x))}(y)\\
&\qquad\ \text{by the size condition of $f$ [see Definition~\ref{testfunction}(i)]}\\
&\sim\frac{1}{(V_\rho)_{r+\rho(x_1,y)}(x_1)}
\left[\frac{r}{\rho(x_1,x)}\right]^\gamma
\mathbf{1}_{B_\rho(x,A^{-2}\rho(x_1,x))}(y)
\quad\text{by Lemma~\ref{A3.2}(vi)}\\
&\sim\frac{1}{V_\rho(x_1,x)}\left[\frac{r}{\rho(x_1,x)}\right]^\gamma,
\quad\text{by \eqref{1533}}
\end{align*}
which completes the proof of \eqref{f1size}.

Next, we show \eqref{f1r}. Note that
${\mathrm{supp\,}} u_1\subset B_\rho(x,A^{-2}\rho(x_1,x))$,
if $|f_1(y)-f_1(y')|\neq0$, then either
$y\in B_\rho(x,A^{-2}\rho(x_1,x))$
or $y'\in B_\rho(x,A^{-2}\rho(x_1,x))$.
Without loss of generality,
we may assume that $y\in B_\rho(x,A^{-2}\rho(x_1,x))$.
Then, by an argument similar to
that used in the proof of \eqref{f1size},
we find that \eqref{1533} still holds.
On the one hand, if $\rho(y,y')>\frac{r+\rho(x_1,y)}{2C_\rho}$,
then \eqref{f1size} and \eqref{1533} imply that
\begin{align}\label{1548}
|f_1(y)-f_1(y')|
&\leq|f_1(y)|+|f_1(y')|
\nonumber\\
&\lesssim\left[\frac{\rho(y,y')}{r+\rho(x_1,y)}\right]^\beta
\frac{1}{V_\rho(x_1,x)}\left[\frac{r}{\rho(x_1,x)}\right]^\gamma
\nonumber\\
&\sim\left[\frac{\rho(y,y')}{\rho(x_1,x)}\right]^\beta
\frac{1}{V_\rho(x_1,x)}\left[\frac{r}{\rho(x_1,x)}\right]^\gamma.
\end{align}
On the other hand, if
$\rho(y,y')\leq\frac{r+\rho(x_1,y)}{2C_\rho}$,
from both the size condition and the regularity condition of $f$,
we infer that
\begin{align*}
|f_1(y)-f_1(y')|
&=|f(y)u_1(y)-f(y)u_1(y')
+f(y)u_1(y')-f(y')u_1(y')|
\\
&\leq|f(y)||u_1(y)-u_1(y')|
+|u_1(y')||f(y)-f(y')|
\\
&\lesssim\frac{1}{(V_\rho)_r(x_1)+V_\rho(x_1,y)}
\left[\frac{r}{r+\rho(x_1,y)}\right]^\gamma
\\
&\quad\times\min\left\{\sup_{y\in X}u(y)-\inf_{y\in X}u(y),\,
\left[\rho(y,y')\right]^\varepsilon
\|u_1\|_{\dot{C}^\varepsilon(X)}\right\}
\\
&\quad+\left[\frac{\rho(y,y')}{r+\rho(x_1,y)}\right]^\beta
\frac{1}{(V_\rho)_r(x_1)+V_\rho(x_1,y)}
\left[\frac{r}{r+\rho(x_1,y)}\right]^\gamma,
\end{align*}
which, combined with $\beta\in(0,\varepsilon)$
and the fact that $r+\rho(x_1,y)\sim\rho(x_1,x)$,
implies that
\begin{align*}
|f_1(y)-f_1(y')|
&\lesssim\frac{1}{(V_\rho)_r(x_1)+V_\rho(x_1,y)}
\left[\frac{r}{r+\rho(x_1,y)}\right]^\gamma
\min\left\{1,\,
\left[\frac{\rho(y,y')}{\rho(x_1,x)}\right]^\varepsilon\right\}
\\
&\quad+\left[\frac{\rho(y,y')}{r+\rho(x_1,y)}\right]^\beta
\frac{1}{(V_\rho)_r(x_1)+V_\rho(x_1,y)}
\left[\frac{r}{r+\rho(x_1,y)}\right]^\gamma
\\
&\lesssim\left[\frac{\rho(y,y')}{\rho(x_1,x)}\right]^\beta
\frac{1}{V_\rho(x_1,x)}\left[\frac{r}{\rho(x_1,x)}\right]^\gamma.
\end{align*}
This, together with \eqref{1548}, finishes the proof of \eqref{f1r}.

Now, we prove \eqref{f3size}.
By $B_\rho(x,A^{-2}\rho(x_1,x))
\cap B_\rho(x_1,A^{-2}\rho(x_1,x))=\varnothing$,
the definitions of $u_1$,
$u_2$, and $u_3$, we find that $$0\leq u_3\leq1$$
and, for any
$y\in B_\rho(x,A^{-3}\rho(x_1,x))\cup B_\rho(x_1,A^{-3}\rho(x_1,x))$,
$u_3(y)=0$.
From this and the size condition of $f$, we deduce that,
for any $y\in X$,
\begin{align*}
|f_3(y)|
\lesssim\frac{1}{(V_\rho)_r(x_1)+V_\rho(x_1,y)}
\left[\frac{r}{r+\rho(x_1,y)}\right]^\gamma
\mathbf{1}_{[B_\rho(x,A^{-3}\rho(x_1,x))]^\complement}(y),
\end{align*}

Next, we show \eqref{intf3}.
By \eqref{f3size} and Lemma~\ref{A3.2}(v), we obtain
\begin{align*}
\int_X|f_3(y)|\,d\mu(y)
&\lesssim\int_{\rho(x_1,y)\ge A^{-3}\rho(x_1,x)}
\frac{1}{(V_\rho)_r(x_1)+V_\rho(x_1,y)}
\left[\frac{r}{r+\rho(x_1,y)}\right]^\gamma\,d\mu(y)
\\
&\lesssim\left[\frac{r}{r+\rho(x_1,x)}\right]^\gamma
\leq\left[\frac{r}{\rho(x_1,x)}\right]^\gamma,
\end{align*}
which is precisely \eqref{intf3}.

Finally, we prove \eqref{intf2}.
Since $f\in\mathring{\mathcal{G}}(x_1,r,\beta,\gamma)$, we infer that
$\int_Xf(y)\,d\mu(y)=0$, which, combined with
the facts that $f=f_1+f_2+f_3$ and
${\mathrm{supp\,}} f_1\subset{\mathrm{supp\,}}
u_1\subset B_\rho(x,A^{-2}\rho(x_1,x))$,
further implies that
\begin{align*}
\left|\int_Xf_2(y)\,d\mu(y)\right|
&=\left|\int_X\left[f_1(y)+f_3(y)\right]\,d\mu(y)\right|
\\
&\leq\int_{B_\rho(x,A^{-2}\rho(x_1,x))}
|f_1(y)|\,d\mu(y)+\int_X|f_3(y)|\,d\mu(y).
\end{align*}
By this, \eqref{f1size}, and \eqref{intf3},
we conclude that
\begin{align*}
\left|\int_Xf_2(y)\,d\mu(y)\right|
&\lesssim\int_{B_\rho(x,A^{-2}\rho(x_1,x))}
\frac{1}{V_\rho(x_1,x)}\left[\frac{r}{\rho(x_1,x)}\right]^\gamma
\,d\mu(y)+\left[\frac{r}{\rho(x_1,x)}\right]^\gamma
\\
&=\left[\frac{r}{\rho(x_1,x)}\right]^\gamma
\left[\frac{(V_\rho)_{A^{-2}\rho(x_1,x)}(x_1)}{V_\rho(x_1,x)}+1\right]
\sim\left[\frac{r}{\rho(x_1,x)}\right]^\gamma,
\end{align*}
which completes the proof of \eqref{intf2}
and hence the above claim.

Now, we show \eqref{Tfs} in the case $\rho(x_1,x)\ge Ar$.
From $\rho(x_1,x)\ge Ar$,
it follows that
\begin{align*}
\frac{r}{r+\rho(x_1,x)}\sim\frac{r}{\rho(x_1,x)}
\end{align*}
and
\begin{align*}
\frac{1}{(V_\rho)_r(x_1)+V_\rho(x_1,x)}\sim\frac{1}{V_\rho(x_1,x)}.
\end{align*}
Thus, to prove \eqref{Tfs}, it suffices to show that
\begin{align}\label{1635}
|Tf(x)|\lesssim\frac{1}{V_\rho(x_1,x)}
\left[\frac{r}{\rho(x_1,x)}\right]^\gamma.
\end{align}
With this goal in mind,
we first use the fact that $f=f_1+f_2+f_3$ to write
\begin{align}\label{Tf}
|Tf(x)|\leq|Tf_1(x)|+|Tf_2(x)|+|Tf_3(x)|.
\end{align}

To deal with $Tf_1(x)$, by Corollary~\ref{Lem3.3}, we find that
there exists a function $u\in{C}^\varepsilon(X)$ such that
$$
\mathbf{1}_{B_\rho(x,A^{-2}\rho(x_1,x))}\leq u
\leq\mathbf{1}_{B_\rho(x,A^{-1}\rho(x_1,x))}
$$
and $\|u\|_{\dot{C}^\varepsilon(X)}
\lesssim[\rho(x_1,x)]^{-\varepsilon}$.
Then, for any $y\in X$, $f_1(y)=u(y)f_1(y)$
and
\begin{align}\label{Tf1}
Tf_1(x)
&=\int_XK(x,y)u(y)[f_1(y)-f_1(x)]\,d\mu(y)
+f_1(x)\int_XK(x,y)u(y)\,d\mu(y)
\nonumber\\
&=:I_1(x)+I_2(x).
\end{align}
From the fact that ${\mathrm{supp\,}} u\subset B_\rho(x,A^{-1}\rho(x_1,x))$,
we infer that
\begin{align*}
|I_1(x)|\leq\int_{B_\rho(x,A^{-1}\rho(x_1,x))}
|K(x,y)||f_1(y)-f_1(x)|\,d\mu(y),
\end{align*}
which, together with \eqref{CZk1}, \eqref{f1r}, and Lemma~\ref{A3.2}(i),
further implies that
\begin{align}\label{Tf1I1}
|I_1(x)|
&\lesssim\frac{1}{V_\rho(x_1,x)}\left[\frac{r}{\rho(x_1,x)}\right]^\gamma
\int_{B_\rho(x,A^{-1}\rho(x_1,x))}\frac{1}{V_\rho(x,y)}
\left[\frac{\rho(x,y)}{\rho(x_1,x)}\right]^\beta\,d\mu(y)
\nonumber\\
&\lesssim\frac{1}{V_\rho(x_1,x)}\left[\frac{r}{\rho(x_1,x)}\right]^\gamma.
\end{align}
By \eqref{f1size}, the definition of $u$, and \eqref{Lem3.4},
we conclude that
\begin{align*}
|I_2(x)|
\lesssim|f_1(x)||Tu(x)|
\lesssim\frac{1}{V_\rho(x_1,x)}
\left[\frac{r}{\rho(x_1,x)}\right]^\gamma,
\end{align*}
which, combined with both \eqref{Tf1} and \eqref{Tf1I1},
further implies that
\begin{align}
\label{1641}
\left|Tf_1(x)\right|\lesssim\frac{1}{V_\rho(x_1,x)}
\left[\frac{r}{\rho(x_1,x)}\right]^\gamma.
\end{align}

To deal with $Tf_2(x)$, note that
$$
{\mathrm{supp\,}} f_2\subset{\mathrm{supp\,}}
u_2\subset B_\rho\left(x_1,A^{-2}\rho(x_1,x)\right)
$$
and
\begin{align}\label{1551}
\rho(x_1,y)<\frac{\rho(x_1,x)}{A^2}
\leq\frac{\rho(x_1,x)}{2C_\rho}.
\end{align}
From this and the fact that $|f_2|\leq|f|$ pointwise on $X$, it follows that
\begin{align*}
|Tf_2(x)|
&=\left|\int_X\left[K(x,y)-K(x,x_1)\right]
f_2(y)\,d\mu(y)+\int_XK(x,x_1)f_2(y)\,d\mu(y)\right|\\
&\leq\int_{\rho(x_1,y)<A^{-2}\rho(x_1,x)}
\left|K(x,y)-K(x,x_1)\right||f(y)|\,d\mu(y)\\
&\quad+|K(x,x_1)|\left|\int_Xf(y)\,d\mu(y)\right|.
\end{align*}
By this, \eqref{CZk2} [here, we need \eqref{1551}],
the size condition of $f$,
\eqref{CZk1}, and \eqref{intf2},
we find that
\begin{align*}
|Tf_2(x)|
&\lesssim\int_{\rho(x_1,y)<A^{-2}\rho(x_1,x)}
\left[\frac{\rho(x_1,y)}{\rho(x_1,x)}\right]^\varepsilon
\frac{1}{V_\rho(x_1,x)}
\frac{1}{(V_\rho)_r(x_1)+V_\rho(x_1,y)}\\
&\quad\times\left[\frac{r}{r+\rho(x_1,y)}\right]^\gamma\,d\mu(y)
+\frac{1}{V_\rho(x_1,x)}\left[\frac{r}{\rho(x_1,x)}\right]^\gamma,
\end{align*}
which, together with $\gamma\in(0,\varepsilon)$
and Lemma~\ref{A3.2}(i),
further implies that
\begin{align}\label{Tf2}
|Tf_2(x)|
&\lesssim\frac{1}{V_\rho(x_1,x)}\int_{\rho(x_1,y)<A^{-2}\rho(x_1,x)}
\left[\frac{\rho(x_1,y)}{\rho(x_1,x)}\right]^{\varepsilon-\gamma}
\frac{1}{V_\rho(x_1,y)}\nonumber\\
&\quad\times\left[\frac{\rho(x_1,y)}{\rho(x_1,x)}
\frac{r}{r+\rho(x_1,y)}\right]^\gamma\,d\mu(y)
+\frac{1}{V_\rho(x_1,x)}\left[\frac{r}{\rho(x_1,x)}\right]^\gamma
\nonumber\\
&\leq\frac{1}{V_\rho(x_1,x)}\left[\frac{r}{\rho(x_1,x)}\right]^\gamma
\left\{\int_{\rho(x_1,y)<A^{-2}\rho(x_1,x)}
\left[\frac{\rho(x_1,y)}{\rho(x_1,x)}\right]^{\varepsilon-\gamma}
\frac{d\mu(y)}{V_\rho(x_1,y)}+1\right\}
\nonumber\\
&\sim\frac{1}{V_\rho(x_1,x)}\left[\frac{r}{\rho(x_1,x)}\right]^\gamma.
\end{align}

To deal with $Tf_3(x)$,
we obtain
\begin{align*}
|Tf_3(x)|
&\leq\int_{\rho(x,y)\ge A^{-3}\rho(x_1,x)}
|K(x,y)||f_3(y)|\,d\mu(y)
\quad\text{by \eqref{f3size}}
\nonumber\\
&\lesssim\int_{\rho(x,y)\ge A^{-3}\rho(x_1,x)}
\frac{1}{V_\rho(x,y)}|f_3(y)|\,d\mu(y)
\quad\text{by \eqref{CZk1}}
\nonumber\\
&\leq\frac{1}{(V_\rho)_{A^{-3}\rho(x_1,x)}(x)}
\int_{\rho(x,y)\ge A^{-3}\rho(x_1,x)}
|f_3(y)|\,d\mu(y)
\nonumber\\
&\lesssim\frac{1}{V_\rho(x,x_1)}\int_{X}
|f_3(y)|\,d\mu(y)
\quad\text{by \eqref{upperdoub}}
\nonumber\\
&\lesssim\frac{1}{V_\rho(x_1,x)}\left[\frac{r}{\rho(x_1,x)}\right]^\gamma
\quad\text{by \eqref{intf3}},
\end{align*}
which, combined with \eqref{1641}, \eqref{Tf}, and
\eqref{Tf2}, further implies \eqref{1635} and hence
\eqref{Tfs} in the case $\rho(x_1,x)\ge Ar$.
This finishes the proof of \eqref{Tfs}.

\emph{Step 2.}
In this step, we prove that,
for any $x,x'\in X$ satisfying
$\rho(x,x')\leq\frac{r+\rho(x_1,x)}{2C_\rho}$,
\begin{align}
\label{Tfr}
|Tf(x)-Tf(x')|\lesssim\left[\frac{\rho(x,x')}
{r+\rho(x_1,x)}\right]^\beta
\frac{1}{(V_\rho)_r(x_1)+V_\rho(x_1,x)}
\left[\frac{r}{r+\rho(x_1,x)}\right]^\gamma.
\end{align}
Assume that $x,x'\in X$
satisfy $\rho(x,x')\leq\frac{r+\rho(x_1,x)}{2C_\rho}$.
We show \eqref{Tfr} by considering the
following three cases.

\emph{Case (1)} $\rho(x,x')\ge\frac{r+\rho(x_1,x)}{A^5}$.
In this case, we have
\begin{align}\label{1943}
r+\rho(x_1,x)\sim\rho(x,x')\sim r+\rho(x_1,x'),
\end{align}
which, together with Lemma~\ref{A3.2}(vi), further implies that
\begin{align}\label{1544}
(V_\rho)_r(x_1)+V_\rho(x_1,x)
&\sim(V_\rho)_{r+\rho(x_1,x)}(x_1)
\sim(V_\rho)_{r+\rho(x_1,x')}(x_1)\nonumber\\
&\sim(V_\rho)_r(x_1)+V_\rho(x_1,x').
\end{align}
From this, we deduce that
\begin{align*}
|Tf(x)-Tf(x')|
&\leq|Tf(x)|+|Tf(x')|\\
&\lesssim\frac{1}{(V_\rho)_r(x_1)+V_\rho(x_1,x)}
\left[\frac{r}{r+\rho(x_1,x)}\right]^\gamma\\
&\quad
+\frac{1}{(V_\rho)_r(x_1)+V_\rho(x_1,x')}
\left[\frac{r}{r+\rho(x_1,x')}\right]^\gamma
\quad\text{by \eqref{Tfs}}\\
&\sim\left[\frac{\rho(x,x')}{r+\rho(x_1,x)}\right]^\beta
\frac{1}{(V_\rho)_r(x_1)+V_\rho(x_1,x)}
\left[\frac{r}{r+\rho(x_1,x)}\right]^\gamma\\
&\qquad\ \text{by \eqref{1544} and \eqref{1943}},
\end{align*}
which completes the proof of \eqref{Tfr}
in the case $\rho(x,x')\ge\frac{r+\rho(x_1,x)}{A^5}$.

\emph{Case (2)} $\rho(x,x')<\frac{r+\rho(x_1,x)}{A^5}$
and $\rho(x_1,x)<Ar$.
In this case, we have
$$
\left[\frac{\rho(x,x')}{r+\rho(x_1,x)}\right]^\beta
\frac{1}{(V_\rho)_r(x_1)+V_\rho(x_1,x)}
\left[\frac{r}{r+\rho(x_1,x)}\right]^\gamma
\sim\left[\frac{\rho(x,x')}{r}\right]^\beta
\frac{1}{(V_\rho)_r(x_1)}.
$$
Thus, to prove \eqref{Tfr},
it suffices to show that
\begin{align}\label{2150}
|Tf(x)-Tf(x')|\lesssim
\left[\frac{\rho(x,x')}{r}\right]^\beta
\frac{1}{(V_\rho)_r(x_1)}.
\end{align}
By Corollary~\ref{Lem3.3},
choose a function $\phi_1\in{C}^\varepsilon(X)$
such that
$$
\mathbf{1}_{B_\rho(x,A^2\rho(x,x'))}
\leq\phi_1
\leq\mathbf{1}_{B_\rho(x,A^3\rho(x,x'))}
$$
and $\|\phi_1\|_{\dot{C}^\varepsilon(X)}
\lesssim[\rho(x,x')]^{-\varepsilon}$.
Let $\phi_2:=1-\phi_1$.
Since $f\in\mathring{\mathcal{G}}(x_1,r,\beta,\gamma)\subset{C}^\beta(X)$,
it follows from (ii) of the present theorem that
\begin{align*}
Tf(x)&=\int_XK(x,y)[f(y)-f(x)]\phi_1(y)\,d\mu(y)
\nonumber\\
&\quad+\left[\int_XK(x,y)f(y)\phi_2(y)\,d\mu(y)
+f(x)\int_XK(x,y)\phi_1(y)\,d\mu(y)\right]
\nonumber\\
&=:Y_1(x)+Y_2(x).
\end{align*}
Since we have
\begin{align}
\label{12.26x1-2}
|Tf(x)-Tf(x')|\leq|Y_1(x)-Y_1(x')|+|Y_2(x)-Y_2(x')|,
\end{align}
we will estimate the two terms on the right-hand side of \eqref{12.26x1-2}.
To estimate the first term
$|Y_1(x)-Y_1(x')|$, noting that
${\mathrm{supp\,}}\phi_1\subset B_\rho(x,A^3\rho(x,x'))$,
we obtain
\begin{align*}
|Y_1(x)|
&\leq\int_{\rho(x,y)<A^3\rho(x,x'),\,
\rho(x,y)\leq\frac{r+\rho(x_1,x)}{2C_\rho}}|K(x,y)||f(y)-f(x)|\,d\mu(y)
\\
&\quad+\int_{\rho(x,y)<A^3\rho(x,x'),\,
\rho(x,y)>\frac{r+\rho(x_1,x)}{2C_\rho}}|K(x,y)|\left[|f(y)|+|f(x)|\right]
\,d\mu(y)
\\
&=:Y_{1,1}(x)+Y_{1,2}(x).
\end{align*}

For $Y_{1,1}(x)$,
by \eqref{CZk1}, the regularity condition of $f$,
and Lemma~\ref{A3.2}(i), we conclude that
\begin{align}\label{Y11}
Y_{1,1}(x)
&\lesssim\int_{\rho(x,y)<A^3\rho(x,x')}
\frac{1}{V_\rho(x,y)}\left[\frac{\rho(x,y)}
{r+\rho(x_1,x)}\right]^\beta
\nonumber\\
&\quad\times\frac{1}{(V_\rho)_r(x_1)+V_\rho(x_1,x)}
\left[\frac{r}{r+\rho(x_1,x)}\right]^\gamma\,d\mu(y)
\nonumber\\
&\leq\left[\frac{\rho(x,x')}{r}\right]^\beta
\frac{1}{(V_\rho)_r(x_1)}
\int_{\rho(x,y)<A^3\rho(x,x')}
\frac{1}{V_\rho(x,y)}
\left[\frac{\rho(x,y)}{\rho(x,x')}\right]^\beta\,d\mu(y)
\nonumber\\
&\lesssim
\left[\frac{\rho(x,x')}{r}\right]^\beta
\frac{1}{(V_\rho)_r(x_1)}.
\end{align}

For $Y_{1,2}(x)$, from the fact that
$\rho(x,y)<A^3\rho(x,x')<\frac{r+\rho(x_1,x)}{A^2}$,
we infer that
$$
\left\{y\in X:\frac{r+\rho(x_1,x)}{2C_\rho}<\rho(x,y)<A^3\rho(x,x')\right\}
=\varnothing
$$
and hence
$Y_{1,2}(x)=0$.
This, combined with \eqref{Y11}, implies that
\begin{align}\label{Y1}
|Y_1(x)|\lesssim\left[\frac{\rho(x,x')}{r}\right]^\beta
\frac{1}{(V_\rho)_r(x_1)}.
\end{align}
Moreover, by Definition~\ref{D2.1}(iii), we find that,
for any $y\in{\mathrm{supp\,}}\phi_1\subset B_\rho(x,A^3\rho(x,x'))$,
$$
\rho(x',y)\leq C_\rho\max\{\rho(x',x),\,\rho(x,y)\}<C_\rho A^3\rho(x,x').
$$
From this and an argument similar to that used in the estimation
of \eqref{Y1}, we deduce that
$$
|Y_1(x')|\lesssim\left[\frac{\rho(x,x')}{r}\right]^\beta
\frac{1}{(V_\rho)_r(x_1)},
$$
which, together with \eqref{Y1}, implies that
\begin{align}\label{2149}
|Y_1(x)-Y_1(x')|
\leq|Y_1(x)|+|Y_1(x')|
\lesssim\left[\frac{\rho(x,x')}{r}\right]^\beta
\frac{1}{(V_\rho)_r(x_1)}.
\end{align}

To estimate $|Y_2(x)-Y_2(x')|$, using (iii),
we write
\begin{align}\label{Y2}
Y_2(x)-Y_2(x')
&=\int_X\left[K(x,y)-K(x',y)\right]f(y)\phi_2(y)\,d\mu(y)
\nonumber\\
&\quad+f(x)\int_XK(x,y)\phi_1(y)\,d\mu(y)
-f(x')\int_XK(x',y)\phi_1(y)\,d\mu(y)
\nonumber\\
&=\int_X\left[K(x,y)-K(x',y)\right]
\left[f(y)-f(x)\right]\phi_2(y)\,d\mu(y)
\nonumber\\
&\quad+\left[f(x)-f(x')\right]\int_XK(x',y)\phi_1(y)\,d\mu(y)
\nonumber\\
&=:Y_{2,1}(x,x')+Y_{2,2}(x,x').
\end{align}

For $Y_{2,1}(x,x')$, by the definitions of $\phi_1$ and $\phi_2$,
we conclude that $\phi_2$ is bounded by 1 and ${\mathrm{supp\,}}\phi_2
\subset [B_\rho(x,A^2\rho(x,x'))]^\complement$.
Therefore, we can write
\begin{align}
\label{Y21}
&|Y_{2,1}(x,x')|\nonumber\\
&\quad\leq\int_{A^2\rho(x,x')\leq
\rho(x,y)<\frac{r+\rho(x_1,x)}{2C_\rho}}
|K(x,y)-K(x',y)||f(y)-f(x)|\,d\mu(y)
\nonumber\\
&\qquad+\int_{\rho(x,y)\ge\max\{A^2\rho(x,x'),\,
\frac{r+\rho(x_1,x)}{2C_\rho}\}}
|K(x,y)-K(x',y)|\left[|f(y)|+|f(x)|\right]\,d\mu(y)
\nonumber\\
&\quad=:Y_{2,1,1}(x,x')+Y_{2,1,2}(x,x').
\end{align}
To deal with $Y_{2,1,1}(x,x')$, from the fact that
$\rho(x,x')\leq\frac{\rho(x,y)}{A^2}\leq\frac{\rho(x,y)}{2C_\rho}$,
\eqref{CZk2}, and
the regularity condition of $f$, it follows that
\begin{align*}
Y_{2,1,1}(x,x')
&\lesssim\int_{A^2\rho(x,x')\leq
\rho(x,y)<\frac{r+\rho(x_1,x)}{2C_\rho}}
\left[\frac{\rho(x,x')}{\rho(x,y)}\right]^\varepsilon\frac{1}{V_\rho(x,y)}
\nonumber\\
&\quad\times\left[\frac{\rho(x,y)}{r+\rho(x_1,x)}\right]^\beta
\frac{1}{(V_\rho)_r(x_1)+V_\rho(x_1,x)}
\left[\frac{r}{r+\rho(x_1,x)}\right]^\gamma
\,d\mu(y)
\nonumber\\
&\leq\frac{1}{(V_\rho)_r(x_1)}\int_{\rho(x,y)\ge A^2\rho(x,x')}
\left[\frac{\rho(x,x')}{\rho(x,y)}\right]^\varepsilon\frac{1}{V_\rho(x,y)}
\left[\frac{\rho(x,y)}{r}\right]^\beta\,d\mu(y),
\end{align*}
which, combined with $\beta\in(0,\varepsilon)$
and Lemma~\ref{A3.2}(i),
further implies that
\begin{align}\label{Y211}
Y_{2,1,1}(x,x')
&\lesssim
\left[\frac{\rho(x,x')}{r}\right]^\beta
\frac{1}{(V_\rho)_r(x_1)}
\int_{\rho(x,y)\ge A^2\rho(x,x')}
\left[\frac{\rho(x,x')}{\rho(x,y)}\right]^{\varepsilon-\beta}\frac{1}{V_\rho(x,y)}
\,d\mu(y)
\nonumber\\
&\lesssim
\left[\frac{\rho(x,x')}{r}\right]^\beta
\frac{1}{(V_\rho)_r(x_1)}.
\end{align}
To deal with $Y_{2,1,2}(x,x')$, by the fact that
$\rho(x,x')\leq\frac{\rho(x,y)}{A^2}
\leq\frac{\rho(x,y)}{2C_\rho}$,
\eqref{CZk2}, and
the size condition of $f$, we obtain
\begin{align*}
Y_{2,1,2}(x,x')
&\lesssim\int_{\rho(x,y)\ge\max\{A^2\rho(x,x'),\,
\frac{r+\rho(x_1,x)}{2C_\rho}\}}
\left[\frac{\rho(x,x')}{\rho(x,y)}\right]^\varepsilon
\frac{1}{V_\rho(x,y)}
\nonumber\\
&\quad\times\Bigg\{\frac{1}{(V_\rho)_r(x_1)+V_\rho(x_1,y)}
\left[\frac{r}{r+\rho(x_1,y)}\right]^\gamma
\nonumber\\
&\quad+\frac{1}{(V_\rho)_r(x_1)+V_\rho(x_1,x)}
\left[\frac{r}{r+\rho(x_1,x)}\right]^\gamma\Bigg\}
\,d\mu(y),
\end{align*}
which, together with $\beta\in(0,\varepsilon)$
and Lemma~\ref{A3.2}(i),
further implies that
\begin{align}\label{Y212}
Y_{2,1,2}(x,x')
&\lesssim\frac{1}{(V_\rho)_r(x_1)}
\int_{\rho(x,y)\ge\max\{A^2\rho(x,x'),\,
\frac{r}{2C_\rho}\}}
\left[\frac{\rho(x,x')}{\rho(x,y)}\right]^{\varepsilon-\beta}
\left[\frac{\rho(x,x')}{\rho(x,y)}\right]^{\beta}
\frac{1}{V_\rho(x,y)}\,d\mu(y)
\nonumber\\
&\lesssim
\left[\frac{\rho(x,x')}{r}\right]^\beta
\frac{1}{(V_\rho)_r(x_1)}
\int_{\rho(x,y)\ge A^2\rho(x,x')}
\left[\frac{\rho(x,x')}{\rho(x,y)}\right]^{\varepsilon-\beta}
\frac{1}{V_\rho(x,y)}\,d\mu(y)\nonumber\\
&\lesssim
\left[\frac{\rho(x,x')}{r}\right]^\beta
\frac{1}{(V_\rho)_r(x_1)}.
\end{align}

For $Y_{2,2}(x,x')$, from \eqref{Lem3.4},
the regularity condition of $f$, and $\rho(x_1,x)<Ar$,
we infer that
\begin{align*}
|Y_{2,2}(x,x')|
&\lesssim\left|f(x)-f(x')\right||T\phi_1(x')|
\\
&\lesssim\left[\frac{\rho(x,x')}{r+\rho(x_1,x)}\right]^\beta
\frac{1}{(V_\rho)_r(x_1)+V_\rho(x_1,x)}
\left[\frac{r}{r+\rho(x_1,x)}\right]^\gamma
\\
&\leq\left[\frac{\rho(x,x')}{r}\right]^\beta
\frac{1}{(V_\rho)_r(x_1)}.
\end{align*}
By this, \eqref{Y2}, \eqref{Y21}, \eqref{Y211}, and \eqref{Y212},
we conclude that
\begin{align*}
|Y_2(x)-Y_2(x')|
&\leq|Y_{2,1}(x,x')|+|Y_{2,2}(x,x')|
\\
&\leq Y_{2,1,1}(x,x')+Y_{2,1,2}(x,x')+|Y_{2,2}(x,x')|
\lesssim\left[\frac{\rho(x,x')}{r}\right]^\beta\frac{1}{(V_\rho)_r(x_1)},
\end{align*}
which, combined with both \eqref{2149} and \eqref{12.26x1-2},
further implies that \eqref{2150} holds and hence
\eqref{Tfr} holds in the case both
$\rho(x,x')<\frac{r+\rho(x_1,x)}{A^5}$ and $\rho(x_1,x)<Ar$.

\emph{Case (3)} $\rho(x,x')<\frac{r+\rho(x_1,x)}{A^5}$
and $\rho(x_1,x)\ge Ar$. In this case, we have
\begin{align*}
r+\rho(x_1,x)\sim\rho(x_1,x)
\ \ \text{and}\ \
(V_\rho)_r(x_1)+V_\rho(x_1,x)\sim V_\rho(x_1,x).
\end{align*}
From this, we deduce that, to prove \eqref{Tfr},
it suffices to show that
\begin{align}\label{TfrCase3}
|Tf(x)-Tf(x')|
\lesssim
\left[\frac{\rho(x,x')}{\rho(x_1,x)}\right]^\beta
\frac{1}{V_\rho(x_1,x)}\left[\frac{r}{\rho(x_1,x)}\right]^\gamma.
\end{align}
For any $i\in\{1,\,2,\,3\}$, let $u_i$ and $f_i$
be defined as in Step 1 of the present proof. Then
\begin{align}\label{qq-pt1}
|Tf(x)-Tf(x')|\leq\sum_{i=1}^3|Tf_i(x)-Tf_i(x')|,
\end{align}
and hence we can estimate each $|Tf_i(x)-Tf_i(x')|$ separately.
To this end, we first estimate $|Tf_1(x)-Tf_1(x')|$.
By Corollary~\ref{Lem3.3}, we find that
there exists a function $u_4\in{C}^\varepsilon(X)$
such that
\begin{align*}
\mathbf{1}_{B_\rho(x,A\rho(x,x'))}\leq
u_4\leq\mathbf{1}_{B_\rho(x,A^2\rho(x,x'))}
\end{align*}
and $\|u_4\|_{\dot{C}^\varepsilon(X)}
\lesssim[\rho(x,x')]^{-\varepsilon}$.
Let $v_4:=1-u_4$.
Using this, we can write
\begin{align}\label{Y3+Y4}
Tf_1(x)
&=\int_XK(x,y)u_4(y)\left[f_1(y)-f_1(x)\right]\,d\mu(y)
\nonumber\\
&\quad+\left[\int_XK(x,y)v_4(y)f_1(y)\,d\mu(y)
+f_1(x)\int_X
K(x,y)u_4(y)\,d\mu(y)\right]
\nonumber\\
&=:Y_3(x)+Y_4(x).
\end{align}

To deal with $Y_3(x)-Y_3(x')$, we first claim that
\begin{align}
\label{1536}
\frac{\rho(x_1,x)}{C_\rho}\leq\rho(x_1,x')\leq C_\rho\rho(x_1,x).
\end{align}
Indeed, from Definition~\ref{D2.1}(iii) and the facts that
\begin{align*}
\rho(x,x')<\frac{r+\rho(x_1,x)}{A^5}
\ \ \text{and}\ \
\rho(x_1,x)\ge Ar,
\end{align*}
we infer that
\begin{align*}
\rho(x_1,x)\leq C_\rho\max\left\{\rho(x_1,x'),\,\rho(x',x)\right\}
=C_\rho\rho(x_1,x')
\end{align*}
and
\begin{align*}
\rho(x_1,x')\leq C_\rho\max\left\{\rho(x_1,x),\,\rho(x,x')\right\}
=C_\rho\rho(x_1,x),
\end{align*}
which completes the proof of \eqref{1536}.
Note that, by $\rho(x,x')<\frac{r+\rho(x_1,x)}{A^5}$,
the choice $A=2C_\rho$, and \eqref{1536},
we conclude that,
for any $y\in B_\rho(x,A^2\rho(x,x'))$,
\begin{align*}
\rho(x,y)<A^2\rho(x,x')<\frac{r+\rho(x_1,x)}{A^3}
<\frac{r+\rho(x_1,x)}{2C_\rho}
\end{align*}
and, for any $y\in B_\rho(x',A^3\rho(x,x'))$,
\begin{align*}
\rho(x',y)<A^3\rho(x,x')<\frac{r+\rho(x_1,x)}{A^2}
\leq\frac{r+C_\rho\rho(x_1,x')}{A^2}
\le\frac{r+\rho(x_1,x')}{2C_\rho}.
\end{align*}
From these, ${\mathrm{supp\,}} u_4\subset B_\rho(x,A^2\rho(x,x'))$,
the fact that
$$
B_\rho\left(x,A^2\rho(x,x')\right)
\subset B_\rho\left(x',A^3\rho(x,x')\right),
$$
\eqref{CZk1}, and \eqref{f1r},
it follows that
\begin{align*}
|Y_3(x)-Y_3(x')|
&\leq\int_{B_\rho(x,A^2\rho(x,x'))}|K(x,y)|
\left|f_1(y)-f_1(x)\right|\,d\mu(y)\nonumber\\
&\quad+\int_{B_\rho(x,A^2\rho(x,x'))}|K(x',y)|
\left|f_1(y)-f_1(x')\right|\,d\mu(y)\nonumber\\
&\lesssim\int_{B_\rho(x,A^2\rho(x,x'))}\frac{1}{V_\rho(x,y)}
\left[\frac{\rho(x,y)}{\rho(x_1,x)}\right]^\beta
\frac{1}{V_\rho(x_1,x)}
\left[\frac{r}{\rho(x_1,x)}\right]^\gamma\,d\mu(y)
\nonumber\\
&\quad+\int_{B_\rho(x',A^3\rho(x,x'))}\frac{1}{V_\rho(x',y)}
\left[\frac{\rho(x',y)}{\rho(x_1,x')}\right]^\beta
\frac{1}{V_\rho(x_1,x')}
\left[\frac{r}{\rho(x_1,x')}\right]^\gamma\,d\mu(y),
\end{align*}
which, together with \eqref{1536}
and Lemma~\ref{A3.2}(i),
further implies that
\begin{align}\label{Y3(x-x')}
|Y_3(x)-Y_3(x')|
&\lesssim\left[\frac{\rho(x,x')}{\rho(x_1,x)}\right]^\beta
\frac{1}{V_\rho(x_1,x)}
\left[\frac{r}{\rho(x_1,x)}\right]^\gamma\nonumber\\
&\quad\times\Bigg\{\int_{B_\rho(x,A^2\rho(x,x'))}\frac{1}{V_\rho(x,y)}
\left[\frac{\rho(x,y)}{\rho(x,x')}\right]^\beta
\,d\mu(y)
\nonumber\\
&\quad+\int_{B_\rho(x',A^3\rho(x,x'))}\frac{1}{V_\rho(x',y)}
\left[\frac{\rho(x',y)}{\rho(x,x')}\right]^\beta
\,d\mu(y)\Bigg\}
\nonumber\\
&\lesssim\left[\frac{\rho(x,x')}{\rho(x_1,x)}\right]^\beta
\frac{1}{V_\rho(x_1,x)}
\left[\frac{r}{\rho(x_1,x)}\right]^\gamma.
\end{align}

For $Y_4$,
by \eqref{CZO2} and \eqref{Y3+Y4},
we write $Y_4(x)-Y_4(x')$ as
\begin{align}\label{Y4(x-x')}
Y_4(x)-Y_4(x')
&=\int_X\left[K(x,y)-K(x',y)\right]v_4(y)
\left[f_1(y)-f_1(x)\right]\,d\mu(y)
\nonumber\\
&\quad+\left[f_1(x)-f_1(x')\right]
\int_XK(x',y)u_4(y)\,d\mu(y)
\nonumber\\
&=:Y_{4,1}(x,x')+Y_{4,2}(x,x').
\end{align}
To deal with $Y_{4,1}(x,x')$, from
${\mathrm{supp\,}} v_4\subset[B_\rho(x,A\rho(x,x'))]^\complement$,
$0\leq v_4\leq1$,
\eqref{CZk2}, and \eqref{f1r},
we deduce that
\begin{align*}
|Y_{4,1}(x,x')|
&\leq\int_{\rho(x,y)\ge A\rho(x,x')}
\left|K(x,y)-K(x',y)\right|
\left|f_1(y)-f_1(x)\right|\,d\mu(y)
\nonumber\\
&\lesssim\int_{\rho(x,y)\ge A\rho(x,x')}
\left[\frac{\rho(x,x')}{\rho(x,y)}\right]^\varepsilon
\frac{1}{V_\rho(x,y)}\\
&\quad\times\left[\frac{\rho(x,y)}{\rho(x_1,x)}\right]^\beta
\frac{1}{V_\rho(x_1,x)}
\left[\frac{r}{\rho(x_1,x)}\right]^\gamma\,d\mu(y).
\end{align*}
Using this, $\beta\in(0,\varepsilon)$,
and Lemma~\ref{A3.2}(i),
we find that
\begin{align}\label{Y41}
|Y_{4,1}(x,x')|
&\lesssim\left[\frac{\rho(x,x')}{\rho(x_1,x)}\right]^\beta
\frac{1}{V_\rho(x_1,x)}
\left[\frac{r}{\rho(x_1,x)}\right]^\gamma\nonumber\\
&\quad\times\int_{\rho(x,y)\ge A\rho(x,x')}
\left[\frac{\rho(x,x')}{\rho(x,y)}\right]^{\varepsilon-\beta}
\frac{1}{V_\rho(x,y)}\,d\mu(y)
\nonumber\\
&\lesssim\left[\frac{\rho(x,x')}{\rho(x_1,x)}\right]^\beta
\frac{1}{V_\rho(x_1,x)}
\left[\frac{r}{\rho(x_1,x)}\right]^\gamma.
\end{align}
To deal with $Y_{4,2}(x,x')$, from \eqref{Lem3.4} and \eqref{f1r},
we infer that
\begin{align*}
|Y_{4,2}(x,x')|
&\lesssim\left|f_1(x)-f_1(x')\right|
\left|Tu_4(x')\right|\nonumber\\
&\lesssim\left[\frac{\rho(x,x')}{\rho(x_1,x)}\right]^\beta
\frac{1}{V_\rho(x_1,x)}
\left[\frac{r}{\rho(x_1,x)}\right]^\gamma.
\end{align*}
By this, \eqref{Y3+Y4}, \eqref{Y3(x-x')}, \eqref{Y4(x-x')},
and \eqref{Y41}, we further conclude that
\begin{align}
\label{126.26x2}
\left|Tf_1(x)-Tf_1(x')\right|
&\leq\left|Y_3(x)-Y_3(x')\right|
+\left|Y_{4,1}(x,x')\right|
+\left|Y_{4,2}(x,x')\right|
\nonumber\\
&\lesssim\left[\frac{\rho(x,x')}{\rho(x_1,x)}\right]^\beta
\frac{1}{V_\rho(x_1,x)}
\left[\frac{r}{\rho(x_1,x)}\right]^\gamma.
\end{align}

For $Tf_2$, we have
\begin{align}\label{Tf2(x-x')}
\left|Tf_2(x)-Tf_2(x')\right|
&=\left|\int_X\left[K(x,y)-K(x',y)\right]
f_2(y)\,d\mu(y)\right|\nonumber\\
&\leq\int_X
\left|K(x,y)-K(x',y)-K(x,x_1)
+K(x',x_1)\right|
|f_2(y)|\,d\mu(y)\nonumber\\
&\quad+\left|K(x,x_1)-K(x',x_1)\right|
\left|\int_Xf_2(y)\,d\mu(y)\right|
\nonumber\\
&=:Y_5(x,x')+Y_6(x,x').
\end{align}
To deal with $Y_5(x,x')$, from
${\mathrm{supp\,}}f_2\subset{\mathrm{supp\,}}u_2
\subset B_\rho(x_1,A^{-2}\rho(x_1,x))$,
\eqref{CZO1}, and the size condition of $f$,
we deduce that
\begin{align*}
Y_5(x,x')
&\leq\int_{\rho(x_1,y)<A^{-2}\rho(x_1,x)}
\left|K(x,y)-K(x',y)-K(x,x_1)+K(x',x_1)\right|
|f(y)|\,d\mu(y)\nonumber\\
&\lesssim
\int_{\rho(x_1,y)<A^{-2}\rho(x_1,x)}
\left\{\frac{\rho(x,x')\rho(x_1,y)}{[\rho(x_1,x)]^2}\right\}^\varepsilon
\frac{1}{V_\rho(x_1,x)}\nonumber\\
&\quad\times\frac{1}{(V_\rho)_r(x_1)+V_\rho(x_1,y)}
\left[\frac{r}{r+\rho(x_1,y)}\right]^\gamma\,d\mu(y),
\end{align*}
which further implies that
\begin{align}\label{Y5}
Y_5(x,x')
&\lesssim
\int_{\rho(x_1,y)<A^{-2}\rho(x_1,x)}
\left[\frac{\rho(x,x')}{\rho(x_1,x)}\right]^\beta
\left[\frac{\rho(x_1,y)}{\rho(x_1,x)}\right]^\varepsilon
\nonumber\\
&\quad\times\frac{1}{V_\rho(x_1,x)}
\frac{1}{V_\rho(x_1,y)}
\left[\frac{r}{\rho(x_1,y)}\right]^\gamma\,d\mu(y)
\quad\text{since $\beta\in(0,\varepsilon)$}\nonumber\\
&\lesssim\left[\frac{\rho(x,x')}{\rho(x_1,x)}\right]^\beta
\frac{1}{V_\rho(x_1,x)}\left[\frac{r}{\rho(x_1,x)}\right]^\gamma\nonumber\\
&\quad\times\int_{\rho(x_1,y)<A^{-2}\rho(x_1,x)}
\frac{1}{V_\rho(x_1,y)}
\left[\frac{\rho(x_1,y)}{\rho(x_1,x)}\right]^{\varepsilon-\gamma}
\,d\mu(y)
\nonumber\\
&\lesssim\left[\frac{\rho(x,x')}{\rho(x_1,x)}\right]^\beta
\frac{1}{V_\rho(x_1,x)}\left[\frac{r}{\rho(x_1,x)}\right]^\gamma
\quad\text{by Lemma~\ref{A3.2}(i)}.
\end{align}
To deal with $Y_6(x,x')$, note that $\rho(x,x')\leq A^{-2}\rho(x_1,x)$,
by this,
\eqref{CZk2}, \eqref{intf2},
and $\beta\in(0,\varepsilon)$, we find that
\begin{align*}
Y_6(x,x')
&\lesssim
\left[\frac{\rho(x,x')}{\rho(x_1,x)}\right]^\varepsilon
\frac{1}{V_\rho(x_1,x)}
\left[\frac{r}{\rho(x_1,x)}\right]^\gamma
\lesssim\left[\frac{\rho(x,x')}{\rho(x_1,x)}\right]^\beta
\frac{1}{V_\rho(x_1,x)}\left[\frac{r}{\rho(x_1,x)}\right]^\gamma.
\end{align*}
From this, \eqref{Tf2(x-x')}, and \eqref{Y5},
we further infer that
\begin{align}
\label{12.26x3}
\left|Tf_2(x)-Tf_2(x')\right|
\lesssim\left[\frac{\rho(x,x')}{\rho(x_1,x)}\right]^\beta
\frac{1}{V_\rho(x_1,x)}\left[\frac{r}{\rho(x_1,x)}\right]^\gamma.
\end{align}

For $Tf_3$, by $\rho(x,x')\leq A^{-2}\rho(x_1,x)$,
we conclude that
\begin{align*}
\left|Tf_3(x)-Tf_3(x')\right|
&\leq\int_{\rho(x,y)\ge A^{-3}\rho(x_1,x)}
\left|K(x,y)-K(x',y)\right|
|f_3(y)|\,d\mu(y)
\quad\text{by \eqref{f3size}}
\\
&\lesssim\int_{\rho(x,y)\ge A^{-3}\rho(x_1,x)}
\left[\frac{\rho(x,x')}{\rho(x,y)}\right]^\varepsilon
\frac{1}{V_\rho(x,y)}|f_3(y)|\,d\mu(y)
\quad\text{by \eqref{CZk2}}
\\
&\lesssim\left[\frac{\rho(x,x')}{\rho(x_1,x)}\right]^\beta
\frac{1}{V_\rho(x_1,x)}\int_X|f_3(y)|\,d\mu(y)
\quad\text{since $\beta\in(0,\varepsilon)$}
\\
&\lesssim\left[\frac{\rho(x,x')}{\rho(x_1,x)}\right]^\beta
\frac{1}{V_\rho(x_1,x)}
\left[\frac{r}{\rho(x_1,x)}\right]^\gamma
\quad\text{by \eqref{intf3}},
\end{align*}
which, combined with the fact that $f=f_1+f_2+f_3$,
\eqref{126.26x2}, \eqref{12.26x3}, and \eqref{qq-pt1},
further implies that
\begin{align*}
\left|Tf(x)-Tf(x')\right|
&\leq\sum_{i=1}^3\left|Tf_i(x)-Tf_i(x')\right|\\
&\lesssim\left[\frac{\rho(x,x')}{\rho(x_1,x)}\right]^\beta
\frac{1}{V_\rho(x_1,x)}
\left[\frac{r}{\rho(x_1,x)}\right]^\gamma,
\end{align*}
which is exactly \eqref{TfrCase3}.
This finishes the proof of \eqref{Tfr} in the case
that $\rho(x,x')<\frac{r+\rho(x_1,x)}{A^5}$ and $\rho(x_1,x)\ge Ar$.

From the above three cases, we finish the proof
of \eqref{Tfr}. This concludes Step 2.
Given what we have shown in Step~1 and Step~2,
we can immediately conclude that
\begin{align*}
\left\|Tf\right\|_{{\mathcal G}(x_1,r,\beta,\gamma)}
\lesssim\|f\|_{{\mathcal G}(x_1,r,\beta,\gamma)},
\end{align*}
which completes the proof of Theorem~\ref{CZOonGO}.
\end{proof}

\begin{remark}
\label{1262rem}
\begin{enumerate}
\item[\textup{(i)}]
We stress here that
the cancellation condition of $f$
is only used in the proof of \eqref{intf2}.

\item[\textup{(ii)}]
It was pointed out in \cite[Remark 3.2]{HLYY}
that,
compared with \cite[Theorem 2.18]{hmy08},
the conditions (ii) and (iii) of Theorem~\ref{CZOonGO}
used here are a little bit stronger
but is enough to establish Calder\'on reproducing
formulae.

\item[\textup{(iii)}]
The proof of Theorem~\ref{CZOonGO} strongly relies on the existence
of smooth cut-off functions as in Corollary~\ref{Lem3.3}.
Let $T$ be an $\varepsilon$-Calder\'on--Zygmund operator
defined as in Theorem~\ref{CZOonGO}.
If $\rho\in\mathbf{q}$ is a metric,
then ${\mathrm{ind\,}}(X,\mathbf{q})\in[1,\infty]$ and,
for any $\beta,\,\gamma\in(0,{\mathrm{ind\,}}(X,\mathbf{q}))$,
$T$ is bounded on $\mathring{{\mathcal G}}(x_1,r,\beta,\gamma)$,
and, if $\rho\in\mathbf{q}$ is an ultrametric,
then, for any $\beta,\gamma\in(0,\infty)$,
$T$ is bounded on $\mathring{{\mathcal G}}(x_1,r,\beta,\gamma)$.
However, in \cite[Theorem 3.1]{HLYY}, He et al. only
proved that $T$ is bounded
on $\mathring{{\mathcal G}}(x_1,r,\beta,\gamma)$
with $0<\beta,\,\gamma<1$. In this sense,
Theorem~\ref{CZOonGO} is sharper than \cite[Theorem 3.1]{HLYY}.
\end{enumerate}
\end{remark}

Finally, we establish the following boundedness
of Calder\'on--Zygmund operators on
${\mathcal G}(x_1,r,\beta,\gamma)$ which sharpens
\cite[Theorem 3.5]{HLYY}; see Remark~\ref{xkh-45} below.

\begin{theorem}
\label{CZOonCG}
Let $(X,\mathbf{q},\mu)$ be a quasi-ultrametric space
of homogeneous type. Let
$\rho\in\mathbf{q}$ have the property that
all $\rho$-balls are $\mu$-measurable and
$0<\varepsilon\preceq{\mathrm{ind\,}}(X,\mathbf{q})$
with ${\mathrm{ind\,}}(X,\mathbf{q})$ as in \eqref{index}.
Assume that $T$ is an
$\varepsilon$-Calder\'on--Zygmund operator (with respect to $\rho$),
where its kernel $K$ is
the same as in Theorem~\ref{CZOonGO}.
Assume that there exist positive
constants $C_T$, $r_0$, and $\sigma$
such that, for any $x,y,x'\in X$,
\begin{enumerate}
\item[\textup{(i)}]
If $\rho(x,y)\ge r_0$, then
\begin{equation}\label{CZO3}
|K(x,y)|\leq C_T\frac{1}{V_\rho(x,y)}
\left[\frac{r_0}{\rho(x,y)}\right]^\sigma.
\end{equation}

\item[\textup{(ii)}]
If $\rho(x,y)\ge r_0$
and $\rho(x,x')\leq\frac{\rho(x,y)}{2C_\rho}$
with $C_\rho\in[1,\infty)$ as in \eqref{Crho}, then
\begin{equation}\label{CZO4}
|K(x,y)-K(x',y)|\leq C_T\frac{1}{V_\rho(x,y)}
\left[\frac{\rho(x,x')}{\rho(x,y)}\right]^\varepsilon
\left[\frac{r_0}{\rho(x,y)}\right]^\sigma.
\end{equation}
\end{enumerate}
Let $\beta\in(0,\varepsilon)$ and
$\gamma\in(0,\varepsilon)\cap(0,\sigma]$.
Then, for any
$f\in{\mathcal G}(x_1,r_0,\beta,\gamma)$ with $x_1\in X$,
\begin{align}\label{A13}
\|Tf\|_{{\mathcal G}(x_1,r_0,\beta,\gamma)}\leq C
\left[C_T+\|T\|_{L^2(X)\to L^2(X)}+c_0\right]
\|f\|_{{\mathcal G}(x_1,r_0,\beta,\gamma)},
\end{align}
where $c_0$ is the same as in \eqref{CZO2} and
$C$ is a positive constant
independent of $x_1$, $r_0$, $c_0$, and $f$.
Moreover, if $T$ further satisfies that, for any $y\in X$,
$$\int_XK(x,y)\,d\mu(x)=0,$$
then, for any $f\in{\mathcal G}(x_1,r_0,\beta,\gamma)$,
$Tf\in\mathring{\mathcal{G}}(x_1,r_0,\beta,\gamma)$.
\end{theorem}

\begin{proof}
Similarly to the proof of Theorem~\ref{CZOonGO},
we only need to show \eqref{Tfs} and \eqref{Tfr}
for any $f\in{\mathcal G}(x_1,r_0,\beta,\gamma)$ with
$\|f\|_{{\mathcal G}(x_1,r_0,\beta,\gamma)}=1$ but without
the cancellation condition on $f$. In addition,
without loss of generality, we may assume that
$\rho\in\mathbf{q}$ is a regularized quasi-ultrametric as in Definition~\ref{D2.3}. According to Remark~\ref{RmkD2.3}(iii),
we may also assume that $\rho\in\mathbf{q}$
satisfies $\varepsilon\leq(\log_2C_\rho)^{-1}$.

We first prove \eqref{Tfs}. Following the
arguments used in the proof of \eqref{Tfs}
in Theorem~\ref{CZOonGO}, we consider the cases
$\rho(x_1,x)<Ar_0$ and $\rho(x_1,x)\ge Ar_0$,
respectively,
where $A=2C_\rho\in[2,\infty)$ is the same
as in the proof of Theorem~\ref{CZOonGO}.
By Remark~\ref{1262rem}(i), we find that
we only need to re-estimate $Tf_2$
in the case $\rho(x_1,x)\ge Ar_0$
without using \eqref{intf2}.
Maintaining the notation used in the proof of Theorem~\ref{CZOonGO},
since we have
$$
{\mathrm{supp\,}} f_2\subset{\mathrm{supp\,}}
u_2\subset B_\rho\left(x_1,A^{-2}\rho(x_1,x)\right),
$$
from Definition~\ref{D2.1}(iii) and $A=2C_\rho$,
it follows that, for any $y\in{\mathrm{supp\,}} f_2$,
\begin{align}
\label{2038}
Ar_0\leq\rho(x_1,x)\leq C_\rho
\max\{\rho(x_1,y),\,\rho(y,x)\}=C_\rho\rho(x,y),
\end{align}
which further implies that
$$
\rho(x,y)\ge r_0
\ \ \text{and}\ \
V_\rho(x,y)\gtrsim V_\rho(x_1,x),
$$
granted \eqref{upperdoub}.
By this, Theorem~\ref{CZOonGO}(ii),
the fact that $|f_2|\leq|f|$,
\eqref{CZO3}, and the size condition of $f$,
we conclude that
\begin{align*}
|Tf_2(x)|
&\leq\int_{\rho(x_1,y)< A^{-2}\rho(x_1,x)}
|K(x,y)||f(y)|\,d\mu(y)
\\
&\lesssim\int_{\rho(x_1,y)< A^{-2}\rho(x_1,x)}
\frac{1}{V_\rho(x,y)}\left[\frac{r_0}{\rho(x,y)}\right]^\sigma
\\
&\quad\times
\frac{1}{(V_\rho)_{r_0}(x_1)+V_\rho(x_1,y)}
\left[\frac{r_0}{r_0+\rho(x_1,y)}\right]^\gamma\,d\mu(y).
\end{align*}
From this, Lemma~\ref{A3.2}(ii),
$\rho(x_1,x)\ge Ar_0$, and $\gamma\leq\sigma$,
we deduce that
\begin{align*}
|Tf_2(x)|
&\lesssim\frac{1}{V_\rho(x,x_1)}
\left[\frac{r_0}{\rho(x_1,x)}\right]^{\sigma-\gamma}
\left[\frac{r_0}{\rho(x_1,x)}\right]^\gamma
\lesssim\frac{1}{V_\rho(x,x_1)}
\left[\frac{r_0}{\rho(x_1,x)}\right]^\gamma,
\end{align*}
which is the desired estimate for $Tf_2$ [see \eqref{Tf2}].
This finishes the proof of \eqref{Tfs}.

Next, we show \eqref{Tfr}.
Let $x,x'\in X$ satisfy $\rho(x,x')\leq\frac{r+\rho(x_1,x)}{2C_\rho}$
and we prove \eqref{Tfr} via following a similar argument to that used
in the proof of \eqref{Tfr} in Theorem~\ref{CZOonGO}.
Note that, in the proof of Theorem~\ref{CZOonGO},
the cancellation condition of $f$
is used in dealing with the regularity
condition of $Tf_2$ (see \eqref{Tf2(x-x')}),
which relies on the estimate of
\eqref{intf2}.
By Remark~\ref{1262rem}(i), we only need to re-estimate $Tf_2$
when $x,x'\in X$ satisfy
$\rho(x,x')<\frac{r_0+\rho(x_1,x)}{A^5}$ and $\rho(x_1,x)\ge Ar_0$,
without using \eqref{intf2}.
To this end, let
$x,x'\in X$ satisfy
$\rho(x,x')<\frac{r_0+\rho(x_1,x)}{A^5}$ and $\rho(x_1,x)\ge Ar_0$.
From \eqref{2038}, it follows that
$$
\rho(x,y)\ge\frac{\rho(x_1,x)}{C_\rho}>r_0,
$$
which further implies that
$$
V_\rho(x,y)\gtrsim V_\rho(x_1,x)
$$
and
$$
\rho(x,x')<\frac{r_0+\rho(x_1,x)}{A^5}\leq
\frac{\rho(x,y)}{2C_\rho}.
$$
By these, we find that
\begin{align*}
&\left|Tf_2(x)-Tf_2(x')\right|\\
&\quad\lesssim\int_{\rho(x_1,y)<A^{-2}\rho(x_1,x)}
\frac{1}{V_\rho(x,y)}
\left[\frac{\rho(x,x')}{\rho(x,y)}\right]^\varepsilon
\left[\frac{r_0}{\rho(x,y)}\right]^\sigma|f(y)|\,d\mu(y)
\\
&\quad\qquad\text{by \eqref{CZO4}}\\
&\quad\lesssim
\frac{1}{V_\rho(x_1,x)}
\left[\frac{\rho(x,x')}{\rho(x_1,x)}\right]^\varepsilon
\left[\frac{r_0}{\rho(x_1,x)}\right]^\sigma
\\
&\qquad\times
\int_{X}
\frac{1}{(V_\rho)_{r_0}(x_1)+V_\rho(x_1,y)}
\left[\frac{r_0}{r_0+\rho(x_1,y)}\right]^\gamma\,d\mu(y)\\
&\quad\qquad\text{by the size condition of $f$ [see Definition~\ref{testfunction}(i)]},
\end{align*}
which, together with Lemma~\ref{A3.2}(ii)
and the assumptions $\beta\in(0,\varepsilon)$ and $\gamma\in(0,\sigma]$,
further implies that
\begin{align*}
&\left|Tf_2(x)-Tf_2(x')\right|\\
&\quad\lesssim
\frac{1}{V_\rho(x_1,x)}
\left[\frac{\rho(x,x')}{\rho(x_1,x)}\right]^{\varepsilon-\beta}
\left[\frac{\rho(x,x')}{\rho(x_1,x)}\right]^\beta
\left[\frac{r_0}{\rho(x_1,x)}\right]^{\sigma-\gamma}
\left[\frac{r_0}{\rho(x_1,x)}\right]^\gamma
\\
&\quad\lesssim\left[\frac{\rho(x,x')}{r_0+\rho(x_1,x)}\right]^{\beta}
\frac{1}{V_\rho(x_1,x)+(V_\rho)_{r_0}(x_1)}
\left[\frac{r_0}{r_0+\rho(x_1,x)}\right]^\gamma.
\end{align*}
This finishes the proof of \eqref{Tfr}
and hence Theorem~\ref{CZOonCG}.
\end{proof}

\begin{remark}\label{xkh-45}
The proof of Theorem~\ref{CZOonCG} strongly relies on the existence
of smooth cut-off functions as in Corollary~\ref{Lem3.3}.
Let $T$ be an $\varepsilon$-Calder\'on--Zygmund operator
defined as in Theorem~\ref{CZOonCG}.
If $\rho\in\mathbf{q}$ is a metric,
then ${\mathrm{ind\,}}(X,\mathbf{q})\in[1,\infty]$ and,
for any $\beta,\,\gamma\in(0,{\mathrm{ind\,}}(X,\mathbf{q}))$,
$T$ is bounded on ${\mathcal G}(x_1,r,\beta,\gamma)$,
and, if $\rho\in\mathbf{q}$ is an ultrametric,
then, for any $\beta,\gamma\in(0,\infty)$,
$T$ is bounded on ${\mathcal G}(x_1,r,\beta,\gamma)$.
However, in \cite[Theorem 3.5]{HLYY}, He et al. only
showed that $T$ is bounded on ${\mathcal G}(x_1,r,\beta,\gamma)$
with $0<\beta,\,\gamma<1$. In this sense,
Theorem~\ref{CZOonCG} is sharper than \cite[Theorem 3.5]{HLYY}.
\end{remark}

\section[Homogeneous Continuous
Calder\'on Reproducing Formulae]{Homogeneous Continuous
Calder\'on Reproducing \\Formulae}
\label{SS3.2}

The main goal of this section is to prove sharp
homogeneous continuous
Calder\'on reproducing formulae
on ultra-RD-spaces $(X,\mathbf{q},\mu)$,
where $\mu$ is assumed to be a Borel-semiregular measure on $X$,
which we present in Theorem~\ref{A3.19}.
For any $\rho\in\mathbf{q}$,
we always assume that $\mathrm{diam}_\rho(X)=\infty$.
From this and Remark~\ref{1031},
it is easy to infer that $\mu(X)=\infty$
and, for any $\rho\in\mathbf{q}$ and $x\in X$,
$R_\rho(x)=\infty$.

We begin with first introducing the following notion of
homogeneous approximations of the identity
on quasi-ultrametric spaces
of homogeneous type.
The reader is reminded that the symbol ``$\preceq$"
was introduced in Remark~\ref{rmk390}(iii).

\begin{definition}\label{DefHATI}
Let $(X,\mathbf{q},\mu)$ be a quasi-ultrametric space
of homogeneous type and
$0<\varepsilon\preceq{\mathrm{ind\,}}(X,\mathbf{q})$
with ${\mathrm{ind\,}}(X,\mathbf{q})$ as in \eqref{index}.
A sequence $\{\mathcal{S}_{2^{-k}}\}_{k\in\mathbb{Z}}$ of
linear operators is called a
\emph{homogeneous approximation of the
identity of order $\varepsilon$}\index[words]{homogeneous approximation of the identity}
if $\{\mathcal{S}_{2^{-k}}\}_{k\in\mathbb{Z}}$
satisfies all the conditions of Definition~\ref{DefATI}
with $t$ replaced by $2^{-k}$.
\end{definition}

By Theorem~\ref{corATI}, we have the following properties
for any homogeneous approximation of the identity.

\begin{proposition}\label{A3.11}
Let $(X,\mathbf{q},\mu)$ be a quasi-ultrametric space
of homogeneous type,
where $\mu$ is assumed to be a Borel-semiregular measure on $X$.
Let $0<\varepsilon\preceq{\mathrm{ind\,}}(X,\mathbf{q})$
with ${\mathrm{ind\,}}(X,\mathbf{q})$ as in \eqref{index}
and
$\{\mathcal{S}_{2^{-k}}\}_{k\in\mathbb{Z}}$
be a homogeneous approximation
of the identity of order $\varepsilon$
as in Definition~\ref{DefHATI}.
Then the following assertions hold:
\begin{enumerate}
\item[\textup{(i)}]
For any $p\in(1,\infty)$ and $f\in L^p(X)$,
$\|\mathcal{S}_{2^{-k}}f\|_{L^p(X)}\to 0$
as $k\to-\infty$.

\item[\textup{(ii)}]
For any $k\in\mathbb{Z}$,
let $\mathcal{D}_k:=\mathcal{S}_{2^{-k}}-\mathcal{S}_{2^{-k+1}}$.
Then, for any $p\in(1,\infty)$,
\begin{align*}
\sum_{k=-\infty}^{\infty}\mathcal{D}_k=\mathcal{I}
\end{align*}
in $L^p(X)$,
where $\mathcal{I}$ is the identity operator on $L^p(X)$.
\end{enumerate}
\end{proposition}

\begin{proof}
Since $\mathrm{diam}_\rho(X)=\infty$, \textup{(i)}
follows immediately
from Theorem~\ref{corATI}(vii),
whereas \textup{(ii)} is a simple
corollary of \textup{(i)} and Theorem~\ref{corATI}(ix).
This finishes the proof of Proposition~\ref{A3.11}.
\end{proof}

To establish the sharp homogeneous continuous
Calder\'on reproducing formulae,
we need several technical lemmas.
The following conclusion follows directly from Definition~\ref{DefATI}
(see also Remark~\ref{Rm} and Lemma~\ref{ATItest});
we omit the details.

\begin{lemma}\label{Pk}
Let $(X,\mathbf{q},\mu)$ be a quasi-ultrametric space
of homogeneous type and
$0<\varepsilon\preceq{\mathrm{ind\,}}(X,\mathbf{q})$.
Assume that
$\{\mathcal{S}_{2^{-k}}\}_{k\in\mathbb{Z}}$
is a homogeneous approximation
of the identity of order $\varepsilon$ as
in Definition~\ref{DefHATI}.
For any $k\in\mathbb{Z}$,
let
$\mathcal{D}_k:=\mathcal{S}_{2^{-k}}-\mathcal{S}_{2^{-k+1}}$.
Then the kernels $\{D_k\}_{k\in\mathbb{Z}}$ of the
operators $\{\mathcal{D}_k\}_{k\in\mathbb{Z}}$ satisfy that
there exist two constants $C,C'\in(0,\infty)$
and $\rho\in\mathbf{q}$, where all $\rho$-balls are $\mu$-measurable,
such that, for any $k\in\mathbb{Z}$,
\begin{enumerate}
\item[\textup{(i)}]
for any $x,y\in X$,
\begin{align*}
\left|D_k(x,y)\right|\leq C\frac{1}{(V_\rho)_{2^{-k}}(x)};
\end{align*}
in particular, if $\rho(x,y)\ge C'2^{-k}$, then $D_k(x,y)=0$;

\item[\textup{(ii)}]
for any $x,y,x'\in X$,
$$
\left|D_k(x,y)-D_k(x',y)\right|
\leq C2^{k\varepsilon}\left[\rho(x,x')\right]^\varepsilon
\left[\frac{1}{(V_\rho)_{2^{-k}}(x)}
+\frac{1}{(V_\rho)_{2^{-k}}(x')}\right];
$$

\item[\textup{(iii)}]
for any $x,y,x',y'\in X$,
\begin{align*}
&\left|\left[D_k(x,y)-D_k(x',y)\right]
-\left[D_k(x,y')-D_k(x',y')\right]\right|\\
&\quad\leq C4^{k\varepsilon}
\left[\rho(x,x')\rho(y,y')\right]^\varepsilon
\left[\frac{1}{(V_\rho)_{2^{-k}}(x)}
+\frac{1}{(V_\rho)_{2^{-k}}(x')}\right];
\end{align*}

\item[\textup{(iv)}]
for any $x,y\in X$,
$D_k(x,y)=D_k(y,x)$;

\item[\textup{(v)}]
for any $x\in X$,
$$
\int_XD_k(x,z)\,d\mu(z)=\int_XD_k(z,x)\,d\mu(z)=0.
$$
\end{enumerate}
Moreover, the following analogous conclusions in Remark~\ref{Rm} also hold
for any $k\in\mathbb{Z}$.
\begin{enumerate}
\item[\textup{(vi)}]
For any $x,y\in X$ with $x\neq y$,
and $\gamma\in\mathbb{R}_+$,
\begin{align*}
|D_k(x,y)|
\lesssim\frac{1}{V_\rho(x,y)+(V_\rho)_{2^{-k}}(x)+(V_\rho)_{2^{-k}}(y)}
\left[\frac{{2^{-k}}}{{2^{-k}}+\rho(x,y)}\right]^\gamma,
\end{align*}
where the implicit positive constant is independent of
$k$, $x$, and $y$, but depends on $\gamma$.

\item[\textup{(vii)}]
For any $x,y,x',y'\in X$,
\begin{align*}
&\left|\left[D_k(x,y)-D_k(x',y)\right]
-\left[D_k(x,y')-D_k(x',y')\right]\right|\\
&\quad\lesssim4^{k\varepsilon}\left[\rho(x,x')
\rho(y,y')\right]^\varepsilon
\left[\frac{1}{(V_\rho)_{2^{-k}}(y)}
+\frac{1}{(V_\rho)_{2^{-k}}(y')}\right],
\end{align*}
where the implicit positive constant is independent of
$k$, $x$, $x'$, $y$, and $y'$.
\end{enumerate}
Finally, for any $x\in X$ and $\gamma\in(0,\infty)$,
\begin{align*}
D_k(\cdot,x),D_k(x,\cdot)\in\mathring{\mathcal{G}}(x,2^{-k},\varepsilon,\gamma)
\end{align*}
with their norms independent of both $x$ and $k$.
\end{lemma}

Now, we establish the following estimates
for the kernel of the compositions of two
homogeneous approximations
of the identity.
The reader is reminded that,
for any $a,b\in\mathbb{R}$,
$a\land b:=\min\{a,\,b\}$.

\begin{lemma}\label{A3.12}
Let $(X,\mathbf{q},\mu)$ be a quasi-ultrametric space
of homogeneous type and
$0<\varepsilon\preceq{\mathrm{ind\,}}(X,\mathbf{q})$
with ${\mathrm{ind\,}}(X,\mathbf{q})$ as in \eqref{index}.
Assume that $\{\mathcal{S}_{2^{-k}}\}_{k\in\mathbb{Z}}$
and $\{\mathcal{E}_{2^{-k}}\}_{k\in\mathbb{Z}}$
are two homogeneous approximations
of the identity of order $\varepsilon$
in Definition~\ref{DefHATI}.
For any $k\in\mathbb{Z}$,
let
$\mathcal{P}_k:=\mathcal{S}_{2^{-k}}-\mathcal{S}_{2^{-k+1}}$
and
$\mathcal{Q}_k:=\mathcal{E}_{2^{-k}}-\mathcal{E}_{2^{-k+1}}$
and denote their kernels by $P_k$ and $Q_k$.
Then, for any given $\varepsilon'\in(0,\varepsilon)$,
$\delta\in(0,\varepsilon-\varepsilon']$,
and $\rho\in\mathbf{q}$ having the property
that all $\rho$-balls are $\mu$-measurable,
there exist $C,C'\in(0,\infty)$,
depending on $\varepsilon'$, $\varepsilon$,
and $\delta$,
such that the kernels $\{P_lQ_k\}_{l,k\in\mathbb{Z}}$
of $\{\mathcal{P}_l\mathcal{Q}_k\}_{l,k\in\mathbb{Z}}$
have the following properties:
\begin{enumerate}
\item [\textup{(i)}]
For any $l,k\in\mathbb{Z}$ and $x,y\in X$,
\begin{equation}\label{3.2}
\left|P_lQ_k(x,y)\right|\leq C2^{-|k-l|\varepsilon}
\frac{1}{(V_\rho)_{2^{-(k\land l)}}(x)}
\end{equation}
and, if $\rho(x,y)>C'2^{-(k\land l)}$, then
$|P_lQ_k(x,y)|=0.$

\item [\textup{(ii)}]
For any $l,k\in\mathbb{Z}$ and $x,y,y'\in X$,
\begin{equation}\label{3.3}
\left|P_lQ_k(x,y)-P_lQ_k(x,y')\right|
\leq C2^{\varepsilon'(l\land k)-|k-l|\delta}
\left[\rho(y,y')\right]^{\varepsilon'}
\frac{1}{(V_\rho)_{2^{-(k\land l)}}(x)}.
\end{equation}

\item [\textup{(iii)}]
For any $l,k\in\mathbb{Z}$ and $x,x',y\in X$,
\begin{equation}\label{3.4}
|P_lQ_k(x,y)-P_lQ_k(x',y)|
\leq C2^{\varepsilon'(l\land k)-|k-l|\delta}
\left[\rho(x,x')\right]^{\varepsilon'}
\frac{1}{(V_\rho)_{2^{-(k\land l)}}(y)}.
\end{equation}

\item [\textup{(iv)}]
For any $l,k\in\mathbb{Z}$ and $x,y,x',y'\in X$,
\begin{align}
\label{3.5}
&\left|\left[P_lQ_k(x,y)-P_lQ_k(x',y)\right]
-\left[P_lQ_k(x,y')-P_lQ_k(x',y')\right]\right|
\nonumber\\
&\quad\leq C2^{2\varepsilon'(l\land k)-|k-l|\delta}
\left[\rho(x,x')\rho(y,y')\right]^{\varepsilon'}
\left[\frac{1}{(V_\rho)_{2^{-(k\land l)}}(x)}
+\frac{1}{(V_\rho)_{2^{-(k\land l)}}(x')}\right].
\end{align}
\end{enumerate}
\end{lemma}

\begin{proof}
Given Lemma~\ref{Pk} and Remark~\ref{Rm}(iii),
without loss of generality,
we may assume that $\rho$
is a regularized quasi-ultrametric as in Definition~\ref{D2.3}.
Observe that, for any $l,k\in\mathbb{Z}$ and $x,y\in X$,
$$
P_lQ_k(x,y)=\int_XP_l(x,z)Q_k(z,y)\,d\mu(z).
$$
By symmetry, we may assume that $l\leq k$.
From Lemma~\ref{Pk}(i),
we deduce that there exists a positive constant $C_1$
such that, for any $x,y,z\in X$
and $l,k\in\mathbb{Z}$,
if $\rho(x,z)\ge C_12^{-l}$, then $P_l(x,z)=0$
and, if $\rho(z,y)\ge C_12^{-k}$, then $Q_k(z,y)=0$.
We first show (i). On the one hand,
we claim that,
for any $x,y,z\in X$ satisfying $\rho(x,y)>C_1C_\rho 2^{-l}$,
$$
P_l(x,z)Q_k(z,y)=0.
$$
Indeed, if there exist $x,y,z\in X$ satisfying that
$\rho(x,y)>C_1C_\rho 2^{-l}$, $P_l(x,z)\neq0$, and $Q_k(z,y)\neq0$,
then, by Lemma~\ref{Pk}, we find that
$\rho(x,z)\leq C_12^{-l}$ and $\rho(z,y)\leq C_12^{-k}$.
From this and Definition~\ref{D2.1}(iii), it follows that
\begin{align*}
\rho(x,y)&\leq C_\rho\max\{\rho(x,z),\,\rho(z,y)\}\\
&\leq C_\rho\max\{C_12^{-l},\,C_12^{-k}\}
=C_1C_\rho2^{-l}<\rho(x,y),
\end{align*}
which is a contradiction. This finishes the proof of the above claim.
Thus, for any $x,y\in X$ with $\rho(x,y)>C_1C_\rho 2^{-l}$,
$$P_lQ_k(x,y)=\int_XP_l(x,z)Q_k(z,y)\,d\mu(z)=0.$$
Hence, the estimate in \eqref{3.2} is trivially satisfied.
On the other hand, by (i), (ii), and (v) of Lemma~\ref{Pk} and
Remark~\ref{RemThmATI}(i), we conclude that,
for any $x,y\in X$ with $\rho(x,y)\leq C_1C_\rho2^{-l}$,
\begin{align*}
|P_lQ_k(x,y)|
&\leq\int_{\rho(z,y)\leq C_12^{-k}}
|P_l(x,z)-P_l(x,y)||Q_k(z,y)|\,d\mu(z)
\\
&\lesssim\int_{\rho(z,y)\leq C_12^{-k}}
\frac{2^{l\varepsilon}[\rho(z,y)]^\varepsilon}
{(V_\rho)_{2^{-l}}(x)}\frac{1}{V_{2^{-k}}(z)}\,d\mu(z),
\end{align*}
which, combined with \eqref{doublingcond}
and Lemma~\ref{A3.2}(i),
further implies that
\begin{align*}
|P_lQ_k(x,y)|
&\lesssim\frac{2^{l\varepsilon}}
{(V_\rho)_{2^{-l}}(x)}
\int_{\rho(z,y)\leq C_12^{-k}}
\frac{[\rho(z,y)]^\varepsilon}{V_\rho(z,y)}\,d\mu(z)\\
&\lesssim\frac{2^{{(l-k)}\varepsilon}}{(V_\rho)_{2^{-l}}(x)}
=\frac{2^{{-|k-l|}\varepsilon}}{(V_\rho)_{2^{-(k\land l)}}(x)}.
\end{align*}
This finishes the proof of \eqref{3.2}.

Next, we prove \eqref{3.3}.
To this end, it suffices to show that,
for any $x,y,y'\in X$,
\begin{equation}\label{Tl}
|P_lQ_k(x,y)-P_lQ_k(x,y')|
\lesssim2^{l\varepsilon}
\left[\rho(y,y')\right]^\varepsilon
\frac{1}{(V_\rho)_{2^{-l}}(x)}.
\end{equation}
Indeed, from \eqref{3.2}, we infer that, for any $x,y,y'\in X$,
\begin{align}\label{T1x}
&\left|P_lQ_k(x,y)-P_lQ_k(x,y')\right|
\nonumber\\
&\quad\leq\left|P_lQ_k(x,y)\right|+\left|P_lQ_k(x,y')\right|
\lesssim 2^{-(k-l)\varepsilon}
\frac{1}{(V_\rho)_{2^{-l}}(x)}.
\end{align}
Using \eqref{Tl} and \eqref{T1x},
we find that,
for any $\sigma\in(0,1)$ and $x,y,y'\in X$,
\begin{align*}
&\left|P_lQ_k(x,y)-P_lQ_k(x,y')\right|\\
&\quad=\left|P_lQ_k(x,y)-P_lQ_k(x,y')\right|^\sigma
\left|P_lQ_k(x,y)-P_lQ_k(x,y')\right|^{1-\sigma}\\
&\quad\lesssim2^{l\varepsilon\sigma-(k-l)\varepsilon(1-\sigma)}
\left[\rho(y,y')\right]^{\varepsilon\sigma}
\frac{1}{(V_\rho)_{2^{-l}}(x)}.
\end{align*}
Choose
$\sigma:=\frac{\varepsilon'}{\varepsilon}$ and
$\delta\in(0,\varepsilon-\varepsilon']$.
Then $\delta\in(0,\varepsilon(1-\sigma)]$
and
\begin{align*}
|P_lQ_k(x,y)-P_lQ_k(x,y')|
&\lesssim2^{l\varepsilon'-(k-l)(\varepsilon-\varepsilon')}
\left[\rho(y,y')\right]^{\varepsilon'}
\frac{1}{(V_\rho)_{2^{-l}}(x)}
\\
&\leq
2^{(k\land l)\varepsilon'-|k-l|\delta}
\left[\rho(y,y')\right]^{\varepsilon'}
\frac{1}{(V_\rho)_{2^{-(k\land l)}}(x)}.
\end{align*}
Thus,
it remains to prove \eqref{Tl}.
From both (i) and (v) of Lemma~\ref{Pk}, we deduce that,
for any $x,y,y'\in X$,
\begin{align}
\label{12.28.x1}
&|P_lQ_k(x,y)-P_lQ_k(x,y')|
\nonumber\\
&\quad=\left|\int_X[P_l(x,z)-P_l(x,y)]
[Q_k(z,y)-Q_k(z,y')]\,d\mu(z)\right|
\nonumber\\
&\quad\leq\int_{\rho(z,y)\leq C_1C_\rho2^{-k}}
\left|P_l(x,z)-P_l(x,y)\right|
\left|Q_k(z,y)-Q_k(z,y')\right|\,d\mu(z)
\nonumber\\
&\qquad+\int_{\rho(z,y)>C_1C_\rho2^{-k}}
\left|P_l(x,z)-P_l(x,y)\right|
\left|Q_k(z,y)-Q_k(z,y')\right|\,d\mu(z)
\nonumber\\
&\quad=:Y_1(x,y,y')+Y_2(x,y,y').
\end{align}
For $Y_1(x,y,y')$, by Remark~\ref{RemThmATI}(i) and
\eqref{upperdoub}, we obtain
\begin{align*}
Y_1(x,y,y')
&\lesssim\int_{\rho(z,y)\leq C_1C_\rho2^{-k}}
\frac{2^{l\varepsilon}[\rho(z,y)]^\varepsilon}{(V_\rho)_{2^{-l}}(x)}
\frac{2^{k\varepsilon}
[\rho(y,y')]^\varepsilon}{(V_\rho)_{2^{-k}}(z)}\,d\mu(z)
\nonumber\\
&\lesssim2^{(l+k)\varepsilon}[\rho(y,y')]^\varepsilon
\frac{1}{(V_\rho)_{2^{-l}}(x)}\int_{\rho(z,y)\leq C_1C_\rho2^{-k}}
\frac{[\rho(z,y)]^\varepsilon}{(V_\rho)_{C_1C_\rho2^{-k}}(z)}\,d\mu(z),
\end{align*}
which, together with Lemma~\ref{A3.2}(i),
further implies that
\begin{align}
\label{12.28.x2}
Y_1(x,y,y')
&\lesssim2^{(l+k)\varepsilon}[\rho(y,y')]^\varepsilon
\frac{1}{(V_\rho)_{2^{-l}}(x)}
\int_{\rho(z,y)\leq C_1C_\rho2^{-k}}
\frac{[\rho(z,y)]^\varepsilon}{V_\rho(z,y)}\,d\mu(z)
\nonumber\\
&\lesssim2^{l\varepsilon}[\rho(y,y')]^\varepsilon
\frac{1}{(V_\rho)_{2^{-l}}(x)}.
\end{align}

To proceed, we write $Q_k=E_{2^{-k}}-E_{2^{-k+1}}$,
where $\{E_{2^{-k}}\}_{k\in\mathbb{Z}}$
satisfies (i)--(v) in Definition~\ref{DefATI}.
For $Y_2(x,y,y')$, on the one hand,
if $\rho(z,y')\ge C_12^{-k}$,
then, from Lemma~\ref{Pk} and $\rho(z,y)>C_1C_\rho2^{-k}$,
it follows that $Q_k(z,y)=0=Q_k(z,y')$ and hence
$$
Y_2(x,y,y')=0.
$$
On the other hand, if $\rho(z,y')<C_12^{-k}$, then,
by Definition~\ref{D2.1}(iii) and $\rho(z,y)>C_1C_\rho2^{-k}$,
we conclude that
\begin{align*}
\rho(z,y)\leq C_\rho\max\{\rho(z,y'),\,\rho(y,y')\}
\leq C_\rho\max\{C_12^{-k},\,\rho(y,y')\}=C_\rho\rho(y,y').
\end{align*}
From this, Remark~\ref{RemThmATI}(i),
and Definition~\ref{DefATI}(v), we infer that
\begin{align*}
Y_2(x,y,y')
&\lesssim\int_{C_\rho\rho(y,y')>\rho(z,y)\ge C_1C_\rho2^{-k}}
\frac{2^{l\varepsilon}[\rho(z,y)]^\varepsilon}{(V_\rho)_{2^{-l}}(x)}
\left|Q_k(z,y')\right|\,d\mu(z)
\\
&\lesssim
\frac{2^{l\varepsilon}[\rho(y,y')]^\varepsilon}{(V_\rho)_{2^{-l}}(x)}
\int_{X}\left[E_{2^{-k}}(z,y')+E_{2^{-k+1}}(z,y')\right]\,d\mu(z)
\sim\frac{2^{l\varepsilon}[\rho(y,y')]^\varepsilon}{(V_\rho)_{2^{-l}}(x)},
\end{align*}
which, combined with both \eqref{12.28.x1} and \eqref{12.28.x2},
further implies \eqref{Tl} and hence \eqref{3.3}.

Since the proof of \eqref{3.4} is similar to that of \eqref{3.3},
we now justify \eqref{3.5}.
To this end, it suffices to verify that,
for any $x,y,x',y'\in X$,
\begin{align}
\label{418x}
&\left|\left[P_lQ_k(x,y)-P_lQ_k(x',y)\right]
-\left[P_lQ_k(x,y')-P_lQ_k(x',y')\right]\right|
\nonumber\\
&\quad\lesssim2^{2l\varepsilon}
\left[\rho(x,x')\rho(y,y')\right]^{\varepsilon}
\left[\frac{1}{(V_\rho)_{2^{-l}}(x)}
+\frac{1}{(V_\rho)_{2^{-l}}(x')}\right].
\end{align}
Indeed, by \eqref{3.2}, we find that,
for any $x,y,x',y'\in X$,
\begin{align*}
&\left|\left[P_lQ_k(x,y)-P_lQ_k(x',y)\right]
-\left[P_lQ_k(x,y')-P_lQ_k(x',y')\right]\right|
\nonumber\\	
&\quad\leq|P_lQ_k(x,y)|+|P_lQ_k(x',y)|
+|P_lQ_k(x,y')|+|P_lQ_k(x',y')|
\nonumber\\
&\quad\lesssim 2^{-(k-l)\varepsilon}
\left[\frac{1}{(V_\rho)_{2^{-l}}(x)}
+\frac{1}{(V_\rho)_{2^{-l}}(x')}\right].
\end{align*}
From this, \eqref{418x},  $\sigma=\frac{\varepsilon'}{\varepsilon}$,
and the fact that $\delta\in(0,\varepsilon-\varepsilon']$,
we deduce that, for any $x,y,x',y'\in X$,
\begin{align*}
&\left|\left[P_lQ_k(x,y)-P_lQ_k(x',y)\right]
-\left[P_lQ_k(x,y')-P_lQ_k(x',y')\right]\right|\\
&\quad=\left|\left[P_lQ_k(x,y)-P_lQ_k(x',y)\right]
-\left[P_lQ_k(x,y')-P_lQ_k(x',y')\right]\right|^{\sigma}\\
&\qquad\times\left|\left[P_lQ_k(x,y)-P_lQ_k(x',y)\right]
-\left[P_lQ_k(x,y')-P_lQ_k(x',y')\right]\right|^{1-\sigma}\\
&\quad\lesssim2^{2l\varepsilon\sigma-(k-l)\varepsilon(1-\sigma)}
\left[\rho(x,x')\rho(y,y')\right]^{\varepsilon\sigma}
\left[\frac{1}{(V_\rho)_{2^{-l}}(x)}
+\frac{1}{(V_\rho)_{2^{-l}}(x')}\right]
\\
&\quad\leq2^{2\varepsilon'(k\land l)-|k-l|\delta}
\left[\rho(x,x')\rho(y,y')\right]^{\varepsilon'}
\left[\frac{1}{(V_\rho)_{2^{-(k\land l)}}(x)}
+\frac{1}{(V_\rho)_{2^{-(k\land l)}}(x')}\right].
\end{align*}
Thus, to complete the proof of \eqref{3.5},
it remains to show \eqref{418x}.
By Lemma~\ref{Pk}(v), we conclude that,
for any $x,y,x',y'\in X$,
\begin{align*}
&\left|\left[P_lQ_k(x,y)-P_lQ_k(x',y)\right]
-\left[P_lQ_k(x,y')-P_lQ_k(x',y')\right]\right|
\\
&\quad\leq\int_{\rho(z,y)\leq C_\rho\rho(y,y')}
\left|[P_l(x,z)-P_l(x',z)]-[P_l(x,y)-P_l(x',y)]\right|
\\
&\qquad\times\left|Q_k(z,y)-Q_k(z,y')\right|\,d\mu(z)
+\int_{\rho(z,y)>C_\rho\rho(y,y')}
\cdots
\\
&\quad=:Z_1(x,y,x',y')+Z_2(x,y,x',y').
\end{align*}
For $Z_1(x,y,x',y')$, from Definition~\ref{DefATI}(v), we infer that,
for any $x,y,x',y'\in X$,
\begin{align*}
&\int_{\rho(z,y)\leq C_\rho\rho(y,y')}
\left|Q_k(z,y)-Q_k(z,y')\right|\,d\mu(z)
\\
&\quad\leq\int_{X}\left[E_{2^{-k}}(z,y)
+E_{2^{-k+1}}(z,y)+E_{2^{-k}}(z,y')+
E_{2^{-k+1}}(z,y')\right]\,d\mu(z)
\lesssim 1.
\end{align*}
By this and Lemma~\ref{Pk}(iii),
we find that,
for any $x,y,x',y'\in X$,
\begin{align*}
Z_1(x,y,x',y')
&\lesssim\int_{\rho(z,y)\leq C_\rho\rho(y,y')}
\left[\frac{\rho(x,x')}{2^{-l}}\right]^\varepsilon
\left[\frac{\rho(z,y)}{2^{-l}}\right]^\varepsilon
\left[\frac{1}{(V_\rho)_{2^{-l}}(x)}
+\frac{1}{(V_\rho)_{2^{-l}}(x')}\right]
\\
&\quad\times
\left|Q_k(z,y)-Q_k(z,y')\right|\,d\mu(z)
\\
&\lesssim2^{2l\varepsilon}
\left[\rho(x,x')\rho(y,y')\right]^\varepsilon
\left[\frac{1}{(V_\rho)_{2^{-l}}(x)}
+\frac{1}{(V_\rho)_{2^{-l}}(x')}\right].
\end{align*}

For $Z_2(x,y,x',y')$,
on the one hand,
if $\rho(y,y')\ge C_12^{-k}$,
then, from Definition~\ref{D2.1}(iii) and
$\rho(z,y)>C_\rho\rho(y,y')$,
it follows that
$$
C_1C_\rho2^{-k}\leq C_\rho\rho(y,y')<\rho(z,y)\leq
C_\rho\max\left\{\rho(z,y'),\,\rho(y',y)\right\}
=C_\rho\rho(z,y'),
$$
which, together with Lemma~\ref{Pk}
and $\rho(z,y)>C_\rho\rho(y,y')\ge C_1C_\rho2^{-k}$,
further implies that $Q_k(z,y)=0=Q_k(z,y')$ and hence,
in this case,
$$
Z_2(x,y,x',y')=0.
$$
On the other hand, if $\rho(y,y')<C_12^{-k}$,
by Definition~\ref{D2.1}(iii)
and $\rho(z,y)>C_\rho\rho(y,y')$, we conclude
\begin{align}
\label{1223-1}
\rho(z,y)\leq C_\rho\max\left\{\rho(z,y'),\,\rho(y',y)\right\}=C_\rho\rho(z,y')
\end{align}
and
\begin{align}
\label{1223-2}
\rho(z,y')\leq C_\rho\max\left\{\rho(z,y),\,\rho(y,y')\right\}=C_\rho\rho(z,y).
\end{align}
If $Z_2(x,y,x',y')\neq0$, then, for any $z\in X$ satisfying
$\rho(z,y)>C_\rho\rho(y,y')$, $$Q_k(z,y)-Q_k(z,y')\neq0,$$
which, combined with \eqref{1223-1}
and \eqref{1223-2}, further implies that
$$
\rho(z,y)\leq C_1C_\rho2^{-k}.
$$
From this, (ii), (iii), and (iv) of Lemma~\ref{Pk},
and Remark~\ref{RemThmATI}(i),
we deduce that
\begin{align*}
Z_2(x,y,x',y')
&\lesssim\int_{C_1C_\rho2^{-k}\ge\rho(z,y)>C_\rho\rho(y,y')}
\left[\frac{\rho(x,x')}{2^{-l}}\right]^\varepsilon
\left[\frac{\rho(z,y)}{2^{-l}}\right]^\varepsilon
\\
&\quad\times\left[\frac{1}{(V_\rho)_{2^{-l}}(x)}
+\frac{1}{(V_\rho)_{2^{-l}}(x')}\right]
\left[\frac{\rho(y,y')}{2^{-k}}\right]^\varepsilon
\frac{1}{(V_\rho)_{2^{-k}}(z)}\,d\mu(z)\\
&\leq\left[\frac{\rho(x,x')}{2^{-l}}\right]^\varepsilon
\left[\frac{\rho(y,y')}{2^{-l}}\right]^\varepsilon
\left[\frac{1}{(V_\rho)_{2^{-l}}(x)}
+\frac{1}{(V_\rho)_{2^{-l}}(x')}\right]
\\
&\quad\times
\int_{C_1C_\rho2^{-k}>\rho(z,y)}
\left[\frac{\rho(z,y)}{2^{-k}}\right]^\varepsilon
\frac{1}{(V_\rho)_{2^{-k}}(z)}\,d\mu(z),
\end{align*}
which, together with \eqref{upperdoub} and Lemma~\ref{A3.2}(i),
further implies that
\begin{align*}
Z_2(x,y,x',y')
&\lesssim\left[\frac{\rho(x,x')}{2^{-l}}\right]^\varepsilon
\left[\frac{\rho(y,y')}{2^{-l}}\right]^\varepsilon
\left[\frac{1}{(V_\rho)_{2^{-l}}(x)}
+\frac{1}{(V_\rho)_{2^{-l}}(x')}\right]
\\
&\quad\times
\int_{C_1C_\rho2^{-k}>\rho(z,y)}
\left[\frac{\rho(z,y)}{2^{-k}}\right]^\varepsilon
\frac{1}{V_\rho(z,y)}\,d\mu(z)
\\
&\lesssim2^{2l\varepsilon}
\left[\rho(x,x')\rho(y,y')\right]^\varepsilon
\left[\frac{1}{(V_\rho)_{2^{-l}}(x)}
+\frac{1}{(V_\rho)_{2^{-l}}(x')}\right].
\end{align*}
This finishes the proof of \eqref{3.5}
and hence Lemma~\ref{A3.12}.
\end{proof}

In what follows,
let $(X,\mathbf{q},\mu)$ be a quasi-ultrametric space
of homogeneous type and
let $0<\varepsilon\preceq{\mathrm{ind\,}}(X,\mathbf{q})$
with ${\mathrm{ind\,}}(X,\mathbf{q})$ as in \eqref{index}.
Assume that $\{\mathcal{S}_{2^{-k}}\}_{k\in\mathbb{Z}}$
is a homogeneous approximation
of the identity of order $\varepsilon$
in Definition~\ref{DefHATI}.
For any $k\in\mathbb{Z}$, define
$$
\mathcal{D}_k:=\mathcal{S}_{2^{-k}}-\mathcal{S}_{2^{-k+1}}.
$$
To establish sharp continuous Calder\'on reproducing formulae,
using Proposition~\ref{A3.11}(ii) and following the idea
of the decomposition of identity from Coifman (see \cite{DaJoSe85} for details),
we rewrite that, for any fixed $N\in\mathbb{N}$ and $p\in(1,\infty)$,
\begin{align}\label{3.6}
\mathcal{I}
&=\sum_{k=-\infty}^{\infty}\mathcal{D}_k
=\sum_{k=-\infty}^\infty\left(\sum_{l=-\infty}^\infty
\mathcal{D}_{k+l}\right)\mathcal{D}_k
\nonumber\\
&=\sum_{k=-\infty}^{\infty}\mathcal{D}_k^N\mathcal{D}_k
+\sum_{k=-\infty}^{\infty}\sum_{|l|>N}\mathcal{D}_{k+l}\mathcal{D}_k
=:\mathcal{T}_N\index[symbols]{T@$\mathcal{T}_N$}
+\mathcal{R}_N\index[symbols]{R@$\mathcal{R}_N$}
\end{align}
in $L^p(X)$, where $\mathcal{I}$ is the identity operator
and
\begin{align}
\label{DkN7}
\mathcal{D}_k^N:=\sum_{|l|\leq N}\mathcal{D}_{k+l}
\end{align}
for any $k\in\mathbb{Z}$.

In Proposition~\ref{A3.16} below, we verify that
$\mathcal{T}_N^{-1}$ exists and is bounded on spaces
of test functions when $N$ is sufficiently large.
In a step towards this goal, we
first prove that $\mathcal{R}_N$ is bounded on $L^2(X)$
with a small operator norm.
To this end, we need the following
Cotlar--Stein lemma\index[words]{Cotlar--Stein lemma}
(see, for instance, \cite[pp.\,279--280]{Pag}
and \cite[p.\,14]{DaJoSe85}).

\begin{lemma}
\label{CotlarStein}
Let $(X,\mathbf{q},\mu)$ be a quasi-ultrametric space
of homogeneous type.
Assume that a sequence $\{\gamma(i)\}_{i\in\mathbb{Z}}$ in $(0,\infty)$
satisfies $\sum_{i\in\mathbb{Z}}\gamma(i)\in(0,\infty)$ and
$\{T_i\}_{i\in\mathbb{Z}}$ is a sequence of bounded linear operators
on $L^2(X)$ such that, for any $k,j\in\mathbb{Z}$,
\begin{align*}
\left\|T_j^*T_k\right\|_{L^2(X)\to L^2(X)}
\leq\left[\gamma(j-k)\right]^2
\ \ \text{and}\ \
\left\|T_jT_k^*\right\|_{L^2(X)\to L^2(X)}
\leq\left[\gamma(j-k)\right]^2,
\end{align*}
where $T^*_j$ and $T^*_k$ denote, respectively, the adjoint
operators of $T_j$ and $T_k$.
Then, for any $f\in L^2(X)$, the series $\sum_{i\in\mathbb{Z}}T_if$
converges to an element in $L^2(X)$, denoted by $Tf$.
Moreover, the operator $T$ is bounded on $L^2(X)$
and
$$
\|T\|_{L^2(X)\to L^2(X)}\leq\sum_{i\in\mathbb{Z}}\gamma(i).
$$
\end{lemma}

We are ready to show the boundedness
of $\mathcal{R}_N$ on $L^2(X)$.

\begin{lemma}\label{A3.13}
Let $(X,\mathbf{q},\mu)$ be a quasi-ultrametric space
of homogeneous type.
Assume that $0<\varepsilon\preceq{\mathrm{ind\,}}(X,\mathbf{q})$
with ${\mathrm{ind\,}}(X,\mathbf{q})$ as in \eqref{index}.
Fix $N\in\mathbb{N}$ and let
$\mathcal{R}_N$ be the same as in \eqref{3.6}.
Then, for any given $\varepsilon'\in(0,\varepsilon)$,
there exists a constant $C\in(0,\infty)$,
independent of $N$, such that, for any $f\in L^2(X)$,
$$
\|\mathcal{R}_Nf\|_{L^2(X)}\leq C2^{-N\varepsilon'}\|f\|_{L^2(X)}.
$$
\end{lemma}

\begin{proof}
Let $\rho\in\mathbf{q}$ be a regularized quasi-ultrametric as in Definition~\ref{D2.3}.
On the one hand, from Lemma~\ref{A3.12}(i), \eqref{upperdoub},
and Lemma~\ref{M} with $p=2$,
we infer that, for any $k,l\in\mathbb{Z}$ and $f\in L^2(X)$,
\begin{align}
\label{RN0}
\|\mathcal{D}_{k+l}\mathcal{D}_kf\|_{L^2(X)}
&\lesssim
2^{-|l|\varepsilon}\left\|
\frac{1}{(V_\rho)_{2^{-[k\land(k+l)]}}(\cdot)}\int_{\rho(\cdot,y)
\leq C'2^{-[k\land(k+l)]}}f(y)\,d\mu(y)\right\|_{L^2(X)}
\nonumber\\
&=2^{-|l|\varepsilon}\left\|
\frac{(V_\rho)_{C'2^{-[k\land(k+l)]}}(\cdot)}
{(V_\rho)_{2^{-[k\land(k+l)]}}(\cdot)}\fint_{\rho(\cdot,y)
\leq C'2^{-[k\land(k+l)]}}f(y)\,d\mu(y)\right\|_{L^2(X)}
\nonumber\\
&\lesssim2^{-|l|\varepsilon}\|\mathcal{M}_\rho f\|_{L^2(X)}
\lesssim2^{-|l|\varepsilon}\|f\|_{L^2(X)},
\end{align}
where $C'\in[1,\infty)$ is the same positive constant as in Lemma~\ref{A3.12}(i),
which further implies that, for any $k_1,k_2,l_1,l_2\in\mathbb{Z}$,
\begin{align}
\label{RN1}
&\left\|\mathcal{D}_{k_1+l_1}\mathcal{D}_{k_1}(\mathcal{D}_{k_2+l_2}
\mathcal{D}_{k_2})^*\right\|_{L^2(X)\to L^2(X)}
\nonumber\\
&\quad\leq\left\|\mathcal{D}_{k_1+l_1}\mathcal{D}_{k_1}\right\|_{L^2(X)\to L^2(X)}
\left\|(\mathcal{D}_{k_2+l_2}\mathcal{D}_{k_2})^*\right\|_{L^2(X)\to L^2(X)}
\lesssim2^{-(|l_1|+|l_2|)\varepsilon}.
\end{align}
On the other hand, by Theorem~\ref{corATI}(i)
(which still holds for both $\{\mathcal{D}_k\}_{k\in\mathbb{Z}}$
and $\{\mathcal{D}^*_k\}_{k\in\mathbb{Z}}$),
Lemma~\ref{Pk}(iv),
and \eqref{RN0},
we find that, for any $k_1,k_2,l_1,l_2\in\mathbb{Z}$,
\begin{align*}
&\left\|\mathcal{D}_{k_1+l_1}\mathcal{D}_{k_1}(\mathcal{D}_{k_2+l_2}
\mathcal{D}_{k_2})^*\right\|_{L^2(X)\to L^2(X)}
\nonumber\\
&\quad=
\left\|\mathcal{D}_{k_1+l_1}\mathcal{D}_{k_1}\mathcal{D}_{k_2}^*
\mathcal{D}_{k_2+l_2}^*\right\|_{L^2(X)\to L^2(X)}
\nonumber\\
&\quad\lesssim\left\|\mathcal{D}_{k_1}\mathcal{D}_{k_2}^*
\right\|_{L^2(X)\to L^2(X)}=\left\|\mathcal{D}_{k_1}\mathcal{D}_{k_2}
\right\|_{L^2(X)\to L^2(X)}
\lesssim2^{-|k_1-k_2|\varepsilon},
\end{align*}
which, combined with \eqref{RN1}, further implies that,
for any given $\theta\in(0,1)$,
$$
\left\|\mathcal{D}_{k_1+l_1}\mathcal{D}_{k_1}
(\mathcal{D}_{k_2+l_2}\mathcal{D}_{k_2})^*
\right\|_{L^2(X)\to L^2(X)}\lesssim
2^{-\varepsilon[|k_1-k_2|(1-\theta)+(|l_1|+|l_2|)\theta]}.
$$
Thus, for any $k_1,k_2\in\mathbb{Z}$,
\begin{align}
\label{RNL2}
&\left\|\left(\sum_{|l_1|>N}\mathcal{D}_{k_1+l_1}\mathcal{D}_{k_1}\right)
\left(\sum_{|l_2|>N}\mathcal{D}_{k_2+l_2}\mathcal{D}_{k_2}\right)^*
\right\|_{L^2(X)\to L^2(X)}\nonumber\\
&\quad\leq
\sum_{|l_1|>N}\sum_{|l_2|>N}\left\|\mathcal{D}_{k_1+l_1}\mathcal{D}_{k_1}
(\mathcal{D}_{k_2+l_2}\mathcal{D}_{k_2})^*\right\|_{L^2(X)\to L^2(X)}
\nonumber\\
&\quad\lesssim2^{-\varepsilon
[|k_1-k_2|(1-\theta)+2N\theta]}.
\end{align}
Similarly, we conclude that,
for any $k_1,k_2\in\mathbb{Z}$,
\begin{align*}
\left\|\left(\sum_{|l_1|>N}\mathcal{D}_{k_1+l_1}\mathcal{D}_{k_1}\right)^*
\left(\sum_{|l_2|>N}\mathcal{D}_{k_2+l_2}\mathcal{D}_{k_2}\right)
\right\|_{L^2(X)\to L^2(X)}
&\lesssim2^{-\varepsilon
[|k_1-k_2|(1-\theta)+2N\theta]}.
\end{align*}
From this, \eqref{RNL2}, and Lemma~\ref{CotlarStein}
with $T:=\mathcal{R}_N$,
$T_k:=\sum_{|l|>N}\mathcal{D}_{k+l}\mathcal{D}_k$ for any $k\in\mathbb{Z}$,
and $\gamma(j)$ for any $j\in\mathbb{Z}$
equaling a harmless constant multiple of
$2^{-\varepsilon[\frac{(1-\theta)|j|}{2}+N\theta]}$, we further deduce that
\begin{align}
\label{2210}
\|\mathcal{R}_N\|_{L^2(X)\to L^2(X)}
\leq\sum_{j\in\mathbb{Z}}\gamma(j)
\lesssim2^{-N\varepsilon\theta}
=2^{-N\varepsilon'},
\end{align}
where we have chosen $\theta:=\frac{\varepsilon'}{\varepsilon}\in(0,1)$.
This finishes the proof of Lemma~\ref{A3.13}.
\end{proof}

Next, to establish the boundedness of
$\mathcal{R}_N$ on spaces of test functions,
we first give some estimates for its kernel
$R_N$
and we begin with the following technical lemma.
Recall that, in this section,
we always assume that $\mu(X)=\infty$ and hence,
for any $x\in X$,
$R_\rho(x)=\infty$,
where $R_\rho(x)$ is the same as in \eqref{Rrho}.

\begin{lemma}\label{A3.14}
Let $(X,\rho,\mu)$ be a $(\kappa,n)$-ultra-RD-space
with $0<\kappa\leq n<\infty$.
Then, for any given $c_0\in(0,\infty)$,
there exists a positive constant $C$,
may depend on $c_0$, such that,
for any $x\in X$ and any constant $a\in\mathbb{R}$
satisfying $2^{-a}\ge c_0c_1r_\rho(x)$,
$$
\sum_{k\in\mathbb{Z}\cap(-\infty,a]}\frac{1}{(V_\rho)_{2^{-k}}(x)}
\leq C\frac{1}{(V_\rho)_{2^{-a}}(x)},
$$
where $c_1$ and $r_\rho(x)$ are
the same as, respectively, in Definition~\ref{RD} and \eqref{rrho}.
\end{lemma}

\begin{proof}
By Lemma~\ref{kst.1} with $\theta_1:=c_0$ and $\theta_2:=1$, we find that,
for any $k\in\mathbb{Z}\cap(-\infty,a]$,
\begin{align*}
(V_\rho)_{2^{-k}}(x)
&=(V_\rho)_{2^{-(k-a)}2^{-a}}(x)
\gtrsim2^{-\kappa(k-a)}(V_\rho)_{2^{-a}}(x).
\end{align*}
Therefore,
\begin{align*}
\sum_{k\in\mathbb{Z}\cap(-\infty,a]}\frac{1}{(V_\rho)_{2^{-k}}(x)}
\lesssim\frac{1}{(V_\rho)_{2^{-a}}(x)}
\sum_{k\in\mathbb{Z}\cap(-\infty,a]}2^{\kappa(k-a)}
\lesssim\frac{1}{(V_\rho)_{2^{-a}}(x)},
\end{align*}
where the implicit positive constants are independent of both $x$ and $a$,
but may depend on $c_0$ due to Lemma~\ref{kst.1}.
This finishes the proof of Lemma~\ref{A3.14}.
\end{proof}

\begin{remark}
Although the approximation of the identity in Definition~\ref{DefATI}
has a bounded support,
we still need the reverse doubling condition to obtain
Lemma~\ref{A3.14}, whereas He et al. \cite[Lemma 4.9]{HLYY}
used the exponential decay factor in their approximation of the identity
to get rid of the dependence on the reverse doubling condition
of the equipped measure.
\end{remark}

\begin{lemma}\label{A3.15}
Let $(X,\mathbf{q},\mu)$ be an ultra-RD-space and
$N\in\mathbb{N}$.
Assume that
$0<\varepsilon\preceq{\mathrm{ind\,}}(X,\mathbf{q})$
with ${\mathrm{ind\,}}(X,\mathbf{q})$ as in \eqref{index}.
Let $\mathcal{R}_N$ be the same as in \eqref{3.6}
and denote its kernel by $R_N$. Then,
for any $\varepsilon'\in(0,\varepsilon)$ and
$\delta\in(0,\varepsilon-\varepsilon']$,
$R_N$ satisfies \eqref{CZk1},
\eqref{CZk2}, and \eqref{CZO1}
with $C_T$ replaced by $\widetilde{C_T}2^{-\delta N}$,
where $\widetilde{C_T}$ is a positive constant independent of $N$.
\end{lemma}

\begin{proof}
Let $\rho\in\mathbf{q}$ be a regularized quasi-ultrametric.
From Lemma~\ref{A3.12}, we infer that there exists a positive constant $C$
(which we assume to be at least 1)
such that, for any $x,y\in X$ satisfying $\rho(x,y)\ge C2^{-(k\land l)}$
with $k,l\in\mathbb{Z}$,
$D_kD_l(x,y)=0$.

We first prove that $R_N$ satisfies \eqref{CZk1}.
Let $x,y\in X$ satisfy $x\neq y$
and choose $k_0\in\mathbb{Z}$ satisfying
$C^{-1}2^{-k_0}<\rho(x,y)\le C2^{-k_0+1}$.
To estimate $R_N(x,y)$, we separately estimate
\begin{align*}
R^1_N(x,y):=\sum_{k=-\infty}^{\infty}\sum_{l=N+1}^{\infty}
D_{k+l}D_k(x,y)
\end{align*}
and
\begin{align*}
R^2_N(x,y):=\sum_{k=-\infty}^{\infty}\sum_{l=-\infty}^{-N-1}
D_{k+l}D_k(x,y),
\end{align*}
which are defined formally.
For $R_N^1(x,y)$,
by $l\in\mathbb{Z}_+$ and \eqref{3.2}, we conclude that,
for any $k\in\mathbb{Z}$ satisfying $\rho(x,y)\leq C2^{-k}$,
$$
\left|D_{k+l}D_k(x,y)\right|
\lesssim2^{-l\varepsilon}
\frac{1}{(V_\rho)_{2^{-k}}(x)}
$$
and, for any $k\in\mathbb{Z}$ satisfying $\rho(x,y)>C2^{-k}$,
$D_{k+l}D_k(x,y)=0$.
From this, Lemma~\ref{A3.14} with both $a:=k_0-1$ and $c_0:=C^{-1}$
[here, we used the fact that
$$
2^{-a}=2^{-k_0+1}\ge\frac{\rho(x,y)}{C}\ge\frac{c_1r_\rho(x)}{C}
$$
since $\rho(x,y)\leq C2^{-k_0+1}$,
where $c_1\in(0,1]$ and $r_\rho(x)$ are
the same as, respectively, in Definition~\ref{RD} and \eqref{rrho}],
and \eqref{doublingcond}, we deduce that
\begin{align}\label{RN-1}
\left|R^1_N(x,y)\right|
&\leq\sum_{k=-\infty}^{\infty}\sum_{l=N+1}^{\infty}
|D_{k+l}D_k(x,y)|
\lesssim\sum_{l=N+1}^{\infty}
\sum_{k=-\infty}^{k_0-1}
2^{-l\varepsilon}\frac{1}{(V_\rho)_{2^{-k}}(x)}
\nonumber\\
&\lesssim\sum_{l=N+1}^{\infty}
2^{-l\varepsilon}\frac{1}{(V_\rho)_{2^{-k_0+1}}(x)}
\lesssim 2^{-\varepsilon N}\frac{1}{V_\rho(x,y)}.
\end{align}
For $R_N^2(x,y)$,
by $l\in\mathbb{Z}
\setminus\mathbb{Z}_+$ and \eqref{3.2},
we find that, for any $k\in\mathbb{Z}$
satisfying
$\rho(x,y)\leq C2^{-(k+l)}$,
$$
|D_{k+l}D_k(x,y)|
\lesssim2^{l\varepsilon}\frac{1}
{(V_\rho)_{2^{-(k+l)}}(x)}
$$
and, for any $k\in\mathbb{Z}$
satisfying
$\rho(x,y)>C2^{-(k+l)}$,
$D_{k+l}D_k(x,y)=0$.
From this, Lemma~\ref{A3.14} with both $a:=k_0-1$ and $c_0:=C^{-1}$,
and \eqref{doublingcond}, we infer that
\begin{align}\label{RN-2}
\left|R^2_N(x,y)\right|
&\leq\sum_{k=-\infty}^{\infty}\sum_{l=-\infty}^{-N-1}
\left|D_{k+l}D_k(x,y)\right|
\lesssim\sum_{l=-\infty}^{-N-1}\sum_{k=-\infty}^{k_0-l-1}
2^{l\varepsilon}\frac{1}{(V_\rho)_{2^{-(k+l)}}(x)}
\nonumber\\
&\lesssim\sum_{l=-\infty}^{-N-1}
2^{l\varepsilon}\frac{1}{(V_\rho)_{2^{-k_0+1}}(x)}
\lesssim 2^{-\varepsilon N}\frac{1}{V_\rho(x,y)}.
\end{align}
In particular, the series defining $R^1_N(x,y)$
and $R^2_N(x,y)$ converge absolutely
and, consequently, we can write
\begin{align}\label{5.3.x1}
R_N(x,y)=R^1_N(x,y)+R^2_N(x,y).
\end{align}
By this, \eqref{RN-1}, \eqref{RN-2},
and $\delta<\varepsilon$,
we conclude that there exists a positive constant $C_1$,
independent of $N$, such that,
for any $x,y\in X$ with $x\neq y$,
\begin{align}\label{RNCZO1}
\left|R_N(x,y)\right|
&\leq\left|R^1_N(x,y)\right|+\left|R^2_N(x,y)\right|\nonumber\\
&\leq C_12^{-\varepsilon N}\frac{1}{V_\rho(x,y)}
\leq C_12^{-\delta N}\frac{1}{V_\rho(x,y)}.
\end{align}
Thus, $R_N$ satisfies \eqref{CZk1}
with $C_T$ replaced by $C_12^{-\delta N}$.

Now, we show that $R_N$ satisfies \eqref{CZk2}.
Note that, if $l\in\mathbb{Z}\cap[N+1,\infty)$ and $k\in\mathbb{Z}$,
then $k\land(k+l)=k$.
Fix $x,x',y\in X$ with
$x\neq y$ and $\rho(x,x')
\leq\frac{\rho(x,y)}{2C_\rho}$, and choose
$k_1\in\mathbb{Z}$ satisfying that
\begin{align}\label{k01}
C^{-1}C_\rho 2^{-k_1}<\rho(x,y)\leq CC_\rho 2^{-k_1+1}.
\end{align}
From Definition~\ref{D2.1}(iii),
it follows that
$$
\rho(x,y)\leq C_\rho
\max\{\rho(x,x'),\,\rho(x',y)\}=C_\rho\rho(x',y).
$$
By this and \eqref{k01}, we obtain $\rho(x',y)>C2^{-k_1}$.
From this and Lemma~\ref{A3.12}(i), we deduce that,
for any $l\in\mathbb{Z}\cap[N+1,\infty)$ and $k\in\mathbb{Z}\cap[k_1,\infty)$,
\begin{align}\label{DklDk0}
D_{k+l}D_k(x,y)=0=D_{k+l}D_k(x',y).
\end{align}
Note that, by \eqref{k01}, we find that,
for any $k\in\mathbb{Z}\cap(-\infty,k_1)$, $\rho(x,y)\leq CC_\rho2^{-k}.$
From this and \eqref{DklDk0}, we infer that
\begin{align}\label{RN1x}
\left|R_N^1(x,y)-R_N^1(x',y)\right|
&\leq\sum_{l=N+1}^\infty\sum_{k=-\infty}^{k_1-1}
\left|D_{k+l}D_k(x,y)-D_{k+l}D_k(x',y)\right|
\nonumber\\
&\lesssim\sum_{l=N+1}^\infty\sum_{k=-\infty}^{k_1-1}2^{-l\delta}
\left[\frac{\rho(x,x')}{2^{-k}}\right]^{\varepsilon'}
\frac{1}{(V_\rho)_{2^{-k}}(y)}
\quad\text{by Lemma~\ref{A3.12}(iii)}
\nonumber\\
&\lesssim\sum_{l=N+1}^\infty\sum_{k=-\infty}^{k_1-1}2^{-l\delta}
\left[\frac{\rho(x,x')}{\rho(x,y)}\right]
^{\varepsilon'}\frac{1}{(V_\rho)_{2^{-k}}(x)}
\quad\text{by Lemma~\ref{A3.2}(vii)}
\nonumber\\
&\lesssim\left[\frac{\rho(x,x')}{\rho(x,y)}\right]
^{\varepsilon'}\frac{1}{(V_\rho)_{2^{-k_1+1}}(x)}
\sum_{l=N+1}^\infty2^{-l\delta}
\nonumber\\
&\qquad\qquad\text{by Lemma~\ref{A3.14} with $a:=k_1-1$ and $c_0:=(CC_\rho)^{-1}$}
\nonumber\\
&\lesssim2^{-\delta N}\left[\frac{\rho(x,x')}{\rho(x,y)}\right]
^{\varepsilon'}\frac{1}{V_\rho(x,y)}
\quad\text{by \eqref{doublingcond}}.
\end{align}
By an argument similar to that used
in the proof of \eqref{RN1x}, we conclude that,
for any $x,x',y\in X$ with
$x\neq y$ and $\rho(x,x')\leq\frac{\rho(x,y)}{2C_\rho}$,
\begin{align*}
\left|R_N^2(x,y)-R_N^2(x',y)\right|\lesssim
2^{-\delta N}\left[\frac{\rho(x,x')}
{\rho(x,y)}\right]^{\varepsilon'}
\frac{1}{V_\rho(x,y)},
\end{align*}
which, together with \eqref{5.3.x1} and \eqref{RN1x},
further implies that,
\begin{align}
\label{RNxx1}
\left|R_N(x,y)-R_N(x',y)\right|&\leq
\left|R_N^1(x,y)-R_N^1(x',y)\right|+
\left|R_N^2(x,y)-R_N^2(x',y)\right|
\nonumber\\
&\lesssim2^{-\delta N}
\left[\frac{\rho(x,x')}{\rho(x,y)}\right]
^{\varepsilon'}\frac{1}{V_\rho(x,y)}.
\end{align}
From an argument similar to that used
in the proof of \eqref{RNxx1},
we deduce that, for any $x,y,y'\in X$ with
$x\neq y$ and $\rho(y,y')\leq\frac{\rho(x,y)}{2C_\rho}$,
$$
|R_N(x,y)-R_N(x,y')|
\lesssim 2^{-\delta N}\left[\frac{\rho(y,y')}
{\rho(x,y)}\right]^\varepsilon
\frac{1}{V_\rho(x,y)},
$$
which, combined with \eqref{RNxx1},
proves that there exists a positive constant $C_2$,
independent of $N$, such that $R_N$ satisfies \eqref{CZk2}
with $C_T$ replaced by $C_22^{-\delta N}$.

Next, we show that $R_N$ satisfies \eqref{CZO1}.
For any $l\in\mathbb{Z}\cap[N+1,\infty)$
and $k\in\mathbb{Z}$, we have $k\land(k+l)=k$.
Fix $x,x',y,y'\in X$ with
$x\neq y$, $\rho(x,x')\leq\frac{\rho(x,y)}{4C_\rho^2}$,
and $\rho(y,y')\leq\frac{\rho(x,y)}{4C_\rho^2}$,
and choose $k_2\in\mathbb{Z}$ satisfying
\begin{align}\label{k0}
C^{-1}C^2_\rho 2^{-k_2}<\rho(x,y)\leq CC^2_\rho 2^{-k_2+1}.
\end{align}
By Definition~\ref{D2.1}(iii),
we find that
\begin{align*}
\rho(x,y)\leq C_\rho\max
\left\{\rho(x,x'),\,\rho(x',y)\right\}=C_\rho\rho(x',y),
\end{align*}
\begin{align*}
\rho(x,y)\leq C_\rho\max
\left\{\rho(x,y'),\,\rho(y',y)\right\}=C_\rho\rho(x,y'),
\end{align*}
and
\begin{align*}
\rho(x,y)\leq C_\rho\rho(x',y)
\leq C^2_\rho\max\left\{\rho(x',y'),\,\rho(y',y)\right\}
=C^2_\rho\rho(x',y').
\end{align*}
From this, \eqref{k0},
and Lemma~\ref{A3.12}(i),
we infer that, for any $l\in\mathbb{Z}
\cap[N+1,\infty)$ and $k\in\mathbb{Z}\cap[k_2,\infty)$,
\begin{align}\label{DklDk1}
D_{k+l}D_k(x,y)=D_{k+l}D_k(x',y)
=D_{k+l}D_k(x,y')=D_{k+l}D_k(x',y')=0.
\end{align}
By \eqref{k0} and $\rho(x,x')
\leq\frac{\rho(x,y)}{4C_\rho^2}$, we conclude that,
for any $k\in\mathbb{N}\cap(-\infty,k_2)$, $\rho(x,y)
\lesssim 2^{-k}$ and $\rho(x,x')\lesssim2^{-k}$,
which, together with Lemma~\ref{A3.2}(vii), further implies that
$$\frac{1}{(V_\rho)_{2^{-k}}(x)}
\sim\frac{1}{(V_\rho)_{2^{-k}}(x')}.$$
From this, \eqref{DklDk1}, and \eqref{3.5},
we deduce that
\begin{align*}
&\left|\left[R_N^1(x,y)-R_N^1(x',y)\right]
-\left[R_N^1(x,y')-R_N^1(x',y')\right]\right|
\nonumber\\
&\quad\leq\sum_{l=N+1}^\infty\sum_{k=-\infty}^{k_2-1}
\left|\left[D_{k+l}D_k(x,y)-D_{k+l}D_k(x',y)\right]\right.\nonumber\\
&\qquad\left.-\left[D_{k+l}D_k(x,y')-D_{k+l}D_k(x',y')\right]\right|
\nonumber\\
&\quad\lesssim\sum_{l=N+1}^\infty
\sum_{k=-\infty}^{k_2-1}2^{-l\delta}
\left[\frac{\rho(x,x')}{2^{-k}}\right]^{\varepsilon'}
\left[\frac{\rho(y,y')}{2^{-k}}\right]^{\varepsilon'}
\left[\frac{1}{(V_\rho)_{2^{-k}}(x)}+
\frac{1}{(V_\rho)_{2^{-k}}(x')}\right]
\nonumber\\
&\quad\lesssim
\left[\frac{\rho(x,x')}{\rho(x,y)}\right]^{\varepsilon'}
\left[\frac{\rho(y,y')}{\rho(x,y)}\right]^{\varepsilon'}
\sum_{l=N+1}^\infty2^{-l\delta}
\sum_{k=-\infty}^{k_2-1}
\frac{1}{(V_\rho)_{2^{-k}}(x)},
\end{align*}
which further implies that
\begin{align}
\label{RN1xx}
&\left|\left[R_N^1(x,y)-R_N^1(x',y)\right]
-\left[R_N^1(x,y')-R_N^1(x',y')\right]\right|
\nonumber\\
&\quad\lesssim
\left[\frac{\rho(x,x')}{\rho(x,y)}\right]
^{\varepsilon'}
\left[\frac{\rho(y,y')}{\rho(x,y)}\right]
^{\varepsilon'}\frac{1}{(V_\rho)_{2^{-k_2+1}}(x)}
\sum_{l=N+1}^\infty2^{-l\delta}
\nonumber\\
&\qquad\qquad\text{by Lemma~\ref{A3.14} with
both $a:=k_2-1$ and $c_0:=(CC_\rho^2)^{-1}$}
\nonumber\\
&\quad\lesssim2^{-\delta N}\left[\frac{\rho(x,x')}{\rho(x,y)}
\right]^{\varepsilon'}
\left[\frac{\rho(y,y')}{\rho(x,y)}\right]
^{\varepsilon'}\frac{1}{V_\rho(x,y)}
\quad\text{by \eqref{doublingcond}}.
\end{align}
By an argument similar to that used
in the estimation of \eqref{RN1xx}, we find that,
for any $x,x',y,y'\in X$ satisfying that
$x\neq y$, $\rho(x,x')\leq\frac{\rho(x,y)}{4C_\rho^2}$,
and $\rho(y,y')\leq\frac{\rho(x,y)}{4C_\rho^2}$,
\begin{align*}
&\left|\left[R_N^2(x,y)-R_N^2(x',y)\right]
-\left[R_N^2(x,y')-R_N^2(x',y')\right]\right|
\\
&\quad\lesssim2^{-\delta N}\left[\frac{\rho(x,x')}
{\rho(x,y)}\right]^{\varepsilon'}
\left[\frac{\rho(y,y')}{\rho(x,y)}\right]
^{\varepsilon'}\frac{1}{V_\rho(x,y)},
\end{align*}
which, combined with \eqref{RN1xx}, further implies that
\begin{align*}
&\left|\left[R_N(x,y)-R_N(x',y)\right]
-\left[R_N(x,y')-R_N(x',y')\right]\right|
\\
&\quad\leq
\left|\left[R_N^1(x,y)-R_N^1(x',y)\right]
-\left[R_N^1(x,y')-R_N^1(x',y')\right]\right|
\\
&\qquad+\left|\left[R_N^2(x,y)-R_N^2(x',y)\right]
-\left[R_N^2(x,y')-R_N^2(x',y')\right]\right|
\\
&\quad\lesssim2^{-\delta N}
\left[\frac{\rho(x,x')}{\rho(x,y)}\right]^{\varepsilon'}
\left[\frac{\rho(y,y')}{\rho(x,y)}\right]
^{\varepsilon'}\frac{1}{V_\rho(x,y)}.
\end{align*}
From this, we infer that there exists a positive constant $C_3$,
independent of $N$, such that $R_N$ satisfies \eqref{CZO1}
with $C_T$ replaced by $C_32^{-\delta N}$.
Let
$$
\widetilde{C_T}:=\max\left\{C_1,\,C_2,\,C_3\right\}.
$$
Then $\widetilde{C_T}$ is independent of $N$
and $R_N$ satisfies \eqref{CZk1}, \eqref{CZk2}, and \eqref{CZO1}
with $C_T$ replaced by $\widetilde{C_T}2^{-\delta N}$.
This finishes the proof of Lemma~\ref{A3.15}.
\end{proof}

As mentioned in \cite[p.\,313]{HLYY}, it is not clear whether
$R_N$, the kernel of $\mathcal{R}_N$,
satisfies both (ii) and (iii) of Theorem~\ref{CZOonGO}.
To overcome this difficulty, we borrow some ideas from \cite{HLYY}.
For any $M\in\mathbb{N}\cap(N,\infty)$,
define\index[symbols]{R@$R_{N,M}$}
\begin{equation}\label{defRNM}
\mathcal{R}_{N,M}:=\sum_{|k|\leq M}\sum_{N<|l|\leq M}
\mathcal{D}_{k+l}\mathcal{D}_k.
\end{equation}
It is easy to infer that $R_{N,M}$, the kernel of $\mathcal{R}_{N,M}$,
satisfies both (ii) and (iii) of Theorem~\ref{CZOonGO} with $c_0=0$.

Indeed, to verify that $\mathcal{R}_{N,M}$ satisfies
Theorem~\ref{CZOonGO}(ii),
from Theorem~\ref{corATI}(iv) and the definition of $\mathcal{D}_k$,
we deduce that, for any $k\in\mathbb{Z}$
and $\beta\in(0,\varepsilon)$ with $\varepsilon$ as in Theorem~\ref{CZOonGO},
\begin{align*}
\left\|\mathcal{D}_k\right\|_{C^\beta(X)\to C^\beta(X)}
&\leq\left\|\mathcal{S}_{2^{-k}}\right\|_{C^\beta(X)\to C^\beta(X)}
+\left\|\mathcal{S}_{2^{-k+1}}\right\|_{C^\beta(X)\to C^\beta(X)}\\
&\lesssim\sup_{t\in(0,\infty)}
\left\|\mathcal{S}_t\right\|_{C^\beta(X)\to C^\beta(X)},
\end{align*}
which further implies that,
for any $g\in C^\beta(X)$,
\begin{align*}
\left\|\mathcal{R}_{N,M}f\right\|_{C^\beta(X)}
\leq\sum_{|k|\leq M}\sum_{N<|l|\leq M}
\left\|\mathcal{D}_{k+l}\mathcal{D}_kf\right\|_{C^\beta(X)}
\lesssim\sum_{|k|\leq M}
\left\|\mathcal{D}_kf\right\|_{C^\beta(X)}
\lesssim\left\|f\right\|_{C^\beta(X)},
\end{align*}
where the implicit positive constants depend only on $N$, $M$,
and $\sup_{t\in(0,\infty)}
\|\mathcal{S}_t\|_{C^\beta(X)\to C^\beta(X)}$.
Thus, for any $g\in C^\beta(X)$ and $x\in X$,
$$
\mathcal{R}_{N,M}g(x)=\int_XR_{N,M}(x,y)g(y)\,d\mu(y)
$$
is well defined and hence $\mathcal{R}_{N,M}$ satisfies
Theorem~\ref{CZOonGO}(ii).

To verify that $R_{N,M}(x,y)$ satisfies
Theorem~\ref{CZOonGO}(iii) with $c_0=0$,
by Fubini's theorem and Lemma~\ref{Pk}(v),
we conclude that,
for any $x\in X$,
\begin{align*}
\int_XR_{N,M}(x,y)\,d\mu(y)
&=\sum_{|k|\leq M}\sum_{N<|l|\leq M}\int_X\int_X
D_{k+l}(x,z)D_k(z,y)\,d\mu(z)\,d\mu(y)\nonumber\\
&=\sum_{|k|\leq M}\sum_{N<|l|\leq M}\int_X\left[\int_X
D_k(z,y)\,d\mu(y)\right]D_{k+l}(x,z)\,d\mu(z)
=0
\end{align*}
and hence
\begin{align*}
\int_XR_{N,M}(y,x)\,d\mu(y)=0.
\end{align*}
From this and Fubini's theorem, we infer that,
for any $f\in\mathring{\mathcal{G}}(x_1,r,\beta,\gamma)$ with
$x_1\in X$, $r\in(0,\infty)$, and $\beta,\gamma\in(0,\varepsilon)$,
$$
\int_XR_{N,M}f(x)\,d\mu(x)
=\int_X\left[\int_XR_{N,M}(x,y)\,d\mu(x)\right]f(y)\,d\mu(y)
=0.
$$
Note that Lemmas~\ref{A3.13} and \ref{A3.15} with $\mathcal{R}_N$ replaced by
$\mathcal{R}_{N,M}$ also hold with all the constants
involved independent of $M$.
By this and Theorem~\ref{CZOonGO},
we find that,
for any $f\in\mathring{\mathcal{G}}(x_1,r,\beta,\gamma)$ with
$x_1\in X$, $r\in(0,\infty)$, and $\beta,\gamma\in(0,\varepsilon)$,
$\mathcal{R}_{N,M}f\in\mathring{\mathcal{G}}(x_1,r,\beta,\gamma)$ and
\begin{align}\label{BRNM}
\left\|\mathcal{R}_{N,M}f\right\|_{\mathring{\mathcal{G}}(x_1,r,\beta,\gamma)}
&\lesssim\left[C_{\mathcal{R}_{N,M}}
+\|\mathcal{R}_{N,M}\|_{L^2(X)\rightarrow L^2(X)}\right]
\|f\|_{\mathring{\mathcal{G}}(x_1,r,\beta,\gamma)}\nonumber\\
&\lesssim 2^{(\varepsilon'-\varepsilon)N}
\|f\|_{\mathring{\mathcal{G}}(x_1,r,\beta,\gamma)},
\end{align}
where $\varepsilon'\in(\max\{\beta,\,\gamma\},\varepsilon)$
and the implicit positive constants are independent of both $N$ and $M$.

To obtain the boundedness on $\mathring{\mathcal{G}}(x_1,r,\beta,\gamma)$ of
$\mathcal{R}_N$ from that of $\mathcal{R}_{N,M}$ in \eqref{BRNM},
we need to consider the relationship between $\mathcal{R}_Nf$ and
$\mathcal{R}_{N,M}f$ with $f\in\mathring{\mathcal{G}}(x_1,r,\beta,\gamma)$.
To this end, we need some technical results which are used
in the proof of sharp homogeneous continuous Calder\'on reproducing formulae
in $L^p(X)$ with $p\in(1,\infty)$;
see Theorem~\ref{A3.19} below. To begin with, using Lemmas~\ref{CGL} and~\ref{density},
we obtain the following conclusion.

\begin{lemma}\label{RNRNM}
Let $(X,\mathbf{q},\mu)$ be an ultra-RD-space,
where $\mu$ is a Borel-semiregular measure on $X$.
Let $0<\beta,\gamma<\varepsilon\preceq{\mathrm{ind\,}}(X,\mathbf{q})$
and
$f\in{\mathcal G}(x_1,r,\beta,\gamma)$ with $x_1\in X$ and $r\in(0,\infty)$.
Fix $N\in\mathbb{N}$ and, for any $M\in\mathbb{N}$, let
$\mathcal{R}_N$ and
$\mathcal{R}_{N,M}$
be the same as, respectively, in \eqref{3.6} and \eqref{defRNM}.
Then the following statements hold:
\begin{enumerate}
\item[\textup{(i)}]
$\lim_{M\to\infty}\|\mathcal{R}_{N,M}f-\mathcal{R}_Nf\|_{L^2(X)}=0$.

\item[\textup{(ii)}]
The function sequence $\{\mathcal{R}_{N,M}f\}_{M=N+1}^\infty$
converges locally uniformly to some element,
denoted by $\widetilde{\mathcal{R}_N}f$,
where $\widetilde{\mathcal{R}_N}f$ differs from
$\mathcal{R}_Nf$ at most on a set of
$\mu$-measure 0.

\item[\textup{(iii)}]
The operator $\widetilde{\mathcal{R}_N}$
can be uniquely extended from ${\mathcal G}(x_1,r,\beta,\gamma)$ to $L^2(X)$,
where the extension operator coincides with $\mathcal{R}_N$
in the sense that, for any
$f\in{\mathcal G}(x_1,r,\beta,\gamma)$ and $\mu$-almost every
$x\in X$,
$$
\lim_{M\to\infty}\mathcal{R}_{N,M}f(x)
=\widetilde{\mathcal{R}_N}f(x)=\mathcal{R}_Nf(x).
$$
\end{enumerate}
\end{lemma}

\begin{proof}
Let $\rho\in\mathbf{q}$ be a regularized quasi-ultrametric
as in Definition~\ref{D2.3}.
We first prove (i). From Lemma~\ref{CGL} with $p=2$,
it suffices to show that, for any $f\in L^2(X)$,
\begin{align}\label{RNMRNL2}
\lim_{M\to\infty}\|\mathcal{R}_{N,M}f-\mathcal{R}_Nf\|_{L^2(X)}=0.
\end{align}
For any $M\in\mathbb{N}\cap(N,\infty)$, we write
\begin{align}\label{RNMRN0}
\mathcal{R}_N-\mathcal{R}_{N,M}
=\sum_{|k|\leq M}\sum_{|l|>M}\mathcal{D}_{k+l}\mathcal{D}_k
+\sum_{|k|>M}\sum_{|l|>N}\mathcal{D}_{k+l}\mathcal{D}_k.
\end{align}
On the one hand,
by an argument similar to that
used in the estimation of \eqref{2210}
and Lemma~\ref{CotlarStein}, we obtain
$$
\left\|\sum_{|k|\leq M}\sum_{|l|>M}
\mathcal{D}_{k+l}\mathcal{D}_k\right\|_{L^2(X)\to L^2(X)}
\lesssim2^{-M\varepsilon'},
$$
where the implicit positive constant is independent of $M$.
Therefore,
\begin{align}\label{RNMRN1}
\lim_{M\to\infty}\left\|\sum_{|k|\leq M}
\sum_{|l|>M}\mathcal{D}_{k+l}\mathcal{D}_kf\right\|_{L^2(X)}=0.
\end{align}
On the other hand, from Lemma~\ref{A3.13}, we deduce that
$$
\mathcal{R}_Nf=\sum_{k=-\infty}^\infty\sum_{|l|>N}
\mathcal{D}_{k+l}\mathcal{D}_kf
=\lim_{M\to\infty}\sum_{|k|\leq M}\sum_{|l|>N}
\mathcal{D}_{k+l}\mathcal{D}_kf
$$
in $L^2(X)$,
which implies that
\begin{align}\label{RNMRN2}
\lim_{M\to\infty}\left\|\sum_{|k|>M}\sum_{|l|>N}
\mathcal{D}_{k+l}\mathcal{D}_kf\right\|_{L^2(X)}=0.
\end{align}
By \eqref{RNMRN0}, \eqref{RNMRN1}, and \eqref{RNMRN2},
we obtain \eqref{RNMRNL2}. This finishes the proof of (i).

Now, we prove (ii). To show that
$\{\mathcal{R}_{N,M}f\}_{M=N+1}^\infty$
is a locally uniformly convergent
sequence, we only need to prove that, for any $x\in X$ and $r\in(0,\infty)$,
there exists a positive sequence $\{c_{k,l}\}_{k,l\in\mathbb{Z}}$
such that, for any $k,l\in\mathbb{Z}$,
\begin{align}\label{ck,l}
\sup_{y\in B_\rho(x,r)}\left|\mathcal{D}_{k+l}
\mathcal{D}_kf(y)\right|\leq c_{k,l}
\ \ \text{and}\ \
\sum_{k\in\mathbb{Z}}\sum_{l\in\mathbb{Z}}c_{k,l}<\infty.
\end{align}
To this end, we first claim that,
for any $x\in X$, $r\in(0,\infty)$,
and $k,l\in\mathbb{Z}$,
\begin{equation}\label{ckl}
\sup_{y\in B_\rho(x,r)}\left|\mathcal{D}_{k+l}
\mathcal{D}_kf(y)\right|\lesssim
\begin{cases}
\displaystyle
2^{-\varepsilon|l|}\frac{1}
{(V_\rho)_{2^{-[k\land(k+l)]}}(x)}\quad
&\text{if }2^{-[k\land(k+l)]}\ge r,\\
\displaystyle
2^{-\varepsilon|l|}\frac{1}{(V_\rho)_r(x_1)}
\left\{\frac{2^{-[k\land(k+l)]}}{r}\right\}^\beta
&\text{if }2^{-[k\land(k+l)]}<r.
\end{cases}
\end{equation}
Let $\{c_{k,l}\}_{k,l\in\mathbb{Z}}$ be as
in the right-hand side of \eqref{ckl}.
From this, Tonelli's theorem, and
Lemma~\ref{A3.14} with
both $a:=\log_{\frac{1}{2}}r$ and $c_0:=1$,
we infer that
\begin{align*}
\sum_{k\in\mathbb{Z}}\sum_{l\in\mathbb{Z}}
c_{k,l}
&=\sum_{l\in\mathbb{Z}}2^{-\varepsilon|l|}
\left[\sum_{k\land(k+l)\leq\log_{\frac{1}{2}}r}
\frac{1}{(V_\rho)_{2^{-[k\land(k+l)]}}(x)}\right.
\\
&\quad\left.+\sum_{k\land(k+l)>\log_{\frac{1}{2}}r}
\frac{1}{(V_\rho)_r(x_1)}
\left\{\frac{2^{-[k\land(k+l)]}}{r}\right\}^\beta\right]
\\
&\sim\sum_{m\in\mathbb{Z}\cap(-\infty,\log_{\frac{1}{2}}r]}
\frac{1}{(V_\rho)_{2^{-m}}(x)}+
\sum_{m\in\mathbb{Z}\cap(\log_{\frac{1}{2}}r,\infty)}
\frac{1}{(V_\rho)_r(x_1)}
\left(\frac{2^{-m}}{r}\right)^\beta
\\
&\lesssim\frac{1}{(V_\rho)_r(x)}+\frac{1}{(V_\rho)_r(x_1)}
<\infty
\end{align*}
which is precisely \eqref{ck,l}.

To show the claim in \eqref{ckl},
by the symmetry, we only need to consider the case
$l\in\mathbb{N}$. When $2^{-[k\land(k+l)]}=2^{-k}\geq r$,
from both (ii) and (vii) of Lemma~\ref{A3.2},
we deduce that, for any $y\in B_\rho(x,r)$,
\begin{align*}
\left|\mathcal{D}_{k+l}\mathcal{D}_kf(y)\right|
&\leq\int_X
\left|D_{k+l}D_k(y,z)f(z)\right|\,d\mu(z)
\\
&\lesssim2^{-l\varepsilon}\frac{1}
{(V_\rho)_{2^{-k}}(y)}\int_X|f(z)|\,d\mu(z)
\quad\text{by \eqref{3.2}}
\\
&\lesssim2^{-l\varepsilon}\frac{1}
{(V_\rho)_{2^{-k}}(y)}\int_X
\frac{1}{(V_\rho)_r(x_1)+V_\rho(x_1,x)}
\left[\frac{r}{r+\rho(x,x_1)}\right]^\gamma\,d\mu(z)
\\
&\qquad\qquad\text{by the size condition of $f$}\\
&\lesssim2^{-l\varepsilon}\frac{1}{(V_\rho)_{2^{-k}}(y)}
\sim2^{-l\varepsilon}\frac{1}{(V_\rho)_{2^{-k}}(x)}.
\end{align*}
When $2^{-[k\land(k+l)]}=2^{-k}<r$,
by the regularity condition of $f$,
we conclude that, for any $y,z\in X$
satisfying $\rho(y,z)\leq\frac{r+\rho(x_1,y)}{2C_\rho}$,
\begin{align}\label{fyfz1}
|f(y)-f(z)|
&\lesssim
\left[\frac{\rho(y,z)}{r+\rho(x_1,y)}\right]^\beta
\frac{1}{(V_\rho)_r(x_1)+V_\rho(x_1,y)}
\left[\frac{r}{r+\rho(x_1,y)}\right]^\gamma
\nonumber\\
&\lesssim\left[\frac{\rho(y,z)}{r}\right]^\beta
\frac{1}{(V_\rho)_r(x_1)}.
\end{align}
When $2^{-[k\land(k+l)]}=2^{-k}<r$,
from the size condition of $f$,
we infer that, for any $y,z\in X$
satisfying $\rho(y,z)>\frac{r+\rho(x_1,y)}{2C_\rho}$,
\begin{align}\label{fyfz2}
|f(y)-f(z)|&\lesssim\frac{1}{(V_\rho)_r(x_1)+V_\rho(x_1,y)}
\left[\frac{r}{r+\rho(x_1,y)}\right]^\gamma
\nonumber\\
&\quad+
\frac{1}{(V_\rho)_r(x_1)+V_\rho(x_1,z)}
\left[\frac{r}{r+\rho(x_1,z)}\right]^\gamma
\nonumber\\
&\lesssim\frac{1}{(V_\rho)_r(x_1)}
\lesssim
\left[\frac{\rho(y,z)}{r}\right]^\beta
\frac{1}{(V_\rho)_r(x_1)}.
\end{align}
Since both Lemma~\ref{Pk}(v) and
Fubini's theorem imply that, for any $y\in X$,
\begin{align}\label{111}
\int_XD_{k+l}D_k(y,z)\,d\mu(z)
&=\int_X\int_XD_{k+l}(y,z)D_k(z,w)\,d\mu(z)d\mu(w)
\nonumber\\
&=\int_X\Bigg(\int_XD_k(z,w)\,d\mu(w)\Bigg)D_{k+l}(y,z)\,d\mu(z)=0.
\end{align}
By this, Lemma~\ref{A3.12}(i), \eqref{fyfz1},
and \eqref{fyfz2},
we find that, for any $y\in X$,
\begin{align*}
\left|\mathcal{D}_{k+l}\mathcal{D}_kf(y)\right|
&=\left|\int_X
D_{k+l}D_k(y,z)[f(z)-f(y)]\,d\mu(z)\right|
\quad\text{by \eqref{111}}
\\
&\lesssim\int_{\rho(y,z)\leq C'2^{-k}}
2^{-l\varepsilon}\frac{1}{(V_\rho)_{2^{-k}}(y)}
\left[\frac{\rho(y,z)}{r}\right]^\beta\frac{1}{(V_\rho)_r(x_1)}
\,d\mu(z),
\end{align*}
which implies that
\begin{align*}
\left|\mathcal{D}_{k+l}\mathcal{D}_kf(y)\right|
&\sim2^{-l\varepsilon}
\left(\frac{2^{-k}}{r}\right)^\beta\frac{1}{(V_\rho)_r(x_1)}
\int_{\rho(y,z)\leq C'2^{-k}}\frac{1}{(V_\rho)_{2^{-k}}(y)}
\left[\frac{\rho(y,z)}{2^{-k}}\right]^\beta\,d\mu(z)\\
&\lesssim2^{-l\varepsilon}
\left(\frac{2^{-k}}{r}\right)^\beta\frac{1}{(V_\rho)_r(x_1)}
\int_{\rho(y,z)\leq C'2^{-k}}\frac{1}{V_\rho(y,z)}
\left[\frac{\rho(y,z)}{2^{-k}}\right]^\beta\,d\mu(z)\\
&\lesssim2^{-l\varepsilon}
\left(\frac{2^{-k}}{r}\right)^\beta\frac{1}{(V_\rho)_r(x_1)}
\quad\text{by Lemma~\ref{A3.2}(i)},
\end{align*}
where $C'$ is the same positive constant as in Lemma~\ref{A3.12}.
This finishes the proof of \eqref{ckl}. Thus,
for any $x\in X$, the sequence
$\{\mathcal{R}_{N,M}f\}_{M=N+1}^\infty$
converges locally uniformly to some element,
which is denoted by $\widetilde{\mathcal{R}_N}f$.
From \eqref{RNMRNL2} and the Riesz theorem
(see, for instance, \cite[p.\,61, Theorem~2.30]{G.B.Folland}),
we deduce that
there exists a (strictly) increasing
sequence $\{M_j\}_{j\in\mathbb{N}}$ in $\mathbb{N}$
which tends to $\infty$ such that,
for $\mu$-almost every $x\in X$,
$$
\lim_{j\to\infty}\mathcal{R}_{N,M_j}f(x)=\mathcal{R}_Nf(x).
$$
Consequently, for $\mu$-almost every $x\in X$,
\begin{align}\label{WRNRN}
\widetilde{\mathcal{R}_N}f(x)
=\lim_{M\to\infty}\mathcal{R}_{N,M}f(x)
=\lim_{j\to\infty}\mathcal{R}_{N,M_j}f(x)
=\mathcal{R}_Nf(x).
\end{align}
This finishes the proof of (ii).

Finally, we prove (iii).
By (ii), we conclude that,
for any $f\in{\mathcal G}(x_1,r,\beta,\gamma)$,
$\widetilde{\mathcal{R}_N}f$ is
well defined. From \eqref{WRNRN},
Lemma~\ref{A3.13}, and Lemma~\ref{density} with $p=2$,
we infer that, $\widetilde{\mathcal{R}_N}$ can be
uniquely extended to a
bounded operator on $L^2(X)$.
This extension operator,
still denoted by $\widetilde{\mathcal{R}_N}$,
satisfies that
$\mathcal{R}_Ng=\widetilde{\mathcal{R}_N}g$
both in $L^2(X)$ and $\mu$-almost
everywhere on $X$ for any $g\in L^2(X)$.
This finishes the proof of (iii)
and hence Lemma~\ref{RNRNM}.
\end{proof}

Next, we show that
$\mathcal{T}_N^{-1}$ exists and is bounded
on spaces of test functions
when $N\in\mathbb{N}$ is large enough,
which plays a key role in
establishing sharp homogeneous continuous
Calder\'on reproducing formulae.

\begin{proposition}\label{A3.16}
Let $(X,\mathbf{q},\mu)$ be an ultra-RD-space,
where $\mu$ is assumed to be a Borel-semiregular measure on $X$ satisfying
the doubling condition \eqref{doublingcond}
with respect to $\rho\in\mathbf{q}$.
Let $0<\varepsilon\preceq{\mathrm{ind\,}}(X,\mathbf{q})$
with ${\mathrm{ind\,}}(X,\mathbf{q})$ as in \eqref{index}
and $\{\mathcal{S}_{2^{-k}}\}_{k\in\mathbb{Z}}$
be a homogeneous approximation
of the identity of order $\varepsilon$
in Definition~\ref{DefHATI}.
For any $N\in\mathbb{N}$,
let $\mathcal{R}_N$ and $\mathcal{T}_N$ be the same as in \eqref{3.6}.
Then there exist constants $C,\delta\in(0,\infty)$,
independent of $N$, such that,
for any $f\in\mathring{\mathcal{G}}(x_1,r,\beta,\gamma)$
with $x_1\in X$, $r\in(0,\infty)$,
and $\beta,\gamma\in(0,\varepsilon)$,
\begin{equation}\label{3.7}
\left\|\mathcal{R}_Nf\right\|_{\mathring{{\mathcal G}}(x_1,r,\beta,\gamma)}
\leq C2^{-N\delta}\|f\|_{\mathring{{\mathcal G}}(x_1,r,\beta,\gamma)}.	
\end{equation}
Moreover, if $N$ is large enough such that
\begin{equation}\label{3.8}
C2^{-N\delta}<1,
\end{equation}
then $\mathcal{T}_N^{-1}$ exists and
maps any space of test functions to itself.
More precisely, there exists a positive constant $\widetilde{C}$
such that, for any $f\in\mathring{\mathcal{G}}(x_1,r,\beta,\gamma)$
with $x_1\in X$, $r\in(0,\infty)$,
and $\beta,\gamma\in(0,\varepsilon)$,
$$
\left\|\mathcal{T}_N^{-1}f\right\|_{\mathring{{\mathcal G}}(x_1,r,\beta,\gamma)}
\leq\widetilde{C}\|f\|_{\mathring{{\mathcal G}}(x_1,r,\beta,\gamma)}.
$$
\end{proposition}

\begin{proof}
Given Lemma~\ref{RNRNM}, it is not necessary to distinguish
$\mathcal{R}_N$ and $\widetilde{\mathcal{R}_N}$.
We begin with showing \eqref{3.7}.
To this end, fix  $f\in\mathring{\mathcal{G}}(x_1,r,\beta,\gamma)$.
We first show that
$\mathcal{R}_Nf\in\mathring{{\mathcal G}}(x_1,r,\beta,\gamma)$.
Indeed, by \eqref{BRNM} and Definition~\ref{testfunction}(i),
we find that, for any $x\in X$ and $\varepsilon'\in(0,\varepsilon)$,
\begin{align}\label{2003}
\left|\mathcal{R}_{N,M}f(x)\right|
\lesssim2^{(\varepsilon'-\varepsilon)N}
\|f\|_{\mathring{\mathcal{G}}(x_1,r,\beta,\gamma)}
\frac{1}{(V_\rho)_r(x_1)+V_\rho(x_1,x)}
\left[\frac{r}{r+\rho(x_1,x)}\right]^\gamma
\end{align}
with the implicit positive constant independent of $M$.
From this, Lemmas~\ref{RNRNM}(iii) and \ref{A3.2}(i), and the Lebesgue dominated convergence theorem,
we deduce that
\begin{align}\label{yj}
\int_X\mathcal{R}_Nf(x)\,d\mu(x)
=\int_X\lim_{M\to\infty}\mathcal{R}_{N,M}f(x)\,d\mu(x)
=\lim_{M\to\infty}\int_X\mathcal{R}_{N,M}f(x)\,d\mu(x)=0,
\end{align}
where the last equality follows from the fact that
$\mathcal{R}_{N,M}f\in\mathring{{\mathcal G}}(x_1,r,\beta,\gamma)$
[see \eqref{BRNM}].
By \eqref{2003} and Lemma~\ref{RNRNM}(iii),
we conclude that, for any $x\in X$,
\begin{align}\label{yj1}
\left|\mathcal{R}_Nf(x)\right|
&\leq\left|\mathcal{R}_Nf(x)-\mathcal{R}_{N,M}f(x)\right|
+\left|\mathcal{R}_{N,M}f(x)\right|\nonumber\\
&\lesssim\left|\mathcal{R}_Nf(x)-\mathcal{R}_{N,M}f(x)\right|\nonumber\\
&\quad+2^{(\varepsilon'-\varepsilon)N}
\|f\|_{\mathring{\mathcal{G}}(x_1,r,\beta,\gamma)}
\frac{1}{(V_\rho)_r(x_1)+V_\rho(x_1,x)}
\left[\frac{r}{r+\rho(x_1,x)}\right]^\gamma\nonumber\\
&\to2^{(\varepsilon'-\varepsilon)N}
\|f\|_{\mathring{\mathcal{G}}(x_1,r,\beta,\gamma)}
\frac{1}{(V_\rho)_r(x_1)+V_\rho(x_1,x)}
\left[\frac{r}{r+\rho(x_1,x)}\right]^\gamma
\end{align}
as $M\to\infty$.
Similarly, from the regularity condition of $\mathcal{R}_{N,M}f$
and Lemma~\ref{RNRNM}(iii) again,
we infer that, for any $x,y\in X$ satisfying
$\rho(x,y)\leq\frac{r+\rho(x_1,x)}{2C_\rho}$,
\begin{align*}
\left|\mathcal{R}_Nf(x)-\mathcal{R}_Nf(y)\right|
&\leq\left|\mathcal{R}_Nf(x)-\mathcal{R}_{N,M}f(x)\right|
+\left|\mathcal{R}_{N,M}f(x)-\mathcal{R}_{N,M}f(y)\right|\\
&\quad+\left|\mathcal{R}_{N,M}f(y)-\mathcal{R}_Nf(y)\right|\\
&\lesssim\left|\mathcal{R}_Nf(x)-\mathcal{R}_{N,M}f(x)\right|
+\left|\mathcal{R}_Nf(y)-\mathcal{R}_{N,M}f(y)\right|\\
&\quad+2^{(\varepsilon'-\varepsilon)N}
\left[\frac{\rho(x,y)}{r+\rho(x_1,x)}\right]^\beta
\frac{1}{(V_\rho)_r(x_1)+V_\rho(x_1,x)}
\left[\frac{r}{r+\rho(x_1,x)}\right]^\gamma\\
&\to2^{(\varepsilon'-\varepsilon)N}
\left[\frac{\rho(x,y)}{r+\rho(x_1,x)}\right]^\beta
\frac{1}{(V_\rho)_r(x_1)+V_\rho(x_1,x)}
\left[\frac{r}{r+\rho(x_1,x)}\right]^\gamma
\end{align*}
as $M\to\infty$, which, together with both
\eqref{yj} and \eqref{yj1}, further implies that
$\mathcal{R}_Nf\in\mathring{{\mathcal G}}(x_1,r,\beta,\gamma)$.
By this and Lemma~\ref{RNRNM}(iii) again, we obtain
\begin{align*}
\left\|\mathcal{R}_Nf\right\|_{\mathring{{\mathcal G}}(x_1,r,\beta,\gamma)}
&\leq\left\|\mathcal{R}_Nf-\mathcal{R}_{N,M}f
\right\|_{\mathring{{\mathcal G}}(x_1,r,\beta,\gamma)}
+\left\|\mathcal{R}_{N,M}f\right\|_{\mathring{{\mathcal G}}(x_1,r,\beta,\gamma)}\\
&\lesssim\left\|\mathcal{R}_Nf-\mathcal{R}_{N,M}f
\right\|_{\mathring{{\mathcal G}}(x_1,r,\beta,\gamma)}
+2^{(\varepsilon'-\varepsilon)N}
\left\|f\right\|_{\mathring{{\mathcal G}}(x_1,r,\beta,\gamma)}\\
&\to2^{(\varepsilon'-\varepsilon)N}
\left\|f\right\|_{\mathring{{\mathcal G}}(x_1,r,\beta,\gamma)}
\end{align*}
as $M\to\infty$, which completes the proof of \eqref{3.7}
with $\delta:=\varepsilon-\varepsilon'$.

Now, we prove that $\mathcal{T}_N^{-1}$ exists and
is bounded on spaces of test functions
under the assumption \eqref{3.8}.
Indeed, if we choose $N\in\mathbb{N}$
such that \eqref{3.8} holds,
then, from \eqref{3.7} and \cite[p.\,209, Theorem]{yosida}, we deduce that
$\mathcal{T}_N=\mathcal{I}-\mathcal{R}_N$
is invertible and that we can write its inverse as
$$
(\mathcal{I}-\mathcal{R}_N)^{-1}=\sum_{l=0}^{\infty}(\mathcal{R}_N)^lf.
$$
Moreover, for any $f\in\mathring{\mathcal{G}}(x_1,r,\beta,\gamma)$,
\begin{align*}
\left\|\mathcal{T}_N^{-1}f\right\|_{\mathring{{\mathcal G}}(x_1,r,\beta,\gamma)}
&=\left\|(\mathcal{I}-\mathcal{R}_N)^{-1}f
\right\|_{\mathring{{\mathcal G}}(x_1,r,\beta,\gamma)}
=\left\|\sum_{l=0}^{\infty}(\mathcal{R}_N)^lf
\right\|_{\mathring{{\mathcal G}}(x_1,r,\beta,\gamma)}
\\
&\leq\sum_{l=0}^{\infty}\left(C2^{-N\delta}\right)^l
\left\|f\right\|_{\mathring{{\mathcal G}}(x_1,r,\beta,\gamma)}
\lesssim\|f\|_{\mathring{{\mathcal G}}(x_1,r,\beta,\gamma)},
\end{align*}
which completes the proof of Proposition~\ref{A3.16}.
\end{proof}

To establish sharp homogeneous continuous Calder\'on reproducing formulae,
we also need the following lemma.

\begin{lemma}\label{GOVV}
Let $(X,\mathbf{q},\mu)$ be a quasi-ultrametric space
of homogeneous type.
Assume that $0<\varepsilon\preceq{\mathrm{ind\,}}(X,\mathbf{q})$
and $\{\mathcal{S}_{2^{-k}}\}_{k\in\mathbb{Z}}$
is a homogeneous approximation
of the identity of order $\varepsilon$.
Then, for any $\beta\in(0,\varepsilon]$, $\gamma\in(0,\infty)$,
and $k\in\mathbb{Z}$, $\mathcal{D}_k:=\mathcal{S}_{2^{-k}}
-\mathcal{S}_{2^{-k+1}}$ is bounded on $\mathring{\mathcal{G}}(\beta,\gamma)$
[resp. $\mathcal{G}(\beta,\gamma)$]; that is, there exists a
constant $C\in(0,\infty)$, independent of $k$, such that,
for any $f\in\mathring{\mathcal{G}}(\beta,\gamma)$
[resp. $f\in\mathcal{G}(\beta,\gamma)$],
$\mathcal{D}_kf\in\mathring{\mathcal{G}}(\beta,\gamma)$
[resp. $\mathcal{D}_kf\in\mathcal{G}(\beta,\gamma)$]
and
\begin{align*}
\left\|\mathcal{D}_kf\right\|_{\mathcal{G}(\beta,\gamma)}
\leq C\max\left\{1,\,2^{-k(n+\gamma)}\right\}
\|f\|_{\mathcal{G}(\beta,\gamma)}.
\end{align*}
\end{lemma}

\begin{proof}
Without loss of generality,
we may assume that $\rho\in\mathbf{q}$ is a regularized quasi-ultrametric.
Let $\beta\in(0,\varepsilon]$, $\gamma\in(0,\infty)$, $k\in\mathbb{Z}$,
and $f\in\mathcal{G}(\beta,\gamma)$.
Without loss of generality, we may assume that
$\|f\|_{\mathcal{G}(\beta,\gamma)}=1$.

Next, we show that $\mathcal{D}_kf$ satisfies Definition~\ref{testfunction}(i).
By Lemma~\ref{Pk}, we find that there
exists a constant $C\in[1,\infty)$ such that,
for any $x,y\in X$ satisfying $\rho(x,y)\ge C2^{-k}$,
$D_k(x,y)=0$.
From Definition~\ref{D2.1}(iii), it follows that,
for any $x,y\in X$ satisfying
$\rho(x,y)<C2^{-k}$,
\begin{align*}
\rho(x_1,x)
\leq C_\rho\max\left\{\rho(x_1,y),\,\rho(y,x)\right\}
\lesssim2^{-k}+\rho(x_1,y),
\end{align*}
which further implies that
\begin{align}\label{6181}
\min\left\{1,\,2^k\right\}\left[1+\rho(x_1,x)\right]\lesssim1+\rho(x_1,y)
\end{align}
and [keeping in mind Lemma~\ref{A3.2}(vi) and \eqref{upperdoub}]
\begin{align}\label{6182}
\frac{1}{(V_\rho)_1(x_1)+V_\rho(x_1,y)}
\lesssim\frac{\max\{1,\,2^{-kn}\}}{(V_\rho)_1(x_1)+V_\rho(x_1,x)}.
\end{align}
By this,
Lemma~\ref{Pk}(i),
the size condition of $f$,
and \eqref{upperdoub}, we conclude that,
for any $x\in X$,
\begin{align}\label{2131}
\left|\mathcal{D}_kf(x)\right|
&\lesssim\int_{\rho(x,y)\leq C2^{-k}}
\frac{1}{(V_\rho)_{2^{-k}}(x)}
\frac{1}{(V_\rho)_1(x_1)+V_\rho(x_1,y)}
\left[\frac{1}{1+\rho(x_1,y)}
\right]^\gamma\,d\mu(y)
\nonumber\\
&\sim\max\left\{1,\,2^{-k(n+\gamma)}\right\}
\nonumber\\
&\quad\times\int_{\rho(x,y)\leq C2^{-k}}
\frac{1}{(V_\rho)_{2^{-k}}(x)}
\frac{1}{(V_\rho)_1(x_1)+V_\rho(x_1,x)}
\left[\frac{1}{1+\rho(x_1,x)}\right]
^\gamma\,d\mu(y)
\nonumber\\
&\lesssim\max\left\{1,\,2^{-k(n+\gamma)}\right\}
\frac{1}{(V_\rho)_1(x_1)+V_\rho(x_1,x)}
\left[\frac{1}{1+\rho(x_1,x)}\right]^\gamma.
\end{align}
Hence, $\mathcal{D}_kf$ satisfies the desired size condition.

Now, we prove that
$\mathcal{D}_kf$ satisfies Definition~\ref{testfunction}(ii).
We first assume that $k\leq0$.
Let $x,x'\in X$ satisfy
$\rho(x,x')\leq\frac{1+\rho(x_1,x)}{2C_\rho}$.
We consider the following four cases.

\emph{Case (1)}
$\rho(x_1,x)\leq  C_12^{-k}$,
where $C_1\in(1,\infty)$ is a constant
specified later. In this
case, we have
$$
\rho(x,x')\le\frac{1+\rho(x_1,x)}{2C_\rho}
\leq\frac{C_12^{-k}}{C_\rho},
$$
which further implies that
\begin{align*}
\left[B_\rho(x,C2^{-k})\cup
B_\rho(x',C2^{-k})\right]
\subset B_\rho\left(x,\max
\left\{CC_\rho,\,C_1\right\}2^{-k}\right).
\end{align*}
From this, Lemma~\ref{Pk}(ii), the size condition of $f$, Lemma~\ref{A3.2}(vii),
$\beta\leq\varepsilon$,
\eqref{6181}, and \eqref{6182},
we infer that
\begin{align*}
&\left|\mathcal{D}_kf(x)-\mathcal{D}_kf(x')\right|\\
&\quad\lesssim\int_{\rho(x,y)\leq\max
\{CC_\rho,\,C_1\}2^{-k}}
\frac{2^{k\varepsilon}[\rho(x,x')]^\varepsilon}{(V_\rho)_{2^{-k}}(y)}
\frac{1}
{(V_\rho)_1(x_1)+V_\rho(x_1,y)}
\left[\frac{1}{1+\rho(x_1,y)}\right]^\gamma\,d\mu(y)
\\
&\quad\lesssim2^{-k(n+\gamma)}\int_{\rho(x,y)\leq\max
\{CC_\rho,\,C_1\}2^{-k}}
\frac{2^{k\beta}[\rho(x,x')]^\beta}{(V_\rho)_{2^{-k}}(x)}
\\
&\qquad\times\frac{1}
{(V_\rho)_1(x_1)+V_\rho(x_1,x)}
\left[\frac{1}{1+\rho(x_1,x)}\right]^\gamma\,d\mu(y).
\end{align*}
This, combined with $\rho(x_1,x)\lesssim2^{-k}$ and \eqref{upperdoub},
implies that
\begin{align*}
&\left|\mathcal{D}_kf(x)-\mathcal{D}_kf(x')\right|\\
&\quad\lesssim2^{-k(n+\gamma)}
\left[\frac{\rho(x,x')}{1+\rho(x_1,x)}\right]^\beta
\frac{1}{(V_\rho)_1(x_1)+V_\rho(x_1,x)}
\left[\frac{1}{1+\rho(x_1,x)}\right]^\gamma,
\end{align*}
which is the desired estimate.

\emph{Case (2)}
$\rho(x_1,x)>C_12^{-k}$ and
$\rho(x,x')\ge C_2\rho(x_1,x)$ for some fixed $C_2\in(0,\frac{1}{C_\rho})$.
In this case,
by $\rho(x,x')\leq\frac{1+\rho(x_1,x)}{2C_\rho}$
and Definition~\ref{D2.1}(iii),
we find that
$$
\rho(x,x')\sim\rho(x_1,x)\sim\rho(x_1,x'),
$$
which, together with the doubling condition \eqref{doublingcond}, gives
$V_\rho(x_1,x)\sim V_\rho(x_1,x')$.
Using this,
we conclude that
\begin{align*}
\left|\mathcal{D}_kf(x)-\mathcal{D}_kf(x')\right|
&\leq\left|\mathcal{D}_kf(x)\right|+\left|\mathcal{D}_kf(x')\right|
\\
&\lesssim2^{-k(n+\gamma)}
\frac{1}{(V_\rho)_1(x_1)+V_\rho(x_1,x)}
\left[\frac{1}{1+\rho(x_1,x)}\right]^\gamma
\\
&\quad+2^{-k(n+\gamma)}
\frac{1}{(V_\rho)_1(x_1)+V_\rho(x_1,x')}
\left[\frac{1}{1+\rho(x_1,x')}\right]^\gamma
\quad\text{by \eqref{2131}}\\
&\sim2^{-k(n+\gamma)}\left[\frac{\rho(x,x')}{1
+\rho(x_1,x)}\right]^\beta
\frac{1}{(V_\rho)_1(x_1)+V_\rho(x_1,x)}
\left[\frac{1}{1+\rho(x_1,x)}\right]^\gamma
\end{align*}
since $\rho(x,x')\ge C_2\rho(x_1,x)\ge\frac{C_2[1+\rho(x_1,x)]}{2}$.
This is also the desired estimate.

\emph{Case (3)}
$\rho(x_1,x)>C_12^{-k}$, $\rho(x,x')<C_2\rho(x_1,x)$,
and $\rho(x,x')<C_32^{-k}$ for some fixed $C_3\in(0,\infty)$.
In this case, we have
\begin{align*}
\left[B_\rho(x,C2^{-k})\cup
B_\rho(x',C2^{-k})\right]
\subset B_\rho\left(x,(C+C_3)C_\rho2^{-k}\right).
\end{align*}
Thus, from Lemma~\ref{Pk}(i), we deduce that,
for any $y\in[B_\rho(x,(C+C_3)C_\rho2^{-k})]^\complement$,
\begin{align}
\label{1604}
D_k(x,y)=D_k(x',y)=0.
\end{align}
Choose $C_1\in(1,\infty)$ such that
\begin{align}
\label{2209}
(C+C_3)C_\rho2^{-k}
\leq\frac{1+C_1}{2C_\rho}2^{-k}
<\frac{1+\rho(x_1,x)}{2C_\rho}.
\end{align}
By this, (ii) and (v) of Lemma~\ref{Pk}, \eqref{1604},
and the regularity condition of $f$,
we obtain
\begin{align*}
&\left|\mathcal{D}_kf(x)-\mathcal{D}_kf(x')\right|
\\
&\quad=\left|\int_X[D_k(x,y)-D_k(x',y)][f(y)-f(x)]\,d\mu(y)\right|
\\
&\quad\lesssim\int_{B_\rho(x,(C+C_3)C_\rho2^{-k})}
\left[\frac{\rho(x,x')}{2^{-k}}\right]^\varepsilon
\frac{1}{(V_\rho)_{2^{-k}}(x)}
\left[\frac{\rho(x,y)}{1+\rho(x_1,x)}\right]^\beta
\\
&\qquad\times\frac{1}{(V_\rho)_1(x_1)+V_\rho(x_1,x)}
\left[\frac{1}{1+\rho(x_1,x)}\right]^\gamma\,d\mu(y),
\end{align*}
which, combined with \eqref{upperdoub},
further implies that
\begin{align*}
&\left|\mathcal{D}_kf(x)-\mathcal{D}_kf(x')\right|
\\
&\quad\lesssim
\left[\frac{\rho(x,x')}{1+\rho(x_1,x)}\right]^\beta
\frac{1}{(V_\rho)_1(x_1)+V_\rho(x_1,x)}
\left[\frac{1}{1+\rho(x_1,x)}\right]^\gamma
\end{align*}
since $\beta\leq\varepsilon$ and
$\rho(x,x')\leq C_32^{-k}$.
This is the desired estimate.

\emph{Case (4)}
$\rho(x_1,x)>C_12^{-k}$, $\rho(x,x')<C_2\rho(x_1,x)$,
and $\rho(x,x')\ge C_32^{-k}$.
In this case, we have
\begin{align}\label{Dkf}
\left|\mathcal{D}_kf(x)-\mathcal{D}_kf(x')\right|\leq
\left|\mathcal{D}_kf(x)\right|+\left|\mathcal{D}_kf(x')\right|.
\end{align}
Note that, in the current case, for any $y\in B_\rho(x,C2^{-k})$,
$\rho(x,y)<C2^{-k}\leq\frac{C\rho(x,x')}{C_3}$.
From this, both (i) and (v) of Lemma~\ref{Pk},
and the regularity condition of $f$ [using \eqref{2209}],
we infer that
\begin{align}\label{6183}
\left|\mathcal{D}_kf(x)\right|
&=\left|\int_XD_k(x,y)[f(y)-f(x)]\,d\mu(y)\right|
\nonumber\\
&\lesssim\int_{\rho(x,y)\leq C2^{-k}}
\frac{1}{(V_\rho)_{2^{-k}}(x)}
\left[\frac{\rho(x,y)}{1+\rho(x_1,x)}\right]^\beta
\nonumber\\
&\quad\times\frac{1}{(V_\rho)_1(x_1)+V_\rho(x_1,x)}
\left[\frac{1}{1+\rho(x_1,x)}\right]^\gamma\,d\mu(y)
\nonumber\\
&\lesssim\left[\frac{\rho(x,x')}{1+\rho(x_1,x)}\right]^\beta
\frac{1}{(V_\rho)_1(x_1)+V_\rho(x_1,x)}
\left[\frac{1}{1+\rho(x_1,x)}\right]^\gamma.
\end{align}
For $\mathcal{D}_kf(x')$, note that,
by Definition~\ref{D2.1}(iii), $\rho(x,x')<C_2\rho(x_1,x)$,
and $C_2<\frac{1}{C_\rho}$, we find that
\begin{align}
\label{6184}
\rho(x_1,x)
&\leq C_\rho\max\left\{\rho(x_1,x'),\,\rho(x',x)\right\}
=C_\rho\rho(x_1,x').
\end{align}
On the other hand, from the assumptions in this case,
we deduce that, for any $y\in B_\rho(x',C2^{-k})$,
\begin{align*}
\rho(x',y)<C2^{-k}\leq\frac{C}{C_3}\rho(x,x').
\end{align*}
By this and an argument similar to that
used in the estimation of \eqref{6183}, we conclude that
\begin{align*}
\left|\mathcal{D}_kf(x')\right|
&\lesssim\left[\frac{\rho(x,x')}{1+\rho(x_1,x')}\right]^\beta
\frac{1}{(V_\rho)_1(x_1)+V_\rho(x_1,x')}
\left[\frac{1}{1+\rho(x_1,x')}\right]^\gamma,
\end{align*}
which, together with \eqref{6184},
further implies that
\begin{align*}
\left|\mathcal{D}_kf(x')\right|
&\lesssim\left[\frac{\rho(x,x')}{1+\rho(x_1,x)}\right]^\beta
\frac{1}{(V_\rho)_1(x_1)+V_\rho(x_1,x)}
\left[\frac{1}{1+\rho(x_1,x)}\right]^\gamma.
\end{align*}
Combining this, \eqref{Dkf}, and \eqref{6183},
we then show that, for any $k\in\mathbb{Z}\cap(-\infty,0]$,
\begin{align}\label{2243}
&\left|\mathcal{D}_kf(x)-\mathcal{D}_kf(x')\right|
\nonumber\\
&\quad\lesssim2^{-k(n+\gamma)}\left[\frac{\rho(x,x')}{1
+\rho(x_1,x)}\right]^\beta
\frac{1}{(V_\rho)_1(x_1)+V_\rho(x_1,x)}
\left[\frac{1}{1+\rho(x_1,x)}\right]^\gamma,
\end{align}
where the implicit positive constant is independent of $k$.
This finishes the proof of the desired regularity of
$\mathcal{D}_kf$ in the case $k\in\mathbb{Z}\cap(-\infty,0]$.

Next, we prove the regularity of $\mathcal{D}_kf$ when
$k\in\mathbb{N}$.
According to Remark~\ref{Rmtest}(i), we only need to consider
$k\in\mathbb{N}\cap[K_0,\infty)$, where $K_0\in\mathbb{N}$ is a fixed integer
satisfying $CC_\rho2^{-k_0}\leq\frac{1}{2C_\rho}$. To this end,
we consider the following two cases for $k$.

\emph{Case (1)} $\rho(x,x')\leq C2^{-k}$. In this case, note that
\begin{align*}
B_\rho(x,C2^{-k})\cup B_\rho(x',C2^{-k})\subset B_\rho(x,CC_\rho2^{-k})
\end{align*}
and the choice of $K_0$ implies that, for any $y\in B_\rho(x,CC_\rho2^{-k})$,
\begin{align*}
\rho(x,y)<CC_\rho2^{-k}\leq\frac{1}{2C_\rho}<\frac{1+\rho(x_1,x)}{2C_\rho}.
\end{align*}
From this, (ii) and (v) of Lemma~\ref{Pk},
the regularity condition of $f$, and $\beta\leq\varepsilon$,
we infer that
\begin{align*}
&\left|\mathcal{D}_kf(x)-\mathcal{D}_kf(x')\right|
\\
&\quad=\left|\int_X[D_k(x,y)-D_k(x',y)][f(y)-f(x)]\,d\mu(y)\right|
\\
&\quad\lesssim\int_{B_\rho(x,CC_\rho2^{-k})}
\left[\frac{\rho(x,x')}{2^{-k}}\right]^\varepsilon
\frac{1}{(V_\rho)_{2^{-k}}(x)}
\left[\frac{\rho(x,y)}{1+\rho(x_1,x)}\right]^\beta
\\
&\qquad\times\frac{1}{(V_\rho)_1(x_1)+V_\rho(x_1,x)}
\left[\frac{1}{1+\rho(x_1,x)}\right]^\gamma\,d\mu(y)
\\
&\quad\lesssim
\left[\frac{\rho(x,x')}{1+\rho(x_1,x)}\right]^\beta
\frac{1}{(V_\rho)_1(x_1)+V_\rho(x_1,x)}
\left[\frac{1}{1+\rho(x_1,x)}\right]^\gamma,
\end{align*}
which is the desired estimate.

\emph{Case (2)} $\rho(x,x')>C2^{-k}$.
In this case, we first claim that
\begin{align}\label{2030}
1+\rho(x_1,x)\lesssim1+\rho(x_1,x').
\end{align}
Indeed, if $\rho(x_1,x)\leq1$, then \eqref{2030} holds automatically.
If $\rho(x_1,x)>1$, then Definition~\ref{D2.1}(iii)
and $\rho(x,x')\leq\frac{1+\rho(x_1,x)}{2C_\rho}$ imply that
\begin{align*}
\rho(x_1,x)&\leq C_\rho\max\left\{\rho(x_1,x'),\,\rho(x,x')\right\}\\
&\leq C_\rho\max\left\{\rho(x_1,x'),\,\frac{1+\rho(x_1,x)}{2C_\rho}\right\}
=C_\rho\rho(x_1,x')
\end{align*}
and hence \eqref{2030} holds,
which then completes the proof of the claim \eqref{2030}.
Note that the choice of $K_0$ implies that, for any $y\in B_\rho(x,C2^{-k})$
[resp. $y\in B_\rho(x',C2^{-k})$],
\begin{align}\label{1129}
\rho(x,y)<C2^{-k}\leq\frac{1+\rho(x_1,x)}{2C_\rho}\ \
\left[\text{resp}\ \ \rho(x',y)<C2^{-k}\leq\frac{1+\rho(x_1,x')}{2C_\rho}\right].
\end{align}
Using both (i) and (v) of Lemma~\ref{Pk},
we find that
\begin{align*}
\left|\mathcal{D}_kf(x)-\mathcal{D}_kf(x')\right|
&\leq\left|\mathcal{D}_kf(x)\right|+\left|\mathcal{D}_kf(x')\right|
\\
&=\left|\int_{B_\rho(x,C2^{-k})}D_k(x,y)[f(y)-f(x)]\,d\mu(y)\right|
\\
&\quad+\left|\int_{B_\rho(x',C2^{-k})}D_k(x',y)[f(y)-f(x')]\,d\mu(y)\right|,
\end{align*}
which, together with Lemma~\ref{Pk}(i),
\eqref{1129}, and the regularity condition of $f$,
further implies that
\begin{align*}
&\left|\mathcal{D}_kf(x)-\mathcal{D}_kf(x')\right|\\
&\quad\lesssim\int_{B_\rho(x,C2^{-k})}\frac{1}{(V_\rho)_{2^{-k}}(x)}
\left[\frac{\rho(x,x')}{2^{-k}}\right]^\beta
\left[\frac{\rho(x,y)}{1+\rho(x_1,x)}\right]^\beta
\\
&\qquad\times\frac{1}{(V_\rho)_1(x_1)+V_\rho(x_1,x)}
\left[\frac{1}{1+\rho(x_1,x)}\right]^\gamma\,d\mu(y)
\\
&\qquad+\int_{B_\rho(x',C2^{-k})}\frac{1}{(V_\rho)_{2^{-k}}(x')}
\left[\frac{\rho(x,x')}{2^{-k}}\right]^\beta
\left[\frac{\rho(x',y)}{1+\rho(x_1,x')}\right]^\beta
\\
&\qquad\times\frac{1}{(V_\rho)_1(x_1)+V_\rho(x_1,x')}
\left[\frac{1}{1+\rho(x_1,x')}\right]^\gamma\,d\mu(y).
\end{align*}
From this, \eqref{upperdoub}, and \eqref{2030},
it follows that
\begin{align*}
&\left|\mathcal{D}_kf(x)-\mathcal{D}_kf(x')\right|\\
&\quad\lesssim\left[\frac{\rho(x,x')}{1+\rho(x_1,x)}\right]^\beta
\frac{1}{(V_\rho)_1(x_1)+V_\rho(x_1,x)}
\left[\frac{1}{1+\rho(x_1,x)}\right]^\gamma
\\
&\qquad+\left[\frac{\rho(x,x')}{1+\rho(x_1,x')}\right]^\beta
\frac{1}{(V_\rho)_1(x_1)+V_\rho(x_1,x')}
\left[\frac{1}{1+\rho(x_1,x')}\right]^\gamma
\\
&\quad\lesssim\left[\frac{\rho(x,x')}{1+\rho(x_1,x)}\right]^\beta
\frac{1}{(V_\rho)_1(x_1)+V_\rho(x_1,x)}
\left[\frac{1}{1+\rho(x_1,x)}\right]^\gamma.
\end{align*}
Combining the above two cases, we then conclude that,
for any $k\in\mathbb{N}$,
\begin{align*}
\left|\mathcal{D}_kf(x)-\mathcal{D}_kf(x')\right|
\lesssim\left[\frac{\rho(x,x')}{1+\rho(x_1,x)}\right]^\beta
\frac{1}{(V_\rho)_1(x_1)+V_\rho(x_1,x)}
\left[\frac{1}{1+\rho(x_1,x)}\right]^\gamma,
\end{align*}
which, together with \eqref{2131} and \eqref{2243},
further implies that, for any $k\in\mathbb{Z}$,
\begin{align}\label{2247}
\|\mathcal{D}_kf\|_{\mathcal{G}(\beta,\gamma)}
\lesssim\max\left\{1,\,2^{-k(n+\gamma)}\right\}\|f\|_{\mathcal{G}(\beta,\gamma)},
\end{align}
where the implicit positive constant is independent of $k$.

Now, let $f\in\mathring{\mathcal{G}}(\beta,\gamma)$.
By Fubini's theorem
and Lemma~\ref{Pk}(v),
we find that
\begin{align*}
\int_X\mathcal{D}_kf(x)\,d\mu(x)
=\int_X\left[\int_XD_k(x,y)\,d\mu(x)
\right]f(y)\,d\mu(y)
=0.
\end{align*}
From this and \eqref{2247},
we deduce that $\mathcal{D}_kf\in\mathring{\mathcal{G}}(\beta,\gamma)$
and $\|\mathcal{D}_kf\|_{\mathcal{G}(\beta,\gamma)}
\lesssim\|f\|_{\mathcal{G}(\beta,\gamma)}$.
This finishes the proof of Lemma~\ref{GOVV}.
\end{proof}

We are ready to
establish the following sharp homogeneous
continuous Calder\'on reproducing formula.
\index[words]{homogeneous continuous Calder\'on reproducing formula}

\begin{theorem}\label{A3.19}
Let $(X,\mathbf{q},\mu)$ be an ultra-RD-space,
where $\mu$ is a Borel-semiregular measure on $X$.
Assume that $0<\varepsilon\preceq{\mathrm{ind\,}}(X,\mathbf{q})$
with ${\mathrm{ind\,}}(X,\mathbf{q})$ as in \eqref{index} and
$\{\mathcal{S}_{2^{-k}}\}_{k\in\mathbb{Z}}$
(with kernels $\{S_{2^{-k}}\}_{k\in\mathbb{Z}}$)
is a homogeneous approximation
of the identity of order $\varepsilon$
in Definition~\ref{DefHATI}.
For any $k\in\mathbb{Z}$, let
$\mathcal{D}_k:=\mathcal{S}_{2^{-k}}-\mathcal{S}_{2^{-k+1}}$.
Then there exist two families
$\{\widetilde{\mathcal{D}}_k\}_{k\in\mathbb{Z}}$ and
$\{\overline{\mathcal{D}}_k\}_{k\in\mathbb{Z}}$
of bounded linear operators on $L^2(X)$,
with their kernels denoted, respectively, by
$\{\widetilde{D}_k\}_{k\in\mathbb{Z}}$ and
$\{\overline{D}_k\}_{k\in\mathbb{Z}}$,
such that, for any $f\in\mathring{{\mathcal G}}_0^\varepsilon(\beta,\gamma)$
with $\beta,\gamma\in(0,\varepsilon)$
[resp. $L^p(X)$ with $p\in(1,\infty)$],
\begin{equation}\label{3.12}
f=\sum_{k=-\infty}^{\infty}\widetilde{\mathcal{D}}_k\mathcal{D}_kf
=\sum_{k=-\infty}^{\infty}\mathcal{D}_k\overline{\mathcal{D}}_kf,
\end{equation}
where the series in \eqref{3.12} converge
in $\mathring{{\mathcal G}}_0^\varepsilon(\beta,\gamma)$
[resp. $L^p(X)$ with $p\in(1,\infty)$]
and also converge pointwise $\mu$-almost everywhere on $X$.
Moreover, for any $k\in\mathbb{Z}$, $x,y\in X$, and
$\widetilde{\beta},\ \widetilde{\gamma}\in(0,\varepsilon)$,
\begin{align}
\label{8223}
\widetilde{D}_k(\cdot,y)\in{\mathcal G}
(y,2^{-k},\widetilde{\beta},\widetilde{\gamma}),
\end{align}
\begin{align}
\label{wDkvanish}
\int_X\widetilde{D}_k(z,y)\,d\mu(z)=0
=\int_X\widetilde{D}_k(x,z)\,d\mu(z),
\end{align}
and
\begin{align}
\label{wDkoDk}
\widetilde{D}_k(x,y)=\overline{D}_k(y,x).
\end{align}
\end{theorem}

\begin{remark}
\begin{enumerate}
\item[\textup{(i)}]
In Theorem~\ref{A3.19}, for any $k\in\mathbb{Z}$,
all the estimates satisfied by
$\widetilde{\mathcal{D}}_k$ and $\overline{\mathcal{D}}_k$
are independent of both $\beta$ and $\gamma$ with
$\beta,\gamma\in(0,\varepsilon)$, but may depend on $\varepsilon$.

\item[\textup{(ii)}]
If $\rho\in\mathbf{q}$ is a metric, then ${\mathrm{ind\,}}(X,\mathbf{q})\in[1,\infty]$
and, for any $0<\beta,\gamma<\varepsilon\preceq{\mathrm{ind\,}}(X,\mathbf{q})$,
\eqref{3.12} holds in $\mathring{{\mathcal G}}^\varepsilon_0(\beta,\gamma)$;
if $\rho\in\mathbf{q}$ is an ultrametric, then,
for any $0<\beta,\gamma<\varepsilon<\infty$,
\eqref{3.12} holds in $\mathring{{\mathcal G}}^\varepsilon_0(\beta,\gamma)$.
However, Han et al. \cite[Theorem 3.10]{hmy08}
and He et al. \cite[Theorem 4.15]{HLYY} only showed
homogeneous continuous Calder\'on reproducing
formulae hold
for any $0<\beta,\gamma<\varepsilon<1$.
In this sense, Theorem~\ref{A3.19} is sharper than
\cite[Theorem 3.10]{hmy08} and \cite[Theorem 4.15]{HLYY}
because of the sharpness of both Theorems~\ref{ThmATI}
and \ref{CZOonGO}.
\end{enumerate}
\end{remark}

\begin{proof}[Proof of Theorem~\ref{A3.19}]
Let $\rho\in\mathbf{q}$ be a regularized quasi-ultrametric.
For any $k\in\mathbb{Z}$ and $N\in\mathbb{N}$,
let $\mathcal{D}^N_k$ be the same as in \eqref{DkN7}
and both $\mathcal{R}_N$ and $\mathcal{T}_N$ the same as in \eqref{3.6}.
Fix $N\in\mathbb{N}$ sufficiently large such that
\eqref{3.8} holds,
and therefore Proposition~\ref{A3.16}
holds. By Lemma~\ref{ATItest} and Definition~\ref{DefATI}(v), we find that,
for any $k\in\mathbb{Z}$, $y\in X$,
and $\widetilde{\beta},\widetilde{\gamma}\in(0,\varepsilon)$,
\begin{align}
\label{2123}
D^N_k(\cdot,y)
=S_{2^{-(k+N)}}(\cdot,y)
-S_{2^{-(k-N)+1}}(\cdot,y)
\in\mathring{\mathcal{G}}(y,2^{-k},\widetilde{\beta},\widetilde{\gamma})
\end{align}
with its norm depending on $N$ but independent of both $y$ and $k$.
For any $k\in\mathbb{Z}$ and $x,y\in X$,
let
\begin{align}\label{1617}
\widetilde{D}_k(x,y):=\mathcal{T}^{-1}_N\left(D^N_k(\cdot,y)\right)(x).
\end{align}

Next, we divide the proof of the present theorem
into the following five steps.

\emph{Step 1.} In this step, we prove that, for any $k\in\mathbb{Z}$,
$\widetilde{D}_k$ satisfies both \eqref{8223} and \eqref{wDkvanish}.
To this end, let $k\in\mathbb{Z}$, $x,y\in X$, and
$\widetilde{\beta},\widetilde{\gamma}\in(0,\varepsilon)$.
We first show \eqref{8223}.
From \eqref{1617}, Proposition~\ref{A3.16}, and
\eqref{2123}, we deduce that
$\widetilde{D}_k(\cdot,y)
\in\mathring{{\mathcal G}}(y,2^{-k},\widetilde{\beta},\widetilde{\gamma})$.
Thus, \eqref{8223} holds.

Then we prove \eqref{wDkvanish}.
Note that $\widetilde{D}_k(\cdot,y)
\in\mathring{{\mathcal G}}(y,2^{-k},\widetilde{\beta},\widetilde{\gamma})$ implies
\begin{align}
\label{wDkvanish1}
\int_X\widetilde{D}_k(z,y)\,d\mu(z)=0,
\end{align}
which is exactly the first equality of \eqref{wDkvanish}.
Now, we turn to show the second equality of \eqref{wDkvanish}.
To this end, we first claim that, for any $j\in\mathbb{Z}_+$,
\begin{align}
\label{919}
\int_X(R_N)^j\left(D_k^N(\cdot,z)\right)(x)\,d\mu(z)=0.
\end{align}
We prove \eqref{919} via a method of induction.
By \eqref{DkN7} and Lemma~\ref{Pk}(v),
we find that \eqref{919} holds when $j=0$.
Assume that \eqref{919} holds for some $j\in\mathbb{Z}_+$.
Next, we show that \eqref{919} with $j$ replaced by $j+1$ holds.
For any $M\in\mathbb{N}\cap(N,\infty)$,
let $\mathcal{R}_{N,M}$ be the same as in \eqref{defRNM}.
From \eqref{BRNM} and \eqref{3.7}, we infer that,
for any given $z\in X$,
\begin{align*}
\left\|(\mathcal{R}_N)^j\left(D_k^N(\cdot,z)\right)\right\|_
{{\mathcal G}(z,2^{-k},\beta,\gamma)}\lesssim
\left(C2^{-N\delta}\right)^j\left\|D_k^N(\cdot,z)\right\|_
{{\mathcal G}(z,2^{-k},\beta,\gamma)}
\lesssim\left(C2^{-N\delta}\right)^j
\end{align*}
and
\begin{align*}
\left\|\mathcal{R}_{N,M}(\mathcal{R}_N)^j\left(D_k^N(\cdot,z)\right)\right\|_
{{\mathcal G}(z,2^{-k},\beta,\gamma)}\lesssim
\left(C2^{-N\delta}\right)^j\left\|D_k^N(\cdot,z)\right\|_
{{\mathcal G}(z,2^{-k},\beta,\gamma)}
\lesssim\left(C2^{-N\delta}\right)^j,
\end{align*}
where $C2^{-N\delta}<1$ is the same as in \eqref{3.7},
which, combined with Definition~\ref{testfunction}(i)
and Lemma~\ref{A3.2}(vi), further implies that
\begin{align}
\label{12182}
\left|(\mathcal{R}_N)^j\left(D_k^N(\cdot,z)\right)(x)\right|
&\lesssim\left(C2^{-N\delta}\right)^j
\frac{1}{(V_\rho)_{2^{-k}}(z)+V_\rho(z,x)}
\left[\frac{2^{-k}}{2^{-k}+\rho(z,x)}\right]^\gamma
\nonumber\\
&\sim\left(C2^{-N\delta}\right)^j
\frac{1}{(V_\rho)_{2^{-k}}(x)+V_\rho(x,z)}
\left[\frac{2^{-k}}{2^{-k}+\rho(x,z)}\right]^\gamma
\end{align}
and
\begin{align}
\label{12181}
\left|\mathcal{R}_{N,M}(\mathcal{R}_N)^j\left(D_k^N(\cdot,z)\right)(x)\right|
&\lesssim\left(C2^{-N\delta}\right)^j
\frac{1}{(V_\rho)_{2^{-k}}(z)+V_\rho(z,x)}
\left[\frac{2^{-k}}{2^{-k}+\rho(z,x)}\right]^\gamma
\nonumber\\
&\sim\left(C2^{-N\delta}\right)^j
\frac{1}{(V_\rho)_{2^{-k}}(x)+V_\rho(x,z)}
\left[\frac{2^{-k}}{2^{-k}+\rho(x,z)}\right]^\gamma.
\end{align}
By Lemma~\ref{RNRNM}(iii),
the Lebesgue dominated convergence theorem together with both
\eqref{12181} and Lemma~\ref{A3.2}(ii), Fubini's theorem combined with
both \eqref{12182}
and the fact that
$$
\int_X|R_{N,M}(x,z)|\,d\mu(z)
\leq\sum_{|k|\leq M}\sum_{N<|l|\leq M}\int_X|D_{k+l}D_k(x,z)|\,d\mu(z)<\infty,
$$
and the above assumption for $j$,
we conclude that
\begin{align*}
&\int_X(\mathcal{R}_N)^{j+1}
\left(D_k^N(\cdot,y)\right)(x)\,d\mu(y)
\\
&\quad=\int_X\lim_{M\to\infty}
\mathcal{R}_{N,M}(\mathcal{R}_N)^j
\left(D_k^N(\cdot,y)\right)(x)\,d\mu(y)
\\
&\quad=\lim_{M\to\infty}\int_X
\mathcal{R}_{N,M}(\mathcal{R}_N)^j
\left(D_k^N(\cdot,y)\right)(x)\,d\mu(y)
\\
&\quad=\lim_{M\to\infty}\int_X\int_X
\mathcal{R}_{N,M}(x,z)(\mathcal{R}_N)^j
\left(D_k^N(\cdot,y)\right)(z)\,d\mu(z)\,d\mu(y)
\\
&\quad=\lim_{M\to\infty}\int_X
\mathcal{R}_{N,M}(x,z)\int_X(\mathcal{R}_N)^j
\left(D_k^N(\cdot,y)\right)(z)\,d\mu(y)\,d\mu(z)
=0.
\end{align*}
This finishes the proof of \eqref{919}
and hence the aforementioned claim.
From this claim, \eqref{1617}, and the Lebesgue dominated convergence theorem
[here, we need \eqref{12182}, Lemma~\ref{A3.2}(ii),
and \eqref{3.8}],
we deduce that
\begin{align*}
\int_X\widetilde{D}_k(x,y)\,d\mu(y)
&=\int_X\mathcal{T}_N^{-1}\left(D_k^N(\cdot,y)\right)(x)\,d\mu(y)
\\
&=\int_X\sum_{j=0}^\infty(\mathcal{R}_N)^j
\left(D_k^N(\cdot,y)\right)(x)\,d\mu(y)
\\
&=\lim_{J\to\infty}\sum_{j=0}^J
\int_X(\mathcal{R}_N)^j\left(D_k^N(\cdot,y)\right)(x)\,d\mu(y)
=0,
\end{align*}
which, together with \eqref{wDkvanish1},
completes the proof of \eqref{wDkvanish}.

\emph{Step 2.} In this step, let $\beta,\gamma\in(0,\varepsilon)$
and we prove that the first series in \eqref{3.12}
converges in $\mathring{{\mathcal G}}_0^\varepsilon(\beta,\gamma)$
for any $f\in\mathring{\mathcal{G}}(\beta',\gamma')$
with any given $\beta'\in(\beta,\varepsilon)$
and $\gamma'\in(\gamma,\varepsilon)$.
Observing that, for any $k\in\mathbb{Z}$ and $x,y\in X$,
since
$$
T_N\widetilde{D}_k(x,y)
=\mathcal{T}_N\left(\widetilde{D}_k(\cdot,y)\right)(x)
=D_k^N(x,y),
$$
where $T_N\widetilde{D}_k(x,y)$ is the kernel of
$\mathcal{T_N}\widetilde{\mathcal{D}_k}$,
it follows that $\mathcal{T}_N^{-1}\mathcal{D}_k^N
=\widetilde{\mathcal{D}}_k$,
which, combined with the definition of $\mathcal{T}_N$
in \eqref{3.6} in the sense of $L^2(X)$,
further implies that, for any $g\in L^2(X)$,
\begin{align}\label{2134}
g&=\mathcal{T}_N^{-1}\mathcal{T}_Ng
=\mathcal{T}_N^{-1}\left(\sum_{k=-\infty}^\infty
\mathcal{D}_k^N\mathcal{D}_k\right)g
\nonumber\\
&=\sum_{k=-\infty}^\infty\left(\mathcal{T}_N^{-1}
\mathcal{D}_k^N\right)\mathcal{D}_kg
=\sum_{k=-\infty}^\infty\widetilde{\mathcal{D}}_k
\mathcal{D}_kg
\end{align}
in $L^2(X)$.
By this, we find that, for any $L\in\mathbb{N}$ and $g\in L^2(X)$,
\begin{align}
\label{jj-4u5}
\sum_{|k|\leq L}\widetilde{\mathcal{D}}_k\mathcal{D}_kg
&=\mathcal{T}^{-1}_N\left(\sum_{|k|\leq L}\mathcal{D}^N_k\mathcal{D}_k\right)g
=\mathcal{T}^{-1}_N\left(\mathcal{T}_N-
\sum_{|k|\ge L+1}\mathcal{D}^N_k\mathcal{D}_k\right)g
\nonumber\\
&=g-\mathcal{T}^{-1}_N\left(\sum_{|k|\ge L+1}\mathcal{D}^N_k\mathcal{D}_k\right)g
\end{align}
in $L^2(X)$. Thus, we need to show that
\begin{equation}\label{3.13}
\lim_{L\to\infty}\left\|f-\sum_{|k|\leq L}
\widetilde{\mathcal{D}}_k\mathcal{D}_kf\right\|_{{\mathcal G}(\beta,\gamma)}
=\lim_{L\to\infty}\left\|\mathcal{T}^{-1}_N\left(\sum_{|k|\ge L+1}
\mathcal{D}^N_k\mathcal{D}_k\right)f\right\|_{{\mathcal G}(\beta,\gamma)}=0.
\end{equation}
To prove \eqref{3.13}, let $L\in\mathbb{N}$.
From Proposition~\ref{A3.16}
and the triangle inequality for
$\|\cdot\|_{{\mathcal G}(\beta,\gamma)}$ [see Remark~\ref{Rmtest}(iii)],
it suffices to show that there exists $\sigma_0\in(0,\infty)$
such that, for any $k\in\mathbb{Z}$ with $|k|\in[L+1,\infty)$,
\begin{equation}\label{3.16}
\left\|\mathcal{D}^N_k\mathcal{D}_kf\right\|_{{\mathcal G}(\beta,\gamma)}
\lesssim 2^{-\sigma_0|k|}\|f\|_{{\mathcal G}(\beta',\gamma')}.
\end{equation}
Indeed, if \eqref{3.16} holds, then,
by this and Proposition~\ref{A3.16}, we conclude that
\begin{align*}
\lim_{L\to\infty}\left\|\mathcal{T}^{-1}_N\left(\sum_{|k|\ge L+1}
\mathcal{D}^N_k\mathcal{D}_k\right)f\right\|_{{\mathcal G}(\beta,\gamma)}
&\lesssim\lim_{L\to\infty}\sum_{|k|\ge L+1}\left\|\mathcal{D}^N_k
\mathcal{D}_kf\right\|_{{\mathcal G}(\beta,\gamma)}
\nonumber\\
&\lesssim\lim_{L\to\infty}\sum_{|k|\ge L+1}2^{-\sigma_0|k|}
\|f\|_{{\mathcal G}(\beta',\gamma')}=0,
\end{align*}
which, together with \eqref{jj-4u5} and the fact that both
\begin{align*}
f-\sum_{|k|\leq L}
\widetilde{\mathcal{D}}_k\mathcal{D}_kf
\ \ \text{and}\ \
\mathcal{T}^{-1}_N\left(\sum_{|k|\ge L+1}
\mathcal{D}^N_k\mathcal{D}_k\right)f
\end{align*}
are continuous functions,
further implies that \eqref{3.13} holds; that is,
the first series in \eqref{3.12}
converges in $\mathring{{\mathcal G}}_0^\varepsilon(\beta,\gamma)$
for any $f\in\mathring{\mathcal{G}}(\beta',\gamma')$
with $\beta'\in(\beta,\varepsilon)$
and $\gamma'\in(\gamma,\varepsilon)$.

Thus, it remains to prove  \eqref{3.16}.
Indeed, we need to show that there exists $\sigma\in(0,\infty)$
such that, for any
$k\in\mathbb{Z}$ with $|k|\in[L+1,\infty)$ and for any $x\in X$,
\begin{equation}\label{3.17}
\left|\mathcal{D}^N_k\mathcal{D}_kf(x)\right|
\lesssim 2^{-\sigma|k|}\|f\|_{{\mathcal G}(\beta',\gamma')}
\frac{1}{(V_\rho)_1(x_1)+V_\rho(x_1,x)}
\left[\frac{1}{1+\rho(x,x_1)}\right]^\gamma
\end{equation}
and, for any $x,x'\in X$
with $\rho(x,x')\leq\frac{1+\rho(x_1,x)}{2C_\rho}$,
\begin{align}
\label{3.18}
&\left|\mathcal{D}^N_k\mathcal{D}_kf(x)
-\mathcal{D}^N_k\mathcal{D}_kf(x')\right|
\nonumber\\
&\quad\lesssim\|f\|_{{\mathcal G}(\beta',\gamma')}
\left[\frac{\rho(x,x')}{1+\rho(x,x_1)}\right]^{\beta'}
\frac{1}{(V_\rho)_1(x_1)+V_\rho(x_1,x)}
\left[\frac{1}{1+\rho(x,x_1)}\right]^\gamma.
\end{align}
Indeed, if both \eqref{3.17} and \eqref{3.18} hold,
then, using the geometric mean,
we conclude that, for any $x,x'\in X$
with $\rho(x,x')\leq\frac{1+\rho(x_1,x)}{2C_\rho}$,
\begin{align*}
&\left|\mathcal{D}^N_k\mathcal{D}_kf(x)
-\mathcal{D}^N_k\mathcal{D}_kf(x')\right|
\\
&\quad\leq\left|\mathcal{D}^N_k\mathcal{D}_kf(x)
-\mathcal{D}^N_k\mathcal{D}_kf(x')\right|^{\frac{\beta}{\beta'}}
\left[\left|\mathcal{D}^N_k\mathcal{D}_kf(x)\right|
+\left|\mathcal{D}^N_k\mathcal{D}_kf(x')
\right|\right]^{1-\frac{\beta}{\beta'}}
\\
&\quad\lesssim2^{-\sigma(1-\frac{\beta}{\beta'})|k|}
\|f\|_{{\mathcal G}(\beta',\gamma')}
\left[\frac{\rho(x,x')}{1+\rho(x,x_1)}\right]^{\beta}
\frac{1}{(V_\rho)_1(x_1)+V_\rho(x_1,x)}
\left[\frac{1}{1+\rho(x,x_1)}\right]^\gamma,
\end{align*}
which, combined with \eqref{3.17}, further implies that
\eqref{3.16} holds
with $\sigma_0:=\sigma(1-\frac{\beta}{\beta'})$.

We first prove \eqref{3.18}. From Fubini's theorem and Lemma~\ref{Pk}(v),
we infer that, for any $k\in\mathbb{Z}$ and $x\in X$,
\begin{align}
\label{921}
\int_XD_k^ND_k(x,y)\,d\mu(y)
=\int_XD_k^N(x,z)\left[\int_XD_k(z,y)\,d\mu(y)\right]\,d\mu(z)
=0.
\end{align}
By the definition of $\mathcal{D}_k^N\mathcal{D}_k$
and Theorem~\ref{corATI}(i) with $p=2$,
we find that, for any $k\in\mathbb{Z}$,
$\mathcal{D}_k^N\mathcal{D}_k$ is bounded on $L^2(X)$ with
its operators norm bounded by a positive constant independent of $k$.
From this, an argument similar to that used in the proof
of Lemma~\ref{A3.15},
Definition~\ref{DefATI},
and \eqref{921},
we deduce that, for any $k\in\mathbb{Z}$,
$\mathcal{D}_k^N\mathcal{D}_k$
satisfies all the conditions of Theorem~\ref{CZOonGO}.
Thus, by Theorem~\ref{CZOonGO},
we find that, for any $k\in\mathbb{Z}$,
$\mathcal{D}_k^N\mathcal{D}_k$ is
bounded on $\mathring{\mathcal{G}}(\beta',\gamma')$,
which, together with Definition~\ref{testfunction}(ii) and $\gamma'>\gamma$,
further implies that \eqref{3.18} holds.

To show \eqref{3.17}, for any $k\in\mathbb{Z}$, let
\begin{align*}
\mathcal{E}_k:=\mathcal{D}^N_k\mathcal{D}_k.
\end{align*}
We claim that the kernels $\{E_k\}_{k\in\mathbb{Z}}$
of $\{\mathcal{E}_k\}_{k\in\mathbb{Z}}$
satisfy (i)--(iii) and (v) of Lemma~\ref{Pk}
with $\varepsilon$ therein replaced by any
given constant $\varepsilon'\in(\gamma',\varepsilon)$.
We first verify that the kernels $\{E_k\}_{k\in\mathbb{Z}}$
satisfy Lemma~\ref{Pk}(i).
On the one hand,
from Lemma~\ref{A3.12}(i), it follows that,
for any $k\in\mathbb{Z}$ and $x,y\in X$,
\begin{align}\label{sizeEk}
\left|E_k(x,y)\right|
&\leq\sum_{|l|<N}\left|(D_{k+l}D_k)(x,y)\right|\nonumber\\
&\lesssim\sum_{0\leq l<N}2^{-l\varepsilon}
\frac{1}{(V_\rho)_{2^{-k}}(x)}
+\sum_{-N<l<0}2^{l\varepsilon}
\frac{1}{(V_\rho)_{2^{-k-l}}(x)}\nonumber\\
&\leq\frac{1}{(V_\rho)_{2^{-k}}(x)}\sum_{|l|<N}2^{-|l|\varepsilon}
\sim\frac{1}{(V_\rho)_{2^{-k}}(x)}.
\end{align}
On the other hand,
by Lemma~\ref{A3.12}(i), we conclude that there exists a constant $C'\in(0,\infty)$
such that, for any $x,y\in X$ with $\rho(x,y)>C'2^{-k+N}$,
\begin{align}\label{Ekbound}
E_k(x,y)=0.
\end{align}
Let $C:=C'2^N$.
Then we conclude that the kernels $\{E_k\}_{k\in\mathbb{Z}}$
satisfy Lemma~\ref{Pk}(i).

We verify that the kernels $\{E_k\}_{k\in\mathbb{Z}}$
satisfy Lemma~\ref{Pk}(ii). From Lemma~\ref{A3.12}(iii),
we infer that, for any $k\in\mathbb{Z}$ and $x,y,x'\in X$,
\begin{align*}
\left|E_k(x,y)-E_k(x',y)\right|
&=\left|\sum_{|l|<N}
\left[(D_{k+l}D_k)(x,y)-(D_{k+l}D_k)(x',y)\right]\right|\nonumber\\
&\lesssim\sum_{0\leq l<N}
2^{-l\delta}2^{k\varepsilon'}[\rho(x,x')]^{\varepsilon'}
\frac{1}{(V_\rho)_{2^{-k}}(y)}\nonumber\\
&\quad+\sum_{-N<l<0}2^{l\delta}2^{(k+l)\varepsilon'}
[\rho(x,x')]^{\varepsilon'}
\frac{1}{(V_\rho)_{2^{-k-l}}(y)},
\end{align*}
which further implies that
\begin{align}\label{1549}
\left|E_k(x,y)-E_k(x',y)\right|
&\lesssim 2^{k\varepsilon'}
[\rho(x,x')]^{\varepsilon'}
\frac{1}{(V_\rho)_{2^{-k}}(y)}
\sum_{0\leq l<N}2^{-l\delta}\nonumber\\
&\quad+2^{k\varepsilon'}[\rho(x,x')]^{\varepsilon'}
\frac{1}{(V_\rho)_{2^{-k}}(y)}
\sum_{-N<l<0}2^{l\delta}2^{l\varepsilon'}\nonumber\\
&\lesssim 2^{k\varepsilon'}
[\rho(x,x')]^{\varepsilon'}\frac{1}{(V_\rho)_{2^{-k}}(y)},
\end{align}
where the constant $\delta\in(0,\varepsilon-\varepsilon']$ is
the same as in Lemma~\ref{A3.12},
which further implies that the kernels $\{E_k\}_{k\in\mathbb{Z}}$
satisfy Lemma~\ref{Pk}(ii) with $\varepsilon$ therein
replaced by $\varepsilon'$.
But, for any given $k\in\mathbb{Z}$, it is unclear
whether $E_k(x,y)=E_k(y,x)$ for any $x,y\in X$,
and hence Lemma~\ref{Pk}(iv) may not hold for $E_k$.
However, we still have the following regularity of $E_k(\cdot,y)$
with respect to the variable $y$. Indeed,
by an argument similar to that used in the estimation
of \eqref{1549} with Lemma~\ref{A3.12}(iii) replaced by Lemma~\ref{A3.12}(ii),
we also have, for any $k\in\mathbb{Z}$ and $x,y,y'\in X$,
\begin{align}\label{yjyj}
\left|E_k(x,y)-E_k(x,y')\right|
\lesssim 2^{k\varepsilon'}
[\rho(y,y')]^{\varepsilon'}\frac{1}{(V_\rho)_{2^{-k}}(x)}.
\end{align}

Then we verify that the kernels $\{E_k\}_{k\in\mathbb{Z}}$
satisfy Lemma~\ref{Pk}(iii). To this end,
from Lemma~\ref{A3.12}(iv), we deduce that,
for any $k\in\mathbb{Z}$ and $x,y,x',y'\in X$,
\begin{align*}
&\left|\left[E_k(x,y)-E_k(x',y)\right]
-\left[E_k(x,y')-E_k(x',y')\right]\right|\\
&\quad\leq\sum_{0\leq l<N}
\left|\left[(D_{k+l}D_k)(x,y)
-(D_{k+l}D_k)(x',y)\right]
-\left[(D_{k+l}D_k)(x,y')-(D_{k+l}D_k)(x',y')\right]\right|\\
&\quad\quad+\sum_{-N<l<0}\left|\left[(D_{k+l}D_k)(x,y)
-(D_{k+l}D_k)(x',y)\right]
-\left[(D_{k+l}D_k)(x,y')-(D_{k+l}D_k)(x',y')\right]\right|\\
&\quad\lesssim\sum_{0\leq l<N}
2^{-l\delta}2^{2k\varepsilon'}
[\rho(x,x')]^{\varepsilon'}[\rho(y,y')]^{\varepsilon'}
\left[\frac{1}{(V_\rho)_{2^{-k}}(x)}
+\frac{1}{(V_\rho)_{2^{-k}}(x')}\right]\\
&\qquad+\sum_{-N<l<0}2^{l\delta}2^{2(k+l)\varepsilon'}
[\rho(x,x')]^{\varepsilon'}[\rho(y,y')]^{\varepsilon'}
\left[\frac{1}{(V_\rho)_{2^{-k-l}}(x)}
+\frac{1}{(V_\rho)_{2^{-k-l}}(x')}\right].
\end{align*}
Using this, we obtain
\begin{align*}
&\left|\left[E_k(x,y)-E_k(x',y)\right]
-\left[E_k(x,y')-E_k(x',y')\right]\right|\\
&\quad\lesssim 2^{2k\varepsilon'}
[\rho(x,x')]^{\varepsilon'}
[\rho(y,y')]^{\varepsilon'}
\left[\frac{1}{(V_\rho)_{2^{-k}}(x)}
+\frac{1}{(V_\rho)_{2^{-k}}(x')}\right]
\sum_{0\leq l<N}2^{-l\delta}\\
&\quad\quad+2^{2k\varepsilon'}
[\rho(x,x')]^{\varepsilon'}
[\rho(y,y')]^{\varepsilon'}
\left[\frac{1}{(V_\rho)_{2^{-k}}(x)}
+\frac{1}{(V_\rho)_{2^{-k}}(x')}\right]
\sum_{-N<l<0}2^{l\delta}2^{2l\varepsilon'}\\
&\quad\lesssim 2^{2k\varepsilon'}
[\rho(x,x')]^{\varepsilon'}
[\rho(y,y')]^{\varepsilon'}
\left[\frac{1}{(V_\rho)_{2^{-k}}(x)}
+\frac{1}{(V_\rho)_{2^{-k}}(x')}\right].
\end{align*}
From this,
it follows that the kernels $\{E_k\}_{k\in\mathbb{Z}}$ satisfy
Lemma~\ref{Pk}(iii) with $\varepsilon$ therein
replaced by $\varepsilon'$.

Finally, we verify that $E_k$ for any given $k\in\mathbb{Z}$
satisfies Lemma~\ref{Pk}(v).
Indeed, by Fubini's theorem and Lemma~\ref{Pk}(v),
we conclude that, for any $k\in\mathbb{Z}$ and $x\in X$,
\begin{align}
\label{tqa-38}
\int_XE_k(x,y)\,d\mu(y)&=\int_X\int_XD_k^N(x,z)
D_k(z,y)\,d\mu(z)\,d\mu(y)\nonumber\\
&=\int_XD_k^N(x,z)\left[\int_XD_k(z,y)\,d\mu(y)\right]\,d\mu(z)=0.
\end{align}
Hence $E_k$ for any given $k\in\mathbb{Z}$ satisfies
Lemma~\ref{Pk}(v), which completes the proof of the above claim.

Having just proven
that the kernels $\{E_k\}_{k\in\mathbb{Z}}$ satisfy
Lemma~\ref{Pk}(i)--(iii) and (v), we now show \eqref{3.17}.
Let $k\in\mathbb{Z}$ with $|k|\in[L+1,\infty)$ and $x\in X$.
Note that we may assume $\|f\|_{{\mathcal G}(\beta',\gamma')}=1$
and $L\in\mathbb{N}$ is large enough so that $C2^{-L-1}<(2C_\rho)^{-1}$.
From \eqref{tqa-38} and the support condition of $E_k$
[see \eqref{Ekbound}], we infer that
\begin{align}
\label{tqa-39}
\left|\mathcal{D}^N_k\mathcal{D}_kf(x)\right|
&=\left|\int_X
E_k(x,y)[f(y)-f(x)]\,d\mu(y)\right|\nonumber\\
&\leq\int_{\rho(x,y)\leq C2^{-k}}
|E_k(x,y)||f(y)-f(x)|\,d\mu(y).
\end{align}
To proceed, we consider two cases beginning
with the case where $k\in\mathbb{Z}\cap[L+1,\infty)$.
Under this assumption we have, for any $y\in B_\rho(x,C2^{-k})$,
$$
\rho(x,y)<C2^{-L-1}<\frac{1+\rho(x_1,x)}{2C_\rho}
\ \ \text{and}\ \
V_\rho(x,y)\lesssim(V_\rho)_{2^{-k}}(x).
$$
By this, \eqref{tqa-39}, \eqref{sizeEk},
the regularity condition of $f\in\mathring{\mathcal{G}}(\beta',\gamma')$,
and Lemma~\ref{A3.2}(i), we obtain
\begin{align}
\label{Ek(f)(x)}
\left|\mathcal{D}^N_k\mathcal{D}_kf(x)\right|
&\lesssim\int_{\rho(x,y)\leq C2^{-k}}
\frac{1}{(V_\rho)_{2^{-k}}(x)}\left[\frac{\rho(x,y)}
{1+\rho(x_1,x)}\right]^{\beta'}\nonumber\\
&\quad\times\frac{1}{(V_\rho)_1(x_1)+V_\rho(x_1,x)}
\left[\frac{1}{1+\rho(x_1,x)}\right]^{\gamma'}\,d\mu(y)\nonumber\\
&\lesssim\int_{\rho(x,y)\leq C2^{-k}}
\frac{[\rho(x,y)]^{\beta'}}{V_\rho(x,y)}\,d\mu(y)
\frac{1}{(V_\rho)_1(x_1)+V_\rho(x_1,x)}
\left[\frac{1}{1+\rho(x_1,x)}\right]^{\gamma'}\nonumber\\
&\lesssim 2^{-k\beta'}\frac{1}{(V_\rho)_1(x_1)+V_\rho(x_1,x)}
\left[\frac{1}{1+\rho(x_1,x)}\right]^{\gamma'}\nonumber\\
&\leq 2^{-k\beta'}\frac{1}{(V_\rho)_1(x_1)+V_\rho(x_1,x)}
\left[\frac{1}{1+\rho(x_1,x)}\right]^{\gamma},
\end{align}
where, in the last inequality, we used the fact that $\gamma'>\gamma$.
This finishes the proof of \eqref{3.17}
in the case where $k\in\mathbb{Z}\cap[L+1,\infty)$.

Next, we consider the case where $k\in\mathbb{Z}\cap(-\infty,-L-1]$.
Since $f\in\mathring{\mathcal{G}}(\beta',\gamma')$, it follows that
$\int_Xf(y)\,d\mu(y)=0$, and hence we can write
\begin{align}\label{2135}
\left|\mathcal{D}^N_k\mathcal{D}_kf(x)\right|
&=\left|\int_X\left[E_k(x,y)-E_k(x,x_1)
\right]f(y)\,d\mu(y)\right|\nonumber\\
&\leq\int_{\rho(x_1,y)\leq C_\rho C2^{-k}}
\left|E_k(x,y)-E_k(x,x_1)\right||f(y)|\,d\mu(y)\nonumber\\
&\quad+\int_{\rho(x_1,y)>C_\rho C2^{-k}}
\left|E_k(x,y)\right||f(y)|\,d\mu(y)\nonumber\\
&\quad+\left|E_k(x,x_1)\right|\int_{\rho(x_1,y)>C_\rho C2^{-k}}
|f(y)|\,d\mu(y)\nonumber\\
&=:\sum_{i=1}^{3}I_{i,k}(x).
\end{align}

For $I_{1,k}(x)$, consider the set
$$
Y:=\left\{y\in X:\rho(x_1,y)\leq C_\rho C2^{-k}\mbox{ and }
\left|E_k(x,y)-E_k(x,x_1)\right|\neq 0\right\}.
$$
If $Y=\varnothing$, then $I_{1,k}(x)=0$ automatically.
Thus, in the reminder of the estimate of $I_{1,k}(x)$,
we may assume $Y\neq\varnothing$.
To estimate $I_{1,k}(x)$ in this case,
we first claim that
\begin{equation}\label{rho(x1,x)}
\rho(x_1,x)\leq\left(C^2_\rho C+1\right)2^{-k}.
\end{equation}
Indeed, by the support condition of $E_k$
[see \eqref{Ekbound}],
we find that
$E_k(x,y)=0$ when $\rho(x,y)>C2^{-k}$,
and $E_k(x,x_1)=0$ when $\rho(x,x_1)>C2^{-k}.$
As such, if $y\in Y$ and
$\rho(x,x_1)>(C^2_\rho C+1)2^{-k}>C2^{-k}$, then
$$
|E_k(x,y)|=\left|E_k(x,y)-E_k(x,x_1)\right|\neq 0,
$$
which further implies that $\rho(x,y)\leq C2^{-k}$.
From this and Definition~\ref{D2.1}(iii), we deduce that
\begin{align*}
\left(C^2_\rho C+1\right)2^{-k}<\rho(x,x_1)
\leq C_\rho\max\left\{\rho(x,y),\,\rho(y,x_1)\right\}
\leq C^2_\rho C2^{-k},
\end{align*}
which is a clear contradiction.
Thus, \eqref{rho(x1,x)} holds whenever $Y\neq\varnothing$,
and hence the proof of this claim is finished.

Now, choose $\widetilde{\gamma}$ such that
$$
0<\gamma<\widetilde{\gamma}<\gamma'<\varepsilon'<\varepsilon.
$$
By this, \eqref{rho(x1,x)}, and
the size condition of $f\in\mathring{{\mathcal G}}(\beta',\gamma')$,
\eqref{yjyj},
we conclude that
\begin{align*}
I_{1,k}(x)
&\lesssim\int_{Y}
\left[\frac{\rho(x_1,y)}{2^{-k}}\right]^{\varepsilon'}
\frac{1}{(V_\rho)_{2^{-k}}(x)}\nonumber\\
&\quad\times\frac{1}{(V_\rho)_1(x_1)+V_\rho(x_1,y)}
\left[\frac{1}{1+\rho(x_1,y)}\right]^{\gamma'}\,d\mu(y)
\nonumber\\
&\lesssim\int_{Y}
\left[\frac{\rho(x_1,y)}{2^{-k}}\right]^{\widetilde{\gamma}}
\frac{1}{(V_\rho)_{2^{-k}}(x)}\nonumber\\
&\quad\times
\frac{1}{(V_\rho)_1(x_1)+V_\rho(x_1,y)}
\left[\frac{1}{1+\rho(x_1,y)}\right]^{\gamma'}\,d\mu(y),
\end{align*}
which, combined with Lemma~\ref{A3.2}(ii),
further implies that
\begin{align}\label{1224}
I_{1,k}(x)
&\lesssim\int_{Y}
\left[\frac{1+\rho(x_1,y)}{2^{-k}}\right]^{\widetilde{\gamma}}
\frac{1}{(V_\rho)_{2^{-k}}(x)}\nonumber\\
&\quad\times
\frac{1}{(V_\rho)_1(x_1)+V_\rho(x_1,y)}
\left[\frac{1}{1+\rho(x_1,y)}\right]^{\gamma'}\,d\mu(y)
\nonumber\\
&\lesssim\frac{2^{k\widetilde{\gamma}}}{(V_\rho)_{2^{-k}}(x)}
\int_{X}\frac{1}{(V_\rho)_1(x_1)+V_\rho(x_1,y)}
\left[\frac{1}{1+\rho(x_1,y)}\right]^
{\gamma'-\widetilde{\gamma}}\,d\mu(y)
\nonumber\\
&\lesssim\frac{2^{k\widetilde{\gamma}}}{(V_\rho)_{2^{-k}}(x)}.
\end{align}
Moreover,
from \eqref{rho(x1,x)} and $2^{-k}\geq1$, we infer that
\begin{align}\label{x1x}
\frac{1}{(V_\rho)_{2^{-k}}(x)}
&\sim\frac{1}{(V_\rho)_{2^{-k}}(x)}
\left[\frac{1}{1+C_\rho^2C+1}\right]^\gamma\nonumber\\
&\leq\frac{1}{(V_\rho)_{2^{-k}}(x)}
\left[\frac{2^{-k}}{2^{-k}+\rho(x_1,x)}\right]^\gamma
\nonumber\\
&\lesssim 2^{-k\gamma}
\frac{1}{(V_\rho)_{C_\rho[1+\rho(x_1,x)]}(x)}
\left[\frac{1}{1+\rho(x_1,x)}\right]^{\gamma}
\quad\text{by \eqref{upperdoub}}
\nonumber\\
&\lesssim2^{-k\gamma}
\frac{1}{(V_\rho)_1(x_1)+V_\rho(x_1,x)}
\left[\frac{1}{1+\rho(x_1,x)}\right]^{\gamma}
\quad\text{by Lemma~\ref{A3.2}(vi)},
\end{align}
which, together with \eqref{1224}, further implies that
\begin{align}\label{I1k}
I_{1,k}(x)
&\lesssim
2^{(\widetilde{\gamma}-\gamma)k}
\frac{1}{(V_\rho)_1(x_1)+V_\rho(x_1,x)}
\left[\frac{1}{1+\rho(x_1,x)}\right]^{\gamma}.
\end{align}
This is the desired estimate.

For $I_{2,k}(x)$, if $E_k(x,y)\neq 0$, then,
by the support condition of $E_k$
[see \eqref{Ekbound}], we find that
$\rho(x,y)\leq C2^{-k}$.
From this and Definition~\ref{D2.1}(iii), it follows that,
if $\rho(x_1,y)>C_\rho C2^{-k}$, then
\begin{align*}
C_\rho C2^{-k}
&<\rho(x_1,y)
\leq C_\rho\max\left\{\rho(x_1,x),\,\rho(x,y)\right\}\\
&\leq C_\rho\max\left\{\rho(x_1,x),\,C2^{-k}\right\}=C_\rho\rho(x_1,x).
\end{align*}
Hence, $\rho(x_1,x)>C2^{-k}$.
By this, Definition~\ref{D2.1}(iii),
and $\rho(x,y)\leq C2^{-k}$,
we conclude that
\begin{align*}
\rho(x_1,y)
&\leq C_\rho\max\left\{\rho(x_1,x),\,\rho(x,y)\right\}
=C_\rho\rho(x_1,x)\\
&\leq C^2_\rho\max\left\{\rho(x_1,y),\,\rho(x,y)\right\}
=C^2_\rho\rho(x_1,y),
\end{align*}
which further implies that $\rho(x_1,x)\sim \rho(x_1,y)$.
From this and the size conditions of both $E_k(x,y)$ and $f$,
we deduce that
\begin{align*}
I_{2,k}(x)
&\lesssim\int_{\{y:\rho(x,y)<C2^{-k}
\text{ and }\rho(x_1,y)>C_\rho C2^{-k}\}}
\frac{1}{(V_\rho)_1(x_1)+V_\rho(x_1,y)}\nonumber\\
&\quad\times\left[\frac{1}{1+\rho(x_1,y)}\right]^{\gamma'}
\frac{d\mu(y)}{(V_\rho)_{2^{-k}}(x)}
\nonumber\\
&\sim \int_{\{y:\rho(x,y)<C2^{-k}
\text{ and }\rho(x_1,y)>C_\rho C2^{-k}\}}
\frac{1}{(V_\rho)_{2^{-k}}(x)}\nonumber\\
&\quad\times\frac{1}{(V_\rho)_1(x_1)+V_\rho(x_1,x)}
\left[\frac{1}{1+\rho(x_1,x)}\right]^{\gamma}
\left[\frac{1}{1+\rho(x_1,y)}\right]^{\gamma'-\gamma}\,d\mu(y),
\end{align*}
which further implies that
\begin{align}
\label{I2k}
I_{2,k}(x)
&\lesssim2^{(\gamma'-\gamma)k}
\frac{1}{(V_\rho)_{2^{-k}}(x)}
\frac{1}{(V_\rho)_1(x_1)+V_\rho(x_1,x)}
\left[\frac{1}{1+\rho(x_1,x)}\right]^{\gamma}
\nonumber\\
&\quad\times
\int_{\rho(x,y)<C2^{-k}}
\left[\frac{\rho(x_1,y)}{1+\rho(x_1,y)}
\right]^{\gamma'-\gamma}\,d\mu(y)
\nonumber\\
&\lesssim2^{(\gamma'-\gamma)k}
\frac{1}{(V_\rho)_1(x_1)+V_\rho(x_1,x)}
\left[\frac{1}{1+\rho(x_1,x)}\right]^{\gamma}
\frac{(V_\rho)_{C2^{-k}}(x)}{(V_\rho)_{2^{-k}}(x)}
\quad\text{since $\gamma'>\gamma$}
\nonumber\\
&\lesssim2^{(\gamma'-\gamma)k}
\frac{1}{(V_\rho)_1(x_1)+V_\rho(x_1,x)}
\left[\frac{1}{1+\rho(x_1,x)}\right]^{\gamma}
\quad\text{by \eqref{upperdoub}}.
\end{align}

For $I_{3,k}(x)$,
if $E_k(x,x_1)\neq 0$,
then, by the support condition of $E_k$ [see \eqref{Ekbound}],
we find that $\rho(x_1,x)\leq C2^{-k}$.
As such, via arguing as in \eqref{x1x}, we obtain
\begin{align*}
\frac{1}{(V_\rho)_{2^{-k}}(x)}
\lesssim 2^{-k\gamma}
\frac{1}{(V_\rho)_1(x_1)+V_\rho(x_1,x)}
\left[\frac{1}{1+\rho(x_1,x)}\right]^{\gamma}.
\end{align*}
From this, the fact that
$\rho(x_1,x)<C2^{-k}<\rho(x_1,y)$, and
the size conditions of both $E_k$ and $f$,
we infer that
\begin{align*}
I_{3,k}(x)
&\lesssim\int_{\rho(x_1,y)>C_\rho C2^{-k}}
\left[\frac{1+\rho(x_1,y)}{2^{-k}}\right]^{\widetilde{\gamma}}
\frac{1}{(V_\rho)_{2^{-k}}(x)}
\nonumber\\
&\quad\times
\frac{1}{(V_\rho)_1(x_1)+V_\rho(x_1,y)}
\left[\frac{1}{1+\rho(x_1,y)}\right]^{\gamma'}\,d\mu(y)
\nonumber\\
&\lesssim2^{k\widetilde{\gamma}}
\frac{1}{(V_\rho)_{2^{-k}}(x)}
\nonumber\\
&\quad\times\int_{\rho(x_1,y)>C_\rho C2^{-k}}
\frac{1}{(V_\rho)_1(x_1)+V_\rho(x_1,y)}
\left[\frac{1}{1+\rho(x_1,y)}\right]
^{\gamma'-\widetilde{\gamma}}\,d\mu(y),
\end{align*}
which, combined with both $\gamma'>\widetilde{\gamma}$
and Lemma~\ref{A3.2}(ii),
further implies that
\begin{align}\label{I3k}
I_{3,k}(x)
&\lesssim2^{(\widetilde{\gamma}-\gamma)k}
\frac{1}{(V_\rho)_1(x_1)+V_\rho(x_1,x)}
\left[\frac{1}{1+\rho(x_1,x)}\right]^{\gamma}
\nonumber\\
&\quad\times\int_X
\frac{1}{(V_\rho)_1(x_1)+V_\rho(x_1,y)}
\left[\frac{1}{1+\rho(x_1,y)}\right]
^{\gamma'-\widetilde{\gamma}}\,d\mu(y)
\nonumber\\
&\lesssim2^{(\widetilde{\gamma}-\gamma)k}
\frac{1}{(V_\rho)_1(x_1)+V_\rho(x_1,x)}
\left[\frac{1}{1+\rho(x_1,x)}\right]^{\gamma}.
\end{align}
Therefore, by \eqref{2135},
\eqref{I1k}, \eqref{I2k}, and \eqref{I3k},
we conclude that
$$
\left|\mathcal{D}^N_k\mathcal{D}_kf(x)\right|
=\sum_{i=1}^{3}
I_{i,k}(x)
\lesssim2^{k\sigma}
\frac{1}{(V_\rho)_1(x_1)+V_\rho(x_1,x)}
\left[\frac{1}{1+\rho(x_1,x)}\right]^{\gamma}
$$
with $\sigma:=\gamma'-\gamma\in(0,\infty)$.
Hence, \eqref{3.17} also holds in the case $k\in\mathbb{Z}\cap(-\infty,-L-1]$,
which, together with \eqref{Ek(f)(x)},
completes the proof of \eqref{3.17}
and hence \eqref{3.16}.

\emph{Step 3.}
Having just proven in Step 2
that the first series in \eqref{3.12} holds
whenever $f\in\mathring{\mathcal{G}}(\beta',\gamma')$
with $\beta'\in(\beta,\varepsilon)$ and $\gamma'\in(\gamma,\varepsilon)$,
in this step,
we return to proving that
the first series in \eqref{3.12} converges in
$\mathring{{\mathcal G}}_0^\varepsilon(\beta,\gamma)$
for any $f\in\mathring{{\mathcal G}}_0^\varepsilon(\beta,\gamma)$.
Let $f\in\mathring{{\mathcal G}}_0^\varepsilon(\beta,\gamma)$.
Then there exists a sequence
$\{f_m\}^\infty_{m=1}$ in $\mathring{\mathcal{G}}(\varepsilon,\varepsilon)$
such that
$\lim_{m\to\infty}\|f-f_m\|_{{\mathcal G}(\beta,\gamma)}=0$.
For any given $\xi\in(0,\infty)$,
choose $m\in\mathbb{N}$ such that
\begin{align}\label{1939}
\|f-f_m\|_{{\mathcal G}(\beta,\gamma)}<\xi.
\end{align}
From arguments similar to those used in the
proofs of Lemmas~\ref{A3.13} and~\ref{A3.15},
\eqref{921}, and Theorem~\ref{CZOonGO},
it is easy to verify that, for any $L\in\mathbb{N}$,
the operator
$$
\widetilde{\mathcal{T}}_L
:=\sum_{|k|\leq L}\mathcal{D}^N_k\mathcal{D}_k
$$
is bounded on $\mathring{\mathcal{G}}(\beta,\gamma)$
with its operators norm bounded by a positive constant independent of $L$,
which, combined with both Remark~\ref{Rmtest}(iii), \eqref{1939},
and Proposition~\ref{A3.16} [namely
the boundedness
of $\mathcal{T}^{-1}_N$ on $\mathring{\mathcal{G}}(\beta,\gamma)$],
further implies that,
for any $L\in\mathbb{N}$,
\begin{align}
\label{9243}
&\left\|f-\sum_{|k|\leq L}
\widetilde{\mathcal{D}}_k\mathcal{D}_kf
\right\|_{{\mathcal G}(\beta,\gamma)}\nonumber\\
&\quad=\left\|f-\mathcal{T}^{-1}_N
\widetilde{\mathcal{T}}_Lf\right\|_{{\mathcal G}(\beta,\gamma)}
\nonumber\\
&\quad\leq\|f-f_m\|_{{\mathcal G}(\beta,\gamma)}
+\left\|f_m-\mathcal{T}^{-1}_N
\widetilde{\mathcal{T}}_Lf_m\right\|_{{\mathcal G}(\beta,\gamma)}
+\left\|\mathcal{T}^{-1}_N\widetilde{\mathcal{T}}_Lf_m
-\mathcal{T}^{-1}_N\widetilde{\mathcal{T}}_Lf\right\|_{{\mathcal G}(\beta,\gamma)}
\nonumber\\
&\quad\lesssim\xi+\left\|f_m-\mathcal{T}^{-1}_N
\widetilde{\mathcal{T}}_Lf_m\right\|_{{\mathcal G}(\beta,\gamma)}.
\end{align}
By \eqref{3.13} and the fact that $f_m\in
\mathring{{\mathcal G}}(\varepsilon,\varepsilon)\subset
\mathring{{\mathcal G}}(\beta',\gamma')$, we find that
there exists $L_0\in\mathbb{N}$
such that, for any $L\in\mathbb{N}\cap(L_0,\infty)$,
\begin{align*}
\left\|f_m-\mathcal{T}^{-1}_N
\widetilde{\mathcal{T}}_Lf_m\right\|_{{\mathcal G}(\beta,\gamma)}
=\left\|f_m-\sum_{|k|\leq L}\widetilde{\mathcal{D}}_k
\mathcal{D}_kf_m\right\|_{{\mathcal G}(\beta,\gamma)}<\xi,
\end{align*}
which, together with \eqref{9243}, further implies that
\begin{equation}\label{3.19}
\left\|f-\sum_{|k|\leq L}
\widetilde{\mathcal{D}}_k\mathcal{D}_kf\right\|_{{\mathcal G}(\beta,\gamma)}
\lesssim\xi.	
\end{equation}

We now show that, for any $L\in\mathbb{N}$,
\begin{align}\label{2225}
f-\sum_{|k|\leq L}
\widetilde{\mathcal{D}}_k\mathcal{D}_kf
\in\mathring{{\mathcal G}}_0^\varepsilon(\beta,\gamma).
\end{align}
With this goal in mind, we first claim that,
for any $m,L\in\mathbb{N}$,
\begin{align}\label{fn1}
f_m-\sum_{|k|\leq L}
\widetilde{\mathcal{D}}_k
\mathcal{D}_kf_m\in\mathring{{\mathcal G}}_0^\varepsilon(\beta,\gamma).
\end{align}
Since $\{f_m\}_{m=1}^\infty\subset
\mathring{{\mathcal G}}(\varepsilon,\varepsilon)\subset
\mathring{{\mathcal G}}_0^\varepsilon(\beta,\gamma)$,
from Remark~\ref{Rmtest}(iii), it follows that,
to prove \eqref{fn1},
it suffices to
show that $\widetilde{\mathcal{D}}_k
\mathcal{D}_kf_m\in\mathring{{\mathcal G}}_0^\varepsilon(\beta,\gamma)$
for any $k\in\mathbb{Z}$ and $m\in\mathbb{N}$.
Suppose that $h\in\mathring{{\mathcal G}}_0^\varepsilon(\beta,\gamma)$.
Then there exists a sequence $\{h_j\}_{j=1}^\infty$ in
$\mathring{\mathcal{G}}(\varepsilon,\varepsilon)$
such that $\|h-h_j\|_{\mathring{\mathcal{G}}(\beta,\gamma)}\to 0$ as $j\to\infty$.
By Lemma~\ref{GOVV}, we conclude that, for any $k\in\mathbb{Z}$ and $j\in\mathbb{N}$,
$\mathcal{D}_kh_j\in\mathring{\mathcal{G}}(\varepsilon,\varepsilon)$.
Thus,
$$
\mathcal{D}_k^N\mathcal{D}_kh_j
=\sum_{|l|\leq N}\mathcal{D}_{k+l}\mathcal{D}_kh_j
\in\mathring{\mathcal{G}}(\varepsilon,\varepsilon).
$$
From this and the fact that
\begin{align*}
\left\|\sum_{|k|\leq L}\mathcal{D}_k^N\mathcal{D}_kh_j
-\sum_{k\in\mathbb{Z}}\mathcal{D}_k^N
\mathcal{D}_kh_j\right\|_{\mathring{\mathcal{G}}(\beta,\gamma)}\to0
\end{align*}
as $L\to\infty$,
we deduce that $\sum_{k\in\mathbb{Z}}\mathcal{D}_k^N\mathcal{D}_kh_j\in
\mathring{{\mathcal G}}_0^\varepsilon(\beta,\gamma)$.
The boundedness of $\mathcal{R}_N$
on $\mathring{\mathcal{G}}(\beta,\gamma)$
(given Proposition~\ref{A3.16} and our choice of $N$) and
\eqref{3.6} imply that
$\mathcal{T}_N$ is bounded on $\mathring{\mathcal{G}}(\beta,\gamma)$.
By this, we obtain
$$
\left\|\mathcal{T}_N\left(h-h_j\right)\right\|_{\mathring{\mathcal{G}}(\beta,\gamma)}
\lesssim\left\|h-h_j\right\|_{\mathring{\mathcal{G}}(\beta,\gamma)}\to0
$$
as $j\to\infty$.
From this and the proven conclusion that
$\mathcal{T}_Nh_j
=\sum_{k\in\mathbb{Z}}\mathcal{D}_k^N\mathcal{D}_kh_j\in
\mathring{{\mathcal G}}_0^\varepsilon(\beta,\gamma)$,
we infer that
$\mathcal{T}_Nh\in\mathring{{\mathcal G}}_0^\varepsilon(\beta,\gamma)$.
Hence, $\mathcal{R}_Nh=h-\mathcal{T}_Nh\in
\mathring{{\mathcal G}}_0^\varepsilon(\beta,\gamma).$
By this, the boundedness of $\mathcal{R}_N$
on $\mathring{\mathcal{G}}(\beta,\gamma)$, and the fact
that $(\mathring{{\mathcal G}}_0^\varepsilon(\beta,\gamma),
\|\cdot\|_{\mathring{\mathcal{G}}(\beta,\gamma)})$
is complete, we conclude that $\mathcal{R}_N$ is bounded on
$\mathring{{\mathcal G}}_0^\varepsilon(\beta,\gamma)$
with the same bound constant as in \eqref{3.7}.
This, combined with
$$
\mathcal{T}_N^{-1}h
=(\mathcal{I}-\mathcal{R}_N)^{-1}h
=\sum_{j=0}^\infty(\mathcal{R}_N)^jh,
$$
implies that $\mathcal{T}_N^{-1}$ is bounded
on $\mathring{{\mathcal G}}_0^\varepsilon(\beta,\gamma)$.
Hence, we obtain
$\mathcal{T}_N^{-1}h\in\mathring{{\mathcal G}}_0^\varepsilon(\beta,\gamma)$.
For any given $m\in\mathbb{N}$,
if we define $h:=\mathcal{D}_k^N\mathcal{D}_kf_m$,
then, from Lemma~\ref{GOVV}, it follows that
$h\in\mathring{\mathcal{G}}(\varepsilon,\varepsilon)
\subset\mathring{{\mathcal G}}_0^\varepsilon(\beta,\gamma)$.
Therefore, $$\mathcal{T}_N^{-1}h=
\mathcal{T}_N^{-1}\mathcal{D}_k^N\mathcal{D}_kf_m
=\widetilde{\mathcal{D}}_k\mathcal{D}_kf_m
\in\mathring{{\mathcal G}}_0^\varepsilon(\beta,\gamma).$$
This finishes the proof of \eqref{fn1}.

Observe that
\begin{align*}
&\left\|\left(f-\sum_{|k|\leq L}
\widetilde{\mathcal{D}}_k\mathcal{D}_kf\right)
-\left(f_m-\sum_{|k|\leq L}
\widetilde{\mathcal{D}}_k\mathcal{D}_kf_m\right)
\right\|_{{\mathcal G}(\beta,\gamma)}\\
&\quad\leq\|f-f_m\|_{{\mathcal G}(\beta,\gamma)}
+\left\|\mathcal{T}^{-1}_N\widetilde{\mathcal{T}}_L(f-f_m)
\right\|_{{\mathcal G}(\beta,\gamma)}
\lesssim\|f-f_m\|_{{\mathcal G}(\beta,\gamma)}\to 0
\end{align*}
as $m\to \infty$, which, together with \eqref{fn1}
and Definition~\ref{2207}, further implies that \eqref{2225} holds.
By this and \eqref{3.19},
we find that, for any $L\in\mathbb{N}\cap(L_0,\infty)$,
$$
\left\|f-\sum_{|k|\leq L}
\widetilde{\mathcal{D}}_k\mathcal{D}_kf\right\|_
{\mathring{{\mathcal G}}_0^\varepsilon(\beta,\gamma)}\lesssim\xi.
$$
Thus, the first series in \eqref{3.12} converges in
$\mathring{{\mathcal G}}_0^\varepsilon(\beta,\gamma)$
for any
$f\in\mathring{{\mathcal G}}_0^\varepsilon(\beta,\gamma)$.
Moreover,
from this, we infer that, for any given $\lambda\in(0,\infty)$,
there exists $N_{(\lambda)}\in\mathbb{N}$ such that,
for any $M\in\mathbb{N}\cap(N_{(\lambda)},\infty)$,
$$
\left\|f-\sum_{|k|\leq M}
\widetilde{\mathcal{D}}_k\mathcal{D}_kf\right\|_{
\mathring{{\mathcal G}}_0^\varepsilon(\beta,\gamma)}
<\lambda\,(V_\rho)_r(x_1),
$$
which, combined with Definition~\ref{testfunction}(i),
further implies that, for any $x\in X$,
\begin{align*}
&\left|f(x)-\sum_{|k|\leq M}
\widetilde{\mathcal{D}}_k\mathcal{D}_kf(x)\right|\\
&\quad<\lambda\,(V_\rho)_r(x_1)
\frac{1}{(V_\rho)_r(x_1)+V_\rho(x_1,x)}
\left[\frac{r}{r+\rho(x_1,x)}\right]^\gamma
<\lambda.
\end{align*}
This implies that the first series in \eqref{3.12} converges
everywhere on $X$ for any
$f\in\mathring{{\mathcal G}}_0^\varepsilon(\beta,\gamma)$.

\emph{Step 4.} In this step, we show that
the first series in \eqref{3.12} converges
in $L^p(X)$ for any $p\in(1,\infty)$.
To this end, for any given $L\in\mathbb{N}$,
we define
$$
\overline{\mathcal{T}}_L:=\sum_{|k|\leq L}
\widetilde{\mathcal{D}}_k\mathcal{D}_k
$$
with its kernel
$\overline{T}_L$ defined by setting,
for any $x,y\in X$,
\begin{align*}
\overline{T}_L(x,y):=\sum_{|k|\leq L}
\int_X\widetilde{D}_k(x,z)D_k(z,y)\,d\mu(z).
\end{align*}
By the facts that $\overline{\mathcal{T}}_L
=\mathcal{T}_N^{-1}\widetilde{\mathcal{T}}_L$
and both $\mathcal{T}_N^{-1}$ and $\widetilde{\mathcal{T}}_L$ are bounded
on $L^2(X)$ with their operator norms bounded by a positive constant
independent of $L$, we conclude that,
for any $L\in\mathbb{N}$,
$\overline{\mathcal{T}}_L$ is bounded
on $L^2(X)$ with its operator norm bounded by a positive constant
independent of $L$.

Next, fix any $\theta\in(0,\beta\land\gamma)$. We prove that
$\overline{T}_L$
is a $\theta$-Calder\'on--Zygmund kernel as in Definition~\ref{CZk}.
To this end, we first show that
$\overline{T}_L$ satisfies \eqref{CZk1}.
From Lemma~\ref{Pk}(i),
we deduce that there exists a constant $\widetilde{C}\in(1,\infty)$ such that,
for any $x,y\in X$ with $\rho(x,y)>\widetilde{C}2^{-k}$, $D_k(x,y)=0$.
Let $x,y\in X$ satisfy $x\neq y$, and define
$$
K_1:=\left\{k\in\mathbb{Z}:|k|\leq L\
\text{and}\ \rho(x,y)>C_12^{-k}\right\}
$$
and
$$
K_2:=\left\{k\in\mathbb{Z}:|k|\leq L\
\text{and}\ \rho(x,y)\leq C_12^{-k}\right\},
$$
where $C_1:=C_\rho\widetilde{C}\in[\widetilde{C},\infty)$. Then,
on the one hand,
for any $k\in K_1$ and $z\in B_\rho(y,\widetilde{C}2^{-k})$,
\begin{align*}
\rho(x,z)\gtrsim\rho(x,y)
\ \ \text{and}\ \
(V_\rho)_{2^{-k}}(z)+V_\rho(x,z)\sim V_\rho(x,y)
\end{align*}
with the implicit positive constants independent of $k$, $x$, $y$, and $z$,
which, together with Lemma~\ref{Pk}(i) and
the size condition of $\widetilde{D}_k(\cdot,z)
\in{\mathcal G}(z,2^{-k},{\beta},{\gamma})$,
further implies that
\begin{align}
\label{922}
&\left|\sum_{k\in K_1}
\int_X\widetilde{D}_k(x,z)D_k(z,y)\,d\mu(z)\right|\nonumber\\
&\quad
\lesssim\sum_{k\in K_1}\int_{B_\rho(y,\widetilde{C}2^{-k})}
\frac{1}{(V_\rho)_{2^{-k}}(z)+V_\rho(x,z)}
\left[\frac{2^{-k}}{2^{-k}+\rho(x,z)}\right]^\gamma
\frac{1}{(V_\rho)_{2^{-k}}(y)}\,d\mu(z)
\nonumber\\
&\quad\lesssim\sum_{k\in K_1}
\frac{(V_\rho)_{\widetilde{C}2^{-k}}(y)}{(V_\rho)_{2^{-k}}(y)}
\frac{1}{V_\rho(x,y)}\left[\frac{2^{-k}}{\rho(x,y)}\right]^\gamma
\nonumber\\
&\quad\lesssim\frac{1}{V_\rho(x,y)}
\frac{1}{\left[\rho(x,y)\right]^\gamma}
\sum_{k\in\mathbb{Z}\cap(\log_2\frac{C_1}{\rho(x,y)},\infty)}2^{-k\gamma}
\quad\text{by \eqref{upperdoub}}
\nonumber\\
&\quad\lesssim\frac{1}{V_\rho(x,y)},
\end{align}
where the implicit positive constants are
independent of $L$, $x$, and $y$.
On the other hand, by \eqref{upperdoub}, we find that,
for any $k\in K_2$ and $z\in B_\rho(y,\widetilde{C}2^{-k})$,
$(V_\rho)_{2^{-k}}(x)\lesssim(V_\rho)_{2^{-k}}(z).$
From this, Lemma~\ref{Pk}(i),
and the size condition of $\widetilde{D}_k(\cdot,z)
\in{\mathcal G}(z,2^{-k},{\beta},{\gamma})$,
it follows that
\begin{align*}
&\left|\sum_{k\in K_2}
\int_X\widetilde{D}_k(x,z)D_k(z,y)\,d\mu(z)\right|\\
&\quad
\lesssim\sum_{k\in K_2}\int_{B_\rho(y,\widetilde{C}2^{-k})}
\frac{1}{(V_\rho)_{2^{-k}}(z)+V_\rho(x,z)}
\left[\frac{2^{-k}}{2^{-k}+\rho(x,z)}\right]^\gamma
\frac{1}{(V_\rho)_{2^{-k}}(y)}\,d\mu(z)
\\
&\quad\leq\sum_{k\in K_2}\frac{1}{(V_\rho)_{2^{-k}}(y)}\int_{X}
\frac{1}{(V_\rho)_{2^{-k}}(x)+V_\rho(x,z)}
\left[\frac{2^{-k}}{2^{-k}+\rho(x,z)}\right]^\gamma\,d\mu(z)\\
&\qquad\qquad\text{by Lemma~\ref{A3.2}(vii)}\\
&\quad
\lesssim\sum_{k\in\mathbb{Z}\cap(-\infty,\log_2\frac{C_1}{\rho(x,y)}]}
\frac{1}{(V_\rho)_{2^{-k}}(x)}
\quad\text{by Lemma~\ref{A3.2}(ii)}\\
&\quad\lesssim\frac{1}{V_\rho(x,y)}
\quad\text{by Lemma~\ref{A3.14}},
\end{align*}
where the implicit positive constants are
independent of $L$, $x$, and $y$.
This, combined with \eqref{922}, further
implies that,
for any $x,y\in X$ with $x\neq y$,
\begin{align*}
\left|\overline{T}_L(x,y)\right|\lesssim\frac{1}{V_\rho(x,y)}.
\end{align*}
Thus, $\overline{T}_L$ satisfies \eqref{CZk1}.

Now, we prove that
$\overline{T}_L$ satisfies \eqref{CZk2}
(with $\varepsilon$ therein replaced by $\theta$).
Fix any $x,x',y\in X$ such that
$\rho(x,x')\leq\frac{\rho(x,y)}{2C_\rho}$
with $x\neq y$, and let
\begin{align*}
K_3:&=\left\{k\in\mathbb{Z}:
|k|\leq L\ \text{and}\ \rho(x,y)\leq C_22^{-k}\right\},
\\
K_4:&=\left\{k\in\mathbb{Z}:
|k|\leq L,\ \rho(x,y)>C_22^{-k},\ \text{and}\
\rho(x,x')\leq\frac{2^{-k}}{2C_\rho}\right\},
\end{align*}
and
\begin{align*}
K_5:&=\left\{k\in\mathbb{Z}:
|k|\leq L,\ \rho(x,y)>C_22^{-k},\ \text{and}\
\rho(x,x')>\frac{2^{-k}}{2C_\rho}\right\},
\end{align*}
where $C_2:=C_\rho^2\widetilde{C}\in(1,\infty)$.
Then, by Lemma~\ref{Pk}(i), we can write
\begin{align}
\label{hat-35}
&\left|\overline{T}_L(x,y)-\overline{T}_L(x',y)\right|
\nonumber\\
&\qquad\leq\sum_{k\in K_3\cup K_4\cup K_5}
\int_{B_\rho(y,\widetilde{C}2^{-k})}
\left|\widetilde{D}_k(x,z)-\widetilde{D}_k(x',z)
\right|\left|D_k(z,y)\right|\,d\mu(z).
\end{align}
Note that, if $k\in\mathbb{Z}$ and
$z\in B_\rho(y,\widetilde{C}2^{-k})$, then
Definition~\ref{D2.1}(iii) implies
\begin{align}
\label{zxp-26}
\rho(x,y)\lesssim\rho(x,z)+2^{-k}.
\end{align}
To estimate the case $k\in K_3$, from Definition~\ref{D2.1}(iii)
and \eqref{upperdoub}, we easily infer that,
for any $k\in K_3$,
\begin{align}\label{2040}
(V_\rho)_{2^{-k}}(x)\sim(V_\rho)_{2^{-k}}(y).
\end{align}
By the regularity
condition of $\widetilde{D}_k(\cdot,z)
\in{\mathcal G}(z,2^{-k},{\beta},{\gamma})$
and Lemma~\ref{Pk}(i),
we find that
\begin{align*}
&\sum_{k\in K_3}\int_{B_\rho(y,\widetilde{C}2^{-k})}
\left|\widetilde{D}_k(x,z)-\widetilde{D}_k(x',z)\right|
\left|D_k(z,y)\right|\,d\mu(z)
\nonumber\\
&\quad\lesssim\sum_{k\in K_3}
\frac{1}{(V_\rho)_{2^{-k}}(y)}
\int_{B_\rho(y,\widetilde{C}2^{-k})}
\left[\frac{\rho(x,x')}{2^{-k}+\rho(x,z)}\right]^\beta
\nonumber\\
&\qquad\times\frac{1}{(V_\rho)_{2^{-k}}(z)+V_\rho(x,z)}
\left[\frac{2^{-k}}{2^{-k}+\rho(x,z)}\right]^\gamma\,d\mu(z).
\end{align*}
From this, \eqref{2040}, \eqref{zxp-26},
and Lemma~\ref{A3.2}(vi),
it follows that
\begin{align*}
&\sum_{k\in K_3}\int_{B_\rho(y,\widetilde{C}2^{-k})}
\left|\widetilde{D}_k(x,z)-\widetilde{D}_k(x',z)\right|
\left|D_k(z,y)\right|\,d\mu(z)
\nonumber\\
&\quad\lesssim\left[\frac{\rho(x,x')}{\rho(x,y)}\right]^\beta\sum_{k\in K_3}
\frac{1}{(V_\rho)_{2^{-k}}(x)}
\nonumber\\
&\qquad\times\int_X\frac{1}{(V_\rho)_{2^{-k}}(x)+V_\rho(x,z)}
\left[\frac{2^{-k}}{2^{-k}+\rho(x,z)}\right]^\gamma\,d\mu(z),
\end{align*}
which further implies that
\begin{align}\label{K3}
&\sum_{k\in K_3}\int_{B_\rho(y,\widetilde{C}2^{-k})}
\left|\widetilde{D}_k(x,z)-\widetilde{D}_k(x',z)\right|
\left|D_k(z,y)\right|\,d\mu(z)
\nonumber\\
&\quad\lesssim\left[\frac{\rho(x,x')}{\rho(x,y)}\right]^\beta
\sum_{k\in\mathbb{Z}\cap(-\infty,\log_2
\frac{C_2}{\rho(x,y)}]}\frac{1}{(V_\rho)_{2^{-k}}(x)}
\quad\text{by Lemma~\ref{A3.2}(ii)}
\nonumber\\
&\quad\lesssim
\left[\frac{\rho(x,x')}{\rho(x,y)}\right]^\beta
\frac{1}{V_\rho(x,y)}
\quad\text{by Lemma~\ref{A3.14}}
\nonumber\\
&\quad\leq\left[\frac{\rho(x,x')}{\rho(x,y)}\right]^\theta
\frac{1}{V_\rho(x,y)}
\end{align}
because $\theta<\beta$ and $\rho(x,x')<\rho(x,y)$,
where the implicit positive constants are
independent of $L$, $x$, $x'$, and $y$.

For $K_4$ and $K_5$, first observe that, by \eqref{zxp-26},
\eqref{upperdoub}, and Lemma~\ref{A3.2}(vi), we conclude that,
for any $k\in\mathbb{Z}$ and $z\in B_\rho(y,\widetilde{C}2^{-k})$,
\begin{align}
\label{923}
V_\rho(x,y)\lesssim(V_\rho)_{2^{-k}+\rho(x,z)}(x)
\sim(V_\rho)_{2^{-k}}(z)+V_\rho(x,z).
\end{align}
If $k\in K_{4}$, then,
for any $z\in B_\rho(y,\widetilde{C}2^{-k})$,
\begin{align*}
\rho(x,x')\leq\frac{2^{-k}}{2C_\rho}\leq\frac{2^{-k}+\rho(x,z)}{2C_\rho}.
\end{align*}
From this, the regularity condition of
$\widetilde{D}_k(\cdot,z)\in{\mathcal G}(z,2^{-k},{\beta},{\gamma})$,
Lemma~\ref{Pk}(i),
\eqref{zxp-26}, and \eqref{923},
we deduce that
\begin{align*}
&\sum_{k\in K_{4}}\int_{B_\rho(y,\widetilde{C}2^{-k})}
\left|\widetilde{D}_k(x,z)-\widetilde{D}_k(x',z)\right|
\left|D_k(z,y)\right|\,d\mu(z)\nonumber\\
&\quad\lesssim\sum_{k\in K_{4}}
\int_{B_\rho(y,\widetilde{C}2^{-k})}
\left[\frac{\rho(x,x')}{2^{-k}+\rho(x,z)}\right]^\beta
\frac{1}{(V_\rho)_{2^{-k}}(z)+V_\rho(x,z)}\nonumber\\
&\qquad\times
\left[\frac{2^{-k}}{2^{-k}+\rho(x,z)}\right]^\gamma
\frac{1}{(V_\rho)_{2^{-k}}(y)}\,d\mu(z)
\nonumber\\
&\quad\lesssim\left[\frac{\rho(x,x')}{\rho(x,y)}\right]^\beta
\frac{1}{V_\rho(x,y)}
\sum_{k\in K_{4}}
\left[\frac{2^{-k}}{\rho(x,y)}\right]^\gamma
\frac{(V_\rho)_{\widetilde{C}2^{-k}}(y)}{(V_\rho)_{2^{-k}}(y)},
\end{align*}
which,
together with \eqref{upperdoub} and $\theta<\beta$,
further implies that
\begin{align}
\label{K41}
&\sum_{k\in K_{4}}\int_{B_\rho(y,\widetilde{C}2^{-k})}
\left|\widetilde{D}_k(x,z)-\widetilde{D}_k(x',z)\right|
\left|D_k(z,y)\right|\,d\mu(z)\nonumber\\
&\quad\lesssim
\left[\frac{\rho(x,x')}{\rho(x,y)}\right]^\beta
\frac{1}{V_\rho(x,y)}\sum_{k\in\mathbb{Z}\cap(\log_2\frac{C_2}{\rho(x,y)},\infty)}
\left[\frac{2^{-k}}{\rho(x,y)}\right]^\gamma\nonumber\\
&\quad\lesssim
\left[\frac{\rho(x,x')}{\rho(x,y)}\right]^\theta
\frac{1}{V_\rho(x,y)},
\end{align}
where the implicit positive constants are
independent of $L$, $x$, $x'$, and $y$.

For $K_{5}$, by Definition~\ref{D2.1}(iii) and the assumption that
$\rho(x,x')\leq\frac{\rho(x,y)}{2C_\rho}$, we find that
\begin{align*}
\rho(x,y)\leq C_\rho\max\left\{\rho(x,x'),\,\rho(x',y)\right\}
=C_\rho\rho(x',y),
\end{align*}
which, combined with Definition~\ref{D2.1}(iii), \eqref{upperdoub}, the assumption
$\rho(x,x')\leq\frac{\rho(x,y)}{2C_\rho}$,
and Lemma~\ref{A3.2}(vi), further implies that,
for any $k\in\mathbb{Z}$ and $z\in B_\rho(y,\widetilde{C}2^{-k})$,
\begin{align*}
\rho(x,y)\leq C_\rho\rho(x',y)\leq C_\rho^2\left[\rho(x',z)+\rho(z,y)\right]
\lesssim\rho(x',z)+2^{-k}
\end{align*}
and
\begin{align*}
V_\rho(x,y)\leq (V_\rho)_{C_\rho\rho(x,y)}(x')
\lesssim(V_\rho)_{2^{-k}+\rho(x',z)}(x')\sim (V_\rho)_{2^{-k}}(z)+V_\rho(x',z).
\end{align*}
From this, the size condition of $\widetilde{D}_k(\cdot,z)
\in{\mathcal G}(z,2^{-k},{\beta},{\gamma})$,
Lemma~\ref{Pk}(i), \eqref{zxp-26}, and \eqref{923},
we infer that
\begin{align*}
&\sum_{k\in K_{5}}\int_{B_\rho(y,\widetilde{C}2^{-k})}
\left[\left|\widetilde{D}_k(x,z)\right|
+\left|\widetilde{D}_k(x',z)\right|\right]
\left|D_k(z,y)\right|\,d\mu(z)\nonumber\\
&\quad\lesssim\sum_{k\in K_{5}}
\int_{B_\rho(y,\widetilde{C}2^{-k})}
\Bigg\{\frac{1}{(V_\rho)_{2^{-k}}(z)+V_\rho(x,z)}\left[
\frac{2^{-k}}{2^{-k}+\rho(x,z)}\right]^\gamma
\\
&\qquad+
\frac{1}{(V_\rho)_{2^{-k}}(z)+V_\rho(x',z)}\left[
\frac{2^{-k}}{2^{-k}+\rho(x',z)}\right]^\gamma\Bigg\}
\frac{1}{(V_\rho)_{2^{-k}}(y)}\,d\mu(z)
\\
&\quad\lesssim\sum_{k\in K_{5}}\frac{1}{V_\rho(x,y)}
\left\{\left[\frac{2^{-k}}{\rho(x,y)}\right]^\gamma
+\left[\frac{2^{-k}}{\rho(x',y)}\right]^\gamma\right\}
\frac{(V_\rho)_{\widetilde{C}2^{-k}}(y)}{(V_\rho)_{2^{-k}}(y)}.
\end{align*}
By this, the fact that $\rho(x,y)\sim\rho(x',y)$ [which can be
deduced from Definition~\ref{D2.1}(iii) and $\rho(x,x')<\rho(x,y)$],
\eqref{upperdoub}, and $\theta<\gamma$,
we conclude that
\begin{align*}
&\sum_{k\in K_{5}}\int_{B_\rho(y,\widetilde{C}2^{-k})}
\left[\left|\widetilde{D}_k(x,z)\right|
+\left|\widetilde{D}_k(x',z)\right|\right]
\left|D_k(z,y)\right|\,d\mu(z)\nonumber\\
&\quad\lesssim\frac{1}{V_\rho(x,y)}\sum_{k\in K_{5}}
\left[\frac{2^{-k}}{\rho(x,y)}\right]^\gamma\\
&\quad\lesssim\frac{1}{V_\rho(x,y)}
\left[\frac{\rho(x,x')}{\rho(x,y)}\right]^\theta
\sum_{k\in\mathbb{Z}\cap(\log_2\frac{C_2}{\rho(x,y)},\infty)}
\left[\frac{2^{-k}}{\rho(x,y)}\right]^{\gamma-\theta}\\
&\quad\sim\frac{1}{V_\rho(x,y)}
\left[\frac{\rho(x,x')}{\rho(x,y)}\right]^\theta,
\end{align*}
where the implicit positive constants are
independent of $L$.
From this, \eqref{K3}, \eqref{K41}, and \eqref{hat-35},
it follows that $\overline{T}_L$ satisfies \eqref{CZk2}.
Thus, $\overline{T}_L$
is a $\theta$-Calder\'on--Zygmund kernel as in Definition~\ref{CZk}.
This, together with the well-known Calder\'on--Zygmund
theory on spaces of homogeneous type
developed in \cite[pp.\,74--75, Theorem 2.4]{CW1},
further implies that $\overline{\mathcal{T}}_L$
is bounded on $L^p(X)$ with $p\in(1,\infty)$
and its operator norm is bounded by a positive constant independent of $L$.
Next, let $p\in(1,\infty)$ and $f\in L^p(X)$.
By Lemma~\ref{density}(ii)
(which is valid since we are currently assuming that ${\rm diam}_\rho(X)=\infty$ or,
equivalently, $\mu(X)=\infty$) and the fact that
$\mathring{{\mathcal G}}(\varepsilon,\varepsilon)\subset
\mathring{{\mathcal G}}^\varepsilon_0(\beta,\gamma)\subset L^p(X)$,
we find that $\mathring{{\mathcal G}}^\varepsilon_0(\beta,\gamma)$
is dense in $L^p(X)$ and hence,
for any given $\xi\in(0,\infty)$, there exists a function
$g\in\mathring{{\mathcal G}}^\varepsilon_0(\beta,\gamma)$ such that
\begin{align}
\label{9231}
\left\|f-g\right\|_{L^p(X)}<\xi.
\end{align}
From the fact that the first series of \eqref{3.12} converges in
$\mathring{{\mathcal G}}^\varepsilon_0(\beta,\gamma)$,
we deduce that there exists $L_0\in\mathbb{N}$
such that, for any $L\in\mathbb{N}\cap(L_0,\infty)$,
\begin{align}
\label{9232}
\left\|g-\overline{\mathcal{T}}_Lg\right\|_
{\mathring{{\mathcal G}}^\varepsilon_0(\beta,\gamma)}<\xi.
\end{align}
By the boundedness of $\overline{\mathcal{T}}_L$ on $L^p(X)$,
we conclude that,
for any $L\in\mathbb{Z}$ and $f\in L^p(X)$,
$\overline{\mathcal{T}}_Lf\in L^p(X)$.
From this, the fact that $\mathring{{\mathcal G}}^\varepsilon_0(\beta,\gamma)$
is continuously embedded into $L^p(X)$ with the
embedding constant depending on both $x_1$ and $p$
[see Lemma~\ref{A3.2}(ii)], the boundedness of $\overline{\mathcal{T}}_L$
on $L^p(X)$ with the operator norm bounded by a positive constant
independent of $L$,
\eqref{9231}, and \eqref{9232}, we infer that,
for any $L\in\mathbb{N}\cap(L_0,\infty)$,
\begin{align*}
\left\|f-\overline{\mathcal{T}}_Lf\right\|_{L^p(X)}
&\leq\|f-g\|_{L^p(X)}+
\left\|g-\overline{\mathcal{T}}_Lg\right\|_{L^p(X)}
+\left\|\overline{\mathcal{T}}_L(g-f)\right\|_{L^p(X)}
\\
&\lesssim\|f-g\|_{L^p(X)}+\left\|g-\overline{\mathcal{T}}_Lg\right\|_
{\mathring{{\mathcal G}}^\varepsilon_0(\beta,\gamma)}
<2\xi.
\end{align*}
By this,
we find that
the first series in \eqref{3.12} converges
in the sense of $L^p(X)$ with $p\in(1,\infty)$,
and hence
the first series in \eqref{3.12} also
converges pointwise $\mu$-almost everywhere on $X$
for any $f\in\bigcup_{p\in(1,\infty)}L^p(X)$.
This finishes the proof of the
convergence of the first series in \eqref{3.12}.

\emph{Step 5.} In this step, we show that
the second series in \eqref{3.12} converges
in the sense of $\mathring{\mathcal{G}}_0^\varepsilon(\beta,\gamma)$,
$L^p(X)$ with $p\in(1,\infty)$, and also $\mu$-almost everywhere.
To this end, we first need to construct $\{\overline{\mathcal{D}}_k\}_{k\in\mathbb{Z}}$.
Applying Proposition~\ref{A3.11}(ii) and some ideas
similar to those used in \eqref{3.6},
we obtain
\begin{align*}
\mathcal{I}&=\sum_{k=-\infty}^\infty\mathcal{D}_{k}
=\sum_{k=-\infty}^\infty\mathcal{D}_{k}
\sum_{l=-\infty}^\infty\mathcal{D}_{k+l}\\
&=\sum_{k=-\infty}^\infty\mathcal{D}_{k}\mathcal{D}_{k}^N
+\sum_{k=-\infty}^\infty\mathcal{D}_{k}\sum_{|l|>N}\mathcal{D}_{k+l}
=\mathcal{T}_N^*+\mathcal{R}_N^*
\end{align*}
in $L^2(X)$, where $\mathcal{T}_N^*$ and $\mathcal{R}_N^*$
are, respectively, the adjoint operators of $\mathcal{T}_N$ and $\mathcal{R}_N$
on $L^2(X)$.
From an argument quite similar to that used in the proof of Proposition~\ref{A3.16}
(since all the estimates in Lemma~\ref{A3.15}
still hold with $R_N(x,y)$ replaced by $R_N(y,x)$,
which is the kernel of $\mathcal{R}_N^*$),
we deduce that, for the same $N\in\mathbb{N}$,
Proposition~\ref{A3.16} with $\mathcal{R}_N$ and $\mathcal{T}_N$
replaced, respectively, by $\mathcal{R}_N^*$ and $\mathcal{T}_N^*$ holds.
Now, we define, for any $k\in\mathbb{Z}$,
\begin{align*}
\overline{\mathcal{D}}_k:=\mathcal{D}_k^N\left(\mathcal{T}_N^*\right)^{-1}
\end{align*}
in the sense of $L^2(X)$,
which, combined with \eqref{1617},
further implies that \eqref{wDkoDk} holds.

It remains to prove that
the second series in \eqref{3.12} converges.
Note that, for any $g\in L^2(X)$,
\begin{align}\label{2023}
g&=\mathcal{T}_N^*\left(\mathcal{T}_N^*\right)^{-1}g
=\sum_{k=-\infty}^\infty\mathcal{D}_k\mathcal{D}_k^N\left(\mathcal{T}_N^*\right)^{-1}g
=\sum_{k=-\infty}^\infty\mathcal{D}_k\overline{\mathcal{D}}_kg
\end{align}
in $L^2(X)$. By this, we conclude that,
for any given $L\in\mathbb{N}$ and for any $g\in L^2(X)$,
\begin{align}\label{2022}
\sum_{|k|\leq L}\mathcal{D}_k\overline{\mathcal{D}}_kf
&=\sum_{|k|\leq L}\mathcal{D}_k\mathcal{D}_k^N\left(\mathcal{T}_N^*\right)^{-1}f
\nonumber\\
&=\left(\mathcal{T}_N^*-\sum_{|k|\ge L+1}\mathcal{D}_k
\mathcal{D}_k^N\right)\left(\mathcal{T}_N^*\right)^{-1}f
\nonumber\\
&=f-\left(\sum_{|k|\ge L+1}\mathcal{D}_k\mathcal{D}_k^N
\right)\left(\mathcal{T}_N^*\right)^{-1}f
\end{align}
in $L^2(X)$. With \eqref{2023} and \eqref{2022}
in mind (which is used to replace \eqref{2134} and \eqref{jj-4u5}, respectively),
from the boundedness of $(\mathcal{T}_N^*)^{-1}$
on $\mathring{\mathcal{G}}(\beta,\gamma)$ and
following the proof of the first series in \eqref{3.12} as above,
we infer that the second series in \eqref{3.12} converges
in the sense of $\mathring{\mathcal{G}}_0^\varepsilon(\beta,\gamma)$,
$L^p(X)$ with $p\in(1,\infty)$, and also $\mu$-almost everywhere.
This finishes the proof of Theorem~\ref{A3.19}.
\end{proof}

By Theorem~\ref{A3.19} and a standard duality argument,
we obtain the following sharp homogeneous continuous Calder\'on
reproducing formula on spaces of distributions;
we omit the details.

\begin{theorem}
\label{HC}
Let $(X,\mathbf{q},\mu)$ be an ultra-RD-space,
where $\mu$ is assumed to be a Borel-semiregular measure on $X$.
Assume that $\varepsilon$, $\beta$, $\gamma$,
and $\{\mathcal{D}_k\}_{k\in\mathbb{Z}}$ are the same as in Theorem~\ref{A3.19}.
Then, for any
$f\in(\mathring{{\mathcal G}}_0^\varepsilon(\beta,\gamma))'$,
\eqref{3.12} holds in
$(\mathring{{\mathcal G}}_0^\varepsilon(\beta,\gamma))'$.
\end{theorem}

\begin{remark}\label{Dkf2}
Let $(X,\mathbf{q},\mu)$ be a quasi-ultrametric space of homogeneous type.
Assume that $0<\beta,\gamma<\varepsilon\preceq\mathrm{ind\,}(X,\mathbf{q})$.
Let $f\in(\mathring{{\mathcal G}}_0^\varepsilon(\beta,\gamma))'$,
$\{\mathcal{D}_k\}_{k\in\mathbb{Z}}$ be the same as in Theorem~\ref{A3.19}
with kernels denoted by $\{D_k\}_{k\in\mathbb{Z}}$,
and $\langle\cdot,\cdot\rangle$ denote the natural distributional pairing.
There are two ways to define
$\mathcal{D}_kf\in(\mathring{{\mathcal G}}_0^\varepsilon(\beta,\gamma))'$
for any given $k\in\mathbb{Z}$:
\begin{enumerate}
\item
We can define $\mathcal{D}_kf$ as the distribution induced by
the following function that is given by setting, for any $x\in X$,
$$
\mathcal{D}_k f(x):=\langle f,D_k(x,\cdot)\rangle.
$$
More precisely, in view of Lemma~\ref{dens},
for any $\phi\in\mathring{{\mathcal G}}_0^\varepsilon(\beta,\gamma)$, we let
$$
\left\langle\mathcal{D}_kf,\phi\right\rangle
:=\lim_{j\to\infty}\int_X\mathcal{D}_kf(y)\phi_j(y)\,d\mu(y),
$$
where $\{\phi_j\}_{j\in\mathbb{N}}$ is any sequence in
$\mathring{\mathcal{G}}_{\mathrm{b}}(\varepsilon,\varepsilon)
\subset\mathring{{\mathcal G}}_0^\varepsilon(\beta,\gamma)$ such that
$$
\lim_{j\to\infty}\|\phi-\phi_j\|_{\mathcal{G}(\beta,\gamma)}=0.
$$

\item
We can also define $\mathcal{D}_k f$
by duality; that is, for any
$\phi\in \mathring{{\mathcal G}}_0^\varepsilon(\beta,\gamma)$,
let
\begin{align*}
\left\langle\mathcal{D}_k f,\phi\right\rangle:=
\left\langle f,\mathcal{D}_k^*\phi\right\rangle,
\end{align*}
where $\mathcal{D}_k^*$ is the integral
operator defined by setting, for any $x\in X$,
$$
\mathcal{D}_k^*\phi(x)
:=\int_XD_k(y,x)\phi(y)\,d\mu(y).
$$
\end{enumerate}
\end{remark}

Next, we show that both of the definitions of $\mathcal{D}_k f$
for any $k\in\mathbb{Z}$ and $f\in(\mathring{{\mathcal G}}_0^\varepsilon(\beta,\gamma))'$
in (i) and (ii) of Remark~\ref{Dkf2} are well defined and actually coincide.

\begin{proposition}\label{1651}
Let the notation be as in Remark~\ref{Dkf2}.
Then, for any $k\in\mathbb{Z}$ and
$f\in(\mathring{{\mathcal G}}_0^\varepsilon(\beta,\gamma))'$,
the definitions of $\mathcal{D}_k f$
in (i) and (ii) of Remark~\ref{Dkf2} are well defined and coincide
in $(\mathring{{\mathcal G}}_0^\varepsilon(\beta,\gamma))'$.
\end{proposition}

Although the proof of Proposition~\ref{1651}
is quite similar to that of \cite[Proposition~2.8]{syy21},
we still include its proof here for
the sake of completeness.
To prove Proposition~\ref{1651},
we need two technical lemmas.

\begin{lemma}\label{welldef}
Let the notation be as in Remark~\ref{Dkf2}.
Then, for any $k\in\mathbb{Z}$ and
$f\in(\mathring{{\mathcal G}}_0^\varepsilon(\beta,\gamma))'$,
the linear functionals  $\mathcal{D}_k f$ defined
in (i) and (ii) of Remark~\ref{Dkf2} are well-defined distributions
in $(\mathring{{\mathcal G}}_0^\varepsilon(\beta,\gamma))'$.
\end{lemma}

\begin{proof}
We only show the case $k=0$ because the proof
of the case $k\neq0$ is similar
and we omit the details.
Let $f\in(\mathring{{\mathcal G}}_0^\varepsilon(\beta,\gamma))'$.
Regarding the linear functional  $\mathcal{D}_0 f$ defined
in (ii) of Remark~\ref{Dkf2}, by Lemma~\ref{GOVV}, we easily find that,
for any $\phi\in\mathring{\mathcal{G}}^\varepsilon_0(\beta,\gamma)$,
\begin{align*}
\mathcal{D}_0^*\phi\in\mathring{\mathcal{G}}^\varepsilon_0(\beta,\gamma)
\ \ \text{and}\ \
\left\|\mathcal{D}_0^*\phi\right\|_{\mathcal{G}(\beta,\gamma)}\lesssim
\|\phi\|_{\mathcal{G}(\beta,\gamma)},
\end{align*}
where the implicit positive constant depends on
$\mathcal{D}_0^*$ but is independent of $\phi$.
Thus,
\begin{align*}
\left|\left\langle f,\mathcal{D}_0^*\phi\right\rangle\right|
\leq\|f\|_{(\mathring{\mathcal{G}}^\varepsilon_0(\beta,\gamma))'}
\left\|\mathcal{D}_0^*\phi\right\|_{\mathcal{G}(\beta,\gamma)}
\lesssim\|f\|_{(\mathring{\mathcal{G}}^\varepsilon_0(\beta,\gamma))'}
\|\phi\|_{\mathcal{G}(\beta,\gamma)},
\end{align*}
and hence the linear functional $\mathcal{D}_0 f$ defined in
(ii) of Remark~\ref{Dkf2} is well-defined distribution
in $(\mathring{{\mathcal G}}_0^\varepsilon(\beta,\gamma))'$.

As concerns the linear functional $\mathcal{D}_0 f$ defined
in (i) of Remark~\ref{Dkf2}, we first prove that, for any
$\phi\in\mathring{\mathcal{G}}_{\mathrm{b}}(\varepsilon,\varepsilon)
\subset\mathring{\mathcal{G}}^\varepsilon_0(\beta,\gamma)$, the integral
\begin{align}\label{art-345}
\int_{X}\left\langle f,D_0(y,\cdot)\right\rangle
\phi(y)\,d\mu(y)
\end{align}
is well defined and finite.
Let $\phi\in\mathring{\mathcal{G}}_{\mathrm{b}}(\varepsilon,\varepsilon)$.
Without loss of generality, we may assume
$\|\varphi\|_{\mathcal{G}(\varepsilon,\varepsilon)}=1$.
First, note that the function $\mathcal{D}_0 f$, defined by setting,
for any $y\in X$,
$$
\mathcal{D}_0 f(y):=\left\langle f,D_0(y,\cdot)\right\rangle,
$$
is $\mu$-measurable (see the proof of Lemma~\ref{1955}).

To show that the integral in \eqref{art-345} is finite, note that,
from Lemma~\ref{Pk}, it follows that, for any $y\in X$,
\begin{align}\label{23112}
\left|\left\langle f,D_0(y,\cdot)\right\rangle\right|
\leq\|f\|_{(\mathring{\mathcal{G}}^\varepsilon_0(\beta,\gamma))'}
\|D_0(y,\cdot)\|_{\mathcal{G}(\beta,\gamma)}.
\end{align}
Now, we estimate the term $\|D_0(y,\cdot)\|_{\mathcal{G}(\beta,\gamma)}$.
Fix $y\in X$. For the size condition of $D_0$, by Lemma~\ref{Pk},
we find that, for any $x\in X$,
\begin{align*}
\left|D_0(y,x)\right|
&\lesssim
\frac{1}{(V_\rho)_1(y)+V_\rho(y,x)}
\left[\frac{1}{1+\rho(y,x)}\right]^\gamma\nonumber\\
&\sim\frac{1}{(V_\rho)_{1+\rho(y,x)}(y)}
\left[\frac{1}{1+\rho(y,x)}\right]^\gamma
\quad\text{by Lemma~\ref{A3.2}(vi)}\nonumber\\
&\lesssim\frac{1}{(V_\rho)_{1+\rho(x_1,y)}(y)}
\max\left\{1,\,\left[\frac{1+\rho(x_1,y)}{1+\rho(y,x)}
\right]^n\right\}\nonumber\\
&\quad\times
\left[\frac{1}{1+\rho(x_1,y)}\right]^\gamma
\left[\frac{1+\rho(x_1,y)}{1+\rho(y,x)}\right]^\gamma
\quad\text{by \eqref{upperdoub}},
\end{align*}
which further implies that
\begin{align}\label{1940}
\left|D_0(y,x)\right|
&\lesssim\frac{1}{(V_\rho)_1(x_1)+V_\rho(x_1,y)}
\max\left\{1,\,\left[1+\rho(x_1,y)\right]^n\right\}\nonumber\\
&\quad\times
\left[\frac{1}{1+\rho(x_1,y)}\right]^\gamma
\left[1+\rho(x_1,y)\right]^\gamma\nonumber\\
&=\frac{1}{(V_\rho)_1(x_1)+V_\rho(x_1,y)}
\left[\frac{1}{1+\rho(x_1,y)}\right]^\gamma
\left[1+\rho(x_1,y)\right]^{n+\gamma}.
\end{align}
For the regularity condition of $D_0(y,\cdot)$,
let $x,x'\in X$ satisfy $\rho(x,x')\leq\frac{1+\rho(x_1,x)}{2C_\rho}$.
We first prove that
\begin{align}\label{1658}
\left|D_0(y,x)-D_0(y,x')\right|
&\lesssim
\left[\frac{\rho(x,x')}{1+\rho(x,y)}\right]^\beta
\left\{
\frac{1}{(V_\rho)_1(y)+V_\rho(y,x)}
\left[\frac{1}{1+\rho(y,x)}\right]^\gamma\right.\nonumber\\
&\quad\left.+\frac{1}{(V_\rho)_1(y)+V_\rho(y,x')}
\left[\frac{1}{1+\rho(y,x')}\right]^\gamma
\right\}.
\end{align}
To this end, on the one hand, if
$\rho(x,x')\leq\frac{1+\rho(y,x)}{2C_\rho}$
and $\rho(y,x)>C_\rho(C'+1)$
with $C'$ constant as in Lemma~\ref{Pk}(i),
then, from Definition~\ref{D2.1}(iii),
we infer that
$$
C_\rho C'<\rho(y,x)\leq C_\rho\max\left\{\rho(y,x'),\,
\rho(x',x)\right\}=C_\rho\rho(y,x'),
$$
which, combined with Lemma~\ref{Pk}(i),
further implies that
\begin{align}\label{1703}
\left|D_0(y,x)-D_0(y,x')\right|=0.
\end{align}
On the other hand,
if $\rho(x,x')\leq\frac{1+\rho(y,x)}{2C_\rho}$
and $\rho(y,x)\leq C_\rho(C'+1)$,
then, by Lemma~\ref{Pk}(ii),
we conclude that
\begin{align}\label{1702}
\left|D_0(y,x)-D_0(y,x')\right|
&\lesssim[\rho(x,x')]^\varepsilon
\left[\frac{1}{(V_\rho)_1(x)}+\frac{1}{(V_\rho)_1(x')}\right]\nonumber\\
&\lesssim\left[\frac{\rho(x,x')}{1+\rho(x,y)}\right]^\beta
\left\{\frac{1}{(V_\rho)_1(y)+V_\rho(y,x)}
\left[\frac{1}{1+\rho(y,x)}\right]^\gamma\right.\nonumber\\
&\quad\left.+\frac{1}{(V_\rho)_1(y)+V_\rho(y,x')}
\left[\frac{1}{1+\rho(y,x')}\right]^\gamma\right\},
\end{align}
where the last inequality follows from
$\rho(y,x),\rho(y,x')\lesssim1$ and Lemma~\ref{A3.2}(vii).
If $\rho(x,x')>\frac{1+\rho(y,x)}{2C_\rho}$,
then, from Lemma~\ref{Pk}(vi),
we deduce that
\begin{align*}
\left|D_0(y,x)-D_0(y,x')\right|
&\leq\left|D_0(y,x)\right|+\left|D_0(y,x')\right|\\
&\lesssim\frac{1}{(V_\rho)_1(y)+V_\rho(y,x)}
\left[\frac{1}{1+\rho(y,x)}\right]^\gamma\\
&\quad+\frac{1}{(V_\rho)_1(y)+V_\rho(y,x')}
\left[\frac{1}{1+\rho(y,x')}\right]^\gamma\\
&\lesssim\left[\frac{\rho(x,x')}{1+\rho(x,y)}\right]^\beta
\left\{\frac{1}{(V_\rho)_1(y)+V_\rho(y,x)}
\left[\frac{1}{1+\rho(y,x)}\right]^\gamma\right.\\
&\quad\left.+\frac{1}{(V_\rho)_1(y)+V_\rho(y,x')}
\left[\frac{1}{1+\rho(y,x')}\right]^\gamma\right\},
\end{align*}
which, together with \eqref{1703} and \eqref{1702},
completes the proof of \eqref{1658}.

By \eqref{1658}, an argument similar to that used
in the estimation of \eqref{1940}, the obvious estimates that
$$
\frac{1+\rho(x_1,x)}{1+\rho(y,x)}
\leq\frac{1+C_\rho\rho(x_1,y)}{1+\rho(y,x)}+C_\rho
\lesssim1+\rho(x_1,y)
$$
(which also holds with $x$ replaced by $x'$),
and the estimate that
$$
\rho(x_1,x)\leq C_\rho\max\left\{\rho(x_1,x'),\,
\rho(x',x)\right\}=C_\rho \rho(x_1,x'),
$$
we find that
\begin{align*}
\left|D_0(y,x)-D_0(y,x')\right|
&\lesssim\left[\frac{\rho(x,x')}{1+\rho(x,y)}\right]^\beta
\left\{\frac{1}{(V_\rho)_1(y)+V_\rho(y,x)}
\left[\frac{1}{1+\rho(y,x)}\right]^\gamma\right.\nonumber\\
&\quad\left.+\frac{1}{(V_\rho)_1(y)+V_\rho(y,x')}
\left[\frac{1}{1+\rho(y,x')}\right]^\gamma\right\}\nonumber\\
&\lesssim\left[\frac{\rho(x,x')}{1+\rho(x,y)}\right]^\beta
\left(\frac{1}{(V_\rho)_1(x_1)+V_\rho(x_1,x)}
\left[\frac{1}{1+\rho(x_1,x)}\right]^\gamma\right.\nonumber\\
&\quad\times\max\left\{\left[\frac{1+\rho(x_1,x)}{1+\rho(y,x)}
\right]^{\gamma},\left[\frac{1+\rho(x_1,x)}{1+\rho(y,x)}
\right]^{n+\gamma}\right\}\nonumber\\
&\quad+\frac{1}{(V_\rho)_1(x_1)+V_\rho(x_1,x')}
\left[\frac{1}{1+\rho(x_1,x')}\right]^\gamma\nonumber\\
&\quad\left.
\times\max\left\{\left[\frac{1+\rho(x_1,x')}{1+\rho(y,x')}
\right]^{\gamma},\left[\frac{1+\rho(x_1,x')}{1+\rho(y,x')}
\right]^{n+\gamma}\right\}\right),
\end{align*}
which further implies that
\begin{align}\label{1115}
&\left|D_0(y,x)-D_0(y,x')\right|\nonumber\\
&\quad\lesssim\left[\frac{\rho(x,x')}{1+\rho(x_1,x)}\right]^\beta
\left\{\frac{1}{(V_\rho)_1(x_1)+V_\rho(x_1,x)}
\left[\frac{1}{1+\rho(x_1,x)}\right]^\gamma\right.\nonumber\\
&\qquad\left.+\frac{1}{(V_\rho)_1(x_1)+V_\rho(x_1,x')}
\left[\frac{1}{1+\rho(x_1,x')}\right]^\gamma\right\}
\left[1+\rho(x_1,y)\right]^{n+\beta+\gamma}\nonumber\\
&\quad\sim\left[\frac{\rho(x,x')}{1+\rho(x_1,x)}\right]^\beta
\frac{1}{(V_\rho)_1(x_1)+V_\rho(x_1,x)}
\left[\frac{1}{1+\rho(x_1,x)}\right]^\gamma\nonumber\\
&\qquad\times\left[1+\rho(x_1,y)\right]^{n+\beta+\gamma}.
\end{align}
From this and the size estimate \eqref{1940},
we infer that
\begin{align}\label{2255}
\left\|D_0(y,\cdot)\right\|_{\mathcal{G}(\beta,\gamma)}
\lesssim\left[1+\rho(x_1,y)\right]^{n+\beta+\gamma}.
\end{align}
By this, \eqref{23112}, and both the size
condition and the support condition of $\phi$,
we obtain
\begin{align*}
\left|\left\langle\mathcal{D}_kf(\cdot),\phi(\cdot)\right\rangle\right|
&\leq\int_{X}\left|\left\langle f,D_0(y,\cdot)\right\rangle\right|
\left|\phi(y)\right|\,d\mu(y)\\
&\lesssim\int_{B_\rho(x_1,M)}\|f\|_{(\mathring{\mathcal{G}}^\varepsilon_0(\beta,\gamma))'}
\left[1+\rho(x_1,y)\right]^{n+\beta+\gamma}\\
&\quad\times\frac{1}{(V_\rho)_1(x_1)+V_\rho(x_1,y)}
\left[\frac{1}{1+\rho(x_1,y)}\right]^{\gamma}\,d\mu(y),
\end{align*}
which, together with Lemma~\ref{A3.2}(ii),
further implies that
\begin{align*}
\left|\left\langle\mathcal{D}_kf(\cdot),\phi(\cdot)\right\rangle\right|
&\lesssim\|f\|_{(\mathring{\mathcal{G}}^\varepsilon_0(\beta,\gamma))'}\int_X
\left[1+\rho(x_1,y)\right]^{n+\beta+\gamma}\\
&\quad\times\frac{1}{(V_\rho)_1(x_1)+V_\rho(x_1,y)}
\left[\frac{1}{1+\rho(x_1,y)}\right]^{n+\beta+\gamma+1}\,d\mu(y)\\
&\lesssim\|f\|_{(\mathring{\mathcal{G}}^\varepsilon_0(\beta,\gamma))'},
\end{align*}
and hence the integral in \eqref{art-345} is finite.

From the above observations, it follows that the linear
functional induced by the function $\mathcal{D}_0 f$ defined
in (i) of Remark~\ref{Dkf2} is well defined and continuous
on the space $(\mathring{\mathcal{G}}_{\mathrm{b}}(\varepsilon,\varepsilon),
\|\cdot\|_{\mathring{\mathcal{G}}^\varepsilon_0(\beta,\gamma)})$.
Using this, Lemma~\ref{dens}, and a standard density argument,
we can extend this linear functional
in a natural fashion to distributions
in $(\mathring{{\mathcal G}}_0^\varepsilon(\beta,\gamma))'$.
That is, for any $\phi\in\mathring{{\mathcal G}}_0^\varepsilon(\beta,\gamma)$,
\begin{align}\label{kgp-375}
\left\langle\mathcal{D}_0f(\cdot),\phi(\cdot)\right\rangle
:=\lim_{j\to\infty}\int_X\left\langle f,
D_0(y,\cdot)\right\rangle\phi_j(y)\,d\mu(y),
\end{align}
where $\{\phi_j\}_{j\in\mathbb{N}}$ is any  a sequence in
$\mathring{\mathcal{G}}_{\mathrm{b}}(\varepsilon,\varepsilon)$ such that
$$
\lim_{j\to\infty}\|\phi-\phi_j\|_{\mathcal{G}(\beta,\gamma)}=0.
$$
This finishes the proof of Lemma~\ref{welldef}.
\end{proof}

\begin{lemma}\label{17001}
Let the notation be as in Remark~\ref{Dkf2}.
Then, for any $j\in\mathbb{Z}$,
$f\in(\mathring{{\mathcal G}}_0^\varepsilon(\beta,\gamma))'$,
and $\phi\in\mathring{\mathcal{G}}_{\mathrm{b}}(\varepsilon,\varepsilon)$,
$$
\left\langle f,\mathcal{D}_j^*\phi\right\rangle
=\int_X\left\langle f,D_j(x,\cdot)\right\rangle
\phi(x)\,d\mu(x).
$$
\end{lemma}

\begin{proof}
We only show the case $j=0$ because the proof
of the case $j\neq0$ is similar
and we omit the details. Fix a point $x_1\in X$.
Let $f\in(\mathring{{\mathcal G}}_0^\varepsilon(\beta,\gamma))'$,
$\phi\in\mathring{\mathcal{G}}_{\mathrm{b}}(\varepsilon,\varepsilon)
\subset\mathring{\mathcal{G}}^\varepsilon_0(\beta,\gamma)$,
and $M\in\mathbb{N}$ be such that $\mathrm{supp\,}(\phi)\subset B_\rho(x_1,M)$,
where $\rho\in\mathbf{q}$ is the regularized quasi-ultrametric
as in Definition~\ref{D2.3}. Without loss of generality, we may assume
that $\|\varphi\|_{\mathcal{G}(\varepsilon,\varepsilon)}=1$.
Suppose that $\{Q_\tau^k\}_{k\in\mathbb{Z},\tau\in I_k}$
(with their centers $\{z_\tau^k\}_{k\in\mathbb{Z},\tau\in I_k}$)
is the system of dyadic cubes of $X$
as in Lemma~\ref{DyadicSys}.
For any $k\in\mathbb{N}$, let
$$
\mathcal{A}_{k,M}:=\left\{\tau\in I_k:
Q_\tau^k\cap B_\rho(x_1,M)\neq\varnothing\right\}.
$$
Using Lemma~\ref{DyadicSys}(v) and \eqref{upperdoub},
it is easy to prove that $\#\mathcal{A}_{k,M}<\infty$
for any given $k\in\mathbb{Z}$.
For any $k\in\mathbb{Z}$ and $x\in X$, let
\begin{align}\label{15589}
\phi_k(x):=\sum_{\tau\in I_k}D_0(z_\tau^k,x)\int_{Q_\tau^k}
\phi(y)\,d\mu(y)=\sum_{\tau\in \mathcal{A}_{k,M}}D_0(z_\tau^k,x)\int_{Q_\tau^k}
\phi(y)\,d\mu(y).
\end{align}

Next, we claim that
\begin{align}\label{2121}
\lim_{k\to\infty}\left\|\mathcal{D}_0^*\phi
-\phi_k\right\|_{\mathcal{G}(\beta,\gamma)}=0.
\end{align}
To this end, let $k\in\mathbb{N}$.
We first show the size estimate of
$\mathcal{D}_0^*\phi-\phi_k$.
We claim that,
for any $k\in\mathbb{N}$, $\tau\in I_k$, $x\in X$, and $y\in Q_\tau^k$,
\begin{align}\label{zrx-358}
\left|D_0(y,x)-D_0(z_\tau^k,x)\right|
\lesssim
\frac{2^{-k\varepsilon}}{(V_\rho)_1(y)+(V_\rho)(x,y)}
\left[\frac{1}{1+\rho(x,y)}\right]^\gamma.
\end{align}
Let $k\in\mathbb{N}$, $\tau\in I_k$, $x\in X$,
and $y\in Q_\tau^k$. By Lemma~\ref{DyadicSys}(iv),
we obtain $\rho(z_\tau^k,y)\lesssim2^{-k}\leq1$.
As such, if $\rho(y,x)\ge C_0$,
where we can choose $C_0\in(0,\infty)$ sufficiently large,
but independent of $x$, $y$, $z_\tau^k$, and $k\in\mathbb{N}$, such that
$\rho(x,y), \rho(z_\tau^k,x)\ge C'$
with $C'$ as in Lemma~\ref{Pk}(i)
[using Definition~\ref{D2.1}(iii)], then
it follows from Lemma~\ref{Pk}(i) that
$$
\left|D_0(y,x)-D_0(z_\tau^k,x)\right|=0.
$$
Thus, \eqref{zrx-358} trivially holds in this case.
If $\rho(y,x)<C_0$, then
$$
1\lesssim\left[\frac{1}{1+\rho(x,y)}\right]^\gamma.
$$
Moreover, by Definition~\ref{D2.1}(iii) and $\rho(z_\tau^k,y)\lesssim1$,
we find that $\rho(z_\tau^k,x)\lesssim1$, which,
combined with $\rho(y,x),\rho(z_\tau^k,y)\lesssim1$,
both (vi) and (vii) of Lemma~\ref{A3.2}, and \eqref{doublingcond}, implies
$$
(V_\rho)_1(y)+(V_\rho)(y,x)\sim(V_\rho)_{1+\rho(y,x)}(y)
\lesssim(V_\rho)_{1}(y)\sim(V_\rho)_{1}(z_\tau^k).
$$
Given these observations, from Lemma~\ref{Pk}(ii), we deduce that
\begin{align*}
\left|D_0(y,x)-D_0(z_\tau^k,x)\right|
&\lesssim[\rho(y,z_\tau^k)]^\varepsilon
\left[\frac{1}{(V_\rho)_1(y)}+\frac{1}{(V_\rho)_1(z_\tau^k)}\right]\\
&\lesssim
\frac{2^{-k\varepsilon}}{(V_\rho)_1(y)+(V_\rho)(x,y)}
\left[\frac{1}{1+\rho(x,y)}\right]^\gamma,
\end{align*}
which completes the proof of \eqref{zrx-358} and hence the aforementioned claim.

Observe that, according to Definition~\ref{D2.1}(iii),
we can write, for any $x,y\in X$,
\begin{align}\label{zmi-34}
X=\left\{y\in X:C_\rho\rho(x_1,y)\geq\rho(x_1,x)\right\}\cup
\left\{y\in X:C_\rho\rho(y,x)\geq\rho(x_1,x)\right\}.
\end{align}
By this, \eqref{zrx-358},
and Lemma~\ref{DyadicSys}(i),
we conclude that, for any $k\in\mathbb{N}$
and $x\in X$,
\begin{align*}
\left|\mathcal{D}_0^*\phi(x)-\phi_k(x)\right|
&\leq\sum_{\tau\in I_k}\int_{Q_\tau^k}\left|D_0(y,x)-D_0(z_\tau^k,x)\right|
|\phi(y)|\,d\mu(y)\nonumber\\
&\lesssim2^{-k\varepsilon}\int_{X}
\frac{1}{(V_\rho)_1(y)+(V_\rho)(x,y)}
\left[\frac{1}{1+\rho(x,y)}\right]^\gamma
|\phi(y)|\,d\mu(y).
\end{align*}
From this and Definition~\ref{testfunction}(i),
we infer that, for any $k\in\mathbb{N}$ and $x\in X$,
\begin{align*}
\left|\mathcal{D}_0^*\phi(x)-\phi_k(x)\right|
&\lesssim2^{-k\varepsilon}\int_{X}
\frac{1}{(V_\rho)_1(y)+(V_\rho)(x,y)}
\left[\frac{1}{1+\rho(x,y)}\right]^\gamma\nonumber\\
&\quad\times\frac{1}{(V_\rho)_1(x_1)+(V_\rho)(x_1,y)}
\left[\frac{1}{1+\rho(x_1,y)}\right]^\varepsilon\,d\mu(y),
\end{align*}
which, together with Lemma~\ref{A3.2}(ii),
further implies that
\begin{align}\label{1547}
\left|\mathcal{D}_0^*\phi(x)-\phi_k(x)\right|
&\lesssim2^{-k\varepsilon}\Bigg\{\int_{\rho(x_1,x)\leq C_\rho\rho(x_1,y)}
\frac{1}{(V_\rho)_1(y)+(V_\rho)(x,y)}
\left[\frac{1}{1+\rho(x,y)}\right]^\gamma\nonumber\\
&\quad\times\frac{1}{(V_\rho)_1(x_1)+(V_\rho)(x_1,y)}
\left[\frac{1}{1+\rho(x_1,y)}\right]^\varepsilon\,d\mu(y)\nonumber\\
&\quad+\int_{\rho(x_1,x)\leq C_\rho\rho(y,x)}
\cdots\Bigg\}\nonumber\\
&\lesssim2^{-k\varepsilon}
\frac{1}{(V_\rho)_1(x_1)+(V_\rho)(x_1,x)}
\left[\frac{1}{1+\rho(x_1,x)}\right]^\gamma.
\end{align}
This is the desired size estimate of $\mathcal{D}_0^*\phi-\phi_k$.

Now, we prove the regularity estimate of
$\mathcal{D}_0^*\phi-\phi_k$.
Assume that $x,x'\in X$ satisfy
$\rho(x,x')\leq\frac{1+\rho(x_1,x)}{2C_\rho}$.
We consider the following two cases for $\rho(x,x')$.

\emph{Case (1)}
$\frac{1+\rho(x_1,x)}{(2C_\rho)^2}\leq
\rho(x,x')\leq\frac{1+\rho(x_1,x)}{2C_\rho}$. In this case,
by the above proven size estimate for
$\mathcal{D}_0^*\phi-\phi_k$, we find that
\begin{align*}
&\left|\left[\mathcal{D}_0^*\phi(x)-\phi_k(x)\right]
-\left[\mathcal{D}_0^*\phi(x')-\phi_k(x')\right]\right|\nonumber\\
&\quad\lesssim2^{-k\varepsilon}
\frac{1}{(V_\rho)_1(x_1)+(V_\rho)(x_1,x)}
\left[\frac{1}{1+\rho(x_1,x)}\right]^\gamma\nonumber\\
&\quad\quad+2^{-k\varepsilon}
\frac{1}{(V_\rho)_1(x_1)+(V_\rho)(x_1,x')}
\left[\frac{1}{1+\rho(x_1,x')}\right]^\gamma,
\end{align*}
which, combined with the fact that
$1+\rho(x_1,x)\sim1+\rho(x_1,x')$
and Lemma~\ref{A3.2}(vi),
further implies that
\begin{align}\label{1546}
&\left|\left[\mathcal{D}_0^*\phi(x)-\phi_k(x)\right]
-\left[\mathcal{D}_0^*\phi(x')-\phi_k(x')\right]\right|\nonumber\\
&\quad\lesssim2^{-k\varepsilon}
\frac{1}{(V_\rho)_1(x_1)+(V_\rho)(x_1,x)}
\left[\frac{1}{1+\rho(x_1,x)}\right]^\gamma\nonumber\\
&\quad\sim2^{-k\varepsilon}
\left[\frac{\rho(x,x')}{1+\rho(x_1,x)}\right]^\beta
\frac{1}{(V_\rho)_1(x_1)+(V_\rho)(x_1,x)}
\left[\frac{1}{1+\rho(x_1,x)}\right]^\gamma.
\end{align}
This is the desired regularity estimate of $\mathcal{D}_0^*\phi-\phi_k$.

\emph{Case (2)}
$\rho(x,x')<\frac{1+\rho(x_1,x)}{(2C_\rho)^2}$. In this case,
we first show that, for any $\tau\in\mathcal{A}_{k,M}$
and $y\in Q_\tau^k$,
\begin{align}\label{1027}
\mathrm{G}:&=\left|\left[D_0(y,x)-D_0(z_\tau^k,x)\right]
-\left[D_0(y,x')-D_0(z_\tau^k,x')\right]\right|\nonumber\\
&\lesssim2^{-k\varepsilon}
\left[\frac{\rho(x,x')}{1+\rho(y,x)}\right]^\varepsilon
\left\{\frac{1}{(V_\rho)_1(y)+(V_\rho)(y,x)}
\left[\frac{1}{1+\rho(y,x)}\right]^\gamma\right.\nonumber\\
&\quad\left.+\frac{1}{(V_\rho)_1(y)+(V_\rho)(y,x')}
\left[\frac{1}{1+\rho(y,x')}\right]^\gamma\right\}.
\end{align}
Indeed, if $\rho(x,x')\leq\frac{1+\rho(x,y)}{2C_\rho}$
and $\rho(x,y)\ge C_0$,
where we can choose $C_0\in(0,\infty)$ sufficiently large,
but independent of $x$, $x'$, $y$, $z_\tau^k$, and $k\in\mathbb{N}$, such that
$$
\rho(z_\tau^k,x),\rho(y,x'),\rho(z_\tau^k,x')\ge C'
$$
with $C'$ as in Lemma~\ref{Pk}(i), then
it follows from Lemma~\ref{Pk}(i) that
$$
\mathrm{G}=0.
$$
If $\rho(x,x')\leq\frac{1+\rho(x,y)}{2C_\rho}$
and $\rho(x,y)<C_0$,
then we can use Lemma~\ref{Pk}(iii), the fact that
$$
\rho(x,y),\rho(x,x'),\rho(y,x')\lesssim1,
$$
Lemma~\ref{DyadicSys}(iv), and an argument similar to that used in the
proof of \eqref{zrx-358} to conclude that,
for any $k\in\mathbb{N}$,
\begin{align}\label{1026}
\mathrm{G}&\lesssim\left[\rho(y,z_\tau^k)
\rho(x,x')\right]^\varepsilon
\left[\frac{1}{(V_\rho)_1(x)}+\frac{1}{(V_\rho)_1(x')}\right]\nonumber\\
&\sim2^{-k\varepsilon}
\left[\frac{\rho(x,x')}{1+\rho(y,x)}\right]^\varepsilon
\left\{\frac{1}{(V_\rho)_1(y)+(V_\rho)(y,x)}
\left[\frac{1}{1+\rho(y,x)}\right]^\gamma\right.\nonumber\\
&\quad\left.+\frac{1}{(V_\rho)_1(y)+(V_\rho)(y,x')}
\left[\frac{1}{1+\rho(y,x')}\right]^\gamma\right\}.
\end{align}
If $\rho(x,x')>\frac{1+\rho(x,y)}{2C_\rho}$,
by both (i) and (ii) of Lemma~\ref{Pk},
Lemma~\ref{DyadicSys}(iv), and an argument similar to the one used in the
proof of \eqref{zrx-358}, we find that,
for any $k\in\mathbb{N}$,
\begin{align*}
\mathrm{G}&\leq\left|D_0(y,x)-D_0(z_\tau^k,x)\right|
+\left|D_0(y,x')-D_0(z_\tau^k,x')\right|\\
&\lesssim\left[\rho(y,z_\tau^k)\right]^\varepsilon
\left\{\frac{1}{(V_\rho)_1(y)+(V_\rho)(y,x)}
\left[\frac{1}{1+\rho(y,x)}\right]^\gamma\right.\\
&\quad\left.+\frac{1}{(V_\rho)_1(y)+(V_\rho)(y,x')}
\left[\frac{1}{1+\rho(y,x')}\right]^\gamma\right\}\\
&\lesssim2^{-k\varepsilon}
\left[\frac{\rho(x,x')}{1+\rho(y,x)}\right]^\varepsilon
\left\{\frac{1}{(V_\rho)_1(y)+(V_\rho)(y,x)}
\left[\frac{1}{1+\rho(y,x)}\right]^\gamma\right.\\
&\quad\left.+\frac{1}{(V_\rho)_1(y)+(V_\rho)(y,x')}
\left[\frac{1}{1+\rho(y,x')}\right]^\gamma\right\}.
\end{align*}
This and \eqref{1026} finish the proof of \eqref{1027}.

Next, we use \eqref{1027} to prove the regularity estimate of
$\mathcal{D}_0^*\phi-\phi_k$
in the case $\rho(x,x')<\frac{1+\rho(x_1,x)}{(2C_\rho)^2}$.
From \eqref{1027}, \eqref{zmi-34}, and the size condition of $\phi$
(here, we also used the fact that $\phi$ has a bounded support),
we infer that
\begin{align*}
&\left|\left[\mathcal{D}_0^*\phi(x)-\phi_k(x)\right]
-\left[\mathcal{D}_0^*\phi(x')-\phi_k(x')\right]\right|\nonumber\\
&\quad\leq\sum_{\tau\in I_k}\int_{Q_\tau^k}
\mathrm{G}|\phi(y)|\,d\mu(y)\nonumber\\
&\quad\lesssim2^{-k\varepsilon}\sum_{\tau\in I_k}\int_{Q_\tau^k}
\left[\frac{\rho(x,x')}{1+\rho(y,x)}\right]^\varepsilon
\left\{\frac{1}{(V_\rho)_1(y)+(V_\rho)(y,x)}
\left[\frac{1}{1+\rho(y,x)}\right]^\gamma\right.\nonumber\\
&\qquad\left.+\frac{1}{(V_\rho)_1(y)+(V_\rho)(y,x')}
\left[\frac{1}{1+\rho(y,x')}\right]^\gamma\right\}
|\phi(y)|\,d\mu(y),
\end{align*}
which further implies that
\begin{align}\label{15431}
&\left|\left[\mathcal{D}_0^*\phi(x)-\phi_k(x)\right]
-\left[\mathcal{D}_0^*\phi(x')-\phi_k(x')\right]\right|\nonumber\\
&\quad\lesssim2^{-k\varepsilon}\int_{X}
\left[\frac{\rho(x,x')}{1+\rho(y,x)}\right]^\varepsilon
\left\{\frac{1}{(V_\rho)_1(y)+(V_\rho)(y,x)}
\left[\frac{1}{1+\rho(y,x)}\right]^\gamma\right.\nonumber\\
&\qquad\left.+\frac{1}{(V_\rho)_1(y)+(V_\rho)(y,x')}
\left[\frac{1}{1+\rho(y,x')}\right]^\gamma\right\}\nonumber\\
&\qquad\times\frac{1}{(V_\rho)_1(x_1)+(V_\rho)(x_1,y)}
\left[\frac{1}{1+\rho(x_1,y)}\right]^{\varepsilon+\gamma}\,d\mu(y)\nonumber\\
&\quad\leq2^{-k\varepsilon}\left[\int_{C_\rho\rho(x,y)\ge\rho(x_1,x)}\cdots+
\int_{C_\rho\rho(x_1,y)\ge\rho(x_1,x)}\cdots\right]
=:2^{-k\varepsilon}(\mathrm{I}+\mathrm{II}).
\end{align}

To deal with $\mathrm{I}$,
by $\rho(x,y)\ge\frac{\rho(x_1,x)}{C_\rho}$ and
$\beta<\varepsilon$,
we find that
\begin{align*}
\mathrm{I}
&\lesssim\left[\frac{\rho(x,x')}{1+\rho(x_1,x)}\right]^\beta
\int_{C_\rho\rho(x,y)\ge\rho(x_1,x)}
\left\{\frac{1}{(V_\rho)_1(y)+(V_\rho)(y,x)}
\left[\frac{1}{1+\rho(y,x)}\right]^\gamma\right.\nonumber\\
&\quad\left.+\frac{1}{(V_\rho)_1(y)+(V_\rho)(y,x')}
\left[\frac{1}{1+\rho(y,x')}\right]^\gamma\right\}\nonumber\\
&\quad\times\frac{1}{(V_\rho)_1(x_1)+(V_\rho)(x_1,y)}
\left[\frac{1}{1+\rho(x_1,y)}\right]^{\varepsilon+\gamma}\,d\mu(y),
\end{align*}
which, together with
the inequality that $1+\rho(x_1,x)\lesssim1+\rho(y,x')$
[which can be deduced from
\begin{align*}
1+\rho(x_1,x)&\lesssim1+\rho(x_1,x)-C_\rho^2\rho(x,x')\\
&\leq1+C_\rho\rho(x,y)-C_\rho^2\rho(x,x')\leq1+C_\rho^2\rho(y,x')]
\end{align*}
and both (ii) and (vi) of Lemma~\ref{A3.2},
further implies that
\begin{align}\label{15432}
\mathrm{I}
&\lesssim\left[\frac{\rho(x,x')}{1+\rho(x_1,x)}\right]^\beta
\int_{C_\rho\rho(x,y)\ge\rho(x_1,x)}
\frac{1}{(V_\rho)_1(x_1)+(V_\rho)(x_1,x)}
\left[\frac{1}{1+\rho(x_1,x)}\right]^\gamma\nonumber\\
&\quad\times\frac{1}{(V_\rho)_1(x_1)+(V_\rho)(x_1,y)}
\left[\frac{1}{1+\rho(x_1,y)}\right]^{\varepsilon+\gamma}\,d\mu(y)\nonumber\\
&\lesssim\left[\frac{\rho(x,x')}{1+\rho(x_1,x)}\right]^\beta
\frac{1}{(V_\rho)_1(x_1)+(V_\rho)(x_1,x)}
\left[\frac{1}{1+\rho(x_1,x)}\right]^\gamma.
\end{align}

To deal with $\mathrm{II}$,
from $\beta<\varepsilon$ and both (ii) and (vi) of Lemma~\ref{A3.2},
we deduce that
\begin{align*}
\mathrm{II}
&\lesssim\frac{1}{(V_\rho)_1(x_1)+(V_\rho)(x_1,x)}
\left[\frac{1}{1+\rho(x_1,x)}\right]^{\gamma}\\
&\quad\times\int_{C_\rho\rho(x_1,y)\ge\rho(x_1,x)}
\left[\frac{\rho(x,x')}{1+\rho(x_1,x)}\right]^\varepsilon
\left\{\frac{1}{(V_\rho)_1(x)+(V_\rho)(y,x)}
\left[\frac{1}{1+\rho(y,x)}\right]^\gamma\right.\\
&\quad\left.+\frac{1}{(V_\rho)_1(x')+(V_\rho)(y,x')}
\left[\frac{1}{1+\rho(y,x')}\right]^\gamma\right\}\,d\mu(y)\\
&\lesssim\left[\frac{\rho(x,x')}{1+\rho(x_1,x)}\right]^\beta
\frac{1}{(V_\rho)_1(x_1)+(V_\rho)(x_1,x)}
\left[\frac{1}{1+\rho(x_1,x)}\right]^\gamma.
\end{align*}
This, together with both \eqref{15431} and \eqref{15432},
further implies that
\begin{align*}
&\left|\left[\mathcal{D}_0^*\phi(x)-\phi_k(x)\right]
-\left[\mathcal{D}_0^*\phi(x')-\phi_k(x')\right]\right|\nonumber\\
&\quad\lesssim2^{-k\varepsilon}
\left[\frac{\rho(x,x')}{1+\rho(x_1,x)}\right]^\beta
\frac{1}{(V_\rho)_1(x_1)+(V_\rho)(x_1,x)}
\left[\frac{1}{1+\rho(x_1,x)}\right]^\gamma.
\end{align*}
By this [the regularity estimate in Case (2)], \eqref{1546}
[the regularity estimate in Case (1)], and \eqref{1547}
(the size estimate),
we then complete the proof of \eqref{2121}.

From Lemma~\ref{GOVV} and the fact that
$\phi\in\mathring{\mathcal{G}}_{\mathrm{b}}(\varepsilon,\varepsilon)
\subset\mathring{\mathcal{G}}^\varepsilon_0(\beta,\gamma)$, we easily infer that
$\mathcal{D}_0^*\phi\in\mathring{\mathcal{G}}^\varepsilon_0(\beta,\gamma)$.
By Lemma~\ref{Pk}, the definition of $\phi_k$,
and the fact that $\phi$ has a bounded support
(and hence the summation in \eqref{15589} only has finitely many terms),
we conclude that $\phi_k\in\mathring{\mathcal{G}}(\varepsilon,\varepsilon)
\subset\mathring{\mathcal{G}}^\varepsilon_0(\beta,\gamma)$
for any $k\in\mathbb{Z}$.
From this, \eqref{2121}, and $\#\mathcal{A}_{k,M}<\infty$,
we deduce that
\begin{align*}
\left\langle f,\mathcal{D}_0^*\phi\right\rangle
&=\lim_{k\to\infty}\left\langle f,\sum_{\tau\in
\mathcal{A}_{k,M}}D_0(z_\tau^k,\cdot)\int_{Q_\tau^k}
\phi(y)\,d\mu(y)\right\rangle\\
&=\lim_{k\to\infty}\sum_{\tau\in \mathcal{A}_{k,M}}\left\langle f,
D_0(z_\tau^k,\cdot)\right\rangle\int_{Q_\tau^k}
\phi(y)\,d\mu(y).
\end{align*}
Thus, to complete the proof of the present lemma,
it remains to show that
\begin{align}\label{1611}
\lim_{k\to\infty}\sum_{\tau\in \mathcal{A}_{k,M}}\left\langle f,
D_0(z_\tau^k,\cdot)\right\rangle\int_{Q_\tau^k}
\phi(y)\,d\mu(y)
=\int_{X}\left\langle f,D_0(y,\cdot)\right\rangle
\phi(y)\,d\mu(y).
\end{align}
To this end, we first prove that,
for any $x,x'\in X$ with
$\rho(x,x')\leq(2C_\rho)^{-2}$,
\begin{align}\label{2235}
\left|\left\langle f,
D_0(x,\cdot)\right\rangle-\left\langle f,
D_0(x',\cdot)\right\rangle\right|
&=\left|\left\langle f,
D_0(x,\cdot)-D_0(x',\cdot)\right\rangle\right|\nonumber\\
&\leq\|f\|_{(\mathring{\mathcal{G}}^\varepsilon_0(\beta,\gamma))'}
\left\|D_0(x,\cdot)-D_0(x',\cdot)\right\|_{\mathcal{G}(\beta,\gamma)}
\nonumber\\
&\lesssim\|f\|_{(\mathring{\mathcal{G}}^\varepsilon_0(\beta,\gamma))'}
\left[\rho(x,x')\right]^\varepsilon
\left[1+\rho(x_1,x)\right]^{n+\beta+\gamma}.
\end{align}
Note that the first inequality follows from Lemma~\ref{Pk}.
To show the last inequality of \eqref{2235},
we only need to estimate
$\|D_0(x,\cdot)-D_0(x',\cdot)\|_{\mathcal{G}(\beta,\gamma)}$.
We first consider the size condition. Assume that
$y\in X$ satisfies $\rho(x,y)\leq C_\rho(C'+1)$;
otherwise $D_k(x,y)=0=D_k(x',y)$ since Lemma~\ref{Pk}(i).
In this case, we have
$\rho(x,y),\rho(x',y)\lesssim1$.
By this, Lemma~\ref{Pk}(ii),
an argument similar to that used in the proof
of \eqref{1940}, and Lemma~\ref{A3.2}(iii),
we find that
\begin{align}\label{2249}
\left|D_0(x,y)-D_0(x',y)\right|
&\lesssim\left[\rho(x,x')\right]^\varepsilon
\left[\frac{1}{(V_\rho)_1(x)}+
\frac{1}{(V_\rho)_1(x')}\right]\nonumber\\
&\sim\left[\rho(x,x')\right]^\varepsilon
\frac{1}{(V_\rho)_1(y)+V_\rho(y,x)}
\left[\frac{1}{1+\rho(y,x)}\right]^\gamma\nonumber\\
&\lesssim\left[\rho(x,x')\right]^\varepsilon
\left[1+\rho(x_1,y)\right]^{n+\beta+\gamma}\nonumber\\
&\quad\times\frac{1}{(V_\rho)_1(x_1)+V_\rho(x_1,y)}
\left[\frac{1}{1+\rho(x_1,y)}\right]^\gamma\nonumber\\
&\sim\left[\rho(x,x')\right]^\varepsilon
\left[1+\rho(x_1,x)\right]^{n+\beta+\gamma}\nonumber\\
&\quad\times\frac{1}{(V_\rho)_1(x_1)+V_\rho(x_1,y)}
\left[\frac{1}{1+\rho(x_1,y)}\right]^\gamma.
\end{align}
For the regularity condition of $D_0(x,\cdot)-D_0(x',\cdot)$,
let $y,y'\in X$ be such that $\rho(y,y')\leq\frac{1+\rho(x_1,y)}{2C_\rho}$.
From arguments similar to those used in the estimations
of both \eqref{1027} and \eqref{1115}, we immediately
infer that
\begin{align*}
&\left|\left[D_0(x,y)-D_0(x',y)\right]-
\left[D_0(x,y')-D_0(x',y')\right]\right|\\
&\quad\lesssim\left[\rho(x,x')\right]^\varepsilon
\left[1+\rho(x_1,x)\right]^{n+\beta+\gamma}\nonumber\\
&\qquad\times\left[\frac{\rho(y,y')}{1+\rho(x_1,y)}\right]^\beta
\frac{1}{(V_\rho)_1(x_1)+V_\rho(x_1,y)}
\left[\frac{1}{1+\rho(x_1,y)}\right]^\gamma,
\end{align*}
which, combined with \eqref{2249},
completes the proof of \eqref{2235}.

By \eqref{2235}, Lemma~\ref{DyadicSys}(iv),
$\mathrm{supp\,}(\phi)\subset B_\rho(x_1,M)$,
and the size condition of $\phi$,
we conclude that
\begin{align*}
&\left|\sum_{\tau\in \mathcal{A}_{k,M}}\left\langle f,
D_0(z_\tau^k,\cdot)\right\rangle\int_{Q_\tau^k}
\phi(y)\,d\mu(y)
-\int_{X}\left\langle f,D_0(y,\cdot)\right\rangle
\phi(y)\,d\mu(y)\right|\\
&\quad\leq\sum_{\tau\in \mathcal{A}_{k,M}}
\int_{Q_\tau^k}\left|\left\langle f,D_0(z_\tau^k,\cdot)\right\rangle
-\left\langle f,D_0(y,\cdot)\right\rangle\right|
|\phi(y)|\,d\mu(y)\\
&\quad\lesssim\sum_{\tau\in \mathcal{A}_{k,M}}
\int_{Q_\tau^k}\|f\|_{(\mathring{\mathcal{G}}^\varepsilon_0(\beta,\gamma))'}
\left[\rho(z_\tau^k,y)\right]^\varepsilon
\left[1+\rho(x_1,y)\right]^{n+\beta+\gamma}\mathbf{1}_{B_\rho(x_1,M)}(y)\\
&\qquad\times\frac{1}{(V_\rho)_1(x_1)+V_\rho(x_1,y)}
\left[\frac{1}{1+\rho(x_1,y)}\right]^\gamma\,d\mu(y),
\end{align*}
which, together with Lemma~\ref{A3.2}(ii),
further implies that
\begin{align*}
&\left|\sum_{\tau\in\mathcal{A}_{k,M}}\left\langle f,
D_0(z_\tau^k,\cdot)\right\rangle\int_{Q_\tau^k}
\phi(y)\,d\mu(y)
-\int_{X}\left\langle f,D_0(y,\cdot)\right\rangle
\phi(y)\,d\mu(y)\right|\\
&\quad\lesssim\|f\|_{(\mathring{\mathcal{G}}^\varepsilon_0(\beta,\gamma))'}
C_r2^{-k\varepsilon}(1+M)^{n+\beta+\gamma}\\
&\qquad\times\int_X
\frac{1}{(V_\rho)_1(x_1)+V_\rho(x_1,y)}
\left[\frac{1}{1+\rho(x_1,y)}\right]^\gamma\,d\mu(y)\\
&\quad\lesssim2^{-k\varepsilon}
\|f\|_{(\mathring{\mathcal{G}}^\varepsilon_0(\beta,\gamma))'}\to0
\end{align*}
as $k\to\infty$. This finishes the proof of \eqref{1611}
and hence Lemma~\ref{17001}.
\end{proof}

Now, we are ready to prove Proposition~\ref{1651}.

\begin{proof}[Proof of Proposition~\ref{1651}]
Let $k\in\mathbb{Z}$ and $f\in(\mathring{{\mathcal G}}_0^\varepsilon(\beta,\gamma))'$.
From Lemma~\ref{welldef}, we deduce that
the two linear functionals $\mathcal{D}_k f$ defined
in (i) and (ii) of Remark~\ref{Dkf2} are well-defined distributions
in $(\mathring{{\mathcal G}}_0^\varepsilon(\beta,\gamma))'$.
Next, by Lemma~\ref{dens}, we find that, for any
$\phi\in\mathring{{\mathcal G}}_0^\varepsilon(\beta,\gamma)$,
there exists a sequence
$\{\phi_j\}_{j\in\mathbb{N}}$ in
$\mathring{\mathcal{G}}_{\mathrm{b}}(\varepsilon,\varepsilon)$ such that
$$
\lim_{j\to\infty}\|\phi-\phi_j\|_{\mathcal{G}(\beta,\gamma)}=0.
$$
From this, the fact that $\mathcal{D}_k^*$ is bounded on
$\mathring{{\mathcal G}}_0^\varepsilon(\beta,\gamma)$
(which follows from Lemma~\ref{GOVV} and a standard density argument),
Lemma~\ref{17001}, and \eqref{kgp-375}, we infer that
\begin{align*}
\left\langle f,\mathcal{D}_k^*\phi\right\rangle
&=\lim_{j\to\infty}
\left\langle f,\mathcal{D}_k^*\phi_j\right\rangle
=\lim_{j\to\infty}\int_X\left\langle f,
D_k(x,\cdot)\right\rangle\phi_j(x)\,d\mu(x)
=\left\langle\mathcal{D}_kf,\phi\right\rangle,
\end{align*}
which completes the proof of Proposition~\ref{1651}.
\end{proof}

\section{Inhomogeneous Continuous
Calder\'on Reproducing Formulae}
\label{S3.3}

The main result of this section is Theorem~\ref{A3.33},
where we establish inhomogeneous continuous
Calder\'on reproducing formulae
on the ultra-RD-space $(X,\mathbf{q},\mu)$,
where $\mu$ is assumed to be a Borel-semiregular measure on $X$.
We point out that, in contrast to the homogeneous case,
we have no restriction on $\mathrm{diam\,}(X)$,
which means $\mathrm{diam\,}(X)<\infty$
or $\mathrm{diam\,}(X)=\infty$.
As pointed out in \cite[Remark~1.7]{hmy08},
an essential difference in the theory occurs when
$\mu(X)<\infty$ [equivalently $\mathrm{diam\,}(X)<\infty$]
versus the case when $\mu(X)=\infty$ [equivalently $\mathrm{diam\,}(X)=\infty$].
This is because, when $\mu(X)<\infty$, it is not true that
$\|\mathcal{S}_{2^{-k}}f\|_{L^p(X)}\to0$ as $k\to-\infty$
for any $f\in L^p(X)$ and
for any given $p\in(1,\infty)$;
see the proof of Theorem~\ref{corATI}(vii).
Thus, when $\mu(X)<\infty$, in the Calder\'on reproducing formula
\eqref{3.12}, we should replace $\mathcal{D}_0$
by $\mathcal{S}_0$ and $\mathcal{D}_k$ for any $k\in\{-1,-2,\ldots\}$ by 0.
Then, with this modification,
we always have an inhomogeneous term related to $\mathcal{S}_0$,
which needs some special care.

We first introduce the following definition of the
inhomogeneous approximation of the identity on $X$.
The reader is reminded that the symbol ``$\preceq$" was introduced
in Remark~\ref{rmk390}(iii) and that $\mathbb{Z}_+:=\{0,1,2,\dots\}$.

\begin{definition}\label{DefIATI}
Let $(X,\mathbf{q},\mu)$ be a quasi-ultrametric space
of homogeneous type.
Assume that $0<\varepsilon\preceq{\mathrm{ind\,}}(X,\mathbf{q})$ with
${\mathrm{ind\,}}(X,\mathbf{q})$ as in \eqref{index}.
A sequence $\{\mathcal{S}_{2^{-k}}\}_{k\in\mathbb{Z}_+}$ of
linear operators is called
an \emph{inhomogeneous approximation of the identity of order $\varepsilon$}
\index[words]{inhomogeneous approximation of the identity}
if $\{\mathcal{S}_{2^{-k}}\}_{k\in\mathbb{Z}_+}$
satisfies all the conditions of Definition~\ref{DefATI}
with $t$ replaced by $2^{-k}$.
\end{definition}

The following proposition for the inhomogeneous approximation
of the identity
is a simple corollary of Theorem~\ref{corATI}(ix);
we omit the details of its proof.

\begin{proposition}\label{A3.24}
Let $(X,\mathbf{q},\mu)$ be a quasi-ultrametric space
of homogeneous type,
where $\mu$ is assumed to be a Borel-semiregular measure on $X$.
Assume that $0<\varepsilon\preceq{\mathrm{ind\,}}(X,\mathbf{q})$ with
${\mathrm{ind\,}}(X,\mathbf{q})$ as in \eqref{index} and
$\{\mathcal{S}_{2^{-k}}\}_{k\in\mathbb{Z}_+}$
is an inhomogeneous approximation
of the identity of order $\varepsilon$ as in Definition~\ref{DefIATI}.
Let
\begin{align*}
\mathcal{D}_0:=\mathcal{S}_1
\ \ \text{and}\ \
\mathcal{D}_k:=\mathcal{S}_{2^{-k}}-\mathcal{S}_{2^{-k+1}}
\end{align*}
for any $k\in\mathbb{N}$.
Then, for any $p\in[1,\infty)$,
$$
\mathcal{I}=\sum_{k=0}^{\infty}\mathcal{D}_k
$$
in $L^p(X)$, where $\mathcal{I}$ is the identity operator on $L^p(X)$.
\end{proposition}

To establish sharp inhomogeneous continuous
Calder\'on reproducing formulae,
we need a technical lemma
which is a variant of Lemma~\ref{A3.12}.

\begin{lemma}\label{A3.25}
Let $(X,\mathbf{q},\mu)$ be a quasi-ultrametric space
of homogeneous type,
where $\mu$ is assumed to be a Borel-semiregular measure on $X$.
Assume that $0<\varepsilon\preceq{\mathrm{ind\,}}(X,\mathbf{q})$ with
${\mathrm{ind\,}}(X,\mathbf{q})$ as in \eqref{index} and
both $\{\mathcal{S}_{2^{-k}}\}_{k\in\mathbb{Z}_+}$
and $\{\mathcal{E}_{2^{-k}}\}_{k\in\mathbb{Z}_+}$ are
two inhomogeneous approximations
of the identity of order $\varepsilon$
defined as in Definition~\ref{DefIATI}.
For any $k\in\mathbb{N}$, let
$\mathcal{P}_k:=\mathcal{S}_{2^{-k}}-\mathcal{S}_{2^{-k+1}}$,
$\mathcal{Q}_k:=\mathcal{E}_{2^{-k}}-\mathcal{E}_{2^{-k+1}}$,
$\mathcal{P}_0:=\mathcal{S}_1$, and $\mathcal{Q}_0:=\mathcal{E}_1$.
Then, for any given
$\varepsilon'\in(0,\varepsilon)$, $\delta\in(0,\varepsilon-\varepsilon']$,
and $\rho\in\mathbf{q}$ having the property
that all $\rho$-balls are $\mu$-measurable,
there exist two constants $C,C'\in(0,\infty)$,
depending on $\varepsilon'$, $\varepsilon$,
and $\delta$, such that, for any $k,l\in\mathbb{Z}_+$,
\eqref{3.2}, \eqref{3.3}, \eqref{3.4},
and \eqref{3.5} still hold.
\end{lemma}

\begin{proof}
We denote, respectively, by $\{S_{2^{-k}}\}_{k\in\mathbb{Z}_+}$
$\{E_{2^{-k}}\}_{k\in\mathbb{Z}_+}$,
$\{P_{2^{-k}}\}_{k\in\mathbb{Z}_+}$,
and $\{Q_{2^{-k}}\}_{k\in\mathbb{Z}_+}$
the kernels of
$\{\mathcal{S}_{2^{-k}}\}_{k\in\mathbb{Z}_+}$,
$\{\mathcal{E}_{2^{-k}}\}_{k\in\mathbb{Z}_+}$,
$\{\mathcal{P}_{2^{-k}}\}_{k\in\mathbb{Z}_+}$,
and $\{\mathcal{Q}_{2^{-k}}\}_{k\in\mathbb{Z}_+}$.
Without loss of generality, we may assume that $\rho$
is a regularized quasi-ultrametric as in Definition~\ref{D2.3}.
The proof of the present lemma is essentially the same
as that of Lemma~\ref{A3.12}.
The only situations that are different are
the cases where $l=0$ or $k=0$.
By symmetry, we only show
\eqref{3.2}, \eqref{3.3}, \eqref{3.4},
and \eqref{3.5} in the case both $k\in\mathbb{Z}_+$
and $l=0$. For the case both $k\in\mathbb{N}$
and $l=0$, the proofs of \eqref{3.2}, \eqref{3.3}, \eqref{3.4},
and \eqref{3.5}
are the same as that of
Lemma~\ref{A3.12} because we only need
the vanish condition of $P_k$ or $Q_k$.
Thus, to complete the proof of the
present lemma, it remains to prove \eqref{3.2}, \eqref{3.3}, \eqref{3.4},
and \eqref{3.5} in the case $k=l=0$.

We first show \eqref{3.2}.
On the one hand, by the support conditions of both
$P_0=S_1$ and $Q_0=E_1$
[see Definition~\ref{DefATI}(i)],
we find that there exists a constant $C\in(1,\infty)$
such that, for any $x,y,z\in X$
satisfying $\rho(x,z)>C$ and $\rho(y,z)>C$,
\begin{align*}
P_0(x,z)=0
\ \ \text{and}\ \
Q_0(z,y)=0.
\end{align*}
This, combined with Definition~\ref{D2.1}(iii),
further implies that,
for any $x,y\in X$ satisfying $\rho(x,y)>CC_\rho$,
\begin{align*}
B_\rho(x,C)\cap B_\rho(y,C)=\varnothing
\end{align*}
and hence
\begin{align}\label{2206}
P_0Q_0(x,y)
&=\int_XP_0(x,z)Q_0(z,y)\,d\mu(z)\nonumber\\
&=\int_{B_\rho(x,C)\cap B_\rho(y,C)}
P_0(x,z)Q_0(z,y)\,d\mu(z)=0.
\end{align}
On the other hand, from Lemmas~\ref{Pk} and~\ref{A3.2}(vii)
and \eqref{upperdoub}, we deduce that,
for any $x,y\in X$ satisfying $\rho(x,y)\leq CC_\rho$,
\begin{align*}
\left|P_0Q_0(x,y)\right|
&\leq\int_X|P_0(x,z)Q_0(z,y)|\,d\mu(z)
\lesssim\int_{B_\rho(x,C)}\frac{1}{(V_\rho)_1(x)}
\frac{1}{(V_\rho)_1(y)}\,d\mu(z)
\\
&\sim\frac{(V_\rho)_C(x)}{(V_\rho)_1(x)}
\frac{1}{(V_\rho)_1(x)}
\lesssim\frac{1}{(V_\rho)_1(x)},
\end{align*}
which, together with \eqref{2206},
completes the proof of \eqref{3.2}.

Now, by an argument similar to that used
in the proof of \eqref{3.3} in Lemma~\ref{A3.12},
to prove \eqref{3.3} in the present lemma,
it suffices to show \eqref{Tl} with both $k=0$
and $l=0$.
From Lemma~\ref{Pk} [see also Remark~\ref{RemThmATI}(i)],
we infer that, for any $x,y,y'\in X$,
\begin{align*}
&\left|P_0Q_0(x,y)-P_0Q_0(x,y')\right|\\
&\quad\leq\int_X|P_0(x,z)|\left|Q_0(z,y)-Q_0(z,y')\right|\,d\mu(z)
\\
&\quad\lesssim\int_{B_\rho(x,C)}
\frac{1}{(V_\rho)_1(x)}
\left[\rho(y,y')\right]^\varepsilon\frac{1}{(V_\rho)_1(z)}\,d\mu(z)
\\
&\quad\sim\int_{B_\rho(x,C)}
\frac{1}{(V_\rho)_1(x)}
\left[\rho(y,y')\right]^\varepsilon\frac{1}{(V_\rho)_1(x)}\,d\mu(z)
\quad\text{by Lemma~\ref{A3.2}(vii)}
\\
&\quad\lesssim\left[\rho(y,y')\right]^\varepsilon\frac{1}{(V_\rho)_1(x)}
\quad\text{by \eqref{upperdoub}}.
\end{align*}
This finishes the proof of \eqref{Tl} and hence \eqref{3.3}.
The proof of \eqref{3.4} is similar to the proof of \eqref{3.3}.

Finally, we prove \eqref{3.5}.
By an argument similar to that used
in the proof of \eqref{3.5} in Lemma~\ref{A3.12},
to show \eqref{3.5} in the present lemma,
it suffices to prove \eqref{418x} with both $k=0$
and $l=0$.
From Lemma~\ref{Pk}(ii),
it follows that, for any $x,y,x',y'\in X$,
\begin{align*}
&\left|\left[P_0Q_0(x,y)-P_0Q_0(x',y)\right]-\left[P_0Q_0(x,y')
-P_0Q_0(x',y')\right]\right|
\\
&\quad\leq\int_{B_\rho(x,C)\cup B_\rho(x',C)}
\left|P_0(x,z)-P_0(x',z)\right|
\left|Q_0(z,y)-Q_0(z,y')\right|\,d\mu(z)
\\
&\quad
\leq\int_{B_\rho(x,C)}\left|P_0(x,z)-P_0(x',z)\right|
\left|Q_0(z,y)-Q_0(z,y')\right|\,d\mu(z)
+\int_{B_\rho(x',C)}\cdots
\\
&\quad
\lesssim\int_{B_\rho(x,C)}\left[\rho(x,x')\right]^\varepsilon
\frac{1}{(V_\rho)_1(z)}\left[\rho(y,y')\right]^\varepsilon
\frac{1}{(V_\rho)_1(z)}\,d\mu(z)
+\int_{B_\rho(x',C)}\cdots,
\end{align*}
which, combined with Lemma~\ref{A3.2}(vii) and \eqref{upperdoub},
further implies that
\begin{align*}
&\left|\left[P_0Q_0(x,y)-P_0Q_0(x',y)\right]-\left[P_0Q_0(x,y')
-P_0Q_0(x',y')\right]\right|
\\
&\quad
\sim\int_{B_\rho(x,C)}\left[\rho(x,x')\right]^\varepsilon
\frac{1}{(V_\rho)_1(x)}\left[\rho(y,y')\right]^\varepsilon
\frac{1}{(V_\rho)_1(x)}\,d\mu(z)
\\
&\qquad+\int_{B_\rho(x',C)}\left[\rho(x,x')\right]^\varepsilon
\frac{1}{(V_\rho)_1(x')}\left[\rho(y,y')\right]^\varepsilon
\frac{1}{(V_\rho)_1(x')}\,d\mu(z)\\
&\quad\lesssim
\left[\rho(x,x')\right]^\varepsilon
\left[\rho(y,y')\right]^\varepsilon
\left[\frac{1}{(V_\rho)_1(x)}+\frac{1}{(V_\rho)_1(x')}\right].
\end{align*}
This finishes the proof of \eqref{418x}.
Thus, we have finished the proof of \eqref{3.5}
and hence Lemma~\ref{A3.25}.
\end{proof}

Let $(X,\mathbf{q},\mu)$ be a quasi-ultrametric space
of homogeneous type,
where $\mu$ is assumed to be a Borel-semiregular measure on $X$.
Assume that $0<\varepsilon\preceq{\mathrm{ind\,}}(X,\mathbf{q})$ with
${\mathrm{ind\,}}(X,\mathbf{q})$ as in \eqref{index}
and $\{\mathcal{S}_{2^{-k}}\}_{k\in\mathbb{Z}_+}$
is an inhomogeneous approximation
of the identity of order $\varepsilon$ as in Definition~\ref{DefIATI}.
Throughout this section, we always
let
$$\mathcal{D}_k:=\mathcal{S}_{2^{-k}}
-\mathcal{S}_{2^{-k+1}}
$$
for any given $k\in\mathbb{N}$,
$\mathcal{D}_0:=\mathcal{S}_1$, and
$\mathcal{D}_k:=0$ for any
given $k\in\mathbb{Z}\setminus\mathbb{Z}_+$.
To establish sharp inhomogeneous continuous
Calder\'on reproducing formulae,
similarly to \eqref{3.6},
by Proposition~\ref{A3.24}, we decompose the identity
operator $\mathcal{I}$ as follows:
for any fixed $N\in\mathbb{N}$ and for any $p\in(1,\infty)$,
\begin{align}\label{3.21}
\mathcal{I}
&=\sum_{k=-\infty}^{\infty}\mathcal{D}_k
=\sum_{k=-\infty}^\infty\left(\sum_{l=-\infty}^\infty
\mathcal{D}_{k+l}\right)\mathcal{D}_k
\nonumber\\
&=\sum_{k=-\infty}^{\infty}\mathcal{D}_k^N\mathcal{D}_k
+\sum_{k=-\infty}^{\infty}\sum_{|l|>N}\mathcal{D}_{k+l}\mathcal{D}_k
\nonumber\\
&=\sum_{k=0}^{\infty}\mathcal{D}_k^N\mathcal{D}_k
+\sum_{k=0}^{\infty}\sum_{|l|>N}\mathcal{D}_{k+l}\mathcal{D}_k
=:\mathcal{T}_N+\mathcal{R}_N
\end{align}
in $L^p(X)$, where, for any $k\in\mathbb{Z}_+$,
\begin{align}\label{DkN}
\mathcal{D}^N_k:=\sum_{|l|\leq N}\mathcal{D}_{k+l}=
\begin{cases}
\displaystyle\sum_{l=0}^{k+N}\mathcal{D}_l&\text{ if }k\in\{0,\ldots,N\},\\
\displaystyle\sum_{l=k-N}^{k+N}\mathcal{D}_l&\text{ if }k\in\{N+1,N+2,\ldots\}.
\end{cases}
\end{align}
From this and both (iv) and (v) of Definition~\ref{DefATI},
we deduce that the kernels $\{D_k^N\}_{k\in\mathbb{Z}_+}$
of operators $\{\mathcal{D}_k^N\}_{k\in\mathbb{Z}_+}$
satisfy that, for any $x\in X$,
\begin{align}\label{IDkN}
\int_XD^N_k(x,y)\,d\mu(y)
&=\int_XD^N_k(y,x)\,d\mu(y)\nonumber\\
&=
\begin{cases}
1\ \text{ if }k\in\{0,\ldots,N\},\\
0\ \text{ if }k\in\{N+1,N+2,\ldots\}.
\end{cases}
\end{align}

Next, we establish the boundedness of
$\mathcal{R}_N$ on $L^2(X)$ and ${\mathcal G}(\beta,\gamma)$.
To this end, we first show the $L^2(X)$ boundedness of $\mathcal{R}_N$
and some estimates related to the kernel of $\mathcal{R}_N$.

\begin{lemma}\label{A3.28}
Let $(X,\mathbf{q},\mu)$ be an ultra-RD-space,
$0<\varepsilon\preceq{\mathrm{ind\,}}(X,\mathbf{q})$
with ${\mathrm{ind\,}}(X,\mathbf{q})$ as in \eqref{index},
$N\in\mathbb{N}$, and $\varepsilon'\in(0,\varepsilon)$.
Then $\mathcal{R}_N$ defined in \eqref{3.21} is
an $\varepsilon'$-Calder\'on--Zygmund operator as in Definition~\ref{CZO} satisfying
$$
\left\|\mathcal{R}_N\right\|_{L^2(X)\to L^2(X)}
\leq C2^{-\varepsilon'N}
$$
with the positive constant $C$ independent of $N$,
and the kernel $R_N$ of $\mathcal{R}_N$ satisfies \eqref{CZO1},
\eqref{CZO3}, and \eqref{CZO4} with $r_0:=1$,
$\sigma\in(0,\infty)$,
and $C_{\mathcal{R}_N}:=C2^{(\varepsilon'-\varepsilon)N}$.
\end{lemma}

\begin{proof}
Let $\rho\in\mathbf{q}$ be a regularized quasi-ultrametric as in Definition~\ref{D2.3}.
In particular, $\rho$ is a symmetric quasi-ultrametric.
Since $D_0=S_1$ has no cancellation
[see Definition~\ref{DefATI}(v)], we write
\begin{align*}
\mathcal{R}_N&=\sum_{k>0,\,k+l>0,\,|l|>N}
\mathcal{D}_{k+l}\mathcal{D}_k
+\sum_{l>N}\mathcal{D}_l\mathcal{D}_0
+\sum_{k>N}\mathcal{D}_0\mathcal{D}_k\\
&=:\mathrm{I}+\mathrm{II}+\mathrm{III}.
\end{align*}

For $\mathrm{I}$,
following the proofs of Lemmas~\ref{A3.13} and~\ref{A3.15},
we conclude that $\mathrm{I}$
is an $\varepsilon'$-Calder\'on--Zygmund operator with its kernel
satisfying \eqref{CZO1} and $\|\mathrm{I}\|_{L^2(X)\to L^2(X)}
\leq C2^{-\varepsilon'N}$ with the positive constant $C$ as in Lemma~\ref{A3.13}.

Now, we prove that $\mathrm{I}$ satisfies \eqref{CZO3} with $r_0:=1$,
$\sigma\in(0,\infty)$,
and $C_{\mathcal{R}_N}:=C2^{(\varepsilon'-\varepsilon)N}$.
Let $x,y\in X$ satisfy $\rho(x,y)\ge1$.
On the one hand, if
$\sum_{k>0,\,k+l>0,\,
|l|>N}D_{k+l}D_k(x,y)=0$,
then $\mathrm{I}$ satisfies \eqref{CZO3}.
On the other hand,
if
$\sum_{k>0,\,k+l>0,\,
|l|>N}D_{k+l}D_k(x,y)\neq 0$,
then we find that there exist $k_0\in\mathbb{N}$
and $l_0\in\mathbb{Z}\cap(-k_0,\infty)$ such that
$D_{k_0+l_0}D_{k_0}(x,y)\neq0$,
which, combined with Lemma~\ref{A3.25},
further implies that
$\rho(x,y)\leq C'2^{-[k_0\land(k_0+l_0)]}<C'$,
where $C'\in[1,\infty)$ is the same positive constant as in Lemma~\ref{A3.25}.
From this, we deduce that, for any $\sigma\in(0,\infty)$,
$1\lesssim[\rho(x,y)]^{-\sigma}$,
where the implicit positive constant depends on $\sigma$.
By this and an argument similar
to that used in the proof of \eqref{RNCZO1},
we conclude that
\begin{align*}
\left|\sum_{k>0,\,k+l>0,\,
|l|>N}D_{k+l}D_k(x,y)\right|
&\lesssim2^{-N\varepsilon}\frac{1}{V_\rho(x,y)}
\left[\frac{1}{\rho(x,y)}\right]^\sigma\\
&\leq2^{(\varepsilon'-\varepsilon)N}\frac{1}{V_\rho(x,y)}
\left[\frac{1}{\rho(x,y)}\right]^\sigma.
\end{align*}
Thus, $\mathrm{I}$ satisfies \eqref{CZO3}.

Next, we show that $\mathrm{I}$ satisfies \eqref{CZO4} with $r_0:=1$,
$\sigma\in(0,\infty)$,
and $C_{\mathcal{R}_N}:=C2^{(\varepsilon'-\varepsilon)N}$.
Let $x,x',y\in X$ satisfy both
$\rho(x,y)\ge1$ and $\rho(x,x')\leq\frac{\rho(x,y)}{2C_\rho}$.
On the one hand, if $|D_{k+l}D_k(x,y)-D_{k+l}D_k(x',y)|=0$,
then $\mathrm{I}$ satisfies \eqref{CZO4}.
On the other hand,
from Definition~\ref{D2.1}(iii),
it follows that
\begin{align}\label{22.1.24x1}
\rho(x,y)\leq C_\rho\max\left\{\rho(x,x'),\,\rho(x',y)\right\}=C_\rho\rho(x',y).
\end{align}
If there exist $k_1\in\mathbb{N}$
and $l_1\in\mathbb{Z}\cap(-k_1,\infty)$ such that
$$
\left|D_{k_1+l_1}D_{k_1}(x,y)-D_{k_1+l_1}D_{k_1}(x',y)\right|\neq0,
$$
then, by the support condition of $D_{k_1+l_1}D_{k_1}$
(see also Lemma~\ref{A3.25})
and \eqref{22.1.24x1},
we find that
$$
1\leq\rho(x,y)\leq C'C_\rho2^{-[(k_1+l_1)\land k_1]}<C'C_\rho.
$$
From this, \eqref{upperdoub}, $\varepsilon'<\varepsilon$,
and Lemma~\ref{A3.12}(iii) (see also Lemma~\ref{A3.25}),
we infer that, for any $\sigma\in(0,\infty)$,
\begin{align*}
&\sum_{k>0,\,k+l>0,\,|l|>N}
\left|D_{k+l}D_k(x,y)-D_{k+l}D_k(x',y)\right|
\\
&\quad=\sum_{\genfrac{}{}{0pt}{}{k>0,k+l>0,|l|>N}{\rho(x,y)
\leq C'C_\rho2^{-[(k+l)\land k]}}}
\left|D_{k+l}D_k(x,y)-D_{k+l}D_k(x',y)\right|
\\
&\quad\lesssim\sum_{\genfrac{}{}{0pt}{}{k>0,k+l>0,|l|>N}{\rho(x,y)
\leq C'C_\rho2^{-[(k+l)\land k]}}}
2^{\varepsilon'[(l+k)\land k]-(\varepsilon-\varepsilon')|l|}
\left[\rho(x,x')\right]^{\varepsilon'}
\frac{1}{(V_\rho)_{2^{-[k\land(k+l)]}}(y)}
\\
&\quad\lesssim\sum_{\genfrac{}{}{0pt}{}{k>0,k+l>0,|l|>N}{\rho(x,y)
\leq C'C_\rho2^{-[(k+l)\land k]}}}
2^{(\varepsilon'-\varepsilon)|l|}\left[\frac{\rho(x,x')}{\rho(x,y)}\right]
^{\varepsilon'}\frac{1}{(V_\rho)_{2^{-[k\land(k+l)]}}(y)},
\end{align*}
which, together with Lemmas~\ref{A3.2}(vii) and
\ref{A3.14} [in the case where $a:=\lceil\log_{\frac{1}{2}}\frac{\rho(x,y)}{C'C_\rho}\rceil$ and $c_0:=(C'C_\rho)^{-1}$],
further implies that
\begin{align*}
&\sum_{k>0,\,k+l>0,\,|l|>N}
\left|D_{k+l}D_k(x,y)-D_{k+l}D_k(x',y)\right|
\\
&\quad\lesssim\sum_{\genfrac{}{}{0pt}{}{k>0,k+l>0,|l|>N}{\rho(x,y)
\leq C'C_\rho2^{-[(k+l)\land k]}}}
2^{(\varepsilon'-\varepsilon)|l|}\left[\frac{\rho(x,x')}{\rho(x,y)}\right]
^{\varepsilon'}\frac{1}{(V_\rho)_{2^{-[k\land(k+l)]}}(x)}
\\
&\quad\leq\sum_{|l|>N}2^{(\varepsilon'-\varepsilon)|l|}
\left[\frac{\rho(x,x')}{\rho(x,y)}\right]
^{\varepsilon'}
\sum_{k'=0}^{\lceil\log_{\frac{1}{2}}\frac{\rho(x,y)}{C'C_\rho}\rceil}
\frac{1}{(V_\rho)_{2^{-k'}}(x)}
\\
&\quad\lesssim2^{(\varepsilon'-\varepsilon)N}
\left[\frac{1}{\rho(x,y)}\right]^\sigma
\left[\frac{\rho(x,x')}{\rho(x,y)}\right]^{\varepsilon'}
\frac{1}{V_\rho(x,y)}.
\end{align*}
This proves that $\mathrm{I}$ satisfies
\eqref{CZO4}.

Now, we show that both $\mathrm{II}$ and $\mathrm{III}$ are
$\varepsilon'$-Calder\'on--Zygmund
operators with their kernels satisfying \eqref{CZO1}
as well as \eqref{CZO3} and \eqref{CZO4} with $r_0:=1$,
$\sigma\in(0,\infty)$,
and $C_{\mathcal{R}_N}:=C2^{(\varepsilon'-\varepsilon)N}$.
We first deal with $\mathrm{II}$.
By an argument similar to
that used in the estimation of \eqref{RN0},
we conclude that, for any $l\in(N,\infty)\cap\mathbb{N}$,
$$
\left\|\mathcal{D}_l\mathcal{D}_0\right\|_{L^2(X)\to L^2(X)}
\lesssim2^{-l\varepsilon},
$$
which, combined with $\varepsilon'<\varepsilon$, further implies
\begin{align*}
\|\mathrm{II}\|_{L^2(X)\to L^2(X)}
\leq\sum_{l>N}\left\|\mathcal{D}_l\mathcal{D}_0\right\|_{L^2(X)\to L^2(X)}
\lesssim2^{-N\varepsilon}\leq2^{-\varepsilon'N}.
\end{align*}
Thus, $\mathrm{II}$ is bounded on $L^2(X)$.

Next, we prove that $\mathrm{II}$ satisfies both
\eqref{CZk1} and \eqref{CZO3}.
Let $x,y\in X$ with $x\neq y$.
On the one hand, if $\sum_{l>N}|D_lD_0(x,y)|=0$,
then $\mathrm{II}$ satisfies \eqref{CZk1} and \eqref{CZO3}.
On the other hand,
if $\sum_{l>N}|D_lD_0(x,y)|\neq 0$, then
we find that there exists $l_2\in(N,\infty)\cap\mathbb{N}$
such that $D_{l_2}D_0(x,y)\neq0$,
which, together with Lemma~\ref{A3.25},
further implies that
$\rho(x,y)\leq C'2^{-(l_2\land 0)}=C'$.
From this, Lemma~\ref{A3.25},
and \eqref{upperdoub}, it follows that
\begin{align*}
\sum_{l>N}|D_lD_0(x,y)|
&\lesssim\sum_{l>N}2^{-l\varepsilon}\frac{1}{(V_\rho)_1(x)}
\sim 2^{-N\varepsilon}\frac{1}{(V_\rho)_1(x)}
\nonumber\\
&\lesssim2^{-N\varepsilon}\frac{1}{(V_\rho)_{C'}(x)}
\leq2^{-N\varepsilon}\frac{1}{V_\rho(x,y)}
\nonumber\\
&\lesssim2^{(\varepsilon'-\varepsilon)N}\frac{1}{V_\rho(x,y)}
\min\left\{1,\,\left[\frac{1}{\rho(x,y)}\right]^\sigma\right\},
\end{align*}
which further implies that $\mathrm{II}$ satisfies
\eqref{CZk1} and \eqref{CZO3}.

Now, we show that $\mathrm{II}$ satisfies \eqref{CZk2} and \eqref{CZO4}.
Let $x,y,x'\in X$ satisfy $\rho(x,x')\leq\frac{\rho(x,y)}{2C_\rho}$ with $x\neq y$.
On the one hand, if
$$\sum_{l>N}|D_lD_0(x,y)-D_lD_0(x',y)|+
\sum_{l>N}|D_lD_0(y,x)-D_lD_0(y,x')|=0,
$$
then $\mathrm{II}$ satisfies \eqref{CZk2} and \eqref{CZO4}.
On the other hand, if
\begin{align}
\label{1511}
\sum_{l>N}|D_lD_0(x,y)-D_lD_0(x',y)|+
\sum_{l>N}|D_lD_0(y,x)-D_lD_0(y,x')|\neq0,
\end{align}
then there exists $l_3\in\mathbb{Z}\cap(N,\infty)$
such that at least one of $|D_{l_3}D_0(y,x)|\neq0$, $|D_{l_3}D_0(y,x')|\neq0$,
$|D_{l_3}D_0(x,y)|\neq0$, and $|D_{l_3}D_0(x',y)|\neq0$ holds,
which, combined with Lemma~\ref{A3.25},
further implies that either $\rho(x,y)\leq C'2^{-l_3\land0}=C'$
or $\rho(x',y)\leq C'2^{-l_3\land0}=C'$.
If $\rho(x',y)\leq C'$, then
Definition~\ref{D2.1}(iii) implies that
$$
\rho(x,y)\leq C_\rho\max\left\{\rho(x,x'),\,\rho(x',y)\right\}
=C_\rho\rho(x',y)\leq C'C_\rho.
$$
Thus, \eqref{1511} implies that $\rho(x,y)\leq C'C_\rho$.
By this, Lemma~\ref{A3.25}
[see also both (ii) and (iii) of Lemma~\ref{A3.12}],
\eqref{upperdoub}, and Lemma~\ref{A3.2}(vii),
we conclude that, for any $\sigma\in(0,\infty)$,
\begin{align*}
&\sum_{l>N}|D_lD_0(x,y)-D_lD_0(x',y)|+
\sum_{l>N}|D_lD_0(y,x)-D_lD_0(y,x')|
\\
&\quad\lesssim\sum_{l>N}2^{\varepsilon'(l\land0)-|l|(\varepsilon-\varepsilon')}
\left[\frac{\rho(x,x')}{\rho(x,y)}\right]^{\varepsilon'}\frac{1}{(V_\rho)_1(y)}
\\
&\quad\sim\sum_{l>N}2^{l(\varepsilon'-\varepsilon)}
\left[\frac{\rho(x,x')}{\rho(x,y)}\right]^{\varepsilon'}\frac{1}{(V_\rho)_1(x)}
\\
&\quad\lesssim 2^{N(\varepsilon'-\varepsilon)}
\left[\frac{\rho(x,x')}{\rho(x,y)}\right]^{\varepsilon'}
\frac{1}{(V_\rho)_{C'C_\rho}(x)}
\\
&\quad
\lesssim2^{N(\varepsilon'-\varepsilon)}
\left[\frac{\rho(x,x')}{\rho(x,y)}\right]^{\varepsilon'}
\frac{1}{V_\rho(x,y)}
\min\left\{1,\,\left[\frac{1}{\rho(x,y)}\right]^\sigma\right\},
\end{align*}
which further implies that $\mathrm{II}$ satisfies \eqref{CZk2} and \eqref{CZO4}.

Finally, we prove that $\mathrm{II}$ satisfies \eqref{CZO1}.
Let $x,y,x',y'\in X$ with $x\neq y$
satisfy both
$\rho(x,x')\leq\frac{\rho(x,y)}{(2C_\rho)^2}$
and
$\rho(y,y')\leq\frac{\rho(x,y)}{(2C_\rho)^2}$.
On the one hand, if
\begin{align*}
\sum_{l>N}\left|\left[D_lD_0(x,y)-D_lD_0(x',y)\right]
-\left[D_lD_0(x,y')-D_lD_0(x',y')\right]\right|=0,
\end{align*}
then $\mathrm{II}$ satisfies \eqref{CZO1}.
On the other hand,
from Definition~\ref{D2.1}(iii),
we infer that
\begin{align*}
\rho(x,y)\leq C_\rho\max\left\{\rho(x,x'),\rho(x',y)\right\}
=C_\rho\rho(x',y),
\end{align*}
\begin{align*}
\rho(x,y)\leq C_\rho\max\left\{\rho(x,y'),\rho(y',y)\right\}=
C_\rho\rho(x,y'),
\end{align*}
and
\begin{align*}
\rho(x,y)&\leq C_\rho\max\left\{\rho(x,x'),\rho(x',y)\right\}
\\
&\leq C_\rho^2\max\left\{\rho(x,x'),\,\rho(x',y'),\,\rho(y',y)\right\}
=C_\rho^2\rho(x',y').
\end{align*}
By this and Lemma~\ref{A3.25}, we find that, if
there exists $l_4\in\mathbb{N}\cap(N,\infty)$ such that
\begin{align*}
\left|\left[D_{l_4}D_0(x,y)-D_{l_4}D_0(x',y)\right]
-\left[D_{l_4}D_0(x,y')-D_{l_4}D_0(x',y')\right]\right|\neq0,
\end{align*}
then $\rho(x,y)\leq C_\rho^2C'2^{-(l_4\land0)}=C_\rho^2C'$.
From this, Lemma~\ref{A3.25} [see also Lemma~\ref{A3.12}(iv)],
\eqref{upperdoub}, and Lemma~\ref{A3.2}(vii)
with the fact that $\rho(x,x')\leq\frac{\rho(x,y)}{(2C_\rho)^2}<C'$,
it follows that
\begin{align*}
&\sum_{l>N}\left|\left[D_lD_0(x,y)-D_lD_0(x',y)\right]
-\left[D_lD_0(x,y')-D_lD_0(x',y')\right]\right|
\\
&\quad\lesssim \sum_{l>N}2^{2\varepsilon'(l\land0)-|l|(\varepsilon-\varepsilon')}
\left[\rho(x,x')\rho(y,y')\right]^{\varepsilon'}
\left[\frac{1}{(V_\rho)_{2^{-(0\land l)}}(x)}
+\frac{1}{(V_\rho)_{2^{-(0\land l)}}(x')}\right]\\
&\quad\lesssim2^{N(\varepsilon'-\varepsilon)}
\left[\frac{\rho(x,x')}{\rho(x,y)}\right]^{\varepsilon'}
\left[\frac{\rho(y,y')}{\rho(x,y)}\right]^{\varepsilon'}
\frac{1}{(V_\rho)_{C_\rho^2C'}(x)}
\\
&\quad\leq2^{N(\varepsilon'-\varepsilon)}
\left[\frac{\rho(x,x')}{\rho(x,y)}\right]^{\varepsilon'}
\left[\frac{\rho(y,y')}{\rho(x,y)}\right]^{\varepsilon'}
\frac{1}{V_\rho(x,y)},
\end{align*}
which further implies that $\mathrm{II}$ satisfies \eqref{CZO1}.
The same estimates also hold for $\mathrm{III}$. Therefore, both
$\mathrm{II}$ and $\mathrm{III}$ are $\varepsilon'$-Calder\'on--Zygmund
operators with the kernel
satisfying \eqref{CZO1},
\eqref{CZO3}, and \eqref{CZO4}.
This finishes the proof of Lemma~\ref{A3.28}.
\end{proof}

Via an argument similar to that used in the proof of Lemma~\ref{RNRNM},
we have the following lemma; we omit the details.

\begin{lemma}\label{IRNM}
Let $(X,\mathbf{q},\mu)$ be an ultra-RD-space and
$0<\varepsilon\preceq{\mathrm{ind\,}}(X,\mathbf{q})$
with ${\mathrm{ind\,}}(X,\mathbf{q})$ as in \eqref{index}.
Let $\{\mathcal{S}_{2^{-k}}\}_{k\in\mathbb{Z}_+}$
be an inhomogeneous approximation
of the identity of order $\varepsilon$ as in Definition~\ref{DefIATI} and, for any fixed
$N\in\mathbb{N}$,
let $\mathcal{R}_N$ be the same as in \eqref{3.21}. For any
$M\in\mathbb{N}\cap(N,\infty)$, let
$$
\mathcal{R}_{N,M}:=\sum_{k=0}^M\sum_{N<|l|\leq M}
\mathcal{D}_{k+l}\mathcal{D}_k.
$$
Then the following assertions hold:
\begin{enumerate}
\item[\textup{(i)}]
All of the conclusions of Lemma~\ref{A3.28}
remain true for $\mathcal{R}_{N,M}$
and all of the involved positive constants are independent of $M$.

\item[\textup{(ii)}]
For any $x,y\in X$,
$$
\int_XR_{N,M}(x,z)\,d\mu(z)=0
=\int_XR_{N,M}(z,y)\,d\mu(z).
$$

\item[\textup{(iii)}]
For any given $p\in(1,\infty)$ and for any
$f\in L^p(X)$ and $x\in X$,
$$
\mathcal{R}_{N,M}f(x)=\int_XR_{N,M}(x,y)f(y)\,d\mu(y).
$$

\item[\textup{(iv)}]
For any $f\in{\mathcal G}(x_1,r,\beta,\gamma)$ [with $x_1\in X$,
$r\in(0,\infty)$, and $\beta,\gamma\in(0,\varepsilon)$],
the function sequence $\{\mathcal{R}_{N,M}f\}_{M=N+1}^\infty$
converges locally uniformly to some element,
denoted by $\widetilde{\mathcal{R}_N}f(x)$,
where $\widetilde{\mathcal{R}_N}f$ differs
from $\mathcal{R}_Nf$ at most on a set of $\mu$-measure 0.

\item[\textup{(v)}]
For any $x_1\in X$, $r\in(0,\infty)$, and $\beta,\gamma\in(0,\varepsilon)$,
the operator $\widetilde{\mathcal{R}_N}$
can be uniquely extended
from ${\mathcal G}(x_1,r,\beta,\gamma)$ to $L^2(X)$
and the extension operator coincides with $\mathcal{R}_N$.
Consequently, for any $f\in L^2(X)$,
$\widetilde{\mathcal{R}_N}f=\mathcal{R}_Nf$ in $L^2(X)$ and
pointwise $\mu$-almost everywhere on $X$.
\end{enumerate}
\end{lemma}

Applying Theorems~\ref{CZOonCG} and
\ref{CZOonGO} to $\mathcal{R}_{N,M}$
[which is valid given Lemma~\ref{IRNM}(i)--(iii)]
and then passing limit to $\mathcal{R}_N$,
we obtain the following conclusion;
we omit the details of its proof.

\begin{proposition}\label{IBRN}
Let $(X,\mathbf{q},\mu)$ be an ultra-RD-space.
Let $x_1\in X$, $r\in(0,\infty)$,
$0<\beta,\gamma<\varepsilon\preceq{\mathrm{ind\,}}(X,\mathbf{q})$,
and $\mathcal{R}_N$ be the same as in \eqref{3.21}.
Then there exist two positive constants
$C$ and $\delta$ such that
\begin{align*}
\left\|\mathcal{R}_N\right\|_{{\mathcal G}(x_1,1,\beta,\gamma)
\to{\mathcal G}(x_1,1,\beta,\gamma)}
+\left\|\mathcal{R}_N\right\|_{\mathring{\mathcal{G}}(x_1,r,\beta,\gamma)
\to\mathring{\mathcal{G}}(x_1,r,\beta,\gamma)}
\leq C2^{N\delta}.
\end{align*}
\end{proposition}

Now, we establish the following sharp
inhomogeneous continuous Calder\'on reproducing
formula on ultra-RD-spaces.
\index[words]{inhomogeneous continuous Calder\'on reproducing formula}

\begin{theorem}\label{A3.33}
Let $(X,\mathbf{q},\mu)$ be an ultra-RD-space,
where $\mu$ is assumed to be a Borel-semiregular measure on $X$.
Assume that $0<\varepsilon\preceq{\mathrm{ind\,}}(X,\mathbf{q})$
with ${\mathrm{ind\,}}(X,\mathbf{q})$ as in \eqref{index}
and $\{\mathcal{S}_{2^{-k}}\}_{k\in\mathbb{Z}_+}$
is an inhomogeneous approximation
of the identity of order $\varepsilon$
as in Definition~\ref{DefIATI} with kernels
$\{S_{2^{-k}}\}_{k\in\mathbb{Z}_+}$.
Let $\mathcal{D}_k:=\mathcal{S}_{2^{-k}}-\mathcal{S}_{2^{-k+1}}$
for any $k\in\mathbb{N}$ and
$\mathcal{D}_0:=\mathcal{S}_1$.
Then there exist $N\in\mathbb{N}$
and two families
$\{\widetilde{\mathcal{D}}_k\}_{k\in\mathbb{Z}_+}$ and
$\{\overline{\mathcal{D}}_k\}_{k\in\mathbb{Z}_+}$
of bounded linear operators on $L^2(X)$
with their kernels
$\{\widetilde{D}_k\}_{k\in\mathbb{Z}_+}$ and
$\{\overline{D}_k\}_{k\in\mathbb{Z}_+}$ such that, for any
$f\in{\mathcal G}^\varepsilon_0(\beta,\gamma)$
with $\beta,\gamma\in(0,\varepsilon)$
[resp. $L^p(X)$ with $p\in(1,\infty)$],
\begin{equation}\label{3.22}
f=\sum_{k=0}^{\infty}\widetilde{\mathcal{D}}_k\mathcal{D}_kf
=\sum_{k=0}^{\infty}\mathcal{D}_k\overline{\mathcal{D}}_kf,
\end{equation}
where all the series in \eqref{3.22} converges in
${\mathcal G}^\varepsilon_0(\beta,\gamma)$
[resp. $L^p(X)$ with $p\in(1,\infty)$]
and also converges pointwise $\mu$-almost everywhere on $X$.
Moreover, for any $k\in\mathbb{Z}_+$, $x,y\in X$,
and $\widetilde{\beta},\widetilde{\gamma}\in(0,\varepsilon)$,
\begin{align}
\label{943-1}
\widetilde{D}_k(\cdot,y)\in
{\mathcal G}(y,2^{-k},\widetilde{\beta},\widetilde{\gamma}),
\end{align}
\begin{align}
\label{943-2}
\int_X\widetilde{D}_k(x,z)\,d\mu(z)
&=\int_X\widetilde{D}_k(z,y)\,d\mu(z)\nonumber\\
&=
\begin{cases}
1\quad
\mathrm{if}\ k\in\{0,\ldots,N\},\\
0\quad
\mathrm{if}\ k\in\{N+1,N+2,\ldots\},
\end{cases}
\end{align}
and
$\widetilde{D}_k(x,y)=\overline{D}_k(y,x)$.
\end{theorem}

\begin{remark}
\begin{enumerate}
\item[\textup{(i)}]
In Theorem~\ref{A3.33}, for any $k\in\mathbb{Z}_+$,
all the estimates satisfied by
$\widetilde{D}_k$ and $\overline{D}_k$
are independent of $\beta$ and $\gamma$ with
$\beta,\gamma\in(0,\varepsilon)$ (but may depend on $\varepsilon$).

\item[\textup{(ii)}]
If $\rho\in\mathbf{q}$ is a metric,
then ${\mathrm{ind\,}}(X,\mathbf{q})\in[1,\infty]$
and,
for any $0<\beta,\gamma<\varepsilon\preceq{\mathrm{ind\,}}(X,\mathbf{q})$,
\eqref{3.22} holds in $\mathring{{\mathcal G}}^\varepsilon_0(\beta,\gamma)$,
and, if $\rho\in\mathbf{q}$ is an ultrametric, then,
for any $0<\beta,\gamma<\varepsilon<\infty$,
\eqref{3.22} holds in $\mathring{{\mathcal G}}^\varepsilon_0(\beta,\gamma)$.
However, Han et al. \cite[Theorem 3.26]{hmy08}
and He et al. \cite[Theorem 6.6]{HLYY} only showed
inhomogeneous continuous Calder\'on
reproducing formulae
for any $0<\beta,\gamma<\varepsilon<1$.
In this sense, Theorem~\ref{A3.33} is sharper than
\cite[Theorem 3.26]{hmy08} and \cite[Theorem 6.6]{HLYY}
because of the sharpness of both Theorems~\ref{ThmATI} and~\ref{CZOonCG}.
\end{enumerate}
\end{remark}

\begin{proof}[Proof of Theorem~\ref{A3.33}]
Let $\rho\in\mathbf{q}$ be a regularized quasi-ultrametric.
The proof of the present theorem
is quite similar to that of Theorem~\ref{A3.19} and hence
we only sketch some of the important steps.
By Lemma~\ref{A3.28} and Proposition~\ref{IBRN},
fix a sufficiently large integer $N\in\mathbb{N}$
such that, for any $x_1\in X$ and $r\in(0,\infty)$,
\begin{align*}
\max\left\{\left\|\mathcal{R}_N\right\|_{L^2(X)\to L^2(X)},\,
\left\|\mathcal{R}_N\right\|_{{\mathcal G}(x_1,1,\beta,\gamma)
\to{\mathcal G}(x_1,1,\beta,\gamma)},\,
\left\|\mathcal{R}_N\right\|_{\mathring{\mathcal{G}}(x_1,r,\beta,\gamma)
\to\mathring{\mathcal{G}}(x_1,r,\beta,\gamma)}\right\}
\leq\frac{1}{2},
\end{align*}
which further implies that $\mathcal{T}_N^{-1}
=(\mathcal{I}-\mathcal{R}_N)^{-1}$ exists
and is bounded, respectively, on $L^2(X)$,
${\mathcal G}(x_1,1,\beta,\gamma)$,
and $\mathring{\mathcal{G}}(x_1,r,\beta,\gamma)$
with their operator norms no more than 2.
For any $k\in\mathbb{Z}_+$,
let $\mathcal{D}_k^N$ be the same as in \eqref{DkN}.

We first prove \eqref{943-1}.
Applying an argument similar to that used
in the proof of Theorem~\ref{A3.19}, we conclude that,
for any $k\in\mathbb{Z}\cap(N,\infty)$,
$D^N_k(\cdot,y)\in\mathring{\mathcal{G}}(y,2^{-k},\beta,\gamma)$ and,
for any $k\in\mathbb{Z}\cap[0,N]$,
$D^N_k(\cdot,y)\in{\mathcal G}(y,1,\beta,\gamma)$.
For any $k\in\mathbb{Z}_+$ and $x,y\in X$, define
$$
\widetilde{D}_k(x,y):=\mathcal{T}^{-1}_N\left(D^N_k(\cdot,y)\right)(x).
$$
This and the above boundedness of $\mathcal{T}_N^{-1}$
on spaces of test functions imply that
\begin{align}\label{2219-1}
\widetilde{D}_k(\cdot,y)\in\mathring{\mathcal{G}}(y,2^{-k},\beta,\gamma)
\text{ if }
k\in\mathbb{Z}\cap(N,\infty)
\end{align}
and
\begin{align}\label{2219-2}
\widetilde{D}_k(\cdot,y)\in{\mathcal G}(y,1,\beta,\gamma)
\text{ if }
k\in\mathbb{Z}\cap[0,N].
\end{align}
Note that, from Remark~\ref{Rmtest}(i),
it follows that,
for any $k\in\{0,\ldots,N\}$,
$${\mathcal G}(y,1,\beta,\gamma)={\mathcal G}(y,2^{-k},\beta,\gamma)$$
with the equivalent norms depending only on $N$,
which, together with both \eqref{2219-1} and \eqref{2219-2},
completes the proof of \eqref{943-1}.

Next, we turn to show \eqref{943-2}.
Similarly to the proof of the homogeneous case in Theorem~\ref{A3.19},
for any $k\in\mathbb{Z}\cap(N,\infty)$
and $x\in X$, we have
$$
\int_X\widetilde{D}_k(y,x)\,d\mu(y)
=\int_X\widetilde{D}_k(x,y)\,d\mu(y)=0.
$$
Thus, for any $k\in\mathbb{Z}\cap(N,\infty)$,
$\widetilde{D}_k$
satisfies all the conditions of Theorem~\ref{A3.33}.
For any $k\in\mathbb{Z}\cap[0,N]$ and $x,y\in X$,
write
\begin{align*}
\widetilde{D}_k(x,y)=\mathcal{T}^{-1}_N\left(D_k^N(\cdot,y)\right)(x)
=\sum_{i=0}^{\infty}(\mathcal{R}_N)^i\left(D_k^N(\cdot,y)\right)(x).
\end{align*}
By Lemma~\ref{IRNM} and
the Lebesgue dominated convergence theorem
[used along with the estimates in Definition~\ref{testfunction}(i)
and Lemma~\ref{A3.2}(ii)],
we find that, for any $y\in X$,
\begin{align*}
\mathcal{R}_N\left(D_k^N(\cdot,y)\right)
\in\mathring{\mathcal{G}}(y,1,\beta,\gamma).
\end{align*}
From this and the boundedness of $\mathcal{R}_N$
on $\mathring{\mathcal{G}}(y,1,\beta,\gamma)$, we deduce that,
for any $j\in\mathbb{N}$
and $y\in X$,
\begin{align}\label{9241}
\int_X(\mathcal{R}_N)^j\left(D_k^N(\cdot,y)\right)(x)\,d\mu(x)=0.
\end{align}
On the other hand, by Fubini's theorem,
\eqref{IDkN}, and Lemma~\ref{IRNM}(ii),
we conclude that, for any $M\in\mathbb{N}$ and $x\in X$,
\begin{align*}
&\int_X\mathcal{R}_{N,M}\left(D_k^N(\cdot,y)\right)(x)\,d\mu(y)\\
&\quad=\int_X\int_XR_{N,M}(x,z)D_k^N(z,y)\,d\mu(y)\,d\mu(z)=0.
\end{align*}
From an argument similar to that used in the proof of \eqref{919},
we infer that, for any $j\in\mathbb{N}$ and $x\in X$,
\begin{align}
\label{9242}
\int_X(\mathcal{R}_N)^j\left(D_k^N(\cdot,y)\right)(x)\,d\mu(y)=0.
\end{align}
By the Lebesgue dominated convergence theorem,
\eqref{9241}, \eqref{9242}, and \eqref{IDkN}, we find that,
for any $k\in\mathbb{Z}\cap[0,N]$ and
$x\in X$,
\begin{align*}
\int_X\widetilde{D}_k(x,y)\,d\mu(y)
&=\int_X\sum_{j=0}^\infty(\mathcal{R}_N)^j
\left(D_k^N(\cdot,y)\right)(x)\,d\mu(y)
\\
&=\sum_{j=0}^\infty\int_X(\mathcal{R}_N)^j
\left(D_k^N(\cdot,y)\right)(x)\,d\mu(y)
\\
&=\int_XD_k^N(x,y)\,d\mu(y)=1
\end{align*}
and
\begin{align*}
\int_X\widetilde{D}_k(y,x)\,d\mu(y)=1.
\end{align*}
Thus, for any $k\in\mathbb{Z}\cap[0,N]$,
$\widetilde{D}_k$
satisfies all the conditions of Theorem~\ref{A3.33}.

Now, it remains to prove that
the series in \eqref{3.22} converges in both
${\mathcal G}^\varepsilon_0(\beta,\gamma)$
with
$\beta,\gamma\in(0,\varepsilon)$
and $L^p(X)$ for $p\in(1,\infty)$.
We only show the convergence of the first series in \eqref{3.22}
because the proof of the convergence of the second
series in \eqref{3.22} is similar to
that in Theorem~\ref{A3.19}.
To verify that the first series in \eqref{3.22} converges
in ${\mathcal G}^\varepsilon_0(\beta,\gamma)$
with $\beta,\gamma\in(0,\varepsilon)$, note that,
for any $L\in\mathbb{N}\cap(N+1,\infty)$ and
$f\in{\mathcal G}(\beta',\gamma')$
with both $\beta'\in(\beta,\varepsilon)$
and $\gamma'\in(\gamma,\varepsilon)$, we have
[by applying an argument similar to that used in \eqref{jj-4u5}]
\begin{equation}\label{3.23}
\left\|\sum_{k=0}^{L}\widetilde{\mathcal{D}}_k\mathcal{D}_kf
-f\right\|_{{\mathcal G}(\beta,\gamma)}
=\left\|\mathcal{T}^{-1}_N\left(\sum_{k=L+1}^{\infty}
\mathcal{D}^N_k\mathcal{D}_k\right)f\right\|_{{\mathcal G}(\beta,\gamma)}.
\end{equation}
To prove \eqref{3.23}
tends to $0$ as $L\to\infty$,
from the boundedness of $\mathcal{T}_N^{-1}$ on ${\mathcal G}(\beta,\gamma)$,
we only need to verify that there exists $\sigma\in(0,\infty)$
such that, for any $k\in\mathbb{Z}\cap(L,\infty)
\subset(N+1,\infty)$,
\begin{equation}\label{3.24}
\left\|\mathcal{D}^N_k\mathcal{D}_kf\right\|_{{\mathcal G}(\beta,\gamma)}
\lesssim 2^{-\sigma k}\|f\|_{{\mathcal G}(\beta',\gamma')}.
\end{equation}
Note that, by \eqref{IDkN}, we conclude that,
for any $k\in\mathbb{Z}\cap(L,\infty)
\subset(N+1,\infty)$, the kernel $D_k^N$
has the cancellation property.
Therefore, from an argument similar to that
used in the proof of \eqref{3.16}, we deduce
that \eqref{3.24} holds,
which, combined with \eqref{3.23},
further implies that
\begin{align*}
\lim_{L\to\infty}
\left\|\sum_{k=0}^{L}
\widetilde{\mathcal{D}}_k\mathcal{D}_kf-f\right\|_{{\mathcal G}(\beta,\gamma)}
\lesssim\lim_{L\to\infty}\sum_{k=L+1}^\infty
\left\|\mathcal{D}^N_k\mathcal{D}_kf\right\|_{{\mathcal G}(\beta,\gamma)}
\lesssim\lim_{L\to\infty}2^{-\sigma L}\|f\|_{{\mathcal G}(\beta',\gamma')}=0.
\end{align*}
This finishes the proof of the convergence of the first series of \eqref{3.22}
in ${\mathcal G}^\varepsilon_0(\beta,\gamma)$ for any
$f\in{\mathcal G}(\beta',\gamma')$
with both $\beta'\in(\beta,\varepsilon)$
and $\gamma'\in(\gamma,\varepsilon)$.

Next, we show the convergence of the first series of \eqref{3.22}
in ${\mathcal G}^\varepsilon_0(\beta,\gamma)$
for any $f\in{\mathcal G}^\varepsilon_0(\beta,\gamma)$.
To this end, for any $L\in\mathbb{N}\cap(N,\infty)$, let
\begin{align*}
\widetilde{\mathcal{T}}_L:=\sum_{k=0}^L\mathcal{D}_k^N
\mathcal{D}_k
=\sum_{k=0}^L\sum_{|l|\leq N}\mathcal{D}_{k+l}\mathcal{D}_k.
\end{align*}
By an argument similar to that used in the
proofs of Theorems~\ref{A3.19} and~\ref{CZOonCG},
we find that, for any $L\in\mathbb{N}\cap(N,\infty)$,
$\widetilde{\mathcal{T}}_L$ is bounded on ${\mathcal G}(\beta,\gamma)$
with its operator norm
bounded by a positive constant independent of $L$.
From this and a standard density argument used in the proof of
Theorem~\ref{A3.19}, we infer that $\mathcal{R}_N$
is bounded on ${\mathcal G}^\varepsilon_0(\beta,\gamma)$ with
its operator norm
\begin{align*}
\left\|\mathcal{R}_N\right\|_{{\mathcal G}^\varepsilon_0(\beta,\gamma)
\to{\mathcal G}^\varepsilon_0(\beta,\gamma)}
\leq\left\|\mathcal{R}_N\right\|_{{\mathcal G}(\beta,\gamma)
\to{\mathcal G}(\beta,\gamma)}\leq\frac{1}{2},
\end{align*}
which, together with \eqref{3.21},
further implies that $\mathcal{T}_N^{-1}$ is bounded on
${\mathcal G}^\varepsilon_0(\beta,\gamma)$. By this
and an argument similar to that used in the proof of
Theorem~\ref{A3.19}, we conclude that
the first series in \eqref{3.22}
converges in ${\mathcal G}^\varepsilon_0(\beta,\gamma)$.

Finally, we prove the convergence of first series in \eqref{3.22}
in $L^p(X)$ for any $p\in(1,\infty)$.
Let $p\in(1,\infty)$.
For any $L\in\mathbb{N}\cap(N,\infty)$, define
$$
\overline{\mathcal{T}}_L:=\sum_{k=0}^L\widetilde{\mathcal{D}}_k
\mathcal{D}_k.
$$
From an argument similar to that used in the proof
of Theorem~\ref{A3.19}, we deduce that, for any
$L\in\mathbb{N}\cap(N,\infty)$,
$\sum_{k=N+1}^L\widetilde{\mathcal{D}}_k\mathcal{D}_k$ is a
$\theta$-Calder\'on--Zygmund operator
for some $\theta\in(0,\gamma\land\beta)$.
So the well-known Calder\'on--Zygmund
theory on spaces of homogeneous type
developed in \cite[pp.\,74--75, Theorem 2.4]{CW1}
further implies that $\sum_{k=N+1}^L\widetilde{\mathcal{D}}_k
\mathcal{D}_k$
is bounded on $L^p(X)$
with its operator norm
bounded by a positive constant independent of $L$.
By this and the boundedness of
$\sum_{k=0}^N\widetilde{\mathcal{D}}_k\mathcal{D}_k$
on $L^p(X)$ [using Theorem~\ref{corATI}(i)], we find that
$\overline{\mathcal{T}}_L$ is bounded on $L^p(X)$
with its operator norm
bounded by a positive constant independent of $L$.
From this, Lemma~\ref{density}(i),
a standard density argument,
and an argument similar to that used in the proof of
Theorem~\ref{A3.19}, we infer that
the convergence of first series in \eqref{3.22}
in $L^p(X)$ holds,
which completes the proof of Theorem~\ref{A3.33}.
\end{proof}

By Theorem~\ref{A3.33} and a standard duality argument,
we obtain the following sharp inhomogeneous continuous Calder\'on
reproducing formula on spaces of distributions;
we omit the details of its proof.

\begin{theorem}
Let $(X,\mathbf{q},\mu)$ be an ultra-RD-space,
where $\mu$ is assumed to be a Borel-semiregular measure on $X$.
Assume that $\varepsilon$, $\beta$, $\gamma$,
and $\{\mathcal{D}_k\}_{k\in\mathbb{Z}_+}$
are the same as in Theorem~\ref{A3.33}.
Then, for any
$f\in({{\mathcal G}}_0^\varepsilon(\beta,\gamma))'$,
\eqref{3.22} also holds in
$({{\mathcal G}}_0^\varepsilon(\beta,\gamma))'$.
\end{theorem}

\begin{remark}\label{Dkf3}
Let $(X,\mathbf{q},\mu)$ be a quasi-ultrametric space
of homogeneous type. Assume that
$0<\beta,\gamma<\varepsilon\preceq\mathrm{ind\,}(X,\mathbf{q})$.
Let $f\in({\mathcal G}_0^\varepsilon(\beta,\gamma))'$,
$\{\mathcal{S}_t\}_{t\in(0,\infty)}$ be an approximation of the identity of order
$\varepsilon$ as in Definition~\ref{DefATI}
with kernels denoted by $\{S_t\}_{t\in(0,\infty)}$,
and $\langle\cdot,\cdot\rangle$ denote the natural
distributional pairing. Similarly to Remark~\ref{Dkf2},
there are two ways to define
$\mathcal{S}_tf\in({\mathcal G}_0^\varepsilon(\beta,\gamma))'$
for any $t\in(0,\infty)$:
\begin{enumerate}
\item[\rm(i)]
We can define $\mathcal{S}_tf$ as the distribution induced by
the following function that is given by setting, for any $x\in X$,
$$
\mathcal{S}_tf(x):=\langle f,S_t(x,\cdot)\rangle.
$$
More precisely, in view of Lemma~\ref{dens},
for any $\phi\in{\mathcal G}_0^\varepsilon(\beta,\gamma)$, we let
$$
\left\langle\mathcal{S}_tf,\phi\right\rangle
:=\lim_{j\to\infty}\int_X\mathcal{S}_tf(y)\phi_j(y)\,d\mu(y),
$$
where $\{\phi_j\}_{j\in\mathbb{N}}$ is any  a sequence in
$\mathcal{G}_{\mathrm{b}}(\varepsilon,\varepsilon)
\subset{\mathcal G}_0^\varepsilon(\beta,\gamma)$ such that
$$
\lim_{j\to\infty}\left\|\phi-\phi_j\right\|_{\mathcal{G}(\beta,\gamma)}=0.
$$

\item[\rm(ii)]
We can also define $\mathcal{S}_tf$
by duality; that is, for any
$\phi\in{\mathcal G}_0^\varepsilon(\beta,\gamma)$,
let
\begin{align*}
\left\langle\mathcal{S}_tf,\phi\right\rangle:=
\left\langle f,\mathcal{S}_t^*\phi\right\rangle,
\end{align*}
where $\mathcal{S}_t^*$ is the integral operator defined by setting, for any $x\in X$,
$$
\mathcal{S}_t^*\phi(x):=\int_XS_t(y,x)\phi(y)\,d\mu(y).
$$
\end{enumerate}
Moreover, we can also define $\mathcal{D}_kf\in({\mathcal G}_0^\varepsilon(\beta,\gamma))'$
for any $k\in\mathbb{Z}$ similarly,
where $\mathcal{D}_k:=\mathcal{S}_{2^{-k}}-\mathcal{S}_{2^{-k+1}}$.
\end{remark}

Applying Lemma~\ref{dens} and
an argument similar to that used in the proof of Proposition~\ref{1651},
we can show that both of the definitions of $\mathcal{S}_tf$
for any $t\in(0,\infty)$ and $f\in({\mathcal G}_0^\varepsilon(\beta,\gamma))'$
in (i) and (ii) of Remark~\ref{Dkf3} are well defined and actually coincide;
we omit the details.

\begin{proposition}
Let the notation be as in Remark~\ref{Dkf3}.
Then, for any $t\in(0,\infty)$ and $f\in({\mathcal G}_0^\varepsilon(\beta,\gamma))'$,
the definitions of $\mathcal{S}_tf$
in (i) and (ii) of Remark~\ref{Dkf3} are well defined and coincide
in $({\mathcal G}_0^\varepsilon(\beta,\gamma))'$.
\end{proposition}

\section{Homogeneous Discrete
Calder\'on Reproducing Formulae}
\label{SS3.4}

In this section,
we aim to establish sharp
homogeneous discrete
Calder\'on reproducing formulae
on the ultra-RD-space $(X,\mathbf{q},\mu)$,
where $\mu$ is assumed to be a Borel-semiregular measure on $X$,
which is accomplished in Theorems~\ref{A3.47} and~\ref{HD}.
Throughout this section,
we always assume
$\mathrm{diam}_\rho(X)=\infty$
for some $\rho\in\mathbf{q}$
and hence for all $\rho\in\mathbf{q}$.
Fix a sufficiently large integer $N\in\mathbb{N}$.
For any $k\in\mathbb{Z}$,
let $\mathcal{D}_k^N$ and $\mathcal{D}_k$
(with their kernels $D_k^N$ and $D_k$) be the same as in \eqref{3.6}.
We introduce the following
\emph{discrete Riemann sum operator}\index[words]{discrete Riemann sum operator}
$\mathfrak{S}_N$ by setting,
for any $f\in L^2(X)$ and $x\in X$,
\begin{equation}\label{DRS}
\mathfrak{S}_Nf(x)\index[symbols]{S@$\mathfrak{S}$}
:=\sum_{k=-\infty}^{\infty}\sum_{\tau\in I_k}
\sum_{v=1}^{N(k,\tau)}
\int_{Q^{k,v}_\tau}
D^N_k(x,y)\,d\mu(y)\mathcal{D}_kf(y^{k,v}_\tau),
\end{equation}
where $y^{k,v}_\tau$
is an arbitrary point in $Q^{k,v}_\tau$.

Our first order of business is to prove that, for any $f\in L^2(X)$,
the function $\mathfrak{S}_Nf$, introduced in \eqref{DRS},
is well defined. Indeed, we show that
the discrete Riemann sum operator $\mathfrak{S}_N$ is bounded on both $L^2(X)$ and
$\mathring{\mathcal{G}}(x_1,r,\beta,\gamma)$
and has a bounded inverse on $\mathring{\mathcal{G}}(x_1,r,\beta,\gamma)$.
To this end, we define $\mathfrak{R}_N:=\mathcal{I}-\mathfrak{S}_N$
and first establish some estimates
on the kernel of the operator $\mathfrak{R}_N$.
Note that, for any $f\in L^2(X)$
and $x\in X$, we write
\begin{align}
\label{RGRN}
\mathfrak{R}_Nf(x)
&=f(x)-
\sum_{k=-\infty}^{\infty}\sum_{\tau\in I_k}
\sum_{v=1}^{N(k,\tau)}
\int_{Q^{k,v}_\tau}D^N_k(x,y)\,d\mu(y)
\mathcal{D}_kf(y^{k,v}_\tau)
\nonumber\\
&=\sum_{k=-\infty}^{\infty}\sum_{|l|>N}
\mathcal{D}_{k+l}\mathcal{D}_kf(x)
+\sum_{k=-\infty}^{\infty}\mathcal{D}^N_k\mathcal{D}_kf(x)
\nonumber\\
&\quad-\sum_{k=-\infty}^{\infty}\sum_{\tau\in I_k}
\sum_{v=1}^{N(k,\tau)}
\int_{Q^{k,v}_\tau}D^N_k(x,y)\,d\mu(y)
\mathcal{D}_kf(y^{k,v}_\tau)
\nonumber\\
&=\sum_{k=-\infty}^{\infty}\sum_{|l|>N}
\mathcal{D}_{k+l}\mathcal{D}_kf(x)
\nonumber\\
&\quad+\sum_{k=-\infty}^{\infty}\sum_{\tau\in I_k}
\sum_{v=1}^{N(k,\tau)}
\int_{Q^{k,v}_\tau}
D^N_k(x,y)\left[\mathcal{D}_kf(y)
-\mathcal{D}_kf(y^{k,v}_\tau)\right]\,d\mu(y)
\nonumber\\
&=:\mathcal{R}_Nf(x)+\sum_{k=-\infty}^{\infty}
\mathfrak{G}_{N,k}f(x)
=:\mathcal{R}_Nf(x)+\mathfrak{G}_Nf(x),
\end{align}
where $\mathcal{R}_N$ is the same as in \eqref{3.6}.
For any $k\in\mathbb{Z}$,
let $G_{N,k}$ be the kernel of $\mathfrak{G}_{N,k}$
and
$G_N:=\sum_{k=-\infty}^{\infty}G_{N,k}$
the kernel of $\mathfrak{G}_N$.
Clearly, for any $x,y\in X$,
\begin{align}
\label{G}
G_N(x,y)&=\sum_{k=-\infty}^{\infty}G_{N,k}(x,y)
\nonumber\\
&=\sum_{k=-\infty}^{\infty}\left\{\sum_{\tau\in I_k}
\sum_{v=1}^{N(k,\tau)}
\int_{Q^{k,v}_\tau}D^N_k(x,z)
\left[D_k(z,y)-D_k(y^{k,v}_\tau,y)\right]\,d\mu(z)\right\}.
\end{align}
From Lemma~\ref{A3.13} and Proposition~\ref{A3.16},
we deduce that $\mathcal{R}_N$ is bounded on both
$L^2(X)$ and $\mathring{{\mathcal G}}(x_1,r,\beta,\gamma)$
with $x_1\in X$, $r\in(0,\infty)$, and $\beta,\gamma\in(0,\varepsilon)$.
Therefore, to prove that, for any $f\in L^2(X)$
[resp. $\mathring{{\mathcal G}}(x_1,r,\beta,\gamma)$],
$\mathfrak{S}_N(f)$ in \eqref{DRS} is a well-defined function in
$L^2(X)$ [resp. $\mathring{{\mathcal G}}(x_1,r,\beta,\gamma)$],
it suffices to show that $\mathfrak{G}_N$ is
bounded on $L^2(X)$ and $\mathring{{\mathcal G}}(x_1,r,\beta,\gamma)$
with $x_1\in X$, $r\in(0,\infty)$, and $\beta,\gamma\in(0,\varepsilon)$.
To this end, we need the following lemma.

\begin{lemma}\label{A3.41}
Let $(X,\mathbf{q},\mu)$ be an ultra-RD-space,
where $\mu$ is assumed to be a Borel-semiregular measure on $X$.
Assume that $0<\varepsilon\preceq{\mathrm{ind\,}}(X,\mathbf{q})$
with ${\mathrm{ind\,}}(X,\mathbf{q})$ as in \eqref{index} and
$\{\mathcal{S}_{2^{-k}}\}_{k\in\mathbb{Z}}$
is a homogeneous approximation
of the identity of order $\varepsilon$
in Definition~\ref{DefHATI}.
Let $j\in\mathbb{N}$ be the same as in \eqref{2cubej} and,
for any $k\in\mathbb{Z}$ and $x,y\in X$,
let $G_{N,k}$ be the kernel of $\mathfrak{G}_{N,k}$ in \eqref{G}.
Then, for any given $\rho\in\mathbf{q}$ having the property that
all $\rho$-balls are $\mu$-measurable
and for any $k,j\in\mathbb{Z}$,
there exist two constants $C,C'\in[1,\infty)$,
independent of $k$, $j$, and $y_\tau^{k,v}$,
where $\{y_\tau^{k,v}\}_{k\in\mathbb{Z},\tau\in I_k,v\in\{1,\ldots,N(k,\tau)\}}$
is the same as in \eqref{DRS}, such that
\begin{enumerate}
\item[\textup{(i)}]
For any $x,y\in X$,
$$
|G_{N,k}(x,y)|\leq
C2^{-j\varepsilon}\frac{1}{(V_\rho)_{2^{-k}}(x)}
$$
and, if $\rho(x,y)\ge C'2^{-k}$, then $G_{N,k}(x,y)=0$.

\item[\textup{(ii)}]
For any $x,y,x'\in X$ with
$\rho(x,x')\leq 2^{-k}$ or
$\rho(x,x')\leq\frac{2^{-k}+\rho(x,y)}{2C_\rho}$
with $C_\rho\in[1,\infty)$ as in \eqref{Crho},
$$
|G_{N,k}(x,y)-G_{N,k}(x',y)|
\leq C2^{(k-j)\varepsilon}
\left[\rho(x,x')\right]^{\varepsilon}
\frac{1}{(V_\rho)_{2^{-k}}(x)}.
$$

\item[\textup{(iii)}]
For any $x,y,x',y'\in X$ with
$\rho(x,x')\leq 2^{-k}$ and $\rho(y,y')\leq 2^{-k}$
or
$\rho(x,x')\leq\frac{2^{-k}+\rho(x,y)}{(2C_\rho)^2}$
and
$\rho(y,y')\leq\frac{2^{-k}+\rho(x,y)}{(2C_\rho)^2}$,
\begin{align*}
&\left|\left[G_{N,k}(x,y)-G_{N,k}(x',y)\right]
-\left[G_{N,k}(x,y')-G_{N,k}(x',y')\right]\right|
\\
&\quad\leq C2^{(2k-j)\varepsilon}
\left[\rho(x,x')\right]^{\varepsilon}
\left[\rho(y,y')\right]^{\varepsilon}
\frac{1}{(V_\rho)_{2^{-k}}(x)}.
\end{align*}
\item[\textup{(iv)}]
For any $x\in X$,
$$
\int_XG_{N,k}(x,z)\,d\mu(z)=0
=\int_XG_{N,k}(z,x)\,d\mu(z).
$$
\end{enumerate}
\end{lemma}

\begin{proof}
Let $N$ be a fixed large integer.
Without loss of generality, we may assume that
$\rho$ is a regularized quasi-ultrametric as in Definition~\ref{D2.3}.

We first prove (iv).
By Lemma~\ref{Pk}(v), we find that, for any $x\in X$,
$$
\int_XD_k(x,y)\,d\mu(y)=0.
$$
From this, the definition of $G_{N,k}(x,y)$, and Fubini's theorem,
it follows that (iv) holds.

Now, we show (i).
Recall that, from Lemma~\ref{Pk}, we deduce that
both $D_k$ and $D_k^N$ satisfy all the conclusions
in Lemma~\ref{Pk} but with the constants in Lemma~\ref{Pk} depending on $N$.
In particular, there exists a positive
constant $\widetilde{C}$, depending on $N$, such that, for any $x,y\in X$
satisfying $\rho(x,y)\ge\widetilde{C}2^{-k}$,
$$
D_k^N(x,y)=0=D_k(x,y).
$$
On the one hand, by both (iv) and (v) of Lemma~\ref{DyadicSys},
we find that, for any given $y^{k,v}_\tau\in Q^{k,v}_\tau$
and for any $z\in Q^{k,v}_\tau$ and $y\in X$,
\begin{align}\label{rhozy}
\rho(z,y^{k,v}_\tau)\leq\mathrm{diam}_\rho(Q^{k,v}_\tau)
\sim2^{-(k+j)},
\end{align}
where $j\in\mathbb{N}$ is the same as in \eqref{2cubej}.
Therefore, using \eqref{upperdoub},
we conclude that
$$
(V_\rho)_{2^{-k}}(z)\lesssim(V_\rho)_{2^{-k}}(y^{k,v}_\tau),
$$
where the implicit positive constant is independent of $j$.
From this, the regularity condition of $D_k$
[see Lemma~\ref{Pk}(ii)],
both the size condition and the support condition of $D^N_k$
[see Lemma~\ref{Pk}(i)],
and Lemma~\ref{DyadicSys}(i),
we infer that, for any $x,y\in X$,
\begin{align*}
|G_{N,k}(x,y)|
&\leq\sum_{\tau\in I_k}\int_{Q^{k}_\tau}
\left|D^N_k(x,z)\right|
\left|D_k(z,y)-D_k(y^{k,v}_\tau,y)\right|\,d\mu(z)\\
&\lesssim\int_{B_\rho(x,\widetilde{C}2^{-k})}
\frac{1}{(V_\rho)_{2^{-k}}(x)}2^{k\varepsilon}
2^{-(k+j)\varepsilon}
\frac{1}{(V_\rho)_{2^{-k}}(z)}\,d\mu(z)\\
&\lesssim 2^{-j\varepsilon}\frac{1}{(V_\rho)_{2^{-k}}(x)}
\frac{(V_\rho)_{\widetilde{C}2^{-k}}(x)}{(V_\rho)_{2^{-k}}(x)}
\quad\text{by Lemma~\ref{A3.2}(vii)}\\
&\lesssim 2^{-j\varepsilon}\frac{1}{(V_\rho)_{2^{-k}}(x)}
\quad\text{by \eqref{upperdoub}}.
\end{align*}
Recall that $\widetilde{C}$ is a positive
constant such that, for any $x,y\in X$
satisfying $\rho(x,y)\ge\widetilde{C}2^{-k}$,
$$
D_k^N(x,y)=0=D_k(x,y).
$$
Note that,
using Definition~\ref{D2.1}(iii), we find that,
for any given $y_\tau^{k,v}\in Q_\tau^{k,v}$
and for any $x,y\in X$ and $z\in Q_\tau^{k,v}$,
\begin{align}
\label{2021}
\rho(x,y)
&\leq C_\rho\max\left\{\rho(x,z),\,\rho(z,y)\right\}
\nonumber\\
&\leq C_\rho^2\max\left\{\rho(x,z),\,\rho(z,y_\tau^{k,v}),\,
\rho(y_\tau^{k,v},y)\right\}.
\end{align}
Choose $C':=C_\rho^3(\max\{C_r2^{-j},\,\widetilde{C}\}+1)$,
where $C_r$ is the same as in Lemma~\ref{DyadicSys}.
We claim that $G_{N,k}(x,y)=0$ if
$\rho(x,y)\ge C'2^{-k}$. Let $x,y\in X$ satisfy
$\rho(x,y)\ge C'2^{-k}$. Then, based on \eqref{2021},
we consider the following
three cases.

\emph{Case (1)} $\rho(x,z)\ge(\max\{C_r2^{-j},\,\widetilde{C}\}+1)C_\rho2^{-k}$.
In this case, the support condition of $D_k^N$
[see also Lemma~\ref{Pk}(i)] implies that that
$D_k^N(x,z)=0$ and hence $G_{N,k}(x,y)=0$.

\emph{Case (2)}
$\rho(z,y_\tau^{k,v})\ge(\max\{C_r2^{-j},\,\widetilde{C}\}+1)C_\rho2^{-k}$.
In this case, from Lemma~\ref{DyadicSys}(iv) and
the fact that $z,y_\tau^{k,v}\in Q_\tau^{k,v}$,
it follows that
$$
\rho(z,y_\tau^{k,v})>\mathrm{diam}_\rho(Q_\tau^{k,v})
\ge\rho(z,y_\tau^{k,v}),
$$
which is a contradiction.

\emph{Case (3)}
$\rho(y_\tau^{k,v},y)\ge(\max\{C_r2^{-j},\,\widetilde{C}\}+1)C_\rho2^{-k}$.
In this case, on the one hand,
by Lemma~\ref{Pk}(i), we conclude that
\begin{align}
\label{925}
D_k(y_\tau^{k,v},y)=0.
\end{align}
On the other hand, from Definition~\ref{D2.1}(iii)
and Lemma~\ref{DyadicSys}(iv),
we deduce that
\begin{align*}
\rho(y_\tau^{k,v},y)
\leq C_\rho\max\left\{\rho(y_\tau^{k,v},z),\,\rho(z,y)\right\}
\leq C_\rho\max\left\{C_r2^{-(j+k)},\,\rho(z,y)\right\}.
\end{align*}
By this and the fact that
$\rho(y_\tau^{k,v},y)\ge(\max\{C_r2^{-j},\,\widetilde{C}\}+1)C_\rho2^{-k}$,
we obtain
\begin{align*}
\rho(z,y)\ge\frac{\rho(y_\tau^{k,v},y)}{C_\rho}
\ge\left(\max\left\{C_r2^{-j},\,\widetilde{C}\right\}+1\right)2^{-k},
\end{align*}
which, combined with the support condition of $D_k$
[see also Lemma~\ref{Pk}(i)],
further implies that $D_k(z,y)=0$.
From this, \eqref{925}, and the definition of $G_{N,k}(x,y)$ in \eqref{G},
we infer that $G_{N,k}(x,y)=0$.

In summary, we have proved that $G_{N,k}(x,y)=0$
when $\rho(x,y)\ge C'2^{-k}$
and hence we have finished the proof of (i).

Next, we show (ii). Let $x,y,x'\in X$ satisfy
$\rho(x,x')\leq2^{-k}$ or
$\rho(x,x')\leq\frac{2^{-k}+\rho(x,y)}{2C_\rho}$.
Suppose first that $\rho(x,x')\leq\frac{2^{-k}+\rho(x,y)}{2C_\rho}$
and assume, for the moment, that $\rho(x,y)>(C'C_\rho+1)2^{-k}$,
where  $C'\in[1,\infty)$ is the constant as in (i)
of the present lemma.
If $\rho(x,x')<\frac{(C'C_\rho+1)2^{-k-1}}{C_\rho}$, then, by
Definition~\ref{D2.1}(iii), we find that
\begin{align*}
(C'C_\rho+1)2^{-k}&<\rho(x,y)\leq
C_\rho\max\left\{\rho(x,x'),\,\rho(x',y)\right\}
\nonumber\\
&<C_\rho\max\left\{\frac{C'C_\rho+1)2^{-k-1}}{C_\rho},\,\rho(x',y)\right\}
=C_\rho\rho(x',y).
\end{align*}
Hence, $\rho(x',y)>C'2^{-k}$, which, together with the support condition of
the kernels $G_{N,k}$ and the fact that $\rho(x,y)>C'2^{-k}$,
gives $G_{N,k}(x,y)=G_{N,k}(x',y)=0$. Thus, the
desired estimate holds under these assumptions.
If $\rho(x,x')>\frac{(C'C_\rho+1)2^{-k-1}}{C_\rho}$, then a rewriting
of the inequality $\rho(x,x')\leq\frac{2^{-k}+\rho(x,y)}{2C_\rho}$ gives
\begin{align*}
\rho(x,y)\geq 2C_\rho\rho(x,x')-2^{-k}
> (C'C_\rho+1)2^{-k}-2^{-k}=C'C_\rho2^{-k}.
\end{align*}
Hence, $\rho(x,y)>C'C_\rho2^{-k}$ which implies
that $G_{N,k}(x,y)=0$, given that $C_\rho\geq1$.
Going further, this also implies that $2^{-k}<\frac{\rho(x,y)}{C'C_\rho}$ and hence
\begin{align*}
C'C_\rho2^{-k}&<\rho(x,y)
\leq C_\rho\max\left\{\rho(x,x'),\,\rho(x',y)\right\}
\nonumber\\
&<C_\rho\max\left\{\frac{2^{-k}+\rho(x,y)}{2C_\rho},\,\rho(x',y)\right\}
\nonumber\\
&<C_\rho\max\left\{\frac{[(C'C_\rho)^{-1}+1]}{2C_\rho}\rho(x,y),\,\rho(x',y)\right\}
=C_\rho\rho(x',y).
\end{align*}
Therefore, $\rho(x',y)>C'2^{-k}$ from which it follows that $G_{N,k}(x',y)=0$.
Thus,
$$
G_{N,k}(x,y)=G_{N,k}(x',y)=0
$$
and hence the
desired estimate holds under these assumptions. Observe that, if
$\rho(x,x')\leq\frac{2^{-k}+\rho(x,y)}{2C_\rho}$
and $\rho(x,y)\leq(C'C_\rho+1)2^{-k}$, then
$\rho(x,x')\leq\frac{(C'C_\rho+2)2^{-k}}{2C_\rho}$ and hence, to finish
the proof of (ii) of the present lemma, we only need to
consider the case when  $\rho(x,x')\leq C_02^{-k}$
for some fixed $C_0\in[1,\infty)$ which is independent of both $j$ and $k$.

Suppose now that  $\rho(x,x')\leq C_02^{-k}$ for some fixed $C_0\in[1,\infty)$
which is independent of both $j$ and $k$.
We first prove that,
for any given $y\in Q^{k,v}_\tau$
with $\tau\in I_k$ and $v\in\{1,\ldots,N(k,\tau)\}$
and for any $z\in Q^{k,v}_\tau$,
\begin{align}
\label{1225-1}
\left|D^N_k(x,z)-D^N_k(x',z)\right|
\lesssim\left[\frac{\rho(x,x')}{2^{-k}}\right]^\varepsilon
\frac{1}{(V_\rho)_{2^{-k}}(x)}
\end{align}
and
\begin{align}
\label{1225-2}
\left|D_k(z,y)-D_k(y^{k,v}_\tau,y)\right|
\lesssim
\left[\frac{\rho(z,y_\tau^{k,v})}{2^{-k}}\right]^\varepsilon
\frac{1}{(V_\rho)_{2^{-k}}(y)}\mathbf{1}_{B_\rho(y,C2^{-k})}(z).
\end{align}
Indeed, if $\rho(x,z)>\max\{C'C_\rho2^{-k},\,C'C_\rho C_02^{-k}\}$,
then, from Definition~\ref{D2.1}(iii), it follows that
\begin{align*}
C'C_\rho2^{-k}<\rho(x,z)\leq
C_\rho\max\left\{\rho(x,x'),\,\rho(x',z)\right\}=C_\rho\rho(x',z),
\end{align*}
which, combined with the support condition of $D_k^N$,
further implies
\begin{align*}
D^N_k(x,z)=0=D^N_k(x',z)
\end{align*}
and hence
$$\left|D^N_k(x,z)-D^N_k(x',z)\right|=0.$$
On the other hand, if $\rho(x,z)\leq
\max\{C'C_\rho2^{-k},\,C'C_\rho C_02^{-k}\}$,
then, by Lemma~\ref{A3.2}(vii), we conclude that
$(V_\rho)_{2^{-k}}(x)\sim (V_\rho)_{2^{-k}}(z)$.
From this, Lemma~\ref{Pk}, and Remark~\ref{RemThmATI}(i),
we deduce that
\begin{align*}
\left|D^N_k(x,z)-D^N_k(x',z)\right|
\lesssim\left[\frac{\rho(x,x')}{2^{-k}}\right]^\varepsilon
\frac{1}{(V_\rho)_{2^{-k}}(z)}
\sim\left[\frac{\rho(x,x')}{2^{-k}}\right]^\varepsilon
\frac{1}{(V_\rho)_{2^{-k}}(x)},
\end{align*}
which completes the proof of \eqref{1225-1}.
To show \eqref{1225-2}, by an argument similar to that used
in the estimation of \eqref{1225-1},
we find that, if $\rho(z,y)\gtrsim2^{-k}$
with the implicit positive constant being large enough and
independent of both $j$ and $k$,
then $|D_k(z,y)-D_k(y^{k,v}_\tau,y)|=0$;
if $\rho(z,y)\lesssim2^{-k}$ with the same positive constant, then
Lemma~\ref{Pk} and Remark~\ref{RemThmATI}(i)
imply that \eqref{1225-2} holds.
From \eqref{1225-1} and \eqref{1225-2},
we infer that
\begin{align*}
&\left|G_{N,k}(x,y)-G_{N,k}(x',y)\right|
\\
&\quad\leq\sum_{\tau\in I_k}\sum_{v=1}^{N(k,\tau)}
\int_{Q^{k,v}_\tau}\left|D^N_k(x,z)-D^N_k(x',z)\right|
\left|D_k(z,y)-D_k(y^{k,v}_\tau,y)\right|
\,d\mu(z)
\\
&\quad\lesssim\sum_{\tau\in I_k}\sum_{v=1}^{N(k,\tau)}
\int_{Q^{k,v}_\tau}
\left[\frac{\rho(x,x')}{2^{-k}}\right]^\varepsilon
\frac{1}{(V_\rho)_{2^{-k}}(z)}
\left[\frac{\rho(z,y_\tau^{k,v})}{2^{-k}}\right]^\varepsilon
\\
&\qquad\times
\frac{1}{(V_\rho)_{2^{-k}}(y)}
\mathbf{1}_{B_\rho(y,C2^{-k})}(z)\,d\mu(z),
\end{align*}
which, together with \eqref{rhozy}
and \eqref{upperdoub},
further implies that
\begin{align*}
&\left|G_{N,k}(x,y)-G_{N,k}(x',y)\right|
\\
&\quad\sim\int_X
\left[\frac{\rho(x,x')}{2^{-k}}\right]^\varepsilon
\frac{1}{(V_\rho)_{2^{-k}}(x)}
\left[\frac{2^{-j}2^{-k}}{2^{-k}}\right]^\varepsilon
\frac{1}{(V_\rho)_{2^{-k}}(y)}
\mathbf{1}_{B_\rho(y,C2^{-k})}(z)
\,d\mu(z)
\\
&\quad\lesssim 2^{(k-j)\varepsilon}
[\rho(x,x')]^\varepsilon\frac{1}{(V_\rho)_{2^{-k}}(x)}.
\end{align*}
This proves that $G_{N,k}(x,y)$ satisfies (ii).

Finally, we show (iii).
By an argument similar to that used in the proof of (ii),
it suffices to prove (iii) in the case
when $\rho(x,x')\leq C_12^{-k}$ and
$\rho(y,y')\leq C_12^{-k}$
for some fixed $C_1\in[1,\infty)$ which is independent of $j$ and $k$.
In this case, from Lemma~\ref{Pk}(iii),
the regularity condition of $D^N_k$ [see Lemma~\ref{Pk}(ii)],
\eqref{rhozy}, and the support conditions of both $D_k$ and $D_N^k$
[see also Lemma~\ref{Pk}(i)], we deduce that
\begin{align*}
&\left|[G_{N,k}(x,y)-G_{N,k}(x',y)]-[G_{N,k}(x,y')-G_{N,k}(x',y')]\right|
\\
&\quad\leq\sum_{\tau\in I_k}\sum_{v=1}^{N(k,\tau)}
\int_{Q^{k,v}_\tau}\left|D^N_k(x,z)-D^N_k(x',z)\right|
\\
&\qquad\times
\left|[D_k(z,y)-D_k(y^{k,v}_\tau,y)]
-[D_k(z,y')-D_k(y^{k,v}_\tau,y')]\right|
\,d\mu(z)
\\
&\quad\lesssim
\sum_{\tau\in I_k}\sum_{v=1}^{N(k,\tau)}
\int_{Q^{k,v}_\tau}\left[\frac{\rho(x,x')}{2^{-k}}\right]^\varepsilon
\frac{1}{(V_\rho)_{2^{-k}}(x)}
\left[\frac{\rho(z,y_\tau^{k,v})}{2^{-k}}\right]^\varepsilon
\left[\frac{\rho(y,y')}{2^{-k}}\right]^\varepsilon
\\
&\qquad\times
\left[\frac{1}{(V_\rho)_{2^{-k}}(y)}+\frac{1}{(V_\rho)_{2^{-k}}(y')}\right]
\mathbf{1}_{B_\rho(y,C2^{-k})}(z)\,d\mu(z),
\end{align*}
which, combined with
an argument similar to that used in the proofs of both \eqref{1225-1}
and \eqref{1225-2} and \eqref{upperdoub},
further implies that
\begin{align*}
&\left|[G_{N,k}(x,y)-G_{N,k}(x',y)]-[G_{N,k}(x,y')-G_{N,k}(x',y')]\right|
\\
&\quad\lesssim
\int_X\left[\frac{\rho(x,x')}{2^{-k}}\right]^\varepsilon
\frac{1}{(V_\rho)_{2^{-k}}(x)}
\left[\frac{2^{-(j+k)}}{2^{-k}}\right]^\varepsilon
\left[\frac{\rho(y,y')}{2^{-k}}\right]^\varepsilon
\\
&\qquad\times
\left[\frac{1}{(V_\rho)_{2^{-k}}(y)}+\frac{1}{(V_\rho)_{2^{-k}}(y')}\right]
\mathbf{1}_{B_\rho(y,C2^{-k})}(z)\,d\mu(z)
\\
&\quad\lesssim2^{(2k-j)\varepsilon}
\left[\rho(x,x')\rho(y,y')\right]^{\varepsilon}
\frac{1}{(V_\rho)_{2^{-k}}(x)}.
\end{align*}
This shows that $G_{N,k}(x,y)$ satisfies (iii),
which completes the proof of Lemma~\ref{A3.41}.
\end{proof}

To establish the boundedness of $\mathfrak{G}_N$ on spaces of test functions,
we use the method similar to that used in
the proof of the boundedness of $\mathcal{R}_N$ in
Proposition~\ref{A3.16}. With this in mind, for any $N,M\in\mathbb{N}$,
$G^{(M)}_N$, the kernel of the integral
operator $\mathfrak{G}^{(M)}_N$,
is defined by setting, for any $x,y\in X$,
\begin{align}
\label{GM}
G^{(M)}_N(x,y)
&:=\sum_{|k|\leq M}G_{N,k}(x,y)
\nonumber\\
&=\sum_{|k|\leq M}
\sum_{\tau\in I_k}
\sum_{v=1}^{N(k,\tau)}
\int_{Q^{k,v}_\tau}D^N_k(x,z)
\left[D_k(z,y)-D_k(y^{k,v}_\tau,y)\right]\,d\mu(z).
\end{align}

\begin{lemma}\label{GL2}
Let $(X,\mathbf{q},\mu)$ be an ultra-RD-space,
where $\mu$ is assumed to be a Borel-semiregular measure on $X$.
Assume that $0<\varepsilon\preceq{\mathrm{ind\,}}(X,\mathbf{q})$
with ${\mathrm{ind\,}}(X,\mathbf{q})$ as in \eqref{index}. Let $x_1\in X$,
$r\in(0,\infty)$, and $\beta,\gamma\in(0,\varepsilon)$.
For any $N,M\in\mathbb{N}$,
let $\mathfrak{G}_N$ (with kernel $G_N$)
and $\mathfrak{G}^{(M)}_N$ (with kernel $G^{(M)}_N$) be
the same as, respectively, in \eqref{G} and \eqref{GM}. Then
there exists a positive constant $C$, depending on $N$,
but independent of $M$, $j$, $x_1$, $r$, and $y_\tau^{k,v}$, such that
\begin{align}\label{G1}
\left\|\mathfrak{G}_N\right\|_{L^2(X)\to L^2(X)}
+\left\|\mathfrak{G}^{(M)}_N\right\|_{L^2(X)\to L^2(X)}
\leq C2^{-j\varepsilon}
\end{align}
and
\begin{align}\label{G2}
\left\|\mathfrak{G}^{(M)}_N\right\|_{\mathring{\mathcal{G}}(x_1,r,\beta,\gamma)
\to\mathring{\mathcal{G}}(x_1,r,\beta,\gamma)}\leq
C2^{-j\varepsilon},
\end{align}
where $j$ is the same as in \eqref{2cubej} and
$\{y_\tau^{k,v}\}_{k\in\mathbb{Z},\tau\in I_k,v\in\{1,\ldots,N(k,\tau)\}}$
the same as in \eqref{DRS}.
\end{lemma}

\begin{proof}
Let $\rho\in\mathbf{q}$ be a regularized quasi-ultrametric
as in Definition~\ref{D2.3}.
We first prove \eqref{G1}.
Using (i) and (ii) of Lemma~\ref{A3.41}
and following the proof of \eqref{3.2} in Lemma~\ref{A3.12},
we find that, for any $k,l\in\mathbb{Z}$
and $x,y\in X$ satisfying $\rho(x,y)\lesssim 2^{-k}$,
$$
\left|G_{N,k}^*G_{N,l}(x,y)\right|+\left|G_{N,k}G_{N,l}^*(x,y)\right|
\lesssim2^{-j\varepsilon}2^{-|k-l|\varepsilon}
\frac{1}{(V_\rho)_{2^{-(k\land l)}}(x)},
$$
where $G_{N,k}^*G_{N,l}$ and $G_{N,k}G_{N,l}^*$ denote the kernels of operators
$\mathfrak{G}_{N,k}^*\mathfrak{G}_{N,l}$ and
$\mathfrak{G}_{N,k}\mathfrak{G}_{N,l}^*$
and where
$\mathfrak{G}_{N,k}^*$ and $\mathfrak{G}_{N,l}^*$ denote the adjoint operators of
$\mathfrak{G}_{N,k}$ and $\mathfrak{G}_{N,l}$.
Otherwise,
$$
\left|G_{N,k}^*G_{N,l}(x,y)\right|
+\left|G_{N,k}G_{N,l}^*(x,y)\right|=0.
$$
From this and an argument similar to that used in the proof of \eqref{RN0},
we infer that, for any $k,l\in\mathbb{Z}$,
$$
\left\|\mathfrak{G}_{N,k}^*\mathfrak{G}_{N,l}\right\|_{L^2(X)\to L^2(X)}
+\left\|\mathfrak{G}_{N,k}\mathfrak{G}_{N,l}^*\right\|_{L^2(X)\to L^2(X)}
\lesssim2^{-j\varepsilon}2^{-|k-l|\varepsilon},
$$
which, together with Lemma~\ref{CotlarStein},
further implies that \eqref{G1} holds.

To show \eqref{G2}, by the definition of $\mathfrak{G}^{(M)}_N$ and Lemma~\ref{A3.41}(i),
for any $f\in C^\beta(X)$, we conclude that
$\mathfrak{G}^{(M)}_N$ satisfies Theorem~\ref{CZOonGO}(ii).
Next, we prove that $\mathfrak{G}^{(M)}_N$ satisfies \eqref{CZk1}, \eqref{CZk2}, and
\eqref{CZO1}.

For \eqref{CZk1}, from Lemma~\ref{A3.41}(i), Lemma~\ref{A3.14},
and an argument similar to that used in the proof of Lemma~\ref{A3.15},
we deduce that,
for any $x,y\in X$ with $x\neq y$,
$$
\left|G^{(M)}_N(x,y)\right|
\leq\sum_{|k|\leq M}\left|G_{N,k}(x,y)\right|
\lesssim2^{-j\varepsilon}
\frac{1}{V_\rho(x,y)}.
$$
Similarly, by both (ii) and (iii) of Lemma~\ref{A3.41} and Lemma~\ref{A3.14},
we obtain \eqref{CZk2} and
\eqref{CZO1}. From \eqref{G1}, we infer that
$\mathfrak{G}^{(M)}_N$ is bounded on $L^2(X)$
with its operator norm bounded by a positive constant independent of $M$. Applying the proof of
Lemma~\ref{A3.41}(iv) to $\mathfrak{G}^{(M)}_N$, we obtain
$G^{(M)}_N$, the kernel of $\mathfrak{G}^{(M)}_N$,
satisfies the cancellation condition.
By this and Theorem~\ref{CZOonGO}, we obtain \eqref{G2}.
This finishes the proof of Lemma~\ref{GL2}.
\end{proof}

By Lemma~\ref{A3.41} and an argument similar to that used in
the proof of Lemma~\ref{RNRNM}, we obtain the
following lemma; we omit the details of its proof.

\begin{lemma}\label{GMG}
Let $(X,\mathbf{q},\mu)$ be a $(\kappa,n)$-ultra-RD-space with
$0<\kappa\leq n<\infty$,
where $\mu$ is assumed to be a Borel-semiregular measure on $X$.
Assume that $0<\varepsilon\preceq{\mathrm{ind\,}}(X,\mathbf{q})$
with ${\mathrm{ind\,}}(X,\mathbf{q})$  as in \eqref{index}.
Let $x_1\in X$, $r\in(0,\infty)$,
and $\beta,\gamma\in(0,\varepsilon)$.
Fix $N\in\mathbb{N}$ and, for any $M\in\mathbb{N}$, let
$\mathfrak{G}^{(M)}_N$ and $\mathfrak{G}_N$ be
the same as, respectively, in \eqref{GM} and \eqref{G}.
Then the following statements hold:
\begin{enumerate}
\item[\textup{(i)}]
For any $f\in{\mathcal G}(x_1,r,\beta,\gamma)$,
\begin{align*}
\lim_{M\to\infty}\left\|\mathfrak{G}^{(M)}_Nf-\mathfrak{G}_Nf\right\|_{L^2(X)}=0.
\end{align*}

\item[\textup{(ii)}]
For any
$f\in{\mathcal G}(x_1,r,\beta,\gamma)$, the sequence
$\{\mathfrak{G}^{(M)}_Nf\}_{M=N+1}^\infty$
converges locally uniformly to some element,
denoted by $\widetilde{\mathfrak{G}}_Nf$,
where $\widetilde{\mathfrak{G}}_Nf$ differs
from $\mathfrak{G}_Nf$ at most on a set of
$\mu$-measure 0.

\item[\textup{(iii)}]
The operator $\widetilde{\mathfrak{G}}_N$
can be uniquely extended from
${\mathcal G}(x_1,r,\beta,\gamma)$ to $L^2(X)$,
where the extension operator coincides with
$\mathfrak{G}_N$ in the sense that, for any
$f\in{\mathcal G}(x_1,r,\beta,\gamma)$ and $\mu$-almost every
$x\in X$,
$$
\lim_{M\to\infty}\mathfrak{G}^{(M)}_Nf(x)
=\widetilde{\mathfrak{G}}_Nf(x)=\mathfrak{G}_Nf(x).
$$
\end{enumerate}
\end{lemma}

By Lemmas~\ref{GL2} and~\ref{GMG},
Definition~\ref{testfunction},
and the Lebesgue dominated convergence theorem,
we obtain the following result;
we omit the details of its proof.

\begin{proposition}\label{A3.44}
Let $(X,\mathbf{q},\mu)$ be an ultra-RD-space,
where $\mu$ is assumed to be a Borel-semiregular measure on $X$.
Assume that $0<\varepsilon\preceq{\mathrm{ind\,}}(X,\mathbf{q})$
with ${\mathrm{ind\,}}(X,\mathbf{q})$ as in \eqref{index}.
Let $x_1\in X$, $r\in(0,\infty)$, $\beta,\gamma\in(0,\varepsilon)$,
and $\mathfrak{G}_N$ and $\mathfrak{R}_N$ be
the same as, respectively, in \eqref{G} and \eqref{RGRN}.
Then, for any given $N\in\mathbb{N}$, there exists a positive constant $C_N$,
depending on $N$, such that
\begin{align}\label{1646124}
\left\|\mathfrak{G}_N\right\|_{L^2(X)\to L^2(X)}+
\left\|\mathfrak{G}_N\right\|_{\mathring{\mathcal{G}}(x_1,r,\beta,\gamma)
\to\mathring{\mathcal{G}}(x_1,r,\beta,\gamma)}
\leq C_N2^{-j\varepsilon}.
\end{align}
Consequently, for any
$\varepsilon'\in(\max\{\beta,\,\gamma\},\varepsilon)$,
there exists $\widetilde{C}\in(0,\infty)$,
independent of $N$,
such that
$$\left\|\mathfrak{R}_N\right\|_{L^2(X)\to L^2(X)}
\leq\widetilde{C}2^{-N\varepsilon'}+C_N2^{-j\varepsilon}$$
and
$$\left\|\mathfrak{R}_N\right\|_{\mathring{\mathcal{G}}(x_1,r,\beta,\gamma)
\to\mathring{\mathcal{G}}(x_1,r,\beta,\gamma)}\leq
\widetilde{C}2^{N(\varepsilon'-\varepsilon)}+C_N2^{-j\varepsilon}.$$
\end{proposition}

By Proposition~\ref{A3.44}
and repeating the proof of
Proposition~\ref{A3.16},
we obtain the boundedness of $\mathfrak{S}_N^{-1}$
on $\mathring{{\mathcal G}}(x_1,r,\beta,\gamma)$;
we omit the details of its proof.

\begin{corollary}\label{A3.45}
Let $(X,\mathbf{q},\mu)$ be a $(\kappa,n)$-ultra-RD-space,
where $\mu$ is assumed to be a Borel-semiregular measure on $X$ satisfying both
the doubling condition \eqref{doublingcond}
and the reverse doubling condition \eqref{RDcondition}
with respect to $\rho\in\mathbf{q}$.
Assume that $0<\varepsilon\preceq{\mathrm{ind\,}}(X,\mathbf{q})$
with ${\mathrm{ind\,}}(X,\mathbf{q})$ as in \eqref{index}.
Let $x_1\in X$, $r\in(0,\infty)$, $\beta,\gamma\in(0,\varepsilon)$, and
$\varepsilon'\in(\max\{\beta,\,\gamma\},\varepsilon)$.
Let $\mathfrak{S}_N$ be the same as in \eqref{DRS},
where $N\in\mathbb{N}$ is fixed, and
let $j\in\mathbb{N}$ be sufficiently large such that
\begin{align}
\label{1838}
\max\left\{\widetilde{C}2^{-N\varepsilon'}+C_N2^{-j\varepsilon},\,
\widetilde{C}2^{-N\varepsilon'}+C_N2^{-j\varepsilon}
\right\}\leq\frac{1}{2},
\end{align}
where $C_N$ and $\widetilde{C}$ are the same as in
Proposition~\ref{A3.44}.
Then $\mathfrak{S}_N$ has a well-defined and bounded inverse
on $\mathring{{\mathcal G}}(x_1,r,\beta,\gamma)$.
That is, there exists a constant $C_1\in(0,\infty)$,
depending only on $\beta$, $\gamma$, and $\rho$,
such that, for any
$f\in\mathring {{\mathcal G}}(x_1,r,\beta,\gamma)$,
\begin{align*}
\left\|\mathfrak{S}_N^{-1}(f)\right\|_{\mathring{{\mathcal G}}(x_1,r,\beta,\gamma)}
\leq C_1\|f\|_{\mathring{{\mathcal G}}(x_1,r,\beta,\gamma)}.
\end{align*}
\end{corollary}

To establish sharp homogeneous discrete Calder\'on
reproducing formulae,
we still need the following technical lemma.

\begin{lemma}\label{A3.46}
Let $(X,\mathbf{q},\mu)$ be an ultra-RD-space,
where $\mu$ is assumed to be a Borel-semiregular measure on $X$.
Assume that $0<\varepsilon\preceq{\mathrm{ind\,}}(X,\mathbf{q})$
with ${\mathrm{ind\,}}(X,\mathbf{q})$ as in \eqref{index}, and
$\{\mathcal{S}_{2^{-k}}\}_{k\in\mathbb{Z}}$
is a homogeneous approximation
of the identity of order $\varepsilon$
in Definition~\ref{DefHATI}.
Fix $N\in\mathbb{N}$.
For any $k\in\mathbb{Z}$,
define
$\mathcal{D}_k:=\mathcal{S}_{2^{-k}}
-\mathcal{S}_{2^{-k+1}}$
and
let $\mathcal{D}_k^N$ be the same as in \eqref{3.6}.
For any fixed $\{y^{k,v}_\tau:y^{k,v}_\tau\in Q^{k,v}_\tau\}_{k\in\mathbb{Z},\,
\tau\in I_k,\,v\in\{1,\ldots,N(k,\tau)\}}$ in $X$
and for any $k\in\mathbb{Z}$ and $x,y\in X$, let
\begin{align}
\label{926}
E_k(x,y)
:=\sum_{\tau\in I_k}\sum_{v=1}^{N(k,\tau)}
\int_{Q^{k,v}_\tau}
D^N_k(x,z)\,d\mu(z)
D_k\left(y^{k,v}_\tau,y\right).
\end{align}
Then, for any given $\rho\in\mathbf{q}$
having the property that all $\rho$-balls are $\mu$-measurable,
there exist two positive constants $C$ and $\widetilde{C}$,
depending on $N$, but independent of $k$, $j$, and $y_\tau^{k,v}$,
such that
\begin{enumerate}
\item[\textup{(i)}]
For any $x,y\in X$ with $\rho(x,y)\leq C2^{-k}$,
$$
\left|E_k(x,y)\right|
\leq \widetilde{C}\frac{1}{(V_\rho)_{2^{-k}}(x)},
$$
and, for any $x,y\in X$ with $\rho(x,y)>C2^{-k}$, $E_k(x,y)=0$.

\item[\textup{(ii)}]
For any $x,y,y'\in X$ with $\rho(y,y')\leq 2^{-k}$
or $\rho(y,y')\leq\frac{2^{-k}+\rho(x,y)}{2C_\rho}$,
$$
\left|E_k(x,y)-E_k(x,y')\right|
\leq\widetilde{C}2^{k\varepsilon}\left[\rho(y,y')\right]^\varepsilon
\frac{1}{(V_\rho)_{2^{-k}}(x)}.
$$

\item[\textup{(iii)}]
For any $x,x',y,y'\in X$ with
$\rho(y,y')\leq 2^{-k}$ and $\rho(x,x')\leq 2^{-k}$, or
$\rho(y,y')\leq\frac{2^{-k}+\rho(x,y)}{2C_\rho}$ and
$\rho(x,x')\leq\frac{2^{-k}+\rho(x,y)}{2C_\rho}$,
\begin{align*}
&\left|\left[E_k(x,y)-E_k(x,y')\right]
-\left[E_k(x,y')-E_k(x',y')\right]\right|\\
&\quad\leq\widetilde{C}2^{2k\varepsilon}[\rho(x,x')\rho(y,y')]^\varepsilon
\frac{1}{(V_\rho)_{2^{-k}}(x)}.
\end{align*}

\item[\textup{(iv)}]
For any $x\in X$,
$$
\int_XE_k(x,z)\,d\mu(z)=\int_XE_k(z,x)\,d\mu(z)=0.
$$
\end{enumerate}
\end{lemma}

\begin{proof}
Without loss of generality, we may assume that
$\rho\in\mathbf{q}$ is a regularized quasi-ultrametric as in Definition~\ref{D2.3}.
Fix $k\in\mathbb{Z}$.
By Lemma~\ref{Pk}(v), we find that, for any $x\in X$,
$$
\int_XD_k(x,y)\,d\mu(y)=0
=\int_XD^N_k(y,x)\,d\mu(y),
$$
which, combined with the definition of $E_k$ and Fubini's theorem,
further implies that (iv) holds.	

Now, we show (i). On the one hand, let $x,y\in X$ satisfy
$\rho(x,y)\leq\widetilde{C}2^{-k}$ for some $\widetilde{C}\in(1,\infty)$.
From Lemma~\ref{Pk}(i), we deduce that
there exists a constant $C_1\in(1,\infty)$
such that, for any $x,y\in X$ satisfying $\rho(x,y)\ge C_12^{-k}$,
$$
D_k(x,y)=D_k^N(x,y)=0.
$$
By this, the size conditions of $D_k$ and $D_k^N$
in Lemma~\ref{Pk}(i),
\eqref{upperdoub}, Lemma~\ref{DyadicSys}(i), \eqref{un-3iu},
and Lemma~\ref{A3.2}(vi), we conclude that
\begin{align}
\label{1102}
|E_k(x,y)|
&\lesssim\int_X
\left|D^N_k(x,z)\right|
\frac{1}{(V_\rho)_{2^{-k}}(y)}\,d\mu(z)
\nonumber\\
&\lesssim\int_{\rho(x,z)\leq C_12^{-k}}
\frac{1}{(V_\rho)_{2^{-k}}(x)}\,d\mu(z)
\frac{1}{(V_\rho)_{2^{-k}}(y)}
\nonumber\\
&\lesssim\frac{1}{(V_\rho)_{2^{-k}}(y)}
\sim\frac{1}{(V_\rho)_{2^{-k}}(x)}.
\end{align}
On the other hand, let $x,y\in X$ satisfy $\widetilde{C}2^{-k}\leq\rho(x,y)$.
From Lemma~\ref{DyadicSys}(iv),
it follows that $\mathrm{diam}_\rho(Q_{\tau}^{k,v})\sim 2^{-k}$,
which, together with Definition~\ref{D2.1}(iii), further implies that,
for any $z,y_\tau^{k,v}\in Q_{\tau}^{k,v}$,
\begin{align}
\label{1053}
\widetilde{C}2^{-k}
&\leq\rho(x,y)\leq C_\rho\max\left\{\rho(x,z),\,\rho(z,y)\right\}
\nonumber\\
&\lesssim\max\left\{\rho(x,z),\,\rho(z,y_\tau^{k,v}),\,
\rho(y_\tau^{k,v},y)\right\}
\nonumber\\
&\lesssim\max\left\{\rho(x,z),\,2^{-k},\,
\rho(y_\tau^{k,v},y)\right\}.
\end{align}
If we choose $\widetilde{C}$ large enough, then, by \eqref{1053},
the support conditions of $D_k$ and $D_k^N$,
and an argument similar to that used in the proof of
Lemma~\ref{A3.41}(i),
it is easy to verify that $E_k(x,y)=0$,
which, combined with \eqref{1102},
completes the proof of (i).

To prove (ii), from an argument similar to
that used in the proof of Lemma~\ref{A3.41}(ii), we infer that
it suffices to show (ii) in the case $\rho(y,y')\leq C_02^{-k}$
for some $C_0\in(0,\infty)$.
By the size condition of $D_k^N$,
the regularity condition of $D_k$,  and the support conditions
of both $D_k^N$ and $D_k$,
we find that, for any $x,y,y'\in X$
with $\rho(y,y')\leq C_02^{-k}$,
\begin{align*}
&\left|E_k(x,y)-E_k(x,y')\right|
\\
&\quad\leq
\sum_{\tau\in I_k}\sum_{v=1}^{N(k,\tau)}
\int_{Q^{k,v}_\tau}
\left|D^N_k(x,z)\right|\,d\mu(z)
\left|D_k\left(y^{k,v}_\tau,y\right)-D_k\left(y^{k,v}_\tau,y'\right)\right|
\\
&\quad\lesssim
\sum_{\tau\in I_k}\sum_{v=1}^{N(k,\tau)}
\int_{Q^{k,v}_\tau\cap
\{z:\rho(z,x)\leq C_12^{-k}\text{ and }\rho(z,y)\leq C_12^{-k}\}}
\frac{1}{(V_\rho)_{2^{-k}}(x)}\,d\mu(z)
\\
&\qquad\times\left[\frac{\rho(y,y')}{2^{-k}}\right]^\varepsilon
\frac{1}{(V_\rho)_{2^{-k}}(y^{k,v}_\tau)},
\end{align*}
which, together with Lemma~\ref{DyadicSys}(i) and \eqref{un-3iu},
further implies that
\begin{align*}
&\left|E_k(x,y)-E_k(x,y')\right|
\\
&\quad\leq
\int_{\{z:\rho(z,y)\leq C_12^{-k}\}}
\frac{1}{(V_\rho)_{2^{-k}}(x)}\,d\mu(z)
\left[\frac{\rho(y,y')}{2^{-k}}\right]^\varepsilon
\frac{1}{(V_\rho)_{2^{-k}}(y)}
\\
&\quad\lesssim
\left[\frac{\rho(y,y')}{2^{-k}}\right]^\varepsilon
\frac{1}{(V_\rho)_{2^{-k}}(x)}.
\end{align*}

To prove (iii), from an argument similar
to that used in the proof of (ii)
in the present lemma, we deduce that
it suffices to show it in the case
that $\rho(x,x')\leq C_22^{-k}$
for some fixed $C_2\in(0,\infty)$.
In this case, Lemma~\ref{A3.2}(vii) implies that,
for such $x'$ and $x$,
\begin{align}
\label{1121}
(V_\rho)_{2^{-k}}(x)\sim (V_\rho)_{2^{-k}}(x').
\end{align}
By the support condition of $D_k^N$
and Definition~\ref{D2.1}(iii), we conclude that
\begin{align}
\label{1429}
&{\mathrm{supp\,}} D_k^N(x,\cdot)\cup
{\mathrm{supp\,}} D_k^N(x',\cdot)\nonumber\\
&\quad\subset
\left[B_\rho(x,C_12^{-k})\cup B_\rho(x',C_12^{-k})\right]
\subset B_\rho(x,C_1C_\rho2^{-k}).
\end{align}
Note that, from Lemma~\ref{DyadicSys}(iv), we infer that,
for any $z\in B_\rho(x,C_1C_\rho2^{-k})
\cap Q_\tau^{k,v}$ and $y_\tau^{k,v}\in Q_\tau^{k,v}$,
$\rho(z,y_\tau^{k,v})\lesssim 2^{-k}$
and
$\rho(z,x)\lesssim 2^{-k}$,
which, combined with Lemma~\ref{A3.2}(vii), further implies that
$V_{2^{-k}}(y_\tau^{k,v})\sim
V_{2^{-k}}(z)\sim
V_{2^{-k}}(x)$.
By this, \eqref{1429}, and both
(i) and (ii) of Lemma~\ref{Pk},
we find that
\begin{align*}
&\left|\left[E_k(x,y)-E_k(x,y')\right]
-\left[E_k(x,y')-E_k(x',y')\right]\right|
\\
&\quad\leq
\sum_{\tau\in I_k}\sum_{v=1}^{N(k,\tau)}
\int_{Q^{k,v}_\tau}
\left|D^N_k(x,z)-D^N_k(x',z)\right|\,d\mu(z)\\
&\qquad\times\left|D_k\left(y^{k,v}_\tau,y\right)
-D_k\left(y^{k,v}_\tau,y'\right)\right|
\\
&\quad\lesssim
\int_{B_\rho(x,C_1C_\rho2^{-k})}
\left[\frac{\rho(x,x')}{2^{-k}}\right]^\varepsilon
\left[\frac{1}{(V_\rho)_{2^{-k}}(x)}+\frac{1}{(V_\rho)_{2^{-k}}(x')}\right]
\,d\mu(z)\\
&\qquad\times\left[\frac{\rho(y,y')}{2^{-k}}\right]^\varepsilon
\left[\frac{1}{(V_\rho)_{2^{-k}}(x)}\right]\\
&\quad\lesssim
\left[\frac{\rho(y,y')}{2^{-k}}\right]^\varepsilon
\left[\frac{\rho(x,x')}{2^{-k}}\right]^\varepsilon
\frac{1}{(V_\rho)_{2^{-k}}(x)}
\quad\text{by \eqref{1121}
and \eqref{upperdoub}}.
\end{align*}
This finishes the proof of (iii)
and hence Lemma~\ref{A3.46}.
\end{proof}

Next, we establish the following sharp homogeneous
discrete Calder\'on reproducing formulae.
\index[words]{homogeneous
discrete Calder\'on reproducing formulae}

\begin{theorem}\label{A3.47}
Let $(X,\mathbf{q},\mu)$ be an ultra-RD-space
with upper dimension $n\in(0,\infty)$,
where $\mu$ is assumed to be a Borel-semiregular measure on $X$.
Assume that $0<\varepsilon\preceq{\mathrm{ind\,}}(X,\mathbf{q})$
with ${\mathrm{ind\,}}(X,\mathbf{q})$ as in \eqref{index}
and $\{\mathcal{S}_{2^{-k}}\}_{k\in\mathbb{Z}}$
is a homogeneous approximation
of the identity of order $\varepsilon$ as in Definition~\ref{DefHATI} with kernels
$\{S_{2^{-k}}\}_{k\in\mathbb{Z}}$.
For any $k\in\mathbb{Z}$, let
$\mathcal{D}_k:=\mathcal{S}_{2^{-k}}
-\mathcal{S}_{2^{-k+1}}$.
Fix
$y^{k,v}_\tau\in Q^{k,v}_\tau$
for any given $k\in\mathbb{Z}$, $\tau\in I_k$, and
$v\in\{1,\ldots,N(k,\tau)\}$.
Then, for any fixed $j,N\in\mathbb{N}$ satisfying \eqref{1838}
with $\beta,\gamma\in(0,\varepsilon)$
and $\varepsilon'\in(\max\{\beta,\,\gamma\},\varepsilon)$,
there exist
$\{\widetilde{\mathcal{D}}^{(i)}_k\}_{k\in\mathbb{Z}}$ and
$\{\overline{\mathcal{D}}^{(i)}_k\}_{k\in\mathbb{Z}}$ for $i\in\{1,\,2,\,3\}$
of bounded linear operators on $L^2(X)$,
with kernels
$\{\widetilde{D}^{(i)}_k\}_{k\in\mathbb{Z}}$ and
$\{\overline{D}^{(i)}_k\}_{k\in\mathbb{Z}}$,
such that, for
any $f\in\mathring{{\mathcal G}}_0^\varepsilon(\beta,\gamma)$
[resp. $L^p(X)$ with $p\in(1,\infty)$],
\begin{align}
\label{3.29}
f(\cdot)
&=\sum_{k=-\infty}^{\infty}\sum_{\tau\in I_k}
\sum_{v=1}^{N(k,\tau)}
\int_{Q^{k,v}_\tau}\widetilde{D}^{(1)}_k(\,\cdot\,,y)\,d\mu(y)
\mathcal{D}_kf(y^{k,v}_\tau)
\nonumber\\
&=\sum_{k=-\infty}^{\infty}\sum_{\tau\in I_k}
\sum_{v=1}^{N(k,\tau)}
\widetilde{D}^{(2)}_k\left(\,\cdot\,,y^{k,v}_\tau\right)
\int_{Q^{k,v}_\tau}\mathcal{D}_kf(y)\,d\mu(y)
\nonumber\\
&=\sum_{k=-\infty}^{\infty}\sum_{\tau\in I_k}
\sum_{v=1}^{N(k,\tau)}
\mu\left(Q^{k,v}_\tau\right)
\widetilde{D}^{(3)}_k\left(\,\cdot\,,y^{k,v}_\tau\right)
\mathcal{D}_kf(y^{k,v}_\tau)
\nonumber\\
&=\sum_{k=-\infty}^{\infty}\sum_{\tau\in I_k}
\sum_{v=1}^{N(k,\tau)}
\int_{Q^{k,v}_\tau}D_k(\,\cdot\,,y)\,d\mu(y)
\overline{\mathcal{D}}^{(1)}_kf(y^{k,v}_\tau)
\nonumber\\
&=\sum_{k=-\infty}^{\infty}\sum_{\tau\in I_k}
\sum_{v=1}^{N(k,\tau)}
D_k\left(\,\cdot\,,y^{k,v}_\tau\right)
\int_{Q^{k,v}_\tau}\overline{\mathcal{D}}^{(2)}_kf(y)\,d\mu(y)
\nonumber\\
&=\sum_{k=-\infty}^{\infty}\sum_{\tau\in I_k}
\sum_{v=1}^{N(k,\tau)}
\mu\left(Q^{k,v}_\tau\right)
D_k\left(\,\cdot\,,y^{k,v}_\tau\right)
\overline{\mathcal{D}}^{(3)}_kf(y^{k,v}_\tau),
\end{align}
where the series converges in
$\mathring{{\mathcal G}}_0^\varepsilon(\beta,\gamma)$
[resp. $L^p(X)$ with $p\in(1,\infty)$]
and also converges pointwise $\mu$-almost everywhere on $X$.
Moreover, for any $i\in\{1,\,2,\,3\}$, $k\in\mathbb{Z}$, $x\in X$,
and $\widetilde{\beta},\widetilde{\gamma}\in(0,\varepsilon)$,
\begin{align*}
\widetilde{D}^{(i)}_k(\cdot,x),\,
\overline{D}^{(i)}_k(x,\cdot)\in{\mathcal G}
(x,2^{-k},\widetilde{\beta},\widetilde{\gamma}),
\end{align*}
\begin{align*}
\int_X\widetilde{D}^{(i)}_k(x,z)\,d\mu(z)=0
=\int_X\widetilde{D}^{(i)}_k(z,x)\,d\mu(z),
\end{align*}
and
\begin{align*}
\int_X\overline{D}^{(i)}_k(x,z)\,d\mu(z)=0
=\int_X\overline{D}^{(i)}_k(z,x)\,d\mu(z).
\end{align*}
\end{theorem}

\begin{remark}
\begin{enumerate}
\item[\textup{(i)}]
By the proof of Theorem~\ref{A3.47},
we find that,
for any $i\in\{1,\,2,\,3\}$
and $k\in\mathbb{Z}$,
all the estimates in Theorem~\ref{A3.47}
satisfied by kernels
$\widetilde{D}_k^{(i)}$
and $\overline{D}_k^{(i)}$
are independent of the choices of
$\{y_\tau^{k,v}\}_{\tau\in I_k,\,
v\in\{1,\ldots,N(k,\tau)\}}$
and also independent of both $\beta$ and $\gamma$ with
$\beta,\gamma\in(0,\varepsilon)$ (but may depend on $\varepsilon$).

\item[\textup{(ii)}]
If $\rho\in\mathbf{q}$ is a metric,
then ${\mathrm{ind\,}}(X,\mathbf{q})\in[1,\infty]$
and, for any $0<\beta,\gamma<\varepsilon\preceq{\mathrm{ind\,}}(X,\mathbf{q})$,
\eqref{3.29} holds in $\mathring{{\mathcal G}}^\varepsilon_0(\beta,\gamma)$,
and,
if $\rho\in\mathbf{q}$ is an ultrametric, then,
for any $0<\beta,\gamma<\varepsilon<\infty$,
\eqref{3.29} holds in $\mathring{{\mathcal G}}^\varepsilon_0(\beta,\gamma)$.
However, Han et al. \cite[Theorem 4.11]{hmy08}
and He et al. \cite[Theorem 5.8]{HLYY}
only proved homogeneous discrete Calder\'on
reproducing formulae
for any $0<\beta,\gamma<\varepsilon<1$.
In this sense, Theorem~\ref{A3.47} is sharper than
\cite[Theorem 4.11]{hmy08} and \cite[Theorem 5.8]{HLYY}
because of the sharpness of both Theorems~\ref{ThmATI} and~\ref{CZOonGO}.
\end{enumerate}
\end{remark}

\begin{proof}[Proof of Theorem~\ref{A3.47}]
Let $\rho\in\mathbf{q}$ be a regularized quasi-ultrametric.
We only show the first equality in \eqref{3.29}
because the proofs of the other equalities
in \eqref{3.29} are similar
(see Remark~\ref{ntt-947} below for additional details).
Fix two sufficiently large constants $N,j\in\mathbb{Z}$
such that Corollary~\ref{A3.45} holds,
where $\beta,\gamma\in(0,\varepsilon)$
and $\varepsilon'\in(\max\{\beta,\,\gamma\},\varepsilon)$ are fixed.
For any $k\in\mathbb{Z}$,
let $\mathcal{D}^N_k$ be the same as in \eqref{3.6} and,
for any $x,y\in X$, define
\begin{align*}
\widetilde{D}^{(1)}_k(x,y):=
\widetilde{D}_k(x,y):=\mathfrak{S}_N^{-1}\left(D^N_k(\cdot,y)\right)(x).
\end{align*}
We can verify that $\widetilde{D}^{(1)}_k$ is well defined,
given \eqref{2123} and Corollary~\ref{A3.45}.
Moreover, from
Corollary~\ref{A3.45} and an argument
similar to that used in the proof of Theorem~\ref{A3.19},
we deduce that, for any $k\in\mathbb{Z}$ and $x,y\in X$,
$\widetilde{D}_k(x,y)$ satisfies all the conclusions
of Theorem~\ref{A3.47} except for
the convergence of the series
in the first equality in \eqref{3.29}.

To prove the convergence of this series, we need to verify that
all the series in the first summation
and the second summation of
the first equality in \eqref{3.29}
converge in $\mathring{{\mathcal G}}^\varepsilon_0(\beta,\gamma)$.
To simplify the presentation, by similarity,
we show this only for the series
in the first summation of
the first equality in \eqref{3.29}.
To this end, we divide
the remaining proof into three steps.

\emph{Step 1}. In this step, we prove that,
for any $f\in\mathring{{\mathcal G}}(\beta',\gamma')$
with $\beta'\in(\beta,\varepsilon)$ and
$\gamma'\in(\gamma,\varepsilon)$,
the first equality in \eqref{3.29}
converges in $\mathring{{\mathcal G}}(\beta,\gamma)$.
To this end,
it suffices to show that
\begin{align}
\label{9260}
\lim_{L\to\infty}
\left\|f-\sum_{|k|\leq L}\sum_{\tau\in I_k}
\sum_{v=1}^{N(k,\tau)}\int_{Q_\tau^{k,v}}
\widetilde{D}_k(\cdot,y)\,d\mu(y)
\mathcal{D}_kf(y_\tau^{k,v})\right\|_
{\mathring{{\mathcal G}}(\beta,\gamma)}=0.
\end{align}
By Corollary~\ref{A3.45},
we find that, for any $L\in\mathbb{N}$,
\begin{align}
\label{9261}
&\left\|f-\sum_{|k|\leq L}\sum_{\tau\in I_k}
\sum_{v=1}^{N(k,\tau)}\int_{Q_\tau^{k,v}}
\widetilde{D}_k(\cdot,y)\,d\mu(y)
\mathcal{D}_kf(y_\tau^{k,v})\right\|_
{\mathring{{\mathcal G}}(\beta,\gamma)}
\nonumber\\
&\quad=\left\|\mathfrak{S}_N^{-1}\left[\sum_{|k|\ge L+1}\sum_{\tau\in I_k}
\sum_{v=1}^{N(k,\tau)}\int_{Q_\tau^{k,v}}
D_k^N(\cdot,y)\,d\mu(y)
\mathcal{D}_kf(y_\tau^{k,v})\right]\right\|_
{\mathring{{\mathcal G}}(\beta,\gamma)}
\nonumber\\
&\quad\lesssim\left\|\sum_{|k|\ge L+1}\sum_{\tau\in I_k}
\sum_{v=1}^{N(k,\tau)}\int_{Q_\tau^{k,v}}
D_k^N(\cdot,y)\,d\mu(y)
\mathcal{D}_kf(y_\tau^{k,v})\right\|_
{\mathring{{\mathcal G}}(\beta,\gamma)}
\nonumber\\
&\quad\leq\sum_{|k|\ge L+1}\left\|\mathcal{E}_kf\right\|_
{\mathring{{\mathcal G}}(\beta,\gamma)},
\end{align}
where, for any $k\in\mathbb{Z}$, $f\in\mathring{{\mathcal G}}(\beta,\gamma)$,
and $x\in X$,
\begin{align*}
\mathcal{E}_kf(x):=\int_XE_k(x,y)f(y)\,d\mu(y)
\end{align*}
with $E_k$ as in \eqref{926}.
From Lemma~\ref{A3.46} and an argument similar to that
used in the proof of \eqref{3.16},
we infer that there exists $\sigma\in(0,\infty)$ such that,
for any $k\in\mathbb{Z}$
and $f\in\mathring{{\mathcal G}}(\beta',\gamma')$
with both $\beta'\in(\beta,\varepsilon)$
and $\gamma'\in(\gamma,\varepsilon)$,
\begin{align}
\label{9262}
\left\|\mathcal{E}_kf\right\|_
{\mathring{{\mathcal G}}(\beta,\gamma)}
\lesssim2^{-|k|\sigma}\|f\|_{\mathring{{\mathcal G}}(\beta',\gamma')},
\end{align}
which, together with \eqref{9261},
further implies that \eqref{9260} holds.
This finishes the proof of Step 1.

\emph{Step 2.} In this step, we prove that,
for any $f\in\mathring{{\mathcal G}}^\varepsilon_0(\beta,\gamma)$,
the first equality in \eqref{3.29}
converges in $\mathring{{\mathcal G}}^\varepsilon_0(\beta,\gamma)$.
To this end, we first claim that,
for any $g\in\mathring{{\mathcal G}}^\varepsilon_0(\beta,\gamma)$,
$\mathfrak{S}_N(g)\in\mathring{{\mathcal G}}^\varepsilon_0(\beta,\gamma)$.
Indeed, for any $g\in\mathring{{\mathcal G}}^\varepsilon_0(\beta,\gamma)$,
there exists a sequence $\{g_i\}_{i\in\mathbb{N}}$ in
$\mathring{{\mathcal G}}(\varepsilon,\varepsilon)$ such that
\begin{align*}
\lim_{i\to\infty}\|g-g_i\|_{\mathring{{\mathcal G}}(\beta,\gamma)}=0.
\end{align*}
By both the last inequality of \eqref{9261} and \eqref{9262}, we conclude that,
for any $i\in\mathbb{N}$,
\begin{align*}
\mathfrak{S}_N(g_i)=\sum_{k=-\infty}^\infty\sum_{\tau\in I_k}
\sum_{v=1}^{N(k,\tau)}\int_{Q_\tau^{k,v}}D_k^N(\cdot,y)\,d\mu(y)
D_kg_i(y_\tau^{k,v})
\end{align*}
converges in $\mathring{{\mathcal G}}(\beta,\gamma)$.
Since
\begin{align*}
\int_{Q_\tau^{k,v}}D_k^N(\cdot,y)\,d\mu(y)\in\mathring{{\mathcal G}}
(\varepsilon,\varepsilon),
\end{align*}
it follows that, for any $i\in\mathbb{N}$,
$\mathcal{E}_kg_i\in\mathring{{\mathcal G}}^\varepsilon_0(\beta,\gamma)$,
which, combined with the definitions of
$\mathcal{E}_k$ and $\mathfrak{S}_N$, further implies that
$\mathfrak{S}_N(g_i)\in\mathring{{\mathcal G}}^\varepsilon_0(\beta,\gamma)$.
Note that, by Proposition~\ref{A3.44} and
Corollary~\ref{A3.45}, we find that
$\|\mathfrak{R}_N\|_{\mathring{{\mathcal G}}(\beta,\gamma)
\to\mathring{{\mathcal G}}(\beta,\gamma)}\leq\frac{1}{2}$.
From this and \eqref{RGRN}, we infer that
\begin{align*}
&\lim_{i\to\infty}\left\|\mathfrak{S}_N(g)
-\mathfrak{S}_N(g_i)\right\|_{\mathring{{\mathcal G}}(\beta,\gamma)}\\
&\quad=\lim_{i\to\infty}\left\|(\mathcal{I}-\mathfrak{R}_N)(g-g_i)
\right\|_{\mathring{{\mathcal G}}(\beta,\gamma)}
\lesssim\lim_{i\to\infty}
\left\|(g-g_i)\right\|_{\mathring{{\mathcal G}}(\beta,\gamma)}
=0,
\end{align*}
which, together with
$\mathfrak{S}_N(g_i)\in\mathring{{\mathcal G}}^\varepsilon_0(\beta,\gamma)$
for any $i\in\mathbb{N}$,
further implies that $\mathfrak{S}_N(g)\in
\mathring{{\mathcal G}}^\varepsilon_0(\beta,\gamma)$.
This finishes the proof of the above claim.

By the above claim and \eqref{RGRN}, we conclude that,
for any $g\in\mathring{{\mathcal G}}^\varepsilon_0(\beta,\gamma)$,
$\mathfrak{R}_Ng\in\mathring{{\mathcal G}}^\varepsilon_0(\beta,\gamma)$.
From this, \eqref{RGRN}, and
$\|\mathfrak{R}_N\|_{\mathring{{\mathcal G}}(\beta,\gamma)
\to\mathring{{\mathcal G}}(\beta,\gamma)}\leq\frac{1}{2}$,
we further deduce that both $\mathfrak{S}_N$
and $\mathfrak{S}_N^{-1}$ are bounded on
$\mathring{{\mathcal G}}^\varepsilon_0(\beta,\gamma)$.

Now, let $f\in\mathring{{\mathcal G}}^\varepsilon_0(\beta,\gamma)$.
Then there exists a sequence $\{f_i\}_{i\in\mathbb{N}}$ in
$\mathring{{\mathcal G}}(\varepsilon,\varepsilon)$ such that
\begin{align}
\label{1728}
\lim_{i\to\infty}\|f-f_i\|_{\mathring{{\mathcal G}}(\beta,\gamma)}=0.
\end{align}
For any $L\in\mathbb{N}$, the kernel $E^{(L)}$ of the
integral operator $\mathcal{E}^{(L)}$
is defined by setting, for any $x,y\in X$,
\begin{align*}
E^{(L)}(x,y):=\sum_{|k|\leq L}
E_k(x,y).
\end{align*}
Then repeating the proof of Lemma~\ref{A3.46}
for $E^{(L)}$, we find that $E^{(L)}$
satisfies all the conclusions of Lemma~\ref{A3.46}
with positive constants independent of $L$.
This, combined with arguments similar to those
used in the proofs of both Lemmas~\ref{A3.41} and~\ref{GL2},
we further conclude that $\mathcal{E}^{(L)}$
is bounded on $\mathring{{\mathcal G}}(\beta,\gamma)$
with its operator norm bounded by a positive constant independent of $L$.
From this and the boundedness of $\mathfrak{S}_N$ on
$\mathring{{\mathcal G}}(\beta,\gamma)$,
it follows that, for any $f\in\mathring{{\mathcal G}}^\varepsilon_0(\beta,\gamma)$,
\begin{align*}
\left\|\mathfrak{S}_N(f)-\mathcal{E}^{(L)}f
\right\|_{\mathring{{\mathcal G}}^\varepsilon_0(\beta,\gamma)}
&\leq\left\|\mathfrak{S}_N(f-f_i)\right\|_{\mathring{{\mathcal G}}(\beta,\gamma)}
+\left\|\mathfrak{S}_N(f_i)-\mathcal{E}^{(L)}
f_i\right\|_{\mathring{{\mathcal G}}(\beta,\gamma)}
\\
&\quad+\left\|\mathcal{E}^{(L)}(f-f_i)
\right\|_{\mathring{{\mathcal G}}(\beta,\gamma)}
\\
&\lesssim\left\|f-f_i\right\|_{\mathring{{\mathcal G}}(\beta,\gamma)}
+\left\|\mathfrak{S}_N(f_i)-\mathcal{E}^{(L)}
f_i\right\|_{\mathring{{\mathcal G}}(\beta,\gamma)}.
\end{align*}
By this and the boundedness of $\mathfrak{S}_N^{-1}$
on $\mathring{{\mathcal G}}^\varepsilon_0(\beta,\gamma)$,
we find that,
for any $f\in\mathring{{\mathcal G}}^\varepsilon_0(\beta,\gamma)$,
\begin{align*}
&\left\|f-\sum_{|k|\leq L}\sum_{\tau\in I_k}
\sum_{v=1}^{N(k,\tau)}\int_{Q_\tau^{k,v}}
\widetilde{D}_k(\cdot,y)\,d\mu(y)
\mathcal{D}_kf(y_\tau^{k,v})
\right\|_{\mathring{{\mathcal G}}^\varepsilon_0(\beta,\gamma)}
\\
&\quad=\left\|\mathfrak{S}_N^{-1}\left[\mathfrak{S}_N(f)
-\mathcal{E}^{(L)}f\right]\right\|_
{\mathring{{\mathcal G}}^\varepsilon_0(\beta,\gamma)}
\lesssim\left\|\mathfrak{S}_N(f)-\mathcal{E}^{(L)}f\right\|_
{\mathring{{\mathcal G}}^\varepsilon_0(\beta,\gamma)}
\\
&\quad\lesssim\left\|f-f_i\right\|_{\mathring{{\mathcal G}}(\beta,\gamma)}
+\left\|\mathfrak{S}_N(f_i)-\mathcal{E}^{(L)}
f_i\right\|_{\mathring{{\mathcal G}}(\beta,\gamma)}.
\end{align*}
This, together with the boundedness of $\mathfrak{S}_N$
on $\mathring{{\mathcal G}}(\beta,\gamma)$,
\eqref{1728},
and the conclusion of Step 1,
implies that
\begin{align*}
&\left\|f-\sum_{|k|\leq L}\sum_{\tau\in I_k}
\sum_{v=1}^{N(k,\tau)}\int_{Q_\tau^{k,v}}
\widetilde{D}_k(\cdot,y)\,d\mu(y)
\mathcal{D}_kf(y_\tau^{k,v})
\right\|_{\mathring{{\mathcal G}}^\varepsilon_0(\beta,\gamma)}
\\
&\quad\lesssim\left\|f-f_i\right\|_{\mathring{{\mathcal G}}(\beta,\gamma)}
+\left\|\mathfrak{S}_N^{-1}\left[\mathfrak{S}_N(f_i)-
\mathcal{E}^{(L)}
f_i\right]\right\|_{\mathring{{\mathcal G}}(\beta,\gamma)}
\\
&\quad=\left\|f-f_i\right\|_{\mathring{{\mathcal G}}(\beta,\gamma)}
+\left\|f_i-\sum_{|k|\leq L}\sum_{\tau\in I_k}
\sum_{v=1}^{N(k,\tau)}\int_{Q_\tau^{k,v}}
\widetilde{D}_k(\cdot,y)\,d\mu(y)\mathcal{D}_kf_i(y_\tau^{k,v})
\right\|_{\mathring{{\mathcal G}}(\beta,\gamma)}\\
&\quad\to0
\end{align*}
via first letting $L\to\infty$ and then letting $i\to\infty$,
which further implies that
the first equality of \eqref{3.29}
converges in $\mathring{{\mathcal G}}^\varepsilon_0(\beta,\gamma)$.
This finishes the proof of Step 2.

\emph{Step 3.} In this step, we show that
the first equality in \eqref{3.29}
converges in $L^p(X)$ for any $p\in(1,\infty)$.
Let $p\in(1,\infty)$.
For any $L\in\mathbb{N}$, $f\in L^p(X)$, and $x\in X$,
let
\begin{align*}
\widetilde{\mathcal{E}}^{(L)}f(x)
:=\mathfrak{S}_N^{-1}\left(\mathcal{E}^{(L)}f\right)
=\sum_{|k|\leq L}\sum_{\tau\in I_k}
\sum_{v=1}^{N(k,\tau)}\int_{Q_\tau^{k,v}}
\widetilde{D}_k(x,y)\,d\mu(y)
\mathcal{D}_kf(y_\tau^{k,v})
\end{align*}
with its kernel $\widetilde{E}^{(L)}$ defined by setting,
for any $x,y\in X$,
\begin{align*}
\widetilde{E}^{(L)}(x,y):=
\sum_{|k|\leq L}\sum_{\tau\in I_k}
\sum_{v=1}^{N(k,\tau)}\int_{Q_\tau^{k,v}}
\widetilde{D}_k(x,z)\,d\mu(z)D_k(y_\tau^{k,v},y).
\end{align*}
From Lemma~\ref{A3.46}
and an argument similar to that used in the proof of \eqref{G1},
we deduce that, for any $L\in\mathbb{N}$,
$\mathcal{E}^{(L)}$ is bounded on $L^2(X)$,
which, combined with the boundedness of $\mathfrak{S}_N^{-1}$
on $L^2(X)$, further implies that
$\widetilde{\mathcal{E}}^{(L)}$ is bounded on $L^2(X)$.
By Lemma~\ref{A3.46}
and an argument similar to that used in the proof of
Theorem~\ref{A3.19}, we find that, for any $L\in\mathbb{N}$,
$\widetilde{E}^{(L)}$,
the kernel of $\widetilde{\mathcal{E}}^{(L)}$,
satisfies \eqref{CZk1} and \eqref{CZk2}.
Therefore, for any $L\in\mathbb{N}$,
$\widetilde{E}^{(L)}$
is a $\theta$-Calder\'on--Zygmund kernel
for some $\theta\in(0,\gamma\land\beta)$ as in Definition~\ref{CZk}.
This, together with the well-known Calder\'on--Zygmund
theory on spaces of
homogeneous type developed in \cite[pp.\,74--75, Theorem 2.4]{CW1},
further implies that $\widetilde{\mathcal{E}}^{(L)}$
is bounded on $L^p(X)$
with its operator norm bounded by a positive constant independent of $L$.
Applying this, Lemma~\ref{density}(ii),
a standard density argument,
and an argument similar to that used in the proof of
Theorem~\ref{A3.19}, we obtain
the convergence of first series in \eqref{3.29}
in $L^p(X)$.
This finishes the proof of Step 3
and hence Theorem~\ref{A3.47}.
\end{proof}

\begin{remark}
\label{ntt-947}
To obtain the other five equalities in \eqref{3.29},
instead of considering $\mathfrak{S}_N$
in \eqref{DRS}, we only need to consider,
for any fixed sufficiently large $N\in\mathbb{N}$
and for any $f\in L^2(X)$ and $x\in X$,
\begin{align*}
\widetilde{\mathfrak{S}}_N^{(2)}f(x):=
\sum_{k=-\infty}^{\infty}\sum_{\tau\in I_k}
\sum_{v=1}^{N(k,\tau)}
D_k^N\left(x,y^{k,v}_\tau\right)
\int_{Q^{k,v}_\tau}\mathcal{D}_kf(y)\,d\mu(y),
\end{align*}
\begin{align*}
\widetilde{\mathfrak{S}}_N^{(3)}f(x):=
\sum_{k=-\infty}^{\infty}\sum_{\tau\in I_k}
\sum_{v=1}^{N(k,\tau)}
\mu\left(Q^{k,v}_\tau\right)
D^N_k\left(x,y^{k,v}_\tau\right)
\mathcal{D}_kf(y^{k,v}_\tau),
\end{align*}
\begin{align*}
\overline{\mathfrak{S}}_N^{(1)}f(x):=
\sum_{k=-\infty}^{\infty}\sum_{\tau\in I_k}
\sum_{v=1}^{N(k,\tau)}
\int_{Q^{k,v}_\tau}D_k(x,y)\,d\mu(y)
\mathcal{D}^N_kf(y^{k,v}_\tau),
\end{align*}
\begin{align*}
\overline{\mathfrak{S}}_N^{(2)}f(x):=
\sum_{k=-\infty}^{\infty}\sum_{\tau\in I_k}
\sum_{v=1}^{N(k,\tau)}
D_k\left(x,y^{k,v}_\tau\right)
\int_{Q^{k,v}_\tau}\mathcal{D}^N_kf(y)\,d\mu(y),
\end{align*}
and
\begin{align*}
\overline{\mathfrak{S}}_N^{(3)}f(x):=
\sum_{k=-\infty}^{\infty}\sum_{\tau\in I_k}
\sum_{v=1}^{N(k,\tau)}
\mu\left(Q^{k,v}_\tau\right)
D_k\left(x,y^{k,v}_\tau\right)
\mathcal{D}^N_kf(y^{k,v}_\tau),
\end{align*}
respectively. The corresponding remainders are defined by setting
\begin{align*}
\widetilde{\mathfrak{R}}_N^{(2)}:=I-\widetilde{\mathfrak{S}}_N^{(2)},\ \
\widetilde{\mathfrak{R}}_N^{(3)}:=I-\widetilde{\mathfrak{S}}_N^{(3)},\\
\overline{\mathfrak{R}}_N^{(1)}:=I-\overline{\mathfrak{S}}_N^{(1)},\ \
\overline{\mathfrak{R}}_N^{(2)}:=I-\overline{\mathfrak{S}}_N^{(2)},
\end{align*}
and
$$
\overline{\mathfrak{R}}_N^{(3)}:=I-\overline{\mathfrak{S}}_N^{(3)},
$$
which satisfy the same estimates of $\mathfrak{R}_N$ as in Proposition~\ref{A3.44}.
The remaining arguments are similar. We omit the details.
\end{remark}

By Theorem~\ref{A3.47} and a standard duality argument,
we obtain the following sharp homogeneous discrete Calder\'on
reproducing formula on spaces of distributions;
we omit the details of its proof.

\begin{theorem}\label{HD}
Let $(X,\mathbf{q},\mu)$ be an ultra-RD-space,
where $\mu$ is assumed to be a Borel-semiregular measure on $X$.
Let $\varepsilon$, $\beta$, $\gamma$,
and $\{\mathcal{D}_k\}_{k\in\mathbb{Z}}$
be the same as in Theorem~\ref{A3.47}. Then, for any
$f\in(\mathring{{\mathcal G}}_0^\varepsilon(\beta,\gamma))'$,
\eqref{3.29} holds in
$(\mathring{{\mathcal G}}_0^\varepsilon(\beta,\gamma))'$.
\end{theorem}

\section{Inhomogeneous Discrete
Calder\'on Reproducing Formulae}
\label{SS3.5}

The target of this section is to establish sharp
inhomogeneous discrete Calder\'on reproducing formulae
on the ultra-RD-space $(X,\mathbf{q},\mu)$,
where $\mu$ is assumed to be a Borel-semiregular measure on $X$.
\index[words]{inhomogeneous discrete Calder\'on reproducing formulae}
In particular, given $\rho\in\mathbf{q}$,
we have no restriction on $\mathrm{diam}_\rho(X)$ and hence we can have
$\mathrm{diam}_\rho(X)<\infty$ or $\mathrm{diam}_\rho(X)=\infty$.
The proofs of the main results of this section are similar to those presented
in Sections~\ref{S3.3} and~\ref{SS3.4}.
As such,
we provide a detailed account of their assertions and leave the
details of the proofs to the reader.

\begin{theorem}\label{ID}
Let $(X,\mathbf{q},\mu)$ be an ultra-RD-space,
where $\mu$ is assumed to be a Borel-semiregular measure on $X$.
Assume that $0<\varepsilon\preceq{\mathrm{ind\,}}(X,\mathbf{q})$
with ${\mathrm{ind\,}}(X,\mathbf{q})$ as in \eqref{index}.
Let $\{\mathcal{S}_{2^{-k}}\}_{k\in\mathbb{Z}_+}$
be an inhomogeneous approximation
of the identity of order $\varepsilon$
as in Definition~\ref{DefIATI} (with kernels
$\{S_{2^{-k}}\}_{k\in\mathbb{Z}_+}$).
For any $k\in\mathbb{N}$, let
\begin{align}\label{DKZ+}
\mathcal{D}_k:=\mathcal{S}_{2^{-k}}
-\mathcal{S}_{2^{-k+1}}
\ \ \text{and}\ \
\mathcal{D}_0:=\mathcal{S}_1.
\end{align}
Fix
$y^{k,v}_\tau\in Q^{k,v}_\tau$
for any $k\in\mathbb{N}$, $\tau\in I_k$, and
$v\in\{1,\ldots,N(k,\tau)\}$.
Then, for any fixed $j,N\in\mathbb{N}$ sufficiently large
and for any $i\in\{1,\,2,\,3\}$,
there exist two families
$\{\widetilde{\mathcal{D}}^{(i)}_k\}_{k\in\mathbb{Z}_+}$
and $\{\overline{\mathcal{D}}^{(i)}_k\}_{k\in\mathbb{Z}_+}$
of bounded integral operators on $L^2(X)$
(with their kernels
$\{\widetilde{D}^{(i)}_k\}_{k\in\mathbb{Z}_+}$
and $\{\overline{D}^{(i)}_k\}_{k\in\mathbb{Z}_+}$)
such that, for any $f\in{\mathcal G}_0^\varepsilon(\beta,\gamma)$
with $\beta,\gamma\in(0,\varepsilon)$
[resp. $L^p(X)$ with $p\in(1,\infty)$],
\begin{align}\label{3.35}
f(\cdot)
&=\sum_{k=0}^{N}\sum_{\tau\in I_k}\sum_{v=1}^{N(k,\tau)}
\int_{Q^{k,v}_\tau }
\widetilde{D}^{(1)}_k(\cdot,y)\,d\mu(y)
\mathcal{D}^{k,v}_{\tau,1}f\nonumber\\
&\quad+\sum_{k=N+1}^{\infty}\sum_{\tau\in I_k}
\sum_{v=1}^{N(k,\tau)}\int_{Q^{k,v}_\tau}
\widetilde{D}^{(1)}_k(\cdot,y)\,d\mu(y)
\mathcal{D}_kf(y^{k,v}_\tau)\nonumber\\
&=\sum_{\tau\in I_0}\sum_{v=1}^{N(0,\tau)}\int_{Q^{0,v}_\tau}
\widetilde{D}^{(2)}_0(\cdot,y)\,d\mu(y)
\mathcal{D}^{0,v}_{\tau,1}f
\nonumber\\
&\quad+\sum_{k=1}^{\infty}\sum_{\tau\in I_k}\sum_{v=1}^{N(k,\tau)}
\mu(Q^{k,v}_\tau)
\widetilde{D}^{(2)}_k(\cdot,y^{k,v}_\tau)
\mathcal{D}^{k,v}_{\tau,1}f
\nonumber\\
&=\sum_{\tau\in I_0}\sum_{v=1}^{N(0,\tau)}\int_{Q^{0,v}_\tau}
\widetilde{D}^{(3)}_0(\cdot,y)\,d\mu(y)
\mathcal{D}^{0,v}_{\tau,1}f
\nonumber\\
&\quad+\sum_{k=1}^{N}\sum_{\tau\in I_k}\sum_{v=1}^{N(k,\tau)}
\mu(Q^{k,v}_\tau)
\widetilde{D}^{(3)}_k(\cdot,y^{k,v}_\tau)
\mathcal{D}^{k,v}_{\tau,1}f
\nonumber\\
&\quad+\sum_{k=N+1}^{\infty}\sum_{\tau\in I_k}
\sum_{v=1}^{N(k,\tau)}
\mu(Q^{k,v}_\tau)
\widetilde{D}^{(3)}_k(\cdot,y^{k,v}_\tau)
\mathcal{D}_kf(y^{k,v}_\tau)
\end{align}
and
\begin{align}\label{3.351}
f(\cdot)&=\sum_{k=0}^{N}\sum_{\tau\in I_k}\sum_{v=1}^{N(k,\tau)}
\int_{Q^{k,v}_\tau }
D_k(\cdot,y)\,d\mu(y)
\overline{\mathcal{D}}_{\tau,1}^{(1),k,v}f
\nonumber\\
&\quad+\sum_{k=N+1}^{\infty}\sum_{\tau\in I_k}
\sum_{v=1}^{N(k,\tau)}
\int_{Q^{k,v}_\tau}
D_k(\cdot,y)\,d\mu(y)
\overline{\mathcal{D}}_k^{(1)}f(y^{k,v}_\tau)
\nonumber\\
&=\sum_{\tau\in I_0}\sum_{v=1}^{N(0,\tau)}\int_{Q^{0,v}_\tau}
D_0(\cdot,y)\,d\mu(y)
\overline{\mathcal{D}}_{\tau,1}^{(2),0,v}f
\nonumber\\
&\quad+\sum_{k=1}^{\infty}\sum_{\tau\in I_k}\sum_{v=1}^{N(k,\tau)}
\mu(Q^{k,v}_\tau)D_k(\cdot,y^{k,v}_\tau)
\overline{\mathcal{D}}_{\tau,1}^{(2),k,v}f
\nonumber\\
&=\sum_{\tau\in I_0}\sum_{v=1}^{N(0,\tau)}\int_{Q^{0,v}_\tau}
D_0(\cdot,y)\,d\mu(y)
\overline{\mathcal{D}}_{\tau,1}^{(3),0,v}f
\nonumber\\
&\quad+\sum_{k=1}^{N}\sum_{\tau\in I_k}\sum_{v=1}^{N(k,\tau)}
\mu(Q^{k,v}_\tau)D_k(\cdot,y^{k,v}_\tau)
\overline{\mathcal{D}}_{\tau,1}^{(3),k,v}f
\nonumber\\
&\quad+\sum_{k=N+1}^{\infty}\sum_{\tau\in I_k}\sum_{v=1}^{N(k,\tau)}
\mu(Q^{k,v}_\tau)D_k(\cdot,y^{k,v}_\tau)
\overline{\mathcal{D}}_{\tau,1}^{(3)}f(y^{k,v}_\tau),
\end{align}
where all the series converge in ${\mathcal G}_0^\varepsilon(\beta,\gamma)$
[resp. $L^p(X)$ with $p\in(1,\infty)$]
and also converge pointwise $\mu$-almost everywhere on $X$
and where, for any $k\in\mathbb{Z}_+$, $\tau\in I_k$,
and $v\in\{1,\ldots,N(k,\tau)\}$, $\mathcal{D}^{k,v}_{\tau,1}$ and
$\overline{\mathcal{D}}_{\tau,1}^{(i),k,v}$ with $i\in\{1,\,2,\,3\}$ are
the corresponding integral operators,
respectively, with their kernels defined by setting,
for any $z\in X$,
\begin{align}
\label{9245}
D^{k,v}_{\tau,1}(z)
:=\frac{1}{\mu(Q^{k,v}_\tau)}
\int_{Q^{k,v}_\tau }D_k(u,z)\,d\mu(u)
\end{align}
and
\begin{align*}
\overline{D}_{\tau,1}^{(i),k,v}(z)
:=\frac{1}{\mu(Q^{k,v}_\tau)}
\int_{Q^{k,v}_\tau }\overline{D}^{(i)}_k(u,z)\,d\mu(u).
\end{align*}
Moreover, for any $i\in\{1,\,2,\,3\}$,
$k\in\mathbb{Z}_+$,
$\widetilde{\beta},\widetilde{\gamma}\in(0,\varepsilon)$,
and $x\in X$,
\begin{align*}
\widetilde{D}^{(i)}_k(\cdot,x),\overline{D}^{(i)}_k(x,\cdot)
\in{\mathcal G}(x,2^{-k},\widetilde{\beta},\widetilde{\gamma}),
\end{align*}
\begin{align*}
\int_X\widetilde{D}^{(i)}_k(x,z)\,d\mu(z)
&=\int_X\widetilde{D}^{(i)}_k(z,x)\,d\mu(z)\\
&=
\begin{cases}
1\quad
\mathrm{if}\ k\in\{0,\ldots,N\},\\
0\quad
\mathrm{if}\ k\in\{N+1,N+2,\ldots\},
\end{cases}
\end{align*}
and
\begin{align*}
\int_X\overline{D}^{(i)}_k(x,z)\,d\mu(z)
&=\int_X\overline{D}^{(i)}_k(z,x)\,d\mu(z)\\
&=
\begin{cases}
1\quad
\mathrm{if}\ k\in\{0,\ldots,N\},\\
0\quad
\mathrm{if}\ k\in\{N+1,N+2,\ldots\}.
\end{cases}
\end{align*}
\end{theorem}

\begin{remark}
\begin{enumerate}
\item[\textup{(i)}]
We explain the decomposition of the identity operator
to derive the first equality in \eqref{3.35}.
Let $\{\mathcal{D}_k\}_{k\in\mathbb{Z}_+}$ be
the same as in \eqref{DKZ+},
$\mathcal{D}_k^N$ the same as in \eqref{DkN}
for any fixed $N\in\mathbb{N}$ and for any $k\in\mathbb{Z}_+$,
$\mathcal{R}_N$ the same as in \eqref{3.21}, and $p\in(1,\infty)$.
By \eqref{3.21}, we find that, for any $f\in L^p(X)$
with $p\in(1,\infty)$,
\begin{align*}
f(\cdot)
&=\sum_{k=0}^\infty\mathcal{D}_k^N\mathcal{D}_kf(\cdot)
+\mathcal{R}_Nf(\cdot)
\\
&=\sum_{k=0}^N\sum_{\tau\in I_k}\sum_{v=1}^{N(k,\tau)}
\int_{Q_\tau^{k,v}}D_k^N(\cdot,y)\mathcal{D}_kf(y)\,d\mu(y)
\\
&\quad+\sum_{k=N+1}^\infty\sum_{\tau\in I_k}\sum_{v=1}^{N(k,\tau)}
\int_{Q_\tau^{k,v}}D_k^N(\cdot,y)\mathcal{D}_kf(y)\,d\mu(y)
+\mathcal{R}_Nf(\cdot)
\\
&=:\mathfrak{T}_Nf(\cdot)+\mathcal{R}^{(1)}_Nf(\cdot)
+\mathcal{R}^{(2)}_Nf(\cdot)+\mathcal{R}_Nf(\cdot)
\end{align*}
in $L^p(X)$ with
\begin{align*}
\mathfrak{T}_Nf(\cdot):&=\sum_{k=0}^N\sum_{\tau\in I_k}\sum_{v=1}^{N(k,\tau)}
\int_{Q_\tau^{k,v}}D_k^N(\cdot,y)\,d\mu(y)\mathcal{D}_{\tau,1}^{k,v}f
\\
&\quad+\sum_{k=N+1}^\infty\sum_{\tau\in I_k}\sum_{v=1}^{N(k,\tau)}
\int_{Q_\tau^{k,v}}D_k^N(\cdot,y)\,d\mu(y)\mathcal{D}_kf(y_\tau^{k,v}),
\end{align*}
\begin{align*}
\mathcal{R}^{(1)}_Nf(\cdot):&=\sum_{k=0}^N\sum_{\tau\in I_k}\sum_{v=1}^{N(k,\tau)}
\frac{1}{\mu(Q_\tau^{k,v})}
\int_{Q_\tau^{k,v}}D_k^N(\cdot,y)\\
&\quad\times
\int_{Q_\tau^{k,v}}\left[\mathcal{D}_kf(y)
-\mathcal{D}_kf(u)\right]\,d\mu(u)\,d\mu(y),
\end{align*}
and
\begin{align*}
\mathcal{R}^{(2)}_Nf(\cdot):=\sum_{k=N+1}^\infty\sum_{\tau\in I_k}\sum_{v=1}^{N(k,\tau)}
\int_{Q_\tau^{k,v}}D_k^N(\cdot,y)
\left[\mathcal{D}_kf(y)-\mathcal{D}_kf(y_\tau^{k,v})\right]\,d\mu(y),
\end{align*}
where
$\mathcal{D}_{\tau,1}^{k,v}f$ is the same as in \eqref{9245}.
Then we can prove that
$\mathcal{R}^{(1)}_N$, $\mathcal{R}^{(2)}_N$, and $\mathcal{R}_N$ are bounded on
both $L^2(X)$ and ${\mathcal G}(\beta,\gamma)$;
see Lemma~\ref{A3.28} for the boundedness
of $\mathcal{R}_N$ on $L^2(X)$ and see Proposition~\ref{IBRN} for the boundedness
of $\mathcal{R}_N$ on ${\mathcal G}(\beta,\gamma)$;
the boundedness of $\mathcal{R}^{(2)}_N$ on
both $L^2(X)$ and ${\mathcal G}(\beta,\gamma)$
can be deduced from the same estimates for
kernels $\{G_{N,k}\}_{k\in\mathbb{N}\cap[N+1,\infty)}$ as in
Lemma~\ref{A3.41} combined with an argument
similar to that used in the proof of \eqref{1646124}
in Proposition~\ref{A3.44}.
In addition, for $\mathcal{R}^{(1)}_N$,
we first need to show similar estimates as in Lemma~\ref{A3.41}
for the kernel of $\mathcal{R}^{(1)}_N$
and then apply an argument similar to that used in the proof of Lemma~\ref{GL2}
together with Theorem~\ref{CZOonCG} to obtain the boundedness of $\mathcal{R}^{(1)}_N$
on both $L^2(X)$ and ${\mathcal G}(\beta,\gamma)$.
With these boundedness conclusions in mind,
if we choose $j,N\in\mathbb{N}$ sufficiently large such that
\begin{align*}
\left\|\mathcal{R}^{(1)}_N\right\|_{L^2(X)\to L^2(X)}
+\left\|\mathcal{R}^{(2)}_N\right\|_{L^2(X)\to L^2(X)}
+\left\|\mathcal{R}_N\right\|_{L^2(X)\to L^2(X)}<1
\end{align*}
and
\begin{align*}
\left\|\mathcal{R}^{(1)}_N
\right\|_{{\mathcal G}(\beta,\gamma)\to{\mathcal G}(\beta,\gamma)}
+\left\|\mathcal{R}^{(2)}_N
\right\|_{{\mathcal G}(\beta,\gamma)\to{\mathcal G}(\beta,\gamma)}
+\left\|\mathcal{R}_N
\right\|_{{\mathcal G}(\beta,\gamma)\to{\mathcal G}(\beta,\gamma)}<1,
\end{align*}
then we conclude that $(\mathfrak{T}_N)^{-1}$ exists and is
also bounded on both $L^2(X)$ and ${\mathcal G}(\beta,\gamma)$
and hence we define, for any $k\in\mathbb{Z}_+$ and $x,y\in X$,
\begin{align*}
\widetilde{D}^{(1)}_k(x,y):=(\mathfrak{T}_N)^{-1}\left(D_k^N(\cdot,y)\right)(x).
\end{align*}
The remaining part of the proof of Theorem~\ref{ID} follows from
an argument similar to that used in the proof of Theorem~\ref{A3.47}
and we omit the details.
\item[\textup{(ii)}]
In Theorem~\ref{ID}, for any $i\in\{1,\,2,\,3\}$
and $k\in\mathbb{Z}_+$,
all the estimates satisfied by
$\widetilde{D}_k^{(i)}$ and $\overline{D}_k^{(i)}$
are independent of $\{y_\tau^{k,v}\}_{\tau\in I_k,\,
v\in\{1,\ldots,N(k,\tau)\}}$
and also independent of both $\beta$ and $\gamma$ with
$\beta,\gamma\in(0,\varepsilon)$ (but may depend on $\varepsilon$).
\item[\textup{(iii)}]
If $\rho\in\mathbf{q}$ is a metric,
then ${\mathrm{ind\,}}(X,\mathbf{q})\in[1,\infty]$
and, for any $0<\beta,\gamma<\varepsilon\preceq{\mathrm{ind\,}}(X,\mathbf{q})$,
\eqref{3.35} and \eqref{3.351} hold in
$\mathring{{\mathcal G}}^\varepsilon_0(\beta,\gamma)$,
and,
if $\rho\in\mathbf{q}$ is an ultrametric, then,
for any $0<\beta,\gamma<\varepsilon<\infty$,
\eqref{3.35} and \eqref{3.351} hold in
$\mathring{{\mathcal G}}^\varepsilon_0(\beta,\gamma)$.
However, Han et al. \cite[Theorem 4.14]{hmy08}
and He et al. \cite[Theorem 6.10]{HLYY} only proved inhomogeneous discrete Calder\'on
reproducing formulae
for any $0<\beta,\gamma<\varepsilon<1$.
In this sense, Theorem~\ref{ID} is sharper than
\cite[Theorem 4.14]{hmy08} and \cite[Theorem 6.10]{HLYY}
because of the sharpness of both Theorems~\ref{ThmATI} and~\ref{CZOonCG}.
\end{enumerate}
\end{remark}

By Theorem~\ref{ID} and a standard duality argument,
we obtain the following sharp inhomogeneous discrete Calder\'on
reproducing formula on spaces of distributions;
we omit the details of its proof.

\begin{theorem}\label{ID2}
Let $(X,\mathbf{q},\mu)$ be an ultra-RD-space,
where $\mu$ is assumed to be a Borel-semiregular measure on $X$.
Let $\varepsilon$, $\beta$, $\gamma$,
and $\{\mathcal{D}_k\}_{k\in\mathbb{Z}_+}$
be the same as in Theorem~\ref{ID}.
Then, for any $f\in({\mathcal G}_0^\varepsilon(\beta,\gamma))'$,
\eqref{3.35} and \eqref{3.351} hold in $({\mathcal G}_0^\varepsilon(\beta,\gamma))'$.
\end{theorem}


\chapter[Littlewood--Paley Function Characterizations
of Hardy Spaces on\\ Ultra-RD-Spaces]
{Littlewood--Paley Function Characterizations
of Hardy Spaces on Ultra-RD-Spaces}
\label{Chapter3}
\markboth{\scriptsize\rm\sc Littlewood--Paley Characterizations
of Hardy Spaces}
{\scriptsize\rm\sc Littlewood--Paley Characterizations
of Hardy Spaces}

In \cite{AlMi2015}, Alvarado and Mitrea introduced
the Hardy spaces $H^p_{{\rm max}}(X)$
defined via the grand maximal function on $n$-Ahlfors-regular
quasi-ultrametric spaces $(X,\mathbf{q},\mu)$,
with $n\in(0,\infty)$, and proved, among other things,
that this brand of Hardy spaces has atomic
and molecular characterizations for
$p\in(\frac{n}{n+\mathrm{ind\,}(X,\mathbf{q})},1]$
by showing that $H^p_{{\rm max}}(X)$ can be naturally identified
with the atomic Hardy space in this setting (which is known to coincide
with the Coifman--Weiss atomic Hardy space from \cite{Cowe77}).

In this chapter, we continue the work initiated in
\cite{AlMi2015} by establishing the Littlewood--Paley characterizations
of the Hardy spaces introduced in \cite{AlMi2015} in
$n$-Ahlfors-regular quasi-ultrametric spaces (see Remark~\ref{2548}).
Indeed, we establish the Littlewood--Paley characterizations
of the Hardy spaces in a more general setting,
namely $(\kappa,n)$-ultra-RD-spaces $(X,\mathbf{q},\mu)$ with
$0<\kappa\leq n<\infty$.
To this end, in Section~\ref{3.1},
we first recall the definitions of both several kinds of Littlewood--Paley
functions, based on the approximation of the identity in Section~\ref{ss2.2},
and the Hardy space $H^p(X)$ defined via the Lusin area function for
the optimal range of
$
p\in(\frac{n}{n+\mathrm{ind\,}(X,\mathbf{q})},\infty).
$
Then, using the
homogeneous discrete Calder\'on reproducing formulae
established in Section~\ref{SS3.4}, we prove that,
for any $p\in(\frac{n}{n+\mathrm{ind\,}(X,\mathbf{q})},\infty)$,
$H^p(X)$ is independent of the choice of
homogeneous approximations
of the identity.
In Section~\ref{Section3.2}, we recall the definition of
the atomic Hardy space $H^{p}_{\mathrm{at}}(X)$ and
show that $H^p(X)$ has an atomic characterization
for $p\in(\frac{n}{n+\mathrm{ind\,}(X,\mathbf{q})},1]$
via the identification $H^p(X)=H^{p}_{\mathrm{at}}(X)$.
We conclude this chapter with Section~\ref{Section3.3}
wherein we also establish
Littlewood--Paley $g$-function and $g_\lambda^*$-function
characterizations of $H^p(X)$
with $p\in(\frac{n}{n+\mathrm{ind\,}(X,\mathbf{q})},1]$.
In this chapter, unless otherwise specified,
we always assume that
$\mu$ is a Borel-semiregular measure on $X$ and $\mathrm{diam}_\rho(X)=\infty$
for some $\rho\in\mathbf{q}$ and hence
$\mathrm{diam}_\rho(X)=\infty$ for all $\rho\in\mathbf{q}$.

\section{Hardy Spaces via Lusin Area Functions}
\label{3.1}

Let $(X,\mathbf{q},\mu)$ be an ultra-RD-space.
In this section, we first
introduce the definitions of several kinds of Littlewood--Paley
functions defined using the approximation of the identity
given in Definition~\ref{DefATI}, then define
the Hardy space $H^p(X)$
by means of the Lusin area function,
and finally prove that this definition is independent of the choice of
approximations of the identity. To this end,
we begin with recalling the definitions of various brands of
various Littlewood--Paley functions;
see, for instance, \cite[(5.1), (5.2), and (5.3)]{ZYJLDW}.

\begin{definition}\label{areafncts}
Let $(X,\mathbf{q},\mu)$ be a quasi-ultrametric space
of homogeneous type
and $\rho\in\mathbf{q}$ be such that all $\rho$-balls
are $\mu$-measurable. Assume that
$0<\beta,\gamma<\varepsilon<\infty$
and $\{\mathcal{S}_{2^{-k}}\}_{k\in\mathbb{Z}}$
is a homogeneous approximation
of the identity of order $\varepsilon$
as in Definition~\ref{DefHATI}.
Let $f\in(\mathring{\mathcal{G}}_0^\varepsilon(\beta,\gamma))'$
and $\mathcal{D}_k:=\mathcal{S}_{2^{-k}}-\mathcal{S}_{2^{-k+1}}$
for any $k\in\mathbb{Z}$.
\begin{enumerate}
\item
For any given $\theta\in(0,\infty)$,
the \emph{Lusin area function $S_\theta(f)$} of $f$ with
\emph{aperture $\theta$}\index[words]{Lusin area function with aperture $\theta$}
is defined by setting, for any $x\in X$,
\begin{equation}\label{Def-S-f}
S_\theta(f)(x)\index[symbols]{S@$S_\theta(f)$}
:=\left[\sum_{k=-\infty}^{\infty}
\int_{B_\rho(x,\theta2^{-k})}
|\mathcal{D}_kf(y)|^2\frac{d\mu(y)}
{(V_\rho)_{\theta 2^{-k}}(x)}\right]^\frac{1}{2}.
\end{equation}
When $\theta:=1$, we write $S_1$ simply as $S$.

\item
The \emph{Littlewood--Paley $g$-function}\index[words]{Littlewood--Paley $g$-function}
$g(f)$ of $f$ is defined by setting,
for any $x\in X$,
\begin{equation}\label{Def-g(f)(x)}
g(f)(x)\index[symbols]{G@$g(f)$}
:=\left[\sum_{k=-\infty}^{\infty}
\left|\mathcal{D}_kf(x)\right|^2\right]^\frac{1}{2}.
\end{equation}

\item
The \emph{Littlewood--Paley $g_\lambda^*$-function} $g_\lambda^*(f)$ of $f$,
\index[words]{Littlewood--Paley $g_\lambda^*$-function}
with any given $\lambda\in(0,\infty)$,
is defined by setting,
for any $x\in X$,
\begin{equation}\label{Def-glamda*}
g_\lambda^*(f)(x)\index[symbols]{G@$g_\lambda^*$}
:=\left\{\sum_{k=-\infty}^{\infty}
\int_X\left|\mathcal{D}_kf(y)\right|^2
\left[\frac{2^{-k}}{2^{-k}+\rho(x,y)}\right]^\lambda
\frac{d\mu(y)}{(V_\rho)_{2^{-k}}(x)+(V_\rho)_{2^{-k}}(y)}\right\}^\frac{1}{2}.
\end{equation}
\end{enumerate}
\end{definition}

\begin{remark}\label{rb-375}
Let the notation be as in Definition~\ref{areafncts}.
If $S_\theta(f)$ is constructed using $\rho$ and $S'_\theta(f)$
is constructed using $\rho'\in\mathbf{q}$, where all $\rho$-balls
and $\rho'$-balls are $\mu$-measurable,
then $S_\theta(f)\sim S'_\theta(f)$ pointwise in $X$.
Similar equivalences for $g(f)$ and $g_\lambda^*(f)$ also hold.
Thus, without loss of generality, we may assume that $\rho\in\mathbf{q}$
in Definition~\ref{areafncts} is a regularized quasi-ultrametric
in the remainder of this monograph.
Indeed, this assumption made on $\rho$
guarantees that the Littlewood--Paley functions in Definition~\ref{areafncts}
are $\mu$-measurable; see Lemma~\ref{1955} below.
\end{remark}

The following lemma is valid if $\mu(X)<\infty$ or $\mu(X)=\infty$,
equivalently, when $X$ is bounded or unbounded.

\begin{lemma}\label{1955}
Let the notation be as in Definition~\ref{areafncts}.
Then, for any $f\in(\mathring{{\mathcal{G}}}_0^\varepsilon(\beta,\gamma))'$,
the functions $S_\theta(f)$, $g(f)$, and $g_\lambda^*(f)$, constructed using
a regularized quasi-ultrametric $\rho\in\mathbf{q}$ as in Definition~\ref{D2.3},
are $\mu$-measurable.
\end{lemma}

\begin{proof}
Let $f\in(\mathring{{\mathcal{G}}}_0^\varepsilon(\beta,\gamma))'$.
If $f=0$ in $(\mathring{{\mathcal{G}}}_0^\varepsilon(\beta,\gamma))'$, then
the claim follows immediately. Thus, we may assume that
$f\neq0$ in $(\mathring{{\mathcal{G}}}_0^\varepsilon(\beta,\gamma))'$; equivalently,
$\left\|f\right\|_{(\mathring{{\mathcal{G}}}_0^\varepsilon(\beta,\gamma))'}>0$.

We first show that $g(f)$ is $\mu$-measurable.
Indeed, it suffices to prove that,
for any $k\in\mathbb{Z}$, $\mathcal{D}_kf$ is $\mu$-measurable.
To this end, according to Remark~\ref{Rmk502}(ii), we only need to show that,
for any $\lambda\in\mathbb{R}$,
\begin{align*}
\Omega_\lambda:=\left\{x\in X:\mathcal{D}_kf(x)>\lambda\right\}
\end{align*}
is open.
For any given $\lambda\in\mathbb{R}$ such that $\Omega_\lambda\neq\varnothing$,
we arbitrarily choose a point $x\in\Omega_\lambda$.
By Remark~\ref{Dkf2}(i) and
an argument similar to that used in the
estimation of \eqref{2235},
we conclude that
there exists a positive constant $C$ such that,
for any $x'\in X$,
\begin{align*}
\left|\mathcal{D}_kf(x)-\mathcal{D}_kf(x')\right|
&=\left|\left\langle f,D_k(x,\cdot)-D_k(x',\cdot)\right\rangle\right|\\
&\leq\left\|f\right\|_{(\mathring{{\mathcal{G}}}_0^\varepsilon(\beta,\gamma))'}
\left\|D_k(x,\cdot)-D_k(x',\cdot)\right\|_{\mathcal{G}(\beta,\gamma)}\\
&\leq C\left\|f\right\|_{(\mathring{{\mathcal{G}}}_0^\varepsilon(\beta,\gamma))'}
\left[\rho(x,x')\right]^\varepsilon
\left[2^{-k}+\rho(x_1,x)\right]^{n+\beta+\gamma}.
\end{align*}
From this, we deduce that,
for any $x'\in B_\rho(x,C_x)$,
\begin{align*}
\mathcal{D}_kf(x')&\ge\mathcal{D}_kf(x)-
\left|\mathcal{D}_kf(x)-\mathcal{D}_kf(x')\right|>\lambda,
\end{align*}
where
\begin{align*}
C_x:=\left\{\frac{\mathcal{D}_kf(x)-\lambda}
{C\left\|f\right\|_{(\mathring{{\mathcal{G}}}_0^\varepsilon(\beta,\gamma))'}
\left[1+\rho(x_1,x)\right]^{n+\beta+\gamma}}\right\}^\frac{1}{\varepsilon},
\end{align*}
which further implies that
$B_\rho(x,C_x)\subset\Omega_\lambda$
for any $x\in\Omega_\lambda$.
Thus, $\Omega_\lambda$ is open and hence we finish the proof of the fact
that $g(f)$ is $\mu$-measurable.

Next, we prove that $S_\theta(f)$
is $\mu$-measurable. Without loss of generality,
we may assume $\theta=1$.
It suffices to show that, for any $k\in\mathbb{Z}$,
the function
$$
\fint_{B_\rho(\cdot,2^{-k})}\left|\mathcal{D}_kf(y)\right|^2\,d\mu(y)
$$
is $\mu$-measurable.
By \cite[p.\,15, Theorem 1.17]{r1987} and
the above proven conclusion that $\mathcal{D}_kf$
for any $k\in\mathbb{Z}$
is $\mu$-measurable,
we find that there exists a sequence
of simple functions $\{h_j\}_{j\in\mathbb{N}}$ such that, for any $y\in X$,
$0\leq h_j(y)\uparrow|\mathcal{D}_kf(y)|^2$
as $j\to\infty$,
which, combined with the monotone convergence theorem,
further implies that
\begin{align*}
\fint_{B_\rho(\cdot,2^{-k})}h_j(y)\,d\mu(y)
\uparrow
\fint_{B_\rho(\cdot,2^{-k})}|\mathcal{D}_kf(y)|^2\,d\mu(y).
\end{align*}
Thus, to show that $\fint_{B_\rho(\cdot,2^{-k})}|\mathcal{D}_kf(y)|^2\,d\mu(y)$
is $\mu$-measurable,
we only need to prove that, for any $j\in\mathbb{N}$,
$\fint_{B_\rho(\cdot,2^{-k})}h_j(y)\,d\mu(y)$ is $\mu$-measurable.
Noting that $h_j$ for any $j\in\mathbb{N}$ is a simple function,
it suffices to show that, for any $\mu$-measurable set $E\subset X$,
$$
\frac{\mu(B_\rho(\cdot,2^{-k})\cap E)}{\mu(B_\rho(\cdot,2^{-k}))}
$$
is $\mu$-measurable. To prove this,
from an argument similar to that used in the proof
Proposition~\ref{1507},
we only need to show that, for any $\mu$-measurable set $E\subset X$,
the function
$$
F(\cdot):=\mu\left(B_\rho(\cdot,2^{-k})\cap E\right)
$$
is lower semicontinuous.
If this holds, then we find that
both $\mu(B_\rho(\cdot,2^{-k})\cap E)$
and $\mu(B_\rho(\cdot,2^{-k}))$
are $\mu$-measurable
and hence
$\mu(B_\rho(\cdot,2^{-k})\cap E)/\mu(B_\rho(\cdot,2^{-k}))$
is $\mu$-measurable.

Now, we prove that $F(\cdot)$
is lower semicontinuous.
To this end, we first claim that,
for any given $x\in X$ and $\{x_j\}_{j\in\mathbb{N}}$
satisfying that $\rho(x,x_j)\to0$ as $j\to\infty$
and for any $y\in X$ and $r\in(0,\infty)$,
\begin{align}\label{11091}
\mathbf{1}_{B_\rho(x,r)}(y)\leq\liminf_{j\to\infty}
\mathbf{1}_{B_\rho(x_j,r)}(y).
\end{align}
Indeed, if $y\in[B_\rho(x,r)]^\complement$,
then \eqref{11091} holds automatically.
Let $y\in B_\rho(x,r)$.
By this and the continuity of $\rho(y,\cdot)$
[here, we used the assumption that $\rho$
is a regularized quasi-ultrametric and Lemma~\ref{L2.3}(iii)],
we conclude that, for any $j\in\mathbb{N}$ sufficiently large,
$\rho(y,x_j)<r$.
Thus, $y\in B_\rho(x_j,r)$ with $j\in\mathbb{N}$ sufficiently large,
which completes the proof of \eqref{11091}.
Consequently,
\begin{align*}
F(x)&=\int_E\mathbf{1}_{B_\rho(x,2^{-k})}(y)\,d\mu(y)\\
&\leq\int_E\liminf_{j\to\infty}
\mathbf{1}_{B_\rho(x_j,2^{-k})}(y)\quad\text{by \eqref{11091}}\\
&\leq\liminf_{j\to\infty}\int_E
\mathbf{1}_{B_\rho(x_j,2^{-k})}(y)\quad\text{by Fatou's lemma}\\
&=\liminf_{j\to\infty}F(x_j)
\end{align*}
and hence $F(\cdot)$ is lower semicontinuous.
This finishes the proof of that
$S_\theta(f)$
is $\mu$-measurable.

The proof of the measurability of $g_\lambda^*(f)$ for any
given $\lambda\in(0,\infty)$
is similar to that of $S_\theta(f)$, we omit the details, which
completes the proof of Lemma~\ref{1955}.
\end{proof}

Next, we recall the definition of the Hardy space
in terms of the Lusin area function as follows.
Recall that, for any $a\in\mathbb{R}$,
$(a)_+:=\max\{a,\,0\}$.

\begin{definition}\label{D3.1.4}
Let $(X,\mathbf{q},\mu)$ be a quasi-ultrametric space
of homogeneous type with upper dimension $n\in(0,\infty)$.
Let $p\in(\frac{n}{n+\mathrm{ind\,}(X,\mathbf{q})},\infty)$ and
$\beta,\gamma,\varepsilon\in(0,\infty)$ satisfy
\begin{align}\label{1945}
n\left(\frac{1}{p}-1\right)_+<\beta,\gamma<\varepsilon
\preceq\mathrm{ind\,}(X,\mathbf{q}),
\end{align}
where $\mathrm{ind\,}(X,\mathbf{q})$ is the same as in \eqref{index}.
The \emph{Hardy space $H^p(X)$
via the Lusin area function}
\index[words]{Hardy space via the Lusin area function}
is defined by setting
\begin{align*}
H^p(X)\index[symbols]{H@$H^p(X)$}
:=\left\{f\in\left(\mathring{{\mathcal{G}}}_0^\varepsilon(\beta,\gamma)\right)':
\|f\|_{H^p(X)}:=\left\|S(f)\right\|_{L^p(X)}<\infty\right\}.
\end{align*}
\end{definition}

\begin{remark}
Let $(X,\mathbf{q},\mu)$ be an ultra-RD-space with upper dimension $n\in(0,\infty)$
and lower smoothness index $\mathrm{ind\,}(X,\mathbf{q})$.
Let $p\in(\frac{n}{n+\mathrm{ind\,}(X,\mathbf{q})},\infty)$ and
$\beta,\gamma,\varepsilon\in(0,\infty)$ satisfy
\eqref{1945}.
Then $H^p(X)$ for any $p\in(\frac{n}{n+\mathrm{ind\,}(X,\mathbf{q})},1)$
is a quasi-Banach space and for any $p\in[1,\infty)$
is a Banach space,
which can be deduced from Theorem~\ref{Thm4.11}
for $p\in(\frac{n}{n+\mathrm{ind\,}(X,\mathbf{q})},1)$
and Theorem~\ref{Hp=Lp} for $p\in[1,\infty)$.
\end{remark}

Now, we verify that $H^p(X)$
for any $p\in(\frac{n}{n+\mathrm{ind\,}(X,\mathbf{q})},\infty)$
is independent of the choice of
homogeneous approximations
of the identity $\{\mathcal{S}_{2^{-k}}\}_{k\in\mathbb{Z}}$.

\begin{theorem}
\label{independentofATI}
Let $(X,\mathbf{q},\mu)$ be an ultra-RD-space
with upper dimension $n\in(0,\infty)$.
Fix
$p\in(\frac{n}{n+\mathrm{ind\,}(X,\mathbf{q})},\infty)$ and
$\beta,\gamma,\varepsilon\in(0,\infty)$ satisfying
\eqref{1945}.
Assume that $\{\mathcal{S}_{2^{-k}}\}_{k\in\mathbb{Z}}$
and $\{\mathcal{P}_{2^{-k}}\}_{k\in\mathbb{Z}}$
are two homogeneous approximations
of the identity of order $\varepsilon$
in Definition~\ref{DefHATI}.
For any $k\in\mathbb{Z}$, let
$\mathcal{E}_k:=\mathcal{S}_{2^{-k}}-\mathcal{S}_{2^{-k+1}}$ and
$\mathcal{Q}_k:=\mathcal{P}_{2^{-k}}-\mathcal{P}_{2^{-k+1}}$.
Let $S_\mathcal{E}$ and $S_\mathcal{Q}$ be
the Lusin area functions  associated, respectively, with
$\{\mathcal{E}_k\}_{k\in\mathbb{Z}}$ and
$\{\mathcal{Q}_k\}_{k\in\mathbb{Z}}$ as in \eqref{Def-S-f} with $\theta:=1$.
Then, for any
$f\in(\mathring{{\mathcal{G}}}_0^\varepsilon(\beta,\gamma))'$,
$$
\left\|S_\mathcal{E}(f)\right\|_{L^p(X)}
\sim\left\|S_\mathcal{Q}(f)\right\|_{L^p(X)},
$$
where the positive equivalence constants are
independent of $f$.
\end{theorem}

\begin{remark}
Let $(X,\mathbf{q},\mu)$ be an ultra-RD-space
with upper dimension $n\in(0,\infty)$,
where $\mu$ is assumed to be Borel-semiregular.
Assume that
$p\in(\frac{n}{n+\mathrm{ind\,}(X,\mathbf{q})},\infty)$ and
$\beta,\gamma,\varepsilon\in(0,\infty)$ satisfy
\eqref{1945}.
Then, by Theorem~\ref{independentofATI}, we conclude that
the Hardy space $H^p(X)$ via the Lusin
area function is independent of the choices of
homogeneous approximations of the identity in
Definition~\ref{DefHATI}. Thus, in the remainder of this monograph,
without loss of generality, we may always assume that
the Hardy space $H^p(X)$ via the Lusin area function
is defined by a fixed
homogeneous approximation of the identity
$\{\mathcal{S}_{2^{-k}}\}_{k\in\mathbb{Z}}$
satisfying Definition~\ref{DefHATI}.
\end{remark}

To show Theorem~\ref{independentofATI},
we need three technical lemmas.
The following lemma, namely the first of three, is valid
if $\mu(X)<\infty$ or $\mu(X)=\infty$, equivalently,
when $X$ is bounded or unbounded.

\begin{lemma}\label{L4.5}
Let $(X,\mathbf{q},\mu)$ be
a quasi-ultrametric space of homogeneous type
with upper dimension $n\in(0,\infty)$ and let
$\varepsilon\in(0,\infty)$.
Assume that $\rho\in\mathbf{q}$ has the property that all $\rho$-balls are
$\mu$-measurable,
$\{a^{k,v}_\tau\}_{k\in\mathbb{Z},\,
\tau\in I_k,\,v\in\{1,...,N(k,\tau)\}}$ in $\mathbb{C}$, and
$r\in(\frac{n}{n+\varepsilon},1]$.
Then there exists a positive constant $C$, independent of
$y^{k,v}_\tau\in Q^{k,v}_\tau$
and $a^{k,v}_\tau$ with $k\in\mathbb{Z}$,
$\tau\in I_k$, and $v\in\{1,...,N(k,\tau)\}$,
such that, for any $k,k'\in\mathbb{Z}$ and $x\in X$,
\begin{align*}
&\sum_{\tau\in I_{k}}\sum_{v=1}^{N(k,\tau)}
\mu(Q^{k,v}_\tau)
\frac{1}{(V_\rho)_{2^{-(k\land k')}}(x)+V_\rho(x,y^{k,v}_\tau)}
\left[\frac{2^{-(k\land k')}}
{2^{-(k\land k')}+\rho(x,y^{k,v}_\tau)}\right]^\varepsilon
\left|a^{k,v}_\tau\right|
\\
&\quad\leq C2^{[(k\land k')-k]n(1-\frac{1}{r})}
\left[\mathcal{M}_\rho\left(\sum_{\tau\in I_{k}}
\sum_{v=1}^{N(k,\tau)}
\left|a^{k,v}_\tau\right|^r
\mathbf{1}_{Q^{k,v}_\tau}\right)(x)\right]^{\frac{1}{r}},
\end{align*}
where $\mathcal{M}_\rho$ is the
Hardy--Littlewood maximal operator defined in \eqref{Mf}.
\end{lemma}

To prove Lemma~\ref{L4.5}, we first recall the
following elementary and useful inequality.

\begin{lemma}\label{discreteinequality}
Let $r\in(0,1]$
and $\{a_j\}_{j\in\mathbb{N}}$ in $\mathbb{C}$.
Then $$
\left(\sum_{j\in\mathbb{N}}\left|a_j\right|\right)^r
\leq\sum_{j\in\mathbb{N}}\left|a_j\right|^r.
$$
\end{lemma}

\begin{proof}[Proof of Lemma~\ref{L4.5}]
Without loss of generality,
we may assume that $\rho\in\mathbf{q}$
is a regularized quasi-ultrametric as in Definition~\ref{D2.3}.
Let $k,k'\in\mathbb{Z}$, $x\in X$, and
\begin{align*}
H(x):=\sum_{\tau\in I_{k}}\sum_{v=1}^{N(k,\tau)}
\mu(Q^{k,v}_\tau)
\frac{1}{(V_\rho)_{2^{-(k\land k')}}(x)+V_\rho(x,y^{k,v}_\tau)}
\left[\frac{2^{-(k\land k')}}
{2^{-(k\land k')}+\rho(x,y^{k,v}_\tau)}\right]^\varepsilon
\left|a^{k,v}_\tau\right|.
\end{align*}
Note that, by both (iv) and (v) of Lemma~\ref{DyadicSys},
Definition~\ref{D2.1}(iii),
and \eqref{upperdoub}, it is easy to show that,
for any given $\tau\in I_k$ and $v\in\{1,\ldots,N(k,\tau)\}$
and for any $z\in Q_\tau^{k,v}$,
\begin{align}\label{1045}
\mu(Q_\tau^{k,v})\sim (V_\rho)_{2^{-k}}(z),
\end{align}
where the positive equivalence constants are independent of $z$.
Observing that, from \eqref{upperdoub},
we deduce that, for any $z\in X$,
$$
(V_\rho)_{2^{-(k\land k')}}(z)\lesssim
2^{[k-(k\land k')]n}(V_\rho)_{2^{-k}}(z),
$$
which, together with Lemma~\ref{discreteinequality}
and \eqref{1045}, implies that
\begin{align*}
[H(x)]^r&\leq\sum_{\tau\in I_{k}}\sum_{v=1}^{N(k,\tau)}
\left[\mu(Q^{k,v}_\tau)\right]^r
\left[\frac{1}{(V_\rho)_{2^{-(k\land k')}}(x)+V_\rho(x,y^{k,v}_\tau)}\right]^r
\left[\frac{2^{-(k\land k')}}
{2^{-(k\land k')}+\rho(x,y^{k,v}_\tau)}\right]^{r\varepsilon}
\left|a^{k,v}_\tau\right|^r\\
&=\sum_{\tau\in I_{k}}\sum_{v=1}^{N(k,\tau)}\int_{X}
\left[\mu(Q^{k,v}_\tau)\right]^{r-1}
\left[\frac{1}{(V_\rho)_{2^{-(k\land k')}}(x)+V_\rho(x,y^{k,v}_\tau)}\right]^r
\left[\frac{2^{-(k\land k')}}
{2^{-(k\land k')}+\rho(x,y^{k,v}_\tau)}\right]^{r\varepsilon}
\\
&\quad\times
\left|a^{k,v}_\tau\right|^r\mathbf{1}_{Q_\tau^{k,v}}(z)\,d\mu(z)\\
&\lesssim2^{[(k\land k')-k]n(r-1)}\sum_{\tau\in I_{k}}\sum_{v=1}^{N(k,\tau)}\int_{X}
\left[(V_\rho)_{2^{-(k\land k')}}(z)\right]^{r-1}
\left[\frac{1}{(V_\rho)_{2^{-(k\land k')}}(x)+V_\rho(x,y^{k,v}_\tau)}\right]^r
\\
&\quad\times
\left[\frac{2^{-(k\land k')}}
{2^{-(k\land k')}+\rho(x,y^{k,v}_\tau)}\right]^{r\varepsilon}
\left|a^{k,v}_\tau\right|^r\mathbf{1}_{Q_\tau^{k,v}}(z)\,d\mu(z).
\end{align*}
By this and Lemma~\ref{A3.2}(iii), we find that
\begin{align*}
[H(x)]^r
&\lesssim2^{[(k\land k')-k]n(r-1)}\sum_{\tau\in I_{k}}\sum_{v=1}^{N(k,\tau)}\int_{X}
\left[(V_\rho)_{2^{-(k\land k')}}(z)\right]^{r-1}
\left[\frac{1}{(V_\rho)_{2^{-(k\land k')}}(x)+V_\rho(x,z)}\right]^r
\\
&\quad\times
\left[\frac{2^{-(k\land k')}}
{2^{-(k\land k')}+\rho(x,z)}\right]^{r\varepsilon}
\left|a^{k,v}_\tau\right|^r\mathbf{1}_{Q_\tau^{k,v}}(z)\,d\mu(z)
\\
&=2^{[(k\land k')-k]n(r-1)}\left\{\int_{\rho(x,z)<2^{-(k\land k')}}+\sum_{l=0}^\infty
\int_{2^l2^{-(k\land k')}\leq\rho(x,z)<2^{l+1}2^{-(k\land k')}}\right\}\\
&\quad\times
\left[(V_\rho)_{2^{-(k\land k')}}(z)\right]^{r-1}
\left[\frac{1}{(V_\rho)_{2^{-(k\land k')}}(x)+V_\rho(x,z)}\right]^r
\left[\frac{2^{-(k\land k')}}
{2^{-(k\land k')}+\rho(x,z)}\right]^{r\varepsilon}
\\
&\quad\times
\sum_{\tau\in I_{k}}\sum_{v=1}^{N(k,\tau)}
\left|a^{k,v}_\tau\right|^r\mathbf{1}_{Q_\tau^{k,v}}(z)\,d\mu(z).
\end{align*}
From this, \eqref{upperdoub},
and $r>\frac{n}{n+\varepsilon}$, we infer that
\begin{align*}
[H(x)]^r
&\lesssim2^{[(k\land k')-k]n(r-1)}
\Bigg\{\mathcal{M}_\rho\left(\sum_{\tau\in I_{k}}
\sum_{v=1}^{N(k,\tau)}
\left|a^{k,v}_\tau\right|^r
\mathbf{1}_{Q^{k,v}_\tau}\right)(x)
\\
&\quad+\sum_{l=0}^\infty2^{l[n-(n+r)\varepsilon]}
\fint_{\rho(x,z)<2^{l+1}2^{-(k\land k')}}
\sum_{\tau\in I_{k}}\sum_{v=1}^{N(k,\tau)}
\left|a^{k,v}_\tau\right|^r\mathbf{1}_{Q_\tau^{k,v}}(z)\,d\mu(z)
\Bigg\}\\
&\lesssim2^{[(k\land k')-k]n(r-1)}
\mathcal{M}_\rho\left(\sum_{\tau\in I_{k}}
\sum_{v=1}^{N(k,\tau)}
\left|a^{k,v}_\tau\right|^r
\mathbf{1}_{Q^{k,v}_\tau}\right)(x).
\end{align*}
This finishes the proof of Lemma~\ref{L4.5}.
\end{proof}

The second lemma is also valid if $\mu(X)<\infty$ or $\mu(X)=\infty$,
equivalently, when $X$ is bounded or unbounded.

\begin{lemma}\label{1122}
Let $(X,\rho,\mu)$ be
a quasi-ultrametric space of homogeneous type and
$\beta,\gamma\in(0,\infty)$.
For any $k,l\in\mathbb{Z}$, let
$\overline{\mathcal{D}}_k$ and $\widetilde{\mathcal{D}}_l$
be two integral operators with their kernels
denoted, respectively, by
$\overline{D}_k$ and $\widetilde{D}_l$ which
satisfy that, for any
$x,z\in X$,
$\overline{D}_k(x,\cdot)\in\mathring{\mathcal{G}}(x,2^{-k},\beta,\gamma)$
and $\widetilde{D}_l(\cdot,z)\in\mathring{\mathcal{G}}(z,2^{-l},\beta,\gamma)$.
Then, for any given $\varepsilon\in(0,\beta)\cap(0,\gamma]$,
there exists a positive constant $C$ such that,
for any $k,l\in\mathbb{Z}$ and $x,z\in X$,
\begin{align}\label{DkDl}
\left|\overline{D}_k\widetilde{D}_l(x,z)\right|
\leq C2^{-|k-l|\varepsilon}
\frac{1}{(V_\rho)_{2^{-(k\land l)}}(x)+V_\rho(x,z)}
\left[\frac{2^{-(k\land l)}}
{2^{-(k\land l)}+\rho(x,z)}\right]^\gamma,
\end{align}
where $\overline{D}_k\widetilde{D}_l$ is the kernel
of the integral operator $\overline{\mathcal{D}}_k\widetilde{\mathcal{D}}_l$.
\end{lemma}

\begin{proof}
Without loss of generality, we may assume that
$\rho$ is a regularized quasi-ultrametric as in Definition~\ref{D2.3}.
Fix $k,l\in\mathbb{Z}$ and $x,z\in X$.
Using symmetry, we find that
it suffices to prove \eqref{DkDl} in the case $k\ge l$.
By the cancellation of $\overline{D}_k(z,\cdot)$
for the second variable,
we conclude that
\begin{align}\label{J123}
\left|\overline{D}_k\widetilde{D}_l(x,z)\right|
&=\left|\int_X\overline{D}_k(z,y)
\left[\widetilde{D}_l(y,x)-\widetilde{D}_l(z,x)
\right]\,d\mu(y)\right|
\nonumber\\
&\leq\int_{\rho(z,y)\leq\frac{2^{-l}+\rho(x,z)}{2C_\rho}}
\left|\overline{D}_k(z,y)\right|
\left|\widetilde{D}_l(y,x)-\widetilde{D}_l(z,x)
\right|\,d\mu(y)
\nonumber\\
&\quad+\int_{\rho(z,y)>\frac{2^{-l}+\rho(x,z)}{2C_\rho}}
\left|\overline{D}_k(z,y)\right|
\left|\widetilde{D}_l(y,x)\right|\,d\mu(y)
\nonumber\\
&\quad+\left|\widetilde{D}_l(z,x)\right|
\int_{\rho(z,y)>\frac{2^{-l}+\rho(x,z)}{2C_\rho}}
\left|\overline{D}_k(z,y)\right|\,d\mu(y)
\nonumber\\
&=:J_1(x,z)+J_2(x,z)+J_3(x,z).
\end{align}
To deal with $J_1(x,z)$,
from the size condition of $\overline{D}_k$,
the regularity condition of $\widetilde{D}_l$,
$\varepsilon\in(0,\beta)\cap(0,\gamma]$,
and Lemma~\ref{A3.2}(i),
we deduce that
\begin{align}\label{J153}
J_1(x,z)
&\lesssim\int_{\rho(z,y)\leq\frac{2^{-l}+\rho(x,z)}{2C_\rho}}
\frac{1}{(V_\rho)_{2^{-k}}(z)+V_\rho(z,y)}
\left[\frac{2^{-k}}{2^{-k}+\rho(z,y)}\right]^\gamma
\nonumber\\
&\quad\times
\left[\frac{\rho(z,y)}{2^{-l}+\rho(x,z)}\right]^\beta
\frac{1}{(V_\rho)_{2^{-l}}(x)+V_\rho(x,z)}
\left[\frac{2^{-l}}{2^{-l}+\rho(x,z)}\right]^\gamma
\,d\mu(y)
\nonumber\\
&\lesssim\int_{\rho(z,y)\leq\frac{2^{-l}+\rho(x,z)}{2C_\rho}}
\frac{1}{V_\rho(z,y)}
\left[\frac{2^{-k}}{\rho(z,y)}\right]^{\varepsilon}
\left[\frac{\rho(z,y)}{2^{-l}}\right]^{\varepsilon}
\left[\frac{\rho(z,y)}{2^{-l}+\rho(x,z)}\right]^{\beta-\varepsilon}
\,d\mu(y)
\nonumber\\
&\quad\times
\frac{1}{(V_\rho)_{2^{-l}}(x)+V_\rho(x,z)}
\left[\frac{2^{-l}}{2^{-l}+\rho(x,z)}\right]^\gamma
\nonumber\\
&\lesssim2^{(l-k)\varepsilon}
\frac{1}{(V_\rho)_{2^{-l}}(x)+V_\rho(x,z)}
\left[\frac{2^{-l}}{2^{-l}+\rho(x,z)}\right]^\gamma.
\end{align}
For $J_2(x,z)$, by
$\rho(z,y)>\frac{2^{-l}+\rho(x,z)}{2C_\rho}$,
\eqref{upperdoub},
Lemma~\ref{A3.2}(vi),
and $k\ge l$, we find that
\begin{align*}
V_\rho(z,y)\ge (V_\rho)_{\frac{2^{-l}+\rho(z,x)}{2C_\rho}}(z)
\gtrsim (V_\rho)_{2^{-l}+\rho(z,x)}(z)\sim
(V_\rho)_{2^{-l}}(x)+V_\rho(x,z)
\end{align*}
and $\rho(z,y)\gtrsim2^{-l}+\rho(z,x)$.
From these, the size conditions of
$\overline{D}_k$ and $\widetilde{D}_l$,
Lemma~\ref{A3.2}(ii),
and the fact that $(l-k)(\varepsilon-\gamma)\ge0$,
we infer that
\begin{align}\label{J253}
J_2(x,z)
&\lesssim\int_{\rho(z,y)>\frac{2^{-l}+\rho(x,z)}{2C_\rho}}
\frac{1}{(V_\rho)_{2^{-k}}(z)+V_\rho(z,y)}
\left[\frac{2^{-k}}{2^{-k}+\rho(z,y)}\right]^\gamma
\nonumber\\
&\quad\times
\frac{1}{(V_\rho)_{2^{-l}}(x)+V_\rho(x,y)}
\left[\frac{2^{-l}}{2^{-l}+\rho(x,y)}\right]^\gamma
\,d\mu(y)
\nonumber\\
&\lesssim2^{(l-k)\gamma}\int_{X}
\frac{1}{(V_\rho)_{2^{-l}}(x)+V_\rho(x,z)}
\left[\frac{2^{-l}}{2^{-l}+\rho(z,x)}\right]^\gamma
\nonumber\\
&\quad\times
\frac{1}{(V_\rho)_{2^{-l}}(x)+V_\rho(x,y)}
\left[\frac{2^{-l}}{2^{-l}+\rho(x,y)}\right]^\gamma
\,d\mu(y)
\nonumber\\
&\lesssim2^{(l-k)\varepsilon}
\frac{1}{(V_\rho)_{2^{-l}}(x)+V_\rho(x,z)}
\left[\frac{2^{-l}}{2^{-l}+\rho(x,z)}\right]^\gamma.
\end{align}
To deal with $J_3(x,z)$, by the size conditions of both
$\overline{D}_k$ and $\widetilde{D}_l$,
Lemma~\ref{A3.2}(i), and the fact that $(l-k)(\varepsilon-\gamma)\ge0$,
we conclude that
\begin{align*}
J_3(x,z)
&\lesssim\frac{1}{(V_\rho)_{2^{-l}}(x)+V_\rho(x,z)}
\left[\frac{2^{-l}}{2^{-l}+\rho(x,z)}\right]^\gamma
\\
&\quad\times\int_{\rho(z,y)>\frac{2^{-l}+\rho(x,z)}{2C_\rho}}
\frac{1}{(V_\rho)_{2^{-k}}(y)+V_\rho(y,z)}
\left[\frac{2^{-k}}{2^{-k}+\rho(y,z)}\right]^\gamma
\,d\mu(y)
\\
&\leq\frac{1}{(V_\rho)_{2^{-l}}(x)+V_\rho(x,z)}
\left[\frac{2^{-l}}{2^{-l}+\rho(x,z)}\right]^\gamma
\int_{\rho(z,y)>\frac{2^{-k}}{2C_\rho}}
\frac{1}{V_\rho(y,z)}\left[\frac{2^{-l}}{\rho(y,z)}\right]^\gamma\,d\mu(y)
\\
&\lesssim2^{(l-k)\gamma}
\frac{1}{(V_\rho)_{2^{-l}}(x)+V_\rho(x,z)}
\left[\frac{2^{-l}}{2^{-l}+\rho(x,z)}\right]^\gamma
\\
&\leq2^{(l-k)\varepsilon}
\frac{1}{(V_\rho)_{2^{-l}}(x)+V_\rho(x,z)}
\left[\frac{2^{-l}}{2^{-l}+\rho(x,z)}\right]^\gamma,
\end{align*}
which, combined with \eqref{J123}, \eqref{J153},
and \eqref{J253}, completes the proof of \eqref{DkDl}
and hence Lemma~\ref{1122}.
\end{proof}

The third and final lemma
is the Fefferman--Stein vector-valued maximal inequality
(see, for instance, \cite[Theorem 1.2]{FS}),
\index[words]{Fefferman--Stein vector-valued maximal inequality}
which is valid if $\mu(X)<\infty$ or $\mu(X)=\infty$, equivalently,
when $X$ is bounded or unbounded.

\begin{lemma}\label{Fefferman--Stein}
Let $(X,\mathbf{q},\mu)$ be a quasi-ultrametric
space of homogeneous type.
Assume that $\rho\in\mathbf{q}$ is a regularized quasi-ultrametric
as in Definition~\ref{D2.3},
$p\in (1,\infty)$, and $u\in(1,\infty]$.
Then there exists a positive constant $C$
such that, for any sequence
of $\mu$-measurable functions $\{f_j\}_{j=1}^\infty$,
$$
\left\|\left[\sum_{j=1}^{\infty}
\left(\mathcal{M}_\rho f_j\right)^u\right]^\frac{1}{u}\right\|_{L^p(X)}
\leq C\left\|\left(\sum_{j=1}^{\infty}
|f_j|^u\right)^\frac{1}{u}\right\|_{L^p(X)}
$$
with the usual modification when $u=\infty$.
\end{lemma}

Next, we show Theorem~\ref{independentofATI}.

\begin{proof}[Proof of Theorem~\ref{independentofATI}]
Let $\rho\in\mathbf{q}$ be a regularized quasi-ultrametric.
Let $f\in(\mathring{{\mathcal{G}}}_0^\varepsilon(\beta,\gamma))'$.
Using symmetry, we only need to prove that
\begin{align}\label{2105}
\left\|S_\mathcal{E}(f)\right\|_{L^p(X)}
\lesssim\left\|S_\mathcal{Q}(f)\right\|_{L^p(X)},
\end{align}
where the implicit positive constant is independent of $f$.
For any $k\in\mathbb{Z}$ and $z\in X$, define
\begin{align}
\label{mkf}
m_k(f)(z)
:=\left[\frac{1}{(V_\rho)_{2^{-k}}(z)}
\int_{B_\rho(z,2^{-k})}
\left|\mathcal{Q}_kf(u)\right|^2\,d\mu(u)\right]^{\frac{1}{2}}.
\end{align}
Suppose that $l\in\mathbb{Z}$,
$x\in X$, and $y\in B_\rho(x,2^{-l})$.
By Theorem~\ref{HD} (keeping in mind
that we are assuming $\mathrm{diam}_\rho(X)=\infty$),
we find that there exists a family of linear operators
$\{\widetilde{\mathcal{Q}}_k\}_{k\in\mathbb{Z}}$ such that the second
equality in \eqref{3.29} holds in
$(\mathring{{\mathcal{G}}}_0^\varepsilon(\beta,\gamma))'$.
As such, we conclude that
\begin{align}\label{Elf}
\mathcal{E}_lf(y)&=\left\langle f,E_l(y,\cdot)\right\rangle
=\left\langle\sum_{k=-\infty}^{\infty}\sum_{\tau\in I_{k}}
\sum_{v=1}^{N(k,\tau)}
\widetilde{Q}_k(\cdot,y_\tau^{k,v})
\int_{Q_\tau^{k,v}}\mathcal{Q}_kf(u)\,d\mu(u),E_l(y,\cdot)\right\rangle
\nonumber\\
&=\sum_{k=-\infty}^{\infty}\sum_{\tau\in I_{k}}
\sum_{v=1}^{N(k,\tau)}\left\langle
\widetilde{Q}_k(\cdot,y_\tau^{k,v})
,E_l(y,\cdot)\right\rangle\int_{Q_\tau^{k,v}}
\mathcal{Q}_kf(u)\,d\mu(u)
\nonumber\\
&=\sum_{k=-\infty}^{\infty}\sum_{\tau\in I_{k}}
\sum_{v=1}^{N(k,\tau)}
E_l\widetilde{Q}_k(y,y_\tau^{k,v})
\int_{Q_\tau^{k,v}}\mathcal{Q}_kf(u)\,d\mu(u),
\end{align}
where the penultimate equality holds because the
series converge in $(\mathring{{\mathcal{G}}}_0^\varepsilon(\beta,\gamma))'$
and where all of the notation are as in Theorem~\ref{A3.47}.

Now, we estimate $\mathcal{E}_lf(y)$.
From Lemma~\ref{DyadicSys}(iv) and \eqref{2cubej},
we infer that, for any given $z\in Q_\tau^{k,v}$ and for any
$u\in Q_\tau^{k,v}$,
\begin{align*}
\rho(z,u)\leq\mathrm{diam}_\rho(Q_\tau^{k,v})
\leq C_r2^{-(k+j)}<2^{-k},
\end{align*}
where $C_r$ and $j\in\mathbb{N}$ are the same positive constants
as, respectively, in Lemma~\ref{DyadicSys}(iv) and \eqref{2cubej},
and	hence
\begin{align}\label{1041}
Q_\tau^{k,v}\subset B_\rho(z,2^{-k}).
\end{align}
By this, Definition~\ref{D2.1}(iii), and Lemma~\ref{DyadicSys}(v),
we find that, for any $z\in Q_\tau^{k,v}$,
\begin{align*}
B_\rho(z,2^{-k})
\subset B_\rho(z_\tau^{k,v},C_\rho2^{-k}),
\end{align*}
where $z_\tau^{k,v}$ denotes the center of $Q_\tau^{k,v}$,
which, together with \eqref{upperdoub}, further implies that
\begin{align*}
\mu\left(B_\rho(z,2^{-k})\right)
\leq\mu\left(B_\rho(z_\tau^{k,v},C_\rho2^{-k})\right)
\lesssim\mu\left(B_\rho(z_\tau^{k,v},C_l2^{-k})\right)
\leq\mu(Q_\tau^{k,v}),
\end{align*}
where $C_l\in(0,\infty)$ is as in
Lemma~\ref{DyadicSys}(v) and the implicit positive constant depends only
on $C_l$, $C_\rho$, and $C_\mu$.
From this, we deduce that, for any $z\in Q_\tau^{k,v}$,
\begin{align*}
&\left|\frac{1}{\mu(Q_\tau^{k,v})}
\int_{Q_\tau^{k,v}}\mathcal{Q}_kf(u)\,d\mu(u)\right|\\
&\quad\leq\left[\frac{1}{\mu(Q_\tau^{k,v})}
\int_{Q_\tau^{k,v}}
\left|\mathcal{Q}_kf(u)\right|^2\,d\mu(u)\right]^{\frac{1}{2}}
\quad\text{by H\"older's inequality}\\
&\quad\lesssim\left[\frac{1}{(V_\rho)_{2^{-k}}(z)}
\int_{B_\rho(z,2^{-k})}
\left|\mathcal{Q}_kf(u)\right|^2\,d\mu(u)\right]^{\frac{1}{2}}
\quad\text{by \eqref{1041}},
\end{align*}
which, combined with the arbitrariness
of $z\in Q_\tau^{k,v}$, further implies that
\begin{align}\label{Elf1}
\left|\frac{1}{\mu(Q_\tau^{k,v})}
\int_{Q_\tau^{k,v}}
\mathcal{Q}_kf(u)\,d\mu(u)\right|
\lesssim\inf_{z\in Q_\tau^{k,v}}m_k(f)(z).
\end{align}
Moreover, by Lemma~\ref{1122}, we conclude that,
for any fixed
$\varepsilon'\in(n(\frac{1}{p}-1)_+,\varepsilon)$
and $y_\tau^{k,v}\in Q_\tau^{k,v}$,
\begin{align*}
\left|E_l\widetilde{Q}_k(y,y_\tau^{k,v})\right|
&\lesssim2^{-|k-l|\varepsilon'}
\frac{1}{(V_\rho)_{2^{-(k\land l)}}(y)
+V_\rho(y,y_\tau^{k,v})}
\left[\frac{2^{-(k\land l)}}
{2^{-(k\land l)}+\rho(y,y_\tau^{k,v})}\right]^{\varepsilon'},
\end{align*}
which, together with \eqref{Elf}, \eqref{Elf1},
and Lemma~\ref{L4.5}, further implies that, for any
\begin{align}\label{2202}
r\in\left(\frac{n}{n+\varepsilon'},1\right]\cap(0,p),
\end{align}
we have
\begin{align*}
\left|\mathcal{E}_lf(y)\right|
&\lesssim\sum_{k=-\infty}^{\infty}
2^{-|k-l|\varepsilon'}\sum_{\tau\in I_k}
\sum_{v=1}^{N(k,\tau)}
\mu(Q_\tau^{k,v})
\\
&\quad\times
\frac{1}{(V_\rho)_{2^{-(k\land l)}}(y)
+V_\rho(y,y_\tau^{k,v})}
\left[\frac{2^{-(k\land l)}}
{2^{-(k\land l)}+\rho(y,y_\tau^{k,v})}\right]^{\varepsilon'}
\inf_{z\in Q_\tau^{k,v}}
m_k(f)(z)
\\
&\lesssim\sum_{k=-\infty}^{\infty}
2^{-|k-l|\varepsilon'+n(k\land l-k)(1-\frac{1}{r})}
\left[\mathcal{M}_\rho
\left(\sum_{\tau\in I_k}\sum_{v=1}^{N(k,\tau)}
\inf_{z\in Q_\tau^{k,v}}
\left[m_k(f)(z)\right]^r
\mathbf{1}_{Q_\tau^{k,v}}\right)(x)\right]^{\frac{1}{r}}.
\end{align*}
From this, H\"older's inequality, and \eqref{2202},
we infer that, for any $x\in X$,
\begin{align*}
\left[S_\mathcal{E}(f)(x)\right]^2
&\lesssim\sum_{l=-\infty}^{\infty}
\left\{\sum_{k=-\infty}^{\infty}
2^{-|k-l|\varepsilon'+n(k\land l-k)(1-\frac{1}{r})}
\left[\mathcal{M}_\rho
\left(\sum_{\tau\in I_k}\sum_{v=1}^{N(k,\tau)}
\inf_{z\in Q_\tau^{k,v}}
\left[m_k(f)(z)\right]^r
\mathbf{1}_{Q_\tau^{k,v}}\right)(x)\right]^{\frac{1}{r}}\right\}^2
\nonumber\\
&\lesssim\sum_{l=-\infty}^{\infty}
\sum_{k=-\infty}^{\infty}
2^{-|k-l|\varepsilon'+n(k\land l-k)(1-\frac{1}{r})}
\left[\mathcal{M}_\rho
\left(\sum_{\tau\in I_k}\sum_{v=1}^{N(k,\tau)}
\inf_{z\in Q_\tau^{k,v}}
\left[m_k(f)(z)\right]^r
\mathbf{1}_{Q_\tau^{k,v}}\right)(x)\right]^{\frac{2}{r}},
\end{align*}
which, combined with Tonelli's theorem,
further implies that
\begin{align}\label{817}
\left[S_\mathcal{E}(f)(x)\right]^2
&\lesssim\sum_{k=-\infty}^{\infty}
\left[\mathcal{M}_\rho
\left(\sum_{\tau\in I_k}\sum_{v=1}^{N(k,\tau)}
\inf_{z\in Q_\tau^{k,v}}
\left[m_k(f)(z)\right]^r
\mathbf{1}_{Q_\tau^{k,v}}\right)(x)\right]^{\frac{2}{r}}
\nonumber\\
&\leq\sum_{k=-\infty}^{\infty}
\left[\mathcal{M}_\rho
\left([m_k(f)]^r\right)(x)\right]^{\frac{2}{r}}.
\end{align}
This, together with Lemma~\ref{Fefferman--Stein}
with $p$ and $u$ therein replaced, respectively,
by $\frac{p}{r}$ and $\frac{2}{r}$,
implies that
\begin{align*}
\left\|S_\mathcal{E}(f)\right\|_{L^p(X)}
&=\left\|\left[S_\mathcal{E}(f)\right]^r\right\|_{L^\frac{p}{r}(X)}^\frac{1}{r}
\lesssim\left\|\left\{\sum_{k=-\infty}^{\infty}
\left[\mathcal{M}_\rho([m_k(f)]^r)\right]^{\frac{2}{r}}
\right\}^{\frac{r}{2}}\right\|^{\frac{1}{r}}_{L^\frac{p}{r}(X)}
\\
&\lesssim\left\|\left\{\sum_{k=-\infty}^{\infty}
[m_k(f)]^2\right\}^{\frac{1}{2}}\right\|_{L^p(X)}
=\left\|S_\mathcal{Q}(f)\right\|_{L^p(X)},
\end{align*}
which completes the proof of \eqref{2105}
and hence Theorem~\ref{independentofATI}.
\end{proof}

The other main theorem of this section reads as follows,
in which we consider the case $p\in(1,\infty)$.
The reader is reminded of Proposition~\ref{distfncttype}.

\begin{theorem}\label{Hp=Lp}
Let $(X,\mathbf{q},\mu)$ be an ultra-RD-space.
Let $p\in(1,\infty)$ and
$\beta,\gamma,\varepsilon\in(0,\infty)$ satisfy
\eqref{1945}.
Then $H^p(X)=L^p(X)$ in the following sense:
the linear map
$\Phi:L^p(X)\to H^p(X)$ that is defined by setting, for any $f\in L^p(X)$,
$\Phi(f):=\Lambda_f$, where
$\Lambda_f\in({\mathcal{G}}^\varepsilon_0(\beta,\gamma))'$
is the distribution defined in Proposition~\ref{distfncttype},
is well defined, linear, bijective,
and has the property that, for any
$f\in L^p(X)$,
$$
\|f\|_{L^p(X)}\sim\left\|\Lambda_f\right\|_{H^p(X)},
$$
where the positive equivalence constants are independent of $f$.
\end{theorem}

\begin{remark}
Similarly to that in Proposition~\ref{distfncttype},
in the context of Theorem~\ref{Hp=Lp}
since the association of $f\in L^p(X)$ to the distribution
$\Lambda_f\in({\mathcal{G}}^\varepsilon_0(\beta,\gamma))'$ is
injective whenever $\mu$ is Borel-semiregular,
we simply write $f$ in place of $\Lambda_f$ in what follows.
\end{remark}

To show Theorem~\ref{Hp=Lp},
we need the following lemma.

\begin{lemma}
\label{gL2L2}
Let $(X,\mathbf{q},\mu)$ be an ultra-RD-space.
Let $p\in(1,\infty)$ and
$\beta,\gamma,\varepsilon\in(0,\infty)$ satisfy
\eqref{1945}.
Assume that $\{\mathcal{S}_{2^{-k}}\}_{k\in\mathbb{Z}}$
is a homogeneous approximation
of the identity of order $\varepsilon$
in Definition~\ref{DefHATI}
and let
$\mathcal{D}_k:=\mathcal{S}_{2^{-k}}-\mathcal{S}_{2^{-k+1}}$
for any $k\in\mathbb{Z}$. Let
$g(f)$ for any $f\in(\mathring{\mathcal{G}}^\varepsilon_0(\beta,\gamma))'$
be as in \eqref{Def-g(f)(x)}.
Then there
exists a constant $C\in[1,\infty)$ such that
\begin{enumerate}
\item[\textup{(i)}]
For any $f\in L^p(X)$,
$\|g(f)\|_{L^p(X)}\leq C\|f\|_{L^p(X)}$.

\item[\textup{(ii)}]
For any $f\in (\mathring{\mathcal{G}}^\varepsilon_0(\beta,\gamma))'$
satisfying $g(f)\in L^p(X)$, there exists $u\in L^p(X)$ such that
$f=u$ in $(\mathring{\mathcal{G}}^\varepsilon_0(\beta,\gamma))'$ and
$\|u\|_{L^p(X)}\leq C\|g(f)\|_{L^p(X)}$.
\end{enumerate}
\end{lemma}

To prove Lemma~\ref{gL2L2},
we need the following well-known Khinchin inequality;\index[words]{Khinchin inequality}
see, for instance, \cite[p.\,141, (2.2)]{m1992}.

\begin{lemma}\label{KhinchinInequality}
There exists $C\in(0,\infty)$
such that, for any $N\in\mathbb{N}$ and $\{a_j\}_{j=1}^N$ in $\mathbb{C}$,
\begin{align*}
\left(\sum_{j=1}^N\left|a_j\right|^2\right)^\frac{1}{2}
\leq C2^{-N}\sum_{\varepsilon_1\in \{1,\,-1\}}\ldots
\sum_{\varepsilon_N\in \{1,\,-1\}}\left|\sum_{j=1}^N\varepsilon_ja_j\right|.
\end{align*}
\end{lemma}

Next, we show Lemma~\ref{gL2L2}.

\begin{proof}[Proof of Lemma~\ref{gL2L2}]
We first prove (i). Let $f\in L^p(X)$.
Using Lemma~\ref{KhinchinInequality} and Minkowski's inequality,
we find that, for any $N\in\mathbb{N}$,
\begin{align}\label{1709}
\left\|\left(\sum_{k=-N}^N
\left|\mathcal{D}_kf\right|^2\right)^\frac{1}{2}\right\|_{L^p(X)}
&\lesssim\left\|2^{-2N}\sum_{\varepsilon_1\in \{1,\,-1\}}\ldots
\sum_{\varepsilon_N\in \{1,\,-1\}}\left|\sum_{k=-N}^N
\varepsilon_k\mathcal{D}_kf\right|\,\right\|_{L^p(X)}\nonumber\\
&\leq2^{-2N}\sum_{\varepsilon_1\in \{1,\,-1\}}\ldots
\sum_{\varepsilon_N\in \{1,\,-1\}}\left\|\sum_{k=-N}^N
\varepsilon_k\mathcal{D}_kf\right\|_{L^p(X)}.
\end{align}
For any fixed $N\in\mathbb{N}$ and $\vec{\varepsilon}:=
\{\varepsilon_k\}_{k=-N}^N$ in $\{1,\,-1\}$,
let
$$
\mathcal{T}_N^{\vec{\varepsilon}}:=\sum_{k=-N}^N\varepsilon_k\mathcal{D}_k
$$
and $K_N^{\vec{\varepsilon}}$ be the kernel of $\mathcal{T}_N^{\vec{\varepsilon}}$.

Now, we claim that $K_N^{\vec{\varepsilon}}$
is a standard Calder\'on--Zygmund kernel as in
Definition~\ref{CZk} with the positive constant $C_T$
in both \eqref{CZk1} and \eqref{CZk2} independent of
both $N$ and $\vec{\varepsilon}$.

To verify the size estimate of $K_N^{\vec{\varepsilon}}(x,y)$,
from Lemmas~\ref{Pk}(i) and~\ref{A3.14}
and an argument similar to that used in the estimation
of \eqref{RN-1},
we deduce that, for any $x,y\in X$ with $x\neq y$,
\begin{align}\label{16411}
\left|K_N^{\vec{\varepsilon}}(x,y)\right|
&\leq\sum_{k=-N}^N\left|D_k(x,y)\right|\lesssim\frac{1}{V_\rho(x,y)}.
\end{align}

To verify the regularity of $K_N^{\vec{\varepsilon}}(x,y)$,
by Lemmas~\ref{Pk}(ii) and~\ref{A3.14}
and an argument similar to that used in the estimation
of \eqref{RN1x},
we conclude that, for any $x,x',y\in X$ with both $x\neq y$
and $\rho(x,x')\leq\frac{\rho(x,y)}{2C_\rho}$,
\begin{align*}
\left|K_N^{\vec{\varepsilon}}(x,y)-K_N^{\vec{\varepsilon}}(x',y)\right|
&\leq\sum_{k=-N}^N\left|D_k(x,y)-D_k(x',y)\right|
\lesssim\left[\frac{\rho(x,x')}{\rho(x,y)}\right]^\varepsilon\frac{1}{V_\rho(x,y)},
\end{align*}
which, combined with \eqref{16411},
completes the proof of the above claim.

From an argument similar to that used in the estimation of \eqref{RN0},
we infer that, for any $k_1,k_2\in\mathbb{Z}$,
\begin{align*}
\left\|\varepsilon_{k_1}\mathcal{D}_{k_1}
\left(\varepsilon_{k_2}\mathcal{D}_{k_2}\right)^*
\right\|_{L^2(X)\to L^2(X)}\lesssim2^{-|k_1-k_2|\varepsilon}
\end{align*}
and
\begin{align*}
\left\|\left(\varepsilon_{k_1}\mathcal{D}_{k_1}\right)^*
\varepsilon_{k_2}\mathcal{D}_{k_2}\right\|_{L^2(X)\to L^2(X)}
\lesssim2^{-|k_1-k_2|\varepsilon},
\end{align*}
which, together with Lemma~\ref{CotlarStein},
further implies that, for any $f\in L^2(X)$,
\begin{align*}
\left\|\mathcal{T}_N^{\vec{\varepsilon}}f\right\|_{L^2(X)}
\lesssim\left\|f\right\|_{L^2(X)},
\end{align*}
where the implicit positive constant is independent of $\vec{\varepsilon}$,
$N$, and $f$.
By this, \cite[p.\,22, Corollary]{Pag}, and
a standard duality argument, we find that
$\mathcal{T}_N^{\vec{\varepsilon}}$ is bounded on $L^p(X)$.
This and \eqref{1709} imply that, for any $f\in L^p(X)$,
\begin{align*}
\left\|\left(\sum_{k=-N}^N\left|\mathcal{D}_kf\right|^2
\right)^\frac{1}{2}\right\|_{L^p(X)}\lesssim\|f\|_{L^p(X)},
\end{align*}
where the implicit positive constant is independent of both
$N$ and $f$, which, combined with Fatou's lemma, further implies that
\begin{align*}
\left\|g(f)\right\|_{L^p(X)}=\left\|\left(\sum_{k=-\infty}^\infty
\left|\mathcal{D}_kf\right|^2\right)^\frac{1}{2}\right\|_{L^p(X)}
\lesssim\|f\|_{L^p(X)}.
\end{align*}
This finishes the proof of Lemma~\ref{gL2L2}(i).

Next, we turn to show Lemma~\ref{gL2L2}(ii). To this end,
suppose that $f\in (\mathring{\mathcal{G}}^\varepsilon_0(\beta,\gamma))'$
satisfies $g(f)\in L^p(X)$. From
Theorem~\ref{HC} [to be precise, the first equality of \eqref{3.12} holds
in $(\mathring{\mathcal{G}}^\varepsilon_0(\beta,\gamma))'$],
we deduce that there exists a family
of linear operators $\{\widetilde{\mathcal{D}}_k\}_{k\in\mathbb{Z}}$ such that,
for any $g\in (\mathring{\mathcal{G}}^\varepsilon_0(\beta,\gamma))'$,
\begin{align}\label{92611}
g=\sum_{k=-\infty}^\infty\widetilde{\mathcal{D}}_k
\mathcal{D}_kg
\end{align}
in $(\mathring{\mathcal{G}}^\varepsilon_0(\beta,\gamma))'$.
Moreover, by Theorem~\ref{A3.19},
we conclude that, if $g\in L^r(X)$ for some $r\in(1,\infty)$, then
the convergence in \eqref{92611} also holds pointwise $\mu$-almost everywhere on $X$.
For any $N\in\mathbb{N}$, let
$$
f_N:=\sum_{k=-N}^N\widetilde{\mathcal{D}}_k
\mathcal{D}_kf.
$$

Now, we claim that, for any $h\in L^{p'}(X)$,
\begin{align}\label{1029}
\left\|\left[\sum_{k=-\infty}^\infty\left|
\left(\widetilde{\mathcal{D}}_k\right)^*h\right|^2
\right]^\frac{1}{2}\right\|_{L^{p'}(X)}
\lesssim\left\|h\right\|_{L^{p'}(X)},
\end{align}
where $p'\in(1,\infty)$ is the conjugate index of $p$.
To this end, let $h\in L^{p'}(X)$.
Lemma~\ref{1122} then gives that, for any
given $\varepsilon'\in(0,\beta)\cap(0,\gamma]$ and for any $k,l\in\mathbb{Z}$,
\begin{align*}
\left|\left(\widetilde{D}_k\right)^*\widetilde{D}_l\left(x,y\right)\right|
\lesssim2^{-|k-l|\varepsilon'}
\frac{1}{(V_\rho)_{2^{-(k\land l)}}(x)+V_\rho(x,y)}
\left[\frac{2^{-(k\land l)}}
{2^{-(k\land l)}+\rho(x,y)}\right]^\gamma,
\end{align*}
which, together with Lemma~\ref{A3.2}(iv),
further implies that, for any $g\in L^{p'}(X)$ and $x\in X$,
\begin{align*}
\left|\left(\widetilde{\mathcal{D}}_k\right)^*\widetilde{\mathcal{D}}_lg(x)\right|
&\lesssim2^{-|k-l|\varepsilon'}
\int_X\frac{1}{(V_\rho)_{2^{-(k\land l)}}(x)+V_\rho(x,y)}
\left[\frac{2^{-(k\land l)}}
{2^{-(k\land l)}+\rho(x,y)}\right]^\gamma|g(y)|\,d\mu(y)\\
&\lesssim2^{-|k-l|\varepsilon'}\fint_{B_\rho(x,2^{-(k\land l)})}|g(y)|\,d\mu(y)\\
&\quad+2^{-|k-l|\varepsilon'}
\int_{(B_\rho(x,2^{-(k\land l)}))^\complement}\frac{1}{V_\rho(x,y)}
\left[\frac{2^{-(k\land l)}}
{\rho(x,y)}\right]^\gamma|g(y)|\,d\mu(y)\\
&\lesssim2^{-|k-l|\varepsilon'}\mathcal{M}_\rho g(x).
\end{align*}
From this [applied to $g=\mathcal{D}_lh\in L^{p'}(X)$; see Theorem~\ref{corATI}(i)] and the pointwise convergence in \eqref{92611}
[applied to $g=h\in L^{p'}(X)$],
we infer that
\begin{align*}
\left\|\left[\sum_{k=-\infty}^\infty\left|
\left(\widetilde{\mathcal{D}}_k\right)^*h\right|^2
\right]^\frac{1}{2}\right\|_{L^{p'}(X)}
&=\left\|\left[\sum_{k=-\infty}^\infty\left|\sum_{l=-\infty}^\infty
\left(\widetilde{\mathcal{D}}_k\right)^*\widetilde{\mathcal{D}}_l\mathcal{D}_lh\right|^2
\right]^\frac{1}{2}\right\|_{L^{p'}(X)}\\
&\lesssim\left\|\left\{\sum_{k=-\infty}^\infty
\left[\sum_{l=-\infty}^\infty2^{-|k-l|\varepsilon'}
\mathcal{M}_\rho(\mathcal{D}_lh)\right]^2
\right\}^\frac{1}{2}\right\|_{L^{p'}(X)},
\end{align*}
which, combined with
Tonelli's theorem, further implies that
\begin{align*}
&\left\|\left[\sum_{k=-\infty}^\infty\left|
\left(\widetilde{\mathcal{D}}_k\right)^*h\right|^2
\right]^\frac{1}{2}\right\|_{L^{p'}(X)}\\
&\quad\lesssim\left\|\left\{\sum_{k=-\infty}^\infty\left[
\left(\sum_{l=-\infty}^\infty2^{-|k-l|\varepsilon'}\right)^\frac{1}{2}\left(
\sum_{l=-\infty}^\infty2^{-|k-l|\varepsilon'}
\left[\mathcal{M}_\rho(\mathcal{D}_lh)\right]^2\right)^\frac{1}{2}\right]^2
\right\}^\frac{1}{2}\right\|_{L^{p'}(X)}\\
&\qquad\qquad\text{by H\"older's inequality}\\
&\quad\sim\left\|\left\{\sum_{l=-\infty}^\infty\left[
\mathcal{M}_\rho(\mathcal{D}_lh)\right]^2
\right\}^\frac{1}{2}\right\|_{L^{p'}(X)}
\lesssim\left\|\left(\sum_{l=-\infty}^\infty\left|
\mathcal{D}_lh\right|^2
\right)^\frac{1}{2}\right\|_{L^{p'}(X)}
\quad\text{by Lemma~\ref{Fefferman--Stein}}\\
&\quad\lesssim
\left\|h\right\|_{L^{p'}(X)}
\quad\text{by Lemma~\ref{gL2L2}(i)}.
\end{align*}
This finishes the proof of \eqref{1029} and hence the above claim.

By H\"older's inequality twice,
we find that, for any $l,m\in\mathbb{N}$ with $l>m$,
\begin{align}\label{2045}
\left\|f_l-f_m\right\|_{L^p(X)}
&=\sup_{\|h\|_{L^{p'}(X)}\leq1}\left|\left\langle
\sum_{m<|k|\leq l}\widetilde{\mathcal{D}}_k
\mathcal{D}_kf,h\right\rangle\right|
=\sup_{\|h\|_{L^{p'}(X)}\leq1}\left|\sum_{m<|k|\leq l}\left\langle
\mathcal{D}_kf,
\left(\widetilde{\mathcal{D}}_k\right)^*h\right\rangle\right|
\nonumber\\
&\leq\sup_{\|h\|_{L^{p'}(X)}\leq1}
\int_X\sum_{m<|k|\leq l}\left|\mathcal{D}_kf(x)\right|
\left|\left(\widetilde{\mathcal{D}}_k\right)^*h(x)\right|\,d\mu(x)
\nonumber\\
&\leq\sup_{\|h\|_{L^{p'}(X)}\leq1}
\int_X\left[\sum_{m<|k|\leq l}\left|\mathcal{D}_kf(x)\right|^2\right]^\frac{1}{2}
\left[\sum_{m<|k|\leq l}\left|\left(\widetilde{\mathcal{D}}_k
\right)^*h(x)\right|^2\right]^\frac{1}{2}\,d\mu(x)
\nonumber\\
&\leq\left\|\left(\sum_{m<|k|\leq l}\left|\mathcal{D}_kf\right|^2
\right)^\frac{1}{2}\right\|_{L^p(X)}\sup_{\|h\|_{L^{p'}(X)}\leq1}
\left\|\left[\sum_{k=-\infty}^\infty\left|
\left(\widetilde{\mathcal{D}}_k\right)^*h\right|^2
\right]^\frac{1}{2}\right\|_{L^{p'}(X)}
\nonumber\\
&\lesssim\left\|\left(\sum_{m<|k|}
\left|\mathcal{D}_kf\right|^2\right)^\frac{1}{2}\right\|_{L^p(X)}
\quad\text{by \eqref{1029}},
\end{align}
which, together with $g(f)\in L^p(X)$,
further implies that, for any given $\delta\in(0,\infty)$,
there exists $N_0\in\mathbb{N}$ such that, for any $l,m>N_0$,
\begin{align*}
\left\|f_l-f_m\right\|_{L^p(X)}
\lesssim\left\|\left(\sum_{m<|k|}\left|\mathcal{D}_kf
\right|^2\right)^\frac{1}{2}\right\|_{L^p(X)}
\leq\left\|\left(\sum_{N_0<|k|}\left|\mathcal{D}_kf
\right|^2\right)^\frac{1}{2}\right\|_{L^p(X)}<\delta.
\end{align*}
From this, we deduce that $\{f_N\}_{N\in\mathbb{N}}$ is a
Cauchy sequence in $L^p(X)$,
which, combined with the completeness of $L^p(X)$,
further implies that there exists $u\in L^p(X)$ such that
$f_N\to u$ in $L^p(X)$ as $N\to\infty$.
By this and \eqref{92611} [applied to
$g=f\in(\mathring{\mathcal{G}}^\varepsilon_0(\beta,\gamma))'$],
we conclude that, for any
$\varphi\in\mathring{\mathcal{G}}^\varepsilon_0(\beta,\gamma)\subset L^{p'}(X)$
[see Lemma~\ref{CGL} for this inclusion],
\begin{align*}
\left\langle f,\varphi\right\rangle=\lim_{N\to\infty}
\left\langle \sum_{k=-N}^N\widetilde{\mathcal{D}}_k\mathcal{D}_kf,
\varphi\right\rangle=\lim_{N\to\infty}\left\langle f_N,
\varphi\right\rangle=\left\langle u,\varphi\right\rangle
\end{align*}
and hence $f=u$ in $(\mathring{\mathcal{G}}^\varepsilon_0(\beta,\gamma))'$.
Moreover, from the proven conclusion that $f_N\to u$
in $L^p(X)$ as $N\to\infty$ and an argument similar to that used in the
estimation of \eqref{2045} with $f_l-f_m$
replaced by $f_N$, we infer that, for any $N\in\mathbb{N}$,
\begin{align*}
\left\|u\right\|_{L^p(X)}
&\leq\left\|u-f_N\right\|_{L^p(X)}+\left\|f_N\right\|_{L^p(X)}
\\
&\lesssim\left\|u-f_N\right\|_{L^p(X)}+\left\|\left(\sum_{|k|\leq N}
\left|\mathcal{D}_kf\right|^2\right)^\frac{1}{2}\right\|_{L^p(X)}
\\
&\leq\left\|u-f_N\right\|_{L^p(X)}+\left\|g(f)\right\|_{L^p(X)},
\end{align*}
where the implicit positive constant is independent of both
$f$ and $N$.
Letting $N\to\infty$, we then
complete the proof of Lemma~\ref{gL2L2}(ii) and hence Lemma~\ref{gL2L2}.
\end{proof}

\begin{remark}\label{2057}
Let the notation be as in Lemma~\ref{gL2L2}.
For any $f\in L^p(X)$ and $x\in X$, let
\begin{align*}
\overline{g}(f)(x):=\left[\sum_{k\in\mathbb{Z}}
\left|\overline{\mathcal{D}}_kf(x)\right|^2\right]^\frac{1}{2},
\end{align*}
where $\{\overline{\mathcal{D}}_k\}_{k\in\mathbb{Z}}$ is the same as
in Theorem~\ref{A3.19}.
Then, by an argument similar to that used in the proof of
Lemma~\ref{gL2L2} [which needs the regularity only for the
second variable of $\overline{D}_k$], we find that there
exists a constant $C\in(0,\infty)$ such that,
for any $f\in L^p(X)$,
\begin{align*}
\left\|\overline{g}(f)\right\|_{L^p(X)}
\leq C\|f\|_{L^p(X)}.
\end{align*}
\end{remark}

We are ready to show Theorem~\ref{Hp=Lp}.

\begin{proof}[Proof of Theorem~\ref{Hp=Lp}]
The conclusions in the present theorem follow immediately from
Theorem~\ref{Lusin} below (with both $s:=0$ and $q:=2$)
and Lemma~\ref{gL2L2}.
\end{proof}

\section{Atomic Characterization}
\label{Section3.2}

Let $(X,\mathbf{q},\mu)$ be an ultra-RD-space with upper
dimension $n\in(0,\infty)$ and $\mathrm{diam\,}(X)=\infty$.
The main target of this section is to establish Theorem~\ref{Thm4.11} below,
namely the atomic characterization
of $H^p(X)$ with $p\in(\frac{n}{n+\mathrm{ind\,}(X,\mathbf{q})},1]$, where
$\mathrm{ind\,}(X,\mathbf{q})$ is the same as in \eqref{index}.
To achieve this, we first recall the notions of atoms
and atomic Hardy spaces $\mathring{H}^{p,q}_{\mathrm{at}}(X)$,
and prove that such atomic Hardy spaces
are independent of the choice of the parameter $q$;
see Theorem~\ref{328} below.
Then we show that $\mathring{H}^{p,2}_{\mathrm{at}}(X)$
coincides with $H^p(X)$; see Theorem~\ref{329} below.

To this end,
we begin with recalling the notion of atoms as follows.

\begin{definition}\label{atom}
Let $(X,\mathbf{q},\mu)$ be a
quasi-ultrametric space of homogeneous type. Assume that
$p\in(0,1]$,
$q\in(p,\infty]\cap[1,\infty]$,
and $\rho\in\mathbf{q}$ is such that all $\rho$-balls are $\mu$-measurable.
A $\mu$-measurable function
$a: X\to\mathbb{C}$ is called
a $(\rho,p,q)$-\emph{atom}\index[words]{atom}
if there exist $x\in X$
and $r\in(0,\infty)$ such that
\begin{enumerate}
\item[\rm(i)]
$\mathrm{\,supp\,} a:=\{x\in X:a(x)\neq0\}\subset B_{\rho}(x,r)$.

\item[\rm(ii)]
$||a||_{L^q(X)}\leq[(V_\rho)_r(x)]^{\frac{1}{q}-\frac{1}{p}}$.

\item[\rm(iii)]
$\int_Xa(x)\,d\mu(x)=0$.
\end{enumerate}
Moreover, when $\mu(X)<\infty$,
it is agreed upon that the constant function $a:=[\mu(X)]^{-\frac{1}{p}}$
is a $(\rho,p,q)$-atom.
\end{definition}

\begin{remark}\label{Rmk-atoms}
\begin{enumerate}
\item[\rm(i)]
It is easy to see from Definition~\ref{atom},
\eqref{1557-xx}, and \eqref{doublingcond}
that, if $a: X\to\mathbb{C}$ is a $(\rho,p,q)$-atom
supported in $B_{\rho}(x,r)$ for some $x\in X$
and $r\in(0,\infty)$, then, for any other $\rho'\in\mathbf{q}$
with the property that all $\rho'$-balls are $\mu$-measurable,
there exist positive constants $C$ and $c$ (depending only on $\rho$,
$\rho'$, and $\mu$) such that $Ca$ is
a $(\rho',p,q)$-atom supported in $B_{\rho'}(x,cr)$.
As such, we sometimes simply refer to $a$ as a
$(p,q)$-\textit{atom} or a
$(p,q)$-\textit{atom modulo a harmless constant multiple}.

\item[\rm(ii)]
Note that a $(p,q)$-atom is in $L^1(X)$. Indeed,
using H\"older's inequality and both (i) and (ii)
of Definition~\ref{atom}, we obtain
\begin{align}\label{atomL1}
\int_X|a(x)|\,d\mu(x)\leq\left\|\mathbf{1}_B\right\|_{L^{q'}(X)}
\left\|a\right\|_{L^q(X)}\leq\left[\mu(B)\right]
^{1-\frac{1}{p}}
\end{align}
for any $(p,q)$-atom $a$ supported in some ball $B\subset X$.

\item[\rm(iii)]
Any $(p,q)$-atom
is normalized in the way that its $L^p(X)$ norm does not exceed 1.
Indeed, from H\"older's inequality and both (i) and (ii)
of Definition~\ref{atom}, we deduce that,
for any $(p,q)$-atom $a$ supported in some ball $B\subset X$,
\begin{align*}
\|a\|_{L^p(X)}\leq\|a\|_{L^q(X)}
\left\|\mathbf{1}_B\right\|_{L^{(\frac{q}{p})'}(X)}^\frac{1}{p}
\leq1.
\end{align*}
\end{enumerate}
\end{remark}

We recall the following concept of atomic Hardy spaces.

\begin{definition}\label{atomichardy}
Let $(X,\mathbf{q},\mu)$ be a quasi-ultrametric space
of homogeneous type with upper dimension $n\in(0,\infty)$
and $\rho\in\mathbf{q}$ have the property that all $\rho$-balls are $\mu$-measurable.
Let $p\in(0,1]$, $q\in(p,\infty]\cap[1,\infty]$,
and $\varepsilon,\beta,\gamma\in(0,\infty)$ satisfy \eqref{1945}.
\begin{enumerate}
\item[\rm(i)]
The \emph{atomic Hardy space}\index[words]{atomic Hardy space}
${H}^{p,q}_{\mathrm{at}}(X,\rho)$\index[symbols]{H@${H}^{p,q}_{\mathrm{at}}(X,\rho)$}
is defined to be the set of all $f\in(\mathcal{G}_0^\varepsilon(\beta,\gamma))'$
such that
$$
f=\sum_{j\in\mathbb{N}}\lambda_ja_j
$$
in $(\mathcal{G}_0^\varepsilon(\beta,\gamma))'$,
where $\{a_j\}_{j\in\mathbb{N}}$ is a sequence of $(\rho,p,q)$-atoms
and $\{\lambda_j\}_{j\in\mathbb{N}}$ in $\mathbb{C}$
satisfies $\sum_{j\in\mathbb{N}}|\lambda_j|^p<\infty$.
Moreover, for any $f\in{H}^{p,q}_{\mathrm{at}}(X,\rho)$,
let
\begin{align*}
\left\|f\right\|_{{H}^{p,q}_{\mathrm{at}}(X,\rho)}:=
\inf\left\{\left(\sum_{j\in\mathbb{N}}|\lambda_j|^p\right)^\frac{1}{p}:
f=\sum_{j\in\mathbb{N}}\lambda_ja_j\text{ in }
(\mathcal{G}_0^\varepsilon(\beta,\gamma))'\right\},
\end{align*}
where the infimum is taken over all decompositions of $f$ as above.

\item[\rm(ii)]
The \emph{atomic Hardy space}
$\mathring{H}^{p,q}_{\mathrm{at}}(X,\rho)$
\index[symbols]{H@$\mathring{H}^{p,q}_{\mathrm{at}}(X,\rho)$}
is defined in the same way as $H^{p,q}_{\mathrm{at}}(X,\rho)$ in (i)
but with $(\mathcal{G}_0^\varepsilon(\beta,\gamma))'$
replaced by $(\mathring{\mathcal{G}}_0^\varepsilon(\beta,\gamma))'$.
\end{enumerate}
\end{definition}

\begin{remark}\label{ebe53}
Let the notation be as in Definition~\ref{atomichardy}.
By Remark~\ref{Rmk-atoms}(i), we find that
the spaces ${H}^{p,q}_{\mathrm{at}}(X,\rho)$
and $\mathring{H}^{p,q}_{\mathrm{at}}(X,\rho)$ are independent of the choice of
$\rho\in\mathbf{q}$. More precisely, if $\rho'\in\mathbf{q}$ is another quasi-ultrametric
with the property that all $\rho'$-balls are $\mu$-measurable, then
${H}^{p,q}_{\mathrm{at}}(X,\rho)={H}^{p,q}_{\mathrm{at}}(X,\rho')$
and $\mathring{H}^{p,q}_{\mathrm{at}}(X,\rho)=\mathring{H}^{p,q}_{\mathrm{at}}(X,\rho')$
with equivalent quasi-norms. Thus, throughout this monograph, we sometimes simply write
${H}^{p,q}_{\mathrm{at}}(X)$ and $\mathring{H}^{p,q}_{\mathrm{at}}(X)$
instead of ${H}^{p,q}_{\mathrm{at}}(X,\rho)$ and $\mathring{H}^{p,q}_{\mathrm{at}}(X,\rho)$.
\end{remark}

For the relationship between ${H}^{p,q}_{\mathrm{at}}(X)$
and $\mathring{H}^{p,q}_{\mathrm{at}}(X)$,
we have the following conclusion,
whose proof is a slight modification of \cite[Theorem~5.4]{LLD}
with ${H}^{p}_{\mathrm{at}}(X)$
and $\mathring{H}^{p}_{\mathrm{at}}(X)$ therein replaced by
${H}^{p,q}_{\mathrm{at}}(X)$
and $\mathring{H}^{p,q}_{\mathrm{at}}(X)$; we omit the details.

\begin{proposition}\label{324}
Let $(X,\mathbf{q},\mu)$ be a quasi-ultrametric space
of homogeneous type with upper dimension $n\in(0,\infty)$.
Let $p\in(0,1]$, $q\in(p,\infty]\cap[1,\infty]$,
and $\varepsilon,\beta,\gamma\in(0,\infty)$ satisfy \eqref{1945}.
Then
\begin{align*}
H^{p,q}_{\mathrm{at}}(X)=\mathring{H}^{p,q}_{\mathrm{at}}(X)
\end{align*}
with equivalent quasi-norms in the following sense:
\begin{enumerate}
\item[\textup{(i)}]
For any $f\in H^{p,q}_{\mathrm{at}}(X,\rho)
\subset({\mathcal{G}}^\varepsilon_0(\beta,\gamma))'$,
the restriction $\widehat{f}$ of $f$ to
$\mathring{{\mathcal{G}}}^\varepsilon_0(\beta,\gamma)$
belongs to $\mathring{H}^{p,q}_{\mathrm{at}}(X,\rho)$ and, moreover,
\begin{align*}
\left\|\widehat{f}\,\right\|_{\mathring{H}^{p,q}_{\mathrm{at}}(X,\rho)}
\leq\left\|f\right\|_{H^{p,q}_{\mathrm{at}}(X,\rho)}.
\end{align*}

\item[\textup{(ii)}]
For any $f\in\mathring{H}^{p,q}_{\mathrm{at}}(X,\rho)\subset
(\mathring{{\mathcal{G}}}^\varepsilon_0(\beta,\gamma))'$,
there exists a unique
$\widetilde{f}\in H^{p,q}_{\mathrm{at}}(X,\rho)$ such that
the restriction of $\widetilde{f}$ to
$\mathring{{\mathcal{G}}}^\varepsilon_0(\beta,\gamma)$ is $f$
and, moreover,
\begin{align*}
\left\|\widetilde{f}\,\right\|_{H^{p,q}_{\mathrm{at}}(X,\rho)}
\leq C\left\|f\right\|_{\mathring{H}^{p,q}_{\mathrm{at}}(X,\rho)},
\end{align*}
where $C$ is a positive constant independent of $f$.
\end{enumerate}
Consequently, the linear map
$\Phi:H^{p,q}_{\mathrm{at}}(X)\to\mathring{H}^{p,q}_{\mathrm{at}}(X)$,
defined by setting $\Phi(f)\in\mathring{H}^{p,q}_{\mathrm{at}}(X)$
to be the restriction of $f\in H^{p,q}_{\mathrm{at}}(X)$
to the space $\mathring{\mathcal{G}}_0^\varepsilon(\beta,\gamma)$,
is well defined, linear, bijective,
and has the property that, for any
$f\in H^{p,q}_{\mathrm{at}}(X)$,
\begin{align*}
\|f\|_{H^{p,q}_{\mathrm{at}}(X)}\sim
\left\|\Phi(f)\right\|_{\mathring{H}^{p,q}_{\mathrm{at}}(X)},
\end{align*}
where the positive equivalence constants are independent of $f$.
\end{proposition}

Let $(X,\mathbf{q},\mu)$ be a quasi-ultrametric space
of homogeneous type.
Recall that, for any $\alpha\in(0,\infty)$,
the \emph{Lipschitz space}\index[words]{Lipschitz space}
$\mathcal{L}_\alpha(X)$
is defined to be the set of all $\mu$-measurable functions $f: X\to\mathbb{C}$ such that
\begin{align*}
\left\|f\right\|_{\mathcal{L}_\alpha(X)}:=\sup_{\genfrac{}{}{0pt}{}{x,y\in X}{x\neq y}}
\frac{|f(x)-f(y)|}{[V_\rho(x,y)]^\alpha}<\infty\index[symbols]{L@$\mathcal{L}_\alpha(X)$}
\end{align*}
for some (hence all) $\rho\in\mathbf{q}$ having the property
that all $\rho$-balls are $\mu$-measurable.
Denote by $(\mathcal{L}_\alpha(X))'$ the \emph{dual space}
of $\mathcal{L}_\alpha(X)$ equipped with the weak-$*$ topology.

Next, we recall the notion of atomic Hardy spaces
on quasi-ultrametric spaces of homogeneous type
introduced by Coifman and Weiss \cite{Cowe77}.

\begin{definition}\label{1539}
Let $(X,\mathbf{q},\mu)$ be a quasi-ultrametric space
of homogeneous type
and $\rho\in\mathbf{q}$ have the property that all $\rho$-balls are $\mu$-measurable.
\begin{enumerate}
\item[\rm(i)]
Let $p\in(0,1)$ and $q\in[1,\infty]$.
The \emph{atomic Hardy space} $H^{p,q}_{\mathrm{CW}}(X,\rho)$
is defined to be the subspace of $(\mathcal{L}_{\frac{1}{p}-1}(X))'$,
which consists of all elements $f$ admitting an atomic decomposition
\begin{align}\label{atmoicdec}
f=\sum_{j\in\mathbb{N}}\lambda_ja_j
\end{align}
in $(\mathcal{L}_{\frac{1}{p}-1}(X))'$,
where $\{a_j\}_{j\in\mathbb{N}}$ is a sequence of $(p,q)$-atoms
and $\{\lambda_j\}_{j\in\mathbb{N}}$ in $\mathbb{C}$
satisfies
$\sum_{j\in\mathbb{N}}|\lambda_j|^p<\infty$.
Moreover, for any $f\in H^{p,q}_{\mathrm{CW}}(X,\rho)$, define
\begin{align*}
\left\|f\right\|_{H^{p,q}_{\mathrm{CW}}(X,\rho)}:=\inf\left\{
\left(\sum_{j\in\mathbb{N}}\left|\lambda_j\right|^p\right)^\frac{1}{p}:
f=\sum_{j\in\mathbb{N}}\lambda_ja_j\text{ in }
\left(\mathcal{L}_{\frac{1}{p}-1}(X)\right)'\right\},
\end{align*}
where the infimum is taken over all decompositions of $f$ as in \eqref{atmoicdec}.
\item[\rm(ii)]
Let $q\in(1,\infty]$.
The \emph{atomic Hardy space} $H^{1,q}_{\mathrm{CW}}(X,\rho)$
is defined to be the subspace of $L^1(X)$,
which consists of all elements $f$ admitting an atomic decomposition
\eqref{atmoicdec} in $L^1(X)$,
where $\{a_j\}_{j\in\mathbb{N}}$ is a sequence of $(1,q)$-atoms
and $\{\lambda_j\}_{j\in\mathbb{N}}$ in $\mathbb{C}$
satisfies
$\sum_{j\in\mathbb{N}}|\lambda_j|<\infty$.
Moreover, for any $f\in H^{1,q}_{\mathrm{CW}}(X,\rho)$, define
\begin{align*}
\left\|f\right\|_{H^{1,q}_{\mathrm{CW}}(X,\rho)}:=\inf\left\{
\sum_{j\in\mathbb{N}}\left|\lambda_j\right|:
f=\sum_{j\in\mathbb{N}}\lambda_ja_j\text{ in }L^1(X)\right\},
\end{align*}
where the infimum is taken over all decompositions of $f$ as above.
\end{enumerate}
\end{definition}

\begin{remark}\label{326}
Let the notation be as in Definition~\ref{1539}.
The following statements can be found in \cite{Cowe77}.
\begin{enumerate}
\item[\rm(i)]
Similarly to Remark~\ref{ebe53}, the space
$H^{p,q}_{\mathrm{CW}}(X,\rho)$ is independent of the choice of
$\rho\in\mathbf{q}$. Thus, throughout this monograph, we sometimes simply write
$H^{p,q}_{\mathrm{CW}}(X)$
instead of $H^{p,q}_{\mathrm{CW}}(X,\rho)$. Moreover,
$H^{p,q}_{\mathrm{CW}}(X)$ is also independent of the choice of
$q\in(p,\infty]\cap[1,\infty]$; see \cite[Theorem~A]{Cowe77}.

\item[\rm(ii)]
$H^{p,q}_{\mathrm{CW}}(X)$ is a quasi-Banach space when $p\in(0,1)$
and a Banach space when $p=1$.

\item[\rm(iii)]
The dual space of $H^{p,q}_{\mathrm{CW}}(X)$ is the Lipschitz space
$\mathcal{L}_{\frac{1}{p}-1}(X)$ when $p\in(0,1)$
and the space of functions with bounded mean oscillation $\mathrm{BMO\,}(X)$
(which was originally introduced
by John and Nirenberg in \cite[pp.\,415--416]{JN})
when $p=1$.
\end{enumerate}
\end{remark}

The following lemma is a slight modification of \cite[Theorem~4.16]{ZYJLDW}.

\begin{lemma}\label{325}
Let $(X,\mathbf{q},\mu)$ be a quasi-ultrametric space
of homogeneous type with upper dimension $n\in(0,\infty)$.
Let $p\in(0,1]$, $q\in(p,\infty]\cap[1,\infty]$,
and $\varepsilon,\beta,\gamma\in(0,\infty)$ satisfy \eqref{1945}.
Then, as subspaces of $(\mathcal{G}_0^\varepsilon(\beta,\gamma))'$,
$H^{p,q}_{\mathrm{at}}(X)=H^{p,q}_{\mathrm{CW}}(X)$
with equivalent (quasi-)norms.
\end{lemma}

As an immediate corollary of Proposition~\ref{324},
Lemma~\ref{325}, and Remark~\ref{326}(i),
we have the following conclusion.

\begin{theorem}\label{328}
Let $(X,\mathbf{q},\mu)$ be a quasi-ultrametric space
of homogeneous type with upper dimension $n\in(0,\infty)$.
Let $p\in(0,1]$, $q\in(p,\infty]\cap[1,\infty]$,
and $\varepsilon,\beta,\gamma\in(0,\infty)$ satisfy \eqref{1945}.
Then $\mathring{H}^{p,q}_{\mathrm{at}}(X)=\mathring{H}^{p,2}_{\mathrm{at}}(X)$.
In particular, $\mathring{H}^{p,q}_{\mathrm{at}}(X)$
is independent of the choice of $q$.
\end{theorem}

To obtain the atomic characterization of $H^p(X)$
[that is, to prove that $H^p(X)=\mathring{H}^{p,2}_{\mathrm{at}}(X)$
and hence $H^p(X)=\mathring{H}^{p,q}_{\mathrm{at}}(X)$
for all $q\in(p,\infty]\cap[1,\infty]$
by Lemma~\ref{325}],
we need several technical lemmas.
First, for any $\theta\in(0,1)$, we have the following relation
between $S_\theta$ and $S$.

\begin{lemma}\label{Sa1}	
Let $(X,\mathbf{q},\mu)$ be a quasi-ultrametric space of homogeneous
type with upper dimension $n\in(0,\infty)$.
Assume that $\theta\in(0,1)$ and $0<\beta,\gamma<\varepsilon<\infty$.
Then there exists a positive constant $C$ such that,
for any $f\in(\mathring{{\mathcal{G}}}_0^\varepsilon(\beta,\gamma))'$ and $x\in X$,
$$
S_\theta(f)(x)
\leq C\theta^{-\frac{n}{2}}S(f)(x),
$$
where $S_\theta(f)$ and $S(f)$ are as in Definition~\ref{areafncts}(i).
\end{lemma}

\begin{proof}
Let $\rho\in\mathbf{q}$ be
a regularized quasi-ultrametric as in Definition~\ref{D2.3}.
Let $f\in(\mathring{{\mathcal{G}}}_0^\varepsilon(\beta,\gamma))'$ and $x\in X$.
By \eqref{upperdoub} [used here with $\lambda=\theta^{-1}\in(1,\infty)$],
we find that, for any $k\in\mathbb{Z}$,
$(V_\rho)_{2^{-k}}(x)\lesssim\theta^{-n}(V_\rho)_{\theta2^{-k}}(x)$,
which further implies that
$$
\int_{B_\rho(x,\theta2^{-k})}
\left|\mathcal{D}_kf(y)\right|^2
\frac{d\mu(y)}{(V_\rho)_{\theta2^{-k}}(x)}
\lesssim\theta^{-n}
\int_{ B_\rho(x,2^{-k})}
\left|\mathcal{D}_kf(y)\right|^2
\frac{d\mu(y)}{(V_\rho)_{2^{-k}}(x)}.
$$
Hence,
\begin{align*}
\left[S_\theta(f)(x)\right]^2&=
\sum_{k=-\infty}^\infty
\int_{B_\rho(x,\theta 2^{-k})}\left|\mathcal{D}_kf(y)\right|^2
\frac{d\mu(y)}{(V_\rho)_{\theta2^{-k}}(x)}\\
&\lesssim\theta^{-n}\sum_{k=-\infty}^\infty
\int_{ B_\rho(x,2^{-k})}
\left|\mathcal{D}_kf(y)\right|^2
\frac{d\mu(y)}{(V_\rho)_{2^{-k}}(x)}
=\theta^{-n}\left[S(f)(x)\right]^2
\end{align*}
and the proof of Lemma~\ref{Sa1} is now complete.
\end{proof}

We first show that $\mathring{H}^{p,2}_{\mathrm{at}}(X)\subset H^p(X)$.
This is accomplished
in Theorem~\ref{Sf}. A key ingredient in its proof is the following lemma which proves
that, for any $\theta\in[1,\infty)$,
there is a uniform bound for the $L^p$-quasi-norm of
the Lusin area function $S_\theta(a)$, where $a$ is an atom.

\begin{lemma}\label{Sa2}
Let $(X,\mathbf{q},\mu)$ be a quasi-ultrametric space of
homogeneous type with upper dimension $n\in(0,\infty)$.
Assume that $\theta\in[1,\infty)$
and $p\in(\frac{n}{n+\mathrm{ind\,}(X,\mathbf{q})},1]$.
Then there exists a positive constant $C$, independent of $\theta$,
such that, for any $(p,2)$-atom $a$,
$$
\left\|S_\theta(a)\right\|_{L^p(X)}\leq C\theta^\frac{n}{p}.
$$
\end{lemma}

\begin{proof}
Let $\varepsilon\in(0,\infty)$ satisfy
$n(\frac{1}{p}-1)<\varepsilon\preceq\mathrm{ind\,}(X,\mathbf{q})$
and let $\rho\in\mathbf{q}$ be a regularized
quasi-ultrametric as in Definition~\ref{D2.3}.
Assume that
$\{\mathcal{S}_{2^{-k}}\}_{k\in\mathbb{Z}}$
is a homogeneous approximation
of the identity of order $\varepsilon$
as in Definition~\ref{DefHATI} and,
for any $k\in\mathbb{Z}$, let
$\mathcal{D}_k:=\mathcal{S}_{2^{-k}}-\mathcal{S}_{2^{-k+1}}$.
Let $S_\theta$ be defined as in \eqref{Def-S-f}
associated with $\{\mathcal{D}_k\}_{k\in\mathbb{Z}}$,
and $a$ be a $(p,2)$-atom
supported in ball $B:=B_\rho(x_0,r_0)$
with $x_0\in X$ and $r_0\in(0,\infty)$.
From Remark~\ref{1032}, we infer that, if $r_0\leq r_{\rho}(x_0)$, then
$$
B_{\rho}(x,r_0)=B_{\rho}\left(x,r_{\rho}(x_0)\right)=\{x_0\}.
$$
Hence, without
loss of generality, we may assume that $r_0>r_{\rho}(x_0)$.
We write
\begin{align}
\label{1904}
\left\|S_\theta(a)\right\|_{L^p(X)}^p
&=\int_{B_\rho(x_0,4C_\rho^2\theta r_0)}
|S_\theta(a)(x)|^p\,d\mu(x)
+\int_{[B_\rho(x_0,4C_\rho^2\theta r_0)]^\complement}
\cdots\nonumber\\
&=:\mathrm{I_1}+\mathrm{I_2}.
\end{align}
To deal with $\mathrm{I_1}$,
by Tonelli's theorem and Lemma~\ref{A3.2}(vii), we find that
\begin{align*}
\left\|S_\theta(a)\right\|_{L^2(X)}
&=\left[\int_X\sum_{k=-\infty}^\infty
\left|\mathcal{D}_ka(y)\right|^2
\int_{B_\rho(y,\theta2^{-k})}
\frac{1}{(V_\rho)_{\theta2^{-k}}(x)}
\,d\mu(x)\,d\mu(y)\right]^\frac{1}{2}
\\
&\sim\left[\int_X\sum_{k=-\infty}^\infty
\left|\mathcal{D}_ka(y)\right|^2
\int_{B_\rho(y,\theta2^{-k})}
\frac{1}{(V_\rho)_{\theta2^{-k}}(y)}
\,d\mu(x)\,d\mu(y)\right]^\frac{1}{2},
\end{align*}
which, together with Lemma~\ref{gL2L2} with $p:=2$,
further implies that
\begin{align}\label{2001}
\left\|S_\theta(a)\right\|_{L^2(X)}
\sim\left[\int_X\sum_{k=-\infty}^\infty
\left|\mathcal{D}_ka(y)\right|^2
\,d\mu(y)\right]^\frac{1}{2}
=||g(a)||_{L^2(X)}
\lesssim||a||_{L^2(X)},
\end{align}
where $g$ is defined as in \eqref{Def-g(f)(x)}
and the implicit positive constants are independent of $\theta$.
From this, H\"older's inequality, Definition~\ref{atom}(ii),
the assumption that $\theta\geq1$,
and \eqref{upperdoub}, we deduce that
\begin{align}
\label{4.5}
\mathrm{I_1}
&\leq\left\|\,\left|S_\theta(a)\right|^p\right\|_{L^{\frac{2}{p}}(X)}
\left\|\mathbf{1}_{B_\rho(x_0,4C_\rho^2\theta r_0)}\right\|_
{L^{\frac{2}{2-p}}(X)}
\nonumber\\
&\lesssim
\left[(V_\rho)_{r_0}(x_0)\right]^{\frac{p}{2}-1}
\theta^{n(1-\frac{p}{2})}\left[(V_\rho)_{r_0}(x_0)\right]^{1-\frac{p}{2}}\leq\theta^n,
\end{align}
where the implicit positive constant is independent of $\theta$.

To deal with $\mathrm{I_2}$,
let $x\in[B_\rho(x_0,4C_\rho^2\theta r_0)]^\complement$ and
$k\in\mathbb{Z}$.
By Lemma~\ref{Pk},
we conclude that there exists a constant $C\in[1,\infty)$
such that, if $\rho(x,y)\ge C2^{-k}$, then $D_k(x,y)=0$.
We consider the following two cases for $k\in\mathbb{Z}$.

\emph{Case (1)}
$k\in\mathbb{Z}$ is such that $\rho(x,x_0)>CC_\rho^2\theta2^{-k}$.
In this case, we claim that, for any $y\in B_\rho(x,\theta2^{-k})$,
\begin{equation}\label{qqw-894}
B_\rho(y,C2^{-k})\cap B_\rho(x_0,r_0)=\varnothing.
\end{equation}
Indeed, let $y\in B_\rho(x,\theta2^{-k})$ and $z\in B_\rho(y,C2^{-k})$.
From $\rho(x,y)<\theta2^{-k}$, the fact that $C,\theta, C_\rho\in[1,\infty)$,
and Definition~\ref{D2.1}(iii),
we infer that
\begin{align*}
CC_\rho^2\theta2^{-k}<\rho(x,x_0)
&\leq C_\rho^2
\max\left\{\rho(x,y),\,\rho(y,z),\,\rho(z,x_0)\right\}\\
&\leq C_\rho^2\max\left\{\theta 2^{-k},\,C2^{-k},\,\rho(z,x_0)\right\}
=C_\rho^2\rho(z,x_0),
\end{align*}
which further implies that
\begin{align*}
\rho(z,x_0)\geq\frac{\rho(x,x_0)}{C_\rho^2}\geq4\theta r_0>r_0.
\end{align*}
Hence, $z\in[B_\rho(x_0,r_0)]^\complement$ and the claim in \eqref{qqw-894} holds.

Next, given the support conditions of $D_k$ and $a$, \eqref{qqw-894} implies that
$\mathcal{D}_ka(y)=0$ for any $y\in B_\rho(x,\theta2^{-k})$ and hence we find that
\begin{align}
\label{1508}
\sum_{\{k\in\mathbb{Z}:\rho(x,x_0)>CC_\rho^2\theta2^{-k}\}}
\int_{B_\rho(x,\theta2^{-k})}
\left|\mathcal{D}_ka(y)\right|^2
\frac{d\mu(y)}{(V_\rho)_{\theta2^{-k}}(x)}=0.
\end{align}

\emph{Case (2)}
$k\in\mathbb{Z}$ is such that $\rho(x,x_0)\leq CC_\rho^2\theta2^{-k}$.
Let $y\in B_\rho(x,\theta2^{-k})$ be arbitrary.
In this case, by Definition~\ref{atom}, we conclude that
\begin{align}
\label{1526}
\left|\mathcal{D}_ka(y)\right|
&\leq\int_{B_\rho(x_0,r_0)}
\left|D_k(y,u)-D_k(y,x_0)\right|\left|a(u)\right|\,d\mu(u)
\nonumber\\
&\lesssim2^{k\varepsilon}
\int_{B_\rho(x_0,r_0)}
\left[\rho(u,x_0)\right]^\varepsilon
\frac{1}{(V_\rho)_{2^{-k}}(y)}
|a(u)|\,d\mu(u)
\quad\text{by Lemma~\ref{Pk}}
\nonumber\\
&\leq2^{k\varepsilon}r_0^\varepsilon
\frac{1}{(V_\rho)_{2^{-k}}(y)}
\int_X|a(u)|\,d\mu(u)
\nonumber\\
&\leq2^{k\varepsilon}
r_0^\varepsilon
\frac{1}{(V_\rho)_{2^{-k}}(y)}
\left[(V_\rho)_{r_0}(x_0)\right]^{1-\frac{1}{p}}
\quad\text{by \eqref{atomL1}}.
\end{align}
From the support conditions of $D_k$ and $a$,
we deduce that, if $\mathcal{D}_ka(y)\neq0$, then
$$B_\rho(y,C2^{-k})\cap B_\rho(x_0,r_0)\neq\varnothing.$$
By this, Definition~\ref{D2.1}(iii),
the fact that $4C_\rho^2\theta r_0\leq
\rho(x,x_0)\leq CC_\rho^2\theta2^{-k}$,
we find that, for any
$u\in B_\rho(y,C2^{-k})\cap B_\rho(x_0,r_0)$,
$$
\rho(y,x_0)\leq C_\rho
\max\left\{\rho(y,u),\,\rho(u,x_0)\right\}
\leq C_\rho\max\left\{C2^{-k},\,r_0\right\}
\leq CC_\rho2^{-k}.
$$
From this, Lemma~\ref{A3.2}(vii),
and \eqref{upperdoub}, we infer that
$(V_\rho)_{2^{-k}}(y)\sim (V_\rho)_{2^{-k}}(x_0)$
and
$$
V_\rho(x,x_0)
\lesssim (V_\rho)_{\theta2^{-k}}(x_0)
\lesssim\theta^{n}(V_\rho)_{2^{-k}}(x_0)
\sim\theta^{n}(V_\rho)_{2^{-k}}(y),
$$
where the implicit positive constants are independent of $\theta$.
By this and \eqref{1526}, we obtain
$$
\left|\mathcal{D}_ka(y)\right|
\lesssim2^{k\varepsilon}
r_0^\varepsilon
\frac{\theta^n}{V_\rho(x,x_0)}
\left[(V_\rho)_{r_0}(x_0)\right]^{1-\frac{1}{p}},
$$
which, given that $y\in B_\rho(x,\theta2^{-k})$ was arbitrary, further implies that
\begin{align*}
&\sum_{\{k\in\mathbb{Z}:
\rho(x,x_0)\leq CC_\rho^2\theta2^{-k}\}}\int_{B_\rho(x,\theta2^{-k})}
\left|\mathcal{D}_ka(y)\right|^2
\frac{d\mu(y)}{(V_\rho)_{\theta2^{-k}}(x)}
\\
&\quad\lesssim\sum_{\{k\in\mathbb{Z}:
\rho(x,x_0)\leq CC_\rho^2\theta2^{-k}\}}
2^{2k\varepsilon}r_0^{2\varepsilon}
\frac{\theta^{2n}}{[V_\rho(x,x_0)]^2}
\left[(V_\rho)_{r_0}(x_0)\right]^{2-\frac{2}{p}}
\\
&\quad\sim r_0^{2\varepsilon}\theta^{2(n+\varepsilon)}
\frac{1}{[V_\rho(x,x_0)]^2}
\left[(V_\rho)_{r_0}(x_0)\right]^{2-\frac{2}{p}}
\left[\frac{1}{\rho(x,x_0)}\right]^{2\varepsilon},
\end{align*}
where the implicit positive constants are independent of $\theta$.
This, combined with \eqref{1508},
implies that, for any
$x\in B_\rho(x_0,4C_\rho^2\theta r_0)^\complement$,
\begin{align}\label{dw-23}
S_\theta(a)(x)
\lesssim\theta^{n+\varepsilon}
\left[(V_\rho)_{r_0}(x_0)\right]^{1-\frac{1}{p}}
\frac{1}{V_\rho(x,x_0)}
\left[\frac{r_0}{\rho(x,x_0)}\right]^\varepsilon,
\end{align}
which, together with \eqref{upperdoub}
and $p\in(\frac{n}{n+\varepsilon},1]$, further
implies that
\begin{align}
\label{4.6}
\mathrm{I_2}
&\lesssim\theta^{p(n+\varepsilon)}
\left[(V_\rho)_{r_0}(x_0)\right]^{p-1}
\sum_{j=2}^{\infty}\int_{(2C_\rho)^j\theta r_0
\leq\rho(x,x_0)<(2C_\rho)^{j+1}\theta r_0}
\left\{\frac{[(2C_\rho)^j\theta]^{-\varepsilon}}{V_\rho(x,x_0)}
\right\}^p\,d\mu(x)
\nonumber\\
&\leq\theta^{p(n+\varepsilon)}
\left[(V_\rho)_{r_0}(x_0)\right]^{p-1}
\sum_{j=2}^{\infty}[(2C_\rho)^j\theta]^{-p\varepsilon}
\frac{\mu(B_\rho(x_0,(2C_\rho)^{j+1}\theta r_0))}
{[\mu(B_\rho(x_0,(2C_\rho)^j\theta r_0))]^p}
\nonumber\\
&\lesssim\theta^{pn}
\left[(V_\rho)_{r_0}(x_0)\right]^{p-1}
\sum_{j=2}^{\infty}(2C_\rho)^{-jp\varepsilon}
\left[\mu(B_\rho(x_0,(2C_\rho)^{j}\theta r_0))\right]^{1-p}
\nonumber\\
&\lesssim\theta^{n}
\sum_{j=2}^{\infty}(2C_\rho)^{j(n-pn-p\varepsilon)}
\sim\theta^{n},
\end{align}
where the implicit positive constants are independent of $\theta$.

Combining \eqref{1904}, \eqref{4.5}, and \eqref{4.6},
we conclude that, for any $\theta\in[1,\infty)$,
$$
\left\|S_\theta(a)\right\|_{L^p(X)}
=\left(\mathrm{I_1}+\mathrm{I_2}\right)^\frac{1}{p}
\lesssim\theta^\frac{n}{p}.
$$
This finishes the proof of Lemma~\ref{Sa2}.
\end{proof}

Using Lemmas~\ref{Sa1} and~\ref{Sa2},
we obtain the following conclusion.

\begin{theorem}\label{Sf}
Let $(X,\mathbf{q},\mu)$ be a quasi-ultrametric spaces
of homogeneous type
with upper dimension $n\in(0,\infty)$.
Assume that $\theta\in(0,\infty)$ and
$p\in(\frac{n}{n+\mathrm{ind\,}(X,\mathbf{q})},1]$.
Then there exists a positive constant $C$,
independent of $\theta$, such that,
for any $f\in\mathring{H}^{p,2}_{\mathrm{at}}(X)$,
\begin{align*}
\left\|S_\theta(f)\right\|_{L^p(X)}
\leq C\max\left\{\theta^{-\frac{n}{2}},\,\theta^\frac{n}{p}\right\}
||f||_{\mathring{H}^{p,2}_{\mathrm{at}}(X)}.
\end{align*}
In particular,
$\mathring{H}^{p,2}_{\mathrm{at}}(X)\subset H^p(X)$.
\end{theorem}

\begin{proof}
Let $\rho\in\mathbf{q}$ be a regularized quasi-ultrametric
and $f\in\mathring{H}^{p,2}_{\mathrm{at}}(X)$.
By the definition of $\mathring{H}^{p,2}_{\mathrm{at}}(X)$,
we find that
there exist $\{a_j\}_{j\in\mathbb{N}}$ of
$(p,2)$-atoms
and $\{\lambda_j\}_{j\in\mathbb{N}}$ in $\mathbb{C}$
such that $f=\sum_{j\in\mathbb{N}}\lambda_ja_j$
in $(\mathring{{\mathcal{G}}}_0^\varepsilon(\beta,\gamma))'$ and
\begin{align}\label{atomsim}
\sum_{j\in\mathbb{N}}|\lambda_j|^p
\sim||f||^p_{\mathring{H}^{p,2}_{\mathrm{at}}(X)}.
\end{align}
Granted this convergence and
the fact that Lemma~\ref{Pk} implies $D_k(y,\cdot)\in
\mathring{{\mathcal{G}}}(\varepsilon,\varepsilon)\subset
\mathring{{\mathcal{G}}}_0^\varepsilon(\beta,\gamma)$, we obtain,
for any $y\in X$,
\begin{align*}
\left|\mathcal{D}_kf(y)\right|&=
\left|\left\langle f,D_k(y,\cdot)\right\rangle\right|\\
&=\left|\lim_{m\to\infty}\sum_{j=1}^m\lambda_j
\langle a_j,D_k(y,\cdot)\rangle\right|
=\left|\lim_{m\to\infty}\sum_{j=1}^m\lambda_j \mathcal{D}_k(a_j)(y)\right|\\
&\leq\sum_{j\in\mathbb{N}}
\left|\lambda_j\right|\cdot\left|\mathcal{D}_k(a_j)(y)\right|.
\end{align*}
From this and Minkowski's inequality,
it follows that, for any $x\in X$,
\begin{align*}
S_\theta(f)(x)
\leq\sum_{j\in\mathbb{N}}
\left|\lambda_j\right|S_\theta(a_j)(x).
\end{align*}
Having this estimate in mind, on the one hand,
for any given $\theta\in(0,1)$,
\begin{align*}
\left\|S_\theta(f)\right\|^p_{L^p(X)}
&\leq\left\|\sum_{j\in\mathbb{N}}
\left|\lambda_j\right|S_\theta(a_j)\right\|^p_{L^p(X)}
=\left\|\left[\sum_{j\in\mathbb{N}}
\left|\lambda_j\right|S_\theta(a_j)\right]^p\right\|_{L^1(X)}\nonumber\\
&\leq\left\|\sum_{j\in\mathbb{N}}
\left|\lambda_j\right|^p\left[S_\theta(a_j)\right]^p\right\|_{L^1(X)}
\quad\text{by Lemma~\ref{discreteinequality}}
\nonumber\\
&=\sum_{j\in\mathbb{N}}
\left|\lambda_j\right|^p\left\|S_\theta(a_j)\right\|^p_{L^p(X)}
\quad\text{by Tonelli's theorem}.
\end{align*}
This then implies that
\begin{align}
\label{17331}
\left\|S_\theta(f)\right\|^p_{L^p(X)}
&\lesssim\theta^{-\frac{np}{2}}\sum_{j\in\mathbb{N}}
\left|\lambda_j\right|^p\left\|S(a_j)\right\|^p_{L^p(X)}
\quad\text{by Lemma~\ref{Sa1}}\nonumber\\
&\lesssim\theta^{-\frac{np}{2}}\sum_{j\in\mathbb{N}}
\left|\lambda_j\right|^p
\quad\text{by Lemma~\ref{Sa2} with $\theta=1$}\nonumber\\
&\sim\theta^{-\frac{np}{2}}
||f||^p_{\mathring{H}^{p,2}_{\mathrm{at}}(X)}
\quad\text{by \eqref{atomsim}}.
\end{align}
On the other hand,
for any given $\theta\in[1,\infty)$,
\begin{align*}
\left\|S_\theta(f)\right\|^p_{L^p(X)}
&\leq\left\|\sum_{j\in\mathbb{N}}
\left|\lambda_j\right|S_\theta(a_j)\right\|^p_{L^p(X)}
=\left\|\left[\sum_{j\in\mathbb{N}}
\left|\lambda_j\right|S_\theta(a_j)\right]^p\right\|_{L^1(X)}\nonumber\\
&\leq\left\|\sum_{j\in\mathbb{N}}
\left|\lambda_j\right|^p\left[S_\theta(a_j)\right]^p\right\|_{L^1(X)}
\quad\text{by Lemma~\ref{discreteinequality}}
\nonumber\\
&=\sum_{j\in\mathbb{N}}
\left|\lambda_j\right|^p\left\|S_\theta(a_j)\right\|^p_{L^p(X)}
\quad\text{by Tonelli's theorem}
\nonumber\\
&\lesssim\theta^{n}\sum_{j\in\mathbb{N}}
\left|\lambda_j\right|^p
\quad\text{by Lemma~\ref{Sa2}}\nonumber\\
&\sim\theta^{n}
||f||^p_{\mathring{H}^{p,2}_{\mathrm{at}}(X)}
\quad\text{by \eqref{atomsim}},
\end{align*}
which, together with \eqref{17331},
further implies that,
for any given $\theta\in(0,\infty)$,
$$
\left\|S_\theta(f)\right\|_{L^p(X)}
\lesssim\max\left\{\theta^{-\frac{n}{2}},\,\theta^\frac{n}{p}\right\}
||f||_{\mathring{H}^{p,2}_{\mathrm{at}}(X)}
$$
and hence $f\in H^p(X)$.
This finishes the proof of Theorem~\ref{Sf}.
\end{proof}

Our next goal is to show that $H^p(X)\subset\mathring{H}^{p,2}_{\mathrm{at}}(X)$.
We begin with the following technical lemma.

\begin{lemma}
\label{2157}
Let $(X,\mathbf{q},\mu)$ be an ultra-RD-space
with upper dimension $n\in(0,\infty)$.
Let	$p\in(\frac{n}{n+\mathrm{ind\,}(X,\mathbf{q})},1]$
and $\beta,\gamma,\varepsilon\in(0,\infty)$ satisfy
\eqref{1945}.
Let $\{\mathcal{S}_{2^{-k}}\}_{k\in\mathbb{Z}}$ be
a homogeneous approximation
of the identity of order $\varepsilon$
in Definition~\ref{DefHATI} and let
$\mathcal{D}_k:=\mathcal{S}_{2^{-k}}-\mathcal{S}_{2^{-k+1}}$
for any $k\in\mathbb{Z}$.
Assume that $S_\mathcal{D}$ is the Lusin area function (with aperture 1) associated with
$\{\mathcal{D}_k\}_{k\in\mathbb{Z}}$.
Let $\{\overline{\mathcal{D}}_k\}_{k\in\mathbb{Z}}$ be the same as in Theorem~\ref{A3.19}
and $S_{\overline{\mathcal{D}}}$ be defined as in \eqref{Def-S-f}
with $\theta:=1$ and $\mathcal{D}_k$ replaced by $\overline{\mathcal{D}}_k$.
Then there exists a positive constant $C$ such that, for any
$f\in(\mathring{{\mathcal{G}}}^\varepsilon_0(\beta,\gamma))'$,
\begin{align}\label{DbarD0}
\left\|S_{\overline{\mathcal{D}}}(f)\right\|_{L^{p}(X)}\leq C
\left\|S_\mathcal{D}(f)\right\|_{L^{p}(X)}.
\end{align}
\end{lemma}

\begin{proof}
Following the proof of Theorem~\ref{independentofATI},
to show \eqref{DbarD0}, it suffices to prove that,
for any given $r\in(\frac{n}{n+\varepsilon},1]\cap(0,p)$ and for
any $f\in(\mathring{{\mathcal{G}}}_0^\varepsilon(\beta,\gamma))'$
and $x\in X$,
\begin{align}\label{SDbarMmf0}
\left[S_{\overline{\mathcal{D}}}(f)(x)\right]^2
\lesssim\sum_{k=-\infty}^{\infty}
\left[\mathcal{M}_\rho
\left([m_k(f)]^r\right)(x)\right]^{\frac{2}{r}},
\end{align}
where, for any $k\in\mathbb{Z}$,
\begin{align}\label{mkfx}
m_k(f)(x):=\left[\frac{1}{(V_\rho)_{2^{-k}}(x)}
\int_{B_\rho(x,2^{-k})}\left|\mathcal{D}_k(f)(u)\right|^2
\,d\mu(u)\right]^\frac{1}{2}.
\end{align}
To show \eqref{SDbarMmf0},
suppose that $l\in\mathbb{Z}$,
$x\in X$, and $y\in B_\rho(x,2^{-l})$.
By Theorem~\ref{HD}, we conclude that
\begin{align*}
\overline{\mathcal{D}}_lf(y)
&=\left\langle f,\overline{D}_l(y,\cdot)\right\rangle\\
&=\left\langle\sum_{k=-\infty}^{\infty}\sum_{\tau\in I_{k}}
\sum_{v=1}^{N(k,\tau)}
\widetilde{D}^{(2)}_k(\cdot,y_\tau^{k,v})
\int_{Q_\tau^{k,v}}\mathcal{D}_kf(u)\,d\mu(u),\overline{D}_l(y,\cdot)\right\rangle
\\
&=\sum_{k=-\infty}^{\infty}\sum_{\tau\in I_{k}}
\sum_{v=1}^{N(k,\tau)}\left\langle
\widetilde{D}^{(2)}_k(\cdot,y_\tau^{k,v})
,\overline{D}_l(y,\cdot)\right\rangle\int_{Q_\tau^{k,v}}
\mathcal{D}_kf(u)\,d\mu(u)
\\
&=\sum_{k=-\infty}^{\infty}\sum_{\tau\in I_{k}}
\sum_{v=1}^{N(k,\tau)}
\overline{D}_l\widetilde{D}^{(2)}_k(y,y_\tau^{k,v})
\int_{Q_\tau^{k,v}}\mathcal{D}_kf(u)\,d\mu(u),	
\end{align*}
where the penultimate equality holds because the
series converge in $(\mathring{{\mathcal{G}}}_0^\varepsilon(\beta,\gamma))'$ and
where $\widetilde{D}^{(2)}_k$ is the same as in Theorem~\ref{A3.47}.
Note that $\overline{D}_l(y,\cdot)$
(the kernel of $\overline{\mathcal{D}}_l$)
and $\widetilde{D}^{(2)}_k(\cdot,y_\tau^{k,v})$
(the kernel of $\widetilde{\mathcal{D}}^{(2)}_k$) satisfy all the conditions of
Lemma~\ref{1122} and hence \eqref{DkDl}
holds for $\overline{D}_l$
and $\widetilde{D}^{(2)}_k$
with $k,l\in\mathbb{Z}$.
From this and an argument similar to
that used in the estimation of \eqref{817},
we deduce that \eqref{SDbarMmf0} holds.
This finishes the proof of Lemma~\ref{2157}.
\end{proof}

Now, we are ready to prove
$H^p(X)\subset\mathring{H}^{p,2}_{\mathrm{at}}(X)$.

\begin{theorem}\label{L4.10}
Let $(X,\mathbf{q},\mu)$ be an ultra-RD-space
with upper dimension $n\in(0,\infty)$.
Let $p\in(\frac{n}{n+\mathrm{ind\,}(X,\mathbf{q})},1]$ and
$\beta,\gamma,\varepsilon\in(0,\infty)$ satisfy \eqref{1945}.
Then, for any $f\in H^p(X)$,
there exist a sequence of $(p,2)$-atoms
$\{a_j\}_{j\in\mathbb{N}}$
and a numerical sequence
$\{\lambda_j\}_{j\in\mathbb{N}}$ in $\mathbb{C}$
such that
$f=\sum_{j\in\mathbb{N}}\lambda_ja_j$
holds in $(\mathring{{\mathcal{G}}}_0^\varepsilon(\beta,\gamma))'$
and there exists a positive constant C,
independent of $f$, such that
$\sum_{j\in\mathbb{N}}|\lambda_j|^p
\leq C ||f||^p_{H^p(X)}$.
In particular,
$H^p(X)\subset\mathring{H}^{p,2}_{\mathrm{at}}(X)$.
\end{theorem}

\begin{proof}
Let $\rho\in\mathbf{q}$ be a regularized quasi-ultrametric
as in Definition~\ref{D2.3}.
Let $\{\mathcal{S}_{2^{-k}}\}_{k\in\mathbb{Z}}$ be
a homogeneous approximation
of the identity of order $\varepsilon$
in Definition~\ref{DefHATI} and let
$\mathcal{D}_k:=\mathcal{S}_{2^{-k}}-\mathcal{S}_{2^{-k+1}}$
for any $k\in\mathbb{Z}$.
Assume that $S_\mathcal{D}$ is the Lusin area
function (with aperture 1) associated with
$\{\mathcal{D}_k\}_{k\in\mathbb{Z}}$.
Let $\{\overline{\mathcal{D}}_k\}_{k\in\mathbb{Z}}$ be the same as in Theorem~\ref{A3.19}
and $S_{\overline{\mathcal{D}}}$ be defined as in \eqref{Def-S-f}
with $\theta:=1$ and $\mathcal{D}_k$ replaced by $\overline{\mathcal{D}}_k$.
Let $\mathfrak{Q}:=\{Q^k_\alpha:k\in\mathbb{Z},\,\alpha\in I_k\}$\index[symbols]{Q@$\mathfrak{Q}$}
be the set of all dyadic cubes (at all levels)
as in Lemma~\ref{DyadicSys} and recall that, by Remark~\ref{remarkcube}, we
may assume that $\delta=\frac{1}{2}$ in Lemma~\ref{DyadicSys}.

Let $f\in H^p(X)$. From this and Theorem~\ref{2157}, we deduce that
$$
\left\|S_{\overline{\mathcal{D}}}(f)\right\|_{L^p(X)}
\lesssim\left\|S_{\mathcal{D}}(f)\right\|_{L^p(X)}
=\left\|f\right\|_{H^p(X)}<\infty.
$$
In particular, $S_{\overline{\mathcal{D}}}(f)$ is finite $\mu$-almost everywhere
on $X$. For any $k\in\mathbb{Z}$, we define
$$
\Omega_k:=\left\{x\in X:S_{\overline{\mathcal{D}}}(f)(x)>2^k\right\},
$$
$$
\mathfrak{Q}_k\index[symbols]{Q@$\mathfrak{Q}_k$}
:=\left\{Q\in\mathfrak{Q}:
\mu(Q\cap\Omega_{k+1})\leq\frac{1}{2}\mu(Q)<\mu(Q\cap \Omega_k)\right\},
$$
$$
\Omega_{-\infty}:=\left\{x\in X:S_{\overline{\mathcal{D}}}(f)(x)=0\right\},
$$
and
$$
\mathfrak{Q}_{-\infty}\index[symbols]{Q@$\mathfrak{Q}_{-\infty}$}
:=\left\{Q\in\mathfrak{Q}:
\mu(Q\cap \Omega_{-\infty})\ge\frac{1}{2}\mu(Q)\right\}.
$$
It is easy to verify that
$\mathfrak{Q}=\bigcup_{k\in\mathbb{Z}\cup\{-\infty\}}\mathfrak{Q}_k$
and, for any dyadic cube $Q\in\mathfrak{Q}$,
there exists a unique $k\in\mathbb{Z}\cup\{-\infty\}$
such that $Q\in\mathfrak{Q}_k$.
A dyadic cube $Q\in\mathfrak{Q}_k$
is called a \emph{maximal dyadic cube}\index[words]{maximal dyadic cube} in $\mathfrak{Q}_k$
if, whenever $Q'\in\mathfrak{Q}$ is such that
$Q\neq Q'$ and $Q\subset Q'$,
it follows that  $Q'\notin\mathfrak{Q}_k$.
For any given $k,j\in\mathbb{Z}$,
denote the set of all maximal dyadic cubes
in $\mathfrak{Q}_k$ at level $j$
by $\{Q^j_{\tau,k}\}_{\tau\in I_{j,k}}$\index[symbols]{Q@$Q^j_{\tau,k}$}
with their centers denoted by
$\{z^j_{\tau,k}\}_{\tau\in I_{j,k}}$,
where $I_{j,k}\subset I_j$
(with $I_j$ as in Lemma~\ref{DyadicSys}) may be empty.
Then we can write
\begin{align}\label{936}
\mathfrak{Q}=\left(\bigcup_{k\in\mathbb{Z}}\mathfrak{Q}_k\right)
\cup\mathfrak{Q}_{-\infty}
=\left(\bigcup_{j,\,k\in\mathbb{Z}}
\bigcup_{\tau\in I_{j,k}}
\left\{Q\in\mathfrak{Q}_k:Q\subset Q^j_{\tau,k}\right\}\right)
\cup\mathfrak{Q}_{-\infty};
\end{align}
see also \cite[(4.12)]{YHYY2}.

Next, choose any $k_0\in\mathbb{Z}$ satisfying $C_\rho C_r\leq 2^{k_0}$,
where $C_\rho$ is the same positive constant as in \eqref{Crho} and
$C_r$ is the same positive constant as in Lemma~\ref{DyadicSys}. Note that we
can assume that $C_r\geq1$ so that $k_0\geq0$. We first claim that,
for any given $Q^{l+k_0}_\tau\in\mathfrak{Q}_{-\infty}$
with both $l\in\mathbb{Z}$ and $\tau\in I_{l+k_0}$
and for $\mu$-almost every $x\in Q^{l+k_0}_\tau$,
\begin{align}\label{934}
\overline{\mathcal{D}}_lf(x)=0.
\end{align}
Indeed, by the definition of $\mathfrak{Q}_{-\infty}$,
we conclude that
there exists $x_0\in Q^{l+k_0}_\tau$ such that
$$
S_{\overline{\mathcal{D}}}(f)(x_0)=0,
$$
which, combined with the definition of $S_{\overline{\mathcal{D}}}(f)$,
further implies that, for $\mu$-almost every $x\in B_\rho(x_0,2^{-l})$,
\begin{align}
\label{1950}
\overline{\mathcal{D}}_lf(x)=0.
\end{align}
Note that, from Lemma~\ref{DyadicSys}(iv),
the fact that $x_0\in Q^{l+k_0}_\tau$,
and $C_\rho C_r\leq 2^{k_0}$, we infer that
\begin{align}\label{1407}
Q^{l+k_0}_\tau&\subset
B_\rho\left(z_\tau^{l+k_0},C_r2^{-(l+k_0)}\right)
\subset B_\rho\left(x_0,C_\rho C_r2^{-(l+k_0)}\right)
\subset B_\rho\left(x_0,2^{-l}\right),
\end{align}
where $z_\tau^{l+k_0}$ is the center of $Q^{l+k_0}_\tau$.
By this and \eqref{1950},
we find that, for any given $Q\in\mathfrak{Q}_{-\infty}$
and for $\mu$-almost every $x\in Q$,
\eqref{934} holds.
This finishes the proof of the above claim.

Now, let $\mathcal{D}_Q:=\mathcal{D}_l$
and $\overline{\mathcal{D}}_Q:=\overline{\mathcal{D}}_l$
whenever $Q=Q^{l+k_0}_\tau$
for some $l\in\mathbb{Z}$
and $\tau\in I_{l+k_0}$.
Then, from Theorem~\ref{HC}, Theorem~\ref{DyadicSys}(i),
and \eqref{936},
we deduce that
\begin{align*}
f(\cdot)&=\sum_{l=-\infty}^{\infty}
\mathcal{D}_l\overline{\mathcal{D}}_lf(\cdot)
=\sum_{l=-\infty}^{\infty}\int_X
D_l(\cdot,y)\overline{\mathcal{D}}_lf(y)\,d\mu(y)
\nonumber\\
&=\sum_{l=-\infty}^{\infty}\sum_{\tau\in I_{l+k_0}}
\int_{Q^{l+k_0}_\tau}
D_l(\cdot,y)\overline{\mathcal{D}}_lf(y)\,d\mu(y)
\nonumber\\
&=\sum_{Q\in\mathfrak{Q}}
\int_QD_Q(\cdot,y)\overline{\mathcal{D}}_Qf(y)\,d\mu(y)
\nonumber\\
&=\sum_{k=-\infty}^{\infty}\sum_{j=-\infty}^{\infty}
\sum_{\tau\in I_{j,k}}
\sum_{Q\in\mathfrak{Q}_k,\,Q\subset Q^j_{\tau,k}}
\int_QD_Q(\cdot,y)\overline{\mathcal{D}}_Qf(y)\,d\mu(y)
\nonumber\\
&\quad+\sum_{Q\in\mathfrak{Q}_{-\infty}}
\int_QD_Q(\cdot,y)\overline{\mathcal{D}}_Qf(y)\,d\mu(y),
\end{align*}
which, together with \eqref{934},
further implies that
\begin{align}\label{lammdab}
f(\cdot)&=\sum_{k=-\infty}^{\infty}\sum_{j=-\infty}^{\infty}
\sum_{\tau\in I_{j,k}}
\sum_{Q\in\mathfrak{Q}_k,\,Q\subset Q^j_{\tau,k}}
\int_QD_Q(\cdot,y)\overline{\mathcal{D}}_Qf(y)\,d\mu(y)
\nonumber\\
&=\sum_{k=-\infty}^{\infty}
\sum_{j=-\infty}^{\infty}
\sum_{\tau\in I_{j,k}}
\lambda^j_{\tau,k}a^j_{\tau,k}(\cdot),
\end{align}
where all the equalities hold
in $(\mathring{{\mathcal{G}}}_0^\varepsilon(\beta,\gamma))'$
and, for any $k,j\in\mathbb{Z}$, $\tau\in I_{j,k}$, and $x\in X$,
\begin{align}\label{2203}
\lambda^j_{\tau,k}
:=\left[\mu(Q^j_{\tau,k})\right]^{\frac{1}{p}-\frac{1}{2}}
\left[\sum_{Q\in\mathfrak{Q}_k,\,Q\subset Q^j_{\tau,k}}
\int_Q\left|\overline{\mathcal{D}}_Qf(y)\right|^2\,d\mu(y)\right]^\frac{1}{2}
\end{align}
and
\begin{align}\label{2201}
a^j_{\tau,k}(x)
:=\frac{1}{\lambda^j_{\tau,k}}
\sum_{Q\in\mathfrak{Q}_k,\,Q\subset Q^j_{\tau,k}}
\int_QD_Q(x,y)\overline{\mathcal{D}}_Qf(y)\,d\mu(y).
\end{align}

Next, we aim to show that
\begin{align}\label{2200}
\sum_{k\in\mathbb{Z}}\sum_{j\in\mathbb{Z}}
\sum_{\tau\in I_{j,k}}\left|\lambda^j_{\tau,k}\right|^p
\lesssim\left\|S_\mathcal{D}(f)\right\|^p_{L^{p}(X)}.
\end{align}
Indeed,
if $Q\in\mathfrak{Q}_k$
with $k\neq-\infty$, then there exists a unique $j\in\mathbb{Z}$ such
that $Q\subset Q^j_{\tau,k}$.
Assume that $Q=Q^{l+k_0}_\alpha$
for some $l\in\mathbb{Z}\cap[j-k_0,\infty)$
and $\alpha\in I_{l+k_0}$.
For any given $y\in Q$, by Lemma~\ref{DyadicSys}(iv) and
an argument similar to that used in the proof of \eqref{1407},
we obtain
$Q=Q^{l+k_0}_\alpha\subset B_\rho(y,2^{-l})$.
From this, the fact that
$\mu(Q\cap\Omega_{k+1})\leq\frac{\mu(Q)}{2}$,
both (iv) and (v) of Lemma~\ref{DyadicSys},
and \eqref{AR},
we infer that
\begin{align}\label{2327}
\mu\left(B_\rho(y,2^{-l})
\cap[Q\setminus\Omega_{k+1}]\right)
=\mu(Q\setminus\Omega_{k+1})
\ge\frac{1}{2}\mu(Q)\sim (V_\rho)_{2^{-l}}(y).
\end{align}
By this,
we conclude that
\begin{align*}
&\sum_{Q\in\mathfrak{Q}_k,\,Q\subset Q^j_{\tau,k}}
\int_Q\left|\overline{\mathcal{D}}_Qf(y)\right|^2\,d\mu(y)
\nonumber\\
&\quad\lesssim\sum_{l=j-k_0}^{\infty}
\sum_{\genfrac{}{}{0pt}{}{\alpha\in I_{l+k_0}}{
\mathfrak{Q}_k\ni Q^{l+k_0}_\alpha\subset Q^j_{\tau,k}}}
\int_{Q^{l+k_0}_\alpha}
\frac{\mu(B_\rho(y,2^{-l})\cap
[Q^{l+k_0}_\alpha\setminus\Omega_{k+1}])}
{(V_\rho)_{2^{-l}}(y)}\left|\overline{\mathcal{D}}_lf(y)\right|^2\,d\mu(y)
\nonumber\\
&\quad\leq\sum_{l=j-k_0}^{\infty}
\sum_{\genfrac{}{}{0pt}{}{\alpha\in I_{l+k_0}}{\mathfrak{Q}_k\ni Q^{l+1}_\alpha\subset Q^j_{\tau,k}}}
\int_{Q^{l+k_0}_\alpha}
\frac{\mu(B_\rho(y,2^{-l})\cap
[Q^j_{\tau,k}\setminus\Omega_{k+1}])}
{(V_\rho)_{2^{-l}}(y)}\left|\overline{\mathcal{D}}_lf(y)\right|^2\,d\mu(y),
\end{align*}
which, combined with Tonelli's theorem and Lemma~\ref{A3.2}(vi),
further implies that
\begin{align*}
&\sum_{Q\in\mathfrak{Q}_k,\,Q\subset Q^j_{\tau,k}}
\int_Q\left|\overline{\mathcal{D}}_Qf(y)\right|^2\,d\mu(y)
\nonumber\\
&\quad
\lesssim\sum_{l=j-k_0}^{\infty}
\int_{Q^j_{\tau,k}}
\frac{\mu(B_\rho(y,2^{-l})\cap
[Q^j_{\tau,k}\setminus\Omega_{k+1}])}
{(V_\rho)_{2^{-l}}(y)}\left|\overline{\mathcal{D}}_lf(y)\right|^2\,d\mu(y)
\nonumber\\
&\quad=
\int_X\sum_{l=j-k_0}^{\infty}
\int_{B_\rho(x,2^{-l})}
\left|\overline{\mathcal{D}}_lf(y)\right|^2
\frac{d\mu(y)}{(V_\rho)_{2^{-l}}(x)}
\mathbf{1}_{[Q^j_{\tau,k}\setminus\Omega_{k+1}]}(x)\,d\mu(x)
\nonumber\\
&\quad
\leq\int_{Q^j_{\tau,k}
\setminus\Omega_{k+1}}\left[S_{\overline{\mathcal{D}}}(f)(x)\right]^2\,d\mu(x)
\lesssim2^{2k}\mu(Q^j_{\tau,k}\setminus\Omega_{k+1})
\leq2^{2k}\mu(Q^j_{\tau,k}).
\end{align*}
This, together with the definition of $\lambda^j_{\tau,k}$,
further implies that
\begin{align*}
\left|\lambda^j_{\tau,k}\right|
\lesssim2^k\left[\mu(Q^j_{\tau,k})\right]^{\frac{1}{p}}.
\end{align*}
From this, the facts that
$\mu(Q^j_{\tau,k})<2\mu(Q^j_{\tau,k}\cap\Omega_k)$
for any $k\in\mathbb{Z}$
and the cubes $\{Q^j_{\tau,k}\}_{j\in\mathbb{Z},\,\tau\in I_{j,k}}$
are mutually disjoint, the definition of $\Omega_k$ with $k\in\mathbb{Z}$,
and Lemma~\ref{2157},
we deduce that
\begin{align}\label{lammdaa}
&\sum_{k\in\mathbb{Z}}\sum_{j\in\mathbb{Z}}
\sum_{\tau\in I_{j,k}}\left|\lambda^j_{\tau,k}\right|^p\nonumber\\
&\quad\lesssim\sum_{k\in\mathbb{Z}}2^{kp}
\sum_{j\in\mathbb{Z}}\sum_{\tau\in I_{j,k}}
\mu(Q^j_{\tau,k})
\lesssim\sum_{k\in\mathbb{Z}}2^{kp}
\sum_{j\in\mathbb{Z}}
\sum_{\tau\in I_{j,k}}\mu(Q^j_{\tau,k}\cap\Omega_k)
\nonumber\\
&\quad\leq\sum_{k\in\mathbb{Z}}2^{kp}
\mu(\Omega_k)
\sim||S_{\overline{\mathcal{D}}}(f)||^p_{L^p(X)}
\lesssim\|S_\mathcal{D}(f)\|^p_{L^{p}(X)},
\end{align}
which completes the proof of \eqref{2200}.

Now, we prove that, for any $k,j\in\mathbb{Z}$
and $\tau\in I_{j,k}$, $a_{\tau,k}^j$ is
a harmless constant multiple of a $(p,2)$-atom.
To this end, let $k,j\in\mathbb{Z}$, $\tau\in I_{j,k}$,
and $B_{\tau,k}^j:=B_\rho(z_{\tau,k}^j,2^{-j+k_0})$.
We first show the size estimate of $a_{\tau,k}^j$ as in Definition~\ref{atom}(ii).
Note that, by both (iv) and (v) of Lemma~\ref{DyadicSys},
we find that
\begin{align}\label{BsubsetQ}
C_l2^{-k_0}B_{\tau,k}^j\subset Q_{\tau,k}^j
\subset B_{\tau,k}^j
\end{align}
with $C_l$ as in Lemma~\ref{DyadicSys}(v),
which, combined with \eqref{upperdoub},
further implies
\begin{align}
\label{BsimQ}
\mu(B_{\tau,k}^j)
\sim\mu(Q_{\tau,k}^j).
\end{align}
From this, \eqref{2201}, \eqref{2203},
and Fubini's theorem,
we infer that,
for any $k,j\in\mathbb{Z}$ and $\tau\in I_{j,k}$,
\begin{align*}
\left\|a_{\tau,k}^j\right\|_{L^2(X)}
&=\sup_{h\in L^2(X),\,\|h\|_{L^2(X)}\leq1}
\left|\int_Xa_{\tau,k}^j(x)h(x)\,d\mu(x)\right|
\nonumber\\
&\leq\frac{1}{\lambda_{\tau,k}^j}
\sup_{h\in L^2(X),\,\|h\|_{L^2(X)}\leq1}
\sum_{Q\in\mathfrak{Q}_k,\,Q\subset Q^j_{\tau,k}}\nonumber\\
&\quad\times\int_Q\left|\int_XD_Q(x,y)h(x)\,d\mu(x)\right|
\left|\overline{\mathcal{D}}_Qf(y)\right|\,d\mu(y)
\nonumber\\
&=\frac{1}{\lambda_{\tau,k}^j}
\sup_{h\in L^2(X),\,\|h\|_{L^2(X)}\leq1}
\sum_{Q\in\mathfrak{Q}_k,\,Q\subset Q^j_{\tau,k}}
\int_Q\left|\mathcal{D}_Qh(y)\right|
\left|\overline{\mathcal{D}}_Qf(y)\right|\,d\mu(y),
\end{align*}
which, together with H\"older's inequality twice, Remark~\ref{2057} with $p=2$,
and \eqref{BsimQ},
further implies that
\begin{align}\label{atkjs}
\left\|a_{\tau,k}^j\right\|_{L^2(X)}
&\leq\frac{1}{\lambda_{\tau,k}^j}
\sup_{h\in L^2(X),\,\|h\|_{L^2(X)}\leq1}
\left[\sum_{Q\in\mathfrak{Q}_k,\,Q\subset
Q^j_{\tau,k}}
\int_Q\left|\mathcal{D}_Qh(y)\right|^2\,d\mu(y)\right]^\frac{1}{2}
\nonumber\\
&\quad\times
\left[\sum_{Q\in\mathfrak{Q}_k,\,Q\subset
Q^j_{\tau,k}}\int_Q\left|\overline{\mathcal{D}}_Qf(y)\right|
^2\,d\mu(y)\right]^\frac{1}{2}
\nonumber\\
&=\left[\mu(Q_{\tau,k}^j)\right]^{\frac{1}{2}-\frac{1}{p}}
\sup_{h\in L^2(X),\,\|h\|_{L^2(X)}\leq1}
\left\|\overline{g}(h)\right\|_{L^2(X)}
\lesssim
\left[\mu(B_{\tau,k}^j)\right]^{\frac{1}{2}-\frac{1}{p}}
\end{align}
with $\overline{g}$ as in Remark~\ref{2057}.
Thus, $a_{\tau,k}^j$ satisfies the desired size estimate as in Definition~\ref{atom}(ii).

Then we prove that,
for any given $k,j\in\mathbb{Z}$
and $\tau\in I_{j,k}$,
$a_{\tau,k}^j$ is
supported in some ball centered in $z_{\tau,k}^j$.
For any $Q\in\mathfrak{Q}_k$ and $Q\subset Q_{\tau,k}^j$,
there exist $l\in[j-k_0,\infty)$
and $\tau'\in I_{l+k_0}$ (both depending on $Q$) such that
$Q=Q^{l+k_0}_{\tau'}$. On the other hand, by Lemma~\ref{Pk},
we conclude that
there exists a positive constant $C$ such that, for any $x,y\in X$ and $l\in\mathbb{Z}$,
if $\rho(x,y)\geq C2^{-l}$, then $D_Q(x,y)=D_l(x,y)=0$.
Note that we may assume that $C\geq C_r$.
Based on these observations, we have, for any given $x\in X$,
\begin{align}
\label{1727}
\left|a_{\tau,k}^j(x)\right|
&\leq\frac{1}{\lambda^j_{\tau,k}}
\sum_{Q\in\mathfrak{Q}_k,\,Q\subset Q^j_{\tau,k}}
\int_Q\left|D_Q(x,y)\overline{\mathcal{D}}_Qf(y)\right|\,d\mu(y)
\nonumber\\
&\leq\frac{1}{\lambda^j_{\tau,k}}
\sum_{Q\in\mathfrak{Q}_k,\,Q\subset Q^j_{\tau,k}}
\int_{Q^j_{\tau,k}\cap B_\rho(x,C2^{-l})}
\left|D_Q(x,y)\overline{\mathcal{D}}_Qf(y)\right|\,d\mu(y).
\end{align}
We consider the following two cases for $x\in X$.

\emph{Case (1)} $x\in X$ satisfies that
$Q_{\tau,k}^j\cap B_\rho(x,C2^{-j+k_0})=\varnothing$.
In this case,
for any $l\in[j-k_0,\infty)$, we have
$Q^j_{\tau,k}\cap B_\rho(x,C2^{-l})=\varnothing$,
which, combined with \eqref{1727},
implies that $a_{\tau,k}^j(x)=0$.

\emph{Case (2)}
$x\in X$ satisfies that $Q_{\tau,k}^j\cap B_\rho(x,C2^{-j+k_0})\neq\varnothing$.
In this case,
we infer that there exists a point $y_0\in Q_{\tau,k}^j\cap B_\rho(x,C2^{-j+k_0})$,
which, together with Definition~\ref{D2.1}(iii)
and $Q_{\tau,k}^j\subset B_\rho(z_{\tau,k}^j,C_r2^{-j})$,
further implies that
\begin{align*}
\rho(x,z_{\tau,k}^j)
&\leq C_\rho\max\left\{\rho(x,y_0),\,\rho(y_0,z_{\tau,k}^j)\right\}\\
&\leq C_\rho\max\left\{C2^{-j+k_0},\,C_r2^{-j}\right\}
=C_\rho C2^{-j+k_0}
\end{align*}
(here we used the choice of $k_0$ that $C_\rho C_r\leq 2^{k_0}$)
and hence
\begin{align}\label{twh-349}
\mathrm{\,supp\,} a_{\tau,k}^j\subset C_\rho CB_{\tau,k}^j.
\end{align}

Finally, we show that $a_{\tau,k}^j$ satisfies Definition~\ref{atom}(iii). To this end,
we first prove that
\begin{align}\label{2133}
&\sum_{Q\in\mathfrak{Q}_k,\,Q\subset Q^j_{\tau,k}}
\int_QD_Q(\cdot,y)\overline{\mathcal{D}}_Qf(y)\,d\mu(y)\nonumber\\
&\quad=\sum_{l=j-k_0}^\infty\sum_{\alpha\in I_l,\,Q_\alpha^l
\in\mathfrak{Q}_k,\,Q_\alpha^l\subset Q^j_{\tau,k}}
\int_{Q_\alpha^l}D_l(\cdot,y)\overline{\mathcal{D}}_lf(y)\,d\mu(y)
\end{align}
converges in $L^q(X)$ for any given $q\in(1,\infty)$
in the sense of weak-$*$ topology. Indeed, by Lemma~\ref{Pk}(i),
Lemma~\ref{A3.2}(vii), Lemma~\ref{DyadicSys}(iv),
and an argument similar to that used in the proof of
\eqref{twh-349},
we find that,
for any $g\in\bigcup_{u\in[1,\infty]}L^u(X)$ and
$Q_\alpha^l\in\mathfrak{Q}_k$ satisfying $Q_\alpha^l\subset Q^j_{\tau,k}$
for some $l\in[j-k_0,\infty)\cap\mathbb{Z}$ and $\alpha\in I_l$,
\begin{align*}
&\int_X\int_{Q_\alpha^l}\left|D_l(x,y)
\overline{\mathcal{D}}_lf(y)\right||g(x)|\,d\mu(y)\,d\mu(x)\nonumber\\
&\quad\lesssim\int_X\int_{Q_\alpha^l}\frac{1}{V_{2^{-l}}(y)}
\mathbf{1}_{B_\rho(y,C2^{-l})}(x)\left|
\overline{\mathcal{D}}_lf(y)\right||g(x)|\,d\mu(y)\,d\mu(x)\nonumber\\
&\quad\lesssim\int_X\int_{Q_\alpha^l}\frac{1}{V_{2^{-l}}(z_\alpha^l)}
\mathbf{1}_{C_\rho CB_{\tau,k}^j}(x)\left|\overline{\mathcal{D}}_lf(y)\right|
|g(x)|\,d\mu(y)\,d\mu(x),
\end{align*}
which, combined with H\"older's inequality and
an argument similar to that used in the estimation of \eqref{2255},
further implies that
\begin{align}\label{220}
&\int_X\int_{Q_\alpha^l}\left|D_l(x,y)
\overline{\mathcal{D}}_lf(y)\right||g(x)|\,d\mu(y)\,d\mu(x)\nonumber\\
&\quad\leq\frac{1}{V_{2^{-l}}(z_\alpha^l)}
\left\|\mathbf{1}_{C_\rho CB_{\tau,k}^j}\right\|_{L^{u'}(X)}\left\|g\right\|_{L^u(X)}
\int_{Q_\alpha^l}\left|\overline{\mathcal{D}}_lf(y)\right|\,d\mu(y)
\nonumber\\
&\quad\leq\frac{1}{V_{2^{-l}}(z_\alpha^l)}
\left\|\mathbf{1}_{C_\rho CB_{\tau,k}^j}\right\|_{L^{u'}(X)}\left\|g\right\|_{L^u(X)}
\left\|f\right\|_{(\mathring{\mathcal{G}}^\varepsilon_0(\beta,\gamma))'}
\int_{Q_\alpha^l}\left\|\overline{D}_l(y,\cdot)
\right\|_{{\mathcal{G}}^\varepsilon_0(\beta,\gamma)}\,d\mu(y)\nonumber\\
&\quad<\infty.
\end{align}
For any $N_1,N_2\in\mathbb{Z}$ with $j-k_0\leq N_1\leq N_2<\infty$
and $h\in L^{q'}(X)$, where $q'\in(1,\infty)$
is the conjugate index of $q$, let
\begin{align*}
A_{N_1,N_2}(h):=\sum_{l=N_1}^{N_2}\sum_{\alpha\in I_l,\,Q_\alpha^l\in
\mathfrak{Q}_k,\,Q_\alpha^l\subset Q^j_{\tau,k}}
\int_X\int_{Q_\alpha^l}D_l(x,y)\overline{\mathcal{D}}_lf(y)\,d\mu(y)h(x)\,d\mu(x).
\end{align*}
From
Fubini's theorem, \eqref{220}, Lemma~\ref{Pk}(iv), \eqref{2327},
and Lemma~\ref{A3.2}(vii),
we deduce that
\begin{align*}
\left|A_{N_1,N_2}(h)\right|
&\leq\sum_{l=N_1}^{N_2}\sum_{\alpha\in I_l,\,Q_\alpha^l\in
\mathfrak{Q}_k,\,Q_\alpha^l\subset Q^j_{\tau,k}}
\int_{Q_\alpha^l}\left|\mathcal{D}_lh(y)\overline{\mathcal{D}}_lf(y)\right|\,d\mu(y)\nonumber\\
&\lesssim\sum_{l=N_1}^{N_2}\sum_{\alpha\in I_l,\,Q_\alpha^l\in
\mathfrak{Q}_k,\,Q_\alpha^l\subset Q^j_{\tau,k}}
\int_{Q_\alpha^l}\frac{\mu(B_\rho(y,2^{-l})\cap(Q_{\tau,k}^j\setminus
\Omega_{k+1}))}{(V_\rho)_{2^{-l}}(y)}\nonumber\\
&\quad\times\left|\mathcal{D}_lh(y)\overline{\mathcal{D}}_lf(y)\right|\,d\mu(y)\nonumber\\
&\leq\sum_{l=N_1}^{N_2}
\int_{Q_{\tau,k}^j}\frac{\mu(B_\rho(y,2^{-l})\cap(Q_{\tau,k}^j
\setminus\Omega_{k+1}))}{(V_\rho)_{2^{-l}}(y)}
\left|\mathcal{D}_lh(y)\overline{\mathcal{D}}_lf(y)\right|\,d\mu(y)\nonumber\\
&\sim\sum_{l=N_1}^{N_2}
\int_{Q_{\tau,k}^j}
\int_{Q_{\tau,k}^j\setminus\Omega_{k+1}}
\frac{\mathbf{1}_{B_\rho(y,2^{-l})}(x)}{(V_\rho)_{2^{-l}}(x)}\,d\mu(x)
\left|\mathcal{D}_lh(y)\overline{\mathcal{D}}_lf(y)\right|\,d\mu(y),
\end{align*}
which, together with Tonelli's theorem, H\"older's inequality,
and Theorem~\ref{Hp=Lp},
further implies that
\begin{align}\label{2026}
\left|A_{N_1,N_2}(h)\right|
&=\int_{Q_{\tau,k}^j\setminus\Omega_{k+1}}\sum_{l=N_1}^{N_2}
\int_{Q_{\tau,k}^j\cap B_\rho(x,2^{-l})}
\left|\mathcal{D}_lh(y)\overline{\mathcal{D}}_lf(y)\right|
\frac{d\mu(y)}{(V_\rho)_{2^{-l}}(x)}\,d\mu(x)\nonumber\\
&\leq\left\{\int_{Q_{\tau,k}^j\setminus\Omega_{k+1}}\left[\sum_{l=N_1}^{N_2}
\int_{Q_{\tau,k}^j\cap B_\rho(x,2^{-l})}
\left|\overline{\mathcal{D}}_lf(y)\right|^2
\frac{d\mu(y)}{(V_\rho)_{2^{-l}}(x)}\right]^\frac{q}{2}\,d\mu(x)\right\}^\frac{1}{q}\nonumber\\
&\quad\times\left\{\int_{Q_{\tau,k}^j\setminus\Omega_{k+1}}\left[\sum_{l=N_1}^{N_2}
\int_{Q_{\tau,k}^j\cap B_\rho(x,2^{-l})}
\left|\mathcal{D}_lh(y)\right|^2\frac{d\mu(y)}{(V_\rho)_{2^{-l}}(x)}
\right]^\frac{q'}{2}\,d\mu(x)\right\}^\frac{1}{q'}\nonumber\\
&\leq\left\{\int_{Q_{\tau,k}^j\setminus\Omega_{k+1}}
\left[S_{\overline{\mathcal{D}}}f(x)\right]^q\,d\mu(x)
\right\}^\frac{1}{q}\left\|S_{\mathcal{D}}(h)\right\|_{L^{q'}(X)}\nonumber\\
&\lesssim2^k\left[\mu(Q_{\tau,k}^j)\right]^\frac{1}{q}\|h\|_{L^{q'}(X)}
\end{align}
and hence $A_{N_1,N_2}$ is a bounded linear functional on $L^{q'}(X)$.
By the fact that
$$\|S_{\mathcal{D}}(h)\|_{L^{q'}(X)}\sim\|h\|_{L^{q'}(X)}<\infty$$
(which can be deduced from Theorem~\ref{Hp=Lp})
and the Lebesgue dominated convergence theorem,
we conclude that, for any given $\delta\in(0,\infty)$,
there exists $N_{\delta}\in\mathbb{Z}\cap[j-k_0,\infty)$
such that, for any $N_1,N_2\in\mathbb{Z}\cap[N_\delta,\infty)$,
\begin{align*}
\left\{\int_{Q_{\tau,k}^j\setminus\Omega_{k+1}}\left[\sum_{l=N_1}^{N_2}
\int_{Q_{\tau,k}^j\cap B_\rho(x,2^{-l})}
\left|\mathcal{D}_lh(y)\right|^2\frac{d\mu(y)}{(V_\rho)_{2^{-l}}(x)}
\right]^\frac{q'}{2}\,d\mu(x)\right\}^\frac{1}{q'}<\delta,
\end{align*}
which, combined with the penultimate inequality in \eqref{2026}, further implies that
$$
\left|A_{N_1,N_2}(h)\right|\lesssim2^k\left[\mu(Q_{\tau,k}^j)\right]^\frac{1}{q}\delta.
$$
Thus, $\{A_{j-k_0,N}(h)\}_{N=j-k_0}^\infty$ is a
Cauchy sequence and $\lim_{N\to\infty}A_{j-k_0,N}(h)$
exists, which, together with the uniform boundedness principle
(see, for instance, \cite[Theorem~2.5]{r91}),
completes the proof of \eqref{2133}. From this, \eqref{twh-349}, \eqref{220}
with $g:=\mathbf{1}_{C_\rho CB_{\tau,k}^j}\in L^{q'}(X)$, and Lemma~\ref{Pk}(v),
we infer that, for any $k,j\in\mathbb{Z}$ and $\tau\in I_{j,k}$,
\begin{align}
\label{atkjc}
&\int_Xa_{\tau,k}^j(x)\,d\mu(x)\nonumber\\
&\quad=\int_X\frac{1}{\lambda^j_{\tau,k}}
\sum_{Q\in\mathfrak{Q}_k,\,Q\subset Q^j_{\tau,k}}
\int_QD_Q(x,y)\overline{\mathcal{D}}_Qf(y)\,d\mu(y)\,d\mu(x)
\nonumber\\
&\quad=\frac{1}{\lambda^j_{\tau,k}}\int_X
\left[\sum_{Q\in\mathfrak{Q}_k,\,Q\subset Q^j_{\tau,k}}
\int_QD_Q(x,y)\overline{\mathcal{D}}_Qf(y)\,d\mu(y)\right]
\mathbf{1}_{C_\rho CB_{\tau,k}^j}(x)\,d\mu(x)
\nonumber\\
&\quad=\frac{1}{\lambda^j_{\tau,k}}
\sum_{Q\in\mathfrak{Q}_k,\,Q\subset Q^j_{\tau,k}}\int_X
\int_QD_Q(x,y)\overline{\mathcal{D}}_Qf(y)\,d\mu(y)
\mathbf{1}_{C_\rho CB_{\tau,k}^j}(x)\,d\mu(x)
\nonumber\\
&\quad=\frac{1}{\lambda^j_{\tau,k}}
\sum_{Q\in\mathfrak{Q}_k,\,Q\subset Q^j_{\tau,k}}
\int_Q\left[\int_XD_Q(x,y)\mathbf{1}_{C_\rho CB_{\tau,k}^j}(x)
\,d\mu(x)\right]\overline{\mathcal{D}}_Qf(y)\,d\mu(y)
=0,
\end{align}
where, in the penultimate inequality, we used the fact
that $\mathrm{supp\,}D_Q(\cdot,y)\subset C_\rho CB_{\tau,k}^j$
[which can be deduced from the proof of \eqref{twh-349}].
Thus, $a_{\tau,k}^j$ satisfies Definition~\ref{atom}(iii).

Note that, by \eqref{atkjs}, \eqref{upperdoub},
and $p<2$, we obtain
$$
\left\|a_{\tau,k}^j\right\|_{L^2(X)}
\lesssim\left[\mu( B_{\tau,k}^j)\right]^{\frac{1}{2}-\frac{1}{p}}
\lesssim\left[\mu\left(C_\rho CB_{\tau,k}^j\right)\right]^{\frac{1}{2}-\frac{1}{p}}.
$$
From this, \eqref{twh-349}, and \eqref{atkjc},
we deduce that, for any $k,j\in\mathbb{Z}$
and $\tau\in I_{j,k}$, $a_{\tau,k}^j$ is
a harmless constant multiple of a $(p,2)$-atom,
which, combined with \eqref{lammdab} and \eqref{lammdaa},
further implies that $f\in\mathring{H}^{p,2}_{\mathrm{at}}(X)$.
This finishes the proof of Theorem~\ref{L4.10}.
\end{proof}

\begin{remark}
In the proof of Theorem~\ref{L4.10},
we decompose the elements of $H^{p}(X)$ into the
linear combination of atoms,
while He et al. \cite[Proposition 5.7]{ZYJLDW}
only decomposed that into a
linear combination of molecules and then used
wavelet reproducing formulae to first obtain the molecular
characterization of $H^p(X)$.
This is because the approximation of the
identity $\{\mathcal{S}_{2^{-k}}\}_{k\in\mathbb{Z}}$
have bounded supports [see Definition~\ref{DefATI}(i) above] in this monograph
while that in \cite[Definition 2.7]{ZYJLDW} may not have bounded supports.
\end{remark}

Using Theorems~\ref{Sf} and~\ref{L4.10},
we obtain the following conclusion.

\begin{theorem}\label{329}
Let $(X,\mathbf{q},\mu)$ be an ultra-RD-space
with upper dimension $n\in(0,\infty)$. Suppose that
$p\in(\frac{n}{n+\mathrm{ind\,}(X,\mathbf{q})},1]$ and
$\varepsilon,\beta,\gamma\in(0,\infty)$ satisfy
\eqref{1945}.
Then, as subspaces of $(\mathring{{\mathcal{G}}}_0^\varepsilon(\beta,\gamma))'$,
$\mathring{H}^{p,2}_{\mathrm{at}}(X)=H^p(X)$
with equivalent (quasi-)norms.
\end{theorem}

According to Theorems~\ref{328} and~\ref{329},
we immediately obtain the following result.

\begin{theorem}\label{Thm4.11}
Let $(X,\mathbf{q},\mu)$ be an ultra-RD-space
with upper dimension $n\in(0,\infty)$. Suppose that
$p\in(\frac{n}{n+\mathrm{ind\,}(X,\mathbf{q})},1]$, $q\in[1,\infty]\cap(p,\infty]$, and
$\varepsilon,\beta,\gamma\in(0,\infty)$ satisfy
\eqref{1945}.
Then, as subspaces of $(\mathring{{\mathcal{G}}}_0^\varepsilon(\beta,\gamma))'$,
$\mathring{H}^{p,q}_{\mathrm{at}}(X)=H^p(X)$
with equivalent (quasi-)norms.
\end{theorem}

\begin{remark}\label{2548}
Alvarado and Mitrea in \cite{AlMi2015} considered Hardy spaces on
$n$-Ahlfors-regular quasi-ultrametric spaces
which were based on a radial maximal function
[denoted by $H^{p}_{\rm rad}(X)$],
a nontangential maximal function (denoted by $H^{p}_{\rm n.t.}(X)$),
and a grand maximal function (denoted by $H^{p}_{\rm max}(X)$),
and also considered two Hardy spaces based on atoms
(denoted by $H^{p,q}_{\mathrm{atom}}(X)$
which is showed to coincide with $H^{p,q}_{\mathrm{CW}}(X)$;
see \cite[Theorem~7.5]{AlMi2015})
and molecules (denoted by $H^{p,q}_{\rm mol}(X)$).
As proved in \cite[Theorem~6.11]{AlMi2015} that
all of these Hardy spaces coincide with each other
for any $p\in(\frac{n}{n+\mathrm{ind\,}(X,\mathbf{q})},1]$.
As a further conclusion of \cite[Theorems~6.11 and~7.5]{AlMi2015}
and Theorem~\ref{Thm4.11}, we find that,
for any given $p\in(\frac{n}{n+\mathrm{ind\,}(X,\mathbf{q})},1]$
and $q\in(p,\infty]\cap[1,\infty]$,
\begin{align*}
H^{p}(X)
&=\mathring{H}^{p,q}_{\mathrm{at}}(X)
=H^{p,q}_{\mathrm{at}}(X)=H^{p,q}_{\mathrm{CW}}(X)
\\
&=H^{p,q}_{\mathrm{atom}}(X)
=H^{p,q}_{\rm mol}(X)=H^{p}_{\rm max}(X)
=H^{p}_{\rm rad}(X)=H^{p}_{\rm n.t.}(X)
\end{align*}
as subspaces of $(\dot{\mathscr{D}}_\varepsilon(X))'$,
with equivalent (quasi-)norms,
where
$$
\varepsilon\in\left(n\left(\frac{1}{p}-1\right),\mathrm{ind\,}(X,\mathbf{q})\right)
\ \text{and}\
\dot{\mathscr{D}}_\varepsilon(X):=
\bigcap_{\alpha\in(0,\varepsilon)}C_{\mathrm{b}}^\alpha(X).
$$
\end{remark}

\begin{remark}\label{20071}
By Theorem~\ref{Thm4.11}, Proposition~\ref{324}, Lemma~\ref{325},
Definition~\ref{1539},
and Theorem~\ref{Hp=Lp},
we conclude that $H^p(X)$
with $p\in(\frac{n}{n+\mathrm{ind\,}(X,\mathbf{q})},\infty)$ is independent of
the choice of the space
$(\mathring{{\mathcal{G}}}_0^\varepsilon(\beta,\gamma))'$ of distributions whenever
$\varepsilon,\beta,\gamma\in(0,\infty)$ satisfy \eqref{1945}.
\end{remark}

\section{Littlewood--Paley Function Characterizations}
\label{Section3.3}

In this section, we separately characterize
$H^p(X)$ by the Littlewood--Paley $g$-function
and $g^*_\lambda$-function
defined, respectively, in \eqref{Def-g(f)(x)} and \eqref{Def-glamda*}.
The main results of this section are
Theorems~\ref{gf} and \ref{Glamda*} below.

We first consider
the Littlewood--Paley $g$-function characterization
of $H^p(X)$.
More specifically,
we have the following conclusion.

\begin{theorem}\label{gf}
Let $(X,\mathbf{q},\mu)$ be an ultra-RD-space
with upper dimension $n\in(0,\infty)$.
Let $p\in(\frac{n}{n+\mathrm{ind\,}(X,\mathbf{q})},1]$,
$\beta,\gamma,\varepsilon\in(0,\infty)$
satisfy \eqref{1945},
and $\{\mathcal{S}_{2^{-k}}\}_{k\in\mathbb{Z}}$
be a homogeneous approximation
of the identity of order $\varepsilon$
in Definition~\ref{DefHATI}.
For any $k\in\mathbb{Z}$, let
$\mathcal{D}_k:=\mathcal{S}_{2^{-k}}-\mathcal{S}_{2^{-k+1}}$.
Let $g$ be the same as in \eqref{Def-g(f)(x)}.
Then $f\in H^p(X)$ if and only if
$f\in(\mathring{{\mathcal{G}}}_0^\varepsilon(\beta,\gamma))'$
and $g(f)\in L^p(X)$.
Moreover,
there exists a constant $C\in[1,\infty)$ such that,
for any $f\in H^p(X)$,
\begin{align*}
C^{-1}||g(f)||_{L^p(X)}
\leq||f||_{H^p(X)}
\leq C||g(f)||_{L^p(X)}.
\end{align*}
\end{theorem}

\begin{remark}
By Remark~\ref{2548}, we find that,
when $(X,\mathbf{q},\mu)$ is an Ahlfors-regular
quasi-ultrametric space as in \cite[Definition~2.11]{AlMi2015},
Theorem~\ref{gf} gives the
Littlewood--Paley $g$-function characterization
of the Hardy space in \cite{AlMi2015}.
\end{remark}

\begin{proof}[Proof of Theorem~\ref{gf}]
Let $\rho\in\mathbf{q}$ be a regularized quasi-ultrametric.
We first show the necessity. To this end,
let $f\in H^p(X)$.
Then Definition~\ref{D3.1.4} and Remark~\ref{20071} imply
that $f\in(\mathring{{\mathcal{G}}}_0^\varepsilon(\beta,\gamma))'$. Moreover,
from arguments similar to those used in the proofs
of both Lemma~\ref{Sa2} and Theorem~\ref{Sf}
and from Theorem~\ref{Thm4.11},
we deduce that, for any $f\in H^p(X)$,
$$
\left\|g(f)\right\|_{L^p(X)}
\lesssim\|f\|_{\mathring{H}^{p,q}_{\mathrm{at}}(X)}
\sim||f||_{H^p(X)},
$$
where $q\in[1,\infty]\cap(p,\infty]$ is the same as in Theorem~\ref{Thm4.11}.

Next, we prove the sufficiency. To this end,
let $f\in(\mathring{{\mathcal{G}}}_0^\varepsilon(\beta,\gamma))'$
and $\|g(f)\|_{L^p(X)}<\infty$.
By Definition~\ref{D2.1}(iii) and Lemma~\ref{DyadicSys}(iv), we find that,
for any fixed $k\in\mathbb{Z}$, $\tau\in I_k$, $v\in\{1,...,N(k,\tau)\}$, and
$x\in Q_{\tau}^{k,v}$
and for any $y\in B_\rho(x,2^{-k})$,
\begin{align*}
\rho(z_{\tau}^{k,v},y)
\leq C_\rho\max\left\{\rho(z_{\tau}^{k,v},x),\,\rho(x,y)\right\}
\leq C_\rho\max\left\{\mathrm{diam}_\rho(Q_{\tau}^{k,v}),\,2^{-k}\right\}
\leq C_\rho2^{-k}
\end{align*}
with both $Q_{\tau}^{k,v}$ and $z_{\tau}^{k,v}$ as in Remark~\ref{jcubes},
which implies that
$$
B_\rho\left(x,2^{-k}\right)\subset
B_\rho\left(z_{\tau}^{k,v},C_\rho2^{-k}\right).
$$
From this, Lemma~\ref{DyadicSys}(i), and \eqref{un-3iu},
we infer that, for $\mu$-almost every $x\in X$,
\begin{align}
\label{1501}
S(f)(x)
&=\left[\sum_{k\in\mathbb{Z}}\sum_{\tau\in I_k}
\sum_{v=1}^{N(k,\tau)}
\int_{B_\rho(x,2^{-k})}
\left|\mathcal{D}_kf(y)\right|^2
\frac{d\mu(y)}{(V_\rho)_{2^{-k}}(x)}
\mathbf{1}_{Q^{k,v}_\tau}(x)\right]^\frac{1}{2}
\nonumber\\
&\leq\left\{\sum_{k\in\mathbb{Z}}\sum_{\tau\in I_k}
\sum_{v=1}^{N(k,\tau)}
\left[\sup_{z\in B_\rho(z^{k,v}_\tau,C_\rho2^{-k})}
\left|\mathcal{D}_kf(z)\right|^2\right]
\mathbf{1}_{Q^{k,v}_\tau}(x)\right\}^\frac{1}{2}.
\end{align}

Now, we first estimate $|\mathcal{D}_kf(z)|$ for any given $k\in\mathbb{Z}$
and for any $z\in B_\rho(z^{k,v}_\tau,C_\rho2^{-k})$.
By Theorem~\ref{HD} (keeping in mind
that we are assuming $\mathrm{diam}_\rho(X)=\infty$),
we conclude that there exists a family of linear operators
$\{\widetilde{\mathcal{D}}_{k'}\}_{k'\in\mathbb{Z}}$
such that the third equality in \eqref{3.29}
holds in $(\mathring{{\mathcal{G}}}_0^\varepsilon(\beta,\gamma))'$.
As such,
for any given $k\in\mathbb{Z}$ and for any $z\in B_\rho(z^{k,v}_\tau,C_\rho2^{-k})$,
\begin{align}
\label{1502}
\mathcal{D}_kf(z)
&=\left\langle f,D_k(z,\cdot)\right\rangle
\nonumber\\
&=\left\langle\sum_{k'=-\infty}^{\infty}\sum_{\tau'\in I_{k'}}
\sum_{v'=1}^{N(k',\tau')}\mu(Q_{\tau'}^{k',v'})
\widetilde{D}_{k'}(\cdot,y_{\tau'}^{k',v'})
\mathcal{D}_{k'}f(y_{\tau'}^{k',v'}),D_k(z,\cdot)\right\rangle
\nonumber\\
&=\sum_{k'=-\infty}^{\infty}\sum_{\tau'\in I_{k'}}
\sum_{v'=1}^{N(k',\tau')}\mu(Q_{\tau'}^{k',v'})\left\langle
\widetilde{D}_{k'}(\cdot,y_{\tau'}^{k',v'})
,D_k(z,\cdot)\right\rangle\mathcal{D}_{k'}f(y_{\tau'}^{k',v'})
\nonumber\\
&=\sum_{k'=-\infty}^\infty\sum_{\tau'\in I_{k'}}
\sum_{v'=1}^{N(k',\tau')}
\mu(Q^{k',v'}_{\tau'})
D_k\widetilde{D}_{k'}(z,y^{k',v'}_{\tau'})
\mathcal{D}_{k'}f(y^{k',v'}_{\tau'}),
\end{align}
where the penultimate equality holds because the
series converge in $(\mathring{{\mathcal{G}}}_0^\varepsilon(\beta,\gamma))'$
and where, for any $k'\in\mathbb{Z}$,
$\tau'\in I_{k'}$, and $v'\in\{1,\ldots,N(k',\tau')\}$,
$y^{k',v'}_{\tau'}\in Q^{k',v'}_{\tau'}$ is the same point as in \eqref{3.29}.
From Lemma~\ref{DyadicSys}(iv),
it follows that, for any $x\in Q_{\tau}^{k,v}$
and $z\in B_\rho(z^{k,v}_\tau,C_\rho2^{-k})$,
$$
\rho(x,z)\lesssim 2^{-k} \leq 2^{-(k\land k')},
$$
which, in conjunction with Definition~\ref{D2.1}(iii),
Lemma~\ref{A3.2}(vi), and \eqref{doublingcond}, implies that
$$
\rho(x,y^{k',v'}_{\tau'})\lesssim\rho(x,z)
+\rho(z,y^{k',v'}_{\tau'})\lesssim 2^{-(k\land k')}+\rho(z,y^{k',v'}_{\tau'})
$$
and
\begin{align*}
&(V_\rho)_{2^{-(k\land k')}}(x)+V_\rho(x,y^{k',v'}_{\tau'})
\\
&\quad\sim(V_\rho)_{2^{-(k\land k')}+\rho(x,y^{k',v'}_{\tau'})}(x)
\lesssim(V_\rho)_{2^{-(k\land k')}+\rho(z,y^{k',v'}_{\tau'})}(z)
\\
&\quad
\sim(V_\rho)_{2^{-(k\land k')}}(z)+V_\rho(z,y^{k',v'}_{\tau'}).
\end{align*}
By this and Lemma~\ref{1122}, we find that,
for any $k\in\mathbb{Z}$,
$x\in Q_{\tau}^{k,v}$, $z\in B_\rho(z^{k,v}_\tau,C_\rho2^{-k})$,
and $\varepsilon'\in(n(\frac{1}{p}-1),\varepsilon)$,
\begin{align*}
\left|D_k\widetilde{D}_{k'}(z,y^{k',v'}_{\tau'})\right|
&\lesssim2^{-|k-k'|\varepsilon'}
\frac{1}{(V_\rho)_{2^{-(k\land k')}}(z)+V_\rho(z,y^{k',v'}_{\tau'})}
\left[\frac{2^{-(k\land k')}}
{2^{-(k\land k')}
+\rho(z,y^{k',v'}_{\tau'})}\right]^{\varepsilon'}
\\
&\lesssim2^{-|k-k'|\varepsilon'}
\frac{1}{(V_\rho)_{2^{-(k\land k')}}(x)+V_\rho(x,y^{k',v'}_{\tau'})}
\left[\frac{2^{-(k\land k')}}
{2^{-(k\land k')}
+\rho(x,y^{k',v'}_{\tau'})}\right]^{\varepsilon'},
\end{align*}
which, together with \eqref{1502}, further implies that
\begin{align}\label{911}
\left|\mathcal{D}_kf(z)\right|
&\lesssim\sum_{k'\in\mathbb{Z}}2^{-|k-k'|\varepsilon'}\sum_{\tau'\in I_{k'}}
\sum_{v'=1}^{N(k',\tau')}
\mu(Q^{k',v'}_{\tau'})
\frac{1}{(V_\rho)_{2^{-(k\land k')}}(x)+V_\rho(x,y^{k',v'}_{\tau'})}
\nonumber\\
&\quad\times\left[\frac{2^{-(k\land k')}}
{2^{-(k\land k')}
+\rho(x,y^{k',v'}_{\tau'})}\right]^{\varepsilon'}
\left|\mathcal{D}_{k'}f(y^{k',v'}_{\tau'})\right|.
\end{align}
Given that $y^{k',v'}_{\tau'}\in Q^{k',v'}_{\tau'}$
was arbitrary, we conclude that
\eqref{911} holds with $|\mathcal{D}_{k'}f(y^{k',v'}_{\tau'})|$ replaced by
$\inf_{u\in Q_{\tau'}^{k',v'}}|\mathcal{D}_{k'}f(u)|$.
From this and Lemma~\ref{L4.5},
we deduce that,
for any fixed $r\in(\frac{n}{n+\varepsilon'},p)$ and $k\in\mathbb{Z}$
and for any $z\in B_\rho(z^{k,v}_\tau,C_\rho2^{-k})$,
\begin{align*}
\left|\mathcal{D}_kf(z)\right|
&\lesssim\sum_{k'\in\mathbb{Z}}2^{-|k-k'|\varepsilon'}
2^{-n(k\land k'-k')(\frac{1}{r}-1)}\\
&\quad\times
\left[\mathcal{M}_\rho\left(\sum_{\tau'\in I_{k'}}
\sum_{v'=1}^{N(k',\tau')}
\inf_{u\in Q_{\tau'}^{k',v'}}\left|\mathcal{D}_{k'}f(u)\right|^r
\mathbf{1}_{Q^{k',v'}_{\tau}}\right)(x)\right]^\frac{1}{r},
\end{align*}
which, combined with \eqref{1501}, \eqref{un-3iu},
Lemma~\ref{DyadicSys}(ii), and H\"older's inequality,
further implies that, for $\mu$-almost every $x\in X$,
\begin{align*}
\left[S(f)(x)\right]^2
&\lesssim\sum_{k\in\mathbb{Z}}\sum_{\tau\in I_{k}}
\sum_{v=1}^{N(k,\tau)}
\Bigg\{\sum_{k'\in\mathbb{Z}}
2^{-n(k\land k'-k')(\frac{1}{r}-1)-|k-k'|\varepsilon'}
\nonumber\\
&\quad\times\left[\mathcal{M}_\rho
\left(\sum_{\tau'\in I_{k'}}\sum_{v'=1}^{N(k',\tau')}
\inf_{u\in Q_{\tau'}^{k',v'}}\left|\mathcal{D}_{k'}f(u)\right|^r
\mathbf{1}_{Q^{k',v'}_{\tau}}\right)(x)\right]^\frac{1}{r}
\mathbf{1}_{Q_\tau^{k,v}}(x)\Bigg\}^2
\nonumber\\
&\leq\sum_{k\in\mathbb{Z}}\sum_{\tau\in I_{k}}
\sum_{v=1}^{N(k,\tau)}
\left\{\sum_{k'\in\mathbb{Z}}
2^{2[-n(k\land k'-k')(\frac{1}{r}-1)-|k-k'|\varepsilon']}\right.
\nonumber\\
&\quad\left.\times\sum_{k'\in\mathbb{Z}}\left[\mathcal{M}_\rho
\left(\sum_{\tau'\in I_{k'}}
\sum_{v'=1}^{N(k',\tau')}
\inf_{u\in Q_{\tau'}^{k',v'}}\left|\mathcal{D}_{k'}f(u)\right|^r
\mathbf{1}_{Q^{k',v'}_{\tau}}\right)(x)\right]^\frac{2}{r}\right\}
\mathbf{1}_{Q_\tau^{k,v}}(x).
\end{align*}
This further implies that
\begin{align}\label{819Sf}
\left[S(f)(x)\right]^2
&=\sum_{k'\in\mathbb{Z}}\left\{\sum_{k\in\mathbb{Z},\,k\ge k'}
2^{2(k-k')[n(\frac{1}{r}-1)-\varepsilon']}
+\sum_{k\in\mathbb{Z},\,k<k'}
2^{2(k-k')\varepsilon'}\right\}
\nonumber\\
&\quad\times
\sum_{k'\in\mathbb{Z}}\left[\mathcal{M}_\rho\left(\sum_{\tau'\in I_{k'}}
\sum_{v'=1}^{N(k',\tau')}
\inf_{u\in Q_{\tau'}^{k',v'}}\left|\mathcal{D}_{k'}f(u)\right|^r
\mathbf{1}_{Q^{k',v'}_{\tau}}\right)(x)\right]^\frac{2}{r}
\nonumber\\
&\lesssim\sum_{k'\in\mathbb{Z}}\left[\mathcal{M}_\rho
\left(\sum_{\tau'\in I_{k'}}
\sum_{v'=1}^{N(k',\tau')}
\inf_{u\in Q_{\tau'}^{k',v'}}\left|\mathcal{D}_{k'}f(u)\right|^r
\mathbf{1}_{Q^{k',v'}_{\tau}}\right)(x)\right]^\frac{2}{r}.
\end{align}
By this, Theorem~\ref{independentofATI},
Lemma~\ref{Fefferman--Stein}
with $p$ and $u$ therein replaced,
respectively, by $\frac{p}{r}$ and $\frac{2}{r}$,
and Lemma~\ref{DyadicSys}(i), we conclude that
\begin{align*}
||f||_{H^p(X)}
&\sim\left\|\left[S(f)\right]^r\right\|^\frac{1}{r}_{L^\frac{p}{r}(X)}\\
&\lesssim\left\|\left\{\sum_{k'\in\mathbb{Z}}
\left[\mathcal{M}_\rho\left(\sum_{\tau'\in I_{k'}}
\sum_{v'=1}^{N(k',\tau')}
\inf_{u\in Q_{\tau'}^{k',v'}}\left|\mathcal{D}_{k'}f(u)\right|^r
\mathbf{1}_{Q^{k',v'}_{\tau}}
\right)\right]^\frac{2}{r}\right\}^{\frac{r}{2}}
\right\|_{L^\frac{p}{r}(X)}^{\frac{1}{r}}\\
&\lesssim\left\|\left\{\sum_{k'\in\mathbb{Z}}
\left[\sum_{\tau'\in I_{k'}}
\sum_{v'=1}^{N(k',\tau')}
\inf_{u\in Q_{\tau'}^{k',v'}}\left|\mathcal{D}_{k'}f(u)\right|^r
\mathbf{1}_{Q^{k',v'}_{\tau}}\right]^{\frac{2}{r}}\right\}^
{\frac{r}{2}}\right\|_{L^\frac{p}{r}(X)}^{\frac{1}{r}}\\
&=\left\|\left[\sum_{k'\in\mathbb{Z}}
\sum_{\tau'\in I_{k'}}\sum_{v'=1}^{N(k',\tau')}
\inf_{u\in Q_{\tau'}^{k',v'}}\left|\mathcal{D}_{k'}f(u)\right|^2
\mathbf{1}_{Q^{k',v'}_{\tau}}\right]
^{\frac{1}{2}}\right\|_{L^p(X)}\leq\left\|g(f)\right\|_{L^p(X)},
\end{align*}
which completes the proof of Theorem~\ref{gf}.	
\end{proof}

Next, we establish the Littlewood--Paley $g^*_\lambda$-function
characterization of $H^p(X)$.
To this end,
we need the following
Littlewood--Paley auxiliary function with apertures.

\begin{definition}\label{auxiliary function}
Let $(X,\mathbf{q},\mu)$ be
a quasi-ultrametric space of homogeneous type.
Let $\rho\in\mathbf{q}$ be such that all $\rho$-balls
are $\mu$-measurable.
Assume that
$0<\beta,\gamma<\varepsilon<\infty$
and $\{\mathcal{S}_{2^{-k}}\}_{k\in\mathbb{Z}}$
is a homogeneous approximation
of the identity of order $\varepsilon$
in Definition~\ref{DefHATI}.
Let $f\in(\mathring{{\mathcal{G}}}_0^\varepsilon(\beta,\gamma))'$
and $\mathcal{D}_k:=\mathcal{S}_{2^{-k}}-\mathcal{S}_{2^{-k+1}}$
for any $k\in\mathbb{Z}$.
The \emph{Littlewood--Paley auxiliary function $S^{(1)}_\theta(f)$} of $f$,
\index[words]{Littlewood--Paley auxiliary function}
with any given $\theta\in(0,\infty)$, is defined
by setting, for any $x\in X$,
\begin{align}
\label{S1thetaf}
S_\theta^{(1)}(f)(x)\index[symbols]{S@$S_\theta^{(1)}(f)$}
:=\left[\sum_{k=-\infty}^{\infty}
\int_{B_\rho(x,\theta2^{-k})}
\left|\mathcal{D}_kf(y)\right|^2
\frac{d\mu(y)}{(V_\rho)_{2^{-k}}(y)}\right]^\frac{1}{2}.
\end{align}
\end{definition}

\begin{remark}
If $S_\theta^{(1)}(f)$ is constructed using $\rho$ and
$\widetilde{S}_\theta^{(1)}(f)$ is constructed using
$\rho'\in\mathbf{q}$ (where all $\rho'$-balls are $\mu$-measurable),
then $S_\theta^{(1)}(f)\sim\widetilde{S}_\theta^{(1)}(f)$ pointwise in $X$.
Thus, without loss of generality, we may assume that $\rho\in\mathbf{q}$
in \eqref{S1thetaf} is a regularized quasi-ultrametric
in the remainder of this monograph.
Then, by this assumption
and an argument similar to that used in the proof of Lemma~\ref{1955},
we find that, for any $\theta\in(0,\infty)$
and $f\in(\mathring{{\mathcal{G}}}_0^\varepsilon(\beta,\gamma))'$,
$S^{(1)}_\theta(f)$ is $\mu$-measurable.
\end{remark}

For the Littlewood--Paley auxiliary function
in Definition~\ref{auxiliary function},
we have the following property.

\begin{lemma}\label{S(1)}
Let $(X,\mathbf{q},\mu)$ be
a quasi-ultrametric space of homogeneous type
with upper dimension $n\in(0,\infty)$.
Let $\theta\in[1,\infty)$, $p\in(0,2]$, and
$0<\beta,\gamma<\varepsilon\preceq\mathrm{ind\,}(X,\mathbf{q})$.
Then, for any
$f\in(\mathring{{\mathcal{G}}}_0^\varepsilon(\beta,\gamma))'$,
\begin{equation}
\label{S(1)eq}
\left\|S_\theta^{(1)}(f)\right\|_{L^p(X)}
\lesssim\theta^\frac{n}{p}\left\|S_1^{(1)}(f)\right\|_{L^p(X)}
\sim\theta^\frac{n}{p}\left\|S(f)\right\|_{L^p(X)},
\end{equation}
where the implicit positive constants are independent of
both $\theta$ and $f$.
\end{lemma}

\begin{proof}
Let $\rho\in\mathbf{q}$ be a regularized
quasi-ultrametric as in Definition~\ref{D2.3}.
Let $f\in(\mathring{{\mathcal{G}}}_0^\varepsilon(\beta,\gamma))'$.
We first prove the second inequality
in \eqref{S(1)eq}.
Observe that, by Definition~\ref{D2.1}(iii),
we find that,
for any $k\in\mathbb{Z}$, $x\in X$, and $y\in B_\rho(x,2^{-k})$,
$$
B_\rho\left(x,2^{-k}\right)
\subset B_\rho\left(y,C_\rho2^{-k}\right)
\subset B_\rho\left(x,C_\rho^22^{-k}\right),
$$
which, together with \eqref{doublingcond}, gives
$(V_\rho)_{2^{-k}}(x)\sim (V_\rho)_{2^{-k}}(y)$,
where the positive equivalence constants are
independent of $x$, $k$, and $y$.
From this and the definitions of $S^{(1)}_1(f)$
and $S(f)$, respectively, in \eqref{S1thetaf} and \eqref{Def-S-f},
we deduce that, for any $x\in X$,
\begin{align}\label{wwm-495}
S^{(1)}_1(f)(x)\sim S(f)(x),
\end{align}
from which the second inequality in \eqref{S(1)eq} follows.

Now, we prove the first inequality in \eqref{S(1)eq}.
If $p=2$, then, by Tonelli's theorem
and \eqref{upperdoub}, we conclude that
\begin{align*}
&\int_{X}
\left[S^{(1)}_\theta(f)(x)\right]^2\,d\mu(x)
\\
&\quad=\sum_{k=-\infty}^\infty\int_{X}
\int_{B_\rho(x,\theta2^{-k})}
\left|\mathcal{D}_kf(y)\right|^2
\frac{d\mu(y)}{(V_\rho)_{2^{-k}}(y)}\,d\mu(x)
\\
&\quad=\sum_{k=-\infty}^\infty
\int_{X}
\left|\mathcal{D}_kf(y)\right|^2
(V_\rho)_{\theta2^{-k}}(y)
\frac{d\mu(y)}{(V_\rho)_{2^{-k}}(y)}
\\
&\quad\lesssim\theta^n
\sum_{k=-\infty}^\infty
\int_{X}
\left|\mathcal{D}_kf(y)\right|^2
(V_\rho)_{2^{-k}}(y)
\frac{d\mu(y)}{(V_\rho)_{2^{-k}}(y)},
\end{align*}
which further implies that
\begin{align*}
&\int_{X}
\left[S^{(1)}_\theta(f)(x)\right]^2\,d\mu(x)
\\
&\quad=\theta^n\sum_{k=-\infty}^\infty\int_{X}
\left|\mathcal{D}_kf(y)\right|^2
\int_{X}
\mathbf{1}_{B_\rho(y,2^{-k})}(x)\,d\mu(x)
\frac{d\mu(y)}{(V_\rho)_{2^{-k}}(y)}\\
&\quad=\theta^n\sum_{k=-\infty}^\infty
\int_{X}\int_{B_\rho(x,2^{-k})}
\left|\mathcal{D}_kf(y)\right|^2
\frac{d\mu(y)}{(V_\rho)_{2^{-k}}(y)}\,d\mu(x)\\
&\quad=\theta^n\int_{X}
\left[S^{(1)}_1(f)(x)\right]^2\,d\mu(x).
\end{align*}
Hence,
\begin{align}\label{2256}
\left\|S_\theta^{(1)}(f)\right\|_{L^2(X)}\lesssim
\theta^\frac{n}{2}\left\|S_1^{(1)}(f)\right\|_{L^2(X)},
\end{align}
which implies that
the first inequality in \eqref{S(1)eq} holds when $p=2$.

To show that the first inequality in \eqref{S(1)eq} holds when $p\in(0,2)$, for
any nonnegative function $g\in L^1(X)$ and any $x\in X$, define
\begin{align}\label{1355}
\widetilde{\mathcal{M}_\rho}(g)(x)\index[symbols]{M@$\widetilde{\mathcal{M}_\rho}(g)$}
:=\sup_{k\in\mathbb{Z}}\sup_{\rho(x,y)<\theta2^{-k}}
\frac{1}{(V_\rho)_{2^{-k}}(y)}
\int_{B_\rho(y,2^{-k})}g(z)\,d\mu(z).
\end{align}
Note that, for any given $g\in L^1(X)$, $\widetilde{\mathcal{M}_\rho}(g)$ is
$\mu$-measurable. Indeed, this follows from an
argument similar to the one used in the proof
of the measurability of the Hardy--Littlewood maximal
function defined in \eqref{Mf}; see, for instance,
\cite[Theorem~3.7]{AlMi2015}.
Note that, from Definition~\ref{D2.1}(iii),
we infer that,
for any $k\in\mathbb{Z}$,
$y\in B_\rho(x,\theta2^{-k})$,
$z\in B_\rho(y,2^{-k})$, and $u\in B_\rho(x,C_\rho\theta2^{-k})$,
$$
\rho(x,z)\leq C_\rho
\max\left\{\rho(x,y),\,\rho(y,z)\right\}
\leq C_\rho\theta2^{-k}
$$
and
$$
\rho(y,u)\leq C_\rho
\max\left\{\rho(y,x),\,\rho(x,u)\right\}
\leq C_\rho^2
\theta2^{-k},
$$
which further implies that
$$
B_\rho(y,2^{-k})\subset B_\rho(x,C_\rho\theta2^{-k})
\subset B_\rho(y,C_\rho^2\theta2^{-k}).
$$
By this and \eqref{upperdoub}, we find that,
for any $k\in\mathbb{Z}$,
$x\in X$, and $y\in B_\rho(x,\theta2^{-k}),$
\begin{align*}
&\frac{1}{(V_\rho)_{2^{-k}}(y)}\int_{B_\rho(y,2^{-k})}g(z)\,d\mu(z)\\
&\quad\leq\frac{(V_\rho)_{C_\rho\theta2^{-k}}(x)}{(V_\rho)_{2^{-k}}(y)}
\frac{1}{(V_\rho)_{C_\rho\theta2^{-k}}(x)}
\int_{B_\rho(x,C_\rho\theta2^{-k})}g(z)\,d\mu(z)
\\
&\quad\leq\frac{(V_\rho)_{C_\rho^2\theta2^{-k}}(y)}{(V_\rho)_{2^{-k}}(y)}
\mathcal{M}_\rho(g)(x)
\lesssim\theta^n\mathcal{M}_\rho(g)(x),
\end{align*}
where $\mathcal{M}_\rho$ denotes the Hardy--Littlewood
maximal function defined in \eqref{Mf}.
This, combined with the boundedness of $\mathcal{M}_\rho$
from $L^1(X)$ to $L^{1,\infty}(X)$
(see Lemma~\ref{M} above),
further implies that, for any $r\in(0,\infty)$,
\begin{align}\label{muM}
&\mu\left(\left\{x\in X:\widetilde{\mathcal{M}_\rho}
(g)(x)>r\right\}\right)\nonumber\\
&\quad\leq\mu\left(\left\{x\in X:
(C_\rho\theta)^n\mathcal{M}_\rho(g)(x)>r\right\}\right)
\lesssim\frac{\theta^n}{r}||g||_{L^1(X)},
\end{align}
where the implicit positive constant is independent of $g$ and $\theta$.
For any $t\in(0,\infty)$,
consider the sets
$$
E_t:=\left\{x\in X:S_1^{(1)}(f)(x)>t\right\}
$$
and
\begin{align}
\label{1403}
\widetilde{E}_t
:=\left\{x\in X:\widetilde{\mathcal{M}_\rho}
(\mathbf{1}_{E_t})(x)>\frac{1}{2}\right\}.
\end{align}
Since the desired inequality is trivially true
if $\|S_1^{(1)}(f)\|_{L^p(X)}=\infty$, in what
follows, we may assume that $\|S_1^{(1)}(f)\|_{L^p(X)}<\infty$.
In particular, for any $t\in(0,\infty)$,
$\mu(E_t)<\infty$ and hence
$\mathbf{1}_{E_t}\in L^1(X)$. We claim that,
for any $t\in(0,\infty)$,
\begin{align}
\label{1403-x}
\int_{\widetilde{E}_t^{\complement}}
\left[S^{(1)}_\theta(f)(x)\right]^2\,d\mu(x)
\lesssim\theta^n\int_{E_t^{\complement}}
\left[S^{(1)}_1(f)(x)\right]^2\,d\mu(x).
\end{align}
Observe that \eqref{1403-x} trivially holds
if $\widetilde{E}_t^{\complement}=\varnothing$ and hence,
in what follows, we may assume that $\widetilde{E}_t^{\complement}\neq\varnothing$.
In this case, for any $y\in X$, let
$\widetilde{\rho}(y)
:=\inf_{x\in\widetilde{E}_t^{\complement}}\rho(x,y)$.
For any $k\in\mathbb{Z}$
and $y\in X$, if there exists $y_0\in X$ such that
$y_0\in\widetilde{E}_t^{\complement}\cap B_\rho(y,\theta2^{-k})$,
then $\widetilde{\rho}(y)<\theta2^{-k}$. In this case,
from the definitions of $\widetilde{E}_t$ in \eqref{1403} and
$\widetilde{\mathcal{M}_\rho}$ in \eqref{1355}
and from $\rho(y,y_0)<\theta2^{-k}$,
we deduce that, for any $t\in(0,\infty)$
and $k\in\mathbb{Z}$,
\begin{align*}
\mu\left(E_t\cap B_\rho(y,2^{-k})\right)
&=\int_{B_\rho(y,2^{-k})}
\mathbf{1}_{E_t}(z)\,d\mu(z)
\\
&=(V_\rho)_{2^{-k}}(y)
\frac{1}{(V_\rho)_{2^{-k}}(y)}
\int_{B_\rho(y,2^{-k})}
\mathbf{1}_{E_t}(z)\,d\mu(z)
\\
&\leq (V_\rho)_{2^{-k}}(y)
\widetilde{\mathcal{M}_\rho}(\mathbf{1}_{E_t})(y_0)
\leq\frac{1}{2}(V_\rho)_{2^{-k}}(y),
\end{align*}
which further implies that
\begin{align}\label{2208}
\mu\left(E^{\complement}_t\cap B_\rho(y,2^{-k})\right)
\ge\frac{1}{2}(V_\rho)_{2^{-k}}(y).
\end{align}
By this, Tonelli's theorem,
and \eqref{upperdoub}, we further conclude that,
for any $t\in(0,\infty)$,
\begin{align*}
&\int_{\widetilde{E}^{\complement}_t}
\left[S^{(1)}_\theta(f)(x)\right]^2\,d\mu(x)\\
&\quad=\int_{\widetilde{E}^{\complement}_t}
\sum_{k=-\infty}^\infty
\int_{B_\rho(x,\theta2^{-k})}
\left|\mathcal{D}_kf(y)\right|^2
\frac{d\mu(y)}{(V_\rho)_{2^{-k}}(y)}\,d\mu(x)
\\
&\quad=\sum_{k=-\infty}^\infty
\int_{\widetilde{\rho}(y)<\theta2^{-k}}
\left|\mathcal{D}_kf(y)\right|^2
\mu\left(\widetilde{E}^{\complement}_t
\cap B_\rho(y,\theta2^{-k})\right)
\frac{d\mu(y)}{(V_\rho)_{2^{-k}}(y)}
\\
&\quad\leq\sum_{k=-\infty}^\infty
\int_{\widetilde{\rho}(y)<\theta2^{-k}}
\left|\mathcal{D}_kf(y)\right|^2
(V_\rho)_{\theta2^{-k}}(y)
\frac{d\mu(y)}{(V_\rho)_{2^{-k}}(y)}
\\
&\quad\lesssim\theta^n
\sum_{k=-\infty}^\infty
\int_{\widetilde{\rho}(y)<\theta2^{-k}}
\left|\mathcal{D}_kf(y)\right|^2
(V_\rho)_{2^{-k}}(y)
\frac{d\mu(y)}{(V_\rho)_{2^{-k}}(y)}.
\end{align*}
This further implies that
\begin{align*}
&\int_{\widetilde{E}^{\complement}_t}
\left[S^{(1)}_\theta(f)(x)\right]^2\,d\mu(x)\\
&\quad\leq\theta^n\sum_{k=-\infty}^\infty
\int_{X}\left|\mathcal{D}_kf(y)\right|^2
\mu\left(E^{\complement}_t\cap B_\rho(y,2^{-k})\right)
\frac{d\mu(y)}{(V_\rho)_{2^{-k}}(y)}\\
&\quad=\theta^n\sum_{k=-\infty}^\infty\int_{X}
\left|\mathcal{D}_kf(y)\right|^2
\int_{B_\rho(y,2^{-k})}
\mathbf{1}_{E^{\complement}_t}(x)\,d\mu(x)
\frac{d\mu(y)}{(V_\rho)_{2^{-k}}(y)}\\
&\quad=\theta^n\sum_{k=-\infty}^\infty
\int_{E^{\complement}_t}\int_{B_\rho(x,2^{-k})}
\left|\mathcal{D}_kf(y)\right|^2
\frac{d\mu(y)}{(V_\rho)_{2^{-k}}(y)}\,d\mu(x)\\
&\quad=\theta^n\int_{E^{\complement}_t}
\left[S^{(1)}_1(f)(x)\right]^2\,d\mu(x),
\end{align*}
which completes the proof of the above claim \eqref{1403-x}.
From \eqref{muM} [used here with $g=\mathbf{1}_{E_t}\in L^1(X)$],
it follows that, for any $t\in(0,\infty)$,
$$
\mu\Big(\widetilde{E}_t\Big)=
\mu\left(\left\{x\in X:\widetilde{\mathcal{M}_\rho}
\mathbf{1}_{E_t}(x)>\frac{1}{2}\right\}\right)
\lesssim\theta^n\left\|\mathbf{1}_{E_t}\right\|_{L^1(X)}
=\theta^n\mu(E_t),
$$
which implies that
\begin{align*}
&\mu\left(\left\{x\in X:
S^{(1)}_\theta(f)(x)>t\right\}\right)
\nonumber\\
&\quad\leq\mu\Big(\widetilde{E}_t\Big)
+\mu\left(\left\{x\in\widetilde{E}_t^{\complement}:
S^{(1)}_\theta(f)(x)>t\right\}\right)
\nonumber\\
&\quad\lesssim\theta^n\mu(E_t)
+t^{-2}\int_{\widetilde{E}_t^{\complement}}
\left[S^{(1)}_\theta(f)(x)\right]^2\,d\mu(x)
\quad\text{by Chebyshev's inequality}
\nonumber\\
&\quad\lesssim\theta^n\left[\mu(E_t)
+t^{-2}\int_{E_t^{\complement}}
\left[S^{(1)}_1(f)(x)\right]^2\,d\mu(x)\right]
\quad\text{by \eqref{1403-x}}.
\end{align*}
This further implies that
\begin{align}\label{8193}
&\mu\left(\left\{x\in X:
S^{(1)}_\theta(f)(x)>t\right\}\right)
\nonumber\\
&\quad=\theta^n\left[\mu(E_t)
+t^{-2}\left\|S^{(1)}_1(f)(x)
\mathbf{1}_{E_t^{\complement}}(x)
\right\|^2_{L^2(X)}\right]
\nonumber\\
&\quad\sim\theta^n\left[\mu(E_t)
+t^{-2}
\int_0^\infty s\mu\left(\left\{x\in X:
S^{(1)}_1(f)(x)
\mathbf{1}_{E_t^{\complement}}
(x)>s\right\}\right)\,ds\right]
\nonumber\\
&\quad=\theta^n\left[\mu(E_t)
+t^{-2}
\int_0^\infty s
\mu\left(E_s\cap\left\{x\in X:S^{(1)}_1(f)(x)
\leq t\right\}\right)\,ds\right]
\nonumber\\	
&\quad\leq\theta^n\left[\mu(E_t)
+t^{-2}\int_0^ts\mu(E_s)\,ds\right].
\end{align}
By this, Tonelli's theorem,
and $p\in(0,2)$,
we obtain
\begin{align*}
\left\|S^{(1)}_\theta(f)\right\|^p_{L^p(X)}
&=p\int_0^\infty t^{p-1}\mu\left(\left\{x\in X:
S^{(1)}_\theta(f)(x)>t\right\}\right)dt
\\
&\lesssim\theta^n
\left[\int_0^\infty t^{p-1}\mu(E_t)dt
+\int_0^\infty t^{p-3}
\int_0^ts\mu(E_s)dsdt\right]
\\
&\sim\theta^n\left\|S^{(1)}_1(f)\right\|^p_{L^p(X)}
+\theta^n\int_0^\infty s\mu(E_s)
\int_s^\infty t^{p-3}dtds
\\
&\sim\theta^n\left\|S^{(1)}_1(f)\right\|^p_{L^p(X)}
+\theta^n\int_0^\infty s^{p-1}\mu(E_s)ds
\sim\theta^n\left\|S^{(1)}_1(f)\right\|^p_{L^p(X)},
\end{align*}
which, together with \eqref{2256}, completes the proof
of the first inequality in \eqref{S(1)eq}
and hence Lemma~\ref{S(1)}.
\end{proof}

Next, we establish the Littlewood--Paley $g_\lambda^*$-function
characterization of $H^p(X)$ in the following theorem.

\begin{theorem}\label{Glamda*}
Let $(X,\mathbf{q},\mu)$ be an ultra-RD-space
with upper dimension $n\in(0,\infty)$.
Let $p\in(0,2]$,
$\lambda\in(\frac{2n}{p},\infty)$,
and $\beta,\gamma,\varepsilon\in(0,\infty)$ be such that
$0<\beta,\,\gamma<\varepsilon\preceq\mathrm{ind\,}(X,\mathbf{q})$.
Assume that $\{\mathcal{S}_{2^{-k}}\}_{k\in\mathbb{Z}}$
is a homogeneous approximation
of the identity of order $\varepsilon$
in Definition~\ref{DefHATI} and, for any $k\in\mathbb{Z}$, let
$\mathcal{D}_k:=\mathcal{S}_{2^{-k}}-\mathcal{S}_{2^{-k+1}}$.
Let $S$ and $g^*_\lambda$ be
the same operators as, respectively, in \eqref{Def-S-f} and \eqref{Def-glamda*}.
Then there exists a constant $C\in[1,\infty)$
such that, for any
$f\in(\mathring{{\mathcal{G}}}_0^\varepsilon(\beta,\gamma))'$,
\begin{equation}\label{A4.22}
C^{-1}\left\|g^*_\lambda(f)\right\|_{L^p(X)}
\leq\left\|S(f)\right\|_{L^{p}(X)}
\leq C\left\|g^*_\lambda(f)\right\|_{L^p(X)}.
\end{equation}
Consequently, if $p\in(\frac{n}{n+\mathrm{ind\,}(X,\mathbf{q})},2]$ and
$\beta,\gamma,\varepsilon\in(0,\infty)$
satisfy \eqref{1945},
then $f\in H^p(X)$ if and only if
$f\in(\mathring{{\mathcal{G}}}_0^\varepsilon(\beta,\gamma))'$
and $g^*_\lambda(f)\in L^p(X)$.
Moreover, for any $f\in H^p(X)$,
$$
\left\|f\right\|_{H^p(X)}\sim\left\|g^*_\lambda(f)\right\|_{L^p(X)},
$$
where the positive equivalence constants are independent of $f$.
\end{theorem}

\begin{proof}
Let $\rho\in\mathbf{q}$ be a regularized quasi-ultrametric
as in Definition~\ref{D2.3}.
Let $f\in(\mathring{{\mathcal{G}}}_0^\varepsilon(\beta,\gamma))'$.
From Definition~\ref{D2.1}(iii),
we infer that,
for any given $k\in\mathbb{Z}$, $x\in X$, and $y\in B_\rho(x,2^{-k})$
and for any $z\in B_\rho(y,2^{-k})$,
$$
\rho(x,z)\leq C_\rho\max\left\{\rho(x,y),\,\rho(y,z)\right\}
<C_\rho2^{-k}
$$
and hence $B_\rho(y,2^{-k})\subset B_\rho(x,C_\rho 2^{-k})$,
which, combined with \eqref{upperdoub}, gives
$$
(V_\rho)_{2^{-k}}(y)\lesssim (V_\rho)_{2^{-k}}(x).
$$
Consequently, on the one hand,
by the definitions of both
$S(f)$ and  $g^*_\lambda(f)$,
we find that,
for any $x\in X$,
\begin{align}
\label{819Sfg*}
\left[S(f)(x)\right]^2
&=2^\lambda\sum_{k=-\infty}^{\infty}
\int_{B_\rho(x,2^{-k})}
\left|\mathcal{D}_kf(y)\right|^2
\left[\frac{2^{-k}}{2^{-k}+2^{-k}}\right]^\lambda
\frac{d\mu(y)}{(V_\rho)_{2^{-k}}(x)}
\nonumber\\
&\leq 2^\lambda\sum_{k=-\infty}^{\infty}
\int_{B_\rho(x,2^{-k})}
\left|\mathcal{D}_kf(y)\right|^2
\left[\frac{2^{-k}}{2^{-k}+\rho(x,y)}\right]^\lambda
\frac{d\mu(y)}{(V_\rho)_{2^{-k}}(x)}
\nonumber\\
&=2^{\lambda+1}\sum_{k=-\infty}^{\infty}
\int_{B_\rho(x,2^{-k})}
\left|\mathcal{D}_kf(y)\right|^2
\left[\frac{2^{-k}}{2^{-k}+\rho(x,y)}\right]^\lambda
\frac{d\mu(y)}{(V_\rho)_{2^{-k}}(x)+(V_\rho)_{2^{-k}}(x)}
\nonumber\\
&\lesssim 2^{\lambda}\sum_{k=-\infty}^{\infty}
\int_X\left|\mathcal{D}_kf(y)\right|^2
\left[\frac{2^{-k}}{2^{-k}+\rho(x,y)}\right]^\lambda
\frac{d\mu(y)}{(V_\rho)_{2^{-k}}(x)+(V_\rho)_{2^{-k}}(y)}
\nonumber\\
&= 2^{\lambda}\left[g^*_\lambda(f)(x)\right]^2,
\end{align}
which, together with Theorem~\ref{independentofATI},
further implies that
\begin{equation}
\label{A4.23}
\left\|f\right\|_{H^p(X)}
\sim\left\|S(f)\right\|_{L^p(X)}
\lesssim\left\|g^*_\lambda(f)\right\|_{L^p(X)}.
\end{equation}
On the other hand,
from the definition of $g^*_\lambda(f)$, we deduce that,
for any $x\in X$,
\begin{align*}
\left[g^*_\lambda(f)(x)\right]^2
&=\sum_{k=-\infty}^{\infty}\sum_{j=1}^{\infty}
\int_{B_\rho(x,2^j2^{-k})\setminus B_\rho(x,2^{j-1}2^{-k})}
\left|\mathcal{D}_kf(y)\right|^2
\nonumber\\
&\quad\times
\left[\frac{2^{-k}}{2^{-k}+\rho(x,y)}\right]^\lambda
\frac{d\mu(y)}{(V_\rho)_{2^{-k}}(x)+(V_\rho)_{2^{-k}}(y)}
\nonumber\\
&\quad+\sum_{k=-\infty}^{\infty}
\int_{B_\rho(x,2^{-k})}\left|\mathcal{D}_kf(y)\right|^2
\left[\frac{2^{-k}}{2^{-k}+\rho(x,y)}\right]^\lambda
\frac{d\mu(y)}{(V_\rho)_{2^{-k}}(x)+(V_\rho)_{2^{-k}}(y)},
\end{align*}
which, combined with the definition of $S_\theta^{(1)}(f)$,
further implies that
\begin{align}
\label{819g*sf}
\left[g^*_\lambda(f)(x)\right]^2
&\leq\sum_{k=-\infty}^{\infty}\sum_{j=1}^{\infty}
\int_{B_\rho(x,2^j2^{-k})}\left|\mathcal{D}_kf(y)\right|^2
\left[\frac{1}{1+2^{j-1}}\right]^\lambda
\frac{d\mu(y)}{(V_\rho)_{2^{-k}}(x)+(V_\rho)_{2^{-k}}(y)}
\nonumber\\
&\quad+\sum_{k=-\infty}^{\infty}
\int_{B_\rho(x,2^{-k})}\left|\mathcal{D}_kf(y)\right|^2
\frac{d\mu(y)}{(V_\rho)_{2^{-k}}(x)+(V_\rho)_{2^{-k}}(y)}
\nonumber\\
&\lesssim \sum_{j=1}^{\infty}2^{-j\lambda}
\sum_{k=-\infty}^{\infty}
\int_{B_\rho(x,2^j2^{-k})}
\left|\mathcal{D}_kf(y)\right|^2
\frac{d\mu(y)}{(V_\rho)_{2^{-k}}(y)}
+\left[S_1(f)(x)\right]^2
\nonumber\\
&=\sum_{j=0}^{\infty}2^{-j\lambda}
\left[S_{2^j}^{(1)}(f)(x)\right]^2.
\end{align}
By this and the triangle inequality for
$||\cdot||_{L^\frac{p}{2}(X)}^\frac{p}{2}$,
we conclude that
\begin{align*}
\left\|g^*_\lambda(f)\right\|_{L^p(X)}^p
&=\left\|\left[g^*_\lambda(f)\right]^2
\right\|_{L^\frac{p}{2}(X)}^\frac{p}{2}\\
&\lesssim\left\|\sum_{j=0}^{\infty}2^{-j\lambda}
\left[S_{2^j}^{(1)}(f)\right]^2\right\|_{L^\frac{p}{2}(X)}^\frac{p}{2}
\quad\text{by \eqref{819g*sf}}\\
&\leq\sum_{j=0}^{\infty}2^{-\frac{j\lambda p}{2}}
\left\|\left[S_{2^j}^{(1)}(f)
\right]^2\right\|_{L^\frac{p}{2}(X)}^\frac{p}{2}
=\sum_{j=0}^{\infty}2^{-\frac{j\lambda p}{2}}
\left\|S_{2^j}^{(1)}(f)\right\|_{L^{p}(X)}^{p}
\\
&\lesssim\sum_{j=0}^{\infty}2^{-j(\frac{\lambda p}{2}-n)}
\left\|S(f)\right\|_{L^{p}(X)}^{p}
\quad\text{by Lemma~\ref{S(1)}}\\
&\sim\left\|S(f)\right\|_{L^{p}(X)}^{p}
\quad\text{since $\lambda\in(\frac{2n}{p},\infty)$},
\end{align*}
which, together with Theorem~\ref{independentofATI},
further implies that
\begin{equation}\label{A4.24}
\left\|g^*_\lambda(f)\right\|_{L^p(X)}
\lesssim\left\|S(f)\right\|_{L^p(X)}
\sim\left\|f\right\|_{H^p(X)}.
\end{equation}
Combining \eqref{A4.23} with \eqref{A4.24},
we obtain \eqref{A4.22},
which completes the proof of Theorem~\ref{Glamda*}.
\end{proof}

\begin{remark}
\label{1629}
\begin{enumerate}
\item Note that we used the same
homogeneous approximation of the identity to define the
Lusin area function and the Littlewood--Paley $g_\lambda^*$ function;
see \eqref{Def-S-f} and \eqref{Def-glamda*}.
From this and the proof of Theorem~\ref{Glamda*},
we can clearly see that,
if $(X,\mathbf{q},\mu)$ is only assumed to be
a quasi-ultrametric space of homogeneous type,
then \eqref{A4.22} also
holds for any given $p\in(0,2]$ and for any
$f\in(\mathring{\mathcal{G}}^\varepsilon_0(\beta,\gamma))'$
with $\varepsilon$, $\beta$, and $\gamma$ as in Theorem~\ref{Glamda*}.
\item   In the proof of Theorem~\ref{Glamda*},
the reverse doubling condition \eqref{RDcondition}
is needed since we use Theorem~\ref{independentofATI}
to guarantee that $H^p(X)$ is independent of the choices
of the homogeneous approximation of the identity
and hence $H^p(X)$ defined by means of different Lusin area functions
coincide with each other with equivalent quasi-norms;
see \eqref{A4.23} and \eqref{A4.24}.
In this case, we also need
$p>\frac{n}{n+\mathrm{ind\,}(X,\mathbf{q})}$
due to Theorem~\ref{independentofATI}.
\item 	Let $X$ be a homogeneous group,
which is also an ultra-RD-space,
with upper dimension $n\in(0,\infty)$.
Note that $\lambda$ in \eqref{Def-glamda*} equals
$2\lambda$ with $\lambda$ as in the Littlewood--Paley
$g_\lambda^*$-function in \cite[p.\,218]{fs1982}.
In this case, Theorem~\ref{Glamda*} coincides with \cite[Corollary 7.4]{fs1982}
which proves that, for any given $p\in(0,2]$ and for any $f$
belonging to the space of distributions on $X$,
\begin{align*}
\left\|g^*_\lambda(f)\right\|_{L^p(X)}\lesssim\left\|S(f)\right\|_{L^p(X)},
\end{align*}
whenever $\lambda\in(\frac{2n}{p},\infty)$,
whose range of $\lambda$ is the known best possible.
Thus, the range of $\lambda$
in Theorem~\ref{Glamda*} is the known best possible.
\item
Let $(X,\mathbf{q},\mu)$ be an Ahlfors-regular
quasi-ultrametric space as in \cite{AlMi2015}.
In this case, by Remark~\ref{2548}, we find that
Theorem~\ref{Glamda*} gives the
Littlewood--Paley $g_\lambda^*$-function characterization
of the Hardy space in \cite{AlMi2015}.
\end{enumerate}
\end{remark}


\chapter[Littlewood--Paley Characterizations of Triebel--Lizorkin Spaces
on\\ Ultra-RD-Spaces]{Littlewood--Paley Characterizations of Triebel--Lizorkin Spaces
on\\ Ultra-RD-Spaces}
\label{Chapter4}
\markboth{\scriptsize\rm\sc Littlewood--Paley
Characterizations of Triebel--Lizorkin Spaces}
{\scriptsize\rm\sc Littlewood--Paley
Characterizations of Triebel--Lizorkin Spaces}

In this chapter, we establish various Littlewood--Paley function
characterizations of the homogeneous (as well as the inhomogeneous)
Triebel--Lizorkin spaces, for the optimal ranges
$s\in(-{\mathrm{ind\,}}(X,\mathbf{q}),{\mathrm{ind\,}}(X,\mathbf{q}))$ and
\begin{align*}
\max\left\{\frac{n}{n+{\mathrm{ind\,}}(X,\mathbf{q})},\,
\frac{n}{n+{\mathrm{ind\,}}(X,\mathbf{q})+s}\right\}<p,q\leq\infty,
\end{align*}
on  standard ultra-RD-spaces
$(X,\mathbf{q},\mu)$ with upper dimension $n\in(0,\infty)$,
where $\mu$ is assumed to be a Borel-semiregular measure on $X$.
Recall that a quasi-ultrametric space $(X,\mathbf{q})$
is said to be \emph{standard}\index[words]{standard}
if, for any $x\in X$ and for some (and hence for all)
$\rho\in\mathbf{q}$, $r_\rho(x)=0$, where
$r_\rho(x)$ is the same as in \eqref{rrho}.
Throughout Sections~\ref{Section4.1}--\ref{1700},
we always assume that $\mathrm{diam\,}(X)=\infty$
and hence $\mu(X)=\infty$,
which is a quite natural assumption because of
the homogeneity of the homogeneous Triebel--Lizorkin spaces.
In Section~\ref{Section4.1}, we recall the definition of
homogeneous Triebel--Lizorkin spaces $\dot{F}^{p,q}_s(X)$.
Then we establish the Lusin area function
and the Littlewood--Paley $g_\lambda^*$-function characterizations of
$\dot{F}^{p,q}_s(X)$ with $p<\infty$
in Sections~\ref{Section4.2} and~\ref{Section4.3}, respectively.
The analogous characterizations of $\dot{F}^{\infty,q}_s(X)$
are recorded in Section~\ref{1700}. We extend all of the
corresponding results to the inhomogeneous case
in Sections~\ref{subsection5.2}
and \ref{Section4.5}, where we have no restriction on
$\mathrm{diam\,}(X)$.

\section{Homogeneous Triebel--Lizorkin Spaces
$\dot{F}^{p,q}_s(X)$}
\label{Section4.1}

In this section, we recall the definitions
of the homogeneous Triebel--Lizorkin spaces,
the homogeneous Lusin area functions,
and the homogeneous Littlewood--Paley functions.
We always assume that, for any $\rho\in\mathbf{q}$,
$\mathrm{diam}_\rho(X)=\infty$
and hence $\mu(X)=\infty$.
We now recall
the definition of the homogeneous Triebel--Lizorkin spaces
on quasi-ultrametric spaces of
homogeneous type;
see \cite[Definitions 9.2 and 9.6]{AlMi2015}
for standard Ahlfors-regular quasi-ultrametric spaces with
a Borel-semiregular measure,
\cite[Definitions~5.8(ii), 5.29(ii), 6.7, and 6.14]{hmy08}
for standard doubling metric spaces with
a Borel regular measure,
and \cite[Definitions~3.1(ii), 4.1(ii), 5.1, and 5.9]{whyy21}
for standard quasi-metric spaces
of homogeneous type with a Borel regular measure.	
Recall that,
for any $r\in\mathbb{R}$, $r_+:=\max\{0,\,r\}$.

\begin{definition}\label{DefT--L}
Let $(X,\mathbf{q},\mu)$
be a quasi-ultrametric space of
homogeneous type with upper dimension $n\in(0,\infty)$
and $\rho\in\mathbf{q}$ be such that all $\rho$-balls are $\mu$-measurable.
Let $0<\varepsilon\preceq\mathrm{ind\,}(X,\mathbf{q})$
with $\mathrm{ind\,}(X,\mathbf{q})$ as in \eqref{index}
and assume that $\{\mathcal{S}_{2^{-k}}\}_{k\in\mathbb{Z}}$ is
a homogeneous approximation of the identity
of order $\varepsilon$
in Definition~\ref{DefHATI}
(whose existence is ensured
by Theorem~\ref{ThmATI}).
For any $k\in\mathbb{Z}$, define
$\mathcal{D}_k:=\mathcal{S}_{2^{-k}}-\mathcal{S}_{2^{-k+1}}$
and let $\{Q^l_\alpha\}_{\alpha\in I_l}$ be the set
of all dyadic cubes (at level $l\in\mathbb{Z}$) given
by Lemma~\ref{DyadicSys}.
Assume that $s\in(-\varepsilon,\varepsilon)$,
$$
p,q\in\left(\max\left\{\frac{n}{n+\varepsilon},\,
\frac{n}{n+\varepsilon+s}\right\},\infty\right],
$$
$$
\gamma\in\left(\max\left\{s-\frac{n}{p},\,
n\left(\frac{1}{p}-1\right)_+,-s+n\left(\frac{1}{p}-1\right)_+
-n\left(1-\frac{1}{p}\right)_+\right\},\varepsilon\right),
$$
and
$$
\beta\in\left(\max\left\{s_+,\,
-s+n\left(\frac{1}{p}-1\right)_+\right\},\varepsilon\right).
$$
The \emph{homogeneous Triebel--Lizorkin
space}\index[words]{homogeneous Triebel--Lizorkin space}
$\dot{F}^{p,q}_s(X)$\index[symbols]{F@$\dot{F}^{p,q}_s(X)$}
is defined to be the set of all distributions
$f\in(\mathring{{\mathcal{G}}}_0^\varepsilon(\beta,\gamma))'$
such that, when $p<\infty$,
\begin{align*}
\|f\|_{\dot{F}^{p,q}_s(X)}
:=\left\|\left(\sum_{k\in\mathbb{Z}}
2^{ksq}\left|\mathcal{D}_kf\right|^q\right)^\frac{1}{q}\right\|_{L^p(X)}<\infty
\end{align*}
and, when $p=\infty$,
\begin{align*}
\|f\|_{\dot{F}^{\infty,q}_s(X)}
:=\sup_{l\in\mathbb{Z}}\sup_{\alpha\in I_l}
\left[\frac{1}{\mu(Q_\alpha^l)}\int_{Q^l_\alpha}\sum_{k=l}^{\infty}
2^{ksq}\left|\mathcal{D}_kf(x)\right|^q\,d\mu(x)\right]^\frac{1}{q}<\infty
\end{align*}
with the usual modification when $q=\infty$.
\end{definition}

\begin{remark}
By the proof of Lemma~\ref{1955}, we conclude that,
for any $k\in\mathbb{Z}$ and
$f\in(\mathring{{\mathcal{G}}}_0^\varepsilon(\beta,\gamma))'$,
$\mathcal{D}_kf$ is actually continuous [when constructed using a regularized
quasi-ultrametric as in Definition~\ref{D2.3}] and hence $\mu$-measurable.
Thus, the integration in Definition~\ref{DefT--L} is meaningful.
\end{remark}

\begin{remark}\label{2101}
The following statements regarding Definition~\ref{DefT--L} are true.
\begin{enumerate}			
\item[\rm(i)]
If $\rho\in\mathbf{q}$ is an ultrametric, then $\mathrm{ind\,}(X,\mathbf{q})=\infty$
and the smoothness parameter $s$  can be chosen arbitrarily from $(-\infty,\infty)$.

\item[\rm(ii)]
Let $(X,\mathbf{q},\mu)$ be an ultra-RD-space
and $\mu$ be a Borel-semiregular measure on $X$.
Then the definition of $\dot{F}^{p,q}_s(X)$ is independent of
the choice of the approximation of the identity
(see Remarks~\ref{zj2037} and~\ref{zj2038} below)
as well as the choice of the indices $\beta,\gamma$; see also
\cite[Propositions~5.6 and~6.5]{hmy08}
and \cite[Propositions~5.7 and~6.6]{hmy08},
respectively.
Moreover, by Remark~\ref{Rm}(iii),
we find that $\dot{F}^{p,q}_s(X)$ is independent of
the choice of the quasi-ultrametric $\rho$
if $\rho\in\mathbf{q}$ is such that
$\mathcal{D}_kf$ for any $k\in\mathbb{Z}$
and $f\in(\mathring{{\mathcal{G}}}_0^\varepsilon(\beta,\gamma))'$
is $\mu$-measurable.

\item[\rm(iii)]
By considering $\dot{F}^{p,q}_s(X)$ with an equivalence relation,
$\sim$, defined by $f\sim g$ if and only if
$f-g$ is a constant on $X$, the space $\dot{F}^{p,q}_s(X)/\sim$ is quasi-Banach
for all $s$, $p$, $q$, $\beta$, and $\gamma$ as in Definition~\ref{DefT--L}, and
genuinely Banach whenever $p,q\geq1$; see also
\cite[Propositions~5.10(vi) and~6.9(v)]{hmy08}.
\end{enumerate}
\end{remark}

Next, we recall the definitions of the homogeneous
Lusin area function
and Littlewood--Paley functions.

\begin{definition}\label{homoareafuncts}
Let $(X,\mathbf{q},\mu)$ be a quasi-ultrametric space of homogeneous type
and $\rho\in\mathbf{q}$ be such that all $\rho$-balls are $\mu$-measurable.
Suppose that
$q\in(0,\infty]$ and
$\varepsilon,\beta,\gamma\in(0,\infty)$
satisfy $s\in(-\varepsilon,\varepsilon)$ and
$0<\beta,\gamma<\varepsilon\preceq\mathrm{ind\,}(X,\mathbf{q})$
with $\mathrm{ind\,}(X,\mathbf{q})$ as in \eqref{index}.
Let $f\in(\mathring{{\mathcal{G}}}_0^\varepsilon(\beta,\gamma))'$.
Assume that $\{\mathcal{S}_{2^{-k}}\}_{k\in\mathbb{Z}}$ is
a homogeneous approximation of the identity
of order $\varepsilon$
in Definition~\ref{DefHATI} and
define $\mathcal{D}_k:=\mathcal{S}_{2^{-k}}-\mathcal{S}_{2^{-k+1}}$
for any given $k\in\mathbb{Z}$.
\begin{enumerate}
\item
The \emph{homogeneous Littlewood--Paley $g$-function
$\dot{g}^s_q(f)$} of $f$\index[symbols]{G@$\dot{g}^s_q(f)$}
is defined
by setting, for any $x\in X$,
$$
\dot{g}^s_q(f)(x)
:=\left[\sum_{k\in\mathbb{Z}}
2^{ksq}|\mathcal{D}_kf(x)|^q\right]^\frac{1}{q}.
$$
\item
The \emph{homogeneous Lusin area function
$\dot{\mathcal{S}}^s_q(f)$}\index[symbols]{S@$\dot{\mathcal{S}}^s_q(f)$} of $f$
is defined
by setting, for any $x\in X$,
$$
\dot{\mathcal{S}}^s_q(f)(x)
:=\left[\sum_{k\in\mathbb{Z}}
2^{ksq}\int_{B_\rho(x,2^{-k})}
|\mathcal{D}_kf(y)|^q\frac{d\mu(y)}{(V_\rho)_{2^{-k}}(x)}\right]^\frac{1}{q}.
$$
\item
The \emph{homogeneous Littlewood--Paley
$g^*_\lambda$-function
$(\dot{g}^*_\lambda)^s_q(f)$} of $f$,\index[symbols]{G@$(\dot{g}^*_\lambda)^s_q(f)$}
with any given $\lambda\in(0,\infty)$,
is defined
by setting, for any $x\in X$,
\begin{align}\label{g*}
(\dot{g}^*_\lambda)^s_q(f)(x)
:&=\left\{\sum_{k\in\mathbb{Z}}2^{ksq}
\int_X|\mathcal{D}_kf(y)|^q
\left[\frac{2^{-k}}{2^{-k}+\rho(x,y)}\right]^\lambda\right.\nonumber\\
&\quad\left.\times\frac{d\mu(y)}{(V_\rho)_{2^{-k}}(x)+
(V_\rho)_{2^{-k}}(y)}\right\}^\frac{1}{q}.
\end{align}
\end{enumerate}
\end{definition}

\begin{remark}
Let the notation be as in Definition~\ref{homoareafuncts}.
\begin{enumerate}
\item
If $s:=0$ and $q:=2$, then
the Lusin area function and the Littlewood--Paley
functions in Definition~\ref{homoareafuncts}
coincide, respectively, with those in Definition~\ref{areafncts}.	
\item
Let $\rho\in\mathbf{q}$ be a regularized quasi-ultrametric as in Definition~\ref{D2.3}.
Then, by an argument similar to that used in the proof of
Lemma~\ref{1955}, we conclude that, for any
$f\in(\mathring{{\mathcal{G}}}_0^\varepsilon(\beta,\gamma))'$,
the functions $\dot{g}^s_q(f)$, $\dot{\mathcal{S}}^s_q(f)$,
and $(\dot{g}^*_\lambda)^s_q(f)$ are $\mu$-measurable.
\item
If $\dot{g}^s_q(f)$ is constructed using $\rho$ and
$\widetilde{\dot{g}^s_q}(f)$ is constructed using
$\widetilde{\rho}\in\mathbf{q}$ (where all
$\rho$-balls and all $\widetilde{\rho}$-balls are $\mu$-measurable),
then $\dot{g}^s_q(f)\sim\widetilde{\dot{g}^s_q}(f)$ pointwise in $X$.
Similar equivalences for $\dot{\mathcal{S}}^s_q(f)$
and $(\dot{g}^*_\lambda)^s_q(f)$ also hold.
Thus, without loss of generality, we may assume that $\rho\in\mathbf{q}$
in Definition~\ref{homoareafuncts} is a regularized quasi-ultrametric
in the remainder of this monograph.
\end{enumerate}
\end{remark}

By the definitions of $\dot{g}^s_q(f)$ and $\dot{F}^s_{p,q}(X)$,
we immediately obtain the following conclusion.

\begin{theorem}
Let $(X,\mathbf{q},\mu)$, $n$, $\varepsilon$, $s$, $p$, $q$, $\beta$, and $\gamma$ be
the same as in Definition~\ref{DefT--L} with $p\neq\infty$.
Then $f\in\dot{F}^{p,q}_s(X)$ if and only if
$f\in(\mathring{{\mathcal{G}}}_0^\varepsilon(\beta,\gamma))'$
and $\dot{g}^s_q(f)\in L^p(X)$.
Moreover,
for any
$f\in(\mathring{{\mathcal{G}}}_0^\varepsilon(\beta,\gamma))'$,
$\|\dot{g}^s_q(f)\|_{L^p(X)}
=\|f\|_{\dot{F}^{p,q}_s(X)}$.
\end{theorem}

\section{Lusin Area Function
Characterizations of
$\dot{F}^{p,q}_s(X)$}
\label{Section4.2}

In this section, we establish the Lusin area function
characterization of $\dot{F}^{p,q}_s(X)$ for $p<\infty$
on the ultra-RD-space $(X,\mathbf{q},\mu)$.
We always assume that, for any $\rho\in\mathbf{q}$,
$\mathrm{diam}_\rho(X)=\infty$
and hence $\mu(X)=\infty$.
The main theorem of this section is as follows
and the key tool used in its proof is the
homogeneous Plancherel--P\'olya
inequality in Lemma~\ref{P-P } below.

\begin{theorem}\label{Lusin}
Let $(X,\mathbf{q},\mu)$ be an ultra-RD-space
with upper dimension $n\in(0,\infty)$,
where $\mu$ is assumed to be a Borel-semiregular measure on $X$.
Let $\varepsilon$, $s$, $p$, $q$, $\beta$, and $\gamma$ be
the same as in Definition~\ref{DefT--L} with $p\neq\infty$.
Then $f\in\dot{F}^{p,q}_s(X)$
if and only if
$f\in(\mathring{{\mathcal{G}}}_0^\varepsilon(\beta,\gamma))'$
and $\dot{\mathcal{S}}^s_q(f)\in L^p(X)$.
Moreover, there exists a constant $C\in[1,\infty)$
such that,
for any $f\in(\mathring{{\mathcal{G}}}_0^\varepsilon(\beta,\gamma))'$,
\begin{align}
\label{1855}
C^{-1}\left\|\dot{\mathcal{S}}^s_q(f)\right\|_{L^p(X)}
\leq\|f\|_{\dot{F}^{p,q}_s(X)}
\leq C\left\|\dot{\mathcal{S}}^s_q(f)\right\|_{L^p(X)}.
\end{align}
\end{theorem}

To show Theorem~\ref{Lusin},
we need the following homogeneous Plancherel--P\'olya
inequality.\index[words]{Plancherel--P\'olya inequality}

\begin{lemma}\label{P-P }
Let $(X,\mathbf{q},\mu)$, $n$, $\varepsilon$, $s$, $p$, $q$, $\beta$, and $\gamma$ be
the same as in Theorem~\ref{Lusin}.
Assume that $\{\mathcal{S}_{2^{-k}}\}_{k\in\mathbb{Z}}$
and $\{\mathcal{Q}_{2^{-k}}\}_{k\in\mathbb{Z}}$
are two homogeneous approximations of the identity
of order $\varepsilon$ as
in Definition~\ref{DefHATI} and, for any $k\in\mathbb{Z}$, let
$\mathcal{D}_k:=\mathcal{S}_{2^{-k}}-\mathcal{S}_{2^{-k+1}}$
and  $\mathcal{P}_k:=\mathcal{Q}_{2^{-k}}-\mathcal{Q}_{2^{-k+1}}$.
Let $\{Q^{k,v}_\tau: k\in\mathbb{Z},\ \tau\in I_k,v=1,\dots,N(k,\tau)\}$
be the collection of dyadic cubes as in Remark~\ref{jcubes}
and suppose that $\rho\in\mathbf{q}$ is a
regularized quasi-ultrametric as in Definition~\ref{D2.3}.
Then, for any $c_0\in[1,\infty)$, there exists a positive constant $C$ such that,
for any $f\in(\mathring{{\mathcal{G}}}_0^\varepsilon(\beta,\gamma))'$,
\begin{align}
\label{P-P ineq}
&\left\|\left[\sum_{k=-\infty}^{\infty}
\sum_{\tau\in I_k}\sum_{v=1}^{N(k,\tau)}
2^{ksq}\sup_{z\in B_\rho(z^{k,v}_\tau,c_02^{-k})}
\left|\mathcal{P}_kf(z)\right|^q
\mathbf{1}_{Q^{k,v}_\tau}\right]^\frac{1}{q}\right\|_{L^p(X)}
\nonumber\\
&\quad\leq C\left\|\left[\sum_{k=-\infty}^{\infty}
\sum_{\tau\in I_k}\sum_{v=1}^{N(k,\tau)}
2^{ksq}\inf_{z\in Q^{k,v}_\tau}
\left|\mathcal{D}_kf(z)\right|^q
\mathbf{1}_{Q^{k,v}_\tau}\right]^\frac{1}{q}\right\|_{L^p(X)}
\end{align}
with the usual modification when $q=\infty$, where,
for any $k\in\mathbb{Z}$, $\tau\in I_k$, and $v\in\{1,\ldots,N(k,\tau)\}$,
$z^{k,v}_\tau$ is the center of $Q^{k,v}_\tau$ as in Lemma~\ref{DyadicSys}.
\end{lemma}

\begin{proof}
From Theorem~\ref{HD}, we deduce that
there exists a family
of linear operators $\{\widetilde{\mathcal{D}}_k\}_{k\in\mathbb{Z}}$
whose kernels satisfy
the same conditions as in Theorem~\ref{A3.47} and are such that,
for any $f\in(\mathring{{\mathcal{G}}}_0^\varepsilon(\beta,\gamma))'$,
the third equality in \eqref{3.29} holds in
$(\mathring{{\mathcal{G}}}_0^\varepsilon(\beta,\gamma))'$.
As such, for any $k\in\mathbb{Z}$,
$f\in(\mathring{{\mathcal{G}}}_0^\varepsilon(\beta,\gamma))'$,
and $z\in X$,
\begin{align}\label{PP1}
\mathcal{P}_kf(z)
&=\left\langle f,P_k(z,\cdot)\right\rangle
\nonumber\\
&=\left\langle\sum_{k'=-\infty}^{\infty}\sum_{\tau'\in I_{k'}}
\sum_{v'=1}^{N(k',\tau')}\mu(Q_{\tau'}^{k',v'})
\widetilde{D}_{k'}(\cdot,y_{\tau'}^{k',v'})
\mathcal{D}_{k'}f(y_{\tau'}^{k',v'}),P_k(z,\cdot)\right\rangle
\nonumber\\
&=\sum_{k'=-\infty}^{\infty}\sum_{\tau'\in I_{k'}}
\sum_{v'=1}^{N(k',\tau')}\mu(Q_{\tau'}^{k',v'})\left\langle
\widetilde{D}_{k'}(\cdot,y_{\tau'}^{k',v'})
,P_k(z,\cdot)\right\rangle\mathcal{D}_{k'}f(y_{\tau'}^{k',v'})
\nonumber\\
&=\sum_{k'\in\mathbb{Z}}
\sum_{\tau'\in I_k}\sum_{v'=1}^{N(k',\tau')}
\mu(Q_{\tau'}^{k',v'})
P_k\widetilde{D}_{k'}(z,y_{\tau'}^{k',v'})
\mathcal{D}_{k'}(f)(y_{\tau'}^{k',v'}),
\end{align}
where $y_{\tau'}^{k',v'}$ is any fixed point in $Q_{\tau'}^{k',v'}$.
By Definition~\ref{D2.1}(iii) and Lemma~\ref{DyadicSys}(iv),
we find that, for any $x\in Q_{\tau}^{k,v}$ and $z\in B_\rho(z^{k,v}_\tau,c_02^{-k})$,
$$
\rho(x,z)\leq C_\rho\max\left\{\rho(x,z^{k,v}_\tau),\,
\rho(z^{k,v}_\tau,z)\right\}\lesssim 2^{-k} \leq 2^{-(k\land k')},
$$
which, in conjunction with Definition~\ref{D2.1}(iii),
Lemma~\ref{A3.2}(vi), and \eqref{doublingcond}, implies that
$$
\rho(x,y^{k',v'}_{\tau'})\leq C_\rho\max\left\{\rho(x,z),\,
\rho(z,y^{k',v'}_{\tau'})\right\}\lesssim 2^{-(k\land k')}
+\rho(z,y^{k',v'}_{\tau'})
$$
and
\begin{align*}
&(V_\rho)_{2^{-(k\land k')}}(x)+V_\rho(x,y^{k',v'}_{\tau'})\\
&\quad\sim(V_\rho)_{2^{-(k\land k')}+\rho(x,y^{k',v'}_{\tau'})}(x)
\lesssim(V_\rho)_{2^{-(k\land k')}+\rho(z,y^{k',v'}_{\tau'})}(z)
\sim(V_\rho)_{2^{-(k\land k')}}(z)+V_\rho(z,y^{k',v'}_{\tau'}).
\end{align*}
From this and Lemma~\ref{1122},
we infer that, for any given $\theta'\in(0,\varepsilon)$ and
for any $x\in Q_{\tau}^{k,v}$ and $z\in B_\rho(z^{k,v}_\tau,c_02^{-k})$,
\begin{align}
\label{PP2}
&\left|P_k\widetilde{D}_{k'}(z,y^{k',v'}_{\tau'})\right|
\nonumber\\
&\quad\lesssim 2^{-|k-k'|\theta'}
\frac{1}{(V_\rho)_{2^{-(k\land k')}}(z)+V_\rho(z,y^{k',v'}_{\tau'})}
\left[\frac{2^{-(k\land k')}}
{2^{-(k\land k')}+\rho(z,y^{k',v'}_{\tau'})}\right]^{\gamma}
\nonumber\\
&\quad\lesssim 2^{-|k-k'|\theta'}
\frac{1}{(V_\rho)_{2^{-(k\land k')}}(x)+V_\rho(x,y^{k',v'}_{\tau'})}
\left[\frac{2^{-(k\land k')}}
{2^{-(k\land k')}+\rho(x,y^{k',v'}_{\tau'})}\right]^{\gamma}.
\end{align}
Since $s\in(-\varepsilon,\varepsilon)$ and
$\max\{\frac{n}{n+\varepsilon},\,
\frac{n}{n+\varepsilon+s}\}<p,q\leq\infty$,
we choose $\theta'\in(|s|,\varepsilon)$ such that
there exists
\begin{align}\label{2143}
r\in\left(\max\left\{\frac{n}{n+\theta'},\,\frac{n}{n+\theta'+s}\right\},
\min\left\{p,\,q,\,1\right\}\right).
\end{align}
By this, \eqref{PP1}, \eqref{PP2}, and Lemma~\ref{L4.5},
we conclude that, for any
$f\in(\mathring{{\mathcal{G}}}_0^\varepsilon(\beta,\gamma))'$
and $x\in X$,
\begin{align}
\label{2147}
&\sum_{k=-\infty}^{\infty}
\sum_{\tau\in I_k}\sum_{v=1}^{N(k,\tau)}
2^{ksq}\sup_{z\in B_\rho(z^{k,v}_\tau,c_02^{-k})}
\left|\mathcal{P}_kf(z)\right|^q
\mathbf{1}_{Q^{k,v}_\tau}(x)
\nonumber\\
&\quad\lesssim\sum_{k=-\infty}^\infty
2^{ksq}\Bigg[\sum_{k'=-\infty}^{\infty}2^{-k's}2^{-|k-k'|\theta'}
\sum_{\tau'\in I_k}\sum_{v'=1}^{N(k',\tau')}
\mu(Q_{\tau'}^{k',v'})
\nonumber\\
&\qquad\times
\frac{1}{(V_\rho)_{2^{-(k\land k')}}(x)+V_\rho(x,y^{k',v'}_{\tau'})}
\left[\frac{2^{-(k\land k')}}
{2^{-(k\land k')}+\rho(x,y^{k',v'}_{\tau'})}\right]^{\gamma}
\left|2^{k's}\mathcal{D}_{k'}f(y_{\tau'}^{k',v'})\right|\Bigg]^q
\nonumber\\
&\quad\lesssim\sum_{k=-\infty}^\infty
\Bigg\{\sum_{k'\in\mathbb{Z}}
2^{(k-k')s-|k-k'|\theta'+n(k\land k'-k')(1-\frac{1}{r})}
\nonumber\\
&\qquad\times\left[\mathcal{M}_\rho\left(\sum_{\tau'\in I_{k'}}
\sum_{v'=1}^{N(k',\tau')}2^{k'sr}
\left|\mathcal{D}_{k'}f(y_{\tau'}^{k',v'})\right|^r
\mathbf{1}_{Q_{\tau'}^{k',v'}}\right)(x)\right]^\frac{1}{r}
\Bigg\}^q.
\end{align}
Given that $y^{k',v'}_{\tau'}\in Q^{k',v'}_{\tau'}$ was arbitrary, the
inequalities in \eqref{2147} hold with
$|\mathcal{D}_{k'}f(y^{k',v'}_{\tau'})|$ replaced by
$\inf_{u\in Q_{\tau'}^{k',v'}}|\mathcal{D}_{k'}f(u)|$.

Now, on
the one hand, when $q\in(0,1]$, from Lemma~\ref{discreteinequality}
with $r$ therein replaced by $q$
and the fact that, for any $k,k'\in\mathbb{Z}$,
\begin{align}
\label{kk'<0}
&(k-k')s-|k-k'|\theta'+(k\land k'-k')n\left(1-\frac{1}{r}\right)\nonumber\\
&\quad=
\begin{cases}
(k-k')(s-\theta')<0&\text{ if }k\ge k',\\
(k-k')\left[s+\theta'+n\left(1-\frac{1}{r}\right)\right]<0&\text{ if }k<k'
\end{cases}
\end{align}
[here, we used $|s|<\theta'$ and \eqref{2143}],
we deduce that
\begin{align}
\label{qleq1}
&\sum_{k=-\infty}^\infty
\Bigg\{\sum_{k'\in\mathbb{Z}}
2^{(k-k')s-|k-k'|\theta'+(k\land k'-k')n(1-\frac{1}{r})}
\nonumber\\
&\qquad\times\left[\mathcal{M}_\rho\left(\sum_{\tau'\in I_{k'}}
\sum_{v'=1}^{N(k',\tau')}2^{k'sr}
\inf_{u\in Q_{\tau'}^{k',v'}}\left|\mathcal{D}_{k'}f(u)\right|^r
\mathbf{1}_{Q_{\tau'}^{k',v'}}\right)(x)\right]^\frac{1}{r}
\Bigg\}^q
\nonumber\\
&\quad
\leq\sum_{k=-\infty}^\infty
\sum_{k'\in\mathbb{Z}}
2^{\{(k-k')s-|k-k'|\theta'+(k\land k'-k')n(1-\frac{1}{r})\}q}
\nonumber\\
&\qquad\times\left[\mathcal{M}_\rho\left(\sum_{\tau'\in I_{k'}}
\sum_{v'=1}^{N(k',\tau')}2^{k'sr}
\inf_{u\in Q_{\tau'}^{k',v'}}\left|\mathcal{D}_{k'}f(u)\right|^r
\mathbf{1}_{Q_{\tau'}^{k',v'}}\right)(x)\right]^\frac{q}{r}
\nonumber\\
&\quad
\lesssim\sum_{k'\in\mathbb{Z}}
\left[\mathcal{M}_\rho\left(\sum_{\tau'\in I_{k'}}
\sum_{v'=1}^{N(k',\tau')}2^{k'sr}
\inf_{u\in Q_{\tau'}^{k',v'}}\left|\mathcal{D}_{k'}f(u)\right|^r
\mathbf{1}_{Q_{\tau'}^{k',v'}}\right)(x)\right]^\frac{q}{r}.
\end{align}
On the other hand, when $q\in(1,\infty]$,
by H\"older's inequality,
we find that
\begin{align*}
&\sum_{k=-\infty}^\infty
\Bigg\{\sum_{k'\in\mathbb{Z}}
2^{(k-k')s-|k-k'|\theta'+(k\land k'-k')n(1-\frac{1}{r})}
\\
&\qquad\times\left[\mathcal{M}_\rho\left(\sum_{\tau'\in I_{k'}}
\sum_{v'=1}^{N(k',\tau')}2^{k'sr}
\inf_{u\in Q_{\tau'}^{k',v'}}\left|\mathcal{D}_{k'}f(u)\right|^r
\mathbf{1}_{Q_{\tau'}^{k',v'}}\right)(x)\right]^\frac{1}{r}
\Bigg\}^q
\\
&\quad
\leq\sum_{k=-\infty}^\infty
\left\{\sum_{k'\in\mathbb{Z}}
2^{(k-k')s-|k-k'|\theta'+(k\land k'-k')n(1-\frac{1}{r})}\right\}^\frac{q}{q'}
\\
&\qquad\times\sum_{k'\in\mathbb{Z}}
2^{(k-k')s-|k-k'|\theta'+(k\land k'-k')n(1-\frac{1}{r})}
\\
&\qquad\times\left[\mathcal{M}_\rho\left(\sum_{\tau'\in I_{k'}}
\sum_{v'=1}^{N(k',\tau')}2^{k'sr}
\inf_{u\in Q_{\tau'}^{k',v'}}\left|\mathcal{D}_{k'}f(u)\right|^r
\mathbf{1}_{Q_{\tau'}^{k',v'}}\right)(x)\right]^\frac{q}{r},
\end{align*}
which, together with \eqref{kk'<0}, further implies that
\begin{align*}
&\sum_{k=-\infty}^\infty
\Bigg\{\sum_{k'\in\mathbb{Z}}
2^{(k-k')s-|k-k'|\theta'+(k\land k'-k')n(1-\frac{1}{r})}
\\
&\qquad\times\left[\mathcal{M}_\rho\left(\sum_{\tau'\in I_{k'}}
\sum_{v'=1}^{N(k',\tau')}2^{k'sr}
\inf_{u\in Q_{\tau'}^{k',v'}}\left|\mathcal{D}_{k'}f(u)\right|^r
\mathbf{1}_{Q_{\tau'}^{k',v'}}\right)(x)\right]^\frac{1}{r}
\Bigg\}^q
\\
&\quad
\lesssim\sum_{k'\in\mathbb{Z}}
\left[\mathcal{M}_\rho\left(\sum_{\tau'\in I_{k'}}
\sum_{v'=1}^{N(k',\tau')}2^{k'sr}
\inf_{u\in Q_{\tau'}^{k',v'}}\left|\mathcal{D}_{k'}f(u)\right|^r
\mathbf{1}_{Q_{\tau'}^{k',v'}}\right)(x)\right]^\frac{q}{r}.
\end{align*}
From this, \eqref{2147}, \eqref{qleq1},
and the fact that the cubes $\{Q_{\tau'}^{k',v'}\}_{\tau',v'}$
are mutually disjoint for any given $k'\in\mathbb{Z}$,
we infer that,
for any $f\in(\mathring{{\mathcal{G}}}_0^\varepsilon(\beta,\gamma))'$,
\begin{align}\label{2145}
&\left\|\left[\sum_{k=-\infty}^{\infty}
\sum_{\tau\in I_k}\sum_{v=1}^{N(k,\tau)}
2^{ksq}\sup_{z\in B_\rho(z^{k,v}_\tau,c_02^{-k})}
\left|\mathcal{P}_kf(z)\right|^q
\mathbf{1}_{Q^{k,v}_\tau}\right]^\frac{1}{q}\right\|_{L^p(X)}
\nonumber\\
&\quad\lesssim\left\|\left\{\sum_{k'=-\infty}^\infty
\left[\mathcal{M}_\rho\left(\sum_{\tau'\in I_{k'}}\sum_{v'=1}^{N(k',\tau')}
2^{k'sr}\inf_{u\in Q_{\tau'}^{k',v'}}
\left|\mathcal{D}_{k'}f(u)\right|^r
\mathbf{1}_{Q_{\tau'}^{k',v'}}\right)\right]^\frac{q}{r}
\right\}^\frac{1}{q}\right\|_{L^{p}(X)}
\nonumber\\
&\quad\lesssim\left\|\left\{\sum_{k'=-\infty}^\infty
\left[\sum_{\tau'\in I_{k'}}\sum_{v'=1}^{N(k',\tau')}
2^{k'sr}\inf_{u\in Q_{\tau'}^{k',v'}}
\left|\mathcal{D}_{k'}f(u)\right|^r
\mathbf{1}_{Q_{\tau'}^{k',v'}}\right]^\frac{q}{r}
\right\}^\frac{r}{q}\right\|_{L^\frac{p}{r}(X)}^\frac{1}{r}\nonumber\\
&\qquad\qquad\text{by Lemma~\ref{Fefferman--Stein}
with $u:=\frac{q}{r}$ and $p:=\frac{p}{r}$}
\nonumber\\
&\quad=\left\|\left[\sum_{k'=-\infty}^{\infty}
\sum_{\tau'\in I_{k'}}\sum_{v'=1}^{N(k',\tau')}
2^{k'sq}\inf_{u\in Q_{\tau'}^{k',v'}}
\left|\mathcal{D}_{k'}f(u)\right|^q
\mathbf{1}_{Q_{\tau'}^{k',v'}}\right]^\frac{1}{q}\right\|_{L^p(X)},
\end{align}
which implies \eqref{P-P ineq}.
This finishes the proof of Lemma~\ref{P-P }.
\end{proof}

\begin{remark}\label{zj2037}
Let the notation be as in Lemma~\ref{P-P }.
By Lemma~\ref{P-P },
we easily conclude that the definition of $\dot{F}^{p,q}_s(X)$
is independent of the choice of
the approximation of the identity;
that is, for any
$f\in(\mathring{{\mathcal{G}}}_0^\varepsilon(\beta,\gamma))'$,
\begin{align*}
\left\|\left(\sum_{k\in\mathbb{Z}}
2^{ksq}\left|\mathcal{D}_kf\right|^q\right)^\frac{1}{q}\right\|_{L^p(X)}
\sim\left\|\left(\sum_{k\in\mathbb{Z}}
2^{ksq}\left|\mathcal{P}_kf\right|^q\right)^\frac{1}{q}\right\|_{L^p(X)},
\end{align*}
where the positive equivalence constants are independent of $f$.
A similar conclusion was also obtained in
\cite[Proposition~5.6(ii)]{hmy08}
in the setting of RD-spaces
and \cite[Proposition~3.13]{WHHY}
in the setting of spaces of homogeneous type.
\end{remark}

Next, we prove Theorem \ref{Lusin}.

\begin{proof}[Proof of Theorem \ref{Lusin}]
Let $\rho\in\mathbf{q}$ be a regularized
quasi-ultrametric.
Note that all of the conclusions of Theorem~\ref{Lusin} will follow
once we establish \eqref{1855}. To this end,
we first show the second inequality of \eqref{1855}.
On the one hand, by Remark~\ref{jcubes} and Lemma~\ref{DyadicSys}(i),
we find that,
for any $f\in(\mathring{{\mathcal{G}}}_0^\varepsilon(\beta,\gamma))'$ and
$\mu$-almost every $x\in X$,
\begin{align}\label{T5.3(1)}
\sum_{k=-\infty}^{\infty}
2^{ksq}\left|\mathcal{D}_kf(x)\right|^q
&=\sum_{k=-\infty}^{\infty}
2^{ksq}\sum_{\tau\in I_k}
\sum_{v=1}^{N(k,\tau)}
\left|\mathcal{D}_kf(x)\right|^q\mathbf{1}_{Q^{k,v}_\tau}(x)
\nonumber\\
&\leq\sum_{k=-\infty}^{\infty}
2^{ksq}\sum_{\tau\in I_k}
\sum_{v=1}^{N(k,\tau)}\sup_{z\in Q^{k,v}_\tau}
\left|\mathcal{D}_kf(z)\right|^q\mathbf{1}_{Q^{k,v}_\tau}(x).
\end{align}
On the other hand,
note that, from Lemma~\ref{DyadicSys}(iv),
it follows that, for any $x\in Q^{k,v}_\tau$
with $k\in\mathbb{Z}$,
$\tau\in I_k$, and $v\in\{1,...,N(k,\tau)\}$, $\rho(z_\tau^{k,v},x)<C_r2^{-k-j}$,
where $z^{k,v}_\tau$ is the center of $Q^{k,v}_\tau$,
which, combined with both Definition~\ref{D2.1}(iii)
and \eqref{2cubej}, further implies that,
for any $y\in B_\rho(x,2^{-k})$,
\begin{align*}
\rho(z_\tau^{k,v},y)\leq C_\rho\max\left\{\rho(z_\tau^{k,v},x),\,\rho(x,y)\right\}
< C_\rho2^{-k}
\end{align*}
and hence
\begin{align}\label{21412}
B_\rho(x,2^{-k})\subset B_\rho(z_\tau^{k,v},C_\rho2^{-k}).
\end{align}
By this, Lemma~\ref{DyadicSys}(v), and \eqref{upperdoub},
we conclude that, for any $x\in Q^{k,v}_\tau$,
\begin{align}\label{muQsimV2-k}
\mu(Q^{k,v}_\tau)&\ge\mu(B_\rho(z_\tau^{k,v},C_l2^{-k}))
\gtrsim\mu(B_\rho(z_\tau^{k,v},C_\rho2^{-k}))
\ge\mu(B_\rho(x,2^{-k})).
\end{align}
Moreover, from Lemma~\ref{DyadicSys}(iv)
and \eqref{2cubej}, we deduce that,
for any given $x\in Q^{k,v}_\tau$
and for any $y\in Q^{k,v}_\tau$,
\begin{align*}
\rho(x,y)\leq\mathrm{diam}_\rho(Q_\tau^{k,v})<C_r2^{-k-j}<2^{-k}
\end{align*}
and hence
\begin{align}\label{21411}
Q^{k,v}_\tau\subset B_\rho(x,2^{-k}).
\end{align}
By this, \eqref{muQsimV2-k}, and Lemma~\ref{DyadicSys}(i), we find that,
for any $f\in(\mathring{{\mathcal{G}}}_0^\varepsilon(\beta,\gamma))'$ and  $x\in X$,
\begin{align*}
&\sum_{k=-\infty}^{\infty}\sum_{\tau\in I_k}
\sum_{v=1}^{N(k,\tau)}2^{ksq}
\inf_{z\in Q^{k,v}_\tau}
\left|\mathcal{D}_kf(z)\right|^q\mathbf{1}_{Q^{k,v}_\tau}(x)
\\
&\quad\leq\sum_{k=-\infty}^{\infty}
\sum_{\tau\in I_k}\sum_{v=1}^{N(k,\tau)}
2^{ksq}\frac{1}{\mu(Q^{k,v}_\tau)}
\int_{Q^{k,v}_\tau}
\left|\mathcal{D}_kf(z)\right|^q\,d\mu(z)\mathbf{1}_{Q^{k,v}_\tau}(x)
\\
&\quad\lesssim\sum_{k=-\infty}^{\infty}
\sum_{\tau\in I_k}\sum_{v=1}^{N(k,\tau)}
2^{ksq}\frac{1}{\mu(B_\rho(x,2^{-k}))}
\int_{B_\rho(x,2^{-k})}
\left|\mathcal{D}_kf(z)\right|^q\,d\mu(z)\mathbf{1}_{Q^{k,v}_\tau}(x)\\
&\quad=\left[\dot{\mathcal{S}}^s_q(f)(x)\right]^q,
\end{align*}
which, together with \eqref{T5.3(1)}
and Lemma~\ref{P-P } with $c_0:=1$,
further implies that
\begin{align*}
\|f\|_{\dot{F}^{p,q}_s(X)}
&\leq\left\|\left[\sum_{k=-\infty}^{\infty}
\sum_{\tau\in I_k}\sum_{v=1}^{N(k,\tau)}
2^{ksq}\sup_{z\in Q^{k,v}_\tau}
\left|\mathcal{D}_kf(z)\right|^q
\mathbf{1}_{Q^{k,v}_\tau}\right]^\frac{1}{q}\right\|_{L^p(X)}
\\
&\leq\left\|\left[\sum_{k=-\infty}^{\infty}
\sum_{\tau\in I_k}\sum_{v=1}^{N(k,\tau)}
2^{ksq}\sup_{z\in B_\rho(z^{k,v}_\tau,2^{-k})}
\left|\mathcal{D}_kf(z)\right|^q
\mathbf{1}_{Q^{k,v}_\tau}\right]^\frac{1}{q}\right\|_{L^p(X)}
\\
&\lesssim\left\|\left[\sum_{k=-\infty}^{\infty}
\sum_{\tau\in I_k}\sum_{v=1}^{N(k,\tau)}
2^{ksq}\inf_{z\in Q^{k,v}_\tau}
\left|\mathcal{D}_kf(z)\right|^q
\mathbf{1}_{Q^{k,v}_\tau}\right]^\frac{1}{q}\right\|_{L^p(X)}
\leq\left\|\dot{\mathcal{S}}^s_q(f)\right\|_{L^p(X)}.
\end{align*}
This finishes the proof of the second inequality of \eqref{1855}.

Now, we prove the first inequality of \eqref{1855}.
Note that, from \eqref{21411} and \eqref{muQsimV2-k}, we infer that,
for any given $x\in Q_\tau^{k,v}$
with $k\in\mathbb{Z}$, $\tau\in I_k$,
and $v\in\{1,...,N(k,\tau)\}$,
\begin{align}\label{1055}
\mu\left(B_\rho(x,2^{-k})\right)\ge\mu\left(Q_\tau^{k,v}\right)
\gtrsim\mu\left(B_\rho(z_\tau^{k,v},C_\rho2^{-k})\right).
\end{align}
By this, \eqref{21412},
Lemma~\ref{P-P } with $c_0:=C_\rho$, and Lemma~\ref{DyadicSys}(i),
we conclude that,
for any $f\in(\mathring{{\mathcal{G}}}_0^\varepsilon(\beta,\gamma))'$,
\begin{align*}
\left\|\dot{\mathcal{S}}^s_q(f)\right\|_{L^p(X)}
&=\left\|\left[\sum_{k=-\infty}^{\infty}
\sum_{\tau\in I_k}\sum_{v=1}^{N(k,\tau)}
2^{ksq}\fint_{B_\rho(\cdot,2^{-k})}
\left|\mathcal{D}_kf(z)\right|^q\,d\mu(z)
\mathbf{1}_{Q_\tau^{k,v}}(\cdot)\right]^\frac{1}{q}
\right\|_{L^p(X)}
\\
&\lesssim\left\|\left[\sum_{k=-\infty}^{\infty}
\sum_{\tau\in I_k}\sum_{v=1}^{N(k,\tau)}
2^{ksq}
\fint_{B_\rho(z_\tau^{k,v},C_\rho2^{-k})}
\left|\mathcal{D}_kf(z)\right|^q\,d\mu(z)
\mathbf{1}_{Q_\tau^{k,v}}(\cdot)\right]^\frac{1}{q}
\right\|_{L^p(X)}
\\
&\leq\left\|\left[\sum_{k=-\infty}^{\infty}
\sum_{\tau\in I_k}\sum_{v=1}^{N(k,\tau)}
2^{ksq}\sup_{z\in B_\rho(z_\tau^{k,v},C_\rho2^{-k})}
\left|\mathcal{D}_kf(z)\right|^q
\mathbf{1}_{Q^{k,v}_\tau}(\cdot)\right]^\frac{1}{q}\right\|_{L^p(X)}
\\
&\lesssim\left\|\left[\sum_{k=-\infty}^{\infty}
\sum_{\tau\in I_k}\sum_{v=1}^{N(k,\tau)}
2^{ksq}\inf_{z\in Q^{k,v}_\tau}
\left|\mathcal{D}_kf(z)\right|^q
\mathbf{1}_{Q^{k,v}_\tau}(\cdot)\right]^\frac{1}{q}\right\|_{L^p(X)}
\leq\|f\|_{\dot{F}^{p,q}_s(X)}.
\end{align*}
This finishes the proof of the first inequality of \eqref{1855}
and hence Theorem~\ref{Lusin}.
\end{proof}

\section[Littlewood--Paley $g^*_\lambda$-Function
Characterization of $\dot{F}^{p,q}_s(X)$]
{Littlewood--Paley $g^*_\lambda$-Function
Characterization of\\ $\dot{F}^{p,q}_s(X)$}
\label{Section4.3}

In this section, we aim to establish
the Littlewood--Paley $g^*_\lambda$-function
characterization of $\dot{F}^{p,q}_s(X)$
on an ultra-RD-space $(X,\mathbf{q},\mu)$.
We always assume that, for any $\rho\in\mathbf{q}$,
$\mathrm{diam}_\rho(X)=\infty$
and hence $\mu(X)=\infty$.
The main results of this section are Propositions~\ref{g*1}
and~\ref{g*lammda} below.
We achieve these by using Theorem~\ref{Lusin}
and the following arguments.
On the one hand, it is straightforward to show that the
Lusin area function is pointwise controlled by the
Littlewood--Paley $g^*_\lambda$-function.
On the other hand,
to estimate the $L^p$ norm of the Littlewood--Paley
$g_\lambda^*$-function via that of the Lusin are function,
we divide the proof into two cases,
where we use a duality argument and Tonelli's theorem
when $p\in[q,\infty)$ and an auxiliary norm estimate
related to both the homogeneous Lusin area function and the
homogeneous Littlewood--Paley auxiliary function
(see Proposition~\ref{Prop5.5} below) when $p\in(0,q)$.

We first consider the case $p\in[q,\infty)$.

\begin{proposition}\label{g*1}
Let $(X,\mathbf{q},\mu)$, $n$, $\varepsilon$, $s$, $p$, $q$, $\beta$, and $\gamma$ be
the same as in Theorem~\ref{Lusin}.
Assume that $p\in[q,\infty)$ and $\lambda\in(n,\infty)$.
Then $f\in\dot{F}^{p,q}_s(X)$
if and only if $f\in(\mathring{{\mathcal{G}}}_0^\varepsilon(\beta,\gamma))'$
and $(\dot{g}^*_\lambda)^s_q(f)\in L^p(X)$.
Moreover, there exists a constant $C\in[1,\infty)$,
independent of $f$, such that, for any
$f\in(\mathring{{\mathcal{G}}}_0^\varepsilon(\beta,\gamma))'$,
\begin{align}
\label{2036}
C^{-1}\left\|(\dot{g}^*_\lambda)^s_q(f)\right\|_{L^p(X)}
\leq\|f\|_{\dot{F}^{p,q}_s(X)}
\leq C\left\|(\dot{g}^*_\lambda)^s_q(f)\right\|_{L^p(X)}.
\end{align}
\end{proposition}

\begin{proof}
Let $\rho\in\mathbf{q}$ be a regularized
quasi-ultrametric as in Definition~\ref{D2.3}.
Note that all of the conclusions of this proposition
will follow once we establish \eqref{2036}.
To this end, we first prove the second inequality of \eqref{2036}.
Observe that, by Lemma~\ref{A3.2}(vii), we conclude that,
for any $f\in(\mathring{{\mathcal{G}}}_0^\varepsilon(\beta,\gamma))'$
and $x\in X$,
\begin{align}\label{ddq-339}
\dot{\mathcal{S}}^s_q(f)(x)
&=2^\frac{\lambda}{q}\left[\sum_{k\in\mathbb{Z}}
2^{ksq}\int_{B_\rho(x,2^{-k})}\left|\mathcal{D}_kf(x)\right|^q
\left(\frac{2^{-k}}{2^{-k}+2^{-k}}\right)^\lambda
\frac{d\mu(y)}{(V_\rho)_{2^{-k}}(x)}\right]^\frac{1}{q}
\nonumber\\
&<2^\frac{\lambda}{q}\left\{\sum_{k\in\mathbb{Z}}
2^{ksq}\int_{B_\rho(x,2^{-k})}\left|\mathcal{D}_kf(x)\right|^q
\left[\frac{2^{-k}}{2^{-k}+\rho(x,y)}\right]^\lambda
\frac{d\mu(y)}{(V_\rho)_{2^{-k}}(x)}\right\}^\frac{1}{q}
\nonumber\\
&\sim2^\frac{\lambda}{q}\left\{\sum_{k\in\mathbb{Z}}
2^{ksq}\int_{B_\rho(x,2^{-k})}\left|\mathcal{D}_kf(x)\right|^q
\left[\frac{2^{-k}}{2^{-k}+\rho(x,y)}\right]^\lambda
\frac{d\mu(y)}{(V_\rho)_{2^{-k}}(x)+(V_\rho)_{2^{-k}}(y)}\right\}^\frac{1}{q}
\nonumber\\
&\leq2^\frac{\lambda}{q}(\dot{g}^*_\lambda)^s_q(f)(x),
\end{align}
where the implicit positive constants depend on $q$.
From this and \eqref{1855}, we deduce that
$$
\|f\|_{\dot{F}^{p,q}_s(X)}
\sim\left\|\dot{\mathcal{S}}^s_q(f)\right\|_{L^p(X)}
\lesssim\left\|(\dot{g}^*_\lambda)^s_q(f)\right\|_{L^p(X)}.
$$

Next, we show the first inequality of \eqref{2036}.
Note that we may assume $f\in\dot{F}^{p,q}_s(X)$ since,
otherwise, the first inequality
of \eqref{2036} trivially holds. To proceed, we
consider the following two cases for $p$.

\emph{Case (1)} $p=q$.
In this case, by Tonelli's theorem and \eqref{upperdoub},
we find that
\begin{align*}
&\left\|(\dot{g}^*_\lambda)^s_q(f)\right\|_{L^q(X)}^q\\
&\quad\leq\int_X\sum_{k\in\mathbb{Z}}2^{ksq}\left|\mathcal{D}_kf(y)\right|^q
\int_X
\left[\frac{2^{-k}}{2^{-k}+\rho(x,y)}\right]^\lambda
\frac{d\mu(x)\,d\mu(y)}{(V_\rho)_{2^{-k}}(x)}
\\
&\quad\lesssim\int_X\sum_{k\in\mathbb{Z}}2^{ksq}\left|\mathcal{D}_kf(y)\right|^q
\int_X\left[\frac{2^{-k}}{2^{-k}+\rho(x,y)}\right]^{\lambda-n}
\frac{d\mu(x)\,d\mu(y)}{(V_\rho)_{2^{-k}+\rho(x,y)}(x)},
\end{align*}
which, combined with
(vi) and (ii) of Lemma~\ref{A3.2},
further implies that
\begin{align*}
&\left\|(\dot{g}^*_\lambda)^s_q(f)\right\|_{L^q(X)}^q\\
&\quad\lesssim\int_X\sum_{k\in\mathbb{Z}}2^{ksq}\left|\mathcal{D}_kf(y)\right|^q
\int_X\left[\frac{2^{-k}}{2^{-k}+\rho(x,y)}\right]^{\lambda-n}
\frac{d\mu(x)\,d\mu(y)}{(V_\rho)_{2^{-k}}(y)+V_\rho(x,y)}
\\
&\quad\lesssim\int_X\sum_{k\in\mathbb{Z}}2^{ksq}
\left|\mathcal{D}_kf(y)\right|^q\,d\mu(y)
=\|f\|_{\dot{F}^{q,q}_s(X)}^q.
\end{align*}
Hence, the first inequality of \eqref{2036} holds when $p=q$.

\emph{Case (2)} $p\in(q,\infty)$.
In this case, note that $(\frac{p}{q})'\in(1,\infty)$ and
$(L^\frac{p}{q}(X))'=L^{(\frac{p}{q})'}(X)$,
where $(\frac{p}{q})'$ is the conjugate index of $\frac{p}{q}\in(1,\infty)$.
Then,
from \eqref{upperdoub} and Tonelli's theorem,
we infer that, for any $\varphi\in L^{(\frac{p}{q})'}(X)$,
\begin{align*}
&\left|\left\langle\left[(\dot{g}^*_\lambda)^s_q(f)\right]^q,
\varphi\right\rangle\right|
\\
&\quad\leq\int_X\sum_{k\in\mathbb{Z}}2^{ksq}
\int_X\left|\mathcal{D}_kf(y)\right|^q
\left[\frac{2^{-k}}{2^{-k}+\rho(x,y)}\right]^\lambda
\frac{|\varphi(x)|\,d\mu(y)\,d\mu(x)}{(V_\rho)_{2^{-k}}(x)+(V_\rho)_{2^{-k}}(y)}
\\
&\quad\lesssim\int_X\sum_{k\in\mathbb{Z}}2^{ksq}
\int_X\left|\mathcal{D}_kf(y)\right|^q
\left[\frac{2^{-k}}{2^{-k}+\rho(x,y)}\right]^{\lambda-n}
\frac{|\varphi(x)|\,d\mu(y)\,d\mu(x)}{(V_\rho)_{2^{-k}+\rho(x,y)}(x)}
\\
&\quad=\int_X\sum_{k\in\mathbb{Z}}2^{ksq}
\left|\mathcal{D}_kf(y)\right|^q
\Bigg\{\int_{B_\rho(y, 2^{-k})}
\left[\frac{2^{-k}}{2^{-k}+\rho(x,y)}\right]^{\lambda-n}
\frac{|\varphi(x)|\,d\mu(x)}{(V_\rho)_{2^{-k}+\rho(x,y)}(x)}
\\
&\quad\quad+\int_{[B_\rho(y, 2^{-k})]^\complement}
\left[\frac{2^{-k}}{2^{-k}+\rho(x,y)}\right]^{\lambda-n}
\frac{|\varphi(x)|\,d\mu(x)}{(V_\rho)_{2^{-k}+\rho(x,y)}(x)}\Bigg\}\,d\mu(y),
\end{align*}
which further implies that
\begin{align*}
&\left|\left\langle\left[(\dot{g}^*_\lambda)^s_q(f)\right]^q,
\varphi\right\rangle\right|
\\
&\quad\lesssim\int_X\sum_{k\in\mathbb{Z}}2^{ksq}
\left|\mathcal{D}_kf(y)\right|^q
\Bigg\{\fint_{B_\rho(y, 2^{-k})}
|\varphi(x)|\,d\mu(x)\\
&\quad\quad+\int_{[B_\rho(y, 2^{-k})]^\complement}
\left[\frac{2^{-k}}{\rho(x,y)}\right]^{\lambda-n}
\frac{1}{V(x,y)}
|\varphi(x)|\,d\mu(x)\Bigg\}\,d\mu(y)\\
&\quad\lesssim\int_X	 \sum_{k\in\mathbb{Z}}
2^{ksq}\left|\mathcal{D}_kf(y)\right|^q
\mathcal{M}_\rho\varphi(y)\,d\mu(y)
\quad\text{by Lemma~\ref{A3.2}(iv)}\\
&\quad\lesssim\left\|\sum_{k\in\mathbb{Z}}
2^{ksq}\left|\mathcal{D}_kf\right|^q\right\|_{L^\frac{p}{q}(X)}
\left\|\mathcal{M}_\rho\varphi\right\|_{L^{(\frac{p}{q})'}(X)}
\quad\text{by H\"older's inequality}
\\
&\lesssim\|f\|_{\dot{F}^{p,q}_s(X)}^q
\left\|\varphi\right\|_{L^{(\frac{p}{q})'}(X)}
\quad\text{by Lemma~\ref{M}}.
\end{align*}
Taking the supremum over all $\varphi\in L^{(\frac{p}{q})'}(X)$
with $\|\varphi\|_{L^{(\frac{p}{q})'}(X)}\leq 1$,
we further conclude that
$$
\left\|(\dot{g}^*_\lambda)^s_q(f)\right\|_{L^p(X)}
=\left\|\left[(\dot{g}^*_\lambda)^s_q(f)\right]^q
\right\|_{L^\frac{p}{q}(X)}^\frac{1}{q}
\lesssim\|f\|_{\dot{F}^{p,q}_s(X)},
$$
which proves that the first inequality
of \eqref{2036} holds and
hence completes the proof of Proposition~\ref{g*1}.
\end{proof}

Now, we consider the case $p\in(0,q)$.
To this end, we first introduce the following homogeneous
Littlewood--Paley auxiliary function.
Let $\varepsilon$, $s$, $q$, $\beta$, and $\gamma$ be
the same as in Definition~\ref{DefT--L}
and $\theta\in[1,\infty)$.
The \emph{homogeneous Littlewood--Paley auxiliary function
$\dot{\mathcal{S}}^{s,(1)}_{q,\theta}(f)$\index[symbols]{S@$\dot{\mathcal{S}}^{s,(1)}_{q,\theta}(f)$}
of $f$ with aperture $\theta$}
is defined by setting, for any
$f\in(\mathring{{\mathcal{G}}}_0^\varepsilon(\beta,\gamma))'$
and $x\in X$,
\begin{equation}\label{auxlusin}
\dot{\mathcal{S}}^{s,(1)}_{q,\theta}(f)(x)
:=\left[\sum_{k\in\mathbb{Z}}2^{ksq}
\int_{ B_\rho(x,\theta2^{-k})}\left|\mathcal{D}_kf(y)\right|^q
\frac{d\mu(y)}{(V_\rho)_{2^{-k}}(y)}\right]^\frac{1}{q}.
\end{equation}

\begin{remark}
\begin{enumerate}
\item
From the definitions of $\dot{\mathcal{S}}^{s,(1)}_{q,\theta}(f)$
and $\dot{\mathcal{S}}^{s}_{q}(f)$,
we deduce that, for any
$f\in(\mathring{{\mathcal{G}}}_0^\varepsilon(\beta,\gamma))'$
and $x\in X$,
$\dot{\mathcal{S}}^{s,(1)}_{q,1}(f)(x)
=\dot{\mathcal{S}}^s_q(f)(x)$.
\item
If $\dot{\mathcal{S}}^{s,(1)}_{q,\theta}(f)$ is constructed using $\rho$
and $\widetilde{\dot{\mathcal{S}}^{s,(1)}_{q,\theta}}(f)$ is constructed
using $\rho'\in\mathbf{q}$ (where all $\rho'$-balls are $\mu$-measurable),
then $\dot{\mathcal{S}}^{s,(1)}_{q,\theta}(f)\sim
\widetilde{\dot{\mathcal{S}}^{s,(1)}_{q,\theta}}(f)$ pointwise in $X$.
Thus, without loss of generality, we may assume that $\rho\in\mathbf{q}$
in \eqref{auxlusin} is a regularized quasi-ultrametric
in the remainder of this monograph.
Then, by this assumption
and an argument similar to that used in the proof of Lemma~\ref{1955},
we find that, for any $\theta\in(0,\infty)$
and $f\in(\mathring{{\mathcal{G}}}_0^\varepsilon(\beta,\gamma))'$,
$\dot{\mathcal{S}}^{s,(1)}_{q,\theta}(f)$ is $\mu$-measurable.
\end{enumerate}
\end{remark}

We have the following relationship between $\dot{\mathcal{S}}^s_q(f)$
and $\dot{\mathcal{S}}^{s,(1)}_{q,\theta}(f)$ with $\theta\in[1,\infty)$.

\begin{proposition}\label{Prop5.5}
Let
$(X,\mathbf{q},\mu)$
be a quasi-ultrametric space of homogeneous type
with upper dimension $n\in(0,\infty)$, $q\in(0,\infty)$, $s\in(0,\infty)$,
and $p\in(0,q)$.
Assume that $0<\varepsilon\preceq\mathrm{ind\,}(X,\mathbf{q})$
with $\mathrm{ind\,}(X,\mathbf{q})$ as in \eqref{index} and
$\beta,\gamma\in(0,\varepsilon)$,
Then there exists a positive constant $C$ such that,
for any $\theta\in[1,\infty)$
and $f\in(\mathring{{\mathcal{G}}}_0^\varepsilon(\beta,\gamma))'$,
$$
\left\|\dot{\mathcal{S}}^
{s,(1)}_{q,\theta}(f)\right\|_{L^p(X)}
\leq C\theta^\frac{n}{p}
\left\|\dot{\mathcal{S}}^s_q(f)\right\|_{L^p(X)}.
$$
\end{proposition}

\begin{proof}
Let $\rho\in\mathbf{q}$ be a regularized quasi-ultrametric
as in Definition~\ref{D2.3}.
For any nonnegative function $g$
and any $x\in X$, let $\widetilde{\mathcal{M}_\rho}(g)(x)$
be the same as in \eqref{1355}.
Let $f\in(\mathring{{\mathcal{G}}}_0^\varepsilon(\beta,\gamma))'$.
For any $t\in(0,\infty)$, consider the sets
$$
E_t:=\left\{x\in X:\dot{\mathcal{S}}^{s}_{q}(f)(x)>t\right\}
$$
and
$$
\widetilde{E}_t:=\left\{x\in X:\widetilde{\mathcal{M}_\rho}
\left(\mathbf{1}_{E_t}\right)(x)>\frac{1}{2}\right\}.
$$
Since the desired inequality is trivially true
if $\|\dot{\mathcal{S}}^{s}_{q}(f)\|_{L^p(X)}=\infty$, in what
follows, we may assume that
$\|\dot{\mathcal{S}}^{s}_{q}(f)\|_{L^p(X)}<\infty$.
In particular, for any $t\in(0,\infty)$,
$\mu(E_t)<\infty$ and hence $\mathbf{1}_{E_t}\in L^1(X)$.
We claim that, for any $t\in(0,\infty)$,
\begin{equation}\label{xip-944}
\int_{\widetilde{E}_t^{\complement}}
\left[\dot{\mathcal{S}}^{s,(1)}_{q,\theta}(f)(x)\right]^q\,d\mu(x)
\lesssim\theta^n\int_{E_t^{\complement}}
\left[\dot{\mathcal{S}}^{s}_{q}(f)(x)\right]^q\,d\mu(x).
\end{equation}
Observe that \eqref{xip-944} trivially holds
if $\widetilde{E}_t^{\complement}=\varnothing$ and hence,
in what follows, we may assume that
$\widetilde{E}_t^{\complement}\neq\varnothing$.
In this case, for any $y\in X$, let
$\rho(y):=\inf_{x\in \widetilde{E}_t^{\complement}}\rho(x,y)$.
Then, for any $k\in\mathbb{Z}$
and $x,y\in X$, if
$x\in\widetilde{E}_t^{\complement}\cap B_\rho(y,\theta2^{-k})$,
then $\rho(y)<\theta2^{-k}$.
By this, we obtain
\begin{align*}
&\int_{\widetilde{E}^{\complement}_t}
\left[\dot{\mathcal{S}}^{s,(1)}_{q,\theta}(f)(x)\right]^q\,d\mu(x)\nonumber\\
&\quad=\int_{\widetilde{E}^{\complement}_t}
\sum_{k=-\infty}^\infty2^{ksq}
\int_{B_\rho(x,\theta2^{-k})}
\left|\mathcal{D}_kf(y)\right|^q
\frac{d\mu(y)}{(V_\rho)_{2^{-k}}(y)}\,d\mu(x)
\nonumber\\
&\quad=\sum_{k=-\infty}^\infty2^{ksq}
\int_{X}
\left|\mathcal{D}_kf(y)\right|^q
\frac{\mu(\widetilde{E}^{\complement}_t
\cap B_\rho(y,\theta2^{-k}))}{(V_\rho)_{2^{-k}}(y)}\,d\mu(y)
\nonumber\\
&\qquad\qquad\text{by Tonelli's theorem}
\nonumber\\
&\quad\lesssim\theta^n\sum_{k=-\infty}^\infty2^{ksq}
\int_{X}
\left|\mathcal{D}_kf(y)\right|^q\,d\mu(y)
\quad\text{by \eqref{upperdoub}}.
\end{align*}
From this, an argument similar to that used
in the estimation of \eqref{2208}, and Tonelli's theorem,
it follows that
\begin{align*}
&\int_{\widetilde{E}^{\complement}_t}
\left[\dot{\mathcal{S}}^{s,(1)}_{q,\theta}(f)(x)\right]^q\,d\mu(x)\\
&\quad\lesssim\theta^n\sum_{k=-\infty}^\infty2^{ksq}
\int_{X}\left|\mathcal{D}_kf(y)\right|^q
\frac{\mu(E^{\complement}_t\cap B_\rho(y,2^{-k}))}{(V_\rho)_{2^{-k}}(y)}
d\mu(y)
\\
&\quad=\theta^n\sum_{k=-\infty}^\infty2^{ksq}
\int_{X}\left|\mathcal{D}_kf(y)\right|^q
\int_{B_\rho(y,2^{-k})}\mathbf{1}_{E^{\complement}_t}(x)\,d\mu(x)
\frac{d\mu(y)}{(V_\rho)_{2^{-k}}(y)}
\\
&\quad=\theta^n\int_{E^{\complement}_t}\sum_{k=-\infty}^\infty2^{ksq}
\int_{B_\rho(x,2^{-k})}\left|\mathcal{D}_kf(y)\right|^q
\frac{d\mu(y)}{(V_\rho)_{2^{-k}}(y)}\,d\mu(x)\\
&\quad=\theta^n\int_{E^{\complement}_t}
\left[\dot{\mathcal{S}}^{s}_{q}(f)(x)\right]^q\,d\mu(x),
\end{align*}
which completes the proof of \eqref{xip-944}
and hence the above claim \eqref{xip-944}.

By \eqref{muM}
with $g:=\mathbf{1}_{E_t}$,
we conclude that, for any $t\in(0,\infty)$,
\begin{align*}
&\mu\left(\left\{x\in X:
\dot{\mathcal{S}}^{s,(1)}_{q,\theta}(f)(x)>t\right\}\right)
\\
&\quad\leq\mu\Big(\widetilde{E}_t\Big)+
\mu\left(\left\{x\in \widetilde{E}_t^{\complement}:
\dot{\mathcal{S}}^{s,(1)}_{q,\theta}(f)(x)>t\right\}\right)
\\
&\quad\lesssim\theta^n\mu(E_t)+
t^{-q}\int_{\widetilde{E}_t^{\complement}}
\left[\dot{\mathcal{S}}^{s,(1)}_{q,\theta}
(f)(x)\right]^q\,d\mu(x)
\quad\text{by Chebyshev's inequality}
\\
&\quad\lesssim\theta^n\mu(E_t)+
t^{-q}\theta^n\int_{E_t^{\complement}}
\left[\dot{\mathcal{S}}^{s}_{q}
(f)(x)\right]^q\,d\mu(x)
\quad\text{by \eqref{xip-944}},
\end{align*}
which, together with the definition of $E_t$,
further implies that
\begin{align*}
&\mu\left(\left\{x\in X:
\dot{\mathcal{S}}^{s,(1)}_{q,\theta}(f)(x)>t\right\}\right)
\\
&\quad\lesssim\theta^n\mu(E_t)+
qt^{-q}\theta^n\int_0^t s^{q-1}\mu\left(\left\{x\in X:
\dot{\mathcal{S}}^{s}_{q}(f)(x)
\mathbf{1}_{E_t^{\complement}}(x)>s\right\}\right)\,ds
\\	
&\quad\lesssim\theta^n\mu(E_t)+
t^{-q}\theta^n\int_0^ts^{q-1}\mu(E_s)\,ds.
\end{align*}
From this, Tonelli's theorem,
and $p\in(0,q)$, we infer that
\begin{align*}
\left\|\dot{\mathcal{S}}^{s,(1)}_{q,\theta}
(f)\right\|^p_{L^p(X)}
&=p\int_0^\infty t^{p-1}
\mu\left(\left\{x\in X:
\dot{\mathcal{S}}^{s,(1)}_{q,\theta}
(f)(x)>t\right\}\right)\,dt\\
&\lesssim\theta^n\left[\int_0^\infty
t^{p-1}\mu(E_t)dt+
\int_0^\infty t^{p-q-1}
\int_0^ts^{q-1}\mu(E_s)\,ds\,dt\right]\\
&\sim\theta^n
\left\|\dot{\mathcal{S}}^{s}_{q}(f)\right\|^p_{L^p(X)}
+\theta^n\int_0^\infty s^{q-1}\mu(E_s)
\int_s^\infty t^{p-q-1}\,dt\,ds\\
&\sim\theta^n
\left\|\dot{\mathcal{S}}^{s}_{q}(f)\right\|^p_{L^p(X)}
+\theta^n\int_0^\infty s^{p-1}\mu(E_s)\,ds
\sim\theta^n
\left\|\dot{\mathcal{S}}^{s}_{q}(f)\right\|^p_{L^p(X)},
\end{align*}
which completes the proof of Proposition~\ref{Prop5.5}.
\end{proof}

With Proposition~\ref{Prop5.5} in mind, we are ready to
show the Littlewood--Paley $g_\lambda^*$-function
characterization of $\dot{F}^{p,q}_s(X)$
in the case $p\in(0,q)$.

\begin{proposition}\label{g*lammda}
Let $(X,\mathbf{q},\mu)$,
$n$, $\varepsilon$, $s$, $p$, $q$, $\beta$, and $\gamma$ be
the same as in Definition~\ref{DefT--L}
with $p\in(0,q)$ and $\lambda\in(\frac{qn}{p},\infty)$.
Then $f\in\dot{F}^{p,q}_s(X)$
if and only if $f\in(\mathring{{\mathcal{G}}}_0^\varepsilon(\beta,\gamma))'$
and $(\dot{g}^*_\lambda)^s_q(f)\in L^p(X)$.
Moreover, there exists a constant $C\in[1,\infty)$ such that,
for any $f\in(\mathring{{\mathcal{G}}}_0^\varepsilon(\beta,\gamma))'$,
\begin{align}
\label{958}
C^{-1}\left\|\dot{\mathcal{S}}^s_q(f)\right\|_{L^p(X)}
\leq\left\|(\dot{g}^*_\lambda)^s_q(f)\right\|_{L^p(X)}
\leq C\left\|\dot{\mathcal{S}}^s_q(f)\right\|_{L^p(X)}.
\end{align}
\end{proposition}

\begin{proof}
Let $\rho\in\mathbf{q}$ be a regularized quasi-ultrametric
as in Definition~\ref{D2.3}.
Note that all of the conclusions of this proposition
will follow once we establish \eqref{958}.
To this end, first note that,
by the same argument as that used in the estimation of \eqref{ddq-339}
(which does not rely on the relation between $p$ and $q$),
we immediately obtain the first inequality of \eqref{958}.

Next, we turn to prove the second inequality of \eqref{958}.
For any $f\in(\mathring{{\mathcal{G}}}_0^\varepsilon(\beta,\gamma))'$
and $x\in X$, we have [keeping in mind the definition
of $\dot{\mathcal{S}}^{s,(1)}_{q,\theta}(f)$
in \eqref{auxlusin}]
\begin{align*}
\left[(\dot{g}^*_\lambda)^s_q(f)(x)\right]^q
&\leq\sum_{k\in\mathbb{Z}}2^{ksq}
\int_{B_\rho(x,2^{-k})}\left|\mathcal{D}_kf(y)\right|^q
\frac{d\mu(y)}{(V_\rho)_{2^{-k}}(y)}
\\
&\quad+\sum_{k\in\mathbb{Z}}2^{ksq}
\int_{[B_\rho(x,2^{-k})]^\complement}
\left|\mathcal{D}_kf(y)\right|^q
\left[\frac{2^{-k}}{2^{-k}+\rho(x,y)}\right]^\lambda
\frac{d\mu(y)}{(V_\rho)_{2^{-k}}(y)}
\\
&\le\left[\dot{\mathcal{S}}^{s,(1)}_{q,1}(f)(x)\right]^q
+\sum_{j=1}^{\infty}\left(\frac{1}{1+2^{j-1}}\right)^\lambda
\\
&\quad\times\sum_{k=-\infty}^{\infty}
2^{ksq}
\int_{2^{j-1}2^{-k}\leq\rho(x,y)<2^j2^{-k}}
\left|\mathcal{D}_kf(y)\right|^q
\frac{d\mu(y)}{(V_\rho)_{2^{-k}}(y)},
\end{align*}
which implies that
\begin{align}\label{1421}
\left[(\dot{g}^*_\lambda)^s_q(f)(x)\right]^q
&\lesssim\left[\dot{\mathcal{S}}^{s,(1)}_{q,1}(f)(x)\right]^q\nonumber\\
&\quad+\sum_{j=1}^{\infty}2^{-j\lambda}
\sum_{k=-\infty}^{\infty}
2^{ksq}
\int_{B_\rho(x,2^{j-k})}
\left|\mathcal{D}_kf(y)\right|^q
\frac{d\mu(y)}{(V_\rho)_{2^{-k}}(y)}\nonumber\\
&\lesssim\sum_{j=0}^{\infty}2^{-j\lambda}
\left[\dot{\mathcal{S}}^{s,(1)}_{q,2^j}(f)(x)\right]^q.
\end{align}
From this and $\lambda\in(\frac{qn}{p},\infty)$,
we deduce that
\begin{align*}
\left\|(\dot{g}^*_\lambda)^s_q(f)\right\|_{L^p(X)}^p
&=\left\|\left[(\dot{g}^*_\lambda)^s_q(f)\right]^q
\right\|_{L^\frac{p}{q}(X)}^\frac{p}{q}\\
&\lesssim\int_X\left\{\sum_{j=0}^{\infty}2^{-j\lambda}
\left[\dot{\mathcal{S}}^{s,(1)}_{q,2^j}
(f)(x)\right]^q\right\}^\frac{p}{q}\,d\mu(x)
\quad\text{by \eqref{1421}}
\\
&\leq\int_X\sum_{j=0}^{\infty}2^{-\frac{j\lambda p}{q}}
\left[\dot{\mathcal{S}}^{s,(1)}_{q,2^j}(f)(x)\right]^p\,d\mu(x)
\quad\text{by Lemma~\ref{discreteinequality}}\\
&=\sum_{j=0}^{\infty}2^{-\frac{j\lambda p}{q}}\int_X
\left[\dot{\mathcal{S}}^{s,(1)}_{q,2^j}(f)(x)\right]^p\,d\mu(x)
\quad\text{by Tonelli's theorem}\\
&\lesssim\sum_{j=0}^{\infty}2^{-j(\frac{\lambda p}{q}-n)}
\left\|\dot{\mathcal{S}}^{s}_{q}(f)\right\|^p_{L^p(X)}
\sim\left\|\dot{\mathcal{S}}^s_q(f)\right\|_{L^p(X)}^p
\quad\text{by Proposition~\ref{Prop5.5}}.
\end{align*}
This finishes the proof of the second inequality of \eqref{958}
and hence Proposition~\ref{g*lammda}.
\end{proof}

Combining Proposition~\ref{g*lammda},
Theorem~\ref{Lusin},
and Proposition~\ref{g*1},
we obtain the following Littlewood--Paley $g^*_\lambda$-function
characterization of $\dot{F}^{p,q}_s(X)$.

\begin{theorem}\label{T--Lg*}
Let $(X,\mathbf{q},\mu)$, $n$, $\varepsilon$,
$s$, $p$, $q$, $\beta$, and $\gamma$ be
the same as in Theorem~\ref{Lusin}.
Assume that
$\lambda\in(\max\{n,\,\frac{qn}{p}\},\infty)$.
Then $f\in\dot{F}^{p,q}_s(X)$
if and only if
$f\in(\mathring{{\mathcal{G}}}_0^\varepsilon(\beta,\gamma))'$
and $(\dot{g}^*_\lambda)^s_q(f)\in L^p(X)$.
Moreover, there exists a constant $C\in[1,\infty)$ such that,
for any $f\in(\mathring{{\mathcal{G}}}_0^\varepsilon(\beta,\gamma))'$,
\begin{align*}
C^{-1}\left\|(\dot{g}^*_\lambda)^s_q(f)\right\|_{L^p(X)}
\leq\|f\|_{\dot{F}^{p,q}_s(X)}
\leq C\left\|(\dot{g}^*_\lambda)^s_q(f)\right\|_{L^p(X)}.
\end{align*}
\end{theorem}

\begin{remark}
\begin{enumerate}
\item
$\lambda$ in \eqref{g*} equals
$q\lambda$ with $\lambda$ as in the Littlewood--Paley
$g_\lambda^*$-function in \cite[p.\,102, Definition 1.1]{p1982} and
\cite[p.\,9, (62)]{h2013}.
In 1982, P\"aiv\"arinta \cite{p1982} showed the Littlewood--Paley
$g_\lambda^*$-function characterization of the Triebel--Lizorkin space
$F^{p,q}_s(\mathbb{R}^n)$ whenever $\lambda\in(n+\frac{nq}{\min\{p,\,q\}},\infty)$.
It is easy to see that the range of $\lambda$ in
Theorem~\ref{T--Lg*} improves upon
\cite[Remark 2.6]{p1982} by regarding $\mathbb{R}^n$
as a standard $n$-Ahlfors-regular quasi-ultrametric space.
Moreover, in 2013, Hu \cite{h2013} established the Littlewood--Paley
$g_\lambda^*$-function characterization of the Triebel--Lizorkin space
on a stratified Lie group whenever $\lambda\in(\max\{\frac{nq}{p},\,n\},\infty)$
and it is easy to see that Theorem~\ref{T--Lg*}
coincides with \cite[Proposition 14]{h2013}
when $X$ is a stratified Lie group which is also a standard
$n$-Ahlfors-regular quasi-ultrametric space.
\item
Let $(X,\rho,\mu)$ be a space of homogeneous type
with upper dimension $n\in(0,\infty)$ and
$\lambda\in(\max\{n,\,\frac{qn}{p}\},\infty)$.
He et al. \cite[Theorem 6.6]{HWYY} proved that
$f\in\dot{F}^{p,q}_s(X)$
if and only if
$f\in(\mathring{{\mathcal{G}}}_0^\varepsilon(\beta,\gamma))'$
and $(\dot{g}^*_\lambda)^s_q(f)\in L^p(X)$ where $|s|<\varepsilon<1$.
Moreover, there exists a constant $C\in[1,\infty)$,
independent of $f$, such that
$$
C^{-1}\left\|(\dot{g}^*_\lambda)^s_q(f)\right\|_{L^p(X)}
\leq\|f\|_{\dot{F}^{p,q}_s(X)}
\leq C\left\|(\dot{g}^*_\lambda)^s_q(f)\right\|_{L^p(X)},
$$
where $\dot{F}^{p,q}_s(X)$ is the homogeneous
Triebel--Lizorkin space introduced in \cite[Definition 3.1(ii)]{WHHY}.
In Theorem~\ref{T--Lg*}, the Littlewood--Paley
$g_\lambda^*$-function characterization of $\dot{F}^{p,q}_s(X)$ holds
for the larger range of $|s|<\varepsilon\preceq\mathrm{ind\,}(X,\mathbf{q})$
with $\mathrm{ind\,}(X,\mathbf{q})$ as in \eqref{index},
but in the more specialized setting of
ultra-RD-spaces.

\item
Note that, by \cite[Proposition 9.5]{AlMi2015},
we find that, for any $p\in(\frac{n}{n+\mathrm{ind\,}(X,\mathbf{q})},1]$,
$\dot{F}^{p,2}_0(X)=H^p_{\mathrm{at}}(X)$ (the atomic Hardy space)
with $\mathrm{ind\,}(X,\mathbf{q})$ as in \eqref{index}.
In this special case, Theorem~\ref{T--Lg*} coincides with
Theorem~\ref{Glamda*} in the case where $(X,\mathbf{q},\mu)$
is standard.
Thus, by Remark~\ref{1629},
the range of $\lambda$ in Theorem~\ref{T--Lg*}
is also the known best possible.
\item
Let $(X,\mathbf{q},\mu)$ be an Ahlfors-regular
quasi-ultrametric space as in \cite{AlMi2015}.
In this setting, Theorems~\ref{Lusin} and~\ref{T--Lg*} give the
Lusin area function characterization and the
Littlewood--Paley $g_\lambda^*$-function characterization
of the homogeneous Triebel--Lizorkin space
$\dot{F}^{p,q}_s(X)$ in \cite[Definition~9.2]{AlMi2015}.
\end{enumerate}
\end{remark}

\section{Littlewood--Paley Characterization of $\dot{F}^{\infty,q}_s(X)$}
\label{1700}

In this section, we aim to establish the
Littlewood--Paley characterization of
homogeneous Triebel--Lizorkin spaces $\dot{F}^{\infty,q}_s(X)$
on the ultra-RD-space $(X,\mathbf{q},\mu)$,
where $\mu$ is assumed to be a Borel-semiregular measure on $X$.
We always assume that, for any $\rho\in\mathbf{q}$,
$\mathrm{diam}_\rho(X)=\infty$
and hence $\mu(X)=\infty$.

The main theorem of this section is as follows.

\begin{theorem}
\label{2158}
Let $(X,\mathbf{q},\mu)$ be an ultra-RD-space with upper
dimension $n\in(0,\infty)$, where $\mu$ is assumed to
be a Borel-semiregular measure on $X$.
Let $\varepsilon$, $s$, $q$, $\beta$, and $\gamma$ be
the same as in Definition~\ref{DefT--L}
with $p=\infty$.
Let $\rho\in\mathbf{q}$ be such that all $\rho$-balls
are $\mu$-measurable.
Assume that $\{\mathcal{S}_{2^{-k}}\}_{k\in\mathbb{Z}}$ is
a homogeneous approximation of the identity
of order $\varepsilon$
as in Definition~\ref{DefHATI} and
$\{Q^{l}_\tau\}_{l\in\mathbb{Z},\,\tau\in I_l}$ the collection of all
dyadic cubes given in Lemma~\ref{DyadicSys}.
For any $k\in\mathbb{Z}$, let
$\mathcal{D}_k:=\mathcal{S}_{2^{-k}}-\mathcal{S}_{2^{-k+1}}$.
Then the following assertions hold.
\begin{enumerate}
\item[\textup{(i)}]
$f\in\dot{F}^{\infty,q}_s(X)$
if and only if
$f\in(\mathring{{\mathcal{G}}}_0^\varepsilon(\beta,\gamma))'$
and
\begin{align*}
\|f\|_{\widetilde{\dot{F}^{\infty,q}_s}(X)}
:&=\sup_{l\in\mathbb{Z}}\sup_{\tau\in I_l}
\left[\frac{1}{\mu(Q_\tau^l)}\int_{Q^l_\tau}\sum_{k=l}^{\infty}
2^{ksq}\right.\\
&\quad\left.\times\int_{B_\rho(x,2^{-k})}
\left|\mathcal{D}_kf(y)\right|^q
\frac{d\mu(y)\,d\mu(x)}{(V_\rho)_{2^{-k}}(x)}
\right]^\frac{1}{q}
\end{align*}
is finite.

\item[\textup{(ii)}]
Let	 $\lambda\in(\max\{2n,2n-sq,n+nq+|s|q\},\infty)$.
Then $f\in\dot{F}^{\infty,q}_s(X)$
if and only if
$f\in(\mathring{{\mathcal{G}}}_0^\varepsilon(\beta,\gamma))'$
and
\begin{align*}
\|f\|_{\overline{\dot{F}^{\infty,q}_s}(X)}
:&=\sup_{l\in\mathbb{Z}}\sup_{\tau\in I_l}
\Bigg\{\frac{1}{\mu(Q_\tau^l)}\int_{Q^l_\tau}\sum_{k=l}^{\infty}
2^{ksq}\int_{X}\left|\mathcal{D}_kf(y)\right|^q
\\
&\quad\times\left[\frac{2^{-k}}{2^{-k}+\rho(x,y)}\right]^\lambda
\frac{d\mu(y)\,d\mu(x)}{(V_\rho)_{2^{-k}}(x)+(V_\rho)_{2^{-k}}(y)}\,
\Bigg\}^\frac{1}{q}
\end{align*}
is finite.
\end{enumerate}
Moreover, for any $f\in(\mathring{{\mathcal{G}}}_0^\varepsilon(\beta,\gamma))'$,
\begin{align*}
\|f\|_{\widetilde{\dot{F}^{\infty,q}_s}(X)}
\sim\|f\|_{\dot{F}^{\infty,q}_s(X)}
\sim\|f\|_{\overline{\dot{F}^{\infty,q}_s}(X)},
\end{align*}
where the positive equivalence constants are independent of $f$.
\end{theorem}

\begin{remark}
Let the notation be as in Theorem~\ref{2158}.
\begin{enumerate}
\item[\rm(i)]
By the proof of Lemma~\ref{1955}, we conclude that, for any
$k\in\mathbb{Z}$ and
$f\in(\mathring{{\mathcal{G}}}_0^\varepsilon(\beta,\gamma))'$,
$\mathcal{D}_kf$ is a continuous function.
Moreover, the proof of Lemma~\ref{1955} implies that,
for any
$k\in\mathbb{Z}$ and $r\in(0,\infty)$, $(V_\rho)_{r}(\cdot)$
is $\mu$-measurable,
where we need that $\rho\in\mathbf{q}$ is a regularized
quasi-ultrametric as in Definition~\ref{D2.3}.
As a conclusion,
the functions being integrated
in both (i) and (ii) of Theorem~\ref{2158}
are $\mu$-measurable.

\item[\rm(ii)]
$\|\cdot\|_{\widetilde{\dot{F}^{\infty,q}_s}(X)}$
defined via different $\rho\in\mathbf{q}$
(where all $\rho$-balls are $\mu$-measurable)
are equivalent.
Similar equivalence for $\|\cdot\|_{\overline{\dot{F}^{\infty,q}_s}(X)}$ also holds.
Thus, without loss of generality, we may assume that $\rho\in\mathbf{q}$
in Theorem~\ref{2158} is a regularized quasi-ultrametric
in the remainder of this monograph.

\item[\rm(iii)]
The Littlewood--Paley $g_\lambda^*$-function characterization
of $\dot{F}^{\infty,q}_s(X)$ when $X=\mathbb{R}^n$
was obtained in \cite[Theorem~3.2]{lsuyy2012}
and when $(X,\rho,\mu)$ is a space of homogeneous type
was obtained in \cite[Theorem~6.1]{wyyTLinfty}.
We point out that, in the setting of ultra-RD-spaces,
Theorem~\ref{2158} gives the Littlewood--Paley
$g_\lambda^*$-function characterization
of $\dot{F}^{\infty,q}_s(X)$ for wider ranges of $s$ and $q$
than those of \cite[Theorem~6.1]{wyyTLinfty};
see Remark~\ref{2101}(i).

\item[\rm(iv)]
Let $(X,\mathbf{q},\mu)$ be an Ahlfors-regular
quasi-ultrametric space as in \cite{AlMi2015}.
In this setting, Theorem~\ref{2158} gives the
Lusin area function characterization and the
Littlewood--Paley $g_\lambda^*$-function characterization
of $\dot{F}^{\infty,q}_s(X)$
in \cite[Definition~9.2]{AlMi2015}.
\end{enumerate}
\end{remark}

To prove Theorem~\ref{2158}, we need the following
homogeneous Plancherel--P\'olya
inequality, which is the endpoint case $p = \infty$ of
Lemma~\ref{P-P }.

\begin{proposition}
\label{ppineq-infty}
Let $(X,\mathbf{q},\mu)$,
$n$, $\varepsilon$, $s$, $q$, $\beta$, and $\gamma$ be
the same as in Theorem~\ref{2158}.
Assume that $\{\mathcal{S}_{2^{-k}}\}_{k\in\mathbb{Z}}$
and $\{\mathcal{Q}_{2^{-k}}\}_{k\in\mathbb{Z}}$
are two homogeneous approximations of the identity
of order $\varepsilon$
in Definition~\ref{DefHATI} and let
$\mathcal{D}_k:=\mathcal{S}_{2^{-k}}-\mathcal{S}_{2^{-k+1}}$
and $\mathcal{P}_k:=\mathcal{Q}_{2^{-k}}-\mathcal{Q}_{2^{-k+1}}$
for any $k\in\mathbb{Z}$.
Let
$\{Q^{l}_\alpha: l\in\mathbb{Z},\ \alpha\in I_l\}$ be the collection of all
dyadic cubes given in Lemma~\ref{DyadicSys}.
Then, for any $c_0\in[1,\infty)$ and
$f\in(\mathring{{\mathcal{G}}}_0^\varepsilon(\beta,\gamma))'$,
\begin{align}\label{1613}
&\sup_{l\in\mathbb{Z}}\sup_{\alpha\in I_l}
\left[\frac{1}{\mu(Q_\alpha^l)}
\sum_{k=l}^{\infty}\sum_{\tau\in I_k}\sum_{v=1}^{N(k,\tau)}
2^{ksq}\mu(Q_\tau^{k,v})\nonumber
\mathbf{1}_{\{(\tau,v):Q_\tau^{k,v}\subset Q_\alpha^l\}}(\tau,v)
\sup_{x\in B_\rho(z^{k,v}_\tau,c_02^{-k})}
\left|\mathcal{D}_kf(x)\right|^q\right]^\frac{1}{q}
\nonumber\\
&\quad\lesssim
\sup_{l\in\mathbb{Z}}\sup_{\alpha\in I_l}
\left[\frac{1}{\mu(Q_\alpha^l)}
\sum_{k=l}^{\infty}\sum_{\tau\in I_k}\sum_{v=1}^{N(k,\tau)}
2^{ksq}\mu(Q_\tau^{k,v})
\mathbf{1}_{\{(\tau,v):Q_\tau^{k,v}\subset Q_\alpha^l\}}(\tau,v)
\inf_{x\in Q_\tau^{k,v}}
\left|\mathcal{P}_kf(x)\right|^q\right]^\frac{1}{q}
\end{align}
with the usual modification when $q=\infty$,
where the implicit positive constant is independent of $f$ and where
$\{Q^{k,v}_\tau: k\in\mathbb{Z},\ \tau\in I_k,\ v=1,\dots,N(k,\tau)\}$
is the collection of dyadic cubes as in Remark~\ref{jcubes}.
\end{proposition}

To show Proposition~\ref{ppineq-infty},
we need the following technical lemma.

\begin{lemma}
\label{444}
Let $(X,\rho,\mu)$ be a quasi-ultrametric space
of homogeneous type with upper dimension $n\in(0,\infty)$.
Assume that $\gamma\in(0,\infty)$ and $p\in(\frac{n}{n+\gamma},1]$.
Then there exists a positive constant $C$ such that, for any
$k,k'\in\mathbb{Z}$, $x\in X$, and $y_\tau^{k,v}\in Q_\tau^{k,v}$
with both $\tau\in I_k$ and $v\in\{1,\ldots,N(k,\tau)\}$,
\begin{align*}
&\sum_{\tau\in I_k}\sum_{v=1}^{N(k,\tau)}\mu(Q_\tau^{k,v})
\left\{\frac{1}{(V_\rho)_{2^{-(k\land k')}}(x)+V_\rho(x,y_\tau^{k,v})}
\left[\frac{2^{-(k\land k')}}{2^{-(k\land k')}+
\rho(x,y_\tau^{k,v})}\right]^\gamma\right\}^p\\
&\quad\sim\left[(V_\rho)_{2^{-(k\land k')}}(x)\right]^{1-p},
\end{align*}
where the positive equivalence constants are independent of
$k$, $k'$, $x$, and $y_\tau^{k,v}$.
\end{lemma}

\begin{proof}
Let $k,k'\in\mathbb{Z}$ and $x\in X$.
For any $\tau\in I_k$, $v\in\{1,\ldots,N(k,\tau)\}$,
and $z\in Q_\tau^{k,v}$,
Lemma~\ref{DyadicSys}(iv) implies that
\begin{align*}
\rho(z,y_\tau^{k,v})\leq\mathrm{diam}_\rho(Q_\tau^{k,v})\lesssim2^{-k}
\leq2^{-(k\land k')},
\end{align*}
which, together with Definition~\ref{D2.1}(iii), further implies that
\begin{align*}
2^{-(k\land k')}+\rho(x,y_\tau^{k,v})\sim
2^{-(k\land k')}+\rho(x,z).
\end{align*}
From this and both (vii) and (vi)
of Lemma~\ref{A3.2}, we infer that,
for any $z\in Q_\tau^{k,v}$,
\begin{align*}
(V_\rho)_{2^{-(k\land k')}}(x)+V_\rho(x,y_\tau^{k,v})
\sim(V_\rho)_{2^{-(k\land k')}}(x)+V_\rho(x,z).
\end{align*}
By this,
we find that
\begin{align}\label{442}
&\sum_{\tau\in I_k}\sum_{v=1}^{N(k,\tau)}\mu(Q_\tau^{k,v})
\left\{\frac{1}{(V_\rho)_{2^{-(k\land k')}}(x)+V_\rho(x,y_\tau^{k,v})}
\left[\frac{2^{-(k\land k')}}{2^{-(k\land k')}+\rho(x,y_\tau^{k,v})}
\right]^\gamma\right\}^p\nonumber\\
&\quad=\sum_{\tau\in I_k}\sum_{v=1}^{N(k,\tau)}
\int_{Q_\tau^{k,v}}
\left\{\frac{1}{(V_\rho)_{2^{-(k\land k')}}(x)+V_\rho(x,y_\tau^{k,v})}
\left[\frac{2^{-(k\land k')}}{2^{-(k\land k')}+\rho(x,y_\tau^{k,v})}
\right]^\gamma\right\}^p\,d\mu(z)\nonumber\\
&\quad\sim\sum_{\tau\in I_k}\sum_{v=1}^{N(k,\tau)}
\int_{Q_\tau^{k,v}}
\left\{\frac{1}{(V_\rho)_{2^{-(k\land k')}}(x)+V_\rho(x,z)}
\left[\frac{2^{-(k\land k')}}{2^{-(k\land k')}+\rho(x,z)}
\right]^\gamma\right\}^p\,d\mu(z)\nonumber\\
&\quad=\int_{B_\rho(x,2^{-(k\land k')})}
\left\{\frac{1}{(V_\rho)_{2^{-(k\land k')}}(x)+V_\rho(x,z)}
\left[\frac{2^{-(k\land k')}}{2^{-(k\land k')}+\rho(x,z)}
\right]^\gamma\right\}^p\,d\mu(z)\nonumber\\
&\qquad+\sum_{l=1}^\infty
\int_{B_\rho(x,2^l2^{-(k\land k')})\setminus
B_\rho(x,2^{l-1}2^{-(k\land k')})}\cdots\nonumber\\
&\quad=:\mathrm{I}_0+\sum_{l=1}^\infty\mathrm{I}_l.
\end{align}
To estimate $\mathrm{I}_0$, we easily conclude that
\begin{align}\label{443}
\mathrm{I}_0\sim\int_{B_\rho(x,2^{-(k\land k')})}
\left[\frac{1}{(V_\rho)_{2^{-(k\land k')}}(x)}\right]^p
\,d\mu(z)=\left[(V_\rho)_{2^{-(k\land k')}}(x)\right]^{1-p}.
\end{align}
To estimate $\mathrm{I}_l$ with $l\in\mathbb{N}$,
from \eqref{upperdoub} and $p\in(\frac{n}{n+\gamma},1]$,
we deduce that
\begin{align*}
\mathrm{I}_l&\sim\int_{B_\rho(x,2^l2^{-(k\land k')})\setminus
B_\rho(x,2^{l-1}2^{-(k\land k')})}\left[\frac{1}{V(x,z)}\right]^p
2^{-lp\gamma}\,d\mu(z)\\
&\lesssim2^{-lp\gamma}(V_\rho)_{2^{l-(k\land k')}}(x)
\left[(V_\rho)_{2^{l-1-(k\land k')}}(x)\right]^p\\
&\lesssim2^{-lp\gamma}\left[(V_\rho)_{2^{l-1-(k\land k')}}(x)\right]^{1-p}
\lesssim2^{l(n-pn-p\gamma)}\left[(V_\rho)_{2^{-(k\land k')}}(x)\right]^{1-p},
\end{align*}
which, combined with \eqref{442}, \eqref{443}, and $p\in(\frac{n}{n+\gamma},1]$,
further implies that
\begin{align*}
&\sum_{\tau\in I_k}\sum_{v=1}^{N(k,\tau)}\mu(Q_\tau^{k,v})
\left\{\frac{1}{(V_\rho)_{2^{-(k\land k')}}(x)+V_\rho(x,y_\tau^{k,v})}
\left[\frac{2^{-(k\land k')}}{2^{-(k\land k')}+
\rho(x,y_\tau^{k,v})}\right]^\gamma\right\}^p\\
&\quad\lesssim\left[1+\sum_{l=1}^\infty2^{l(n-pn-p\gamma)}\right]
\left[(V_\rho)_{2^{-(k\land k')}}(x)\right]^{1-p}
\sim\left[(V_\rho)_{2^{-(k\land k')}}(x)\right]^{1-p}\sim\mathrm{I}_0\\
&\quad\lesssim\sum_{\tau\in I_k}\sum_{v=1}^{N(k,\tau)}\mu(Q_\tau^{k,v})
\left\{\frac{1}{(V_\rho)_{2^{-(k\land k')}}(x)+V_\rho(x,y_\tau^{k,v})}
\left[\frac{2^{-(k\land k')}}{2^{-(k\land k')}+
\rho(x,y_\tau^{k,v})}\right]^\gamma\right\}^p.
\end{align*}
This finishes the proof of Lemma~\ref{444}.
\end{proof}

Now, we prove Proposition~\ref{ppineq-infty}.

\begin{proof}[Proof of Proposition~\ref{ppineq-infty}]
We only show the present proposition in the case $q<\infty$
because the proof of the case $q=\infty$
is similar and we omit the details.
By arguments similar to those used in the estimations
of both \eqref{PP1} and \eqref{PP2},
we find that, for any fixed $\varepsilon'\in(0,\varepsilon)$ and for any
$l\in\mathbb{Z}$, $\alpha\in I_l$, and
$f\in(\mathring{{\mathcal{G}}}_0^\varepsilon(\beta,\gamma))'$,
\begin{align}\label{2115}
&\frac{1}{\mu(Q_\alpha^l)}
\sum_{k=l}^{\infty}\sum_{\tau\in I_k}\sum_{v=1}^{N(k,\tau)}
2^{ksq}\mu(Q_\tau^{k,v})\mathbf{1}_{\{(\tau,v):Q_\tau^{k,v}\subset Q_\alpha^l\}}(\tau,v)
\sup_{x\in B_\rho(z^{k,v}_\tau,c_02^{-k})}
\left|\mathcal{D}_kf(x)\right|^q
\nonumber\\
&\quad=\frac{1}{\mu(Q_\alpha^l)}
\sum_{k=l}^{\infty}\sum_{\tau\in I_k}\sum_{v=1}^{N(k,\tau)}
2^{ksq}\mu(Q_\tau^{k,v})\mathbf{1}_{\{(\tau,v):Q_\tau^{k,v}\subset Q_\alpha^l\}}(\tau,v)
\nonumber\\
&\qquad\times
\sup_{x\in B_\rho(z^{k,v}_\tau,c_02^{-k})}
\left|\sum_{k'\in\mathbb{Z}}
\sum_{\tau'\in I_k}\sum_{v'=1}^{N(k',\tau')}
\mu(Q_{\tau'}^{k',v'})
D_k\widetilde{P}_{k'}(x,y_{\tau'}^{k',v'})
\mathcal{P}_{k'}(f)(y_{\tau'}^{k',v'})\right|^q
\nonumber\\
&\quad\lesssim\mathrm{I}+\mathrm{II},
\end{align}
where
\begin{align*}
\mathrm{I}:&=\frac{1}{\mu(Q_\alpha^l)}
\sum_{k=l}^{\infty}\sum_{\tau\in I_k}\sum_{v=1}^{N(k,\tau)}
2^{ksq}\mu(Q_\tau^{k,v})\mathbf{1}_{\{(\tau,v):Q_\tau^{k,v}\subset Q_\alpha^l\}}(\tau,v)
\nonumber\\
&\quad\times
\left\{\sum_{k'=l}^\infty
\sum_{\tau'\in I_{k'}}\sum_{v'=1}^{N(k',\tau')}2^{-|k-k'|\varepsilon'}
\mu(Q_{\tau'}^{k',v'})
\left|\mathcal{P}_{k'}(f)(y_{\tau'}^{k',v'})\right|
\right.\nonumber\\
&\quad\left.\times\frac{1}{(V_\rho)_{2^{-(k\land k')}}
(y_\tau^{k,v})+V_\rho(y_\tau^{k,v},y^{k',v'}_{\tau'})}
\left[\frac{2^{-(k\land k')}}
{2^{-(k\land k')}+\rho(y_\tau^{k,v},y^{k',v'}_{\tau'})}\right]^{\varepsilon'}\right\}^q
\end{align*}
and
\begin{align*}
\mathrm{II}:&=\frac{1}{\mu(Q_\alpha^l)}
\sum_{k=l}^{\infty}\sum_{\tau\in I_k}\sum_{v=1}^{N(k,\tau)}
2^{ksq}\mu(Q_\tau^{k,v})\mathbf{1}_{\{(\tau,v):Q_\tau^{k,v}\subset Q_\alpha^l\}}(\tau,v)
\nonumber\\
&\quad\times\left\{\sum_{k'=-\infty}^{l-1}
\sum_{\tau'\in I_{k'}}\sum_{v'=1}^{N(k',\tau')}2^{-|k-k'|\varepsilon'}
\mu(Q_{\tau'}^{k',v'})
\left|\mathcal{P}_{k'}(f)(y_{\tau'}^{k',v'})\right|
\right.\nonumber\\
&\quad\left.\times\frac{1}{(V_\rho)_{2^{-(k\land k')}}
(y_\tau^{k,v})+V_\rho(y_\tau^{k,v},y^{k',v'}_{\tau'})}
\left[\frac{2^{-(k\land k')}}
{2^{-(k\land k')}+\rho(y_\tau^{k,v},y^{k',v'}_{\tau'})}\right]^{\varepsilon'}\right\}^q.
\end{align*}
We first estimate $\mathrm{I}$.
To this end, for any $l\in\mathbb{Z}$ and $\alpha\in I_l$,
let
\begin{align*}
\mathcal{B}_\alpha^{l,0}:=\left\{\widetilde{\alpha}\in I_l:
\rho(z_{\widetilde{\alpha}}^l,z_\alpha^l)<A2^{-l}\right\},
\end{align*}
where $A\in(0,\infty)$ is a constant specified later.
From this and Definition~\ref{D2.1}(iii), we infer that,
for any $\widetilde{\alpha}\in I_l$,
$$
B_\rho(z_{\widetilde{\alpha}}^l,C_l2^{-l})
\subset B_\rho(z_\alpha^l,C_\rho(A+C_r)2^{-l}),
$$
which, together with both (i) and (v) of Lemma~\ref{DyadicSys}
and with Remark~\ref{940}(iii),
further implies that $\sup_{l\in\mathbb{Z},\,\alpha\in I_l}\#\mathcal{B}_\alpha^{l,0}<\infty$.
Moreover, by Lemmas~\ref{DyadicSys} and~\ref{A3.2}(vii),
we conclude that, for any $\widetilde{\alpha}\in\mathcal{B}_\alpha^{l,0}$,
\begin{align}\label{9252}
\mu(Q_{\widetilde{\alpha}}^l)\sim\mu(Q_\alpha^l).
\end{align}
From this, we deduce that
\begin{align}\label{2116}
\mathrm{I}\lesssim\mathrm{I_{1,1}}+\mathrm{I_{1,2}},
\end{align}
where
\begin{align*}
\mathrm{I_{1,1}}:&=\frac{1}{\mu(Q_\alpha^l)}
\sum_{k=l}^{\infty}\sum_{\tau\in I_k}\sum_{v=1}^{N(k,\tau)}
2^{ksq}\mu(Q_\tau^{k,v})\mathbf{1}_{\{(\tau,v):Q_\tau^{k,v}\subset Q_\alpha^l\}}(\tau,v)
\nonumber\\
&\quad\times
\left\{\sum_{k'=l}^\infty
\sum_{\tau'\in I_{k'}}\sum_{v'=1}^{N(k',\tau')}2^{-|k-k'|\varepsilon'}
\mu(Q_{\tau'}^{k',v'})
\left|\mathcal{P}_{k'}(f)(y_{\tau'}^{k',v'})\right|
\right.\nonumber\\
&\quad\times\frac{1}{(V_\rho)_{2^{-(k\land k')}}(y_\tau^{k,v})+
V_\rho(y_\tau^{k,v},y^{k',v'}_{\tau'})}
\left[\frac{2^{-(k\land k')}}
{2^{-(k\land k')}+\rho(y_\tau^{k,v},y^{k',v'}_{\tau'})}\right]^{\varepsilon'}
\nonumber\\
&\quad\times\left.\mathbf{1}_{\{(\tau',v'):
Q_{\tau'}^{k',v'}\subset\bigcup_{\widetilde{\alpha}\in
\mathcal{B}_\alpha^{l,0}}Q_{\widetilde{\alpha}}^l\}}(\tau',v')
\right\}^q
\end{align*}
and
\begin{align*}
\mathrm{I_{1,2}}:&=\frac{1}{\mu(Q_\alpha^l)}
\sum_{k=l}^{\infty}\sum_{\tau\in I_k}\sum_{v=1}^{N(k,\tau)}
2^{ksq}\mu(Q_\tau^{k,v})\mathbf{1}_{\{(\tau,v):Q_\tau^{k,v}\subset Q_\alpha^l\}}(\tau,v)
\nonumber\\
&\quad\times
\left\{\sum_{k'=l}^\infty
\sum_{\tau'\in I_{k'}}\sum_{v'=1}^{N(k',\tau')}2^{-|k-k'|\varepsilon'}
\mu(Q_{\tau'}^{k',v'})
\left|\mathcal{P}_{k'}(f)(y_{\tau'}^{k',v'})\right|
\right.\nonumber\\
&\quad\times\frac{1}{(V_\rho)_{2^{-(k\land k')}}(y_\tau^{k,v})+
V_\rho(y_\tau^{k,v},y^{k',v'}_{\tau'})}
\left[\frac{2^{-(k\land k')}}
{2^{-(k\land k')}+\rho(y_\tau^{k,v},y^{k',v'}_{\tau'})}\right]^{\varepsilon'}
\nonumber\\
&\quad\times\left.\mathbf{1}_{\{(\tau',v'):
Q_{\tau'}^{k',v'}\cap(\bigcup_{\widetilde{\alpha}\in
\mathcal{B}_\alpha^{l,0}}Q_{\widetilde{\alpha}}^l)=\varnothing\}}(\tau',v')
\right\}^q.
\end{align*}
To deal with $\mathrm{I_{1,1}}$,
we consider the following two cases for $q$.

\emph{Case (1)} $q\in(\max\{\frac{n}{n+\varepsilon},\,
\frac{n}{n+\varepsilon+s}\},1]$. In this case,
by Lemma~\ref{discreteinequality},
we find that
\begin{align*}
\mathrm{I_{1,1}}
&\leq\frac{1}{\mu(Q_\alpha^l)}
\sum_{k=l}^{\infty}\sum_{\tau\in I_k}\sum_{v=1}^{N(k,\tau)}
2^{ksq}\mu(Q_\tau^{k,v})\mathbf{1}_{\{(\tau,v):Q_\tau^{k,v}\subset Q_\alpha^l\}}(\tau,v)
\\
&\quad\times\sum_{\widetilde{\alpha}\in\mathcal{B}_\alpha^{l,0}}
\sum_{k'=l}^\infty
\sum_{\tau'\in I_{k'}}\sum_{v'=1}^{N(k',\tau')}2^{-|k-k'|q\varepsilon'}
\left[\mu(Q_{\tau'}^{k',v'})\right]^q
\left|\mathcal{P}_{k'}(f)(y_{\tau'}^{k',v'})\right|^q
\\
&\quad\times\left\{\frac{1}{(V_\rho)_{2^{-(k\land k')}}
(y_\tau^{k,v})+V_\rho(y_\tau^{k,v},y^{k',v'}_{\tau'})}
\left[\frac{2^{-(k\land k')}}
{2^{-(k\land k')}+\rho(y_\tau^{k,v},y^{k',v'}_{\tau'})}\right]^{\varepsilon'}\right\}^q
\\
&\quad\times\mathbf{1}_{\{(\tau',v'):
Q_{\tau'}^{k',v'}\subset Q_{\widetilde{\alpha}}^l\}}(\tau',v'),
\end{align*}
which, combined with Tonelli's theorem for sequences,
further implies that
\begin{align}\label{1442}
\mathrm{I_{1,1}}
&=\frac{1}{\mu(Q_\alpha^l)}
\sum_{\widetilde{\alpha}\in\mathcal{B}_\alpha^{l,0}}
\sum_{k'=l}^\infty
\sum_{\tau'\in I_{k'}}\sum_{v'=1}^{N(k',\tau')}2^{k'sq}
\left[\mu(Q_{\tau'}^{k',v'})\right]^q\left|\mathcal{P}_{k'}(f)(y_{\tau'}^{k',v'})\right|^q
\mathbf{1}_{\{(\tau',v'):
Q_{\tau'}^{k',v'}\subset Q_{\widetilde{\alpha}}^l\}}(\tau',v')
\nonumber\\
&\quad\times\sum_{k=l}^{\infty}
2^{ksq-k'sq-|k-k'|q\varepsilon'}
\sum_{\tau\in I_k}\sum_{v=1}^{N(k,\tau)}
\mu(Q_\tau^{k,v})\mathbf{1}_{\{(\tau,v):Q_\tau^{k,v}\subset Q_\alpha^l\}}(\tau,v)
\nonumber\\
&\quad\times\left\{\frac{1}{(V_\rho)_{2^{-(k\land k')}}(y_\tau^{k,v})
+V_\rho(y_\tau^{k,v},y^{k',v'}_{\tau'})}
\left[\frac{2^{-(k\land k')}}
{2^{-(k\land k')}+\rho(y_\tau^{k,v},y^{k',v'}_{\tau'})}\right]^{\varepsilon'}\right\}^q.
\end{align}
Note that, from \eqref{upperdoub} and Lemma~\ref{DyadicSys},
it follows that
\begin{align*}
(V_\rho)_{2^{-(k\land k')}}(y_{\tau'}^{k',v'})\lesssim
2^{-(k\land k'-k')n}(V_\rho)_{2^{-k'}}(y_{\tau'}^{k',v'})
\sim2^{-(k\land k'-k')n}
\mu(Q_{\tau'}^{k',v'}).
\end{align*}
By this, \eqref{1442}, and \eqref{9252},
we conclude that
\begin{align*}
\mathrm{I_{1,1}}
&\lesssim\frac{1}{\mu(Q_\alpha^l)}
\sum_{\widetilde{\alpha}\in\mathcal{B}_\alpha^{l,0}}
\sum_{k'=l}^\infty
\sum_{\tau'\in I_{k'}}\sum_{v'=1}^{N(k',\tau')}2^{k'sq}
\left[\mu(Q_{\tau'}^{k',v'})\right]^q\left|\mathcal{P}_{k'}(f)(y_{\tau'}^{k',v'})\right|^q
\mathbf{1}_{\{(\tau',v'):
Q_{\tau'}^{k',v'}\subset Q_{\widetilde{\alpha}}^l\}}(\tau',v')
\\
&\quad\times\sum_{k=l}^{\infty}
2^{ksq-k'sq-|k-k'|q\varepsilon'}
\left[(V_\rho)_{2^{-(k\land k')}}(y^{k',v'}_{\tau'})\right]^{1-q}\\
&\lesssim\frac{1}{\mu(Q_\alpha^l)}
\sum_{\widetilde{\alpha}\in\mathcal{B}_\alpha^{l,0}}
\sum_{k'=l}^\infty
\sum_{\tau'\in I_{k'}}\sum_{v'=1}^{N(k',\tau')}2^{k'sq}
\mu(Q_{\tau'}^{k',v'})\left|\mathcal{P}_{k'}(f)(y_{\tau'}^{k',v'})\right|^q
\mathbf{1}_{\{(\tau',v'):
Q_{\tau'}^{k',v'}\subset Q_{\widetilde{\alpha}}^l\}}(\tau',v')
\\
&\quad\times\sum_{k=l}^{\infty}
2^{ksq-k'sq-|k-k'|q\varepsilon'-(k\land k'-k')n(1-q)}\\
&\sim\sum_{\widetilde{\alpha}\in\mathcal{B}_\alpha^{l,0}}
\frac{1}{\mu(Q_{\widetilde{\alpha}}^l)}
\sum_{k'=l}^\infty
\sum_{\tau'\in I_{k'}}\sum_{v'=1}^{N(k',\tau')}2^{k'sq}
\mu(Q_{\tau'}^{k',v'})\left|\mathcal{P}_{k'}(f)(y_{\tau'}^{k',v'})
\right|^q\mathbf{1}_{\{(\tau',v'):
Q_{\tau'}^{k',v'}\subset Q_{\widetilde{\alpha}}^l\}}(\tau',v')
\end{align*}
with choosing $\varepsilon'\in(s_+,\varepsilon)$
such that $q\in(\frac{n}{n+\varepsilon'+s},1]$,
which, together with the arbitrary choice of $y_{\tau'}^{k',v'}\in Q_{\tau'}^{k',v'}$
and the proven conclusion that $\sup_{l\in\mathbb{Z},\,
\alpha\in I_l}\#\mathcal{B}_\alpha^{l,0}<\infty$,
further implies that
\begin{align}\label{2228}
\mathrm{I_{1,1}}&\lesssim\sup_{l\in\mathbb{Z}}\sup_{\alpha\in I_l}
\frac{1}{\mu(Q_\alpha^l)}
\sum_{k=l}^{\infty}\sum_{\tau\in I_k}\sum_{v=1}^{N(k,\tau)}
2^{ksq}\mu(Q_\tau^{k,v})\nonumber\\
&\quad\times\mathbf{1}_{\{(\tau,v):Q_\tau^{k,v}\subset Q_\alpha^l\}}(\tau,v)
\inf_{x\in Q_\tau^{k,v}}
\left|\mathcal{P}_kf(x)\right|^q.
\end{align}

\emph{Case (2)} $q\in(1,\infty)$. In this case, from
$\sup_{l\in\mathbb{Z},\,
\alpha\in I_l}\#\mathcal{B}_\alpha^{l,0}<\infty$,
we infer that, for any given $\varepsilon'\in
(|s|,\varepsilon)$,
\begin{align*}
\mathrm{I_{1,1}}
&\sim\frac{1}{\mu(Q_\alpha^l)}\sum_{\widetilde{\alpha}\in\mathcal{B}_\alpha^{l,0}}
\sum_{k=l}^{\infty}\sum_{\tau\in I_k}\sum_{v=1}^{N(k,\tau)}
\mu(Q_\tau^{k,v})\mathbf{1}_{\{(\tau,v):Q_\tau^{k,v}\subset Q_\alpha^l\}}(\tau,v)
\\
&\quad\times
\left\{\sum_{k'=l}^\infty
\sum_{\tau'\in I_{k'}}\sum_{v'=1}^{N(k',\tau')}2^{ks-|k-k'|\varepsilon'}
\mu(Q_{\tau'}^{k',v'})
\left|\mathcal{P}_{k'}(f)(y_{\tau'}^{k',v'})\right|
\right.\\
&\quad\times\frac{1}{(V_\rho)_{2^{-(k\land k')}}(y_\tau^{k,v})
+V_\rho(y_\tau^{k,v},y^{k',v'}_{\tau'})}
\left[\frac{2^{-(k\land k')}}
{2^{-(k\land k')}+\rho(y_\tau^{k,v},y^{k',v'}_{\tau'})}\right]^{\varepsilon'}
\\
&\quad\times\left.\mathbf{1}_{\{(\tau',v'):
Q_{\tau'}^{k',v'}\subset Q_{\widetilde{\alpha}}^l\}}(\tau',v')
\right\}^q.
\end{align*}
By this and H\"older's inequality,
we find that
\begin{align*}
\mathrm{I_{1,1}}
&\lesssim\frac{1}{\mu(Q_\alpha^l)}\sum_{\widetilde{\alpha}\in\mathcal{B}_\alpha^{l,0}}
\sum_{k=l}^{\infty}\sum_{\tau\in I_k}\sum_{v=1}^{N(k,\tau)}
\mu(Q_\tau^{k,v})\mathbf{1}_{\{(\tau,v):Q_\tau^{k,v}\subset Q_\alpha^l\}}(\tau,v)
\\
&\quad\times\left[
\left\{\sum_{k'=l}^\infty
\sum_{\tau'\in I_{k'}}\sum_{v'=1}^{N(k',\tau')}2^{k'sq+(k-k')s-|k-k'|\varepsilon'}
\mu(Q_{\tau'}^{k',v'})
\left|\mathcal{P}_{k'}(f)(y_{\tau'}^{k',v'})\right|^q
\right.\right.\\
&\quad\times\frac{1}{(V_\rho)_{2^{-(k\land k')}}(y_\tau^{k,v})
+V_\rho(y_\tau^{k,v},y^{k',v'}_{\tau'})}
\left[\frac{2^{-(k\land k')}}
{2^{-(k\land k')}+\rho(y_\tau^{k,v},y^{k',v'}_{\tau'})}\right]^{\varepsilon'}
\\
&\quad\times\left.\mathbf{1}_{\{(\tau',v'):
Q_{\tau'}^{k',v'}\subset Q_{\widetilde{\alpha}}^l\}}(\tau',v')
\right\}^\frac{1}{q}
\\
&\quad\times
\left\{\sum_{k'=l}^\infty
\sum_{\tau'\in I_{k'}}\sum_{v'=1}^{N(k',\tau')}2^{(k-k')s-|k-k'|\varepsilon'}
\mu(Q_{\tau'}^{k',v'})
\right.\\
&\quad\left.\left.\times\frac{1}{(V_\rho)_{2^{-(k\land k')}}
(y_\tau^{k,v})+V_\rho(y_\tau^{k,v},y^{k',v'}_{\tau'})}
\left[\frac{2^{-(k\land k')}}
{2^{-(k\land k')}+\rho(y_\tau^{k,v},y^{k',v'}_{\tau'})}\right]^{\varepsilon'}
\right\}^\frac{1}{q'}\right]^q,
\end{align*}
which, combined with Lemma~\ref{444},
further implies that
\begin{align*}
\mathrm{I_{1,1}}
&\lesssim\frac{1}{\mu(Q_\alpha^l)}\sum_{\widetilde{\alpha}\in\mathcal{B}_\alpha^{l,0}}
\sum_{k=l}^{\infty}\sum_{\tau\in I_k}\sum_{v=1}^{N(k,\tau)}
\mu(Q_\tau^{k,v})\mathbf{1}_{\{(\tau,v):Q_\tau^{k,v}\subset Q_\alpha^l\}}(\tau,v)
\\
&\quad\times\sum_{k'=l}^\infty
\sum_{\tau'\in I_{k'}}\sum_{v'=1}^{N(k',\tau')}2^{k'sq+(k-k')s-|k-k'|\varepsilon'}
\mu(Q_{\tau'}^{k',v'})
\left|\mathcal{P}_{k'}(f)(y_{\tau'}^{k',v'})\right|^q
\\
&\quad\times\frac{1}{(V_\rho)_{2^{-(k\land k')}}(y_\tau^{k,v})+
V_\rho(y_\tau^{k,v},y^{k',v'}_{\tau'})}
\left[\frac{2^{-(k\land k')}}
{2^{-(k\land k')}+\rho(y_\tau^{k,v},y^{k',v'}_{\tau'})}\right]^{\varepsilon'}
\\
&\quad\times\mathbf{1}_{\{(\tau',v'):
Q_{\tau'}^{k',v'}\subset Q_{\widetilde{\alpha}}^l\}}(\tau',v')
\left\{\sum_{k'=l}^\infty2^{(k-k')s-|k-k'|\varepsilon'}\right\}^\frac{q}{q'}.
\end{align*}
From this, Tonelli's theorem for sequences, Lemma~\ref{444}, and \eqref{9252},
we deduce that
\begin{align*}
\mathrm{I_{1,1}}
&\lesssim\frac{1}{\mu(Q_\alpha^l)}\sum_{\widetilde{\alpha}\in\mathcal{B}_\alpha^{l,0}}
\sum_{k'=l}^\infty
\sum_{\tau'\in I_{k'}}\sum_{v'=1}^{N(k',\tau')}2^{k'sq}
\mu(Q_{\tau'}^{k',v'})
\left|\mathcal{P}_{k'}(f)(y_{\tau'}^{k',v'})\right|^q
\\
&\quad\times\mathbf{1}_{\{(\tau',v'):
Q_{\tau'}^{k',v'}\subset Q_{\widetilde{\alpha}}^l\}}(\tau',v')
\sum_{k=l}^{\infty}2^{(k-k')s-|k-k'|\varepsilon'}\sum_{\tau\in I_k}\sum_{v=1}^{N(k,\tau)}
\mu(Q_\tau^{k,v})
\\
&\quad\times\frac{1}{(V_\rho)_{2^{-(k\land k')}}(y_\tau^{k,v})
+V_\rho(y_\tau^{k,v},y^{k',v'}_{\tau'})}
\left[\frac{2^{-(k\land k')}}
{2^{-(k\land k')}+\rho(y_\tau^{k,v},y^{k',v'}_{\tau'})}\right]^{\varepsilon'}
\\
&\lesssim\sum_{\widetilde{\alpha}\in
\mathcal{B}_\alpha^{l,0}}\frac{1}{\mu(Q_{\widetilde{\alpha}}^l)}
\sum_{k'=l}^\infty
\sum_{\tau'\in I_{k'}}\sum_{v'=1}^{N(k',\tau')}2^{k'sq}
\mu(Q_{\tau'}^{k',v'})
\left|\mathcal{P}_{k'}(f)(y_{\tau'}^{k',v'})\right|^q
\mathbf{1}_{\{(\tau',v'):
Q_{\tau'}^{k',v'}\subset Q_{\widetilde{\alpha}}^l\}}(\tau',v'),
\end{align*}
which, together with the arbitrary choice of $y_{\tau'}^{k',v'}\in Q_{\tau'}^{k',v'}$
and the proven conclusion that $\sup_{l\in\mathbb{Z},\,
\alpha\in I_l}\#\mathcal{B}_\alpha^{l,0}<\infty$,
further implies that
\begin{align}\label{2229}
\mathrm{I_{1,1}}&\lesssim\sup_{l\in\mathbb{Z}}\sup_{\alpha\in I_l}
\frac{1}{\mu(Q_\alpha^l)}
\sum_{k=l}^{\infty}\sum_{\tau\in I_k}\sum_{v=1}^{N(k,\tau)}
2^{ksq}\mu(Q_\tau^{k,v})\nonumber\\
&\quad\times\mathbf{1}_{\{(\tau,v):Q_\tau^{k,v}\subset Q_\alpha^l\}}(\tau,v)
\inf_{x\in Q_\tau^{k,v}}
\left|\mathcal{P}_kf(x)\right|^q.
\end{align}
Next, we turn to deal with $\mathrm{I_{1,2}}$.
For any $j\in\mathbb{N}$, let
\begin{align*}
\mathcal{A}_l^j:=\left\{\widetilde{\alpha}\in I_l:
A2^{j-l-1}\leq\rho(z_{\widetilde{\alpha}}^l,z_\alpha^l)
<A2^{j-l}\right\}
\end{align*}
and
\begin{align*}
\mathcal{B}_\alpha^{l,j}:=\left\{\widetilde{\tau}\in I_{l-j}:\text{there exists }
\widetilde{\alpha}\in\mathcal{A}_l^j\text{ such that }
Q_{\widetilde{\alpha}}^l\subset Q_{\widetilde{\tau}}^{l-j}\right\}.
\end{align*}
Then, by both (i) and (ii) of Lemma~\ref{DyadicSys},
we obtain, for any $j\in\mathbb{N}$,
\begin{align}\label{2233}
\bigcup_{\widetilde{\alpha}\in\mathcal{A}_l^j}Q_{\widetilde{\alpha}}^l
\subset\bigcup_{\widetilde{\tau}\in\mathcal{B}_\alpha^{l,j}}Q_{\widetilde{\tau}}^{l-j}.
\end{align}
Now, we claim that there exists $M\in(0,\infty)$, independent of $\alpha$, $l$, and $j$,
such that
\begin{align}\label{m2}
\#\mathcal{B}_\alpha^{l,j}\leq M.
\end{align}
Indeed, from Definition~\ref{D2.1}(iii), Lemma~\ref{DyadicSys}(iv),
and the definitions of both $\mathcal{A}_l^j$ and $\mathcal{B}_\alpha^{l,j}$,
we infer that, for any given $j\in\mathbb{N}$
and for any $\widetilde{\tau}\in\mathcal{B}_\alpha^{l,j}$ satisfying
$Q_{\widetilde{\alpha}}^l\subset Q_{\widetilde{\tau}}^{l-j}$
for some $\widetilde{\alpha}\in\mathcal{A}_l^j$,
\begin{align}\label{945}
\rho(z_{\widetilde{\tau}}^{l-j},z_\alpha^l)
\leq C_\rho\max\left\{\rho(z_{\widetilde{\tau}}^{l-j},z_{\widetilde{\alpha}}^l),\,
\rho(z_{\widetilde{\alpha}}^l,z_\alpha^l)\right\}
\leq C_\rho\max\left\{C_r,\,A\right\}2^{j-l}.
\end{align}
By this, Lemma~\ref{DyadicSys}(v) and (i),
and the doubling property of $(X,\rho,\mu)$,
we conclude that \eqref{m2} holds. This finishes the proof of
\eqref{m2} and hence the above claim.

Moreover, from \eqref{945}, Lemma~\ref{DyadicSys}(iv),
and \eqref{upperdoub}, we deduce that, for any given
$j\in\mathbb{N}$ and $\widetilde{\alpha}\in\mathcal{A}_l^j$ and
for any $y'\in Q_{\widetilde{\alpha}}^l$,
\begin{align}\label{2252}
\mu\left(\bigcup_{\widetilde{\tau}\in\mathcal{B}_\alpha^{l,j}}
Q_{\widetilde{\tau}}^{l-j}\right)
\leq\mu(B_\rho(z_\alpha^l,M_12^{j-l}))\leq
\mu(B_\rho(y',M_22^{j-l}))
\lesssim(V_\rho)_{2^{j-l}}(y'),
\end{align}
where $M_1,M_2\in(0,\infty)$ are independent of
$j$, $\widetilde{\alpha}$, and
$y'$.
By Definition~\ref{D2.1}(iii),
we find that, for any given $j\in\mathbb{N}$ and
$\widetilde{\alpha}\in\mathcal{A}_l^j$
and for any $y\in Q_\tau^{k,v}\subset Q_\alpha^l$
and $y'\in Q_{\tau'}^{k',v'}\subset Q_{\widetilde{\alpha}}^l$,
\begin{align}\label{2245}
A2^{j-l-1}\leq\rho(z_\alpha^l,z_{\widetilde{\alpha}}^l)\leq (C_\rho)^2
\max\left\{\rho(z_\alpha^l,y),\,\rho(y,y'),\,\rho(y',z_{\widetilde{\alpha}}^l)\right\}.
\end{align}
From Lemma~\ref{DyadicSys}(iv), we infer that
$\max\{\rho(z_\alpha^l,y),\,\rho(y',z_{\widetilde{\alpha}}^l)\}<C_r2^{-l}$.
Choosing $A\in(0,\infty)$ such that $A\ge C_r(C_\rho)^2$,
then \eqref{2245} implies that, for any $y\in Q_\tau^{k,v}\subset Q_\alpha^l$
and $y'\in Q_{\tau'}^{k',v'}\subset Q_{\widetilde{\alpha}}^l$,
\begin{align}\label{2155}
\rho(y,y')\gtrsim2^{j-l}.
\end{align}
By this and \eqref{2252}, we conclude that,
for any $j\in\mathbb{N}$, $\widetilde{\tau}\in\mathcal{B}_\alpha^{l,j}$,
$\widetilde{\alpha}\in\mathcal{A}_l^j$,
$y\in Q_\alpha^l$, and $y'\in Q_{\widetilde{\alpha}}^l$,
\begin{align}\label{939}
\mu(Q_{\widetilde{\tau}}^{l-j})\leq\mu\left(\bigcup_{\widetilde{\tau}
\in\mathcal{B}_\alpha^{l,j}}Q_{\widetilde{\tau}}^{l-j}\right)
\lesssim(V_\rho)_{2^{j-l}}(y')\lesssim V_\rho(y,y').
\end{align}
To estimate $\mathrm{I_{1,2}}$,
we consider the following two cases for $q$.

\emph{Case (1)} $q\in(\max\{\frac{n}{n+\varepsilon},\,
\frac{n}{n+\varepsilon+s}\},1]$. In this case,
from Lemma~\ref{discreteinequality},
we deduce that, for any fixed $\varepsilon'\in(|s|,\varepsilon)$,
\begin{align*}
\mathrm{I_{1,2}}
&\leq\frac{1}{\mu(Q_\alpha^l)}
\sum_{k=l}^{\infty}\sum_{\tau\in I_k}\sum_{v=1}^{N(k,\tau)}
2^{ksq}\mu(Q_\tau^{k,v})\mathbf{1}_{\{(\tau,v):Q_\tau^{k,v}\subset Q_\alpha^l\}}(\tau,v)
\nonumber\\
&\quad\times
\sum_{k'=l}^\infty
\sum_{\tau'\in I_{k'}}\sum_{v'=1}^{N(k',\tau')}2^{-|k-k'|q\varepsilon'}
\left[\mu(Q_{\tau'}^{k',v'})\right]^q
\left|\mathcal{P}_{k'}(f)(y_{\tau'}^{k',v'})\right|^q\nonumber\\
&\quad\times\sum_{j=1}^\infty
\sum_{\widetilde{\alpha}\in\mathcal{A}_j^l}\mathbf{1}_{\{(\tau',v'):
Q_{\tau'}^{k',v'}\subset Q_{\widetilde{\alpha}}^l\}}(\tau',v')
\nonumber\\
&\quad\times\left\{\frac{1}{(V_\rho)_{2^{-(k\land k')}}
(y_\tau^{k,v})+V_\rho(y_\tau^{k,v},y^{k',v'}_{\tau'})}
\left[\frac{2^{-(k\land k')}}
{2^{-(k\land k')}+\rho(y_\tau^{k,v},y^{k',v'}_{\tau'})}\right]^{\varepsilon'}\right\}^q.
\end{align*}
By this and the following two facts that
\begin{align*}
\bigcup_{\widetilde{\alpha}\in\mathcal{B}_\alpha^{l,0}}Q_{\widetilde{\alpha}}^l
\cup\left(\bigcup_{j\in\mathbb{N}}\bigcup_{\widetilde{\beta}\in\mathcal{A}_j^l}
Q_{\widetilde{\beta}}^l\right)=X
\end{align*}
[which can be deduced from Lemma~\ref{DyadicSys}(i)
and the definitions of both $\mathcal{B}_\alpha^{l,0}$
and $\mathcal{A}_l^j$]
and
\begin{align*}
\rho(y_\tau^{k,v},y_{\tau'}^{k',v'})\gtrsim2^{j-l}>2^{-l}\ge2^{-(k\land k')}
\end{align*}
for any $y_\tau^{k,v}\in Q_\tau^{k,v}\subset Q_\alpha^l$
and $y_{\tau'}^{k',v'}\in Q_{\tau'}^{k',v'}\subset Q_{\widetilde{\alpha}}^l$
[which can be deduced from \eqref{2155}],
we obtain
\begin{align*}
\mathrm{I_{1,2}}
&\lesssim\frac{1}{\mu(Q_\alpha^l)}
\sum_{k=l}^{\infty}\sum_{\tau\in I_k}\sum_{v=1}^{N(k,\tau)}
2^{ksq}\mu(Q_\tau^{k,v})\mathbf{1}_{\{(\tau,v):Q_\tau^{k,v}\subset Q_\alpha^l\}}(\tau,v)
\nonumber\\
&\quad\times
\sum_{k'=l}^\infty
\sum_{\tau'\in I_{k'}}\sum_{v'=1}^{N(k',\tau')}2^{-|k-k'|q\varepsilon'}
\left[\mu(Q_{\tau'}^{k',v'})\right]^q
\left|\mathcal{P}_{k'}(f)(y_{\tau'}^{k',v'})\right|^q
\nonumber\\
&\quad\times\sum_{j=1}^\infty
\sum_{\widetilde{\tau}\in\mathcal{A}_j^l}\mathbf{1}_{\{(\tau',v'):
Q_{\tau'}^{k',v'}\subset Q_{\widetilde{\alpha}}^l\}}(\tau',v')
\left\{\frac{1}{V_\rho(y_\tau^{k,v},y^{k',v'}_{\tau'})}
\left[\frac{2^{-(k\land k')}}
{\rho(y_\tau^{k,v},y^{k',v'}_{\tau'})}\right]^{\varepsilon'}\right\}^q.
\end{align*}
From this and \eqref{2233},
we infer that
\begin{align}\label{2236}
\mathrm{I_{1,2}}
&\lesssim\frac{1}{\mu(Q_\alpha^l)}
\sum_{k=l}^{\infty}\sum_{\tau\in I_k}\sum_{v=1}^{N(k,\tau)}
2^{ksq}\mu(Q_\tau^{k,v})\mathbf{1}_{\{(\tau,v):Q_\tau^{k,v}\subset Q_\alpha^l\}}(\tau,v)
\nonumber\\
&\quad\times
\sum_{k'=l}^\infty
\sum_{\tau'\in I_{k'}}\sum_{v'=1}^{N(k',\tau')}2^{-|k-k'|q\varepsilon'}
\left[\mu(Q_{\tau'}^{k',v'})\right]^q
\left|\mathcal{P}_{k'}(f)(y_{\tau'}^{k',v'})\right|^q
\nonumber\\
&\quad\times\sum_{j=1}^\infty
\sum_{\widetilde{\tau}\in\mathcal{B}_\alpha^{l,j}}\mathbf{1}_{\{(\tau',v'):
Q_{\tau'}^{k',v'}\subset Q_{\widetilde{\tau}}^{l-j}\}}(\tau',v')
\left\{\frac{1}{V_\rho(y_\tau^{k,v},y^{k',v'}_{\tau'})}
\left[\frac{2^{-(k\land k')}}
{\rho(y_\tau^{k,v},y^{k',v'}_{\tau'})}\right]^{\varepsilon'}\right\}^q.
\end{align}
By Lemmas~\ref{DyadicSys} and~\ref{A3.2}(vii),
we find that, for any $j\in\mathbb{N}$,
$\widetilde{\tau}\in\mathcal{B}_\alpha^{l,j}$, and
$y_{\tau'}^{k',v'}\in Q_{\tau'}^{k',v'}\subset Q_{\widetilde{\tau}}^{l-j}$,
\begin{align}\label{1512}
(V_\rho)_{2^{j-l}}(y_{\tau'}^{k',v'})\sim\mu(Q^{l-j}_{\widetilde{\tau}}).
\end{align}
Thus,
\begin{align}\label{2204}
&\mathbf{1}_{\{(\tau',v'):
Q_{\tau'}^{k',v'}\subset Q_{\widetilde{\tau}}^{l-j}\}}(\tau',v')
\left\{\frac{1}{V_\rho(y_\tau^{k,v},y^{k',v'}_{\tau'})}
\left[\frac{2^{-(k\land k')}}
{\rho(y_\tau^{k,v},y^{k',v'}_{\tau'})}\right]^{\varepsilon'}\right\}^q
\nonumber\\
&\quad\lesssim\mathbf{1}_{\{(\tau',v'):
Q_{\tau'}^{k',v'}\subset Q_{\widetilde{\tau}}^{l-j}\}}(\tau',v')
\left[\frac{1}{V_\rho(y_\tau^{k,v},y^{k',v'}_{\tau'})}\right]^q
2^{-(k\land k'-l+j)\varepsilon' q}
\quad\text{by  \eqref{2155}}
\nonumber\\
&\quad\lesssim\mathbf{1}_{\{(\tau',v'):
Q_{\tau'}^{k',v'}\subset Q_{\widetilde{\tau}}^{l-j}\}}(\tau',v')
\left[\frac{1}{\mu(Q_{\widetilde{\tau}}^{j-l})}\right]^q
2^{-(k\land k'-l+j)\varepsilon' q}
\quad\text{by \eqref{939}}
\nonumber\\
&\quad\sim\mathbf{1}_{\{(\tau',v'):
Q_{\tau'}^{k',v'}\subset Q_{\widetilde{\tau}}^{l-j}\}}(\tau',v')
\frac{1}{\mu(Q_{\widetilde{\tau}}^{l-j})}\nonumber\\
&\qquad\times\left[(V_\rho)_{2^{j-l}}(y_{\tau'}^{k',v'})\right]^{1-q}
2^{-(k\land k'-l+j)\varepsilon' q}
\quad\text{by \eqref{1512}}.
\end{align}
Choose $\varepsilon'\in(|s|,\varepsilon)$ such that
$q\in(\max\{\frac{n}{n+\varepsilon'},\,
\frac{n}{n+\varepsilon'+s}\},1]$.
Then, from \eqref{2236} and \eqref{2204},
we deduce that
\begin{align*}
\mathrm{I_{1,2}}
&\lesssim\frac{1}{\mu(Q_\alpha^l)}
\sum_{k=l}^{\infty}\sum_{\tau\in I_k}\sum_{v=1}^{N(k,\tau)}
2^{ksq}\mu(Q_\tau^{k,v})\mathbf{1}_{\{(\tau,v):Q_\tau^{k,v}\subset Q_\alpha^l\}}(\tau,v)
\nonumber\\
&\quad\times
\sum_{k'=l}^\infty
\sum_{\tau'\in I_{k'}}\sum_{v'=1}^{N(k',\tau')}2^{-|k-k'|q\varepsilon'}
\left[\mu(Q_{\tau'}^{k',v'})\right]^q
\left|\mathcal{P}_{k'}(f)(y_{\tau'}^{k',v'})\right|^q
\nonumber\\
&\quad\times\sum_{j=1}^\infty
\sum_{\widetilde{\tau}\in\mathcal{B}_\alpha^{l,j}}\mathbf{1}_{\{(\tau',v'):
Q_{\tau'}^{k',v'}\subset Q_{\widetilde{\tau}}^{l-j}\}}(\tau',v')
\frac{1}{\mu(Q_{\widetilde{\tau}}^{l-j})}
\left[(V_\rho)_{2^{j-l}}(y_{\tau'}^{k',v'})\right]^{1-q}
2^{-(k\land k'-l+j)\varepsilon' q}.
\end{align*}
By this and \eqref{upperdoub},
we conclude that
\begin{align*}
\mathrm{I_{1,2}}
&\lesssim\frac{1}{\mu(Q_\alpha^l)}
\sum_{k=l}^{\infty}\sum_{\tau\in I_k}\sum_{v=1}^{N(k,\tau)}
2^{ksq}\mu(Q_\tau^{k,v})\mathbf{1}_{\{(\tau,v):Q_\tau^{k,v}\subset Q_\alpha^l\}}(\tau,v)
\nonumber\\
&\quad\times
\sum_{k'=l}^\infty
\sum_{\tau'\in I_{k'}}\sum_{v'=1}^{N(k',\tau')}2^{-|k-k'|q\varepsilon'}
\left[\mu(Q_{\tau'}^{k',v'})\right]^q
\left|\mathcal{P}_{k'}(f)(y_{\tau'}^{k',v'})\right|^q
\nonumber\\
&\quad\times\sum_{j=1}^\infty
\sum_{\widetilde{\tau}\in\mathcal{B}_\alpha^{l,j}}\mathbf{1}_{\{(\tau',v'):
Q_{\tau'}^{k',v'}\subset Q_{\widetilde{\tau}}^{l-j}\}}(\tau',v')
\frac{1}{\mu(Q_{\widetilde{\tau}}^{l-j})}
\left[(V_\rho)_{2^{-k'}}(y_{\tau'}^{k',v'})\right]^{1-q}
\nonumber\\
&\quad\times2^{(j-l+k')n(1-q)}2^{-(k\land k'-l+j)\varepsilon' q}.
\end{align*}
This and the fact that
$\mu(Q_{\tau'}^{k',v'})\sim(V_\rho)_{2^{-k'}}(y_{\tau'}^{k',v'})$
imply that
\begin{align*}
\mathrm{I_{1,2}}
&\lesssim\sum_{j=1}^\infty2^{j[n(1-q)-\varepsilon' q]}
\sum_{\widetilde{\tau}\in\mathcal{B}_\alpha^{l,j}}
\frac{1}{\mu(Q_{\widetilde{\tau}}^{l-j})}
\sum_{k'=l}^\infty
\sum_{\tau'\in I_{k'}}\sum_{v'=1}^{N(k',\tau')}2^{k'sq}
\mu(Q_{\tau'}^{k',v'})
\nonumber\\
&\quad\times\left|\mathcal{P}_{k'}(f)(y_{\tau'}^{k',v'})\right|^q
\mathbf{1}_{\{(\tau',v'):
Q_{\tau'}^{k',v'}\subset Q_{\widetilde{\tau}}^{l-j}\}}(\tau',v')
\nonumber\\
&\quad\times\sum_{k=l}^\infty
2^{(k-k')sq-|k-k'|q\varepsilon'-(k\land k')\varepsilon' q+k'n(1-q)+l[\varepsilon' q-n(1-q)]},
\end{align*}
which, together with \eqref{m2},
further implies that
\begin{align}\label{2244}
\mathrm{I_{1,2}}
&\leq\sum_{j=1}^\infty2^{j[n(1-q)-\varepsilon' q]}
\sum_{\widetilde{\tau}\in\mathcal{B}_\alpha^{l,j}}
\frac{1}{\mu(Q_{\widetilde{\tau}}^{l-j})}
\sum_{k'=l-j}^\infty
\sum_{\tau'\in I_{k'}}\sum_{v'=1}^{N(k',\tau')}2^{k'sq}
\mu(Q_{\tau'}^{k',v'})
\nonumber\\
&\quad\times\left|\mathcal{P}_{k'}(f)(y_{\tau'}^{k',v'})\right|^q
\mathbf{1}_{\{(\tau',v'):
Q_{\tau'}^{k',v'}\subset Q_{\widetilde{\tau}}^{l-j}\}}(\tau',v')
\nonumber\\
&\lesssim\sup_{l\in\mathbb{Z}}\sup_{\widetilde{\tau}\in I_l}
\frac{1}{\mu(Q_{\widetilde{\tau}}^l)}
\sum_{k'=l}^{\infty}\sum_{\tau'\in I_{k'}}\sum_{v'=1}^{N(k',\tau')}
2^{k'sq}\mu(Q_{\tau'}^{k',v'})\nonumber\\
&\quad\times\left|\mathcal{P}_{k'}f(y_{\tau'}^{k',v'})\right|^q
\mathbf{1}_{\{(\tau',v'):Q_{\tau'}^{k',v'}\subset Q_{\widetilde{\tau}}^l\}}(\tau',v').
\end{align}

\emph{Case (2)} $q\in(1,\infty)$. In this case,
from H\"older's inequality,
it follows that, for any fixed $\varepsilon'\in(|s|,\varepsilon)$,
\begin{align*}
\mathrm{I_{1,2}}
&\lesssim\frac{1}{\mu(Q_\alpha^l)}\sum_{k=l}^\infty
\sum_{\tau\in I_k}\sum_{v=1}^{N(k,\tau)}\mu(Q_\tau^{k,v})
\mathbf{1}_{\{(\tau,v):Q_\tau^{k,v}\subset Q_\alpha^l\}}(\tau,v)
\nonumber\\
&\quad\times\sum_{k'=l}^\infty\sum_{\tau'\in I_{k'}}\sum_{v'=1}^{N(k',\tau')}
2^{-|k-k'|\varepsilon'}2^{-k'sq}2^{(k-k')s}\mu(Q_{\tau'}^{k',v'})
\left|\mathcal{P}_{k'}(f)(y_{\tau'}^{k',v'})\right|^q
\nonumber\\
&\quad\times\mathbf{1}_{\{(\tau',v'):Q_{\tau'}^{k',v'}\cap\bigcup_{\widetilde{\alpha}\in
\mathcal{B}_\alpha^{l,0}}Q_{\widetilde{\alpha}}^l=\varnothing\}}(\tau',v')
\nonumber\\
&\quad\times\frac{1}{(V_\rho)_{2^{-(k\land k')}}(y_\tau^{k,v})
+V_\rho(y_\tau^{k,v},y^{k',v'}_{\tau'})}
\left[\frac{2^{-(k\land k')}}
{2^{-(k\land k')}+\rho(y_\tau^{k,v},y^{k',v'}_{\tau'})}\right]^{\varepsilon'},
\end{align*}
which, combined with Lemma~\ref{444}, \eqref{2155}, \eqref{m2},
\eqref{2252},
and an argument similar to that used in the estimation of \eqref{2244},
further implies that
\begin{align}\label{2246}
\mathrm{I_{1,2}}
&\lesssim\sum_{j=0}^\infty2^{-j\varepsilon'}
\sum_{\widetilde{\tau}\in\mathcal{B}_\alpha^{l,j}}
\frac{1}{\mu(Q_{\widetilde{\tau}}^{l-j})}
\sum_{k'=l-j}^\infty
\sum_{\tau'\in I_{k'}}\sum_{v'=1}^{N(k',\tau')}2^{k'sq}
\mu(Q_{\tau'}^{k',v'})
\nonumber\\
&\quad\times\left|\mathcal{P}_{k'}(f)(y_{\tau'}^{k',v'})\right|^q
\mathbf{1}_{\{(\tau',v'):
Q_{\tau'}^{k',v'}\subset Q_{\widetilde{\tau}}^{l-j}\}}(\tau',v')
\nonumber\\
&\lesssim\sup_{l\in\mathbb{Z}}\sup_{\widetilde{\tau}\in I_l}
\frac{1}{\mu(Q_{\widetilde{\tau}}^l)}
\sum_{k'=l}^{\infty}\sum_{\tau'\in I_{k'}}\sum_{v'=1}^{N(k',\tau')}
2^{k'sq}\mu(Q_{\tau'}^{k',v'})\nonumber\\
&\quad\times\left|\mathcal{P}_{k'}f(y_{\tau'}^{k',v'})\right|^q
\mathbf{1}_{\{(\tau',v'):Q_{\tau'}^{k',v'}\subset Q_{\widetilde{\tau}}^l\}}(\tau',v').
\end{align}
Using \eqref{2116}, \eqref{2228},
\eqref{2229}, \eqref{2244}, \eqref{2246},
and the arbitrary choice of $y_{\tau'}^{k',v'}$
in both \eqref{2244} and \eqref{2246},
we find that
\begin{align}\label{24126}
\mathrm{I}&\lesssim\sup_{l\in\mathbb{Z}}\sup_{\alpha\in I_l}
\frac{1}{\mu(Q_\alpha^l)}
\sum_{k=l}^{\infty}\sum_{\tau\in I_k}\sum_{v=1}^{N(k,\tau)}
2^{ksq}\mu(Q_\tau^{k,v})\nonumber\\
&\quad\times\mathbf{1}_{\{(\tau,v):Q_\tau^{k,v}\subset Q_\alpha^l\}}(\tau,v)
\inf_{x\in Q_\tau^{k,v}}
\left|\mathcal{P}_kf(x)\right|^q.
\end{align}
This finishes the estimate of $\mathrm{I}$.

Next, we turn to estimate $\mathrm{II}$.
Note that Lemma~\ref{DyadicSys} and \eqref{upperdoub}
imply that, for any $k'\in\mathbb{Z}$,
$\tau'\in I_{k'}$, and $v'\in\{1,\ldots,N(k',\alpha')\}$,
$$
\mu(Q_{\tau'}^{k'})\lesssim\mu(Q_{\tau'}^{k',v'}),
$$
which further implies that
\begin{align}\label{1601}
&\sum_{\tau'\in I_{k'}}\sum_{m'=1}^{N(k',\tau')}2^{k'sq}
\left|\mathcal{P}_{k'}f(y_{\tau'}^{k',v'})\right|^q\nonumber\\
&\quad\lesssim\sum_{\tau'\in I_{k'}}\sum_{m'=1}^{N(k',\tau')}2^{k'sq}
\frac{\mu(Q_{\tau'}^{k',v'})}{\mu(Q_{\tau'}^{k'})}
\left|\mathcal{P}_{k'}f(y_{\tau'}^{k',v'})\right|^q\nonumber\\
&\quad=\sum_{\tau'\in I_{k'}}\sum_{m'=1}^{N(k',\tau')}2^{k'sq}
\frac{\mu(Q_{\tau'}^{k',v'})}{\mu(Q_{\tau'}^{k'})}
\left|\mathcal{P}_{k'}f(y_{\tau'}^{k',v'})\right|^q
\mathbf{1}_{\{(\tau',v'):Q_{\tau'}^{k',v'}\subset Q_{\tau'}^{k'}\}}(\tau',v')\nonumber\\
&\quad\leq\sup_{l\in\mathbb{Z}}\sup_{\alpha\in I_l}
\frac{1}{\mu(Q_\alpha^l)}
\sum_{k=l}^{\infty}\sum_{\tau\in I_k}\sum_{v=1}^{N(k,\tau)}
2^{ksq}\mu(Q_\tau^{k,v})\nonumber\\
&\qquad\times\mathbf{1}_{\{(\tau,v):Q_\tau^{k,v}\subset Q_\alpha^l\}}(\tau,v)
\left|\mathcal{P}_kf(y_\tau^{k,v})\right|^q.
\end{align}
Now, we consider the following two cases for $q$.

\emph{Case (1)} $q\in(\max\{\frac{n}{n+\varepsilon},\,\frac{n}{n+\varepsilon+s}\},1]$.
In this case, choosing $\varepsilon'\in(s_+,\varepsilon)$,
from Lemma~\ref{discreteinequality} and
the fact that \begin{align*}
\mu(Q_{\tau'}^{k',v'})&\sim
(V_\rho)_{2^{-k'}}(y_{\tau'}^{k',v'})
\leq(V_\rho)_{2^{-(k\land k')}}(y_{\tau'}^{k',v'})\\
&\lesssim(V_\rho)_{2^{-(k\land k')}}(y_\tau^{k,v})+V_\rho(y_\tau^{k,v},y^{k',v'}_{\tau'}),
\end{align*}
we infer that
\begin{align*}
\mathrm{II}&\leq\frac{1}{\mu(Q_\alpha^l)}
\sum_{k=l}^{\infty}\sum_{\tau\in I_k}\sum_{v=1}^{N(k,\tau)}
2^{ksq}\mu(Q_\tau^{k,v})\mathbf{1}_{\{(\tau,v):Q_\tau^{k,v}\subset Q_\alpha^l\}}(\tau,v)
\nonumber\\
&\quad\times\sum_{k'=-\infty}^{l-1}
\sum_{\tau'\in I_{k'}}\sum_{v'=1}^{N(k',\tau')}2^{-(k-k')q\varepsilon'}
\left[\mu(Q_{\tau'}^{k',v'})\right]^q
\left|\mathcal{P}_{k'}(f)(y_{\tau'}^{k',v'})\right|^q
\nonumber\\
&\quad\times\left\{\frac{1}{(V_\rho)_{2^{-(k\land k')}}(y_\tau^{k,v})
+V_\rho(y_\tau^{k,v},y^{k',v'}_{\tau'})}
\left[\frac{2^{-(k\land k')}}
{2^{-(k\land k')}+\rho(y_\tau^{k,v},y^{k',v'}_{\tau'})}\right]^{\varepsilon'}\right\}^q
\nonumber\\
&\lesssim\frac{1}{\mu(Q_\alpha^l)}
\sum_{k=l}^{\infty}\sum_{\tau\in I_k}\sum_{v=1}^{N(k,\tau)}
2^{ksq}\mu(Q_\tau^{k,v})\mathbf{1}_{\{(\tau,v):Q_\tau^{k,v}\subset Q_\alpha^l\}}(\tau,v)
\nonumber\\
&\quad\times\sum_{k'=-\infty}^{l-1}
\sum_{\tau'\in I_{k'}}\sum_{v'=1}^{N(k',\tau')}2^{-(k-k')q\varepsilon'}
\left|\mathcal{P}_{k'}(f)(y_{\tau'}^{k',v'})\right|^q,
\end{align*}
which, together with \eqref{1601},
further implies that
\begin{align}\label{2231}
\mathrm{II}
&\lesssim\sum_{k=l}^{\infty}2^{kq(s-\varepsilon')}\frac{1}{\mu(Q_\alpha^l)}
\sum_{\tau\in I_k}\sum_{v=1}^{N(k,\tau)}
\mu(Q_\tau^{k,v})\mathbf{1}_{\{(\tau,v):Q_\tau^{k,v}\subset Q_\alpha^l\}}(\tau,v)
\nonumber\\
&\quad\times\sum_{k'=-\infty}^{l-1}2^{k'q(\varepsilon'-s)}
\sum_{\tau'\in I_{k'}}\sum_{v'=1}^{N(k',\tau')}2^{k'sq}
\left|\mathcal{P}_{k'}(f)(y_{\tau'}^{k',v'})\right|^q
\nonumber\\
&\lesssim\sup_{l\in\mathbb{Z}}\sup_{\alpha\in I_l}
\frac{1}{\mu(Q_\alpha^l)}
\sum_{k=l}^{\infty}\sum_{\tau\in I_k}\sum_{v=1}^{N(k,\tau)}
2^{ksq}\mu(Q_\tau^{k,v})\nonumber\\
&\quad\times\mathbf{1}_{\{(\tau,v):Q_\tau^{k,v}\subset Q_\alpha^l\}}(\tau,v)
\left|\mathcal{P}_kf(y_\tau^{k,v})\right|^q.
\end{align}

\emph{Case (2)} $q\in(1,\infty)$. In this case,
choose $\varepsilon'\in(|s|,\varepsilon)$
and $\widetilde{q}\in(1,q)$ such that
\begin{align}\label{1607}
sq<\varepsilon'\widetilde{q}<\varepsilon'q.
\end{align}
By this and H\"older's inequality,
we conclude that
\begin{align*}
\mathrm{II}
&\leq\frac{1}{\mu(Q_\alpha^l)}
\sum_{k=l}^{\infty}\sum_{\tau\in I_k}\sum_{v=1}^{N(k,\tau)}
2^{ksq}\mu(Q_\tau^{k,v})\mathbf{1}_{\{(\tau,v):Q_\tau^{k,v}\subset Q_\alpha^l\}}(\tau,v)
\nonumber\\
&\quad\times\left[\sum_{k'=-\infty}^{l-1}
\sum_{\tau'\in I_{k'}}\sum_{v'=1}^{N(k',\tau')}
2^{-(k-k')\varepsilon'\widetilde{q}}
\left|\mathcal{P}_{k'}(f)(y_{\tau'}^{k',v'})\right|^q
\right]\nonumber\\
&\quad\times\left\{\sum_{k'=-\infty}^{l-1}2^{(k-k')(\frac{\widetilde{q}}{q}-1)\varepsilon'q'}
\sum_{\tau'\in I_{k'}}\sum_{v'=1}^{N(k',\tau')}
\left[\frac{\mu(Q_{\tau'}^{k',v'})}{(V_\rho)_{2^{-(k\land k')}}(y_\tau^{k,v})+
V_\rho(y_\tau^{k,v},y^{k',v'}_{\tau'})}\right]^{q'}
\right.\nonumber\\
&\quad\left.\times
\left[\frac{2^{-(k\land k')}}
{2^{-(k\land k')}+\rho(y_\tau^{k,v},y^{k',v'}_{\tau'})}
\right]^{q'\varepsilon'}\right\}^\frac{q}{q'}.
\end{align*}
From this and Lemma~\ref{discreteinequality},
it follows that
\begin{align*}
\mathrm{II}
&\leq\frac{1}{\mu(Q_\alpha^l)}
\sum_{k=l}^{\infty}\sum_{\tau\in I_k}\sum_{v=1}^{N(k,\tau)}
2^{ksq}\mu(Q_\tau^{k,v})\mathbf{1}_{\{(\tau,v):Q_\tau^{k,v}\subset Q_\alpha^l\}}(\tau,v)
\nonumber\\
&\quad\times\left[\sum_{k'=-\infty}^{l-1}
\sum_{\tau'\in I_{k'}}\sum_{v'=1}^{N(k',\tau')}
2^{-(k-k')\varepsilon'\widetilde{q}}
\left|\mathcal{P}_{k'}(f)(y_{\tau'}^{k',v'})\right|^q
\right]\nonumber\\
&\quad\times\left\{\sum_{k'=-\infty}^{l-1}2^{(k-k')
(\frac{\widetilde{q}}{q}-1)\varepsilon'}
\sum_{\tau'\in I_{k'}}\sum_{v'=1}^{N(k',\tau')}
\frac{\mu(Q_{\tau'}^{k',v'})}{(V_\rho)_{2^{-(k\land k')}}(y_\tau^{k,v})+
V_\rho(y_\tau^{k,v},y^{k',v'}_{\tau'})}
\right.\nonumber\\
&\quad\left.\times
\left[\frac{2^{-(k\land k')}}
{2^{-(k\land k')}+\rho(y_\tau^{k,v},y^{k',v'}_{\tau'})}
\right]^{\varepsilon'}\right\}^q.
\end{align*}
By this, Lemma~\ref{444}, and $\widetilde{q}<q$,
we obtain
\begin{align*}
\mathrm{II}&\lesssim\frac{1}{\mu(Q_\alpha^l)}
\sum_{k=l}^{\infty}\sum_{\tau\in I_k}\sum_{v=1}^{N(k,\tau)}
2^{ksq}\mu(Q_\tau^{k,v})\mathbf{1}_{\{(\tau,v):Q_\tau^{k,v}\subset Q_\alpha^l\}}(\tau,v)
\nonumber\\
&\quad\times\sum_{k'=-\infty}^{l-1}
\sum_{\tau'\in I_{k'}}\sum_{v'=1}^{N(k',\tau')}
2^{-(k-k')\varepsilon'\widetilde{q}}
\left|\mathcal{P}_{k'}(f)(y_{\tau'}^{k',v'})\right|^q,
\end{align*}
which, combined with \eqref{1601} and \eqref{1607},
further implies that
\begin{align}\label{2232}
\mathrm{II}
&\lesssim\sum_{k=l}^{\infty}2^{k(sq-\varepsilon'\widetilde{q})}
\frac{1}{\mu(Q_\alpha^l)}
\sum_{\tau\in I_k}\sum_{v=1}^{N(k,\tau)}
\mu(Q_\tau^{k,v})\mathbf{1}_{\{(\tau,v):Q_\tau^{k,v}\subset Q_\alpha^l\}}(\tau,v)
\nonumber\\
&\quad\times\sum_{k'=-\infty}^{l-1}2^{k'(\varepsilon'\widetilde{q}-sq)}
\sum_{\tau'\in I_{k'}}\sum_{v'=1}^{N(k',\tau')}2^{k'sq}
\left|\mathcal{P}_{k'}(f)(y_{\tau'}^{k',v'})\right|^q
\nonumber\\
&\lesssim\sup_{l\in\mathbb{Z}}\sup_{\alpha\in I_l}
\frac{1}{\mu(Q_\alpha^l)}
\sum_{k=l}^{\infty}\sum_{\tau\in I_k}\sum_{v=1}^{N(k,\tau)}
2^{ksq}\mu(Q_\tau^{k,v})\nonumber\\
&\quad\times\mathbf{1}_{\{(\tau,v):Q_\tau^{k,v}\subset Q_\alpha^l\}}(\tau,v)
\left|\mathcal{P}_kf(y_\tau^{k,v})\right|^q.
\end{align}
From \eqref{2231}, \eqref{2232}, and
the arbitrary choice of $y_\tau^{k,v}$ in both
\eqref{2231} and \eqref{2232}, we further deduce that
\begin{align}\label{2212}
\mathrm{II}&\lesssim\sup_{l\in\mathbb{Z}}\sup_{\alpha\in I_l}
\frac{1}{\mu(Q_\alpha^l)}
\sum_{k=l}^{\infty}\sum_{\tau\in I_k}\sum_{v=1}^{N(k,\tau)}
2^{ksq}\mu(Q_\tau^{k,v})\nonumber\\
&\quad\times\mathbf{1}_{\{(\tau,v):Q_\tau^{k,v}\subset Q_\alpha^l\}}(\tau,v)
\inf_{x\in Q_\tau^{k,v}}
\left|\mathcal{P}_kf(x)\right|^q.
\end{align}
This finishes the estimate of $\mathrm{II}$.

Finally, using \eqref{2115}, \eqref{24126},
and \eqref{2212}
and taking the supremum over all $\alpha\in I_l$ and $l\in\mathbb{Z}$,
we conclude that \eqref{1613} holds.
This finishes the proof of Proposition~\ref{ppineq-infty}.
\end{proof}

\begin{remark}\label{zj2038}
Let the notation be as in Proposition~\ref{ppineq-infty}.
Then, by Proposition~\ref{ppineq-infty},
we easily conclude that the definition of $\dot{F}^{\infty,q}_s(X)$
is independent of the choice of
the approximation of the identity;
that is, for any
$f\in(\mathring{{\mathcal{G}}}_0^\varepsilon(\beta,\gamma))'$,
\begin{align*}
&\sup_{l\in\mathbb{Z}}\sup_{\alpha\in I_l}
\left[\frac{1}{\mu(Q_\alpha^l)}\int_{Q^l_\alpha}\sum_{k=l}^{\infty}
2^{ksq}\left|\mathcal{D}_kf(x)\right|^q\,d\mu(x)\right]^\frac{1}{q}\\
&\quad\sim\sup_{l\in\mathbb{Z}}\sup_{\alpha\in I_l}
\left[\frac{1}{\mu(Q_\alpha^l)}\int_{Q^l_\alpha}\sum_{k=l}^{\infty}
2^{ksq}\left|\mathcal{P}_kf(x)\right|^q\,d\mu(x)\right]^\frac{1}{q},
\end{align*}
where the positive equivalence constants are independent of $f$.
A similar conclusion was also obtained in
\cite[Proposition~6.5]{hmy08}
on RD-spaces
and \cite[Proposition~5.4]{WHHY}
on spaces of homogeneous type.
\end{remark}

We are ready to show Theorem~\ref{2158}.

\begin{proof}[Proof of Theorem~\ref{2158}]
We first prove the sufficiency of (i).
Let
$f\in(\mathring{{\mathcal{G}}}_0^\varepsilon(\beta,\gamma))'$
and $\|f\|_{\widetilde{\dot{F}^{\infty,q}_s}(X)}<\infty$.
By Tonelli's theorem, both (i) and (ii) of Lemma~\ref{DyadicSys},
and Proposition~\ref{ppineq-infty}, we obtain
\begin{align*}
&\sup_{l\in\mathbb{Z}}\sup_{\alpha\in I_l}\left[
\frac{1}{\mu(Q_\alpha^l)}\int_{Q^l_\alpha}\sum_{k=l}^{\infty}
2^{ksq}\left|\mathcal{D}_kf(x)\right|^q\,d\mu(x)\right]^\frac{1}{q}
\nonumber\\
&\quad=\sup_{l\in\mathbb{Z}}\sup_{\alpha\in I_l}
\left[\frac{1}{\mu(Q_\alpha^l)}\sum_{k=l}^{\infty}
\sum_{\tau\in I_k}\sum_{v=1}^{N(k,\tau)}2^{ksq}\int_{Q_\tau^{k,v}}
\mathbf{1}_{Q_\alpha^l}(x)
\left|\mathcal{D}_kf(x)\right|^q\,d\mu(x)\right]^\frac{1}{q}
\nonumber\\
&\quad\leq\sup_{l\in\mathbb{Z}}\sup_{\alpha\in I_l}
\left[\frac{1}{\mu(Q_\alpha^l)}
\sum_{k=l}^{\infty}\sum_{\tau\in I_k}\sum_{v=1}^{N(k,\tau)}
2^{ksq}\mu(Q_\tau^{k,v})\mathbf{1}_{\{(\tau,v):Q_\tau^{k,v}\subset Q_\alpha^l\}}(\tau,v)
\sup_{z\in Q_\tau^{k,v}}
\left|\mathcal{D}_kf(z)\right|^q\right]^\frac{1}{q}
\nonumber\\
&\quad\sim\sup_{l\in\mathbb{Z}}\sup_{\alpha\in I_l}
\left[\frac{1}{\mu(Q_\alpha^l)}
\sum_{k=l}^{\infty}\sum_{\tau\in I_k}\sum_{v=1}^{N(k,\tau)}
2^{ksq}\mu(Q_\tau^{k,v})\mathbf{1}_{\{(\tau,v):Q_\tau^{k,v}\subset Q_\alpha^l\}}(\tau,v)
\inf_{z\in Q_\tau^{k,v}}
\left|\mathcal{D}_kf(z)\right|^q\right]^\frac{1}{q}.
\end{align*}
From this and an argument similar to
that used in the proof of the second inequality of
\eqref{1855}, it follows that
\begin{align*}
&\sup_{l\in\mathbb{Z}}\sup_{\alpha\in I_l}\left[
\frac{1}{\mu(Q_\alpha^l)}\int_{Q^l_\alpha}\sum_{k=l}^{\infty}
2^{ksq}\left|\mathcal{D}_kf(x)\right|^q\,d\mu(x)\right]^\frac{1}{q}
\nonumber\\
&\quad\lesssim\sup_{l\in\mathbb{Z}}\sup_{\alpha\in I_l}\left[\frac{1}{\mu(Q_\alpha^l)}
\sum_{k=l}^{\infty}\sum_{\tau\in I_k}\sum_{v=1}^{N(k,\tau)}
2^{ksq}\mu(Q_\tau^{k,v})\mathbf{1}_{\{(\tau,v):Q_\tau^{k,v}\subset Q_\alpha^l\}}(\tau,v)
\right.\nonumber\\
&\qquad\left.\times
\frac{1}{\mu(Q_\tau^{k,v})}\int_{Q_\tau^{k,v}}
\left|\mathcal{D}_kf(y)\right|^q\,d\mu(y)\right]^\frac{1}{q}
\nonumber\\
&\quad\sim\sup_{l\in\mathbb{Z}}\sup_{\alpha\in I_l}
\left[\frac{1}{\mu(Q_\alpha^l)}
\sum_{k=l}^{\infty}\sum_{\tau\in I_k}\sum_{v=1}^{N(k,\tau)}
2^{ksq}\mu(Q_\tau^{k,v})\mathbf{1}_{\{(\tau,v):Q_\tau^{k,v}\subset Q_\alpha^l\}}(\tau,v)
\right.\nonumber\\
&\qquad\left.\times
\inf_{x\in Q_\tau^{k,v}}\frac{1}{V_{2^{-k}}(x)}
\int_{Q_\tau^{k,v}}
\left|\mathcal{D}_kf(y)\right|^q\,d\mu(y)\right]^\frac{1}{q}.
\end{align*}
By this and \eqref{21411},
we find that
\begin{align}\label{2009}
&\sup_{l\in\mathbb{Z}}\sup_{\alpha\in I_l}\left[
\frac{1}{\mu(Q_\alpha^l)}\int_{Q^l_\alpha}\sum_{k=l}^{\infty}
2^{ksq}\left|\mathcal{D}_kf(x)\right|^q\,d\mu(x)\right]^\frac{1}{q}
\nonumber\\
&\quad\lesssim\sup_{l\in\mathbb{Z}}\sup_{\alpha\in I_l}
\left[\frac{1}{\mu(Q_\alpha^l)}
\sum_{k=l}^{\infty}\sum_{\tau\in I_k}\sum_{v=1}^{N(k,\tau)}
2^{ksq}\mu(Q_\tau^{k,v})\mathbf{1}_{\{(\tau,v):Q_\tau^{k,v}\subset Q_\alpha^l\}}(\tau,v)
\right.\nonumber\\
&\qquad\left.\times
\inf_{x\in Q_\tau^{k,v}}\frac{1}{(V_\rho)_{2^{-k}}(x)}\int_{B_\rho(x,2^{-k})}
\left|\mathcal{D}_kf(y)\right|^q\,d\mu(y)\right]^\frac{1}{q}
\nonumber\\
&\quad\leq
\sup_{l\in\mathbb{Z}}\sup_{\alpha\in I_l}
\left[\frac{1}{\mu(Q_\alpha^l)}\int_{Q^l_\alpha}\sum_{k=l}^{\infty}
2^{ksq}\int_{B_\rho(x,2^{-k})}\left|\mathcal{D}_kf(y)\right|^q
\frac{d\mu(y)\,d\mu(x)}{(V_\rho)_{2^{-k}}(x)}\right]^\frac{1}{q}
\end{align}
with the usual modification when $q=\infty$.
From this, we infer that
\begin{align*}
\left\|f\right\|_{\dot{F}^{\infty,q}_s(X)}
\lesssim\left\|f\right\|_{\widetilde{\dot{F}^{\infty,q}_s}(X)}<\infty
\end{align*}
and hence $f\in\dot{F}^{\infty,q}_s(X)$,
which completes the proof of the sufficiency of (i).

Next, we show the necessity of (i).
To this end, let $f\in\dot{F}^{\infty,q}_s(X)
\subset(\mathring{{\mathcal{G}}}_0^\varepsilon(\beta,\gamma))'$.
By Tonelli's theorem, both (i) and (ii) of Lemma~\ref{DyadicSys},
\eqref{1055}, and \eqref{21412},
we conclude that
\begin{align*}
&\sup_{l\in\mathbb{Z}}\sup_{\alpha\in I_l}
\left[\frac{1}{\mu(Q_\alpha^l)}\int_{Q^l_\alpha}\sum_{k=l}^{\infty}
2^{ksq}\int_{B_\rho(x,2^{-k})}\left|\mathcal{D}_kf(y)\right|^q
\frac{d\mu(y)\,d\mu(x)}{(V_\rho)_{2^{-k}}(x)}\right]^\frac{1}{q}
\\
&\quad=\sup_{l\in\mathbb{Z}}\sup_{\alpha\in I_l}
\left[\frac{1}{\mu(Q_\alpha^l)}\sum_{k=l}^{\infty}
\sum_{\tau\in I_k}\sum_{v=1}^{N(k,\tau)}2^{ksq}\int_{Q_\tau^{k,v}}\mathbf{1}_{Q_\alpha^l}(x)
\right.\\
&\qquad\times
\left.\int_{B_\rho(x,2^{-k})}\left|\mathcal{D}_kf(y)\right|^q
\frac{d\mu(y)\,d\mu(x)}{(V_\rho)_{2^{-k}}(x)}\right]^\frac{1}{q}
\\
&\quad\lesssim\sup_{l\in\mathbb{Z}}\sup_{\alpha\in I_l}
\left[\frac{1}{\mu(Q_\alpha^l)}\sum_{k=l}^{\infty}\sum_{\tau\in I_k}\sum_{v=1}^{N(k,\tau)}
2^{ksq}\mu(Q_\tau^{k,v})\mathbf{1}_{\{(\tau,v):Q_\tau^{k,v}\subset Q_\alpha^l\}}(\tau,v)
\right.
\\
&\qquad\times
\left.
\frac{1}{(V_\rho)_{C_\rho2^{-k}}(z_\tau^{k,v})}\int_{B_\rho(z_\tau^{k,v},C_\rho2^{-k})}
\left|\mathcal{D}_kf(y)\right|^q\,d\mu(y)\right]^\frac{1}{q}.
\end{align*}
Using this and Proposition~\ref{ppineq-infty} with $c_0:=C_\rho$,
we obtain
\begin{align*}
&\sup_{l\in\mathbb{Z}}\sup_{\alpha\in I_l}
\left[\frac{1}{\mu(Q_\alpha^l)}\int_{Q^l_\alpha}\sum_{k=l}^{\infty}
2^{ksq}\int_{B_\rho(x,2^{-k})}\left|\mathcal{D}_kf(y)\right|^q
\frac{d\mu(y)\,d\mu(x)}{(V_\rho)_{2^{-k}}(x)}\right]^\frac{1}{q}
\\
&\quad\lesssim\sup_{l\in\mathbb{Z}}\sup_{\alpha\in I_l}
\left\{\frac{1}{\mu(Q_\alpha^l)}\sum_{k=l}^{\infty}\sum_{\tau\in I_k}\sum_{v=1}^{N(k,\tau)}
2^{ksq}\mu(Q_\tau^{k,v})\mathbf{1}_{\{(\tau,v):Q_\tau^{k,v}\subset Q_\alpha^l\}}(\tau,v)
\right.
\\
&\qquad\times
\left.
\sup_{y\in B_\rho(z_\tau^{k,v},C_\rho2^{-k})}
\left|\mathcal{D}_kf(y)\right|^q\right\}^\frac{1}{q}
\\
&\quad\sim
\sup_{l\in\mathbb{Z}}\sup_{\alpha\in I_l}
\left\{\frac{1}{\mu(Q_\alpha^l)}
\sum_{k=l}^{\infty}\sum_{\tau\in I_k}\sum_{v=1}^{N(k,\tau)}
2^{ksq}\mu(Q_\tau^{k,v})\mathbf{1}_{\{(\tau,v):Q_\tau^{k,v}\subset Q_\alpha^l\}}(\tau,v)\inf_{x\in Q_\tau^{k,v}}
\left|\mathcal{D}_kf(x)\right|^q\right\}^\frac{1}{q}
\\
&\quad\leq
\sup_{l\in\mathbb{Z}}\sup_{\alpha\in I_l}
\left[\frac{1}{\mu(Q_\alpha^l)}\int_{Q^l_\alpha}\sum_{k=l}^{\infty}
2^{ksq}\left|\mathcal{D}_kf(x)\right|^q\,d\mu(x)\right]^\frac{1}{q},
\end{align*}
which further implies that
\begin{align*}
\|f\|_{\widetilde{\dot{F}^{\infty,q}_s}(X)}
\lesssim\|f\|_{\dot{F}^{\infty,q}_s(X)}<\infty.
\end{align*}
This finishes the proof of the necessity of (i)
and hence the proof of (i).

Now, we turn to prove the sufficiency of (ii).
To this end, let
$f\in(\mathring{{\mathcal{G}}}_0^\varepsilon(\beta,\gamma))'$
and $\|f\|_{\overline{\dot{F}^{\infty,q}_s}(X)}<\infty$.
Note that, from Lemma~\ref{A3.2}(vii), it follows that,
for any $k\in\mathbb{Z}$ and $x\in X$,
\begin{align*}
&\int_{B_\rho(x,2^{-k})}\left|\mathcal{D}_kf(y)\right|^q
\frac{d\mu(y)}{(V_\rho)_{2^{-k}}(x)}\\
&\quad\leq2^\lambda \int_{B_\rho(x,2^{-k})}\left|\mathcal{D}_kf(y)\right|^q
\left[\frac{2^{-k}}{2^{-k}+\rho(x,y)}\right]^\lambda
\frac{d\mu(y)}{(V_\rho)_{2^{-k}}(x)}\\
&\quad\sim2^\lambda\int_{B_\rho(x,2^{-k})}
\left|\mathcal{D}_kf(y)\right|^q
\left[\frac{2^{-k}}{2^{-k}+\rho(x,y)}\right]^\lambda
\frac{d\mu(y)}{(V_\rho)_{2^{-k}}(x)+(V_\rho)_{2^{-k}}(y)}\\
&\quad\leq2^\lambda\int_{X}\left|\mathcal{D}_kf(y)\right|^q
\left[\frac{2^{-k}}{2^{-k}+\rho(x,y)}\right]^\lambda
\frac{d\mu(y)}{(V_\rho)_{2^{-k}}(x)+(V_\rho)_{2^{-k}}(y)},
\end{align*}
which, together with \eqref{2009}, further implies that
\begin{align*}
\|f\|_{\dot{F}^{\infty,q}_s(X)}
\lesssim\|f\|_{\widetilde{\dot{F}^{\infty,q}_s}(X)}
\lesssim\|f\|_{\overline{\dot{F}^{\infty,q}_s}(X)}<\infty
\end{align*}
and hence $f\in\dot{F}^{\infty,q}_s(X)$.
This finishes the proof of the sufficiency of (ii).

Next, we show the necessity of (ii).
To this end, let $f\in\dot{F}^{\infty,q}_s(X)
\subset(\mathring{{\mathcal{G}}}_0^\varepsilon(\beta,\gamma))'$.
By Lemma~\ref{A3.2}(vi) and \eqref{upperdoub},
we find that, for any $k\in\mathbb{Z}$ and $x,y\in X$,
\begin{align*}
(V_\rho)_{2^{-k}}(x)+V_\rho(x,y)
&\sim(V_\rho)_{2^{-k}+\rho(x,y)}(x)+(V_\rho)_{2^{-k}+\rho(x,y)}(y)
\\
&\lesssim\left[\frac{2^{-k}+\rho(x,y)}{2^{-k}}\right]^n
\left[(V_\rho)_{2^{-k}}(x)+(V_\rho)_{2^{-k}}(y)\right],
\end{align*}
which further implies that
\begin{align*}
&\int_X\left|\mathcal{D}_kf(y)\right|^q
\left[\frac{2^{-k}}{2^{-k}+\rho(x,y)}\right]^\lambda
\frac{d\mu(y)}{(V_\rho)_{2^{-k}}(x)+(V_\rho)_{2^{-k}}(y)}
\nonumber\\
&\quad=\int_{B_\rho(x,2^{-k})}\left|\mathcal{D}_kf(y)\right|^q
\left[\frac{2^{-k}}{2^{-k}+\rho(x,y)}\right]^\lambda
\frac{d\mu(y)}{(V_\rho)_{2^{-k}}(x)+(V_\rho)_{2^{-k}}(y)}
\nonumber\\
&\qquad+\sum_{m=1}^\infty\int_{B_\rho(x,2^{-k+m})\setminus
B_\rho(x,2^{-k+m-1})}\ldots
\nonumber\\
&\quad\lesssim\int_{B_\rho(x,2^{-k})}\left|\mathcal{D}_kf(y)
\right|^q\frac{d\mu(y)}{(V_\rho)_{2^{-k}}(x)}
+\sum_{m=1}^\infty\int_{B_\rho(x,2^{-k+m})\setminus
B_\rho(x,2^{-k+m-1})}\left|\mathcal{D}_kf(y)\right|^q
\nonumber\\
&\qquad\times\left[\frac{2^{-k}}{2^{-k}+\rho(x,y)}
\right]^{\lambda-n}\frac{d\mu(y)}{(V_\rho)_{2^{-k}}(x)+V_\rho(x,y)}
\end{align*}
and hence
\begin{align}\label{2309}
&\int_X\left|\mathcal{D}_kf(y)\right|^q
\left[\frac{2^{-k}}{2^{-k}+\rho(x,y)}\right]^\lambda
\frac{d\mu(y)}{(V_\rho)_{2^{-k}}(x)+(V_\rho)_{2^{-k}}(y)}
\nonumber\\
&\quad\lesssim\int_{B_\rho(x,2^{-k})}\left|\mathcal{D}_kf(y)
\right|^q\frac{d\mu(y)}{(V_\rho)_{2^{-k}}(x)}
\nonumber\\
&\qquad+\sum_{m=1}^\infty2^{-m(\lambda-n)}\int_{B_\rho(x,2^{-k+m})\setminus
B_\rho(x,2^{-k+m-1})}\left|\mathcal{D}_kf(y)\right|^q\frac{d\mu(y)}{(V_\rho)_{2^{-k+m}}(x)}.
\end{align}
Moreover, from Definition~\ref{D2.1}(iii), we deduce that, for any
$k,k'\in\mathbb{Z}$, $\tau'\in I_{k'}$,
$v'\in\{1,...,N(k',\tau')\}$,
$m\in\mathbb{N}$,
$x\in X$, and $y\in B_\rho(x,2^{-k+m})\setminus
B_\rho(x,2^{-k+m-1})$,
\begin{align*}
2^{-(k\land k')}+\rho(x,y_{\tau'}^{k',v'})
&\leq2^{-(k\land k')}+C_\rho\max\left\{\rho(x,y),\,\rho(y,y_{\tau'}^{k',v'})\right\}
\\
&<2^{-(k\land k')}+C_\rho\max\left\{2^{-k+m},\,\rho(y,y_{\tau'}^{k',v'})\right\}
\\
&\lesssim2^m\left[2^{-(k\land k')}+\rho(y,y_{\tau'}^{k',v'})\right],
\end{align*}
where $y_{\tau'}^{k',v'}$ is any fixed point in $Q_{\tau'}^{k',v'}$,
which implies that, for any given $\gamma\in(0,\infty)$,
\begin{align}\label{16292}
&\frac{1}{(V_\rho)_{2^{-(k\land k')}}(y)+V_\rho(y,y^{k',v'}_{\tau'})}
\left[\frac{2^{-(k\land k')}}{2^{-(k\land k')}+\rho(y,y_{\tau'}^{k',v'})}\right]^{\gamma}
\nonumber\\
&\quad\sim\frac{1}{(V_\rho)_{2^{-(k\land k')}+\rho(y,y_{\tau'}^{k',v'})}(y_{\tau'}^{k',v'})}
\left[\frac{2^{-(k\land k')}}{2^{-(k\land k')}+\rho(y,y_{\tau'}^{k',v'})}\right]^{\gamma}
\quad\text{by Lemma~\ref{A3.2}(vi)}
\nonumber\\
&\quad\lesssim\frac{2^{m(n+\gamma)}}{(V_\rho)_{2^{-(k\land k')}
+\rho(x,y_{\tau'}^{k',v'})}(y_{\tau'}^{k',v'})}
\left[\frac{2^{-(k\land k')}}{2^{-(k\land k')}+\rho(x,y_{\tau'}^{k',v'})}\right]^{\gamma}
\quad\text{by \eqref{upperdoub}}
\nonumber\\
&\quad\sim\frac{2^{m(n+\gamma)}}{(V_\rho)_{2^{-(k\land k')}}(x)+V_\rho(x,y^{k',v'}_{\tau'})}
\left[\frac{2^{-(k\land k')}}{2^{-(k\land k')}+\rho(x,y_{\tau'}^{k',v'})}\right]^{\gamma}.
\end{align}
By Theorem~\ref{HD},
we conclude that, for any given $\varepsilon'\in(0,\varepsilon)$
and for any $k\in\mathbb{Z}$, $x\in X$, and $y\in B_\rho(x,2^{-k+m})\setminus
B_\rho(x,2^{-k+m-1})$ if $m\in\mathbb{N}$
or $y\in B_\rho(x,2^{-k})$ if $m=0$,
\begin{align*}
\left|\mathcal{D}_kf(y)\right|
&=\left|\left\langle f,D_k(y,\cdot)\right\rangle\right|
\nonumber\\
&=\left|\left\langle\sum_{k'=-\infty}^{\infty}\sum_{\tau'\in I_{k'}}
\sum_{v'=1}^{N(k',\tau')}\mu(Q_{\tau'}^{k',v'})
\widetilde{D}_{k'}(\cdot,y_{\tau'}^{k',v'})
\mathcal{D}_{k'}f(y_{\tau'}^{k',v'}),D_k(y,\cdot)\right\rangle
\right|\nonumber\\
&=\left|\sum_{k'=-\infty}^{\infty}\sum_{\tau'\in I_{k'}}
\sum_{v'=1}^{N(k',\tau')}\mu(Q_{\tau'}^{k',v'})\left\langle
\widetilde{D}_{k'}(\cdot,y_{\tau'}^{k',v'})
,D_k(y,\cdot)\right\rangle\mathcal{D}_{k'}f(y_{\tau'}^{k',v'})
\right|\nonumber\\
&=\left|\sum_{k'=-\infty}^\infty
\sum_{\tau'\in I_{k'}}\sum_{v'=1}^{N(k',\tau')}
\mu(Q_{\tau'}^{k',v'})
D_k\widetilde{D}_{k'}(y,y_{\tau'}^{k',v'})
\mathcal{D}_{k'}f(Q_{\tau'}^{k',v'})\right|,
\end{align*}
which, combined with Lemma~\ref{1122}
and \eqref{16292}, implies that
\begin{align*}
\left|\mathcal{D}_kf(y)\right|
&\lesssim\sum_{k'=-\infty}^\infty2^{-|k-k'|\varepsilon'}
\sum_{\tau'\in I_{k'}}\sum_{v'=1}^{N(k',\tau')}
\mu(Q_{\tau'}^{k',v'})
\nonumber\\
&\quad\times
\frac{1}{(V_\rho)_{2^{-(k\land k')}}(y)+V_\rho(y,y^{k',v'}_{\tau'})}
\left[\frac{2^{-(k\land k')}}{2^{-(k\land k')}+
\rho(y,y_{\tau'}^{k',v'})}\right]^{\varepsilon'}
\left|\mathcal{D}_{k'}f(y_{\tau'}^{k',v'})\right|
\nonumber\\
&\lesssim2^{m(n+\varepsilon')}\sum_{k'=-\infty}^\infty2^{-|k-k'|\varepsilon'}
\sum_{\tau'\in I_{k'}}\sum_{v'=1}^{N(k',\tau')}
\mu(Q_{\tau'}^{k',v'})
\nonumber\\
&\quad\times
\frac{1}{(V_\rho)_{2^{-(k\land k')}}(x)+V_\rho(x,y^{k',v'}_{\tau'})}
\left[\frac{2^{-(k\land k')}}{2^{-(k\land k')}+
\rho(x,y_{\tau'}^{k',v'})}\right]^{\varepsilon'}
\left|\mathcal{D}_{k'}f(y_{\tau'}^{k',v'})\right|.
\end{align*}
This, together with \eqref{2309} and $\lambda\in(\max\{2n,2n-sq,n+nq+|s|q\},\infty)$,
further implies that, for any fixed $\varepsilon'\in(|s|,\frac{\lambda-n-nq}{q})$,
\begin{align*}
&\int_X\left|\mathcal{D}_kf(y)\right|^q
\left[\frac{2^{-k}}{2^{-k}+\rho(x,y)}\right]^\lambda
\frac{d\mu(y)}{(V_\rho)_{2^{-k}}(x)+(V_\rho)_{2^{-k}}(y)}
\\
&\quad\lesssim\sum_{m=0}^\infty2^{m(n+\varepsilon'q+nq-\lambda)}
\left\{\sum_{k'=-\infty}^\infty2^{-|k-k'|\varepsilon'}
\sum_{\tau'\in I_{k'}}\sum_{v'=1}^{N(k',\tau')}
\mu(Q_{\tau'}^{k',v'})\right.
\nonumber\\
&\qquad\times\left.
\frac{1}{(V_\rho)_{2^{-(k\land k')}}(x)+V_\rho(x,y^{k',v'}_{\tau'})}
\left[\frac{2^{-(k\land k')}}{2^{-(k\land k')}+
\rho(x,y_{\tau'}^{k',v'})}\right]^{\varepsilon'}
\left|\mathcal{D}_{k'}f(y_{\tau'}^{k',v'})\right|\right\}^q
\\
&\quad\sim\left\{\sum_{k'=-\infty}^\infty2^{-|k-k'|\varepsilon'}
\sum_{\tau'\in I_{k'}}\sum_{v'=1}^{N(k',\tau')}
\mu(Q_{\tau'}^{k',v'})\right.
\nonumber\\
&\qquad\times\left.
\frac{1}{(V_\rho)_{2^{-(k\land k')}}(x)+V_\rho(x,y^{k',v'}_{\tau'})}
\left[\frac{2^{-(k\land k')}}{2^{-(k\land k')}+
\rho(x,y_{\tau'}^{k',v'})}\right]^{\varepsilon'}
\left|\mathcal{D}_{k'}f(y_{\tau'}^{k',v'})\right|\right\}^q.
\end{align*}
From this, Tonelli's theorem,
and both (i) and (ii) of Lemma~\ref{DyadicSys},
we infer that, for any $l\in\mathbb{Z}$ and $\alpha\in I_l$,
\begin{align*}
&\int_{Q^l_\alpha}\sum_{k=l}^{\infty}
2^{ksq}\int_{X}\left|\mathcal{D}_kf(y)\right|^q
\left[\frac{2^{-k}}{2^{-k}+\rho(x,y)}\right]^\lambda\frac{d\mu(y)
\,d\mu(x)}{(V_\rho)_{2^{-k}}(x)+(V_\rho)_{2^{-k}}(y)}
\\
&\quad\lesssim\int_{Q^l_\alpha}\sum_{k=l}^{\infty}
2^{ksq}\left\{\sum_{k'=-\infty}^\infty2^{-|k-k'|\varepsilon'}
\sum_{\tau'\in I_{k'}}\sum_{v'=1}^{N(k',\tau')}
\mu(Q_{\tau'}^{k',v'})\right.
\nonumber\\
&\qquad\times\left.
\frac{1}{(V_\rho)_{2^{-(k\land k')}}(x)+V_\rho(x,y^{k',v'}_{\tau'})}
\left[\frac{2^{-(k\land k')}}{2^{-(k\land k')}+
\rho(x,y_{\tau'}^{k',v'})}\right]^{\varepsilon'}
\left|\mathcal{D}_{k'}f(y_{\tau'}^{k',v'})\right|\right\}^q\,d\mu(x)
\\
&\quad\lesssim\sum_{k=l}^{\infty}2^{ksq}\sum_{\tau\in I_k}\sum_{v=1}^{N(k,\tau)}
\mathbf{1}_{\{(\tau,v):Q_{\tau}^{k,v}\subset Q_\alpha^l\}}(\tau,v)\int_{Q_{\tau}^{k,v}}
\left\{\sum_{k'=l}^\infty2^{-|k-k'|\varepsilon'}
\sum_{\tau'\in I_{k'}}\sum_{v'=1}^{N(k',\tau')}
\mu(Q_{\tau'}^{k',v'})\right.
\\
&\qquad\times\left.
\frac{1}{(V_\rho)_{2^{-(k\land k')}}(x)+V_\rho(x,y^{k',v'}_{\tau'})}
\left[\frac{2^{-(k\land k')}}{2^{-(k\land k')}+\rho(x,y_{\tau'}^{k',v'})}\right]^{\varepsilon'}
\left|\mathcal{D}_{k'}f(y_{\tau'}^{k',v'})\right|\right\}^q\,d\mu(x),
\end{align*}
which, combined with the fact that
\begin{align*}
2^{-(k\land k')}+\rho(x,y_{\tau'}^{k',v'})
\sim2^{-(k\land k')}+\rho(y_{\tau}^{k,v},y_{\tau'}^{k',v'})
\end{align*}
for any $x,y_{\tau}^{k,v}\in Q_{\tau}^{k,v}$,
further implies that
\begin{align*}
&\int_{Q^l_\alpha}\sum_{k=l}^{\infty}
2^{ksq}\int_{X}\left|\mathcal{D}_kf(y)\right|^q
\left[\frac{2^{-k}}{2^{-k}+\rho(x,y)}\right]^\lambda\frac{d\mu(y)
\,d\mu(x)}{(V_\rho)_{2^{-k}}(x)+(V_\rho)_{2^{-k}}(y)}\\
&\quad\lesssim\sum_{k=l}^{\infty}2^{ksq}\sum_{\tau\in I_k}\sum_{v=1}^{N(k,\tau)}\mu(Q_{\tau}^{k,v})
\mathbf{1}_{\{(\tau,v):Q_{\tau}^{k,v}\subset Q_\alpha^l\}}(\tau,v)
\left\{\sum_{k'=l}^\infty2^{-|k-k'|\varepsilon'}
\sum_{\tau'\in I_{k'}}\sum_{v'=1}^{N(k',\tau')}
\mu(Q_{\tau'}^{k',v'})\right.
\\
&\qquad\times\left.
\frac{1}{(V_\rho)_{2^{-(k\land k')}}(y_\tau^{k,v})+V_\rho(y_\tau^{k,v},y^{k',v'}_{\tau'})}
\left[\frac{2^{-(k\land k')}}{2^{-(k\land k')}+
\rho(y_\tau^{k,v},y_{\tau'}^{k',v'})}\right]^{\varepsilon'}
\left|\mathcal{D}_{k'}f(y_{\tau'}^{k',v'})\right|\right\}^q.
\end{align*}
This, together with the same argument as that used in the
estimations of both $\mathrm{I}$ and $\mathrm{II}$ in the proof of
Proposition~\ref{ppineq-infty}
[see \eqref{2116}, \eqref{2228}, \eqref{2229}, \eqref{2244},
and \eqref{2246} for $\mathrm{I}$ and see \eqref{2231} and \eqref{2232}
for $\mathrm{II}$], implies that, for any fixed
$\varepsilon'\in(\max\{\frac{n}{q}-n,\frac{n}{q}-n-s,|s|\},\varepsilon)$,
\begin{align*}
&\int_{Q^l_\alpha}\sum_{k=l}^{\infty}
2^{ksq}\int_{X}\left|\mathcal{D}_kf(y)\right|^q
\left[\frac{2^{-k}}{2^{-k}+\rho(x,y)}\right]^\lambda
\frac{d\mu(y)\,d\mu(x)}{(V_\rho)_{2^{-k}}(x)+(V_\rho)_{2^{-k}}(y)}
\\
&\quad\lesssim\sup_{l\in\mathbb{Z}}\sup_{\alpha\in I_l}
\sum_{k=l}^\infty\sum_{\tau\in I_k}\sum_{v=1}^{N(k,\tau)}2^{ksq}
\mu(Q_\tau^{k,v})\mathbf{1}_{\{(\tau,v):Q_\tau^{k,v}\subset Q_\alpha^l\}}(\tau,v)
\left|\mathcal{D}_{k}f(y_{\tau}^{k,v})\right|^q.
\end{align*}
Using this, the arbitrariness of $y_\tau^{k,v}\in Q_\tau^{k,v}$,
and Lemma~\ref{DyadicSys} and
taking the supremum over all $\alpha\in I_l$
and $l\in\mathbb{Z}$,
we find that
\begin{align*}
&\sup_{l\in\mathbb{Z}}\sup_{\alpha\in I_l}
\left\{\frac{1}{\mu(Q_\alpha^l)}\int_{Q^l_\alpha}\sum_{k=l}^{\infty}
2^{ksq}\int_{X}\left|\mathcal{D}_kf(y)\right|^q\right.
\\	 &\qquad\left.\times\left[\frac{2^{-k}}{2^{-k}+\rho(x,y)}\right]^\lambda
\frac{d\mu(y)\,d\mu(x)}{(V_\rho)_{2^{-k}}(x)+(V_\rho)_{2^{-k}}(y)}\right\}^\frac{1}{q}
\\
&\quad\lesssim\sup_{l\in\mathbb{Z}}\sup_{\alpha\in I_l}
\left[\frac{1}{\mu(Q_\alpha^l)}\int_{Q^l_\alpha}\sum_{k=l}^{\infty}
2^{ksq}\left|\mathcal{D}_kf(x)\right|^q\,d\mu(x)\right]^\frac{1}{q}
\end{align*}
and hence
\begin{align*}
\|f\|_{\overline{\dot{F}^{\infty,q}_s}(X)}
\lesssim\|f\|_{\dot{F}^{\infty,q}_s(X)}<\infty,
\end{align*}
which completes the proof of the necessity of (ii).
This finishes the proof of (ii)
and hence Theorem~\ref{2158}.
\end{proof}

\section[Inhomogeneous Triebel--Lizorkin Spaces
$F^{p,q}_s(X)$ and $F^{\infty,q}_s(X)$]{Inhomogeneous Triebel--Lizorkin Spaces
$F^{p,q}_s(X)$ and\\ $F^{\infty,q}_s(X)$}
\label{subsection5.2}

In this section,
we recall the definition
of the inhomogeneous Triebel--Lizorkin space
$F^{p,q}_s(X)$ on the quasi-ultrametric space of homogeneous type
$(X,\mathbf{q},\mu)$; see also \cite[Definition 9.6]{AlMi2015} on
Ahlfors-regular	 quasi-ultrametric spaces,
\cite[Definitions~5.29(ii) and~6.19]{hmy08} on RD-spaces,
and \cite[Definitions~4.1(ii) and~5.9]{WHHY}
on spaces of homogeneous type.
From now on, we do not need to assume that
$\mu(X)=\infty$; that is, $\mu(X)$ can be finite or infinite.

Let $\rho\in\mathbf{q}$ and
$E$ be a $\rho$-bounded and $\mu$-measurable
subset of $X$ satisfying that $\mu(E)\in(0,\infty)$.
For any $f\in L^1_{\mathrm{loc}}(X)$,
the \emph{integral average of $f$ over $E$}, denoted by $m_E(f)$\index[symbols]{E@$m_E(f)$},
is defined by setting
\begin{align*}
m_E(f):=\fint_Ef\,d\mu:=\frac{1}{\mu(E)}\int_Ef(x)\,d\mu(x).
\end{align*}
Additionally, recall that, for any $r\in\mathbb{R}$, $r_+:=\max\{0,\,r\}$.
With these in mind, we now record the definition of inhomogeneous
Triebel--Lizorkin spaces.

\begin{definition}\label{iTL}
Let $(X,\mathbf{q},\mu)$
be a quasi-ultrametric space of homogeneous type
with upper dimension $n\in(0,\infty)$ and
$\rho\in\mathbf{q}$ be such that all $\rho$-balls
are $\mu$-measurable.
Let $0<\varepsilon\preceq\mathrm{ind\,}(X,\mathbf{q})$
with $\mathrm{ind\,}(X,\mathbf{q})$ as in \eqref{index}.
Assume that $\{\mathcal{S}_{2^{-k}}\}_{k\in\mathbb{Z}_+}$ is
an inhomogeneous approximation of the identity
of order $\varepsilon$ as in Definition~\ref{DefIATI}
(whose existence is ensured
by Theorem~\ref{ThmATI}).
Let $\mathcal{D}_{0}:=\mathcal{S}_{1}$ and
$\mathcal{D}_k:=\mathcal{S}_{2^{-k}}-\mathcal{S}_{2^{-k+1}}$
for any $k\in\mathbb{N}$. Let $N\in\mathbb{N}$ be the
same as in Theorem~\ref{ID}.
Assume that $s\in(-\varepsilon,\varepsilon)$,
$$
p,q\in\left(\max\left\{\frac{n}{n+\varepsilon},\,
\frac{n}{n+\varepsilon+s}\right\},\infty\right],\ \
\gamma\in\left(n\left(\frac{1}{p}-1\right)_+,\varepsilon\right),
$$
and
$$
\beta\in\left(\max\left\{s_+,\,
-s+n\left(\frac{1}{p}-1\right)_+\right\},\varepsilon\right).
$$
The \emph{inhomogeneous Triebel--Lizorkin
space}\index[words]{inhomogeneous Triebel--Lizorkin space}
$F^{p,q}_s(X)$\index[symbols]{F@$F^{p,q}_s(X)$}
is defined to be the set of all distributions
$f\in({\mathcal{G}}_0^\varepsilon(\beta,\gamma))'$
such that, when $p<\infty$,
\begin{align*}
\|f\|_{F^{p,q}_s(X)}
:&=\left\{\sum_{k=0}^N\sum_{\tau\in I_{k}}
\sum_{v=1}^{N(k,\tau)}
\mu(Q_\tau^{k,v})
\left[m_{Q_\tau^{k,v}}\left(\left|\mathcal{D}_{k}f\right|\right)\right]^p\right\}^\frac{1}{p}\\
&\quad+\left\|\left(\sum_{k=N+1}^\infty
2^{ksq}\left|\mathcal{D}_kf\right|^q\right)^\frac{1}{q}\right\|_{L^p(X)}
\end{align*}
is finite and, when $p=\infty$,
\begin{align*}
\|f\|_{F^{\infty,q}_s(X)}
:&=\max\Bigg\{\sup_{k\in\{0,\ldots,N\}}
\sup_{\tau\in I_{k}}\sup_{v\in\{1,...,N(k,\tau)\}}m_{Q_\tau^{k,v}}
\left(\left|\mathcal{D}_{0}f\right|\right),\,\\
&\qquad\sup_{l\in\mathbb{N},\,l>N}
\sup_{\alpha\in I_l}\left[\frac{1}{\mu(Q_\alpha^l)}\int_{Q_\alpha^l}\sum_{k=l}^\infty
2^{ksq}\left|\mathcal{D}_kf(x)\right|^q\,d\mu(x)\right]^\frac{1}{q}\Bigg\}
\end{align*}
(with the usual modification when $q=\infty$) is finite,
where, for any $k\in\mathbb{Z}$ and $\tau\in I_k$,
$Q_\tau^{k,v}$ is the dyadic cube at level $j+k$
such that $Q_\tau^{k,v}\subset Q_\tau^k$ as
in Remark~\ref{jcubes}.
\end{definition}

\begin{remark}
By the proof of Lemma~\ref{1955}, we conclude that, for any
$k\in\mathbb{Z_+}$ and
$f\in({{\mathcal{G}}}_0^\varepsilon(\beta,\gamma))'$,
$\mathcal{D}_kf$ is a continuous function [when constructed using a regularized
quasi-ultrametric as in Definition~\ref{D2.3}] and hence
the functions being integrated
in Definition~\ref{iTL}
are $\mu$-measurable.
\end{remark}

\begin{remark}\label{22013}
The following statements regarding Definition~\ref{iTL} are true.
\begin{enumerate}
\item[\rm(i)]
If $\mathbf{q}$ contains an ultrametric,
then $\mathrm{ind\,}(X,\mathbf{q})=\infty$
and $s$ in Definition~\ref{iTL} can be
chosen arbitrarily from $(-\infty,\infty)$.

\item[\rm(ii)]
Let $(X,\mathbf{q},\mu)$ be an ultra-RD-space
and $\mu$ be a Borel-semiregular measure on $X$.
Then the definition of $F^{p,q}_s(X)$ is independent
of the choice of the approximation of the identity
(see Remarks~\ref{2205}(i) and~\ref{22012}(i) below)
as well as the choice of the indices $\beta,\gamma$; see also
\cite[Propositions~5.27 and~6.17]{hmy08} and \cite[Propositions~5.28 and~6.18]{hmy08},
respectively.
Moreover, by Remark~\ref{Rm}(iii),
we find that ${F}^{p,q}_s(X)$ is independent of
the choice of the quasi-ultrametric $\rho$
if $\rho\in\mathbf{q}$ is such that
$\mathcal{D}_kf$ for any $k\in\mathbb{Z_+}$
and $f\in({\mathcal{G}}_0^\varepsilon(\beta,\gamma))'$ is $\mu$-measurable.

\item[\rm(iii)]
The space $F^{p,q}_s(X)$ is quasi-Banach
for all $s$, $p$, $q$, $\beta$, and $\gamma$ as in Definition~\ref{iTL}, and
genuinely Banach whenever $p,q\geq1$; see also
\cite[Propositions~5.31(vii) and~6.21(vi)]{hmy08}.
\end{enumerate}
\end{remark}

\section[Lusin Area Function
and Littlewood--Paley $g^*_\lambda$-Function
Characterizations of $F^{p,q}_s(X)$ and $F^{\infty,q}_s(X)$]
{Lusin Area Function
and Littlewood--Paley $g^*_\lambda$-\\Function
Characterizations of $F^{p,q}_s(X)$ and $F^{\infty,q}_s(X)$}
\label{Section4.5}

In this section, we record
the inhomogeneous Lusin area function and the
inhomogeneous Littlewood--Paley $g^*_\lambda$-function
characterizations of the inhomogeneous Triebel--Lizorkin spaces
on the ultra-RD-space $(X,\mathbf{q},\mu)$.
In what follows, we have no restriction on $\mu(X)$,
which means $\mu(X)<\infty$
or $\mu(X)=\infty$ (equivalently, $X$ can be bounded or unbounded).
To this end, we first recall the following
inhomogeneous Littlewood--Paley functions
on quasi-ultrametric spaces of homogeneous type.

\begin{definition}\label{inhomoareafuncts}
Let $(X,\mathbf{q},\mu)$
be a quasi-ultrametric space of homogeneous type.
Let $\rho\in\mathbf{q}$ be such that all $\rho$-balls
are $\mu$-measurable and
$\varepsilon,\beta,\gamma\in(0,\infty)$
be such that $0<\beta,\gamma<\varepsilon\preceq\mathrm{ind\,}(X,\mathbf{q})$
with $\mathrm{ind\,}(X,\mathbf{q})$ as in \eqref{index}.
Assume that $\{\mathcal{S}_{2^{-k}}\}_{k\in\mathbb{Z}_+}$ is
an inhomogeneous approximation of the identity
of order $\varepsilon$ as in Definition~\ref{DefIATI}.
Let $\mathcal{D}_{0}:=\mathcal{S}_{1}$ and
$\mathcal{D}_k:=\mathcal{S}_{2^{-k}}-\mathcal{S}_{2^{-k+1}}$
for any $k\in\mathbb{N}$.
Let $N\in\mathbb{N}$ be the same as in Theorem~\ref{ID}.
Assume that $\{Q^{k,v}_\tau\}_{k\in\{0,\dots,N\},\,
\tau\in I_k,\,v\in\{1,\dots,N(0,\tau)\}}$
is the collection of dyadic cubes as in Remark~\ref{jcubes}.
Let $s\in(-\varepsilon,\varepsilon)$, $q\in(0,\infty]$,
and $f\in(\mathring{{\mathcal{G}}}_0^\varepsilon(\beta,\gamma))'$.
\begin{enumerate}
\item
The \emph{inhomogeneous Littlewood--Paley
$g$-function}\index[words]{inhomogeneous Littlewood--Paley
$g$-function} $g^s_q(f)$ of $f$\index[symbols]{G@$g^s_q(f)$}
is defined by setting, for any $x\in X$,
\begin{align*}
g^s_q(f)(x):&=\sum_{k=0}^N\sum_{\tau\in I_k}
\sum_{v=1}^{N(k,\tau)}
m_{Q_\tau^{k,v}}\left(\left|\mathcal{D}_{k}f\right|\right)
\mathbf{1}_{Q_\tau^{k,v}}(x)
+\left[\sum_{k=N+1}^\infty2^{ksq}\left|\mathcal{D}_kf(x)\right|^q\right]^\frac{1}{q}.
\end{align*}
\item
The \emph{inhomogeneous Lusin area function}
\index[words]{inhomogeneous Lusin area function}
$\mathcal{S}^s_q(f)$ of $f$\index[symbols]{S@$\mathcal{S}^s_q(f)$}
is defined by setting, for any $x\in X$,
\begin{align*}
\mathcal{S}^s_q(f)(x):=\left[\sum_{k={0}}^\infty
2^{ksq}
\int_{B_\rho(x,2^{-k})}
\left|\mathcal{D}_kf(y)\right|^q
\frac{d\mu(y)}{(V_\rho)_{2^{-k}}(x)}\right]^\frac{1}{q}.
\end{align*}
\item
The \emph{inhomogeneous Littlewood--Paley $g_\lambda^*$-function}
\index[words]{inhomogeneous Littlewood--Paley $g_\lambda^*$-function}
$(g_\lambda^*)^s_q(f)$ of $f$,\index[symbols]{G@$(g_\lambda^*)^s_q(f)$}
with any given $\lambda\in(0,\infty)$,
is defined by setting, for any $x\in X$,
\begin{align*}
(g_\lambda^*)^s_q(f)(x):
&=\left\{\sum_{k={0}}^\infty2^{ksq}
\int_X\left|\mathcal{D}_kf(y)\right|^q\left[\frac{2^{-k}}{2^{-k}+
\rho(x,y)}\right]^\lambda
\frac{d\mu(y)}{(V_\rho)_{2^{-k}}(x)+(V_\rho)_{2^{-k}}(y)}
\right\}^\frac{1}{q}.
\end{align*}
\end{enumerate}
\end{definition}

\begin{remark}
Let the notation be as in Definition~\ref{inhomoareafuncts}.
\begin{enumerate}
\item[\rm(i)]
From an argument similar to that used in the proof of
Lemma~\ref{1955}, we infer that, for any
$f\in(\mathring{{\mathcal{G}}}_0^\varepsilon(\beta,\gamma))'$,
the functions $g^s_q(f)$, $\mathcal{S}^s_q(f)$,
and $(g_\lambda^*)^s_q(f)$ are $\mu$-measurable,
where we need that $\rho\in\mathbf{q}$ is a regularized
quasi-ultrametric as in Definition~\ref{D2.3}.

\item[\rm(ii)]
If $g^s_q(f)$ is constructed using $\rho$ and $\widetilde{g^s_q}(f)$
is constructed using $\widetilde{\rho}\in\mathbf{q}$ (where all $\rho$-balls
and all $\widetilde{\rho}$-balls are $\mu$-measurable),
then $g^s_q(f)\sim\widetilde{g^s_q}(f)$ pointwise in $X$.
Similar equivalences for $\mathcal{S}^s_q(f)$ and $(g_\lambda^*)^s_q(f)$ also hold.
Thus, without loss of generality, we may assume that $\rho\in\mathbf{q}$
in Definition~\ref{inhomoareafuncts} is a regularized quasi-ultrametric
in the remainder of this monograph.
\end{enumerate}
\end{remark}

Using the definitions of $g^s_q(f)$ and $F^{p,q}_s(X)$,
Lemma~\ref{DyadicSys}(i), and Remark~\ref{jcubes},
we immediately obtain the following theorem.

\begin{theorem}
Let $(X,\mathbf{q},\mu)$, $n$, $\varepsilon$, $s$, $p$, $q$, $\beta$, and $\gamma$ be
the same as in Definition~\ref{iTL} with $p\neq\infty$.
Then $f\in{F}^{p,q}_s(X)$ if and only if
$f\in({{\mathcal{G}}}_0^\varepsilon(\beta,\gamma))'$
and ${g}^s_q(f)\in L^p(X)$;
moreover, for such $f$,
$\|g^s_q(f)\|_{L^p(X)}\sim\|f\|_{F^{p,q}_s(X)}$,
where the positive equivalence constants are independent of $f$.
\end{theorem}

In order to establish the Lusin area function and
Littlewood--Paley $g_\lambda^*$-function characterizations
$F^{p,q}_s(X)$ with $p<\infty$, we need the following inhomogeneous analogue of the
Plancherel--P\'olya inequality as in Lemma~\ref{P-P }.
The two key tools in establishing this is the
inhomogeneous discrete Calder\'on reproducing formula
(Theorem~\ref{ID2}) and the almost orthogonal estimate (Lemma~\ref{1122}).

\begin{lemma}
\label{ppin}
Let $(X,\mathbf{q},\mu)$ be an ultra-RD-space with upper dimension $n\in(0,\infty)$,
where $\mu$ is assumed to be a Borel-semiregular measure on $X$.
Let $\varepsilon$, $s$, $p$, $q$, $\beta$, and $\gamma$ be
the same as in Definition~\ref{iTL} with $p\neq\infty$.
Assume that $\{\mathcal{S}_{2^{-k}}\}_{k\in\mathbb{Z}_+}$
and $\{\mathcal{Q}_{2^{-k}}\}_{k\in\mathbb{Z}_+}$
are two inhomogeneous approximations of the identity
of order $\varepsilon$ as in Definition~\ref{DefIATI}.
For any $k\in\mathbb{N}$, let
$\mathcal{D}_k:=\mathcal{S}_{2^{-k}}-\mathcal{S}_{2^{-k+1}}$,
$\mathcal{P}_k:=\mathcal{Q}_{2^{-k}}-\mathcal{Q}_{2^{-k+1}}$,
$\mathcal{D}_0:=\mathcal{S}_1$, and $\mathcal{P}_0:=\mathcal{Q}_1$.
Let $N,N'\in\mathbb{N}$ be the same as in Theorem~\ref{ID}
associated, respectively, with $\{\mathcal{D}_k\}_{k\in\mathbb{Z_+}}$
and $\{\mathcal{P}_k\}_{k\in\mathbb{Z_+}}$.
Then, for any given $c_0\in[1,\infty)$, there exists a positive constant $C$ such that,
for any $f\in({{\mathcal{G}}}_0^\varepsilon(\beta,\gamma))'$,
\begin{align}\label{2127}
&\left\{\sum_{k=0}^N\sum_{\tau\in I_k}
\sum_{v=1}^{N(k,\tau)}
\mu(Q_\tau^{k,v})
\left[m_{Q_\tau^{k,v}}\left(\left|\mathcal{D}_kf\right|\right)
\right]^p\right\}^\frac{1}{p}
\nonumber\\
&\quad+\left\|\left[\sum_{k=N+1}^\infty\sum_{\tau\in I_k}\sum_{v=1}^{N(k,\tau)}
2^{ksq}\sup_{z\in B_\rho(z_\tau^{k,v},c_02^{-k})}\left|\mathcal{D}_kf(z)\right|^q
\mathbf{1}_{Q_\tau^{k,v}}\right]^\frac{1}{q}\right\|_{L^p(X)}
\nonumber\\
&\quad\leq C\left\{\sum_{k=0}^{N'}\sum_{\tau\in I_k}
\sum_{v=1}^{N(k,\tau)}
\mu(Q_\tau^{k,v})
\left[m_{Q_\tau^{k,v}}\left(\left|\mathcal{P}_kf\right|\right)
\right]^p\right\}^\frac{1}{p}
\nonumber\\
&\qquad
+C\left\|\left[\sum_{k=N'+1}^\infty\sum_{\tau\in I_k}\sum_{v=1}^{N(k,\tau)}
2^{ksq}\inf_{z\in Q_\tau^{k,v}}\left|\mathcal{P}_kf(z)\right|^q
\mathbf{1}_{Q_\tau^{k,v}}\right]^\frac{1}{q}\right\|_{L^p(X)},
\end{align}
where $\{Q^{k,v}_\tau\}_{k\in\mathbb{Z}_+,\,\tau\in I_k,\,v\in\{1,\dots,N(k,\tau)\}}$
(with centers $\{z^{k,v}_\tau\}_{k\in\mathbb{Z}_+,\,\tau\in I_k,\,v\in\{1,\dots,N(k,\tau)\}}$)
is the collection of dyadic cubes given in Remark~\ref{jcubes}.
\end{lemma}

\begin{proof}
Let $f\in({{\mathcal{G}}}_0^\varepsilon(\beta,\gamma))'$.
We first estimate the first term of the left-hand side of \eqref{2127}.
From the inhomogeneous discrete Calder\'on reproducing formula in Theorem~\ref{ID2},
we deduce that, for any $k\in\mathbb{Z}_+$ and $x\in X$,
\begin{align}\label{462}
\left|\mathcal{D}_kf(x)\right|
&\leq\sum_{\tau'\in I_0}\sum_{v'=1}^{N(0,\tau')}
\int_{Q_{\tau'}^{0,v'}}\left|D_k\widetilde{P}_0(x,y)\right|\,d\mu(y)
\left|\mathcal{P}^{0,v'}_{\tau',1}(f)\right|
\nonumber\\
&\quad+\sum_{k'=1}^{N'}\sum_{\tau'\in I_{k'}}\sum_{v'=1}^{N(k',\tau')}
\mu(Q_{\tau'}^{k',v'})\left|D_k\widetilde{P}_{k'}(x,y_{\tau'}^{k',v'})\right|
\left|\mathcal{P}_{\tau',1}^{k',v'}f\right|
\nonumber\\
&\quad+\sum_{k'=N'+1}^{\infty}\sum_{\tau'\in I_{k'}}\sum_{v'=1}^{N(k',\tau')}
\mu(Q_{\tau'}^{k',v'})\left|D_k\widetilde{P}_{k'}(x,y_{\tau'}^{k',v'})\right|
\left|\mathcal{P}_{k'}f(y_{\tau'}^{k',v'})\right|.
\end{align}

Next, we claim that, for any $k'\in\{0,\dots,N'\}$,
$\tau'\in I_{k'}$, and $v'\in\{1,\dots,N(k',\tau')\}$,
\begin{align}\label{464}
\left|\mathcal{P}_{\tau',1}^{k',v'}f\right|
\leq m_{Q_{\tau'}^{k',v'}}(|\mathcal{P}_{k'}f|).
\end{align}
To this end, we first prove that,
for any $k'\in\{0,\dots,N'\}$,
$\tau'\in I_{k'}$, and $v'\in\{1,\dots,N(k',\tau')\}$,
\begin{align}\label{1958}
\left|\left\langle f,\fint_{Q_{\tau'}^{k',v'}}P_{k'}(x,\cdot)\,d\mu(x)\right\rangle\right|
\leq\fint_{Q_{\tau'}^{k',v'}}\left|\left\langle f,P_{k'}(x,\cdot)\right\rangle\right|\,d\mu(x).
\end{align}
Indeed, if $f\in L^2(X)$, then, by Fubini's theorem,
we immediately conclude that
\eqref{1958} holds.
From Theorem~\ref{ID2}, it follows that there exists a sequence $\{f_n\}_{n=1}^\infty$ in $L^2(X)$
such that $f_n\to f$ in $({{\mathcal{G}}}_0^\varepsilon(\beta,\gamma))'$ as $n\to\infty$.
By this and the Banach--Steinhaus theorem (see, for instance, \cite[Theorem~2.2]{b2011}),
we find that $\{\|f_n\|_{({{\mathcal{G}}}_0^\varepsilon(\beta,\gamma))'}\}_{n\in\mathbb{N}}$
is bounded, which, combined with \eqref{1958} for $\{f_n\}_{n\in\mathbb{N}}$ and
the Lebesgue dominated convergence theorem
(here, we used the fact that
$|\langle f_n,P_{k'}(x,\cdot)\rangle|\leq\sup_{n\in\mathbb{N}}
\|f_n\|_{({{\mathcal{G}}}_0^\varepsilon(\beta,\gamma))'}
\|P_{k'}(x,\cdot)\|_{{{\mathcal{G}}}_0^\varepsilon(\beta,\gamma)}$
for any $x\in Q_{\tau'}^{k',v'}$),
further implies that
\begin{align*}
&\left|\left\langle f,\fint_{Q_{\tau'}^{k',v'}}P_{k'}(x,\cdot)
\,d\mu(x)\right\rangle\right|
\\
&\quad=\lim_{n\to\infty}\left|\left\langle f_n,
\fint_{Q_{\tau'}^{k',v'}}P_{k'}(x,\cdot)\,d\mu(x)\right\rangle\right|
\leq\lim_{n\to\infty}
\fint_{Q_{\tau'}^{k',v'}}\left|\left\langle f_n,P_{k'}(x,\cdot)
\right\rangle\right|\,d\mu(x)
\\
&\quad=\fint_{Q_{\tau'}^{k',v'}}\lim_{n\to\infty}
\left|\left\langle f_n,P_{k'}(x,\cdot)\right\rangle\right|\,d\mu(x)
=\fint_{Q_{\tau'}^{k',v'}}\left|\left\langle f,P_{k'}(x,\cdot)
\right\rangle\right|\,d\mu(x).
\end{align*}
This finishes the proof of \eqref{1958}. With this in mind,
we easily conclude that,
for any $k'\in\{0,\dots,N'\}$,
$\tau'\in I_{k'}$, and $v'\in\{1,\dots,N(k',\tau')\}$,
\begin{align*}
\left|\mathcal{P}_{\tau',1}^{k',v'}f\right|
\leq\fint_{Q_{\tau'}^{k',v'}}\left|\left\langle f,P_{k'}(x,\cdot)
\right\rangle\right|\,d\mu(x)
=m_{Q_{\tau'}^{k',v'}}(|\mathcal{P}_{k'}f|),
\end{align*}
which completes the proof of
\eqref{464} and hence the above claim.

From \eqref{462} and \eqref{464}, we infer that, for any $k\in\{0,\ldots,N\}$,
$\tau\in I_k$, and $v\in\{1,\ldots,N(k,\tau)\}$,
\begin{align}\label{946}
m_{Q_\tau^{k,v}}(|\mathcal{D}_kf|)
&=\frac{1}{\mu(Q_\tau^{k,v})}\int_{Q_\tau^{k,v}}\left|\mathcal{D}_kf(z)\right|\,d\mu(z)\leq\mathrm{I}+\mathrm{II}+\mathrm{III},
\end{align}
where
\begin{align*}
\mathrm{I}:=\sum_{\tau'\in I_0}\sum_{v'=1}^{N(0,\tau')}
\frac{1}{\mu(Q_\tau^{k,v})}\int_{Q_\tau^{k,v}}
\int_{Q_{\tau'}^{0,v'}}\left|D_k\widetilde{P}_0(z,y)\right|\,d\mu(y)\,d\mu(z)
m_{Q_{\tau'}^{0,v'}}(|\mathcal{P}_{0}f|),
\end{align*}
\begin{align*}
\mathrm{II}:&=\sum_{k'=1}^{N'}\sum_{\tau'\in I_{k'}}\sum_{v'=1}^{N(k',\tau')}
\mu(Q_{\tau'}^{k',v'})\\
&\quad\times\frac{1}{\mu(Q_\tau^{k,v})}\int_{Q_\tau^{k,v}}
\left|D_k\widetilde{P}_{k'}(z,y_{\tau'}^{k',v'})\right|\,d\mu(z)
m_{Q_{\tau'}^{k',v'}}(|\mathcal{P}_{k'}f|),
\end{align*}
and
\begin{align*}
\mathrm{III}:&=\sum_{k'=N'+1}^{\infty}\sum_{\tau'\in I_{k'}}\sum_{v'=1}^{N(k',\tau')}
\mu(Q_{\tau'}^{k',v'})\\
&\quad\times\frac{1}{\mu(Q_\tau^{k,v})}\int_{Q_\tau^{k,v}}
\left|D_k\widetilde{P}_{k'}(z,y_{\tau'}^{k',v'})\right|\,d\mu(z)
\left|\mathcal{P}_{k'}f(y_{\tau'}^{k',v'})\right|.
\end{align*}
Note that, for any $k\in\{0,\ldots,N\}$ and
$k'\in\{0,\ldots,N'\}$, we have $2^{-(k\land k')}\sim1\sim2^{-|k-k'|}$
with the positive equivalence constants depending only on both $N$ and $N'$,
which, together with Lemma~\ref{1122},
further implies that, for any fixed $\varepsilon'\in(0,\varepsilon)$ and for any $y,z\in X$,
\begin{align}\label{2039}
\left|D_k\widetilde{P}_{k'}(z,y)\right|
&\lesssim2^{-|k-k'|\varepsilon'}
\frac{1}{(V_\rho)_{2^{-(k\land k')}}(y)
+V_\rho(z,y)}\left[\frac{2^{-(k\land k')}}{2^{-(k\land k')}+
\rho(z,y)}\right]^{\varepsilon'}
\nonumber\\
&\sim\frac{1}{(V_\rho)_{1}(y)
+V_\rho(z,y)}\left[\frac{1}{1+\rho(z,y)}\right]^{\varepsilon'}.
\end{align}
Observe that, for any given $k\in\{0,\ldots,N\}$ and
$k'\in\{0,\ldots,N'\}$, by Definition~\ref{2.1}(iii) and Lemma~\ref{DyadicSys}(iv),
we find that, for any $y_1,y_2\in Q_{\tau'}^{k',v'}$
and $z_1,z_2\in B_\rho(z_\tau^{k,v},c_02^{-k})$,
\begin{align}\label{1706}
1+\rho(z_1,y_1)
&\leq1+C_\rho\max\left\{\rho(z_1,y_2),\,\rho(y_2,y_1)\right\}
\nonumber\\
&\lesssim1+\max\left\{\rho(z_1,y_2),\,2^{-k'}\right\}
\sim1+\rho(z_1,y_2)
\nonumber\\
&\leq1+C_\rho\max\left\{\rho(z_1,z_2),\,\rho(z_2,y_2)\right\}
\nonumber\\
&\lesssim1+\max\left\{2^{-k},\,\rho(z_2,y_2)\right\}
\sim1+\rho(z_2,y_2),
\end{align}
which, combined with both \eqref{2039} and Lemma~\ref{A3.2}(vi),
further implies that
\begin{align*}
&\sup_{y\in Q_{\tau'}^{k',v'}}\sup_{z\in B_\rho(z_\tau^{k,v},c_02^{-k})}
\left|D_k\widetilde{P}_{k'}(z,y)\right|
\nonumber\\
&\quad\lesssim\inf_{y\in Q_{\tau'}^{k',v'}}\inf_{z\in B_\rho(z_\tau^{k,v},c_02^{-k})}
\frac{1}{(V_\rho)_{1}(y)
+V_\rho(z,y)}\left[\frac{1}{1+\rho(z,y)}\right]^{\varepsilon'}
\\
&\quad\leq\inf_{y\in Q_{\tau'}^{k',v'}}\inf_{z\in Q_\tau^{k,v}}
\frac{1}{(V_\rho)_{1}(y)
+V_\rho(z,y)}\left[\frac{1}{1+\rho(z,y)}\right]^{\varepsilon'}.
\end{align*}
From this, we deduce that
\begin{align}\label{2048}
\mathrm{I}+\mathrm{II}&\lesssim
\sum_{k'=0}^{N'}\sum_{\tau'\in I_{k'}}\sum_{v'=1}^{N(k',\tau')}
\mu(Q_{\tau'}^{k',v'})m_{Q_{\tau'}^{k',v'}}(|\mathcal{P}_{k'}f|)
\nonumber\\
&\quad\times\inf_{y\in Q_{\tau'}^{k',v'}}\inf_{z\in Q_\tau^{k,v}}
\frac{1}{(V_\rho)_{1}(y)
+V_\rho(z,y)}\left[\frac{1}{1+\rho(z,y)}\right]^{\varepsilon'}.
\end{align}
Now, we consider the following two cases for $p$.

\emph{Case (1)} $p\in(\max\{\frac{n}{n+\varepsilon},\,\frac{n}{n+\varepsilon+s}\},1]$.
In this case, by \eqref{2048} and
Lemmas~\ref{discreteinequality} and~\ref{444},
we conclude that
\begin{align}\label{469}
&\sum_{k=0}^N\sum_{\tau\in I_k}\sum_{v=1}^{N(k,\tau)}
\mu(Q_\tau^{k,v})\left(\mathrm{I}+\mathrm{II}\right)^p
\nonumber\\
&\quad\lesssim\sum_{k=0}^N\sum_{\tau\in I_k}\sum_{v=1}^{N(k,\tau)}
\mu(Q_\tau^{k,v})\sum_{k'=0}^{N'}\sum_{\tau'\in I_{k'}}\sum_{v'=1}^{N(k',\tau')}
\left[\mu(Q_{\tau'}^{k',v'})\right]^p\left[m_{Q_{\tau'}^{k',v'}}(|\mathcal{P}_{k'}f|)\right]^p
\nonumber\\
&\qquad\times\inf_{y\in Q_{\tau'}^{k',v'}}\inf_{z\in Q_\tau^{k,v}}
\left(\frac{1}{(V_\rho)_{1}(y)
+V_\rho(z,y)}\left[\frac{1}{1+\rho(z,y)}\right]^{\varepsilon'}\right)^p
\nonumber\\
&\quad\lesssim\sum_{k'=0}^{N'}\sum_{\tau'\in I_{k'}}\sum_{v'=1}^{N(k',\tau')}
\left[\mu(Q_{\tau'}^{k',v'})\right]^p\left[m_{Q_{\tau'}^{k',v'}}(|\mathcal{P}_{k'}f|)\right]^p
\left[(V_\rho)_1(y_{\tau'}^{k',v'})\right]^{1-p}
\nonumber\\
&\quad\sim\sum_{k'=0}^{N'}\sum_{\tau'\in I_{k'}}\sum_{v'=1}^{N(k',\tau')}
\mu(Q_{\tau'}^{k',v'})\left[m_{Q_{\tau'}^{k',v'}}(|\mathcal{P}_{k'}f|)\right]^p,
\end{align}
where we chose $\varepsilon'\in(0,\varepsilon)$ such that
$p\in(\max\{\frac{n}{n+\varepsilon'},\,\frac{n}{n+\varepsilon'+s}\},1]$
and hence Lemma~\ref{444} is valid.

\emph{Case (2)} $p\in(1,\infty)$. In this case, from \eqref{2048},
H\"older's inequality, and Lemma~\ref{444} with
$p$ therein replaced by 1,
we infer that
\begin{align*}
\left(\mathrm{I}+\mathrm{II}\right)^p
&=\left\{\sum_{k'=0}^{N'}\sum_{\tau'\in I_{k'}}\sum_{v'=1}^{N(k',\tau')}
\left[\mu(Q_{\tau'}^{k',v'})\right]^\frac{1}{p}m_{Q_{\tau'}^{k',v'}}(|\mathcal{P}_{k'}f|)
\right.\nonumber\\
&\quad\times\left[\inf_{y\in Q_{\tau'}^{k',v'}}\inf_{z\in Q_\tau^{k,v}}
\frac{1}{(V_\rho)_{1}(y)
+V_\rho(z,y)}\left\{\frac{1}{1+\rho(z,y)}\right\}^{\varepsilon'}\right]^\frac{1}{p}
\nonumber\\
&\quad\left.\times\left[\mu(Q_{\tau'}^{k',v'})\right]^{1-\frac{1}{p}}
\left[\inf_{y\in Q_{\tau'}^{k',v'}}\inf_{z\in Q_\tau^{k,v}}
\frac{1}{(V_\rho)_{1}(y)
+V_\rho(z,y)}\left\{\frac{1}{1+\rho(z,y)}\right\}^{\varepsilon'}\right]^{1-\frac{1}{p}}\right\}^p
\nonumber\\
&\lesssim\sum_{k'=0}^{N'}\sum_{\tau'\in I_{k'}}\sum_{v'=1}^{N(k',\tau')}
\mu(Q_{\tau'}^{k',v'})\left[m_{Q_{\tau'}^{k',v'}}(|\mathcal{P}_{k'}f|)\right]^p
\nonumber\\
&\quad\times\inf_{y\in Q_{\tau'}^{k',v'}}\inf_{z\in Q_\tau^{k,v}}
\frac{1}{(V_\rho)_{1}(y)
+V_\rho(z,y)}\left[\frac{1}{1+\rho(z,y)}\right]^{\varepsilon'},
\end{align*}
which, together with Tonelli's theorem and Lemma~\ref{444} with $p$ therein replaced by 1,
implies that
\begin{align}\label{4610}
&\sum_{k=0}^N\sum_{\tau\in I_k}\sum_{v=1}^{N(k,\tau)}
\mu(Q_\tau^{k,v})\left(\mathrm{I}+\mathrm{II}\right)^p
\nonumber\\
&\quad\lesssim\sum_{k=0}^N\sum_{\tau\in I_k}\sum_{v=1}^{N(k,\tau)}
\mu(Q_\tau^{k,v})\sum_{k'=0}^{N'}\sum_{\tau'\in I_{k'}}\sum_{v'=1}^{N(k',\tau')}
\mu(Q_{\tau'}^{k',v'})\left[m_{Q_{\tau'}^{k',v'}}(|\mathcal{P}_{k'}f|)\right]^p
\nonumber\\
&\qquad\times\inf_{y\in Q_{\tau'}^{k',v'}}\inf_{z\in Q_\tau^{k,v}}
\frac{1}{(V_\rho)_{1}(y)
+V_\rho(z,y)}\left[\frac{1}{1+\rho(z,y)}\right]^{\varepsilon'}
\nonumber\\
&\quad\lesssim\sum_{k'=0}^{N'}\sum_{\tau'\in I_{k'}}\sum_{v'=1}^{N(k',\tau')}
\mu(Q_{\tau'}^{k',v'})\left[m_{Q_{\tau'}^{k',v'}}(|\mathcal{P}_{k'}f|)\right]^p.
\end{align}

Next, we turn to estimate $\mathrm{III}$.
By Lemma~\ref{1122}
[here, we need the cancellation of $\widetilde{P}_{k'}(x,y)$ for the first variable $x$]
and \eqref{1706},
we obtain, for any fixed $\varepsilon'\in(0,\varepsilon)$
and $y_{\tau'}^{k',v'}\in Q_{\tau'}^{k',v'}$
and for any
$k\in\{0,\ldots,N\}$ and $k'\in\{N'+1,N'+2,\ldots\}$,
\begin{align*}
&\sup_{z\in Q_\tau^{k,v}}\left|D_k\widetilde{P}_{k'}
(z,y_{\tau'}^{k',v'})\right|\\
&\quad\lesssim2^{-|k-k'|\varepsilon'}\sup_{z\in Q_\tau^{k,v}}
\frac{1}{(V_\rho)_{2^{-(k\land k')}}(y_{\tau'}^{k',v'})
+V_\rho(z,y_{\tau'}^{k',v'})}\left[\frac{2^{-(k\land k')}}
{2^{-(k\land k')}+\rho(z,y_{\tau'}^{k',v'})}\right]^{\varepsilon'}\\
&\quad\sim2^{-k'\varepsilon'}\sup_{z\in Q_\tau^{k,v}}
\frac{1}{(V_\rho)_{1}(y_{\tau'}^{k',v'})
+V_\rho(z,y_{\tau'}^{k',v'})}\left[\frac{1}{1+\rho(z,y_{\tau'}^{k',v'})}\right]^{\varepsilon'}\\
&\quad\sim2^{-k'\varepsilon'}\inf_{z\in Q_\tau^{k,v}}
\frac{1}{(V_\rho)_{1}(y_{\tau'}^{k',v'})
+V_\rho(z,y_{\tau'}^{k',v'})}\left[\frac{1}{1+\rho(z,y_{\tau'}^{k',v'})}\right]^{\varepsilon'},
\end{align*}
where, in the penultimate step, we used the facts that
$$
|k-k'|\sim k'
\ \ \text{and}\ \
2^{-(k\land k')}\sim1
$$
with the positive equivalence
constants depending on both $N$ and $N'$.
This implies that
\begin{align}\label{2211}
\mathrm{III}&\lesssim\sum_{k'=N'+1}^\infty2^{-k'\varepsilon'}\sum_{\tau'\in I_{k'}}
\sum_{v'=1}^{N(k',\tau')}\mu(Q_{\tau'}^{k',v'})
\left|\mathcal{P}_{k'}f(y_{\tau'}^{k',v'})\right|
\nonumber\\
&\quad\times\inf_{z\in Q_\tau^{k,v}}
\frac{1}{(V_\rho)_{1}(y_{\tau'}^{k',v'})
+V_\rho(z,y_{\tau'}^{k',v'})}\left[\frac{1}{1+\rho(z,y_{\tau'}^{k',v'})}\right]^{\varepsilon'}.
\end{align}
Note that, from \eqref{upperdoub} and Lemma~\ref{DyadicSys},
we deduce that, for any $y_{\tau'}^{k',v'}\in Q_{\tau'}^{k',v'}$,
$$
(V_\rho)_1(y_{\tau'}^{k',v'})
\lesssim2^{k'n}(V_\rho)_{2^{-k'}}(y_{\tau'}^{k',v'})
\sim2^{k'n}\mu(Q_{\tau'}^{k',v'}).
$$
By this and \eqref{2211},
we find that, for any fixed $r\in(\frac{n}{n+\varepsilon'+s},\min\{1,p,q\})$,
\begin{align*}
&\left[\sum_{k=0}^N\sum_{\tau\in I_k}\sum_{v=1}^{N(k,\tau)}
\mu(Q_\tau^{k,v})\left(\mathrm{III}\right)^p\right]^\frac{1}{p}
\\
&\quad\lesssim
\left\{\sum_{k=0}^N\sum_{\tau\in I_k}\sum_{v=1}^{N(k,\tau)}
\mu(Q_\tau^{k,v})\left[\sum_{k'=N'+1}^\infty2^{-k'\varepsilon'}\sum_{\tau'\in I_{k'}}
\sum_{v'=1}^{N(k',\tau')}\mu(Q_{\tau'}^{k',v'})
\left|\mathcal{P}_{k'}f(y_{\tau'}^{k',v'})\right|
\right.\right.\\
&\qquad\left.\left.\times\inf_{z\in Q_\tau^{k,v}}
\frac{1}{(V_\rho)_{1}(y_{\tau'}^{k',v'})
+V_\rho(z,y_{\tau'}^{k',v'})}\left\{\frac{1}{1+\rho(z,y_{\tau'}^{k',v'})}
\right\}^{\varepsilon'}\right]^p\right\}^\frac{1}{p}
\\
&\quad\lesssim\left\|\sum_{k'=N'+1}^\infty2^{-k'\varepsilon'}\sum_{\tau'\in I_{k'}}
\sum_{v'=1}^{N(k',\tau')}\mu(Q_{\tau'}^{k',v'})
\left|\mathcal{P}_{k'}f(y_{\tau'}^{k',v'})\right|
\right.\\
&\qquad\left.\times\frac{1}{(V_\rho)_{1}(y_{\tau'}^{k',v'})
+V_\rho(\cdot,y_{\tau'}^{k',v'})}\left[\frac{1}{1+
\rho(\cdot,y_{\tau'}^{k',v'})}\right]^{\varepsilon'}
\right\|_{L^p(X)},
\end{align*}
which, combined with Lemmas~\ref{L4.5}
and~\ref{Fefferman--Stein}, further implies that
\begin{align*}
&\left[\sum_{k=0}^N\sum_{\tau\in I_k}\sum_{v=1}^{N(k,\tau)}
\mu(Q_\tau^{k,v})\left(\mathrm{III}\right)^p\right]^\frac{1}{p}
\\
&\quad\lesssim\left\|\sum_{k'=N'+1}^\infty2^{-k'[\varepsilon'+s+n(1-\frac{1}{r})]}\left[
\mathcal{M}_\rho\left(\sum_{\tau'\in I_{k'}}\sum_{v'=1}^{N(k',\tau')}
2^{k'sr}\left|\mathcal{P}_{k'}f(y_{\tau'}^{k',v'})\right|^r\mathbf{1}_{Q_{\tau'}^{k',v'}}\right)
\right]^\frac{1}{r}\right\|_{L^p(X)}
\\
&\quad\lesssim\left\|\left\{\sum_{k'=N'+1}^\infty\left[
\mathcal{M}_\rho\left(\sum_{\tau'\in I_{k'}}\sum_{v'=1}^{N(k',\tau')}
2^{k'sr}\left|\mathcal{P}_{k'}f(y_{\tau'}^{k',v'})\right|^r\mathbf{1}_{Q_{\tau'}^{k',v'}}\right)
\right]^\frac{q}{r}\right\}^\frac{1}{q}
\right\|_{L^p(X)}
\\
&\quad=\left\|\left\{\sum_{k'=N'+1}^\infty\left[
\mathcal{M}_\rho\left(\sum_{\tau'\in I_{k'}}\sum_{v'=1}^{N(k',\tau')}
2^{k'sr}\left|\mathcal{P}_{k'}f(y_{\tau'}^{k',v'})\right|^r\mathbf{1}_{Q_{\tau'}^{k',v'}}\right)
\right]^\frac{q}{r}\right\}^\frac{r}{q}
\right\|_{L^\frac{p}{r}(X)}^\frac{1}{r}
\\
&\quad\lesssim\left\|\left[\sum_{k'=N'+1}^\infty
\sum_{\tau'\in I_{k'}}\sum_{v'=1}^{N(k',\tau')}
2^{k'sq}\left|\mathcal{P}_{k'}f(y_{\tau'}^{k',v'})\right|^q\mathbf{1}_{Q_{\tau'}^{k',v'}}
\right]^\frac{1}{q}
\right\|_{L^p(X)},
\end{align*}
where, in the antepenultimate step, we used
Lemma~\ref{discreteinequality} when $q\leq1$
or H\"older's inequality when $q>1$.
From this, the arbitrariness of
$y_{\tau'}^{k',v'}\in Q_{\tau'}^{k',v'}$, \eqref{946},
\eqref{469}, and \eqref{4610}, we infer that
\begin{align}\label{2253}
&\left\{\sum_{k=0}^N\sum_{\tau\in I_k}\sum_{v=1}^{N(k,\tau)}
\mu(Q_\tau^{k,v})\left[m_{Q_\tau^{k,v}}(|\mathcal{D}_kf|)\right]^p\right\}^\frac{1}{p}
\nonumber\\
&\quad\lesssim\left\{\sum_{k'=0}^{N'}\sum_{\tau'\in I_{k'}}\sum_{v'=1}^{N(k',\tau')}
\mu(Q_{\tau'}^{k',v'})\left[m_{Q_{\tau'}^{k',v'}}(|\mathcal{P}_{k'}f|)\right]^p\right\}^\frac{1}{p}
\nonumber\\
&\qquad+\left\|\left[\sum_{k'=N'+1}^\infty
\sum_{\tau'\in I_{k'}}\sum_{v'=1}^{N(k',\tau')}
2^{k'sq}\inf_{z\in Q_{\tau'}^{k',v'}}\left|\mathcal{P}_{k'}f(z)
\right|^q\mathbf{1}_{Q_{\tau'}^{k',v'}}
\right]^\frac{1}{q}
\right\|_{L^p(X)}.
\end{align}
This finishes the proof of the desired estimate of
the first term of the left-hand side of \eqref{2127}.

Now, we estimate the second term of the left-hand side of \eqref{2127}.
On the one hand, \eqref{462} and \eqref{464}
imply that, for any $k\in\{N+1,N+2,\ldots\}$
and $z\in X$,
\begin{align}\label{4615}
\left|\mathcal{D}_kf(z)\right|
&\leq\sum_{\tau'\in I_0}\sum_{v'=1}^{N(0,\tau')}
\int_{Q_{\tau'}^{0,v'}}\left|D_k\widetilde{P}_0(z,y)\right|\,d\mu(y)
m_{Q_{\tau'}^{0,v'}}(|\mathcal{P}_{0}f|)
\nonumber\\
&\quad+\sum_{k'=1}^{N'}\sum_{\tau'\in I_{k'}}\sum_{v'=1}^{N(k',\tau')}
\mu(Q_{\tau'}^{k',v'})\left|D_k\widetilde{P}_{k'}(z,y_{\tau'}^{k',v'})\right|
m_{Q_{\tau'}^{k',v'}}(|\mathcal{P}_{k'}f|)
\nonumber\\
&\quad+\sum_{k'=N'+1}^{\infty}\sum_{\tau'\in I_{k'}}\sum_{v'=1}^{N(k',\tau')}
\mu(Q_{\tau'}^{k',v'})\left|D_k\widetilde{P}_{k'}(z,y_{\tau'}^{k',v'})\right|
\left|\mathcal{P}_{k'}f(y_{\tau'}^{k',v'})\right|
\nonumber\\
&=:\mathrm{IV}+\mathrm{V}+\mathrm{VI}.
\end{align}
By Lemma~\ref{1122}, we conclude that,
for any fixed $\varepsilon'\in(0,\varepsilon)$ and for any
$k\in\{N+1,N+2,\ldots\}$, $k'\in\{0,\ldots,N'\}$, and $z,y\in X$,
\begin{align*}
\left|D_k\widetilde{P}_{k'}(z,y)\right|
&\lesssim2^{-|k-k'|\varepsilon'}
\frac{1}{(V_\rho)_{2^{-(k'\land k)}}(z)+V_\rho(z,y)}\left[\frac{2^{-(k'\land k)}}
{2^{-(k'\land k)}+\rho(z,y)}\right]^{\varepsilon'}\\
&\sim2^{-k\varepsilon'}
\frac{1}{(V_\rho)_1(z)+V_\rho(z,y)}\left[\frac{1}{1+\rho(z,y)}\right]^{\varepsilon'}
\end{align*}
since $|k-k'|\sim k$ and
$2^{-(k'\land k)}\sim1$
with the implicit positive constants depending on $N'$
but independent of both $k$ and $k'$,
which, together with an argument similar to that used in the estimation of \eqref{2048},
further implies that
\begin{align}\label{1531}
&\sup_{z\in B_\rho(z_\tau^{k,v},c_02^{-k})}\left(\mathrm{IV}+\mathrm{V}\right)
\nonumber\\
&\quad\lesssim2^{-k\varepsilon'}\sum_{k'=0}^{N'}\sum_{\tau'\in I_{k'}}
\sum_{v'=1}^{N(k',\tau')}\mu(Q_{\tau'}^{k',v'})m_{Q_{\tau'}^{k',v'}}(|\mathcal{P}_{k'}f|)
\nonumber\\
&\qquad\times\inf_{z\in Q_\tau^{k,v}}\inf_{y\in Q_{\tau'}^{k',v'}}
\frac{1}{(V_\rho)_1(z)+V_\rho(z,y)}\left[\frac{1}{1+\rho(z,y)}\right]^{\varepsilon'}.
\end{align}
From this and Lemma~\ref{discreteinequality} when $\frac{p}{q}\in(0,1]$
or H\"older's inequality when $\frac{p}{q}\in(1,\infty)$,
we deduce that, for any $\varepsilon'\in(0,\varepsilon)$,
\begin{align*}
&\left\|\left\{\sum_{k=N+1}^\infty\sum_{\tau\in I_k}
\sum_{v=1}^{N(k,\tau)}2^{ksq}\left[\sup_{z\in
B_\rho(z_\tau^{k,v},c_02^{-k})}\mathrm{IV}+\mathrm{V}\right]^q
\mathbf{1}_{Q_\tau^{k,v}}\right\}^\frac{1}{q}\right\|_{L^p(X)}^p\\
&\quad\lesssim\sum_{k=N+1}^\infty2^{k(s-\varepsilon')(p\land q)}
\sum_{\tau\in I_k}\sum_{v=1}^{N(k,\tau)}\mu(Q_\tau^{k,v})
\left\{\sum_{k'=0}^{N'}\sum_{\tau'\in I_{k'}}\sum_{v'=1}^{N(k',\tau')}
\mu(Q_{\tau'}^{k',v'})
\right.\\
&\qquad\left.\times m_{Q_{\tau'}^{k',v'}}(|\mathcal{P}_{k'}f|)
\inf_{z\in Q_\tau^{k,v}}\inf_{y\in Q_{\tau'}^{k',v'}}
\frac{1}{(V_\rho)_1(z)+V_\rho(z,y)}\left[\frac{1}{1+\rho(z,y)}\right]^{\varepsilon'}
\right\}^p\\
&\quad\leq\sum_{k=N+1}^\infty2^{k(s-\varepsilon')(p\land q)}
\sum_{\tau\in I_k}\sum_{v=1}^{N(k,\tau)}\mu(Q_\tau^{k,v})
\left\{\sum_{k'=0}^{N'}\sum_{\tau'\in I_{k'}}\sum_{v'=1}^{N(k',\tau')}
\mu(Q_{\tau'}^{k',v'})
\right.\\
&\qquad\left.\times m_{Q_{\tau'}^{k',v'}}(|\mathcal{P}_{k'}f|)
\frac{1}{(V_\rho)_1(y_\tau^{k,v})+V_\rho(y_\tau^{k,v},y_{\tau'}^{k',v'})}
\left[\frac{1}{1+\rho(y_\tau^{k,v},y_{\tau'}^{k',v'})}\right]^{\varepsilon'}
\right\}^p.
\end{align*}
Choose $\varepsilon'\in(|s|,\varepsilon)$ such
that $p\in(\frac{n}{n+\varepsilon'},\infty)$.
By this, Lemma~\ref{discreteinequality} when $p\leq1$
or H\"older's inequality when $p\in(1,\infty)$,
Lemma~\ref{444}, and the fact that
$(V_\rho)_1(y_{\tau'}^{k',v'})\sim\mu(Q_{\tau'}^{k',v'})$
with the positive equivalence constants independent of $k'$
but depending on $N'$,
we find that
\begin{align}\label{2251}
&\left\|\left\{\sum_{k=N+1}^\infty\sum_{\tau\in I_k}
\sum_{v=1}^{N(k,\tau)}2^{ksq}\left[\sup_{z\in
B_\rho(z_\tau^{k,v},c_02^{-k})}\mathrm{IV}+\mathrm{V}\right]^q
\mathbf{1}_{Q_\tau^{k,v}}\right\}^\frac{1}{q}\right\|_{L^p(X)}^p
\nonumber\\
&\quad\lesssim\sum_{k'=0}^{N'}\sum_{\tau'\in I_{k'}}\sum_{v'=1}^{N(k',\tau')}
\mu(Q_{\tau'}^{k',v'})\left[m_{Q_{\tau'}^{k',v'}}(|\mathcal{P}_{k'}f|)\right]^p.
\end{align}
On the other hand,
from an argument similar to that used in the estimation of \eqref{2145},
we infer that
\begin{align*}
&\left\|\left[\sum_{k=N+1}^\infty\sum_{\tau\in I_k}\sum_{v=1}^{N(k,\tau)}
2^{ksq}\sup_{z\in B_\rho(z_\tau^{k,v},c_02^{-k})}\left|\mathrm{VI}\right|^q
\mathbf{1}_{Q_\tau^{k,v}}\right]^\frac{1}{q}\right\|_{L^p(X)}
\nonumber\\
&\quad\lesssim\left\|\left[\sum_{k'=N'+1}^\infty
\sum_{\tau'\in I_{k'}}\sum_{v'=1}^{N(k',\tau')}
2^{k'sq}\left|\mathcal{P}_{k'}f(y_{\tau'}^{k',v'})\right|^q
\mathbf{1}_{Q_{\tau'}^{k',v'}}\right]^\frac{1}{q}\right\|_{L^p(X)},
\end{align*}
which, combined with both \eqref{2251} and \eqref{2253},
further implies \eqref{2127} and hence completes the proof of Lemma~\ref{ppin}.
\end{proof}

\begin{remark}\label{2205}
Let the notation be as in Lemma~\ref{ppin}.
\begin{enumerate}
\item[(i)]
Then, by Lemma~\ref{ppin},
we easily conclude that the definition of ${F}^{p,q}_s(X)$
is independent of the choice of
the approximation of the identity;
that is, for any
$f\in({{\mathcal{G}}}_0^\varepsilon(\beta,\gamma))'$,
\begin{align*}
&\left\{\sum_{k=0}^N\sum_{\tau\in I_{k}}
\sum_{v=1}^{N(k,\tau)}
\mu(Q_\tau^{k,v})
\left[m_{Q_\tau^{k,v}}\left(\left|\mathcal{D}_{k}f
\right|\right)\right]^p\right\}^\frac{1}{p}
+\left\|\left(\sum_{k=N+1}^\infty
2^{ksq}\left|\mathcal{D}_kf\right|^q\right)^\frac{1}{q}\right\|_{L^p(X)}\\
&\quad\sim\left\{\sum_{k=0}^{N'}\sum_{\tau\in I_{k}}
\sum_{v=1}^{N(k,\tau)}
\mu(Q_\tau^{k,v})
\left[m_{Q_\tau^{k,v}}\left(\left|\mathcal{P}_{k}f
\right|\right)\right]^p\right\}^\frac{1}{p}
+\left\|\left(\sum_{k=N'+1}^\infty
2^{ksq}\left|\mathcal{P}_kf\right|^q\right)^\frac{1}{q}\right\|_{L^p(X)},
\end{align*}
where the positive equivalence constants are independent of $f$.
A similar conclusion was also obtained in
\cite[Proposition~5.27(ii)]{hmy08}
on RD-spaces
and \cite[Proposition~4.4(ii)]{WHHY}
on spaces of homogeneous type.
\item[(ii)]
By an argument similar to that used in the proof of Lemma~\ref{ppin},
we obtain the following flexible estimate:
for any given $c_0\in[1,\infty)$ and for any
$f\in({{\mathcal{G}}}_0^\varepsilon(\beta,\gamma))'$,
\begin{align*}
&\left[\sum_{k=0}^N\sum_{\tau\in I_k}
\sum_{v=1}^{N(k,\tau)}
\mu(Q_\tau^{k,v})\sup_{z\in B_\rho(z_\tau^{k,v},c_02^{-k})}
\left|\mathcal{D}_kf(z)\right|^p\right]^\frac{1}{p}
\nonumber\\
&\qquad+\left\|\left[\sum_{k=N+1}^\infty\sum_{\tau\in I_k}\sum_{v=1}^{N(k,\tau)}
2^{ksq}\sup_{z\in B_\rho(z_\tau^{k,v},c_02^{-k})}\left|\mathcal{D}_kf(z)\right|^q
\mathbf{1}_{Q_\tau^{k,v}}\right]^\frac{1}{q}\right\|_{L^p(X)}
\nonumber\\
&\quad\lesssim\left\{\sum_{k=0}^{N'}\sum_{\tau\in I_k}
\sum_{v=1}^{N(k,\tau)}
\mu(Q_\tau^{k,v})
\left[m_{Q_\tau^{k,v}}\left(\left|\mathcal{P}_kf\right|\right)
\right]^p\right\}^\frac{1}{p}
\nonumber\\
&\qquad
+\left\|\left[\sum_{k=N'+1}^\infty\sum_{\tau\in I_k}\sum_{v=1}^{N(k,\tau)}
2^{ksq}\inf_{z\in Q_\tau^{k,v}}\left|\mathcal{P}_kf(z)\right|^q
\mathbf{1}_{Q_\tau^{k,v}}\right]^\frac{1}{q}\right\|_{L^p(X)}.
\end{align*}
\item[(iii)]
Recall that \cite[Definition~7.11]{HWYY} introduce the following version of the
\emph{inhomogeneous Littlewood--Paley
$g$-function} $G^s_q(f)$ \emph{of $f\in({{\mathcal{G}}}_0^\varepsilon(\beta,\gamma))'$}, which
is defined by setting, for any $x\in X$,
\begin{align*}
G^s_q(f)(x):&=\left[\sum_{k=0}^N\sum_{\tau\in I_k}
\sum_{v=1}^{N(k,\tau)}
m_{Q_\tau^{k,v}}\left(\left|\mathcal{D}_{k}f\right|^q\right)
\mathbf{1}_{Q_\tau^{k,v}}(x) +\sum_{k=N+1}^\infty2^{ksq}
\left|\mathcal{D}_kf(x)\right|^q\right]^\frac{1}{q}.
\end{align*}
If $q\in[1,\infty]$, then,
using (ii) of the present remark, we easily conclude that,
for any $f\in({{\mathcal{G}}}_0^\varepsilon(\beta,\gamma))'$,
\begin{align*}
\left\|G^s_q(f)\right\|_{L^p(X)}\sim\left\|g^s_q(f)\right\|_{L^p(X)}.
\end{align*}
\end{enumerate}
\end{remark}

Using Remark~\ref{2205}(ii)
and proceeding as in the proof of the homogeneous case,
properly speaking, the proofs of Theorems~\ref{Lusin} and~\ref{T--Lg*},
we obtain the following Lusin area function and
Littlewood--Paley $g_\lambda^*$-function
characterizations of $F^{p,q}_s(X)$;
we omit the details.

\begin{theorem}\label{1602}
Let $(X,\mathbf{q},\mu)$, $n$, $\varepsilon$, $s$, $p$, $\beta$, and $\gamma$ be
the same as in Lemma~\ref{ppin} and $q\in[1,\infty]$.
Then $f\in F^{p,q}_s(X)$
if and only if
$f\in({\mathcal{G}}_0^\varepsilon(\beta,\gamma))'$
and $\mathcal{S}^s_q(f)\in L^p(X)$;
moreover, for such $f$,
$$
\left\|\mathcal{S}^s_q(f)\right\|_{L^p(X)}
\sim\|f\|_{F^{p,q}_s(X)},
$$
where the positive equivalence constants are independent of $f$.
\end{theorem}

\begin{theorem}\label{T--Lg*-inhom}
Let $(X,\mathbf{q},\mu)$, $n$, $\varepsilon$, $s$, $p$, $\beta$, and $\gamma$ be
the same as in Lemma~\ref{ppin}. Assume that $q\in[1,\infty)$
and
$\lambda\in(\max\{n,\,\frac{qn}{p}\},\infty)$.
Then $f\in F^{p,q}_s(X)$
if and only if
$f\in({\mathcal{G}}_0^\varepsilon(\beta,\gamma))'$
and $(g^*_\lambda)^s_q(f)\in L^p(X)$;
moreover, for such $f$,
$$
\left\|(g^*_\lambda)^s_q(f)\right\|_{L^p(X)}
\sim\|f\|_{F^{p,q}_s(X)},
$$
where the positive equivalence constants are independent of $f$.
\end{theorem}

\begin{remark}
\begin{enumerate}
\item
Let $(X,\rho,\mu)$ be a space of homogeneous type
with upper dimension $n\in(0,\infty)$.
Let $\eta\in(0,1)$ be the H\"older regularity of the
wavelets constructed in \cite[Theorem~7.1]{AH1}.
Assume that $s\in(-\eta,\eta)$, $p\in(\max\{\frac{n}{n+\eta},\,\frac{n}{n+\eta+s}\},\infty)$,
$q\in(\max\{\frac{n}{n+\eta},\,\frac{n}{n+\eta+s}\},\infty]$,
and $\lambda\in(\max\{n,\,\frac{qn}{p}\},\infty)$.
He et al. \cite[Theorem 7.13]{HWYY} showed that
$f\in F^{p,q}_s(X)$
if and only if
$f\in({\mathcal{G}}_0^\eta(\beta,\gamma))'$
and $(g^*_\lambda)^s_q(f)\in L^p(X)$; moreover, for such $f$,
$$
\left\|(g^*_\lambda)^s_q(f)\right\|_{L^p(X)}
\sim\|f\|_{F^{p,q}_s(X)}
$$
with the positive equivalence constants independent of $f$,
where $F^{p,q}_s(X)$ is the inhomogeneous
Triebel--Lizorkin space in \cite[Definition 4.1(ii)]{WHHY}.
However, in more specialized setting of
ultra-RD-spaces,
Theorem~\ref{T--Lg*-inhom} gives the Littlewood--Paley
$g_\lambda^*$-function characterization of inhomogeneous
Triebel--Lizorkin spaces
for a larger range of
\begin{align*}
s\in\left(-\mathrm{ind\,}(X,\mathbf{q}),\mathrm{ind\,}(X,\mathbf{q})\right)
\end{align*}
(and hence a larger range of both $p$ and $q$)
with $\mathrm{ind\,}(X,\mathbf{q})$ as in \eqref{index};
see Remark~\ref{22013}(i).
\item
Let $(X,\mathbf{q},\mu)$ be an Ahlfors-regular
quasi-ultrametric space as in \cite{AlMi2015}.
Then, in this setting,
Theorems~\ref{1602} and~\ref{T--Lg*-inhom} give the
Lusin area function characterization and the
Littlewood--Paley $g_\lambda^*$-function characterization
of the inhomogeneous Triebel--Lizorkin space ${F}^{p,q}_s(X)$
in \cite[Definition~9.6]{AlMi2015}.
\end{enumerate}
\end{remark}

In order to establish Littlewood--Paley characterizations of
$F^{\infty,q}_s(X)$, we need the following inhomogeneous Plancherel--P\'olya
inequality in the case $p=\infty$.

\begin{proposition}
\label{proin}
Let $(X,\mathbf{q},\mu)$ be an ultra-RD-space with upper dimension
$n\in(0,\infty)$, where $\mu$ is assumed to be
a Borel-semiregular measure on $X$.
Let $\varepsilon$, $s$, $q$, $\beta$, and $\gamma$ be
the same as in Definition~\ref{iTL}.
Assume that $\{\mathcal{S}_{2^{-k}}\}_{k\in\mathbb{Z}_+}$
and $\{\mathcal{Q}_{2^{-k}}\}_{k\in\mathbb{Z}_+}$
are two inhomogeneous approximations of the identity
of order $\varepsilon$ as in Definition~\ref{DefIATI}. For
any $k\in\mathbb{N}$, let
$\mathcal{D}_k:=\mathcal{S}_{2^{-k}}-\mathcal{S}_{2^{-k+1}}$,
$\mathcal{P}_k:=\mathcal{Q}_{2^{-k}}-\mathcal{Q}_{2^{-k+1}}$,
$\mathcal{D}_0:=\mathcal{S}_1$, and $\mathcal{P}_0:=\mathcal{Q}_1$.
Let $N,N'\in\mathbb{N}$ be the same as in Theorem~\ref{ID}
associated, respectively, with $\{\mathcal{D}_k\}_{k\in\mathbb{Z_+}}$
and $\{\mathcal{P}_k\}_{k\in\mathbb{Z_+}}$.
Let
$\{Q^{l}_\alpha\}_{l\in\mathbb{Z},\,\alpha\in I_l}$ be the collection of all
dyadic cubes given in Lemma~\ref{DyadicSys}.
Then, for any $c_0\in[1,\infty)$
and $f\in(\mathring{{\mathcal{G}}}_0^\varepsilon(\beta,\gamma))'$,
\begin{align}\label{1528}
&\max\Bigg\{\sup_{k\in\{0,\ldots,N\}}
\sup_{\tau\in I_k}\sup_{v\in\{1,\ldots,N(k,\tau)\}}m_{Q_\tau^{k,v}}
\left(\left|\mathcal{D}_kf\right|\right),\,
\nonumber\\
&\qquad\sup_{l\in\mathbb{N},\,l>N}
\sup_{\alpha\in I_l}\left[
\frac{1}{\mu(Q^l_\alpha)}\sum_{k=l}^\infty\sum_{\tau\in I_k}
\sum_{v=1}^{N(k,\tau)}
2^{ksq}\mu(Q_\tau^{k,v})
\right.\nonumber\\
&\qquad\left.\times\mathbf{1}_{\{(\tau,v):Q_\tau^{k,v}\subset Q^l_\alpha\}}(\tau,v)
\sup_{x\in B_\rho(z^{k,v}_\tau,c_02^{-k})}
\left|\mathcal{D}_kf(x)\right|^q\right]^\frac{1}{q}\Bigg\}
\nonumber\\
&\quad\lesssim\max\Bigg\{\sup_{k\in\{0,\ldots,N'\}}
\sup_{\tau\in I_k}\sup_{v\in\{1,\ldots,N(k,\tau)\}}m_{Q_\tau^{k,v}}
\left(\left|\mathcal{P}_kf\right|\right),\,
\nonumber\\
&\qquad\sup_{l\in\mathbb{N},\,l>N'}
\sup_{\alpha\in I_l}\left[
\frac{1}{\mu(Q^l_\alpha)}\sum_{k=l}^\infty\sum_{\tau\in I_k}
\sum_{v=1}^{N(k,\tau)}
2^{ksq}\mu(Q_\tau^{k,v})
\right.\nonumber\\
&\qquad\left.\times\mathbf{1}_{\{(\tau,v):Q_\tau^{k,v}\subset Q^l_\alpha\}}(\tau,v)
\inf_{x\in Q^{k,v}_\tau}\left|\mathcal{P}_kf(x)\right|^q\right]^\frac{1}{q}\Bigg\}
\end{align}
with the usual modification when $q=\infty$,
where the implicit positive constant is independent of $f$
and $\{Q^{k,v}_\tau\}_{k\in\mathbb{Z},\,\tau\in I_k,\,v=1,\dots,N(k,\tau)}$
(with centers $\{z^{k,v}_\tau\}_{k\in\mathbb{Z},\,\tau\in I_k,\,v=1,\dots,N(k,\tau)}$)
is the collection of dyadic cubes given in Remark~\ref{jcubes}.
\end{proposition}

\begin{proof}
We only prove the present proposition in the case $q<\infty$ because the proof
in the case $q=\infty$ is similar and hence we omit the details.
Let $f\in(\mathring{{\mathcal{G}}}_0^\varepsilon(\beta,\gamma))'$.
From \eqref{4615}, we infer that, for any $k\in\mathbb{Z_+}$ znd $z\in X$,
\begin{align}\label{ut2149}
\left|\mathcal{D}_kf(z)\right|
&\leq\sum_{\tau'\in I_0}\sum_{v'=1}^{N(0,\tau')}
\int_{Q_{\tau'}^{0,v'}}\left|D_k\widetilde{P}_0(z,y)\right|\,d\mu(y)
m_{Q_{\tau'}^{0,v'}}(|\mathcal{P}_{0}f|)
\nonumber\\
&\quad+\sum_{k'=1}^{N'}\sum_{\tau'\in I_{k'}}\sum_{v'=1}^{N(k',\tau')}
\mu(Q_{\tau'}^{k',v'})\left|D_k\widetilde{P}_{k'}(z,y_{\tau'}^{k',v'})\right|
m_{Q_{\tau'}^{k',v'}}(|\mathcal{P}_{k'}f|)
\nonumber\\
&\quad+\sum_{k'=N'+1}^{\infty}\sum_{\tau'\in I_{k'}}\sum_{v'=1}^{N(k',\tau')}
\mu(Q_{\tau'}^{k',v'})\left|D_k\widetilde{P}_{k'}(z,y_{\tau'}^{k',v'})\right|
\left|\mathcal{P}_{k'}f(y_{\tau'}^{k',v'})\right|
\nonumber\\
&=:\mathrm{I}+\mathrm{II}+\mathrm{III}.
\end{align}
Then, to estimate the left-hand side of \eqref{1528},
we consider the following two cases for $k\in\mathbb{Z_+}$.

\emph{Case (1)}
$k\in\{0,\ldots,N\}$. In this case, by \eqref{2048},
we find that,
for any fixed $\varepsilon'\in(0,\varepsilon)$ and
for any $\tau\in I_k$ and $v\in\{1,\dots,N(k,\tau)\}$,
\begin{align*}
m_{Q_\tau^{k,v}}\left(\mathrm{I}+\mathrm{II}\right)
&\lesssim\sum_{k'=0}^{N'}\sum_{\tau'\in I_{k'}}
\sum_{v'=1}^{N(k',\tau')}\mu(Q_{\tau'}^{k',v'})
m_{Q_{\tau'}^{k',v'}}(|\mathcal{P}_{k'}f|)
\nonumber\\
&\quad\times\inf_{y\in Q_{\tau'}^{k',v'}}
\frac{1}{(V_\rho)_1(z_\tau^{k,v})+
V_\rho(z_\tau^{k,v},y)}\left[\frac{1}{1+\rho(z_\tau^{k,v},y)}\right]^{\varepsilon'},
\end{align*}
which, together with Lemmas~\ref{DyadicSys}(i)
and~\ref{A3.2}(ii), further implies that
\begin{align}\label{4618}
m_{Q_\tau^{k,v}}\left(\mathrm{I}+\mathrm{II}\right)
&\lesssim\sup_{k'\in\{0,\ldots,N'\}}
\sup_{\tau'\in I_{k'}}\sup_{v'\in\{1,\ldots,N(k',\tau')\}}
m_{Q_{\tau'}^{k',v'}}(|\mathcal{P}_{k'}f|)
\nonumber\\
&\quad\times\int_X\frac{1}{(V_\rho)_1(z_\tau^{k,v})+
V_\rho(z_\tau^{k,v},y)}\left[\frac{1}{1+
\rho(z_\tau^{k,v},y)}\right]^{\varepsilon'}\,d\mu(y)
\nonumber\\
&\lesssim\sup_{k'\in\{0,\ldots,N'\}}
\sup_{\tau'\in I_{k'}}\sup_{v'\in\{1,\ldots,N(k',\tau')\}}
m_{Q_{\tau'}^{k',v'}}(|\mathcal{P}_{k'}f|).
\end{align}
Next, we turn to estimate $m_{Q_\tau^{k,v}}(\mathrm{III})$
in the case $k\in\{0,\ldots,N\}$.
To this end, we consider the following two cases for $q$.

\emph{Subcase (1)} $q\in(\max\{\frac{n}{n+\varepsilon},\,
\frac{n}{n+\varepsilon+s}\},1]$. In this case,
from \eqref{2211} and Lemma~\ref{discreteinequality} with
$r$ therein replaced by $q$,
we deduce that, for any fixed $\varepsilon'\in((-s)_+,\varepsilon)$,
\begin{align*}
m_{Q_\tau^{k,v}}(\mathrm{III})
&\lesssim
\sum_{k'=N'+1}^\infty2^{-k'\varepsilon'}\sum_{\tau'\in I_{k'}}
\sum_{v'=1}^{N(k',\tau')}\mu(Q_{\tau'}^{k',v'})
\left|\mathcal{P}_{k'}f(y_{\tau'}^{k',v'})\right|
\nonumber\\
&\quad\times\inf_{z\in Q_\tau^{k,v}}
\frac{1}{(V_\rho)_{1}(y_{\tau'}^{k',v'})
+V_\rho(z,y_{\tau'}^{k',v'})}\left[\frac{1}{1+
\rho(z,y_{\tau'}^{k',v'})}\right]^{\varepsilon'}
\nonumber\\
&\leq
\sum_{k'=N'+1}^\infty2^{-k'\varepsilon'-k's}\left\{\sum_{\tau'\in I_{k'}}
\sum_{v'=1}^{N(k',\tau')}
2^{k'sq}\left|\mathcal{P}_{k'}f(y_{\tau'}^{k',v'})\right|^q
\right.\nonumber\\
&\quad\left.\times\inf_{z\in Q_\tau^{k,v}}\left[
\frac{\mu(Q_{\tau'}^{k',v'})}{(V_\rho)_{1}(y_{\tau'}^{k',v'})
+V_\rho(z,y_{\tau'}^{k',v'})}\left\{\frac{1}{1+\rho(z,y_{\tau'}^{k',v'})}
\right\}^{\varepsilon'}\right]^q\right\}^\frac{1}{q},
\end{align*}
which, combined with the fact that $\mu(Q_{\tau'}^{k',v'})\sim
(V_\rho)_{2^{-k'}}(y_{\tau'}^{k',v'})
\leq(V_\rho)_1(y_{\tau'}^{k',v'})$ and an argument similar to
that used in the estimation of \eqref{1601},
further implies that
\begin{align}\label{1506}
m_{Q_\tau^{k,v}}(\mathrm{III})
&\lesssim\sum_{k'=N'+1}^\infty2^{-k'(\varepsilon'+s)}
\left[\sum_{\tau'\in I_{k'}}\sum_{v'=1}^{N(k',\tau')}
2^{k'sq}\left|\mathcal{P}_{k'}f(y_{\tau'}^{k',v'})\right|^q\right]^\frac{1}{q}
\nonumber\\
&\lesssim\sup_{l'\in\mathbb{N},\,l'>N'}
\sup_{\alpha'\in I_{l'}}
\left\{\frac{1}{\mu(Q_{\alpha'}^{l'})}
\sum_{k'=l'}^\infty\sum_{\tau'\in I_{k'}}\sum_{v'=1}^{N(k',\tau')}
2^{k'sq}\mu(Q_{\tau'}^{k',v'})
\right.\nonumber\\
&\quad\left.\times\mathbf{1}_{\{(\tau',v'):Q_{\tau'}^{k',v'}
\subset Q_{\alpha'}^{l'}\}}(\tau',v')
\left|\mathcal{P}_{k'}f(y_{\tau'}^{k',v'})\right|^q\right\}^\frac{1}{q}.
\end{align}

\emph{Subcase (2)} $q\in(1,\infty)$. In this case,
by \eqref{2211} and H\"older's inequality,
we conclude that, for any fixed $\varepsilon'\in((-s)_+,\varepsilon)$,
\begin{align*}
m_{Q_\tau^{k,v}}(\mathrm{III})
&\lesssim
\sum_{k'=N'+1}^\infty2^{-k'(\varepsilon'+s)}\sum_{\tau'\in I_{k'}}
\sum_{v'=1}^{N(k',\tau')}\mu(Q_{\tau'}^{k',v'})2^{k's}
\left|\mathcal{P}_{k'}f(y_{\tau'}^{k',v'})\right|
\nonumber\\
&\quad\times\inf_{z\in Q_\tau^{k,v}}
\inf_{y\in Q_{\tau'}^{k',v'}}
\frac{1}{(V_\rho)_{1}(y)
+V_\rho(z,y)}\left[\frac{1}{1+\rho(z,y)}\right]^{\varepsilon'}
\nonumber\\
&\leq
\sum_{k'=N'+1}^\infty2^{-k'(\varepsilon'+s)}\left[\sum_{\tau'\in I_{k'}}
\sum_{v'=1}^{N(k',\tau')}2^{k'sq}
\left|\mathcal{P}_{k'}f(y_{\tau'}^{k',v'})\right|^q\right]^\frac{1}{q}
\nonumber\\
&\quad\times\left[\sum_{\tau'\in I_{k'}}
\sum_{v'=1}^{N(k',\tau')}\left\{
\inf_{y\in Q_{\tau'}^{k',v'}}
\frac{\mu(Q_{\tau'}^{k',v'})}{(V_\rho)_{1}(y)
+V_\rho(y_\tau^{k,v},y)}\left[\frac{1}{1+\rho(y_\tau^{k,v},y)}
\right]^{\varepsilon'}\right\}^{q'}\right]^\frac{1}{q'}.
\end{align*}
From this, the fact that $\mu(Q_{\tau'}^{k',v'})\sim
(V_\rho)_{2^{-k'}}(y)
\leq(V_\rho)_1(y)$ for any $y\in Q_{\tau'}^{k',v'}$,
and $q'>1$,
it follows that
\begin{align*}
m_{Q_\tau^{k,v}}(\mathrm{III})
&\lesssim
\sum_{k'=N'+1}^\infty2^{-k'(\varepsilon'+s)}\left[\sum_{\tau'\in I_{k'}}
\sum_{v'=1}^{N(k',\tau')}2^{k'sq}
\left|\mathcal{P}_{k'}f(y_{\tau'}^{k',v'})\right|^q\right]^\frac{1}{q}
\nonumber\\
&\quad\times\left\{\sum_{\tau'\in I_{k'}}
\sum_{v'=1}^{N(k',\tau')}
\inf_{y\in Q_{\tau'}^{k',v'}}
\frac{\mu(Q_{\tau'}^{k',v'})}{(V_\rho)_{1}(y)
+V_\rho(y_\tau^{k,v},y)}\left[\frac{1}{1+\rho(y_\tau^{k,v},y)}
\right]^{\varepsilon'}\right\}^\frac{1}{q'}.
\end{align*}
By this, Lemma~\ref{A3.2}(ii),
and an argument similar to that used in the estimation of \eqref{1601},
we find that
\begin{align*}
m_{Q_\tau^{k,v}}(\mathrm{III})
&\lesssim\sum_{k'=N'+1}^\infty2^{-k'(\varepsilon'+s)}
\left[\sum_{\tau'\in I_{k'}}\sum_{v'=1}^{N(k',\tau')}
2^{k'sq}\left|\mathcal{P}_{k'}f(y_{\tau'}^{k',v'})\right|^q\right]^\frac{1}{q}
\\
&\quad\times\left\{\int_X\frac{1}{(V_\rho)_1(y)+V_\rho(z_\tau^{k,v},y)}
\left[\frac{1}{1+\rho(z_\tau^{k,v},y)}\right]^{\varepsilon'}
\,d\mu(y)\right\}^\frac{1}{q'}\\
&\lesssim\sup_{l'\in\mathbb{N},\,l'>N'}
\sup_{\alpha'\in I_{l'}}
\left\{\frac{1}{\mu(Q_{\alpha'}^{l'})}
\sum_{k'=l'}^\infty\sum_{\tau'\in I_{k'}}\sum_{v'=1}^{N(k',\tau')}
2^{k'sq}\mu(Q_{\tau'}^{k',v'})
\right.\\
&\quad\left.\times\mathbf{1}_{\{(\tau',v'):Q_{\tau'}^{k',v'}
\subset Q_{\alpha'}^{l'}\}}(\tau',v')
\left|\mathcal{P}_{k'}f(y_{\tau'}^{k',v'})\right|^q\right\}^\frac{1}{q},
\end{align*}
which, together with \eqref{1506}
and the arbitrariness of $y_\tau^{k,v}\in Q_\tau^{k,v}$,
further implies that,
for any $q\in(\max\{\frac{n}{n+\varepsilon},\,\frac{n}{n+\varepsilon+s}\},\infty)$,
\begin{align}\label{1110}
m_{Q_\tau^{k,v}}(\mathrm{III})
&\lesssim\sup_{l'\in\mathbb{N},\,l'>N'}
\sup_{\alpha'\in I_{l'}}
\left\{\frac{1}{\mu(Q_{\alpha'}^{l'})}
\sum_{k'=l'}^\infty\sum_{\tau'\in I_{k'}}\sum_{v'=1}^{N(k',\tau')}
2^{k'sq}\mu(Q_{\tau'}^{k',v'})
\right.\nonumber\\
&\quad\left.\times\mathbf{1}_{\{(\tau',v'):Q_{\tau'}^{k',v'}
\subset Q_{\alpha'}^{l'}\}}(\tau',v')
\left|\mathcal{P}_{k'}f(y_{\tau'}^{k',v'})\right|^q\right\}^\frac{1}{q}.
\end{align}
This finishes the
estimation of $m_{Q_\tau^{k,v}}(\mathrm{III})$
in the case $k\in\{0,\ldots,N\}$.

From \eqref{ut2149}, \eqref{4618}, \eqref{1110},
and taking the supremum over all $k\in\{0,\dots,N\}$,
$\tau\in I_k$, and $v\in\{1,\dots,N(k,\tau)\}$,
we deduce that
\begin{align}\label{2019}
&\sup_{k\in\{0,\ldots,N\}}
\sup_{\tau\in I_k}\sup_{v\in\{1,\ldots,N(k,\tau)\}}m_{Q_\tau^{k,v}}
\left(\left|\mathcal{D}_kf\right|\right)
\nonumber\\
&\quad\leq\sup_{k\in\{0,\ldots,N\}}
\sup_{\tau\in I_k}\sup_{v\in\{1,\ldots,N(k,\tau)\}}\left[m_{Q_\tau^{k,v}}
\left(\mathrm{I}+\mathrm{II}\right)
+m_{Q_\tau^{k,v}}(\mathrm{III})\right]
\nonumber\\
&\quad\lesssim\sup_{k'\in\{0,\ldots,N'\}}
\sup_{\tau'\in I_{k'}}\sup_{v'\in\{1,\ldots,N(k',\tau')\}}
m_{Q_{\tau'}^{k',v'}}(|\mathcal{P}_{k'}f|)
\nonumber\\
&\qquad+\sup_{l'\in\mathbb{N},\,l'>N'}
\sup_{\alpha'\in I_{l'}}
\left\{\frac{1}{\mu(Q_{\alpha'}^{l'})}
\sum_{k'=l'}^\infty\sum_{\tau'\in I_{k'}}\sum_{v'=1}^{N(k',\tau')}
2^{k'sq}\mu(Q_{\tau'}^{k',v'})
\right.\nonumber\\
&\qquad\left.\times\mathbf{1}_{\{(\tau',v'):Q_{\tau'}^{k',v'}
\subset Q_{\alpha'}^{l'}\}}(\tau',v')
\inf_{x\in Q_{\tau'}^{k',v'}}\left|\mathcal{P}_{k'}f(x)\right|^q\right\}^\frac{1}{q},
\end{align}
which is the desired estimate for the
first term in the left-hand side of \eqref{1528}.

\emph{Case (2)} $k\in\{N+1,N+2,\ldots\}$. In this case,
we first estimate $\sup_{z\in Q_\tau^{k,v}}\mathrm{I}+\mathrm{II}$.
By \eqref{1531},
we conclude that, for any fixed $\varepsilon'\in(0,\varepsilon)$,
\begin{align*}
\sup_{z\in B_\rho(z_\tau^{k,v},c_02^{-k})}\left(\mathrm{I}+\mathrm{II}\right)
&\lesssim2^{-k\varepsilon'}\sum_{k'=0}^{N'}
\sum_{\tau'\in I_{k'}}\sum_{v'=1}^{N(k',v')}
\mu(Q_{\tau'}^{k',v'})m_{Q_{\tau'}^{k',v'}}(|\mathcal{P}_{k'}f|)\\
&\quad\times\inf_{z\in Q_\tau^{k,v}}\inf_{y\in Q_{\tau'}^{k',v'}}
\frac{1}{(V_\rho)_1(z)+V_\rho(z,y)}\left[
\frac{1}{1+\rho(z,y)}\right]^{\varepsilon'}\\
&\leq2^{-k\varepsilon'}\sup_{k'\in\{0,\ldots,N'\}}
\sup_{\tau'\in I_{k'}}\sup_{v'\in\{1,\ldots,N(k',\tau')\}}
m_{Q_{\tau'}^{k',v'}}(|\mathcal{P}_{k'}f|)\\
&\quad\times\sum_{k'=0}^{N'}\sum_{\tau'\in I_{k'}}
\sum_{v'=1}^{N(k',\tau')}\mu(Q_{\tau'}^{k',v'})\\
&\quad\times\inf_{y\in Q_{\tau'}^{k',v'}}
\frac{1}{(V_\rho)_1(y_\tau^{k,v})+V_\rho(y_\tau^{k,v},y)}\left[
\frac{1}{1+\rho(y_\tau^{k,v},y)}\right]^{\varepsilon'}.
\end{align*}
This and Lemma~\ref{A3.2}(ii) imply that
\begin{align*}
\sup_{z\in B_\rho(z_\tau^{k,v},c_02^{-k})}\left(\mathrm{I}+\mathrm{II}\right)
&\lesssim(N'+1)2^{-k\varepsilon'}\sup_{k'\in\{0,\ldots,N'\}}
\sup_{\tau'\in I_{k'}}\sup_{v'\in\{1,\ldots,N(k',\tau')\}}
m_{Q_{\tau'}^{k',v'}}(|\mathcal{P}_{k'}f|)\\
&\quad\times\int_X
\frac{1}{(V_\rho)_1(y_\tau^{k,v})+V_\rho(y_\tau^{k,v},y)}\left[
\frac{1}{1+\rho(y_\tau^{k,v},y)}\right]^{\varepsilon'}\,d\mu(y)\\
&\lesssim2^{-k\varepsilon'}\sup_{k'\in\{0,\ldots,N'\}}
\sup_{\tau'\in I_{k'}}\sup_{v'\in\{1,\ldots,N(k',\tau')\}}
m_{Q_{\tau'}^{k',v'}}(|\mathcal{P}_{k'}f|),
\end{align*}
which further implies that,
for any fixed $\varepsilon'\in(s_+,\varepsilon)$,
\begin{align}\label{153}
&\sup_{l\in\mathbb{N},l>N}\sup_{\alpha\in I_l}
\left\{\frac{1}{\mu(Q_\alpha^l)}\sum_{k=l}^\infty
\sum_{\tau\in I_k}\sum_{v=1}^{N(k,\tau)}
2^{ksq}\mu(Q_\tau^{k,v})
\right.\nonumber\\
&\qquad\left.\times\mathbf{1}_{\{(\tau,v):Q_\tau^{k,v}
\subset Q_\alpha^l\}}(\tau,v)
\sup_{x\in B_\rho(z_\tau^{k,v},c_02^{-k})}\left(\mathrm{I}+
\mathrm{II}\right)^q\right\}^\frac{1}{q}
\nonumber\\
&\quad\lesssim\sup_{k'\in\{0,\ldots,N'\}}
\sup_{\tau'\in I_{k'}}\sup_{v'\in\{1,\ldots,N(k',\tau')\}}
m_{Q_{\tau'}^{k',v'}}(|\mathcal{P}_{k'}f|)
\nonumber\\
&\qquad\times\sup_{l\in\mathbb{N},l>N}
\sup_{\alpha\in I_l}\left[\sum_{k=l}^\infty
2^{kq(s-\varepsilon')}
\sum_{\tau\in I_k}
\sum_{v=1}^{N(k,\tau)}\frac{\mu(Q_\tau^{k,v})}{\mu(Q_\alpha^l)}
\mathbf{1}_{\{(\tau,v):Q_\tau^{k,v}\subset Q_\alpha^l\}}(\tau,v)\right]^\frac{1}{q}
\nonumber\\
&\quad\lesssim\sup_{k'\in\{0,\ldots,N'\}}
\sup_{\tau'\in I_{k'}}\sup_{v'\in\{1,\ldots,N(k',\tau')\}}
m_{Q_{\tau'}^{k',v'}}(|\mathcal{P}_{k'}f|).
\end{align}
To estimate $\mathrm{III}$,
from an argument similar to that used in the proof of
Proposition~\ref{ppineq-infty}, we infer that
\begin{align}\label{154}
&\sup_{l\in\mathbb{N},l>N}\sup_{\alpha\in I_l}
\left\{\frac{1}{\mu(Q_\alpha^l)}\sum_{k=l}^\infty
\sum_{\tau\in I_k}\sum_{v=1}^{N(k,\tau)}
2^{ksq}\mu(Q_\tau^{k,v})\right.\nonumber\\
&\qquad\left.\times\mathbf{1}_{(\tau,v):
Q_\tau^{k,v}\subset Q_\alpha^l}(\tau,v)
\sup_{x\in B_\rho(z_\tau^{k,v},c_02^{-k})}\mathrm{III}^q\right\}^\frac{1}{q}
\nonumber\\
&\quad\lesssim\sup_{l'\in\mathbb{N},\,l'>N'}
\sup_{\alpha'\in I_{l'}}
\left\{\frac{1}{\mu(Q_{\alpha'}^{l'})}
\sum_{k'=l'}^\infty\sum_{\tau'\in I_{k'}}\sum_{v'=1}^{N(k',\tau')}
2^{k'sq}\mu(Q_{\tau'}^{k',v'})
\right.\nonumber\\
&\qquad\left.\times\mathbf{1}_{\{(\tau',v'):Q_{\tau'}^{k',v'}\subset Q_{\alpha'}^{l'}\}}(\tau',v')
\inf_{x\in Q_{\tau'}^{k',v'}}\left|\mathcal{P}_{k'}f(x)\right|^q\right\}^\frac{1}{q}.
\end{align}
Consequently,
\begin{align*}
&\sup_{l\in\mathbb{N},\,l>N}
\sup_{\alpha\in I_l}\left[
\frac{1}{\mu(Q^l_\alpha)}\sum_{k=l}^\infty\sum_{\tau\in I_k}
\sum_{v=1}^{N(k,\tau)}
2^{ksq}\mu(Q_\tau^{k,v})
\right.\\
&\qquad\left.\times\mathbf{1}_{\{(\tau,v):Q_\tau^{k,v}\subset Q^l_\alpha\}}(\tau,v)
\sup_{x\in B_\rho(z^{k,v}_\tau,c_02^{-k})}\left|\mathcal{D}_kf(x)\right|^q\right]^\frac{1}{q}\\
&\quad\lesssim
\sup_{l\in\mathbb{N},l>N}\sup_{\alpha\in I_l}
\left\{\frac{1}{\mu(Q_\alpha^l)}\sum_{k=l}^\infty
\sum_{\tau\in I_k}\sum_{v=1}^{N(k,\tau)}
2^{ksq}\mu(Q_\tau^{k,v})
\right.\nonumber\\
&\qquad\left.\times\mathbf{1}_{\{(\tau,v):Q_\tau^{k,v}\subset Q_\alpha^l\}}(\tau,v)
\sup_{x\in Q_\tau^{k,v}}\left(\mathrm{I}+\mathrm{II}\right)^q\right\}^\frac{1}{q}
\nonumber\\
&\qquad+\sup_{l\in\mathbb{N},l>N}\sup_{\alpha\in I_l}
\left\{\frac{1}{\mu(Q_\alpha^l)}\sum_{k=l}^\infty
\sum_{\tau\in I_k}\sum_{v=1}^{N(k,\tau)}
2^{ksq}\mu(Q_\tau^{k,v})\mathbf{1}_{\{(\tau,v):
Q_\tau^{k,v}\subset Q_\alpha^l\}}(\tau,v)
\sup_{x\in Q_\tau^{k,v}}\mathrm{III}^q\right\}^\frac{1}{q},
\end{align*}
which, combined with \eqref{153} and \eqref{154},
further implies that
\begin{align*}
&\sup_{l\in\mathbb{N},\,l>N}
\sup_{\alpha\in I_l}\left[
\frac{1}{\mu(Q^l_\alpha)}\sum_{k=l}^\infty\sum_{\tau\in I_k}
\sum_{v=1}^{N(k,\tau)}
2^{ksq}\mu(Q_\tau^{k,v})
\right.\\
&\qquad\left.\times\mathbf{1}_{\{(\tau,v):Q_\tau^{k,v}\subset Q^l_\alpha\}}(\tau,v)
\sup_{x\in B_\rho(z^{k,v}_\tau,c_02^{-k})}\left|\mathcal{D}_kf(x)\right|^q\right]^\frac{1}{q}\\
&\quad\lesssim\sup_{k'\in\{0,\ldots,N'\}}
\sup_{\tau'\in I_{k'}}\sup_{v'\in\{1,\ldots,N(k',\tau')\}}
m_{Q_{\tau'}^{k',v'}}(|\mathcal{P}_{k'}f|)
\nonumber\\
&\qquad+\sup_{l'\in\mathbb{N},\,l'>N'}
\sup_{\alpha'\in I_{l'}}
\left\{\frac{1}{\mu(Q_{\alpha'}^{l'})}
\sum_{k'=l'}^\infty\sum_{\tau'\in I_{k'}}\sum_{v'=1}^{N(k',\tau')}
2^{k'sq}\mu(Q_{\tau'}^{k',v'})
\right.\nonumber\\
&\qquad\left.\times\mathbf{1}_{\{(\tau',v'):Q_{\tau'}^{k',v'}\subset Q_{\alpha'}^{l'}\}}(\tau',v')
\inf_{x\in Q_{\tau'}^{k',v'}}\left|\mathcal{P}_{k'}f(x)\right|^q\right\}^\frac{1}{q}.
\end{align*}
This is the desired estimate for the
second term in the left-hand side of \eqref{1528},
which, together with \eqref{2019},
completes the proof of \eqref{1528} and hence Proposition~\ref{proin}.
\end{proof}

\begin{remark}\label{22012}
Let the notation be as in Proposition~\ref{proin}.
\begin{enumerate}
\item[(i)]
Then, by Proposition~\ref{proin},
we easily conclude that the definition of ${F}^{\infty,q}_s(X)$
is independent of the choice of
the approximation of the identity;
that is, for any
$f\in({{\mathcal{G}}}_0^\varepsilon(\beta,\gamma))'$,
\begin{align*}
&\max\Bigg\{\sup_{k\in\{0,\ldots,N\}}
\sup_{\tau\in I_{k}}\sup_{v\in\{1,...,N(k,\tau)\}}m_{Q_\tau^{k,v}}
\left(\left|\mathcal{D}_{0}f\right|\right),\,\\
&\qquad\sup_{l\in\mathbb{N},\,l>N}
\sup_{\alpha\in I_l}\left[\frac{1}{\mu(Q_\alpha^l)}\int_{Q_\alpha^l}\sum_{k=l}^\infty
2^{ksq}\left|\mathcal{D}_kf(x)\right|^q\,d\mu(x)\right]^\frac{1}{q}\Bigg\}\\
&\quad\sim\max\Bigg\{\sup_{k\in\{0,\ldots,N'\}}
\sup_{\tau\in I_{k}}\sup_{v\in\{1,...,N(k,\tau)\}}m_{Q_\tau^{k,v}}
\left(\left|\mathcal{P}_{0}f\right|\right),\,\\
&\qquad\sup_{l\in\mathbb{N},\,l>N'}
\sup_{\alpha\in I_l}\left[\frac{1}{\mu(Q_\alpha^l)}\int_{Q_\alpha^l}\sum_{k=l}^\infty
2^{ksq}\left|\mathcal{P}_kf(x)\right|^q\,d\mu(x)\right]^\frac{1}{q}\Bigg\},
\end{align*}
where the positive equivalence constants are independent of $f$.
A similar conclusion was also obtained in
\cite[Proposition~6.17]{hmy08}
on RD-spaces
and \cite[Proposition~5.11]{WHHY}
on spaces of homogeneous type.
\item[(ii)]
By an argument similar to that used in the proof of Proposition~\ref{proin},
we obtain the following flexible estimate:
for any given $c_0\in[1,\infty)$ and for any
$f\in({{\mathcal{G}}}_0^\varepsilon(\beta,\gamma))'$,
\begin{align*}
&\max\Bigg\{\sup_{k\in\{0,\ldots,N\}}
\sup_{\tau\in I_k}\sup_{v\in\{1,\ldots,N(k,\tau)\}}
\sup_{z\in B_\rho(z_\tau^{k,v},c_02^{-k})}\left|\mathcal{D}_kf(z)\right|,\,
\nonumber\\
&\qquad\sup_{l\in\mathbb{N},\,l>N}
\sup_{\alpha\in I_l}\left[
\frac{1}{\mu(Q^l_\alpha)}\sum_{k=l}^\infty\sum_{\tau\in I_k}
\sum_{v=1}^{N(k,\tau)}
2^{ksq}\mu(Q_\tau^{k,v})
\right.\nonumber\\
&\qquad\left.\times\mathbf{1}_{\{(\tau,v):Q_\tau^{k,v}\subset Q^l_\alpha\}}(\tau,v)
\sup_{x\in B_\rho(z^{k,v}_\tau,c_02^{-k})}
\left|\mathcal{D}_kf(x)\right|^q\right]^\frac{1}{q}\Bigg\}
\nonumber\\
&\quad\lesssim\max\Bigg\{\sup_{k\in\{0,\ldots,N'\}}
\sup_{\tau\in I_k}\sup_{v\in\{1,\ldots,N(k,\tau)\}}m_{Q_\tau^{k,v}}
\left(\left|\mathcal{P}_kf\right|\right),\,
\nonumber\\
&\qquad\sup_{l\in\mathbb{N},\,l>N'}
\sup_{\alpha\in I_l}\left[
\frac{1}{\mu(Q^l_\alpha)}\sum_{k=l}^\infty\sum_{\tau\in I_k}
\sum_{v=1}^{N(k,\tau)}
2^{ksq}\mu(Q_\tau^{k,v})
\right.\nonumber\\
&\qquad\left.\times\mathbf{1}_{\{(\tau,v):Q_\tau^{k,v}\subset Q^l_\alpha\}}(\tau,v)
\inf_{x\in Q^{k,v}_\tau}\left|\mathcal{P}_kf(x)\right|^q\right]^\frac{1}{q}\Bigg\}.
\end{align*}
\end{enumerate}
\end{remark}

Applying Remark~\ref{22012}(ii)
and an argument similar to that used in the proof
of Theorem~\ref{2158},
we obtain the following conclusion;
we omit the details.

\begin{theorem}\label{1605}
Let $(X,\mathbf{q},\mu)$,
$n$, $\varepsilon$, $s$, $\beta$, and $\gamma$ be
the same as in Proposition~\ref{proin} and $q\in[1,\infty]$.
Assume that $\{\mathcal{S}_{2^{-k}}\}_{k\in\mathbb{Z}_+}$ is
an inhomogeneous approximation of the identity
of order $\varepsilon$
as in Definition~\ref{DefIATI}. Let
$\mathcal{D}_{0}:=\mathcal{S}_{1}$ and
$\mathcal{D}_k:=\mathcal{S}_{2^{-k}}-\mathcal{S}_{2^{-k+1}}$ for any $k\in\mathbb{N}$.
Let $\{Q^{l}_\alpha\}_{l\in\mathbb{Z},\,\alpha\in I_l}$ be the collection of all
dyadic cubes given in Lemma~\ref{DyadicSys}. Then
the following assertions hold.
\begin{enumerate}
\item[\textup{(i)}]
$f\in{F}^{\infty,q}_s(X)$
if and only if
$f\in({{\mathcal{G}}}_0^\varepsilon(\beta,\gamma))'$
and \begin{align*}
\|f\|_{\widetilde{{F}^{\infty,q}_s}(X)}
:&=\sup_{l\in\mathbb{Z}_+}\sup_{\alpha\in I_l}
\left[\frac{1}{\mu(Q_\alpha^l)}\int_{Q^l_\alpha}\sum_{k=l}^{\infty}
2^{ksq}\right.\\
&\quad\left.\times\int_{B_\rho(x,2^{-k})}\left|\mathcal{D}_kf(y)\right|^q
\frac{d\mu(y)\,d\mu(x)}{(V_\rho)_{2^{-k}}(x)}
\right]^\frac{1}{q}
\end{align*}
is finite.

\item[\textup{(ii)}]
Let $\lambda\in(n+nq+|s|q,\infty)$. Then
$f\in{F}^{\infty,q}_s(X)$
if and only if
$f\in({{\mathcal{G}}}_0^\varepsilon(\beta,\gamma))'$
and
\begin{align*}
\|f\|_{\overline{{F}^{\infty,q}_s}(X)}
:&=\sup_{l\in\mathbb{Z}_+}\sup_{\alpha\in I_l}
\Bigg\{\frac{1}{\mu(Q_\alpha^l)}\int_{Q^l_\alpha}\sum_{k=l}^{\infty}
2^{ksq}\int_{X}\left|\mathcal{D}_kf(y)\right|^q
\\
&\quad\times\left[\frac{2^{-k}}{2^{-k}+\rho(x,y)}\right]^\lambda
\frac{d\mu(y)\,d\mu(x)}{(V_\rho)_{2^{-k}}(x)
+(V_\rho)_{2^{-k}}(y)}
\Bigg\}^\frac{1}{q}
\end{align*}
is finite.
\end{enumerate}
Moreover, for any $f\in(\mathring{{\mathcal{G}}}_0^\varepsilon(\beta,\gamma))'$,
\begin{align*}
\|f\|_{\widetilde{{F}^{\infty,q}_s}(X)}
\sim\|f\|_{{F}^{\infty,q}_s(X)}
\sim\|f\|_{\overline{{F}^{\infty,q}_s}(X)},
\end{align*}
where the positive equivalence constants are independent of $f$.
\end{theorem}

\begin{remark}
\begin{enumerate}
\item
The Littlewood--Paley $g_\lambda^*$-function characterization
of ${F}^{\infty,q}_s(X)$ on spaces of homogeneous type
was obtained in \cite[Theorem~7.14]{wyyTLinfty},
but with the lower smoothness index $\mathrm{ind\,}(X,\mathbf{q})$
replaced by $\eta\in(0,1)$ which is the
smoothness index of the wavelets on $X$
constructed by Auscher and Hyt\"onen \cite[Theorem~7.1]{AH1}.
We point out that, in more specialized setting of ultra-RD-spaces,
Theorem~\ref{2158} gives the Littlewood--Paley
$g_\lambda^*$-function characterization
of ${F}^{\infty,q}_s(X)$ in
wider ranges of $s$ and $q$ than the corresponding ones of \cite[Theorem~7.14]{wyyTLinfty};
see Remark~\ref{22013}(i).
\item
Let $(X,\mathbf{q},\mu)$ be an Ahlfors-regular
quasi-ultrametric space as in \cite{AlMi2015}.
Then, in this setting, Theorems~\ref{1605} gives the
Lusin area function characterization and the
Littlewood--Paley $g_\lambda^*$-function characterization
of the inhomogeneous Triebel--Lizorkin space ${F}^{\infty,q}_s(X)$
in \cite[Definition~9.6]{AlMi2015}.
\end{enumerate}
\end{remark}


\chapter{Real-Variable Theory
of Hardy--Lorentz Spaces on Ultra-RD-Spaces}
\label{ChapterHL}
\markboth{\scriptsize\rm\sc Real-Variable Theory
of Hardy--Lorentz Spaces on Ultra-RD-Spaces}
{\scriptsize\rm\sc Real-Variable Theory
of Hardy--Lorentz Spaces on Ultra-RD-Spaces}

The main target of this chapter is to establish various real-variable
characterizations of Hardy--Lorentz spaces $H^{p,q}(X)$
for the optimal ranges
$$
p\in\left(\frac{n}{n+{\mathrm{ind\,}}(X,\mathbf{q})},1\right]
\ \ \text{and}\ \
q\in(0,\infty]
$$
on ultra-RD-spaces $(X,\mathbf{q},\mu)$,
where $\mu$ is
a Borel-semiregular measure and $n\in(0,\infty)$
is the upper dimension of $(X,\mathbf{q},\mu)$.
In this chapter, unless otherwise specified,
we allow $\mu(X)$ to be finite or infinite (equivalently,
by Remark~\ref{1031}, $X$ can be bounded or unbounded).
Specifically,
in Section~\ref{section5.1},
we first recall the definition and some basic properties of
Lorentz spaces.
Then, by means of different types of maximal functions,
we introduce several brands of Hardy--Lorentz spaces
and further show that all these Hardy--Lorentz spaces are equivalent.
In Section~\ref{section5.3}, we
establish the atomic and the finite atomic characterizations of
Hardy--Lorentz spaces.
In Section~\ref{section5.5}, we
establish the molecular characterization of Hardy--Lorentz spaces.
In Section~\ref{section5.6}, we first introduce the
Hardy--Lorentz space by means of the Lusin area function and then
prove that it is independent of the
choice of approximations of the identity.
Next, we establish its atomic
and its Littlewood--Paley function characterizations.
In Section~\ref{section5.7},
we obtain the interpolation properties of Hardy--Lorentz spaces via the real method.
In Section~\ref{Section5.5}, we establish the duality theorem
for Hardy--Lorentz spaces.
In Section~\ref{section5.8},
we show the boundedness of
Calder\'on--Zygmund operators on Hardy--Lorentz spaces.
It is important to emphasize that many of the main results in
this chapter are established in the more general setting of
quasi-ultrametric spaces of homogeneous type.

\section{Maximal Theory of Hardy--Lorentz Spaces}
\label{section5.1}

In this section, we first recall the definition of
Lorentz spaces and collect some its properties.
Then we introduce several brands of
Hardy--Lorentz spaces based on
different maximal functions
and prove that all of these spaces coincide.

Let $(X,\mu)$ be an arbitrary measure space.
To define Lorentz spaces, we need the following notion.
The \emph{non-increasing rearrangement}\index[words]{non-increasing
rearrangement} $R_f$\index[symbols]{R@$R_f$} of any $\mu$-measurable function
$f: X\to\overline{\mathbb{R}}:=\mathbb{R}\cup\{\pm\infty\}$
\index[symbols]{R@$\overline{\mathbb{R}}$} is defined by setting,
for any $t\in\mathbb{R}_+$,
$$
R_f(t):=\inf\left\{\alpha\in\mathbb{R}_+:d_f(\alpha)\leq t\right\}.\index[symbols]{R@$R_f(t)$}
$$
Here, and thereafter,
for any $\alpha\in\mathbb{R}_+$,
$$
d_f(\alpha):=\mu\left(\left\{x\in X:|f(x)|>\alpha\right\}\right)\index[symbols]{D@$d_f$}
$$
denotes the \emph{distribution function}
\index[words]{distribution function} of $f$.
Note that we adopt the convention $\inf\varnothing=\infty$.
The function $R_f:\mathbb{R}_+\to[0,\infty]$ is right continuous and hence measurable on
$\mathbb{R}_+$ with respect to the 1-dimensional Lebesgue measure; see, for instance,
\cite[Proposition~1.4.5]{G2014C}.

In the above context, for any $p,q\in(0,\infty]$,
the \emph{Lorentz space}\index[words]{Lorentz space}
$L^{p,q}(X)$\index[symbols]{L@$L^{p,q}(X)$}
is defined to be the set of all $\mu$-measurable
functions $f: X\to\overline{\mathbb{R}}$
such that
\begin{align*}
\|f\|_{L^{p,q}(X)}:=
\begin{cases}
\displaystyle
\left\{\int_0^\infty\left[t^{\frac{1}{p}}
R_f(t)\right]^q\,\frac{dt}{t}\right\}^{\frac{1}{q}}
\ &\text{if }q\in(0,\infty),\\
\displaystyle
\displaystyle\sup_{t\in(0,\infty)}t^{\frac{1}{p}}R_f(t)
&\text{if }q=\infty
\end{cases}
\end{align*}
is finite.

We take a moment to make a few comments and collect a few basic properties of $L^{p,q}(X)$ spaces.
The reader is referred to \cite[Section~1.4]{G2014C} for a more detailed account of these spaces.
As is the case with the Lebesgue space $L^p(X)$, two functions
in $L^{p,q}(X)$ are regarded as equal if they coincide pointwise $\mu$-almost
everywhere on $X$. Moreover, the spaces $L^{p,q}(X)$ are quasi-Banach
spaces with respect to the quasi-norm
$\|\cdot\|_{L^{p,q}(X)}$ for all $p,q\in(0,\infty]$.
The next proposition highlights the fact that $\|\cdot\|_{L^{p,q}(X)}$ scales naturally
under power dilations; see, for instance, \cite[Remark 1.4.7]{G2014C}.

\begin{proposition}
\label{r;p,q=pr,qr;r}
Let $(X,\mu)$ be a measure space, $p,q\in(0,\infty]$,
and $r\in(0,\infty)$.
Then, for any $\mu$-measurable function $f: X\to\overline{\mathbb{R}}$,
$$
\left\|\,|f|^r\right\|_{L^{p,q}(X)}=\|f\|^r_{L^{pr,qr}(X)}.
$$
\end{proposition}

The following proposition is precisely \cite[Proposition 1.4.9]{G2014C} and it
provides an important description of the $\|\cdot\|_{L^{p,q}(X)}$ quasi-norm in terms
of the distribution function $d_f$.

\begin{proposition}
\label{Lpqnorm=c,d}
Let $(X,\mu)$ be a measure space, $p\in(0,\infty)$,
and $q\in(0,\infty]$.
Then, for any $\mu$-measurable function $f$,
\begin{align*}
\|f\|_{L^{p,q}(X)}
=p^{\frac{1}{q}}\left\{\int_0^\infty\alpha^{q-1}
\left[d_f(\alpha)\right]^{\frac{q}{p}}
\,d\alpha\right\}^{\frac{1}{q}}
\sim p^{\frac{1}{q}}
\left\{\sum_{k\in\mathbb{Z}}2^{kq}
\left[d_f(2^k)\right]^{\frac{q}{p}}\right\}^{\frac{1}{q}}
\end{align*}
with the usual modification when $q=\infty$,
where the positive equivalence constants
are independent of $p$, $q$, and $f$.
\end{proposition}

The scale of Lorentz spaces naturally contains Lebesgue spaces and
weak Lebesgue spaces. Indeed, it is easy to see from
Proposition~\ref{Lpqnorm=c,d} and the definition of $\|\cdot\|_{L^{p,q}(X)}$
that $L^{p,p}(X)=L^p(X)$ for all $p\in(0,\infty]$ and,
when $q=\infty$, a simple power rescaling of $t$ in
the definition of $\|\cdot\|_{L^{p,q}(X)}$ gives
$L^{p,\infty}(X)=$ weak-$L^p(X)$ for all $p\in(0,\infty)$.
It is also worth noting that, when $p=\infty$,
$L^{\infty,q}(X)=\{0\}$ for all $q\in(0,\infty)$;
see, for instance, \cite[Example~1.4.8]{G2014C}.

As illustrated by the next proposition, the spaces $L^{p,q}(X)$ are nested
in the second exponent $q$; see, for instance, \cite[Proposition 1.4.10]{G2014C}.

\begin{proposition}
\label{p,q<p,s}
Let $(X,\mu)$ be a measure space, $p,q\in(0,\infty]$, and $s\in(0,q]$.
Then $L^{p,s}(X)\subset L^{p,q}(X)$ and there exists a positive constant $C_{p,q,s}$,
depending on $p$, $q$, and $s$,
such that, for any $\mu$-measurable function $f: X\to\overline{\mathbb{R}}$,
\begin{align*}
\|f\|_{L^{p,q}(X)}\leq C_{p,q,s}\|f\|_{L^{p,s}(X)}.
\end{align*}
\end{proposition}
Collecting these observations, we have the following embeddings.
That is, for any
$0<q\leq p\leq\infty$,
$$
L^{p,q}(X)\hookrightarrow L^{p}(X)\hookrightarrow L^{p,\infty}(X).
$$

We conclude our brief survey of Lorentz spaces with the following proposition which
collects a few additional properties of $L^{p,q}(X)$ spaces;
see, for instance, \cite[Proposition~1.4.5 and Example~1.4.8]{G2014C}.

\begin{proposition}
\label{inc-norm}
Let $(X,\mu)$ be a measure space and $p,q\in(0,\infty]$.
Then the following statements hold.
\begin{enumerate}
\item[\emph{(i)}]
If $f\in L^{p,q}(X)$, then, for $\mu$-almost every $x\in X$,
$|f(x)|<\infty$.

\item[\emph{(ii)}]
If $f,g:X\to\overline{\mathbb{R}}$ are two $\mu$-measurable functions such that
$|f|\leq|g|$ pointwise $\mu$-almost everywhere on $X$, then
$\|f\|_{L^{p,q}(X)}\leq\|g\|_{L^{p,q}(X)}$.

\item[\emph{(iii)}]
If $p\in(0,\infty)$ and $E\subset X$ is a $\mu$-measurable set, then
$\|{\mathbf{1}}_E\|_{L^{p,q}(X)}\sim[\mu(E)]^\frac{1}{p}$,
where the positive equivalence constants are independent of $E$.
\end{enumerate}
\end{proposition}

Moving on, we introduce the following three maximal
functions that will in turn be used to define
different brands of maximal Hardy--Lorentz spaces.

\begin{definition}\label{maxfunc}
Let $(X,\mathbf{q},\mu)$ be a quasi-ultrametric
space of homogeneous type,
$\rho\in\mathbf{q}$
have the property that all $\rho$-balls are $\mu$-measurable,
and $\varepsilon,\beta,\gamma\in(0,\infty)$ be such that
$0<\beta,\gamma<\varepsilon
\preceq\mathrm{ind\,}(X,\mathbf{q})$.
Assume that $\{\mathcal{S}_{2^{-k}}\}_{k\in\mathbb{Z}}$
is a homogeneous approximation of the identity
of order $\varepsilon$ as in Definition~\ref{DefHATI}
and let $f\in({\mathcal{G}}^\varepsilon_0(\beta,\gamma))'$.
\begin{enumerate}
\item[\rm(i)]
The \emph{radial maximal function}\index[words]{radial maximal function}
$\mathcal{M}_\rho^+f$\index[symbols]{M@$\mathcal{M}_\rho^+$} of $f$
is defined by setting, for any $x\in X$,
$$
\mathcal{M}_\rho^+f(x):=\sup_{k\in\mathbb{Z}}
\left|\mathcal{S}_{2^{-k}}f(x)\right|.
$$

\item[\rm(ii)]
The \emph{non-tangential maximal
function}\index[words]{non-tangential maximal function}
$(\mathcal{M}_\rho)_\theta f$\index[symbols]{M@$(\mathcal{M}_\rho)_\theta$} of $f$,
with $\theta\in(0,\infty)$,
is defined by setting, for any $x\in X$,
\begin{align*}
(\mathcal{M}_\rho)_\theta f(x):=\sup_{k\in\mathbb{Z}}
\sup_{y\in B_\rho(x,\theta2^{-k})}
\left|\mathcal{S}_{2^{-k}}f(y)\right|.
\end{align*}

\item[\rm(iii)]
The \emph{grand maximal function}\index[words]{grand maximal function}
$f^*_\rho$\index[symbols]{F@$f^*_\rho$} of $f$
is defined by setting, for any $x\in X$,
$$
f^*_\rho(x):=\sup\left\{|\langle f,\phi\rangle|:
\phi\in{\mathcal{G}}^\varepsilon_0(\beta,\gamma)\text{ and }
\|\phi\|_{G(x,r^\alpha,\frac{\beta}{{\alpha}},\frac{\gamma}{{\alpha}})}\leq 1
\text{ for some }r\in(0,\infty)\right\},
$$
where $\alpha:=\min\{1,\,(\log_2C_\rho)^{-1}\}$ and
the \emph{space $G(x,r,\beta,\gamma)$},\index[symbols]{G@$G(x,r,\beta,\gamma)$}
with $x\in X$, $r\in(0,\infty)$, and $0<{\beta},{\gamma}<\infty$,
is defined to be the set of all
$\mu$-measurable functions $g:X\to\mathbb{C}$
satisfying that there exists $C\in(0,\infty)$ such that,
for any $y\in X$,
\begin{align}\label{1G1}
|g(y)|\leq C\frac{1}{\mu(B_{\varrho}(y,r+\varrho(x,y)))}
\left[\frac{r}{r+\varrho(x,y)}\right]^{{\gamma}}
\end{align}
and, for any $y,y'\in X$ with $\varrho(y,y')\leq\frac{r+\varrho(x,y)}{2C_\varrho}$,
\begin{align}\label{2G2}
\left|g(y)-g(y')\right|\leq C
\left[\frac{\varrho(y,y')}{r+\varrho(x,y)}
\right]^{{\beta}}
\frac{1}{\mu(B_{\varrho}(y,r+\varrho(x,y)))}
\left[\frac{r}{r+\varrho(x,y)}\right]^{{\gamma}},
\end{align}
where
$\varrho:=(\rho_{\#})^{\alpha}$ and
$\rho_{\#}\in\mathbf{q}$ is the regularized
quasi-ultrametric of $\rho$ as in Definition~\ref{D2.3};
moreover, for any $g\in G(x,r,{\beta},{\gamma})$, define
$$
\|g\|_{G(x,r,{\beta},{\gamma})}:=\inf\left\{C\in(0,\infty):
\eqref{1G1}\ \mathrm{and}\ \eqref{2G2}\ \mathrm{hold}\right\}.
$$
\end{enumerate}
\end{definition}

\begin{remark}\label{RemMf}
Let the notation be as in Definition~\ref{maxfunc}.
\begin{enumerate}
\item[\rm(i)]
From Remark~\ref{Rm}(iii),
we deduce that, for any $\rho'\in\mathbf{q}$ such
that all $\rho'$-balls are $\mu$-measurable,
$f\in({\mathcal{G}}^\varepsilon_0(\beta,\gamma))'$, and $x\in X$,
\begin{align*}
\mathcal{M}_\rho^+f(x)\approx\mathcal{M}_{\rho'}^+f(x)
\end{align*}
and hence we simply write $\mathcal{M}_\rho^+f$ by $\mathcal{M}^+f$
in the remainder of this chapter.

\item[\rm(ii)]
It is easy to see that, for any $0<\theta_1\leq\theta_2<\infty$,
$f\in({\mathcal{G}}^\varepsilon_0(\beta,\gamma))'$, and $x\in X$,
\begin{align*}
(\mathcal{M}_\rho)_{\theta_1}f(x)\leq(\mathcal{M}_\rho)_{\theta_2}f(x).
\end{align*}
If $\rho'\in\mathbf{q}$ has the property that
all $\rho'$-balls are $\mu$-measurable
and $C\in[1,\infty)$ is such that $C^{-1}\rho\leq\rho'\leq C\rho$,
then, for any $\theta\in(0,\infty)$,
$f\in({\mathcal{G}}^\varepsilon_0(\beta,\gamma))'$, and $x\in X$,
\begin{align*}
(\mathcal{M}_\rho)_{C^{-1}\theta}f(x)
\leq(\mathcal{M}_{\rho'})_{\theta}f(x)
\leq(\mathcal{M}_\rho)_{C\theta}f(x).
\end{align*}

\item[\rm(iii)]
Let $x\in X$, $r\in(0,\infty)$,
and $0<\beta,\gamma<\infty$. Then,
by Lemmas~\ref{A3.2}(vi) and~\ref{L2.3}(ii), we find that
$$
\mathcal{G}(x,r,\beta,\gamma)
=G\left(x,r^{\alpha},\frac{\beta}{\alpha},\frac{\gamma}{\alpha}\right)
$$
with equivalent norms,
where the positive equivalence constants are independent of both $x$ and $r$.
Note that, if $\rho$ is a metric, then we have $\alpha=1$.

\item[\rm(iv)]
From (i) of the present remark and Remark~\ref{Rmtest}(ii),
we infer that,
if $\rho,\rho'\in\mathbf{q}$ have the property that all $\rho$-balls
and all $\rho'$-balls are $\mu$-measurable,
then, for any $f\in({\mathcal{G}}^\varepsilon_0(\beta,\gamma))'$
and $x\in X$,
\begin{align*}
f^*_\rho(x)\sim f^*_{\rho'}(x),
\end{align*}
where the positive equivalence constants are independent of both
$f$ and $x$.
\end{enumerate}
\end{remark}

The following proposition gives the measurability of
functions $\mathcal{M}^+f$, $(\mathcal{M}_\rho)_\theta f$,
and $f^*_\rho$;
see \cite[Lemma 4.7 and Proposition 4.10]{AlMi2015}
in the setting of Ahlfors-regular
quasi-ultrametric spaces
and see also \cite[Remark~2.9(iii)]{LLD}
for the proof of the measurability of $f^*_\rho$
on doubling metric measure spaces.

\begin{proposition}\label{measurability}
Let $(X,\mathbf{q},\mu)$ be a quasi-ultrametric
space of homogeneous type
with upper dimension $n\in(0,\infty)$,
$\rho\in\mathbf{q}$
have the property that all $\rho$-balls are $\mu$-measurable,
and $\varepsilon,\beta,\gamma\in(0,\infty)$ be such that
$0<\beta,\gamma<\varepsilon
\preceq\mathrm{ind\,}(X,\mathbf{q})$.
Assume that $\{\mathcal{S}_{2^{-k}}\}_{k\in\mathbb{Z}}$
is a homogeneous approximation of the identity
of order $\varepsilon$ as in Definition~\ref{DefHATI}.
Then, for any $f\in({\mathcal{G}}^\varepsilon_0(\beta,\gamma))'$,
both $\mathcal{M}^+f$ and $f^*_\rho$ are $\mu$-measurable.
Moreover, if $\rho$ is a regularized quasi-ultrametric,
then, for any $\theta\in(0,\infty)$,
$(\mathcal{M}_\rho)_\theta f$ is also $\mu$-measurable.
\end{proposition}

\begin{proof}
Let $f\in({\mathcal{G}}^\varepsilon_0(\beta,\gamma))'$.
We first show that $\mathcal{M}^+f$ is $\mu$-measurable.
An argument similar to that used in the estimation
of \eqref{2235} implies that, for any $k\in\mathbb{Z}$
and $x,y\in X$,
\begin{align*}
\left|\mathcal{S}_{2^{-k}}f(x)-\mathcal{S}_{2^{-k}}f(y)\right|
&=\left|\left\langle f,S_{2^{-k}}(x,\cdot)-
S_{2^{-k}}(y,\cdot)\right\rangle\right|
\nonumber\\
&\leq\left\|f\right\|_{({\mathcal{G}}^\varepsilon_0(\beta,\gamma))'}
\left\|S_{2^{-k}}(x,\cdot)-
S_{2^{-k}}(y,\cdot)\right\|_{{\mathcal{G}}^\varepsilon_0(\beta,\gamma)}
\nonumber\\
&\leq C\left\|f\right\|_{({\mathcal{G}}^\varepsilon_0(\beta,\gamma))'}
\left[\rho(x,y)\right]^\varepsilon
\left[2^{-k}+\rho(x_1,x)\right]^{n+\beta+\gamma},
\end{align*}
where the positive constant $C$ is independent of $k$, $f$, $x$, and $y$,
which, further implies that, for any given $\delta\in(0,1)$
and for any $z\in B_\rho(x,r_x(\delta))$,
\begin{align*}
\mathcal{M}^+f(z)=\sup_{k\in\mathbb{Z}}\left|\mathcal{S}_{2^{-k}}f(z)
\right|\ge\sup_{k\in\mathbb{Z}}\left|\mathcal{S}_{2^{-k}}f(x)
\right|-\delta=\mathcal{M}^+f(x)-\delta,
\end{align*}
where
$$
r_x(\delta):=\left\{\frac{\delta}{C\|f\|_{({\mathcal{G}}^\varepsilon_0(\beta,\gamma))'}
[2^{-k}+\rho(x_1,x)]^{n+\beta+\gamma}}\right\}^\frac{1}{\varepsilon}.
$$
Thus, $\mathcal{M}^+f$ is lower semicontinuous.
This, combined with an argument similar to
that used in the proof of Proposition~\ref{1507},
implies that $\mathcal{M}^+f$ is $\mu$-measurable.

Now, we prove that
$f^*_\rho$ is $\mu$-measurable.
It suffices to show that, for any $\lambda\in(0,\infty)$,
$$
\Gamma_\lambda:=\left\{x\in X:f^*_\rho(x)>\lambda\right\}
$$
is open.
To this end, let $x\in\Gamma_\lambda$.
Then there exists $\phi\in\mathcal{G}^\varepsilon_0(\beta,\gamma)$
such that $\|\phi\|_{G(x,r^{\alpha},\frac{\beta}{\alpha},
\frac{\gamma}{{\alpha}})}\leq1$ for some $r\in(0,\infty)$
and $|\langle f,\phi\rangle|>\lambda$,
where $\alpha:=\min\{1,\,(\log_2C_\rho)^{-1}\}$.
Let $\eta\in(0,\infty)$ be such that
$|\langle f,\phi\rangle|>(1+\eta)\lambda$.
For any $y\in B_\varrho(x,r^\alpha-r^\alpha(1+\eta)^{-\frac{\alpha}{\gamma}})$,
choose $s\in(0,\infty)$ such that
$$
s^\alpha\in\left(r^\alpha(1+\eta)^{-\frac{\alpha}{\gamma}},r^\alpha-\varrho(x,y)\right),
$$
where $\varrho:=(\rho_{\#})^\alpha$.
By this, we conclude that, for any $z\in X$,
\begin{align*}
|\phi(z)|
&\leq\frac{1}{\mu(B_{\varrho}(z,r^\alpha+\varrho(x,z)))}
\left[\frac{r^\alpha}{r^\alpha+\varrho(x,z)}\right]^\frac{\gamma}{\alpha}
\quad\text{by \eqref{1G1}}
\\
&\leq\frac{1}{\mu(B_{\varrho}(z,r^\alpha+\varrho(y,z)-\varrho(x,y)))}
\left[\frac{r^\alpha}{r^\alpha+\varrho(y,z)-\varrho(x,y)}\right]^\frac{\gamma}{\alpha}
\\
&\qquad\ \text{by Lemma~\ref{L2.3}(iv)}\\
&\leq(1+\eta)\frac{1}{\mu(B_{\varrho}(z,s^\alpha+\varrho(y,z)))}
\left[\frac{s^\alpha}{s^\alpha+\varrho(y,z)}\right]^\frac{\gamma}{\alpha}
\end{align*}
and, for any $z,z'\in X$ with
$\varrho(z,z')\leq\frac{s^\alpha+\varrho(y,z)}{2C_\varrho}
\leq\frac{s^\alpha+\varrho(y,x)+\varrho(x,z)}{2C_\varrho}
<\frac{r^\alpha+\varrho(x,z)}{2C_\varrho}$,
\begin{align*}
\left|\phi(z)-\phi(z')\right|
&\leq\left[\frac{\varrho(z,z')}{r^\alpha+\varrho(x,z)}
\right]^\frac{\beta}{\alpha}
\frac{1}{\mu(B_{\varrho}(z,r^\alpha+\varrho(x,z)))}
\left[\frac{r^\alpha}{r^\alpha+\varrho(x,z)}\right]^\frac{\gamma}{\alpha}
\\
&\qquad\ \text{by \eqref{2G2}}\\
&\leq(1+\eta)
\left[\frac{\varrho(z,z')}{s^\alpha+\varrho(y,z)}
\right]^\frac{\beta}{\alpha}
\frac{1}{\mu(B_{\varrho}(z,s^\alpha+\varrho(y,z)))}
\left[\frac{s^\alpha}{s^\alpha+\varrho(y,z)}\right]^\frac{\gamma}{\alpha},
\end{align*}
which further implies that $\|\phi\|_{G(y,s^\alpha,
\frac{\beta}{\alpha},\frac{\gamma}{\alpha})}\leq1+\eta$.
Thus, for any $y\in B_\varrho(x,r-r^\alpha(1+\eta)^{-\frac{\alpha}{\gamma}})$,
\begin{align*}
f^*_\rho(y)\ge\left|\left\langle f,\frac{\phi}{1+\eta}\right\rangle\right|>\lambda,
\end{align*}
from which we deduce that
$B_\varrho(x,r-r^\alpha(1+\eta)^{-\frac{\alpha}{\gamma}})\subset\Gamma_\lambda$.
Hence, $\Gamma_\lambda$ is open in $\tau_\varrho=\tau_{\mathbf{q}}$
and hence we can conclude that $f^*_\rho$ is $\mu$-measurable.

Finally, we prove that, for any $\theta\in(0,\infty)$,
$(\mathcal{M}_\rho)_\theta f$
is $\mu$-measurable, assuming that $\rho$ is
a regularized quasi-ultrametric.
Let $\theta\in(0,\infty)$.
It suffices to show that, for any $\lambda\in(0,\infty)$,
the level set
\begin{align*}
\Omega_\lambda:=\left\{x\in X:(\mathcal{M}_\rho)_\theta f(x)>\lambda\right\}
\end{align*}
is open. Let $\lambda\in(0,\infty)$ and fix $x\in\Omega_\lambda$.
Then there exist $k\in\mathbb{Z}$ and $y\in B_\rho(x,\theta2^{-k})$
such that $|\mathcal{S}_{2^{-k}}f(y)|>\lambda$.
Since we are assuming that $\rho$ is a regularized quasi-ultrametric,
we infer that
there exists $\tau\in(0,\infty)$ such that $\rho^\tau$ is a metric on $X$.
Let
$$
r_x:=\left(\left(\theta2^{-k}\right)^{\tau}-\left[\rho(x,y)
\right]^{\tau}\right)^\frac{1}{\tau}\in(0,\infty).
$$
Since $\rho^\tau$ is a metric on $X$, we deduce that,
for any given $x'\in B_\rho(x,r_x)$,
\begin{align*}
\rho(x',y)&=\left(\left[\rho(x',y)\right]^\tau\right)^\frac{1}{\tau}
\leq\left(\left[\rho(x',x)\right]^\tau
+\left[\rho(x,y)\right]^\tau\right)^\frac{1}{\tau}<\theta 2^{-k},
\end{align*}
which further implies that
$y\in B_\rho(x',\theta2^{-k})$ and
\begin{align*}
(\mathcal{M}_\rho)_\theta f(x')\ge
\sup_{z\in B_\rho(x',\theta2^{-k})}
\left|\mathcal{S}_{2^{-k}}f(z)\right|
\ge\left|\mathcal{S}_{2^{-k}}f(y)\right|
>\lambda.
\end{align*}
Thus, $B_\rho(x,r_x)\subset\Omega_\lambda$ and hence $\Omega_\lambda$ is open,
from which we infer that $(\mathcal{M}_\rho)_\theta f$
is $\mu$-measurable.
This finishes the proof of Proposition~\ref{measurability}.
\end{proof}

Next, based on the maximal functions in Definition~\ref{maxfunc},
we introduce three different maximal-function-based Hardy--Lorentz spaces as follows.

\begin{definition}
\label{DefHardyL}
Let $(X,\mathbf{q},\mu)$ be a quasi-ultrametric
space of homogeneous type with upper dimension $n\in(0,\infty)$,
$\rho\in\mathbf{q}$ have
the property that all $\rho$-balls are $\mu$-measurable,
$p,q\in(0,\infty]$, $\theta\in(0,\infty)$,
and
$n(\frac{1}{p}-1)_+<\beta,\gamma<\varepsilon<\infty$.
The \emph{Hardy--Lorentz spaces}\index[words]{Hardy--Lorentz space}
$H^{p,q}_+(X,\rho)$\index[symbols]{H@$H^{p,q}_+(X)$}
and $H^{p,q}_*(X,\rho)$\index[symbols]{H@$H^{p,q}_*(X)$}
are defined, respectively, by setting
$$
H^{p,q}_+(X,\rho):=\left\{f\in({\mathcal{G}}^\varepsilon_0(\beta,\gamma))':
\|f\|_{H^{p,q}_+(X)}:=\left\|\mathcal{M}^+f\right\|_{L^{p,q}(X)}<\infty\right\}
$$
and
$$
H^{p,q}_*(X,\rho):=\left\{f\in({\mathcal{G}}^\varepsilon_0(\beta,\gamma))':
\|f\|_{H^{p,q}_*(X)}:=\left\|f^*_\rho\right\|_{L^{p,q}(X)}<\infty\right\}.
$$
If $\rho\in\mathbf{q}$ is a regularized
quasi-ultrametric as in Definition~\ref{D2.3},
then the \emph{Hardy--Lorentz space}
$H^{p,q}_\theta(X,\rho)$\index[symbols]{H@$H^{p,q}_\theta(X)$}
is defined by setting
$$
H^{p,q}_\theta(X,\rho):=\left\{f\in({\mathcal{G}}^\varepsilon_0(\beta,\gamma))':
\|f\|_{H^{p,q}_\theta(X)}:=
\left\|(\mathcal{M}_\rho)_\theta f\right\|_{L^{p,q}(X)}<\infty\right\}.
$$
\end{definition}

\begin{remark}\label{R6.6}
Let the notation be as in Definition~\ref{DefHardyL}.
\begin{enumerate}
\item Then,
for any $p\in(0,\frac{n}{n+{\mathrm{ind}}_H(X,\mathbf{q})}]$,
the Hardy--Lorentz spaces
in Definition~\ref{DefHardyL} are trivial
because the space ${\mathcal{G}}^\varepsilon_0(\beta,\gamma)$
only contains constant functions.
\item
Using the fact that $L^{\infty,q}(X)=\{0\}$
for any $q\in(0,\infty)$ and Proposition~\ref{2250},
we obtain, for any $q\in(0,\infty)$,
\begin{align*}
H^{\infty,q}_+(X,\rho)
=H^{\infty,q}_\theta(X,\rho)
=H^{\infty,q}_*(X,\rho)
=\{0\}.
\end{align*}
\item
Assume that $p\in(\frac{n}{n+\mathrm{ind\,}(X,\mathbf{q})},\infty]$,
$q\in(0,\infty]$, $\theta\in(0,\infty)$,
and $n(\frac{1}{p}-1)_+<\beta,\gamma<\varepsilon\preceq\mathrm{ind\,}(X,\mathbf{q})$.
Then $H^{p,q}_+(X,\rho)$, $H^{p,q}_\theta(X,\rho)$,
and $H^{p,q}_*(X,\rho)$ are quasi-Banach spaces,
which can be seen from both Theorem~\ref{HL-M}
(see also Theorem~\ref{HL-Mfinite}) and Proposition~\ref{5.1.11} below.
\item
If $\rho,\varrho\in\mathbf{q}$ are
such that both $f^*_\rho$ and $f^*_\varrho$
are $\mu$-measurable functions, then,
by Remarks~\ref{Rmtest}(ii) and~\ref{RemMf}(iv),
we find that $H^{p,q}_*(X,\rho)$
and $H^{p,q}_*(X,\varrho)$
defined, respectively, in terms of $f^*_\rho$
and $f^*_\varrho$ coincide with equivalent quasi-norms.
Thus, in this sense, $H^{p,q}_*(X,\rho)$ is
independent of the choice of $\rho\in\mathbf{q}$.
Obviously, $H^{p,q}_+(X,\rho)$ is also
independent of the choice of $\rho\in\mathbf{q}$.
Moreover, we prove that
$H^{p,q}_\theta(X,\rho)$ is also
independent of the choice of $\rho\in\mathbf{q}$. Similar considerations
apply to the parameter $\theta$
in $H^{p,q}_\theta(X,\rho)$;
see Theorem~\ref{HL-M} and Theorem~\ref{HL-Mfinite} below.
Thus, throughout this monograph,
we sometimes simply write $H^{p,q}_\theta(X)$,
$H^{p,q}_+(X)$, and $H^{p,q}_*(X)$
instead of $H^{p,q}_\theta(X,\rho)$,
$H^{p,q}_+(X,\rho)$, and $H^{p,q}_*(X,\rho)$, where,
without loss of generality, we may assume that $H^{p,q}_\theta(X)$,
$H^{p,q}_+(X)$, and $H^{p,q}_*(X)$ are defined using a regularized quasi-ultrametric.

\item
Recall that, for any $p\in(0,\infty]$, $L^{p,p}(X)=L^p(X)$, which implies that
$$
H^{p,p}_+(X)=H^{p}_+(X):=\left\{f\in({\mathcal{G}}^\varepsilon_0(\beta,\gamma))':
\|f\|_{H^{p}_+(X)}:=\left\|\mathcal{M}^+f\right\|_{L^{p}(X)}<\infty\right\},
$$
$$
H^{p,p}_\theta(X)=H^{p}_\theta(X):=\left\{f\in({\mathcal{G}}^\varepsilon_0(\beta,\gamma))':
\|f\|_{H^{p}_\theta(X)}:=
\left\|(\mathcal{M}_\rho)_\theta f\right\|_{L^{p}(X)}<\infty\right\},
$$
and
$$
H^{p,p}_*(X)=H^{p}_*(X):=\left\{f\in({\mathcal{G}}^\varepsilon_0(\beta,\gamma))':
\|f\|_{H^{p}_*(X)}:=\left\|f^*_\rho\right\|_{L^{p}(X)}<\infty\right\}.
$$
\item We show
that each of the maximal-function-based Hardy--Lorentz spaces
in Definition~\ref{DefHardyL} is independent of
the choice of parameters $\beta$ and $\gamma$
within certain specified ranges.
In light of this fact, we choose not
to include these parameters in the notation of
$H^{p,q}_+(X)$, $H^{p,q}_\theta(X)$, and
$H^{p,q}_*(X)$;
see Theorem~\ref{indep2} below.
\end{enumerate}
\end{remark}

We now aim to prove that
the three maximal-function-based Hardy--Lorentz spaces
in Definition~\ref{DefHardyL} coincide with
equivalent quasi-norms. This is accomplished in
Theorem~\ref{HL-M} and Theorem~\ref{HL-Mfinite}; however,
we first need to collect a number of needed lemmas.

\begin{lemma}\label{MMM}
Let $(X,\mathbf{q},\mu)$ be a quasi-ultrametric space
of homogeneous type, $\rho\in\mathbf{q}$
have the property that all $\rho$-balls are
$\mu$-measurable,
and
$0<\beta,\gamma<\varepsilon
\preceq\mathrm{ind\,}(X,\mathbf{q})$.
Then, for any given $\theta\in(0,\infty)$,
there exists a positive constant $C$
such that,
for any $f\in({\mathcal{G}}^\varepsilon_0(\beta,\gamma))'$ and $x\in X$,
\begin{align}
\label{2114}
\mathcal{M}^+f(x)\leq(\mathcal{M}_\rho)_\theta f(x)
\leq Cf^*_\rho(x).
\end{align}
\end{lemma}

\begin{proof}
Without loss of generality, we may assume that
$\rho\in\mathbf{q}$ is a regularized quasi-ultrametric as in Definition~\ref{D2.3}.
Let $f\in({\mathcal{G}}^\varepsilon_0(\beta,\gamma))'$ and $x\in X$.
The first inequality of \eqref{2114} directly follows from
the definitions of $\mathcal{M}^+f(x)$
and $(\mathcal{M}_\rho)_\theta f(x)$.

To show the second inequality of \eqref{2114},
we first claim that
\begin{align}\label{wza-47}
\sup_{k\in\mathbb{Z}}
\sup_{y\in B_\rho(x,\theta2^{-k})}
\left\|S_{2^{-k}}(y,\cdot)\right\|_{{\mathcal{G}}(x,2^{-k},\beta,\gamma)}<\infty.
\end{align}
Let $k\in\mathbb{Z}$
and $y\in B_\rho(x,\theta2^{-k})$.
We first prove the size condition.
Let $\tau\in X$.
From Definition~\ref{DefATI}(i),
we deduce that there exists a constant $C'\in(1,\infty)$
such that,
if $\rho(y,\tau)\ge C'2^{-k}$,
then $S_{2^{-k}}(y,\tau)=0$. Thus, we may assume that
$\rho(y,\tau)<C'2^{-k}$.
By Definition~\ref{D2.1}(iii), we find that
\begin{align*}
\rho(x,\tau)
\leq C_\rho\max\left\{\rho(x,y),\,\rho(y,\tau)\right\}
<C_\rho(\theta+C')2^{-k},
\end{align*}
which, together with \eqref{doublingcond},
further implies that
$V_\rho(x,\tau)\lesssim (V_\rho)_{2^{-k}}(x)$,
where the implicit positive constant depends
on $\theta$, $C'$, $\rho$, and $\mu$.
From this, Definition~\ref{DefATI}(i),
and Lemma~\ref{A3.2}(vii), it follows that
\begin{align*}
|S_{2^{-k}}(y,\tau)|
&\lesssim\frac{1}{(V_\rho)_{2^{-k}}(y)}
\sim\frac{1}{(V_\rho)_{2^{-k}}(x)}
\\
&\lesssim\frac{1}{(V_\rho)_{2^{-k}}(x)+V_\rho(x,\tau)}
\left[\frac{2^{-k}}{2^{-k}+\rho(x,\tau)}\right]^{\gamma},
\end{align*}
where, in the last step, we used the fact that
$[\frac{2^{-k}}{2^{-k}+\rho(x,\tau)}]^{\gamma}
\ge[\frac{1}{1+C_\rho(C'+\theta)}]^{\gamma}\sim1$. Note that the
implicit positive constants in this calculation
depend on $\theta$, $C'$, $\rho$, $\mu$,
and the constant $C$ in Definition~\ref{DefATI}.

Next, we show the regularity condition.
Let $\tau_1,\tau_2\in X$ satisfy
$\rho(\tau_1,\tau_2)\leq\frac{2^{-k}+\rho(x,\tau_1)}{2C_\rho}$.
We consider the following two cases for $\rho(x,\tau_1)$.

\emph{Case (1)}
$\rho(x,\tau_1)>C_\rho^2(C'+\theta)2^{-k}$.
In this case,
by Definition~\ref{D2.1}(iii) and
the fact that $y\in B_\rho(x,\theta2^{-k})$,
we conclude that
\begin{align*}
C_\rho^2(C'+\theta)2^{-k}<\rho(x,\tau_1)\leq
C_\rho\max\left\{\rho(x,y),\,\rho(y,\tau_1)\right\}
=C_\rho\rho(y,\tau_1),
\end{align*}
which further implies that $\rho(y,\tau_1)>C_\rho(C'+\theta)2^{-k}>C'2^{-k}$.
From this and the support condition of $S_{2^{-k}}(y,\cdot)$,
we infer that $S_{2^{-k}}(y,\tau_1)=0$.
If $\rho(y,\tau_2)\ge C'2^{-k}$,
then
$$
S_{2^{-k}}(y,\tau_2)=0.
$$
If $\rho(y,\tau_2)<C'2^{-k}$,
by the fact that $y\in B_\rho(x,\theta2^{-k})$,
we find that
\begin{align*}\rho(x,\tau_2)\leq
C_\rho\max\left\{\rho(x,y),\,\rho(y,\tau_2)\right\}
<C_\rho(C'+\theta)2^{-k},
\end{align*}
which, combined with Definition~\ref{D2.1}(iii),
the assumption $\rho(\tau_1,\tau_2)\leq\frac{2^{-k}+\rho(x,\tau_1)}{2C_\rho}$,
and the fact that $\rho(x,\tau_1)
>C_\rho^2(C'+\theta)2^{-k}>2^{-k}$, further implies that
\begin{align*}
C_\rho^2(C'+\theta)2^{-k}&<\rho(x,\tau_1)
\leq C_\rho\max\left\{\rho(x,\tau_2),\,\rho(\tau_2,\tau_1)\right\}
\\
&\leq\max\left\{C_\rho^2(C'+\theta)2^{-k},\,2^{-k-1}
+\frac{\rho(x,\tau_1)}{2}\right\}
\\
&=2^{-k-1}
+\frac{\rho(x,\tau_1)}{2}<\rho(x,\tau_1),
\end{align*}
which is a contradiction.
Thus, $|S_{2^{-k}}(y,\tau_1)-S_{2^{-k}}(y,\tau_2)|=0$
in this case.

\emph{Case (2)}
$\rho(x,\tau_1)\leq C_\rho^2(C'+\theta)2^{-k}$.
In this case, we have
$$
\rho(\tau_1,\tau_2)\leq\frac{2^{-k}+\rho(x,\tau_1)}{2C_\rho}\lesssim2^{-k},
$$
which, together with Lemma~\ref{A3.2}(vii), further implies that
$$
{V_\rho(x,\tau_1)}
\lesssim{(V_\rho)_{2^{-k}}(x)}
\sim{(V_\rho)_{2^{-k}}(\tau_1)}
\sim{(V_\rho)_{2^{-k}}(\tau_2)}.
$$
From this, (ii) and (iv) of Definition~\ref{DefATI},
the assumption that $\rho(x,\tau_1)\leq C_\rho^2(C'+\theta)2^{-k}$,
and $\beta<\varepsilon$,
we deduce that
\begin{align*}
&\left|S_{2^{-k}}(y,\tau_1)-S_{2^{-k}}(y,\tau_2)\right|\\
&\quad\lesssim 2^{k\varepsilon}[\rho(\tau_1,\tau_2)]^\varepsilon
\left[\frac{1}{(V_\rho)_{2^{-k}}(\tau_1)}
+\frac{1}{(V_\rho)_{2^{-k}}(\tau_2)}\right]
\\
&\quad\lesssim\left[\frac{\rho(\tau_1,\tau_2)}
{2^{-k}+\rho(x,\tau_1)}\right]^\varepsilon
\frac{1}{(V_\rho)_{2^{-k}}(x)+V_\rho(x,\tau_1)}
\left[\frac{2^{-k}}{2^{-k}+\rho(x,\tau_1)}\right]^\gamma
\\
&\quad\lesssim\left[\frac{\rho(\tau_1,\tau_2)}{2^{-k}+\rho(x,\tau_1)}\right]^\beta
\frac{1}{(V_\rho)_{2^{-k}}(x)+V_\rho(x,\tau_1)}
\left[\frac{2^{-k}}{2^{-k}+\rho(x,\tau_1)}\right]^\gamma,
\end{align*}
which completes the proof of \eqref{wza-47}. By \eqref{wza-47}
and Remark~\ref{RemMf}(iii),
we conclude that, for any
$f\in({\mathcal{G}}^\varepsilon_0(\beta,\gamma))'$ and $x\in X$,
\begin{align*}
(\mathcal{M}_\rho)_\theta f(x)=\sup_{k\in\mathbb{Z}}
\sup_{y\in B_\rho(x,\theta2^{-k})}
\left|\left\langle f,S_{2^{-k}}(y,\cdot)\right\rangle\right|
\lesssim f^*_\rho(x).
\end{align*}
This finishes the proof of the second inequality of
\eqref{2114} and hence Lemma~\ref{MMM}.
\end{proof}

Using Lemma~\ref{MMM}, we obtain the following proposition.

\begin{proposition}\label{2250}
Let $(X,\mathbf{q},\mu)$ be a quasi-ultrametric
space of homogeneous type,
$\rho\in\mathbf{q}$
have the property that all $\rho$-balls are $\mu$-measurable,
and $\varepsilon,\beta,\gamma\in(0,\infty)$ be such that
$0<\beta,\gamma<\varepsilon\preceq\mathrm{ind\,}(X,\mathbf{q})$.
Let $f\in({\mathcal{G}}^\varepsilon_0(\beta,\gamma))'$.
Then $f^*_\rho=0$ pointwise in $X$ if and
only if $f=0$ in $({\mathcal{G}}^\varepsilon_0(\beta,\gamma))'$.
\end{proposition}

To prove Proposition~\ref{2250},
we need the following lemma.

\begin{lemma}\label{96}
Let $(X,\mathbf{q},\mu)$ be a quasi-ultrametric space of homogeneous type.
Assume that $0<\varepsilon\preceq\mathrm{ind\,}(X,\mathbf{q})$
and $\{\mathcal{S}_{2^{-k}}\}_{k\in\mathbb{Z}_+}$
is a homogeneous approximation of the identity
of order $\varepsilon$ as in Definition~\ref{DefHATI}.
Then, for any $\beta\in(0,\varepsilon)$,
$\gamma\in(0,\varepsilon]$, and
$f\in\mathcal{G}^\varepsilon_0(\beta,\gamma)$
[resp. $f\in(\mathcal{G}^\varepsilon_0(\beta,\gamma))'$],
\begin{align*}
\lim_{k\to\infty}\mathcal{S}_{2^{-k}}f=f
\end{align*}
in $\mathcal{G}^\varepsilon_0(\beta,\gamma)$
[resp. $(\mathcal{G}^\varepsilon_0(\beta,\gamma))'$].
\end{lemma}

\begin{proof}
Let $\beta\in(0,\varepsilon)$,
$\gamma\in(0,\varepsilon]$, $f\in\mathcal{G}^\varepsilon_0(\beta,\gamma)$,
and $\rho\in\mathbf{q}$ be
a regularized quasi-ultrametric.
We first claim that, for any $\beta'\in(\beta,\varepsilon)$
and $\phi\in\mathcal{G}(\beta',\gamma)$,
\begin{align*}
\lim_{k\to\infty}\left\|\mathcal{S}_{2^{-k}}
\phi-\phi\right\|_{\mathcal{G}(\beta,\gamma)}=0.
\end{align*}
Let $\beta'\in(\beta,\varepsilon)$ and
$\phi\in\mathcal{G}(\beta',\gamma)\subset\mathcal{G}(\beta,\gamma)$.
Note that Definition~\ref{DefATI} and Theorem~\ref{corATI}(i)
imply that $\{\mathcal{S}_{2^{-k}}\}_{k\in\mathbb{Z_+}}$
satisfies all the conditions in Theorem~\ref{CZOonCG}
with the constant in \eqref{A13} independent of $k\in\mathbb{Z_+}$.
By this and Theorem~\ref{CZOonCG},
we find that, for any $k\in\mathbb{Z}_+$,
\begin{align}\label{1718}
\left\|\mathcal{S}_{2^{-k}}\phi\right\|_{\mathcal{G}(\beta,\gamma)}
\lesssim\left\|\phi\right\|_{\mathcal{G}(\beta,\gamma)},
\end{align}
where the implicit positive constant is independent of both $k$ and $\phi$.
Therefore, to show the above claim, it remains to prove that
\begin{align}\label{900}
\lim_{k\to\infty}\left\|\mathcal{S}_{2^{-k}}
\phi-\phi\right\|_{\mathcal{G}(\beta,\gamma)}=0.
\end{align}
To this end,
let $k\in\mathbb{N}$ be sufficiently large such that
$C'2^{-k}<(2C_\rho)^{-1}$,
where $C'$ is the same positive constant as in Definition~\ref{DefATI}.

We first show the size condition of
$\mathcal{S}_{2^{-k}}\phi-\phi$.
From (i) and (v) in Definition~\ref{DefATI}
and the choice of $k$, we deduce that, for any $x\in X$,
\begin{align*}
\left|\mathcal{S}_{2^{-k}}\phi(x)-\phi(x)\right|
&\leq\int_{X}
S_{2^{-k}}(x,y)|\phi(y)-\phi(x)|\,d\mu(y)
\\
&=\int_{\rho(x,y)\leq\frac{\rho(x_1,x)+1}{2C_\rho}}
S_{2^{-k}}(x,y)|\phi(y)-\phi(x)|\,d\mu(y)
\\
&\quad+\int_{\rho(x,y)>\frac{\rho(x_1,x)+1}{2C_\rho}}
S_{2^{-k}}(x,y)|\phi(y)-\phi(x)|\,d\mu(y)
\\
&=\int_{\rho(x,y)\leq\frac{\rho(x_1,x)+1}{2C_\rho}}
S_{2^{-k}}(x,y)|\phi(y)-\phi(x)|\,d\mu(y).
\end{align*}
Then, by the size and the support conditions of $S_{2^{-k}}$,
the regularity condition of $\phi$, and \eqref{upperdoub},
we conclude that, for any $x\in X$,
\begin{align}\label{825}
&\left|\mathcal{S}_{2^{-k}}\phi(x)-\phi(x)\right|
\nonumber\\
&\quad\lesssim\int_{\rho(x,y)\leq C'2^{-k}}
\frac{1}{V_{2^{-k}}(x)}\left[\frac{\rho(x,y)}{1+\rho(x_1,x)}\right]^{\beta'}
\frac{1}{(V_\rho)_1(x_1)+V_\rho(x_1,x)}
\left[\frac{1}{1+\rho(x_1,x)}\right]^\gamma\,d\mu(y)
\nonumber\\
&\quad\lesssim2^{-k\beta'}\int_{\rho(x,y)
\leq C'2^{-k}}
\frac{1}{V_{2^{-k}}(x)}
\frac{1}{(V_\rho)_1(x_1)+V_\rho(x_1,x)}
\left[\frac{1}{1+\rho(x_1,x)}\right]^\gamma\,d\mu(y)
\nonumber\\
&\quad\lesssim2^{-k\beta'}
\frac{1}{(V_\rho)_1(x_1)+V_\rho(x_1,x)}
\left[\frac{1}{1+\rho(x_1,x)}\right]^\gamma.
\end{align}

Now, we show the regularity condition of
$\mathcal{S}_{2^{-k}}\phi-\phi$.
Let $x,x'\in X$ satisfy $\rho(x,x')\leq\frac{\rho(x_1,x)+1}{2C_\rho}$.
On the one hand, from \eqref{1718} and
the regularity conditions of
both $\phi$ and $\mathcal{S}_{2^{-k}}\phi$,
we infer that
\begin{align}\label{833}
&\left|\left[\mathcal{S}_{2^{-k}}\phi(x)-\phi(x)\right]
-\left[\mathcal{S}_{2^{-k}}\phi(x')-\phi(x')\right]\right|
\nonumber\\
&\quad\leq\left|\mathcal{S}_{2^{-k}}\phi(x)
-\mathcal{S}_{2^{-k}}\phi(x')\right|
+|\phi(x)-\phi(x')|
\nonumber\\
&\quad\lesssim
\left[\frac{\rho(x,x')}{1+\rho(x_1,x)}\right]^{\beta'}
\frac{1}{(V_\rho)_1(x_1)+V_\rho(x_1,x)}
\left[\frac{1}{1+\rho(x_1,x)}\right]^\gamma.
\end{align}
On the other hand, by $\rho(x,x')\leq\frac{\rho(x_1,x)+1}{2C_\rho}$,
we find that
\begin{align*}
1+\rho(x_1,x)\sim 1+\rho(x_1,x')
\ \ \text{and}\ \
(V_\rho)_1(x_1)+V_\rho(x_1,x)
\sim
(V_\rho)_1(x_1)+V_\rho(x_1,x'),
\end{align*}
which, combined with \eqref{825}, further implies that
\begin{align*}
\left|\left[\mathcal{S}_{2^{-k}}\phi(x)-\phi(x)\right]
-\left[\mathcal{S}_{2^{-k}}\phi(x')-\phi(x')\right]\right|
\lesssim2^{-k\beta'}
\frac{1}{(V_\rho)_1(x_1)+V_\rho(x_1,x)}
\left[\frac{1}{1+\rho(x_1,x)}\right]^\gamma.
\end{align*}
From this and \eqref{833},
we deduce that, for any $\sigma\in(0,1)$,
\begin{align*}
&\left|\left[\mathcal{S}_{2^{-k}}\phi(x)-\phi(x)\right]
-\left[\mathcal{S}_{2^{-k}}\phi(x')-\phi(x')\right]\right|
\\
&\quad\lesssim2^{-k(\beta'-\beta)}
\left[\frac{\rho(x,x')}{1+\rho(x_1,x)}\right]^\beta
\frac{1}{(V_\rho)_1(x_1)+V_\rho(x_1,x)}
\left[\frac{1}{1+\rho(x_1,x)}\right]^\gamma,
\end{align*}
which completes the proof of the regularity condition
of $\mathcal{S}_{2^{-k}}\phi-\phi$. These estimates imply
$$
\left\|\mathcal{S}_{2^{-k}}
\phi-\phi\right\|_{\mathcal{G}(\beta,\gamma)}\lesssim2^{-k(\beta'-\beta)}
\rightarrow0
$$
as $k\to\infty$, which finishes the proof of \eqref{900},
and hence the claim.

Let $f\in\mathcal{G}^\varepsilon_0(\beta,\gamma)$. Then
there exists a sequence $\{f_j\}_{j\in\mathbb{N}}$ in
$\mathcal{G}(\varepsilon,\varepsilon)\subset\mathcal{G}(\beta',\gamma)$ such that
$$
\lim_{j\to\infty}\left\|f-f_j\right\|_{\mathcal{G}(\beta,\gamma)}=0.
$$
By this and \eqref{1718}, we conclude that
\begin{align*}
\|\mathcal{S}_{2^{-k}}f-f\|_{\mathcal{G}(\beta,\gamma)}
&\leq\|\mathcal{S}_{2^{-k}}f-\mathcal{S}_{2^{-k}}f_j\|_{\mathcal{G}(\beta,\gamma)}
+\|\mathcal{S}_{2^{-k}}f_j-f_j\|_{\mathcal{G}(\beta,\gamma)}
+\|f_j-f\|_{\mathcal{G}(\beta,\gamma)}
\\
&\lesssim\|f_j-f\|_{\mathcal{G}(\beta,\gamma)}
+\|\mathcal{S}_{2^{-k}}f_j-f_j\|_{\mathcal{G}(\beta,\gamma)}\to0
\end{align*}
via first letting $k\to\infty$ and then letting $j\to\infty$,
which, together with a standard duality argument,
further implies that, for any $g\in(\mathcal{G}^\varepsilon_0(\beta,\gamma))'$,
\begin{align*}
\lim_{k\to\infty}\left\|\mathcal{S}_{2^{-k}}g-g
\right\|_{(\mathcal{G}^\varepsilon_0(\beta,\gamma))'}=0.
\end{align*}
This finishes the proof of Lemma~\ref{96}.
\end{proof}

\begin{proof}[Proof of Proposition~\ref{2250}]
The sufficiency is trivial and we only prove the necessity.
Let $f\in({\mathcal{G}}^\varepsilon_0(\beta,\gamma))'$ and
assume $f^*_\rho(x)=0$ for any $x\in X$.
From Lemma~\ref{MMM}, we infer that, for any $x\in X$,
$$
\mathcal{M}^+f(x)\lesssim f^*_\rho(x)=0,
$$
which further implies that, for any $k\in\mathbb{Z}$,
$\mathcal{S}_{2^{-k}}f(x)=0$. This and Lemma~\ref{96} imply that
$$
f=\lim_{k\to\infty}\mathcal{S}_{2^{-k}}f=0
$$
in $({\mathcal{G}}^\varepsilon_0(\beta,\gamma))'$,
which completes the proof of Proposition~\ref{2250}.
\end{proof}

The following lemma details the action of the
Hardy--Littlewood maximal operator on Lorentz spaces and
it is a slight generalization of \cite[Lemma 3.5]{zhy}
wherein Zhou et al. \cite[Lemma 3.5]{zhy} considered symmetric quasi-metrics in
unbounded spaces of homogeneous type where $\mu(\{x\})=0$ for all $x\in X$.

\begin{lemma}
\label{MLpq}
Let $(X,\mathbf{q},\mu)$ be a quasi-ultrametric
space of homogeneous type and $\rho\in\mathbf{q}$
a regularized quasi-ultrametric.
Then, for any $p\in(1,\infty]$ and $q\in(0,\infty]$,
$$
\mathcal{M}_{\rho}: L^{p,q}(X)\to L^{p,q}(X)
$$
is well defined and bounded.
\end{lemma}

\begin{proof}
Let $E\subset X$ be an arbitrary
$\mu$-measurable set satisfying $\mu(E)<\infty$.
By Lemma~\ref{M}, we find that
$\mathcal{M}_{\rho}$ is bounded on $L^\infty(X)$ and from $L^1(X)$
to $L^{1,\infty}(X)$, which implies
$$
\left\|\mathcal{M}_\rho{\mathbf{1}}_E\right\|_{L^{1,\infty}(X)}
\lesssim \left\|{\mathbf{1}}_E\right\|_{L^{1}(X)}=\mu(E)
$$
and
$$
\left\|\mathcal{M}_\rho{\mathbf{1}}_E\right\|_{L^{\infty,\infty}(X)}
=\left\|\mathcal{M}_\rho{\mathbf{1}}_E\right\|_{L^{\infty}(X)}
\lesssim\|{\mathbf{1}}_E\|_{L^{\infty}(X)}=1.
$$
Therefore, it follows from \cite[Theorem~1.1 and Remark~1.4]{YLD}
(used here with $p_0=q_0=1$ and $p_1=q_1=\infty$)
that there exists a positive constant $C$ such that,
for any $f\in L^{p,q}(X)$,
$$
\|\mathcal{M}_\rho f\|_{L^{p,q}(X)}\leq C\|f\|_{L^{p,q}(X)}.
$$
This finishes the proof of Lemma~\ref{MLpq}.
\end{proof}

\begin{lemma}
\label{fstarMf}
Let $(X,\mathbf{q},\mu)$ be an ultra-RD-space with
upper dimension $n\in(0,\infty)$,
where $\mu$ is assumed to be a Borel-semiregular
measure. Assume that $\mathrm{diam\,}(X)=\infty$.
Let $\rho\in\mathbf{q}$ be a regularized
quasi-ultrametric.
Let $p\in(\frac{n}{n+\mathrm{ind\,}(X,\mathbf{q})},\infty)$ and
$\beta,\widetilde{\beta},\gamma,\varepsilon\in(0,\infty)$ satisfy
$n(\frac{1}{p}-1)_+<\widetilde{\beta}<\beta,\gamma<
\varepsilon\preceq\mathrm{ind\,}(X,\mathbf{q})$.
Then there exists a positive constant $C$ such that,
for any given $r\in(\frac{n}{n+\widetilde{\beta}},p)\cap(0,1]$
and for any
$f\in({\mathcal{G}}^\varepsilon_0(\beta,\gamma))'$ and $x\in X$,
$$
f^*_\rho(x)\leq C\mathcal{M}^+f(x)
+C\left[\mathcal{M}_\rho\left((\mathcal{M}^+f)^r\right)(x)\right]^\frac{1}{r}.
$$
\end{lemma}

\begin{proof}
Assume that $\{\mathcal{S}_{2^{-k}}\}_{k\in\mathbb{Z}}$ is
a homogeneous approximation of the identity
of order $\varepsilon$ as in Definition~\ref{DefHATI}.
Let $f\in({\mathcal{G}}^\varepsilon_0(\beta,\gamma))'$
and $x\in X$.
For any $k\in\mathbb{Z}$, define
$$
\mathcal{D}_k:=\mathcal{S}_{2^{-k}}-\mathcal{S}_{2^{-k+1}}.
$$
We first claim that,
for any
$k\in\mathbb{Z}$ and
$\phi\in\mathring{{\mathcal{G}}}^\varepsilon_0(\beta,\gamma)$
with $\|\phi\|_{{\mathcal{G}}(x,2^{-k},\beta,\gamma)}\leq1$,
\begin{equation}\label{fphi}
|\langle f,\phi\rangle|^r\lesssim
\mathcal{M}_\rho\left(\left[\mathcal{M}^+
f\right]^r\right)(x).
\end{equation}
Suppose $k\in\mathbb{Z}$ and
$\phi\in\mathring{{\mathcal{G}}}^\varepsilon_0(\beta,\gamma)$
with $\|\phi\|_{{\mathcal{G}}(x,2^{-k},\beta,\gamma)}\leq1$.
Then, for any given $\delta\in(0,\infty)$ and
for any $\tau\in I_k$ and $v\in\{1,\ldots,N(k,\tau)\}$,
we can choose $y_\tau^{k,v}\in Q_\tau^{k,v}$ such that
\begin{align*}
\left|\mathcal{D}_kf(y_\tau^{k,v})\right|
&\leq\inf_{z\in Q_\tau^{k,v}}\left|\mathcal{D}_kf(z)\right|+\delta\\
&\leq\inf_{z\in Q_\tau^{k,v}}\left[\left|\mathcal{S}_{2^{-k}}f(z)\right|
+\left|\mathcal{S}_{2^{-k+1}}f(z)\right|\right]+\delta
\\
&\leq2\inf_{z\in Q_\tau^{k,v}}\mathcal{M}^+f(z)+\delta,
\end{align*}
where $Q_\tau^{k,v}$ is as in Remark~\ref{jcubes}.
Let $g:=f|_{\mathring{{\mathcal{G}}}^\varepsilon_0(\beta,\gamma)}$
be the restriction of $f$ to the set
$\mathring{{\mathcal{G}}}^\varepsilon_0(\beta,\gamma)$.
Then $g\in(\mathring{{\mathcal{G}}}^\varepsilon_0(\beta,\gamma))'$
and $\|g\|_{(\mathring{{\mathcal{G}}}^\varepsilon_0(\beta,\gamma))'}
\leq\|f\|_{({\mathcal{G}}^\varepsilon_0(\beta,\gamma))'}$.
From Lemma~\ref{Pk}
and Definition~\ref{DefATI}(v), we deduce that
$D_k(x,\cdot)\in\mathring{{\mathcal{G}}}(x,2^{-k},\varepsilon,\varepsilon)$.
Moreover, by Theorem~\ref{HD} [keeping in mind
that we are assuming $\mathrm{diam}_\rho(X)=\infty$], we find that
there exists a family of linear operators
$\{\widetilde{\mathcal{D}}_{l}\}_{l\in\mathbb{Z}}$ such that
\begin{align}
\label{930}
\langle f,\phi\rangle
&=\langle g,\phi\rangle
=\left\langle\sum_{l=-\infty}^\infty\sum_{\tau\in I_l}
\sum_{v=1}^{N(\tau,l)}\mu(Q_\tau^{l,v})
\widetilde{D}_l(\cdot,y_{\tau}^{l,v})
\mathcal{D}_lg(y_\tau^{l,v}),\phi\right\rangle
\nonumber\\
&=\sum_{l=-\infty}^\infty\sum_{\tau\in I_l}
\sum_{v=1}^{N(\tau,l)}\mu(Q_\tau^{l,v})\left\langle
\widetilde{D}_l(\cdot,y_{\tau}^{l,v})
,\phi\right\rangle\mathcal{D}_lg(y_\tau^{l,v})
\nonumber\\
&=\sum_{l=-\infty}^\infty\sum_{\tau\in I_l}
\sum_{v=1}^{N(\tau,l)}\mu(Q_\tau^{l,v})
\widetilde{\mathcal{D}}_l^*\phi(y_\tau^{l,v})
\mathcal{D}_lg(y_\tau^{l,v})
\nonumber\\
&=\sum_{l=-\infty}^\infty\sum_{\tau\in I_l}
\sum_{v=1}^{N(\tau,l)}\mu(Q_\tau^{l,v})
\widetilde{\mathcal{D}}_l^*\phi(y_\tau^{l,v})
\mathcal{D}_lf(y_\tau^{l,v}),
\end{align}
where $\widetilde{\mathcal{D}}_l^*$ denotes the adjoint operator
of $\widetilde{\mathcal{D}}_l$. From Lemma~\ref{1122} with
both $\overline{D}_k$ and $\widetilde{D}_l$
replaced, respectively, by
$\widetilde{D}_{l}^*(y_\tau^{l,v},\cdot)\in
\mathring{\mathcal{G}}(y_\tau^{l,v},2^{-l},\beta,\gamma)$ and
$\phi\in\mathring{\mathcal{G}}(x,2^{-k},\beta,\gamma)$,
we infer that, for any fixed
$\widetilde{\beta}\in(0,\beta\wedge\gamma)$ and for any $l\in\mathbb{Z}$,
\begin{align*}
\left|\widetilde{\mathcal{D}}_l^*\phi(y_\tau^{l,v})\right|
\lesssim2^{-|k-l|\widetilde{\beta}}\frac{1}{(V_\rho)_{2^{-(k\wedge l)}}(x)+
V_\rho(x,y_\tau^{l,v})}\left[\frac{2^{-(k\wedge l)}}
{2^{-(k\wedge l)}+\rho(x,y_\tau^{l,v})}\right]^\gamma.
\end{align*}
By this, \eqref{930}, and Lemma~\ref{L4.5},
we obtain
\begin{align*}
\left|\langle f,\phi\rangle\right|
&\lesssim\sum_{l=-\infty}^\infty2^{-|k-l|\widetilde{\beta}}
\sum_{\tau\in I_l}\sum_{v=1}^{N(l,\tau)}\mu(Q_\tau^{l,v})
\frac{\inf_{z\in Q_\tau^{l,v}}\mathcal{M}^+f(z)
+\delta}{(V_\rho)_{2^{-(k\wedge l)}}(x)+
V_\rho(x,y_\tau^{l,v})}\left[\frac{2^{-(k\wedge l)}}
{2^{-(k\wedge l)}+\rho(x,y_\tau^{l,v})}\right]^\gamma
\nonumber\\
&\lesssim\sum_{l=-\infty}^\infty2^{-|k-l|\widetilde{\beta}+
n(k\wedge l-l)(1-\frac{1}{r})}
\left[\mathcal{M}_\rho\left(\sum_{\tau\in I_l}
\sum_{v=1}^{N(l,\tau)}\left[\inf_{z\in Q_\tau^{l,v}}\mathcal{M}^+f(z)
+\delta\right]^r\mathbf{1}_{Q_\tau^{l,v}}\right)(x)\right]^\frac{1}{r},
\end{align*}
which, combined with Lemma~\ref{DyadicSys}(i)
and $r>\frac{n}{n+\widetilde{\beta}}$,
further implies that
\begin{align}\label{1999}
\left|\langle f,\phi\rangle\right|
&\lesssim\sum_{l=-\infty}^\infty
2^{-|k-l|\widetilde{\beta}+n(k\wedge l-l)(1-\frac{1}{r})}
\left[\mathcal{M}_\rho\left(\left(\mathcal{M}^+f
+\delta\right)^r\right)(x)\right]^\frac{1}{r}
\nonumber\\
&\sim\left[\mathcal{M}_\rho\left(\left(\mathcal{M}^+f
\right)^r\right)(x)\right]^\frac{1}{r}+\delta
\to\left[\mathcal{M}_\rho\left(\left(\mathcal{M}^+f
\right)^r\right)(x)\right]^\frac{1}{r}
\end{align}
as $\delta\to0^+$.
This finishes the the proof of \eqref{fphi}.

Next, let $\psi\in{\mathcal{G}}^\varepsilon_0(\beta,\gamma)$
satisfy $\|\psi\|_{{\mathcal{G}}(x,r_0,\beta,\gamma)}\leq 1$
for some $r_0\in(0,\infty)$ and let $k_0\in\mathbb{Z}$ be
such that $2^{-(k_0+1)}\leq r_0<2^{-k_0}$.
Then it follows from Remark~\ref{Rmtest}(i)
that $\|\psi\|_{{\mathcal{G}}(x,2^{-k_0},\beta,\gamma)}\lesssim 1$,
where the implicit positive constant is independent of $x$ and $k_0$.
Let
$$
\sigma:=\int_X\psi(z)\,d\mu(z)
$$
(which is a well-defined
number in light of Lemma~\ref{CGL}) and, for any $y\in X$, let
$$
\phi(y):=\psi(y)-\sigma S_{2^{-k_0}}(x,y).
$$
Using Lemma~\ref{A3.2}(ii) and Definition~\ref{testfunction}(i),
we conclude that $|\sigma|\lesssim1$.
From this, Lemma~\ref{ATItest},
${\mathcal{G}}(x,2^{-k_0},\varepsilon,\varepsilon)
\subset{\mathcal{G}}(x,2^{-k_0},\beta,\gamma)$,
and Definition~\ref{DefATI}(v),
we deduce that $\phi\in\mathring{{\mathcal{G}}}^\varepsilon_0(\beta,\gamma)$
and
\begin{align*}
\|\phi\|_{{\mathcal{G}}(x,2^{-k_0},\beta,\gamma)}
\leq\|\psi\|_{{\mathcal{G}}(x,2^{-k_0},\beta,\gamma)}
+|\sigma|\,\|S_{2^{-k_0}}(x,\cdot)\|_
{{\mathcal{G}}(x,2^{-k_0},\beta,\gamma)}\lesssim 1,
\end{align*}
where the implicit positive constant is independent of
$x$ and $k_0$.
By this and \eqref{fphi},
we find that
\begin{align*}
|\langle f,\psi\rangle|
&\leq|\langle f,\phi\rangle|
+|\langle f,\sigma S_{2^{-k}}(x,\cdot)\rangle|
\\
&\lesssim\left[\mathcal{M}_\rho\left(
(\mathcal{M}^+f)^r\right)(x)\right]
^\frac{1}{r}+|\sigma|\left|\mathcal{S}_{2^{-k}}f(x)\right|
\\
&\leq\left[\mathcal{M}_\rho\left(
(\mathcal{M}^+f)^r\right)(x)\right]
^\frac{1}{r}+|\sigma|\mathcal{M}^+f(x).
\end{align*}
Taking the supremum of all $\psi$
as in Definition~\ref{maxfunc}(iii)
and using Remark~\ref{RemMf}(iii) then complete the proof of Lemma~\ref{fstarMf}.
\end{proof}

The stage has now been set for us to
present the main theorem of this section,
which shows that the three brands of Hardy--Lorentz spaces
defined in Definition~\ref{DefHardyL}
coincide with equivalent quasi-norms. The following theorem
focuses on the case where $\mathrm{diam\,}(X)=\infty$
which is equivalent to $\mu(X)=\infty$;
see Theorem~\ref{HL-Mfinite} below for the case
where $\mathrm{diam\,}(X)<\infty$.

\begin{theorem}
\label{HL-M}
Let $(X,\mathbf{q},\mu)$ be an ultra-RD-space with
upper dimension $n\in(0,\infty)$,
where $\mu$ is assumed to be a Borel-semiregular
measure. Assume that $\mathrm{diam\,}(X)=\infty$.
Let $p\in(\frac{n}{n+\mathrm{ind\,}(X,\mathbf{q})},\infty)$,
$q\in(0,\infty]$, $\theta\in(0,\infty)$,
and $\beta,\gamma,\varepsilon\in(0,\infty)$ satisfy
$n(\frac{1}{p}-1)_+<\beta,\gamma<\varepsilon\preceq\mathrm{ind\,}(X,\mathbf{q})$.
Then, as subspaces of $({\mathcal{G}}^\varepsilon_0(\beta,\gamma))'$,
$$
H^{p,q}_+(X)
=H^{p,q}_\theta(X)
=H^{p,q}_*(X)
$$
with equivalent quasi-norms.
\end{theorem}

\begin{proof}
Without loss of generality, we may assume that
$\rho\in\mathbf{q}$ is a regularized quasi-ultrametric.
By Lemma~\ref{MMM}, Proposition~\ref{inc-norm}(ii),
and Definition~\ref{DefHardyL}, we conclude that
$$
H^{p,q}_*(X)\subset H^{p,q}_\theta(X)\subset H^{p,q}_+(X),
$$
where, for any $f\in ({\mathcal{G}}^\varepsilon_0(\beta,\gamma))'$,
\begin{align}\label{dqp-38}
\left\|\mathcal{M}^+f\right\|_{L^{p,q}(X)}
\lesssim\left\|(\mathcal{M}_\rho)_\theta f\right\|_{L^{p,q}(X)}
\lesssim\left\|f^*_\rho\right\|_{L^{p,q}(X)}.
\end{align}
In light of this, to see that these three spaces coincide (with equivalent
quasi-norms), it suffice to prove
that $H^{p,q}_+(X)\subset H^{p,q}_*(X)$. To this end, let
$f\in H^{p,q}_+(X)\subset ({\mathcal{G}}^\varepsilon_0(\beta,\gamma))'$,
$\widetilde{\beta}\in(n(\frac{1}{p}-1)_+,\beta\wedge\gamma)$,
and $r\in(\frac{n}{n+\widetilde{\beta}},p)\cap(0,1]$.
Then $\mathcal{M}^+f\in L^{p,q}(X)$ and hence
\begin{align}\label{221111}
\left\|f^*_\rho\right\|_{L^{p,q}(X)}
&\lesssim\left\|\mathcal{M}^+f+\left[\mathcal{M}_\rho
\left((\mathcal{M}^+f)^r\right)
\right]^\frac{1}{r}\right\|_{L^{p,q}(X)}
\quad\text{by Lemma~\ref{fstarMf}}
\nonumber\\
&\sim\left\|\mathcal{M}^+f\right\|_{L^{p,q}(X)}
+\left\|\mathcal{M}_\rho\left((\mathcal{M}^+f
)^r\right)\right\|^\frac{1}{r}_{L^{\frac{p}{r},\frac{q}{r}}(X)}
\quad\text{by Proposition~\ref{r;p,q=pr,qr;r}}
\nonumber\\
&\lesssim\left\|\mathcal{M}^+f\right\|_{L^{p,q}(X)}
+\left\|(\mathcal{M}^+f
)^r\right\|^\frac{1}{r}_{L^{\frac{p}{r},\frac{q}{r}}(X)}
\quad\text{by Lemma~\ref{MLpq}}\nonumber\\
&\sim\left\|\mathcal{M}^+f\right\|_{L^{p,q}(X)}
\quad\text{by Proposition~\ref{r;p,q=pr,qr;r}},
\end{align}
which implies $f\in H^{p,q}_*(X)$. Moreover, this estimate,
together with \eqref{dqp-38},
further implies that, for any $f\in({\mathcal{G}}^\varepsilon_0(\beta,\gamma))'$,
\begin{align*}
\left\|\mathcal{M}^+f\right\|_{L^{p,q}(X)}
\sim\left\|(\mathcal{M}_\rho)_\theta f\right\|_{L^{p,q}(X)}
\sim\left\|f^*_\rho\right\|_{L^{p,q}(X)}.
\end{align*}
This finishes the proof of Theorem~\ref{HL-M}.
\end{proof}

We conclude this section by recording a few
results that is useful in the sequel. First,
with Proposition~\ref{distfncttype} in hand, we easily obtain the
following pointwise relation between the grand maximal
and Hardy--Littlewood maximal functions.

\begin{proposition}\label{HL-grandmax}
Let $(X,\mathbf{q},\mu)$ be a
quasi-ultrametric spaces of homogeneous type, $\rho\in\mathbf{q}$
be a regularized quasi-ultrametric, $r\in[1,\infty]$,
and $\varepsilon,\beta,\gamma\in(0,\infty)$ be such that
$0<\beta,\gamma<\varepsilon
\preceq\mathrm{ind\,}(X,\mathbf{q})$.
Then there exists a positive constant $C$ such that,
for any $f\in L^r(X)$ and $x\in X$,
\begin{align}
\label{*<M}
f^*_\rho(x)\leq C\mathcal{M}_\rho f(x).
\end{align}
\end{proposition}

\begin{proof}
Assume that $\phi\in{\mathcal{G}}^\varepsilon_0(\beta,\gamma)$ satisfies
$\|\phi\|_{G(x,r_0^\alpha,\frac{\beta}{\alpha},\frac{\gamma}{\alpha})}\leq 1$
for some $r_0\in(0,\infty)$, where $\alpha:=\min\{1,\,(\log_2C_\rho)^{-1}\}$.
From Remark~\ref{RemMf}(iii), we infer that
$$
\|\phi\|_{{\mathcal{G}}(x,r_0,\beta,\gamma)}\lesssim1.
$$
By this and an argument similar to that used in the
estimation of \eqref{alk-34}, we find that
\begin{align*}
\left|\left\langle f,\varphi\right\rangle\right|
&\leq\int_{\rho(x,y)<r_0}|f(y)\phi(y)|\,d\mu(y)
+\int_{\rho(x,y)\geq r_0}|f(y)\phi(y)|\,d\mu(y)
\nonumber\\
&\lesssim\frac{1}{(V_\rho)_{r_0}(x)}\int_{\rho(x,y)<r_0}|f(y)|\,d\mu(y)
\nonumber\\
&\quad+\int_{\rho(x,y)\geq r_0}|f(y)|\frac{1}{V_\rho(x,y)}
\left[\frac{1}{\rho(x,y)}\right]^\gamma\,d\mu(y)
\nonumber\\
&\lesssim\mathcal{M}_\rho f(x),
\end{align*}
where the implicit positive constants are independent of $f$, $x$, and $\phi$.
Taking the supremum of all $\phi$
as above then completes the proof of \eqref{*<M}
and hence Proposition~\ref{HL-grandmax}.
\end{proof}

Next, we show that
$H_*^{p,q}(X)$ is independent of the
choices of parameters $\beta$ and $\gamma$
that intervene in the definition of the spaces of distributions
$({\mathcal{G}}^\varepsilon_0(\beta,\gamma))'$;
see Theorem~\ref{indep2} for the case $\mu(X)=\infty$
and Theorem~\ref{indep3} for the case $\mu(X)<\infty$.

\begin{theorem}\label{indep2}
Let $(X,\mathbf{q},\mu)$ be an ultra-RD-space with upper dimension $n\in(0,\infty)$,
where $\mu$ is assumed to be a Borel-semiregular
measure. Assume that $\mu(X)=\infty$. Let
$p\in(\frac{n}{n+\mathrm{ind\,}(X,\mathbf{q})},\infty)$,
$q\in(0,\infty]$, and $n(\frac{1}{p}-1)_+<
\varepsilon\preceq\mathrm{ind\,}(X,\mathbf{q})$.
Then $H_*^{p,q}(X)$ is independent of the
choice of parameters $\beta$ and $\gamma$
that intervene in the definition of the spaces of distributions
$({\mathcal{G}}^\varepsilon_0(\beta,\gamma))'$ whenever
$\beta,\gamma\in(0,\infty)$ satisfy
$n(\frac{1}{p}-1)_+<\beta,\gamma<\varepsilon$.
\end{theorem}

\begin{proof}
Let $\rho\in\mathbf{q}$ be a regularized quasi-ultrametric.
Let $\beta_1,\gamma_1,\beta_2,\gamma_2\in(0,\infty)$ satisfy
$$
n\left(\frac{1}{p}-1\right)_+
<\beta_1,\gamma_1,\beta_2,\gamma_2<\varepsilon.
$$
Assume that
$f\in H^{p,q}_*(X)\subset(\mathcal{G}^\varepsilon_0(\beta_1,\gamma_1))'$.
We first claim that, for any $\phi\in\mathcal{G}(\varepsilon,\varepsilon)$ with
$\|\phi\|_{\mathcal{G}(\beta_2,\gamma_2)}\leq1$,
\begin{align}\label{1515}
\left|\left\langle f,\phi\right\rangle\right|
\lesssim\left\|(\mathcal{M}_\rho)_\theta f\right\|_{L^{p,q}(X)},
\end{align}
where $\theta:=C_\rho C_r2^{-j}$ with
$C_\rho\in[1,\infty)$ as in \eqref{Crho},
$C_r\in(0,\infty)$ as in Lemma~\ref{DyadicSys}(iv),
and $j\in\mathbb{N}$ as in Remark~\ref{jcubes}.
From Theorem~\ref{ID}, we deduce that
\begin{align}\label{1612}
\langle f,\phi\rangle
&=\left\langle\sum_{k=0}^{N}\sum_{\tau\in I_k}\sum_{v=1}^{N(k,\tau)}
\int_{Q^{k,v}_\tau }
\widetilde{D}^{(1)}_k(\cdot,y)\,d\mu(y)
\mathcal{D}^{k,v}_{\tau,1}f,\phi\right\rangle
\nonumber\\
&\quad+\left\langle\sum_{k=N+1}^{\infty}\sum_{\tau\in I_k}
\sum_{v=1}^{N(k,\tau)}\int_{Q^{k,v}_\tau}
\widetilde{D}^{(1)}_k(\cdot,y)\,d\mu(y)
\mathcal{D}_kf(y^{k,v}_\tau),\phi\right\rangle
\nonumber\\
&=\sum_{k=0}^{N}\sum_{\tau\in I_k}\sum_{v=1}^{N(k,\tau)}
\int_{Q^{k,v}_\tau}
\left(\widetilde{\mathcal{D}}^{(1)}_k\right)^*\phi(y)\,d\mu(y)
\mathcal{D}^{k,v}_{\tau,1}f
\nonumber\\
&\quad+\sum_{k=N+1}^{\infty}\sum_{\tau\in I_k}
\sum_{v=1}^{N(k,\tau)}\int_{Q^{k,v}_\tau}
\left(\widetilde{\mathcal{D}}^{(1)}_k\right)^*\phi(y)\,d\mu(y)
\mathcal{D}_kf(y^{k,v}_\tau)
\nonumber\\
&=:\mathrm{Z}_1+\mathrm{Z}_2.
\end{align}
By Lemma~\ref{DyadicSys}(iv) and the fact that the cube $Q_\tau^{k,v}$ is
at level $k+j$ (see Remark~\ref{jcubes}), we conclude that,
for any $z\in Q_\tau^{k,v}$,
\begin{align}\label{1556}
Q_\tau^{k,v}\subset B_\rho(z_\tau^{k,v},C_r2^{-k-j})
\subset B_\rho(z,C_\rho C_r2^{-k-j})
=B_\rho(z,\theta2^{-k}).
\end{align}
Fix $x\in B_\rho(x_1,1)$.
Then
$$
\|\phi\|_{\mathcal{G}(x,1,\beta_2,\gamma_2)}
\sim\|\phi\|_{\mathcal{G}(x_1,1,\beta_2,\gamma_2)}\leq1.
$$
If $k\in\mathbb{Z_+}\cap[0,N]$,
then
$$
\|\phi\|_{\mathcal{G}(x,2^{-k},\beta_2,\gamma_2)}
\sim\|\phi\|_{\mathcal{G}(x,1,\beta_2,\gamma_2)}\lesssim1,
$$
where the implicit positive constants are independent of $x$
but may depend on $N$.
Let $\beta_-:=\min\{\beta_1,\,\gamma_1,\,\beta_2,\,\gamma_2\}$.
Since $\widetilde{D}_k^{(1)}(\cdot,y)\in\mathcal{G}(y,2^{-k},\beta_1,\gamma_1)$
and $\phi\in\mathcal{G}(x,1,\beta_2,\gamma_2)$,
from Definition~\ref{testfunction}(i) and
the fact that $2^{-k}\sim1$ with the positive
equivalence constants depending only on $N$,
we infer that, for any $y\in Q_\tau^{k,v}$,
\begin{align*}
\left|\left(\widetilde{\mathcal{D}}^{(1)}_k\right)^*\phi(y)\right|
&=\left|\int_X\widetilde{D}_k^{(1)}(z,y)\phi(z)\,d\mu(z)\right|
\\
&\lesssim\int_X\frac{1}{(V_\rho)_{2^{-k}+\rho(z,y)}(y)}
\left[\frac{2^{-k}}{2^{-k}+\rho(z,y)}\right]^{\gamma_1}
\frac{1}{(V_\rho)_{1+\rho(x,z)}(x)}\left[\frac{1}{1+
\rho(x,z)}\right]^{\gamma_2}\,d\mu(z)
\\
&\sim\int_X\frac{1}{(V_\rho)_{1+\rho(z,y)}(y)}
\left[\frac{1}{1+\rho(z,y)}\right]^{\gamma_1}\frac{1}{(V_\rho)_{1+\rho(x,z)}(x)}
\left[\frac{1}{1+\rho(x,z)}\right]^{\gamma_2}\,d\mu(z).
\end{align*}
By this, Definition~\ref{D2.1}(iii),
and both (vii) and (ii) of Lemma~\ref{A3.2},
we find that, for any $y\in Q_\tau^{k,v}$,
\begin{align*}
\left|\left(\widetilde{\mathcal{D}}^{(1)}_k\right)^*\phi(y)\right|
&\lesssim\frac{1}{(V_\rho)_{1+\frac{\rho(x,y)}{C_\rho}}(x)}
\left[\frac{1}{1+\frac{\rho(x,y)}{C_\rho}}\right]^{\gamma_2}\\
&\quad\times\int_{\{z\in X:\rho(x,z)\ge\rho(z,y)\}}
\frac{1}{(V_\rho)_{1+\rho(z,y)}(y)}
\left[\frac{1}{1+\rho(z,y)}\right]^{\gamma_1}\,d\mu(z)\\
&\quad+\frac{1}{(V_\rho)_{1+\frac{\rho(x,y)}{C_\rho}}(y)}
\left[\frac{1}{1+\frac{\rho(x,y)}{C_\rho}}\right]^{\gamma_1}\\
&\quad\times\int_{\{z\in X:\rho(x,z)<\rho(z,y)\}}
\frac{1}{(V_\rho)_{1+\rho(x,z)}(x)}
\left[\frac{1}{1+\rho(x,z)}\right]^{\gamma_2}\,d\mu(z)\\
&\lesssim\frac{1}{(V_\rho)_1(x)+V_\rho(x,y)}
\left[\frac{1}{1+\rho(x,y)}\right]^{\beta_-},
\end{align*}
which, combined with the fact that
\begin{align}\label{949}
1+\rho(x,y_\tau^{k,v})&\leq1+C_\rho^2\left[\rho(x,y)
+\rho(y,z_\tau^{k,v})+\rho(z_\tau^{k,v},y_\tau^{k,v})\right]
\nonumber\\
&\lesssim1+\rho(x,y)+2^{-k}\sim1+\rho(x,y),
\end{align}
further implies that
\begin{align}\label{15}
\left|\left(\widetilde{\mathcal{D}}^{(1)}_k\right)^*\phi(y)\right|
\lesssim\frac{1}{(V_\rho)_1(x)+V_\rho(x,y_\tau^{k,v})}
\left[\frac{1}{1+\rho(x,y_\tau^{k,v})}\right]^{\beta_-}.
\end{align}
In addition, note that, for any $k\in\mathbb{Z_+}\cap[0,N]$
and $z\in Q_\tau^{k,v}$,
\begin{align}\label{2031}
\left|\mathcal{D}^{k,v}_{\tau,1}f\right|
&=\left|\fint_{Q_\tau^{k,v}}\int_XD_k(x,y)f(y)\,d\mu(y)\,d\mu(x)\right|
\quad\text{by Fubini's theorem}
\nonumber\\
&\leq\fint_{Q_\tau^{k,v}}\left|\mathcal{S}_{2^{-k}}f(x)\right|
+\left|\mathcal{S}_{2^{-k+1}}f(x)\right|\,d\mu(x)
\quad\text{since $\mathcal{D}_k=\mathcal{S}_{2^{-k}}-\mathcal{S}_{2^{-k+1}}$}
\nonumber\\
&\leq 2(\mathcal{M}_\rho)_\theta f(z)
\quad\text{using \eqref{1556} and Definition~\ref{maxfunc}(ii)},
\end{align}
which, together with \eqref{15},
further implies that
\begin{align}\label{15333}
\left|\mathrm{Z}_1\right|
&\lesssim\sum_{k=0}^{N}\sum_{\tau\in I_k}
\sum_{v=1}^{N(k,\tau)}\mu(Q_\tau^{k,v})
\frac{1}{(V_\rho)_1(x)+V_\rho(x,y_\tau^{k,v})}
\nonumber\\
&\quad\times\left[\frac{1}{1+\rho(x,y_\tau^{k,v})}\right]^{\beta_-}
\inf_{z\in Q_\tau^{k,v}}(\mathcal{M}_\rho)_\theta f(z).
\end{align}
If $k\in\mathbb{Z_+}\cap[N+1,\infty)$,
then, from \eqref{1556}, it follows that,
for any $y^{k,v}_\tau\in Q_\tau^{k,v}$,
\begin{align}\label{161}
\left|\mathcal{D}_kf(y^{k,v}_\tau)\right|
\leq2\inf_{z\in Q_\tau^{k,v}}(\mathcal{M}_\rho)_\theta f(z).
\end{align}
By $\|\phi\|_{\mathcal{G}(x,1,\beta_2,\gamma_2)}\lesssim1$,
the cancellation of $(\widetilde{D}^{(1)}_k)^*$,
Lemma~\ref{1122} with both $\overline{D}_k$ and $\widetilde{D}_l$
replaced, respectively, by
$\widetilde{D}_k^{(1)}(\cdot,y)\in
\mathring{\mathcal{G}}(y,2^{-k},\beta_1,\gamma_1)$ and
$\phi\in\mathcal{G}(x,1,\beta_2,\gamma_2)$,
and \eqref{949}, we obtain, for any given
$\beta'\in(0,\beta_-)$ and for any $y\in Q_\tau^{k,v}$,
\begin{align*}
\left|\left(\widetilde{\mathcal{D}}^{(1)}_k\right)^*\phi(y)\right|
&\lesssim2^{-k\beta'}\frac{1}{(V_\rho)_1(x)+V_\rho(x,y)}
\left[\frac{1}{1+\rho(x,y)}\right]^{\beta_-}
\nonumber\\
&\sim2^{-k\beta'}\frac{1}{(V_\rho)_1(x)+V_\rho(x,y_\tau^{k,v})}
\left[\frac{1}{1+\rho(x,y_\tau^{k,v})}\right]^{\beta_-},
\end{align*}
which, combined with \eqref{161},
further implies that
\begin{align*}
\left|\mathrm{Z}_2\right|
&\lesssim\sum_{k=N+1}^\infty2^{-k\beta'}
\sum_{\tau\in I_k}\sum_{v=1}^{N(k,\tau)}
\mu(Q_\tau^{k,v})\frac{1}{(V_\rho)_1(x)+V_\rho(x,y_\tau^{k,v})}
\\
&\quad\times\left[\frac{1}{1+\rho(x,y_\tau^{k,v})}\right]^{\beta_-}
\inf_{z\in Q_\tau^{k,v}}(\mathcal{M}_\rho)_\theta f(z).
\end{align*}
From this, \eqref{1612}, \eqref{15333},
and an argument similar to that used in
the estimation of \eqref{1999},
we deduce that, for any given $r\in(\frac{n}{n+\varepsilon},p)$,
\begin{align*}
\left|\langle f,\phi\rangle\right|\lesssim
\left\{\mathcal{M}_\rho\left(\left[(\mathcal{M}_\rho)_\theta
f\right]^r\right)(x)\right\}^\frac{1}{r},
\end{align*}
where $x\in B_\rho(x_1,1)$ is arbitrary,
which, together with Lemma~\ref{MLpq}, further implies
\begin{align*}
\left|\langle f,\phi\rangle\right|
&\lesssim\frac{\|\{\mathcal{M}_\rho([(\mathcal{M}_\rho)_\theta f]^r)
\}^\frac{1}{r}\|_{L^{p,q}(X)}}{\|
\mathbf{1}_{B_\rho(x_1,1)}\|_{L^{p,q}(X)}}\\
&\sim\left\|\mathcal{M}_\rho\left(\left[\left(\mathcal{M}_\rho
\right)_\theta f\right]^r\right)\right\|_{L^{\frac{p}{r},\frac{q}{r}}(X)}
\lesssim\left\|(\mathcal{M}_\rho)_\theta f\right\|_{L^{p,q}(X)}.
\end{align*}
This finishes the proof of \eqref{1515}
and hence the aforementioned claim.

By \eqref{1515}, Lemma~\ref{MMM}, and Definition~\ref{DefHardyL},
we conclude that, for any $\phi\in\mathcal{G}(\varepsilon,\varepsilon)$,
\begin{align}\label{163}
\left|\left\langle f,\phi\right\rangle\right|
\lesssim\left\|(\mathcal{M}_\rho)_\theta f\right\|_{L^{p,q}(X)}
\|\phi\|_{\mathcal{G}(\beta_2,\gamma_2)}
\lesssim\|f\|_{H^{p,q}_*(X)}
\|\phi\|_{\mathcal{G}(\beta_2,\gamma_2)}.
\end{align}
Now, for any given $h\in\mathcal{G}^\varepsilon_0(\beta_2,\gamma_2)$,
there exists a sequence $\{\phi_m\}_{m\in\mathbb{N}}$ in
$\mathcal{G}(\varepsilon,\varepsilon)$ such that
$$
\lim_{m\to\infty}\left\|h-\phi_m\right\|_{\mathcal{G}(\beta_2,\gamma_2)}
\to0
$$
as $m\to0$, which, combined with \eqref{163},
implies that
$$
\left|\left\langle f,\phi_m-\phi_k\right\rangle\right|
\lesssim\|f\|_{H^{p,q}_*(X)}
\left\|\phi_m-\phi_k\right\|_{\mathcal{G}(\beta_2,\gamma_2)}\to0
$$
as $m,k\to\infty$. Thus, the limit of
$\langle f,\phi_m\rangle$
exists and is independent of the choice of $\{\phi_m\}_{m\in\mathbb{N}}$
and we let
$$
\langle f,h\rangle:=\lim_{m\to\infty}
\left\langle f,\phi_m\right\rangle.
$$
From this and \eqref{163} again, we infer that
\begin{align*}
\left|\left\langle f,h\right\rangle\right|
&=\lim_{m\to\infty}\left|\left\langle f,\phi_m\right\rangle\right|
\lesssim\|f\|_{H^{p,q}_*(X)}
\lim_{m\to\infty}\|\phi_m\|_{\mathcal{G}(\beta_2,\gamma_2)}\\
&=\|f\|_{H^{p,q}_*(X)}\|h\|_{\mathcal{G}(\beta_2,\gamma_2)}.
\end{align*}
Thus, $f\in(\mathcal{G}^\varepsilon_0(\beta_2,\gamma_2))'$
and $\|f\|_{(\mathcal{G}^\varepsilon_0(\beta_2,\gamma_2))'}
\lesssim\|f\|_{H^{p,q}_*(X)}$.

On the other hand, by Theorem~\ref{HL-M} and the assumption
that $f\in H^{p,q}_*(X)$,
we find $\mathcal{M}^+f\in L^{p,q}(X)$,
which, together with the proven conclusion that
$f\in(\mathcal{G}^\varepsilon_0(\beta_2,\gamma_2))'$
and Theorem~\ref{HL-M} again,
further implies that $f\in H^{p,q}_*(X)
\subset(\mathcal{G}^\varepsilon_0(\beta_2,\gamma_2))'$.
This finishes the proof of Theorem~\ref{indep2}.
\end{proof}

Applying an argument similar to that used
in the proof of Theorem~\ref{indep2},
we obtain the following equivalence in the case
$\mathrm{diam\,}(X)<\infty$,
which is equivalent to $\mu(X)<\infty$.

\begin{theorem}\label{HL-Mfinite}
Let $(X,\mathbf{q},\mu)$ be an ultra-RD-space with
upper dimension $n\in(0,\infty)$,
where $\mu$ is assumed to be a Borel-semiregular
measure. Assume that $\mathrm{diam\,}(X)<\infty$.
Let $p\in(\frac{n}{n+\mathrm{ind\,}(X,\mathbf{q})},\infty)$,
$q\in(0,\infty]$, $\theta\in(0,\infty)$,
and $\beta,\gamma,\varepsilon\in(0,\infty)$ satisfy
$n(\frac{1}{p}-1)_+<\beta,\gamma<\varepsilon\preceq\mathrm{ind\,}(X,\mathbf{q})$.
Then, as subspaces of $({\mathcal{G}}^\varepsilon_0(\beta,\gamma))'$,
$$
H^{p,q}_\theta(X)=H^{p,q}_*(X)
$$
with equivalent quasi-norms.
\end{theorem}

\begin{proof}
Fix $\theta\in(0,\infty)$. By \eqref{2114} and Proposition~\ref{inc-norm}(ii),
we find that
\begin{align}\label{2213}
H^{p,q}_*(X)\hookrightarrow H^{p,q}_\theta(X).
\end{align}

On the other hand, choose $j\in\mathbb{N}$ sufficiently large
so that both \eqref{2cubej}
and
\begin{align*}
C_\rho C_r2^{-j}\leq\theta
\end{align*}
hold.
Let $f\in H^{p,q}_\theta(X)$.
Similarly to the proof of Lemma~\ref{fstarMf},
we first claim that, for any given $r\in(\frac{n}{n+\widetilde{\beta}},p)\cap(0,1]$
with $\widetilde{\beta}\in(n(\frac{1}{p}-1)_+,\beta\wedge\gamma)$ and for
any $l\in\mathbb{Z}$ and
$\phi\in\mathring{{\mathcal{G}}}^\varepsilon_0(\beta,\gamma)$
with $\|\phi\|_{{\mathcal{G}}(x,2^{-l},\beta,\gamma)}\leq1$,
\begin{align}\label{202}
|\langle f,\phi\rangle|\lesssim
\left[\mathcal{M}_\rho\left(\left[(\mathcal{M}_\rho)_\theta
f\right]^r\right)(x)\right]^\frac{1}{r}.
\end{align}
To this end, from Theorem~\ref{ID}, we deduce that
\begin{align}\label{16134}
\langle f,\phi\rangle
&=\sum_{k=0}^{N}\sum_{\tau\in I_k}\sum_{v=1}^{N(k,\tau)}
\int_{Q^{k,v}_\tau}
\left(\widetilde{\mathcal{D}}_k\right)^*\phi(y)\,d\mu(y)
\mathcal{D}^{k,v}_{\tau,1}f
\nonumber\\
&\quad+\sum_{k=N+1}^{\infty}\sum_{\tau\in I_k}
\sum_{v=1}^{N(k,\tau)}\int_{Q^{k,v}_\tau}
\left(\widetilde{\mathcal{D}}_k\right)^*\phi(y)\,d\mu(y)
\mathcal{D}_kf(y^{k,v}_\tau)
\nonumber\\
&=:\mathrm{Z}_1+\mathrm{Z}_2.
\end{align}
Note that, by \eqref{2031} and \eqref{161} [here, the
above choice of $j$ guarantees
that \eqref{1556} and hence \eqref{161} hold], we conclude that,
for any $k\in\mathbb{N}\cap[0,N]$, $\tau\in I_k$,
and $v\in\mathbb{N}\cap[1,N(k,\tau)]$,
\begin{align}\label{2037}
\left|\mathcal{D}^{k,v}_{\tau,1}f\right|
\leq 2\inf_{z\in Q_\tau^{k,v}}(\mathcal{M}_\rho)_\theta f(z)
\end{align}
and, for any $k\in\mathbb{N}\cap[N+1,\infty]$, $\tau\in I_k$,
and $v\in\mathbb{N}\cap[1,N(k,\tau)]$,
\begin{align}\label{2032}
\left|\mathcal{D}_kf(y^{k,v}_\tau)\right|
\leq2\inf_{z\in Q_\tau^{k,v}}(\mathcal{M}_\rho)_\theta f(z).
\end{align}
Next, from Lemma~\ref{1122} with
both $\overline{D}_k$ and $\widetilde{D}_l$
replaced, respectively, by
$(\widetilde{D}_{k})^*(y_\tau^{l,v},\cdot)\in
\mathring{\mathcal{G}}(y_\tau^{l,v},2^{-k},\beta,\gamma)$ and
$\phi\in\mathring{\mathcal{G}}(x,2^{-l},\beta,\gamma)$
and from \eqref{949} with $1$ therein replaced by $2^{-(k\wedge l)}$,
we infer that,
for any $k\in\mathbb{N}\cap[N+1,\infty]$, $\tau\in I_k$,
$v\in\mathbb{N}\cap[1,N(k,\tau)]$, and $y\in Q_\tau^{k,v}$,
\begin{align}\label{1113}
\left|\left(\widetilde{\mathcal{D}}_k\right)^*\phi(y)\right|
&\lesssim2^{-|k-l|\widetilde{\beta}}\frac{1}{(V_\rho)_{2^{-(k\wedge l)}}(x)+
V_\rho(x,y)}\left[\frac{2^{-(k\wedge l)}}
{2^{-(k\wedge l)}+\rho(x,y)}\right]^\gamma
\nonumber\\
&\lesssim2^{-|k-l|\widetilde{\beta}}\frac{1}{(V_\rho)_{2^{-(k\wedge l)}}(x)+
V_\rho(x,y_\tau^{k,v})}\left[\frac{2^{-(k\wedge l)}}
{2^{-(k\wedge l)}+\rho(x,y_\tau^{k,v})}\right]^\gamma.
\end{align}

Now, let $k\in\mathbb{N}\cap[0,N]$, $\tau\in I_k$,
$v\in\mathbb{N}\cap[1,N(k,\tau)]$, and $y\in Q_\tau^{k,v}$.
We turn to estimate $(\widetilde{\mathcal{D}}_k)^*\phi(y)$.
Note that $(\widetilde{D}_k)^*(y,\cdot)\in
\mathcal{G}(y,2^{-k},\beta,\gamma)$, $2^{-k}\sim1$,
with the positive equivalence constants depending only on $N$,
and $\phi\in\mathring{\mathcal{G}}(x,2^{-l},\beta,\gamma)$.
We consider the following two cases for $k\wedge l$.

\emph{Case (1)} $k\wedge l=k$. In this case,
by the size conditions of both $(\widetilde{D}_k)^*(y,\cdot)$
and $\phi$ and the fact that $2^{-(k\wedge l)}=2^{-k}\sim1$,
we find that
\begin{align*}
\left|(\widetilde{\mathcal{D}}_k)^*\phi(y)\right|
&\leq\int_X\left|\widetilde{D}_k(z,y)\right|
\left|\phi(z)\right|\,d\mu(z)
\\
&\lesssim\int_X\frac{1}{(V_\rho)_{1+\rho(z,y)}(y)}
\left[\frac{1}{1+\rho(z,y)}\right]^\gamma
\\
&\quad\times\frac{1}{(V_\rho)_{2^{-l}+\rho(z,x)}(x)}
\left[\frac{2^{-l}}{2^{-l}+\rho(z,x)}\right]^\gamma\,d\mu(z)
\\
&=\int_{\rho(x,z)\leq\rho(z,y)+1}\cdots+\int_{\rho(x,z)>\rho(z,y)+1}\cdots.
\end{align*}
From this and the facts that
\begin{align*}
1+\rho(x,y)&\leq1+C_\rho\max\left\{\rho(x,z),\,\rho(z,y)
\right\}\\
&\leq1+C_\rho\left[\rho(z,y)+1\right]
\sim1+\rho(z,y)
\end{align*}
for any $z\in X$ satisfying $\rho(x,z)\leq\rho(z,y)+1$
and that
\begin{align*}
1+\rho(x,y)&\leq1+C_\rho\max\left\{\rho(x,z),\,\rho(z,y)\right\}\\
&=1+C_\rho\rho(x,z)
\sim\rho(x,z)\leq2^{-l}+\rho(x,z)
\end{align*}
for any $z\in X$ satisfying $\rho(x,z)>\rho(z,y)+1$,
we deduce that
\begin{align*}
\left|(\widetilde{\mathcal{D}}_k)^*\phi(y)\right|
&\lesssim\frac{1}{(V_\rho)_{1+\rho(x,y)}(y)}\left[\frac{1}{1+\rho(x,y)}\right]^\gamma\\
&\quad\times\int_{\rho(x,z)\leq\rho(z,y)+1}\frac{1}{(V_\rho)_{2^{-l}+\rho(z,x)}(x)}
\left[\frac{2^{-l}}{2^{-l}+\rho(z,x)}\right]^\gamma\,d\mu(z)\\
&\quad+\frac{1}{(V_\rho)_{1+\rho(x,y)}(y)}\left[\frac{1}{1+\rho(x,y)}\right]^\gamma\\
&\quad\times\int_{\rho(x,z)>\rho(z,y)+1}\frac{1}{(V_\rho)_{1+\rho(z,y)}(y)}
\left[\frac{1}{1+\rho(z,y)}\right]^\gamma\,d\mu(z),
\end{align*}
which, combined with Lemma~\ref{A3.2}(ii),
further implies that
\begin{align*}
\left|(\widetilde{\mathcal{D}}_k)^*\phi(y)\right|
&\lesssim\frac{1}{(V_\rho)_{1+\rho(x,y)}(y)}\left[\frac{1}{1+\rho(x,y)}\right]^\gamma\\
&\sim\frac{1}{(V_\rho)_{2^{-(k\wedge l)}+\rho(x,y)}(y)}
\left[\frac{2^{-(k\wedge l)}}{2^{-(k\wedge l)}+\rho(x,y)}\right]^\gamma
\end{align*}
because $1\sim2^{-k}=2^{-(k\wedge l)}$ in this case.

\emph{Case (2)} $k\wedge l=l$. In this case, using
Lemma~\ref{1122} with both $\overline{D}_k$ and $\widetilde{D}_l$
replaced, respectively, by
$(\widetilde{D}_k)^*(y,\cdot)\in\mathcal{G}(y,2^{-k},\beta,\gamma)$ and
$\phi\in\mathring{\mathcal{G}}(x,2^{-l},\beta,\gamma)$
and using \eqref{949} with $1$ replaced by $2^{-(k\wedge l)}$, we conclude that,
for any $y\in Q_\tau^{k,v}$, \eqref{1113} still holds in this case.

Combining the above two cases, we obtain,
for any $k\in\mathbb{N}\cap[0,N]$, $\tau\in I_k$,
$v\in\mathbb{N}\cap[1,N(k,\tau)]$, and $y\in Q_\tau^{k,v}$,
\begin{align*}
\left|\left(\widetilde{\mathcal{D}}_k\right)^*\phi(y)\right|
\lesssim2^{-|k-l|\widetilde{\beta}}\frac{1}{(V_\rho)_{2^{-(k\wedge l)}}(x)+
V_\rho(x,y_\tau^{k,v})}\left[\frac{2^{-(k\wedge l)}}
{2^{-(k\wedge l)}+\rho(x,y_\tau^{k,v})}\right]^\gamma.
\end{align*}
This, together with \eqref{16134}, \eqref{2037}, \eqref{1113}, \eqref{2032},
and Lemma~\ref{L4.5}, further implies that
\begin{align*}
\left|\langle f,\phi\rangle\right|
&\lesssim\sum_{k=0}^\infty2^{-|k-l|\widetilde{\beta}}
\sum_{\tau\in I_k}\sum_{v=1}^{N(k,\tau)}
\mu(Q_\tau^{k,v})\frac{\inf_{z\in Q_\tau^{k,v}}(\mathcal{M}_\rho)_\theta f(z)}
{(V_\rho)_{2^{-(k\wedge l)}}(x)+
V_\rho(x,y_\tau^{k,v})}\left[\frac{2^{-(k\wedge l)}}
{2^{-(k\wedge l)}+\rho(x,y_\tau^{k,v})}\right]^\gamma\\
&\lesssim\sum_{k=0}^\infty2^{-|k-l|\widetilde{\beta}+n(k\wedge l-k)(1-\frac{1}{r})}
\left[\mathcal{M}_\rho\left(\sum_{\tau\in I_k}\sum_{v=1}^{N(k,\tau)}
\left[\inf_{z\in Q_\tau^{k,v}}(\mathcal{M}_\rho)_\theta f(z)
\right]^r\mathbf{1}_{Q_\tau^{k,v}}\right)(x)\right]^\frac{1}{r}
\\
&\lesssim\left[\mathcal{M}_\rho\left(\left[(\mathcal{M}_\rho)_\theta
f\right]^r\right)(x)\right]^\frac{1}{r},
\end{align*}
which completes the proof of \eqref{202} and hence the
aforementioned claim.
From this and an argument similar to that used in the proof of Lemma~\ref{fstarMf}
with \eqref{fphi} replaced by \eqref{202},
we infer that, for any $\psi\in{\mathcal{G}}^\varepsilon_0(\beta,\gamma)$
with $\|\psi\|_{{\mathcal{G}}(x,r_0,\beta,\gamma)}\leq 1$
for some $r_0\in(0,\infty)$,
\begin{align*}
\left|\langle f,\psi\rangle\right|
\lesssim\left[\mathcal{M}_\rho\left(\left[(\mathcal{M}_\rho)_\theta
f\right]^r\right)(x)\right]^\frac{1}{r}+(\mathcal{M}_\rho)_\theta
f(x),
\end{align*}
which, combined with an argument similar to that
used in the estimation of \eqref{221111} with $\mathcal{M}^+f$
replaced by $(\mathcal{M}_\rho)_\theta$,
further implies that
\begin{align*}
\left\|f^*_\rho\right\|_{L^{p,q}(X)}
\lesssim\left\|(\mathcal{M}_\rho)_\theta\right\|_{L^{p,q}(X)}<\infty.
\end{align*}
Thus, $H^{p,q}_\theta(X)\hookrightarrow H^{p,q}_*(X)$.
This, together with \eqref{2213},
then finishes the proof of Theorem~\ref{HL-Mfinite}.
\end{proof}

By the proof of Theorem~\ref{HL-M} and applying
an argument similar to that used in the proof of
Theorem~\ref{indep2} with
Theorem~\ref{HL-M} therein
replaced by Theorem~\ref{HL-Mfinite}, we obtain the following
conclusion  in the case $\mathrm{diam\,}(X)<\infty$; we omit the details.

\begin{theorem}\label{indep3}
Let $(X,\mathbf{q},\mu)$ be an ultra-RD-space with upper dimension $n\in(0,\infty)$,
where $\mu$ is assumed to be a Borel-semiregular measure.
Assume that $\mathrm{diam\,}(X)<\infty$. Let
$p\in(\frac{n}{n+\mathrm{ind\,}(X,\mathbf{q})},\infty)$,
$q\in(0,\infty]$, and $n(\frac{1}{p}-1)_+<
\varepsilon\preceq\mathrm{ind\,}(X,\mathbf{q})$.
Then $H_*^{p,q}(X)$ is independent of the
choice of parameters $\beta$ and $\gamma$
that intervene in the definition of the spaces of distributions
$({\mathcal{G}}^\varepsilon_0(\beta,\gamma))'$ whenever
$\beta,\gamma\in(0,\infty)$ satisfy
$n(\frac{1}{p}-1)_+<\beta,\gamma<\varepsilon$.
\end{theorem}

The following two propositions show that Hardy--Lorentz spaces
can be continuously embedded into the space
of distributions and are complete.

\begin{proposition}\label{5.1.10}
Let $(X,\mathbf{q},\mu)$ be a
quasi-ultrametric spaces of homogeneous type,
$p\in(0,\infty)$,
$q\in(0,\infty]$,
and $0<\beta,\gamma<\varepsilon\preceq\mathrm{ind\,}(X,\mathbf{q})$.
Then $H^{p,q}_*(X)\subset({\mathcal{G}}^\varepsilon_0(\beta,\gamma))'$
is continuously embedded into $({\mathcal{G}}^\varepsilon_0(\beta,\gamma))'$.
\end{proposition}

\begin{proof}
We only prove in the case $q\in(0,\infty)$ because the
proof in the case $q=\infty$ is similar and we omit the details.
Let $f\in H^{p,q}_*(X)$ and
$\phi\in\mathcal{G}^\varepsilon_0(\beta,\gamma)$
with $\|\phi\|_{\mathcal{G}(x_1,1,\beta,\gamma)}=1$.
Recall that we may assume, without loss of generality, that
$H^{p,q}_*(X)$ is defined via
a regularized quasi-ultrametric $\rho\in\mathbf{q}$
as in Definition~\ref{D2.3}.
By Remarks~\ref{Rmtest}(i) and~\ref{RemMf}(iii), we find that there exists
$C\in(0,\infty)$ such that, for any $x\in B_\rho(x_1,1)$,
$$
\|\phi\|_{G(x,1,\frac{\beta}{\alpha},\frac{\gamma}{\alpha})}
\leq C\|\phi\|_{\mathcal{G}(x_1,1,\beta,\gamma)}=C,
$$
where $\alpha:=\min\{1,\,(\log_2C_\rho)^{-1}\}$,
which further implies that
$$
|\langle f,\phi\rangle|\leq Cf^*_\rho(x).
$$
From this and Proposition~\ref{Lpqnorm=c,d},
we infer that
\begin{align*}
\left|\langle f,\phi\rangle\right|^q
&\sim\sum_{k\in\mathbb{Z}\cap(-\infty,\log_2
|\langle f,\phi\rangle|)}2^{kq}
\sim\sum_{k\in\mathbb{Z}\cap(-\infty,\log_2
|\langle f,\phi\rangle|)}2^{kq}\left[(V_\rho)_1(x_1)\right]^\frac{q}{p}\\
&\sim\sum_{k\in\mathbb{Z}\cap(-\infty,\log_2
|\langle f,\phi\rangle|)}2^{kq}\left[\mu\left(\left\{x\in B_\rho(x_1,1):
\left|\langle f,\phi\rangle\right|>2^k\right\}\right)
\right]^\frac{q}{p}\\
&\leq\sum_{k\in\mathbb{Z}\cap(-\infty,\log_2
|\langle f,\phi\rangle|)}2^{kq}\left[\mu\left(\left\{x\in B_\rho(x_1,1):
Cf^*_\rho(x)>2^k\right\}\right)
\right]^\frac{q}{p}\\
&\leq\sum_{k\in\mathbb{Z}}2^{kq}\left[\mu\left(
\left\{x\in B_\rho(x_1,1):Cf^*_\rho(x)>2^k\right\}\right)
\right]^\frac{q}{p}\lesssim\|f\|_{H^{p,q}_*(X)}^q.
\end{align*}
Taking the supremum of all
$\phi\in\mathcal{G}^\varepsilon_0(\beta,\gamma)$
with $\|\phi\|_{\mathcal{G}(x_1,1,\beta,\gamma)}=1$,
we conclude that
$$
f\in({\mathcal{G}}^\varepsilon_0(\beta,\gamma))'
\ \ \text{and}\ \
\|f\|_{({\mathcal{G}}^\varepsilon_0(\beta,\gamma))'}\lesssim
\|f\|_{H^{p,q}_*(X)}.
$$
This finishes the proof of Proposition~\ref{5.1.10}.
\end{proof}

\begin{proposition}
\label{5.1.11}
Let $(X,\mathbf{q},\mu)$ be a
quasi-ultrametric spaces of homogeneous type, $p\in(0,\infty)$,
$q\in(0,\infty]$.
Then $H^{p,q}_*(X)$ is complete.
\end{proposition}

\begin{proof}
Recall that we may assume, without loss of generality, that
$H^{p,q}_*(X)$ is defined by means of
a regularized quasi-ultrametric, $\rho\in\mathbf{q}$,
as in Definition~\ref{D2.3}.
Let $\{f_k\}_{k\in\mathbb{N}}$ be a Cauchy sequence
in $H^{p,q}_*(X)$. By this and Proposition~\ref{5.1.10},
we find that $\{f_k\}_{k\in\mathbb{N}}$
is also a Cauchy sequence
in $({\mathcal{G}}^\varepsilon_0(\beta,\gamma))'$,
which, combined with the completeness of
$({\mathcal{G}}^\varepsilon_0(\beta,\gamma))'$,
further implies that there exists $f\in({\mathcal{G}}^\varepsilon_0(\beta,\gamma))'$
such that $f_k\to f$ as $k\to\infty$
in $({\mathcal{G}}^\varepsilon_0(\beta,\gamma))'$.
Let $x\in X$ and $\phi\in{\mathcal{G}}^\varepsilon_0(\beta,\gamma)$
be such that $\|\phi\|_{G(x,r_0^\alpha,\frac{\beta}{\alpha},
\frac{\gamma}{\alpha})}\leq 1$ for some $r_0\in(0,\infty)$, where
$\alpha:=\min\{1,\,(\log_2C_\rho)^{-1}\}$.
Then, for any $k,l\in\mathbb{N}$,
$$
\left|\left\langle f_{k+l}-f_k,\phi\right\rangle\right|
\leq(f_{k+l}-f_k)_\rho^*(x).
$$
Letting $l\to\infty$, we then obtain
$$
\left|\langle f-f_k,\phi\rangle\right|
\leq\liminf_{l\to\infty}(f_{k+l}-f_k)_\rho^*(x).
$$
Taking the supremum of all $\phi\in{\mathcal{G}}^\varepsilon_0(\beta,\gamma)$
satisfying that there exists $r_0\in(0,\infty)$ such that
$\|\phi\|_{G(x,r_0^\alpha,\frac{\beta}{\alpha},\frac{\gamma}{\alpha})}\leq 1$,
we conclude that
$$
(f-f_k)_\rho^*(x)\leq\liminf_{l\to\infty}(f_{k+l}-f_k)_\rho^*(x),
$$
which, together with Proposition~\ref{inc-norm}(ii),
Fatou's lemma, and the fact that $\{f_k\}_{k\in\mathbb{N}}$ is a Cauchy sequence
in $H^{p,q}_*(X)$, further implies that
\begin{align*}
\left\|(f-f_k)_\rho^*\right\|_{L^{p,q}(X)}
&\leq\left\|\liminf_{l\to\infty}(f_{k+l}-f_k)_\rho^*\right\|_{L^{p,q}(X)}\\
&\leq\liminf_{l\to\infty}\left\|(f_{k+l}-f_k)_\rho^*\right\|_{L^{p,q}(X)}
\to0
\end{align*}
as $k\to\infty$.
From this, we deduce that $f\in H^{p,q}_*(X)$ and
$$
\lim_{k\to\infty}\|f-f_k\|_{H^{p,q}_*(X)}=0.
$$
Thus, $H^{p,q}_*(X)$ is complete,
which completes the proof of
Proposition~\ref{5.1.11}.
\end{proof}

\section[Atomic and Finite Atomic Characterizations]
{Atomic and Finite Atomic Characterizations}
\label{section5.3}

In this section, we establish
the atomic and the finite atomic characterizations of
the grand maximal-function-based Hardy--Lorentz space
$H^{p,q}_*(X)$,
respectively, in Subsection~\ref{5.2.1}
and in Subsection~\ref{5.2.2}
for the optimal ranges
\begin{align*}
p\in\left(\frac{n}{n+{\mathrm{ind\,}}(X,\mathbf{q})},1\right]
\ \ \text{and}\ \
q\in(0,\infty]
\end{align*}
on quasi-ultrametric spaces of homogeneous type
with upper dimension $n\in(0,\infty)$.
As corollaries
of these and Theorems~\ref{HL-M} and~\ref{HL-Mfinite},
we also obtain such atomic
characterization for the Hardy--Lorentz spaces $H^{p,q}_+(X)$
and $H^{p,q}_\theta(X)$.

\subsection{Atomic Characterization}\label{5.2.1}

In this subsection, we establish
the atomic characterization of $H^{p,q}_*(X)$.
To this end,
based on the $(p,r)$-atom in Definition~\ref{atom},
we first introduce the definition of
atomic Hardy--Lorentz spaces.

\begin{definition}\label{6.3}
Let $(X,\mathbf{q},\mu)$ be a quasi-ultrametric space
of homogeneous type with upper dimension $n\in(0,\infty)$,
$\rho\in\mathbf{q}$ have the property that
all $\rho$-balls are $\mu$-measurable,
$p\in(0,1]$,
$q\in(0,\infty]$, $r\in(p,\infty]\cap[1,\infty]$,
and $n(\frac{1}{p}-1)<\beta,\gamma<\varepsilon<\infty$.
The $\emph{atomic Hardy--Lorentz space}$
\index[words]{atomic Hardy--Lorentz space}
$H^{p,q,r}_{\mathrm{at}}(X,\rho)$
\index[symbols]{H@$H^{p,q,r}_{\mathrm{at}}(X,\rho)$}
is defined to be the set of
all $f\in({\mathcal{G}}^\varepsilon_0(\beta,\gamma))'$
such that there exist a sequence
$\{\lambda_{i,k}\}_{i\in\mathbb{N},k\in\mathbb{Z}}$
in $\mathbb{R}_+$ and a sequence
of $(\rho,p,r)$-atoms
$\{a_{i,k}\}_{i\in\mathbb{N},k\in\mathbb{Z}}$
supported, respectively,
in $\rho$-balls $\{B_{i,k}\}_{i\in\mathbb{N},k\in\mathbb{Z}}$
such that
\begin{align}\label{213}
f=\sum_{k\in\mathbb{Z}}
\sum_{i\in\mathbb{N}}
\lambda_{i,k}a_{i,k}
\end{align}
in $({\mathcal{G}}^\varepsilon_0(\beta,\gamma))'$ and
the following conditions hold:
\begin{enumerate}
\item There exist constants
$C_0\in\mathbb{N}$ and $j_0\in\mathbb{Z}\setminus\mathbb{N}$,
both of which are independent of $f$, such that, for any $k\in\mathbb{Z}$ and $x\in X$,
$$
\sum_{i\in\mathbb{N}}\mathbf{1}_{(2C_\rho)^{j_0}B_{i,k}}(x)\leq C_0.
$$

\item For any $i\in\mathbb{N}$ and $k\in\mathbb{Z}$,
$\lambda_{i,k}\sim2^k[\mu(B_{i,k})]^\frac{1}{p}$
with the positive equivalence constants independent of $f$
and
\begin{align*}
\left[\sum_{k\in\mathbb{Z}}\left(\sum_{i\in\mathbb{N}}
\lambda_{i,k}^p\right)^\frac{q}{p}\right]^\frac{1}{q}<\infty
\end{align*}
with the usual modification
when $q=\infty$.
\end{enumerate}
Moreover, for any $f\in H^{p,q,r}_{\mathrm{at}}(X,\rho)$, define
$$
\|f\|_{H^{p,q,r}_{\mathrm{at}}(X,\rho)}:
=\inf\left\{\left[\sum_{k\in\mathbb{Z}}\left(\sum_{i\in\mathbb{N}}
\lambda_{i,k}^p\right)^\frac{q}{p}\right]^\frac{1}{q}:
f=\sum_{k\in\mathbb{Z}}
\sum_{i\in\mathbb{N}}
\lambda_{i,k}a_{i,k}\right\}
$$
with the usual modification
when $q=\infty$, where the infimum is taken
over all decompositions of $f$ as in \eqref{213}.
\end{definition}

\begin{remark}\label{5216}
Let the notation be as in Definition~\ref{6.3}.
\begin{enumerate}
\item[\rm(i)]
By $n(\frac{1}{p}-1)<\beta,\gamma<\infty$,
\eqref{1550}, and Remark~\ref{Rmtest}(v), we find that,
for any $p\in(0,\frac{n}{n+\mathrm{ind}_H(X,\mathbf{q})}]$,
$H^{p,q,r}_{\mathrm{at}}(X,\rho)$ is trivial because the space
${\mathcal{G}}^\varepsilon_0(\beta,\gamma)$ only contains constant functions.

\item[\rm(ii)]
If $p\in(\frac{n}{n+\mathrm{ind\,}(X,\mathbf{q})},1]$,
$q\in(0,\infty]$, $r\in(p,\infty]\cap[1,\infty]$,
and $n(\frac{1}{p}-1)<\beta,\gamma<\varepsilon\preceq\mathrm{ind\,}(X,\mathbf{q})$,
then $H^{p,q,r}_{\mathrm{at}}(X,\rho)$
is a quasi-Banach space in light of
Proposition~\ref{5.1.11} and Theorem~\ref{6.14}.

\item[\rm(iii)]
Suppose $f\in H^{p,q,r}_{\mathrm{at}}(X,\rho)$ and let
$\{\lambda_{i,k}\}_{i\in\mathbb{N},k\in\mathbb{Z}}$ in $\mathbb{R}_+$ and
$\{a_{i,k}\}_{i\in\mathbb{N},k\in\mathbb{Z}}$ be the same as in Definition~\ref{6.3}.
If $\rho'\in\mathbf{q}$ has the property that all $\rho'$-balls are $\mu$-measurable, then,
according to Remark~\ref{Rmk-atoms}(i), we can write
$$
f=\sum_{k\in\mathbb{Z}}
\sum_{i\in\mathbb{N}}
(C\lambda_{i,k}) (C^{-1}a_{i,k})
$$
in $({\mathcal{G}}^\varepsilon_0(\beta,\gamma))'$,
where $C$ is a positive constant depending only on $\rho$, $\rho'$, and $\mu$
and where $\{C^{-1}a_{i,k}\}_{i\in\mathbb{N},k\in\mathbb{Z}}$ is a sequence
of $(\rho',p,r)$-atoms supported, respectively,
in $\rho'$-balls $\{\overline{B}_{i,k}\}_{i\in\mathbb{N},k\in\mathbb{Z}}$.
Moreover, by \eqref{1557-xx} and \eqref{doublingcond}, we conclude that
the sequences $\{\overline{B}_{i,k}\}_{i\in\mathbb{N},k\in\mathbb{Z}}$
and $\{C\lambda_{i,k}\}_{i\in\mathbb{N},k\in\mathbb{Z}}$ in $\mathbb{R}_+$ satisfy
conditions (i) and (ii) in Definition~\ref{6.3} for some constants
depending on $\rho$, $\rho'$, $\mu$, $C_0$, $C_1$, and $j_0$.
It follows from these observations that the
space $H^{p,q,r}_{\mathrm{at}}(X,\rho)$ is independent of the choice of $\rho\in\mathbf{q}$.
Thus, in what follows, we sometimes simply write
$H^{p,q,r}_{\mathrm{at}}(X)$ instead of $H^{p,q,r}_{\mathrm{at}}(X,\rho)$.

\item[\rm(iv)] By Definition~\ref{6.3}, the proof of Theorem~\ref{6.16},
and Lemma~\ref{325},
we find that, when $p=q\in(\frac{n}{n+\mathrm{ind\,}(X,\mathbf{q})},1]$
and $r\in(p,\infty]\cap[1,\infty]$,
as subspaces of $(\mathcal{G}_0^\varepsilon(\beta,\gamma))'$,
\begin{align*}
H^{p,p,r}_{\rm at}(X)=H^{p,r}_{\rm at}(X)=H^{p,r}_{\rm CW}(X)
\end{align*}
with equivalent quasi-norms,
where $\varepsilon,\beta,\gamma\in(0,\infty)$ satisfy
$n(\frac{1}{p}-1)<\beta,\gamma<\varepsilon\preceq\mathrm{ind\,}(X,\mathbf{q})$ and
$H^{p,r}_{\rm at}(X)$
and $H^{p,r}_{\rm CW}(X)$ are the atomic Hardy spaces as in
Definitions~\ref{atomichardy}(i) and~\ref{1539}, respectively.
\end{enumerate}
\end{remark}

Next, we present the atomic characterization
of $H^{p,q}_*(X)$ as follows.

\begin{theorem}\label{6.14}
Let $(X,\mathbf{q},\mu)$ be a quasi-ultrametric space
of homogeneous type with upper dimension $n\in(0,\infty)$,
$\rho\in\mathbf{q}$ have the property that
all $\rho$-balls are $\mu$-measurable,
$p\in(\frac{n}{n+\mathrm{ind\,}(X,\mathbf{q})},1]$,
$q\in(0,\infty]$, $r\in(p,\infty]\cap[1,\infty]$,
and $\varepsilon,\beta,\gamma\in(0,\infty)$ be such that
$n(\frac{1}{p}-1)<\beta,\gamma<\varepsilon\preceq\mathrm{ind\,}(X,\mathbf{q})$.
Then, as subspaces of $({\mathcal{G}}^\varepsilon_0(\beta,\gamma))'$,
$$
H^{p,q,r}_{\mathrm{at}}(X,\rho)=H^{p,q}_*(X)
$$
with equivalent (quasi-)norms.
\end{theorem}

\begin{remark}
Alvarado and Mitrea \cite[Theorem~5.27]{AlMi2015}
established the atomic characterization of the maximal Hardy space $H^p(X)$
on $n$-Ahlfors regular quasi-ultrametric spaces $(X,\mathbf{q},\mu)$
for the optimal ranges
\begin{align}\label{1659}
p\in\left(\frac{n}{n+\mathrm{ind\,}(X,\mathbf{q})},1\right]
\ \ \text{and}\ \
r\in(p,\infty]\cap[1,\infty].
\end{align}
From Remarks~\ref{2548} and~\ref{5216}(iv),
we deduce that Theorem~\ref{6.14}
is a generalization of \cite[Theorem~5.27]{AlMi2015}
from maximal Hardy spaces on $n$-Ahlfors regular quasi-ultrametric spaces
to maximal Hardy--Lorentz spaces on quasi-ultrametric spaces of
homogeneous type, both of which hold for the optimal ranges for $p$ and $r$
as in \eqref{1659}.
\end{remark}

To show Theorem~\ref{6.14},
we need to prove the following three theorems.

\begin{theorem}\label{6.15}
Let $(X,\mathbf{q},\mu)$ be a quasi-ultrametric space
of homogeneous type with upper dimension $n\in(0,\infty)$,
$\rho\in\mathbf{q}$ have the property that
all $\rho$-balls are $\mu$-measurable,
$p\in(\frac{n}{n+\mathrm{ind\,}(X,\mathbf{q})},1]$,
$q\in(0,\infty]$, $r\in(p,\infty]\cap[1,\infty]$,
and $\varepsilon,\beta,\gamma\in(0,\infty)$ be such that
$n(\frac{1}{p}-1)<\beta,\gamma<\varepsilon\preceq\mathrm{ind\,}(X,\mathbf{q})$.
Suppose that a sequence
$\{\lambda_{i,k}\}_{i\in\mathbb{N},k\in\mathbb{Z}}$ in $\mathbb{R}_+$ and a sequence
of $(\rho,p,r)$-atoms $\{a_{i,k}\}_{i\in\mathbb{N},k\in\mathbb{Z}}$ supported,
respectively, in $\rho$-balls $
\{B_{i,k}\}_{i\in\mathbb{N},k\in\mathbb{Z}}
:=\{B_\rho(x_{i,k},r_{i,k})\}_{i\in\mathbb{N},k\in\mathbb{Z}}
$
with $\{x_{i,k}\}_{i\in\mathbb{N},k\in\mathbb{Z}}$ in $X$
and $\{r_{i,k}\}_{i\in\mathbb{N},k\in\mathbb{Z}}$ in $(0,\infty)$
satisfy that
\begin{enumerate}
\item[\emph{(a)}]
There exist constants $j_0\in\mathbb{Z}
\setminus\mathbb{N}$ and $C_0\in\mathbb{N}$
such that, for any $x\in X$ and $k\in\mathbb{Z}$,
\begin{align}\label{4ee0}
\sum_{i\in\mathbb{N}}\mathbf{1}_{(2C_\rho)^{j_0}B_{i,k}}(x)\leq C_0.
\end{align}

\item[\emph{(b)}]
There exists a positive constant $C_1\in[1,\infty)$ such that, for any $k\in\mathbb{Z}$
and $i\in\mathbb{N}$,
\begin{align}\label{lambdaik}
\frac{\lambda_{i,k}}{C_1}\leq2^k\left[\mu(B_{i,k})\right]^\frac{1}{p}
\leq C_1\lambda_{i,k}
\end{align}
and
\begin{align}\label{4ee1}
\sum_{k\in\mathbb{Z}}\left(\sum_{i\in\mathbb{N}}\lambda_{i,k}^p\right)
^{\frac{q}{p}}<\infty
\end{align}
with the usual modification when $q=\infty$.
\end{enumerate}
Then the following assertions hold.
\begin{enumerate}
\item[\rm(i)]
The series
$\sum_{k\in\mathbb{Z}}\sum_{i\in\mathbb{N}}
\lambda_{i,k}a_{i,k}$ converges
in $({\mathcal{G}}^\varepsilon_0(\beta,\gamma))'$.

\item[\rm(ii)]
$f:=\sum_{k\in\mathbb{Z}}\sum_{i\in\mathbb{N}}
\lambda_{i,k}a_{i,k}\in H^{p,q}_*(X)$ and
there exists a positive constant $C$,
independent of $f$,
such that
$$
\|f\|_{H^{p,q}_*(X)}\leq C\left[\sum_{k\in\mathbb{Z}}
\left(\sum_{i\in\mathbb{N}}\lambda_{i,k}^p\right)
^{\frac{q}{p}}\right]^{\frac{1}{q}}
$$
with the usual modification when $q=\infty$.

\item[\rm(iii)]
The series
$\sum_{k\in\mathbb{Z}}\sum_{i\in\mathbb{N}}
\lambda_{i,k}a_{i,k}$ converges
in $H^{p,q}_*(X)$.
\end{enumerate}
\end{theorem}

\begin{theorem}\label{6.16}
Let $(X,\mathbf{q},\mu)$ be a quasi-ultrametric space
of homogeneous type with upper dimension $n\in(0,\infty)$,
$\rho\in\mathbf{q}$ have the property that
all $\rho$-balls are $\mu$-measurable,
$p\in(\frac{n}{n+\mathrm{ind\,}(X,\mathbf{q})},1]$,
$q\in(0,\infty]$,
and $\varepsilon,\beta,\gamma\in(0,\infty)$ be such that
$n(\frac{1}{p}-1)<\beta,\gamma<\varepsilon\preceq\mathrm{ind\,}(X,\mathbf{q})$.
Then, as subspaces of $({\mathcal{G}}^\varepsilon_0(\beta,\gamma))'$,
$$
H^{p,q}_*(X)\subset H^{p,q,\infty}_{\mathrm{at}}(X,\rho).
$$
Moreover, there exists a positive constant $C$ such that,
for any $f\in H^{p,q}_*(X)$,
$$
\|f\|_{H^{p,q,\infty}_{\mathrm{at}}(X,\rho)}\leq C\|f\|_{H^{p,q}_*(X)}.
$$
\end{theorem}

\begin{theorem}\label{6.17}
Let $(X,\mathbf{q},\mu)$ be a quasi-ultrametric space
of homogeneous type with upper dimension $n\in(0,\infty)$,
$\rho\in\mathbf{q}$ have the property that
all $\rho$-balls are $\mu$-measurable,
$p\in(0,1]$,
$q\in(0,\infty]$, $r\in(p,\infty]\cap[1,\infty]$,
and $\varepsilon,\beta,\gamma\in(0,\infty)$ be such that
$n(\frac{1}{p}-1)<\beta,\gamma<\varepsilon<\infty$.
Then, as subspaces of $({\mathcal{G}}^\varepsilon_0(\beta,\gamma))'$,
$$
H^{p,q,\infty}_{\mathrm{at}}(X,\rho)\subset H^{p,q,r}_{\mathrm{at}}(X,\rho).
$$
Moreover, for any $f\in H^{p,q,\infty}_{\mathrm{at}}(X,\rho)$,
$$
\|f\|_{H^{p,q,r}_{\mathrm{at}}(X,\rho)}\leq \|f\|_{H^{p,q,\infty}_{\mathrm{at}}(X,\rho)}.
$$
\end{theorem}

Theorem~\ref{6.17} is a direct consequence
of Definitions~\ref{6.3} and~\ref{atom} because any $(p,\infty)$-atom is also a
$(p,r)$-atom whenever $p\in(0,1]$ and
$r\in(p,\infty]\cap[1,\infty]$.
To show Theorem~\ref{6.15}, we need several technical lemmas.
The following two lemmas are a generalization of
\cite[Propositions~4.8 and~4.9]{zhy}, respectively,
from unbounded spaces of homogeneous type satisfying $\mu(\{x\})=0$
for all $x\in X$ to quasi-ultrametric
spaces of homogeneous type.

\begin{lemma}\label{4r2}
Let $(X,\mathbf{q},\mu)$ be a quasi-ultrametric space
of homogeneous type with upper dimension $n\in(0,\infty)$ and
$\rho\in\mathbf{q}$ have the property that all $\rho$-balls are $\mu$-measurable.
Let $p\in (0,\infty)$,
$s\in (\max\{1,\,\frac{1}{p}\},\infty)$,
$\tau\in (1,\infty)$, and
$\{B_k\}_{k\in\mathbb{N}}$ be a sequence of $\rho$-balls.
Then there exists a positive constant
$C$, independent of $\{B_k\}_{k\in\mathbb{N}}$ and $\tau$, such that
\begin{equation*}
\left\|\sum_{k\in\mathbb{N}}\mathbf{1}_{\tau B_k}\right\|_{L^p(X)}\leq C\tau
^{ns}\left\|\sum_{k\in\mathbb{N}}\mathbf{1}_{B_k}\right\|_{L^p(X)}.
\end{equation*}
\end{lemma}

\begin{proof}
Let $\rho_\#\in\mathbf{q}$ be the regularized quasi-ultrametric of $\rho$
(see Definition~\ref{D2.3}).
Using \eqref{eq-538}, \eqref{HLmaxequiv}, and \eqref{upperdoub}, we find that,
for any $k\in\mathbb{N}$ and $x\in\tau B_k$,
\begin{align*}
\mathcal{M}_{\rho_\#}(\mathbf{1}_{B_k})(x)
\gtrsim\widetilde{\mathcal{M}}_\rho(\mathbf{1}_{B_k})(x)
\ge\frac{1}{\mu(\tau B_k)}\int_{\tau B_k}\mathbf{1}_{B_k}(y)\,d\mu(y)
=\frac{\mu(B_k)}{\mu(\tau B_k)}
\ge\frac{1}{C_\mu\tau^n},
\end{align*}
which, combined with Lemma~\ref{Fefferman--Stein}
[here, we need $s\in(1,\infty]$ and $ps\in(1,\infty)$],
further implies that
\begin{align*}
\left\|\sum_{k\in\mathbb{N}}\mathbf{1}_{\tau B_k}\right\|_{L^p(X)}
&\lesssim\tau^{ns}
\left\|\sum_{k\in\mathbb{N}}\left[\mathcal{M}_{\rho_\#}
(\mathbf{1}_{B_k})\right]^s\right\|_{L^p(X)}
=\tau^{ns}
\left\|\left\{\sum_{k\in\mathbb{N}}\left[\mathcal{M}_{\rho_\#}(\mathbf{1}_{B_k})
\right]^s\right\}^\frac{1}{s}\right\|_{L^{ps}(X)}^s\\
&\lesssim\tau^{ns}
\left\|\left(\sum_{k\in\mathbb{N}}\mathbf{1}_{B_k}
\right)^\frac{1}{s}\right\|_{L^{ps}(X)}^s=\tau^{ns}
\left\|\sum_{k\in\mathbb{N}}\mathbf{1}_{B_k}\right\|_{L^p(X)}.
\end{align*}
This finishes the proof of Lemma~\ref{4r2}.
\end{proof}

\begin{lemma}\label{4l4-1}
Let $(X,\mathbf{q},\mu)$ be a quasi-ultrametric space
of homogeneous type with upper dimension $n\in(0,\infty)$  and
$\rho\in\mathbf{q}$ have the property that all $\rho$-balls are $\mu$-measurable.
Let $p\in (0,\infty)$, $t\in (0,p]$, and $r\in (p,\infty]\cap[1,\infty]$.
Then there exists a positive constant $C$ such that,
for any $\rho$-balls $\{B_k\}_{k\in\mathbb{N}}$,
any sequence $\{\eta_k\}_{k\in\mathbb{N}}$ in $\mathbb{C}$,
and any sequence of
measurable functions $\{\xi_k\}_{k\in\mathbb{N}}$
satisfying$\mathrm{\,supp\,} \xi_k\subset B_k$ and
$\|\xi_k\|_{L^r(X)}\leq [\mu(B_k)]^\frac{1}{r}$
for any $k\in\mathbb{N}$,
\begin{align*}
\left\|\left(\sum_{k\in\mathbb{N}}|\eta_k\xi_k|^t\right)^\frac{1}{t}\right\|_{L^p(X)}
\leq C\left\|\left(\sum_{k\in\mathbb{N}}
|\eta_k|^t\mathbf{1}_{B_k}\right)^\frac{1}{t}\right\|_{L^p(X)}.
\end{align*}
\end{lemma}

\begin{proof}
Let $\rho_\#\in\mathbf{q}$ be the regularized quasi-ultrametric of $\rho$.
From $(L^\frac{p}{t}(X))'=L^{(\frac{p}{t})'}(X)$,
Tonelli's theorem,
H\"older's inequality,
and the assumptions that
$$
\mathrm{supp\,} \xi_k\subset B_k
\ \ \text{and}\ \
\|\xi_k\|_{L^r(X)}\leq [\mu(B_k)]^\frac{1}{r},
$$
we infer that
\begin{align*}
\left\|\left(\sum_{k\in\mathbb{N}}|\eta_k\xi_k|^t\right)^\frac{1}{t}
\right\|_{L^p(X)}^t
&=\left\|\sum_{k\in\mathbb{N}}
|\eta_k\xi_k|^t\right\|_{L^\frac{p}{t}(X)}\\
&=\sup_{\|g\|_{L^{(\frac{p}{t})'}(X)}=1}\sum_{k\in\mathbb{N}}
|\eta_k|^t\int_X|\xi_k(x)|^tg(x)\,d\mu(x)\\
&\leq\sup_{\|g\|_{L^{(\frac{p}{t})'}(X)}=1}\sum_{k\in\mathbb{N}}
|\eta_k|^t\left\|\xi_k\right\|_{L^r(X)}^t
\left\|g\mathbf{1}_{B_k}\right\|_{L^{(\frac{r}{t})'}(X)}\\
&\leq\sup_{\|g\|_{L^{(\frac{p}{t})'}(X)}=1}\sum_{k\in\mathbb{N}}
|\eta_k|^t\left[\mu(B_k)\right]^\frac{t}{r}
\left\|g\mathbf{1}_{B_k}\right\|_{L^{(\frac{r}{t})'}(X)},
\end{align*}
which, together with \eqref{HLmaxequiv} and \eqref{eq-538},
further implies that
\begin{align*}
\left\|\left(\sum_{k\in\mathbb{N}}|\eta_k\xi_k|^t\right)^\frac{1}{t}
\right\|_{L^p(X)}^t
&\lesssim\sup_{\|g\|_{L^{(\frac{p}{t})'}(X)}=1}\sum_{k\in\mathbb{N}}
|\eta_k|^t\int_X\mathbf{1}_{B_k}(x)
\left[\mathcal{M}_{\rho_\#}(g^{(\frac{r}{t})'})(x)
\right]^\frac{1}{(\frac{r}{t})'}\,d\mu(x)\\
&\leq\sup_{\|g\|_{L^{(\frac{p}{t})'}(X)}=1}
\left\|\sum_{k\in\mathbb{N}}
|\eta_k|^t\mathbf{1}_{B_k}\right\|_{L^\frac{p}{t}(X)}
\left\|\left[\mathcal{M}_{\rho_\#}(g^{(\frac{r}{t})'})
\right]^\frac{1}{(\frac{r}{t})'}\right\|_{L^{(\frac{p}{t})'}(X)}\\
&\qquad\ \text{by H\"older's inequality}\\
&\lesssim\left\|\left(\sum_{k\in\mathbb{N}}
|\eta_k|^t\mathbf{1}_{B_k}\right)^\frac{1}{t}\right\|_{L^p(X)}^t
\quad\text{by Lemma~\ref{M}}.
\end{align*}
This finishes the proof of Lemma~\ref{4l4-1}.
\end{proof}

The following lemma is a generalization of \cite[Lemma 4.3]{ZYJLDW}
from unbounded spaces of homogeneous type satisfying $\mu(\{x\})=0$
for all $x\in X$ to quasi-ultrametric
spaces of homogeneous type.

\begin{lemma}\label{4l3}
Let $(X,\mathbf{q},\mu)$ be a quasi-ultrametric space
of homogeneous type with upper dimension $n\in(0,\infty)$ and
$\rho\in\mathbf{q}$ be a regularized quasi-ultrametric.
Assume that $p\in(\frac{n}{n+\mathrm{ind\,}(X,\mathbf{q})},1]$
and $q\in(p,\infty]\cap[1,\infty]$.
Then there exists a positive constant $C$ such that,
for any $(p,q)$-atom $a$
supported in $B:=B_\rho(x_B,r_B)$ with $x_B\in X$
and $r_B\in (0,\infty)$,
\begin{equation}
\label{4l3.1}
a^*_\rho(x)\leq
\begin{cases}
C\mathcal{M}_\rho(a)
\quad&\text{if}\ x\in2C_\rho B,\\
\displaystyle
C\left[\frac{r_B}{\rho(x_B,x)}\right]^{\beta}
\frac{[\mu(B)]^{1-\frac{1}{p}}}
{V_\rho(x_B,x)}
&\text{if}\ x\in[2C_\rho B]^\complement,
\end{cases}
\end{equation}
where the atom $a$ is viewed as
a distribution on ${\mathcal{G}}^\varepsilon_0(\beta,\gamma)$
in the sense described in Proposition~\ref{distfncttype},
where $\varepsilon,\beta,\gamma\in(0,\infty)$ satisfy
$n(\frac{1}{p}-1)<\beta,\gamma<\varepsilon\preceq\mathrm{ind\,}(X,\mathbf{q})$.
Consequently, there exists $C_1\in(0,\infty)$ such that
\begin{equation}\label{4l3.2}
\left\|a^*_\rho\right\|_{L^p(X)}\leq C_1.
\end{equation}
Moreover, if $\mu(X)<\infty$ and $a:=[\mu(X)]^{-\frac{1}{p}}$ is a constant function
on $X$, then there exists a positive constant $C_2$ such that,
for any $x\in X$,
$a^*_\rho(x)\leq C_2[\mu(X)]^{-\frac{1}{p}}$ and,
for any $r\in(0,\infty]$,
\begin{align}\label{20281}
\left\|a^*_\rho\right\|_{L^{p,r}(X)}\leq C_2.
\end{align}
\end{lemma}

\begin{proof}
We first prove \eqref{4l3.1}.
Let $x\in X$ and
$\phi\in\mathcal{G}^\varepsilon_0(\beta,\gamma)$ with
$\|\phi\|_{G(x,r^\alpha,\frac{\beta}{\alpha},\frac{\gamma}{\alpha})}\leq1$
for some $r\in(0,\infty)$,
where $\alpha:=\min\{1,\,(\log_2C_\rho)^{-1}\}$.
Then, from Remark~\ref{RemMf}(iii), we deduce that
$\|\phi\|_{\mathcal{G}(x,r,\beta,\gamma)}\lesssim1$.
On the one hand, if $x\in 2C_\rho B$,
then, by the size condition of $\phi$
and Lemma~\ref{A3.2}(iv), we find that
\begin{align}\label{2124}
\left|\langle a,\phi\rangle\right|
&\lesssim\int_X|a(y)|\frac{1}{V_r(x)+V_\rho(x,y)}
\left[\frac{r}{r+\rho(x,y)}\right]^\gamma\,d\mu(y)\nonumber\\
&=\int_{B_\rho(x,r)}|a(y)|\frac{1}{V_r(x)+V_\rho(x,y)}
\left[\frac{r}{r+\rho(x,y)}\right]^\gamma\,d\mu(y)\nonumber\\
&\quad+\int_{[B_\rho(x,r)]^\complement}|a(y)|\frac{1}{V_r(x)+V_\rho(x,y)}
\left[\frac{r}{r+\rho(x,y)}\right]^\gamma\,d\mu(y)\nonumber\\
&\leq\frac{1}{V_r(x)}\int_{B_\rho(x,r)}|a(y)|\,d\mu(y)
+\int_{[B_\rho(x,r)]^\complement}\frac{|a(y)|}{V_\rho(x,y)}
\left[\frac{r}{\rho(x,y)}\right]^\gamma\,d\mu(y)\nonumber\\
&\lesssim\mathcal{M}_\rho(a)(x).
\end{align}
On the other hand, if $x\notin 2C_\rho B$,
then, from (i) and (iii) of Definition~\ref{atom},
the regularity condition of $\phi$
[note that, for any $y\in\mathrm{supp\,}a=B$,
$\rho(x_B,y)<r_B\leq\frac{\rho(x_B,x)}{2C_\rho}
<\frac{\rho(x_B,x)+r}{2C_\rho}$], and
Definition~\ref{atom}(ii), it follows that
\begin{align*}
\left|\langle a,\phi\rangle\right|
&\leq\int_X|a(y)|\left|\phi(y)-\phi(x_B)\right|\,d\mu(y)
=\int_{B}|a(y)|\left|\phi(y)-\phi(x_B)\right|\,d\mu(y)
\\
&\lesssim\int_{B}|a(y)|
\left[\frac{\rho(x_B,y)}{r+\rho(x_B,x)}\right]^\beta
\frac{1}{(V_\rho)_r(x)+V_\rho(x_B,x)}
\left[\frac{r}{r+\rho(x_B,x)}\right]^\gamma\,d\mu(y)
\\
&<\left[\frac{r_B}{\rho(x_B,x)}\right]^\beta
\frac{1}{V_\rho(x_B,x)}\|a\|_{L^1(X)}
\leq\left[\frac{r_B}{\rho(x_B,x)}\right]^\beta
\frac{[\mu(B)]^{1-\frac{1}{p}}}{V_\rho(x_B,x)},
\end{align*}
which, combined with \eqref{2124}
and taking the supremum of $\phi$ as above,
then completes the proof of \eqref{4l3.1}.

Now, we show \eqref{4l3.2}.
We first estimate
$\|a_\rho^*\mathbf{1}_{2C_\rho B}\|_{L^p(X)}$.
To this end,
we consider the following two cases for $p$ and $q$.

\emph{Case (1)}
$p\in(\frac{n}{n+{\mathrm{ind\,}}(X,\mathbf{q})},1)$ and $q=1$.
In this case, we conclude that
\begin{align*}
\left\|a_\rho^*\mathbf{1}_{2C_\rho B}\right\|_{L^p(X)}^p
&\lesssim\int_{\rho(x_B,x)<2C_\rho r_B}
\left[\mathcal{M}_\rho(a)(x)\right]^p\,d\mu(x)
\quad\text{by \eqref{4l3.1}}
\nonumber\\
&=\int_0^\infty\lambda^{p-1}
\mu\left(\left\{x\in 2C_\rho B:
\mathcal{M}_\rho(a)(x)>\lambda\right\}\right)\,d\lambda
\nonumber\\
&\qquad\ \text{by Cavalieri's principle (see
\cite[Corollary~1.2.4]{G2024})}\\
&\lesssim\int_0^\infty\lambda^{p-1}
\min\left\{\mu(2C_\rho B),\,\frac{\|a\|_{L^1(X)}}{\lambda}\right\}\,d\lambda
\quad\text{by Lemma~\ref{M}},
\end{align*}
which, together with Definition~\ref{atom}(ii)
and \eqref{upperdoub},
further implies that
\begin{align}\label{1858}
\left\|a_\rho^*\mathbf{1}_{2C_\rho B}\right\|_{L^p(X)}^p
&\lesssim\int_0^\frac{\|a\|_{L^1(X)}}{\mu(2C_\rho B)}
\lambda^{p-1}\mu(2C_\rho B)\,d\lambda
+\int_\frac{\|a\|_{L^1(X)}}{\mu(2C_\rho B)}^\infty
\lambda^{p-1}\frac{\|a\|_{L^1(X)}}{\lambda}\,d\lambda
\nonumber\\
&\sim\|a\|_{L^1(X)}^p\left[\mu(2C_\rho B)\right]^{1-p}
\leq\left[\frac{\mu(2C_\rho B)}{\mu(B)}\right]^{1-p}
\sim1.
\end{align}

\emph{Case (2)}
$p\in(\frac{n}{n+{\mathrm{ind\,}}(X,\mathbf{q})},1]$
and $q\in(\max\{1,\,p\},\infty]$.
In this case, we infer that
\begin{align}\label{1905}
\left\|a_\rho^*\mathbf{1}_{2C_\rho B}\right\|_{L^p(X)}
&\lesssim\left\|\mathcal{M}_\rho(a)\mathbf{1}_{2C_\rho B}\right\|_{L^p(X)}
\quad\text{by \eqref{4l3.1}}\nonumber\\
&\lesssim\left\|\mathcal{M}_\rho(a)\right\|_{L^q(X)}
\left[\mu(2C_\rho B)\right]^{\frac{1}{p(\frac{q}{p})'}}
\quad\text{by H\"older's inequality}\nonumber\\
&\lesssim\left\|a\right\|_{L^q(X)}
\left[\mu(2C_\rho B)\right]^{\frac{1}{p(\frac{q}{p})'}}
\quad\text{by Lemma~\ref{M}}
\nonumber\\
&\lesssim\left\|a\right\|_{L^q(X)}
\left[\mu(B)\right]^{\frac{1}{p(\frac{q}{p})'}}
\quad\text{by \eqref{upperdoub}}
\nonumber\\
&\leq\left[\mu(B)\right]^{\frac{1}{q}-\frac{1}{p}+\frac{1}{p(\frac{q}{p})'}}=1\quad\text{by Definition~\ref{atom}(ii)}.
\end{align}
On the other hand, to estimate
$\|a_\rho^*\mathbf{1}_{(2C_\rho B)^\complement}\|_{L^p(X)}$,
we find that
\begin{align*}
\left\|a_\rho^*\mathbf{1}_{(2C_\rho B)^\complement}\right\|_{L^p(X)}^p
&\lesssim\int_{(2C_\rho B)^\complement}\left[\frac{r_B}
{\rho(x_B,x)}\right]^{(\beta\land\gamma)p}
\frac{[\mu(B)]^{p-1}}{[V_\rho(x_B,x)]^p}\,d\mu(x)
\quad\text{by \eqref{4l3.1}}
\\
&=\sum_{k=1}^\infty
\int_{2^kC_\rho r_B\leq\rho(x_B,x)<2^{k+1}C_\rho r_B }\left[\frac{r_B}
{\rho(x_B,x)}\right]^{(\beta\land\gamma)p}
\frac{[\mu(B)]^{p-1}}{[V_\rho(x_B,x)]^p}\,d\mu(x)
\\
&\lesssim\sum_{k=1}^\infty
2^{-(\beta\land\gamma)kp}
\int_{2^kC_\rho r_B\leq\rho(x_B,x)<2^{k+1}C_\rho r_B}
\frac{[\mu(B)]^{p-1}}
{[\mu(2^kC_\rho B)]^p}\,d\mu(x)
\\
&\lesssim\sum_{k=1}^\infty
2^{-(\beta\land\gamma)kp}
\left[\frac{\mu(B)}
{\mu(2^kC_\rho B)}\right]^{p-1}
\quad\text{by \eqref{upperdoub}}\\
&\lesssim\sum_{k=1}^\infty
2^{-kp[\beta\land\gamma-n(\frac{1}{p}-1)]}
\sim1
\end{align*}
since $n(\frac{1}{p}-1)<\beta,\gamma$.
From this, \eqref{1858}, and \eqref{1905},
we deduce that
$$
\left\|a^*_\rho\right\|_{L^p(X)}
\lesssim\left\|a^*_\rho\mathbf{1}_{2C_\rho B}\right\|_{L^p(X)}
+\left\|a_\rho^*\mathbf{1}_{(2C_\rho B)^\complement}\right\|_{L^p(X)}
\lesssim1,
$$
which completes the proof of \eqref{4l3.2}.

Finally, assume that $\mu(X)<\infty$ and
$a:=[\mu(X)]^{-\frac{1}{p}}$ is a constant function
on $X$. Using Remark~\ref{RemMf}(iii),
Definition~\ref{testfunction}(i),
and Lemma~\ref{A3.2}(ii), we conclude that,
for any given $x\in X$ and
$\phi\in\mathcal{G}^\varepsilon_0(\beta,\gamma)$ with
$\|\phi\|_{G(x,r^\alpha,\frac{\beta}{\alpha},\frac{\gamma}{\alpha})}\leq1$
for some $r\in(0,\infty)$,
where $\alpha:=\min\{1,\,(\log_2C_\rho)^{-1}\}$,
\begin{align*}
\left|\langle a,\phi\rangle\right|
\lesssim\left[\mu(X)\right]^{-\frac{1}{p}}\int_X\frac{1}{V_r(x)+V_\rho(x,y)}
\left[\frac{r}{r+\rho(x,y)}\right]^\gamma\,d\mu(y)
\lesssim\left[\mu(X)\right]^{-\frac{1}{p}}.
\end{align*}
Taking the supremum over all $\phi$ as above,
we then obtain, for any $x\in X$,
$$
a^*_\rho(x)\lesssim[\mu(X)]^{-\frac{1}{p}}.
$$
From this and both (ii) and (iii) of Proposition~\ref{inc-norm},
it follows that, for any $r\in(0,\infty]$,
$$
\left\|a^*_\rho\right\|_{L^{p,r}(X)}
\lesssim[\mu(X)]^{-\frac{1}{p}}
\left\|\mathbf{1}_X\right\|_{L^{p,r}(X)}\sim1,
$$
which completes the proof of \eqref{20281}
and hence Lemma~\ref{4l3}.
\end{proof}

Next, we prove Theorem~\ref{6.15}.

\begin{proof}[Proof of Theorem~\ref{6.15}]
We only show the present theorem in the case where both
$q\in(0,\infty)$ and $\mu(X)=\infty$
because the proof in the case where $q=\infty$ or $\mu(X)<\infty$
[when $\mu(X)<\infty$, we need to take advantage of \eqref{20281}]
is similar and hence we omit the details.
Let $\{a_{i,k}\}_{i\in\mathbb{N},k\in\mathbb{Z}}$,
$\{B_{i,k}\}_{i\in\mathbb{N},k\in\mathbb{Z}}
:=\{B_\rho(x_{i,k},r_{i,k})\}_{i\in\mathbb{N},k\in\mathbb{Z}}$
with $\{x_{i,k}\}_{i\in\mathbb{N},k\in\mathbb{Z}}$ in $X$
and $\{r_{i,k}\}_{i\in\mathbb{N},k\in\mathbb{Z}}$ in $(0,\infty)$,
and
$\{\lambda_{i,k}\}_{i\in\mathbb{N},k\in\mathbb{Z}}$ in
$\mathbb{R}_+$ be the same as in the present theorem.
According to Remark~\ref{5216}, we may assume that
$\rho$ is a regularized quasi-ultrametric as in Definition~\ref{D2.3}.

We first prove (i).
To this end, let $\phi\in{\mathcal{G}}^\varepsilon_0(\beta,\gamma)$
with $\|\phi\|_{{\mathcal{G}}^\varepsilon_0(\beta,\gamma)}\leq1$.
Let $x_1\in X$ be the same as in Definition~\ref{testfunction}.
Then, by an argument similar to that used in the proof of
Proposition~\ref{5.1.10}, we find that
\begin{align}\label{16177}
\sum_{k\in\mathbb{Z}}\sum_{i\in\mathbb{N}}
\lambda_{i,k}\left|\langle a_{i,k},\phi\rangle\right|
&\lesssim\sum_{k\in\mathbb{Z}}\sum_{i\in\mathbb{N}}
\lambda_{i,k}\inf_{x\in B_\rho(x_1,1)}(a_{i,k})^*_\rho(x)\nonumber\\
&\leq\inf_{x\in B_\rho(x_1,1)}\sum_{k\in\mathbb{Z}}\sum_{i\in\mathbb{N}}
\lambda_{i,k}(a_{i,k})^*_\rho(x)\nonumber\\
&\leq\frac{\|\sum_{k\in\mathbb{Z}}\sum_{i\in\mathbb{N}}
\lambda_{i,k}(a_{i,k})^*_\rho\|_{L^{p,q}(X)}}
{\|\mathbf{1}_{B_\rho(x_1,1)}\|_{L^{p,q}(X)}}.
\end{align}
It remains to estimate
$\|\sum_{k\in\mathbb{Z}}\sum_{i\in\mathbb{N}}
\lambda_{i,k}(a_{i,k})^*_\rho\|_{L^{p,q}(X)}$.
From Proposition~\ref{Lpqnorm=c,d},
we infer that
\begin{align}\label{1059}
\left\|\sum_{k\in\mathbb{Z}}\sum_{i\in\mathbb{N}}
\lambda_{i,k}(a_{i,k})^*_\rho\right\|_{L^{p,q}(X)}
\sim\left\{\sum_{k_0\in\mathbb{Z}}2^{k_0q}
\left[d_g(2^{k_0})\right]^\frac{q}{p}\right\}^\frac{1}{q},
\end{align}
where, for any given $k_0\in\mathbb{Z}$, we write
$$
g:=\sum_{k\in\mathbb{Z}}\sum_{i\in\mathbb{N}}
\lambda_{i,k}(a_{i,k})^*_\rho
=\sum_{k=-\infty}^{k_0-1}\sum_{i\in\mathbb{N}}
\lambda_{i,k}(a_{i,k})^*_\rho
+\sum_{k=k_0}^{\infty}
\sum_{i\in\mathbb{N}}\ldots=:g_1+g_2.
$$
Fix $k_0\in\mathbb{Z}$
and let
$$
E_{k_0}:=\bigcup_{k=k_0}^{\infty}\bigcup_{i\in\mathbb{N}}
2C_\rho B_{i,k}.
$$
To estimate \eqref{1059},
we first write
\begin{align}
\label{730}
d_{g}(2^{k_0})
&\leq\mu\left(\left\{x\in X:
g_1(x)>2^{k_0-1}\right\}\right)
\nonumber\\
&\quad+\mu\left(\left\{x\in E_{k_0}:g_2(x)>2^{k_0-2}\right\}\right)
\nonumber\\
&\quad+\mu\left(\left\{x\in(E_{k_0})^\complement:g_2(x)>2^{k_0-2}\right\}\right)
\nonumber\\
&=:(\mathrm J_1)^p+(\mathrm J_2)^p+(\mathrm J_3)^p.
\end{align}

Now, we estimate each of $\mathrm J_1$, $\mathrm J_2$, and $\mathrm J_3$.	
To estimate $\mathrm J_1$, by the sublinearity of the grand maximal
function, we obtain,
for any given $\delta\in(0,1)$,
\begin{align*}
\mathrm J_1
&\lesssim\left[\mu\left(\left\{x\in X:
\sum_{k=-\infty}^{k_0-1}\sum_{i\in\mathbb{N}}\lambda_{i,k}(a_{i,k})^*_\rho(x)
\mathbf{1}_{2C_\rho B_{i,k}}(x)>2^{k_0-2}\right\}\right)\right]^\frac{1}{p}
\\
&\quad+\left[\mu\left(\left\{x\in X:
\sum_{k=-\infty}^{k_0-1}
\sum_{i\in\mathbb{N}}\lambda_{i,k}(a_{i,k})^*_\rho(x)
\mathbf{1}_{(2C_\rho B_{i,k})^{\complement}}(x)
>2^{k_0-2}\right\}\right)\right]^\frac{1}{p}
\\
&=:\mathrm J_{1,1}+\mathrm J_{1,2}.
\end{align*}

We first estimate $\mathrm J_{1,1}$. Let
$t\in(0,\min\{p,\,q\})$,
$\widetilde{q}\in(1,\min\{r,\,\frac{1}{t}\})$,
$\delta\in(0,1-\frac{1}{\widetilde{q}})$, and
$$
C_{k_0,\delta,\widetilde{q}}:=\frac{2^{k_0\delta}}{[2^\frac{\delta\widetilde{q}}
{\widetilde{q}-1}-1]^\frac{\widetilde{q}-1}{\widetilde{q}}}.
$$
Then
\begin{align*}
\mathrm J_{1,1}
&\leq\left[\mu\left(\left\{x\in X:
C_{k_0,\delta,\widetilde{q}}
\left\{\sum_{k=-\infty}^{k_0-1}2^{-k\delta\widetilde{q}}
\left[\sum_{i\in\mathbb{N}}\lambda_{i,k}(a_{i,k})^*_\rho(x)
\mathbf{1}_{2C_\rho B_{i,k}}(x)\right]^{\widetilde{q}}
\right\}^\frac{1}{\widetilde{q}}>2^{k_0-2}\right\}
\right)\right]^\frac{1}{p}\nonumber\\
&\qquad\ \text{by H\"older's inequality}\nonumber\\
&\leq2^{k_0\widetilde{q}(\delta-1)}
\left\|\sum_{k=-\infty}^{k_0-1}2^{-k\delta\widetilde{q}}
\left[\sum_{i\in\mathbb{N}}\lambda_{i,k}(a_{i,k})^*_\rho
\mathbf{1}_{2C_\rho B_{i,k}}\right]^{\widetilde{q}}
\right\|_{L^p(X)}
\quad\text{by Chebyshev's inequality}\nonumber\\
&=2^{k_0\widetilde{q}(\delta-1)}
\left\|\left\{\sum_{k=-\infty}^{k_0-1}2^{-k\delta\widetilde{q}}
\left[\sum_{i\in\mathbb{N}}\lambda_{i,k}(a_{i,k})^*_\rho
\mathbf{1}_{2C_\rho B_{i,k}}\right]^{\widetilde{q}}\right\}
^t\right\|_{L^\frac{p}{t}(X)}^\frac{1}{t}
\quad\text{by Proposition~\ref{r;p,q=pr,qr;r}}.
\end{align*}
From this, Lemma~\ref{discreteinequality},
\eqref{lambdaik}, \eqref{*<M},
and Minkowski's inequality,
we deduce that
\begin{align}
\label{xxk-749}
\mathrm J_{1,1}
&\lesssim2^{k_0\widetilde{q}(\delta-1)}
\left\|\sum_{k=-\infty}^{k_0-1}2^{-k\delta\widetilde{q}t}
\sum_{i\in\mathbb{N}}\left[\lambda_{i,k}(a_{i,k})^*_\rho
\mathbf{1}_{2C_\rho B_{i,k}}\right]^{\widetilde{q}t}
\right\|_{L^\frac{p}{t}(X)}^\frac{1}{t}\nonumber\\
&\lesssim2^{k_0\widetilde{q}(\delta-1)}\left\|\sum_{k=-\infty}^{k_0-1}
2^{(1-\delta)k\widetilde{q} t}\sum_{i\in\mathbb{N}}
\left\{[\mu(B_{i,k})]^\frac{1}{p}\mathcal{M}_\rho(a_{i,k})
\mathbf{1}_{2C_\rho B_{i,k}}\right\}^{\widetilde{q}t}
\right\|^\frac{1}{t}_{L^\frac{p}{t}(X)}\nonumber\\
&\lesssim2^{k_0\widetilde{q}(\delta-1)}
\left\{\sum_{k=-\infty}^{k_0-1}2^{(1-\delta)k\widetilde{q} t}\left\|
\left[\sum_{i\in\mathbb{N}}\left\{[\mu(B_{i,k})]^\frac{1}{p}
\mathcal{M}_\rho(a_{i,k})\mathbf{1}_{2C_\rho B_{i,k}}
\right\}^{\widetilde{q} t}\right]
^\frac{1}{t}\right\|^t_{L^p(X)}\right\}^\frac{1}{t}.
\end{align}
With the goal of employing Lemma~\ref{4l4-1},
by Lemma~\ref{M} and Definition~\ref{atom}(ii),
we conclude that,
for any $k\in\mathbb{Z}$ and $i\in\mathbb{N}$,
\begin{align*}
&\left\|\left\{\left[\mu(B_{i,k})\right]^\frac{1}{p}
\mathcal{M}_\rho(a_{i,k})\mathbf{1}_{2C_\rho B_{i,k}}
\right\}^{\widetilde{q}}\right\|
_{L^{\frac{r}{\widetilde{q}}}(X)}
\\
&\quad\leq\left[\mu(B_{i,k})\right]^{\frac{\widetilde{q}}{p}}
\left\|\mathcal{M}_\rho(a_{i,k})\right\|^{\widetilde{q}}_{L^r(X)}
\lesssim\left[\mu(B_{i,k})\right]^{\frac{\widetilde{q}}{p}}
\left\|a_{i,k}\right\|^{\widetilde{q}}_{L^r(X)}
\leq[\mu(2C_\rho B_{i,k})]^\frac{\widetilde{q}}{r}.
\end{align*}
From this, \eqref{xxk-749}, Lemma~\ref{4l4-1}
[used here with $r$, $\eta_k$, and $\xi_k$ replaced,
respectively, by $\frac{r}{\widetilde{q}}$,
$1$, and $\{[\mu(B_{i,k})]^\frac{1}{p}
\mathcal{M}_\rho(a_{i,k})\mathbf{1}_{2C_\rho B_{i,k}}
\}^{\widetilde{q}}$ (modulo a harmless constant)],
and Lemma~\ref{4r2} with $s\in(0,\min\{1,\,\frac{p}{t}\}\}$ and
$\tau:=(2C_\rho)^{1-j_0}>1$,
we infer that,
for any fixed $\eta\in(1,(1-\delta)\widetilde{q})$,
\begin{align*}
\mathrm{J}_{1,1}
&\lesssim2^{k_0\widetilde{q}(\delta-1)}
\left[\sum_{k=-\infty}^{k_0-1}2^{(1-\delta)k\widetilde{q}t}
\left\|\left(\sum_{i\in\mathbb{N}}\mathbf{1}_{2C_\rho B_{i,k}}
\right)^\frac{1}{t}\right\|^t_{L^p(X)}\right]^\frac{1}{t}
\\
&\sim2^{k_0\widetilde{q}(\delta-1)}
\left[\sum_{k=-\infty}^{k_0-1}2^{(1-\delta)k\widetilde{q}t}
\left\|\sum_{i\in\mathbb{N}}\mathbf{1}_{\tau(2C_\rho)^{j_0}
B_{i,k}}\right\|_{L^\frac{p}{t}(X)}\right]^\frac{1}{t}
\\
&\lesssim2^{k_0\widetilde{q}(\delta-1)}
\left\{\sum_{k=-\infty}^{k_0-1}2^{(1-\delta)k\widetilde{q}t}
\left\|\sum_{i\in\mathbb{N}}\mathbf{1}_{(2C_\rho)^{j_0}
B_{i,k}}\right\|_{L^\frac{p}{t}(X)}\right\}^\frac{1}{t},
\end{align*}
which, combined with \eqref{4ee0},
Lemma~\ref{discreteinequality} with $r$
therein replaced by $p$,
H\"older's inequality, and \eqref{lambdaik},
further implies that
\begin{align*}
\mathrm{J}_{1,1}
&\lesssim2^{k_0\widetilde{q}(\delta-1)}\left\{\sum_{k=-\infty}^{k_0-1}
2^{(1-\delta)k\widetilde{q}t}\left\|\sum_{i\in\mathbb{N}}
\mathbf{1}_{(2C_\rho)^{j_0}B_{i,k}}\right\|^t_{L^p(X)}\right\}^\frac{1}{t}
\\
&\lesssim2^{k_0\widetilde{q}(\delta-1)}
\left\{\sum_{k=-\infty}^{k_0-1}2^{[(1-\delta)\widetilde{q}-\eta]kt}
2^{\eta kt}\left[\sum_{i\in\mathbb{N}}\mu((2C_\rho)^{j_0}
B_{i,k})\right]^\frac{t}{p}\right\}^\frac{1}{t}
\\
&\lesssim2^{k_0\widetilde{q}(\delta-1)}
\left\{\sum_{k=-\infty}^{k_0-1}2^{[(1-\delta)\widetilde{q}-\eta]kt
\frac{q}{q-t}}\right\}^{\frac{q-t}{qt}}
\left\{\sum_{k=-\infty}^{k_0-1}2^{\eta kq}
\left[\sum_{i\in\mathbb{N}}\mu(B_{i,k})\right]^\frac{q}{p}\right\}^\frac{1}{q}
\\
&\sim 2^{-\eta k_0}
\left\{\sum_{k=-\infty}^{k_0-1}2^{\eta kq}
\left[\sum_{i\in\mathbb{N}}\mu(B_{i,k})\right]^\frac{q}{p}\right\}^\frac{1}{q}
\sim2^{-\eta k_0}\left[\sum_{k=-\infty}^{k_0-1}
2^{(\eta-1) kq}\left(\sum_{i\in\mathbb{N}}
\lambda_{i,k}^p\right)^\frac{q}{p}\right]^\frac{1}{q},
\end{align*}
where the implicit positive constants are independent of $k_0$.
By this and Tonelli's theorem, we find that
\begin{align}
\label{4ee5}
\sum_{k_0\in\mathbb{Z}}2^{k_0q}(\mathrm{J}_{1,1})^q
&\lesssim\sum_{k_0\in\mathbb{Z}}2^{k_0q}
2^{-\eta k_0q}\sum_{k=-\infty}^{k_0-1}
2^{(\eta-1) kq}\left(\sum_{i\in\mathbb{N}}
\lambda_{i,k}^p\right)^\frac{q}{p}
\nonumber\\
&=\sum_{k\in\mathbb{Z}}
\sum_{k_0=k+1}^{\infty}2^{(k_0-k)(1-\eta)q}
\left(\sum_{i\in\mathbb{N}}
\lambda_{i,k}^p\right)^\frac{q}{p}
\sim\sum_{k\in\mathbb{Z}}
\left(\sum_{i\in\mathbb{N}}\lambda_{i,k}^p\right)
^{\frac{q}{p}}.
\end{align}

Next, we estimate $\mathrm J_{1,2}$.
From Definition~\ref{D2.1}(iii), we deduce that,
for any $i\in\mathbb{N}$, $k\in\mathbb{Z}$, $y\in B_{i,k}$, and
$x\in(2C_\rho B_{i,k})^{\complement}$,
$$
\rho(x,y)\leq C_\rho\max\left\{\rho(x,x_{i,k}),\,\rho(x_{i,k},y)\right\}
=C_\rho\rho(x,x_{i,k})
$$
and hence
$B_{i,k}\subset B_\rho(x,C_\rho \rho(x,x_{i,k}))$,
which, together with \eqref{upperdoub},
further implies that
\begin{align}\label{4ee6-1}
\mathcal{M}_\rho(\mathbf{1}_{B_{i,k}})(x)
&\geq\fint_{B_\rho(x,C_\rho\rho(x,x_{i,k}))}
\mathbf{1}_{B_{i,k}}(y)\,d\mu(y)
\nonumber\\
&=\frac{\mu(B_{i,k})}{\mu(B_\rho(x,C_\rho \rho(x,x_{i,k})))}
\gtrsim\frac{\mu(B_{i,k})}{V_\rho(x,x_{i,k})}.
\end{align}
By Lemma~\ref{A3.2}(vii) and \eqref{upperdoub},
we conclude that,
for any $i\in\mathbb{N}$, $k\in\mathbb{Z}$,
and $x\in (2C_\rho B_{i,k})^{\complement}$,
\begin{equation}
\label{4ee6-2}
\frac{\mu(B_{i,k})}{V_\rho(x,x_{i,k})}\sim \frac{\mu(B_{i,k})}{V_\rho(x_{i,k},x)}
\gtrsim\left[\frac{r_{i,k}}{\rho(x_{i,k},x)}\right]^n.
\end{equation}
From this, \eqref{4l3.1},
and \eqref{4ee6-1},
we infer that, for any
$i\in\mathbb{N}$, $k\in\mathbb{Z}$, and $x\in (2C_\rho B_{i,k})^{\complement}$,
\begin{align}
\label{4ee6}
({a_{i,k}})^*_\rho(x)
&\lesssim\left[\frac{r_{i,k}}{\rho(x_{i,k},x)}\right]
^{\beta}\frac{\left[\mu(B_{i,k})\right]^{1-\frac{1}{p}}}
{V_\rho(x_{i,k},x)}
\lesssim\frac{\left[\mu(B_{i,k})\right]^{1+\frac{\beta}{n}-\frac{1}{p}}}
{[V_\rho(x_{i,k},x)]^{1+\frac{\beta}{n}}}
\nonumber\\
&\lesssim\left[\mu(B_{i,k})\right]^{-\frac{1}{p}}
\left[\mathcal{M}_\rho(\mathbf{1}_{B_{i,k}})(x)\right]^{1+\frac{\beta}{n}}.
\end{align}

Note that $1+\frac{\beta}{n}>\frac{1}{p}$.
Choose $t\in (0,\min\{p,\,q,\,\frac{1}{1+\frac{\beta}{n}}\})$,
$q_1\in(\frac{1}{(1+\frac{\beta}{n})t},\frac{1}{t})\subset(1,\infty)$,
and $\delta\in(0,1-\frac{1}{q_1})$.
By Chebyshev's inequality,
\eqref{4ee6}, \eqref{lambdaik},
the fact that $t<1$,
Minkowski's inequality with the norm $\|\cdot\|_{L^\frac{p}{t}(X)}$,
and Lemma~\ref{discreteinequality} with $r$ therein replaced by $q_1t$,
we find that,
for any fixed $\eta\in (1,(1-\delta)q_1)$,
\begin{align*}
\mathrm{J}_{1,2}
&\lesssim2^{k_0q_1(\delta-1)}\left\|\sum_{k=-\infty}^{k_0-1}2^{-k\delta q_1}
\left[\sum_{i\in\mathbb{N}}\lambda_{i,k}(a_{i,k})^*_\rho
\mathbf{1}_{(2C_\rho B_{i,k})^{\complement}}\right]^{q_1}\right\|_{L^p(X)}
\\
&\lesssim2^{k_0q_1(\delta-1)}\left\|\sum_{k=-\infty}^{k_0-1}2^{(1-\delta)kq_1}
\left\{\sum_{i\in\mathbb{N}}\left[\mathcal{M}_\rho
(\mathbf{1}_{B_{i,k}})(x)\right]^{1+\frac{\beta}{n}}
\right\}^{q_1}\right\|_{L^p(X)}
\\
&\leq2^{k_0q_1(\delta-1)}\left\|
\left(\sum_{k=-\infty}^{k_0-1}2^{(1-\delta)kq_1t}
\left\{\sum_{i\in\mathbb{N}}\left[\mathcal{M}_\rho
(\mathbf{1}_{B_{i,k}})(x)\right]^{1+\frac{\beta}{n}}
\right\}^{q_1t}\right)^\frac{1}{t}\right\|_{L^p(X)}
\\
&\lesssim2^{k_0q_1(\delta-1)}\left\{\sum_{k=-\infty}^{k_0-1}2^{(1-\delta)kq_1t}
\left\|\sum_{i\in\mathbb{N}}\left[\mathcal{M}_\rho(\mathbf{1}_{B_{i,k}})\right]
^{(1+\frac{\beta}{n})tq_1}\right\|_{L^\frac{p}{t}(X)}\right\}^\frac{1}{t}
\\
&\lesssim2^{k_0q_1(\delta-1)}\left\{\sum_{k=-\infty}^{k_0-1}2^{(1-\delta)kq_1t}
\left\|\left(\sum_{i\in\mathbb{N}}\left[\mathcal{M}_\rho(\mathbf{1}_{B_{i,k}})\right]
^{u}\right)^\frac{1}{u}\right\|_{L^\frac{pu}{t}(X)}^u\right\}^\frac{1}{t},
\end{align*}
which, combined with Lemma \ref{Fefferman--Stein} with $u:=(1+\frac{\beta}{n})tq_1>1$ and
$p$ replaced by $\frac{pu}{t}>1$,
Lemma~\ref{4r2} with $\tau:=(2C_\rho)^{-j_0}>1$, \eqref{4ee0},
and H\"older's inequality with $\frac{q}{t}>1$,
further implies that
\begin{align*}
\mathrm{J}_{1,2}
&\lesssim2^{k_0q_1(\delta-1)}\left[\sum_{k=-\infty}^{k_0-1}2^{(1-\delta)kq_1t}
\left\|\left(\sum_{i\in\mathbb{N}}\mathbf{1}_{B_{i,k}}
\right)^\frac{1}{u}\right\|_{L^\frac{pu}{t}(X)}^u\right]^\frac{1}{t}
\\
&\lesssim2^{k_0q_1(\delta-1)}
\left\{\sum_{k=-\infty}^{k_0-1}2^{(1-\delta)kq_1t}
\left\|\left[\sum_{i\in\mathbb{N}}
\mathbf{1}_{(2C_\rho)^{j_0}B_{i,k}}\right]^\frac{1}{t}\right\|
^t_{L^p(X)}\right\}^\frac{1}{t}
\\
&\lesssim2^{k_0q_1(\delta-1)}
\left\{\sum_{k=-\infty}^{k_0-1}2^{(1-\delta)kq_1t}
\left\|\sum_{i\in\mathbb{N}}\mathbf{1}_{(2C_\rho)^{j_0}B_{i,k}}\right\|
^t_{L^p(X)}\right\}^\frac{1}{t}
\\
&\lesssim2^{k_0q_1(\delta-1)}
\left\{\sum_{k=-\infty}^{k_0-1}2^{[(1-\delta)q_1-\eta]kt}2^{\eta kt}
\left[\sum_{i\in\mathbb{N}}\mu\big((2C_\rho)^{j_0}B_{i,k}\big)\right]
^\frac{t}{p}\right\}^\frac{1}{t}
\\
&\lesssim2^{-\eta k_0}
\left\{\sum_{k=-\infty}^{k_0-1}2^{\eta kq}
\left[\sum_{i\in\mathbb{N}}\mu(B_{i,k})\right]^\frac{q}{p}\right\}
^\frac{1}{q}
\sim2^{-\eta k_0}
\left[\sum_{k=-\infty}^{k_0-1}2^{(\eta-1)kq}
\left(\sum_{i\in\mathbb{N}}\lambda_{i,k}^p\right)^\frac{q}{p}\right]^\frac{1}{q},
\end{align*}
where the implicit positive constants are independent of $k_0$.
From this and an argument similar to
that used in the estimation of \eqref{4ee5},
we deduce that
\begin{align}
\label{2007}
\left[\sum_{k_0\in\mathbb{Z}}2^{k_0q}
(\mathrm{J}_{1,2})^q\right]^\frac{1}{q}
\lesssim\left[\sum_{k\in\mathbb{Z}}
\left(\sum_{i\in\mathbb{N}}\lambda_{i,k}^p\right)
^{\frac{q}{p}}\right]^{\frac{1}{q}},
\end{align}
which, together with \eqref{4ee5},
further implies that
\begin{align}\label{J1}
\left[\sum_{k_0\in\mathbb{Z}}2^{k_0q}(\mathrm J_{1})^{q}\right]^\frac{1}{q}
&\lesssim\left[\sum_{k_0\in\mathbb{Z}}2^{k_0q}(\mathrm{J}_{1,1})^q\right]^\frac{1}{q}
+\left[\sum_{k_0\in\mathbb{Z}}2^{k_0q}(\mathrm{J}_{1,2})^q\right]^\frac{1}{q}
\nonumber\\
&\lesssim\left[\sum_{k\in\mathbb{Z}}\left(\sum_{i\in\mathbb{N}}\lambda_{i,k}^p\right)
^{\frac{q}{p}}\right]^{\frac{1}{q}}.
\end{align}

To estimate $\mathrm{J}_2$, by
\eqref{upperdoub}, Lemma~\ref{discreteinequality} with
$r$ therein replaced by $\frac{t}{p}$,
H\"older's inequality, and \eqref{lambdaik},
we conclude that, for any fixed $t\in(0,\min\{p,\,q\})$ and $\delta\in(0,1)$,
\begin{align*}
\mathrm{J}_2
&\leq\left[\mu(E_{k_0})\right]^\frac{1}{p}
\leq\left[\sum_{k=k_0}^{\infty}
\sum_{i\in\mathbb{N}}\mu(2C_\rho B_{i,k})\right]^\frac{1}{p}
\lesssim\left[\sum_{k=k_0}^{\infty}
\sum_{i\in\mathbb{N}}\mu(B_{i,k})\right]^\frac{1}{p}
\\
&\leq\left\{\sum_{k=k_0}^{\infty}
\left[\sum_{i\in\mathbb{N}}\mu(B_{i,k})
\right]^\frac{t}{p}\right\}^\frac{1}{t}
\lesssim\left(\sum_{k=k_0}^{\infty}
2^{-k\delta t\frac{q}{q-t}}\right)^{\frac{q-t}{qt}}
\left\{\sum_{k=k_0}^{\infty}2^{k\delta q}
\left[\sum_{i\in\mathbb{N}}\mu(B_{i,k})\right]^\frac{q}{p}\right\}^\frac{1}{q}
\\
&\sim2^{-k_0\delta}
\left\{\sum_{k=k_0}^{\infty}2^{k\delta q}
\left[\sum_{i\in\mathbb{N}}\mu(B_{i,k})\right]^\frac{q}{p}\right\}^\frac{1}{q}
\sim 2^{-k_0\delta}
\left[\sum_{k=k_0}^{\infty}2^{kq(\delta -1)}
\left(\sum_{i\in\mathbb{N}}\lambda_{i,k}^p\right)^\frac{q}{p}\right]^\frac{1}{q},
\end{align*}
where the implicit positive constants are independent of $k_0$.
This, combined with an argument similar to
that used in the estimation of \eqref{4ee5},
implies that
\begin{equation}
\label{4ee8}
\left[\sum_{k_0\in\mathbb{Z}}2^{k_0q}
(\mathrm J_{2})^{q}\right]^\frac{1}{q}
\lesssim\left[\sum_{k\in\mathbb{Z}}\left(\sum_{i\in\mathbb{N}}\lambda_{i,k}^p\right)
^{\frac{q}{p}}\right]^{\frac{1}{q}}.
\end{equation}

To estimate $\mathrm J_3$,
from the fact that $1+\frac{\beta}{n}>\frac{1}{p}$, we infer
that there exist $t\in(0,1)$ and $L\in (1,\infty)$
such that
$$
\left(1+\frac{\beta}{n}\right)t>\frac{1}{p}\geq1
\ \ \text{and}\ \
\left(1+\frac{\beta}{n}\right)qtL>1.
$$
By this, Chebyshev's inequality,
Lemma~\ref{discreteinequality} (used three times with $t<1$,
$\frac{1}{(1+\frac{\beta}{n})t}\leq1$, and $\frac{1}{L}<1$),
and Minkowski's inequality with the norm $\|\cdot\|_{L^{(1+\frac{\beta}{n})pt}}$,
we find that, for any fixed $\delta\in(t,1)$,
\begin{align*}
\mathrm J_3
&\lesssim 2^{-tk_0}\left\|\left\{\sum_{k=k_0}^{\infty}\sum_{i\in\mathbb{N}}
|\lambda_{i,k}|(a_{i,k})^*_\rho
\mathbf{1}_{E_{k_0}^{\complement}}\right\}^t\right\|_{L^p(X)}
\\
&\lesssim 2^{-tk_0}\left\|\sum_{k=k_0}^{\infty}\sum_{i\in\mathbb{N}}
\left[|\lambda_{i,k}|(a_{i,k})^*_\rho\right]^t
\mathbf{1}_{(2C_\rho B_{i,k})^{\complement}}\right\|_{L^p(X)}
\\
&\sim2^{-tk_0}
\left\|\left\{\sum_{k=k_0}^{\infty}
\sum_{i\in\mathbb{N}}\left[|\lambda_{i,k}|(a_{i,k})^*_\rho\right]^t
\mathbf{1}_{(2C_\rho B_{i,k})^{\complement}}\right\}^\frac{1}{(1+\frac{\beta}{n})t}
\right\|_{L^{(1+\frac{\beta}{n})tp}(X)}^{(1+\frac{\beta}{n})t}
\\
&\lesssim2^{-tk_0}
\left\|\sum_{k=k_0}^{\infty}\left\{
\sum_{i\in\mathbb{N}}\left[|\lambda_{i,k}|(a_{i,k})^*_\rho\right]^t
\mathbf{1}_{(2C_\rho B_{i,k})^{\complement}}\right\}^\frac{1}{(1+\frac{\beta}{n})t}
\right\|_{L^{(1+\frac{\beta}{n})tp}(X)}^{(1+\frac{\beta}{n})t}
\\
&\sim2^{-tk_0}\left[\sum_{k=k_0}^{\infty}
\left\|\left\{\sum_{i\in\mathbb{N}}\left[
|\lambda_{i,k}|(a_{i,k})^*_\rho\right]^t
\mathbf{1}_{(2C_\rho B_{i,k})^{\complement}}
\right\}^\frac{1}{(1+\frac{\beta}{n})t}\right\|
_{L^{(1+\frac{\beta}{n})tp}(X)}\right]^{(1+\frac{\beta}{n})t},
\end{align*}
which, together with \eqref{4ee6}, \eqref{lambdaik},
Lemma \ref{Fefferman--Stein}
with $u=(1+\frac{\beta}{n})t\ge(1+\frac{\beta}{n})tp>1$
and $p$ replaced by $(1+\frac{\beta}{n})tp>1$,
and H\"older's inequality with $(1+\frac{\beta}{n})qtL>1$,
further implies that
\begin{align*}
\mathrm J_3
&\lesssim2^{-tk_0}\left[\sum_{k=k_0}^{\infty}2^{\frac{k}{1+\frac{\beta}{n}}}
\left\|\left\{\sum_{i\in\mathbb{N}}
\left[\mathcal{M}_\rho(\mathbf{1}_{B_{i,k}})\right]^
{(1+\frac{\beta}{n})t}\right\}^\frac{1}{(1+\frac{\beta}{n})t}\right\|
_{L^{(1+\frac{\beta}{n})tp}(X)}\right]^{(1+\frac{\beta}{n})t}
\\
&\lesssim 2^{-tk_0}\left[\sum_{k=k_0}^{\infty}2^\frac{k}{1+\frac{\beta}{n}}
\left\|\sum_{i\in\mathbb{N}}\mathbf{1}_{B_{i,k}}\right\|
^\frac{1}{(1+\frac{\beta}{n})t}_{L^p(X)}\right]^{(1+\frac{\beta}{n})t}\\
&\leq 2^{-tk_0}\left[\sum_{k=k_0}^{\infty}2^\frac{k}{(1+\frac{\beta}{n})L}
\left\|\sum_{i\in\mathbb{N}}\mathbf{1}_{B_{i,k}}\right\|
^\frac{1}{(1+\frac{\beta}{n})tL}_{L^p(X)}\right]^{(1+\frac{\beta}{n})tL}
\\
&\lesssim 2^{-tk_0}2^{(t-\delta)k_0}\left[\sum_{k=k_0}^{\infty}2^{k\delta q}
\left\|\sum_{i\in\mathbb{N}}\mathbf{1}_{B_{i,k}}\right\|^q_{L^p(X)}\right]^\frac{1}{q}
\\
&\lesssim 2^{-\delta k_0}\left\{\sum_{k=k_0}^{\infty}2^{k\delta q}
\left[\sum_{i\in\mathbb{N}}\mu(B_{i,k})\right]^\frac{q}{p}\right\}^\frac{1}{q}
\sim 2^{-k_0\delta}
\left[\sum_{k=k_0}^{\infty}2^{kq(\delta -1)}
\left(\sum_{i\in\mathbb{N}}\lambda_{i,k}^p\right)^\frac{q}{p}\right]^\frac{1}{q},
\end{align*}
where the implicit positive constants are independent of $k_0$.
Applying this and an argument similar to
that used in the estimation of \eqref{4ee5},
we obtain
\begin{equation}
\label{4ee9}
\left[\sum_{k_0\in\mathbb{Z}}2^{k_0q}(\mathrm J_{3})^{q}\right]^\frac{1}{q}
\lesssim\left[\sum_{k\in\mathbb{Z}}\left(\sum_{i\in\mathbb{N}}\lambda_{i,k}^p\right)
^{\frac{q}{p}}\right]^{\frac{1}{q}}.
\end{equation}

From \eqref{J1}, \eqref{4ee8}, \eqref{4ee9},
\eqref{730}, and \eqref{1059},
we deduce that
\begin{align}\label{1608}
\sum_{k\in\mathbb{Z}}\sum_{i\in\mathbb{N}}
\lambda_{i,k}\left|\langle a_{i,k},\phi\rangle\right|
&\lesssim\frac{\|\sum_{k\in\mathbb{Z}}\sum_{i\in\mathbb{N}}
\lambda_{i,k}(a_{i,k})^*_\rho\|_{L^{p,q}(X)}}
{\|\mathbf{1}_{B_\rho(x_1,1)}\|_{L^{p,q}(X)}}\nonumber\\
&\lesssim\frac{[\sum_{k\in\mathbb{Z}}(\sum_{i\in\mathbb{N}}\lambda_{i,k}^p)
^{\frac{q}{p}}]^{\frac{1}{q}}}
{\|\mathbf{1}_{B_\rho(x_1,1)}\|_{L^{p,q}(X)}},
\end{align}
which further implies that
the series
$\sum_{k\in\mathbb{Z}}\sum_{i\in\mathbb{N}}
\lambda_{i,k}a_{i,k}$ converges
in $({\mathcal{G}}^\varepsilon_0(\beta,\gamma))'$.
This finishes the proof of (i).

Now, we show (ii).
Let $f:=\sum_{k\in\mathbb{Z}}\sum_{i\in\mathbb{N}}
\lambda_{i,k}a_{i,k}$ in $({\mathcal{G}}^\varepsilon_0(\beta,\gamma))'$.
Then, by an argument similar to that used
in the estimation of \eqref{16177}, we conclude that,
for any $x\in X$,
$$
f^*_\rho(x)\lesssim\sum_{k\in\mathbb{N}}\sum_{i\in\mathbb{N}}
\lambda_{i,k}(a_{i,k})^*_\rho(x).
$$
From this and \eqref{1608},
we infer that $f\in H^{p,q}_*(X)$ and
$$
\|f\|_{H^{p,q}_*(X)}
=\left\|f^*_\rho\right\|_{L^{p,q}(X)}
\lesssim\left[\sum_{k\in\mathbb{Z}}
\left(\sum_{i\in\mathbb{N}}\lambda_{i,k}^p\right)
^{\frac{q}{p}}\right]^{\frac{1}{q}},
$$
which completes the proof (ii).

Finally, we prove (iii).
Let
\begin{align*}
f:=\sum_{k\in\mathbb{Z}}\sum_{i\in\mathbb{N}}
\lambda_{i,k}a_{i,k}
\end{align*}
in $({\mathcal{G}}^\varepsilon_0(\beta,\gamma))'$
and, for any $N\in\mathbb{N}$ and $x\in X$, let
\begin{align*}
f_N(x):=\sum_{k=-N}^N\sum_{i=1}^N
\lambda_{i,k}a_{i,k}(x).
\end{align*}
Then, by the proof of (ii) and \eqref{4ee1}, we find that
$f,f_N\in H^{p,q}_*(X)$ for any $N\in\mathbb{N}$ and
\begin{align*}
\left\|f-f_N\right\|_{H^{p,q}_*(X)}
&=\left\|(f-f_N)_\rho^*\right\|_{L^{p,q}(X)}\\
&\lesssim\left[\sum_{|k|\leq N}\left(\sum_{i=N+1}^\infty\lambda_{i,k}^p
\right)^{\frac{q}{p}}+\sum_{|k|>N}\left(\sum_{i\in\mathbb{N}}^\infty
\lambda_{i,k}^p\right)^{\frac{q}{p}}\right]^{\frac{1}{q}}\to0
\end{align*}
as $N\to\infty$.
This finishes the proof of (iii) and hence Theorem~\ref{6.15}.
\end{proof}

To show Theorem~\ref{6.16}, we need several
technical lemmas. The following lemma is a Whitney decomposition theorem of
open sets  proven in  \cite[Theorem~6.2]{AlMiMi13}
(see also \cite[Theorem~2.4]{AlMi2015})
which is a slight generalization of \cite[Theorem~4.21]{MiMiMiMo13}.

\begin{lemma}
\label{Wtd}
Let $(X,\mathbf{q})$ be a geometrically doubling
quasi-ultrametric space and fix $\rho\in\mathbf{q}$.
Then, for any $\lambda\in(1,\infty)$ and
$\Lambda\in(4C_\rho^2\widetilde{C}_\rho\lambda,\infty)$,
with $C_\rho$ and $\widetilde{C}_\rho$ as, respectively, in \eqref{Crho} and \eqref{Crhotilde},
there exists a constant $L_0\in\mathbb{N}$, depending only on $C_\rho$,
$\widetilde{C}_\rho$,
$\lambda$, and the geometric doubling constant of $(X,\mathbf{q})$,
such that, for any proper, nonempty, open subset $\Omega$ of
the topological space $(X,\tau_{\mathbf{q}})$,
there exist a sequence $\{x_k\}_{k\in\mathbb{N}}$ in $\Omega$ of points
and a sequence $\{r_k\}_{k\in\mathbb{N}}$ in $(0,\infty)$
having the following properties:
\begin{enumerate}
\item[\textup{(i)}]
$\Omega=\bigcup_{k\in\mathbb{N}}B_\rho(x_k,r_k)$.

\item[\textup{(ii)}]
For any $x\in\Omega$,
$$
\sum_{k\in\mathbb{N}}\mathbf{1}_{B_\rho(x_k,\lambda r_k)}(x)
\leq L_0.
$$
Indeed, there exists $\delta\in(0,1)$,
which depends only on $C_\rho$, $\widetilde{C}_\rho$,
$\lambda$, and the geometric doubling constant of $(X,\mathbf{q})$, such that,
for any $x_0\in\Omega$,
\begin{align*}
\#\left\{k\in\mathbb{N}:
B_\rho\left(x_0,\delta\mathrm{\,dist}_\rho\left(x_0,\Omega^\complement\right)\right)
\cap B_\rho(x_k,\lambda r_k)\neq\varnothing\right\}\leq L_0,
\end{align*}
where $\#E$\index[symbols]{E@$\#E$} denotes the \emph{cardinality} of a set $E$.

\item[\textup{(iii)}]
For any $k\in\mathbb{N}$,
$B_\rho(x_k,\lambda r_k)\subset\Omega$ and there exists
$y_k\in B_\rho(x_k,\Lambda r_k)\cap\Omega^\complement$;
\item[\textup{(iv)}]
for any $i,j\in\mathbb{N}$ satisfying
$B_\rho(x_i,\lambda r_i)\cap B_\rho(x_j,\lambda r_j)\neq\varnothing$,
$r_i\sim r_j$,
where the positive equivalence constants
are independent of both $i$ and $j$.
\end{enumerate}
\end{lemma}

\begin{remark}
In Lemma~\ref{Wtd}, the flexibility in choosing $\Lambda$
is implicit in the construction
in the proofs of \cite[Theorem~6.2]{AlMiMi13} and
\cite[Theorem 4.21]{MiMiMiMo13}.
\end{remark}

The following lemma shows the existence of a partition of unity
on geometrically doubling quasi-ultrametric spaces,
which is exactly \cite[Theorem 2.5]{AlMi2015}
and also a direct result of \cite[Theorem 5.1]{AlMiMi13}.

\begin{lemma}
\label{pou}
Let $(X,\mathbf{q})$ be a geometrically doubling quasi-ultrametric space and fix
$\rho\in\mathbf{q}$.
Assume that $\Omega$ is a proper nonempty subset of $X$.
For any given $\lambda\in(C_\rho^2,\infty)$
and $\lambda'\in(C_\rho,\frac{\lambda}{C_\rho})$,
consider the decomposition of $\Omega$ into the family
$\{B_\rho(x_k,r_k)\}_{k\in\mathbb{N}}$ as given in Lemma~\ref{Wtd}.
Then, for any $\varepsilon\in(0,(\log_2C_\rho)^{-1}]\cap(0,\infty)$,
there exist a constant $C\in[1,\infty)$, depending only on
$\rho$, $\varepsilon$, $L_0$, and the proportionality constants in
Lemma~\ref{Wtd}(iv), and a family
of real-valued functions $\{\phi_k\}_{k\in\mathbb{N}}$
defined on $X$ such that
\begin{enumerate}
\item[\textup{(i)}]
For any $k\in\mathbb{N}$ and $\beta\in(0,\varepsilon]$,
$\phi_k\in\dot{C}^\beta(X,\mathbf{q})$ and
$\|\phi_k\|_{\dot{C}^\beta(X,\rho)}
\leq Cr_k^{-\beta}$.

\item[\textup{(ii)}]
For any $k\in\mathbb{N}$,
$0\leq\phi_k\leq1$ on $X$,
$\phi_k\equiv0$ on $X\setminus B_\rho(x_k,\lambda'r_k)$,
and $\phi_k\ge C^{-1}$ on $B_\rho(x_k,r_k)$.

\item[\textup{(iii)}]
$\sum_{k\in\mathbb{N}}\phi_k
=\mathbf{1}_{\bigcup_{k\in\mathbb{N}}B_\rho(x_k,r_k)}
=\mathbf{1}_{\bigcup_{k\in\mathbb{N}}B_\rho(x_k,\lambda'r_k)}
=\mathbf{1}_{\bigcup_{k\in\mathbb{N}}B_\rho(x_k,\lambda r_k)}
=\mathbf{1}_{\Omega}$.
\end{enumerate}
\end{lemma}

\begin{remark}
\label{lemma5.14iv}
From the proof of Lemma~\ref{CGL} and both (i) and (ii) of Lemma~\ref{pou},
we infer that there exists a positive constant $C$,
independent of both $x_k$ and $r_k$, such that,
for any $k\in\mathbb{N}$, $\beta\in(0,\varepsilon]$,
and $\gamma\in(0,\infty)$,
$\phi_k\in{\mathcal{G}}(x_k,r_k,\beta,\gamma)$ and
$\|\phi_k\|_{{\mathcal{G}}(x_k,r_k,\beta,\gamma)}
\leq C(V_\rho)_{r_k}(x_k)$.
\end{remark}

For a homogeneous approximation of the identity
$\{\mathcal{S}_{2^{-k}}\}_{k\in\mathbb{Z}}$ as in
Definition~\ref{DefHATI}, we have the following lemmas.

\begin{lemma}\label{4l4.1}
Let $(X,\mathbf{q},\mu)$ be a quasi-ultrametric space
of homogeneous type and $0<\beta,\gamma<\varepsilon
\preceq\mathrm{ind\,}(X,\mathbf{q})$.
Let $\{\mathcal{S}_{2^{-k}}\}_{k\in\mathbb{Z}}$
be a homogeneous approximation
of the identity of order $\varepsilon$ as in Definition~\ref{DefHATI}
and $f\in {\mathcal{G}}(x_0,r_0,\beta,\gamma)$ for some
$x_0\in X$ and $r_0\in (0,\infty)$.
Then there exists a positive constant $C$, independent
of $x_0,\ r_0$, and $f$, such that
\begin{enumerate}
\item[\textup{(i)}]
For any $k\in\mathbb{Z}$ satisfying $r_0\geq2^{-k}$,
\begin{align*}
\left\|\mathcal{S}_{2^{-k}}f\right\|_{{\mathcal{G}}(x_0,r_0,\beta,\gamma)}
\leq C\|f\|_{{\mathcal{G}}(x_0,r_0,\beta,\gamma)}.
\end{align*}

\item[\textup{(ii)}]
For any $k\in\mathbb{Z}$ satisfying $r_0<2^{-k}$,
\begin{align*}
\left\|\mathcal{S}_{2^{-k}}f\right\|_{{\mathcal{G}}(x_0,2^{-k},\beta,\gamma)}
\leq C\|f\|_{{\mathcal{G}}(x_0,r_0,\beta,\gamma)}.
\end{align*}
\end{enumerate}
\end{lemma}

\begin{proof}
(i) follows directly from Theorem~\ref{CZOonCG}.
A similar argument can be found in the estimation of \eqref{2247}.

Next, we prove (ii). Let $\rho\in\mathbf{q}$.
Without loss of generality,
we may assume that $\rho$ is a regularized quasi-ultrametric
and $f\in{\mathcal{G}}(x_0,r_0,\beta,\gamma)$
with $\|f\|_{{\mathcal{G}}(x_0,r_0,\beta,\gamma)}\leq1$.
To show $\|\mathcal{S}_{2^{-k}}f\|_{{\mathcal{G}}(x_0,2^{-k},\beta,\gamma)}
\lesssim1$, we first estimate the size condition of
$\mathcal{S}_{2^{-k}}f$.
To this end, let $x\in X$. We
consider the following two cases for $\rho(x_0,x)$.
Let $C'$ be the same positive constant
as in Definition~\ref{DefATI}(i).

\emph{Case (1)} $\rho(x_0,x)\leq C'C_\rho2^{-k}$.
In this case, by Remark~\ref{Rm}(i),
the size condition of $f$,
and (iii) and (ii) of Lemma~\ref{A3.2},
we find that
\begin{align*}
\left|\mathcal{S}_{2^{-k}}f(x)\right|
&\lesssim\int_{B_\rho(x,C'2^{-k})}\frac{1}{(V_\rho)_{2^{-k}}(x)+V_\rho(x,y)}
\left[\frac{2^{-k}}{2^{-k}+\rho(x,y)}\right]^\gamma\\
&\quad\times
\frac{1}{(V_\rho)_{r_0}(x_0)+V_\rho(x_0,y)}
\left[\frac{r_0}{r_0+\rho(x_0,y)}\right]^\gamma\,d\mu(y)\\
&\sim\frac{1}{(V_\rho)_{2^{-k}}(x)+V_\rho(x,x_0)}
\left[\frac{2^{-k}}{2^{-k}+\rho(x,x_0)}\right]^\gamma\\
&\quad\times\int_{B_\rho(x,C'2^{-k})}
\frac{1}{(V_\rho)_{r_0}(x_0)+V_\rho(x_0,y)}
\left[\frac{r_0}{r_0+\rho(x_0,y)}\right]^\gamma\,d\mu(y)\\
&\lesssim\frac{1}{(V_\rho)_{2^{-k}}(x)+V_\rho(x,x_0)}
\left[\frac{2^{-k}}{2^{-k}+\rho(x,x_0)}\right]^\gamma.
\end{align*}

\emph{Case (2)} $\rho(x_0,x)>C'C_\rho2^{-k}$.
In this case, from Definition~\ref{D2.1}(iii),
we deduce that, for any $y\in B_\rho(x,C'2^{-k})$,
$$
\rho(x_0,y)>C'2^{-k}
$$
and hence, using Definition~\ref{D2.1}(iii) again, we obtain
$$
\rho(x_0,y)\sim\rho(x_0,x).
$$
By this,
Remark~\ref{Rm}(i),
the size condition of $f$,
and the assumption that $r_0<2^{-k}$,
we conclude that
\begin{align*}
\left|\mathcal{S}_{2^{-k}}f(x)\right|
&\lesssim\int_{B_\rho(x,C'2^{-k})}\frac{1}{(V_\rho)_{2^{-k}}(x)
+V_\rho(x,y)}\left[\frac{2^{-k}}{2^{-k}+\rho(x,y)}\right]^\gamma\\
&\quad\times
\frac{1}{(V_\rho)_{r_0}(x_0)+V_\rho(x_0,y)}
\left[\frac{r_0}{r_0+\rho(x_0,y)}\right]^\gamma\,d\mu(y)\\
&\lesssim\int_{B_\rho(x,C'2^{-k})}\frac{1}{(V_\rho)_{2^{-k}}(x)
+V_\rho(x,y)}\left[\frac{2^{-k}}{2^{-k}+\rho(x,y)}\right]^\gamma\,d\mu(y)\\
&\quad\times
\frac{1}{(V_\rho)_{r_0}(x_0)+V_\rho(x_0,x)}
\left[\frac{2^{-k}}{r_0+\rho(x_0,x)}\right]^\gamma,
\end{align*}
which, combined with both (ii) and (iii) of Lemma~\ref{A3.2},
further implies that
\begin{align*}
\left|\mathcal{S}_{2^{-k}}f(x)\right|
&\lesssim\frac{1}{(V_\rho)_{r_0}(x_0)+V_\rho(x_0,x)}
\left[\frac{2^{-k}}{r_0+\rho(x_0,x)}\right]^\gamma\\
&\sim\frac{1}{(V_\rho)_{2^{-k}}(x_0)+V_\rho(x_0,x)}
\left[\frac{2^{-k}}{2^{-k}+\rho(x_0,x)}\right]^\gamma.
\end{align*}
This is the desired estimate.

Now, we estimate the regularity condition of
$\mathcal{S}_{2^{-k}}f$.
To this end, let $x,x'\in X$ satisfy
$\rho(x,x')\leq\frac{2^{-k}+\rho(x_0,x)}{2C_\rho}$.
We consider the following two cases.

\emph{Case (1)} $\rho(x_0,x)\leq2^{-k}$.
In this case, we have
$$
\rho(x,x')\leq\frac{2^{-k}}{C_\rho}.
$$
From this, Definition~\ref{DefATI}(ii), Lemma~\ref{A3.2}(vii),
$\beta\in(0,\varepsilon]$, the size condition of $f$,
and both (ii) and (vi) of Lemma~\ref{A3.2},
we infer that
\begin{align*}
&\left|\mathcal{S}_{2^{-k}}f(x)-\mathcal{S}_{2^{-k}}f(x')\right|\\
&\quad\lesssim\int_X\left[\frac{\rho(x,x')}{2^{-k}}
\right]^\varepsilon\left[\frac{1}{(V_\rho)_{2^{-k}}(x)}
+\frac{1}{(V_\rho)_{2^{-k}}(x')}\right]|f(y)|\,d\mu(y)\\
&\quad\sim\left[\frac{\rho(x,x')}{2^{-k}+\rho(x_0,x)}
\right]^\varepsilon\frac{1}{(V_\rho)_{2^{-k}}(x)}
\|f\|_{L^1(X)}\\
&\quad\lesssim\left[\frac{\rho(x,x')}{2^{-k}+\rho(x_0,x)}
\right]^\beta\frac{1}{(V_\rho)_{2^{-k}}(x)}\\
&\quad\sim\left[\frac{\rho(x,x')}{2^{-k}+\rho(x_0,x)}
\right]^\beta\frac{1}{(V_\rho)_{2^{-k}}(x)+V_\rho(x,x_0)}
\left[\frac{2^{-k}}{2^{-k}+\rho(x,x_0)}\right]^\gamma.
\end{align*}

\emph{Case (2)} $\rho(x_0,x)>2^{-k}$
and $\frac{2^{-k}+\rho(x_0,x)}{8C_\rho^2}\leq
\rho(x,x')\leq\frac{2^{-k}+\rho(x_0,x)}{2C_\rho}$.
In this case, by Definition~\ref{D2.1}(iii),
we find that
\begin{align*}
\rho(x,x_0)
&\leq C_\rho\max\left\{\rho(x,x'),\,\rho(x',x_0)\right\}
\\
&\leq\max\left\{\frac{2^{-k}+\rho(x_0,x)}{2},\,
C_\rho\rho(x',x_0)\right\}=C_\rho\rho(x',x_0).
\end{align*}
Similarly, we can prove that $\rho(x',x_0)\leq C_\rho\rho(x_0,x)$.
From this, the size estimate of $\mathcal{S}_{2^{-k}}f$,
Lemma~\ref{A3.2}(iii),
and the assumption that $\frac{2^{-k}+\rho(x_0,x)}{8C_\rho^2}\leq\rho(x,x')$,
we deduce that
\begin{align*}
&\left|\mathcal{S}_{2^{-k}}f(x)-\mathcal{S}_{2^{-k}}f(x')
\right|\\
&\quad\leq\left|\mathcal{S}_{2^{-k}}f(x)\right|
+\left|\mathcal{S}_{2^{-k}}f(x')\right|\\
&\quad\lesssim\frac{1}{(V_\rho)_{2^{-k}}(x)+V_\rho(x,x_0)}
\left[\frac{2^{-k}}{2^{-k}+\rho(x,x_0)}\right]^\gamma\\
&\qquad+
\frac{1}{(V_\rho)_{2^{-k}}(x')+V_\rho(x',x_0)}
\left[\frac{2^{-k}}{2^{-k}+\rho(x',x_0)}\right]^\gamma\\
&\quad\sim\left[\frac{\rho(x,x')}{2^{-k}+\rho(x_0,x)}
\right]^\beta\frac{1}{(V_\rho)_{2^{-k}}(x)+V_\rho(x,x_0)}
\left[\frac{2^{-k}}{2^{-k}+\rho(x,x_0)}\right]^\gamma.
\end{align*}

\emph{Case (3)} $\rho(x_0,x)>2^{-k}$
and $\rho(x,x')<\frac{2^{-k}+\rho(x_0,x)}{8C_\rho^2}$.
In this case,
\begin{align*}
\rho(x,x')<\frac{\rho(x_0,x)}{(2C_\rho)^2}
<\frac{r_0+\rho(x_0,x)}{(2C_\rho)^2}.
\end{align*}
By Definition~\ref{DefATI}(v), we conclude that
\begin{align*}
\left|\mathcal{S}_{2^{-k}}f(x)-\mathcal{S}_{2^{-k}}f(x')\right|
&=\left|\int_X\left[S_{2^{-k}}(x,y)-S_{2^{-k}}(x',y)\right]
\left[f(y)-f(x)\right]\,d\mu(y)\right|\\
&\leq\int_X\left|S_{2^{-k}}(x,y)-S_{2^{-k}}(x',y)\right|
\left|f(y)-f(x)\right|\,d\mu(y)\\
&=\sum_{i=1}^3\int_{D_i}\left|S_{2^{-k}}(x,y)-S_{2^{-k}}(x',y)\right|
\left|f(y)-f(x)\right|\,d\mu(y)\\
&=:\mathrm{I}+\mathrm{II}+\mathrm{III},
\end{align*}
where
$$
D_1:=\left\{y\in X:
\rho(x,x')\leq\frac{2^{-k}+\rho(x,y)}{2C_\rho}\leq
\frac{r_0+\rho(x_0,x)}{(2C_\rho)^2}\right\},
$$
$$
D_2:=\left\{y\in X:
\rho(x,x')\leq\frac{r_0+\rho(x_0,x)}{(2C_\rho)^2}
\leq\frac{2^{-k}+\rho(x,y)}{2C_\rho}\right\},
$$
and
$$
D_3:=\left\{y\in X:
\rho(x,x')>\frac{2^{-k}+\rho(x,y)}{2C_\rho}\right\}.
$$
To estimate $\mathrm{I}$,
note that $y\in D_1$ implies that
$\rho(x,y)\leq\frac{r_0+\rho(x_0,x)}{2C_\rho}$.
From this, the regularity
conditions of both $S_{2^{-k}}(\cdot,y)$ and $f$
[see Lemma~\ref{ATItest} and Definition~\ref{testfunction}(ii)],
$\beta\in(0,\varepsilon]$,
$\rho(x_0,x)>2^{-k}>r_0$, and Lemma~\ref{A3.2}(ii),
we infer that
\begin{align*}
\mathrm{I}&\lesssim\int_{D_1}
\left[\frac{\rho(x,x')}{2^{-k}+\rho(y,x)}\right]^\varepsilon
\frac{1}{(V_\rho)_{2^{-k}}(y)+V_\rho(y,x)}
\left[\frac{2^{-k}}{2^{-k}+\rho(y,x)}\right]^\gamma\\
&\quad\times\left[\frac{\rho(x,y)}{r_0+\rho(x_0,x)}\right]^\beta
\frac{1}{(V_\rho)_{r_0}(x_0)+V_\rho(x_0,x)}
\left[\frac{r_0}{r_0+\rho(x_0,x)}\right]^\gamma\,d\mu(y)\\
&\lesssim\left[\frac{\rho(x,x')}{r_0+\rho(x_0,x)}\right]^\beta
\frac{1}{(V_\rho)_{r_0}(x_0)+V_\rho(x_0,x)}
\left[\frac{r_0}{r_0+\rho(x_0,x)}\right]^\gamma\\
&\quad\times\int_{D_1}
\left[\frac{\rho(x,y)}{2^{-k}+\rho(y,x)}\right]^\beta
\frac{1}{(V_\rho)_{2^{-k}}(y)+V_\rho(y,x)}
\left[\frac{2^{-k}}{2^{-k}+\rho(y,x)}\right]^\gamma
\,d\mu(y)\\
&\lesssim\left[\frac{\rho(x,x')}{2^{-k}+\rho(x_0,x)}
\right]^\beta\frac{1}{(V_\rho)_{2^{-k}}(x)+V_\rho(x,x_0)}
\left[\frac{2^{-k}}{2^{-k}+\rho(x,x_0)}\right]^\gamma,
\end{align*}
where, in the last inequality, we used (ii), (vi),
and (vii) of Lemma~\ref{A3.2} to estimate
\begin{align*}
&\int_{D_1}
\left[\frac{\rho(x,y)}{2^{-k}+\rho(y,x)}\right]^\beta
\frac{1}{(V_\rho)_{2^{-k}}(y)+V_\rho(y,x)}
\left[\frac{2^{-k}}{2^{-k}+\rho(y,x)}\right]^\gamma\,d\mu(y)
\\
&\quad\leq\int_{X}
\frac{1}{(V_\rho)_{2^{-k}}(y)+V_\rho(y,x)}
\left[\frac{2^{-k}}{2^{-k}+\rho(y,x)}\right]^\gamma\,d\mu(y)
\\
&\quad\sim\int_{X}
\frac{1}{(V_\rho)_{2^{-k}}(x)+V_\rho(x,y)}
\left[\frac{2^{-k}}{2^{-k}+\rho(x,y)}\right]^\gamma\,d\mu(y)
\lesssim1.
\end{align*}
To estimate $\mathrm{II}$, using the size condition of $f$,
Lemma~\ref{ATItest}, and Definition~\ref{testfunction}(ii) for $S_{2^{-k}}(\cdot,y)$,
we find that
\begin{align*}
\mathrm{II}
&\leq|f(x)|\int_{D_2}
\left|S_{2^{-k}}(x,y)-S_{2^{-k}}(x',y)\right|\,d\mu(y)
\\
&\quad+\int_{D_2}
\left|S_{2^{-k}}(x,y)-S_{2^{-k}}(x',y)\right|
|f(y)|\,d\mu(y)
\\
&\lesssim\frac{1}{(V_\rho)_{r_0}(x_0)
+V_\rho(x_0,x)}
\left[\frac{r_0}{r_0+\rho(x_0,x)}\right]^\gamma
\\
&\quad\times\int_{D_2}
\left[\frac{\rho(x,x')}{2^{-k}+\rho(y,x)}\right]^\varepsilon
\frac{1}{(V_\rho)_{2^{-k}}(y)
+V_\rho(y,x)}
\left[\frac{2^{-k}}{2^{-k}+\rho(y,x)}\right]^\gamma\,d\mu(y)
\\
&\quad+\int_{D_2}\frac{1}{(V_\rho)_{r_0}(x_0)
+V_\rho(x_0,y)}
\left[\frac{r_0}{r_0+\rho(x_0,y)}\right]^\gamma
\\
&\quad\times
\left[\frac{\rho(x,x')}{2^{-k}+\rho(y,x)}\right]^\varepsilon
\frac{1}{(V_\rho)_{2^{-k}}(y)
+V_\rho(y,x)}
\left[\frac{2^{-k}}{2^{-k}+\rho(y,x)}\right]^\gamma\,d\mu(y),
\end{align*}
which, together with $\beta\in(0,\varepsilon]$,
the fact that $$
2^{-k}+\rho(x,y)\gtrsim2^{-k}+\rho(x_0,x)>2^{-k}>r_0,
$$
and (ii), (vi), and (vii) of Lemma~\ref{A3.2},
further implies that
\begin{align*}
\mathrm{II}	
&\lesssim\left[\frac{\rho(x,x')}{2^{-k}+\rho(x_0,x)}
\right]^\beta\frac{1}{(V_\rho)_{2^{-k}}(x_0)
+V_\rho(x_0,x)}
\left[\frac{2^{-k}}{2^{-k}+\rho(x_0,x)}\right]^\gamma
\\
&\quad\times\int_{X}
\frac{1}{(V_\rho)_{2^{-k}}(x)
+V_\rho(x,y)}
\left[\frac{2^{-k}}{2^{-k}+\rho(x,y)}\right]^\gamma\,d\mu(y)
\\
&\quad+\left[\frac{\rho(x,x')}{2^{-k}+\rho(x_0,x)}
\right]^\beta\frac{1}{(V_\rho)_{2^{-k}}(x_0)
+V_\rho(x_0,x)}
\left[\frac{2^{-k}}{2^{-k}+\rho(x_0,x)}\right]^\gamma
\\
&\quad\times\int_{X}
\frac{1}{(V_\rho)_{r_0}(x_0)
+V_\rho(x_0,y)}
\left[\frac{r_0}{r_0+\rho(x_0,y)}\right]^\gamma\,d\mu(y)
\\
&\lesssim\left[\frac{\rho(x,x')}{2^{-k}+
\rho(x_0,x)}\right]^\beta\frac{1}{(V_\rho)_{2^{-k}}(x_0)
+V_\rho(x_0,x)}
\left[\frac{2^{-k}}{2^{-k}+\rho(x_0,x)}\right]^\gamma.
\end{align*}
This implies the desired regularity estimate of $\mathrm{II}$.

To estimate $\mathrm{III}$, note that, for any $y\in D_3$,
we have $\frac{\rho(x,y)}{2C_\rho}<\rho(x,x')<\frac{\rho(x_0,x)}{(2C_\rho)^2}$,
which further implies that $\rho(x,y)<\frac{r_0+\rho(x_0,x)}{2C_\rho}$.
From this, Definition~\ref{DefATI}(i),
the regularity condition of $f$,
$\rho(x_0,x)>2^{-k}>r_0$,
and Definition~\ref{DefATI}(v),
we deduce that
\begin{align*}
\mathrm{III}&\lesssim\int_X
\left[S_{2^{-k}}(x,y)
+S_{2^{-k}}(x',y)\right]\\
&\quad\times\left[\frac{\rho(x,y)}{r_0+\rho(x_0,x)}\right]^\beta
\frac{1}{(V_\rho)_{r_0}(x_0)+V_\rho(x_0,x)}
\left[\frac{r_0}{r_0+\rho(x_0,x)}\right]^\gamma\,d\mu(y)\\
&\lesssim\left[\frac{\rho(x,x')}{2^{-k}+\rho(x_0,x)}\right]^\beta
\frac{1}{(V_\rho)_{2^{-k}}(x_0)+V_\rho(x_0,x)}
\left[\frac{2^{-k}}{2^{-k}+\rho(x_0,x)}\right]^\gamma\\
&\quad\times\int_X\left[S_{2^{-k}}(x,y)
+S_{2^{-k}}(x',y)\right]\,d\mu(y)\\
&=2\left[\frac{\rho(x,x')}{2^{-k}+\rho(x_0,x)}\right]^\beta
\frac{1}{(V_\rho)_{2^{-k}}(x_0)+V_\rho(x_0,x)}
\left[\frac{2^{-k}}{2^{-k}+\rho(x_0,x)}\right]^\gamma.
\end{align*}

To summarize, we obtain the desired regularity estimate
for $\mathcal{S}_{2^{-k}}f$ and hence we conclude that
$$
\|\mathcal{S}_{2^{-k}}f\|_{{\mathcal{G}}(x_0,2^{-k},\beta,\gamma)}
\lesssim1.
$$
This finishes the proof of Lemma~\ref{4l4.1}.
\end{proof}

The following lemma when $(X,\mathbf{q},\mu)$
is an RD-space reduces to \cite[Lemma~4.10]{LLD}.

\begin{lemma}\label{2051}
Let $(X,\mathbf{q},\mu)$ be a quasi-ultrametric space of
homogeneous type with upper dimension $n\in(0,\infty)$ and $\rho\in\mathbf{q}$
a regularized quasi-ultrametric.
Let $\beta\in(0,\infty)$,
$q\in(\frac{n}{n+\beta},\infty)$,
and $M\in\mathbb{N}$.
Then there exists a positive constant $C$
such that, for any given sequences $\{x_k\}_{k\in\mathbb{N}}$ in $X$
and $\{r_k\}_{k\in\mathbb{N}}$ in $(0,\infty)$
satisfying $\sum_{k\in\mathbb{N}}\mathbf{1}_{B_\rho(x_k,r_k)}\leq M$ pointwise in $X$,
\begin{align*}
&\int_X\left\{\sum_{k\in\mathbb{N}}
\frac{(V_\rho)_{r_k}(x_k)}{(V_\rho)_{r_k}(x_k)+V_\rho(x_k,x)}
\left[\frac{r_k}{r_k+\rho(x_k,x)}\right]^\beta\right\}^q\,d\mu(x)\\
&\quad\leq C\mu\left(\bigcup_{k\in\mathbb{N}}B_\rho(x_k,r_k)\right).
\end{align*}
\end{lemma}

\begin{proof}
From \eqref{upperdoub},
Lemma~\ref{A3.2}(vi), and \eqref{HLmaxequiv},
we infer that, for any $k\in\mathbb{N}$ and $x\in X$,
\begin{align*}
\frac{r_k}{r_k+\rho(x_k,x)}
&\lesssim\left[\frac{(V_\rho)_{r_k}(x_k)}
{(V_\rho)_{r_k+\rho(x_k,x)}(x_k)}\right]^\frac{1}{n}
\leq\left[\widetilde{\mathcal{M}}_\rho
\left(\mathbf{1}_{B_\rho(x_k,r_k)}\right)(x)\right]^\frac{1}{n}\\
&\sim\left[\mathcal{M}_{\rho}
\left(\mathbf{1}_{B_\rho(x_k,r_k)}\right)(x)\right]^\frac{1}{n}
\end{align*}
and
\begin{align*}
\frac{(V_\rho)_{r_k}(x_k)}{(V_\rho)_{r_k}(x_k)+V_\rho(x_k,x)}
&\sim\frac{(V_\rho)_{r_k}(x_k)}{(V_\rho)_{r_k+\rho(x_k,x)}(x_k)}
\leq\widetilde{\mathcal{M}}_\rho\left(\mathbf{1}_{B_\rho(x_k,r_k)}\right)(x)\\
&\sim\mathcal{M}_{\rho}\left(\mathbf{1}_{B_\rho(x_k,r_k)}\right)(x),
\end{align*}
which, combined with $q>\frac{n}{n+\beta}$,
Lemma~\ref{Fefferman--Stein}
with $u=\frac{n+\beta}{n}$ and $p=qu$,
and the finite intersection property of the $\rho$-balls
$\{B_\rho(x_k,r_k)\}_{k\in\mathbb{N}}$,
further implies that
\begin{align*}
&\int_X\left\{\sum_{k\in\mathbb{N}}
\frac{(V_\rho)_{r_k}(x_k)}{(V_\rho)_{r_k}(x_k)+V_\rho(x_k,x)}
\left[\frac{r_k}{r_k+\rho(x_k,x)}\right]^\beta\right\}^q\,d\mu(x)
\\
&\quad\lesssim\int_X\left\{\sum_{k\in\mathbb{N}}
\left[\mathcal{M}_\rho\left(\mathbf{1}_{B_\rho(x_k,r_k)}\right)(x)
\right]^\frac{n+\beta}{n}\right\}^q\,d\mu(x)
\\
&\quad\lesssim\int_X\left[\sum_{k\in\mathbb{N}}
\mathbf{1}_{B_\rho(x_k,r_k)}(x)\right]^q\,d\mu(x)
\\
&\quad\lesssim\int_X\sum_{k\in\mathbb{N}}
\mathbf{1}_{B_\rho(x_k,r_k)}(x)\,d\mu(x)
\lesssim\mu\left(\bigcup_{k\in\mathbb{N}}B_\rho(x_k,r_k)\right).
\end{align*}
This finishes the proof of Lemma~\ref{2051}.
\end{proof}

To show Theorem~\ref{6.16} for the case $\mu(X)<\infty$,
we also need a local atom characterization of $H^{p,q,\infty}_{\rm at}(X)$
in the case $\mu(X)<\infty$
(see Theorem~\ref{1314}).
We begin with
recalling the following concept of local atoms.

\begin{definition}\label{localatom}
Let $(X,\mathbf{q},\mu)$ be a
quasi-ultrametric space of homogeneous type. Assume that
$p\in(0,1]$,
$q\in(p,\infty]\cap[1,\infty]$,
and $\rho\in\mathbf{q}$ is such that all $\rho$-balls are $\mu$-measurable.
A $\mu$-measurable function
$a: X\to\mathbb{C}$ is called
a \emph{local $(\rho,p,q)$-atom}\index[words]{local $(\rho,p,q)$-atom}
if there exist $x\in X$
and $r\in(0,\infty)$ such that
\begin{enumerate}
\item[\rm(i)]
$\mathrm{\,supp\,} a:=\{x\in X:a(x)\neq0\}\subset B_{\rho}(x,r)$;

\item[\rm(ii)]
$||a||_{L^q(X)}\leq[(V_\rho)_r(x)]^{\frac{1}{q}-\frac{1}{p}}$;

\item[\rm(iii)]
$\int_Xa(x)\,d\mu(x)=0$ if $r\in(0,1]$.
\end{enumerate}
\end{definition}

\begin{remark}
\begin{enumerate}
\item[\rm(i)]
It is easy to see from Definition~\ref{localatom},
\eqref{1557-xx}, and \eqref{doublingcond}
that, if $a: X\to\mathbb{C}$ is a local $(\rho,p,q)$-atom
supported in $B_{\rho}(x,r)$ for some $x\in X$
and $r\in(0,\infty)$, then, for any other $\rho'\in\mathbf{q}$
with the property that all $\rho'$-balls are $\mu$-measurable,
there exist positive constants $C$ and $c$ (depending only on $\rho$,
$\rho'$, and $\mu$) such that $Ca$ is a
local $(\rho',p,q)$-atom supported in $B_{\rho'}(x,cr)$.
As such, we sometimes simply refer to $a$ as a
$(p,q)$-\textit{atom} or a $(p,q)$-\textit{atom modulo a harmless constant multiple}.

\item[\rm(ii)]
Note that a $(p,q)$-atom is also a local $(p,q)$-atom.
Indeed, when $\mu(X)=\infty$, it holds automatically.
When $\mu(X)<\infty$, we only need to prove that $a:=[\mu(X)]^{-\frac{1}{p}}$
is a local $(p,q)$-atom. Note that $\mu(X)<\infty$
implies $\mathrm{diam\,}(X)<\infty$ (see, for instance, \cite[Lemma~5.1]{ny1997}).
Then, for any fixed $x_1\in X$, we have
$$
\mathrm{supp\,}a=X=B_\rho(x_1,\mathrm{diam\,}(X)+1)=:B,
$$
$$
\|a\|_{L^q(X)}=[\mu(X)]^{\frac{1}{q}-\frac{1}{p}}=[\mu(B)]^{\frac{1}{q}-\frac{1}{p}},
$$
and the radius of $B$ is greater than $1$.
Thus, $a$ is a local $(p,q)$-atom.

\item[\rm(iii)]
For any local $(p,q)$-atom $a$ supported in some ball $B\subset X$,
we have
\begin{align*}
\|a\|_{L^1(X)}\leq[\mu(B)]^{1-\frac{1}{p}}
\ \ \text{and}\ \
\|a\|_{L^p(X)}\leq1.
\end{align*}
\end{enumerate}
\end{remark}

Based on the local atoms in Definition~\ref{localatom},
we introduce the following concept of
local atomic Hardy--Lorentz spaces.

\begin{definition}
Let $(X,\mathbf{q},\mu)$ be a quasi-ultrametric space
of homogeneous type,
$\rho\in\mathbf{q}$ have the property that
all $\rho$-balls are $\mu$-measurable,
$p\in(0,1]$,
$q\in(0,\infty]$, $r\in(p,\infty]\cap[1,\infty]$,
and $0<\beta,\gamma<\varepsilon<\infty$.
The $\emph{local atomic Hardy--Lorentz space}$
\index[words]{local atomic Hardy--Lorentz space}
$h^{p,q,r}_{\mathrm{at}}(X,\rho)$
\index[symbols]{H@$h^{p,q,r}_{\mathrm{at}}(X,\rho)$}
is defined to be the same as in Definition~\ref{6.3}
but only with $(p,r)$-atoms therein replaced by local $(p,r)$-atoms.
\end{definition}

To show the equivalence between $H^{p,q,r}_{\mathrm{at}}(X,\rho)$
and $h^{p,q,r}_{\mathrm{at}}(X,\rho)$,
we need the following proposition.

\begin{proposition}\label{2304}
Let $(X,\mathbf{q},\mu)$ be a quasi-ultrametric space
of homogeneous type with upper dimension $n\in(0,\infty)$,
$\rho\in\mathbf{q}$ have the property that
all $\rho$-balls are $\mu$-measurable,
and $\mu(X)<\infty$. Let
$p\in(\frac{n}{n+\mathrm{ind\,}(X,\mathbf{q})},1]$,
$q\in(0,\infty]$, $r\in(p,\infty]\cap[1,\infty]$,
and $\varepsilon,\beta,\gamma\in(0,\infty)$ be such that
$n(\frac{1}{p}-1)<\beta,\gamma<\varepsilon\preceq\mathrm{ind\,}(X,\mathbf{q})$.
Then there exists a positive constant $C$ such that, for any local $(p,r)$-atom $a$,
$$
\|a\|_{H^{p,q,r}_{\mathrm{at}}(X,\rho)}\leq C.
$$
\end{proposition}

\begin{proof}
Let $a$ be a local $(p,r)$-atom supported in $B:=B_\rho(x,r)$
for some $x\in X$ and $r\in(0,\infty)$.
If $r\in(0,1]$, then $a$ is a $(p,r)$-atom,
which further implies that $\|a\|_{H^{p,q,r}_{\mathrm{at}}(X,\rho)}\leq1$.
Assume that $r>1$. Since $\mu(X)<\infty$,
it follows that $\mathrm{diam\,}(X)<\infty$ and
$$
X=B_\rho(x,\mathrm{diam\,}(X)+1).
$$
Using this and \eqref{upperdoub}, we conclude that
\begin{align}\label{2345}
\mu(B)\leq\mu(X)\leq C_\mu\left[\frac{\mathrm{diam\,}(X)+1}{r}\right]^n\mu(B)
\leq C_\mu\left[\mathrm{diam\,}(X)+1\right]^n\mu(B).
\end{align}
We write
$$
a=\left[a-\frac{1}{\mu(X)}\int_Xa(x)\,d\mu(x)\right]
+\frac{1}{\mu(X)}\int_Xa(x)\,d\mu(x)=:a_1+a_2.
$$
For $a_2$, from H\"older's inequality, Definition~\ref{localatom}(ii),
and \eqref{2345}, we infer that
\begin{align*}
\left|a_2\right|\leq\left[\frac{1}{\mu(X)}\int_X|a(x)|^q\,d\mu(x)\right]^\frac{1}{q}
\leq\left[\mu(X)\right]^{-\frac{1}{q}}\left[\mu(B)\right]^{\frac{1}{q}-\frac{1}{p}}
\sim\left[\mu(X)\right]^{-\frac{1}{p}}.
\end{align*}
Thus, $a_2$ is a harmless constant multiple of the $(p,r)$-atom $[\mu(X)]^{-\frac{1}{p}}$.
For $a_1$, it is easy to verify that
$\mathrm{supp\,}a_1=X$, $\int_Xa_1(x)\,d\mu(x)=0$,
and
\begin{align*}
\left\|a_1\right\|_{L^q(X)}\leq\|a\|_{L^q(X)}+\|a_2\|_{L^q(X)}
\lesssim\left[\mu(B)\right]^{\frac{1}{q}-\frac{1}{p}}+
\left[\mu(X)\right]^{\frac{1}{q}-\frac{1}{p}}
\sim\left[\mu(X)\right]^{\frac{1}{q}-\frac{1}{p}}
\end{align*}
and hence $a_1$ is a harmless constant multiple of a $(p,r)$-atom
supported in $B_\rho(x,\mathrm{diam\,}(X)+1)$.
Thus,
$$
\|a\|_{H^{p,q,r}_{\mathrm{at}}(X,\rho)}\lesssim
\|a_1\|_{H^{p,q,r}_{\mathrm{at}}(X,\rho)}
+\|a_2\|_{H^{p,q,r}_{\mathrm{at}}(X,\rho)}
\lesssim1.
$$
This finishes the proof of Proposition~\ref{2304}.
\end{proof}

Using Proposition~\ref{2304}, we immediately obtain
the following conclusion.

\begin{theorem}\label{1314}
Let $(X,\mathbf{q},\mu)$ be a quasi-ultrametric space
of homogeneous type with upper dimension $n\in(0,\infty)$,
$\rho\in\mathbf{q}$ have the property that
all $\rho$-balls are $\mu$-measurable,
and $\mu(X)<\infty$. Let
$p\in(\frac{n}{n+\mathrm{ind\,}(X,\mathbf{q})},1]$,
$q\in(0,\infty]$, $r\in(p,\infty]\cap[1,\infty]$,
and $\varepsilon,\beta,\gamma\in(0,\infty)$ be such that
$n(\frac{1}{p}-1)<\beta,\gamma<\varepsilon\preceq\mathrm{ind\,}(X,\mathbf{q})$.
Then, as subspaces of $({\mathcal{G}}^\varepsilon_0(\beta,\gamma))'$,
$$
H^{p,q,r}_{\mathrm{at}}(X,\rho)=h^{p,q,r}_{\mathrm{at}}(X,\rho)
$$
with equivalent (quasi-)norms.
\end{theorem}

Next, we are ready to prove Theorem~\ref{6.16}.

\begin{proof}[Proof of Theorem~\ref{6.16}]
By Remarks~\ref{5216}(iii) and \ref{RmkD2.3}(iii),
we may assume that $\rho$ is a regularized quasi-ultrametric
satisfying $\varepsilon\leq(\log_2C_\rho)^{-1}$.
Let $\{\mathcal{S}_{2^{-k}}\}_{k\in\mathbb{Z}}$
be a homogeneous approximation
of the identity of order $\varepsilon$ as in Definition~\ref{DefHATI}.
Let $f\in H^{p,q}_*(X)\subset(\mathcal{G}^\varepsilon_0(\beta,\gamma))'$.

If $\|f\|_{H^{p,q}_*(X)}=0$,
then $f^*_\rho=0$ $\mu$-almost everywhere,
which, together with Proposition~\ref{2250},
further implies $f=0$ in
$(\mathcal{G}^\varepsilon_0(\beta,\gamma))'$.
In this case,
for any $i\in\mathbb{N}$ and $k\in\mathbb{Z}$,
let
$$
\lambda_{i,k}:=0
\ \ \text{and}\ \
a_{i,k}:=0.
$$
Then
it is easy to verify that
\begin{align*}
f\in H^{p,q,\infty}_{\mathrm{at}}(X)
\ \ \text{and}\ \
\|f\|_{H^{p,q,\infty}_{\mathrm{at}}(X)}
=0=\|f\|_{H^{p,q}_*(X)}.
\end{align*}
This finishes the proof of Theorem~\ref{6.16}
in the case $\|f\|_{H^{p,q}_*(X)}=0$.

In the remainder of the proof,
we may assume that $\|f\|_{H^{p,q}_*(X)}\in(0,\infty)$.
We first consider the case where $\mu(X)=\infty$.
For any
$m\in\mathbb{N}$, let
\begin{align*}
f^{(m)}:=\mathcal{S}_{2^{-m}}f.
\end{align*}
Note that, from Remark~\ref{Dkf3},
it follows that, for any $m\in\mathbb{N}$,
$f^{(m)}\in(\mathcal{G}^\varepsilon_0(\beta,\gamma))'$.
To obtain the atomic decomposition of $f$ via utilizing $f^{(m)}$,
we divide the proof into the following three steps.

\emph{Step 1.}
In this step, we show that, for any $m\in\mathbb{N}$,
\begin{equation}\label{4e2}
f^{(m)}=\sum_{k\in\mathbb{Z}}\sum_{i\in\mathbb{N}}h^{(m)}_{i,k}
\end{equation}
in $({\mathcal{G}}^\varepsilon_0(\beta,\gamma))'$,
where, for any $k\in\mathbb{Z}$ and $i\in\mathbb{N}$,
$h^{(m)}_{i,k}$ is a $(p,\infty)$-atom multiplied
by a constant depending on $k$ and
$i$, but independent of both $f$ and $m$.
To this end, let $m\in\mathbb{N}$.
We first prove that
\begin{align}\label{1711}
f^{(m)}\in H^{p,q}_*(X)
\ \ \text{and}\ \
\left\|f^{(m)}\right\|_{H^{p,q}_*(X)}
\lesssim\|f\|_{H^{p,q}_*(X)},
\end{align}
where the implicit
positive constant is independent of both $f$ and $m$.
For any given $x\in X$,
let $\varphi\in{\mathcal{G}}^\varepsilon_0(\beta,\gamma)$
satisfy $\|\varphi\|_{G(x,r_0^\alpha,\frac{\beta}{\alpha},
\frac{\gamma}{\alpha})}\leq 1$ for some $r_0\in(0,\infty)$,
where $\alpha:=\min\{1,\,(\log_2 C_\rho)^{-1}\}$.
Then Remark~\ref{RemMf}(iii) implies that
$\|\phi\|_{\mathcal{G}(x,r_0,\beta,\gamma)}\lesssim1$.
We consider the following two cases for $m$.

\emph{Case (1)} $2^{-m}\leq r_0$. In this case,
by Lemma~\ref{4l4.1}(i), we conclude that
$\mathcal{S}_{2^{-m}}\varphi\in{\mathcal{G}}(x,r_0,\beta,\gamma)$
and
$\|\mathcal{S}_{2^{-m}}\varphi\|_{{\mathcal{G}}(x,r_0,\beta,\gamma)}
\lesssim\|\varphi\|_{{\mathcal{G}}(x,r_0,\beta,\gamma)}\lesssim 1$,
where the implicit positive constant
is independent of $\varphi$, $m$, $x$, and $r_0$.
From this, Remark~\ref{Dkf3}(ii), and Definition~\ref{DefATI}(iv),
we infer that
\begin{equation}\label{4e1-1}
\left|\left\langle f^{(m)},\varphi\right\rangle\right|
=\left|\left\langle\mathcal{S}_{2^{-m}}f,\varphi\right\rangle\right|
=\left|\left\langle f,\mathcal{S}_{2^{-m}}
\varphi\right\rangle\right|\lesssim f^*_\rho(x).
\end{equation}

\emph{Case (2)} $2^{-m}>r_0$.
In this case,
according to Lemma~\ref{4l4.1}(ii),
we find that
$\mathcal{S}_{2^{-m}}\varphi\in{\mathcal{G}}(x,2^{-m},\beta,\gamma)$
and
\begin{align*}
\left\|\mathcal{S}_{2^{-m}}\varphi\right\|_{{\mathcal{G}}(x,2^{-m},\beta,\gamma)}
\lesssim\|\varphi\|_{{\mathcal{G}}(x,r_0,\beta,\gamma)}
\lesssim1,
\end{align*}
which, combined with Remark~\ref{Dkf3}(ii)
and Definition~\ref{DefATI}(iv), implies that
\begin{equation}\label{4e1-2}
\left|\left\langle f^{(m)},\varphi\right\rangle\right|
=\left|\left\langle f,\mathcal{S}_{2^{-m}}\varphi\right\rangle\right|
\lesssim f^*_\rho(x).
\end{equation}
Using \eqref{4e1-1} and \eqref{4e1-2} and taking the supremum
over all
$$
\phi\in\mathcal{G}^\varepsilon_0(\beta,\gamma)
\ \ \text{with}\ \
\|\varphi\|_{G(x,r_0^\alpha,\frac{\beta}{\alpha},
\frac{\gamma}{\alpha})}\leq 1
$$
for some $r_0\in(0,\infty)$,
we conclude that,
for any $x\in X$,
\begin{align*}
(f^{(m)})^*_\rho(x)\lesssim f^*_\rho(x),
\end{align*}
which, together with Proposition~\ref{inc-norm}(ii),
then completes the proof of \eqref{1711}.

To establish the atomic decomposition of $f^{(m)}$,
let $k\in\mathbb{Z}$ and define
\begin{align*}
\Omega_k:=\left\{x\in X:f_\rho^*(x)>2^k\right\}.
\end{align*}
From the proof of Proposition~\ref{measurability},
we deduce that $\Omega_k$ is an open set.
Since $f\in H^{p,q}_*(X)$ and $\mu(X)=\infty$,
it immediately follows that $\mu(\Omega_k)<\infty$
and hence $\Omega_k$ is a proper subset of $X$.
Moreover, if $\Omega_k=\varnothing$, then, for any $i\in\mathbb{N}$,
we let $h_{i,k}^{(m)}:=0$ directly.
Since the desired claims still hold true in this case,
without loss of generality, we may assume that
$\Omega_k\neq\varnothing$.
Therefore,
$\Omega_k$ is a proper, nonempty,
and open subset of $X$.
Then, by Lemma~\ref{Wtd},
we find that, for any given $\lambda\in(4C_\rho^6,\infty)$
and $\Lambda\in(4C_\rho^2\lambda,
\frac{\lambda^2}{C_\rho^4})$,
there exist $L_0\in\mathbb{N}$,
a sequence $\{x_{i,k}\}_{i\in\mathbb{N}}$ in  $X$,
a sequence $\{r_{i,k}\}_{i\in \mathbb{N}}$ in $(0,\infty)$,
and a sequence $\{y_{i,k}\}_{i\in \mathbb{N}}$ in $X$
such that (i)--(iv) of Lemma~\ref{Wtd} hold,
which, combined with Lemma~\ref{pou},
further implies that, for any given
$C_0\in(\frac{C_\rho\Lambda}{\lambda},\frac{\lambda}{C_\rho^3})$ and
$\lambda'\in(C_\rho,\frac{\lambda}{C_0C_\rho^2})$,
there exists a sequence
$\{\xi_{i,k}\}_{i\in \mathbb{N}}$ of nonnegative
functions satisfying (i)--(iii) of Lemma~\ref{pou};
in particular, Remark~\ref{lemma5.14iv} holds for
$\{\xi_{i,k}\}_{i\in \mathbb{N}}$.
For any $i\in\mathbb{N}$, let
$B_{i,k}:=B_\rho(x_{i,k},r_{i,k})$,
\begin{align*}
P^{(m)}_{i,k}:=\frac{1}{\|\xi_{i,k}\|_{L^1(X)}}
\int_{X}f^{(m)}(y)\xi_{i,k}(y)\,d\mu(y)
\end{align*}
and, for any $x\in X$,
\begin{align}
\label{1911}
b^{(m)}_{i,k}(x):=\left[f^{(m)}(x)-P^{(m)}_{i,k}\right]\xi_{i,k}(x),
\end{align}
where $P^{(m)}_{i,k}$ is well defined,
$b^{(m)}_{i,k}$ is a $\mu$-measurable function, and $b^{(m)}_{i,k}\in L^1(X)$.
Indeed, from Lemma~\ref{pou}
and an argument similar to that used in the
estimation of \eqref{2235} with $D_0$ therein replaced
by $S_{2^{-m}}$,
we infer that $f^{(m)}$ is continuous and
$f^{(m)}\xi_{i,k}$
is bounded with bounded support and hence is integrable.
Thus, $P^{(m)}_{i,k}$ is well defined and $b^{(m)}_{i,k}\in L^1(X)$.
Then, for any $i\in\mathbb{N}$,
\begin{equation}\label{4e3-1}
\int_{X}b^{(m)}_{i,k}(x)\,d\mu(x)=0.
\end{equation}
By the bounded overlap property of the balls $\{B_{i,k}\}_{i,k}$
in Lemma~\ref{Wtd}(ii) and the fact that
$\mathrm{supp\,}\xi_{i,k}\subset\lambda B_{i,k}$,
we may conclude that, for any $x\in X$,
the series $\sum_{i\in\mathbb{N}}b^{(m)}_{i,k}(x)$
contains only finitely many nonzero terms.

Now, we show that the series
$b^{(m)}_k:=\sum_{i\in\mathbb{N}}b^{(m)}_{i,k}$ converges
to $0$ in $({\mathcal{G}}^\varepsilon_0(\beta,\gamma))'$
uniformly in $m$ as $k\to\infty$.
Indeed, let $i\in\mathbb{N}$ and $x\in X$.
Assume that
$\phi\in\mathcal{G}^\varepsilon_0(\beta,\gamma)$
with $\|\phi\|_{G(x,r_0^\alpha,\frac{\beta}{\alpha},\frac{\gamma}{\alpha})}\leq1$
for some $r_0\in(0,\infty)$.
Then, from Fubini's theorem, we deduce that
\begin{align}\label{215}
\left\langle b_{i,k}^{(m)},\phi\right\rangle
&=\int_X\left[f^{(m)}(y)-P^{(m)}_{i,k}\right]\xi_{i,k}(y)\phi(y)\,d\mu(y)\nonumber\\
&=\int_Xf^{(m)}(y)\Phi_{i,k}(\phi)(y)\,d\mu(y)
=\left\langle f,\mathcal{S}_{2^{-m}}\left(\Phi_{i,k}(\phi)\right)\right\rangle,
\end{align}
where, for any $y\in X$,
$$
\Phi_{i,k}(\phi)(y):=
\frac{\xi_{i,k}(y)}{\|\xi_{i,k}\|_{L^1(X)}}
\int_X\left[\phi(y)-\phi(z)\right]\xi_{i,k}(z)\,d\mu(z).
$$
We claim that, for any $x\in(\lambda B_{i,k})^\complement$,
\begin{align}\label{claim11}
\left\|\Phi_{i,k}(\phi)\right\|_{\mathcal{G}(y_{i,k},r_{i,k},\beta,\gamma)}
\lesssim\frac{(V_\rho)_{r_{i,k}}(x_{i,k})}{(V_\rho)_{r_{i,k}}(x_{i,k})
+V_\rho(x_{i,k},x)}\left[\frac{r_{i,k}}{r_{i,k}+\rho(x_{i,k},x)}\right]^\beta
\end{align}
and, for any $x\in\lambda B_{i,k}$,
\begin{align}\label{claim12}
\left\|\Phi_{i,k}(\phi)\right\|_{\mathcal{G}(x,r_{i,k},\beta,\gamma)}
\lesssim1\ \ \text{if}\ \ r_0\ge r_{i,k}
\end{align}
and
\begin{align}\label{claim13}
\Phi_{i,k}(\phi)=h_1-h_2\ \ \text{if}\ \ r_0<r_{i,k},
\end{align}
where $\|h_1\|_{\mathcal{G}(x,r_0,\beta,\gamma)}\lesssim1$
and $\|h_2\|_{\mathcal{G}(x,r_{i,k},\beta,\gamma)}\lesssim1$.
To prove the size estimate of \eqref{claim11},
observing that $\mathrm{supp\,}\xi_{i,k}\in B_\rho(x_{i,k},\lambda' r_{i,k})$,
we only need to consider the case where $y\in B_\rho(x_{i,k},\lambda' r_{i,k})$.
Let $y\in B_\rho(x_{i,k},\lambda' r_{i,k})$
and $x\in(\lambda B_{i,k})^\complement$.
Note that, for any $z\in B_\rho(x_{i,k},\lambda' r_{i,k})$, it holds that
\begin{align}\label{010}
\rho(z,x_{i,k})<\lambda'r_{i,k}\leq\frac{\lambda'}{\lambda}\rho(x,x_{i,k})
<\frac{1}{C_0C_\rho^2}\rho(x,x_{i,k})
\leq\frac{\rho(x,x_{i,k})}{2C_\rho}<\frac{r_0+\rho(x,x_{i,k})}{2C_\rho},
\end{align}
where we used the fact that $C_0\in(\frac{C_\rho\Lambda}{\lambda},\frac{\lambda}{C_\rho^3})$
ensures $C_0\ge2$ [we also need this in the estimation of \eqref{048}],
which also implies that
\begin{align*}
\lambda r_{i,k}\leq\rho(x,x_{i,k})\leq C_\rho\max\left\{\rho(x,z),\,\rho(z,x_{i,k})\right\}
=C_\rho\rho(x,z).
\end{align*}
Using this, \eqref{010}, the regularity condition of $\phi$,
and Remark~\ref{lemma5.14iv}, we obtain
\begin{align*}
\left|\Phi_{i,k}(\phi)(y)\right|
&\leq\xi_{i,k}(y)\left|\phi(y)-\phi(x_{i,k})\right|\\
&\quad+
\frac{\xi_{i,k}(y)}{\|\xi_{i,k}\|_{L^1(X)}}
\int_X\left|\phi(x_{i,k})-\phi(z)\right|\xi_{i,k}(z)\,d\mu(z)
\nonumber\\
&\lesssim\frac{(V_\rho)_{r_{i,k}}(x_{i,k})}{(V_\rho)_{r_{i,k}}(x_{i,k})
+V_\rho(x_{i,k},y)}\left[\frac{r_{i,k}}{r_{i,k}+\rho(x_{i,k},y)}\right]^\gamma
\nonumber\\
&\quad\times
\left[\frac{r_{i,k}}{r_0+\rho(x,x_{i,k})}\right]^\beta
\frac{1}{(V_\rho)_{r_0}(x)
+V_\rho(x,x_{i,k})}\left[\frac{r_0}{r_0+\rho(x,x_{i,k})}\right]^\gamma,
\end{align*}
which, together with the fact that
$$
r_{i,k}+\rho(y_{i,k},y)\lesssim r_{i,k}+\rho(y_{i,k},x_{i,k})
+\rho(x_{i,k},y)\sim r_{i,k}+\rho(x_{i,k},y),
$$
further implies that
\begin{align}\label{038}
\left|\Phi_{i,k}(\phi)(y)\right|
&\lesssim\frac{(V_\rho)_{r_{i,k}}(x_{i,k})}{(V_\rho)_{r_{i,k}}(x_{i,k})
+V_\rho(x_{i,k},x)}\left[\frac{r_{i,k}}{r_{i,k}+\rho(x_{i,k},x)}\right]^\beta
\nonumber\\
&\quad\times\frac{1}{(V_\rho)_{r_{i,k}}(x_{i,k})
+V_\rho(x_{i,k},y)}\left[\frac{r_{i,k}}{r_{i,k}+\rho(x_{i,k},y)}\right]^\gamma
\nonumber\\
&\lesssim\frac{(V_\rho)_{r_{i,k}}(x_{i,k})}{(V_\rho)_{r_{i,k}}(x_{i,k})
+V_\rho(x_{i,k},x)}\left[\frac{r_{i,k}}{r_{i,k}+\rho(x_{i,k},x)}\right]^\beta
\nonumber\\
&\quad\times\frac{1}{(V_\rho)_{r_{i,k}}(y_{i,k})
+V_\rho(y_{i,k},y)}\left[\frac{r_{i,k}}{r_{i,k}+\rho(y_{i,k},y)}\right]^\gamma.
\end{align}
This is the desired size estimate of \eqref{claim11}.

Next, we show the regularity estimate of \eqref{claim11}.
Let $y,y'\in X$ be such that $\rho(y,y')\leq\frac{r_{i,k}+\rho(y,y_{i,k})}{2C_\rho}$.
Without loss of generality, we may assume that $y\in B_\rho(x_{i,k},\lambda' r_{i,k})$;
otherwise, the desired regularity estimate follows immediately from
the size estimate of \eqref{claim11}.
We write
\begin{align}\label{046}
&\Phi_{i,k}(\phi)(y)-\Phi_{i,k}(\phi)(y')\nonumber\\
&\quad=\frac{\xi_{i,k}(y)}{\|\xi_{i,k}\|_{L^1(X)}}
\int_X\left[\phi(y)-\phi(z)\right]\xi_{i,k}(z)\,d\mu(z)
-\xi_{i,k}(y')\left[\phi(y')-\phi(y)\right]\nonumber\\
&\qquad-\frac{\xi_{i,k}(y')}{\|\xi_{i,k}\|_{L^1(X)}}
\int_X\left[\phi(y)-\phi(z)\right]\xi_{i,k}(z)\,d\mu(z)
\nonumber\\
&\quad=\mathrm{I}+\mathrm{II},
\end{align}
where
\begin{align*}
\mathrm{I}:=\xi_{i,k}(y')\left[\phi(y)-\phi(y')\right]
\end{align*}
and
\begin{align*}
\mathrm{II}:=\frac{\xi_{i,k}(y)-\xi_{i,k}(y')}{\|\xi_{i,k}\|_{L^1(X)}}
\int_X\left[\phi(y)-\phi(z)\right]\xi_{i,k}(z)\,d\mu(z).
\end{align*}
To estimate $\mathrm{I}$, applying the regularity condition of $\phi$,
an argument similar to that used in the estimation of \eqref{038},
and $\rho(x,y)\sim\rho(x_{i,k},x)\gtrsim r_{i,k}$,
we find that
\begin{align}\label{047}
|\mathrm{I}|
&\lesssim\frac{(V_\rho)_{r_{i,k}}(x_{i,k})}{(V_\rho)_{r_{i,k}}(x_{i,k})
+V_\rho(x_{i,k},y)}\left[\frac{r_{i,k}}{r_{i,k}+\rho(x_{i,k},y')}\right]^\gamma
\nonumber\\
&\quad\times
\left[\frac{\rho(y,y')}{r_0+\rho(x,y)}\right]^\beta
\frac{1}{(V_\rho)_{r_0}(x)
+V_\rho(x,y)}\left[\frac{r_0}{r_0+\rho(x,y)}\right]^\gamma
\nonumber\\
&\lesssim\frac{(V_\rho)_{r_{i,k}}(x_{i,k})}{(V_\rho)_{r_{i,k}}(x_{i,k})
+V_\rho(x_{i,k},x)}\left[\frac{r_{i,k}}{r_{i,k}+\rho(x_{i,k},x)}\right]^\beta
\nonumber\\
&\quad\times\left[\frac{\rho(y,y')}{r_{i,k}+\rho(y_{i,k},y)}\right]^\beta
\frac{1}{(V_\rho)_{r_{i,k}}(y_{i,k})
+V_\rho(y_{i,k},y)}\left[\frac{r_{i,k}}{r_{i,k}+\rho(y_{i,k},y)}\right]^\gamma.
\end{align}
To estimate $\mathrm{II}$, we consider the following two cases.

\emph{Case (1)} $\rho(y,y')\leq\frac{r_{i,k}+\rho(x_{i,k},y)}{2C_\rho}$.
In this case, from the regularity condition of both $\xi_{i,k}$ and $\phi$
and an argument similar to that used in the estimation of \eqref{038},
we infer that
\begin{align}\label{048}
|\mathrm{II}|
&\lesssim
\left[\frac{\rho(y,y')}{r_{i,k}+\rho(x_{i,k},y)}\right]^\beta
\frac{(V_\rho)_{r_{i,k}}(x_{i,k})}{(V_\rho)_{r_{i,k}}(x_{i,k})
+V_\rho(x_{i,k},y)}\left[\frac{r_{i,k}}{r_{i,k}+\rho(x_{i,k},y)}\right]^\gamma
\nonumber\\
&\quad\times\left[\frac{r_{i,k}}{r_0+\rho(x,y)}\right]^\beta
\frac{1}{(V_\rho)_{r_0}(x)
+V_\rho(x,y)}\left[\frac{r_0}{r_0+\rho(x,y)}\right]^\gamma
\nonumber\\
&\lesssim\frac{(V_\rho)_{r_{i,k}}(x_{i,k})}{(V_\rho)_{r_{i,k}}(x_{i,k})
+V_\rho(x_{i,k},x)}\left[\frac{r_{i,k}}{r_{i,k}+\rho(x_{i,k},x)}\right]^\beta
\nonumber\\
&\quad\times\left[\frac{\rho(y,y')}{r_{i,k}+\rho(y_{i,k},y)}\right]^\beta
\frac{1}{(V_\rho)_{r_{i,k}}(y_{i,k})
+V_\rho(y_{i,k},y)}\left[\frac{r_{i,k}}{r_{i,k}+\rho(y_{i,k},y)}\right]^\gamma.
\end{align}

\emph{Case (2)} $\rho(y,y')>\frac{r_{i,k}+\rho(x_{i,k},y)}{2C_\rho}$.
In this case, applying the size condition of $\xi_{i,k}$, the
regularity condition of $\phi$,
and an argument similar to that used in the estimation of \eqref{048},
we conclude that
\begin{align*}
|\mathrm{II}|
&\lesssim\frac{(V_\rho)_{r_{i,k}}(x_{i,k})}{(V_\rho)_{r_{i,k}}(x_{i,k})
+V_\rho(x_{i,k},x)}\left[\frac{r_{i,k}}{r_{i,k}+\rho(x_{i,k},x)}\right]^\beta
\nonumber\\
&\quad\times\left[\frac{\rho(y,y')}{r_{i,k}+\rho(y_{i,k},y)}\right]^\beta
\frac{1}{(V_\rho)_{r_{i,k}}(y_{i,k})
+V_\rho(y_{i,k},y)}\left[\frac{r_{i,k}}{r_{i,k}+\rho(y_{i,k},y)}\right]^\gamma.
\end{align*}
From this, \eqref{046}, \eqref{047}, and \eqref{048},
we deduce that the desired regularity estimate of \eqref{claim11} holds,
which completes the proof of \eqref{claim11}.
The proof of \eqref{claim12} is similar to that of \eqref{claim11}
and hence we omit the details. To prove \eqref{claim13},
note that, for any $y\in X$, we can write
\begin{align*}
\Phi_{i,k}(\phi)(y)=\xi_{i,k}(y)\phi(y)-\frac{\xi_{i,k}(y)}{\|\xi_{i,k}\|_{L^1(X)}}
\int_X\phi(z)\xi_{i,k}(z)\,d\mu(z)=:h_1-h_2.
\end{align*}
Then it is easy to show that $\|h_1\|_{\mathcal{G}(x,r_0,\beta,\gamma)}\lesssim1$
and $\|h_2\|_{\mathcal{G}(x,r_{i,k},\beta,\gamma)}\lesssim1$
via an argument similar to that used in the proof of \eqref{claim11}.
This finishes the proof of \eqref{claim13} and hence the
proof of the aforementioned claim.

By \eqref{claim11} and Lemma~\ref{4l4.1},
we find that, for any $x\in(\lambda B_{i,k})^\complement$,
\begin{align}\label{21456}
&\left\|\mathcal{S}_{2^{-m}}\left(\Phi_{i,k}(\phi)\right)
\right\|_{\mathcal{G}(y_{i,k},\max\{r_{i,k},\,2^{-m}\},\beta,\gamma)}
\nonumber\\
&\quad\lesssim\frac{(V_\rho)_{r_{i,k}}(x_{i,k})}{(V_\rho)_{r_{i,k}}(x_{i,k})
+V_\rho(x_{i,k},x)}\left[\frac{r_{i,k}}{r_{i,k}+\rho(x_{i,k},x)}\right]^\beta.
\end{align}
From \eqref{claim12}, \eqref{claim13},
and Lemma~\ref{4l4.1},
we infer that, for any $x\in\lambda B_{i,k}$,
\begin{align*}
\left\|\mathcal{S}_{2^{-m}}\left(\Phi_{i,k}(\phi)\right)
\right\|_{\mathcal{G}(x,\max\{r_{i,k},\,2^{-m}\},\beta,\gamma)}\lesssim1.
\end{align*}
By this, \eqref{215}, \eqref{21456}, Remark~\ref{RemMf}(iii),
and taking the supremum over all $\phi\in\mathcal{G}^\varepsilon_0(\beta,\gamma)$
with $\|\phi\|_{G(x,r_0^\alpha,\frac{\beta}{\alpha},\frac{\gamma}{\alpha})}\leq1$
for $r_0\in(0,\infty)$,
we conclude that
\begin{align*}
\left(b^{(m)}_{i,k}\right)^*_\rho(x)
&\lesssim f^*_\rho(y_{i,k})
\frac{(V_\rho)_{r_{i,k}}(x_{i,k})}{(V_\rho)_{r_{i,k}}(x_{i,k})+V_\rho(x_{i,k},x)}
\left[\frac{r_{i,k}}{r_{i,k}+\rho(x_{i,k},x)}
\right]^\beta\mathbf{1}_{(\lambda B_{i,k})^\complement}(x)\\
&\quad+f^*_\rho(x)\mathbf{1}_{\lambda B_{i,k}}(x)\\
&\leq 2^k\frac{(V_\rho)_{r_{i,k}}(x_{i,k})}{(V_\rho)_{r_{i,k}}(x_{i,k})+V_\rho(x_{i,k},x)}
\left[\frac{r_{i,k}}{r_{i,k}+\rho(x_{i,k},x)}
\right]^\beta\mathbf{1}_{(\lambda B_{i,k})^\complement}(x)\\
&\quad+f^*_\rho(x)\mathbf{1}_{\lambda B_{i,k}}(x),
\end{align*}
which, together with both (i) and (iii) of Lemma~\ref{Wtd} and Lemma~\ref{2051},
further implies that, for
any given $p_0\in(\frac{n}{n+\beta},\min\{1,\,p\})$,
\begin{align}\label{4e5}
\left\|\left(\sum_{i\in\mathbb{N}}b^{(m)}_{i,k}
\right)^*_\rho\right\|_{L^{p_0}(X)}^{p_0}
&\leq\int_X\left[\sum_{i\in\mathbb{N}}\left(b^{(m)}_{i,k}
\right)^*_\rho(x)\right]^{p_0}\,d\mu(x)
\nonumber\\
&\lesssim2^{kp_0}\int_X\left\{\sum_{i\in\mathbb{N}}
\frac{(V_\rho)_{r_{i,k}}(x_{i,k})}{(V_\rho)_{r_{i,k}}(x_{i,k})+V_\rho(x_{i,k},x)}
\left[\frac{r_{i,k}}{r_{i,k}
+\rho(x_{i,k},x)}\right]^\beta\right\}^{p_0}\,d\mu(x)
\nonumber\\
&\quad+\int_{\bigcup_{i\in\mathbb{N}}\lambda B_{i,k}}
\left[f^*_\rho(x)\right]^{p_0}\,d\mu(x)
\nonumber\\
&\lesssim2^{kp_0}\mu(\Omega_k)+\left\|f^*_\rho
\mathbf{1}_{\Omega_k}\right\|_{L^{p_0}(X)}^{p_0}
\sim\left\|f^*_\rho\mathbf{1}_{\Omega_k}\right\|_{L^{p_0}(X)}^{p_0}.
\end{align}
Note that, for any $\alpha\in (0,\infty)$,
\begin{equation*}
\mu\left(\Omega_k\cap \left\{x\in X:f^*_\rho(x)>\alpha\right\}\right)
\leq \min\left\{\mu(\Omega_k),\,
\left[\frac{\|f^*_\rho\|_{L^{p,\infty}(X)}}{\alpha}\right]^p
\right\},
\end{equation*}
which, combined with Proposition~\ref{p,q<p,s},
further implies that
\begin{align*}
\left\|f^*_\rho\mathbf{1}_{\Omega_k}\right\|_{L^{p_0}(X)}^{p_0}
&=\int^{\infty}_0p_0\alpha^{p_0-1}
\mu\left(\left\{x\in\Omega_k:
f^*_\rho(x)>\alpha\right\}\right)\,d\alpha
\nonumber\\
&\leq \int^{[\mu(\Omega_k)]^{-\frac{1}{p}}
\|f^*_\rho\|_{L^{p,\infty}(X)}}_0p_0\alpha^{p_0-1}\mu(\Omega_k)\,d\alpha
\nonumber\\
&\quad+
\int^{\infty}_{[\mu(\Omega_k)]^{-\frac{1}{p}}\|f^*_\rho\|_{L^{p,\infty}(X)}}
p_0\alpha^{p_0-p-1}\|f^*_\rho\|^p_{L^{p,\infty}(X)}\,d\alpha
\nonumber\\
&=\frac{p}{p-p_0}\left[\mu(\Omega_k)\right]^{1-\frac{p_0}{p}}
\left\|f^*_\rho\right\|^{p_0}_{L^{p,\infty}(X)}
\lesssim\left[\mu(\Omega_k)\right]^{1-\frac{p_0}{p}}
\left\|f^*_\rho\right\|^{p_0}_{L^{p,q}(X)}.
\end{align*}
From this, \eqref{4e5}, $p_0<p$, and the fact that
$\mu(\Omega_k)\to 0$ as $k\to \infty$, it follows that
\begin{align}
\label{x}
\left\|\sum_{i\in\mathbb{N}}b^{(m)}_{i,k}\right\|_{H_*^{p_0}(X)}
&=\left\|\left(\sum_{i\in\mathbb{N}}b^{(m)}_{i,k}\right)^*_\rho\right\|_{L^{p_0}(X)}
\lesssim\left\|f^*_\rho\mathbf{1}_{\Omega_k}\right\|_{L^{p_0}(X)}
\nonumber\\
&\lesssim \left[\mu(\Omega_k)\right]^{\frac{1}{p_0}-\frac{1}{p}}
\left\|f^*_\rho\right\|_{L^{p,q}(X)}
\to 0
\end{align}
as $k\to \infty$,
where the implicit positive constants are independent of $m$,
which, together with Proposition~\ref{5.1.10},
further implies that
$$
\lim_{k\to\infty}\left\|\sum_{i\in\mathbb{N}}b^{(m)}_{i,k}
\right\|_{(\mathcal{G}^\varepsilon_0(\beta,\gamma))'}
\lesssim\lim_{k\to\infty}\left\|\sum_{i\in\mathbb{N}}
b^{(m)}_{i,k}\right\|_{H_*^{p_0}(X)}=0
$$
and hence the series
$b^{(m)}_k:=\sum_{i\in\mathbb{N}}b^{(m)}_{i,k}$ converges to $0$
in $({\mathcal{G}}^\varepsilon_0(\beta,\gamma))'$ uniformly in $m$ as $k\to\infty$.

Using \eqref{1911} and Lemma \ref{pou}(iii),
we find that, for any $x\in X$,
\begin{align}\label{defg}
g^{(m)}_k(x):&=f^{(m)}(x)-b^{(m)}_k(x)
\nonumber\\
&=f^{(m)}(x)-\sum_{i\in\mathbb{N}}\left[f^{(m)}(x)-P^{(m)}_{i,k}\right]\xi_{i,k}(x)
\nonumber\\
&=f^{(m)}(x)\mathbf{1}_{\left(\Omega_k\right)^{\complement}}(x)
+\sum_{i\in\mathbb{N}}P^{(m)}_{i,k}\xi_{i,k}(x).
\end{align}
Now, we prove that
\begin{equation}\label{4e7}
f^{(m)}=\sum^{\infty}_{k=-\infty}\left[g^{(m)}_{k+1}-g^{(m)}_k\right]
\end{equation}
in $({\mathcal{G}}^\varepsilon_0(\beta,\gamma))'$.
Indeed, Lemma \ref{Wtd}(iii) implies that, for any $i\in\mathbb{N}$,
there exists $y_{i,k}\in(\Omega_k)^{\complement}$
such that $\rho(x_{i,k},y_{i,k})<\Lambda r_{i,k}$,
which, combined with Lemma~\ref{A3.2}(iii),
further implies that, for any $x\in X$,
$r_{i,k}+\rho(y_{i,k},x)\sim r_{i,k}+\rho(x_{i,k},x)$.
From this, Definition~\ref{testfunction},
$\beta,\gamma\in(0,\varepsilon)$,
and Remark~\ref{lemma5.14iv}, we deduce that,
for any $i\in\mathbb{N}$,
\begin{align*}
\|\xi_{i,k}\|_{{\mathcal{G}}(y_{i,k},r_{i,k},\beta,\gamma)}
\sim\|\xi_{i,k}\|_{{\mathcal{G}}(x_{i,k},r_{i,k},\beta,\gamma)}
\lesssim(V_\rho)_{r_{i,k}}(x_{i,k}),
\end{align*}
which, together with Lemma~\ref{4l4.1},
further implies that
$$
\left\|\mathcal{S}_{2^{-m}}\xi_{i,k}\right\|_{\mathcal{G}(y_{i,k},
r_{i,k},\beta,\gamma)}\lesssim(V_\rho)_{r_{i,k}}(x_{i,k})
$$
if $r_{i,k}\ge2^{-m}$
and
$$
\left\|\mathcal{S}_{2^{-m}}\xi_{i,k}\right\|_{\mathcal{G}(y_{i,k},
2^{-m},\beta,\gamma)}\lesssim(V_\rho)_{r_{i,k}}(x_{i,k})
$$
if $r_{i,k}<2^{-m}$.
By this, Remark~\ref{Dkf3}(ii),
Definition~\ref{DefATI}(iv), and
Lemma~\ref{pou}(ii),
we conclude that, for any $i\in\mathbb{N}$,
\begin{align}\label{4e4-2}
\left|P^{(m)}_{i,k}\right|
&=\left|\left\langle\mathcal{S}_{2^{-m}}f,\frac{\xi_{i,k}}
{\|\xi_{i,k}\|_{L^1(X)}}\right\rangle\right|
=\left|\left\langle f,\frac{\mathcal{S}_{2^{-m}}\xi_{i,k}}
{\|\xi_{i,k}\|_{L^1(X)}}\right\rangle\right|\nonumber\\
&\lesssim
\left|\left\langle f,\frac{\mathcal{S}_{2^{-m}}\xi_{i,k}}
{(V_\rho)_{r_{i,k}}(x_{i,k})}\right\rangle\right|
\lesssim f^*_\rho(y_{i,k})\leq2^k.
\end{align}
From the definitions of both $\mathcal{M}^+f$ and $\Omega_k$
and Lemma~\ref{MMM},
it follows that, for any $x\in\Omega_k^\complement$,
\begin{align}
\label{1453}
\left|f^{(m)}(x)\right|
\leq\mathcal{M}^+f(x)
\lesssim f^*_\rho(x)\leq2^k.
\end{align}
Moreover, by Lemmas~\ref{pou}(ii)
and~\ref{Wtd}(ii) and $\lambda'<\lambda$,
we find that, for any $x\in X$,
\begin{align*}
0\leq\sum_{i\in\mathbb{N}}\xi_{i,k}(x)\leq
\sum_{i\in\mathbb{N}}\mathbf{1}_{\lambda'B_{i,k}}(x)\leq
\sum_{i\in\mathbb{N}}\mathbf{1}_{\lambda B_{i,k}}(x)\leq L_0,
\end{align*}
which, combined with \eqref{defg}, \eqref{1453},
and \eqref{4e4-2},
further implies that
\begin{align*}
\left\|g^{(m)}_k\right\|_{L^{\infty}(X)}\leq \left\|f^{(m)}
\mathbf{1}_{\left(\Omega_k\right)^{\complement}}\right\|_{L^{\infty}(X)}
+\left\|\sum_{i\in\mathbb{N}}P^{(m)}_{i,k}\xi_{i,k}\right\|_{L^{\infty}(X)}
\lesssim 2^k.
\end{align*}
From this,  Proposition~\ref{5.1.10},
the fact that $L^\infty(X)\hookrightarrow({\mathcal{G}}^\varepsilon_0(\beta,\gamma))'$
(see Proposition~\ref{distfncttype}),
and \eqref{x},
we infer that, for any
$\phi\in{\mathcal{G}}^\varepsilon_0(\beta,\gamma)$ with
$\|\phi\|_{{\mathcal{G}}^\varepsilon_0(\beta,\gamma)}\leq1$,
\begin{align*}
&\left|\left\langle f^{(m)}-\sum^N_{k=-N}\left[g^{(m)}_{k+1}
-g^{(m)}_k\right],\phi\right\rangle\right|\\
&\quad\leq
\left|\left\langle\sum_{i\in\mathbb{N}}b^{(m)}_{i,N+1},\phi\right\rangle\right|
+\left|\left\langle g^{(m)}_{-N},\phi\right\rangle\right|\\
&\quad\leq\left\|\sum_{i\in\mathbb{N}}b^{(m)}_{i,N+1}
\right\|_{({\mathcal{G}}^\varepsilon_0(\beta,\gamma))'}
+\left\|g^{(m)}_{-N}\right\|_{({\mathcal{G}}^\varepsilon_0(\beta,\gamma))'}\\
&\quad\lesssim
\left\|\sum_{i\in\mathbb{N}}b^{(m)}_{i,N+1}\right\|_{H^{p_0}_*(X)}
+\left\|g^{(m)}_{-N}\right\|_{L^\infty(X)}
\to0
\end{align*}
as $N\to\infty$,
which further implies that \eqref{4e7} holds.

For any $i,j\in\mathbb{N}$, let
\begin{align}
\label{1916}
L^{m,k+1}_{i,j}:&=\frac{1}{\|\xi_{j,k+1}\|_{L^1(X)}}
\int_{X}b^{(m)}_{j,k+1}(y)\xi_{i,k}(y)\,d\mu(y)
\nonumber\\
&=\frac{1}{\|\xi_{j,k+1}\|_{L^1(X)}}
\int_{X}\left[f^{(m)}(y)-
P^{(m)}_{j,k+1}\right]\xi_{i,k}(y)\xi_{j,k+1}(y)\,d\mu(y).
\end{align}
Next, we show that, for any $i,j\in\mathbb{N}$,
\begin{equation}\label{4e12}
\left|L^{m,k+1}_{i,j}\right|
\lesssim 2^k.
\end{equation}
To this end, let $i,j\in\mathbb{N}$
and we first prove that
\begin{align}\label{1600}
\left\|\frac{\xi_{i,k}\xi_{j,k+1}}{\|\xi_{j,k+1}\|_{L^1(X)}}\right\|
_{{\mathcal{G}}(y_{j,k+1},r_{j,k+1},\varepsilon,\varepsilon)}\lesssim 1.
\end{align}
If $\xi_{i,k}\xi_{j,k+1}=0$ pointwise everywhere in $X$, then
\eqref{1600} holds automatically.
Assume that $\xi_{i,k}\xi_{j,k+1}$
is a non-zero function;
in this case, by Lemma~\ref{pou}(ii),
we conclude that
$$
\varnothing\neq\left[\mathrm{\,supp\,}\xi_{i,k}
\cap\mathrm{\,supp\,}\xi_{j,k+1}\right]\subset
\left[B_\rho(x_{i,k},\lambda'r_{i,k})\cap
B_\rho(x_{j,k+1},\lambda'r_{j,k+1})\right],
$$
and hence we can choose a point
$z_0\in B_\rho(x_{i,k},\lambda'r_{i,k})\cap
B_\rho(x_{j,k+1},\lambda'r_{j,k+1})$. To proceed,
we turn to show that
\begin{align}
\label{1615}
r_{j,k+1}\leq C_0 r_{i,k}.
\end{align}
Indeed, if $r_{j,k+1}>C_0 r_{i,k}$,
then, from Definition~\ref{D2.1}(iii) twice, $C_0>\frac{C_\rho\Lambda}{\lambda}>1$,
and $\lambda'<\frac{\lambda}{C_0C_\rho^2}$,
we deduce that, for any $y\in B_\rho(x_{i,k},\Lambda r_{i,k})$,
\begin{align*}
\rho(y,x_{j,k+1})
&\leq C_\rho\max\left\{\rho(y,x_{i,k}),\,C_\rho\max\left\{\rho(x_{i,k},z_0),\,
\rho(z_0,x_{j,k+1})\right\}\right\}
\\
&\leq\max\left\{C_\rho\Lambda r_{i,k},\,C_\rho^2\lambda'r_{i,k},\,
C_\rho^2\lambda'r_{j,k+1}\right\}
\\
&\leq\max\left\{\frac{C_\rho\Lambda}{C_0} r_{j,k+1},\,
\frac{C_\rho^2\lambda'}{C_0}r_{j,k+1},\,
C_\rho^2\lambda'r_{j,k+1}\right\}
\\
&\leq\max\left\{\lambda,\,\frac{\lambda}{C_0^2},\,
\frac{\lambda}{C_0}\right\} r_{j,k+1}
\leq\lambda r_{j,k+1},
\end{align*}
which further implies that
\begin{align*}
B_\rho(x_{i,k},\Lambda r_{i,k})\subset B_\rho(x_{j,k+1},\lambda r_{j,k+1})
\subset\Omega_{k+1}\subset\Omega_k.
\end{align*}
This clearly contradicts Lemma~\ref{Wtd}(iii) and hence \eqref{1615} holds.

To prove the size condition
of $\frac{\xi_{i,k}\xi_{j,k+1}}{\|\xi_{j,k+1}\|_{L^1(X)}}$,
by Lemma~\ref{pou}(ii), we find that
$$
\mathrm{supp\,}\frac{\xi_{i,k}\xi_{j,k+1}}{\|\xi_{j,k+1}\|_{L^1(X)}}\subset\lambda' B_{j,k+1}.
$$
From Lemmas~\ref{Wtd}(iii) and~\ref{pou}(ii) and Remark~\ref{lemma5.14iv},
we infer that, for any $x\in\lambda' B_{j,k+1}$,
\begin{align}\label{Xisize}
\left|\frac{\xi_{i,k}(x)\xi_{j,k+1}(x)}{\|\xi_{j,k+1}\|_{L^1(X)}}\right|
&\lesssim\frac{|\xi_{j,k+1}(x)|}{(V_\rho)_{r_{j,k+1}}(x_{j,k+1})}
\nonumber\\
&\lesssim\frac{1}{(V_\rho)_{r_{j,k+1}}(x_{j,k+1})+V_\rho(x_{j,k+1},x)}
\left[\frac{r_{j,k+1}}{r_{j,k+1}+\rho(x_{j,k+1},x)}\right]^\varepsilon
\nonumber\\
&\sim\frac{1}{(V_\rho)_{r_{j,k+1}}(y_{j,k+1})+V_\rho(y_{j,k+1},x)}
\left[\frac{r_{j,k+1}}{r_{j,k+1}+\rho(y_{j,k+1},x)}\right]^\varepsilon.
\end{align}

Now, we show the regularity condition of
$\frac{\xi_{i,k}\xi_{j,k+1}}{\|\xi_{j,k+1}\|_{L^1(X)}}$.
Let $x,x'\in X$ satisfy $\rho(x,x')\leq\frac{r_{j,k+1}+\rho(y_{j,k+1},x)}{2C_\rho}$.
We consider the following two cases for $\rho(x,x')$.

\emph{Case (1)} $\rho(x,x')\leq\frac{r_{j,k+1}+\rho(y_{j,k+1},x)}{(2C_\rho)^2}$.
In this case,
we may assume that
\begin{align*}
\xi_{i,k}(x)\xi_{j,k+1}(x)-\xi_{i,k}(x')\xi_{j,k+1}(x')\neq0
\end{align*}
since otherwise the desired estimate is trivially satisfied. We claim that
\begin{align}\label{155}
\rho(x,x_{j,k+1})\leq C_\rho\Lambda r_{j,k+1}.
\end{align}
Indeed, if $\rho(x,x_{j,k+1})>C_\rho\Lambda r_{j,k+1}$,
then, by $\Lambda>\lambda'$ and $\mathrm{supp\,}\frac{\xi_{i,k}
\xi_{j,k+1}}{\|\xi_{j,k+1}\|_{L^1(X)}}\subset\lambda' B_{j,k+1}$, we conclude that
$\xi_{i,k}(x)\xi_{j,k+1}(x)=0$.
From Definition~\ref{D2.1}(iii)
and $\rho(x_{j,k+1},y_{j,k+1})<\Lambda r_{j,k+1}$,
we deduce that
$(C_\rho)^{-1}\rho(x,y_{j,k+1})\leq\rho(x,x_{j,k+1})
\leq C_\rho\rho(x,y_{j,k+1})$,
which, together with $\rho(x,x_{j,k+1})>C_\rho r_{j,k+1}$,
implies that
$$
\rho(x,x')\leq\frac{\rho(x,y_{j,k+1})}{2C_\rho^2}
\leq\frac{\rho(x,x_{j,k+1})}{2C_\rho}.
$$
Using this and Definition~\ref{D2.1}(iii) again,
we obtain
$$
C_\rho\Lambda r_{j,k+1}<\rho(x,x_{j,k+1})
\leq C_\rho\max\left\{\rho(x,x'),\,
\rho(x',x_{j,k+1})\right\}=C_\rho\rho(x',x_{j,k+1}),
$$
which, combined with $\Lambda>\lambda\ge\lambda'$,
implies that $\xi_{i,k}(x')\xi_{j,k+1}(x')=0$.
Thus,
$$
\xi_{i,k}(x)\xi_{j,k+1}(x)-\xi_{i,k}(x')\xi_{j,k+1}(x')=0
$$
which contradicts
$\xi_{i,k}(x)\xi_{j,k+1}(x)-\xi_{i,k}(x')\xi_{j,k+1}(x')\neq0$.
Hence, \eqref{155} holds in this case.
By Remark~\ref{lemma5.14iv} and
both (i) and (ii) of Lemma~\ref{pou},
we find that
\begin{align*}
&\left|\frac{\xi_{i,k}(x)\xi_{j,k+1}(x)}{\|\xi_{j,k+1}\|_{L^1(X)}}
-\frac{\xi_{i,k}(x')\xi_{j,k+1}(x')}{\|\xi_{j,k+1}\|_{L^1(X)}}\right|\\
&\quad\lesssim\frac{\xi_{i,k}(x)|\xi_{j,k+1}(x)-\xi_{j,k+1}(x')|
+|\xi_{i,k}(x)-\xi_{i,k}(x')|\xi_{j,k+1}(x')}{(V_\rho)_{r_{j,k+1}}(x_{j,k+1})}\\
&\quad\lesssim\frac{1}{(V_\rho)_{r_{j,k+1}}(x_{j,k+1})}
\left\{\left[\frac{\rho(x,x')}{r_{j,k+1}}\right]^\varepsilon
+\left[\frac{\rho(x,x')}{r_{i,k}}\right]^\varepsilon\right\}.
\end{align*}
From this, \eqref{155}, \eqref{1615}, and the fact that
$\rho(y_{j,k+1},x_{j,k+1})\sim r_{j,k+1}$,
it follows that
\begin{align*}
&\left|\frac{\xi_{i,k}(x)\xi_{j,k+1}(x)}{\|\xi_{j,k+1}\|_{L^1(X)}}
-\frac{\xi_{i,k}(x')\xi_{j,k+1}(x')}{\|\xi_{j,k+1}\|_{L^1(X)}}\right|\\
&\quad\sim\left[\frac{\rho(x,x')}{r_{j,k+1}+\rho(x_{j,k+1},x)}\right]^\varepsilon
\frac{1}{(V_\rho)_{r_{j,k+1}}(x_{j,k+1})+V_\rho(x_{j,k+1},x)}\\
&\quad\sim\left[\frac{\rho(x,x')}{r_{j,k+1}+\rho(y_{j,k+1},x)}\right]^\varepsilon\\
&\qquad\times\frac{1}{(V_\rho)_{r_{j,k+1}}(y_{j,k+1})+V_\rho(y_{j,k+1},x)}
\left[\frac{r_{j,k+1}}{r_{j,k+1}+\rho(y_{j,k+1},x)}\right]^\varepsilon.
\end{align*}
Thus,  the desired regularity estimate holds in this case.

\emph{Case (2)}
$\frac{r_{j,k+1}+\rho(y_{j,k+1},x)}{(2C_\rho)^2}
\leq\rho(x,x')\leq\frac{r_{j,k+1}+\rho(y_{j,k+1},x)}{2C_\rho}$.
In this case, by Definition~\ref{D2.1}(iii), we conclude that
\begin{align*}
\rho(y_{j,k+1},x)
&\leq C_\rho\max\left\{\rho(y_{j,k+1},x'),\,\rho(x',x)\right\}
\\
&\leq C_\rho\max\left\{\rho(y_{j,k+1},x'),\,
\frac{r_{j,k+1}+\rho(y_{j,k+1},x)}{2C_\rho}\right\}
\\
&\leq C_\rho\rho(y_{j,k+1},x')+\frac{r_{j,k+1}+\rho(y_{j,k+1},x)}{2},
\end{align*}
from which we can infer that $\rho(y_{j,k+1},x)\lesssim r_{j,k+1}+\rho(y_{j,k+1},x')$.
Similarly, we have
\begin{align*}
\rho(y_{j,k+1},x')
&\leq C_\rho\max\left\{\rho(y_{j,k+1},x),\,\rho(x,x')\right\}
\\
&\leq C_\rho\max\left\{\rho(y_{j,k+1},x),\,
\frac{r_{j,k+1}+\rho(y_{j,k+1},x)}{2C_\rho}\right\}
\\
&\lesssim r_{j,k+1}+\rho(y_{j,k+1},x).
\end{align*}
Hence,
\begin{align}\label{2238}
r_{j,k+1}+\rho(y_{j,k+1},x)\sim r_{j,k+1}+\rho(y_{j,k+1},x').
\end{align}
By the proven size estimate, we find that
\begin{align*}
&\left|\frac{\xi_{i,k}(x)\xi_{j,k+1}(x)}{\|\xi_{j,k+1}\|_{L^1(X)}}
-\frac{\xi_{i,k}(x')\xi_{j,k+1}(x')}{\|\xi_{j,k+1}\|_{L^1(X)}}\right|\\
&\quad\leq\left|\frac{\xi_{i,k}(x)\xi_{j,k+1}(x)}{\|\xi_{j,k+1}\|_{L^1(X)}}
\right|+\left|\frac{\xi_{i,k}(x')\xi_{j,k+1}(x')}{\|\xi_{j,k+1}\|_{L^1(X)}}\right|\\
&\quad\lesssim\frac{1}{(V_\rho)_{r_{j,k+1}}(y_{j,k+1})+V_\rho(y_{j,k+1},x)}
\left[\frac{r_{j,k+1}}{r_{j,k+1}+\rho(y_{j,k+1},x)}\right]^\varepsilon\\
&\qquad+\frac{1}{(V_\rho)_{r_{j,k+1}}(y_{j,k+1})+V_\rho(y_{j,k+1},x')}
\left[\frac{r_{j,k+1}}{r_{j,k+1}+\rho(y_{j,k+1},x')}\right]^\varepsilon.
\end{align*}
From this, \eqref{2238}, and Lemma~\ref{A3.2}(vi),
we deduce that
\begin{align*}
&\left|\frac{\xi_{i,k}(x)\xi_{j,k+1}(x)}{\|\xi_{j,k+1}\|_{L^1(X)}}
-\frac{\xi_{i,k}(x')\xi_{j,k+1}(x')}{\|\xi_{j,k+1}\|_{L^1(X)}}\right|\\
&\quad\sim\frac{1}{(V_\rho)_{r_{j,k+1}}(y_{j,k+1})+V_\rho(y_{j,k+1},x)}
\left[\frac{r_{j,k+1}}{r_{j,k+1}+\rho(y_{j,k+1},x)}\right]^\varepsilon\\
&\quad\sim\left[\frac{\rho(x,x')}{r_{j,k+1}+\rho(y_{j,k+1},x)}\right]^\varepsilon\\
&\qquad\times
\frac{1}{(V_\rho)_{r_{j,k+1}}(y_{j,k+1})+V_\rho(y_{j,k+1},x)}
\left[\frac{r_{j,k+1}}{r_{j,k+1}+\rho(y_{j,k+1},x)}\right]^\varepsilon.
\end{align*}
In summary, we have shown that the desired regularity estimate holds.
By this and the above proven size estimate \eqref{Xisize},
we conclude that \eqref{1600} holds
and hence
\begin{align*}
\left\|\frac{\xi_{i,k}\xi_{j,k+1}}{\|\xi_{j,k+1}\|_{L^1(X)}}\right\|
_{{\mathcal{G}}(y_{j,k+1},r_{j,k+1},\beta,\gamma)}\leq\left\|
\frac{\xi_{i,k}\xi_{j,k+1}}{\|\xi_{j,k+1}\|_{L^1(X)}}\right\|
_{{\mathcal{G}}(y_{j,k+1},r_{j,k+1},\varepsilon,\varepsilon)}\lesssim1.
\end{align*}
This, together with Lemma~\ref{4l4.1},
further implies that
\begin{align}\label{1427}
\left\|\mathcal{S}_{2^{-m}}
\left(\frac{\xi_{i,k}\xi_{j,k+1}}
{\|\xi_{j,k+1}\|_{L^1(X)}}\right)\right\|
_{{\mathcal{G}}(y_{j,k+1},\max\{2^{-m},\,r_{j,k+1}\},\beta,\gamma)}\lesssim1.
\end{align}
Note that
\begin{align*}
\left|L^{m,k+1}_{i,j}\right|
&\leq\frac{|P^{(m)}_{j,k+1}|}{\|\xi_{j,k+1}\|_{L^1(X)}}
\int_{X}\xi_{i,k}(y)\xi_{j,k+1}(y)\,d\mu(y)
\nonumber\\
&\quad+\left|\int_{X}f^{(m)}(y)\frac{\xi_{i,k}(y)
\xi_{j,k+1}(y)}{\|\xi_{j,k+1}\|_{L^1(X)}}\,d\mu(y)\right|
\nonumber\\
&\leq\left|P^{(m)}_{j,k+1}\right|
+\left|\left\langle f^{(m)},\frac{\xi_{i,k}
\xi_{j,k+1}}{\|\xi_{j,k+1}\|_{L^1(X)}}\right\rangle\right|.
\end{align*}
From this, \eqref{4e4-2}, Remark~\ref{Dkf3}(ii), \eqref{1427},
and $y_{j,k+1}\in(\Omega_{k+1})^\complement$,
we infer that
\begin{align*}
\left|L^{m,k+1}_{i,j}\right|
\lesssim2^k+\left|\left\langle f,\mathcal{S}_{2^{-m}}
\left(\frac{\xi_{i,k}\xi_{j,k+1}}
{\|\xi_{j,k+1}\|_{L^1(X)}}\right)\right\rangle\right|
\lesssim 2^k+f^*_\rho(y_{j,k+1})
\sim 2^k.
\end{align*}
This finishes the proof of \eqref{4e12}.

Next, we show that, for any $i,j\in\mathbb{N}$ and $x\in X$,
\begin{equation}\label{4e10}
\sum_{j\in\mathbb{N}}\left|L^{m,k+1}_{i,j}\xi_{j,k+1}(x)\right
|\lesssim 2^k\mathbf{1}_{\lambda B_{i,k}}(x).
\end{equation}
Let $i,j\in\mathbb{N}$ and $x\in X$.
We first claim that,
if $L^{m,k+1}_{i,j}\neq 0$, then
\begin{align}
\label{1245}
B_\rho(x_{j,k+1},\lambda'r_{j,k+1})\subset
B_\rho(x_{i,k},\lambda r_{i,k})\subset\Omega_k.
\end{align}
Indeed, by \eqref{1615}, Definition~\ref{D2.1}(iii), and the facts that $C_0>1$ and
$\lambda'<\frac{\lambda}{C_0C_\rho^2}$,
we find that,
for any $y\in B_\rho(x_{j,k+1},\lambda'r_{j,k+1})$,
\begin{align*}
\rho(y,x_{i,k})
&\leq C_\rho\max\left\{\rho(y,x_{j,k+1}),\,\rho(x_{j,k+1},x_{i,k})\right\}
\\
&\leq C_\rho\max\left\{\rho(y,x_{j,k+1}),\,C_\rho\max\left\{\rho(x_{j,k+1},z_0),\,
\rho(z_0,x_{i,k})\right\}\right\}
\\
&\leq\max\left\{C_\rho\lambda'r_{j,k+1},\,C_\rho^2\lambda'r_{j,k+1},\,
C_\rho^2\lambda'r_{i,k}\right\}
\leq C_\rho^2\lambda' C_0r_{i,k}<\lambda r_{i,k},
\end{align*}
which, combined with Lemma~\ref{Wtd}(iii), implies \eqref{1245}.
From \eqref{1245}, (ii) and (iii) in Lemma~\ref{pou}, and \eqref{4e12},
we deduce that
\begin{align}
\label{z9112}
\sum_{j\in\mathbb{N}}\left|L^{m,k+1}_{i,j}\xi_{j,k+1}(x)\right|
&=\sum_{j\in\mathbb{N}}\left|L^{m,k+1}_{i,j}\xi_{j,k+1}(x)
\mathbf{1}_{\lambda'B_{j,k+1}}(x)\right|
\nonumber\\
&\lesssim2^k\sum_{j\in\mathbb{N}}
\xi_{j,k+1}(x)\mathbf{1}_{\lambda B_{i,k}}(x)
\nonumber\\
&\sim2^k\mathbf{1}_{\Omega_{k+1}}(x)\mathbf{1}_{\lambda B_{i,k}}(x)
\sim2^k\mathbf{1}_{\lambda B_{i,k}}(x),
\end{align}
which further implies that \eqref{4e10} holds.
Moreover, by \eqref{z9112} and (ii) and (iii) of Lemma~\ref{Wtd},
we conclude that
\begin{align}
\label{z911}
\sum_{i\in\mathbb{N}}\sum_{j\in\mathbb{N}}\left|L^{m,k+1}_{i,j}
\xi_{j,k+1}(x)\right|\lesssim 2^k
\sum_{i\in\mathbb{N}}\mathbf{1}_{\lambda B_{i,k}}(x)
\lesssim 2^k\mathbf{1}_{\Omega_k}(x).
\end{align}

Now, we prove that
\begin{equation}\label{4e9}
\sum_{i\in\mathbb{N}}\sum_{j\in\mathbb{N}}L^{m,k+1}_{i,j}\xi_{j,k+1}=0\quad
\end{equation}
pointwise and in $({\mathcal{G}}^\varepsilon_0(\beta,\gamma))'$.
From \eqref{1911},
Lemmas~\ref{pou}(ii) and~\ref{Wtd}(ii), $\lambda'<\lambda$,
and $\Omega_{k+1}\subset\Omega_k$,
we infer that
\begin{align}
\label{zvp-23}
\mathrm{\,supp\,} b^{(m)}_{j,k+1}\subset\mathrm{\,supp\,}\xi_{j,k+1}\subset
B_\rho(x_{j,k+1},\lambda'r_{j,k+1})\subset\Omega_{k+1}
\subset\Omega_{k}.
\end{align}
By \eqref{z911}, Fubini's theorem, \eqref{1916},
(ii) and (iii) of Lemma \ref{pou}, $b^{(m)}_{j,k+1}\in L^1(X)$, the Lebesgue
dominated convergence theorem, \eqref{zvp-23}, and
\eqref{4e3-1}, we find that, for any $x\in X$,
\begin{align*}
\sum_{i\in\mathbb{N}}\sum_{j\in\mathbb{N}}
L^{m,k+1}_{i,j}\xi_{j,k+1}(x)
&=\sum_{j\in\mathbb{N}}\frac{\xi_{j,k+1}(x)}{\|\xi_{j,k+1}\|_{L^1(X)}}
\sum_{i\in\mathbb{N}}\int_{X}b^{(m)}_{j,k+1}(y)\,\xi_{i,k}(y)\,d\mu(y)
\nonumber\\
&=\sum_{j\in\mathbb{N}}\frac{\xi_{j,k+1}(x)}{\|\xi_{j,k+1}\|_{L^1(X)}}
\int_{X}b^{(m)}_{j,k+1}(y)\sum_{i\in\mathbb{N}}\xi_{i,k}(y)\,d\mu(y)
\notag\\
&=\sum_{j\in\mathbb{N}}\frac{\xi_{j,k+1}(x)}{\|\xi_{j,k+1}\|_{L^1(X)}}
\int_{X}b^{(m)}_{j,k+1}(y)\,\mathbf{1}_{\Omega_k}(y)\,d\mu(y)
\notag\\
&=\sum_{j\in\mathbb{N}}\frac{\xi_{j,k+1}(x)}{\|\xi_{j,k+1}\|_{L^1(X)}}\int_{X}
b^{(m)}_{j,k+1}(y)\,d\mu(y)
=0,
\end{align*}
which further implies that \eqref{4e9} holds pointwise.
From this, \eqref{z911}, $\mu(\Omega_k)<\infty$,
and the Lebesgue dominated convergence theorem,
we deduce that \eqref{4e9} also holds
in $({\mathcal{G}}^\varepsilon_0(\beta,\gamma))'$,
which completes the proof of \eqref{4e9}.

By \eqref{zvp-23}, Lemma~\ref{lemma5.14iv}(iii), and \eqref{4e9},
we conclude that,
for any $x\in X$,
\begin{align}
\label{4e11}
g^{(m)}_{k+1}(x)-g^{(m)}_k(x)
&=\left[f^{(m)}(x)-\sum_{j\in\mathbb{N}}b^{(m)}_{j,k+1}(x)\right]
-\left[f^{(m)}(x)-\sum_{i\in\mathbb{N}}b^{(m)}_{i,k}(x)\right]
\nonumber\\
&=\sum_{i\in\mathbb{N}}b^{(m)}_{i,k}(x)
-\sum_{j\in\mathbb{N}}b^{(m)}_{j,k+1}(x)
\nonumber\\
&=\sum_{i\in\mathbb{N}}\left\{b^{(m)}_{i,k}(x)
-\sum_{j\in\mathbb{N}}\left[b^{(m)}_{j,k+1}(x)\xi_{i,k}(x)
-L^{m,k+1}_{i,j}\xi_{j,k+1}(x)\right]\right\}
\nonumber\\
&=:\sum_{i\in\mathbb{N}}h^{(m)}_{i,k}(x),
\end{align}
where all the series converge in
$({\mathcal{G}}^\varepsilon_0(\beta,\gamma))'$ and pointwise.

Next, we claim that,
for any $i\in\mathbb{N}$,
$[\mu(\lambda B_{i,k})]^{-\frac{1}{p}}2^{-k}h^{(m)}_{i,k}$ is
a harmless constant multiple of a $(p,\infty)$-atom.
To this end, let $i\in\mathbb{N}$ and
we first show that
$[\mu(\lambda B_{i,k})]^{-\frac{1}{p}}2^{-k}h^{(m)}_{i,k}$
satisfies Definition~\ref{atom}(i).
Note that Lemma~\ref{pou}(ii) and \eqref{1911} imply
$$
\mathrm{\,supp\,} b^{(m)}_{i,k}\subset\mathrm{\,supp\,}\xi_{i,k}\subset \lambda'B_{i,k},
$$
which, together with $\lambda'<\lambda$
and \eqref{4e10}, implies that
\begin{equation}\label{4e13}
\mathrm{\,supp\,} h^{(m)}_{i,k}\subset\lambda B_{i,k}.
\end{equation}
This further implies that $[\mu(\lambda B_{i,k})]^{-\frac{1}{p}}2^{-k}h^{(m)}_{i,k}$
satisfies Definition~\ref{atom}(i).

Now, we prove that
$[\mu(\lambda B_{i,k})]^{-\frac{1}{p}}2^{-k}h^{(m)}_{i,k}$
satisfies Definition~\ref{atom}(ii).
From \eqref{1911} and Lemma~\ref{pou}(iii),
it follows that,
for any $x\in X$,
\begin{align}
\label{dyw-284}
h^{(m)}_{i,k}(x)
&=b^{(m)}_{i,k}(x)-\sum_{j\in\mathbb{N}}
\left[b^{(m)}_{j,k+1}(x)\xi_{i,k}(x)
-L^{m,k+1}_{i,j}\xi_{j,k+1}(x)\right]
\nonumber\\
&=\left[f^{(m)}(x)\xi_{i,k}(x)
-P_{i,k}^{(m)}\xi_{i,k}(x)\right]
-\sum_{j\in\mathbb{N}}f^{(m)}(x)\xi_{j,k+1}(x)\xi_{i,k}(x)
\nonumber\\
&\quad+\sum_{j\in\mathbb{N}}
\left[P^{(m)}_{j,k+1}\xi_{j,k+1}(x)\xi_{i,k}(x)
+L^{m,k+1}_{i,j}\xi_{j,k+1}(x)\right]
\nonumber\\
&=f^{(m)}(x)\xi_{i,k}(x)\mathbf{1}_{
\left(\Omega_{k+1}\right)^{\complement}}(x)
-P^{(m)}_{i,k}\xi_{i,k}(x)
\nonumber\\
&\quad+\xi_{i,k}(x)\sum_{j\in\mathbb{N}}P^{(m)}_{j,k+1}\xi_{j,k+1}(x)
+\sum_{j\in\mathbb{N}}L^{m,k+1}_{i,j}\xi_{j,k+1}(x),
\end{align}
which, combined with \eqref{1453}, \eqref{4e4-2}, \eqref{4e12},
and both (ii) and (iii) of Lemma~\ref{pou},
further implies that
\begin{align}\label{4e15}
\left\|h^{(m)}_{i,k}\right\|_{L^{\infty}(X)}
\lesssim2^k\left\|\xi_{i,k}\right\|_{L^{\infty}(X)}
\left[2+\left\|\sum_{j\in\mathbb{N}}\xi_{j,k+1}\right\|_{L^{\infty}(X)}\right]
+2^k\left\|\sum_{j\in\mathbb{N}}\xi_{j,k+1}\right\|_{L^{\infty}(X)}
\lesssim2^k,
\end{align}
where the implicit positive constants are independent of $m$, $i$, and $k$.
Thus, a harmless constant multiple of
$[\mu(\lambda B_{i,k})]^{-\frac{1}{p}}2^{-k}h^{(m)}_{i,k}$
satisfies Definition~\ref{atom}(ii).

Next, we show that
$[\mu(\lambda B_{i,k})]^{-\frac{1}{p}}2^{-k}h^{(m)}_{i,k}$
satisfies Definition~\ref{atom}(iii).
Using \eqref{4e3-1}, the Lebesgue dominated convergence theorem, \eqref{1916},
and \eqref{1911},
we find that
\begin{align}
\label{912}
&\int_{X}h^{(m)}_{i,k}(x)\, d\mu(x)
\nonumber\\
&\quad=\int_X\sum_{j\in\mathbb{N}}
\left[L^{m,k+1}_{i,j}\xi_{j,k+1}(x)
-b^{(m)}_{j,k+1}(x)\xi_{i,k}(x)\right]\,d\mu(x)
\nonumber\\
&\quad=\sum_{j\in\mathbb{N}}\Bigg\{\int_{X}\left[f^{(m)}(x)-
P^{(m)}_{j,k+1}\right]\xi_{i,k}(x)\xi_{j,k+1}(x)\,d\mu(x)
\nonumber\\
&\qquad-\int_Xb^{(m)}_{j,k+1}(x)\xi_{i,k}(x)\,d\mu(x)\Bigg\}
\nonumber\\
&\quad=0,
\end{align}
which, together with \eqref{4e13} and \eqref{4e15},
completes the proof of the claim that a harmless
constant multiple of $[\mu(\lambda B_{i,k})]^{-\frac{1}{p}}2^{-k}h^{(m)}_{i,k}$
is a $(p,\infty)$-atom supported in $\lambda B_{i,k}$.
From this claim, \eqref{4e7}, and \eqref{4e11}, we deduce that \eqref{4e2} holds.
This finishes the proof of Step 1.

\emph{Step 2.}
In this step, we prove that there exist a sequence
$\{m_l\}_{l\in\mathbb{N}}$ in $\mathbb{N}$ and a
sequence of functions
$\{h_{i,k}\}_{k\in\mathbb{Z},\,i\in\mathbb{N}}$ such that,
for any $k\in\mathbb{Z}$, $i\in\mathbb{N}$, and $g\in L^1(X)$,
\begin{align*}
\lim_{l\to\infty}\left\langle h_{i,k}^{(m_l)},g\right\rangle
=\left\langle h_{i,k},g\right\rangle,
\end{align*}
where, for any $m\in\mathbb{N}$,
$\{h^{(m)}_{i,k}\}_{k\in\mathbb{Z},\,i\in\mathbb{N}}$
is the family of functions from Step 1.
Moreover, we further show that we can
write $h_{i,k}=\lambda_{i,k}a_{i,k}$,
where each $a_{i,k}$ is a $(p,\infty)$-atom and
$\{\lambda_{i,k}\}_{k\in\mathbb{Z},\,i\in\mathbb{N}}$ in $\mathbb{R}_+$
satisfies
\begin{align*}
\sum_{k=-\infty}^{\infty}
\left(\sum_{i\in\mathbb{N}}\lambda_{i,k}^p\right)^{\frac{q}{p}}
\lesssim\|f\|^q_{H^{p,q}_*(X)}.
\end{align*}
To this end, let $k\in\mathbb{Z}$ and $i\in\mathbb{N}$.
By the boundedness of $\{\|h_{i,k}^{(m)}\|_{L^\infty(X)}\}_{m\in\mathbb{N}}$
[see \eqref{4e15}] and the Banach--Alaoglu theorem
(see, for instance, \cite[Theorem 3.17]{r91}),
we conclude that
there exist a sequence $\{n_{i,k,l}\}_{l\in\mathbb{N}}$ in $\mathbb{N}$
satisfying that $n_{i,k,l}\to\infty$ as $l\to\infty$
and a function ${h}_{i,k}\in L^\infty(X)$ such that,
for any $g\in L^1(X)$,
\begin{align*}
\lim_{l\to\infty}\left\langle h_{i,k}^{(n_{i,k,l})},g\right\rangle
=\left\langle{h}_{i,k},g\right\rangle.
\end{align*}
From this and the diagonal rule,
respectively,
for $k\in\mathbb{Z}$ and $i\in\mathbb{N}$,
we infer that there exists a sequence
$\{m_l\}_{l\in\mathbb{N}}$ in $\mathbb{N}$
satisfying $\lim_{l\to\infty}m_l=\infty$
and, for any $g\in L^1(X)$,
\begin{equation}\label{4e15.5}
\lim_{l\to\infty}\left\langle h^{(m_l)}_{i,k},g
\right\rangle=\langle h_{i,k},g\rangle.
\end{equation}
Now, we prove that
$[\mu(\lambda B_{i,k})]^{-\frac{1}{p}}2^{-k}h_{i,k}$
is a harmless constant multiple of a $(p,\infty)$-atom.
Indeed, using \eqref{4e15.5} and \eqref{4e13},
we obtain
\begin{equation}\label{hik-support}
\mathrm{\,supp\,} h_{i,k}\subset\lambda B_{i,k}
\end{equation}
and, according to Lemma~\ref{Wtd}(ii), for any
$j_0\in\mathbb{Z}\setminus\mathbb{N}$, $k\in\mathbb{Z}$,
and $x\in X$,
$$
\sum_{i\in\mathbb{N}}\mathbf{1}_{(2C_\rho)^{j_0}\lambda B_{i,k}}(x)\leq
\sum_{i\in\mathbb{N}}\mathbf{1}_{\lambda B_{i,k}}(x)\leq L_0.
$$
By \eqref{4e15.5} and \eqref{4e15},
we find that
\begin{align}\label{hik-size}
\left\|h_{i,k}\right\|_{L^\infty(X)}\lesssim 2^k.
\end{align}
From \eqref{hik-support}, \eqref{4e15.5},
\eqref{4e13}, and \eqref{912},
we deduce that
\begin{align*}
\int_{X}h_{i,k}(x)\,d\mu(x)
&=\int_{X}h_{i,k}(x)\mathbf{1}_{\lambda B_{i,k}}(x)\,d\mu(x)
\\
&=\lim_{l\to\infty}\int_{X}h^{(m_l)}_{i,k}(x)
\mathbf{1}_{\lambda B_{i,k}}(x)\,d\mu(x)=0,
\end{align*}
which, combined with \eqref{hik-support} and \eqref{hik-size},
further implies that
$[\mu(\lambda B_{i,k})]^{-\frac{1}{p}}2^{-k}h_{i,k}$
is a harmless constant multiple of $(p,\infty)$-atom.
As such, we can
write
\begin{align*}
h_{i,k}=\lambda_{i,k}a_{i,k},
\end{align*}
where the function $a_{i,k}$ is a $(p,\infty)$-atom supported in $\lambda B_{i,k}$,
the coefficient $\lambda_{i,k}\in\mathbb{R}_+$ satisfies
\begin{align}\label{9121}
\lambda_{i,k}\sim 2^k[\mu(\lambda B_{i,k})]^\frac{1}{p},
\end{align}
and, for any $j_0\in\mathbb{Z}\setminus\mathbb{N}$, $k\in\mathbb{Z}$,
and $x\in X$,
$$
\sum_{i\in\mathbb{N}}\mathbf{1}_{(2C_\rho)^{j_0}\lambda B_{i,k}}(x)\leq L_0.
$$
By this,
we conclude that,
for any given $q\in(0,\infty)$,
\begin{align}
\label{4e16}
\sum_{k=-\infty}^{\infty}
\left(\sum_{i\in\mathbb{N}}\lambda_{i,k}^p\right)^{\frac{q}{p}}
&\sim\sum_{k=-\infty}^{\infty}2^{kq}
\left[\sum_{i\in\mathbb{N}}\mu(\lambda B_{i,k})\right]^{\frac{q}{p}}
\quad\text{by \eqref{9121}}
\nonumber\\
&=\sum_{k=-\infty}^{\infty}2^{kq}
\left[\int_X\sum_{i\in\mathbb{N}}\mathbf{1}_{
\lambda B_{i,k}}(x)\,d\mu(x)\right]^{\frac{q}{p}}
\quad\text{by Tonelli's theorem}
\nonumber\\
&\lesssim\sum_{k=-\infty}^{\infty}2^{kq}
\left[\mu(\Omega_k)\right]^{\frac{q}{p}}
\quad\text{by (ii) and (iii) of Lemma~\ref{Wtd}}
\nonumber\\
&\sim\left\|f^*_\rho\right\|^q_{L^{p,q}(X)}
=\|f\|^q_{H^{p,q}_*{(X)}}
\quad\text{by Proposition~\ref{Lpqnorm=c,d}},
\end{align}
which is the desired estimate when $q\in (0,\infty)$.
When $q=\infty$, the proof is similar to
that used in the estimation of \eqref{4e16}
with the natural modifications;
we omit the details here.
This finishes the proof of Step 2.

\emph{Step 3.}
In this step, we show that
\begin{equation}\label{4e16.5}
f=\sum_{k\in\mathbb{Z}}\sum_{i\in\mathbb{N}}h_{i,k}
\end{equation}
in $({\mathcal{G}}^\varepsilon_0(\beta,\gamma))'$.
To this end, from \eqref{4e13}, \eqref{4e15},
and Lemma~\ref{Wtd}(ii),
we infer that, for any $m\in\mathbb{N}$, $k\in\mathbb{Z}$, and $x\in X$,
\begin{align}
\label{qf1-3}
\sum_{i\in\mathbb{N}}\left|h^{(m)}_{i,k}(x)\right|
\lesssim
2^k\sum_{i\in\mathbb{N}}\mathbf{1}_{\lambda B_{i,k}}(x)<\infty,
\end{align}
which, together with Lemma~\ref{CGL}, implies
that, for any $\phi\in{\mathcal{G}}^\varepsilon_0(\beta,\gamma)$,
\begin{align}
\label{qf1}
\sum_{i\in\mathbb{N}}\left|\left\langle h^{(m)}_{i,k},\phi\right\rangle\right|
\lesssim
2^k\int_{X}\sum_{i\in\mathbb{N}}\mathbf{1}_{\lambda B_{i,k}}(x)|\phi(x)|\,d\mu(x)
\lesssim 2^k\int_{X}|\phi(x)|\,d\mu(x)<\infty.
\end{align}
Similarly, by \eqref{hik-support}, \eqref{hik-size},
Lemmas~\ref{Wtd}(ii) and~\ref{CGL},
and an argument similar to that used
in the estimations of \eqref{qf1-3} and \eqref{qf1},
we find that, for any $k\in\mathbb{Z}$, $x\in X$,
and $\phi\in{\mathcal{G}}^\varepsilon_0(\beta,\gamma)$,
\begin{align}
\label{qf2-3}
\sum_{i\in\mathbb{N}}\left|h_{i,k}(x)\right|
\lesssim
2^k\sum_{i\in\mathbb{N}}\mathbf{1}_{\lambda B_{i,k}}(x)<\infty
\end{align}
and
\begin{align}
\label{qf2}
\sum_{i\in\mathbb{N}}\left|\left\langle h_{i,k},\phi\right\rangle\right|
\lesssim 2^k\int_{X}|\phi(x)|\,d\mu(x)<\infty.
\end{align}
As such, for any $m\in\mathbb{N}$ and $k\in\mathbb{Z}$, we can define
\begin{align}
\label{Def-fmk}
f_k:=\sum_{i\in\mathbb{N}}h_{i,k}
\ \ \text{and}\ \
f^{(m)}_k:=\sum_{i\in\mathbb{N}}h^{(m)}_{i,k},
\end{align}
where the convergences hold pointwise and in
$({\mathcal{G}}^\varepsilon_0(\beta,\gamma))'$.

We first prove that, for any $k\in\mathbb{Z}$
and $\phi\in{\mathcal{G}}^\varepsilon_0(\beta,\gamma)$,
\begin{align}
\label{qf3}
\lim_{l\to\infty}\left\langle f^{(m_l)}_k,\phi\right\rangle
=\left\langle f_k,\phi\right\rangle.
\end{align}
Indeed, from \eqref{qf1-3}, \eqref{qf1}, \eqref{qf2-3}, and \eqref{qf2},
the Lebesgue dominated convergence theorem,
and \eqref{4e15.5},
we deduce that, for any $k\in\mathbb{Z}$ and
$\phi\in{\mathcal{G}}^\varepsilon_0(\beta,\gamma)$,
\begin{equation*}
\langle f^{(m_l)}_k,\phi\rangle=
\left\langle\sum_{i\in\mathbb{N}}h^{(m_l)}_{i,k},\phi\right\rangle
=\sum_{i\in\mathbb{N}}\left\langle h^{(m_l)}_{i,k},\phi\right\rangle
\to\sum_{i\in\mathbb{N}}\left\langle h_{i,k},\phi\right\rangle
=\left\langle\sum_{i\in\mathbb{N}}h_{i,k},\phi\right\rangle
=\langle f_k,\phi\rangle
\end{equation*}
as $l\to\infty$,
which further implies that \eqref{qf3} holds.

Next, we show that, for any
$\phi\in{\mathcal{G}}^\varepsilon_0(\beta,\gamma)$,
\begin{equation}\label{4e17}
\lim_{N\to\infty}\left\langle\sum_{|k|\geq N}f^{(m)}_k,\phi\right\rangle=0
\end{equation}
uniformly in $m\in\mathbb{N}$.
To this end, we consider the following two cases for $p$.

\emph{Case (1)} $p\in(\frac{n}{n+\mathrm{ind\,}(X,\mathbf{q})},1)$.
In this case, to prove \eqref{4e17}, it suffices to show that
\begin{equation}\label{4e18}
\lim_{N_1\to-\infty}\sup_{m\in\mathbb{N}}
\left\|\sum_{k\leq N_1}f^{(m)}_k \right\|_{H^1_*(X)}
=0\ \ \mbox{and}\ \
\lim_{N_2\to\infty}\sup_{m\in\mathbb{N}}
\left\|\sum_{k\geq N_2}f^{(m)}_k\right\|_{H^{p_0}_*(X)}
=0,
\end{equation}
where $p_0\in(\frac{n}{n+\mathrm{ind\,}(X,\mathbf{q})},p)$
is any fixed number. Indeed, if \eqref{4e18} holds, then,
by Proposition~\ref{5.1.10}, we conclude that
\eqref{4e17} holds.

Now, we prove the first equality of \eqref{4e18}.
From \eqref{4e13}, \eqref{4e15}, and \eqref{912},
we infer that,
for any $k\in\mathbb{Z}$
and $i,m\in\mathbb{N}$,
$2^{-k}[\mu(\lambda B_{i,k})]^{-1}h^{(m)}_{i,k}$
is a harmless constant multiple of
a $(1,\infty)$-atom,
which
further implies that, for any given $N_1\in\mathbb{N}$
and for any $M,N\in\mathbb{N}$,
\begin{align}\label{4120}
&\left\|\sum_{k=N_1-M}^{N_1}f^{(m)}_k
-\sum_{k=N_1-M-N}^{N_1}f^{(m)}_k\right\|_{H^1_*(X)}\nonumber\\
&\quad\leq\sum_{k=-\infty}^{k=N_1-M-1}\sum_{i\in\mathbb{N}}
\left\|h_{i,k}^{(m)}\right\|_{H^1_*(X)}
\quad\text{by \eqref{Def-fmk}}\nonumber\\
&\quad\lesssim\sum_{K=-\infty}^{k=N_1-M-1}
2^{k}\sum_{i\in\mathbb{N}}\mu(\lambda B_{i,k})
\quad\text{by \eqref{4l3.2}}\nonumber\\
&\quad\lesssim\sum_{K=-\infty}^{k=N_1-M-1}2^{k}\mu(\Omega_k)
\quad\text{by Lemma~\ref{Wtd}}
\nonumber\\
&\quad\lesssim\sum_{K=-\infty}^{k=N_1-M-1}
2^{k(1-p)}
\left\|f^*_\rho\right\|^p_{L^{p,\infty}(X)}
\quad\text{by Proposition~\ref{Lpqnorm=c,d}}\nonumber\\
&\quad\lesssim2^{(N_1-M)(1-p)}\|f\|^p_{H^{p,q}_*(X)}\to0
\quad\text{by $p<1$ and Proposition~\ref{p,q<p,s}}
\end{align}
as $M\to\infty$.
Thus,
$\{\sum_{k=N_1-M}^{N_1}f_k^{(m)}\}_{M\in\mathbb{N}}$
is a Cauchy sequence in $H^1_*(X)$. Hence,
by Proposition~\ref{5.1.11}, we find that
the series $\sum_{k=-\infty}^{N_1}f_k^{(m)}$
for any given $N_1\in\mathbb{N}$
converges in $H^1_*(X)$.
From this and an argument similar to
that used in the estimation of \eqref{4120},
we deduce that
\begin{align}\label{22062}
\left\|\sum_{k\leq N_1}f^{(m)}_k \right\|_{H^1_*(X)}
\lesssim2^{N_1(1-p)}\|f\|^p_{H^{p,q}_*(X)}\to0
\end{align}
(uniformly in $m\in\mathbb{N}$)
as $N_1\to-\infty$.

Next, we show the second equality of \eqref{4e18}.
Similarly, note that, for any $k\in\mathbb{Z}$ and $m,i\in\mathbb{N}$,
$2^{-k}[\mu(\lambda B_{i,k})]^{-\frac{1}{p_0}}h^{(m)}_{i,k}$
is a harmless constant multiple of
a $(p_0,\infty)$-atom,
which further implies that, for any given $N_2\in\mathbb{N}$ and for any $M,N\in\mathbb{N}$,
\begin{align}\label{21572}
&\left\|\sum_{k=N_2}^{N_2+M}f^{(m)}_k
-\sum_{k=N_2}^{N_2+M+N}f^{(m)}_k\right\|^{p_0}_{H^{p_0}_*(X)}\nonumber\\
&\quad\leq\sum_{k=N_2+M+1}^\infty\sum_{i\in\mathbb{N}}
\left\|h^{(m)}_k\right\|^{p_0}_{H^{p_0}_*(X)}
\quad\text{by \eqref{Def-fmk}}\nonumber\\
&\quad\lesssim\sum_{k=N_2+M+1}^\infty
\sum_{i\in\mathbb{N}}2^{kp_0}\mu(\lambda B_{i,k})
\quad\text{by \eqref{4l3.2}}\nonumber\\
&\quad\lesssim\sum_{k=N_2+M+1}^\infty2^{kp_0}\mu(\Omega_k)
\quad\text{by Lemma~\ref{Wtd}}
\nonumber\\
&\quad\lesssim\sum_{k=N_2+M+1}^\infty
2^{k(p_0-p)}\left\|f^*_\rho\right\|^p_{L^{p,\infty}(X)}
\quad\text{by Proposition~\ref{Lpqnorm=c,d}}\nonumber\\
&\quad\lesssim2^{(N_2+M)(p_0-p)}\|f\|^p_{H^{p,q}_*(X)}\to0
\quad\text{by $p_0<p$ and Proposition~\ref{p,q<p,s}}
\end{align}
as $M\to\infty$.
Thus,
$\{\sum_{k=N_2}^{N_2+M}f_k^{(m)}\}_{M\in\mathbb{N}}$
is a Cauchy sequence in $H^{p_0}_*(X)$.
Hence, by Proposition~\ref{5.1.11}, we conclude that
the series $\sum_{k=N_2}^\infty f_k^{(m)}$ for any given $N_2\in\mathbb{N}$
converges in $H^{p_0}_*(X)$.
From this and an argument similar to
that used in the estimation of \eqref{21572},
we infer that
\begin{align*}
\left\|\sum_{k\ge N_2}f^{(m)}_k \right\|_{H^1_*(X)}
\lesssim2^{N_2(p_0-p)}\|f\|^p_{H^{p,q}_*(X)}\to0
\end{align*}
(uniformly in $m\in\mathbb{N}$)
as $N_2\to\infty$,
which, combined with \eqref{22062}, then
completes the proof of \eqref{4e18}, and hence \eqref{4e17} holds in this case.

\emph{Case (2)} $p=1$.
In this case, to prove \eqref{4e17}, it suffices to show that
\begin{equation}\label{4e18-5}
\lim_{N_1\to-\infty}\sup_{m\in\mathbb{N}}
\left\|\sum_{k\leq N_1}f^{(m)}_k \right\|
_{L^{\widetilde{p}}(X)}=0
\ \ \mbox{and}\ \
\lim_{N_2\to\infty}\sup_{m\in\mathbb{N}}
\left\|\sum_{k\geq N_2}f^{(m)}_k\right\|_{H^{p_0}_*(X)}=0,
\end{equation}
where $\widetilde{p}\in(p,\infty)$ and
$p_0\in(\frac{n}{n+\mathrm{ind\,}(X,\mathbf{q})},1)$ are any fixed numbers.
By an argument similar to that used in the proof of Case (1),
we find that, for any given $N_1\in\mathbb{N}$,
the series $\sum_{k\leq N_1}f^{(m)}_k$
converges in $L^{\widetilde{p}}(X)$.
Then, for any $m\in\mathbb{N}$,
\begin{align*}
\left\|\sum_{k\leq N_1}f^{(m)}_k\right\|_{L^{\widetilde{p}}(X)}
&\leq\sum_{k\leq N_1}
\left\|\sum_{i\in\mathbb{N}}h_{i,k}^{(m)}\right\|_{L^{\widetilde{p}}(X)}
\quad\text{by \eqref{Def-fmk}}\\
&\leq\sum_{k\leq N_1}2^k
\left\|\sum_{i\in\mathbb{N}}\mathbf{1}_{\lambda
B_{i,k}}\right\|_{L^{\widetilde{p}}(X)}
\quad\text{by \eqref{qf1-3}}\\
&\lesssim\sum_{k\leq N_1}2^k\left[\mu(\Omega_k)
\right]^\frac{1}{\widetilde{p}}
\quad\text{by Lemma~\ref{Wtd}}\\
&\lesssim\sum_{k\leq N_1}2^{k(1-\frac{p}{\widetilde{p}})}
\left\|f^*_\rho\right\|^{\frac{p}{\widetilde{p}}}_{L^{p,\infty}(X)}
\quad\text{by Proposition~\ref{Lpqnorm=c,d}}\\
&\lesssim2^{N_1(1-\frac{p}{\widetilde{p}})}
\|f\|^{\frac{p}{\widetilde{p}}}_{H^{p,q}_*(X)}\to0
\quad\text{by $\widetilde{p}>p$ and Propositions~\ref{p,q<p,s}}
\end{align*}
as $N_1\to-\infty$,
where the implicit positive constants are independent of $m$.
This finishes the proof of
the first equality in \eqref{4e18-5}. The proof of
the second one in \eqref{4e18-5} is similar
to that in the case
$p\in(\frac{n}{n+\mathrm{ind\,}(X,\mathbf{q})},1)$.
Thus, \eqref{4e17} also holds in this case,
which finally completes the proof of \eqref{4e17}.

From an argument similar to that used in the proof
of \eqref{4e17} with $f_k^{(m)}$ and $h_{i,k}^{(m)}$
therein replaced, respectively, by
$f_k$ and $h_{i,k}$, we deduce that \eqref{4e17}
with $f_k^{(m)}$ replaced by $f_k$
still holds;
moreover, the proof of \eqref{4e18} implies that,
if $p\in(\frac{n}{n+\mathrm{ind\,}(X,\mathbf{q})},1)$, then
$$
f=\sum_{k=-\infty}^0f_k+\sum_{k=1}^\infty f_k
\in H^1_*(X)+H^{p_0}_*(X)\subset(\mathcal{G}^\varepsilon_0(\beta,\gamma))'
$$
and, if $p=1$, then
$$
f=\sum_{k=-\infty}^0f_k+\sum_{k=1}^\infty f_k
\in L^{\widetilde{p}}(X)+H^{p_0}_*(X)
\subset(\mathcal{G}^\varepsilon_0(\beta,\gamma))'.
$$
Thus, $f=\sum_{k\in\mathbb{Z}}f_k$ converges in
$(\mathcal{G}^\varepsilon_0(\beta,\gamma))'$.
These, together with \eqref{4e17}, imply that,
for any $\delta\in (0,\infty)$ and
$\phi\in{\mathcal{G}}^\varepsilon_0(\beta,\gamma)$,
there exists $K_0\in\mathbb{N}$,
depending on $\phi$ and $\varepsilon$ but
independent of $m$, such that
\begin{equation}\label{4e19}
\left|\left\langle\sum_{|k|>K_0}
f^{(m_l)}_k,\phi\right\rangle\right|
<\frac{\delta}{3}\ \mathrm{and}\
\left|\left\langle\sum_{|k|>K_0}f_k,\phi\right
\rangle\right|<\frac{\delta}{3}.
\end{equation}
Note that, by the fact that
$f^{(m_l)}_k\to f_k$ in
$({\mathcal{G}}^\varepsilon_0(\beta,\gamma))'$ as $l\to\infty$
for any $k\in\mathbb{Z}$, we conclude that there
exists $l_0\in\mathbb{N}$ such that,
for any $l\in\mathbb{N}\cap(l_0,\infty)$
and $k\in\mathbb{Z}$ with $|k|\in[0,K_0]$,
\begin{equation*}
|\langle f^{(m_l)}_k-f_k,\phi\rangle|<\frac{\delta}{6K_0+3},
\end{equation*}
which, combined with \eqref{4e19},
further implies that, for any $l\in\mathbb{N}\cap(l_0,\infty)$,
\begin{align*}
&\left|\left\langle\sum_{k\in\mathbb{Z}}f^{(m_l)}_k,\phi\right\rangle
-\left\langle\sum_{k\in\mathbb{Z}}f_k,\phi\right\rangle\right|
\\
&\quad\leq\left|\left\langle\sum_{|k|>K_0}
f^{(m_l)}_k,\phi\right\rangle\right|+
\left|\left\langle\sum_{|k|>K_0}f_k,\phi\right\rangle\right|+
\left|\sum_{|k|\leq K_0}\left\langle
f^{(m_l)}_k-f_k,\phi\right\rangle\right|<\delta.
\end{align*}
From this and the arbitrariness of $\delta$,
we infer that, for any $\phi\in{\mathcal{G}}^\varepsilon_0(\beta,\gamma)$,
\begin{align}\label{2034}
\lim_{l\to\infty}\left\langle\sum_{k\in\mathbb{Z}}f^{(m_l)}_k,\phi\right\rangle
=\left\langle\sum_{k\in\mathbb{Z}}f_k,\phi\right\rangle.
\end{align}
By this,
we find that, for any $\phi\in{\mathcal{G}}^\varepsilon_0(\beta,\gamma)$,
\begin{align*}
\langle f,\phi\rangle
&=\lim_{l\to\infty}\left\langle f^{(m_l)},\phi\right\rangle
\quad\text{by Lemma~\ref{96}}\\
&=\lim_{l\to\infty}\left\langle\sum_{k\in\mathbb{Z}}
\sum_{i\in\mathbb{N}}h_{i,k}^{(m_l)},\phi\right\rangle
\quad\text{by \eqref{4e2}}\\
&=\lim_{l\to\infty}\left\langle\sum_{k\in\mathbb{Z}}
\left[g_{k+1}^{(m_l)}-g_{k}^{(m_l)}\right],\phi\right\rangle
\quad\text{by \eqref{4e11}}\\
&=\lim_{l\to\infty}\left\langle\sum_{k\in\mathbb{Z}}f^{(m_l)}_k,\phi\right\rangle
\quad\text{by \eqref{Def-fmk}}\\
&=\left\langle\sum_{k\in\mathbb{Z}}f_k,\phi\right\rangle
=\left\langle\sum_{k\in\mathbb{Z}}\sum_{i\in\mathbb{N}}h_{i,k},\phi\right\rangle
\quad\text{by \eqref{2034} and \eqref{Def-fmk}}.
\end{align*}
Using this and the arbitrariness of
$\phi\in{\mathcal{G}}^\varepsilon_0(\beta,\gamma)$,
we complete the proof of \eqref{4e16.5}
and hence Step 3.

According to the above three steps, we have proven that
$$
f=\sum_{k\in\mathbb{Z}}\sum_{i\in\mathbb{N}}\lambda_{i,k}a_{i,k}
$$
in $({\mathcal{G}}^\varepsilon_0(\beta,\gamma))'$,
where, for any $k\in\mathbb{Z}$ and $i\in\mathbb{N}$,
$a_{i,k}$ is a $(p,\infty)$-atom and the sequence
$\{\lambda_{i,k}\}_{k\in\mathbb{Z},\,i\in\mathbb{N}}$ in $\mathbb{R}_+$
satisfies
\begin{align*}
\sum_{k\in\mathbb{Z}}
\left(\sum_{i\in\mathbb{N}}\lambda_{i,k}^p\right)^{\frac{q}{p}}
\lesssim\|f\|^q_{H^{p,q}_*{(X)}}.
\end{align*}
Thus, $f\in H^{p,q,\infty}_{\mathrm{at}}(X)$,
which completes the proof of Theorem~\ref{6.16} in the case where $\mu(X)=\infty$.

Now, we are turn to the case $\mu(X)<\infty$.
In this case, we may assume that $\mu(X)=1$.
Let $R:=\mathrm{diam\,}(X)$. Then $R\in(0,\infty)$ by Remark~\ref{rm1112}(ii).
Let $j_0\in\mathbb{Z}$ be the smallest integer such that
$2^{\frac{j_0p}{q}}>\|f\|_{H^{p,q}_*(X)}$.
Then we have, for any $j\in[j_0,\infty)\cap\mathbb{N}$,
$\Omega_j\subsetneqq X$; otherwise, it holds that $\Omega_j=X$ and hence
\begin{align*}
2^{\frac{j_0p}{q}}&>\|f\|_{H^{p,q}_*(X)}
\ge p^{\frac{1}{q}}\left\{\int_0^{2^j}\alpha^{q-1}
\left[d_{f^*}(\alpha)\right]^{\frac{q}{p}}
\,d\alpha\right\}^{\frac{1}{q}}
=2^{\frac{jp}{q}}\ge2^{\frac{j_0p}{q}},
\end{align*}
which contradicts the definition of $j_0$. Then, by an argument similar
to that used in the proof of the present theorem when
$\mu(X)=\infty$ for $j\in[j_0,\infty)\cap\mathbb{N}$, we obtain,
for any given $m\in\mathbb{Z}$,
there exist $\{h_{j,k}^{(m)}\}_{j\in[j_0,\infty)\cap\mathbb{N},
k\in\mathbb{Z}}$ and $g_{j_0}^{(m)}$ such that
\begin{align*}
f^{(m)}=\sum_{j=j_0}^\infty\sum_{k\in\mathbb{Z}}h_{j,k}^{(m)}+g_{j_0}^{(m)}
\end{align*}
in $({\mathcal{G}}^\varepsilon_0(\beta,\gamma))'$,
where $\|g_{j_0}^{(m)}\|_{L^\infty(X)}\leq C2^{j_0}$
and, for any $j\in[j_0,\infty)\cap\mathbb{N}$ and $k\in\mathbb{Z}$,
$$
h_{j,k}^{(m)}=\lambda_{j,k}^{(m)}a_{j,k}^{(m)}
$$
with both $a_{j,k}^{(m)}$ being a $(p,\infty)$-atom and
$\{\lambda_{j,k}^{(m)}\}_{j\in[j_0,\infty)\cap\mathbb{N},
k\in\mathbb{Z}}$ in $\mathbb{R}_+$
satisfying
\begin{align*}
\sum_{k\in\mathbb{Z}}
\left[\sum_{j=j_0}^\infty\left(\lambda_{j,k}^{(m)}\right)^p\right]^{\frac{q}{p}}
\lesssim\|f\|^q_{H^{p,q}_*{(X)}}.
\end{align*}
This, together with an argument similar to that used in the proof
of Cases (2) and (3) when $\mu(X)=\infty$, further implies that
there exist $\{h_{j,k}\}_{j\in[j_0,\infty)\cap\mathbb{N},
k\in\mathbb{Z}}$ and $g_{j_0}$ such that
\begin{align}\label{1312}
f=\sum_{j=j_0}^\infty\sum_{k\in\mathbb{Z}}h_{j,k}+g_{j_0}
\end{align}
in $({\mathcal{G}}^\varepsilon_0(\beta,\gamma))'$,
where $\|g_{j_0}\|_{L^\infty(X)}\leq C2^{j_0}$
and, for any $j\in[j_0,\infty)\cap\mathbb{N}$ and $k\in\mathbb{Z}$,
$h_{j,k}=\lambda_{j,k}a_{j,k}$
with both $a_{j,k}$ being a $(p,\infty)$-atom and
$\{\lambda_{j,k}\}_{j\in[j_0,\infty)\cap\mathbb{N},
k\in\mathbb{Z}}$ in $\mathbb{R}_+$
satisfying
\begin{align}\label{1313}
\sum_{k\in\mathbb{Z}}
\left(\sum_{j=j_0}^\infty\lambda_{j,k}^p\right)^{\frac{q}{p}}
\lesssim\|f\|^q_{H^{p,q}_*{(X)}}.
\end{align}
We write
\begin{align*}
g_{j_0}=C2^{j_0}\left[\mu(B_\rho(x_0,R+1))\right]^\frac{1}{p}
\frac{g_{j_0}}{C2^{j_0}[\mu(B_\rho(x_0,R+1))]^\frac{1}{p}}=:\lambda_0a_0,
\end{align*}
where $x_0\in X$ is any fixed point.
Then it is easy to show that $\lambda_0=C2^{j_0}\sim\|f\|_{H^{p,q}_*{(X)}}$
[here, we used $\mu(B_\rho(x_0,R+1))=\mu(X)=1$ and the choice of $j_0$]
and that $a_0$ is a local $(p,q)$-atom as in Definition~\ref{localatom}.
Combining this, \eqref{1312}, and \eqref{1313},
we conclude that
\begin{align*}
f\in h^{p,q,\infty}_{\mathrm{at}}(X)
\ \ \text{and}\ \
\|f\|_{h^{p,q,\infty}_{\mathrm{at}}(X)}\lesssim\|f\|_{H^{p,q}_*{(X)}},
\end{align*}
which, together with Theorem~\ref{1314},
further implies that
\begin{align*}
f\in H^{p,q,\infty}_{\mathrm{at}}(X)
\ \ \text{and}\ \
\|f\|_{H^{p,q,\infty}_{\mathrm{at}}(X)}
\sim\|f\|_{h^{p,q,\infty}_{\mathrm{at}}(X)}\lesssim\|f\|_{H^{p,q}_*{(X)}}.
\end{align*}
This finishes the proof of Theorem~\ref{6.16} in the case where $\mu(X)<\infty$,
which completes the proof of Theorem~\ref{6.16}.
\end{proof}

\subsection{Finite Atomic characterization}\label{5.2.2}

In this subsection, we establish
the finite atomic characterization of $H^{p,q}_*(X)$.
To be precise, we prove that,
for any given finite linear
combination of $(p,r)$-atoms
[or continuous $(p,\infty)$-atoms],
its $H^{p,q}_*(X)$ (quasi-)norm
can be achieved
via all its finite atomic decompositions (see Theorem~\ref{6.29}).
When establishing the boundedness of linear operators on Hardy spaces, one
typically resorts to its (infinite) atomic characterization since estimates on the
linear operator can, in principle, be reduced to studying the action of the
operator on individual atoms. However, in light of an example presented in
\cite{bow05} (which is based on a result due to Y.~Meyer in \cite{MTW85})
one must exercise great care when interchanging a linear operator and an
infinite linear combination of atoms.
The advantage of these finite decompositions is
that one does not need to be concerned about
the convergence of the summation. This perspective
will be useful in Section~\ref{section5.8}.

We begin with
the following definition of finite
atomic Hardy--Lorentz spaces.

\begin{definition}
\label{5d1}
Let $(X,\mathbf{q},\mu)$ be a quasi-ultrametric space
of homogeneous type with upper dimension $n\in(0,\infty)$,
$\rho\in\mathbf{q}$ have the property that
all $\rho$-balls are $\mu$-measurable,
$p\in(0,1]$,
$q\in(0,\infty]$, $r\in(p,\infty]\cap[1,\infty]$,
and $\varepsilon,\beta,\gamma\in(0,\infty)$ be such that
$n(\frac{1}{p}-1)<\beta,\gamma<\varepsilon<\infty$.
The \emph{finite atomic Hardy--Lorentz
space}\index[words]{finite atomic Hardy--Lorentz space}
$H^{p,q,r}_{\mathrm{at},\,\mathrm{fin}}(X,\rho)$
\index[symbols]{H@$H^{p,q,r}_{\mathrm{at},\,\mathrm{fin}}(X,\rho)$}
is defined to be
the set of all
$f\in({\mathcal{G}}^\varepsilon_0(\beta,\gamma))'$
such that
there exist $W,K\in\mathbb{N}$, a sequence
$\{\lambda_{i,k}\}_{i\in\mathbb{N}\cap[1,W],\,k\in\mathbb{N}\cap[1,K]}$ in $\mathbb{R}_+$,
and a sequence of $(\rho,p,r)$-atoms
$\{a_{i,k}\}_{i\in\mathbb{N}\cap[1,W],\,k\in\mathbb{N}\cap[1,K]}$ supported,
respectively, in $\rho$-balls
$\{B_{i,k}\}_{i\in\mathbb{N}\cap[1,W],\,k\in\mathbb{N}\cap[1,K]}$
such that
\begin{align}\label{finATdec}
f=\sum_{k=1}^K\sum_{i=1}^W\lambda_{i,k}a_{i,k}
\end{align}
in $({\mathcal{G}}^\varepsilon_0(\beta,\gamma))'$
and the following conditions hold:
\begin{enumerate}
\item[\rm(i)]
There exist constants $C_0\in\mathbb{N}$ and $j_0\in\mathbb{Z}\setminus\mathbb{N}$,
both of which are independent of $f$,
such that, for any $k\in\mathbb{N}\cap[1,K]$ and $x\in X$,
$$
\sum_{i=1}^W\mathbf{1}_{(2C_\rho)^{j_0}B_{i,k}}(x)\leq C_0.
$$

\item[\rm(ii)]
For any $i\in\mathbb{N}\cap[1,W]$ and $k\in\mathbb{N}\cap[1,K]$,
$\lambda_{i,k}\sim2^k[\mu(B_{i,k})]^\frac{1}{p}$
with the positive equivalence constants independent of $f$ and
\begin{align*}
\sum_{k=1}^K
\left(\sum_{i=1}^W\lambda_{i,k}^p\right)^\frac{q}{p}<\infty
\end{align*}
with the usual modification
when $q=\infty$.
\end{enumerate}
Moreover, for any $f\in H^{p,q,r}_{\mathrm{at},\,\mathrm{fin}}(X,\rho)$,
define
\begin{equation*}
\|f\|_{H^{p,q,r}_{\mathrm{at},\,\mathrm{fin}}(X,\rho)}:=\inf\left\{\left[\sum_{k=1}^K
\left(\sum_{i=1}^W\lambda_{i,k}^p\right)^\frac{q}{p}\right]^\frac{1}{q}:
\ f=\sum_{k=1}^K\sum_{i=1}^W\lambda_{i,k}a_{i,k}\right\}
\end{equation*}
with the usual modification
when $q=\infty$, where the infimum is
taken over all decompositions of $f$ as in \eqref{finATdec}.
\end{definition}

\begin{remark}
\label{5216-fin}
\begin{enumerate}
\item [\rm(i)]
Let the notation be as in Definition~\ref{5d1}.
By the assumption $n(\frac{1}{p}-1)<\beta,\gamma<\infty$,
\eqref{1550}, and Remark~\ref{Rmtest}(v), we find that,
for any $p\in(0,\frac{n}{n+\mathrm{ind}_H(X,\mathbf{q})}]$,
$H^{p,q,r}_{\mathrm{at,\,fin}}(X,\rho)$ is trivial because the space
${\mathcal{G}}^\varepsilon_0(\beta,\gamma)$ only contains constant functions.

\item[\rm(ii)]
Following Remark~\ref{5216}(iii), if $\rho'\in\mathbf{q}$ has the property that
all $\rho'$-balls are $\mu$-measurable, then
$f\in H^{p,q,r}_{\mathrm{at},\,\mathrm{fin}}(X,\rho)$ if and only if
$f\in H^{p,q,r}_{\mathrm{at},\,\mathrm{fin}}(X,\rho')$,
with equivalent quasi-norms, where
the choice of constants in conditions
(i) and (ii) in Definition~\ref{6.3} depend on
$\rho$, $\rho'$, and $\mu$. In this sense,
$H^{p,q,r}_{\mathrm{at},\,\mathrm{fin}}(X,\rho)$
is independent of the choice of $\rho\in\mathbf{q}$,
and hence, in what follows, we sometimes simply write
$H^{p,q,r}_{\mathrm{at},\,\mathrm{fin}}(X)$ instead of
$H^{p,q,r}_{\mathrm{at},\,\mathrm{fin}}(X,\rho)$.
\end{enumerate}
\end{remark}

Using Theorem~\ref{6.14}, we obtain
the following density theorem.

\begin{theorem}
\label{5r1}
Let $(X,\mathbf{q},\mu)$ be a quasi-ultrametric space
of homogeneous type with upper dimension $n\in(0,\infty)$,
$\rho\in\mathbf{q}$ have the property that
all $\rho$-balls are $\mu$-measurable,
$p\in(\frac{n}{n+\mathrm{ind\,}(X,\mathbf{q})},1]$,
$q\in(0,\infty)$, $r\in(p,\infty]\cap[1,\infty]$,
and
$n(\frac{1}{p}-1)<\beta,\gamma<\varepsilon\preceq\mathrm{ind\,}(X,\mathbf{q})$.
Then the space
$H^{p,q,r}_{\mathrm{at},\,\mathrm{fin}}(X,\rho)$ is dense in
$H^{p,q}_{*}(X)$ with respect to the quasi-norm
$\|\cdot\|_{H^{p,q}_{*}(X)}$.
\end{theorem}

\begin{proof}
Let $f\in H^{p,q}_{*}(X)$.
By Theorems~\ref{6.16} and~\ref{6.17},
we find that there exist a sequence
$\{\lambda_{i,k}\}_{i\in\mathbb{N},k\in\mathbb{Z}}$ in $\mathbb{R}_+$
and a sequence of $(\rho,p,r)$-atoms
$\{a_{i,k}\}_{i\in\mathbb{N},k\in\mathbb{Z}}$
supported, respectively, in
$\{B_{i,k}\}_{i\in\mathbb{N},k\in\mathbb{Z}}$
such that $f=\sum_{i\in\mathbb{N}}\sum_{k\in\mathbb{Z}}
\lambda_{i,k}a_{i,k}$ in $(\mathcal{G}^\varepsilon_0(\beta,\gamma))'$
and
\begin{align}\label{915}
\left[\sum_{k\in\mathbb{Z}}\left(\sum_{i\in\mathbb{N}}
\lambda_{i,k}^p\right)^\frac{q}{p}\right]^\frac{1}{q}
\lesssim\|f\|_{H^{p,q}_{*}(X)}<\infty,
\end{align}
where the positive equivalence constants are
independent of $f$.
For any $N\in\mathbb{N}$, let
$$
f_N:=\sum_{i=1}^N\sum_{k=-N}^N\lambda_{i,k}a_{i,k}.
$$
Then $f_N\in H^{p,q,r}_{\mathrm{at,fin}}(X,\rho)$
for any given $N\in\mathbb{N}$. From \eqref{915},
we deduce that, for any given $\delta\in(0,\infty)$,
there exists $N_{(\delta)}\in\mathbb{N}$ such that
$$
\left[\sum_{k=1}^{N_{\delta}}\left(
\sum_{i>N_{(\delta)}}\lambda_{i,k}^p\right)^\frac{q}{p}
+\sum_{|k|>N_{(\delta)}}\left(\sum_{i\in\mathbb{Z}}
\lambda_{i,k}^p\right)^\frac{q}{p}\right]^\frac{1}{q}<\delta,
$$
which, combined with Theorem~\ref{6.15},
further implies that
\begin{align*}
\left\|f-f_{N_{(\delta)}}\right\|_{H^{p,q}_*(X)}
&\lesssim\left\|f-f_{N_{(\delta)}}\right\|_{H^{p,q,r}_{\mathrm{at}}(X)}\\
&\leq\left[\sum_{k=1}^{N_{\delta}}\left(
\sum_{i>N_{(\delta)}}\lambda_{i,k}^p\right)^\frac{q}{p}
+\sum_{|k|>N_{(\delta)}}\left(\sum_{i\in\mathbb{Z}}
\lambda_{i,k}^p\right)^\frac{q}{p}\right]^\frac{1}{q}<\delta.
\end{align*}
This finishes the proof of Theorem~\ref{5r1}.
\end{proof}

\begin{definition}\label{unifcont}
Let $(X,\mathbf{q})$ be a quasi-ultrametric space and
$\rho\in\mathbf{q}$.
The \emph{space
of uniformly continuous functions}
$\mathrm{UC}(X,\rho)$
on $X$\index[symbols]{U@$\mathrm{UC}(X)$}
is defined to be the set of
all functions $f:X\to\mathbb{C}$ satisfying that,
for any given $\eta\in(0,\infty)$,
there exists $\delta\in(0,\infty)$
such that,
for any $x,y\in X$ with $\rho(x,y)<\delta$,
$$
|f(x)-f(y)|<\eta.
$$
\end{definition}

\begin{remark}\label{9362}
In the context of Definition~\ref{unifcont}, it
is easy to see that, for any $\beta\in(0,\infty)$,
$$
\dot{C}^\beta(X)\subset\mathrm{UC}(X,\rho)
$$
and, for any $\rho'\in\mathbf{q}$,
$$
\mathrm{UC}(X,\rho)=\mathrm{UC}(X,\rho').
$$
Thus, we sometimes simply write
$\mathrm{UC}(X)$ instead of $\mathrm{UC}(X,\rho)$.
\end{remark}

Next, we show the following density lemma.

\begin{lemma}
\label{lemma81}
Let $(X,\mathbf{q},\mu)$ be a quasi-ultrametric space
of homogeneous type with upper dimension $n\in(0,\infty)$.
Assume that $p\in(\frac{n}{n+\mathrm{ind\,}(X,\mathbf{q})},1]$,
$q\in(0,\infty)$,
and $n(\frac{1}{p}-1)<\beta,\gamma<\varepsilon\preceq\mathrm{ind\,}(X,\mathbf{q})$.
Then the following statements hold:
\begin{enumerate}
\item [\textup{(i)}]
For any $s\in[1,\infty]$, $H^{p,q}_*(X)\cap L^s(X)$
is dense in $H^{p,q}_*(X)$ with respect to the quasi-norm
$\|\cdot\|_{H^{p,q}_{*}(X)}$.

\item [\textup{(ii)}]
$H^{p,q}_*(X)\cap\mathrm{UC}(X)$
is dense in $H^{p,q}_*(X)$ with respect to the quasi-norm
$\|\cdot\|_{H^{p,q}_{*}(X)}$.
\end{enumerate}
\end{lemma}

\begin{proof}
Let $\rho\in\mathbf{q}$ be a regularized quasi-ultrametric as in Definition~\ref{D2.3}.
We first prove (i).
Using Definition~\ref{5d1} and a straightforward
calculation, we conclude that,
for any $s\in[1,\infty]$,
$H^{p,q,\infty}_{\mathrm{at},\,\mathrm{fin}}(X,\rho)\subset L^s(X)$.
From this and Theorem~\ref{5r1}, we deduce that
$$
H^{p,q,\infty}_{\mathrm{at},\,\mathrm{fin}}(X,\rho)
=H^{p,q,\infty}_{\mathrm{at},\,\mathrm{fin}}(X,\rho)\cap L^s(X)
\subset H^{p,q}_*(X)\cap L^s(X)
\subset H^{p,q}_*(X),
$$
which, together with Theorem~\ref{5r1},
further implies that
$$
H^{p,q}_*(X)=\overline{H^{p,q,\infty}_{\mathrm{at},\,
\mathrm{fin}}(X,\rho)}^{\|\cdot\|_{H^{p,q}_*(X)}}
\subset\overline{H^{p,q}_*(X)\cap L^s(X)}^{\|\cdot\|_{H^{p,q}_*(X)}}
\subset H^{p,q}_*(X).
$$
This finishes the proof of (i).

Now, we show (ii).
We first prove that every $(p,\infty)$-atom
can be approximated arbitrarily well in the quasi-norm
$\|\cdot\|_{H^{p,q}_{*}(X)}$ by functions in $\mathrm{UC}(X)$.
With this goal in mind,
let $a\in L^\infty(X)$ be a $(p,\infty)$-atom
supported in $B_\rho(x_0,r_0)$
with $x_0\in X$ and $r_0\in(0,\infty)$,
and let $\{\mathcal{S}_{2^{-k}}\}_{k\in\mathbb{Z}}$ be
a homogeneous approximation of the identity
of order $\varepsilon$ as in Definition~\ref{DefHATI}.
By Theorem~\ref{corATI}(viii) and $\mathrm{\,supp\,} a\subset B_\rho(x_0,r_0)$,
we find that there exists a constant $C\in[1,\infty)$ such that,
for any $k\in\mathbb{Z}$,
$\mathrm{\,supp\,} \mathcal{S}_{2^{-k}}a\subset
B_\rho(x_0,C(r_0+2^{-k}))$.
Choose $k_0\in\mathbb{Z}$ large enough such that,
for any $k\in[k_0,\infty)\cap\mathbb{Z}$,
\begin{align}\label{810}
\mathrm{\,supp\,}\mathcal{S}_{2^{-k}}a\subset B_\rho(x_0,2C r_0).
\end{align}
From Theorem~\ref{corATI}(i),
Definition~\ref{atom}(ii), and \eqref{doublingcond},
we infer that, for any $k\in[k_0,\infty)\cap\mathbb{Z}$,
\begin{align}
\label{811}
\left\|\mathcal{S}_{2^{-k}}a\right\|_{L^\infty(X)}
\lesssim\left\|a\right\|_{L^\infty(X)}
\lesssim
\left[(V_\rho)_{r_0}(x_0)\right]^{-\frac{1}{p}}
\lesssim\left[(V_\rho)_{2Cr_0}(x_0)\right]^{-\frac{1}{p}}.
\end{align}
In particular, for any $k\in[k_0,\infty)\cap\mathbb{Z}$,
$\mathcal{S}_{2^{-k}}a\in L^1(X)$.
By this, Fubini's theorem, and Definitions~\ref{DefATI}(v)
and~\ref{atom}(iii),
we conclude that, for any $k\in[k_0,\infty)\cap\mathbb{Z}$,
\begin{align}
\label{811-x}
\int_X\left(\mathcal{S}_{2^{-k}}a\right)(x)\,d\mu(x)
&=\int_X\int_XS_{2^{-k}}(x,y)\,d\mu(x)a(y)\,d\mu(y)
\nonumber\\
&=\int_Xa(x)\,d\mu(x)=0,
\end{align}
which, combined with \eqref{810} and \eqref{811}, implies that,
for any $k\in[k_0,\infty)\cap\mathbb{Z}$,
$\mathcal{S}_{2^{-k}}a$ is  a harmless constant multiple of a
$(p,\infty)$-atom
supported in $B_\rho(x_0,2C r_0)$, where the harmless constant is independent
of $k$. From this, Theorem~\ref{6.14}, Theorem~\ref{corATI}(vi),
and Remark~\ref{9362}, we deduce that, for any  $k\in[k_0,\infty)\cap\mathbb{Z}$,
$\mathcal{S}_{2^{-k}}a\in H^{p,q}_*(X)\cap\mathrm{UC}(X)$.

In a step towards proving (ii), we first show that,
for any $k\in[k_0,\infty)\cap\mathbb{Z}$,
\begin{align}
\label{813-X}
\left\|\mathcal{S}_{2^{-k}}a-a\right\|_{H^{p,q}_*(X)}
\lesssim\frac{\|\mathcal{S}_{2^{-k}}a-a\|_{L^2(X)}}
{[(V_\rho)_{2Cr_0}(x_0)]^{\frac{1}{2}-\frac{1}{p}}}.
\end{align}
If $\|\mathcal{S}_{2^{-k}}a-a\|_{L^2(X)}=0$,
then \eqref{813-X} is trivial and hence we can
assume $\|\mathcal{S}_{2^{-k}}a-a\|_{L^2(X)}>0$.
Observe that \eqref{810}, $\mathrm{supp\,}a\subset B_\rho(x_0,r_0)$,
and \eqref{811-x}
immediately imply that, for any  $k\in[k_0,\infty)\cap\mathbb{Z}$,
$\mathrm{\,supp\,} (\mathcal{S}_{2^{-k}}a-a)\subset B_\rho(x_0,2Cr_0)$
and
\begin{align*}
\int_X\left[\mathcal{S}_{2^{-k}}a(x)-a(x)\right]\,d\mu(x)=0.
\end{align*}
By these observations, it is easy to verify that,
for any $k\in[k_0,\infty)\cap\mathbb{Z}$,
$$
\left[(V_\rho)_{2Cr_0}(x_0)\right]^{\frac{1}{2}-\frac{1}{p}}
\frac{\mathcal{S}_{2^{-k}}a-a}{\|\mathcal{S}_{2^{-k}}a-a\|_{L^2(X)}}
$$
is a $(p,2)$-atom
supported in $B_\rho(x_0,2Cr_0)$.
From this, Theorem~\ref{6.14}, and Definition~\ref{6.3},
we infer that, for any $k\in[k_0,\infty)\cap\mathbb{Z}$,
\begin{align*}
\left\|\mathcal{S}_{2^{-k}}a-a\right\|_{H^{p,q}_*(X)}
\sim\left\|\mathcal{S}_{2^{-k}}a-a\right\|_{H^{p,q,2}_{\mathrm{at}}(X)}
\lesssim\frac{\|\mathcal{S}_{2^{-k}}a-a\|_{L^2(X)}}
{[(V_\rho)_{2C_\rho r_0}(x_0)]^{\frac{1}{2}-\frac{1}{p}}},
\end{align*}
which proves \eqref{813-X}.
Next, by \eqref{813-X} and Theorem~\ref{corATI}(ix), keeping
in mind that we are assuming $\mu$ is Borel-semiregular, we find that
\begin{align}
\label{813}
\left\|\mathcal{S}_{2^{-k}}a-a\right\|_{H^{p,q}_*(X)}
\lesssim\frac{\|\mathcal{S}_{2^{-k}}a-a\|_{L^2(X)}}
{[(V_\rho)_{2C_\rho r_0}(x_0)]^{\frac{1}{2}-\frac{1}{p}}}\to 0
\end{align}
as $k\to\infty$.

To finish the proof of (ii), let $f\in H^{p,q}_*(X)$. Then Theorem~\ref{5r1}
implies that there exists a sequence
$\{f_j\}_{j\in\mathbb{N}}$ in
$H^{p,q,\infty}_{\mathrm{at},\,\mathrm{fin}}(X,\rho)$
such that
\begin{align}\label{814}
\lim_{j\to\infty}\left\|f-f_j\right\|_{H^{p,q}_*(X)}=0.
\end{align}
Note that, for any $j\in\mathbb{N}$,
$f_j$ is a finite linear combination of $(p,\infty)$-atoms.
Consequently, from what we have already established above,
we deduce that,
for any $j\in\mathbb{N}$ and $k\in\mathbb{Z}$,
$\mathcal{S}_{2^{-k}}f_j\in H^{p,q}_*(X)\cap\mathrm{UC}(X)$.
Moreover, \eqref{813} implies that
$$
\lim_{k\to\infty}\left\|
\mathcal{S}_{2^{-k}}f_j-f_j\right\|_{H^{p,q}_*(X)}=0.
$$
By this and \eqref{814}, we conclude that
$$
\left\|
\mathcal{S}_{2^{-k}}f_j-f\right\|_{H^{p,q}_*(X)}
\leq\left\|
\mathcal{S}_{2^{-k}}f_j-f_j\right\|_{H^{p,q}_*(X)}
+\left\|f_j-f\right\|_{H^{p,q}_*(X)}\to0
$$
as $k,j\to\infty$,
which completes the proof of Lemma~\ref{lemma81}.
\end{proof}

The following lemma is a generalization of \cite[Lemma 5.6]{zhy}
on unbounded spaces of homogeneous type satisfying $\mu(\{x\})=0$
for all $x\in X$ to  quasi-ultrametric spaces of homogeneous type.

\begin{lemma}
\label{815}
Let $(X,\mathbf{q},\mu)$ be a quasi-ultrametric space
of homogeneous type with upper dimension $n\in(0,\infty)$,
$\rho\in\mathbf{q}$ have the property that
all $\rho$-balls are $\mu$-measurable,
$p\in(\frac{n}{n+\mathrm{ind\,}(X,\mathbf{q})},1]$,
$q\in(0,\infty]$, $r\in(p,\infty]\cap[1,\infty]$,
and $\varepsilon,\beta,\gamma\in(0,\infty)$ be such that
$n(\frac{1}{p}-1)<\beta,\gamma<
\varepsilon\preceq\mathrm{ind\,}(X,\mathbf{q})$.
Then, for any $f\in H^{p,q,r}_{\mathrm{at},\,\mathrm{fin}}(X,\rho)$,
there exist $x_1\in X$ and $R\in(0,\infty)$ such that
$\mathrm{supp\,} f\subset B_\rho(x_1,R)$.
Moreover, for any $\lambda\in(4C_\rho^6,\infty)$,
there exist positive constants $C$ and $C'$, which are
independent of $f$, $x_1$, and $R$, such that,
for any $x\in[B_\rho(x_1,\lambda R)]^\complement$,
\begin{equation}
\label{xzy.478}
f^*_\rho(x)\leq C\inf_{y\in B_\rho(x,C'\rho(x_1,x))}f^*_\rho(y).
\end{equation}
\end{lemma}

\begin{proof}
Let $f\in H^{p,q,r}_{\mathrm{at},\,\mathrm{fin}}(X,\rho)$.
Note that $f$ is a finite linear combination of $(p,r)$-atoms,
which, together with Definition~\ref{atom}(i),
implies that there exist a point $x_1\in X$ and a radius $R\in(0,\infty)$
such that $\mathrm{supp\,}f\subset B_\rho(x_1,R)$.

To show \eqref{xzy.478}, we split the proof into two parts.
In the first part, we assume that $\rho\in\mathbf{q}$
is a regularized quasi-ultrametric as in Definition~\ref{D2.3}
satisfying $\varepsilon\leq(\log_2C_\rho)^{-1}$.
Let $\lambda\in(4C_\rho^6,\infty)$ and $x\in[B_\rho(x_1,\lambda R)]^\complement$,
and assume that $\phi\in\mathcal{G}^\varepsilon_0(\beta,\gamma)$
with $\|\phi\|_{G(x,r^\alpha,\frac{\beta}{\alpha},
\frac{\gamma}{\alpha})}\leq1$ for some $r\in(0,\infty)$,
where $\alpha:=\min\{1,\,(\log_2C_\rho)^{-1}\}$.
Then Remark~\ref{RemMf}(iii) implies that
$\|\phi\|_{\mathcal{G}(x,r,\beta,\gamma)}\lesssim1$.
Now, we consider the following two cases for $r$.

\emph{Case (1)} $r\ge\rho(x_1,x)$.
In this case, for any given $y\in B_\rho(x,\rho(x_1,x))$
and for any $z\in X$,
we have $r+\rho(x,z)\sim r+\rho(y,z)$,
which further implies that
$\|\phi\|_{\mathcal{G}(y,r,\beta,\gamma)}\sim
\|\phi\|_{\mathcal{G}(x,r,\beta,\gamma)}\lesssim1$.
By this, we find that
\begin{align}\label{1540}
\left|\left\langle f,\phi\right\rangle\right|\lesssim f^*_\rho(y),
\end{align}
where the implicit positive constant is independent of $y$.

\emph{Case (2)} $r<\rho(x_1,x)$.
In this case, from Corollary~\ref{Lem3.3},
we deduce that, for any given $A\in[C_\rho,4C_\rho^4]$,
there exists a function $g\in\dot{C}^\varepsilon(X)$ satisfying that
\begin{align}\label{3b4-sds}
\mathbf{1}_{B_\rho(x_1,\frac{\rho(x_1,x)}{AC_\rho^2})}
\leq g
\leq\mathbf{1}_{B_\rho(x_1,\frac{\rho(x_1,x)}{C_\rho^2})}
\ \ \text{and}\ \
\|g\|_{\dot{C}^\varepsilon(X)}\lesssim\left[\rho(x_1,x)\right]^{-\varepsilon}.
\end{align}
By this,  $\mathrm{supp\,}f\subset B_\rho(x_1,R)$,
and the fact that $\frac{\rho(x_1,x)}{AC_\rho^2}\ge\frac{\lambda R}{4C_\rho^6}>R$,
we conclude that $f=fg$. Let $\widetilde{\phi}:=\phi g$.
We claim that, for any $y\in B_\rho(x,\rho(x_1,x))$,
\begin{align}\label{2517}
\left\|\widetilde{\phi}\right\|_{\mathcal{G}(y,r,\beta,\gamma)}\lesssim1.
\end{align}
Fix $y\in B_\rho(x,\rho(x_1,x))$.
We first consider the size condition of $\widetilde{\phi}$.
If $z\in X$ is such that $\widetilde{\phi}(z)\neq0$, then
$\rho(z,x_1)<\frac{\rho(x_1,x)}{C_\rho^2}$,
which, combined with Definition~\ref{D2.1}(iii),
further implies that
$$
\rho(x_1,x)\leq C_\rho\max\left\{\rho(x_1,z),\,\rho(x,z)\right\}=C_\rho\rho(x,z).
$$
This, together with Definition~\ref{D2.1}(iii)
and $\rho(y,x)<\rho(x_1,x)$,
further implies that
$$
\rho(y,z)\leq C_\rho\max\left\{\rho(y,x),\,\rho(x,z)\right\}
\leq C_\rho\max\left\{\rho(x_1,x),\,\rho(x,z)\right\}
\lesssim\rho(x,z).
$$
From this, the size condition of $\phi$,
and Lemma~\ref{A3.2}(iii), we infer that,
for any $z\in X$ such that $\widetilde{\phi}(z)\neq0$,
\begin{align*}
\left|\widetilde{\phi}(z)\right|
&\leq\left|\phi(z)\right|
\lesssim\frac{1}{(V_\rho)_r(x)+V_\rho(x,z)}
\left[\frac{r}{r+\rho(x,z)}\right]^\gamma\\
&\sim\frac{1}{(V_\rho)_r(y)+V_\rho(y,z)}
\left[\frac{r}{r+\rho(y,z)}\right]^\gamma,
\end{align*}
which is the desired size estimate.
Next, we consider the regularity condition of $\widetilde{\phi}$. By the above
size condition of $\widetilde{\phi}$ and an argument similar to the one
used in the proof of \eqref{1600},
we only need to consider the case where
$z,z'\in X$ such that $\rho(z,z')\leq\frac{r+\rho(y,z)}{4C_\rho^4}$.
We consider the following two subcases for $\rho(z,x_1)$.

\emph{Subcase (1)} $\rho(z,x_1)\ge\frac{\rho(x_1,x)}{C_\rho}$.
In this case, we have
$\rho(z,x_1)\geq\frac{\rho(x_1,x)}{C_\rho^2}$, which, combined with
\eqref{3b4-sds}, implies
$\widetilde{\phi}(z)=\phi(z)g(z)=0$.
Moreover, by Definition~\ref{D2.1}(iii) and the assumptions
$\rho(x,y)<\rho(x_1,x)$
and $r<\rho(x_1,x)$,
we find that
\begin{align*}
\rho(z,z')
&\leq\frac{r+\rho(y,z)}{4C_\rho^4}<
\frac{\rho(x_1,x)+\max\{C_\rho\rho(y,x_1),\,C_\rho\rho(z,x_1)\}}{4C_\rho^4}
\\
&\leq\frac{C_\rho\rho(z,x_1)+
\max\{C_\rho^2\rho(x_1,x),\,C_\rho\rho(z,x_1)\}}{4C_\rho^4}
\leq\frac{\rho(z,x_1)}{C_\rho},
\end{align*}
which, together with Definition~\ref{D2.1}(iii),
further implies
$$
\frac{\rho(x_1,x)}{C_\rho}\leq\rho(z,x_1)\leq C_\rho\max\left\{
\rho(z,z'),\,\rho(z',x_1)\right\}=C_\rho\rho(z',x_1)
$$
and hence
$\widetilde{\phi}(z')=0$.
Thus, $\widetilde{\phi}(z)-\widetilde{\phi}(z')=0$
in this case.

\emph{Subcase (2)} $\rho(z,x_1)<\frac{\rho(x_1,x)}{C_\rho}$.
In this case, Definition~\ref{D2.1}(iii) implies that
\begin{align*}
\rho(x_1,x)
&\leq C_\rho\max\left\{\rho(x_1,z),\,\rho(z,x)\right\}
=C_\rho\rho(z,x)\\
&\leq C_\rho^2\max\left\{\rho(z,x_1),\,\rho(x_1,x)\right\}
=C_\rho^2\rho(x_1,x),
\end{align*}
which, combined with Definition~\ref{D2.1}(iii) and $\rho(y,x_1)<\rho(x_1,x)$,
implies that
\begin{align*}
\rho(y,z)
&\leq C_\rho^2\max\left\{\rho(y,x_1),\,\rho(x_1,x),\,
\rho(x,z)\right\}\\
&=C_\rho^2\max\left\{\rho(x_1,x),\,\rho(x,z)\right\}
\leq C_\rho^3\rho(x,z).
\end{align*}
From this and Definition~\ref{D2.1}(iii) again,
we deduce that
$$
\rho(z,z')\leq\frac{r+\rho(y,z)}{4C_\rho^4}
\leq\frac{r+\rho(x,z)}{2C_\rho}
$$
and
$$
r+\rho(y,z)\lesssim\max\left\{r+\rho(x,z),\,r+\rho(x,z'),\,\rho(x_1,x)\right\}.
$$
By this, both the size and the regularity conditions of $\phi$,
\eqref{3b4-sds}, and $\beta<\varepsilon$,
we conclude that
\begin{align*}
\left|\widetilde{\phi}(z)-\widetilde{\phi}(z')\right|
&\leq\left|\phi(z)-\phi(z')\right||g(z)|
+\left|g(z)-g(z')\right||\phi(z')|\\
&\lesssim\left[\frac{\rho(z,z')}{r+\rho(x,z)}\right]^\beta
\frac{1}{(V_\rho)_r(x)+V_\rho(x,z)}
\left[\frac{r}{r+\rho(x,z)}\right]^\gamma\\
&\quad+\left[\frac{\rho(z,z')}{\rho(x_1,x)}
\right]^\varepsilon\frac{1}{(V_\rho)_r(x)+V_\rho(x,z')}
\left[\frac{r}{r+\rho(x,z')}\right]^\gamma\\
&\lesssim\left[\frac{\rho(z,z')}{r+\rho(y,z)}\right]^\beta
\frac{1}{(V_\rho)_r(y)+V_\rho(y,z)}
\left[\frac{r}{r+\rho(y,z)}\right]^\gamma,
\end{align*}
which is the desired regularity estimate of $\widetilde{\phi}$.
This finishes the proof of \eqref{2517}.

From \eqref{2517} and Proposition~\ref{distfncttype},
we infer that,
for any $y\in B_\rho(x,\rho(x_1,x))$,
\begin{align*}
|\langle f,\phi\rangle|=\left|\int_Xf(x)\phi(x)\,d\mu(x)\right|
=\left|\int_Xf(x)g(x)\phi(x)\,d\mu(x)\right|
=|\langle f,\widetilde{\phi}\rangle|\lesssim f^*_\rho(y).
\end{align*}
Using this and \eqref{1540}
and taking the supremum over all $\phi$ as above
give the desired conclusion.
This finishes the proof of \eqref{xzy.478} (with $C'=1$)
assuming that $\rho\in\mathbf{q}$
is a regularized quasi-ultrametric
satisfying $\varepsilon\leq(\log_2C_\rho)^{-1}$.

There remains to prove \eqref{xzy.478} assuming that
$\rho\in\mathbf{q}$ only has the property that
all $\rho$-balls are $\mu$-measurable.
To this end, by Remark~\ref{RmkD2.3}
and the assumption $\varepsilon\preceq\mathrm{ind\,}(X,\mathbf{q})$,
we find that there exists a
regularized quasi-ultrametric $\varrho\in\mathbf{q}$
satisfying $\varepsilon\leq(\log_2C_{\varrho})^{-1}$.
Now, from what we have already proven to $\varrho$
(with $\lambda=5C_\varrho^6>4C_\varrho^6$),
we deduce that, for any $x\in[B_{\varrho}(x_1,5C_\varrho^6 R_0)]^\complement$,
\begin{equation}
\label{xzy.478-X}
f^*_{\varrho}(x)\lesssim\inf_{y\in B_{\varrho}(x,\varrho(x_1,x))}f^*_{\varrho}(y),
\end{equation}
where $x_1\in X$ and $R_0\in(0,\infty)$ are
such that $\mathrm{supp\,}f\subset B_{\varrho}(x_1,R_0)$.
By $\varrho\in\mathbf{q}$ and \eqref{1557-xx},
we conclude that there exists $R\in(0,\infty)$ such that
$\mathrm{supp\,} f\subset B_{\rho}(x_1,R)$.
Via increasing $R$, if necessary, we may assume
$c^{-1}4C_\rho^6R\geq5C_\varrho^6R_0$,
where $c\in[1,\infty)$ satisfies, for any $x,y\in X$,
\begin{equation}
\label{xzy.478-XX}
c^{-1}\varrho(x,y)\leq\rho(x,y)\leq c\varrho(x,y).
\end{equation}
Let $\lambda\in(4C_\rho^6,\infty)$
and $x\in[B_\rho(x_1,\lambda R)]^\complement$.
Then, from \eqref{xzy.478-XX} and the choices of $\lambda$ and
$R$, we infer that
$$
\varrho(x_1,x)\geq c^{-1}\rho(x_1,x)\geq c^{-1}\lambda R
\geq c^{-1}4C_\rho^6R\geq5C_\varrho^6R_0,
$$
which implies $x\in[B_\varrho(x_1,5C_\varrho^6R_0)]^\complement$.
This, combined with
Remark~\ref{RemMf}(iv), \eqref{xzy.478-X},
\eqref{xzy.478-XX}, and \eqref{1557-xx}, gives
\begin{align*}
f^*_{\rho}(x)\sim
f^*_{\varrho}(x)\lesssim\inf_{y\in B_{\varrho}(x,\varrho(x_1,x))}f^*_{\varrho}(y)
\lesssim\inf_{y\in B_{\rho}(x,C'\rho(x_1,x))}f^*_{\rho}(y),
\end{align*}
where $C':=c^{-2}\in(0,\infty)$. Hence, \eqref{xzy.478}
holds and the proof  of Lemma~\ref{815}
is now complete.
\end{proof}

Next, we list a few conclusions
that follow from the arguments presented in the proof of
Theorem~\ref{6.16}.

\begin{lemma}\label{5.4}
Let $(X,\mathbf{q},\mu)$ be a quasi-ultrametric
space of homogeneous type with upper dimension $n\in(0,\infty)$,
$\rho\in\mathbf{q}$
have the property that all $\rho$-balls are $\mu$-measurable,
$\lambda\in(4C_\rho^6,\infty)$,
$p\in(\frac{n}{n+\mathrm{ind\,}(X,\mathbf{q})},1]$,
$q\in(0,\infty]$, $r\in(p,\infty]\cap[1,\infty]$,
and $\varepsilon,\beta,\gamma\in(0,\infty)$ be such that
$n(\frac{1}{p}-1)<\beta,\gamma<\varepsilon
\preceq\mathrm{ind\,}(X,\mathbf{q})$.
Then, for any $f\in H^{p,q}_*(X)\cap L^r(X)$,
there exist non-negative real numbers $\{\lambda_{i,k}\}_{i\in\mathbb{N},\,
k\in\mathbb{Z}}$,
$\rho$-balls $\{B_{i,k}:=B_\rho(x_{i,k},r_{i,k})\}_
{i\in\mathbb{N},k\in\mathbb{Z}}$
with $\{x_{i,k}\}_
{i\in\mathbb{N},k\in\mathbb{Z}}$ in $X$ and $\{r_{i,k}\}_
{i\in\mathbb{N},k\in\mathbb{Z}}$ in $(0,\infty)$,
and $(\rho,p,\infty)$-atoms $\{a_{i,k}\}_
{i\in\mathbb{N},k\in\mathbb{Z}}$
such that
$$
f=\sum_{k\in\mathbb{Z}}\sum_{i\in\mathbb{N}}
\lambda_{i,k}a_{i,k}
$$
pointwise $\mu$-almost everywhere, in $L^r(X)$,
in $({\mathcal{G}}^\varepsilon_0(\beta,\gamma))'$,
and in $H^{p,q}_\ast(X)$,
where, for any $k\in\mathbb{Z}$ and $i\in\mathbb{N}$,
\begin{align}
\label{suppaik}
\mathrm{\,supp\,} a_{i,k}\subset \lambda B_{i,k}
\subset\Omega_k=\bigcup_{i\in\mathbb{N}}B_{i,k}
\end{align}
with $\Omega_k:=\{x\in X:f^*_\rho(x)>2^{k}\}$.
Moreover, there exists $L_0\in\mathbb{N}$ such that,
for any $k\in\mathbb{Z}$,
\begin{align}\label{9581}
\sum_{i\in\mathbb{N}}\mathbf{1}_{\lambda B_{i,k}}(x)\leq L_0.
\end{align}
Finally, there exists a positive constant $C$,
independent of $f$, such that, for any $k\in\mathbb{Z}$
and $i\in\mathbb{N}$,
\begin{align}
\label{lamdikaik}
\left|\lambda_{i,k}a_{i,k}\right|\leq C2^k
\end{align}
pointwise $\mu$-almost everywhere and
\begin{align*}
\left[\sum_{k\in\mathbb{Z}}
\left(\sum_{i\in\mathbb{N}}
\lambda_{i,k}^p\right)^\frac{q}{p}\right]^\frac{1}{q}
\leq C\|f\|_{H^{p,q}_*(X)}.
\end{align*}
\end{lemma}

\begin{proof}
With the exception of the convergence in $L^r(X)$, all of the conclusions
of the present lemma follow from arguing
as in the proof of Theorem~\ref{6.16}
with $f^{(m)}$ replaced by $f$. To see that
$$
f=\sum_{k\in\mathbb{Z}}\sum_{i\in\mathbb{N}}
\lambda_{i,k}a_{i,k}
$$
converges in $L^r(X)$,
we estimate
\begin{align}
\label{kfp-39}
\sum\limits_{k\in\mathbb{Z}}\sum_{i\in\mathbb{N}}|\lambda_{i,k}a_{i,k}|
&\lesssim\sum\limits_{k\in\mathbb{Z}}
\sum_{i\in\mathbb{N}}2^k{\bf 1}_{B_{\rho}(x_{i,k},\lambda r_{i,k})}
\quad\text{by \eqref{lamdikaik} and \eqref{suppaik}}
\nonumber\\
&\lesssim
\sum\limits_{k\in\mathbb{Z}}2^k{\bf 1}_{\Omega_k}
\quad\text{by \eqref{9581}}\nonumber\\
&=\sum\limits_{k\in\mathbb{Z}}2^k\sum
\limits_{\ell=k}^{\infty}{\bf 1}_{\Omega_\ell\setminus \Omega_{\ell+1}}
=\sum\limits_{\ell\in\mathbb{Z}}2^{\ell}\sum\limits_{k=-\infty}^{\ell}
2^{k-\ell}{\bf 1}_{\Omega_\ell\setminus \Omega_{\ell+1}}
\nonumber\\
&\qquad\ \text{by Tonelli's theorem for sequences}
\nonumber\\
&\sim\sum\limits_{\ell\in\mathbb{Z}}2^{\ell}{\bf 1}_{\Omega_\ell\setminus \Omega_{\ell+1}}
\leq f^*_\rho\sum\limits_{\ell\in\mathbb{Z}}{\bf 1}_{\Omega_\ell\setminus \Omega_{\ell+1}}
=f^*_\rho
\end{align}
pointwise $\mu$-almost everywhere on $X$.
On the other hand, since $f\in L^r(X)$, by
Proposition~\ref{HL-grandmax},
the fact that $L^r(X)=L^{r,r}(X)$, and Lemma~\ref{M},
we find that $f^*_\rho\in L^r(X)$,
which, together with \eqref{kfp-39}, implies that
$$
\sum_{k\in\mathbb{Z}}\sum_{i\in\mathbb{N}}
|\lambda_{i,k}a_{i,k}|\in L^r(X).
$$
As such,
it follows from the Lebesgue dominated convergence theorem that
$$
f=\sum_{k\in\mathbb{Z}}\sum_{i\in\mathbb{N}}
\lambda_{i,k}a_{i,k}
$$
also holds in $L^r(X)$,
which completes the proof of Lemma~\ref{5.4}.
\end{proof}

\begin{remark}\label{5.5}
Let the notation be as in Lemma~\ref{5.4}.
For any $i,j\in\mathbb{N}$
and $k\in\mathbb{Z}$, define
$$
P_{i,k}:=\frac{1}{\|\xi_{i,k}\|_{L^1(X)}}
\int_Xf(y)\xi_{i,k}(y)\,d\mu(y)
$$
and
$$
L^{(k+1)}_{i,j}:=\frac{1}{\|\xi_{i,k}\|_{L^1(X)}}
\int_X\left[f(y)-P_{j,k+1}\right]
\xi_{i,k}(y)\xi_{j,k+1}(y)\,d\mu(y),
$$
where $\xi_{i,k}$ is defined as in Lemma~\ref{pou}.
By Lemma~\ref{pou}(ii),
we find that $\mathrm{supp\,}\xi_{i,k}\subset\lambda' B_{i,k}$,
where $\lambda'\in(C_\rho,\frac{\lambda}{C_\rho^2})$
is the same as in the proof of Theorem~\ref{6.16}.
Replacing $f^{(m)}$, $P^{(m)}_{i,k}$,
and $L^{m,k+1}_{i,j}$,
respectively, by $f$, $P_{i,k}$, and
$L^{(k+1)}_{i,j}$ in the proof of Theorem~\ref{6.16}, it follows from
\eqref{4e7}, \eqref{4e11}, \eqref{dyw-284}, and \eqref{1911} that,
for any $i\in\mathbb{N}$ and $k\in\mathbb{Z}$,
$$
\lambda_{i,k}a_{i,k}=\left(f-P_{i,k}\right)\xi_{i,k}
-\sum_{j\in\mathbb{N}}\left[
\left(f-P_{j,k+1}\right)\xi_{i,k}-L^{(k+1)}_{i,k}\right]\xi_{j,k+1}
$$
$\mu$-almost everywhere.
\end{remark}

Now, we present the main theorem of this section as follows.

\begin{theorem}
\label{6.29}
Let $(X,\mathbf{q},\mu)$ be a quasi-ultrametric space
of homogeneous type with upper dimension $n\in(0,\infty)$,
$\rho\in\mathbf{q}$ have the property that
all $\rho$-balls are $\mu$-measurable,
$p\in(\frac{n}{n+\mathrm{ind\,}(X,\mathbf{q})},1]$,
$q\in(0,\infty]$,
and $\varepsilon,\beta,\gamma\in(0,\infty)$ be such that
$n(\frac{1}{p}-1)<\beta,\gamma<\varepsilon
\preceq\mathrm{ind\,}(X,\mathbf{q})$.
Then the following assertions hold:
\begin{enumerate}
\item [\textup{(i)}]
If $r\in(p,\infty)\cap[1,\infty)$,
then $\|\cdot\|_{H^{p,q,r}_{\mathrm{at},\,\mathrm{fin}}(X,\rho)}$
and $\|\cdot\|_{H^{p,q}_*(X)}$
are equivalent quasi-norms on $H^{p,q,r}_{\mathrm{at},\,\mathrm{fin}}(X,\rho)$.

\item [\textup{(ii)}]
$\|\cdot\|_{H^{p,q,\infty}_{\mathrm{at},\,\mathrm{fin}}(X,\rho)}$
and $\|\cdot\|_{H^{p,q}_*(X)}$
are equivalent quasi-norms on
$H^{p,q,\infty}_{\mathrm{at},\,\mathrm{fin}}(X,\rho)\cap\mathrm{UC}(X)$.
\end{enumerate}
\end{theorem}

\begin{proof}
Based on Remarks~\ref{5216-fin} and \ref{RmkD2.3}(iii),
without loss of generality,
we may assume that $\rho\in\mathbf{q}$ is a regularized quasi-ultrametric
satisfying $\varepsilon\leq(\log_2C_{\rho})^{-1}$.
By the fact that $H^{p,q,r}_{\mathrm{at},\,\mathrm{fin}}(X,\rho)$
continuously embeds into $H^{p,q,r}_{\mathrm{at}}(X)$
and Theorem~\ref{6.14}, we conclude that,
for any $r\in(p,\infty]\cap[1,\infty]$ and
$f\in H^{p,q,r}_{\mathrm{at},\,\mathrm{fin}}(X,\rho)$,
$$
\left\|f\right\|_{H^{p,q}_*(X)}
\sim\left\|f\right\|_{H^{p,q,r}_{\mathrm{at}}(X)}\leq
\left\|f\right\|_{H^{p,q,r}_{\mathrm{at},\,\mathrm{fin}}(X,\rho)}.
$$
Thus, it remains to prove that,
for any $f\in H^{p,q,r}_{\mathrm{at},\,\mathrm{fin}}(X,\rho)$
when $r\in(p,\infty)\cap[1,\infty)$
and for any $f\in H^{p,q,r}_{\mathrm{at},\,\mathrm{fin}}(X,\rho)
\cap\mathrm{UC}(X)$ when $r=\infty$,
$$
\left\|f\right\|_{H^{p,q,r}_{\mathrm{at},\,\mathrm{fin}}(X,\rho)}
\lesssim\left\|f\right\|_{H^{p,q}_*(X)}.
$$

Before proceeding, we make a number of observations and fix some notation.
Let $f\in H^{p,q,r}_{\mathrm{at},\,\mathrm{fin}}(X,\rho)$.
Without loss of generality,
we may assume that $\|f\|_{H^{p,q}_*(X)}=1$.
Let
\begin{align*}
\widetilde{r}:=
\begin{cases}
r\ &\text{if}\ r\in(p,\infty)\cap[1,\infty),\\
2\ &\text{if}\ r=\infty.
\end{cases}
\end{align*}
Then, from Theorem~\ref{5r1} and the bounded support of $f$ (see Lemma~\ref{815}),
we deduce that
$f\in H^{p,q}_*(X)\cap L^{\widetilde{r}}(X)$.
By this and Lemma~\ref{5.4}, we find that, for any $\lambda\in(4C_\rho^6,\infty)$,
there exist a sequence
$\{\lambda_{i,k}\}_{i\in\mathbb{N},k\in\mathbb{Z}}$
in $\mathbb{R}_+$ and a sequence
of $(p,\infty)$-atoms $\{a_{i,k}\}_{i\in\mathbb{N},k\in\mathbb{Z}}$
supported, respectively,
in $\rho$-balls $\{\lambda B_{i,k}\}_{i\in\mathbb{N},k\in\mathbb{Z}}$
such that
\begin{align}
\label{djp-329}
f=\sum_{k\in\mathbb{Z}}\sum_{i\in\mathbb{N}}
\lambda_{i,k}a_{i,k}
\end{align}
both $\mu$-almost everywhere and in
$({\mathcal{G}}^\varepsilon_0(\beta,\gamma))'$. Moreover,
we have
\begin{align}\label{djp-329-2}
\left[\sum_{k\in\mathbb{Z}}
\left(\sum_{i\in\mathbb{N}}
\lambda_{i,k}^p\right)^\frac{q}{p}\right]^\frac{1}{q}
\lesssim\|f\|_{H^{p,q}_*(X)}.
\end{align}
Our goal is to show that the sum in \eqref{djp-329} can be written as a finite
linear combination of $(p,r)$-atoms.
In a step towards this goal,
there exist $x_1\in X$, $R\in(0,\infty)$, and constants $c,C\in(0,\infty)$ (which is
independent of $f$, $x_1$, and $R$) such that
$$
\mathrm{supp\,} f\subset B_\rho(x_1,R)
$$
and, for any $x\in[B_\rho(x_1,\lambda R)]^\complement$ satisfying $f^*_\rho(x)>0$,
\begin{align}\label{816}
1
&=\|f\|_{H^{p,q}_*(X)}
\gtrsim\|f\|_{H^{p,\infty}_*(X)}
=\|f^*_\rho\|_{L^{p,\infty}(X)}
\quad\text{by Proposition~\ref{p,q<p,s}}
\nonumber\\
&=\sup_{t\in(0,\infty)}t\left[
\mu\left(\left\{z\in X:f^*_\rho(z)>t\right\}
\right)\right]^\frac{1}{p}
\quad\text{by Proposition~\ref{Lpqnorm=c,d}}
\nonumber\\
&\ge\frac{f^*_\rho(x)}{2C}\left[
\mu\left(\left\{z\in X:f^*_\rho(z)>\frac{f^*_\rho(x)}{2C}
\right\}\right)\right]^\frac{1}{p}
\nonumber\\
&\ge c^{-1}f^*_\rho(x)\left[
V_\rho(x_1,x)\right]^\frac{1}{p}
\quad\text{by Lemma~\ref{815}}
\nonumber\\
&\ge c^{-1}f^*_\rho(x)\left[
(V_\rho)_R(x_1)\right]^\frac{1}{p}
\quad\text{since $\rho(x_1,x)>R$}.
\end{align}
Let $\widetilde{k}\in\mathbb{Z}$ be such that
\begin{align}
\label{qbk-38}
2^{\widetilde{k}}<c\left[
(V_\rho)_R(x_1)\right]^{-\frac{1}{p}}
\leq 2^{\widetilde{k}+1}.
\end{align}
From this and \eqref{816},
it follows that,
for any $k\in\mathbb{Z}\cap(\widetilde{k},\infty]$,
\begin{align}
\label{omegak}
\Omega_{k}
:=&\left\{x\in X:f^*_\rho(x)>2^{k}\right\}
\subset\left\{x\in X:f^*_\rho(x)>2^{\widetilde{k}+1}\right\}
\nonumber\\
\subset&\left\{x\in X:f^*_\rho(x)>c\left[
(V_\rho)_R(x_1)\right]^{-\frac{1}{p}}\right\}
\subset B_\rho(x_1,\lambda R).
\end{align}
Finally, in the above context, we define
$$
h:=\sum_{k=-\infty}^{\widetilde{k}}\sum_{i\in\mathbb{N}}
\lambda_{i,k}a_{i,k}
\ \ \text{and}\ \
l:=\sum_{k=\widetilde{k}+1}^{\infty}\sum_{i\in\mathbb{N}}
\lambda_{i,k}a_{i,k},
$$
where both sums converge $\mu$-almost everywhere
and in $({\mathcal{G}}^\varepsilon_0(\beta,\gamma))'$.

We are ready to prove (i).
Let $r\in(p,\infty]\cap[1,\infty]$ and
$f\in H^{p,q,r}_{\mathrm{at},\,\mathrm{fin}}(X,\rho)$.
We first claim that
$h$ is a harmless constant multiple of a
$(p,r)$-atom supported in $B_\rho(x_1,\lambda R)$.
Indeed, by Lemma~\ref{5.4},
we conclude that, for any
$k\in\mathbb{Z}$ and $i\in\mathbb{N}$,
$$
\mathrm{\,supp\,} a_{i,k}\subset\lambda B_{i,k}\subset\Omega_k.
$$
From this, the definition of $l$, and \eqref{omegak},
we infer that
\begin{align}
\label{89.00}
\mathrm{\,supp\,} l\subset\bigcup_{k=\widetilde{k}+1}^{\infty}\Omega_k
\subset B_\rho(x_1,\lambda R),
\end{align}
which, combined with the fact that $f=h+l$
pointwise $\mu$-almost everywhere and
$\mathrm{\,supp\,} f\subset B_\rho(x_1,R)\subset B_\rho(x_1,\lambda R)$,
further implies that
\begin{align}
\label{821}
\mathrm{\,supp\,} h\subset B_\rho(x_1,\lambda R).
\end{align}
Note that
\begin{align}
\label{822}
\left\|h\right\|_{L^\infty(X)}
&\leq\sum_{k=-\infty}^{\widetilde{k}}\left\|
\sum_{i\in\mathbb{N}}\left|\lambda_{i,k}a_{i,k}
\right|\right\|_{L^\infty(X)}
\nonumber\\
&\lesssim\sum_{k=-\infty}^{\widetilde{k}}2^k\left\|
\sum_{i\in\mathbb{N}}\mathbf{1}_{\lambda B_{i,k}}\right\|_{L^\infty(X)}
\quad\text{by \eqref{suppaik} and \eqref{lamdikaik}}
\nonumber\\
&\lesssim\sum_{k=-\infty}^{\widetilde{k}}2^k\left\|
\mathbf{1}_{\Omega_{k}}\right\|_{L^\infty(X)}
=2^{\widetilde{k}+1}
\quad\text{by \eqref{suppaik} and \eqref{9581}}
\nonumber\\
&\sim\left[(V_\rho)_R(x_1)\right]^{-\frac{1}{p}}
\quad\text{by \eqref{qbk-38}}
\nonumber\\
&\lesssim\left[(V_\rho)_{\lambda R}(x_1)\right]^{-\frac{1}{p}}
\quad\text{by \eqref{upperdoub}}.
\end{align}
To see that $\int_Xh(x)\,d\mu(x)=0$, by the fact that
$f$ is a finite linear combination of $(p,r)$-atoms, each of which
satisfy Definition~\ref{atom}(iii),
we find that $\int_Xf(x)\,d\mu(x)=0$.
From this and the fact that $h=f-l$,
we deduce that,
to show that $h$ satisfies Definition~\ref{atom}(iii),
it suffices to prove that $l\in L^1(X)$ and satisfies
$\int_Xl(x)\,d\mu(x)=0$. To this end,
we first show
\begin{align*}
\sum_{k=\widetilde{k}+1}^{\infty}\sum_{i\in\mathbb{N}}
\lambda_{i,k}\left|a_{i,k}\right|\in L^1(X).
\end{align*}
Indeed, by H\"older's inequality and the facts that
$f\in H^{p,q,r}_{\mathrm{at},\,\mathrm{fin}}(X,\rho)\subset L^r(X)$ and
$\mathrm{\,supp\,} f\subset B_\rho(x_1,R)$,
we obtain
$$
\|f\|_{L^{1}(X)}\leq
\left\|\mathbf{1}_{B_\rho(x_1,R)}\right\|_{L^{\frac{r}{r-1}}(X)}
\|f\|_{L^{r}(X)}<\infty,
$$
which, together with
an argument similar to that used in the proof of \eqref{kfp-39},
implies that
\begin{align*}
\int_X\sum_{k=\widetilde{k}+1}^{\infty}\sum_{i\in\mathbb{N}}
\lambda_{i,k}\left|a_{i,k}(x)\right|\,d\mu(x)
\lesssim\left\|f^*_\rho\right\|_{L^1(X)}
<\infty.
\end{align*}
From this, the Lebesgue dominated convergence theorem,
and Definition~\ref{atom}(iii), we deduce that
$l\in L^1(X)$ and
\begin{align}
\label{89.0}
\int_Xl(x)\,d\mu(x)
=\sum_{k=\widetilde{k}+1}^{\infty}
\sum_{i\in\mathbb{N}}\lambda_{i,k}\int_X
a_{i,k}(x)\,d\mu(x)=0.
\end{align}
Hence,
\begin{align*}
\int_Xh(x)\,d\mu(x)=\int_X\left[f(x)-l(x)\right]\,d\mu(x)=0.
\end{align*}
By this, \eqref{821}, and \eqref{822}, we conclude that
$h$ is a harmless constant multiple of a
$(p,\infty)$-atom supported in $B_\rho(x_1,\lambda R)$
and hence a harmless constant multiple of a
$(p,r)$-atom supported in $B_\rho(x_1,\lambda R)$.
This finishes the proof of the claim.

Next, we prove that
$l\in L^r(X)$. On the one hand,
since $f\in L^r(X)$, from \eqref{*<M} and Lemma~\ref{M},
we deduce that $f^*_\rho\in L^r(X)$.
On the other hand,
we can argue as in \eqref{kfp-39} to estimate
$$
\sum_{k=\widetilde{k}+1}^\infty\sum_{i\in\mathbb{N}}
\lambda_{i,k}\left|a_{i,k}\right|\lesssim f^*_\rho.
$$
These observations, combined with the
Lebesgue dominated convergence theorem,
imply that $l\in L^r(X)$ and
\begin{align}
\label{apt-235}
\sum_{k=\widetilde{k}+1}^{N}\sum_{i\in\mathbb{N}}
\lambda_{i,k}a_{i,k}\to l
\end{align}
in $L^r(X)$ as $N\to\infty$.
Now, for any given $K\in\mathbb{Z}\cap[\widetilde{k}+1,\infty)$
and for any $k\in\mathbb{Z}\cap[\widetilde{k}+1,K]$, define
$$
I_{(K,k)}:=\left\{i\in\mathbb{N}:
|i|+|k|\leq K\right\}
\ \ \text{and}\ \
l_{(K)}:=\sum_{k=\widetilde{k}+1}^{K}
\sum_{i\in I_{(K,k)}}\lambda_{i,k}a_{i,k}.
$$
For any given $\delta\in(0,1)$,
since $l\in L^r(X)$,
it follows that
\begin{align*}
\frac{l-l_{(K)}}{\delta}\in L^r(X),
\end{align*}
which, together with the definition of $l_{(K)}$ and \eqref{apt-235},
further implies that
\begin{align}\label{89.01}
\lim_{K\to\infty}\left\|\frac{l-l_{(K)}}
{\delta}\right\|_{L^r(X)}=0.
\end{align}
Therefore, by choosing $\delta\in(0,1)$ with $\delta
\leq[(V_\rho)_{\lambda R}(x_1)]^{\frac{1}{r}-\frac{1}{p}}$, it follows from
\eqref{89.01} that there exists $K_0\in[\widetilde{k}+1,\infty)
\cap\mathbb{Z}$ satisfying
\begin{align}
\label{89.01-2}
\left\|l-l_{(K_0)}\right\|_{L^r(X)}\leq\left[(V_\rho)_{\lambda R}(x_1)
\right]^{\frac{1}{r}-\frac{1}{p}}.
\end{align}
On the other hand, by \eqref{89.0} and the fact that $l_{(K_0)}$
is a finite linear combination of atoms,
we find that
\begin{align}
\label{89.2}
\int_X\left[l(x)-l_{(K_0)}(x)\right]\,d\mu(x)=0.
\end{align}
Finally, from \eqref{suppaik}, \eqref{omegak},
\eqref{89.00}, and the definition of $l_{(K_0)}$,
we infer that $l-l_{(K_0)}$
is supported in $B_\rho(x_1,\lambda R)$.
By this, \eqref{89.01-2}, and \eqref{89.2},
we conclude that $l-l_{(K_0)}$ is a $(p,r)$-atom.
Thus,
$$f=h+l_{(K_0)}+\left[l-l_{(K_0)}\right]$$
is a finite linear combination of $(p,r)$-atoms,
which, together with \eqref{djp-329-2},
implies that
\begin{align*}
\|f\|^q_{H^{p,q,r}_{\mathrm{at},\,\mathrm{fin}}(X,\rho)}
\lesssim1+\sum_{k=\widetilde{k}+1}^\infty
\left(\sum_{i\in\mathbb{N}}\lambda_{i,k}^p\right)^\frac{q}{p}
\lesssim\|f\|^q_{H^{p,q}_*(X)}.
\end{align*}
This finishes the proof of (i).

Next, we show (ii).
Let $f\in\mathrm{UC}(X)\cap
H^{p,q,\infty}_{\mathrm{at},\,\mathrm{fin}}(X,\rho)$.
As in the proof of (i), we may assume that
$\|f\|_{H^{p,q}_*(X)}=1$.
From \eqref{*<M} and Lemma~\ref{M},
we deduce that there exists a positive constant $C_1$
such that
\begin{align}\label{8910}
\left\|f^*_\rho\right\|_{L^\infty(X)}\leq
C_1\left\|\mathcal{M}_\rho f\right\|_{L^\infty(X)}
=C_1\|f\|_{L^\infty(X)}<\infty.
\end{align}
Let $K_1\in\mathbb{Z}$ be such that
$$
2^{K_1-1}\leq C_1\|f\|_{L^\infty(X)}<2^{K_1}.
$$
By increasing the value of $C_1$, we may assume that
$K_1>\widetilde{k}$, where $\widetilde{k}\in\mathbb{Z}$
is the same as in \eqref{qbk-38}.
From this, \eqref{suppaik}, and \eqref{8910},
it follows that, for any $k\in\mathbb{Z}\cap[K_1,\infty)$,
$\mu(\Omega_k)=0$, where $\Omega_k$ is the same as in \eqref{omegak}.
Then
\begin{align}
\label{f=h+l}
f
&=\sum_{k=-\infty}^{\infty}
\sum_{i\in\mathbb{N}}\lambda_{i,k}a_{i,k}
=\sum_{k=-\infty}^{K_1}
\sum_{i\in\mathbb{N}}\lambda_{i,k}a_{i,k}
\nonumber\\
&=\sum_{k=-\infty}^{\widetilde{k}}
\sum_{i\in\mathbb{N}}\lambda_{i,k}a_{i,k}
+\sum_{k=\widetilde{k}+1}^{K_1}
\sum_{i\in\mathbb{N}}\lambda_{i,k}a_{i,k}
=:h+l,
\end{align}
where all of the series converge pointwise $\mu$-almost
everywhere and in $({\mathcal{G}}_0^\varepsilon(\beta,\gamma))'$.
By reasoning as in the proof of part (i)
in the present theorem with $r$ replaced by $\infty$, we find that
$h$ is a harmless constant multiple of a
$(p,\infty)$-atom.

There remains to deal with $l$.
Since $f\in\mathrm{UC}(X)$,
it follows that,
for any given $\eta\in(0,\infty)$,
there exists $\delta\in(0,\infty)$ such that,
for any $x,y\in X$ satisfying $\rho(x,y)<\delta$,
$|f(x)-f(y)|<\eta$.
Let
\begin{align}
\label{l=l1+l2}
l=\sum_{k=\widetilde{k}+1}^{K_1}\sum_{i\in G_1(k)}
\lambda_{i,k}a_{i,k}
+\sum_{k=\widetilde{k}+1}^{K_1}\sum_{i\in G_2(k)}
\lambda_{i,k}a_{i,k}=:l_{1,\delta}+l_{2,\delta},
\end{align}
where, for any $k\in\mathbb{Z}$,
$$
G_1(k):=\left\{i\in\mathbb{N}:
\lambda^2 r_{i,k}\ge\delta\right\}
\ \ \text{and}\ \
G_2(k):=\left\{i\in\mathbb{N}:
\lambda^2 r_{i,k}<\delta\right\}.
$$
To deal with $l_{1,\delta}$,
we infer that, for any $k\in[\widetilde{k}+1,K_1)\cap\mathbb{Z}$,
\begin{align*}
\mu\left(B_\rho(x_1,\lambda R)\right)
&\geq\mu(\Omega_k)
\quad\text{by \eqref{omegak}}\\
&\ge
\sum_{i\in G_1(k)}\mu(\lambda B_{i,k})
\quad\text{by \eqref{9581}}\\
&\geq\sum_{i\in G_1(k)}\mu
\left(B_\rho\left(x_{i,k},\frac{\delta}{\lambda^2}\right)\right)
\quad\text{by the definition of $G_1(k)$}\\
&\gtrsim\sum_{i\in G_1(k)}\mu\left(B_\rho(x_{i,k},\lambda R)\right)
\quad\text{by \eqref{doublingcond}}\\
&\sim\#G_1(k)\,\mu\left(B_\rho(x_1,\lambda R)\right)
\quad\text{by Lemma~\ref{A3.2}(vii)},
\end{align*}
which, combined with the fact that
$\mu(B_\rho(x_1,\lambda R))\in(0,\infty)$, further
implies that $G_1(k)$ is a finite set. Hence, $l_{1,\delta}$ is
a finite linear combination of $(p,\infty)$-atoms
and, moreover, by \eqref{djp-329-2}, we obtain
\begin{align}
\label{l1delta}
\left\|l_{1,\delta}\right\|^q_{H^{p,q,\infty}_{\mathrm{at},\,\mathrm{fin}}(X,\rho)}
\leq\sum_{k=\widetilde{k}+1}^{K_1}
\left(\sum_{i\in G_1(k)}\lambda_{i,k}^p\right)^\frac{q}{p}
\lesssim\|f\|_{H^{p,q}_*(X)}^q.
\end{align}

Now, we turn to prove that $l_{2,\delta}$ is a $(p,\infty)$-atom.
From \eqref{suppaik}, \eqref{omegak},
\eqref{l=l1+l2}, \eqref{89.0}, and the fact that
$l_{1,\delta}$ is a finite linear
combination of $(p,\infty)$-atoms, we infer that
\begin{align*}
\mathrm{\,supp\,} l_{2,\delta}\subset
\bigcup_{k=\widetilde{k}+1}^{K_1}\Omega_k
\subset B_\rho(x_1,\lambda R)
\end{align*}
and
\begin{align*}\int_Xl_{2,\delta}(x)\,d\mu(x)
=\int_Xl(x)\,d\mu(x)-\int_Xl_{1,\delta}(x)\,d\mu(x)=0.
\end{align*}
Thus, to finish showing that $l_{2,\delta}$ is
a $(p,\infty)$-atom, it remains to
prove that
$$
\left\|l_{2,\delta}\right\|_{L^\infty(X)}
\leq\left[(V_\rho)_{\lambda R}(x_1)\right]^{-\frac{1}{p}}.
$$
To this end, by Remark~\ref{5.5}, we conclude that,
for any $k\in\mathbb{Z}$ and
$i\in\mathbb{N}$ and for $\mu$-almost every $x\in X$,
\begin{align}
\label{89-0}
\left|\lambda_{i,k}a_{i,k}(x)\right|
&\leq\left|\left[f(x)-P_{i,k}\right]\xi_{i,k}(x)\right|
+\sum_{j\in\mathbb{N}}\left|\left[f(x)-P_{j,k+1}\right]\xi_{i,k}
\xi_{j,k+1}(x)\right|
\nonumber\\
&\quad
+\sum_{j\in\mathbb{N}}
\left|L^{(k+1)}_{i,j}\xi_{j,k+1}(x)\right|
\nonumber\\
&=\left|b_{i,k}(x)\right|
+\sum_{j\in\mathbb{N}}\left|b_{j,k+1}(x)\xi_{i,k}(x)\right|
+\sum_{j\in\mathbb{N}}\left|L^{(k+1)}_{i,j}\xi_{j,k+1}(x)\right|,
\end{align}
where $b_{i,k}(x):=[f(x)-P_{i,k}(x)]\xi_{i,k}(x)$.
From Remark~\ref{5.5},
it follows that,
for any $k\in\mathbb{Z}\cap[\widetilde{k}+1,K_1]$
and $i\in G_2(k)$,
$$
\mathrm{supp\,} b_{i,k}\subset
\lambda'B_{i,k},
$$
where $\lambda'\in(C_\rho,\frac{\lambda}{C_\rho^2})$
is the same positive constant as in the proof of Theorem~\ref{6.16}.
Then, by the definition of $G_2(k)$,
we find that,
for any $x\in\lambda'B_{i,k}$,
$$
\rho(x,x_{i,k})<\lambda'r_{i,k}<\frac{\lambda}{C_\rho^2} r_{i,k}<
\lambda^2 r_{i,k}<\delta.
$$
From this and the fact that $f\in\mathrm{UC}(X)$,
we deduce that, for any
$k\in\mathbb{Z}\cap[\widetilde{k}+1,K_1]$,
$i\in G_2(k)$, and
$x\in\mathrm{\,supp\,} b_{i,k}\subset\lambda'B_{i,k}$,
\begin{align}
\label{89-1}
\left|b_{i,k}(x)\right|
&\leq\left|f(x)-\frac{1}{\|\xi_{i,k}\|_{L^1(X)}}
\int_{B_\rho(x_{i,k},\lambda'r_{i,k})}
f(y)\xi_{i,k}(y)\,d\mu(y)\right|
\nonumber\\
&\leq\left|f(x)-f(x_{i,k})\right|
\nonumber\\
&\quad+\left|\frac{1}{\|\xi_{i,k}\|_{L^1(X)}}
\int_{B_\rho(x_{i,k},\lambda'r_{i,k})}
\left[f(y)-f(x_{i,k})\right]
\xi_{i,k}(y)\,d\mu(y)\right|
\nonumber\\
&\leq2\eta.
\end{align}
Note that, for any
$k\in\mathbb{Z}\cap[\widetilde{k}+1,K_1]$,
$i\in G_2(k)$, and $j\in\mathbb{N}$,
if there exists $x\in X$ such that $b_{j,k+1}(x)\xi_{i,k}(x)\neq0$,
then
\begin{align}\label{89-2}
r_{j,k+1}\leq\frac{\Lambda}{C_\rho\lambda'}r_{i,k},
\end{align}
where $\Lambda\in(4C_\rho^2\lambda,\frac{\lambda^2}{C_\rho^4})$
is the same as in the proof of Theorem~\ref{6.16} and
satisfies the condition listed in Lemma~\ref{Wtd}(iii).
To the contrary, assume that there exist
$k\in\mathbb{Z}\cap[\widetilde{k}+1,K_1]$,
$i\in G_2(k)$, $j\in\mathbb{N}$, and $x\in X$
such that $b_{j,k+1}(x)\xi_{i,k}(x)\neq0$ and
$r_{j,k+1}>\frac{\Lambda r_{i,k}}{C_\rho\lambda'}$.
Since $b_{j,k+1}(x)\xi_{i,k}(x)\neq0$,
it follows from the support conditions listed
in Remark~\ref{5.5} that
$$
x\in B_\rho(x_{j,k+1},\lambda'r_{j,k+1})
\cap B_\rho(x_{i,k},\lambda'r_{i,k}).
$$
Then, by Definition~\ref{D2.1}(iii) and the fact that
$C_\rho\lambda'<\lambda<\Lambda$, we conclude that
$$
\rho(x_{j,k+1},x_{i,k})
\leq C_\rho\max\{\rho(x_{j,k+1},x),\,\rho(x,x_{i,k})\}
\leq C_\rho\lambda'\max\{r_{j,k+1},\,r_{i,k}\}
=C_\rho\lambda'r_{j,k+1}.
$$
From this and $r_{j,k+1}>\frac{\Lambda r_{i,k}}{C_\rho\lambda'}$,
we infer that, for any
$y\in B_\rho(x_{i,k},\Lambda r_{i,k})$,
\begin{align*}
\rho(y,x_{j,k+1})
&\leq C_\rho\max\{\rho(y,x_{i,k}),\,\rho(x_{i,k},x_{j,k+1})\}
\\
&\leq C_\rho\max\left\{\Lambda r_{i,k},\,C_\rho\lambda'r_{j,k+1}\right\}
\leq C_\rho^2\lambda'r_{j,k+1}<\lambda r_{j,k+1},
\end{align*}
which implies that
$B_\rho(x_{i,k},\Lambda r_{i,k})
\subset B_\rho(x_{j,k+1},\lambda r_{j,k+1})$.
By this, Lemma~\ref{Wtd}(iii),
and $\Omega_{k+1}\subset\Omega_k$,
we find that,
for any $k\in\mathbb{Z}\cap[\widetilde{k}+1,K_1]$,
$i\in G_2(k)$, and $j\in\mathbb{N}$,
$$
\varnothing\neq\left[B_\rho(x_{i,k},\Lambda r_{i,k})
\cap\left(\Omega_k\right)^\complement\right]
\subset\left[B_\rho(x_{j,k+1},\lambda r_{j,k+1})
\cap\left(\Omega_{k+1}\right)^\complement\right]=\varnothing,
$$
which is a contradiction. Thus, \eqref{89-2} holds.
From this, it follows that,
for any $y\in B_\rho(x_{j,k+1},\lambda'r_{j,k+1})$,
$$
\rho(y,x_{j,k+1})<\lambda'r_{j,k+1}
\leq\frac{\Lambda}{C_\rho}r_{i,k}<
\lambda^2r_{i,k}<\delta.
$$
By this and both (ii) and (iii) of Lemma~\ref{pou}, we obtain,
for any $x\in\mathrm{supp\,}b_{j,k+1}\xi_{i,k}
\subset[\lambda'B_{j,k+1}\cap\lambda'B_{i,k}]$,
\begin{align}
\label{89-3}
\left|b_{j,k+1}(x)\xi_{i,k}(x)\right|
&\leq\left|f(x)-f(x_{j,k+1})\right|
+\frac{1}{\|\xi_{j,k+1}\|_{L^1(X)}}
\nonumber\\
&\quad\times
\int_{B_\rho(x_{j,k+1},\lambda'r_{j,k+1})}
\left|f(y)-f(x_{j,k+1})\right|
\xi_{j,k+1}(y)\,d\mu(y)\nonumber
\\
&<2\eta
\end{align}
and, for any $x\in X$,
\begin{align*}
\sum_{j\in\mathbb{N}}\left|L_{i,j}^{(k+1)}
\xi_{j,k+1}(x)\right|
&=\sum_{j\in\mathbb{N}}
\frac{1}{\|\xi_{j,k+1}\|_{L^1(X)}}
\left|\int_Xb_{j,k+1}(y)\xi_{i,k}(y)
\,d\mu(y)\right|\xi_{j,k+1}(x)
\nonumber\\
&<2\eta\sum_{j\in\mathbb{N}}
\frac{1}{\|\xi_{j,k+1}\|_{L^1(X)}}
\left|\int_X\xi_{j,k+1}(y)
\,d\mu(y)\right|\xi_{j,k+1}(x)
\nonumber\\
&\leq2\eta\sum_{j\in\mathbb{N}}
\xi_{j,k+1}(x)\leq2\eta.
\end{align*}
From this, \eqref{89-0}, \eqref{89-1},
\eqref{89-3}, and \eqref{suppaik}, we deduce that,
for any $k\in\mathbb{Z}\cap[\widetilde{k}+1,K_1]$,
$i\in G_2(k)$, and $\mu$-almost every $x\in X$,
$$
\left|\lambda_{i,k}a_{i,k}(x)\right|\leq
6\eta\mathbf{1}_{\lambda B_{i,k}}(x),
$$
which, together with \eqref{9581},
further implies that
$$
\left\|l_{2,\delta}\right\|_{L^\infty(X)}
\leq\sum_{k=\widetilde{k}+1}^{K_1}\left\|\sum_{i\in G_2(k)}6\eta
\mathbf{1}_{\lambda B_{i,k}}\right\|_{L^\infty(X)}
\leq6L_0\eta(K_1-\widetilde{k}).
$$
Choose $\eta\in(0,\infty)$ small enough such that
$\left\|l_{2,\delta}\right\|_{L^\infty(X)}
\leq[(V_\rho)_{\lambda R}(x_1)]^{-\frac{1}{p}}$.
Then $l_{2,\delta}$ becomes a $(p,\infty)$-atom.
By this, \eqref{f=h+l}, the fact that $h$ is
a harmless constant multiple of $(p,\infty)$-atom,
\eqref{l=l1+l2}, and \eqref{l1delta},
we find that
\begin{align*}
\|f\|_{H^{p,q,\infty}_{\mathrm{at},\,\mathrm{fin}}(X,\rho)}^q
\lesssim
\|h\|_{H^{p,q,\infty}_{\mathrm{at},\,\mathrm{fin}}(X,\rho)}^q
+\|l_{1,\delta}\|_{H^{p,q,\infty}_{\mathrm{at},\,\mathrm{fin}}(X,\rho)}^q
+\|l_{2,\delta}\|_{H^{p,q,\infty}_{\mathrm{at},\,\mathrm{fin}}(X,\rho)}^q
\lesssim\|f\|_{H^{p,q}_*(X)}^q.
\end{align*}
This finishes the proof of (ii)
and hence Theorem~\ref{6.29}.
\end{proof}

\section{Molecular Characterization}
\label{section5.5}

In this section, we establish
the molecular characterization of the maximal-function-based
Hardy space $H^{p,q}_*(X)$
for the optimal ranges
$$
p\in\left(\frac{n}{n+{\mathrm{ind\,}}(X,\mathbf{q})},1\right]
\ \ \text{and}\ \
q\in(0,\infty]
$$
on ultra-RD spaces with upper dimension $n\in(0,\infty)$.
Indeed, the molecular characterization established in this section
also holds in the more general setting of quasi-ultrametric spaces
of homogeneous type for the same range of $p$ and $q$ above, but
with a different range of parameters for the molecular
Hardy--Lorentz space [see Remark~\ref{Rmk-mol}(i)].
As we see in Section~\ref{section5.8},
the importance of such a characterization naturally
arises when one is seeking to
extend a bounded linear operator on $L^2(X)$ to a
bounded operator on $H^{p,q}_*(X)$.

We begin with
the following notion of $(\rho,p,r,\eta)$-molecules.

\begin{definition}\label{mol}
Let $(X,\mathbf{q},\mu)$ be a quasi-ultrametric space
of homogeneous type,
$p\in(0,1]$, $r\in(p,\infty]\cap[1,\infty]$,
$\eta\in(0,\infty)$, and $\rho\in\mathbf{q}$ have
the property that all $\rho$-balls are $\mu$-measurable.
A function
$m\in L^r(X)$ is called a $(\rho,p,r,\eta)$-\emph{molecule}
\index[words]{molecule}
if there exist $x\in X$
and $r\in(0,\infty)$ such that
\begin{enumerate}
\item[\textup{(i)}]
$\|m\mathbf{1}_{B_\rho(x,r)}\|_{L^r(X)}
\leq[(V_\rho)_r(x)]
^{\frac{1}{r}-\frac{1}{p}}$.

\item[\textup{(ii)}]
For any $j\in\mathbb{N}$,
\begin{align*}
\left\|m\mathbf{1}_{B_{\rho}(x,2^jr)\setminus
B_{\rho}(x,2^{j-1}r)}\right\|_{L^r(X)}\leq2^{-j\eta}
\left[(V_\rho)_{2^jr}(x)\right]^{\frac{1}{r}-\frac{1}{p}}.
\end{align*}

\item[\textup{(iii)}]
$\int_Xm(x)\,d\mu(x)=0$.
\end{enumerate}
In this case, the $(\rho,p,r,\eta)$-molecule $m$ is said
to be \emph{centered at the ball} $B_\rho(x_0,r_0)$.
Moreover, when $\mu(X)<\infty$,
it is agreed upon that the constant function $m:=[\mu(X)]^{-\frac{1}{p}}$
is a $(\rho,p,r,\eta)$-molecule.
\end{definition}

\begin{remark}
\begin{enumerate}
\item[\rm(i)]
It is easy to see from Definition~\ref{mol},
\eqref{1557-xx}, and \eqref{doublingcond}
that, if $m: X\to\mathbb{C}$ is a $(\rho,p,r,\eta)$-molecule
centered at $B_{\rho}(x,r)$ for some $x\in X$
and $r\in(0,\infty)$, then, for any other $\rho'\in\mathbf{q}$
with the property that all $\rho'$-balls are $\mu$-measurable,
there exist positive constants $C$ and $c$ (depending only on $\rho$,
$\rho'$, $\mu$, and $\eta$) such that $Cm$ is a
$(\rho',p,r,\eta)$-molecule centered at $B_{\rho'}(x,cr)$.
As such, we sometimes simply refer to $m$ as a
$(p,r,\eta)$-\textit{molecule} or a
$(p,r,\eta)$-\textit{molecule modulo a harmless constant multiple}.

\item[\rm(ii)]
Conditions (i) and (ii) in Definition~\ref{mol} ensure that any
$(\rho,p,r,\eta)$-molecule is in $L^1(X)$, and hence the integral in
Definition~\ref{mol}(iii) is well defined.
Indeed, using H\"older's inequality, we find that,
for any $(\rho,p,r,\eta)$-molecule $m$ centered at some ball $B\subset X$,
\begin{align*}
\|m\|_{L^1(X)}
&=\int_B|m(x)|\,d\mu(x)
+\sum_{j=1}^\infty\int_{2^jB\setminus2^{j-1}B}
|m(x)|\,d\mu(x)\\
&\leq\left\|m\mathbf{1}_B\right\|_{L^r(X)}
\left\|\mathbf{1}_B\right\|_{L^{r'}(X)}\\
&\quad+\sum_{j=1}^\infty\left\|m\mathbf{1}_{2^jB
\setminus2^{j-1}B}\right\|_{L^r(X)}
\left\|\mathbf{1}_{2^jB}\right\|_{L^{r'}(X)}\\
&\leq\left[\mu(B)\right]^{1-\frac{1}{p}}
+\sum_{j=1}^\infty2^{-j\eta}
\left[\mu(2^jB)\right]^{1-\frac{1}{p}}
\lesssim\left[\mu(B)\right]^{1-\frac{1}{p}}.
\end{align*}

\item[\rm(iii)]
Any $(\rho,p,r,\eta)$-molecule
is normalized in the way that its $L^p(X)$ norm does not exceed
$(\frac{2^{\eta p}}{2^{\eta p}-1})^\frac{1}{p}$.
Indeed, from H\"older's inequality and both (i) and (ii)
of Definition~\ref{mol}, we deduce that,
for any $(\rho,p,r,\eta)$-molecule $m$ centered at some ball $B\subset X$,
\begin{align*}
\|m\|_{L^p(X)}
&=\left[\int_B|m(x)|^p\,d\mu(x)+
\sum_{j=1}^\infty\int_{2^jB\setminus2^{j-1}B}|m(x)|^p\,d\mu(x)\right]^\frac{1}{p}\\
&\leq\left[\left\|\,|m|^p\mathbf{1}_B\right\|_{L^\frac{r}{p}(X)}
\left\|\mathbf{1}_B\right\|_{L^{(\frac{r}{p})'}(X)}\right.\\
&\quad\left.+\sum_{j=1}^\infty\left\||m|^p
\mathbf{1}_{2^jB\setminus2^{j-1}B}\right\|_{L^\frac{r}{p}(X)}
\left\|\mathbf{1}_{2^jB}\right\|_{L^{(\frac{r}{p})'}(X)}\right]^\frac{1}{p}\\
&\leq\left(1+\sum_{j=1}^\infty2^{-j\eta p}\right)^\frac{1}{p}
=\left(\frac{2^{\eta p}}{2^{\eta p}-1}\right)^\frac{1}{p}.
\end{align*}
\end{enumerate}
\end{remark}

It is easy to see from definitions that every $(\rho,p,r)$-atom
is a $(\rho,p,r,\eta)$-molecule for any $\eta\in(0,\infty)$.
While the converse is not true,
in general, the following lemma shows that every $(\rho,p,r,\eta)$-molecule
can be written as a linear combination of $(\rho,p,r)$-atoms.

\begin{lemma}\label{L6.2}
Let $(X,\mathbf{q},\mu)$ be a quasi-ultrametric space
of homogeneous type,
$\rho\in\mathbf{q}$ have the property that
all $\rho$-balls are $\mu$-measurable,
$p\in(0,1]$, $r\in(p,\infty]\cap[1,\infty]$, and $\eta\in(0,\infty)$.
If $m$ is a $(\rho,p,r,\eta)$-molecule
centered at a ball $B:=B_\rho(x_0,r_0)$
with $x_0\in X$ and $r_0\in(0,\infty)$,
then, for any $x\in X$,
\begin{align*}
m(x)=\sum_{k=0}^{\infty}s_ka_k(x)
+\sum_{i=1}^{\infty}
\sum_{j=i}^{\infty}s_jb_{j,i}(x),
\end{align*}
where, for any $k\in\mathbb{Z}_+$,
$i\in\mathbb{N}$,
and $j\in\mathbb{N}\cap[i,\infty)$,
$a_k$ and $b_{j,i}$ are
$(\rho,p,r)$-atoms supported, respectively, in $2^{k}B$
and $2^iB$,
and $|s_k|\leq2^{1-k\eta}C_\mu^{1-\frac{1}{r}}$,
where $C_\mu\in(1,\infty)$ is the same as in \eqref{upperdoub}.
\end{lemma}

\begin{proof}
Let $m$ be a $(\rho,p,r,\eta)$-molecule
centered at a ball $B:=B_\rho(x_0,r_0)$
with $x_0\in X$ and $r_0\in(0,\infty)$.
Let $R_0:=B$ and, for any $k\in\mathbb{N}$,
$R_k:=2^kB\setminus 2^{k-1}B$.
For any $k\in\mathbb{Z}_+$, let
\begin{align*}
m_k:=\int_Xm(y)\mathbf{1}_{R_k}(y)\,d\mu(y)
\end{align*}
and
\begin{align*}
A_k:=m\mathbf{1}_{R_{k}}-\frac{\mathbf{1}_{2^kB}}{\mu(2^kB)}
m_k.
\end{align*}
By this, we find that
\begin{align*}
m&=\sum_{k\in\mathbb{Z}_+}A_k
+\sum_{k\in\mathbb{Z}_+}\frac{\mathbf{1}_{2^kB}}{\mu(2^kB)}m_k
=\sum_{k\in\mathbb{Z}_+}A_k
+\sum_{k\in\mathbb{Z}_+}\frac{\mathbf{1}_{2^kB}}{\mu(2^kB)}
\left(\sum_{j=k}^\infty m_j-\sum_{j=k+1}^\infty m_j\right)\\
&=\sum_{k\in\mathbb{Z}_+}A_k
+\sum_{k\in\mathbb{N}}
\sum_{j=k}^\infty\left[\frac{\mathbf{1}_{2^kB}}{\mu(2^kB)}-
\frac{\mathbf{1}_{2^{k-1}B}}{\mu(2^{k-1}B)}\right]
m_j=:\sum_{k\in\mathbb{Z}_+}A_k
+\sum_{k\in\mathbb{N}}\sum_{j=k}^\infty B_{j,k},
\end{align*}
where all series converge pointwise everywhere in $X$.

Next, we claim that, for any $k\in\mathbb{Z}_+$,
$A_k$ is a harmless constant multiple a $(\rho,p,r)$-atom
supported in $2^{k}B$ and, for any $k\in\mathbb{N}$
and $j\in\mathbb{Z}_+\cap[k,\infty)$,
$B_{j,k}$ is also a harmless constant multiple of a $(\rho,p,r)$-atom
supported in $2^{k}B$.
On the one hand,
from H\"older's inequality, Definition~\ref{mol}(ii),
and $C_\mu^{1-\frac{1}{r}}\geq1$,
we deduce that, for any $k\in\mathbb{Z}_+$,
\begin{align*}
\left\|A_k\right\|_{L^r(X)}
&\leq\left\|m\mathbf{1}_{R_k}\right\|_{L^r(X)}
+\left[\mu(2^kB)\right]^{\frac{1}{r}-1}
\left\|m\mathbf{1}_{R_k}\right\|_{L^1(X)}\\
&\leq2\left\|m\mathbf{1}_{R_k}\right\|_{L^r(X)}
\leq2^{1-k\eta}\left[\mu(2^kB)\right]^{\frac{1}{r}-\frac{1}{p}}
\leq2^{1-k\eta}C_\mu^{1-\frac{1}{r}}
\left[\mu(2^{k}B)\right]^{\frac{1}{r}-\frac{1}{p}},
\end{align*}
which, combined with the facts that
$\mathrm{supp\,}A_k\subset2^kB$ and $\int_XA_k(x)\,d\mu(x)=0$,
further implies that $a_k:=\frac{A_k}{2^{1-k\eta}C_\mu^{1-\frac{1}{r}}}$
is a $(\rho,p,r)$-atom supported in $2^kB$.
On the other hand, by $p\in(0,1]$,
we conclude that,
for any $k\in\mathbb{N}$
and $j\in\mathbb{N}\cap[k,\infty)$,
\begin{align*}
\left\|B_{j,k}\right\|_{L^r(X)}
&\leq\left\{\left[\mu(2^kB)\right]^{\frac{1}{r}-1}
+\left[\mu(2^{k-1}B)\right]^{\frac{1}{r}-1}\right\}
\left\|m\mathbf{1}_{R_j}\right\|_{L^1(X)}\\
&\leq2C_\mu^{1-\frac{1}{r}}\left[\mu(2^{k}B)\right]^{\frac{1}{r}-1}
\left\|m\mathbf{1}_{R_j}\right\|_{L^1(X)}
\quad\text{by \eqref{doublingcond}}\\
&\leq2C_\mu^{1-\frac{1}{r}}\left[\mu(2^{k}B)\right]^{\frac{1}{r}-1}
\left[\mu(2^jB)\right]^\frac{1}{r'}
\left\|m\mathbf{1}_{R_j}\right\|_{L^r(X)}
\quad\text{by H\"older's inequality}\\
&\leq2^{1-j\eta}C_\mu^{1-\frac{1}{r}}\left[\mu(2^{k}B)\right]^{\frac{1}{r}-1}
\left[\mu(2^jB)\right]^{1-\frac{1}{p}}
\quad\text{by Definition~\ref{mol}(ii)}\\
&\leq2^{1-j\eta}C_\mu^{1-\frac{1}{r}}
\left[\mu(2^{k}B)\right]^{\frac{1}{r}-\frac{1}{p}}
\quad\text{use $j\ge k$ and $p\leq1$},
\end{align*}
which, together with the facts that
$\mathrm{supp\,}B_{j,k}\subset2^{k}B$ and $\int_XB_{j,k}(x)\,d\mu(x)=0$,
further implies that $b_{j,k}:=
\frac{B_{j,k}}{2^{1-j\eta}C_\mu^{1-\frac{1}{r}}}$
is a $(\rho,p,r)$-atom supported in $2^{k}B$.
Note that all the desired claims follow from these observations
with $s_k:=2^{1-k\eta}C_\mu^{1-\frac{1}{r}}$ for
$k\in\mathbb{Z}_+$. This finishes the proof of Lemma~\ref{L6.2}.
\end{proof}

Based on the notion of molecules in Definition~\ref{mol},
we introduce molecular Hardy--Lorentz spaces as follows.

\begin{definition}
\label{Hpqm}
Let $(X,\mathbf{q},\mu)$ be a quasi-ultrametric space
of homogeneous type,
$\rho\in\mathbf{q}$ have the property that
all $\rho$-balls are $\mu$-measurable,
$p\in(0,1]$,
$q\in(0,\infty]$, $r\in(p,\infty]\cap[1,\infty]$,
$\eta\in(0,\infty)$,
and $\varepsilon,\beta,\gamma\in(0,\infty)$ be such that
$n(\frac{1}{p}-1)<\beta,\gamma<\varepsilon<\infty$.
The \emph{molecular Hardy--Lorentz
space}\index[words]{molecular Hardy--Lorentz space}
$H^{p,q,r,\eta}_{\mathrm{mol}}(X,\rho)$
\index[symbols]{H@$H^{p,q,r,\eta}_{\mathrm{mol}}(X,\rho)$}
is defined to be the set of all
$f\in({\mathcal{G}}^\varepsilon_0(\beta,\gamma))'$ satisfying that
there exist a sequence
$\{\lambda_{i,k}\}_{i\in\mathbb{N},k\in\mathbb{Z}}$ in $\mathbb{R}_+$
and a sequence
of $(\rho,p,r,\eta)$-molecules
$\{m_{i,k}\}_{i\in\mathbb{N},k\in\mathbb{Z}}$
supported, respectively,
in $\rho$-balls $\{B_{i,k}\}_{i\in\mathbb{N},k\in\mathbb{Z}}$
such that
\begin{align}\label{decmol}
f=\sum_{k\in\mathbb{Z}}\sum_{i\in\mathbb{N}}
\lambda_{i,k}m_{i,k}
\end{align}
in $({\mathcal{G}}^\varepsilon_0(\beta,\gamma))'$
and the following conditions hold:
\begin{enumerate}
\item[\rm(i)]
There exist constants $C_0\in\mathbb{N}$
and $j_0\in\mathbb{Z}\setminus\mathbb{N}$,
independent of $f$,
such that, for any $x\in X$ and $k\in\mathbb{Z}$,
\begin{align}
\label{830}
\sum_{i\in\mathbb{N}}
\mathbf{1}_{(2C_\rho)^{j_0}B_{i,k}}(x)\leq C_0.
\end{align}

\item[\rm(ii)]
For any $i\in\mathbb{N}$ and $k\in\mathbb{Z}$,
$\lambda_{i,k}\sim2^k[\mu(B_{i,k})]^\frac{1}{p}$
with the positive equivalence constants independent of $f$ and
$$
\sum_{k\in\mathbb{Z}}
\left(\sum_{i\in\mathbb{N}}
\lambda_{i,k}^p\right)^\frac{q}{p}<\infty
$$
with the usual modification when $q=\infty$.
\end{enumerate}
Moreover, for any $f\in H^{p,q,r,\eta}_{\mathrm{mol}}(X,\rho)$, define
$$
\|f\|_{H^{p,q,r,\eta}_{\mathrm{mol}}(X,\rho)}:=
\inf\left\{\left[\sum_{k\in\mathbb{Z}}
\left(\sum_{i\in\mathbb{N}}
\lambda_{i,k}^p\right)^\frac{q}{p}\right]^\frac{1}{q}:
f=\sum_{k\in\mathbb{Z}}\sum_{i\in\mathbb{N}}
\lambda_{i,k}m_{i,k}\right\}
$$
with the usual modification when $q=\infty$,
where the infimum is taken over all
decompositions of $f$ as in \eqref{decmol}.
\end{definition}

\begin{remark}
Let the notation be as in Definition~\ref{Hpqm}.
\begin{enumerate}
\item[\rm(i)]
From the assumption $n(\frac{1}{p}-1)<\beta,\gamma<\infty$,
\eqref{1550}, and Remark~\ref{Rmtest}(v), we infer that,
for any $p\in(0,\frac{n}{n+\mathrm{ind}_H(X,\mathbf{q})}]$,
$H^{p,q,r,\eta}_{\mathrm{mol}}(X,\rho)$ is trivial because the space
${\mathcal{G}}^\varepsilon_0(\beta,\gamma)$
only contains constant functions.

\item[\rm(ii)]
According to Proposition~\ref{5.1.11},
Theorem~\ref{thmHpqm}, and Remark~\ref{Rmk-mol}(i),
we find that
$H^{p,q,r,\eta}_{\mathrm{mol}}(X,\rho)$
is a quasi-Banach space for certain ranges of
$p$, $q$, $r$, $\eta$, $\varepsilon$, $\beta$, and $\gamma$.
\end{enumerate}
\end{remark}

Now, we present the main theorem of this
section, which establishes the molecular characterization
of the maximal-function-based Hardy space $H^{p,q}_*(X)$.

\begin{theorem}\label{thmHpqm}
Let $(X,\mathbf{q},\mu)$ be a $(\kappa,n)$-ultra-RD-space with
$0<\kappa\leq n<\infty$,
$p\in(\frac{n}{n+\mathrm{ind\,}(X,\mathbf{q})},1]$,
$\rho\in\mathbf{q}$ have the property that all $\rho$-balls
are $\mu$-measurable,
$q\in(0,\infty]$, $r\in(p,\infty]\cap[1,\infty]$,
$\varepsilon,\beta,\gamma\in(0,\infty)$ be such that
$n(\frac{1}{p}-1)<\beta,\gamma<
\varepsilon\preceq\mathrm{ind\,}(X,\mathbf{q})$,
and
\begin{align}
\label{etarange}
\eta\in\left(\max\left\{\frac{n-\kappa}{p},\,\beta+
\frac{\kappa(p-1)}{p}\right\},\infty\right).
\end{align}
Then, as subspaces of
$({\mathcal{G}}^\varepsilon_0(\beta,\gamma))'$,
$$
H^{p,q,r,\eta}_{\mathrm{mol}}(X,\rho)=H^{p,q}_*(X)
$$
with equivalent (quasi-)norms.
Moreover, if a sequence
$\{\lambda_{i,k}\}_{i\in\mathbb{N},k\in
\mathbb{Z}}$ in $\mathbb{R}_+$ and a sequence
of $(\rho,p,r,\eta)$-molecules
$\{m_{i,k}\}_{i\in\mathbb{N},k\in\mathbb{Z}}$
supported, respectively,
in $\rho$-balls $\{B_{i,k}\}_{i\in\mathbb{N},k\in\mathbb{Z}}$
satisfy both (i) and (ii) of Definition~\ref{Hpqm},
then
$\sum_{k\in\mathbb{Z}}\sum_{i\in\mathbb{N}}\lambda_{i,k}m_{i,k}$
converges both in $(\mathcal{G}^\varepsilon_0(\beta,\gamma))'$
and in $H^{p,q}_*(X)$.
\end{theorem}

\begin{remark}\label{Rmk-mol}
Let the notation be as in Theorem~\ref{thmHpqm}.
\begin{enumerate}
\item[\rm(i)]
Let $(X,\mathbf{q},\mu)$
be a quasi-ultrametric space of homogeneous type with
upper dimension $n\in(0,\infty)$,
but without assuming that $(X,\mathbf{q},\mu)$
satisfies the reverse doubling condition \eqref{RDcondition}.
Then, following an argument similar to that used in
the proof of Theorem~\ref{thmHpqm},
we conclude that Theorem~\ref{thmHpqm} in this case
holds, but with \eqref{etarange}
replaced by
$$
\eta\in\left(\max\left\{\frac{n}{p},\,\beta\right\},\infty\right).
$$

\item[\rm(ii)]
Let $(X,\mathbf{q},\mu)$ be an $n$-Ahlfors-regular
quasi-ultrametric space
with $n\in(0,\infty)$, which means that
$(X,\mathbf{q},\mu)$ is an $(n,n)$-ultra-RD-space.
Let $\rho\in\mathbf{q}$ be such that all $\rho$-balls
are $\mu$-measurable.
Then, by specializing Theorem~\ref{thmHpqm} to this setting,
\eqref{etarange} becomes
$$
\eta\in\left(\beta+\frac{n(p-1)}{p},\infty\right).
$$

\item[\rm(iii)]
The range of $\eta$ in Theorem~\ref{thmHpqm}
is greater than \cite[Theorem 6.4]{zhy}
(wherein Zhou et al. worked in general
spaces of homogeneous type) because we establish
Theorem~\ref{thmHpqm} in the more specialized setting of ultra-RD-spaces.
Even in the setting of $\mathbb{R}^n$, this
range of $\eta$ in Theorem~\ref{thmHpqm}
is greater than \cite[Theorem 3.9]{lyy2016}
because we decompose molecules into
linear combinations of atoms (see Lemma~\ref{L6.2}).

\item[\rm(iv)]
From Theorems~\ref{thmHpqm} and \ref{6.14},
combined with Remark~\ref{5216}(iii), we deduce that
$H^{p,q,r,\eta}_{\mathrm{mol}}(X,\rho)$ is
independent of the choice of $\rho\in\mathbf{q}$.
Thus, in what follows, we sometimes simply write
$H^{p,q,r,\eta}_{\mathrm{mol}}(X)$ instead
of $H^{p,q,r,\eta}_{\mathrm{mol}}(X,\rho)$.

\item[\rm(v)]
Alvarado and Mitrea \cite{AlMi2015}
introduced the molecular Hardy space $H^{p,r,A,\epsilon}_{\rm mol}(X)$
(see \cite[pp.\,272--273]{AlMi2015}) and
established the molecular characterization of the maximal Hardy space $H^p(X)$
on $n$-Ahlfors regular quasi-ultrametric spaces $(X,\mathbf{q},\mu)$
for the optimal ranges
\begin{align*}
p\in\left(\frac{n}{n+\mathrm{ind\,}(X,\mathbf{q})},1\right]
\ \ \text{and}\ \
r\in(p,\infty]\cap[1,\infty];
\end{align*}
see \cite[Theorem~6.11]{AlMi2015}.
By Remarks~\ref{2548} and~\ref{5216}(iv),
we find that Theorem~\ref{thmHpqm}
is a generalization of the molecular
characterization in \cite[Theorem~6.11]{AlMi2015}
from maximal Hardy spaces on $n$-Ahlfors regular quasi-ultrametric spaces
to maximal Hardy--Lorentz spaces on ultra-RD-spaces
(or quasi-ultrametric spaces of homogeneous type) for the same and
optimal ranges of $p$ and $r$.
\end{enumerate}
\end{remark}

\begin{proof}[Proof of Theorem~\ref{thmHpqm}]
From Definitions~\ref{atom} and \ref{mol},
we infer that any
$(\rho,p,\infty)$-atom is a $(\rho,p,r,\eta)$-molecule.
By this and Theorem~\ref{6.14}, we conclude that,
for any $f\in H^{p,q}_*(X)=H^{p,q,\infty}_{\mathrm{at}}(X,\rho)$,
$f\in H^{p,q,r,\eta}_{\mathrm{mol}}(X,\rho)$ and
$$
\|f\|_{H^{p,q,r,\eta}_{\mathrm{mol}}(X,\rho)}
\leq\|f\|_{H^{p,q,\infty}_{\mathrm{at}}(X,\rho)}
\sim\|f\|_{H^{p,q}_*(X)}.
$$

It remains to show that
$H^{p,q,r,\eta}_{\mathrm{mol}}(X,\rho)\subset H^{p,q}_*(X)$.
We only prove this in the case where both
$q\in(0,\infty)$ and $\mu(X)=\infty$
because the proof in the case where $q=\infty$ or $\mu(X)<\infty$
[when $\mu(X)<\infty$, we need to take advantage of \eqref{20281}]
is similar and hence we omit the details.
Let $\{\lambda_{i,k}\}_{i\in\mathbb{N},k\in\mathbb{Z}}$ be in $\mathbb{R}_+$
and $\{m_{i,k}\}_{i\in\mathbb{N},k\in\mathbb{Z}}$ be a sequence
of $(\rho,p,r,\eta)$-molecules supported, respectively,
in $\rho$-balls $\{B_{i,k}\}_{i\in\mathbb{N},k\in\mathbb{Z}}$
such that there exist constants
$C_0\in\mathbb{N}$,
$j_0\in\mathbb{Z}\setminus\mathbb{N}$,
and $C_1\in(0,\infty)$ such that
for any $k\in\mathbb{Z}$ and $x\in X$
\begin{align}\label{L0}
\sum_{i\in\mathbb{N}}
\mathbf{1}_{(2C_\rho)^{j_0}B_{i,k}}(x)\leq C_0,
\end{align}
for any $i\in\mathbb{N}$
\begin{align}
\label{lamdik}
\lambda_{i,k}\sim2^k\left[\mu(B_{i,k})\right]^\frac{1}{p},
\end{align}	
and
\begin{align*}
\sum_{k\in\mathbb{Z}}
\left(\sum_{i\in\mathbb{N}}
\lambda_{i,k}^p\right)^\frac{q}{p}<\infty.
\end{align*}	

To proceed, fix an arbitrary $k_0\in\mathbb{Z}$. We write
\begin{align}\label{10311}
g:=\sum_{k\in\mathbb{Z}}\sum_{i\in\mathbb{N}}
\lambda_{i,k}(m_{i,k})_\rho^*
=\sum_{k=-\infty}^{k_0-1}\sum_{i\in\mathbb{N}}
\lambda_{i,k}(m_{i,k})_\rho^*
+\sum_{k=k_0}^{\infty}\sum_{i\in\mathbb{N}}
\cdots=:g_1+g_2.
\end{align}
Since $\rho\sim\rho_{\#}$, where $\rho_{\#}$ is the
regularized quasi-ultrametric defined in
Lemma~\ref{L2.3}, by Remarks~\ref{RemMf}(iv) and~\ref{Rmk-atoms}(i), \eqref{eq-538},
and \eqref{1557-xx},
without loss of generality,
we may assume that
$\rho$ is a regularized quasi-ultrametric.
From this, we deduce that
\begin{align}
\label{muf1f2}
d_{g}(2^{k_0})
&\leq\mu\left(\left\{x\in X:g_1(x)
>2^{k_0-1}\right\}\right)
\nonumber\\
&\quad+\mu\left(\left\{x\in X:g_2(x)
>2^{k_0-1}\right\}\right)
\nonumber\\
&=:(\mathrm{I_1})^p+(\mathrm{I_2})^p.
\end{align}
For any $k\in\mathbb{Z}$ and $i\in\mathbb{N}$,
by Lemma~\ref{L6.2}, we can write,
for any $x\in X$,
$$
m_{i,k}(x)=\sum_{l=0}^{\infty}s_la^{i,k}_l(x)
+\sum_{l=0}^{\infty}
\sum_{u=l}^{\infty}s_ub_{u,l}^{i,k}(x),
$$
where, for any $l\in\mathbb{Z}_+$
and $u\in\mathbb{Z}_+\cap[l,\infty)$,
$a^{i,k}_l$ and $b_{u,l}^{i,k}$ are two
$(\rho,p,r)$-atoms both supported in $2^{l+1}B_{i,k}$
and $s_l:=2^{1-l\eta}(C_\mu2^n)^{\frac{1}{p}-\frac{1}{r}}$
with $ C_\mu\in[1,\infty)$ as in \eqref{doublingcond}.
By this, we find that
\begin{align}\label{I1}
\mathrm{I_1}&\leq
\left[\mu\left(\left\{x\in X:\sum_{k=-\infty}^{k_0-1}
\sum_{i\in\mathbb{N}}\lambda_{i,k}(m_{i,k})^*_\rho(x)
>2^{k_0-1}\right\}\right)\right]^\frac{1}{p}
\nonumber\\
&\leq2^\frac{1}{p}\left[\mu\left(\left\{x\in X:\sum_{k=-\infty}^{k_0-1}
\sum_{i\in\mathbb{N}}\sum_{l=0}^{\infty}
\lambda_{i,k}\mathrm{J}^{i,k}_{1,l}(x)
>2^{k_0-2}\right\}\right)\right]^\frac{1}{p}
\nonumber\\
&\quad+2^\frac{1}{p}\left[\mu\left(\left\{x\in X:\sum_{k=-\infty}^{k_0-1}
\sum_{i\in\mathbb{N}}\sum_{l=0}^{\infty}
\lambda_{i,k}\mathrm{J}^{i,k}_{2,l}(x)
>2^{k_0-2}\right\}\right)\right]^\frac{1}{p}
\nonumber\\
&=:2^\frac{1}{p}\left(\mathrm{I_{1,1}}+\mathrm{I_{1,2}}\right),
\end{align}
where, for any $k\in\mathbb{Z}$, $i\in\mathbb{N}$,
$l\in\mathbb{Z}_+$, and $x\in X$,
\begin{align*}
\mathrm{J}_{1,l}^{i,k}(x):=
\left[s_l\left(a^{i,k}_l\right)^*_\rho(x)
+\sum_{u=l}^{\infty}s_u
\left(b^{i,k}_{u,l}\right)^*_\rho(x)\right]
\mathbf{1}_{C_\rho2^{l+2}B_{i,k}}(x)
\end{align*}
and
\begin{align*}
\mathrm{J}_{2,l}^{i,k}(x):=
\left[s_l\left(a^{i,k}_l\right)^*_\rho(x)
+\sum_{u=l}^{\infty}s_u
\left(b^{i,k}_{u,l}\right)^*_\rho(x)\right]
\mathbf{1}_{[C_\rho2^{l+2}B_{i,k}]^\complement}(x).
\end{align*}
Similarly to \eqref{I1}, we have
\begin{align}
\label{I2}
\mathrm{I_2}
&\leq
\left[\mu\left(\left\{x\in X:\sum_{k=k_0}^{\infty}
\sum_{i\in\mathbb{N}}\lambda_{i,k}(m_{i,k})^*_\rho(x)
>2^{k_0-1}\right\}\right)\right]^\frac{1}{p}
\nonumber\\
&\leq2^\frac{1}{p}\left[\mu\left(\left\{x\in X:\sum_{k=k_0}^{\infty}
\sum_{i\in\mathbb{N}}\sum_{l=0}^{\infty}
\lambda_{i,k}\mathrm{J}^{i,k}_{1,l}(x)
>2^{k_0-2}\right\}\right)\right]^\frac{1}{p}
\nonumber\\
&\quad+2^\frac{1}{p}\left[\mu\left(\left\{x\in X:\sum_{k=k_0}^{\infty}
\sum_{i\in\mathbb{N}}\sum_{l=0}^{\infty}
\lambda_{i,k}\mathrm{J}^{i,k}_{2,l}(x)
>2^{k_0-2}\right\}\right)\right]^\frac{1}{p}
\nonumber\\
&=:2^\frac{1}{p}\left(\mathrm{I_{2,1}}+\mathrm{I_{2,2}}\right).
\end{align}
We claim that there exist constants
$\eta_1,\eta_2\in(1,\infty)$
and $\delta_1,\delta_2\in(0,1)$ (independent of $k_0$) such that,
for any $j\in\{1,\,2\}$,
\begin{align}
\label{82I1}
\mathrm{I_{1,j}}\lesssim
2^{-k_0\eta_j}
\left\{\sum_{k=-\infty}^{k_0-1}
2^{\eta_jkq}
\left[\sum_{i\in\mathbb{N}}
\mu(B_{i,k})\right]^\frac{q}{p}\right\}^\frac{1}{q}
\end{align}
and
\begin{align}
\label{82I2}
\mathrm{I_{2,j}}\lesssim
2^{-k_0\delta_j}
\left\{\sum_{k=k_0}^{\infty}
2^{\delta_jkq}
\left[\sum_{i\in\mathbb{N}}
\mu(B_{i,k})\right]^\frac{q}{p}\right\}^\frac{1}{q}.
\end{align}
Assume the moment that both
\eqref{82I1} and \eqref{82I2} hold. Then,
from \eqref{10311},
Proposition~\ref{Lpqnorm=c,d}, \eqref{muf1f2},
\eqref{I1}, and \eqref{I2},
it follows that
\begin{align*}
\|g\|_{L^{p,q}(X)}
&\sim\left\{\sum_{k_0\in\mathbb{Z}}2^{k_0q}
\left[d_g(2^{k_0})\right]^\frac{q}{p}\right\}^\frac{1}{q}\lesssim
\left\{\sum_{k_0\in\mathbb{Z}}2^{k_0q}
\left[(\mathrm{I_1})^p+(\mathrm{I_2})^p\right]^\frac{q}{p}\right\}^\frac{1}{q}
\\
&\lesssim
\left\{\sum_{k_0\in\mathbb{Z}}2^{k_0q}
\left[(\mathrm{I_{1,1}}+\mathrm{I_{1,2}})^p
+(\mathrm{I_{2,1}}+\mathrm{I_{2,2}})^p\right]^\frac{q}{p}\right\}^\frac{1}{q}
\\
&\lesssim
\left\{\sum_{k_0\in\mathbb{Z}}2^{k_0q}
(\mathrm{I_{1,1}})^q\right\}^\frac{1}{q}+
\left\{\sum_{k_0\in\mathbb{Z}}2^{k_0q}
(\mathrm{I_{1,2}})^q\right\}^\frac{1}{q}
\\
&\quad+\left\{\sum_{k_0\in\mathbb{Z}}2^{k_0q}
(\mathrm{I_{2,1}})^q\right\}^\frac{1}{q}+
\left\{\sum_{k_0\in\mathbb{Z}}2^{k_0q}
(\mathrm{I_{2,2}})^q\right\}^\frac{1}{q},
\end{align*}
which, together with \eqref{82I1}, \eqref{82I2},
\eqref{lamdik},
and an argument similar to that used in the estimation of \eqref{4ee5},
further implies that
\begin{align*}
\|g\|_{L^{p,q}(X)}
&\lesssim\sum_{j=1}^{2}
\left\{\sum_{k_0\in\mathbb{Z}}2^{k_0q}
2^{-k_0\eta_jq}
\sum_{k=-\infty}^{k_0-1}
2^{\eta_jkq}
\left[\sum_{i\in\mathbb{N}}
\mu(B_{i,k})\right]^\frac{q}{p}\right\}^\frac{1}{q}
\\
&\quad
+\sum_{j=1}^{2}
\left\{\sum_{k_0\in\mathbb{Z}}2^{k_0q}
2^{-k_0\delta_jq}
\sum_{k=k_0}^{\infty}
2^{\delta_jkq}
\left[\sum_{i\in\mathbb{N}}
\mu(B_{i,k})\right]^\frac{q}{p}\right\}^\frac{1}{q}
\\
&\lesssim\sum_{j=1}^{2}
\left[\sum_{k_0\in\mathbb{Z}}2^{(1-\eta_j)k_0q}
\sum_{k=-\infty}^{k_0-1}
2^{(\eta_j-1)kq}
\left(\sum_{i\in\mathbb{N}}
\lambda_{i,k}^p\right)^\frac{q}{p}\right]^\frac{1}{q}
\\
&\quad
+\sum_{j=1}^{2}
\left[\sum_{k_0\in\mathbb{Z}}2^{(1-\delta_j)k_0q}
\sum_{k=k_0}^{\infty}
2^{(\delta_j-1)kq}
\left(\sum_{i\in\mathbb{N}}
\lambda_{i,k}^p\right)^\frac{q}{p}\right]^\frac{1}{q}
\\
&\lesssim\left[\sum_{k\in\mathbb{Z}}
\left(\sum_{i\in\mathbb{N}}
\lambda_{i,k}^p\right)^\frac{q}{p}\right]^\frac{1}{q}.
\end{align*}
This, combined with an argument similar to that
used in the estimation of \eqref{16177},
further implies that, for any
$\phi\in\mathcal{G}^\varepsilon_0(\beta,\gamma)$
with $\|\phi\|_{\mathcal{G}^\varepsilon_0(\beta,\gamma)}\leq1$,
$$
\sum_{k\in\mathbb{Z}}
\sum_{i\in\mathbb{N}}\lambda_{i,k}\left|\left\langle m_{i,k},
\phi\right\rangle\right|<\infty
$$
and hence the series $f:=\sum_{k\in\mathbb{Z}}
\sum_{i\in\mathbb{N}}\lambda_{i,k}m_{i,k}$
converges in $(\mathcal{G}^\varepsilon_0(\beta,\gamma))'$.
By this and an argument similar to that used in the proof
of both (ii) and (iii) of Theorem~\ref{6.15},
we conclude that
the series $\sum_{k\in\mathbb{Z}}
\sum_{i\in\mathbb{N}}\lambda_{i,k}m_{i,k}$
converges in $H^{p,q}_*(X)$,
$f\in H^{p,q}_*(X)$, and
$$
\|f\|_{H^{p,q}_*(X)}=\left\|f^*_\rho\right\|_{L^{p,q}(X)}
\leq\|g\|_{L^{p,q}(X)}\lesssim\left[\sum_{k\in\mathbb{Z}}
\left(\sum_{i\in\mathbb{N}}
\lambda_{i,k}^p\right)^\frac{q}{p}\right]^\frac{1}{q}.
$$
Thus,
$H^{p,q,r,\eta}_{\mathrm{mol}}(X,\rho)\subset H^{p,q}_*(X)$
(with appropriate control of the quasi-norms),
which completes the proof
of Theorem~\ref{thmHpqm} when $q\in(0,\infty)$
under the assumption that \eqref{82I1} and \eqref{82I2} are true.

It remains to show the estimates in \eqref{82I1} and \eqref{82I2}.
We first deal with $\mathrm{I_{1,1}}$.
For any $k\in\mathbb{Z}$, $i\in\mathbb{N}$,
and $l\in\mathbb{Z_+}$, define
\begin{align}
\label{oikl}
O_{l}^{i,k}:=
s_l\mathcal{M}_\rho\left(a^{i,k}_l\right)
+\sum_{u=l}^{\infty}s_u
\mathcal{M}_\rho\left(b^{i,k}_{u,l}\right),
\end{align}
where $\mathcal{M}_{\rho}$ is the Hardy--Littlewood
maximal operator as in \eqref{Mf}.
Then, from Proposition~\ref{HL-grandmax}, we infer that,
for any $k\in\mathbb{Z}$, $i\in\mathbb{N}$, and $l\in\mathbb{Z_+}$,
\begin{align}
\label{oikl-ineq}
\mathrm{J}_{1,l}^{i,k}
&=\left[s_l\left(a^{i,k}_l\right)^*_\rho
+\sum_{u=l}^{\infty}s_u
\left(b^{i,k}_{u,l}\right)^*_\rho\right]
\mathbf{1}_{C_\rho2^{l+2}B_{i,k}}
\lesssim \left(O_{l}^{i,k}\right)\mathbf{1}_{C_\rho2^{l+2}B_{i,k}}.
\end{align}
Next, let $t\in(\frac{n}{(\eta+\frac{\kappa}{p})r},p)$,
$\widetilde{q}\in(\max\{1,\,\frac{n}{(\eta+\frac{\kappa}{p})t}\},
\min\{r,\,\frac{1}{t}\})$,
$\delta\in(0,1-\frac{1}{\widetilde{q}})$,
and
$$
s\in\left(\frac{n}{(\eta+\frac{\kappa}{p})t\widetilde{q}},
\min\left\{1,\,\frac{p}{t}\right\}\right).
$$
Using H\"older's inequality and $\widetilde{q}\in(1,\infty)$,
we find that
\begin{align*}
\sum_{k=-\infty}^{k_0-1}
\sum_{i\in\mathbb{N}}\sum_{l=0}^{\infty}
\lambda_{i,k}\mathrm{J}^{i,k}_{1,l}
&=\sum_{k=-\infty}^{k_0-1}2^{k\delta}\left(2^{-k\delta}
\sum_{i\in\mathbb{N}}\sum_{l=0}^{\infty}
\lambda_{i,k}\mathrm{J}^{i,k}_{1,l}\right)
\\
&\leq\left(\sum_{k=-\infty}^{k_0-1}
2^{k\delta\widetilde{q}'}\right)^\frac{1}{\widetilde{q}'}
\left[\sum_{k=-\infty}^{k_0-1}2^{-k\delta\widetilde{q}}\left(
\sum_{i\in\mathbb{N}}\sum_{l=0}^{\infty}
\lambda_{i,k}\mathrm{J}^{i,k}_{1,l}
\right)^{\widetilde{q}}\right]^\frac{1}{\widetilde{q}}
\\
&\sim 2^{k_0\delta}\left[\sum_{k=-\infty}^{k_0-1}
2^{-k\delta\widetilde{q}}\left(
\sum_{i\in\mathbb{N}}\sum_{l=0}^{\infty}
\lambda_{i,k}\mathrm{J}^{i,k}_{1,l}
\right)^{\widetilde{q}}\right]^\frac{1}{\widetilde{q}},
\end{align*}
where $\widetilde{q}'\in(1,\infty)$ satisfies
$\frac{1}{\widetilde{q}}+\frac{1}{\widetilde{q}'}=1$.
From this, Chebyshev's inequality,
\eqref{lamdik}, \eqref{oikl-ineq}, and
Lemma \ref{discreteinequality} with $r$ therein replaced by $t$,
we deduce that
\begin{align*}
\mathrm{I_{1,1}}
&\lesssim
\left[\mu\left(\left\{x\in X:2^{k_0\delta\widetilde{q}}
\sum_{k=-\infty}^{k_0-1}2^{-k\delta\widetilde{q}}\left(
\sum_{i\in\mathbb{N}}\sum_{l=0}^{\infty}
\lambda_{i,k}\mathrm{J}^{i,k}_{1,l}(x)
\right)^{\widetilde{q}}
>2^{(k_0-2)\widetilde{q}}\right\}\right)\right]^\frac{1}{p}
\nonumber\\
&\lesssim
2^{(\delta-1)k_0\widetilde{q}}\left\|
\sum_{k=-\infty}^{k_0-1}2^{-k\delta\widetilde{q}}\left(
\sum_{i\in\mathbb{N}}\sum_{l=0}^{\infty}
\lambda_{i,k}\mathrm{J}^{i,k}_{1,l}
\right)^{\widetilde{q}}\right\|_{L^p(X)}
\nonumber\\
&\lesssim
2^{(\delta-1)k_0\widetilde{q}}\left\|
\sum_{k=-\infty}^{k_0-1}2^{(1-\delta)k\widetilde{q}}\left(
\sum_{i\in\mathbb{N}}\sum_{l=0}^{\infty}
\left[\mu(B_{i,k})\right]^\frac{1}{p}
\mathrm{J}^{i,k}_{1,l}
\right)^{\widetilde{q}}\right\|_{L^p(X)}\nonumber
\\
&\lesssim
2^{(\delta-1)k_0\widetilde{q}}\left\|
\sum_{k=-\infty}^{k_0-1}2^{(1-\delta)k\widetilde{q}}\left(
\sum_{i\in\mathbb{N}}\sum_{l=0}^{\infty}
\left[\mu(B_{i,k})\right]^\frac{1}{p}
\left(O_{l}^{i,k}\right)
\mathbf{1}_{C_\rho2^{l+2}B_{i,k}}
\right)^{\widetilde{q}}\right\|_{L^p(X)}\nonumber
\\
&\leq
2^{(\delta-1)k_0\widetilde{q}}\left\|
\sum_{k=-\infty}^{k_0-1}2^{(1-\delta)k\widetilde{q}t}\left(
\sum_{i\in\mathbb{N}}\sum_{l=0}^{\infty}
\left[\mu(B_{i,k})\right]^\frac{1}{p}
\left(O_{l}^{i,k}\right)
\mathbf{1}_{C_\rho2^{l+2}B_{i,k}}
\right)^{\widetilde{q}t}\right\|_{L^\frac{p}{t}(X)}^\frac{1}{t}\nonumber,
\end{align*}
which, together with Lemma \ref{discreteinequality} with $r$ therein replaced by $\widetilde{q}t\in(0,1)$,
Minkowski's inequality for the norm $\|\cdot\|_{L^\frac{p}{t}(X)}$,
and Tonelli's theorem,
further implies that
\begin{align}\label{83I11}
\mathrm{I_{1,1}}
&\lesssim
2^{(\delta-1)k_0\widetilde{q}}\left\{
\sum_{k=-\infty}^{k_0-1}2^{(1-\delta)k\widetilde{q}t}\left\|\left(
\sum_{i\in\mathbb{N}}\sum_{l=0}^{\infty}
\left[\mu(B_{i,k})\right]^\frac{1}{p}
\left(O_{l}^{i,k}\right)
\mathbf{1}_{C_\rho2^{l+2}B_{i,k}}
\right)^{\widetilde{q}t}
\right\|_{L^\frac{p}{t}(X)}\right\}^\frac{1}{t}\nonumber
\\
&\leq
2^{(\delta-1)k_0\widetilde{q}}\left\{
\sum_{k=-\infty}^{k_0-1}2^{(1-\delta)k\widetilde{q}t}\left\|
\sum_{i\in\mathbb{N}}\sum_{l=0}^{\infty}\left[
\left[\mu(B_{i,k})\right]^\frac{1}{p}
\left(O_{l}^{i,k}\right)
\mathbf{1}_{C_\rho2^{l+2}B_{i,k}}\right]^{\widetilde{q}t}
\right\|_{L^\frac{p}{t}(X)}\right\}^\frac{1}{t}\nonumber
\\
&=
2^{(\delta-1)k_0\widetilde{q}}\left\{
\sum_{k=-\infty}^{k_0-1}2^{(1-\delta)k\widetilde{q}t}\left\|
\sum_{l=0}^{\infty}\sum_{i\in\mathbb{N}}\left[
\left[\mu(B_{i,k})\right]^\frac{1}{p}
\left(O_{l}^{i,k}\right)
\mathbf{1}_{C_\rho2^{l+2}B_{i,k}}\right]^{\widetilde{q}t}
\right\|_{L^\frac{p}{t}(X)}\right\}^\frac{1}{t}
\nonumber
\\
&\leq
2^{(\delta-1)k_0\widetilde{q}}\left\{
\sum_{k=-\infty}^{k_0-1}2^{(1-\delta)k\widetilde{q}t}\right.\nonumber\\
&\quad\left.\times\sum_{l=0}^{\infty}\left\|\left\{
\sum_{i\in\mathbb{N}}
\left[\left[\mu(B_{i,k})\right]^\frac{1}{p}
\left(O_{l}^{i,k}\right)
\mathbf{1}_{C_\rho2^{l+2}B_{i,k}}\right]^{\widetilde{q}t}
\right\}^\frac{1}{t}
\right\|_{L^{p}(X)}^t\right\}^\frac{1}{t}.
\end{align}
With the goal of employing Lemma~\ref{4l4-1},
observe that, by \eqref{oikl},
we conclude that,
for any $k\in\mathbb{Z}$, $i\in\mathbb{N}$,
and $l\in\mathbb{Z}_+$,
\begin{align*}
&\left\|\left\{\left[\mu(B_{i,k})\right]^\frac{1}{p}
\left(O_{l}^{i,k}\right)\mathbf{1}_{C_\rho2^{l+2}B_{i,k}}\right\}^
{\widetilde{q}}\right\|_{L^\frac{r}{\widetilde{q}}(X)}
\\
&\quad=\left\|\left\{\left[\mu(B_{i,k})\right]^\frac{1}{p}
\left[s_l\mathcal{M}_\rho\left(a^{i,k}_l\right)
+\sum_{u=l}^{\infty}s_u
\mathcal{M}_\rho\left(b^{i,k}_{u,l}\right)
\right]\mathbf{1}_{C_\rho2^{l+2}B_{i,k}}\right\}^
{\widetilde{q}}\right\|_{L^\frac{r}{\widetilde{q}}(X)}
\\
&\quad\leq\left[\mu(B_{i,k})\right]^\frac{\widetilde{q}}{p}
\left[s_l\left\|\mathcal{M}_\rho\left(a^{i,k}_l\right)
\right\|_{L^r(X)}
+\sum_{u=l}^{\infty}s_u
\left\|\mathcal{M}_\rho\left(b^{i,k}_{u,l}\right)
\right\|_{L^r(X)}\right]^{\widetilde{q}},
\end{align*}
which, combined with Lemma~\ref{M},
the fact that each $a^{i,k}_l$
and $b_{u,l}^{i,k}$ are supported
in $2^{l+1}B_{i,k}$, Definition~\ref{atom}(ii),
and the reverse doubling condition \eqref{RDcondition},
further implies that
\begin{align*}
&\left\|\left\{\left[\mu(B_{i,k})\right]^\frac{1}{p}
\left(O_{l}^{i,k}\right)\mathbf{1}_{C_\rho2^{l+2}B_{i,k}}\right\}^
{\widetilde{q}}\right\|_{L^\frac{r}{\widetilde{q}}(X)}
\\
&\quad\lesssim\left[\mu(B_{i,k})\right]^\frac{\widetilde{q}}{p}
\left[s_l\left\|a^{i,k}_l
\right\|_{L^r(X)}
+\sum_{u=l}^{\infty}s_u
\left\|b^{i,k}_{u,l}
\right\|_{L^r(X)}\right]^{\widetilde{q}}\\
&\quad\lesssim\left[\mu(B_{i,k})\right]^\frac{\widetilde{q}}{p}
\left[\mu(2^{l+1}B_{i,k})\right]^{\frac{\widetilde{q}}{r}-\frac{\widetilde{q}}{p}}
\left(s_l+\sum_{u=l}^{\infty}s_u\right)^{\widetilde{q}}
\\
&\quad\lesssim2^{-l\widetilde{q}(\eta+\frac{\kappa}{p})}
\left[\mu(2^{l+1}B_{i,k})\right]^\frac{\widetilde{q}}{r}
\lesssim2^{-l\widetilde{q}(\eta+\frac{\kappa}{p})}
\left[\mu(C_\rho2^{l+2}B_{i,k})\right]^\frac{\widetilde{q}}{r}.
\end{align*}
From this, \eqref{83I11}, and Lemma~\ref{4l4-1}
[used here with $r$ replaced by $\frac{r}{\widetilde{q}}\in(1,\infty]$,
$\eta_i:=1$, and $\xi_i:=2^{l\widetilde{q}
(\eta+\frac{\kappa}{p})}\{[\mu(B_{i,k})]^\frac{1}{p}
(O_{l}^{i,k})\mathbf{1}_{C_\rho2^{l+2}B_{i,k}}\}^
{\widetilde{q}}$ (modulo a harmless constant)],
it follows that
\begin{align*}
\mathrm{I_{1,1}}
&\lesssim2^{(\delta-1)k_0\widetilde{q}}\left\{
\sum_{k=-\infty}^{k_0-1}2^{(1-\delta)k\widetilde{q}t}
\sum_{l=0}^{\infty}\left\|\left(
\sum_{i\in\mathbb{N}}
2^{-l t\widetilde{q}(\eta+\frac{\kappa}{p})}
\mathbf{1}_{C_\rho2^{l+2}B_{i,k}}
\right)^\frac{1}{t}
\right\|_{L^{p}(X)}^t\right\}^\frac{1}{t}
\\
&\sim2^{(\delta-1)k_0\widetilde{q}}\left(
\sum_{k=-\infty}^{k_0-1}2^{(1-\delta)k\widetilde{q}t}\sum_{l=0}^{\infty}\left\|
\sum_{i\in\mathbb{N}}
2^{-l t\widetilde{q}(\eta+\frac{\kappa}{p})}
\mathbf{1}_{C_\rho2^{l+2}B_{i,k}}
\right\|_{L^\frac{p}{t}(X)}\right)^\frac{1}{t}
\\
&\sim2^{(\delta-1)k_0\widetilde{q}}\left(
\sum_{k=-\infty}^{k_0-1}2^{(1-\delta)k\widetilde{q}t}
\sum_{l=0}^{\infty}
2^{-l t\widetilde{q}(\eta+\frac{\kappa}{p})}\left\|\sum_{i\in\mathbb{N}}
\mathbf{1}_{\tau (2C_\rho)^{j_0}B_{i,k}}
\right\|_{L^\frac{p}{t}(X)}\right)^\frac{1}{t}.
\end{align*}
By this, Lemma~\ref{4r2} with $s\in(\frac{n}{(\eta+\frac{\kappa}{p})t\widetilde{q}},
\min\{1,\,\frac{p}{t}\})$ and $\tau:=(2C_\rho)^{1-j_0}2^{l+1}>1$
[where $j_0$ is the same as in \eqref{L0}],
the fact that $-t\widetilde{q}(\eta+\frac{\kappa}{p})+\frac{n}{s}<0$,
\eqref{L0}, and \eqref{upperdoub},
we obtain
\begin{align*}
\mathrm{I_{1,1}}
&\lesssim2^{(\delta-1)k_0\widetilde{q}}\left(
\sum_{k=-\infty}^{k_0-1}2^{(1-\delta)k\widetilde{q}t}
\sum_{l=0}^{\infty}
2^{-l[t\widetilde{q}(\eta+\frac{\kappa}{p})-\frac{n}{s}]}
\left\|\sum_{i\in\mathbb{N}}
\mathbf{1}_{(2C_\rho)^{j_0}B_{i,k}}
\right\|_{L^\frac{p}{t}(X)}\right)^\frac{1}{t}
\\
&\sim2^{(\delta-1)k_0\widetilde{q}}\left\{
\sum_{k=-\infty}^{k_0-1}2^{(1-\delta)k\widetilde{q}t}
\left[\sum_{i\in\mathbb{N}}
\mu\left((2C_\rho)^{j_0}B_{i,k}\right)\right]^\frac{t}{p}\right\}^\frac{1}{t}
\\
&\lesssim2^{(\delta-1)k_0\widetilde{q}}\left\{
\sum_{k=-\infty}^{k_0-1}2^{(1-\delta)k\widetilde{q}t}
\left[\sum_{i\in\mathbb{N}}
\mu\left(B_{i,k}\right)\right]^\frac{t}{p}\right\}^\frac{1}{t}.
\end{align*}
Let $\eta_1\in(1,(1-\delta)\widetilde{q})$.
On the one hand, if $\frac{q}{t}\in(1,\infty)$,
then, from H\"older's inequality,
we infer that
\begin{align}
\label{1012}
\mathrm{I}_{1,1}
&\sim2^{(\delta-1)k_0\widetilde{q}}\left\{
\sum_{k=-\infty}^{k_0-1}2^{(1-\delta)k\widetilde{q}t}
2^{-kt\eta_1}2^{kt\eta_1}
\left[\sum_{i\in\mathbb{N}}
\mu\left(B_{i,k}\right)\right]^\frac{t}{p}\right\}^\frac{1}{t}
\nonumber\\
&\leq2^{(\delta-1)k_0\widetilde{q}}\left\{
\sum_{k=-\infty}^{k_0-1}
2^{\frac{[(1-\delta)\widetilde{q}-\eta_1]
ktq}{q-t}}\right\}^\frac{q-t}{qt}
\left\{\sum_{k=-\infty}^{k_0-1}
2^{kq\eta_1}
\left[\sum_{i\in\mathbb{N}}
\mu\left(B_{i,k}\right)\right]^\frac{q}{p}\right\}^\frac{1}{q}
\nonumber\\
&\sim2^{-k_0\eta_1}
\left\{\sum_{k=-\infty}^{k_0-1}
2^{kq\eta_1}
\left[\sum_{i\in\mathbb{N}}
\mu\left(B_{i,k}\right)\right]^\frac{q}{p}\right\}^\frac{1}{q}.
\end{align}
On the other hand, if $\frac{q}{t}\in(0,1]$,
then, by Lemma \ref{discreteinequality}
with $r$ therein replaced by $\frac{q}{t}$, we find that
\begin{align}
\label{10121}
\mathrm{I}_{1,1}
&\lesssim2^{(\delta-1)k_0\widetilde{q}}\left\{
\sum_{k=-\infty}^{k_0-1}2^{(1-\delta)k\widetilde{q}t}
\left[\sum_{i\in\mathbb{N}}
\mu\left(B_{i,k}\right)\right]^\frac{t}{p}\right\}^{\frac{q}{t}\times\frac{1}{q}}
\nonumber\\
&\leq
2^{(\delta-1)k_0\widetilde{q}}\left(
\sum_{k=-\infty}^{k_0-1}2^{(1-\delta)k\widetilde{q}q}
\left[\sum_{i\in\mathbb{N}}
\mu\left(B_{i,k}\right)\right]^\frac{q}{p}\right)^\frac{1}{q}
\nonumber\\
&=2^{(\delta-1)k_0\widetilde{q}}\left(
\sum_{k=-\infty}^{k_0-1}2^{kq\eta_1}2^{kq[(1-\delta)\widetilde{q}-\eta_1]}
\left[\sum_{i\in\mathbb{N}}
\mu\left(B_{i,k}\right)\right]^\frac{q}{p}\right)^\frac{1}{q}
\nonumber\\
&\leq2^{-k_0\eta_1}
\left\{\sum_{k=-\infty}^{k_0-1}
2^{kq\eta_1}
\left[\sum_{i\in\mathbb{N}}
\mu\left(B_{i,k}\right)\right]^\frac{q}{p}\right\}^\frac{1}{q},
\end{align}
which, together with \eqref{1012}, completes
the proof of \eqref{82I1} when $j=1$.

To estimate $\mathrm{I_{1,2}}$,
from \eqref{4l3.1},
the definition of $s_l$ in Lemma~\ref{L6.2},
an argument similar to that used in the estimation of
\eqref{4ee6}, and the reverse doubling condition \eqref{RDcondition},
we deduce that, for any $k\in\mathbb{Z}$, $i\in\mathbb{N}$, and
$x\in X$,
\begin{align}
\label{88}
\mathrm{J}^{i,k}_{2,l}(x)
&\lesssim\left(s_l+\sum_{u=l}^{\infty}s_u\right)
\left[\frac{2^{l+1}r_{i,k}}{\rho(x_{i,k},x)}\right]^\beta
\frac{[\mu(2^{l+1}B_{i,k})]^{1-\frac{1}{p}}}{V_\rho(x_{i,k},x)}
\mathbf{1}_{(C_\rho2^{l+2}B_{i,k})^\complement}(x)
\nonumber\\
&\lesssim2^{-l\eta}
\left[\frac{2^{l+1}r_{i,k}}{\rho(x_{i,k},x)}\right]^\beta
\frac{[\mu(2^{l+1}B_{i,k})]^{1-\frac{1}{p}}}{V_\rho(x_{i,k},x)}
\mathbf{1}_{(C_\rho2^{l+2}B_{i,k})^\complement}(x)
\nonumber\\
&\lesssim2^{-l\eta}
\left[\mu(2^{l+1}B_{i,k})\right]^{-\frac{1}{p}}
\left[\mathcal{M}_\rho(\mathbf{1}_{2^{l+1}
B_{i,k}})(x)\right]^{1+\frac{\beta}{n}}
\mathbf{1}_{(C_\rho2^{l+2}B_{i,k})^\complement}(x)
\nonumber\\
&\lesssim2^{-l(\eta+\frac{\kappa}{p})}
\left[\mu(B_{i,k})\right]^{-\frac{1}{p}}
\left[\mathcal{M}_\rho(\mathbf{1}_{2^{l+1}B_{i,k}})(x)
\right]^{1+\frac{\beta}{n}}.
\end{align}
Suppose that
$t\in(0,\min\{p,\,\frac{1}{1+\frac{\beta}{n}}\})$,
$q_1\in(\max\{\frac{n}{(\eta+\frac{\kappa}{p})t},\,
\frac{1}{t(1+\frac{\beta}{n})}\},\frac{1}{t})$,
$\delta\in(0,1-\frac{1}{q_1})$, and
$s\in(\frac{n}{(\eta+\frac{\kappa}{p})tq_1},\min\{1,\,\frac{p}{t}\})$.
Then, by Chebyshev's inequality,
H\"older's inequality with exponent $q_1>1$,
\eqref{lamdik}, and \eqref{88},
we conclude that
\begin{align*}
\mathrm{I_{1,2}}
&\lesssim2^{-k_0q_1}\left\|\left(\sum_{k=-\infty}^{k_0-1}
2^{k\delta}2^{-k\delta}\sum_{i\in\mathbb{N}}
\sum_{l=0}^{\infty}\lambda_{i,k}\mathrm{J}^{i,k}_{2,l}
\right)^{q_1}\right\|_{L^p(X)}
\\
&\leq2^{-k_0q_1}\left\|\left(\sum_{k=-\infty}^{k_0-1}
2^{k\delta q_1'}\right)^\frac{q_1}{q_1'}
\left[\sum_{k=-\infty}^{k_0-1}
2^{-k\delta q_1}\left(\sum_{i\in\mathbb{N}}
\sum_{l=0}^{\infty}\lambda_{i,k}\mathrm{J}^{i,k}_{2,l}
\right)^{q_1}\right]\right\|_{L^p(X)}
\\
&\sim2^{k_0q_1(\delta-1)}\left\|
\sum_{k=-\infty}^{k_0-1}
2^{-k\delta q_1}\left(\sum_{i\in\mathbb{N}}
\sum_{l=0}^{\infty}\lambda_{i,k}\mathrm{J}^{i,k}_{2,l}
\right)^{q_1}\right\|_{L^p(X)}
\\
&\lesssim2^{k_0q_1(\delta-1)}\left\|
\sum_{k=-\infty}^{k_0-1}
2^{kq_1(1-\delta)}\left(\sum_{i\in\mathbb{N}}
\sum_{l=0}^{\infty}2^{-l(\eta+\frac{\kappa}{p})}
\left[\mathcal{M}_\rho(\mathbf{1}_{2^{l+1}B_{i,k}})
\right]^{1+\frac{\beta}{n}}
\right)^{q_1}\right\|_{L^p(X)}
\\
&=2^{k_0q_1(\delta-1)}\left\|\left\{
\sum_{k=-\infty}^{k_0-1}
2^{kq_1(1-\delta)}\left(\sum_{i\in\mathbb{N}}
\sum_{l=0}^{\infty}2^{-l(\eta+\frac{\kappa}{p})}
\left[\mathcal{M}_\rho(\mathbf{1}_{2^{l+1}B_{i,k}})
\right]^{1+\frac{\beta}{n}}
\right)^{q_1}\right\}^t\right\|_{L^\frac{p}{t}(X)}^\frac{1}{t}.
\end{align*}
From this, Lemma \ref{discreteinequality}
(used once with $t<1$ and again with $q_1t<1$), and
Minkowski's inequality with the norm $\|\cdot\|_{L^\frac{p}{t}(X)}$,
we infer that
\begin{align*}
\mathrm{I_{1,2}}
&\lesssim2^{k_0q_1(\delta-1)}\left\|
\sum_{k=-\infty}^{k_0-1}
2^{kq_1t(1-\delta)}\left(\sum_{i\in\mathbb{N}}
\sum_{l=0}^{\infty}2^{-l(\eta+\frac{\kappa}{p})}
\left[\mathcal{M}_\rho(\mathbf{1}_{2^{l+1}B_{i,k}})
\right]^{1+\frac{\beta}{n}}
\right)^{q_1t}\right\|_{L^\frac{p}{t}(X)}^\frac{1}{t}
\\
&\lesssim2^{k_0q_1(\delta-1)}\left\{
\sum_{k=-\infty}^{k_0-1}
2^{kq_1t(1-\delta)}\left\|\left(\sum_{i\in\mathbb{N}}
\sum_{l=0}^{\infty}2^{-l(\eta+\frac{\kappa}{p})}
\left[\mathcal{M}_\rho(\mathbf{1}_{2^{l+1}B_{i,k}})
\right]^{1+\frac{\beta}{n}}
\right)^{q_1t}\right\|_{L^\frac{p}{t}(X)}\right\}^\frac{1}{t}
\\
&\lesssim2^{k_0q_1(\delta-1)}\left\{
\sum_{k=-\infty}^{k_0-1}
2^{kq_1t(1-\delta)}\left\|\sum_{i\in\mathbb{N}}
\sum_{l=0}^{\infty}2^{-lq_1t(\eta+\frac{\kappa}{p})}
\left[\mathcal{M}_\rho(\mathbf{1}_{2^{l+1}B_{i,k}})
\right]^{q_1t(1+\frac{\beta}{n})}
\right\|_{L^\frac{p}{t}(X)}\right\}^\frac{1}{t}
\\
&\lesssim2^{k_0q_1(\delta-1)}\left\{
\sum_{k=-\infty}^{k_0-1}
2^{kq_1t(1-\delta)}\sum_{l=0}^{\infty}2^{-lq_1t(\eta+\frac{\kappa}{p})}
\left\|\sum_{i\in\mathbb{N}}
\left[\mathcal{M}_\rho(\mathbf{1}_{2^{l+1}B_{i,k}})
\right]^{q_1t(1+\frac{\beta}{n})}
\right\|_{L^\frac{p}{t}(X)}\right\}^\frac{1}{t}.
\end{align*}
Using this, Lemma~\ref{Fefferman--Stein} with $u=q_1t(1+\frac{\beta}{n})>1$
and $p$ replaced by $\frac{pu}{t}$,
Lemma~\ref{4r2} with $s\in(\frac{n}{(\eta+
\frac{\kappa}{p})tq_1},\min\{1,\,\frac{p}{t}\})$
and $\tau=(2C_\rho)^{-j_0}2^{l+1}>1$
[where $j_0$ is the same as in \eqref{L0}],
\eqref{L0}, and \eqref{doublingcond},
we find that
\begin{align*}
\mathrm{I_{1,2}}
&\lesssim2^{k_0q_1(\delta-1)}\left\{
\sum_{k=-\infty}^{k_0-1}
2^{kq_1t(1-\delta)}\sum_{l=0}^{\infty}2^{-lq_1t(\eta+\frac{\kappa}{p})}
\left\|\sum_{i\in\mathbb{N}}
\mathbf{1}_{2^{l+1}B_{i,k}}
\right\|_{L^\frac{p}{t}(X)}\right\}^\frac{1}{t}
\\
&\sim2^{k_0q_1(\delta-1)}\left\{
\sum_{k=-\infty}^{k_0-1}
2^{kq_1t(1-\delta)}\sum_{l=0}^{\infty}2^{-lq_1t(\eta+\frac{\kappa}{p})}
\left\|\sum_{i\in\mathbb{N}}
\mathbf{1}_{\tau (2C_\rho)^{j_0}B_{i,k}}
\right\|_{L^\frac{p}{t}(X)}\right\}^\frac{1}{t}
\\
&\lesssim2^{k_0q_1(\delta-1)}\left\{
\sum_{k=-\infty}^{k_0-1}
2^{kq_1t(1-\delta)}\left[\sum_{l=0}^{\infty}
2^{-l[q_1t(\eta+\frac{\kappa}{p})-\frac{n}{s}]}
\right]
\left\|\sum_{i\in\mathbb{N}}
\mathbf{1}_{(2C_\rho)^{j_0}B_{i,k}}
\right\|_{L^\frac{p}{t}(X)}\right\}^\frac{1}{t}
\\
&\lesssim2^{k_0q_1(\delta-1)}\left\{
\sum_{k=-\infty}^{k_0-1}
2^{kq_1t(1-\delta)}
\left[\sum_{i\in\mathbb{N}}
\mu\left((2C_\rho)^{j_0}B_{i,k}\right)\right]^\frac{t}{p}\right\}^\frac{1}{t}
\\
&\lesssim2^{k_0q_1(\delta-1)}\left\{
\sum_{k=-\infty}^{k_0-1}
2^{kq_1t(1-\delta)}
\left[\sum_{i\in\mathbb{N}}
\mu\left(B_{i,k}\right)\right]^\frac{t}{p}\right\}^\frac{1}{t},
\end{align*}
which, combined with arguments similar to those used in the estimations
appearing in \eqref{1012} when $\frac{q}{t}\in(1,\infty)$ and
\eqref{10121} when $\frac{q}{t}\in(0,1]$,
further implies that \eqref{82I1}
holds when $j=2$ with $\eta_2\in(1,q_1(1-\delta))$.

Now, we deal with estimating $\mathrm{I_{2,1}}$.
From Chebyshev's inequality, Lemma \ref{discreteinequality},
\eqref{lamdik}, \eqref{oikl-ineq},
and Minkowski's inequality,
it follows that,
for any
$t\in(\frac{n}{\eta+\frac{\kappa}{p}},p),$
$\widetilde{q}\in(\frac{n}{(\eta+\frac{\kappa}{p})t},1)$,
$s\in(\frac{n}{(\eta+\frac{\kappa}{p})t
\widetilde{q}},\min\{1,\,\frac{p}{t}\})$,
and $M\in\mathbb{Z}$ such that
$2^{-M}\leq(2C_\rho)^{j_0}<2^{-M+1}$,
\begin{align*}
\mathrm{I_{2,1}}
&\lesssim2^{-k_0\widetilde{q}}
\left\|\left(\sum_{k=k_0}^{\infty}
\sum_{i\in\mathbb{N}}
\sum_{l=0}^{\infty}\lambda_{i,k}
\mathrm{J}^{i,k}_{1,l}\right)^{\widetilde{q}}\right\|_{L^p(X)}
\leq2^{-k_0\widetilde{q}}
\left\|\sum_{k=k_0}^{\infty}
\sum_{i\in\mathbb{N}}
\sum_{l=0}^{\infty}\left(\lambda_{i,k}
\mathrm{J}^{i,k}_{1,l}\right)^{\widetilde{q}}\right\|_{L^p(X)}
\nonumber\\
&\lesssim2^{-k_0\widetilde{q}}
\left\|\sum_{k=k_0}^{\infty}2^{k\widetilde{q}}
\sum_{i\in\mathbb{N}}
\sum_{l=0}^{\infty}
\left\{\left[\mu(B_{i,k})\right]^\frac{1}{p}
\left(O_{l}^{i,k}\right)
\mathbf{1}_{C_\rho2^{l+2}B_{i,k}}
\right\}^{\widetilde{q}}\right\|_{L^p(X)}
\\
&=2^{-k_0\widetilde{q}}
\left\|\left(\sum_{k=k_0}^{\infty}2^{k\widetilde{q}}
\sum_{i\in\mathbb{N}}
\sum_{l=0}^{\infty}
\left\{\left[\mu(B_{i,k})\right]^\frac{1}{p}
\left(O_{l}^{i,k}\right)
\mathbf{1}_{C_\rho2^{l+2}B_{i,k}}
\right\}^{\widetilde{q}}\right)^{t}
\right\|_{L^\frac{p}{t}(X)}^\frac{1}{t}
\\
&\lesssim2^{-k_0\widetilde{q}}
\left\|\sum_{k=k_0}^{\infty}2^{k\widetilde{q}t}
\sum_{i\in\mathbb{N}}
\sum_{l=0}^{\infty}
\left\{\left[\mu(B_{i,k})\right]^\frac{1}{p}
\left(O_{l}^{i,k}\right)
\mathbf{1}_{C_\rho2^{l+2}B_{i,k}}
\right\}^{\widetilde{q}t}\right\|_{L^\frac{p}{t}(X)}^\frac{1}{t}
\\
&\leq2^{-k_0\widetilde{q}}\left(
\sum_{k=k_0}^{\infty}2^{k\widetilde{q}t}
\left\|\sum_{i\in\mathbb{N}}
\sum_{l=0}^{\infty}
\left\{\left[\mu(B_{i,k})\right]^\frac{1}{p}
\left(O_{l}^{i,k}\right)
\mathbf{1}_{C_\rho2^{l+2}B_{i,k}}
\right\}^{\widetilde{q}t}\right\|_{L^\frac{p}{t}(X)}\right)^\frac{1}{t},
\end{align*}
which, together with an argument similar to that
used in the estimation of $\mathrm{I_{1,1}}$,
Lemmas~\ref{4l4-1} and~\ref{4r2},
and \eqref{L0},
further implies that
\begin{align*}
\mathrm{I_{2,1}}
&\lesssim2^{-k_0\widetilde{q}}\left(
\sum_{k=k_0}^{\infty}2^{k\widetilde{q}t}
\left\|\sum_{i\in\mathbb{N}}
\sum_{l=0}^{\infty}
2^{-lt\widetilde{q}(\eta+\frac{\kappa}{p})}
\mathbf{1}_{C_\rho2^{l+2}B_{i,k}}
\right\|_{L^\frac{p}{t}(X)}\right)^\frac{1}{t}
\\
&\leq2^{-k_0\widetilde{q}}\left(
\sum_{k=k_0}^{\infty}2^{k\widetilde{q}t}
\sum_{l=0}^{\infty}
2^{-lt\widetilde{q}(\eta+\frac{\kappa}{p})}
\left\|\sum_{i\in\mathbb{N}}
\mathbf{1}_{C_\rho2^{l+2}B_{i,k}}
\right\|_{L^\frac{p}{t}(X)}\right)^\frac{1}{t}
\\
&\lesssim2^{-k_0\widetilde{q}}\left(
\sum_{k=k_0}^{\infty}2^{k\widetilde{q}t}
\sum_{l=0}^{\infty}
2^{-l[t\widetilde{q}(\eta+\frac{\kappa}{p})-\frac{n}{s}]}
\left\|\sum_{i\in\mathbb{N}}
\mathbf{1}_{2^{-M}B_{i,k}}
\right\|_{L^\frac{p}{t}(X)}\right)^\frac{1}{t}
\\
&\sim2^{-k_0\widetilde{q}}
\left\{\sum_{k=k_0}^{\infty}2^{k\widetilde{q}t}
\left[\sum_{i\in\mathbb{N}}\mu(B_{i,k})\right]^\frac{t}{p}
\right\}^\frac{1}{t}
\end{align*}
with $O_{l}^{i,k}$ as in \eqref{oikl}.
This, combined with arguments similar to those in the estimations
of \eqref{1012} when $\frac{q}{t}\in(1,\infty)$ and
\eqref{10121} when $\frac{q}{t}\in(0,1]$, further implies
\eqref{82I2} when $j=1$ with
$\delta_1\in(\widetilde{q},1)$.

Next, we deal with $\mathrm{I_{2,2}}$.
Applying an argument similar to that
used in the estimation of \eqref{88},
\eqref{RDcondition}, and \eqref{4ee6},
we obtain,
for any $x\in X$,
\begin{align*}
\mathrm{J}^{i,k}_{2,l}(x)
&\lesssim2^{-l\eta}\left[
\frac{2^{l+1}r_{i,k}}{\rho(x_{i,k},x)}\right]^\beta
\frac{[\mu(2^{l+1}B_{i,k})]^{1-\frac{1}{p}}}{V_\rho(x_{i,k},x)}
\mathbf{1}_{(C_\rho2^{l+2}B_{i,k})^\complement}(x)
\\
&\lesssim2^{l(\beta-\eta)}\left[
\frac{r_{i,k}}{\rho(x_{i,k},x)}\right]^\beta
\frac{[2^{l\kappa}\mu(B_{i,k})]^{1-\frac{1}{p}}}{V_\rho(x_{i,k},x)}
\mathbf{1}_{(C_\rho2^{l+2}B_{i,k})^\complement}(x)
\\
&\lesssim2^{l[\beta-\eta+\kappa(1-\frac{1}{p})]}
\left[\mu(B_{i,k})\right]^{-\frac{1}{p}}
\left[\mathcal{M}_\rho(
\mathbf{1}_{B_{i,k}})(x)\right]^{1+\frac{\beta}{n}}
\mathbf{1}_{(C_\rho2^{l+2}B_{i,k})^\complement}(x).
\end{align*}
From this, Chebyshev's inequality,
Lemma \ref{discreteinequality} twice,
\eqref{lamdik}, and
$\eta\in(\beta+\kappa(1-\frac{1}{p}),\infty)$,
we deduce that,
for any given $t\in(\frac{1}{(1+\frac{\beta}{n})p},1)$,
$\delta_2\in(t,1)$, and
$L\in(\max\{1,\,\frac{1}{(1+\frac{\beta}{n})qt}\},\infty)$,
\begin{align*}
\mathrm{I_{2,2}}
&\lesssim2^{-k_0t}\left\|\left(
\sum_{k=k_0}^\infty\sum_{i\in\mathbb{N}}
\sum_{l=0}^\infty\lambda_{i,k}\mathrm{J}^{i,k}_{2,l}
\right)^t\right\|_{L^{p}(X)}
\leq2^{-k_0t}\left\|
\sum_{k=k_0}^\infty\sum_{i\in\mathbb{N}}
\sum_{l=0}^\infty\left(\lambda_{i,k}\mathrm{J}^{i,k}_{2,l}
\right)^t\right\|_{L^{p}(X)}
\nonumber\\
&\lesssim2^{-k_0t}\left\|
\sum_{k=k_0}^\infty\sum_{i\in\mathbb{N}}
\sum_{l=0}^\infty2^{kt}
2^{lt[\beta-\eta+\kappa(1-\frac{1}{p})]}\left[\mathcal{M}_\rho
(\mathbf{1}_{B_{i,k}})
\right]^{t(1+\frac{\beta}{n})}\right\|_{L^{p}(X)}
\\
&\sim2^{-k_0t}\left\|
\sum_{k=k_0}^\infty2^{kt}\sum_{i\in\mathbb{N}}
\left[\mathcal{M}_\rho
(\mathbf{1}_{B_{i,k}})
\right]^{t(1+\frac{\beta}{n})}\right\|_{L^{p}(X)}
\\
&=2^{-k_0t}\left\|\left\{
\sum_{k=k_0}^\infty2^{kt}\sum_{i\in\mathbb{N}}
\left[\mathcal{M}_\rho
(\mathbf{1}_{B_{i,k}})
\right]^{t(1+\frac{\beta}{n})}\right\}^\frac{1}{t(1+\frac{\beta}{n})}
\right\|_{L^{pt(1+\frac{\beta}{n})}(X)}^{t(1+\frac{\beta}{n})}
\\
&\leq2^{-k_0t}\left\|
\sum_{k=k_0}^\infty2^\frac{k}{1+\frac{\beta}{n}}\left\{\sum_{i\in\mathbb{N}}
\left[\mathcal{M}_\rho
(\mathbf{1}_{B_{i,k}})
\right]^{t(1+\frac{\beta}{n})}\right\}^\frac{1}{t(1+\frac{\beta}{n})}
\right\|_{L^{pt(1+\frac{\beta}{n})}(X)}^{t(1+\frac{\beta}{n})}.
\end{align*}
By this, Minkowski's inequality,
Lemma~\ref{Fefferman--Stein},
and H\"older's inequality with exponent $(1+\frac{\beta}{n})qtL$,
we further conclude that
\begin{align*}
\mathrm{I_{2,2}}
&\lesssim2^{-k_0t}\left(
\sum_{k=k_0}^\infty2^\frac{k}{1+\frac{\beta}{n}}
\left\|\left\{\sum_{i\in\mathbb{N}}
\left[\mathcal{M}_\rho
(\mathbf{1}_{B_{i,k}})
\right]^{t(1+\frac{\beta}{n})}\right\}^\frac{1}{t(1+\frac{\beta}{n})}
\right\|_{L^{pt(1+\frac{\beta}{n})}(X)}\right)^{t(1+\frac{\beta}{n})}
\\
&\lesssim2^{-k_0t}\left\{
\sum_{k=k_0}^\infty2^\frac{k}{1+\frac{\beta}{n}}
\left\|\left(\sum_{i\in\mathbb{N}}
\mathbf{1}_{B_{i,k}}\right)^\frac{1}{t(1+\frac{\beta}{n})}
\right\|_{L^{pt(1+\frac{\beta}{n})}(X)}\right\}^{t(1+\frac{\beta}{n})}
\\
&\leq2^{-tk_0}\left\{\sum_{k=k_0}^{\infty}
2^\frac{(1-\frac{\delta_2}{t})k}{(1+\frac{\beta}{n})L}
2^\frac{\delta_2k}{(1+\frac{\beta}{n})tL}
\left\|\sum_{i\in\mathbb{N}}\mathbf{1}_{B_{i,k}}\right\|
^\frac{1}{(1+\frac{\beta}{n})tL}_{L^p(X)}
\right\}^{(1+\frac{\beta}{n})tL}
\\
&\lesssim 2^{-tk_0}2^{(t-\delta_2)k_0}
\left[\sum_{k=k_0}^{\infty}2^{k\delta_2q}
\left\|\sum_{i\in\mathbb{N}}\mathbf{1}_{B_{i,k}}
\right\|^q_{L^p(X)}\right]^\frac{1}{q}
\\
&\lesssim2^{-\delta_2k_0}
\left\{\sum_{k=k_0}^{\infty}2^{k\delta_2q}
\left[\sum_{i\in\mathbb{N}}
\mu(B_{i,k})\right]^\frac{q}{p}\right\}^\frac{1}{q},
\end{align*}
which is precisely \eqref{82I2} when $j=2$.
This finishes the proof of the estimates in \eqref{82I1} and \eqref{82I2}
and hence Theorem~\ref{thmHpqm}.
\end{proof}

\section{Littlewood--Paley Function
Characterizations}
\label{section5.6}

In Chapter~\ref{Chapter3}, we introduced the Hardy space
$H^p(X)$ defined via the Lusin area function and established its
atomic and Littlewood--Paley function characterizations.
In this section, we generalize these results by establishing
the atomic and the Littlewood--paley function
characterizations of Hardy--Lorentz spaces
on ultra-RD-spaces, $(X,\mathbf{q},\mu)$,
for the optimal ranges
$$
p\in\left(\frac{n}{n+{\mathrm{ind\,}}(X,\mathbf{q})},1\right]
\ \ \text{and}\ \
q\in(0,\infty].
$$
In contrast to Sections~\ref{section5.1}--\ref{section5.5},
we use $(\mathring{{\mathcal{G}}}^\varepsilon_0(\beta,\gamma))'$
as the considered space of distributions
(see Definition~\ref{1558} for the definition of this
space). Throughout this section, we assume that,
for some (hence all) $\rho\in\mathbf{q}$,
$\mathrm{diam}_\rho(X)=\infty$, which,
by Remark~\ref{1031}, is equivalent to
the demand that $\mu(X)=\infty$.

\subsection{Hardy--Lorentz Spaces
via Lusin Area Functions}

Let $(X,\mathbf{q},\mu)$ be a quasi-ultrametric space
of homogeneous type
and
$0<\beta,\gamma<\varepsilon\preceq\mathrm{ind\,}(X,\mathbf{q})$.
Assume that $\{\mathcal{S}_{2^{-k}}\}_{k\in\mathbb{Z}}$
is a homogeneous approximation
of the identity of order $\varepsilon$ as in Definition~\ref{DefHATI} and,
for any $k\in\mathbb{Z}$,
$\mathcal{D}_k:=\mathcal{S}_{2^{-k}}-\mathcal{S}_{2^{-k+1}}$.
For any $f\in(\mathring{{\mathcal{G}}}^\varepsilon_0(\beta,\gamma))'$,
let the Lusin area function $S_\theta(f)$
with aperture $\theta\in(0,\infty)$,
the Littlewood--Paley $g$-function $g(f)$,
and the Littlewood--Paley $g_\lambda^*$-function
$g_\lambda^*(f)$ with $\lambda\in(0,\infty)$
of $f$ be the same as, respectively, in
\eqref{Def-S-f}, \eqref{Def-g(f)(x)},
and \eqref{Def-glamda*}.
As before, we denote $S_1$ simply by $S$.
Recall from Remark~\ref{rb-375} that we may assume
$S_\theta(f)$, $g(f)$, and $g_\lambda^*(f)$ are constructed
using a regularized quasi-ultrametric belonging to $\mathbf{q}$
as in Definition~\ref{D2.3}.
In particular, the use of a regularized quasi-ultrametric
guarantees that $S_\theta(f)$, $g(f)$, and $g_\lambda^*(f)$
are $\mu$-measurable; see Lemma~\ref{1955}.

Now, we first introduce the following concept of
Lusin-area-function-based
Hardy--Lorentz spaces.

\begin{definition}
\label{HLviaS}
Let $(X,\mathbf{q},\mu)$ be a quasi-ultrametric space
of homogeneous type with upper dimension $n\in(0,\infty)$,
$p\in(0,\infty)$,
$q\in(0,\infty]$,
and
$n(\frac{1}{p}-1)_+<\beta,\gamma<\varepsilon\preceq\mathrm{ind\,}(X,\mathbf{q})$.
The \emph{Lusin-area-function-based Hardy--Lorentz space $H^{p,q}(X)$}
\index[words]{Lusin-area-function-based Hardy--Lorentz space}
\index[symbols]{H@$H^{p,q}(X)$} is defined by setting
\begin{align*}
H^{p,q}(X):=\left\{f\in(\mathring{{\mathcal{G}}}^\varepsilon_0(\beta,\gamma))'
:\|f\|_{H^{p,q}(X)}:=\|S(f)\|_{L^{p,q}(X)}<\infty\right\}.
\end{align*}
\end{definition}

\begin{remark}\label{2214}
Let the notation be as in Definition~\ref{HLviaS}.
\begin{enumerate}
\item[\rm(i)]
Let $p\in(\frac{n}{n+\mathrm{ind\,}(X,\mathbf{q})},\infty)$.
Then $H^{p,q}(X)$ is a quasi-Banach space,
which can be seen from Theorem~\ref{atHL} and Remark~\ref{221}(ii) below.

\item[\rm(ii)]
Recall that, for any $p\in(0,\infty]$,
$L^{p,p}(X)=L^p(X)$, which implies that
$$
H^{p,p}(X)=H^{p}(X):=\left\{f\in
(\mathring{{\mathcal{G}}}^\varepsilon_0(\beta,\gamma))'
:\|f\|_{H^{p}(X)}:=\|S(f)\|_{L^{p}(X)}<\infty\right\}.
$$

\item[\rm(iii)]
From the assumption that $n(\frac{1}{p}-1)_+<\beta,\gamma<\infty$,
\eqref{1550}, and Remark~\ref{Rmtest}(v), we infer that,
for any $p\in(0,\frac{n}{n+\mathrm{ind}_H(X,\mathbf{q})}]$,
$H^{p,q}(X)$ is trivial
because the space
$\mathring{\mathcal{G}}^\varepsilon_0(\beta,\gamma)$
only contains constant functions.
\end{enumerate}
\end{remark}

Next,
we prove that $H^{p,q}(X)$
defined in Definition~\ref{HLviaS}
is independent of the
choice of homogeneous approximations
of the identity $\{\mathcal{S}_{2^{-k}}\}_{k\in\mathbb{Z}}$
in Definition~\ref{DefHATI}.
The main ingredient in this vein is the following lemma
which is a generalization of the Fefferman--Stein
vector-valued inequality in Lemma~\ref{Fefferman--Stein}
from Lebesgue spaces to Lorentz spaces.
This result extends \cite[Theorem 7.3]{zhy}
from nonatomic spaces of homogeneous type to
quasi-ultrametric spaces of homogeneous type.

\begin{lemma}
\label{HLf-s}
Let $(X,\mathbf{q},\mu)$ be a quasi-ultrametric space of homogeneous type and
assume $\rho\in\mathbf{q}$ is a regularized quasi-ultrametric
as in Definition~\ref{D2.3}.
\begin{enumerate}
\item[\textup{(i)}]
Let $p\in(1,\infty)$,
$u\in(1,\infty]$, and $q\in(0,\infty]$.
Then there exists a positive constant $C_1$ such that,
for any sequence of
$\mu$-measurable functions $\{f_j\}_{j=1}^\infty$,
\begin{align*}
\left\|\left\{\sum_{j=1}^\infty
\left[\mathcal{M}_\rho(f_j)\right]^{u}
\right\}^\frac{1}{u}\right\|_{L^{p,q}(X)}
\leq C_1\left\|\left(\sum_{j=1}^\infty
|f_j|^{u}\right)^\frac{1}{u}\right\|_{L^{p,q}(X)}
\end{align*}
with the usual modification when $q=\infty$ or $u=\infty$.
\item[\textup{(ii)}]
Let $p\in(0,\infty)$,
$u\in(0,\infty]$,
$q\in(0,\infty]$,
and $s\in(0,\min\{p,\,u\})$.
Then there exists a positive constant $C_2$
such that, for any sequence
of $\mu$-measurable functions $\{f_j\}_{j=1}^\infty$,
\begin{align*}
\left\|\left\{\sum_{j=1}^\infty
\left[\mathcal{M}_\rho(f_j)\right]^\frac{u}{s}
\right\}^\frac{1}{u}\right\|_{L^{p,q}(X)}
\leq C_2\left\|\left(\sum_{j=1}^\infty
|f_j|^\frac{u}{s}\right)^\frac{1}{u}\right\|_{L^{p,q}(X)}
\end{align*}
with the usual modification when $q=\infty$ or $u=\infty$.
\end{enumerate}
\end{lemma}

To show Lemma~\ref{HLf-s},
we need the following interpolation lemma,
whose proof is a slightly modification of \cite[Theorem~3.1]{jzzw2019}
with $p(\cdot)$ and $\mathbb{R}^n$
therein replaced, respectively,
by $p$ and $X$.

\begin{lemma}\label{160}
Let $(X,\mu)$ be a measure space.
Let $p\in(0,\infty)$, $p_1\in(0,1)$, $p_2\in(1,\infty)$,
and $q\in(0,\infty]$. Assume that $T$ is a sublinear operator defined on
$$
L^{p_1p}(X)+L^{p_2p}(X)
:=\left\{f:f=f_1+f_2,\ f_1\in L^{p_1p}(X),\ f_2\in L^{p_2p}(X)\right\}
$$
satisfying that there exist positive constants
$C_1$ and $C_2$ such that, for any $i\in\{1,2\}$ and
$f\in L^{p_ip}(X)$,
$$
\left\|T(f)\right\|_{L^{p_ip,\infty}(X)}\leq C_i
\left\|f\right\|_{L^{p_ip}(X)}.
$$
Then $T$ is bounded on $L^{p,q}(X)$.
\end{lemma}

\begin{proof}[Proof of Lemma~\ref{HLf-s}]
Let $\{f_j\}_{j\in\mathbb{N}}$
be a sequence of $\mu$-measurable functions on $X$.
We first prove (i).
Define the sublinear operator $T$ by setting,
for any $\mu$-measurable function $g$ and for any $x\in X$,
$$
T(g)(x):=\left\{\sum_{j=1}^\infty\left[
\mathcal{M}_\rho(g\eta_j)(x)\right]^u\right\}^\frac{1}{u},
$$
where, for any $j\in\mathbb{N}$ and $y\in X$,
\begin{align*}
\eta_j(y):=
\begin{cases}
\displaystyle
\frac{f_j(y)}{(\sum_{j=1}^\infty|f_j(y)|^u)^\frac{1}{u}}
\ &\text{if}\ \sum_{j=1}^\infty|f_j(y)|^u\neq0,\\
0&\text{if}\ \sum_{j=1}^\infty|f_j(y)|^u=0.
\end{cases}
\end{align*}
By Lemma~\ref{Fefferman--Stein},
we find that, for any $p_1\in(\frac{1}{p},1)$ and
$p_2\in(1,\infty)$ and for any $\mu$-measurable function $g$,
we have
$\|T(g)\|_{L^{p_ip}(X)}\lesssim\|g\|_{L^{p_ip}(X)}$,
where $i\in\{1,2\}$,
which, together with Lemma~\ref{160},
further implies that $T$ is bounded on $L^{p,q}(X)$.
From this and letting
$g:=(\sum_{j=1}^\infty|f_j|^u)^\frac{1}{u}$,
we deduce that
\begin{align*}
\left\|\left\{\sum_{j=1}^\infty
\left[\mathcal{M}_\rho(f_j)\right]^{u}
\right\}^\frac{1}{u}\right\|_{L^{p,q}(X)}=
\left\|T(g)\right\|_{L^{p,q}(X)}
\lesssim\left\|g\right\|_{L^{p,q}(X)}
=\left\|\left(\sum_{j=1}^\infty
|f_j|^{u}\right)^\frac{1}{u}\right\|_{L^{p,q}(X)}.
\end{align*}
This finishes the proof of (i).

We turn to show (ii).
By Proposition~\ref{r;p,q=pr,qr;r},
$s\in(0,\min\{p,\,u\})$,
and (i) with $p$, $u$, and $q$ therein
replaced, respectively, by $\frac{p}{s}$,
$\frac{u}{s}$, and $\frac{q}{s}$,
we conclude that (ii) holds,
which completes the proof of Lemma~\ref{160}.
\end{proof}

We are now ready to prove that  $H^{p,q}(X)$,
defined in Definition~\ref{HLviaS},
is independent of the
choice of homogeneous approximations
of the identity
$\{\mathcal{S}_{2^{-k}}\}_{k\in\mathbb{Z}}$
in Definition~\ref{DefHATI}.
The analogous result for
$H^p(X)=H^{p,p}(X)$ spaces was shown in Theorem~\ref{independentofATI}.

\begin{theorem}
\label{thm638}
Let $(X,\mathbf{q},\mu)$ be an ultra-RD-space with
upper dimension $n\in(0,\infty)$,
where $\mu$ is assumed to be a Borel-semiregular
measure on $X$ and $\mathrm{diam\,}(X)=\infty$.
Let $p\in(\frac{n}{n+\mathrm{ind\,}(X,\mathbf{q})},\infty)$,
$q\in(0,\infty]$,
and $\varepsilon,\beta,\gamma\in(0,\infty)$ be such that
$n(\frac{1}{p}-1)_+<\beta,\gamma<
\varepsilon\preceq\mathrm{ind\,}(X,\mathbf{q})$.
Assume that $\{\mathcal{S}_{2^{-k}}\}_{k\in\mathbb{Z}}$
and $\{\mathcal{P}_{2^{-k}}\}_{k\in\mathbb{Z}}$
are two homogeneous approximations
of the identity of order $\varepsilon$
in Definition~\ref{DefHATI}.
For any $k\in\mathbb{Z}$, let
$\mathcal{E}_k:=\mathcal{S}_{2^{-k}}-\mathcal{S}_{2^{-k+1}}$ and
$\mathcal{Q}_k:=\mathcal{P}_{2^{-k}}-\mathcal{P}_{2^{-k+1}}$.
Let $S_\mathcal{E}$ and $S_\mathcal{Q}$ be
the Lusin area functions (with aperture 1) as in \eqref{Def-S-f}
associated, respectively, with $\{\mathcal{E}_k\}_{k\in\mathbb{Z}}$ and
$\{\mathcal{Q}_k\}_{k\in\mathbb{Z}}$.
Then there exists a constant $C\in[1,\infty)$
such that, for any
$f\in(\mathring{{\mathcal{G}}}_0^\varepsilon(\beta,\gamma))'$,
$$
C^{-1}\left\|S_\mathcal{E}(f)\right\|_{L^{p,q}(X)}
\leq\left\|S_\mathcal{Q}(f)\right\|_{L^{p,q}(X)}
\leq C\left\|S_\mathcal{E}(f)\right\|_{L^{p,q}(X)}.
$$
\end{theorem}

\begin{proof}
Let
$\rho\in\mathbf{q}$ be a regularized quasi-ultrametric
as in Definition~\ref{D2.3}.
Also, let
$\varepsilon'\in(n(\frac{1}{p}-1)_+,\varepsilon)$
and
$r\in(\frac{n}{n+\varepsilon'},1]\cap(0,p)$.
By symmetry, we only need to prove that,
for any
$f\in(\mathring{{\mathcal{G}}}_0^\varepsilon(\beta,\gamma))'$,
\begin{align*}
\left\|S_\mathcal{E}(f)\right\|_{L^{p,q}(X)}
\lesssim\left\|S_\mathcal{Q}(f)\right\|_{L^{p,q}(X)}.
\end{align*}
From \eqref{817}, Propositions~\ref{inc-norm}(ii)
and~\ref{r;p,q=pr,qr;r},
$r<\min\{p,\,2\}$,
and Lemma~\ref{HLf-s}(i) with $u$ therein replaced by $\frac{2}{r}$,
we infer that, for any
$f\in(\mathring{{\mathcal{G}}}^\varepsilon_0(\beta,\gamma))'$,
\begin{align*}
\left\|S_\mathcal{E}(f)\right\|_{L^{p,q}(X)}
&\lesssim\left\|\left\{\sum_{k=-\infty}^\infty
\left[\mathcal{M}_\rho\left(\left(m_kf
\right)^r\right)\right]^\frac{2}{r}\right\}
^\frac{1}{2}\right\|_{L^{p,q}(X)}
\\
&=\left\|\left\{\sum_{k=-\infty}^\infty
\left[\mathcal{M}_\rho\left(\left(m_kf
\right)^r\right)\right]^\frac{2}{r}\right\}
^\frac{r}{2}\right\|_{L^{\frac{p}{r},\frac{q}{r}}(X)}^\frac{1}{r}
\lesssim\left\|\left[\sum_{k=-\infty}^\infty
\left(m_kf\right)^2\right]
^\frac{r}{2}\right\|_{L^{\frac{p}{r},\frac{q}{r}}(X)}^\frac{1}{r}
\\
&=\left\|\left[\sum_{k=-\infty}^\infty
\left(m_kf\right)^2\right]
^\frac{1}{2}\right\|_{L^{p,q}(X)}
=\left\|S_\mathcal{Q}(f)\right\|_{L^{p,q}(X)},
\end{align*}
where, for any $k\in\mathbb{Z}$,
$m_kf$ is defined as in \eqref{mkf}.
This finishes the proof of Theorem~\ref{thm638}.
\end{proof}

\subsection{Atomic Characterization
of $H^{p,q}(X)$}

In this subsection, we aim to establish the atomic characterization
of $H^{p,q}(X)$ when $p\in(\frac{n}{n+\mathrm{ind\,}(X,\mathbf{q})},1]$
and $q\in(0,\infty]$ with $n$ and $\mathrm{ind\,}(X,\mathbf{q})$ as,
respectively, in \eqref{upperdoub} and \eqref{index}.
To this end,
we first introduce the other
atomic Hardy--Lorentz space
$\mathring{H}^{p,q,r}_{\mathrm{at}}(X,\rho)$, which is closely related
to the atomic Hardy--Lorentz space
$H^{p,q,r}_{\mathrm{at}}(X,\rho)$ in Definition~\ref{6.3},
but with the distribution space
$({\mathcal{G}}^\varepsilon_0(\beta,\gamma))'$ replaced by
$(\mathring{{\mathcal{G}}}^\varepsilon_0(\beta,\gamma))'$.

\begin{definition}\label{1953}
Let $(X,\mathbf{q},\mu)$ be a quasi-ultrametric space
of homogeneous type with upper dimension $n\in(0,\infty)$,
$\rho\in\mathbf{q}$ have the property that
all $\rho$-balls are $\mu$-measurable,
$p\in(0,1]$,
$q\in(0,\infty]$, $r\in(p,\infty]\cap[1,\infty]$,
and $\varepsilon,\beta,\gamma\in(0,\infty)$ be such that
$n(\frac{1}{p}-1)<\beta,\gamma<\varepsilon<\infty$.
The $\emph{atomic Hardy--Lorentz space}$
$\mathring{H}^{p,q,r}_{\mathrm{at}}(X,\rho)$
\index[symbols]{H@$\mathring{H}^{p,q,r}_{\mathrm{at}}(X,\rho)$}
is defined in the same way as $H^{p,q,r}_{\mathrm{at}}(X,\rho)$
in Definition~\ref{6.3} but with the distribution space
$({\mathcal{G}}^\varepsilon_0(\beta,\gamma))'$ replaced by
$(\mathring{{\mathcal{G}}}^\varepsilon_0(\beta,\gamma))'$.
\end{definition}

\begin{remark}\label{221}
Let the notation be as in Definition~\ref{1953}.
\begin{enumerate}
\item
By the assumption $n(\frac{1}{p}-1)<\beta,\gamma<\infty$,
\eqref{1550}, and Remark~\ref{Rmtest}(v), we find that,
for any $p\in(0,\frac{n}{n+\mathrm{ind}_H(X,\mathbf{q})}]$,
$\mathring{H}^{p,q,r}_{\mathrm{at}}(X,\rho)$ is trivial
because the space
$\mathring{\mathcal{G}}^\varepsilon_0(\beta,\gamma)$
only contains constant functions.
\item
Let $p\in(\frac{n}{n+\mathrm{ind\,}(X,\mathbf{q})},1]$.
Then $\mathring{H}^{p,q,r}_{\mathrm{at}}(X,\rho)$ is a quasi-Banach space,
which can be deduced from Theorems~\ref{6.14} and~\ref{HLatring}
and Proposition~\ref{5.1.11}.
\end{enumerate}
\end{remark}

We first establish the following equivalence
between $\mathring{H}^{p,q,r}_{\mathrm{at}}(X,\rho)$
and ${H}^{p,q,r}_{\mathrm{at}}(X,\rho)$.
The analogous result for
Hardy spaces was shown in Proposition~\ref{324}.

\begin{theorem}
\label{HLatring}
Let $(X,\mathbf{q},\mu)$ be a quasi-ultrametric space
of homogeneous type with upper dimension $n\in(0,\infty)$
and $\mathrm{diam\,}(X)=\infty$,
$\rho\in\mathbf{q}$ be such that all $\rho$-balls
are $\mu$-measurable,
$p\in(\frac{n}{n+\mathrm{ind\,}(X,\mathbf{q})},1]$,
$q\in(0,\infty]$,
$r\in(p,\infty]\cap[1,\infty]$,
and $\varepsilon,\beta,\gamma\in(0,\infty)$ be such that
$n(\frac{1}{p}-1)<\beta,\gamma<\varepsilon\preceq\mathrm{ind\,}(X,\mathbf{q})$.
Then
$\mathring{H}^{p,q,r}_{\mathrm{at}}(X,\rho)
=H^{p,q,r}_{\mathrm{at}}(X,\rho)$
with equivalent quasi-norms in the following sense:
\begin{enumerate}
\item[\textup{(i)}]
For any $f\in H^{p,q,r}_{\mathrm{at}}(X,\rho)
\subset({\mathcal{G}}^\varepsilon_0(\beta,\gamma))'$,
the restriction $\widehat{f}$ of $f$ to
$\mathring{{\mathcal{G}}}^\varepsilon_0(\beta,\gamma)$
belongs to $\mathring{H}^{p,q,r}_{\mathrm{at}}(X,\rho)$ and
\begin{align*}
\left\|\widehat{f}\,\right\|_{\mathring{H}^{p,q,r}_{\mathrm{at}}(X,\rho)}
\leq\left\|f\right\|_{H^{p,q,r}_{\mathrm{at}}(X,\rho)}.
\end{align*}

\item[\textup{(ii)}]
For any $f\in\mathring{H}^{p,q,r}_{\mathrm{at}}(X,\rho)\subset
(\mathring{{\mathcal{G}}}^\varepsilon_0(\beta,\gamma))'$,
there exists a unique
$\widetilde{f}\in H^{p,q,r}_{\mathrm{at}}(X,\rho)$ such that
the restriction of $\widetilde{f}$ to
$\mathring{{\mathcal{G}}}^\varepsilon_0(\beta,\gamma)$ is $f$
and moreover,
\begin{align*}
\left\|\widetilde{f}\,\right\|_{H^{p,q,r}_{\mathrm{at}}(X,\rho)}
\leq C\left\|f\right\|_{\mathring{H}^{p,q,r}_{\mathrm{at}}(X,\rho)},
\end{align*}
where $C$ is a positive constant independent of $f$.
\end{enumerate}
Consequently, the linear map
$\Phi:H^{p,q,r}_{\mathrm{at}}(X)\to\mathring{H}^{p,q,r}_{\mathrm{at}}(X)$,
defined by setting $\Phi(f)\in\mathring{H}^{p,q,r}_{\mathrm{at}}(X)$
to be the restriction of $f\in H^{p,q,r}_{\mathrm{at}}(X)$ to
the space $\mathring{\mathcal{G}}_0^\varepsilon(\beta,\gamma)$,
is well defined, linear, bijective,
and has the property that, for any
$f\in H^{p,q,r}_{\mathrm{at}}(X)$,
\begin{align*}
\|f\|_{H^{p,q,r}_{\mathrm{at}}(X)}\sim
\left\|\Phi(f)\right\|_{\mathring{H}^{p,q,r}_{\mathrm{at}}(X)},
\end{align*}
where the positive equivalence constants are independent of $f$.
\end{theorem}

\begin{remark}
Let the notation be as in Theorem~\ref{HLatring}.
Then, from Theorem~\ref{HLatring} and Remark~\ref{5216}(iii),
we deduce that $\mathring{H}^{p,q,r}_{\mathrm{at}}(X,\rho)$
is independent of the choice of $\rho\in\mathbf{q}$.
Thus, in what follows, we sometimes simply write
$\mathring{H}^{p,q,r}_{\mathrm{at}}(X)$ instead
of $\mathring{H}^{p,q,r}_{\mathrm{at}}(X,\rho)$.
\end{remark}

\begin{proof}[Proof of Theorem~\ref{HLatring}]
We first prove (i).
Since $\mathring{\mathcal{G}}^\varepsilon_0(\beta,\gamma)
\subset\mathcal{G}^\varepsilon_0(\beta,\gamma)$,
it follows that
$(\mathcal{G}^\varepsilon_0(\beta,\gamma))'\subset
(\mathring{\mathcal{G}}^\varepsilon_0(\beta,\gamma))'$.
By this, we conclude that, for any $f\in H^{p,q,r}_{\mathrm{at}}(X,\rho)
\subset({\mathcal{G}}^\varepsilon_0(\beta,\gamma))'$,
$\widehat{f}\in \mathring{H}^{p,q,r}_{\mathrm{at}}(X,\rho)$
and
$$\left\|\widehat{f}\,\right\|_{\mathring{H}^{p,q,r}_{\mathrm{at}}(X,\rho)}
\leq\left\|f\right\|_{H^{p,q,r}_{\mathrm{at}}(X,\rho)}.$$
This finishes the proof of (i).

Next, we show (ii).
Let $f\in\mathring{H}^{p,q,r}_{\mathrm{at}}(X,\rho)\subset
(\mathring{{\mathcal{G}}}^\varepsilon_0(\beta,\gamma))'$.
Then there exist a sequence
$\{\lambda_{i,k}\}_{i\in\mathbb{N},k\in\mathbb{Z}}$ in $\mathbb{R}_+$
and a sequence $\{a_{i,k}\}_{i\in\mathbb{N},k\in\mathbb{Z}}$
of $(\rho,p,r)$-atoms satisfying all the conditions of
Theorem~\ref{6.15} such that
$f=\sum_{i\in\mathbb{N}}\sum_{k\in\mathbb{Z}}\lambda_{i,k}a_{i,k}$
in $(\mathring{\mathcal{G}}^\varepsilon_0(\beta,\gamma))'$
and
$$
\|f\|_{\mathring{H}^{p,q,r}_{\mathrm{at}}(X,\rho)}
\sim\left[\sum_{k\in\mathbb{Z}}\left(\sum_{i\in\mathbb{N}}
\lambda_{i,k}^p\right)^\frac{q}{p}\right]^\frac{1}{q}.
$$
From this and Theorem~\ref{6.15}, we infer that
there exists $\widetilde{f}\in
(\mathcal{G}^\varepsilon_0(\beta,\gamma))'$
such that
$$
\widetilde{f}=\sum_{i\in\mathbb{N}}\sum_{k\in\mathbb{Z}}\lambda_{i,k}a_{i,k}
$$ in $({\mathcal{G}}^\varepsilon_0(\beta,\gamma))'$
and
$$
\|f\|_{H^{p,q,r}_{\mathrm{at}}(X,\rho)}\lesssim
\left[\sum_{k\in\mathbb{Z}}\left(\sum_{i\in\mathbb{N}}
\lambda_{i,k}^p\right)^\frac{q}{p}\right]^\frac{1}{q}
\sim\|f\|_{\mathring{H}^{p,q,r}_{\mathrm{at}}(X,\rho)}.
$$

Suppose that there exists another extension
$\widetilde{g}\in H^{p,q,r}_{\mathrm{at}}(X,\rho)$ of $f$
satisfies (ii).
Then $\widetilde{g}\in({\mathcal{G}}^\varepsilon_0(\beta,\gamma))'$
and $\widetilde{g}=f$ on $\mathring{{\mathcal{G}}}^\varepsilon_0(\beta,\gamma)$.
By this and the trivial observation that,
if $h\in({\mathcal{G}}^\varepsilon_0(\beta,\gamma))'$
and $\langle h,\phi\rangle=0$ for any $\phi\in
\mathring{{\mathcal{G}}}^\varepsilon_0(\beta,\gamma)$,
then $h$ is a constant function,
we find that $\widetilde{g}-\widetilde{f}$
is a constant function.
Note that there exists no non-zero constant
function belonging to $H^{p,q}_*(X)$.
Thus, $\widetilde{g}-\widetilde{f}=0$.
This finishes the proof of the uniqueness of $\widetilde{f}$
and hence (ii),
which completes the proof of Theorem~\ref{HLatring}.
\end{proof}

Combining Theorems~\ref{6.14} and~\ref{HLatring},
we immediately obtain the following corollary.

\begin{corollary}
\label{c42}
Let the notation be as in
Theorem~\ref{HLatring} and let $s\in(1,\infty]$.
Then, as subspaces of $(\mathring{{\mathcal{G}}}^\varepsilon_0(\beta,\gamma))'$,
$\mathring{H}^{p,q,r}_{\mathrm{at}}(X)
=\mathring{H}^{{p,q,s}}_{\mathrm{at}}(X)$
with equivalent quasi-norms.
\end{corollary}

The following theorem is the main result of this subsection,
which establishes the atomic characterization of $H^{p,q}(X)$.
The analogue of this result for $H^{p}(X)=H^{p,p}(X)$ was
proven in Theorem~\ref{Thm4.11}.

\begin{theorem}
\label{atHL}
Let $(X,\mathbf{q},\mu)$ be an ultra-RD-space with
upper dimension $n\in(0,\infty)$,
where $\mu$ is assumed to be a Borel-semiregular
measure on $X$ and $\mathrm{diam\,}(X)=\infty$.
Let $p\in(\frac{n}{n+\mathrm{ind\,}(X,\mathbf{q})},1]$,
$q\in(0,\infty]$,
$r\in(p,\infty]\cap[1,\infty]$,
and $\varepsilon,\beta,\gamma\in(0,\infty)$ be such that
$n(\frac{1}{p}-1)<\beta,\gamma<\varepsilon\preceq\mathrm{ind\,}(X,\mathbf{q})$.
Then, as subspaces of $(\mathring{{\mathcal{G}}}^\varepsilon_0(\beta,\gamma))'$,
$$
H^{p,q}(X)
=\mathring{H}^{p,q,r}_{\mathrm{at}}(X)
$$
with equivalent quasi-norms.
\end{theorem}

To show Theorem~\ref{atHL},
we employ the following
essential estimates for the Lusin area function and the
Littlewood--Paley $g$-function. Since these estimates
can be obtained via an argument
similar to that used in the proof of \eqref{dw-23},
we omit the details.

\begin{lemma}
\label{ga+Sa}
Let $(X,\mathbf{q},\mu)$ be a space of homogeneous type
with upper dimension $n\in(0,\infty)$,
$\rho\in\mathbf{q}$ have the property that
all $\rho$-balls are $\mu$-measurable,
and $p\in(\frac{n}{n+\mathrm{ind\,}(X,\mathbf{q})},1]$.
Then there exists a positive constant $C$
such that, for any $(p,\infty)$-atom $a$ supported
in $B:=B_\rho(x_B,r_B)$ for some $x_B\in X$
and $r_B\in(0,\infty)$ and for any
$x\in(4C_\rho^2B)^\complement$,
\begin{align*}
g(a)(x)+S(a)(x)\leq C\left[\mu(B)\right]^{1-\frac{1}{p}}
\left[\frac{r_B}{\rho(x_B,x)}\right]^\varepsilon
\frac{1}{V_\rho(x_B,x)},
\end{align*}
where $\varepsilon\in(0,\infty)$ satisfies
$n(\frac{1}{p}-1)<\varepsilon\preceq\mathrm{ind\,}(X,\mathbf{q})$.
\end{lemma}

\begin{proof}[Proof of Theorem~\ref{atHL}]
Let $\rho\in\mathbf{q}$ be a regularized quasi-ultrametric.
We first prove that, for any
$f\in\mathring{H}^{{p,q,\infty}}_{\mathrm{at}}(X)\subset
(\mathring{{\mathcal{G}}}^\varepsilon_0(\beta,\gamma))'$,
\begin{align}
\label{819--}
f\in H^{p,q}(X)
\ \ \text{and}\ \
\|f\|_{H^{p,q}(X)}\lesssim
\|f\|_{\mathring{H}^{{p,q,\infty}}_{\mathrm{at}}(X)},
\end{align}
where the implicit positive constant is independent of $f$.
To this end,  let $f\in\mathring{H}^{{p,q,\infty}}_{\mathrm{at}}(X)$.
We first consider the case that $q\in(0,\infty)$.
From the definition
of $\mathring{H}^{{p,q,\infty}}_{\mathrm{at}}(X)$,
we deduce that there exist a sequence of
$\{\lambda_{i,k}\}_{i\in\mathbb{N},
k\in\mathbb{Z}}$ in $\mathbb{R}_+$
and a sequence
$\{a_{i,k}\}_{i\in\mathbb{N},k\in\mathbb{Z}}$
of $(p,\infty)$-atoms supported,
respectively, in $\rho$-balls
$\{B_{i,k}\}_{i\in\mathbb{N},k\in\mathbb{Z}}$ such that
$f=\sum_{k\in\mathbb{Z}}\sum_{i\in\mathbb{N}}
\lambda_{i,k}a_{i,k}$
in $(\mathring{{\mathcal{G}}}^\varepsilon_0(\beta,\gamma))'$
and
\begin{align}
\label{819---}
\left[\sum_{k\in\mathbb{Z}}\left(\sum_{i\in\mathbb{N}}
\lambda_{i,k}^p\right)
^\frac{q}{p}\right]^\frac{1}{q}\sim
\|f\|_{\mathring{H}^{p,q,r}_{\mathrm{at}}(X)}.
\end{align}
Moreover, there exist constants $C_0\in\mathbb{N}$,
$j_0\in\mathbb{Z}\setminus\mathbb{N}$, and
$C_1\in (0,\infty)$ (independent of $f$) such that,
for any $x\in X$ and $k\in\mathbb{Z}$,
\begin{align*}
\sum_{i\in\mathbb{N}}\mathbf{1}_{(2C_\rho)^{j_0}B_{i,k}}(x)\leq C_0
\end{align*}
and, for any $k\in\mathbb{Z}$ and $i\in\mathbb{N}$,
$\lambda_{i,k}\sim2^k[\mu(B_{i,k})]^\frac{1}{p}$.
Fix an arbitrary $k_0\in\mathbb{Z}$. We write
\begin{equation*}
f=\sum_{k=-\infty}^{k_0-1}\sum_{i\in\mathbb{N}}\lambda_{i,k}a_{i,k}
+\sum_{k=k_0}^{\infty}
\sum_{i\in\mathbb{N}}\lambda_{i,k}a_{i,k}=:f_1+f_2,
\end{equation*}
where the series converge in $({\mathcal{G}}^\varepsilon_0(\beta,\gamma))'$
by Theorem~\ref{6.15}.
Let $E_{k_0}:=\bigcup_{k=k_0}^{\infty}\bigcup_{i\in\mathbb{N}}
2C_\rho B_{i,k}$,
\begin{align*}
\mathrm{J_1}:=\left[\mu\left(\left\{x\in X:
\sum_{k=-\infty}^{k_0-1}\sum_{i\in\mathbb{N}}
\lambda_{i,k}S(a_{i,k})\mathbf{1}_{4C_\rho^2B_{i,k}}
(x)>2^{k_0-2}\right\}\right)\right]^\frac{1}{p},
\end{align*}
\begin{align*}
\mathrm{J_2}:=\left[\mu\left(\left\{x\in X:
\sum_{k=-\infty}^{k_0-1}\sum_{i\in\mathbb{N}}
\lambda_{i,k}S(a_{i,k})\mathbf{1}_{(4C_\rho^2B_{i,k})
^\complement}(x)>2^{k_0-2}\right\}\right)\right]^\frac{1}{p},
\end{align*}
\begin{align*}
\mathrm{J_3}:=\left[\mu\left(\left\{x\in E_{k_0}:
S(f_2)(x)>2^{k_0-2}\right\}\right)\right]^\frac{1}{p},
\end{align*}
and
\begin{align*}
\mathrm{J_4}:=\left[\mu\left(\left\{x\in (E_{k_0})^\complement:
S(f_2)(x)>2^{k_0-2}\right\}\right)\right]^\frac{1}{p}.
\end{align*}
Then
\begin{align}
\label{23hg}
\left[d_{S(f)}(2^{k_0})\right]^\frac{1}{p}
\lesssim\mathrm{J_1}+\mathrm{J_2}+\mathrm{J_3}+\mathrm{J_4}.
\end{align}
By an argument similar to that used in the proof of
\eqref{4ee5} with both $(a_{i,k})^*_\rho$ and
$\mathcal{M}_\rho(a_{i,k})$ replaced by $S(a_{i,k})$
for any $i\in\mathbb{N}$ and $k\in\mathbb{Z}$
and by \eqref{2001}, we conclude that
\begin{align}
\label{819yy1}
\left[\sum_{k_0\in\mathbb{Z}}2^{k_0q}(\mathrm{J}_{1})^q\right]^\frac{1}{q}
\lesssim\left[\sum_{k\in\mathbb{Z}}
\left(\sum_{i\in\mathbb{N}}\lambda_{i,k}^p\right)
^{\frac{q}{p}}\right]
^{\frac{1}{q}}.
\end{align}
From Lemma~\ref{ga+Sa}, \eqref{4ee6-1}, \eqref{4ee6-2},
and an argument similar to that used in the proof of \eqref{4ee6},
it follows that, for any
$i\in\mathbb{N}$, $k\in\mathbb{Z}$, and $x\in (4C_\rho^2B_{i,k})^{\complement}$,
\begin{align}
\label{S*}
S(a_{i,k})(x)
\lesssim\left[\mu(B_{i,k})\right]^{-\frac{1}{p}}
\left[\mathcal{M}_\rho(\mathbf{1}_{B_{i,k}})(x)\right]^{1+\frac{\beta}{n}},
\end{align}
which, combined with an argument similar to that used in the estimation
of \eqref{2007},
further implies that
\begin{align}
\label{819yy2}
\left[\sum_{k_0\in\mathbb{Z}}2^{k_0q}(\mathrm{J}_{2})^q\right]^\frac{1}{q}
\lesssim\left[\sum_{k\in\mathbb{Z}}
\left(\sum_{i\in\mathbb{N}}\lambda_{i,k}^p\right)
^{\frac{q}{p}}\right]
^{\frac{1}{q}}.
\end{align}
By an argument similar to that used in the estimation
of \eqref{4ee8},
we find that
\begin{align}
\label{819yy3}
\left[\sum_{k_0\in\mathbb{Z}}2^{k_0q}(\mathrm{J}_{3})^q\right]^\frac{1}{q}
\lesssim\left[\sum_{k\in\mathbb{Z}}
\left(\sum_{i\in\mathbb{N}}\lambda_{i,k}^p\right)
^{\frac{q}{p}}\right]
^{\frac{1}{q}}.
\end{align}
An argument similar to that used in the estimation
of \eqref{4ee9}
with \eqref{4ee6} replaced by \eqref{S*} implies that
\begin{align*}
\left[\sum_{k_0\in\mathbb{Z}}2^{k_0q}(\mathrm{J}_{4})^q\right]^\frac{1}{q}
\lesssim\left[\sum_{k\in\mathbb{Z}}
\left(\sum_{i\in\mathbb{N}}\lambda_{i,k}^p\right)
^{\frac{q}{p}}\right]
^{\frac{1}{q}}.
\end{align*}
From this, Proposition~\ref{Lpqnorm=c,d}, \eqref{23hg}, \eqref{819yy1},
\eqref{819yy2}, \eqref{819yy3},
and \eqref{819---}, we infer that
\begin{align*}
\|f\|_{H^{p,q}(X)}
&=\left\|S(f)\right\|_{L^{p,q}(X)}
\sim\left\{\sum_{k_0\in\mathbb{Z}}
2^{k_0q}\left[d_{S(f)}(2^{k_0})\right]^\frac{q}{p}\right\}^\frac{1}{q}
\\
&\lesssim
\left[\sum_{k_0\in\mathbb{Z}}2^{k_0q}(\mathrm{J}_{1})^q\right]^\frac{1}{q}
+\left[\sum_{k_0\in\mathbb{Z}}2^{k_0q}(\mathrm{J}_{2})^q\right]^\frac{1}{q}
\\
&\quad+\left[\sum_{k_0\in\mathbb{Z}}2^{k_0q}(\mathrm{J}_{3})^q\right]^\frac{1}{q}
+\left[\sum_{k_0\in\mathbb{Z}}2^{k_0q}(\mathrm{J}_{4})^q\right]^\frac{1}{q}
\\
&\lesssim\left[\sum_{k\in\mathbb{Z}}
\left(\sum_{i\in\mathbb{N}}\lambda_{i,k}^p\right)
^{\frac{q}{p}}\right]
^{\frac{1}{q}}
\sim\|f\|_{\mathring{H}^{{p,q,2}}_{\mathrm{at}}(X)}.
\end{align*}
This finishes the proof of \eqref{819--}
when $q\in(0,\infty)$. When $q=\infty$, the proof is similar to
the above arguments with some slight modifications;
we omit the details. Having established \eqref{819--},
we immediately obtain the inclusion
$\mathring{H}^{{p,q,\infty}}_{\mathrm{at}}(X)\subset
H^{p,q}(X)$ with appropriate control of the (quasi-)norms.

Now, let $f\in H^{p,q}(X)\subset
(\mathring{{\mathcal{G}}}^\varepsilon_0(\beta,\gamma))'$
and we show that
\begin{align}\label{9}
f\in\mathring{H}^{{p,q,r}}_{\mathrm{at}}(X)
\ \ \text{and}\ \
\|f\|_{\mathring{H}^{{p,q,r}}_{\mathrm{at}}(X)}
\lesssim\|f\|_{H^{p,q}(X)},
\end{align}
where $r\in(p,\infty)\cap[1,\infty)$ and where
the implicit positive constant is independent of
$f$. Let $\{\mathcal{S}_{2^{-k}}\}_{k\in\mathbb{Z}}$
be a homogeneous approximation
of the identity of order $\varepsilon$ as in Definition~\ref{DefHATI}.
For any $k\in\mathbb{Z}$,
define $\mathcal{D}_k:=\mathcal{S}_{2^{-k}}-\mathcal{S}_{2^{-k+1}}$
and let $S_\mathcal{D}$ and $S_{\overline{\mathcal{D}}}$ be
the Lusin area functions (with aperture 1) as in \eqref{Def-S-f}
associated, respectively, with
$\{\mathcal{D}_k\}_{k\in\mathbb{Z}}$ and
$\{\overline{\mathcal{D}}_k\}_{k\in\mathbb{Z}}$,
where $\{\overline{\mathcal{D}}_k\}_{k\in\mathbb{Z}}$
is the same as in Theorem~\ref{A3.19} related to
$\{\mathcal{D}_k\}_{k\in\mathbb{Z}}$.
We first prove that
\begin{align}\label{DbarD}
\left\|S_{\overline{\mathcal{D}}}f\right\|_{L^{p,q}(X)}\lesssim
\left\|S_\mathcal{D}f\right\|_{L^{p,q}(X)}.
\end{align}
Following the proof of Theorem~\ref{thm638},
to show \eqref{DbarD}, it suffices to prove that,
for any given $s\in(\frac{n}{n+\varepsilon},1]\cap(0,p)$
and for any $x\in X$,
\begin{align}\label{SDbarMmf}
\left[S_{\overline{\mathcal{D}}}(f)(x)\right]^2
\lesssim\sum_{k=-\infty}^{\infty}
\left[\mathcal{M}_\rho
\left((m_kf)^s\right)(x)\right]^\frac{2}{s},
\end{align}
where, for any $k\in\mathbb{Z}$,
$m_kf(x)$ is defined as in \eqref{mkfx}.
To show \eqref{SDbarMmf},
suppose that $l\in\mathbb{Z}$,
$x\in X$, and $y\in B_\rho(x,2^{-l})$.
By Theorem~\ref{A3.47},
we conclude that, for $\mu$-almost every $y\in X$,
\begin{align*}
\overline{\mathcal{D}}_lf(y)
&=\left\langle f,\overline{D}_l(y,\cdot)\right\rangle\\
&=\left\langle\sum_{k=-\infty}^{\infty}\sum_{\tau\in I_{k}}
\sum_{v=1}^{N(k,\tau)}
\widetilde{D}_k(\cdot,y_\tau^{k,v})
\int_{Q_\tau^{k,v}}\mathcal{D}_kf(u)\,d\mu(u),
\overline{D}_l(y,\cdot)\right\rangle\\
&=\sum_{k=-\infty}^{\infty}\sum_{\tau\in I_{k}}
\sum_{v=1}^{N(k,\tau)}\left\langle
\widetilde{D}_k(\cdot,y_\tau^{k,v})
,\overline{D}_l(y,\cdot)\right\rangle\int_{Q_\tau^{k,v}}
\mathcal{D}_kf(u)\,d\mu(u)\\
&=\sum_{k=-\infty}^{\infty}\sum_{\tau\in I_{k}}
\sum_{v=1}^{N(k,\tau)}
\overline{D}_l\widetilde{D}_k(y,y_\tau^{k,v})
\int_{Q_\tau^{k,v}}\mathcal{D}_kf(u)\,d\mu(u),	
\end{align*}
where $\widetilde{D}_k$ is the same as in
Theorem~\ref{A3.47} and $y_\tau^{k,v}$ and $Q_\tau^{k,v}$
are the same as in Remark~\ref{jcubes}.
Note that the kernels, $\overline{D}_k$
and $\widetilde{D}_l$, satisfy all the conditions of
Lemma~\ref{1122} and hence \eqref{DkDl}
holds for $\overline{D}_k$
and $\widetilde{D}_l$
with $k,l\in\mathbb{Z}$.
From this and an argument similar to
that used in the estimation of \eqref{817},
it follows that \eqref{SDbarMmf} holds,
which completes the proof of \eqref{DbarD}.

Next, we fix some notation. Let
$\mathcal{D}:=\{Q^k_\alpha:k\in\mathbb{Z},\ \alpha\in I_k\}$
be the set of all dyadic cubes (at all levels) as in Lemma~\ref{DyadicSys},
and recall that, by Remark~\ref{remarkcube}, we
may assume that $\delta=\frac{1}{2}$ in Lemma~\ref{DyadicSys}.
For any $k\in\mathbb{Z}\cup\{-\infty\}$,
let $\Omega_k$ and $\mathcal{D}_k$ be
the same as in the proof of Theorem~\ref{L4.10}.
For any $k,j\in\mathbb{Z}$,
let $I_{j,k}$, $Q_{\tau,k}^j$, and $\tau\in I_{j,k}$,
be the same as in the proof of Theorem~\ref{L4.10}
with the center of $Q_{\tau,k}^j$ denoted by $z_{\tau,k}^j$,
and define $B_{\tau,k}^j:=B_\rho(z_{\tau,k}^j,2^{-j+k_0})$,
where $k_0\in\mathbb{Z}_+$ satisfies $C_\rho C_r\leq 2^{k_0}$
with both $C_\rho$ as in \eqref{Crho} and
$C_r$ as in Lemma~\ref{DyadicSys}.
From \eqref{BsubsetQ} and \eqref{BsimQ},
we deduce that,
for any $k,j\in\mathbb{Z}$ and $\tau\in I_{j,k}$,
\begin{align}\label{btaukj}
C_l2^{-k_0} B_{\tau,k}^j\subset Q_{\tau,k}^j
\subset B_{\tau,k}^j
\ \ \text{and}\ \
\mu\left(B_{\tau,k}^j\right)
\sim\mu\left(Q_{\tau,k}^j\right),
\end{align}
where $C_l$ is the same as in Lemma~\ref{DyadicSys}(v).
Let $\mathcal{D}_Q:=\mathcal{D}_l$
and $\overline{\mathcal{D}}_Q:=\overline{\mathcal{D}}_l$
whenever $Q=Q^{l+{k_0}}_\tau$
for some $l\in\mathbb{Z}$
and $\tau\in I_{l+{k_0}}$ with
$Q^{l+{k_0}}_\tau$ and $I_{l+{k_0}}$ as in Lemma~\ref{DyadicSys}.
Then, by \eqref{3.12} and an argument similar to
that used in the estimation of \eqref{lammdab},
we find that
\begin{align}
\label{HLlammdab}
f(\cdot)
&=\sum_{k=-\infty}^{\infty}\sum_{j=-\infty}^{\infty}
\sum_{\tau\in I_{j,k}}
\sum_{Q\in\mathcal{D}_k,Q\subset Q^j_{\tau,k}}
\int_QD_Q(\cdot,y)\overline{\mathcal{D}}_Qf(y)\,d\mu(y)
\nonumber\\
&=\sum_{k=-\infty}^{\infty}
\sum_{j=-\infty}^{\infty}
\sum_{\tau\in I_{j,k}}
\lambda^j_{\tau,k}a^j_{\tau,k}(\cdot)
\end{align}
in $(\mathring{{\mathcal{G}}}_0^\varepsilon(\beta,\gamma))'$,
where, for any $k,j\in\mathbb{Z}$ and $\tau\in I_{j,k}$,
$$
\lambda_{\tau,k}^j:=2^k
\left[\mu(B_{\tau,k}^j)\right]^\frac{1}{p}
$$
and
\begin{align*}
a^j_{\tau,k}(\cdot)
:=\frac{1}{\lambda^j_{\tau,k}}
\sum_{Q\in\mathcal{D}_k,Q\subset Q^j_{\tau,k}}
\int_QD_Q(\cdot,y)\overline{\mathcal{D}}_Qf(y)\,d\mu(y).
\end{align*}
From this and an argument similar
to that used in the proof of \eqref{lammdaa},
we infer that
\begin{align}
&\label{Hllambda}
\sum_{k\in\mathbb{Z}}
\left(\sum_{j\in\mathbb{Z}}
\sum_{\tau\in I_{j,k}}
\left|\lambda_{\tau,k}^j\right|^p\right)
^\frac{q}{p}\nonumber\\
&\quad\sim\sum_{k\in\mathbb{Z}}
\left[2^{kp}\sum_{j\in\mathbb{Z}}
\sum_{\tau\in I_{j,k}}
\mu\left(Q_{\tau,k}^j\right)\right]
^\frac{q}{p}
\quad\text{by \eqref{btaukj}}
\nonumber\\
&\quad\lesssim\sum_{k\in\mathbb{Z}}
\left[2^{kp}\sum_{j\in\mathbb{Z}}
\sum_{\tau\in I_{j,k}}
\mu\left(Q_{\tau,k}^j\cap\Omega_k\right)\right]
^\frac{q}{p}
\lesssim\sum_{k\in\mathbb{Z}}2^{kp}
\left[\mu\left(\Omega_k\right)\right]^\frac{q}{p}
\nonumber\\
&\quad\sim\left\|S_{\overline{\mathcal{D}}}(f)\right\|^q_{L^{p,q}(X)}
\quad\text{by Lemma~\ref{Lpqnorm=c,d}}
\nonumber\\
&\quad\lesssim\left\|S_{\mathcal{D}}(f)\right\|^q_{L^{p,q}(X)}
=\|f\|_{H^{p,q}(X)}^q
\quad\text{by \eqref{DbarD}}
\end{align}
with the usual modification when $q=\infty$.
By an argument similar to that used in the
proofs of \eqref{twh-349}, \eqref{2026},
and \eqref{atkjc}, we conclude that,
for any $k,j\in\mathbb{Z}$ and $\tau\in I_{j,k}$,
$a_{\tau,k}^j$ is a harmless constant multiple
of a $(p,r)$-atom supported in $C_\rho CB_{\tau,k}^j$,
where $C\in(0,\infty)$ is the same as
in the proof of Theorem~\ref{L4.10}.
Moreover, note that, from \eqref{btaukj}
and the definition of the sequence of
maximal cubes
$\{Q_{\tau,k}^j\}_{j\in\mathbb{Z},\tau\in I_{j,k}}$
of $\mathcal{D}_k$ [which satisfy Lemma~\ref{DyadicSys}(i)], it follows that,
for any given $j,k\in\mathbb{Z}$ and $\tau_1,\tau_2\in I_{j,k}$,
\begin{align*}
\left(C_l2^{-k_0}B_{\tau_1,k}^j\cap
C_l2^{-k_0}B_{\tau_2,k}^j\right)\subset
\left(Q_{\tau_1,k}^j\cap Q_{\tau_2,k}^j\right)=\varnothing,
\end{align*}
which further implies that, for any $k\in\mathbb{Z}$
and $x\in X$,
\begin{align*}
\sum_{j\in\mathbb{Z}}
\sum_{\tau\in I_{j,k}}
\mathbf{1}_{C_l2^{-k_0}B_{\tau,k}^j}(x)
\leq1.
\end{align*}
Thus, \eqref{830} holds with any
$j_0\in(-\infty,\log_{2C_\rho}C_l2^{-k_0}]\cap\mathbb{Z}$ and $C_0=1$,
which, together with the proven conclusion that,
for any $k,j\in\mathbb{Z}$ and
$\tau\in I_{j,k}$,
$a_{\tau,k}^j$ is a harmless constant multiple
of a $(p,r)$-atom supported in the ball $C_\rho CB_{\tau,k}^j$,
\eqref{HLlammdab}, and \eqref{Hllambda},
further implies that \eqref{9} holds,
which, combined with both \eqref{819--}
and Corollary~\ref{c42},
then completes the proof of Theorem~\ref{atHL}.
\end{proof}

\begin{remark}
Zhou et al. \cite[Theorem 7.10]{zhy}
established the atomic characterization of the Lusin-area-function-based
Hardy--Lorentz space $H^{p,q}(X)$
by decomposing the elements of $H^{p,q}(X)$ into the
linear combination of molecules and then using a
wavelet reproducing formula to first obtain its molecule
characterization. In contrast, in the proof of Theorem~\ref{atHL},
we are able to directly decompose the elements of $H^{p,q}(X)$ into
linear combination of atoms since our approximation of the identity
has bounded support (see Definition~\ref{DefATI}(i))
while the approximation of the identity in \cite[Definition 2.8]{zhy}
only has exponential decay.
\end{remark}

\subsection{Littlewood--Paley Function Characterizations}

In this subsection, we establish the
Littlewood--Paley $g$-function
and the Littlewood--Paley $g_\lambda^*$-function
characterizations of
$H^{p,q}(X)$, where $H^{p,q}(X)$ is the
Lusin-area-function-based Hardy--Lorentz
space in Definition~\ref{HLviaS}. The reader is reminded
that $g$ and $g_\lambda^*$ were defined,
respectively, in
\eqref{Def-g(f)(x)} and \eqref{Def-glamda*}.

For Littlewood--Paley $g$-function,
we have the following conclusion,
which generalizes Theorem~\ref{gf} from Hardy spaces
to Hardy--Lorentz spaces.

\begin{theorem}\label{HLg}
Let $(X,\mathbf{q},\mu)$ be an ultra-RD-space with
upper dimension $n\in(0,\infty)$,
where $\mu$ is assumed to be a Borel-semiregular
measure on $X$ and $\mathrm{diam\,}(X)=\infty$.
Let
$p\in(\frac{n}{n+\mathrm{ind\,}(X,\mathbf{q})},1]$,
$q\in(0,\infty]$,
and $\varepsilon,\beta,\gamma\in(0,\infty)$ be such that
$n(\frac{1}{p}-1)<\beta,\gamma<\varepsilon\preceq\mathrm{ind\,}(X,\mathbf{q})$.
Then $f\in H^{p,q}(X)
\subset(\mathring{{\mathcal{G}}}^\varepsilon_0(\beta,\gamma))'$
if and only if
$f\in(\mathring{{\mathcal{G}}}^\varepsilon_0(\beta,\gamma))'$
and $g(f)\in L^{p,q}(X)$;
moreover, for such $f$,
\begin{align*}
\|f\|_{H^{p,q}(X)}
\sim\left\|g(f)\right\|_{L^{p,q}(X)}
\end{align*}
with the positive equivalence constants independent of $f$.
\end{theorem}

\begin{proof}
Let $\rho\in\mathbf{q}$ be a regularized quasi-ultrametric
as in Definition~\ref{D2.3}.
By Theorems~\ref{thm638} and~\ref{atHL}
with $r=\infty$, we find that, for any
$f\in(\mathring{{\mathcal{G}}}^\varepsilon_0(\beta,\gamma))'$,
\begin{align*}
\|f\|_{H^{p,q}(X)}\sim\|S(f)\|_{L^{p,q}(X)}
\sim\|f\|_{\mathring{H}^{{p,q,\infty}}_{\mathrm{at}}(X)}.
\end{align*}
From this, we deduce that,
to show the present theorem, it suffices to
prove that, for any
$f\in(\mathring{{\mathcal{G}}}^\varepsilon_0(\beta,\gamma))'$,
\begin{align}
\label{gflessf}
\|g(f)\|_{L^{p,q}(X)}\lesssim
\|f\|_{\mathring{H}^{{p,q,\infty}}_{\mathrm{at}}(X)}
\end{align}
and
\begin{align}
\label{sflesgf}
\|S(f)\|_{L^{p,q}(X)}\lesssim
\|g(f)\|_{L^{p,q}(X)}.
\end{align}
Let $f\in(\mathring{{\mathcal{G}}}^\varepsilon_0(\beta,\gamma))'$.
By Lemma~\ref{ga+Sa} and an argument
similar to that used in the proof
of Theorem~\ref{atHL} with
\begin{align*}
S(a)\lesssim\left[\mu(B)\right]^{1-\frac{1}{p}}
\left[\frac{r_B}{\rho(x_B,x)}\right]^\varepsilon
\frac{1}{V_\rho(x_B,x)}
\end{align*}
therein replaced by
\begin{align*}
g(a)\lesssim\left[\mu(B)\right]^{1-\frac{1}{p}}
\left[\frac{r_B}{\rho(x_B,x)}\right]^\varepsilon
\frac{1}{V_\rho(x_B,x)},
\end{align*}
we obtain \eqref{gflessf}.

To show \eqref{sflesgf},
let $\varepsilon'\in(n(\frac{1}{p}-1),\varepsilon)$.
From
Lemma~\ref{HLf-s}(i) with $u$, $p$, and $q$ replaced,
respectively,
by $\frac{2}{r}$, $\frac{p}{r}$, and $\frac{q}{r}$,
we infer that, for any given
$r\in(\frac{n}{n+\varepsilon'},p)$
and for any
$f\in(\mathring{{\mathcal{G}}}^\varepsilon_0(\beta,\gamma))'$,
\begin{align*}
\left\|S(f)\right\|_{L^{p,q}(X)}
&=\left\|\left[S(f)\right]^r\right\|_{L^{\frac{p}{r},
\frac{q}{r}}(X)}^\frac{1}{r}
\quad\text{by Proposition~\ref{r;p,q=pr,qr;r}}
\\
&\lesssim\left\|\left\{\sum_{k\in\mathbb{Z}}
\left[\mathcal{M}_\rho\left(\sum_{\tau\in I_k}
\sum_{v=1}^{N(k,\tau)}
\inf_{z\in Q_{\tau}^{k,v}}\left|\mathcal{D}_{k}f(z)\right|^r
\mathbf{1}_{Q_\tau^{k,v}}
\right)\right]^\frac{2}{r}\right\}^\frac{r}{2}\right\|_{
L^{\frac{p}{r},\frac{q}{r}}(X)}^\frac{1}{r}\\
&\qquad\ \text{by \eqref{819Sf} and Proposition~\ref{inc-norm}(ii)}
\\
&\lesssim\left\|\left\{\sum_{k\in\mathbb{Z}}
\left[\sum_{\tau\in I_k}
\sum_{v=1}^{N(k,\tau)}
\inf_{z\in Q_{\tau}^{k,v}}\left|\mathcal{D}_{k}f(z)\right|^r
\mathbf{1}_{Q_\tau^{k,v}}
\right]^\frac{2}{r}\right\}^\frac{r}{2}\right\|_{L^{\frac{p}{r},
\frac{q}{r}}(X)}^\frac{1}{r}
\quad\text{by Lemma~\ref{HLf-s}(i)}
\\
&=\left\|\left[\sum_{k\in\mathbb{Z}}
\sum_{\tau\in I_k}
\sum_{v=1}^{N(k,\tau)}
\inf_{z\in Q_{\tau}^{k,v}}\left|\mathcal{D}_{k}f(z)\right|^2
\mathbf{1}_{Q_\tau^{k,v}}
\right]^\frac{1}{2}\right\|_{L^{p,q}(X)}
\quad\text{by Lemma~\ref{DyadicSys}(i)},
\end{align*}
which, together with Lemma~\ref{DyadicSys}(i) and Remark~\ref{jcubes},
further implies that
\begin{align*}
\left\|S(f)\right\|_{L^{p,q}(X)}
&\lesssim
\left\|\left[\sum_{k\in\mathbb{Z}}
\sum_{\tau\in I_k}
\sum_{v=1}^{N(k,\tau)}
\inf_{z\in Q_\tau^{k,v}}
\left|\mathcal{D}_kf(z)\right|^2
\mathbf{1}_{Q_\tau^{k,v}}
\right]^\frac{1}{2}\right\|_{L^{p,q}(X)}
\\
&\leq
\left\|\left(\sum_{k\in\mathbb{Z}}
\left|\mathcal{D}_kf\right|^2
\right)^\frac{1}{2}\right\|_{L^{p,q}(X)}
=\left\|g(f)\right\|_{L^{p,q}(X)}.
\end{align*}
This finishes the proof of \eqref{sflesgf}
and hence Theorem~\ref{HLg}.
\end{proof}

Now,
we establish the
Littlewood--Paley $g_\lambda^*$-function
characterization of $H^{p,q}(X)$,
which generalizes Theorem~\ref{Glamda*}.

\begin{theorem}
\label{HLg^*}
Let $(X,\mathbf{q},\mu)$ be an ultra-RD-space with
upper dimension $n\in(0,\infty)$,
where $\mu$ is assumed to be a Borel-semiregular
measure on $X$ and $\mathrm{diam\,}(X)=\infty$.
Let $p\in(0,1]$, $q\in(0,\infty]$,
$\lambda\in(\frac{2n}{p},\infty)$,
and $0<\beta,\gamma<\varepsilon\preceq\mathrm{ind\,}(X,\mathbf{q})$.
Then there exists a constant $C\in[1,\infty)$
such that, for any
$f\in(\mathring{{\mathcal{G}}}^\varepsilon_0(\beta,\gamma))'$
,
\begin{align}\label{2242}
C^{-1}\left\|g^*_\lambda(f)\right\|_{L^{p,q}(X)}
\leq\left\|S(f)\right\|_{H^{p,q}(X)}
\leq C\left\|g^*_\lambda(f)\right\|_{L^{p,q}(X)}.
\end{align}
Consequently, if $p\in(\frac{n}{n+\mathrm{ind\,}(X,\mathbf{q})},1]$
and $n(\frac{1}{p}-1)<\beta,\gamma<\varepsilon
\preceq\mathrm{ind\,}(X,\mathbf{q})$,
then $f\in H^{p,q}(X)
\subset(\mathring{{\mathcal{G}}}^\varepsilon_0(\beta,\gamma))'$
if and only if
$f\in(\mathring{{\mathcal{G}}}^\varepsilon_0(\beta,\gamma))'$
and $g^*_\lambda(f)\in L^{p,q}(X)$;
moreover, for such $f$,
$$
\|f\|_{H^{p,q}(X)}\sim
\|g^*_\lambda(f)\|_{L^{p,q}(X)}
$$
with the positive equivalence constants
independent of $f$.
\end{theorem}

To prove Theorem~\ref{HLg^*}, we need the following lemma
which is an analogue of Lemma~\ref{S(1)} for Lorentz spaces.
Before stating it, we remind the reader of
the following piece of notation:
with the notation as in
Theorem~\ref{HLg^*}, recall that,
for any $\theta\in[1,\infty)$,
$f\in(\mathring{{\mathcal{G}}}^\varepsilon_0(\beta,\gamma))'$,
and $x\in X$, $S_\theta^{(1)}(f)(x)$ denotes the
Littlewood--Paley auxiliary function defined in \eqref{S1thetaf}.

\begin{lemma}
\label{fujiao}
Let $(X,\mathbf{q},\mu)$ be a quasi-ultrametric space
of homogeneous type with upper dimension $n\in(0,\infty)$,
$p\in(\frac{n}{n+\mathrm{ind\,}(X,\mathbf{q})},1]$,
$q\in(0,\infty]$, $\theta\in[1,\infty)$,
and $\varepsilon,\beta,\gamma\in(0,\infty)$ be such that
$0<\beta,\gamma<\varepsilon\preceq\mathrm{ind\,}(X,\mathbf{q})$,
Then there exist positive constants $C_1,C_2\in[1,\infty)$,
independent of $\theta$,
such that, for
any $f\in(\mathring{{\mathcal{G}}}^\varepsilon_0(\beta,\gamma))'$,
\begin{align}\label{8192}
\left\|S_\theta^{(1)}(f)\right\|_{L^{p,q}(X)}
\leq C_1\theta^\frac{n}{p}
\left\|S_1^{(1)}(f)\right\|_{L^{p,q}(X)}
\end{align}
and
\begin{align}\label{8191}
C_2^{-1}\left\|S_1^{(1)}(f)\right\|_{L^{p,q}(X)}
\leq\left\|S(f)\right\|_{L^{p,q}(X)}
\leq C_2\left\|S_1^{(1)}(f)\right\|_{L^{p,q}(X)}.
\end{align}
\end{lemma}

\begin{proof}
Let $\rho\in\mathbf{q}$ be a regularized quasi-ultrametric
and $f\in(\mathring{{\mathcal{G}}}^\varepsilon_0(\beta,\gamma))'$.
By \eqref{wwm-495} and Proposition~\ref{inc-norm}(ii),
we find that \eqref{8191} holds.

Next, we show \eqref{8192}. From an argument similar to
that used in the proof of \eqref{8193}, it follows that,
for any $f\in(\mathring{{\mathcal{G}}}^\varepsilon_0(\beta,\gamma))'$
and $t\in(0,\infty)$,
\begin{align}
\label{claim1}
\mu\left(\left\{x\in X:
S^{(1)}_\theta(f)(x)>t\right\}\right)
\lesssim\theta^n\left[\mu\left(E_t\right)
+t^{-2}\int_0^ts\mu\left(E_s\right)\,ds\right],
\end{align}
where, for any $s\in(0,\infty)$,
\begin{align*}
E_s:=\left\{x\in X:S^{(1)}_1(f)(x)>s\right\},
\end{align*}
which further implies that, for any $l\in\mathbb{Z}$,
\begin{align}
\label{claim2}
\mu\left(\left\{x\in X:
S^{(1)}_\theta(f)(x)>2^l\right\}\right)
&\lesssim\theta^n\left[\mu\left(E_{2^l}\right)
+{2}^{-2l}\int_0^{2^l}s\mu\left(E_s\right)\,ds\right]
\nonumber\\
&=\theta^n\left[\mu\left(E_{2^l}\right)
+{2}^{-2l}\sum_{m=-\infty}^l\int_{2^{m-1}}^{2^m}
s\mu\left(E_s\right)\,ds\right]
\nonumber\\
&\lesssim\theta^n\left[\mu\left(E_{2^l}\right)
+{2}^{-2l}\sum_{m=-\infty}^l2^{2m}
\mu\left(E_{2^{m-1}}\right)\right].
\end{align}
Then we consider the following three cases for $q$ and $p$.

\emph{Case (1)} $q\in(0,\infty)$ and $\frac{q}{p}\in(1,\infty)$.
In this case, we conclude that
\begin{align*}
&\left\|S_\theta^{(1)}(f)\right\|_{L^{p,q}(X)}^q
\\
&\quad\sim\int_0^\infty t^{q-1}\left[\mu
\left(\left\{x\in X:S_\theta^{(1)}(f)(x)>t
\right\}\right)\right]^\frac{q}{p}\,dt
\quad\text{by Proposition~\ref{Lpqnorm=c,d}}
\nonumber\\
&\quad\lesssim\theta^\frac{nq}{p}\int_0^\infty t^{q-1}
\left[\mu\left(E_t\right)
+t^{-2}\int_0^ts\mu\left(E_s\right)\,ds\right]^\frac{q}{p}\,dt
\quad\text{by \eqref{claim1}}
\nonumber\\
&\quad\sim\theta^\frac{nq}{p}\left\{\int_0^\infty t^{q-1}
\left[\mu\left(E_t\right)\right]^\frac{q}{p}\,dt
+\int_0^\infty t^{q-1-\frac{2q}{p}}\left[
\int_0^ts\mu\left(E_s\right)\,ds\right]^\frac{q}{p}\,dt\right\},
\end{align*}
which, combined with H\"older's inequality, $p\leq1$,
and Fubini's theorem,
further implies that
\begin{align}
\label{case1}
\left\|S_\theta^{(1)}(f)\right\|_{L^{p,q}(X)}^q
&\lesssim\theta^\frac{nq}{p}\Bigg\{
\left\|S_1^{(1)}(f)\right\|_{L^{p,q}(X)}^q
+\int_0^\infty t^{q-1-\frac{2q}{p}}
\left[\int_0^ts^\frac{q(1-p)}{q-p}\,ds\right]^{\frac{q-p}{p}}
\nonumber\\
&\quad\times
\int_0^t\left[s^p\mu\left(E_s\right)\right]^\frac{q}{p}\,ds
\,dt\Bigg\}\nonumber\\
&\sim\theta^\frac{nq}{p}\left\{
\left\|S_1^{(1)}(f)\right\|_{L^{p,q}(X)}^q
+\int_0^\infty s^{q-1}\left[\mu\left(E_s\right)
\right]^\frac{q}{p}\,ds
\right\}
\nonumber\\
&\sim\theta^\frac{nq}{p}
\left\|S_1^{(1)}(f)\right\|_{L^{p,q}(X)}^q.
\end{align}

\emph{Case (2)} $q\in(0,\infty)$ and $\frac{q}{p}\in(0,1]$.
In this case,
we deduce that
\begin{align*}
\left\|S_\theta^{(1)}(f)\right\|_{L^{p,q}(X)}^q
&\sim\sum_{l\in\mathbb{Z}}2^{lq}
\left[\mu\left(\left\{x\in X:
S_\theta^{(1)}f(x)>2^l\right\}\right)\right]^\frac{q}{p}
\quad\text{by Proposition~\ref{Lpqnorm=c,d}}
\nonumber\\
&\lesssim\theta^\frac{nq}{p}\sum_{l\in\mathbb{Z}}2^{lq}
\left[\mu\left(E_{2^l}\right)
+{2}^{-2l}\sum_{m=-\infty}^l2^{2m}
\mu\left(E_{2^{m-1}}\right)\right]^\frac{q}{p}
\quad\text{by \eqref{claim2}}
\nonumber\\
&\leq\theta^\frac{nq}{p}\left\{\sum_{l\in\mathbb{Z}}
2^{lq}\left[\mu\left(E_{2^l}\right)\right]^\frac{q}{p}
+\sum_{l\in\mathbb{Z}}2^{lq(1-\frac{2}{p})}
\sum_{m=-\infty}^l2^\frac{2mq}{p}
\left[\mu\left(E_{2^{m-1}}\right)\right]^\frac{q}{p}\right\}\\
&\qquad\ \text{by Lemma~\ref{discreteinequality} with
$r$ therein replaced by $\frac{q}{p}$,}
\end{align*}
which, together with Fubini's theorem and $p<2$,
further implies that
\begin{align}
\label{case2}
\left\|S_\theta^{(1)}(f)\right\|_{L^{p,q}(X)}^q
&\lesssim\theta^\frac{nq}{p}
\left\{\left\|S_1^{(1)}(f)\right\|_{L^{p,q}(X)}^q
+\sum_{m\in\mathbb{Z}}\left[\sum_{l=m}^\infty2^{lq(1-\frac{2}{p})}
\right]2^\frac{2mq}{p}\left[\mu
\left(E_{2^{m-1}}\right)\right]^\frac{q}{p}\right\}
\nonumber\\
&\sim\theta^\frac{nq}{p}
\left\{\left\|S_1^{(1)}(f)\right\|_{L^{p,q}(X)}^q+
\sum_{m\in\mathbb{Z}}2^{mq}\left[\mu
\left(E_{2^{m-1}}\right)\right]^\frac{q}{p}\right\}
\nonumber\\
&\sim\theta^\frac{nq}{p}\left\|S_1^{(1)}(f)\right\|_{L^{p,q}(X)}^q.
\end{align}

\emph{Case (3)} $q=\infty$.
In this case, by \eqref{claim1},
we obtain
\begin{align*}
\left\|S_\theta^{(1)}(f)\right\|_{L^{p,\infty}(X)}
&=\sup_{\alpha\in(0,\infty)}
\alpha\mu\left(\left\{x\in X:
\left|S_\theta^{(1)}f(x)\right|>\alpha\right\}\right)^\frac{1}{p}
\\
&\lesssim\sup_{\alpha\in(0,\infty)}
\alpha\theta^\frac{n}{p}\left[\mu\left(E_\alpha\right)
+\alpha^{-2}\int_0^\alpha
s\mu\left(E_s\right)\,ds\right]^\frac{1}{p}
\\
&\sim\sup_{\alpha\in(0,\infty)}
\alpha\theta^\frac{n}{p}
\left\{\left[\mu\left(E_\alpha\right)\right]^\frac{1}{p}
+\left[\alpha^{-2}\int_0^\alpha
s\mu\left(E_s\right)\,ds\right]^\frac{1}{p}\right\},
\end{align*}
which, combined with $p<2$,
further implies that
\begin{align*}
\left\|S_\theta^{(1)}(f)\right\|_{L^{p,\infty}(X)}
&\lesssim\theta^\frac{n}{p}\left\|S^{(1)}_1(f)\right\|_{L^{p,\infty}(X)}+
\sup_{\alpha\in(0,\infty)}
\alpha^{1-\frac{2}{p}}\theta^\frac{n}{p}\left[\int_0^\alpha
s^{1-p}s^p\mu\left(E_s\right)\,ds\right]^\frac{1}{p}
\\
&\lesssim\theta^\frac{n}{p}\left\|S^{(1)}_1(f)\right\|_{L^{p,\infty}(X)}
\left[1+\sup_{\alpha\in(0,\infty)}
\alpha^{1-\frac{2}{p}}\left(\int_0^\alpha
s^{1-p}\,ds\right)^\frac{1}{p}\right]
\\
&\sim\theta^\frac{n}{p}\left\|S^{(1)}_1(f)\right\|_{L^{p,\infty}(X)}.
\end{align*}
Using this, \eqref{case1}, and \eqref{case2} then
completes the proof of \eqref{8192}
and hence Lemma \ref{fujiao}.
\end{proof}

Now, we prove Theorem~\ref{HLg^*}.

\begin{proof}[Proof of Theorem~\ref{HLg^*}]
Let $\rho\in\mathbf{q}$ be a regularized quasi-ultrametric.
On the one hand, by Theorem~\ref{thm638},
\eqref{819Sfg*}, and Proposition~\ref{inc-norm}(ii), we find that,
for any $f\in(\mathring{{\mathcal{G}}}^\varepsilon_0(\beta,\gamma))'$,
\begin{align*}
\|f\|_{H^{p,q}(X)}\sim
\left\|S(f)\right\|_{L^{p,q}(X)}\lesssim
\left\|g^*_\lambda(f)\right\|_{L^{p,q}(X)},
\end{align*}
where the implicit positive constant only depends on $\lambda$.

On the other hand, from the
Aoki--Rolewicz theorem
(see, for instance, \cite[Theorem 5]{a42}) applied to the quasi-norm
$\|\cdot\|_{L^{\frac{p}{2},\frac{q}{2}}(X)}$,
we infer that there exists a constant $\tau\in(0,1]$
such that,
for any $f\in(\mathring{{\mathcal{G}}}^\varepsilon_0(\beta,\gamma))'$,
\begin{align*}
\left\|g^*_\lambda(f)\right\|_{L^{p,q}(X)}^{2\tau}
&=\left\|\left[g^*_\lambda(f)\right]^2\right\|
_{L^{\frac{p}{2},\frac{q}{2}}(X)}^{\tau}
\quad\text{by Proposition~\ref{r;p,q=pr,qr;r}}\\
&\lesssim\left\|\sum_{j=0}^\infty2^{-j\lambda}
\left[S^{(1)}_{2^j}(f)\right]^2
\right\|_{L^{\frac{p}{2},\frac{q}{2}}(X)}^{\tau}
\quad\text{by \eqref{819g*sf} and Proposition~\ref{inc-norm}(ii)}
\\
&\lesssim\sum_{j=0}^\infty2^{-\tau j\lambda}
\left\|\left[S^{(1)}_{2^j}(f)\right]^2
\right\|_{L^{\frac{p}{2},\frac{q}{2}}(X)}^{\tau}
\quad\text{by the
Aoki--Rolewicz theorem}\\
&=\sum_{j=0}^\infty2^{-\tau j\lambda}
\left\|S^{(1)}_{2^j}(f)\right\|_{L^{p,q}(X)}^{2\tau}
\quad\text{by Proposition~\ref{r;p,q=pr,qr;r}}.
\end{align*}
By this, Lemma~\ref{fujiao},
and $\lambda\in(\frac{2n}{p},\infty)$,
we then conclude that
\begin{align*}
\left\|g^*_\lambda(f)\right\|_{L^{p,q}(X)}^{2\tau}
&\lesssim\sum_{j=0}^\infty
2^{-\tau j(\lambda-\frac{2n}{p})}
\left\|S^{(1)}_1(f)\right\|_{L^{p,q}(X)}^{2\tau}
\sim\left\|S(f)\right\|_{L^{p,q}(X)}^{2\tau}
\sim\|f\|_{H^{p,q}(X)}^{2\tau},
\end{align*}
which completes the proof of Theorem~\ref{HLg^*}.
\end{proof}

\begin{remark}
\begin{enumerate}
\item
Note that we use the same
homogeneous approximation of the identity to define the
Lusin area function and the Littlewood--Paley $g_\lambda^*$ function;
see \eqref{Def-S-f} and \eqref{Def-glamda*}.
From this and the proof of Theorem~\ref{HLg^*},
we deduce that,
if $(X,\mathbf{q},\mu)$ is only assumed to be
a quasi-ultrametric space of homogeneous type,
then \eqref{2242} also
holds for any given $p\in(0,\infty)$ and
$f\in(\mathring{\mathcal{G}}^\varepsilon_0(\beta,\gamma))'$
with $\varepsilon$, $\beta$, and $\gamma$ as in Theorem~\ref{HLg^*}.

\item
In the proof of Theorem~\ref{HLg^*},
the reverse doubling condition \eqref{RDcondition}
is needed since we use Theorem~\ref{thm638}
to guarantee that $H^{p,q}(X)$ is independent of the choices
of the homogeneous approximation of the identity
and hence $H^{p,q}(X)$ via different Lusin area functions
coincide with each other with equivalent quasi-norms.
In this case, we also need
$p>\frac{n}{n+\mathrm{ind\,}(X,\mathbf{q})}$
due to Theorem~\ref{thm638}.

\item
Let $X:=\mathbb{R}^n$ be an anisotropic Euclidean space,
which is also a special ultra-RD-space,
and let $A$ be an expansive matrix.
Note that $\lambda$ in \eqref{Def-glamda*} equals
$n\lambda$ with $\lambda$ as in the Littlewood--Paley
$g_\lambda^*$-function in \cite[Definition 2.6]{lyy2018}.
In 2018, Liu et al. \cite[Theorem 2.9]{lyy2018} showed that,
for any $f$ belonging to the classical space of tempered distributions on
$\mathbb{R}^n$ and for any $p\in(0,1]$, $q\in(0,\infty)$,
and $\lambda\in(\frac{2n}{p},\infty)$,
\begin{align*}
\left\|g^*_\lambda(f)\right\|_{L^{p,q}(\mathbb{R}^n)}
\lesssim\left\|f\right\|_{H^{p,q}_A(\mathbb{R}^n)},
\end{align*}
where $H^{p,q}_A(\mathbb{R}^n)$ is the anisotropic Hardy--Lorentz space
as in \cite[Definition 2.4]{lyy2018}.
In this sense, Theorem~\ref{HLg^*} coincides with
\cite[Theorem 2.9]{lyy2018}.

\item
For any
$p\in(0,\infty)$, it holds $H^{p,p}(X)=H^p(X)$,
where $H^{p,p}(X)$ and $H^p(X)$ are
the same as, respectively, in Definitions~\ref{HLviaS}
and~\ref{D3.1.4}. Therefore, in this case,
Theorem~\ref{HLg^*} coincides with Theorem~\ref{Glamda*}
and hence, by Remark~\ref{1629},
the range of $\lambda$ in Theorem~\ref{HLg^*}
is also the known best possible.
\end{enumerate}
\end{remark}

By Theorems~\ref{HL-M},~\ref{6.14},~\ref{thmHpqm},~\ref{HLatring},~\ref{atHL},
\ref{HLg}, and~\ref{HLg^*},
we obtain the following conclusion
of the mutual equivalencies of
all these Hardy--Lorentz spaces.

\begin{theorem}\label{equivelent}
Let $(X,\mathbf{q},\mu)$ be an ultra-RD-space with
upper dimension $n\in(0,\infty)$,
where $\mu$ is assumed to be a Borel-semiregular
measure on $X$.
Let $p\in(\frac{n}{n+\mathrm{ind\,}(X,\mathbf{q})},1]$,
$q\in(0,\infty]$,
$r\in(p,\infty]\cap[1,\infty]$,
$\varepsilon,\beta,\gamma\in(0,\infty)$ be such that
$n(\frac{1}{p}-1)<\beta,\gamma<\varepsilon\preceq
\mathrm{ind\,}(X,\mathbf{q})$, and
$\lambda\in(\frac{2n}{p},\infty)$.
\begin{enumerate}
\item If $\mu(X)=\infty$, then
$H^{p,q}_*(X)$,
$H^{p,q}_{\theta}(X)$ with $\theta\in(0,\infty)$,
$H^{p,q}_+(X)$,
$H^{p,q,r}_{\mathrm{at}}(X)$,
$H^{p,q,r,\eta}_{\mathrm{mol}}(X)$ with $\eta$ as in \eqref{etarange},
$H^{p,q}(X)$, and $\mathring{H}^{p,q,r}_{\mathrm{at}}(X)$
are naturally identified with one another as subspaces of
$(\mathring{\mathcal{G}}^\varepsilon_0(\beta,\gamma))'$;
moreover, for any $f\in(\mathring{\mathcal{G}}^\varepsilon_0(\beta,\gamma))'$,
\begin{align*}
\|f\|_{H^{p,q}_*(X)}&\sim\|f\|_{H^{p,q}_{\theta}(X)}
\sim\|f\|_{H^{p,q}_+(X)}\sim\|f\|_{H^{p,q,r}_{\mathrm{at}}(X)}
\sim\|f\|_{H^{p,q,r,\eta}_{\mathrm{mol}}(X)}
\\
&\sim\|f\|_{H^{p,q}(X)}\sim\left\|g(f)\right\|_{L^{p,q}(X)}\sim
\left\|g^*_\lambda(f)\right\|_{L^{p,q}(X)}
\end{align*}
with the positive equivalence constants independent of $f$.
\item If $\mu(X)<\infty$,
then $H^{p,q}_*(X)$,
$H^{p,q}_{\theta}(X)$ with $\theta\in(0,\infty)$,
$H^{p,q,r}_{\mathrm{at}}(X)$, and
$H^{p,q,r,\eta}_{\mathrm{mol}}(X)$ with $\eta$ as in \eqref{etarange}
are naturally identified with one another as subspaces of
$({\mathcal{G}}^\varepsilon_0(\beta,\gamma))'$;
moreover, for any $f\in({\mathcal{G}}^\varepsilon_0(\beta,\gamma))'$,
\begin{align*}
\|f\|_{H^{p,q}_*(X)}\sim\|f\|_{H^{p,q}_{\theta}(X)}
\sim\|f\|_{H^{p,q,r}_{\mathrm{at}}(X)}
\sim\|f\|_{H^{p,q,r,\eta}_{\mathrm{mol}}(X)}
\end{align*}
with the positive equivalence constants independent of $f$.
\end{enumerate}
\end{theorem}

\begin{remark}
From Remarks~\ref{2548} and~\ref{2214}(ii),
we infer that Theorem~\ref{equivelent} is a generalization of
the maximal function characterizations in \cite[Theorem~4.11]{AlMi2015} and
the atomic and the molecular characterizations in \cite[Theorem~6.11]{AlMi2015}
from the maximal Hardy space $H^p(X)$ on $n$-Ahlfors
regular quasi-ultrametric spaces
to the maximal Hardy--Lorentz space $H^{p,q}_*(X)$ on ultra-RD-spaces
for the optimal ranges
\begin{align*}
p\in\left(\frac{n}{n+\mathrm{ind\,}(X,\mathbf{q})},1\right]
\ \ \text{and}\ \
q\in(0,\infty].
\end{align*}
\end{remark}

\section{Real Interpolation of Hardy--Lorentz Spaces}
\label{section5.7}

In this section, we consider the interpolation
properties of Hardy--Lorentz spaces
via the real method on general quasi-ultrametric spaces of homogeneous type.
We first recall
some basic notions about the real interpolation.
A pair of quasi-normed spaces
$(X_1,X_2)$\index[symbols]{X@$(X_1,X_2)$} is called a \emph{compatible couple
of quasi-normed spaces}\index[words]{compatible couple
of quasi-normed spaces} if $X_1$ and $X_2$ are
two quasi-normed linear spaces which are continuously
embedded in some larger topological vector space.
Let
\begin{align*}
X_1+X_2:=\left\{f_1+f_2:f_1\in X_1,\ f_2\in X_2\right\}.
\index[symbols]{X@$X_1+X_2$}
\end{align*}
The \emph{Peetre $K$-functional}
\index[words]{Peetre $K$-functional}
on $(X_1,X_2)$
is defined by setting, for any $t\in(0,\infty]$
and $f\in X_1+X_2$,
$$
K(t,f;X_1,X_2):=\inf\left\{\|f_1\|_{X_1}
+t\|f_2\|_{X_2}:f=f_1+f_2,\ f_1\in X_1,\ f_2\in X_2\right\}.
$$
For any $\theta\in(0,1)$ and $q\in(0,\infty]$,
the \emph{interpolation space}\index[words]{interpolation space}
$(X_1,X_2)_{\theta,q}$\index[symbols]{X@$(X_1,X_2)_{\theta,q}$}
is defined by setting
$$
(X_1,X_2)_{\theta,q}:=\left\{f\in X_1+X_2:
\|f\|_{(X_1,X_2)_{\theta,q}}<\infty\right\},
$$
where, for any $f\in X_1+X_2$,
$$
\|f\|_{(X_1,X_2)_{\theta,q}}:=\left\{\int_0^\infty
\left[t^{-\theta}K(t,f;X_1,X_2)\right]^q
\frac{dt}{t}\right\}^\frac{1}{q}.
$$

The following interpolation lemma
is exactly \cite[Theorem 5.3.1]{bl76}.

\begin{lemma}
\label{LczLpq}
Let $(X,\mu)$ be a measure space.
Suppose that $p_1,p_2,q_1,q_2,q\in(0,\infty]$ and
$\theta\in(0,1)$. If $p_1\neq p_2$
and
$\frac{1}{p}=\frac{1-\theta}{p_1}+\frac{\theta}{p_2}$,
then
\begin{align*}
\left(L^{p_1,q_1}(X),L^{p_2,q_2}(X)\right)_{\theta,q}
=L^{p,q}(X).
\end{align*}
This identification also holds when $p=p_1=p_2$ and
$\frac{1}{q}=\frac{1-\theta}{q_1}+\frac{\theta}{q_2}$.
\end{lemma}

\begin{remark}
\label{LczLpq-x}
Specializing Lemma~\ref{LczLpq} to the case when $q_1=p_1$ and $q_2=p_2$
and using the fact that $L^{p_1,p_1}(X)=L^{p_1}(X)$
and $L^{p_2,p_2}(X)=L^{p_2}(X)$,
it follows that, for any
$p,p_1,p_2,q\in(0,\infty]$ and $\theta\in(0,1)$ satisfying $p_1\neq p_2$
and $\frac{1}{p}=\frac{1-\theta}{p_1}+\frac{\theta}{p_2}$,
\begin{align*}
\left(L^{p_1}(X),L^{p_2}(X)\right)_{\theta,q}
=L^{p,q}(X).
\end{align*}
\end{remark}

Recall that, for any $r\in\mathbb{R}$, $r_+:=\max\{0,\,r\}$.
With this in mind,
the main result of this section is as follows.

\begin{theorem}\label{chazhi}
Let $(X,\mathbf{q},\mu)$ be a quasi-ultrametric space
of homogeneous type with upper dimension $n\in(0,\infty)$.
Suppose that
$p,p_1,p_2\in(\frac{n}{n+\mathrm{ind\,}(X,\mathbf{q})},\infty)$,
$q,q_1,q_2\in(0,\infty]$, $\theta\in(0,1)$,
and $\varepsilon,\beta,\gamma\in(0,\infty)$ satisfy
$$
n\left(\frac{1}{\min\{p_1,\,p_2\}}-1\right)_+<\beta,
\gamma<\varepsilon\preceq\mathrm{ind\,}(X,\mathbf{q}).
$$
\begin{enumerate}
\item[{\rm(i)}]
If $p_1\neq p_2$ and $\frac{1}{p}
=\frac{1-\theta}{p_1}+\frac{\theta}{p_2}$,
then, as subspaces of
$({\mathcal{G}}^\varepsilon_0(\beta,\gamma))'$,
\begin{align*}
\left(H^{p_1,q_1}_\ast(X),H^{p_2,q_2}_\ast(X)\right)_{\theta,\,q}
=H^{p,q}_\ast(X)
\end{align*}
with equivalent quasi-norms.

\item[{\rm(ii)}]
Assume further that $(X,\mathbf{q},\mu)$ satisfies \eqref{RDcondition}.
If $\frac{1}{q}=
\frac{1-\theta}{q_1}+\frac{\theta}{q_2}$, then, as subspaces of
$({\mathcal{G}}^\varepsilon_0(\beta,\gamma))'$,
\begin{align*}
\left(H^{p,q_1}_\ast(X),H^{p,q_2}_\ast(X)\right)_{\theta,\,q}
=H^{p,q}_\ast(X)
\end{align*}
with equivalent quasi-norms.
\end{enumerate}
\end{theorem}

\begin{remark}
\label{chazhi-x}
Specializing Theorem~\ref{chazhi}(i) to the case when $q_1=p_1$ and $q_2=p_2$
and using the fact that $H^{p_1,p_1}_\ast(X)=H^{p_1}_\ast(X)$
and $H^{p_2,p_2}_\ast(X)=H^{p_2}_\ast(X)$,
it follows that, for any $p,p_1,p_2
\in(\frac{n}{n+\mathrm{ind\,}(X,\mathbf{q})},\infty)$,
$q\in(0,\infty]$, and $\theta\in(0,1)$ satisfying $p_1\neq p_2$
and $\frac{1}{p}=\frac{1-\theta}{p_1}+\frac{\theta}{p_2}$,
\begin{align*}
\left(H^{p_1}_\ast(X),H^{p_2}_\ast(X)\right)_{\theta,q}
=H^{p,q}_\ast(X)
\end{align*}
as subspaces of
$({\mathcal{G}}^\varepsilon_0(\beta,\gamma))'$.
In particular, for any $p,p_1,p_2\in(\frac{n}{n+
\mathrm{ind\,}(X,\mathbf{q})},\infty)$
and $\theta\in(0,1)$ satisfying $p_1\neq p_2$
and $\frac{1}{p}=\frac{1-\theta}{p_1}+\frac{\theta}{p_2}$,
\begin{align*}
\left(H^{p_1}_\ast(X),H^{p_2}_\ast(X)\right)_{\theta,p}
=H^p_\ast(X)
\end{align*}
as subspaces of
$({\mathcal{G}}^\varepsilon_0(\beta,\gamma))'$.
\end{remark}

Before proving Theorem~\ref{chazhi}, we collect
a few auxiliary results. First,
by an argument similar to that used in the proof of Theorem~\ref{6.16},
we immediately obtain the following lemma;
we omit the details of its proof here.

\begin{lemma}
\label{fenjieHpq}
Let $(X,\mathbf{q},\mu)$ be a quasi-ultrametric space
of homogeneous type with upper dimension $n\in(0,\infty)$. Suppose that
$p\in(\frac{n}{n+\mathrm{ind\,}(X,\mathbf{q})},\infty)$,
$p_0\in(\frac{n}{n+\mathrm{ind\,}(X,\mathbf{q})},p)\cap(0,1]$,
$q\in(0,\infty]$, and
$\beta, \gamma,\varepsilon\in(0,\infty)$ satisfy
$n(\frac{1}{p_0}-1)<\beta,\gamma<\varepsilon
\preceq\mathrm{ind\,}(X,\mathbf{q})$.
Then, for any $f\in H^{p,q}_\ast(X)$
and $k\in\mathbb{Z}$,
there exist $g_{k}\in L^\infty(X)$
and $b_{k}\in H^{p_0}_\ast(X)$ such that
$f=g_{k}+b_{k}$,
\begin{align*}
\left\|g_{k}\right\|_{L^\infty(X)}\lesssim2^{k},
\end{align*}
and
\begin{align*}
\left\|b_{k}\right\|_{H^{p_0}_\ast(X)}^{p_0}
\lesssim\int_{\Omega_{k}}\left[f^*_\rho(x)\right]^{p_0}\,d\mu(x)
\end{align*}
with the implicit positive constants
independent of $f$ and $k$,
where $f^*_\rho$ and $\Omega_{k}$ are the same as,
respectively, in Definition \ref{maxfunc}(iii)
and \eqref{omegak} and where $\rho\in\mathbf{q}$ is
any regularized quasi-ultrametric as in Definition~\ref{D2.3}.
\end{lemma}

\begin{proof}
Let $f\in H^{p,q}_\ast(X)$ and $k_0\in\mathbb{Z}$.
By the proof of Theorem~\ref{6.16} [which also holds for $p\in(1,\infty)$],
we find that
\begin{align}\label{2215}
f=\sum_{k\in\mathbb{Z}}\sum_{i\in\mathbb{N}}
h_{i,k}
=\sum_{k=-\infty}^{k_0}\sum_{i\in\mathbb{N}}h_{i,k}
+\sum_{k=k_0+1}^\infty\sum_{i\in\mathbb{N}}h_{i,k}
=:g_{k_0}+b_{k_0}
\end{align}
in $(\mathcal{G}^\varepsilon_0(\beta,\gamma))'$,
where, for any $k\in\mathbb{Z}$ and $i\in\mathbb{N}$,
$h_{i,k}$ is a harmless constant multiple of
a $(p,\infty)$-atom supported in a $\rho$-ball $B_{i,k}$
and satisfies
\begin{align}\label{15581}
\mathrm{supp\,}h_{i,k}\subset B_{i,k},
\end{align}
\begin{align}\label{15582}
\|h_{i,k}\|_{L^\infty(X)}\lesssim2^{k},
\end{align}
\begin{align}\label{15583}
\sum_{i\in\mathbb{N}}\mathbf{1}_{ B_{i,k}}(x)\lesssim1,
\end{align}
and
\begin{align}\label{15584}
\bigcup_{i\in\mathbb{N}}B_{i,k}=\Omega_k.
\end{align}
From \eqref{15581}, \eqref{15582}, and \eqref{15583},
we deduce that
\begin{align}\label{2216}
\left\|g_{k_0}\right\|_{L^\infty(X)}
\lesssim\sum_{k=-\infty}^{k_0}\left\|\sum_{i\in\mathbb{N}}
2^k\mathbf{1}_{B_{i,k}}\right\|_{L^\infty(X)}
\lesssim\sum_{k=-\infty}^{k_0}2^{k}=2^{k_0+1}.
\end{align}

On the other hand, the proof of Theorem~\ref{6.16}
implies that, for any $k\in\mathbb{Z}$
and $i\in\mathbb{N}$,
we can write $h_{i,k}$ as
$$
h_{i,k}=\lambda_{i,k}a_{i,k},
$$
where $0\leq\lambda_{i,k}\lesssim2^k[\mu(B_{i,k})]^\frac{1}{p}$
and $a_{i,k}$ is a $(p,\infty)$-atom.
Then, for any $k\in\mathbb{Z}$ and $i\in\mathbb{N}$,
$[\mu(B_{i,k})]^{\frac{1}{p}-\frac{1}{p_0}}a_{i,k}$
is a $(p_0,\infty)$-atom.
By this,
we conclude that
\begin{align*}
&\sum_{k=k_0+1}^\infty\sum_{i\in\mathbb{N}}
\left\{\frac{\lambda_{i,k}}{[\mu(B_{i,k})]^{
\frac{1}{p}-\frac{1}{p_0}}}\right\}^{p_0}
\\
&\quad\lesssim\sum_{k=k_0+1}^\infty
\sum_{i\in\mathbb{N}}2^{kp_0}\mu(B_{i,k})
\quad\text{since $\lambda_{i,k}\lesssim2^k[\mu(B_{i,k})]^\frac{1}{p}$}\\
&\quad\sim\sum_{k=k_0+1}^\infty2^{kp_0}\mu(\Omega_k)
\quad\text{by \eqref{15583} and \eqref{15584}}\\
&\quad\lesssim\int_{2^{k_0+1}}^\infty
\alpha^{p_0-1}\mu\left(\left\{x\in X:
f^*_\rho(x)>\alpha\right\}\right)\,d\alpha
\\
&\quad=\int_0^\infty\int_X\alpha^{p_0-1}
\mathbf{1}_{\{\alpha\in(0,\infty):\alpha\ge2^{k_0+1}\}}(\alpha)
\mathbf{1}_{\{x\in X:f^*_\rho(x)>\alpha\}}(y)\,d\mu(y)\,d\alpha,
\end{align*}
which, together with Cavalieri's principle (see, for instance,
\cite[Corollary~1.2.4]{G2024}),
further implies that
\begin{align*}
&\sum_{k=k_0+1}^\infty\sum_{i\in\mathbb{N}}
\left\{\frac{\lambda_{i,k}}{[\mu(B_{i,k})]^{
\frac{1}{p}-\frac{1}{p_0}}}\right\}^{p_0}
\\
&\quad\lesssim\int_0^\infty\int_X\alpha^{p_0-1}
\mathbf{1}_{\{x\in X:f^*_\rho(x)>2^{k_0+1}\}}(y)
\mathbf{1}_{\{x\in X:f^*_\rho(x)>\alpha\}}(y)\,d\mu(y)\,d\alpha
\\
&\quad\sim\left\|f^*_\rho
\mathbf{1}_{\Omega_{k_0}}\right\|_{L^{p_0}(X)}^{p_0}.
\end{align*}
From this and Theorem~\ref{6.15}
with $p$, $q$, and $r$ therein replaced,
respectively, by $p_0$, $p_0$, and $\infty$,
we infer that
\begin{align*}
\left\|b_{k_0}\right\|_{H^{p_0}_\ast(X)}^{p_0}
\lesssim\sum_{k=k_0+1}^\infty\sum_{i\in\mathbb{N}}
\left\{\frac{\lambda_{i,k}}{[\mu(B_{i,k})]^{
\frac{1}{p}-\frac{1}{p_0}}}\right\}^{p_0}
\lesssim\left\|f^*_\rho
\mathbf{1}_{\Omega_{k_0}}\right\|_{L^{p_0}(X)}^{p_0},
\end{align*}
which, combined with \eqref{2215} and \eqref{2216},
completes the proof of Lemma~\ref{fenjieHpq}.
\end{proof}

Using Lemma~\ref{fenjieHpq},
we obtain the following interpolation result.

\begin{lemma}
\label{L6.40}
Let $(X,\mathbf{q},\mu)$ be a quasi-ultrametric space
of homogeneous type with upper dimension $n\in(0,\infty)$.
Suppose that
$p_0\in(\frac{n}{n+\mathrm{ind\,}(X,\mathbf{q})},1]$,
$q\in(0,\infty]$,
$\theta\in(0,1)$, $p:=\frac{p_0}{1-\theta}$, and
$\beta,\gamma,\varepsilon\in(0,\infty)$
satisfy
$n(\frac{1}{p}-1)<\beta,\gamma<\varepsilon
\preceq\mathrm{ind\,}(X,\mathbf{q})$.
Then, as subspaces of
$({\mathcal{G}}^\varepsilon_0(\beta,\gamma))'$,
$$
\left(H^{p_0}_\ast(X),L^\infty(X)\right)_{\theta,q}=H^{p,q}_\ast(X)
$$
with equivalent quasi-norms.
\end{lemma}

\begin{remark}
Note that Lemma~\ref{L6.40} also holds when $p_0\in(1,\infty)$;
see Remark~\ref{rml95}.
\end{remark}

\begin{proof}[Proof of Lemma~\ref{L6.40}]
We first prove that
\begin{align}\label{22}
H^{p,q}_\ast(X)\subset
\left(H^{p_0}_\ast(X),L^\infty(X)\right)_{\theta,q}.
\end{align}
Let $f\in H^{p,q}_\ast(X)$. Then, by
Lemma~\ref{fenjieHpq},
we find that, for any $t\in(0,\infty)$
and $k\in\mathbb{Z}$,
\begin{align*}
K(t,f;H^{p_0}_\ast(X),L^\infty(X))
&\leq\left\|b_k\right\|_{H^{p_0}_\ast(X)}
+t\left\|g_k\right\|_{L^\infty(X)}
\nonumber\\
&\lesssim\left\|f^*_\rho\mathbf{1}_{\{x\in X:
f^*_\rho(x)>2^k\}}\right\|_{L^{p_0}(X)}+t2^k
\nonumber\\
&=\left\|f^*_\rho\sum_{j=0}^\infty\mathbf{1}_{\{x\in X:2^{k+j+1}
\ge f^*_\rho(x)>2^{k+j}\}}\right\|_{L^{p_0}(X)}+t2^k,
\end{align*}
which further implies that
\begin{align}\label{1315}
K(t,f;H^{p_0}_\ast(X),L^\infty(X))
&\leq\left[\sum_{j=0}^\infty\left\|f^*_\rho\mathbf{1}_{\{x\in X:2^{k+j+1}\ge
f^*_\rho(x)>2^{k+j}\}}\right\|_{L^{p_0}(X)}^{p_0}\right]^\frac{1}{p_0}+t2^k
\nonumber\\
&\lesssim\left[\sum_{j=0}^\infty2^{(k+j)p_0}\left\|\mathbf{1}_{\{x\in X:
f^*_\rho(x)>2^{k+j}\}}\right\|_{L^{p_0}(X)}^{p_0}\right]^\frac{1}{p_0}+t2^k
\nonumber\\
&=\left[\sum_{j=0}^\infty2^{(k+j)p_0}\mu(\Omega_{k+j})\right]^\frac{1}{p_0}+t2^k,
\end{align}
where $f^*_\rho$ and $\Omega_{k}$ are the same as,
respectively, in Definition \ref{maxfunc}(iii)
and \eqref{omegak}
and where $b_k$ and $g_k$ are the same as in Lemma~\ref{fenjieHpq}.
For any $t\in(0,\infty)$, let
$$
k(t):=\inf\left\{\tau\in\mathbb{Z}:\left[\sum_{j=0}^\infty
2^{jp_0}\mu(\Omega_{\tau+j})\right]^\frac{1}{p_0}\leq t\right\}.
$$
From this and \eqref{1315},
we deduce that, for any $t\in(0,\infty)$,
\begin{align*}
K\left(t,f;H^{p_0}_\ast(X),L^\infty(X)\right)\lesssim t2^{k(t)},
\end{align*}
which further implies that
\begin{align}\label{1516}
&\int_0^\infty\left[t^{-\theta}K(t,f;H^{p_0}_\ast(X),
L^\infty(X))\right]^q\frac{dt}{t}
\nonumber\\
&\quad\lesssim\sum_{l\in\mathbb{Z}}
2^{lq}\int_{\{t\in(0,\infty):
2^l<2^{k(t)}\leq 2^{l+1}\}}t^{(1-\theta)q}\frac{dt}{t}\nonumber\\
&\quad\lesssim\sum_{l\in\mathbb{Z}}2^{lq}\int_0^{[\sum_{j=0}^\infty2^{jp_0}
\mu(\Omega_{l+j})]^\frac{1}{p_0}}t^{(1-\theta)q}\frac{dt}{t}
\nonumber\\
&\quad\sim\sum_{l\in\mathbb{Z}}2^{lq}
\left[\sum_{j=0}^\infty2^{jp_0}\mu(\Omega_{l+j})\right]^\frac{(1-\theta)q}{p_0}.
\end{align}

Next, we consider the following two cases for $\frac{(1-\theta)q}{p_0}$.

\emph{Case (1)} $\frac{(1-\theta)q}{p_0}\in(0,1]$.
In this case, we obtain
\begin{align*}
&\sum_{l\in\mathbb{Z}}2^{lq}
\left[\sum_{j=0}^\infty2^{jp_0}\mu(\Omega_{l+j})
\right]^\frac{(1-\theta)q}{p_0}\\
&\quad\leq\sum_{l\in\mathbb{Z}}2^{lq}
\sum_{j=0}^\infty2^{(1-\theta)jq}
\left[\mu(\Omega_{l+j})\right]^\frac{(1-\theta)q}{p_0}
\quad\text{by Lemma~\ref{discreteinequality}}\\
&\quad=\sum_{j=0}^\infty\sum_{l\in\mathbb{Z}}2^{lq}
2^{(1-\theta)jq}\left[\mu(\Omega_{l+j})\right]^\frac{(1-\theta)q}{p_0}
\quad\text{by Tonelli's theorem for series}\\
&\quad=\sum_{j=0}^\infty\sum_{l'\in\mathbb{Z}}2^{(l'-j)q}
2^{(1-\theta)jq}\left[\mu(\Omega_{l'})\right]^\frac{(1-\theta)q}{p_0}
\quad\text{by a change of variable}\\
&\quad=\sum_{l'\in\mathbb{Z}}2^{l'q}\left[\mu(\Omega_{l'})
\right]^\frac{(1-\theta)q}{p_0}\sum_{j=0}^\infty
2^{-\theta jq}
\quad\text{by Tonelli's theorem for series}\\
&\quad\sim\sum_{l'\in\mathbb{Z}}2^{l'q}
\left[\mu(\Omega_{l'})\right]^\frac{(1-\theta)q}{p_0}.
\end{align*}

\emph{Case (2)}
$\frac{(1-\theta)q}{p_0}\in(1,\infty)$.
In this case, let $\eta:=\frac{(1-\theta)q}{p_0}$.
Then, from H\"older's inequality,
we infer that, for any $\delta\in(0,\frac{\theta}{1-\theta})$,
\begin{align*}
&\sum_{l\in\mathbb{Z}}2^{lq}
\left[\sum_{j=0}^\infty2^{jp_0}\mu(\Omega_{l+j})\right]^\frac{(1-\theta)q}{p_0}\\
&\quad\leq\sum_{l\in\mathbb{Z}}2^{lq}
\left(\sum_{j=0}^\infty2^{-\delta p_0j\eta'}\right)^\frac{\eta}{\eta'}
\sum_{j=0}^\infty2^{(1+\delta)p_0j\eta}\left[\mu(\Omega_{l+j})\right]^{\eta}\\
&\quad\sim\sum_{l\in\mathbb{Z}}2^{lq}
\sum_{j=0}^\infty2^{(1+\delta)p_0j\eta}\left[\mu(\Omega_{l+j})\right]^{\eta},
\end{align*}
which, together with Tonelli's theorem for series twice
and a change of variable,
further implies that
\begin{align*}
&\sum_{l\in\mathbb{Z}}2^{lq}
\left[\sum_{j=0}^\infty2^{jp_0}\mu(\Omega_{l+j})\right]^\frac{(1-\theta)q}{p_0}\\
&\quad\lesssim\sum_{j=0}^\infty\sum_{l'\in\mathbb{Z}}2^{(l'-j)q}
2^{(1+\delta)p_0j\eta}\left[\mu(\Omega_{l'})\right]^{\eta}\\
&\quad=\sum_{l'\in\mathbb{Z}}2^{l'q}\left[\mu(\Omega_{l'})
\right]^{\frac{(1-\theta)q}{p_0}}\sum_{j=0}^\infty
2^{[(1+\delta)(1-\theta)-1]jq}\sim\sum_{l'\in\mathbb{Z}}2^{l'q}
\left[\mu(\Omega_{l'})\right]^{\frac{(1-\theta)q}{p_0}}.
\end{align*}
Combining the above two cases,
by \eqref{1516}, $p=\frac{p_0}{1-\theta}$, and Proposition~\ref{Lpqnorm=c,d},
we conclude that
\begin{align*}
&\int_0^\infty\left[t^{-\theta}K(t,f;H^{p_0}_\ast(X)
L^\infty(X))\right]^q\frac{dt}{t}
\\
&\quad\lesssim\sum_{l'\in\mathbb{Z}}2^{l'q}
\left[\mu(\Omega_{l'})\right]^{\frac{(1-\theta)q}{p_0}}
\sim\left\|f^*_\rho\right\|_{L^{p,q}(X)}^q
=\|f\|_{H^{p,q}_\ast(X)}^q<\infty,
\end{align*}
which further implies that
$f\in (H^{p_0}_\ast(X),L^\infty(X))_{\theta,q}$.
This finishes the proof of \eqref{22}.

It remains to show that
\begin{align}\label{33}
\left(H^{p_0}_\ast(X),L^\infty(X)\right)_{\theta,q}
\subset H^{p,q}_\ast(X).
\end{align}
Let $f\in(H^{p_0}_\ast(X),L^\infty(X))_{\theta,q}$.
Then there exist $f_1\in H^{p_0}_\ast(X)$
and $f_2\in L^\infty(X)$ such that $f=f_1+f_2$
and
\begin{align*}
\int_0^\infty t^{-\theta q}
\left[\|f_1\|_{H^{p_0}_\ast(X)}
+\|f_2\|_{L^\infty(X)}\right]^q\frac{dt}{t}
\sim\|f\|_{(H^{p_0}_\ast(X),L^\infty(X))_{\theta,q}}^q.
\end{align*}
From Proposition~\ref{HL-grandmax},
we deduce that the grand maximal operator
is bounded on $L^\infty(X)$.
Now, we let
$$
E_1:=\left\{x\in X:f^*_\rho(x)\leq2(f_1)_\rho^*(x)\right\}
$$
and
$$
E_2:=\left\{x\in X:f^*_\rho(x)\leq2(f_2)_\rho^*(x)\right\}.
$$
Then we can write
$$
f^*_\rho=f^*_\rho\mathbf{1}_{E_1}
+f^*_\rho\mathbf{1}_{E_2\setminus E_1}\in L^{p_0}(X)+L^\infty(X).
$$
This, together with Remark~\ref{LczLpq-x},
implies that
\begin{align*}
\left\|f^*_\rho\right\|_{L^{p,q}(X)}
&=\left\|f^*_\rho\right\|_{(L^{p_0}(X),L^\infty(X))_{\theta,q}}
\\
&\leq\left\{\int_0^\infty t^{-\theta q}
\left[\left\|f^*_\rho\mathbf{1}_{E_1}\right\|_{L^{p_0}(X)}
+\left\|f^*_\rho\mathbf{1}_{E_2\setminus E_1}\right\|_{L^\infty(X)}
\right]^q\frac{dt}{t}\right\}^\frac{1}{q}
\\
&\lesssim\left\{\int_0^\infty t^{-\theta q}
\left[\left\|(f_1)^*_\rho\right\|_{L^{p_0}(X)}
+\left\|(f_2)^*_\rho\right\|_{L^\infty(X)}\right]^q
\frac{dt}{t}\right\}^\frac{1}{q}\\
&\lesssim\left\{\int_0^\infty t^{-\theta q}
\left[\|f_1\|_{H^{p_0}_\ast(X)}
+\|f_2\|_{L^\infty(X)}\right]^q\frac{dt}{t}\right\}^\frac{1}{q}\\
&\sim\|f\|_{(H^{p_0}_\ast(X),L^\infty(X))_{\theta,q}}
\end{align*}
and hence $f\in H^{p,q}_\ast(X)$,
which completes the proof of \eqref{33}.
This, combined with \eqref{22}, then
finishes the proof of Lemma~\ref{L6.40}.
\end{proof}

We are ready to prove Theorem~\ref{chazhi}(i).

\begin{proof}[Proof of Theorem~\ref{chazhi}(i)]
Suppose first $p_1\neq p_2$ and $\frac{1}{p}=
\frac{1-\theta}{p_1}+\frac{\theta}{p_2}$.
Let
$$
r\in\left(\max\left\{\frac{n}{n+\beta},\,
\frac{n}{n+\gamma}\right\},\min\{p_1,\,p_2\}\right)\cap(0,1].
$$
Then there exist constants $s_1,s_2\in(0,1)$ such that
\begin{align}\label{811-1}
\frac{1}{p_1}=\frac{1-s_1}{r}
\ \ \text{and}\ \
\frac{1}{p_2}=\frac{1-s_2}{r}.
\end{align}
Let $s:=(1-\theta)s_1+\theta s_2$.
By this and \eqref{811-1}, we find that
$$
\frac{1}{p}=\frac{1-\theta}{p_1}
+\frac{\theta}{p_2}=\frac{(1-s_1)(1-\theta)+\theta(1-s_2)}{r}
=\frac{1-s}{r},
$$
which, together with Lemma~\ref{L6.40}
and \cite[Theorem 3.5.3]{bl76},
further implies that
\begin{align*}
\left(H^{p_1,q_1}_\ast(X),H^{p_2,q_2}_\ast(X)\right)_{\theta,q}
&=\left(\left(H^r_\ast(X),L^\infty(X)\right)_{s_1,q_1},
\left(H^r_\ast(X),L^\infty(X)\right)_{s_2,q_2}\right)_{\theta,q}
\\
&=\left(H^r_\ast(X),L^\infty(X)\right)_{s,q}
=H^{p,q}_\ast(X).
\end{align*}
This finishes the proof of Theorem~\ref{chazhi}(i).
\end{proof}

Using Theorem~\ref{chazhi}(i), we obtain the following
conclusion.

\begin{theorem}\label{Hpq=Lpq}
Let $(X,\mathbf{q},\mu)$ be a quasi-ultrametric space
of homogeneous type. Assume that
$p\in(1,\infty)$, $q\in(0,\infty]$,
and $\beta,\gamma,\varepsilon\in(0,\infty)$ satisfy
$0<\beta,\gamma<\varepsilon\preceq\mathrm{ind\,}(X,\mathbf{q})$.
Then
$$
H^{p,q}_\ast(X)=L^{p,q}(X)
$$
in the following sense:
the linear map
$\Phi:L^{p,q}(X)\to H^{p,q}_\ast(X)$ that is
defined by setting, for any $f\in L^{p,q}(X)$,
$\Phi(f):=\Lambda_f$, where
$\Lambda_f\in({\mathcal{G}}^\varepsilon_0(\beta,\gamma))'$
is the distribution on $\mathcal{G}^\varepsilon_0(\beta,\gamma)$
defined in \eqref{1503},
is well defined, linear, bijective,
and has the property that, for any
$f\in L^{p,q}(X)$, $\|f\|_{L^{p,q}(X)}\sim\|\Lambda_f\|_{H^{p,q}_\ast(X)}$,
where the positive equivalence constants are independent of $f$.
\end{theorem}

To show Theorem~\ref{Hpq=Lpq},
we also need the following lemma.

\begin{lemma}\label{1419}
Let $(X,\mathbf{q},\mu)$ be a quasi-ultrametric space
of homogeneous type.
Assume that $p\in(1,\infty]$ and $0<\beta,\gamma<
\varepsilon\preceq\mathrm{ind\,}(X,\mathbf{q})$.
Then
$$
H_*^p(X)=H^p_+(X)=H^p_\theta(X)=L^p(X)
$$
in the following sense:
\begin{enumerate}
\item[\rm(i)]
If $f\in(\mathcal{G}^\varepsilon_0(\beta,\gamma))'$
belongs to $H_*^p(X)$ [resp. $H^p_+(X)$ or $H^p_\theta(X)$],
then there exists $\widetilde{f}\in L^p(X)$ such that,
for any $\phi\in\mathcal{G}^\varepsilon_0(\beta,\gamma)$,
\begin{align}\label{1434}
\left\langle f,\phi\right\rangle=\int_X\widetilde{f}(x)\phi(x)\,d\mu(x)
\end{align}
and $\|\widetilde{f}\|_{L^p(X)}\lesssim\|f\|_{H_*^p(X)}$
[resp. $\|\widetilde{f}\|_{L^p(X)}\lesssim\|f\|_{H_\theta^p(X)}$
or $\|\widetilde{f}\|_{L^p(X)}\lesssim\|f\|_{H_+^p(X)}$],
where the implicit positive constants are independent of $f$.

\item[\rm(ii)]
Any $f\in L^p(X)$ induces a distribution on
$\mathcal{G}^\varepsilon_0(\beta,\gamma)$
as in \eqref{1434}, still denoted by $f$,
such that $f\in H_*^p(X)$
[resp. $f\in H_\theta^p(X)$ or $f\in H_+^p(X)$]
and $\|f\|_{H_*^p(X)}\lesssim\|f\|_{L^p(X)}$
[resp. $\|f\|_{H_\theta^p(X)}\lesssim\|f\|_{L^p(X)}$
or $\|f\|_{H_+^p(X)}\lesssim\|f\|_{L^p(X)}$],
where the implicit positive constants are independent of $f$.
\end{enumerate}
\end{lemma}
\begin{proof}
We first prove (i). From Lemma~\ref{MMM}, it follows that
\begin{align*}
H^{p,q}_*(X)\subset H^{p,q}_\theta(X)\subset H^{p,q}_+(X),
\end{align*}
and hence it suffices to show (i)
for any $f\in H_+^p(X)$.
Let $f\in(\mathcal{G}^\varepsilon_0(\beta,\gamma))'$
satisfy $\mathcal{M}^+f\in L^p(X)$.
By Definition~\ref{maxfunc}(i),
we conclude that
$\{\mathcal{S}_{2^{-k}}f\}_{k\in\mathbb{Z}}$
is a bounded set in $L^p(X)$. Note that $p'\in[1,\infty)$
and $L^{p'}(X)$ is separable.
From this and Alaoglu's theorem
(see, for instance, \cite[Theorem 3.17]{r91}),
we infer that there exist $\widetilde{f}\in L^{p'}(X)$
and a sequence $\{k_j\}_{j\in\mathbb{N}}$ satisfying that
$k_j\to\infty$ as $j\to\infty$ such that
$\mathcal{S}_{2^{-k_j}}f\to\widetilde{f}$ in the weak-$*$
topology of $L^p(X)$,
which, together with H\"older's inequality,
further implies that, for any $g\in L^{p'}(X)$
with $\|g\|_{L^{p'}(X)}\leq1$,
\begin{align*}
\left|\int_X\widetilde{f}(x)g(x)\,d\mu(x)\right|
&=\lim_{j\to\infty}
\left|\int_X\mathcal{S}_{2^{-k_j}}f(x)g(x)\,d\mu(x)\right|\\
&\leq\lim_{j\to\infty}
\left\|\mathcal{S}_{2^{-k_j}}f\right\|_{L^p(X)}\|g\|_{L^{p'}(X)}
\leq\left\|\mathcal{M}^+f\right\|_{L^p(X)}=\|f\|_{H_+^p(X)}.
\end{align*}
Taking the supremum over all $g\in L^{p'}(X)$
with $\|g\|_{L^{p'}(X)}\leq1$,
we obtain $\|\widetilde{f}\|_{L^p(X)}\lesssim\|f\|_{H_+^p(X)}$.
This finishes the proof of (i) for any $f\in H_+^p(X)$.

Next, we prove that (ii) holds for $H_*^p(X)$.
Let $f\in L^p(X)$. By Proposition~\ref{distfncttype},
we find that $f$ induces a distribution on
$\mathcal{G}^\varepsilon_0(\beta,\gamma)$ as in \eqref{1434}. On the other hand,
from Proposition~\ref{HL-grandmax}
and Lemma~\ref{M},
we deduce that
$$
\|f\|_{H^p_*(X)}=\left\|f_\rho^*\right\|_{L^p(X)}
\lesssim\left\|\mathcal{M}_\rho f\right\|_{L^p(X)}
\lesssim\left\|f\right\|_{L^p(X)},
$$
which completes the proof of (ii) for $H_*^p(X)$.
This, together with the conclusion that (i)
holds for any $f\in H_+^p(X)$, and Lemma~\ref{MMM},
then finishes the proof of Lemma~\ref{1419}.
\end{proof}

\begin{proof}[Proof of Theorem~\ref{Hpq=Lpq}]
From Remark~\ref{R6.6} and Lemma~\ref{1419},
we infer that \begin{align}\label{1513}
H^{p,p}_\ast(X)=H_*^p(X)=L^p(X).
\end{align}
Note that there exist $p_1,p_2\in(1,\infty)$
and $\theta\in(0,1)$
satisfying
$$
p_1\neq p_2
\ \ \text{and}\ \
\frac{1}{p}=\frac{1-\theta}{p_1}+\frac{\theta}{p_2}.
$$
By this, Theorem~\ref{chazhi}(i) with $q_1$ and $q_2$
therein replaced, respectively, by $p_1$ and $p_2$,
\eqref{1513},
and Lemma~\ref{LczLpq}, we conclude that
\begin{align*}
H^{p,q}_\ast(X)
=\left(H^{p_1,p_1}_\ast(X),H^{p_2,p_2}_\ast(X)\right)_{\theta,q}
=\left(L^{p_1}(X),L^{p_2}(X)\right)_{\theta,q}=L^{p,q}(X),
\end{align*}
which completes the proof of Theorem~\ref{Hpq=Lpq}.
\end{proof}

Using Theorem~\ref{Hpq=Lpq},
we finally show Theorem~\ref{chazhi}(ii).

\begin{proof}[Proof of Theorem~\ref{chazhi}(ii)]
We consider the following two cases for $p$.

\emph{Case (1)} $p\in(\frac{n}{n+\mathrm{ind\,}(X,\mathbf{q})},1]$.
In this case, we first prove
\begin{align}\label{1654}
\left(H^{p,q_1}_\ast(X),H^{p,q_2}_\ast(X)\right)_{\theta,q}
\hookrightarrow H^{p,q}_\ast(X).
\end{align}
From Lemma~\ref{LczLpq} and
the sublinearity of the grand maximal operator,
we deduce that,
for any $f\in(H^{p,q_1}_\ast(X),H^{p,q_2}_\ast(X))_{\theta,q}$,
\begin{align*}
\|f\|_{H^{p,q}_\ast(X)}=\|f_\rho^*\|_{L^{p,q}(X)}
\sim\|f_\rho^*\|_{(L^{p,q_1}(X),L^{p,q_2}(X))_{\theta,q}}
\leq\|f\|_{(H^{p,q_1}_\ast(X),H^{p,q_2}_\ast(X))_{\theta,q}},
\end{align*}
which completes the proof of \eqref{1654}.

Now, we show that
\begin{align}\label{1655}
H^{p,q}_\ast(X)\hookrightarrow
\left(H^{p,q_1}_\ast(X),H^{p,q_2}_\ast(X)\right)_{\theta,q}.
\end{align}
Without of generality, we may assume that $q_1\leq q_2$.
Let $f\in H^{p,q}_\ast(X)$.
Then, by Theorems~\ref{6.16} and~\ref{6.17}, we find that
there exist a sequence $\{a_{i,k}\}_{i\in\mathbb{N},k\in\mathbb{Z}}$
of $(p,2)$-atoms and a sequence
$\{\lambda_{i,k}\}_{i\in\mathbb{N},k\in\mathbb{Z}}$ in
$\mathbb{R}_+$ satisfying all conditions in Definition~\ref{6.3};
in particular, we have
$f=\sum_{k\in\mathbb{Z}}\sum_{i\in\mathbb{N}}
\lambda_{i,k}a_{i,k}$ in $(\mathcal{G}^\varepsilon_0(\beta,\gamma))'$
and
\begin{align}\label{22481}
\left[\sum_{k\in\mathbb{Z}}\left(\sum_{i\in\mathbb{N}}
\lambda_{i,k}^p\right)^\frac{q}{p}\right]^\frac{1}{q}
\sim\|f\|_{H^{p,q}_\ast(X)},
\end{align}
where the positive equivalence constants are
independent of $f$.
Fix $k_0\in\mathbb{N}$.
Let
$$
f_{1,k_0}:=\sum_{|k|\leq k_0}\sum_{i\in\mathbb{N}}\lambda_{i,k}a_{i,k}
\ \ \text{and}\ \
f_{2,k_0}:=\sum_{|k|>k_0}\sum_{i\in\mathbb{N}}
\lambda_{i,k}a_{i,k},
$$
both of which converge in $(\mathcal{G}^\varepsilon_0(\beta,\gamma))'$
by Theorem~\ref{6.15}.
From Theorem~\ref{6.15} and $q_2\ge q$, we further infer that
$$
f_{1,k_0}\in H^{p,q_1}_\ast(X)
\ \ \text{with}\ \
\|f_{1,k_0}\|_{H^{p,q_1}_\ast(X)}
\lesssim\left[\sum_{|k|\leq k_0}\left(\sum_{i\in\mathbb{N}}
\lambda_{i,k}^p\right)^\frac{q_1}{p}\right]^\frac{1}{q_1}
$$
and
$$
f_{2,k_0}\in H^{p,q_2}_\ast(X)
\ \ \text{with}\ \
\|f_{2,k_0}\|_{H^{p,q_2}_\ast(X)}
\lesssim\left[\sum_{|k|> k_0}\left(\sum_{i\in\mathbb{N}}
\lambda_{i,k}^p\right)^\frac{q_2}{p}\right]^\frac{1}{q_2},
$$
which, combined with Holmstedt's formula
(see, for instance, \cite[Chapter~5, Theorem~2.1]{bs1988}),
further implies that,
for any given $t\in(0,\infty)$,
there exists $K_t\in\mathbb{N}$
such that
\begin{align*}
&K\left(t,f;H^{p,q_1}_\ast(X),H^{p,q_2}_\ast(X)\right)\\
&\quad\lesssim\left[\sum_{|k|\leq K_t}\left(\sum_{i\in\mathbb{N}}
\lambda_{i,k}^p\right)^\frac{q_1}{p}\right]^\frac{1}{q_1}
+t\left[\sum_{|k|> K_t}\left(\sum_{i\in\mathbb{N}}
\lambda_{i,k}^p\right)^\frac{q_2}{p}\right]^\frac{1}{q_2}\\
&\quad\sim K\left(t,\left\{\left(\sum_{i\in\mathbb{N}}
\lambda_{i,k}^p\right)^\frac{1}{p}
\right\}_{k\in\mathbb{Z}};\ell^{q_1},\ell^{q_2}\right),
\end{align*}
where the implicit positive constants are independent of $t$.
By this, Lemma~\ref{LczLpq}, and \eqref{22481},
we conclude that
$f\in(H^{p,q_1}_\ast(X),H^{p,q_2}_\ast(X))_{\theta,q}$ and
\begin{align*}
\|f\|_{(H^{p,q_1}_\ast(X),H^{p,q_2}_\ast(X))
_{\theta,q}}
&\lesssim\left\|\left\{\left(\sum_{i\in\mathbb{N}}
\lambda_{i,k}^p\right)^\frac{1}{p}\right\}_{k\in\mathbb{Z}}
\right\|_{(\ell^{q_1},\ell^{q_2})_{\theta,q}}\\
&\sim
\left\|\left\{\left(\sum_{i\in\mathbb{N}}\lambda_{i,k}^p
\right)^\frac{1}{p}\right\}_{k\in\mathbb{Z}}\right\|_{\ell^q}
\sim\|f\|_{H^{p,q}_\ast(X)},
\end{align*}
which completes the proof of \eqref{1655}
and hence (ii) in this case.

\emph{Case (2)} $p\in(1,\infty)$. In this case,
from Theorem~\ref{Hpq=Lpq} and Lemma~\ref{LczLpq},
we deduce that
\begin{align*}
H^{p,q}_\ast(X)=L^{p,q}(X)
=\left(L^{p,q_1}(X),L^{p,q_2}(X)\right)_{\theta,q}
=\left(H^{p,q_1}_\ast(X),H^{p,q_2}_\ast(X)\right)_{\theta,q},
\end{align*}
which completes the proof of Theorem~\ref{chazhi}(ii) 
and hence Theorem~\ref{chazhi}.
\end{proof}

\begin{remark}
\label{rml95}
Lemma~\ref{L6.40} also holds
for any given $p_0\in(1,\infty)$, $q\in(0,\infty]$,
and $p:=\frac{p_0}{1-\theta}$.
Indeed, by the fact that $H^{p_0}_\ast(X)=L^{p_0}(X)$
for any $p_0\in(1,\infty)$, Lemma~\ref{LczLpq},
and Theorem~\ref{Hpq=Lpq},
we find that
\begin{align*}
\left(H^{p_0}_\ast(X),L^\infty(X)\right)_{\theta,q}
=\left(L^{p_0}(X),L^\infty(X)\right)_{\theta,q}
=L^{p,q}(X)=H^{p,q}_\ast(X).
\end{align*}
\end{remark}

\section{Dual Spaces of Hardy--Lorentz Spaces}\label{Section5.5}

In this section, we aim to explore the topological dual of the maximal
Hardy--Lorentz spaces $H^{p,q}_*(X)$
for the optimal ranges
$$
p\in\left(\frac{n}{n+{\mathrm{ind\,}}(X,\mathbf{q})},\infty\right)
\ \ \text{and}\ \
q\in(0,\infty)
$$
on quasi-ultrametric spaces of homogeneous type
$(X,\mathbf{q},\mu)$
with upper dimension $n\in(0,\infty)$.
Since, when $p\in(1,\infty)$ and $q\in(0,\infty)$,
Theorem~\ref{Hpq=Lpq} implies that $H^{p,q}_*(X)$
can be identified with $L^{p,q}(X)$,
it follows from \cite[(v) and (vi) of Theorem~1.4.16]{G2014C} that
\begin{align}\label{DHpq}
\left(H^{p,q}_*(X)\right)'=\left(L^{p,q}(X)\right)'=
\begin{cases}
L^{p',\infty}(X)&\text{if }p\in(1,\infty)\text{ and }q\in(0,1],\\
L^{p',q'}(X)&\text{if }p\in(1,\infty)\text{ and }q\in(1,\infty).
\end{cases}
\end{align}

Next, we turn to the case where $p\leq1$.
To this end, we first introduce the concept
of Campanato--Lorentz spaces on quasi-ultrametric spaces
of homogeneous type.
Recall that, for any given $p,q\in(0,\infty]$,
the \emph{space} $\ell^q(\ell^p)$\index[symbols]{L@$\ell^q(\ell^p)$}
is defined to be
the set of all sequences
$\{\lambda_{i,j}\}_{i,j\in\mathbb{N}}$ in $\mathbb{C}$
such that
$$
\left\|\{\lambda_{i,j}\}_{i,j\in\mathbb{N}}\right\|_{\ell^q(\ell^p)}
:=\left[\sum_{i\in\mathbb{N}}\left(\sum_{j\in\mathbb{N}}
\left|\lambda_{i,j}\right|^p\right)^\frac{q}{p}\right]^\frac{1}{q}<\infty
$$
with the usual modification when $p=\infty$ or $q=\infty$.

\begin{definition}\label{20191}
Let $(X,\mathbf{q},\mu)$ be a quasi-ultrametric space
of homogeneous type and $\rho\in\mathbf{q}$ be a regularized
quasi-ultrametric.
Let $p,q\in(0,\infty]$ and $r\in[1,\infty)$.
The \emph{Campanato--Lorentz space}\index[words]{Campanato--Lorentz space}
$\mathcal{C}^{p,q,r}(X,\rho)$\index[symbols]{C@$\mathcal{C}^{p,q,r}(X,\rho)$}
is defined to be the set of
all $f\in L^r_{\mathrm{loc}}(X)$ such that
\begin{align*}
\|f\|_{\mathcal{C}^{p,q,r}(X)}:&=\sup
\left[\sum_{i=1}^m\left(\sum_{j=1}^v
\left|\lambda_{i,j}\right|^p\right)^\frac{q}{p}
\right]^{-\frac{1}{q}}
\sum_{i=1}^m\sum_{j=1}^v
\left|\lambda_{i,j}\right|
\left[\mu(B_{i,j})\right]^{1-\frac{1}{p}}\\
&\qquad\times\left[\fint_{B_{i,j}}
\left|f(x)-m_{B_{i,j}}(f)\right|^r\,d\mu(x)\right]^\frac{1}{r}\\
&<\infty,
\end{align*}
where the supremum is taken over all
$m,v\in\mathbb{N}$ and all
sequences $\{B_{i,j}\}_{i\in\mathbb{N}\cap[1,m],j\in\mathbb{N}\cap[1,v]}$
of $\rho$-balls and $\{\lambda_{i,j}\}_{i,j\in\mathbb{N}}\in\ell^q(\ell^p)$
satisfying both (i) and (ii) in Definition~\ref{5d1}
and where, for any $\rho$-ball $B\subset X$, $m_B(f)$ is the same as in \eqref{mBf}.
\end{definition}

\begin{remark}\label{1545}
Let the notation be as in Definition~\ref{20191}.
\begin{enumerate}
\item[\rm(i)]
The Campanato--Lorentz space
$\mathcal{C}^{p,q,r}(X,\rho)$ is a Banach space.

\item[\rm(ii)]
Note that the Campanato--Lorentz space
$\mathcal{C}^{p,q,r}(X,\rho)$ when $p,q\in(0,1]$
reduces to the \emph{Campanato space}
\index[words]{Campanato space} $\mathcal{C}^{p,r}(X,\rho)$,
\index[symbols]{C@$\mathcal{C}^{p,r}(X,\rho)$}
which is defined to be the set of all
$f\in L^r_{\mathrm{loc}}(X)$ such that
\begin{align*}
\|f\|_{\mathcal{C}^{p,r}(X)}:=\sup_{{\rm\rho-balls}\,
B\subset X}\left[\mu(B)\right]^{1-\frac{1}{p}}
\left[\fint_{B}
\left|f(x)-m_{B}(f)\right|^r\,d\mu(x)\right]^\frac{1}{r}<\infty.
\end{align*}
Indeed, it is easy to see that, for any $f\in\mathcal{C}^{p,q,r}(X,\rho)$,
\begin{align}\label{1042}
\|f\|_{\mathcal{C}^{p,r}(X,\rho)}\leq\|f\|_{\mathcal{C}^{p,q,r}(X,\rho)}.
\end{align}
On the other hand, assume that $f\in\mathcal{C}^{p,r}(X,\rho)$. Then,
for any $m,v\in\mathbb{N}$ and sequences
$\{B_{i,j}\}_{i\in\mathbb{N}\cap[1,m],j\in\mathbb{N}\cap[1,v]}$
and $\{\lambda_{i,j}\}_{i\in\mathbb{N}\cap[1,m],j\in\mathbb{N}\cap[1,v]}$
satisfying all the conditions in Definition~\ref{20191},
\begin{align*}
&\left[\sum_{i=1}^m\left(\sum_{j=1}^v
\left|\lambda_{i,j}\right|^p\right)^\frac{q}{p}
\right]^{-\frac{1}{q}}
\sum_{i=1}^m\sum_{j=1}^v
\left|\lambda_{i,j}\right|
\left[\mu(B_{i,j})\right]^{1-\frac{1}{p}}\\
&\qquad\times\left[\fint_{B_{i,j}}
\left|f(x)-m_{B_{i,j}}(f)\right|^r\,d\mu(x)\right]^\frac{1}{r}\\
&\quad\leq\left(\sum_{i=1}^m\sum_{j=1}^v
\left|\lambda_{i,j}\right|
\right)^{-1}
\sum_{i=1}^m\sum_{j=1}^v
\left|\lambda_{i,j}\right|
\left[\mu(B_{i,j})\right]^{1-\frac{1}{p}}\\
&\qquad\times\left[\fint_{B_{i,j}}
\left|f(x)-m_{B_{i,j}}(f)\right|^r\,d\mu(x)\right]^\frac{1}{r}\\
&\quad\leq\left(\sum_{i=1}^m\sum_{j=1}^v
\left|\lambda_{i,j}\right|
\right)^{-1}
\sum_{i=1}^m\sum_{j=1}^v
\left|\lambda_{i,j}\right|\|f\|_{\mathcal{C}^{p,r}(X,\rho)}
=\|f\|_{\mathcal{C}^{p,r}(X,\rho)}.
\end{align*}
Taking the supremum over all
$m,v\in\mathbb{N}$ and all
sequences $\{B_{i,j}\}_{i\in\mathbb{N}\cap[1,m],j\in\mathbb{N}\cap[1,v]}$
of $\rho$-balls and $\{\lambda_{i,j}\}_{i,j\in\mathbb{N}}\in\ell^q(\ell^p)$
satisfying both (i) and (ii) in Definition~\ref{5d1}
and using \eqref{1042},
we then conclude that, for any $p,q\in(0,1]$ and $r\in[1,\infty)$,
$$
\mathcal{C}^{p,q,r}(X,\rho)=\mathcal{C}^{p,r}(X,\rho).
$$
Moreover, from \cite[Theorem~5]{MaSe79i}
(a similar geometric characterization can also be found in \cite[Theorem~3.2
and Corollary~3.11]{wyy2024}),
it follows that,
if $(X,\rho,\mu)$ is an $n$-Ahlfors regular
quasi-ultrametric space for some $n\in(0,\infty)$,
then, for any $p\in(0,1)$, $q\in(0,1]$, and $r\in[1,\infty)$,
\begin{align}\label{1620}
\mathcal{C}^{p,q,r}(X,\rho)=\mathcal{C}^{p,r}(X,\rho)=\dot{C}^{n(\frac{1}{p}-1)}(X,\rho).
\end{align}
In addition, for any $q\in(0,1]$ and $r\in[1,\infty)$,
$$
\mathcal{C}^{1,q,r}(X,\rho)=\mathrm{BMO}_r(X,\rho),
$$
where
$$
\mathrm{BMO}_r(X,\rho):=\left\{f\in L^r_{\mathrm{loc}}(X):
\sup_{{\rm \rho-balls}\,B\subset X}
\fint_B\left|f(x)-m_B(f)\right|^r\,d\mu(x)<\infty\right\}.
\index[symbols]{B@$\mathrm{BMO}_r(X)$}
$$
In particular, for any $q\in(0,1]$,
$\mathcal{C}^{1,q,1}(X,\rho)$ coincides with
$\mathrm{BMO}(X,\rho)$, namely the space of functions
with bounded mean oscillation introduced by
John and Nirenberg \cite{JN}.

\item[\rm(iii)]
By the proof of Theorem~\ref{duality},
we find that $\mathcal{C}^{p,q,r}(X,\rho)$,
with $p\in(\frac{n}{n+\mathrm{ind\,}(X,\mathbf{q})},1]$,
$q\in(0,\infty)$, and $r\in(1,\infty]$,
is independent of the choice of regularized
quasi-ultrametrics $\rho\in\mathbf{q}$.
In this sense, we may simply write $\mathcal{C}^{p,q,r}(X)$
instead of $\mathcal{C}^{p,q,r}(X,\rho)$.
\end{enumerate}
\end{remark}

The main theorem of this section reads as follows.

\begin{theorem}\label{duality}
Let $(X,\mathbf{q},\mu)$ be a quasi-ultrametric space
of homogeneous type. Assume that
$p\in(\frac{n}{n+\mathrm{ind\,}(X,\mathbf{q})},1]$,
$q\in(0,\infty)$, and $r\in(1,\infty]$.
Then the dual space of $H^{p,q}_*(X)$,
denoted by $(H^{p,q}_*(X))'$,
is $\mathcal{C}^{p,q,r'}(X)$
in the following sense:
\begin{enumerate}
\item[\rm(i)]
For any $g\in\mathcal{C}^{p,q,r'}(X)$,
the linear functional
\begin{align}\label{049}
G_g:f\mapsto G_g(f):=\int_Xf(x)g(x)\,d\mu(x),
\end{align}
initially defined for any
$f\in H^{p,q,r}_{\mathrm{at,\,fin}}(X)$ when $r\neq\infty$
and for any $f\in H^{p,q,\infty}_{\mathrm{at,\,fin}}(X)\cap\mathrm{UC}(X)$
when $r=\infty$,
has a bounded extension to $H^{p,q}_*(X)$.

\item[\rm(ii)]
For any $G\in(H^{p,q}_*(X))'$,
there exists a unique $g\in\mathcal{C}^{p,q,r'}(X)$
such that $G=G_g$ in $H^{p,q}_*(X)$,
where $G_g$ is the same as in \eqref{049}.
\end{enumerate}
Consequently, the linear map
$\widetilde{G}:\mathcal{C}^{p,q,r'}(X)
\to(H^{p,q}_*(X))'$, defined by setting
$\widetilde{G}(g)$ to be the bounded extension
of $G_g$ as in \eqref{049} into $H^{p,q}_*(X)$
for each  $g\in\mathcal{C}^{p,q,r'}(X)$,
is well defined, linear, and bijective,
and has the property that, for any
$g\in\mathcal{C}^{p,q,r'}(X)$,
\begin{align*}
\|g\|_{\mathcal{C}^{p,q,r'}(X)}
\sim\left\|\widetilde{G}_g\right\|_{(H^{p,q}_*(X))'},
\end{align*}
where the positive equivalence constants are independent of $g$.
\end{theorem}

\begin{remark}
\begin{enumerate}
\item[\rm(i)]
Theorem~\ref{duality} implies that,
for any $p\in(\frac{n}{n+\mathrm{ind\,}(X,\mathbf{q})},1]$,
$q\in(0,\infty)$, and $r_1,r_2\in[1,\infty)$,
$$
\mathcal{C}^{p,q,r_1}(X)=\mathcal{C}^{p,q,r_2}(X)
$$
with equivalent norms.

\item[\rm(ii)]
From Theorem~\ref{duality} and Remark~\ref{1545}(ii),
we infer that, for any $p\in(\frac{n}{n+\mathrm{ind\,}(X,\mathbf{q})},1]$,
$q\in(0,1]$, and $r\in[1,\infty)$,
\begin{align*}
\left(H^{p,q}_*(X)\right)'=\mathcal{C}^{p,r}(X).
\end{align*}
Hence, for any $p\in(\frac{n}{n+\mathrm{ind\,}(X,\mathbf{q})},1]$
and $r_1,r_2\in[1,\infty)$,
$$
\mathcal{C}^{p,r_1}(X)=\mathcal{C}^{p,r_2}(X)
$$
with equivalent norms.

\item[\rm(iii)]
When $p=q\in(\frac{n}{n+\mathrm{ind\,}(X,\mathbf{q})},1]$,
Theorem~\ref{duality} proves that the dual space
of the maximal Hardy space $H^p_*(X)$ is
the Campanato space $\mathcal{C}^{p,r}(X)$
which, when the underlying space is an $n$-Ahlfors regular quasi-ultrametric
space, is just the H\"older space $\dot{C}^{n(\frac{1}{p}-1)}(X)$
[see \eqref{1620}],
where $r\in[1,\infty)$;
in particular, Remark~\ref{1545}(ii) implies that
the dual space of $H^1_*(X)$ is $\mathrm{BMO\,}(X)$.
Thus, Theorem~\ref{duality} and \eqref{DHpq}
are generalizations of \cite[Theorem~7.22]{AlMi2015}
from the maximal Hardy space $H^p_*(X)$ with
$p\in(\frac{n}{n+\mathrm{ind\,}(X,\rho)},\infty)$
on $n$-Ahlfors regular quasi-ultrametric spaces
to the maximal Hardy--Lorentz space $H^{p,q}_*(X)$
with $p\in(\frac{n}{n+\mathrm{ind\,}(X,\rho)},\infty)$ and $q\in(0,\infty)$
on quasi-ultrametric spaces of homogeneous type.
\end{enumerate}
\end{remark}

\begin{proof}[Proof of Theorem~\ref{duality}]
Let $\rho\in\mathbf{q}$ be a
regularized quasi-ultrametric.
Without loss of generality,
we may assume that
$\mathcal{C}^{p,q,r}(X)$, $H^{p,q}_*(X)$,
$H^{p,q,r}_{\mathrm{at,\,fin}}(X)$,
and $H^{p,q,r}_{\mathrm{at}}(X)$ are constructed in terms of $\rho$.

We first show (i). To this end,
we consider the following two cases for $r$.

\emph{Case (1)} $r\in(1,\infty)$. In this case,
let $g\in\mathcal{C}^{p,q,r'}(X)$
and $f\in H^{p,q,r}_{\mathrm{at,\,fin}}(X)$.
Then there exist $m,v\in\mathbb{N}$,
a sequence of $(p,r)$-atoms
$\{a_{i,j}\}_{i\in\mathbb{N}\cap[1,m],j\in\mathbb{N}\cap[1,v]}$
supported, respectively, in $\rho$-balls
$\{B_{i,j}\}_{i\in\mathbb{N}\cap[1,m],j\in\mathbb{N}\cap[1,v]}$,
and a sequence $\{\lambda_{i,j}\}_{i\in\mathbb{N}\cap[1,m],j\in
\mathbb{N}\cap[1,v]}$ in $\mathbb{C}$ satisfying that
$$
\left[\sum_{i=1}^m\left(\sum_{j=1}^v
\left|\lambda_{i,j}\right|^p\right)^\frac{q}{p}
\right]^\frac{1}{q}\sim\|f\|_{H^{p,q,r}_{\mathrm{at,\,fin}}(X)}
$$
and
$f=\sum_{i=1}^m\sum_{j=1}^v
\lambda_{i,j}a_{i,j}$ $\mu$-almost everywhere,
where the positive equivalence constants are
independent of $f$.
By this, Definition~\ref{atom}(iii),
and H\"older's inequality,
we conclude that
\begin{align*}
\left|G_g(f)\right|
&\leq\sum_{i=1}^m\sum_{j=1}^v
\left|\lambda_{i,j}\right|\left|\int_Xa_{i,j}(x)g(x)\,d\mu(x)\right|\\
&=\sum_{i=1}^m\sum_{j=1}^v
\left|\lambda_{i,j}\right|\left|\int_X
a_{i,j}(x)\left[g(x)-m_{B_{i,j}}(g)\right]\,d\mu(x)\right|\\
&\leq\sum_{i=1}^m\sum_{j=1}^v
\left|\lambda_{i,j}\right|
\left\|a_{i,j}\right\|_{L^r(X)}
\left[\int_{B_{i,j}}\left|g(x)-m_{B_{i,j}}(g)
\right|^{r'}\,d\mu(x)\right]^\frac{1}{r'},
\end{align*}
which, together with Definition~\ref{atom}(ii)
and Theorem~\ref{6.29}(i),
further implies that
\begin{align*}
\left|G_g(f)\right|
&\leq\sum_{i=1}^m\sum_{j=1}^v
\left|\lambda_{i,j}\right|
\left[\mu(B_{i,j})\right]^{1-\frac{1}{p}}
\left[\fint_{B_{i,j}}\left|g(x)-m_{B_{i,j}}(g)
\right|^{r'}\,d\mu(x)\right]^\frac{1}{r'}\\
&\leq\left\|g\right\|_{\mathcal{C}^{p,q,r'}(X)}
\left[\sum_{i=1}^m\left(\sum_{j=1}^v
\left|\lambda_{i,j}\right|^p\right)^\frac{q}{p}
\right]^\frac{1}{q}\sim\left\|g\right\|_{\mathcal{C}^{p,q,r'}(X)}
\|f\|_{H^{p,q,r}_{\mathrm{at,\,fin}}(X)}\\
&\sim\left\|g\right\|_{\mathcal{C}^{p,q,r'}(X)}
\|f\|_{H^{p,q}_{*}(X)}.
\end{align*}
From this, Theorem~\ref{5r1},
and a standard density argument,
we deduce that (i) holds in this case.

\emph{Case (2)} $r=\infty$. Let $g\in\mathcal{C}^{p,q,1}(X)$.
In this case,
by Theorem~\ref{6.29}(ii),
Lemma~\ref{lemma81}, and repeating the proof
in Case (1), we find that
$g$ induces a bounded linear functional $G_g$ on $H^{p,q}_*(X)$,
which is initially defined on
$H^{p,q,\infty}_{\mathrm{at,\,fin}}(X)\cap\mathrm{UC}(X)$
by setting, for any
$f\in H^{p,q,\infty}_{\mathrm{at,\,fin}}(X)\cap\mathrm{UC}(X)$,
\begin{align}\label{22022}
G_g:f\mapsto G_g(f):=\int_Xf(x)g(x)\,d\mu(x),
\end{align}
has a bounded extension to $H^{p,q}_*(X)$.
Thus, to complete the proof of (i) in this case,
it remains to prove that, for any
$f\in H^{p,q,\infty}_{\mathrm{at,\,fin}}(X)$,
\begin{align}\label{2041}
G_g(f)=\int_Xf(x)g(x)\,d\mu(x).
\end{align}
Indeed,
let $f\in H^{p,q,\infty}_{\mathrm{at,\,fin}}(X)$.
From Lemma~\ref{815}, we infer that
$\mathrm{supp\,}(f)\subset B_\rho(x_1,R)$
for some $x_1\in X$ and $R\in(0,\infty)$.
Let $\{\mathcal{S}_t\}_{t\in(0,1)}$ be an approximation
of the identity as in Definition~\ref{DefATI}.
Then, applying Theorem~\ref{corATI}(vi) and an argument
similar to that used in the proof of Lemma~\ref{lemma81}(ii),
it is easy to see that, for any $t\in(0,1)$,
$\mathcal{S}_tf\in H^{p,q,\infty}_{\mathrm{at,\,fin}}(X)
\cap\mathrm{UC}(X)$.
Observing that $f\in L^2(X)$, from Theorem~\ref{corATI}(ix),
we deduce that
\begin{align*}
\lim_{t\to0^+}\left\|\mathcal{S}_tf-f\right\|_{L^2(X)}=0,
\end{align*}
which further implies that there exists
a sequence $\{t_k\}_{k\in\mathbb{N}}$ in $(0,1)$
satisfying $\lim_{k\to\infty}t_k=0$
such that
$\mathcal{S}_{t_k}f\to f$ pointwise $\mu$-almost everywhere on $X$
as $k\to\infty$.

Now, we claim that, for any given $(p,\infty)$-atom $a$
supported in $B_\rho(x_a,r_a)$ with $x_a\in X$ and $r_a\in(0,\infty)$,
\begin{align}\label{20455}
\lim_{t\to0^+}\left\|a-\mathcal{S}_{t}a\right\|_{H^{p,q}_*(X)}=0.
\end{align}
Indeed, by Theorem~\ref{corATI}(ix),
we obtain
\begin{align}\label{2011}
\lim_{t\to0^+}\left\|a-\mathcal{S}_{t}a\right\|_{L^2(X)}=0.
\end{align}
From Definitions~\ref{atom}(iii)
and~\ref{DefATI}(v), Fubini's theorem, Theorem~\ref{corATI}(viii),
and an argument similar to that used in
the proof of \eqref{811-x},
we infer that, for any $t\in(0,1)$,
$$
\int_X\left[a(x)-\mathcal{S}_ta(x)\right]\,d\mu(x)=0
$$
and there exists $C\in(0,\infty)$ such that
$\mathrm{supp\,}(a-\mathcal{S}_{t}a)\subset
B_\rho(x_a,C(r_a+1))$. Thus, for any $t\in(0,1)$
satisfying $\|a-\mathcal{S}_ta\|_{L^2(X)}\neq0$,
$$\frac{[\mu(B_\rho(x_a,C(r_a+1)))]^{\frac{1}{2}-
\frac{1}{p}}}{\|a-\mathcal{S}_ta\|_{L^2(X)}}\left(a-\mathcal{S}_ta\right)$$
is a $(p,2)$-atom supported in
$B_\rho(x_a,C(r_a+1))$. On the other hand,
if $\|a-\mathcal{S}_ta\|_{L^2(X)}=0$,
then $\left\|a-\mathcal{S}_ta\right\|_{H^{p,q}_*(X)}=0$
by Propositions~\ref{distfncttype} and \ref{2250}.
Using these observations, Theorems~\ref{6.15} and~\ref{6.14},
and \eqref{2011},
we conclude that
\begin{align*}
\left\|a-\mathcal{S}_ta\right\|_{H^{p,q}_*(X)}
\lesssim\frac{\|a-\mathcal{S}_ta\|_{L^2(X)}}
{[\mu(B_\rho(x_a,C(r_a+1)))]^{\frac{1}{2}-\frac{1}{p}}}\to0
\end{align*}
as $t\to0^+$, which completes the proof of \eqref{20455}.
From this, we deduce that
$$
\lim_{t\to0^+}\left\|f-\mathcal{S}_{t}f\right\|_{H^{p,q}_*(X)}=0.
$$
By this, the proven conclusion that $G_g$ is
a bounded linear functional on $H^{p,q}_*(X)$,
\eqref{22022}, the fact that
$|(\mathcal{S}_tf)g|\leq\|f\|_{L^\infty(X)}
\mathbf{1}_{B_\rho(x_1,C(R+1))}|g|\in L^1(X)$
[since $g\in L^r_{\mathrm{loc}}(X)$],
the proven conclusion that
$\mathcal{S}_{t_k}f\to f$ $\mu$ almost everywhere
as $k\to\infty$,
and the Lebesgue dominated convergence theorem,
we find that
\begin{align*}
G_g(f)
=\lim_{k\to\infty}G_g(\mathcal{S}_{t_k}f)
=\lim_{k\to\infty}\int_X\mathcal{S}_{t_k}f(x)
g(x)\,d\mu(x)=\int_Xf(x)g(x)\,d\mu(x),
\end{align*}
which completes the proof of \eqref{2041}
and hence (i) in this case.
This finishes the proof of (i).

Next, we show (ii).
Let $G\in(H^{p,q}_*(X))'=(H^{p,q,r}_{\mathrm{at}}(X))'$
(see Theorem~\ref{6.14} for this equality).
For any $\rho$-ball $B\subset X$,
let
$$
L^r(B):=\left\{f\in L^r(X):f(x)=
0\text{ for any }x\in B^\complement\right\}
$$
and
$$
\mathring{L}^{r}(B)
:=\left\{f\in L^{r}(B):\int_{X}f(x)\,d\mu(x)=0\right\}.
$$
Then it is easy to verify that,
for any $h\in\mathring{L}^{r}(B)$
with $\|h\|_{L^r(X)}\neq0$,
$\frac{[\mu(B)]^{\frac{1}{r}-\frac{1}{p}}}{\|h\|_{L^r(X)}}h$
is a $(p,r)$-atom,
which, together with $G\in(H^{p,q}_{\mathrm{at}}(X))'$,
further implies that
\begin{align*}
\|G\|_{(\mathring{L}^{r}(B))'}
&=\sup_{\|h\|_{\mathring{L}^{r}(B)=1}}G(h)
=\left[\mu(B)\right]^{\frac{1}{p}-\frac{1}{r}}
\sup_{\|h\|_{\mathring{L}^{r}(B)=1}}G
\left([\mu(B)]^{\frac{1}{r}-\frac{1}{p}}h\right)\\
&\leq\left[\mu(B)\right]^{\frac{1}{p}-\frac{1}{r}}
\|G\|_{(H^{p,q}_{\mathrm{at}}(X))'}.
\end{align*}
Thus, $G$ is a bounded linear functional on
$\mathring{L}^{r}(B)$.
This, combined with the Hahn--Banach theorem,
implies that
$G$ can be extended to a bounded linear functional on $L^{r}(B)$
without increasing its norm.
Now, we consider the following two cases for $r$.

\emph{Case (1)} $r\in(1,\infty)$.
In this case, since $[L^r(B)]'=L^{r'}(B)$ by
Riesz's representation theorem,
it follows that there exists a function $g_B\in L^{r'}(B)$
such that, for any $h\in\mathring{L}^{r}(B)$,
\begin{align}\label{19401}
G(h)=\int_Xg_B(x)h(x)\,d\mu(x)
\end{align}
and
\begin{align}\label{19402}
\left\|g_B\right\|_{L^{r'}(B)}
\leq\left\|G\right\|_{(\mathring{L}^{r}(B))'}
\leq\left[\mu(B)\right]^{\frac{1}{p}-\frac{1}{r}}
\left\|G\right\|_{(H^{p,q}_{\mathrm{at}}(X))'}.
\end{align}

\emph{Case (2)} $r=\infty$.
In this case, by Theorem~\ref{6.14},
we conclude that
\begin{align*}
G\in(H^{p,q,\infty}_{\mathrm{at}}(X))'
=(H^{p,q,\widetilde{r}}_{\mathrm{at}}(X))',
\end{align*}
where $\widetilde{r}\in(1,\infty)$.
Thus, there exists a function
$g_B\in L^{\widetilde{r}}(B)\subset L^1(B)$ such that,
for any $h\in\mathring{L}^{\infty}(B)$,
$$
G(h)=\int_Xg_B(x)h(x)\,d\mu(x)
$$
and
$$
\left\|g_B\right\|_{L^{1}(B)}
\leq\left\|g_B\right\|_{L^{\widetilde{r}'}(B)}
\left[\mu(B)\right]^\frac{1}{\widetilde{r}}
\leq\left[\mu(B)\right]^{\frac{1}{p}}
\left\|G\right\|_{(H^{p,q}_{\mathrm{at}}(X))'}.
$$
Thus, we have proven that, for any $r\in(1,\infty]$,
\eqref{19401} and \eqref{19402} hold
for some $g_B\in L^{r'}(B)$.

Next, we claim that, if $\widetilde{g}_B\in L^{r'}(B)$ is
such that, for any $h\in\mathring{L}^{r}(B)$,
\begin{align}\label{kbr-58}
G(h)=\int_X\widetilde{g}_B(x)h(x)\,d\mu(x),
\end{align}
then $g_B-\widetilde{g}_B=m_B(g_B-\widetilde{g}_B)$
$\mu$-almost everywhere.
Indeed, \eqref{19401} and \eqref{kbr-58} imply
that, for any $h\in\mathring{L}^{r}(B)$,
$$
\int_X\left[g_B(x)-\widetilde{g}_B(x)\right]
h(x)\,d\mu(x)=G(h)-G(h)=0,
$$
which, together with the fact that
$u-m_B(u)\in\mathring{L}^{r}(B)$
for any $u\in L^r(B)$,
further implies that, for any $u\in L^r(B)$,
\begin{align*}
\int_X\left[g_B(x)-\widetilde{g}_B(x)\right]
\left[u(x)-m_B(u)\right]\,d\mu(x)=0,
\end{align*}
equivalently,
\begin{align*}
\int_X\left[g_B(x)-\widetilde{g}_B(x)\right]
u(x)\,d\mu(x)=m_B(u)\int_X\left[g_B(x)-\widetilde{g}_B(x)\right]\,d\mu(x).
\end{align*}
From this, we infer that,
for any $u\in L^r(B)$,
\begin{align*}
&\int_Xu(x)\left[g_B(x)-\widetilde{g}_B(x)
-m_B(g_B-\widetilde{g}_B)\right]\,d\mu(x)\\
&\quad=\int_X\left[g_B(x)-\widetilde{g}_B(x)\right]
u(x)\,d\mu(x)-m_B(g_B-\widetilde{g}_B)\int_Xu(x)\,d\mu(x)\\
&\quad=m_B(u)\int_X\left[g_B(x)-\widetilde{g}_B(x)
\right]\,d\mu(x)-m_B(g_B-\widetilde{g}_B)\int_Xu(x)\,d\mu(x)=0,
\end{align*}
which, combined with the arbitrariness
of $u\in L^r(B)$ and an argument similar to the
one used in Proposition~\ref{distfncttype}, further implies that
$g_B-\widetilde{g}_B=m_B(g_B-\widetilde{g}_B)$
$\mu$-almost everywhere. This finishes the proof
of the above claim.

Fix $x_0\in X$ and, for any $N\in\mathbb{N}$,
let $B_N:=B_\rho(x_0,N)$. Then, by
what have just proven, we find that there exists
a sequence
$\{g_N\}_{N\in\mathbb{N}}$ of $\mu$-measurable
functions such that, for any $N\in\mathbb{N}$,
$g_N\in L^{r'}(B_N)$ and, for any $h\in \mathring{L^r}(B_N)$,
\begin{align*}
G(h)=\int_Xg_N(x)h(x)\,d\mu(x).
\end{align*}
Let $\widetilde{g}_1:=g_1$ and,
for any $N\in\mathbb{N}\cap[2,\infty)$,
$$
\widetilde{g}_N:=g_N+m_{B_N}\left(\widetilde{g}_{N-1}-g_N\right).
$$
From this and the above claim,
we deduce that, for any given $N\in\mathbb{N}$
and for $\mu$-almost every $x\in B_N$,
$\widetilde{g}_{N+1}(x)=\widetilde{g}_N(x)$
and $\widetilde{g}_N\in L^{r'}(B_N)$.
Let $f\in H^{p,q,r}_{\mathrm{at,\,fin}}(X)$
and $g$ be a function such that, for any $N\in\mathbb{N}$,
$g(x)=\widetilde{g}_N(x)$ $\mu$-almost everywhere on $B_N$.
Then, by Lemma~\ref{815}, we can choose $N_{(f)}\in\mathbb{N}$
large enough such that
$f\in\mathring{L^r}(B_{N_{(f)}})$ and hence
\begin{align*}
G(f)=\int_Xg_{N_{(f)}}(x)f(x)\,d\mu(x)
=\int_X\widetilde{g}_{N_{(f)}}(x)f(x)\,d\mu(x)
=\int_Xg(x)f(x)\,d\mu(x).
\end{align*}

Finally, we show that $g\in\mathcal{C}^{p,q,r'}(X)$.
Indeed, let $m,v\in\mathbb{N}$,
$\{B_{i,j}\}_{i\in\mathbb{N}\cap[1,m],j\in\mathbb{N}\cap[1,v]}$
be a sequence of $\rho$-balls,
and $\{\lambda_{i,j}\}_{i\in\mathbb{N}
\cap[1,m],j\in\mathbb{N}\cap[1,v]}$ in $\mathbb{R}_+$
satisfying all the conditions in Definition~\ref{20191}.
For any $i\in\mathbb{N}\cap[1,m]$ and $j\in\mathbb{N}\cap[1,v]$,
let $h_{i,j}\in L^r(B_{i,j})$ with $\|h_{i,j}\|_{L^r(B_{i,j})}=1$
such that
\begin{align}\label{21501}
\left\|g-m_{B_{i,j}}(g)\right\|_{L^{r'}(B_{i,j})}
=\int_{B_{i,j}}\left[g(x)-m_{B_{i,j}}(g)\right]h_{i,j}(x)\,d\mu(x).
\end{align}
In addition, for any $x\in X$,
if $h_{i,j}$ is a constant function,
then let $a_{i,j}(x):=0$,
and, if $h_{i,j}$ is not a constant function,
then let
\begin{align*}
a_{i,j}(x):=
\frac{[\mu(B_{i,j})]^{\frac{1}{r}-\frac{1}{p}}
[h_{i,j}(x)-m_{B_{i,j}}(h_{i,j})(x)]\mathbf{1}_{B_{i,j}}(x)}
{\|h_{i,j}-m_{B_{i,j}}(h_{i,j})\|_{L^r(B_{i,j})}}.
\end{align*}
Then it is easy to verify that,
for any $i\in\mathbb{N}\cap[1,m]$ and $j\in\mathbb{N}\cap[1,v]$,
$a_{i,j}$ is a $(p,r)$-atom,
which, together with Theorem~\ref{6.15},
further implies that
\begin{align}\label{1439}
\sum_{i=1}^m\sum_{j=1}^v\lambda_{i,j}a_{i,j}\in H^{p,q}_*(X)
\end{align}
and
\begin{align}\label{1440}
\left\|\sum_{i=1}^m
\sum_{j=1}^v\lambda_{i,j}a_{i,j}\right\|_{H^{p,q}_*(X)}
\lesssim\left[\sum_{i=1}^m\left(\sum_{j=1}^v
\left|\lambda_{i,j}\right|^p\right)^\frac{q}{p}\right]^\frac{1}{q}.
\end{align}
From \eqref{21501},
we infer that
\begin{align*}
&\sum_{i=1}^m\sum_{j=1}^v
\left|\lambda_{i,j}\right|
\left[\mu(B_{i,j})\right]^{1-\frac{1}{p}}
\left[\fint_{B_{i,j}}
\left|g(x)-m_{B_{i,j}}(g)\right|^{r'}\,d\mu(x)\right]^\frac{1}{r'}\\
&\quad=\sum_{i=1}^m\sum_{j=1}^v
\left|\lambda_{i,j}\right|
\left[\mu(B_{i,j})\right]^{\frac{1}{r}-\frac{1}{p}}
\int_{B_{i,j}}\left[g(x)-m_{B_{i,j}}(g)\right]h_{i,j}(x)\,d\mu(x)
\\
&\quad=\sum_{i=1}^m\sum_{j=1}^v
\left|\lambda_{i,j}\right|
\left[\mu(B_{i,j})\right]^{\frac{1}{r}-\frac{1}{p}}
\int_{B_{i,j}}\left[h_{i,j}(x)-m_{B_{i,j}}(h_{i,j})\right]g(x)
\,d\mu(x),
\end{align*}
which, combined with the facts that
$\|h_{i,j}-m_{B_{i,j}}(h_{i,j})\|_{L^r(B_{i,j})}\leq2$
and $G\in(H^{p,q}_*(X))'$,
\eqref{1439}, and \eqref{1440},
further implies that
\begin{align*}
&\sum_{i=1}^m\sum_{j=1}^v
\left|\lambda_{i,j}\right|
\left[\mu(B_{i,j})\right]^{1-\frac{1}{p}}
\left[\fint_{B_{i,j}}
\left|g(x)-m_{B_{i,j}}(g)\right|^{r'}\,d\mu(x)\right]^\frac{1}{r'}\\
&\quad\lesssim\sum_{i=1}^m\sum_{j=1}^v\lambda_{i,j}
\int_{B_{i,j}}a_{i,j}(x)g(x)\,d\mu(x)
=G\left(\sum_{i=1}^m\sum_{j=1}^v\lambda_{i,j}a_{i,j}\right)\\
&\quad\leq\|G\|_{(H^{p,q}_*(X))'}\left\|\sum_{i=1}^m
\sum_{j=1}^v\lambda_{i,j}a_{i,j}\right\|_{H^{p,q}_*(X)}
\lesssim\|G\|_{(H^{p,q}_*(X))'}\left[\sum_{i=1}^m\left(\sum_{j=1}^v
\left|\lambda_{i,j}\right|^p\right)^\frac{q}{p}\right]^\frac{1}{q}.
\end{align*}
Taking the supremum over all
$m,v\in\mathbb{N}$, $\{B_{i,j}\}_{i\in\mathbb{N}
\cap[1,m],j\in\mathbb{N}\cap[1,v]}$,
and $\{\lambda_{i,j}\}_{i,j\in\mathbb{N}}\in\ell^q(\ell^p)$
as above,
we find that $g\in\mathcal{C}^{p,q,r'}(X)$.

Finally, we prove the uniqueness of $g$.
Assume that $h\in\mathcal{C}^{p,q,r'}(X)$ is such that
$G=G_h$ in $(H^{p,q}_*(X))'$.
Then, for any atom $a$,
\begin{align*}
0=\left\langle G,a\right\rangle-\left\langle G,a\right\rangle
=\int_X\left[g(x)-h(x)\right]a(x)\,d\mu(x).
\end{align*}
By this and the arbitrariness of atoms,
we conclude that $g-h\equiv C$ $\mu$-almost everywhere,
where $C$ is a constant.
This finishes the proof of (ii) and hence
Theorem~\ref{duality}.
\end{proof}

\section{Boundedness of Calder\'on--Zygmund Operators}
\label{section5.8}

In this section, we aim to show the boundedness on $H^{p,q}_\ast(X)$ of
Calder\'on--Zygmund operators defined in Definition~\ref{CZO}
for optimal ranges of $p$ and $q$
on general quasi-ultrametric spaces of homogeneous type.
To this end, we first establish the boundedness of
Calder\'on--Zygmund operators
from $H^p_\ast(X)$ to $L^p(X)$ and
from $H^p_\ast(X)$ to $H^p_\ast(X)$.
Recall that, for any $p\in(0,\infty)$,
$L^{p,p}(X)=L^p(X)$,
which further implies
$H^p_\ast(X)=H^{p,p}_\ast(X)$.
In what follows, for any $p\in(0,\infty)$ and
$r\in[1,\infty)$, let
$$
H^{p,r}_{\mathrm{at}}(X):=H^{p,p,r}_{\mathrm{at}}(X)
\ \ \text{and}\ \
H^{p,r}_{\mathrm{at},\,\mathrm{fin}}(X)
:=H^{p,p,r}_{\mathrm{at},\,\mathrm{fin}}(X).
$$

We first present the molecular characterization of
$H^p_\ast(X)$.
To this end, we introduce the following
notion of molecules, which is a slight variant of
Hu et al. \cite[Definition~1.2]{hyz2009}.

\begin{definition}
\label{m109}
Let $(X,\mathbf{q},\mu)$ be a quasi-ultrametric space
of homogeneous type. Assume that
$\rho\in\mathbf{q}$ has the property
that all $\rho$-balls are $\mu$-measurable,
$p\in(0,1]$, $r\in(p,\infty]\cap[1,\infty]$,
and $\overrightarrow{\eta}:
=\{\eta_k\}_{k\in\mathbb{N}}$ in $\mathbb{R}_+$
satisfies
\begin{align}
\label{1091}
\sum_{k=1}^\infty k(\eta_k)^p<\infty.
\end{align}
A function $m\in L^r(X)$ is called a
\emph{$(\rho,p,r,\overrightarrow{\eta})$-molecule}
\index[words]{molecule}
centered at a ball $B:=B_\rho(x,r)$
with $x\in X$ and $r\in(0,\infty)$ if
\begin{enumerate}
\item[\textup{(i)}]
$\|m\mathbf{1}_B\|_{L^r(X)}
\leq[\mu(B)]^{\frac{1}{r}-\frac{1}{p}}$,

\item[\textup{(ii)}]
for any $k\in\mathbb{N}$,
\begin{align*}
\left\|m\mathbf{1}_{2^kB\setminus2^{k-1}B}\right\|_{L^r(X)}
\leq\eta_k\left[\mu(2^kB)\right]^{\frac{1}{r}-\frac{1}{p}},
\end{align*}

\item[\textup{(iii)}]
$\int_Xm(x)\,d\mu(x)=0$.
\end{enumerate}
Moreover, when $\mu(X)<\infty$,
it is agreed upon that the constant function $m:=[\mu(X)]^{-\frac{1}{p}}$
is a $(\rho,p,r,\overrightarrow{\eta})$-molecule.
\end{definition}

\begin{remark}
\begin{enumerate}
\item[\rm(i)]
Let $\eta\in(0,\infty)$
and $\vec{\eta}:=\{2^{-k\eta}\}_{k\in\mathbb{N}}$.
Then the $(\rho,p,r,\overrightarrow{\eta})$-molecule in Definition~\ref{m109}
reduces to the $(\rho,p,r,\eta)$-molecule in Definition~\ref{mol}.

\item[\rm(ii)]
It is easy to see from Definition~\ref{m109},
\eqref{1557-xx}, and \eqref{doublingcond}
that, if $m: X\to\mathbb{C}$ is a $(\rho,p,r,\vec{\eta})$-molecule
centered at $B_{\rho}(x,r)$ for some $x\in X$
and $r\in(0,\infty)$, then, for any other $\rho'\in\mathbf{q}$
with the property that all $\rho'$-balls are $\mu$-measurable,
there exist positive constants $C$ and $c$ (depending only on $\rho$,
$\rho'$, $\mu$, and $\vec{\eta}$) such that $Cm$ is a
$(\rho',p,r,\vec{\eta})$-molecule centered at $B_{\rho'}(x,cr)$.
As such, we sometimes simply refer to $m$ as a
$(p,r,\vec{\eta})$-\textit{molecule} or a
$(p,r,\vec{\eta})$-\textit{molecule modulo a harmless constant multiple}.

\item[\rm(iii)]
Conditions (i) and (ii) in Definition~\ref{m109}
and \eqref{1091} ensure that any
$(\rho,p,r,\vec{\eta})$-molecule is in $L^1(X)$, and hence the integral in
Definition~\ref{m109}(iii) is well defined.
Moreover, we also find that
any $(\rho,p,r,\eta)$-molecule
is normalized in the way that its $L^p(X)$ norm does not exceed
$(1+\sum_{j=1}^\infty\eta_k^p)^\frac{1}{p}$.
\end{enumerate}
\end{remark}

Based on the above molecules, by the atomic
characterization of $H^p_\ast(X)$
(see Theorem~\ref{6.14} with $p=q$)
and Lemma~\ref{L6.2},
we obtain the following molecular characterization of
$H^p(X)$,
which is different from Theorem~\ref{thmHpqm} with $p=q$;
see also
\cite[Theorem~3.3]{LDC2} for the case $p\in(0,1]$
and \cite[Theorem~3.2]{LDC1} for the case $p=1$.

\begin{theorem}\label{109}
Let $(X,\mathbf{q},\mu)$ be a quasi-ultrametric space
of homogeneous type with upper dimension $n\in(0,\infty)$.
Assume that
$p\in(\frac{n}{n+\mathrm{ind\,}(X,\mathbf{q})},1]$,
$r\in(p,\infty]\cap[1,\infty]$,
$\overrightarrow{\eta}:=\{\eta_k\}_{k\in\mathbb{N}}$ in $\mathbb{R}_+$
satisfies \eqref{1091}, and $\beta,\gamma,\varepsilon\in(0,\infty)$ satisfy
$n(\frac{1}{p}-1)<\beta,\gamma<\varepsilon\preceq\mathrm{ind\,}(X,\mathbf{q})$.
Then $f\in H^p_\ast(X)$ if and only if there exist
a sequence of $(p,r,\overrightarrow{\eta})$-molecules
$\{m_j\}_{j\in\mathbb{N}}$
and $\{\lambda_{j}\}_{j\in\mathbb{N}}$ in $\mathbb{R}_+$
such that
$f=\sum_{j=1}^\infty\lambda_j m_j$
in $({\mathcal{G}}^\varepsilon_0(\beta,\gamma))'$.
Moreover, there exists a constant $C\in[1,\infty)$ such that,
for any $f\in H^p_\ast(X)$,
\begin{align*}
C^{-1}\|f\|_{H^p_\ast(X)}
\leq\inf\left\{\left(\sum_{j=1}^\infty
\lambda_j^p\right)^\frac{1}{p}:
f=\sum_{j=1}^\infty\lambda_j m_j\right\}
\leq C\|f\|_{H^p_\ast(X)},
\end{align*}
where the infimum is taken over all
decompositions of $f$ as above.
\end{theorem}

Granted this molecular characterization of
$H^p_\ast(X)$,
we are now ready to prove the following conclusion.

\begin{theorem}
\label{CZOextend}
Let $(X,\mathbf{q},\mu)$ be a quasi-ultrametric space
of homogeneous type with upper dimension $n\in(0,\infty)$,
$0<\varepsilon\preceq\mathrm{ind\,}(X,\mathbf{q})$,
$p\in(\frac{n}{n+\varepsilon},1]$,
and $T$ be an $\varepsilon$-Calder\'on--Zygmund
operator as in Definition~\ref{CZO}.
\begin{enumerate}
\item[\textup{(i)}]
Then $T$ can be extended to a unique bounded
linear operator from $H^p_\ast(X)$ to $L^p(X)$.
\item[\textup{(ii)}]
If $T$ satisfies that $T^*1=0$
in the sense that
$$\int_XTf(x)\,d\mu(x)=0$$
for any $f\in L^2(X)$ with bounded support and
$\int_Xf(x)\,d\mu(x)=0$,
then $T$ can be extended to a unique bounded
linear operator from $H^p_\ast(X)$ to $H^p_\ast(X)$.
\end{enumerate}
\end{theorem}

\begin{proof}
Let $\rho\in\mathbf{q}$ be a regularized quasi-ultrametric
as in Definition~\ref{D2.3}.
To show (i),
we first claim that there exists a positive constant $C$
such that, for any $(p,2)$-atom $a\in L^2(X)$,
\begin{align}\label{814-}
\left\|Ta\right\|_{L^p(X)}\leq C.
\end{align}
To prove this claim, assume that $a\in L^2(X)$ is a $(p,2)$-atom
supported in $B:=B_\rho(x_B,r_B)$ with $x_B\in X$
and $r_B\in(0,\infty)$.
Observe that
\begin{align}
\label{8140}
\left\|Ta\right\|_{L^p(X)}^p
&=\int_{2C_\rho B}\left|Ta(x)\right|^p\,d\mu(x)
\nonumber\\
&\quad+\sum_{k=1}^\infty
\int_{(2C_\rho)^{k+1}B\setminus(2C_\rho)^kB}
\left|\int_BK(x,y)a(y)\,d\mu(y)\right|^p\,d\mu(x)
\nonumber\\
&=:\mathrm{I_1}+\mathrm{I_2}.
\end{align}
To estimate $\mathrm{I_1}$, using H\"older's inequality,
\eqref{upperdoub},
the boundedness of $T$ on $L^2(X)$,
and Definition~\ref{atom}(ii),
we find that
\begin{align}
\label{TaC-I1}
\mathrm{I_1}\leq\left[\mu(2C_\rho B)\right]^\frac{2-p}{2}
\left\|Ta\right\|_{L^2(X)}^p
\lesssim\left[\mu(B)\right]^\frac{2-p}{2}
\left\|a\right\|_{L^2(X)}^p
\lesssim1.
\end{align}
To estimate $\mathrm{I_2}$, note that,
for any $x\in(2C_\rho B)^\complement$
and $y\in B$,
\begin{align}\label{1627}
\rho(y,x_B)<\frac{\rho(x,x_B)}{2C_\rho}.
\end{align}
From this, we deduce that
\begin{align*}
\mathrm{I_2}
&\leq\sum_{k=1}^\infty
\int_{(2C_\rho)^{k+1}\setminus(2C_\rho)^kB}
\left[\int_B\left|K(x,y)-K(x,x_B)\right|
|a(y)|\,d\mu(y)\right]^p\,d\mu(x)
\\
&\qquad\ \text{by Definition~\ref{atom}(iii)}
\\
&\lesssim\sum_{k=1}^\infty
\int_{(2C_\rho)^{k+1}\setminus(2C_\rho)^kB}
\left\{\int_B\left[\frac{\rho(x_B,y)}
{\rho(x,x_B)}\right]^\varepsilon
\frac{1}{V_\rho(x,x_B)}|a(y)|
\,d\mu(y)\right\}^p\,d\mu(x)
\\
&\qquad\ \text{by \eqref{CZk2} and \eqref{1627}}\\
&\leq\sum_{k=1}^\infty
\int_{(2C_\rho)^{k+1}\setminus(2C_\rho)^kB}
\|a\|_{L^2(X)}^p
\left\{\int_B\left[\frac{\rho(x_B,y)}
{\rho(x,x_B)}\right]^{2\varepsilon}
\frac{1}{[V_\rho(x,x_B)]^2}
\,d\mu(y)\right\}^\frac{p}{2}\,d\mu(x)\\
&\qquad\ \text{by H\"older's inequality}.
\end{align*}
By this,
we further conclude that
\begin{align*}
\mathrm{I_2}
&\lesssim\sum_{k=1}^\infty
\int_{(2C_\rho)^{k+1}\setminus(2C_\rho)^kB}
[\mu(B)]^{\frac{p}{2}-1}
\left[\frac{r_B}{\rho(x,x_B)}\right]^{p\varepsilon}
\frac{1}{[V_\rho(x,x_B)]^p}[\mu(B)]^\frac{p}{2}\,d\mu(x)
\\
&\qquad\ \text{by Definition~\ref{atom}(ii)}\\
&\leq\sum_{k=1}^\infty
[\mu(B)]^{p-1}
(2C_\rho)^{-kp\varepsilon}
\frac{\mu((2C_\rho)^{k+1}B)}{[\mu((2C_\rho)^kB)]^p}
\\
&\lesssim\sum_{k=1}^\infty
(2C_\rho)^{k[n(1-p)-p\varepsilon]}
\quad\text{by \eqref{upperdoub}}\\
&\sim 1
\quad\text{since $p\in\left(\frac{n}{n+\varepsilon},1\right]$}.
\end{align*}
From this, \eqref{8140}, and \eqref{TaC-I1},
it follows that the claim in \eqref{814-} holds.

Next, we use \eqref{814-} to show that,
for any $f\in H^{p,2}_{\mathrm{at},\,\mathrm{fin}}(X)$,
$$
\left\|Tf\right\|_{L^p(X)}\lesssim\|f\|_{H^p_\ast(X)}.
$$
Let $f\in H^{p,2}_{\mathrm{at},\,\mathrm{fin}}(X)$.
Then there exist $N\in\mathbb{N}$,
a sequence
$\{\lambda_i\}_{i=1}^N$ in $\mathbb{R}_+$,
and a sequence $\{a_i\}_{i=1}^N$
of $(p,2)$-atoms such that
\begin{align}\label{yy1}
f=\sum_{i=1}^N\lambda_ia_i
\end{align}
in $({\mathcal{G}}_0^\varepsilon(\beta,\gamma))'$
and
$$
\|f\|_{H^{p,2}_{\mathrm{at},\,\mathrm{fin}}(X)}
\sim\left(\sum_{i=1}^N\lambda_i^p\right)^\frac{1}{p}.
$$
By this and Lemma~\ref{6.29}(i), we find that
\begin{align}\label{yy2}
\|f\|_{H^p_\ast(X)}\sim\|f\|_{H^{p,2}_{\mathrm{at},\,\mathrm{fin}}(X)}
\sim\left(\sum_{i=1}^N\lambda_i^p\right)^\frac{1}{p}.
\end{align}
From this, $p\in(\frac{n}{n+\varepsilon},1]$,
and \eqref{814-}, we infer that
\begin{align*}
\left\|Tf\right\|_{L^p(X)}^p
\leq\sum_{i=1}^N\lambda_i^p
\left\|Ta_i\right\|_{L^p(X)}^p
\lesssim\sum_{i=1}^N\lambda_i^p\sim\|f\|_{H^p_\ast(X)}^p.
\end{align*}
Thus, $T$ is a bounded operator from
$H^{p,2}_{\mathrm{at},\,\mathrm{fin}}(X)$
[equipped with $\|\cdot\|_{H^p_\ast(X)}$] to $L^p(X)$.
By this, the density of $H^{p,2}_{\mathrm{at},\,\mathrm{fin}}(X)$
in $H^p_\ast(X)$ (see Theorem~\ref{5r1} with $q=p$),
and a standard density argument,
we conclude that $T$ can be extended to a unique bounded
linear operator from $H^p_\ast(X)$ to $L^p(X)$.
This finishes the proof of (i).

Now, we prove (ii).
Assume that $T^*1=0$. To show that
$T$ can be extended to a unique bounded
linear operator from $H^p_\ast(X)$ to $H^p_\ast(X)$,
we first prove that, for any $(p,2)$-atom $a$ supported
in $B:=B_\rho(x_B,r_B)$ with $x_B\in X$ and $r_B\in(0,\infty)$,
$Ta$ is a harmless constant multiple of a
$(p,2,\overrightarrow{\eta})$-molecule as
in Definition~\ref{m109},
where $\overrightarrow{\eta}
:=\{\eta_k\}_{k\in\mathbb{N}}$ in $\mathbb{R}_+$
will be specified later.
From the assumption that $T^*1=0$
and the fact that $a$ is a $(p,2)$-atom,
we deduce that
\begin{align}
\label{814m1}
\int_XTa(x)\,d\mu(x)=0.
\end{align}
Let $\widetilde{B}:=B_\rho(x_B,2C_\rho r_B)$.
Then, by the boundedness of $T$ on $L^2(X)$,
Definition~\ref{atom}(ii), the fact that $p<2$,
and \eqref{upperdoub},
we find that
\begin{align}
\label{814m2}
\left\|Ta\mathbf{1}_{\widetilde{B}}\right\|_{L^2(X)}
\leq\left\|Ta\right\|_{L^2(X)}
\lesssim\left\|a\right\|_{L^2(X)}
\lesssim\left[\mu(B)\right]^{\frac{1}{2}-\frac{1}{p}}
\lesssim\left[\mu\left(\widetilde{B}\right)
\right]^{\frac{1}{2}-\frac{1}{p}}.
\end{align}
Observe that, for any
$x\in(2C_\rho B)^\complement$
and $y\in B$,
$$
\rho(y,x_B)<\frac{\rho(x,x_B)}{2C_\rho}.
$$
From this, Definition~\ref{atom}(iii),
\eqref{CZk2}, and H\"older's inequality,
we infer that, for any $m\in\mathbb{N}$,
\begin{align*}
&\left\|Ta\mathbf{1}_{2^m\widetilde{B}
\setminus2^{(m-1)}\widetilde{B}}\right\|_{L^2(X)}^2
\nonumber\\
&\quad=\int_{2^m\widetilde{B}
\setminus2^{(m-1)}\widetilde{B}}
\left|\int_B\left[K(x,y)-K(x,x_B)\right]
a(y)\,d\mu(y)\right|^2\,d\mu(x)
\nonumber\\
&\quad\lesssim\int_{2^m\widetilde{B}
\setminus2^{(m-1)}\widetilde{B}}
\left\{\int_B\left[\frac{\rho(y,x_B)}
{\rho(x,x_B)}\right]^\varepsilon
\frac{1}{V_\rho(x,x_B)}|a(y)|\,d\mu(y)\right\}^2
\,d\mu(x)
\nonumber\\
&\quad\leq\int_{2^m\widetilde{B}
\setminus2^{(m-1)}\widetilde{B}}
\left[\frac{r_B}{\rho(x,x_B)}\right]^{2\varepsilon}
\frac{1}{[V_\rho(x,x_B)]^2}\mu(B)
\|a\|_{L^2(X)}^2\,d\mu(x),
\end{align*}
which, together with Definition~\ref{atom}(ii)
and \eqref{upperdoub}, further implies that
\begin{align}
\label{814m3}
&\left\|Ta\mathbf{1}_{2^m\widetilde{B}
\setminus2^{(m-1)}\widetilde{B}}\right\|_{L^2(X)}
\nonumber\\
&\quad\lesssim2^{-m\varepsilon}
\left[\mu(2^{m}C_\rho\widetilde{B})\right]^{-1}
\left[\mu(B)\right]^{1-\frac{1}{p}}
\left[\mu(2^{m+1}C_\rho\widetilde{B})\right]^\frac{1}{2}
\nonumber\\
&\quad\lesssim2^{-m[\varepsilon+n(1-\frac{1}{p})]}
\left[\mu(2^m\widetilde{B})\right]^{\frac{1}{2}-\frac{1}{p}}.
\end{align}
For any $m\in\mathbb{N}$,
let $\eta_m:=2^{-m[\varepsilon+n(1-\frac{1}{p})]}$.
By $p\in(\frac{n}{n+\varepsilon},1]$,
we conclude that $\overrightarrow{\eta}
:=\{\eta_k\}_{k\in\mathbb{N}}$ satisfies \eqref{1091}.
From this, \eqref{814m1}, \eqref{814m2},
and \eqref{814m3}, we deduce that
$Ta$ is a harmless constant multiple of a
$(p,2,\overrightarrow{\eta})$-molecule as in Definition~\ref{m109},
which further implies that, for any
$f\in H^{p,2}_{\mathrm{at},\,\mathrm{fin}}(X)$ satisfying
both \eqref{yy1} and \eqref{yy2},
\begin{align*}
Tf=\sum_{i=1}^N\lambda_i Ta_i
=:\sum_{i=1}^N\lambda_i' m_i,
\end{align*}
where, for any $i\in\mathbb{N}\cap[1,N]$,
$\lambda_i'\sim\lambda_i$ and
$m_i$ is a $(p,2,\overrightarrow{\eta})$-molecule,
By this and Theorem~\ref{109},
we find that,
for any $f\in H^{p,2}_{\mathrm{at},\,\mathrm{fin}}(X)$,
\begin{align*}
\left\|Tf\right\|_{H^p_\ast(X)}
&\sim\inf\left\{\left(\sum_{i=1}^\infty\lambda_i^p\right)^\frac{1}{p}
:Tf=\sum_{i=1}^\infty\lambda_i m_i\right\}\\
&\lesssim\left[\sum_{i=1}^N(\lambda_i')^p\right]^\frac{1}{p}
\sim\left(\sum_{i=1}^N\lambda_i^p\right)^\frac{1}{p}
\sim\|f\|_{H^p_\ast(X)}.
\end{align*}
From this, the density of $H^{p,2}_{\mathrm{at},\,\mathrm{fin}}(X)$
in $H^p_\ast(X)$ (see Theorem~\ref{5r1} with $q=p$)
and a standard density argument,
we infer that $T$ can be extended to a unique bounded
linear operator from $H^p_\ast(X)$ to $H^p_\ast(X)$,
which completes the proof of (ii)
and hence Theorem~\ref{CZOextend}.
\end{proof}

By Theorem~\ref{CZOextend} and a classical interpolation argument,
we obtain the following conclusion.

\begin{theorem}
\label{CZOBHL}
Let $(X,\mathbf{q},\mu)$ be a quasi-ultrametric space
of homogeneous type with upper dimension $n\in(0,\infty)$.
Assume that
$0<\varepsilon\preceq\mathrm{ind\,}(X,\mathbf{q})$,
$p\in(\frac{n}{n+\varepsilon},1]$,
$q\in(0,\infty]$,
and $T$ is an $\varepsilon$-Calder\'on--Zygmund
operator as in Definition~\ref{CZO}.
\begin{enumerate}
\item[\textup{(i)}]
Then $T$ can be extended to a unique bounded linear operator
from $H^{p,q}_\ast(X)$ to $L^{p,q}(X)$.

\item[\textup{(ii)}]
Assume further that $T$ satisfies $T^*1=0$;
that is, for any $f\in L^2(X)$ with both bounded support and
$\int_Xf(x)\,d\mu(x)=0$,
$$
\int_XTf(x)\,d\mu(x)=0.
$$
Then $T$ can be extended to a unique bounded linear operator from
$H^{p,q}_\ast(X)$ to $H^{p,q}_\ast(X)$.
\end{enumerate}
\end{theorem}

\begin{proof}
We first show (i).
Let $\theta\in(0,1)$ and
$p_1,p_2\in(\frac{n}{n+\varepsilon},1]$ be
such that $p_1\neq p_2$ and
$\frac{1}{p}=\frac{1-\theta}{p_1}+\frac{\theta}{p_2}$.
By Theorem~\ref{CZOextend}, we conclude
that $T$ can be extended to a unique bounded linear operator
from $H^{p_1}_\ast(X)$ to $L^{p_1}(X)$ and from $H^{p_2}_\ast(X)$
to $L^{p_2}(X)$. Combining this,
the definition of the interpolation space,
Theorem~\ref{chazhi-x}, and Lemma~\ref{LczLpq-x},
we find that $T$ is bounded from
$(H^{p_1}_\ast(X),H^{p_2}_\ast(X))_{\theta,q}$
to $(L^{p_1}(X),L^{p_2}(X))_{\theta,q}$
and hence from $H^{p,q}_\ast(X)$ to $L^{p,q}(X)$.
This finishes the proof of (i).

Next, we prove (ii).
Assume that $T$ satisfies the additional property that
$T^*1=0$. From Theorem~\ref{CZOextend} and an argument similar to
that used in the proof of (i) with
$L^p(X)$ and $L^{p,q}(X)$ replaced, respectively,
by $H^p_\ast(X)$ and $H^{p,q}_\ast(X)$,
we deduce that $T$ is bounded from
$(H^{p_1}_\ast(X),H^{p_2}_\ast(X))_{\theta,q}$
to itself
and hence from $H^{p,q}_\ast(X)$ to itself.
This finishes the proof of (ii)
and hence Theorem~\ref{CZOBHL}.
\end{proof}

In the critical case $p=\frac{n}{n+\varepsilon}$,
we have the following conclusion.

\begin{theorem}
\label{CZOcritical}
Let $(X,\mathbf{q},\mu)$ be a quasi-ultrametric space
of homogeneous type with upper dimension $n\in(0,\infty)$,
$\varepsilon\in(0,\mathrm{ind\,}(X,\mathbf{q}))$,
and $T$ be an $\varepsilon$-Calder\'on--Zygmund
operator as in Definition~\ref{CZO}.
\begin{enumerate}
\item[\textup{(i)}]
Then $T$ can be extended to a unique bounded linear operator
from $H^{\frac{n}{n+\varepsilon}}_\ast(X)$
to $L^{\frac{n}{n+\varepsilon},\infty}(X)$.

\item[\textup{(ii)}]
If $T$ satisfies that $T^*1=0$,
then $T$ can be extended to a unique bounded linear operator
from $H^{\frac{n}{n+\varepsilon}}_\ast(X)$
to $H^{\frac{n}{n+\varepsilon},\infty}_\ast(X)$.
\end{enumerate}
\end{theorem}

To show Theorem~\ref{CZOcritical},
we need some technical lemmas.
The following lemma is a weak-type
summability principle;
see, for instance, \cite[p.\,9, Lemma]{fs86}.

\begin{lemma}
\label{814L1}
Let $(X,\mu)$ be a measure space and
$p\in(0,1)$. Assume that
$\{c_i\}_{i\in\mathbb{N}}$ in $\mathbb{C}$
is a sequence satisfying
$\sum_{i\in\mathbb{N}}|c_i|^p<\infty$
and $\{f_i\}_{i\in\mathbb{N}}$ is a sequence of
measurable functions with the property that there exists
a positive constant $C$ such that, for any $i\in\mathbb{N}$
and $\lambda\in(0,\infty)$,
$$\mu(\{x\in X:
\left|f_i(x)\right|>\lambda\})
\leq C\lambda^{-p}.$$
Then the series $\sum_{i\in\mathbb{N}}c_if_i(x)$
is absolutely convergent for $\mu$-almost every $x\in X$
and there exists a positive constant $C_p$,
independent of $\{f_i\}_{i\in\mathbb{N}}$,
such that, for any $\lambda\in(0,\infty)$,
\begin{align*}
\mu\left(\left\{x\in X:
\left|\sum_{i\in\mathbb{N}}c_if_i(x)\right|
>\lambda\right\}\right)
\leq C_p
\frac{\sum_{i\in\mathbb{N}}|c_i|^p}{\lambda^p}.
\end{align*}
\end{lemma}

The following lemma is a generalization of \cite[Lemma 8.15]{zhy}
by allowing the measure of a singleton
to be strictly positive and $\rho\in\mathbf{q}$
to be a quasi-ultrametric.
Note that, when
$\rho\in\mathbf{q}$
is a regularized quasi-ultrametric,
we have, for any given $x_0\in X$,
the mapping $x\mapsto V_\rho(x_0,x)$ is $\mu$-measurable on $X$.

\begin{lemma}
\label{814L2}
Let $(X,\mathbf{q},\mu)$ be a quasi-ultrametric space
of homogeneous type and $\rho\in\mathbf{q}$
be a regularized quasi-ultrametric as in Definition~\ref{D2.3}.
Then there exists a constant $C\in[1,\infty)$
such that, for any $x_0\in X$ and $s\in(0,\mu(X)]\cap(0,\infty)$,
\begin{align}\label{1111}
\mu\left(\left\{x\in X:
V_\rho(x_0,x)<s\right\}\right)\leq
C\max\left\{s,\,\mu\left(\left\{x_0\right\}\right)\right\},
\end{align}
where $V_\rho(x_0,x_0):=0$.
\end{lemma}

\begin{proof}
Let $x_0\in X$ and $s\in(0,\mu(X)]\cap(0,\infty)$.
We consider the following two cases for $s$.

\emph{Case (1)} $s\in(0,\mu(\{x_0\})]$ with
$\mu(\{x_0\})\in(0,\infty)$. In this case,
noting that $V_\rho(x_0,x_0)=0$ and
$V_\rho(x_0,x)\ge\mu(\{x_0\})\ge s$
for any $x\in X$ with $x\neq x_0$, we have
$$
\mu\left(\left\{x\in X:
V_\rho(x_0,x)<s\right\}\right)=\mu(\{x_0\})
$$
and hence \eqref{1111} holds.

\emph{Case (2)} $s\in(\mu(\{x_0\}),\infty)$.
Let
$$
r_s:=\sup\left\{r\in(0,\mathrm{diam}_\rho(X)]
\cap(0,\infty):
(V_\rho)_r(x_0)<s\right\}\in(0,\infty).
$$
Then, for any $x\in X$ satisfying $V_\rho(x_0,x)<s$, we have
$$
\rho(x_0,x)\leq r_s<2r_s,
$$
which, together with the definition of $r_s$ and \eqref{upperdoub},
further implies that
$$
\mu\left(\left\{x\in X:
V_\rho(x,x_0)<s\right\}\right)
\leq\mu\left(B_\rho\left(x_0,2r_s\right)\right)
\lesssim\mu\left(B_\rho\left(x_0,\frac{r_s}{2}\right)\right)
<s.
$$
Combining the above two cases then completes the proof of Lemma~\ref{814L2}.
\end{proof}

\begin{remark}\label{1982}
By Lemma~\ref{A3.2}(vii),
we find that \eqref{1111} also holds with $V_\rho(x_0,x)$ therein
replaced by $V_\rho(x,x_0)$.
\end{remark}

\begin{proof}[Proof of Theorem~\ref{CZOcritical}]
Let $\rho\in\mathbf{q}$ be a regularized quasi-ultrametric and $p:=\frac{n}{n+\varepsilon}$.
To prove (i),
we first claim that, for any $(p,2)$-atom $a$,
\begin{align}
\label{qkw-34}
\sup_{k\in\mathbb{Z}}2^{kp}
\mu\left(\left\{x\in X:
\left|Ta(x)\right|>2^k\right\}\right)
\lesssim1.
\end{align}
To see this, suppose that $a\in L^2(X)$
is $(p,2)$-atom supported
in $B:=B_\rho(x_B,r_B)$, where $x_B\in X$ and $r_B\in(0,\infty)$.
We write
\begin{align}
\label{815yy0}
&\sup_{k\in\mathbb{Z}}2^{kp}
\mu\left(\left\{x\in X:
\left|Ta(x)\right|>2^k\right\}\right)
\nonumber\\
&\quad\leq
\sup_{k\in\mathbb{Z}}2^{kp}
\mu\left(\left\{x\in 2C_\rho B:
\left|Ta(x)\right|>2^k\right\}\right)
\nonumber\\
&\qquad+\sup_{k\in\mathbb{Z}}2^{kp}
\mu\left(\left\{x\in(2C_\rho B)^\complement:
\left|Ta(x)\right|>2^k\right\}\right).
\end{align}
On the one hand,
\begin{align}
\label{815yy1}
&\sup_{k\in\mathbb{Z}}2^{kp}
\mu\left(\left\{x\in 2C_\rho B:
\left|Ta(x)\right|>2^k\right\}\right)\nonumber\\
&\quad\leq\int_{2C_\rho B}\left|Ta(x)\right|^p\,d\mu(x)
\quad\text{by Chebyshev's inequality}\nonumber\\
&\quad\leq\left[\mu(2C_\rho B)\right]^{1-\frac{p}{2}}
\|Ta\|_{L^2(X)}^p
\quad\text{by H\"older's inequality}
\nonumber\\
&\quad\lesssim\left[\mu(2C_\rho B)
\right]^{1-\frac{p}{2}}\|a\|_{L^2(X)}^p
\quad\text{by the boundedness of $T$ on $L^2(X)$}\nonumber\\
&\quad\leq\left[\mu(2C_\rho B)\right]^{1-\frac{p}{2}}
\left[\mu(B)\right]^{\frac{p}{2}-1}
\quad\text{by Definition~\ref{atom}(ii)}\nonumber\\
&\quad\sim1\quad\text{by \eqref{upperdoub}}.
\end{align}
On the other hand, note that,
for any $x\in (2C_\rho B)^\complement$
and $y\in B$,
\begin{align}\label{1559}
\rho(y,x_B)<\frac{\rho(x,x_B)}{2C_\rho}.
\end{align}
From this, it follows that,
for any $x\in(2C_\rho B)^\complement$,
\begin{align*}
\left|Ta(x)\right|
&\leq\int_B\left|K(x,y)-K(x,x_B)\right|
\left|a(y)\right|\,d\mu(y)
\quad\text{by Definition~\ref{atom}(iii)}\\
&\lesssim\int_B\left[\frac{\rho(y,x_B)}
{\rho(x,x_B)}\right]^\varepsilon\frac{1}{V_\rho(x,x_B)}
|a(y)|\,d\mu(y)
\quad\text{by \eqref{1559} and \eqref{CZk2}}\\
&\leq\left[\frac{r_B}
{\rho(x,x_B)}\right]^\varepsilon\frac{1}{V_\rho(x,x_B)}
\int_B|a(y)|\,d\mu(y),
\end{align*}
which further implies that
\begin{align*}
\left|Ta(x)\right|
&\lesssim\left[\frac{r_B}
{\rho(x,x_B)}\right]^\varepsilon\frac{1}{V_\rho(x,x_B)}
\left[\mu(B)\right]^{1-\frac{1}{p}}
\quad\text{by \eqref{atomL1}}\\
&\lesssim\left[\frac{r_B}
{\rho(x,x_B)}\right]^{\varepsilon+
n(1-\frac{1}{p})}\frac{1}{V_\rho(x,x_B)}
\left[V_\rho(x,x_B)\right]^{1-\frac{1}{p}}
\quad\text{by \eqref{upperdoub}}\\
&\sim\left[V_\rho(x,x_B)\right]^{-\frac{1}{p}}
\quad\text{since $p=\frac{n}{n+\varepsilon}$},
\end{align*}
where the implicit positive constant $C$ is
independent of $a$.
By this and Lemma~\ref{1982},
we conclude that there exists a positive constant $C$,
independent of $x$ and $k$, such that
\begin{align}
\label{815yy6}
&\sup_{k\in\mathbb{Z}}2^{kp}
\mu\left(\left\{x\in(2C_\rho B)^\complement:
\left|Ta(x)\right|>2^k\right\}\right)
\nonumber\\
&\quad\leq\sup_{k\in\mathbb{Z}}2^{kp}
\mu\left(\left\{x\in X:
V_\rho(x,x_B)<(C2^k)^{-p}\right\}\cap
(2C_\rho B)^\complement\right)
\nonumber\\
&\quad=\sup_{\{k\in\mathbb{Z}:
(C2^{-k})^p\ge\mu(\{x_B\})\}}2^{kp}
\mu\left(\left\{x\in X:
V_\rho(x,x_B)<(C2^k)^{-p}\right\}\cap
(2C_\rho B)^\complement\right)
\nonumber\\
&\qquad+\sup_{\{k\in\mathbb{Z}:
(C2^{-k})^p<\mu(\{x_B\})\}}\cdots
\nonumber\\
&\quad\leq\sup_{\{k\in\mathbb{Z}:
(C2^{-k})^p\ge\mu(\{x_B\})\}}2^{kp}
\mu\left(\left\{x\in X:
V_\rho(x,x_B)<(C2^k)^{-p}\right\}\right)
\lesssim1,
\end{align}
which, together with \eqref{815yy1} and \eqref{815yy0},
completes the proof of the claim in \eqref{qkw-34}.
In light of \eqref{qkw-34}, we find that, for any given
$\lambda\in(0,\infty)$, there exists
$k_0\in\mathbb{Z}$
such that $2^{k_0}<\lambda\leq2^{k_0+1}$ and hence
\begin{align}
\label{814yy}
\lambda^p\mu\left(\left\{
x\in X:\left|Ta(x)\right|>\lambda\right\}\right)
\leq2^{(k_0+1)p}\mu\left(\left\{
x\in X:\left|Ta(x)\right|>2^{k_0}\right\}\right)
\lesssim1,
\end{align}
where the implicit positive constant
is independent of $k_0$ (hence also $\lambda$)
and the atom $a$.

Now, for any $f\in H^{p,2}_{\mathrm{at},\,\mathrm{fin}}(X)$,
there exist $N\in\mathbb{N}$,
a sequence
$\{\lambda_i\}_{i=1}^N$ in $\mathbb{R}_+$,
and a sequence of $(p,2)$-atoms $\{a_i\}_{i=1}^N$ such that
$f=\sum_{i=1}^N\lambda_ia_i$
in $({\mathcal{G}}_0^\varepsilon(\beta,\gamma))'$ and
$$
\|f\|_{H^{p,2}_{\mathrm{at},\,\mathrm{fin}}(X)}
\sim\left(\sum_{i=1}^N\lambda_i^p\right)^\frac{1}{p},
$$
which, combined with Lemma~\ref{6.29}(i), further implies that
\begin{align*}
\|f\|_{H^p_\ast(X)}\sim\|f\|_{H^{p,2}_{\mathrm{at},\,\mathrm{fin}}(X)}
\sim\left(\sum_{i=1}^N\lambda_i^p\right)^\frac{1}{p}.
\end{align*}
From this,
the linearity of $T$, and Lemma~\ref{814L1}
[here, we need \eqref{814yy}],
we infer that
\begin{align*}
&\sup_{k\in\mathbb{Z}}2^{kp}\mu
\left(\left\{x\in X:
\left|Tf(x)\right|>2^k\right\}\right)\\
&\quad=\sup_{k\in\mathbb{Z}}2^{kp}\mu
\left(\left\{x\in X:
\left|\sum_{i=1}^N\lambda_i
Ta_i(x)\right|>2^k\right\}\right)
\\
&\quad\leq\sup_{\lambda\in(0,\infty)}\lambda^p\mu
\left(\left\{x\in X:
\left|\sum_{i=1}^N\lambda_i
Ta_i(x)\right|>\lambda\right\}\right)
\\
&\quad\lesssim\sum_{i=1}^N\lambda_i^p
\sim\|f\|_{H^p_\ast(X)}^p,
\end{align*}
which further implies that
\begin{align*}
\left\|Tf\right\|_{L^{p,\infty}(X)}\lesssim\|f\|_{H^p_\ast(X)}.
\end{align*}
By this, the density of $H^{p,2}_{\mathrm{at},\,\mathrm{fin}}(X)$
in $H^p_\ast(X)$ (see Theorem~\ref{5r1} with $q=p$),
and a standard density argument,
we conclude that $T$ can be extended to a unique bounded
linear operator from $H^p_\ast(X)$ to $L^{p,\infty}(X)$.
This finishes the proof of (i).

Next, we show (ii).
Assume that $T$ satisfies $T^*1=0$.
Similarly to the proof above,
it suffices to prove that,
for any $(p,2)$-atom $a$ supported
in $B:=B_\rho(x_B,r_B)$, where $x_B\in X$ and $r_B\in(0,\infty)$,
\begin{align}
\label{815yy5}
\sup_{k\in\mathbb{Z}}2^{kp}
\mu\left(\left\{x\in X:
(Ta)^*_\rho(x)>2^k\right\}\right)
\lesssim1.
\end{align}
To this end,
from \eqref{*<M}, the boundedness of the Hardy--Littlewood maximal
operator on $L^2(X)$ (see Lemma~\ref{M}),
and an argument similar to that used
in the proof of \eqref{815yy1},
we deduce that
\begin{align}
\label{815yy2}
\sup_{k\in\mathbb{Z}}2^{kp}
\mu\left(\left\{x\in4C_\rho^2B:
(Ta)^*_\rho(x)>2^k\right\}\right)
\lesssim1.
\end{align}
Now, we consider the case
where $x\in(4C_\rho^2B)^\complement$.
Note that $\varepsilon<\mathrm{ind\,}(X,\mathbf{q})$
implies that there exist
$\beta,\gamma,\varepsilon_1\in(0,\infty)$
such that
$\varepsilon<\beta,\gamma<
\varepsilon_1\preceq\mathrm{ind\,}(X,\mathbf{q})$.
By $T^*1=0$, we find that,
for any $\phi\in{\mathcal{G}}^{\varepsilon_1}_0(\beta,\gamma)$
satisfying $\|\phi\|_{{\mathcal{G}}(x,r,\beta,\gamma)}\leq1$
for some $r\in(0,\infty)$,
\begin{align}
\label{815yy4}
\langle Ta,\phi\rangle
&=\int_XTa(y)\left[\phi(y)-\phi(x_B)\right]
\,d\mu(y)
\nonumber\\
&=\int_{\rho(y,x_B)<2C_\rho r_B}Ta(y)\left[\phi(y)-\phi(x_B)\right]
\,d\mu(y)
\nonumber\\
&\quad+\int_{2C_\rho r_B\leq\rho(y,x_B)<\frac{r+\rho(x,x_B)}{2C_\rho}}
\cdots
+\int_{\rho(y,x_B)\ge\frac{r+\rho(x,x_B)}{2C_\rho}}\cdots
\nonumber\\
&=:\mathrm{I_1}+\mathrm{I_2}+\mathrm{I_3}.
\end{align}
We claim that, for any $i\in\{1,2,3\}$,
\begin{align}
\label{815yy3}
\left|\mathrm{I}_i\right|\lesssim\left[\frac{r_B}{\rho(x,x_B)}
\right]^\varepsilon
\frac{1}{V_\rho(x,x_B)}\left[\mu(B)\right]^{1-\frac{1}{p}}.
\end{align}
Assume this for the moment.
From this, \eqref{815yy4}, \eqref{upperdoub},
and $p=\frac{n}{n+\varepsilon}$, we deduce that,
for any $x\in(4C_\rho^2B)^\complement$,
\begin{align*}
\left(Ta\right)^*_\rho(x)&\lesssim\left[\frac{r_B}{\rho(x,x_B)}
\right]^\varepsilon
\frac{1}{V_\rho(x,x_B)}\left[\mu(B)\right]^{1-\frac{1}{p}}
\\
&\lesssim\left[\frac{r_B}{\rho(x,x_B)}\right]
^{\varepsilon+n(1-\frac{1}{p})}\frac{1}{V_\rho(x,x_B)}
\left[V_\rho(x,x_B)\right]^{1-\frac{1}{p}}
\\
&\sim\left[V_\rho(x,x_B)\right]^{-\frac{1}{p}},
\end{align*}
which, together with arguments similar to those
used in the proofs of \eqref{815yy6}
and \eqref{815yy0} and with
\eqref{815yy2}, further implies \eqref{815yy5} holds.

Thus, to finish the proof of \eqref{815yy5}, we only need to show \eqref{815yy3}.
To estimate $\mathrm{I_1}$, note that $x\in(4C_\rho^2B)^\complement$ implies
$$
2C_\rho r_B\leq\frac{\rho(x,x_B)}{2C_\rho}
<\frac{r+\rho(x,x_B)}{2C_\rho},
$$
which further implies that, for any $y\in B$,
$$
\rho(y,x_B)<2C_\rho r_B<\frac{r+\rho(x,x_B)}{2C_\rho}.
$$
By this,
we conclude that
\begin{align*}
|\mathrm{I_1}|
&\lesssim\int_{\rho(y,x_B)<2C_\rho r_B}|Ta(y)|
\left[\frac{\rho(y,x_B)}{r+\rho(x,x_B)}\right]^\beta
\frac{1}{(V_\rho)_r(x)+V_\rho(x,x_B)}
\nonumber\\
&\quad\times\left[\frac{r}{r+\rho(x,x_B)}\right]^\gamma
\,d\mu(y)
\quad\text{by the regularity of
$\phi\in{\mathcal{G}}^{\varepsilon_1}_0(\beta,\gamma)$}
\nonumber\\
&\lesssim\left[\frac{r_B}{r+\rho(x,x_B)}\right]^\beta
\frac{1}{(V_\rho)_r(x)+V_\rho(x,x_B)}
\left[\frac{r}{r+\rho(x,x_B)}\right]^\gamma
\nonumber\\
&\quad\times\left[\mu(B)\right]^\frac{1}{2}
\|Ta\|_{L^2(X)}
\quad\text{by H\"older's inequality},
\end{align*}
which, combined with the boundedness of $T$ on $L^2(X)$,
$\beta>\varepsilon$, and Definition~\ref{atom}(ii),
further implies that
\begin{align}\label{815yyI1}
|\mathrm{I_1}|
&\lesssim\left[\frac{r_B}{\rho(x,x_B)}\right]^\beta
\frac{1}{V_\rho(x,x_B)}
\left[\mu(B)\right]^\frac{1}{2}
\|a\|_{L^2(X)}\nonumber\\
&\leq\left[\frac{r_B}{\rho(x,x_B)}\right]^\varepsilon
\frac{1}{V_\rho(x,x_B)}
\left[\mu(B)\right]^{1-\frac{1}{p}}.
\end{align}
To estimate $\mathrm{I_2}$, note that, for any $y\in X$ satisfying
$2C_\rho r_B\leq\rho(y,x_B)<\frac{r+\rho(x,x_B)}{2C_\rho}$
and for any $z\in B$, we have
$$
\rho(z,x_B)<r_B\leq\frac{\rho(y,x_B)}{2C_\rho}.
$$
From this, both (i) and (iii) of Definition~\ref{atom},
\eqref{CZk2}, and
the regularity of  $\phi\in{\mathcal{G}}^{\varepsilon_1}_0(\beta,\gamma)$,
we infer that
\begin{align*}
|\mathrm{I_2}|
&\leq\int_{2C_\rho r_B\leq\rho(y,x_B)<\frac{r+\rho(x,x_B)}{2C_\rho}}
\int_{B}\left|K(y,z)-K(y,x_B)\right||a(z)|\,d\mu(z)
\nonumber\\
&\quad\times
\left|\phi(y)-\phi(x_B)\right|\,d\mu(y)
\nonumber\\
&\lesssim\int_{2C_\rho r_B\leq\rho(y,x_B)<\frac{r+\rho(x,x_B)}{2C_\rho}}
\int_B\left[\frac{\rho(z,x_B)}{\rho(y,x_B)}\right]^\varepsilon
\frac{1}{V_\rho(y,x_B)}
|a(z)|\,d\mu(z)
\nonumber\\
&\quad\times
\left[\frac{\rho(y,x_B)}{r+\rho(x,x_B)}\right]^\beta
\frac{1}{(V_\rho)_r(x)+V_\rho(x,x_B)}
\left[\frac{r}{r+\rho(x,x_B)}\right]^\gamma\,d\mu(y),
\end{align*}
which, together with $\beta>\varepsilon$,
further implies that
\begin{align}\label{815yyI2}
|\mathrm{I_2}|
&\lesssim\int_{2C_\rho r_B\leq\rho(y,x_B)<\frac{r+\rho(x,x_B)}{2C_\rho}}
\left[\frac{r_B}{\rho(y,x_B)}\right]^\varepsilon
\frac{1}{V_\rho(y,x_B)}
\left[\mu(B)\right]^{1-\frac{1}{p}}
\nonumber\\
&\quad\times
\left[\frac{\rho(y,x_B)}{r+\rho(x,x_B)}\right]^\beta
\frac{1}{V_\rho(x,x_B)}\,d\mu(y)
\quad\text{by \eqref{atomL1}}
\nonumber\\
&\lesssim\left[\frac{r_B}{r+\rho(x,x_B)}\right]^\varepsilon
\frac{1}{V_\rho(x,x_B)}\left[\mu(B)\right]^{1-\frac{1}{p}}
\nonumber\\
&\quad\times
\int_{\rho(y,x_B)<\frac{r+\rho(x,x_B)}{2C_\rho}}
\left[\frac{\rho(y,x_B)}{r+\rho(x,x_B)}\right]^{\beta-\varepsilon}
\frac{1}{V_\rho(y,x_B)}\,d\mu(y)
\nonumber\\
&\lesssim\left[\frac{r_B}{\rho(x,x_B)}\right]^\varepsilon
\frac{1}{V_\rho(x,x_B)}\left[\mu(B)\right]^{1-\frac{1}{p}}
\quad\text{by Lemma~\ref{A3.2}(i)}.
\end{align}
For $\mathrm{I_3}$, we have
\begin{align}
\label{815yyI3}
|\mathrm{I_3}|&\leq
\int_{\rho(y,x_B)\ge\frac{r+\rho(x,x_B)}{2C_\rho}}
\left|Ta(y)\right||\phi(y)|
\,d\mu(y)
\nonumber\\
&\quad+|\phi(x_B)|
\int_{\rho(y,x_B)\ge\frac{r+\rho(x,x_B)}{2C_\rho}}
\left|Ta(y)\right|
\,d\mu(y)
\nonumber\\
&=:\mathrm{I_{3,1}}+\mathrm{I_{3,2}}.
\end{align}
To deal with $\mathrm{I_{3,1}}$,
note that, for any $z\in B$ and $y\in X$
satisfying
$\rho(y,x_B)\ge\frac{r+\rho(x,x_B)}{2C_\rho}$,
\begin{align*}
\rho(z,x_B)<r_B\leq\frac{\rho(x,x_B)}{4C_\rho^2}
<\frac{\rho(y,x_B)}{2C_\rho}.
\end{align*}
By this, \eqref{CZk2},
and the size condition of
$\phi\in{\mathcal{G}}^{\varepsilon_1}_0(\beta,\gamma)$,
we obtain
\begin{align*}
\mathrm{I_{3,1}}
&\leq
\int_{\rho(y,x_B)\ge\frac{r+\rho(x,x_B)}{2C_\rho}}
\int_B\left|K(y,z)-K(y,x_B)\right|
|a(z)|\,d\mu(z)|\phi(y)|\,d\mu(y)
\nonumber\\
&\qquad\ \text{by (i) and (iii) of Definition~\ref{atom}}
\\
&\lesssim\int_{\rho(y,x_B)\ge\frac{r+\rho(x,x_B)}{2C_\rho}}
\int_B\left[\frac{\rho(z,x_B)}{\rho(y,x_B)}\right]^{\varepsilon}
\frac{1}{V_\rho(y,x_B)}|a(z)|\,d\mu(z)
\nonumber\\
&\quad\times
\frac{1}{(V_\rho)_r(x)+V_\rho(x,y)}
\left[\frac{r}{r+\rho(x,y)}\right]^\gamma
\,d\mu(y).
\end{align*}
From this, it follows that
\begin{align}\label{815yyI31}
\mathrm{I_{3,1}}
&\lesssim\left[\frac{r_B}{\rho(x,x_B)}\right]^\varepsilon
\frac{1}{V_\rho(x,x_B)}\left[\mu(B)\right]^{1-\frac{1}{p}}\nonumber\\
&\quad\times\int_X\frac{1}{(V_\rho)_r(x)+V_\rho(x,y)}
\left[\frac{r}{r+\rho(x,y)}\right]^\gamma\,d\mu(y)
\quad\text{by \eqref{atomL1}}
\nonumber\\
&\lesssim\left[\frac{r_B}{\rho(x,x_B)}\right]^\varepsilon
\frac{1}{V_\rho(x,x_B)}\left[\mu(B)\right]^{1-\frac{1}{p}}
\quad\text{by Lemma~\ref{A3.2}(ii)}.
\end{align}
To deal with $\mathrm{I_{3,2}}$,
by the size condition of
$\phi\in{\mathcal{G}}^{\varepsilon_1}_0(\beta,\gamma)$
and \eqref{CZk2},
we find that
\begin{align*}
\mathrm{I_{3,2}}
&\leq|\phi(x_B)|
\int_{\rho(y,x_B)\ge\frac{r+\rho(x,x_B)}{2C_\rho}}
\int_B\left|K(y,z)-K(y,x_B)\right|
|a(z)|\,d\mu(z)\,d\mu(y)
\\
&\qquad\ \text{by (i) and (iii) of Definition~\ref{atom}}
\\
&\lesssim\frac{1}{(V_\rho)_r(x)+V_\rho(x,x_B)}
\left[\frac{r}{r+\rho(x,x_B)}\right]^\gamma
\int_{\rho(y,x_B)\ge\frac{r+\rho(x,x_B)}{2C_\rho}}
\\
&\quad\times\int_B\left[\frac{\rho(z,x_B)}
{\rho(y,x_B)}\right]^\varepsilon
\frac{1}{V_\rho(y,x_B)}
|a(z)|\,d\mu(z)\,d\mu(y),
\end{align*}
which further implies that
\begin{align*}
\mathrm{I_{3,2}}
&\lesssim\left[\frac{r}{\rho(x,x_B)}\right]^\varepsilon
\frac{1}{V_\rho(x,x_B)}\|a\|_{L^1(X)}\\
&\quad\times\int_{\rho(y,x_B)\ge\frac{r+\rho(x,x_B)}{2C_\rho}}\left[\frac{r_B}
{\rho(y,x_B)}\right]^\varepsilon
\frac{1}{V_\rho(y,x_B)}\,d\mu(y)
\quad\text{since $\gamma>\varepsilon$}\\
&\leq\left[\frac{r_B}{\rho(x,x_B)}\right]^\varepsilon
\frac{1}{V_\rho(x,x_B)}\left[\mu(B)\right]^{1-\frac{1}{p}}
\\
&\quad\times
\int_{\rho(y,x_B)\ge\frac{r+\rho(x,x_B)}{2C_\rho}}
\left[\frac{r+\rho(x,x_B)}
{\rho(y,x_B)}\right]^\varepsilon
\frac{1}{V_\rho(y,x_B)}\,d\mu(y)
\quad\text{by \eqref{atomL1}}
\\
&\lesssim\left[\frac{r_B}{\rho(x,x_B)}\right]^\varepsilon
\frac{1}{V_\rho(x,x_B)}\left[\mu(B)\right]^{1-\frac{1}{p}}
\quad\text{by Lemma~\ref{A3.2}(i)}.
\end{align*}
From this, \eqref{815yyI1}, \eqref{815yyI2},
\eqref{815yyI3}, and \eqref{815yyI31},
we finally deduce that \eqref{815yy3} holds.
This finishes the proof of (ii)
and hence Theorem~\ref{CZOcritical}.
\end{proof}

\begin{remark}
Let $(X,\rho,\mu)$ be a space of homogeneous type
with upper dimension $n\in(0,\infty)$, and assume that
$\varepsilon\in(0,\eta)$, where $\eta\in(0,1)$ is the
smoothness index of the wavelets on $X$
constructed by Auscher and Hyt\"onen in \cite{AH1},
and that $T$ is an $\varepsilon$-Calder\'on--Zygmund
operator as in Definition~\ref{CZO}.
Zhou et al. \cite[Theorems 8.12 and 8.13]{zhy}
and Sun et al. \cite[Theorem 8.1(iii)]{syy21} both proved
that $T$ can be extended to a unique bounded linear operator
from $H^{\frac{n}{n+\varepsilon}}_\ast(X)$
to $L^{\frac{n}{n+\varepsilon},\infty}(X)$
and, if $T$ satisfies that $T^*1=0$,
then $T$ can be extended to a unique bounded linear operator
from $H^{\frac{n}{n+\varepsilon}}_\ast(X)$
to $H^{\frac{n}{n+\varepsilon},\infty}_\ast(X)$.
Theorem~\ref{CZOcritical} improves on this work
since it is possible for $\varepsilon$ to be greater than 1 in
certain quasi-ultrametric spaces. For example, if
$(X,\mathbf{q})$ is an ultrametric space then
$\mathrm{ind\,}(X,\mathbf{q})=\infty$
and hence Theorem~\ref{CZOcritical}
is valid for all $\varepsilon\in(0,\infty)$.
\end{remark}




\printindex[symbols]
\printindex[words]


\newpage

\markboth{\scriptsize\rm\sc Abstract}{\scriptsize\rm\sc Abstract}
\addcontentsline{toc}{chapter}{Abstract}

\chapter*{Abstract}

Let $(X,\mathbf{q},\mu)$ be an ultra-RD-space
with upper dimension $n\in(0,\infty)$;
i.e., it is a quasi-ultrametric space of homogeneous type whose
measure $\mu$ satisfies an additional reverse doubling property.
Let $\mathrm{ind\,}(X,\mathbf{q})\in(0,\infty]$ denote
its lower smoothness index, as
introduced by Mitrea et al.
In this monograph, the authors first construct a new
approximation of the identity
on quasi-ultrametric
spaces of homogeneous type,
achieving a maximal degree of smoothness
$0<\varepsilon\preceq\mathrm{ind\,}(X,\mathbf{q})$.
This fundamental tool is then used to
derive sharp homogeneous (as well as inhomogeneous)
continuous/discrete Calder\'on reproducing formulae
on ultra-RD-spaces. As applications,
the authors establish Littlewood--Paley function
characterizations for both Hardy spaces and Triebel--Lizorkin
spaces on ultra-RD-spaces.
The authors further introduce
Hardy--Lorentz spaces $H^{p,q}_\ast(X)$ via
the grand maximal function, with the sharp range
$p\in(\frac{n}{n+\mathrm{ind\,}(X,\mathbf{q})},\infty)$
and $q\in(0,\infty]$, and provide their
real-variable characterizations
using radial/non-tangential
maximal functions, (finite) atoms, molecules,
and various Littlewood--Paley functions.
Based on these characterizations, the authors prove
a duality theorem between Hardy--Lorentz spaces and
Campanato--Lorentz spaces, establish a real interpolation
theorem for Hardy--Lorentz spaces, and
derive boundedness results for Calder\'on--Zygmund operators
on them.
It should be emphasized
that many of the main results in
this monograph are indeed established in the more general setting of
quasi-ultrametric spaces of homogeneous type.

The novelty of this monograph resides in two key aspects.
First, the assumptions on the underlying space are minimal:
the authors only assume that $\mu$ is Borel-semiregular,
that singletons may have strictly positive measure,
and that $\rho$ is a quasi-ultrametric,
which provides a more general theoretic setting
than the Coifman and Weiss notion of spaces of homogeneous type.
Second, the maximal order of smoothness achieved
by the approximation of the identity and
the regularity exponent in the Calder\'on reproducing
formulae are both sharp,
as they are intrinsically based on the sharp
regularity index $\mathrm{ind\,}(X,\mathbf{q})$
of the quasi-ultrametric.
Given the sharpness of these tools,
the authors are able to develop the real-variable
theory of the Hardy--Lorentz space
$H^{p,q}_\ast(X)$ for optimal
ranges of $p$ and $q$ (which cannot
be obtained in any of the related results
currently available in the literature),
under minimal assumptions on $X$.

All results presented here are new and to appear
in print for the first time.
Designed to be largely self-contained,
this monograph is an invaluable resource for
graduate students and researchers interested
in the frontier of the function space theory
on ultra-RD spaces.
\end{document}